\documentclass[12pt]{book}
\usepackage{fancyhdr, array,calc,graphicx,url,tabularx}
\usepackage{geometry}
\usepackage{latexsym}
\usepackage{amssymb}
\usepackage{amsmath}
\usepackage{thmtools}

\usepackage{makeidx}
\usepackage{enumerate}
\usepackage{verbatim}
\usepackage{tikz}
\usepackage[colorlinks=true,linkcolor=blue]{hyperref}
\usepackage{ntheorem,lipsum}

\usepackage{changepage}
\usepackage{color}

\usepackage[most]{tcolorbox}

\usepackage{fancybox}

\usepackage{xcolor}
\usepackage{blindtext}
\usepackage{tikz}

\makeindex

\newenvironment{dedication}
{
   \cleardoublepage
   \thispagestyle{empty}
   \vspace*{\stretch{1}}
   \hfill\begin{minipage}[t]{0.66\textwidth}
   \raggedright
}%
{
   \end{minipage}
   \vspace*{\stretch{3}}
   \clearpage
}

\def\undertilde#1{{\baselineskip=0pt\vtop
  {\hbox{$#1$}\hbox{$\scriptscriptstyle\sim$}}}{}}

\newcommand{\gTheta}{\mathsf{G}}
\newcommand{\utilde}{\undertilde}

\renewcommand{\gg}{\gamma}

\newcommand{\gG}{\Gamma}

\newcommand{\J}{{\mathcal{J}}}
\newcommand{\C}{{\sf{core}}}

\newcommand{\layer}{{\sf{layer}}}

\newcommand{\rnc}{{\sf{rnc}}}

\newcommand{\bP}{\bold \Pi}
\newcommand{\bR}{{\mathbb{R}}}

\newcommand{\rest}{\restriction}
\newcommand{\la}{\langle}
\newcommand{\ra}{\rangle}

\newcommand{\ind}{{\rm ind}}
\newcommand{\rge}{{\rm rge}}
\newcommand{\ml}{{\rm ml}}

\newcommand{\card}[1]{{\vert #1 \vert} }

\newcommand{\forces}{\Vdash}

\renewcommand{\models}{\vDash}
\newcommand{\powerset}{{\wp}}
\newcommand{\univ}[1]{{\lfloor #1 \rfloor}}

\newcommand{\dom}{{\rm dom}}
\newcommand{\rng}{{\rm rng}}
\newcommand{\cp}{{\rm crit }}

\newcommand{\cf}{{\rm cf}}
\newcommand{\lh}{{\rm lh}}

\newcommand{\ord}{{\rm ord}}

\newcommand{\sfc}{{\sf{Code}}}

\newcommand{\m}{{\rm m}}

\newcommand{\myqedhere}{\hfill \dashv}

\newtheorem{theorem}{Theorem}[section]
\newtheorem{proposition}[theorem]{Proposition}

\theorembodyfont{\normalfont}
\newtheorem{definition}[theorem]{Definition}
\newtheorem{rem}[theorem]{Remark}

\newtheorem{lemma}[theorem]{Lemma}
\newtheorem{corollary}[theorem]{Corollary}
\newtheorem{claim}[theorem]{Claim}

\newtheorem{conjecture}[theorem]{Conjecture}

\newtheorem{sublemma}[theorem]{Sublemma}

\newtheorem{remark}[theorem]{Remark}
\newtheorem{notation}[theorem]{Notation}
\newtheorem{terminology}[theorem]{Terminology}
\numberwithin{figure}{section}
\newtheorem{problem}[theorem]{Problem}

\newenvironment{proof}{{\it{
Proof.}}}{\nopagebreak\mbox{}{\hfill$\square$}
\par\bigskip}
\newenvironment{Claim}{\noindent{\it
Claim} {}}{}
\newenvironment{Case}{\noindent{\it
Case}}{}

\newcommand{\rcon}[1]{Conjecture~\ref{#1}}

\newcommand{\rprop}[1]{Proposition~\ref{#1}}
\newcommand{\rthm}[1]{Theorem~\ref{#1}}
\newcommand{\rlem}[1]{Lemma~\ref{#1}}
\newcommand{\rsublem}[1]{Sublemma~\ref{#1}}

\newcommand{\rcor}[1]{Corollary~\ref{#1}}
\newcommand{\rdef}[1]{Definition~\ref{#1}}

\newcommand{\rsec}[1]{Section~\ref{#1}}
\newcommand{\rchap}[1]{Chapter~\ref{#1}}
\newcommand{\rsubsec}[1]{Section~\ref{#1}}
\newcommand{\rrem}[1]{Remark~\ref{#1}}
\newcommand{\rnot}[1]{Notation~\ref{#1}}
\newcommand{\rter}[1]{Terminology~\ref{#1}}

\def\so{{\sf{short}}}
\def\ma{{\sf{max}}}

\def\inseg{\trianglelefteq}

\def\k{\kappa}
\def\a{\alpha}
\def\b{\beta}
\def\d{\delta}

\def\l{\lambda}

\def\P{{\mathcal{P} }}
\def\W{{\mathcal{W} }}
\def\V{{\mathcal{V} }}
\def\Q{{\mathcal{ Q}}}
\def\mH{{\mathcal{ H}}}
\def\K{{\mathcal{ K}}}
\def\J{{\mathcal{ J}}}

\def\R{{\mathcal R}}
\def\X{{\mathcal X}}
\def\Y{{\mathcal Y}}
\def\H{{\rm{HOD}}}
\def\M{{\mathcal{M}}}
\def\N{{\mathcal{N}}}
\def\F{{\mathcal{F}}}
\def\T {{\mathcal{T}}}
\def\U{{\mathcal{U}}}
\def\S{{\mathcal{S}}}
\def\F{{\mathcal{F}}}

\def\VT{{\vec{\mathcal{T}}}}
\def\VU{{\vec{\mathcal{U}}}}

\def\card#1{\left|#1\right|}

\def\iff{\mathrel{\leftrightarrow}}

\def\and{\mathrel{\kern1pt\&\kern1pt}}

\def\inseg{\triangleleft}
\def\insegeq{\trianglelefteq}

\def\<#1>{\langle\,#1\,\rangle}

 \input xy
 \xyoption{all}

\title{The Largest Suslin Axiom 
}

\author{Grigor Sargsyan \\
        Institute of Mathematics\\
        Polish Academy of Sciences\\
        https://www.impan.pl/$\sim$gsargsyan/\\
        gsargsyan@impan.pl\\\\
        Nam Trang\\
Department of Mathematics\\ 
University of North Texas, Denton\\
http://math.uci.edu/$\sim$ntrang\\
Email: ntrang@math.unt.edu.}

\date{\today}
\begin{document}

\maketitle

\begin{dedication}
The first author dedicates the book to his father, mother, wife, son and daughter: Edvard Sargsyan, Gayane Shahinyan, Aleksandra Sawa, Tigran Sargsyan and Nairi Sargsyan.\vspace{1cm}

The second author dedicates the book to his father, mother, wife and son: Phil Trang, Tram Pham, Tran Tran and Shawn Trang. 
\end{dedication}

\tableofcontents 


\chapter{Introduction}

This manuscript is a contribution to \textit{descriptive inner model theory}, which is the area of set theory that lies between descriptive set theory as developed in \cite{Moschovakis} and inner model theory. The main goal of this manuscript is to advance the descriptive inner model theoretic methods to the level of the $\sf{Largest\ Suslin\ Axiom}$ ($\sf{LSA}$), which is a strong determinacy axiom asserting that there is a largest Suslin cardinal and that the largest Suslin cardinal is a member of the Solovay sequence. In more concrete terms, our goal is twofold: Firstly develop methods for analyzing the minimal model of $\sf{LSA}$, and secondly, develop methods for building the minimal model of ${\sf{LSA}}$ under various hypotheses such as the $\sf{Proper\ Forcing\ Axiom}$ or $\sf{Large\ Cardinals}$. Since the introduction of Steel's recent manuscript  \cite{SteelComp}, the expository paper \cite{BSL}  and the introduction of \cite{SargTrang2021} contain all the introductory information we need, here we will not introduce the subject matter of this book and instead will hope that the reader has consulted these sources. 

The first problem is an instance of the problem Steel mentions on page xii of \cite{SteelComp} where he writes: ``The most important of the remaining open problems is whether, assuming determinacy, there actually are mouse pairs at every appropriate level of logical complexity". \rthm{the generation of mouse full pointclasses} shows that the aformentioned problem has a positive solution in the minimal model of ${\sf{LSA}}$. As explained in any of the sources cited above, the goal for doing this is to show that letting $\Theta$ be the least ordinal that is not a surjective image of the reals, $V_\Theta^\H$ as computed inside a determinacy model is a hod premouse. The above sources explain the importance of having a hod premouse representation of  $V_\Theta^\H$.

The second problem amounts to advancing the $\sf{Core\ Model\ Induction}$ to the level of ${\sf{LSA}}$. \rcor{cor:pfa_lsa} constructs the minimal model of ${\sf{LSA}}$ assuming ${\sf{PFA}}$. More dramatically, the paper \cite{SargTrang2021}, which extends the methods of this manuscript, demonstrates that the $\sf{Core\ Model\ Induction}$, in its current form, cannot be used to go much further than ${\sf{LSA}}$. 

\rcor{cor:pfa_lsa} also builds the minimal model of $\sf{LSA}$ directly from large cardinals, namely strongly compact cardinals. However, \rthm{lst from wlw} shows that ${\sf{LSA}}$ is weaker than a Woodin cardinal that is a limit of Woodin cardinals, and so strong compactness seems to be much more than needed. Nevertheless, while it is widely believed that strongly compact cardinals are consistency wise stronger than a Woodin cardinal that is a limit of Woodin cardinals, this is not yet known. 
Still we strongly believe that the methods developed in this manuscript, the methods of \cite{SargAdolf} and the main theorem of \cite{Neeman} can be used to show that assuming the existence of a Woodin cardinal that is a limit of Woodin cardinals, the minimal model of ${\sf{LSA}}$ exists (cf. Definition \ref{min model of lsa}). 

Historically, ${\sf{LSA}}$ was introduced by Woodin in \cite[Remark 9.28]{Woodin}, and it features prominently  in Woodin's $\sf{Ultimate\ L}$ framework (see \cite[Definition 7.14]{UltimateL} and Axiom I and Axiom II on page 97 of \cite{UltimateL}\footnote{The requirement in these axioms that there is a strong cardinal which is a limit of Woodin cardinals is only possible if $L(A, \bR)\models \sf{LSA}$.}). \rthm{lst from wlw} is historically the first proof of the consistency of ${\sf{LSA}}$ relative to large cardinals. Cramer and Woodin established the consistency of ${\sf{LSA}}$ from large cardinals in the region of $I_0$ (see \cite[Theorem 65]{CramerIVLC}). \\\\
\textbf{The technical content of the manuscript}\\\\
\textbf{1. The Largest Suslin Axiom}\\
 ${\sf{LSA}}$ is a determinacy theory whose underlying theory is Woodin's $\sf{AD}^+$. Chapter 9.1 of \cite{Woodin} provides a quick overview of ${\sf{AD^+}}$, and Larson's recent manuscript \cite{LarsonAD} provides more details. Perhaps the most important consequence of ${\sf{AD^+}}$ is the fact that assuming $V=L(\powerset(\bR))$, the fragment of $V$ coded by the Suslin, co-Suslin sets of reals is $\Sigma_1$ elementary in $V$ (see Theorem 9.7 of \cite{Woodin}). 

 We will need the following concepts to introduce ${\sf{LSA}}$. A cardinal $\k$ is $\sf{OD}$-inaccessible if for every $\a<\k$ there is no surjection $f: \powerset(\a)\rightarrow \k$ that is definable from ordinal parameters. A set of reals $A\subseteq \bR$ is $\k$-\textit{Suslin} if for some tree $T$ on $\k$, $A=p[T]$\footnote{Given a cardinal $\k$, we say $T\subseteq \bigcup_{n<\omega} \omega^n \times \k^n$ is a \textit{tree} on $\k$ if $T$ is closed under initial segments. Given a tree $T$ on $\k$, we let $[T]$ be the set of its branches, i.e., $b\in [T]$ if $b\in \omega^\omega\times \k^\omega$ and letting $b=(b_0, b_1)$, for each $n\in \omega$, $(b_0\rest n, b_1\rest n)\in T$. We then let $p[T]=\{ x\in \bR: \exists f((x, f)\in [T])\}$.}. A set $A$ is \textit{Suslin} if it is $\kappa$-Suslin for some $\kappa$; $A$ is \textit{co-Suslin} if its complement $\mathbb{R}\backslash A$ is Suslin.  A set $A$ is \textit{Suslin, co-Suslin} if both $A$ and its complement are Suslin. A cardinal $\k$ is a \textit{Suslin cardinal} if there is a set of reals $A$ such that $A$ is $\k$-Suslin but $A$ is not $\l$-Suslin for any $\l<\k$.\index{Suslin cardinal}\index{$\kappa$-Suslin}\index{$\sf{OD}$-inaccessible}\index{co-Suslin} Suslin cardinals play an important role in the study of models of determinacy as can be seen by just flipping through the Cabal Seminar Volumes (\cite{Cabal2008}, \cite{Cabal2012}, \cite{Cabal2016}, \cite{Cabal2021}). $\sf{LSA}$ is then the following theory.
\begin{definition}\label{dfn:lsa} \index{$\sf{Largest\ Suslin\ Axiom}$}
The $\sf{Largest\ Suslin\ Axiom}$, abbreviated as $\sf{LSA}$, is the conjunction of the following statements:
\begin{enumerate}
\item ${\sf{ZF}}+\sf{AD}^+$.
\item There is a largest Suslin cardinal.
\item The largest Suslin cardinal is $\sf{OD}$-inaccessible.
\end{enumerate}
$\myqedhere$
 \end{definition}
 
 $\sf{LSA}$ can also be defined in terms of the \textit{Solovay sequence}.

\begin{definition}\label{solovay sequence}\index{Solovay sequence}
The Solovay sequence is a sequence $(\theta_\a: \a\leq \Omega)$ such that 
\begin{enumerate}
\item $\theta_0=\sup\{\b:\exists f:\powerset(\omega)\rightarrow \b (f$ is an $OD$ surjection$)\}$,
\item if $\theta_\a<\Theta$ then $\theta_{\a+1}=\sup\{\b:\exists f:\powerset(\theta_\a)\rightarrow \b (f$ is an $OD$ surjection$)\}$,
\item for limit $\l\leq\Omega$, $\theta_\lambda=\sup_{\a<\l}\theta_\a$. 
\item $\theta_\Omega = \Theta$.
\end{enumerate}
$\myqedhere$
\end{definition}

\begin{remark}
$\sf{LSA}$ is then equivalent to the conjunction of the following axioms:
\begin{enumerate}
\item ${\sf{ZF}}+\sf{AD}^+$.
\item For some ordinal $\a$, $\Theta=\theta_{\a+1}$ and $\theta_\a$ is the largest Suslin cardinal $<\Theta$.  
\end{enumerate}
$\myqedhere$
\end{remark}

The above equivalence can be shown using the material of Chapter 9.1 of \cite{Woodin}. We note that it follows from \cite[Theorem 9.12]{Woodin} that $\sf{LSA}$ implies $\neg \sf{AD}_{\bR}$. \\\\
\textbf{2. The minimal model of LSA}\\
Suppose $V$ is a model of ${\sf{LSA}}$. Let $\k$ be the largest Suslin cardinal and suppose $A\subseteq \bR$ has Wadge rank $\k$. It then follows that $L(A, \bR)\models {\sf{LSA}}$. Keeping this fact in mind, we make the following definition. 
\begin{definition}\label{min model of lsa} Suppose $T$ is a first order theory extending ${\sf{AD^+}}$. We say that $M$ is \textbf{a minimal model} of $T$ if 
\begin{itemize}
\item $M$ is transitive and $M\models T$,
\item $\bR, Ord\subseteq M$, and 
\item for every $N$ that is a (definable) class of $M$ and contains all the reals and ordinals, either $N=M$ or $N\models \neg {\sf{LSA}}$.
\end{itemize}
$\myqedhere$
\end{definition}
It follows that all minimal models of ${\sf{LSA}}$ have the form $L(A, \bR)$. A natural question is whether there is a unique minimal model of ${\sf{LSA}}$. We will show (see the proof of \rthm{lst from wlw}) that in fact there is a unique minimal model of ${\sf{LSA}}$ which is naturally \textit{the} minimal model of ${\sf{LSA}}$. Woodin's proof of the existence of divergent models of ${\sf{AD^+}}$ also shows that not all extensions of ${\sf{AD^+}}$ have a unique minimal model (see \cite[Theorem 6.1]{EA}). 

The minimal model of ${\sf{LSA}}$ may not actually be big. For example, if $N$ is a transitive model of $\sf{AD^+}$ that contains the minimal model $M$ of ${\sf{LSA}}$ and has a Suslin cardinal $>\Theta^M$ then $\Theta^M<\theta_0^N$. In particular, every set of reals in $M$ is ordinal definable from a real in $N$. Motivated by this fact, we make the following definition.

\begin{definition} Suppose $M$ is a transitive model containing all the reals and ordinals and such that $M\models {\sf{AD^+}}+V=L(\powerset(\bR))$. We say $M$ is \textbf{full} if for all transitive $N$ such that
\begin{itemize}
\item $M\subseteq N$ and
\item $N\models ``{\sf{AD^+}}+V=L(\powerset(\bR))$",
\end{itemize}
$\Theta^M$ is a member of the Solovay sequence of $N$. $\myqedhere$
\end{definition}

The following interesting problem seems central to our understanding of the models of ${\sf{AD^+}}$ that we build from large cardinals or from other hypothesis.
\begin{problem} Do large cardinals or forcing axioms such as ${\sf{PFA}}$ imply that there is a full model of ${\sf{LSA}}$?
\end{problem}
\textit{In particular, whether the models of determinacy obtained as derived models of $V$ contain full models of ${\sf{LSA}}$ or not is a major open problem of the area.} Here we make the following conjecture which is motivated by Woodin's \textit{Sealing Theorem} (see \cite{Stationarytower}). Below ${\sf{uB}}$ stands for the set of universally Baire sets and for a generic $g$, ${\sf{uB}}_g=({\sf{uB}})^{V[g]}$ and $\bR_g=\bR^{V[g]}$.
\begin{conjecture}\label{lsa conjecture} Suppose $\k$ is a supercompact cardinal and there is a proper class of Woodin cardinals. Let $g\subseteq Coll(\omega, 2^{2^\k})$ be generic. Then in $L({\sf{uB}}_g, \bR_g)$, for each $\xi<\Theta$ there is $\a$ such that $\theta_{\a}\in (\xi, \Theta)$ and $\theta_\a$ is the largest Suslin cardinal below $\theta_{\a+1}$.
\end{conjecture} 
Thus, in the set up of the conjecture, $L({\sf{uB}}_g, \bR_g)$ has full models of ${\sf{LSA}}$ that are cofinal in its Wadge hierarchy. The following is what is known on \rcon{lsa conjecture}.  Woodin (unpublished) has shown that $L({\sf{uB}}_g, \bR_g)\models ``{\sf{AD}}_{\mathbb{R}}+\Theta$ is a regular cardinal". Sandra M\"uller and the first author recently showed that $L({\sf{uB}}_g, \bR_g)$ can be represented as a derived model of some iterate of $V$. They also found a stationary-tower-free proof of Woodin's Sealing Theorem. These results are unpublished. \cite{Steel2007} presents a stationary-tower-free proof of the derived model theorem.

The question of whether the Cramer-Woodin model of ${\sf{LSA}}$ from \cite[Theorem 65]{CramerIVLC} is a full model of ${\sf{LSA}}$ or not seems not only interesting but also important for understanding the relationship between large cardinals and models of ${\sf{AD^+}}$. \\\\
\textbf{3. The content of this manuscript}\\
In this manuscript, we establish three kinds of results that can be stated without mentioning the  technology developed to prove them. The first set of results deals with the minimal model of $\sf{LSA}$. Assume $V$ is the minimal model of $\sf{LSA}$. Then the following holds.\\\\
(A) (\rthm{computation of hod}) $\H\models \sf{GCH}$.\\
(B)  (\rthm{msc}) The Mouse Set Conjecture holds. \\

The second set of results contains a single result which shows the consistency of $\sf{LSA}$ relative to large cardinals. We will show the following.\\\\
(C)  (\rcor{lst from wlw}) Suppose the theory $\sf{ZFC}+``$there is a Woodin cardinal that is a limit of Woodin cardinals" is consistent. Then so is $\sf{LSA}$.\\\\
The third type of result establishes the existence of the minimal model of $\sf{LSA}$ assuming combinatorial principles or forcing axioms. The following belongs to this group.\\\\
(D)  (\rcor{cor:pfa_lsa}) Assume $\sf{PFA}$. Then the minimal model of ${\sf{LSA}}$ exists.\\

The precursors of these results already exist in print. The first author demonstrated versions of (A), (B), and (C) for the theory $\sf{AD}_{\mathbb{R}}+``\Theta$ is a regular cardinal". The second author proved the version of (D) for the same theory. The interested reader may consult \cite{ATHM}, \cite{CuBF} and \cite{Trang2015PFA}. The reason to prove such results is to demonstrate that the underlying technical theory is robust and can be used in a  wide range of situations. 

Recently the authors of \cite{JAMS} showed that the theory ${\sf{CH}}+``$there is an $\omega_1$-dense ideal on $\omega_1$" implies that the minimal model of ${\sf{AD}}_{\bR}+``\Theta$ is a regular cardinal" exists. This, along with an earlier result of Woodin, show that these two theories are equiconsistent. This solved part of Problem 12 of \cite{Woodin}. Whether there is a natural hypothesis asserting the existence of an ideal on a small cardinal that is equiconsistent with ${\sf{LSA}}$ is an interesting problem. In particular, letting $M'$ be the minimal model of ${\sf{LSA}}$, $\k$ be the largest Suslin cardinal of $M'$ and $M=L(\Gamma, \bR)$ where $\Gamma=\{ A\in \powerset(\bR)\cap M': w(A)<\k\}$\footnote{$w(A)$ is the Wadge rank of $A$.}, the model $M[G*H]$ where $G*H\subseteq Coll(\omega_1, \bR)*Add(1, \omega_2)$ is $M$-generic has not be studied at all. The model $M[G*H]$ where $G*H\subseteq \mathbb{P}_{max}*Add(1, \omega_3)$ has been investigated in \cite{CLSSSZ}, but not much is known beyond \cite{CLSSSZ}\footnote{But see also \cite{LarsonSarg}.}. \\\\
\textbf{4. The necessity of the short-tree-strategy mice}\\
 We do not know how to prove (B)-(D) using the methods of \cite{SteelComp}, and whether this is possible or not is a very interesting question\footnote{\cite{SteelComp} does show that $\mH\models {\sf{GCH}}$ but only assuming ${\sf{HPC}}$.}. The main issue seems to be the absence of an analysis of the ${\sf{LSA}}$ stages of the Solovay sequence using the least-branch hierarchy. The main technical concept we use to analyze such levels is the notion of a \textit{short-tree-strategy mice}, which is developed in \rchap{chap: shorttreestrategymice}. Thus, the question is whether it is necessary to develop this theory in order to prove results like (B)-(D). 

The main issue is the following. Assume ${\sf{AD}^+}$. Suppose $\theta_{\a+1}<\Theta$ and $\theta_\a$ is the largest Suslin cardinal below $\theta_{\a+1}$. Then if $(\P, \Sigma)$ is the hod pair generating the pointclass $\Gamma_1=\{A\subseteq \bR: w(A)<\theta_{\a+1}\}$ then letting $\d$ be the largest Woodin cardinal of $\P$, $((\P|\d)^\#, \Sigma^{stc})$ is the pair generating the pointclass $\Gamma_0=\{ A\subseteq \bR: w(A)<\theta_{\a}\}$. If one's goal is to show that assuming ${\sf{AD}_{\mathbb{R}}}+{\sf{DC}}+V=L(\powerset(\bR))$, $\H\models {\sf{GCH}}$ then it maybe possible to \textit{skip} $\Gamma_0$ and build the generator of $\Gamma_1$. The problem with skipping $\Gamma_0$ and moving to $\Gamma_1$ is exactly the fact that it is then unclear how to prove theorems like (A)-(D). What one would have liked is some sort of hybrid method that doesn't skip $\Gamma_0$ but also incorporates ideas from \cite{SteelComp} to avoid the theory of short-tree-strategy mice. It seems to us that this may not be possible.

Suppose then we decide not to skip over $\Gamma_0$, and suppose we have succeeded in building a generator $((\P|\d)^\#, \Sigma^{stc})$ for $\Gamma_0$. At this stage, we do not know what $(\P, \Sigma)$ must be and can only see $((\P|\d)^\#, \Sigma^{stc})$. Set then $\Q=(\P|\d)^\#$ and $\Lambda=\Sigma^{stc}$. What we need to show next is that we can extend $\Q$ to $\P$ in such a way that the following hold\footnote{Below $H_\d^\P$ is the set of all $X\in \P$ whose hereditary cardinality is $<\delta$.}:
\begin{enumerate}
\item $\d$ is the largest cardinal of $\P$ and $H_\d^\P=\Q|\d$,
\item for all $A\subseteq \d$, $A\in \P$ if and only if $A$ is ordinal definable from $(\Q, \Lambda)$,
\item $\P\models ``\d$ is a Woodin cardinal".
\end{enumerate}
The main issue seems to be with proving clause 2. It is a version of $\sf{MSC}$ for $\Lambda$, and the only way we know how to prove it is by building a $\Lambda$-mouse over $\Q$ whose derived model contains the set $\{(x, y)\in \bR^2: x$ is ordinal definable from $y$ and $(\Q, \Lambda)\}$. This requires a certain level of uniformity: $\Q$ and what we build on the top of $\Q$ have to be the same kind of objects, as otherwise the construction over $\Q$ can project across $\d$ violating clause 3 above. \\\\
\textbf{5. Some historical remarks on the large cardinal structure of hod mice}\\
The large cardinal structure of hod mice has been somewhat of a mystery. While originally it seemed hod mice must have very limited large cardinal structure, nowadays the prevailing belief is that they in fact can have any large cardinal whatsoever\footnote{At least in the short-extender region.}. First we make the following definition.
\begin{definition} \index{$\Theta_{reg}$}$\Theta_{reg}$ is the theory ${\sf{ZF}}+{\sf{AD}}_{\mathbb{R}}+``\Theta$ is a regular cardinal". $\myqedhere$
\end{definition}
Prior to \cite{ATHM},  the theory $\Theta_{reg}$ was believed to be beyond the short extender region and was believed to be at the complexity level of supercompact cardinals. Because Woodin was able to force strong combinatorial statements over a model of $\Theta_{reg}$ that would normally require large cardinals at the level of supercompact cardinals or beyond\footnote{E.g., ${\sf{MM^{++}}}(c)$ (see \cite{Woodin}) and ${\sf{CH}}+``$there is an $\omega_1$-dense ideal on $\omega_1$" (see \cite{JAMS}).}, the aforementioned belief seemed to be very reasonable. 

The main goal of \cite{ATHM}, which is based on the first author's PhD thesis, was to analyze $\H$ of the minimal model of $\Theta_{reg}$\footnote{Prior \cite{ATHM}, it was not know that there is a unique minimal model of $\Theta_{reg}$.}. While any model of determinacy has a rich large cardinal structure below its $\Theta$\footnote{E.g., $\Theta$ is a limit of strong partition cardinals, see \cite{KKMW}.}, $V_\Theta^\H$ of the minimal model of $\Theta_{reg}$ is very simple in the following sense (see \rthm{exp: woodin cardinals}).
 
Suppose $V\models {\sf{AD^+}}$. The Solovay pointclasses are exactly the stages of the Wadge hierarchy where a ``new" \textit{non-definable from below} set appears. For $\a$ such that $\theta_\a\leq \Theta$ let ${\sf{SP}}_\a=\{B\subseteq \bR: w(B)<\theta_\a\}$\index{$\sf{SP}_\a$}.  If $\theta_\a<\Theta$ and $A\subseteq \bR$ has Wadge rank $\theta_\a$ then $A$ is not ordinal definable from any set of reals $B\in  {\sf{SP}}_\a$ and moreover, every set in ${\sf{SP}}_{\a+1}$ is ordinal definable from $A$ and a real. Thus, in a sense, once we perceive a set of reals of Wadge rank $\theta_\a$, we know everything about ${\sf{SP}}_{\a+1}$. Putting it differently,\\\\
$\dagger:$ in the Wadge hierarchy, nothing of any interest happens among sets whose Wadge rank belongs to the interval $(\theta_\a, \theta_{\a+1})$.\\\\
In general, $\dagger$ is not true. All sorts of structures: Suslin cardinals, large cardinals with complicated partition properties etc, exist in that Wadge interval. However, the hod mice that are below the theory  $\Theta_{reg}$ cannot have regular limits of Woodin cardinals, and moreover, the Woodin cardinals and their limits of such a hod mouse exactly correspond to the Solovay sequence\footnote{By a theorem of Woodin, each $\theta_{\a+1}$ is a Woodin cardinal of $\H$. See \cite{KoelWoodin}.} in the following sense. 
\begin{theorem}[\cite{ATHM} and \rthm{computation of hod}]\label{exp: woodin cardinals}
 In the minimal model of $\Theta_{reg}$, and in fact of $\sf{LSA}$, $\d$ is a Woodin cardinal of $\H$ or a limit of Woodin cardinals of $\H$ if and only if $\d$ is a member of the Solovay sequence. 
 \end{theorem}
\rthm{exp: woodin cardinals} implies that $\H$ of the minimal model of $\Theta_{reg}$ has no Woodin cardinals in the interval $(\theta_\a, \theta_{\a+1})$, and in this sense, $\dagger$ is true below $\Theta_{reg}$\footnote{It is a well-known fact from inner model theory dating back to \cite{IT} that iteration strategies of mice or hod mice acquire complexity only because of Woodin cardinals.}. Therefore, to represent $V_\Theta^\H$ of the minimal model of $\Theta_{reg}$ as a hod mouse, we do not need to understand exactly what happens between $(\theta_\a, \theta_{\a+1})$ in $V$ as none of what happens there makes $\H$ look complicated in that interval\footnote{This was the original motivation of the so-called ``layering" used both in \cite{ATHM} and in this manuscript.}. 

The world of determinacy might have been a simpler place if $\dagger$ was always true, but \cite{ATHM} shows that the theory $\Theta_{reg}$ is much weaker than a Woodin cardinal that is a limit of Woodin cardinals. $\sf{LSA}$, the main topic of this manuscript, is the next natural determinacy theory that is consistency wise stronger than $\Theta_{reg}$, and while the hod mice of this manuscript do have inaccessible limit of Woodin cardinals, \rthm{exp: woodin cardinals} is still true. This once again implies that  the large cardinal structure of hod mice at the level of ${\sf{LSA}}$ is limited and in fact, in such hod mice\\\\
($\ddagger$) there is no Woodin cardinal $\d$ and a $\k<\d$ such that $\k$ is $\d$-strong.\\\\
Moreover, prior to the current work, it was believed that $\ddagger$ and \rthm{exp: woodin cardinals} are just consequences of ${\sf{AD^+}}$. This belief was based on various arguments due to Woodin that showed that if $\d$ is a member of the Solovay sequence then there  cannot be $\k<\d$ whose Mitchell order was much bigger than $\d$. However, \rthm{lst from wlw} shows that ${\sf{LSA}}$ is weaker than a Woodin cardinal that is a limit of Woodin cardinals, and further unpublished work of the first author showed that the large cardinal structure of hod mice, at least in the short extender region, may not be limited. In particular, neither $\dagger$ nor \rthm{exp: woodin cardinals} are consequences of ${\sf{AD}}^+$. The first author then made the following conjecture.
\begin{conjecture}\label{woodin cardinals of hod} Assume ${\sf{AD^+}}+V=L(\powerset(\bR))$. Define the sequence $(\eta_\a: \a\leq \Omega)$ as follows:
\begin{enumerate}
\item $\eta_0=\theta_0$.
\item Assuming $\eta_\a<\Theta$ and setting $\k=(\eta_\a^+)^\H$, $\eta_{\a+1}$ is the supremum of all $\b$ such that there is an ordinal definable surjection $f: \powerset_{\omega_1}(\k)\rightarrow \b$.\footnote{Recall that $\powerset_{\omega_1}(\kappa)$ is the set of countable subsets of $\kappa$.}
\item For a limit ordinal $\xi$, $\eta_\xi=\sup_{\a<\xi}\eta_\a$. 
\end{enumerate}
Then $\d$ is a Woodin cardinal or a limit of Woodin cardinals of $\H$ if and only if $\d=\eta_\a$ for some $\a$. 
\end{conjecture} 
Using the methods of \cite{SteelComp}, Steel verified \rcon{woodin cardinals of hod} assuming ${\sf{HPC}}+{\sf{NLE}}$ (see \cite[Theorem 11.5.7]{SteelComp}). More recently, the first author, using ideas from \cite{SteelComp}, constructed a hod mouse that has a Woodin cardinal that is a  limit of Woodin cardinals. This result confirms the belief that hod mice may have a complicated large cardinal structure. \\\\ 
\textbf{6. Organization.}\\
Chapters 2-8 develop the basic theory of hod mice for $\sf{AD}^+$ models up to the minimal model of $\sf{LSA}$; a consequence of this analysis is (A). The last four chapters focus on applications. Chapter \ref{chap:square} proves that  $\square_{\kappa,2}$ holds in HOD of $\sf{AD}^+$ models up to the minimal model of $\sf{LSA}$ for all HOD-cardinals $\kappa$. Our main use of this chapter is Chapter \ref{chap:lsa_from_pfa}, where a proof of (D) is given. Chapter \ref{chap:condensing_sets} develops the basic theory of \textit{condensing sets}, which is needed in constructions of hod mice in various situations. Chapter \ref{applications} uses the material developed in the previous chapters to prove (B) and (C). The last chapter (Chapter 12) proves (D) by constructing a hybrid version of $K^c$. This chapter uses methods developed in the previous chapters, \cite{SargTrang2021}, and \cite{Trang2015PFA}.

$\sf{Acknowledgments:}$ We are indebted to Dominik Adolf for a long list of very useful corrections. We are also indebted to the referees for a very thorough list of corrections. In particular, the referee for the earlier chapters and Takehiko Gappo have done a tremendous amount of work correcting numerous typos and identifying many hidden mistakes. We thank Paul Larson for the incredible work he has done as the managing editor.

Both authors would like to thank the National Science Foundation for providing financial assistance through Career Award DMS-1352034. Also, the first author's work is funded by the National Science Center, Poland under the Weave-UNISONO call in the Weave program, registration number UMO-2021/03/Y/ST1/00281. The second author would like to thank the National Science Foundation for partial support through NSF Award DMS-1565808 and Career Award DMS-1945592. The first author would like to thank Mathematics Forshungsinstitute Oberwolfach, Germany for hosting him during the winter of 2012 for 7 weeks as a Leibniz Fellow. The first draft of this manuscript was written there. Both authors would like to thank Mathematics Forshungsinstitute Oberwolfach, Germany for hosting them for a week in the winter of 2012. Several of their joint projects go back to this one week. Both authors would like to thank the Newton Institute at Cambridge, UK for hosting them during the Fall of 2015. Several chapters of this manuscript were written during this period. 

Finally, as it must be clear to anyone flipping through the pages of this manuscript, the inspiration behind this work comes from seminal contributions to descriptive inner model theory made by John Steel and Hugh Woodin. We thank them for the monumental work they have done during the past four decades.

\chapter{Hybrid $\mathcal{J}$-structures}

The main goal of this chapter is to prepare some terminology to be used in the rest of this manuscript. One important notion introduced in this chapter is that of the \textit{un-dropping game} (see \rdef{the un-dropping iteration game}). We will use it to prove a comparison theorem for hod mice (see \rcor{comparison holds}). None of the results stated in this chapter are originally due to the authors, though some of them do not appear in literature in exactly the same form that we state here. 

Throughout this book, the reader is assumed to know the basics of inner model theory. Starting from the beginning would have added many more pages to this book, and moreover, the basics of the theory have been developed in several places. The reader is encouraged to review the basic fine structural terminology as presented for example in \cite{ANS}, \cite{SZ} or in \cite{OIMT}. 

\section{$\J$-structures}\label{sec: jstructures}
We say $M=(\univ{M}, Q, \in....)$ is a transitive structure if $\univ{M}$ a transitive set. If $M$ is just a set\footnote{All mathematical objects are sets; here we just mean that $M$ doesn't have any extra structure defined on it.} then we let $\univ{M}=M$\footnote{It seems that this notation is due to Farmer Schlutzenberg.}. In what follows, given a transitive set or a structure $M$  we set $\ord(M)=Ord\cap \univ{M}$. Also, given a set $X$, we let $trc(X)$ be the transitive closure of $X$. We also let $trc^{X}=(trc(X\cup\{X\}), \in)$. 

Recall the inductive definition of $\J_{\omega\alpha}^A(X)$ (for example see \cite[Definition 1.6]{SZ}). In this manuscript, we will also use the round bracket notation while \cite[Definition 1.6]{SZ} only introduces the square bracket notation. We give the definition below, which uses the $rud_A$ function defined in \cite[Definition 1.1]{SZ}. 

\begin{definition}\label{j-structures def}
Suppose $\vec{A}=(A_0, A_1,..., A_n)$ is a finite sequence such that for each $i\leq n$, $A_i$ is a partial set or  class function, and suppose $X$ is a set or a transitive structure. Then set
\begin{align*}
\J_0^{\vec{A}}(X)&=trc(\{X\})\ \text{if $X$ is a set,}\\
\J_0^{\vec{A}}(X)&=trc(\{\univ{X}, Y_0,..., Y_n\})\ \text{if $X=(\univ{X}, Y_0, Y_1, ..., Y_n)$ is a structure,}\\
\J_{\omega\a+\omega}^{\vec{A}}(X)&=rud_{\vec{A}}(\J_{\omega\a}^{\vec{A}}(X)\cup \{\J_{\omega\a}^{\vec{A}}(X)\}),\\
\J_{\omega\lambda}^{\vec{A}}(X)&=\bigcup_{\a<\lambda}\J_{\omega\alpha}^{\vec{A}}(X)\  \text{for limit ordinals $\l$},\\
\J^{\vec{A}}(X)&=\bigcup_{\a\in Ord} \J_{\omega\alpha}^{\vec{A}}(X). 
\end{align*}
$\myqedhere$
\end{definition}

 Recall that a transitive structure $\M=(M, A_1,.., A_k, \in)$ is called \textit{amenable} if for every $X\in M$ and $1\leq i \leq n$, $A_i\cap X\in M$. Following \cite{Zeman}, we say $\M$ is a \textit{$\mathcal{J}$-structure over $X$} if $\M$ is an amenable structure, and 
 \begin{center}
 $\M=(\univ{\mathcal{J}_{\omega \a}^{\vec{A}}(X)}, \vec{A}\cap \mathcal{J}_{\omega \alpha}^{\vec{A}}(X), B_0,..., B_m, X, \in)$
 \end{center}
 where for any set $M$, $\vec{A}\cap M=(A_0\cap M, ..., A_n\cap M)$.

 We think of $\vec{B}=(B_0, \dots, B_m)$ as a sequence of predicates. We  will usually just need three such predicates, one for the last extender, one for the last branch and one for the set of layers to be defined later. At most one of $B_0$ and $B_1$ will be non-empty. $X$ and its predicates (if there are any) are treated as constants. Thus, the language of $\J$-structures is the language of set theory augmented by infinitely many relation symbols and infinitely many constant symbols\footnote{We do not need infintely many such symbols but a large finite number of them.}. As we said above, most cases that will come up in this book will only have three predicates. $X$ usually will itself be a $\J$-structure. 
 
 It is often convenient to think of $A_i$ as a partial function $A_i: \univ{\M}\rightarrow \univ \M$ rather than some larger external function. Notice that for any $\vec{A}$, $\J_{\omega\a}^{\vec{A}}(X)=\J_{\omega \a}^{\vec{A}\cap \J_{\omega \a}^{\vec{A}}(X)}(X)$.
 
 \begin{definition}\label{hierarchical}
 Suppose $\M=(\J_{\omega\a}^{\vec{A}}(X), \vec{A}, B_0, ...B_n, X, \in)$ is a $\J$-structure with $\vec{A}=(A_0, ..., A_n)$. We say $\M$ is \textbf{hierarchical} if the following clauses hold:
  \begin{enumerate}
  \item $\M$ is amenable.
 \item For every $i\leq n$, $\dom(A_i)\subseteq \{\omega\b: \b <\a\}$.
 \item For every $i\leq n$ and for every $\b<\a$ such that $\omega\b\in \dom(A_i)$, $A_i(\omega\b)\subseteq \univ{\mathcal{J}_{\omega\b}^{\vec{A}}(X)}$.
 \item The structure $(\univ{\mathcal{J}_{\omega\b}^{\vec{A}}(X)}, A_0(\omega\b), ..., A_n(\omega))$ is amenable.
\end{enumerate}
$\myqedhere$
%
%
 \end{definition}
 
 The intuition behind a hierarchical structure is that the objects indexed at \textit{active stages}\footnote{Here we say that $\omega\b$ is an active stage of for $A_i$ if $\omega\b\in \dom(A_i)$.} are amenable subsets of the model up to that stage.
 Often hierarchical structures are not represented in this fashion. For example, if $\vec{E}$ is a fine extender sequence (see \cite[Definition 2.4]{OIMT}) then intuitively $\mathcal{J}^{\vec{E}}$ is a hierarchical structures in the above sense, but in reality one needs to use the amenable code (see  \cite[Lemma 2.9]{OIMT}) of each of the extenders in $\vec{E}$ in order to obtain a hierarchical structure. In this book, to avoid making things even more technical than they are, we will simply let the strategy predicates index the iterations and their branches. Thus, if $A_i$ corresponds to the strategy predicate then according to our definition (see \rdef{amenable function}) $A_i$ will be represented as a strategy rather than a function whose domain consists of ordinals. However, it is a simple matter to re-design our hybrid structures so that they fit into our definition of hierarchical. In this book, all our $\mathcal{J}$ structures can be easily represented as hierarchical $\mathcal{J}$-structures. 
 
 Suppose now that $\M=(\J_{\omega\a}^{\vec{A}}(X), \vec{A}, B_0, ...B_n, X, \in)$ is a hierarchical $\J$-structure and $\omega\b< \ord(\M)$. We then set
\begin{center}
 $\M|\omega\b=(\J_{\omega\b}^{\vec{A}}(X), \vec{A}\cap \J_{\omega\b}^{\vec{A}}(X), X, \in)$.
 \end{center}
 and 
 \begin{center}
 $\M||\omega\b=(\J_{\omega\b}^{\vec{A}}(X), \vec{A}\cap \J_{\omega\b}^{\vec{A}}(X), A_0(\omega\b), A_1(\omega\b),..., A_n(\omega\b), X, \in)$.
 \end{center}
 Thus, $\M||\ord(\M)=\M$, and $\M|\omega\b$ is $\M$ ``up to" $\omega\b$ and $\M||\omega\b$ is $\M$ ``up to and including" $\omega\b$. Below we will say that $A_i(\omega\b)$ is \textit{indexed} at $\omega\b$. 
 \begin{remark}\label{omegaalpha} Thus, $\M|\gg$ and $\M||\gg$ are defined only when $\gg=\omega\alpha$ for some $\a$. $\myqedhere$
 \end{remark}
 
We say $X$ is \textit{self-well-ordered}\index{self-well-ordered} if there is a wellordering of $\univ X$ in $\mathcal{J}_1(X)$ definable over $\J_0(X)$ using only the predicates of $X$ as parameters. For example, if $X$ is a premouse then $\vec{E}^X$ is allowed to be used. Unless indicated otherwise, all our $\J$-structures will be over self-well-ordered sets. If $\M$ is a $\J$-structure then we let $X^\M$ be the $X$ above. It follows that each $\J$-structure has a canonical well-ordering given we fix a recursive enumeration of formulas. And so in what follows, we will assume that such an enumeration has been fixed and hence, every $\J$-structure, unless otherwise indicated, has a canonical well-ordering which we will denote by $<_\M$\footnote{$<_\M$ depends on the well-ordering of $X$ that is definable over $\J_0(X)$, and there can be many such well-orderings. It doesn't matter for us which of them is chosen, but one could simply take the well-ordering that is definable via the least formula that defines a well-ordering of $X$ over $\J_0(X)$.}.  We then must have that for $\b<\a$, $<_{\M|\omega\b}=<_{\M|\omega\a}\rest \univ{\M|\omega\b}$.

 \section{Some fine structure}\label{fine structure: sec}

The goal of this section is to review some fine structural ideas. It is not our goal to develop fine structure, but only import some of the standard terminology that is developed in the literature. It is important to note that while new ideas and concepts do appear in the definition of a short-tree-strategy mouse, no new fine structural issues arise. All such fine structural issues have been handled elsewhere, and so we will not dwell on them. The reader unfamiliar with fine structural issues is advised to review some of the following sources:  \cite{ANS}, \cite{FSIT}, \cite{NeemanSteelSub}, \cite{SZ}, \cite{SchSteelZee}, \cite{trang2013}, \cite{FarmerRes}, \cite{OIMT}, \cite{NormIter},   \cite{Zeman}. Our fine structural set up will follow \cite{NeemanSteelSub} and \cite{SZ}. 

We say $\M$ is an \textit{acceptable $\mathcal{J}$-structure} if $\M$ is a $\mathcal{J}$-structure and for all $\tau$ and for all $\b$ such that
 $\ord(\M|0)\leq \tau$ and $\tau< \omega\b$, if $\powerset(\tau)\cap \M|\omega(\b+1)\not \subseteq \M|\omega\b$ then there is a surjection $f:\tau\rightarrow \omega\b$ in $\M|\omega(\b+1)$\footnote{For now, we will need this concept only for $\M$ with $X^\M$ self-well-ordered.}.
 
 \begin{remark}
  From now on all $\J$-structures we will consider will be assumed to be acceptable and hierarchical.$\myqedhere$
  \end{remark}

Suppose $\M$ is a $\J$-structure (over a self-well-ordered set $X$). We then let $\rho_1(\M)$, the $\Sigma_1$ projectum of $\M$, be the least $\rho\leq \ord(\M)$ such that for some $p\in (\ord(\M)^{<\omega})$ and some $\Sigma_1$ formula $\phi$\footnote{In the language of $\J$-structures.} the set $A=\{\xi<\rho: \M\models \phi[\xi, p]\}$ is not in $\M$. The $\Sigma_1$ standard parameter of $\M$, $p_1(\M)$, is the least\footnote{With respect to the lexicographic order on decreasing sequences of ordinals.} $p$ as above. The $\Sigma_1$-reduct of $\M$ is the $J$-structure $(\M||\rho, T)$ where $T$ codes the $\Sigma_1$ theory of $\M$ with parameters in $\rho_1(\M)\cup\{p_1(\M)\}$. The $\Sigma_1$ core of $\M$, $\C_1(\M)$, is the transitive collapse of the $\Sigma_1$ Skolem hull in $\M$ of 
\begin{center}
$\rho_1(\M)\cup \{p_1(\M)\}\cup X^\M\cup\{X^\M\}$. 
\end{center}We say $\M$ is 1-sound if $\C_1(\M)=\M$ and $\M$ is 1-solid. The definition of solidity appears in \cite[Definition 2.15]{OIMT} or in \cite[Definition 7.5]{SZ}. 

\begin{definition}\label{fine structural j-structure} We say $\M$ is a \textbf{fine structural $\mathcal{J}$-structure} (f.s. $\J$-structure) if $\M=(\M', k)$ where $\M'$ is a $\J$-structure, $k\leq \omega$ and letting $(\M_i: i\leq k)$ be given by 
\begin{enumerate}
\item $\M_0=\M'$ and 
\item for $i+1\leq  k$, $\M_{i+1}$ is the $\Sigma_1$ reduct of $\M_i$,
\end{enumerate}
then for all $i<k$, $\M_i$ is 1-sound. We say that $(\M_i: i\leq k)$ is the reduct sequence (r-sequence) of $\M$, and set 
\begin{enumerate}
\item $\rho_0(\M)=\ord(\M)$ and $p_0(\M)=\emptyset$,
\item for $i\leq k$, $\rho_{i+1}(\M)=\rho_1(\M_{i})$ and $p_{i+1}(\M)=p_i(\M)^\frown p_1(\M_i)$,
\item $\rho(\M)=\rho_{k+1}(\M)$ and $p(\M)=p_{k+1}(\M)$.
\end{enumerate}
We also say that $\M'=_{def}j(\M)$ is the $\J$-component of $\M$ and $k=_{def}k(\M)$ is the f.s.-component of $\M$, and set  $l(\M)=(\ord(\M), k)$. Finally, we say $\M$ is sound if $\M_k$ is 1-sound. 

We also say $\M=(\M', \omega)$ is a f.s. $\J$-structure in case $(\M', k)$ is a f.s. $\J$-structure for all $k<\omega$. In this case, $\rho(\M)$ is the eventual value of $\ord(\M_i)$ for $i<\omega$.$\myqedhere$
\end{definition}

Suppose now that $\M=(\M', k)$ is a f.s. $\J$-structure and $E$ is an $\M$ extender such that $\cp(E)<\rho_k(\M)$. Let $(\M_i: i\leq k)$ be the r-sequence of $\M$. We then let $Ult(\M, E)$ be the f.s. $\J$-structure whose $\J$-component is obtained by decoding $Ult_0(\M_k, E)$. We also have a map $\pi_E: \M\rightarrow Ult(\M, E)$ which is a $k$-embedding. The reader can review the relevant notions by consulting \cite[Chapter 2]{ANS}, \cite[Chapter 3 and 4]{SZ}, \cite[Definition 2.8]{NormIter} and \cite[Section 2.5]{NormIter}. 

Suppose $\M=(\M', k)$ is a f.s. $\J$-structure and $(\omega\a, m)\leq l(\M)$ (here $\leq$ is the lexicographical order). We then let $\M|(\omega\a, m)=(\M'|\omega\a, m)$ and $\M||(\omega\a, m)=(\M'||\omega\a, m)$. Also, we write $\N\insegeq \M$ if for some $(\omega\a, m)\leq l(\M)$, $\N=\M|(\omega\a, m)$ or $\N=\M||(\omega\a, m)$. We will often write $\M|\gg$ or $\M||\gg$ for $\M'|\gg$ and $\M'||\gg$.\footnote{Notice that our definitions do not guarantee that $\M|\gg$ or $\M||\gg$ are f.s $\J$-structures. However, the structures that we will eventually consider will have this property.} 

The next definition defines the core of a $\mathcal{J}$-structure. One way of defining it is by doing what is described after \cite[Definition 7.13]{SZ}. Here is an outline of essentially that same construction.

\begin{definition}\label{the core} Suppose $\M$ is a $\J$-structure. We define $(\C_k(\M): k< \omega)$, $(\rho_k(\M): k< \omega)$ and $(p_k(\M): k<\omega)$ by induction as follows. 
\begin{enumerate}
\item Set $\C_0(\M)=\M$.
\item If $\M$ is not 1-solid then let $\C_k(\M)$ for $k\geq 1$ be undefined. Otherwise, $\C_1(\M)$ is defined as above.
\item Suppose $\C_k(\M)$ has been defined and that $\N=(\C_k(\M), k)$ is a f.s. $\J$-structure\footnote{Notice that this condition is simply part of the induction. Above, we have that $(\C_1(\M), 1)$ is an f.s. $\J$-structure.}. Let $(\N_j: j\leq k)$ be the $r$-sequence of $\N$. If $\N_k$ is not 1-solid  then let $\C_{i}(\M)$ for $i\geq k+1$ be undefined. Otherwise, letting $\pi: \C_1(\N_k)\rightarrow \N_k$ be the core map, we let $\C_{k+1}(\M)$ be the decoding of $\C_1(\N_k)$\footnote{The decoding process is similar to the Downward Extension of Embeddings Lemma (see \cite[Lemma 3.3]{SZ}). The decoding gives a $\Sigma_1$-map $\pi': \C_{k+1}(\M)\rightarrow \C_k(\M)$ extending $\pi$.}.
\end{enumerate}
If $\C_k(\M)$ is defined for all $k<\omega$, then let $\C(\M)$ be the eventual value of $\C_k(\M)$. 

Suppose for some $k<\omega$, $\C_k(\M)$ and $p_k(\M)$ have been defined. Then letting $(\N_j: j\leq k)$ be the $r$-sequence of $\C_k(\M)$, set $\rho_{k+1}(\M)=\rho_1(\N_k)$ and $p_{k+1}(\M)=p_k(\M)^\frown p_{1}(\N_k)$. Let $\rho(\M)$ be the eventual value of the sequence $(\rho_k(\M): k<\omega)$ and let ${\sf{ep}}(\M)$ be the least $k$ such that for all $i\geq k+1$, $\rho_i(\M)=\rho_{k+1}(\M)$. $\myqedhere$
\end{definition}

Thinking of $\mathcal{J}$ structures as f.s. $\mathcal{J}$-structures is useful in introducing iteration trees and in the proof of convergence of $K^c$-constructions (for example see \cite{ANS}). 

 \section{Layered hybrid $\mathcal{J}$-structures}\label{layered sec}
 
We say $w$ is a \textit{sequential structure} if $w=(\J_\omega(s), s, \in)$ where $s$ is a sequence $(u_\a: \a<\gg^w)$. 
 
 \begin{definition}[Definition 1.1 of \cite{ATHM}]\label{amenable function}
Given a function $f$, we say $f$ is \textbf{amenable} if for all $w\in \dom(f)$
\begin{enumerate}
\item $w$ is a sequential structure,
\item $f(w)\subseteq \ord(w)$,
\item $\sup f(w)=\gg^w$ and $0\in f(w)$, 
\item  whenever $\eta<\gg^w$, $f(w)\cap \eta\in w$. 
\end{enumerate}
$\myqedhere$
\end{definition}

We say $f$ is a \textit{shift of an amenable function}\index{shifted amenable function} or a \textit{shifted amenable function} if there is an amenable function $g$ with $\dom(g)=\dom(f)$ and such that for all $w\in \dom(f)$, 
\begin{enumerate}
\item $f(w)\subseteq Ord$,
\item $f(w)\subseteq [\min(f(w)), \min(f(w))+\gg^w)$, and 
\item $f(w)=\{ \min(f(w))+\omega\gg : \gg\in g(w)\}$. 
\end{enumerate}
Notice that if $f$ is a shift of an amenable function then it uniquely determines $g$. We say that $g$ is the \textit{amenable component} of $f$.

Jumping ahead, we remark that iteration strategies and mouse operators provide an ample source of amenable functions. For instance, let $\M=\M_1^{\#}$ and let $\Sigma$ be its canonical iteration strategy. We define $f$ as follows. Let first $\dom(f)$ be the set of structures of the form $w=(\mathcal{J}_{\omega}(\T^w), \T^w, \in)$ where $\T^w$\footnote{One could think of $\T$ as a sequence $(\T\rest \a: \a<\lh(\T))$.} is a normal iteration tree on $\M$ of limit length and is according to $\Sigma$. Next, define $f(w)=b$ where $b=\Sigma(\T^w)$. Then $f$ is amenable. We will refer to such an $f$ as \textit{an amenable function given by an iteration strategy}. 

The definitions that follow explain how our indexing schemes work. We first isolate those iterations whose branches will be indexed. The reader may think of the formula $\phi$ appearing in \rnot{smsphi notation} as the formula that defines the set of iterations whose branches need to be indexed. However, $\phi$ alone does not define such iterations as we need to add clause 2 for technical reasons. The ordinal $\b$ essentially identifies the location where the branch of the iteration tree defined by $\phi$ should be indexed.  

In general, to develop a reasonable theory of hybrid $\mathcal{J}$-structures, we need to use indexing schemes to index branches of stacks. The reason for this is that if no particular coherent method of indexing is used to organize such structures then one cannot in general hope to develop a comparison theory for the resulting structures. Indeed, if $\M$ and $\N$ are unindexed hybrid $\mathcal{J}$-structures then it is possible that some $b$ is indexed at $\a$ in $\M$ but nothing is indexed at $\a$ in $\N$, causing a fatal breakdown of the comparison argument. Nevertheless, in \rsec{sec short tree strategy mice}, it will be convenient to work with unindexed  hybrid $\mathcal{J}$-structures (as defined in \rdef{unindexed ses}).

Another important remark is that it might be convenient to think of the indexing scheme as a parameter of the hybrid $\mathcal{J}$-structures, in the same way we internalize the fine structural parameter. Thus, instead of $\M$ we could consider $(\M, \phi)$ where $\phi$ is the indexing scheme. However, in most cases, isolating $\phi$ won't matter so much, and so we will not take this path.
 We suspect that $\M$ may even recognize many different $\phi$s as its own indexing scheme\footnote{Of course, this can be achieved trivially; $\phi$ and $0=0\wedge \phi$ are equivalent. But there could be two different indexing schemes $\phi$ and $\phi'$ such that ${\sf{ZFC}}$ or any natural extension of it, does not prove $\phi\iff \phi'$ yet there is $\M$ which is both $\phi$ indexed and $\phi'$-indexed.}. Thus we may have two $\phi$-indexed $\M$ and $\M'$ and an indexing scheme $\phi'$ such that $\M'$ is $\phi'$-indexed while $\M$ is not. However, such issues will not come up in the sequel.

\begin{notation}\label{smsphi notation}
 Suppose now that $\M=\mathcal{J}_{\iota}^{A_0, A_1,..., A_n}(X)$ is a $\J$-structure or an f.s. $\J$-structure and $\phi(\vec{x}, u)$ is a formula in the language of $\J$-structures such that it implies that ``$u$ is a sequential structure". Suppose $\vec{s}\in \univ{\M}^{<\omega}$. Let $S^{\M}_{\vec{s}, \phi}$ be the set of pairs $(\b, w)$ such that
 \begin{enumerate}
 \item $\omega\beta+\omega\gg^w\leq \ord(\M)$, 
 \item $\M|\omega\b\models ``\cf(\gg^w)$ is not a measurable cardinal as witnessed by extenders in $A_0$"(see \rrem{measurable cofinality issue}), and 
 \item $\M|\omega\b\models {\sf{ZFC}}+\phi[\vec{s}, w]$. 
 \end{enumerate}
 Let $\sf{nmc}(\a)$ be the statement ``$\cf(\a)$ is not a measurable cardinal as witnessed by the extenders in $A_0$". $\myqedhere$
 \end{notation}

 \begin{definition}\label{important notation}  Suppose that $(\M, \vec{s}, \phi)$ are as in \rnot{smsphi notation}. 
 Suppose further that $f$ is a shifted amenable function with amenable component $g$ such that $\dom(f)\subseteq \univ{M}$ and for all $w\in \dom(f)$, $\min(f(w))+\gg^w\leq \ord(\M)$\footnote{Recall our convention that $X^\M$ is self-well-ordered.}.  We say $w$ is \textbf{weakly $(f, \vec{s}, \phi)$-\textit{minimal}} if there is $\b$ such that
 \begin{enumerate}
 \item $(\b, w)\in S^\M_{\vec{s}, \phi}$ (in particular, because $\M|\omega\b\models {\sf{ZFC}}$, $\omega\b=\b$),
 \item $w\not \in \dom(f\cap \univ{\M|\b})$, 
 \item  $\{ u\in \univ{\M|\b}: u<_{\M|\b}w$ and there is $\xi< \b$ such that $(\xi, u)\in S^\M_{\vec{s}, \phi}\}\subseteq \dom(f\cap \univ{\M|\b})$.
 \end{enumerate}
 We say $w$ is \textbf{$(f, \vec{s}, \phi)$-minimal} if there is $\b$ witnessing that $w$ is weakly $(f, \vec{s}, \phi)$-\textit{minimal} and such that $w$ is the $<_{\M|\omega\b}$-minimal $w'$ which is weakly $(f, \vec{s}, \phi)$-\textit{minimal} as witnessed by $\b$.

 

If $w$ is $(f, \vec{s}, \phi)$-\textit{minimal} then we let $\b^{\M, f, \vec{s}, \phi}_w$ be the least $\b$ witnessing that $w$ is $(f, \vec{s}, \phi)$-\textit{minimal}. In many cases, $(\M, f, \phi)$ will be clear from context and so we will drop it from our notation. If $\vec{s}=\emptyset$ then we drop it from our notation.$\myqedhere$
\end{definition}

 \begin{remark}[The measurable cofinality issue]\label{measurable cofinality issue} The reader unfamiliar with\\ strategic mice may find clause 2 of \rdef{passive hybrid j-structure} somewhat odd. This clause has to do with an issue known to experts and was first discovered in earlier versions of \cite{CMI}\footnote{The authors were unable to locate the discussion involving the measurable cofinality issue in \cite{CMI}.}. The problem was fully treated in \cite{trang2013}, and the discussion appears in \cite[Remark 2.47]{trang2013}. Without getting too much into the technical details, the issue is simply that if $\M=\mathcal{J}^{A, f}$ is a $\J$-structure such that the $f$ predicate codes a strategy for some $\N\in \M$, $w=_{def}(\J_\omega(\T), \T, \in) \in \dom(f)$, $\k=_{def}\cf^\M(\lh(\T))$ is a measurable cardinal in $\M$, $f(w)$ is indexed at $\l$ and $E\in \vec{E}^\M$ is an extender with $\cp(E)=\k$ then
\begin{center}
 $\sup(\pi^{\M||(\l, 0)}_E[\lh(\T)])<\pi_E(\lh(\T))$ while $\l=\ord(Ult(\M||(\l, 0), E))$.
 \end{center}
  The issue is hiding in the fact that in most of the natural attempts to organize strategic mice, $\cf^\M(\l)=\k$, while in the above situation this fails in $Ult(\M||(\l, 0), E)$ for $\pi_E^{\M||(\l, 0)}(\T)$. $\myqedhere$
\end{remark}

In general, there may not be a unique $w$ which is $(f, \vec{s}, \phi)$-\textit{minimal}. However, the following holds.

\begin{lemma}\label{order of min objects} Suppose $(f, \vec{s}, \phi)$ and $\M$ are as in \rdef{important notation}. Suppose $w\not =w'$ are two $(f, \vec{s}, \phi)$-\textit{minimal} sets. Set $\b=\b^{\M, f, \vec{s}, \phi}_w$ and $\b'=\b^{\M, f, \vec{s}, \phi}_{w'}$, and suppose that $\b\leq \b'$. Then $\b<\b'$.
\end{lemma}


\begin{remark}\label{idea behind indexing} The $f$ of \rdef{important notation} is designed to code an iteration strategy, and it will be the strategy predicate of a  hybrid $\mathcal{J}$-structures which indexes an iteration strategy. The iterations that will get indexed are exactly the $(f, \vec{s}, \phi)$-minimal ones, and \rlem{order of min objects} implies that $(f, \vec{s}, \phi)$-minimal $w$'s are well-ordered. We will then use the function $w\mapsto \b_w$ to index the branch of $w$ at $\b_w+\omega\gg^w$. $\myqedhere$
\end{remark}

\begin{remark} It is perhaps illuminating to figure out the least iteration tree whose branch will be indexed by the predicate $f$ of \rdef{important notation}. We assume $\vec{s}=\emptyset$. First we pick the least $\b$ such that for some $\T$ of length $\omega$, $\M|\omega\b\models {\rm{ZFC}}+\phi[\T]$. Then we take the $<_{\M|\omega\b}$-least $\T$ as above and index its branch at $\omega\b+\omega ^2$. $\myqedhere$
\end{remark}

 We are now in a position to introduce the  \textit{passive hybrid $\mathcal{J}$-structures}. 

\begin{definition}[Passive Hybrid $\mathcal{J}$-structures]\label{passive hybrid j-structure}\index{passive hybrid $\mathcal{J}$-structures} We say $\M$ is a \textbf{passive hybrid $\mathcal{J}$-structure} over a self-well-ordered set $X$ with indexing scheme $\phi(x)$\footnote{$\phi$ is in the language of $\J$-structures.} if $\M=(\M', k)$ is an f.s. $\mathcal{J}$-structure such that the following conditions hold.
\begin{enumerate}
\item  For some $\a$, $A \subseteq \univ{\M'}$ and $f\subseteq \univ{\M'}$, \begin{center}
$\M'=(\mathcal{J}_{\omega\a}^{A, f}(X), A, f, X, \in)$\footnote{We would like to emphasize that $\M'$ has only the displayed predicates. Also, below $(\M', f, \phi)$ are omitted from $\b_w$ notation.},
\end{center}
\item $f$ is a shift of an amenable function.
\item For all $w\in \univ{\M'}$, $w\in \dom(f)$ if and only if  $w$ is $(f, \phi)$-minimal and $\b_w+\omega\gg^w<\ord(\M)$\footnote{Here $\b_w$ is defined in \rdef{important notation}.}.
\item For all $w\in \dom(f)$, 
\begin{enumerate}
\item $\b_w=\min(f(w))$,
\item $\univ{\M'|(\b_w+\omega\gg^w)}=\J_{\b_w+\omega\gg^w}(\M'||\omega\b_w)$ and $A\cap \univ{\M'|(\b_w+\omega\gg^w)}=A\cap \univ{\M'|\omega\b_w}$\footnote{It also follows that $f\cap \univ{\M'|(\b_w+\gg^w)}=f\cap \univ{\M'|\b_w}$.}.
\end{enumerate}
\end{enumerate}
$\myqedhere$
\end{definition}

\begin{remark}\label{eventual indexing of w} \rdef{passive hybrid j-structure} leaves open one important question. Does it follow that $S^\M_\phi= \dom(f^\M)$? The answer is of course that none of the conditions we have imposed on $f^\M$ guarantees that $S^\M_\phi=\dom(f^\M)$. It could be that some $w\in S^\M_\phi$ but it is not $(f, \phi)$-minimal. However, if $f^\M$ is supposed to code an iteration strategy $\Sigma$ of some $P$ then the fact that $w\in S^\M_\phi$ implies that $w$ is an iteration according to $\Sigma$ and that we must have that $w\in \dom(f^\M)$. What will in fact happen, in intuitive terms, is that while $w$ may not be in $\dom(f^\M)$, it will be in $\dom(f^\N)$ for some $\N$ extending $\M$. This may not be possible to arrange if for example $\N=_{def}\J_\omega(\M)$ projects in a way that say $w$ is no longer in $\C_1(\N)$, but if $\M$ is the final model of some reasonable fully backgrounded construction that produces hybrid premouse then we will indeed have that $S^\M_\phi=\dom(f^\M)$. This is because it can be shown that any $w\in S^\M_\phi$ is $(f, \phi)$-minimal.

To see this in intuitive terms, suppose towards a contradiction that some $w\in S^\M_\phi$ doesn't belong to $\dom(f^\M)$. We can assume $w$ is $<_\M$-minimal. Now, and this depends on our choice of $\phi$, the indexing scheme $\phi$ will be $\Sigma_1$ and hence, it will have the following upward absoluteness property: if for some $\nu$, $\M|\nu\models \phi[w]$ then for all $\nu'\geq \nu$, $\M|\nu'\models \phi[w]$. Let $\xi<\ord(\M)$ be such that whenever $u\in S^\M_\phi$ and $u<_{\M} w$ then $\sup(f(u))<\xi$. Such $\xi$ will exist because, in concrete applications, $\ord(\M)$ will be a Woodin  cardinal of the universe. It follows that to show that $w\in \dom(f^\M)$ it is enough to show that there is $\b\in (\xi, \ord(\M))$ such that $\M|\b\models {\rm{ZFC}}$. Since $\ord(\M)$ is a Woodin cardinal of the background universe, there are plenty of such $\b$. $\myqedhere$
\end{remark}

\begin{definition}[Hybrid $\mathcal{J}$-structures]\label{hybrid j-structure}\index{hybrid $\mathcal{J}$-structures} We say $\M$ is a \textbf{hybrid $\mathcal{J}$-structure} over a self-well-ordered set $X$ with indexing scheme $\phi(x)$ if $\M=(\M', k)$ is an f.s.  $\mathcal{J}$-structure such that 
\begin{enumerate}
\item  for some $\a$,  $A \subseteq \univ{\M'}$ and $f\subseteq \univ{\M'}$, \begin{center}
$\M'=(\mathcal{J}_{\omega\a}^{A, f}(X), A, f, B, F, X, \in)$\footnote{Below $(\M', f, \phi)$ are omitted from $\b_w$ notation.},
\end{center}
\item $(\mathcal{J}_{\omega\a}^{A, f}(X), A, f, X, \in)$ is a passive hybrid $\J$-structure, 
\item at most one of $B$ and $F$ is not empty,
\item if $F\not=\emptyset$ then $F$ is an ordered pair $(w, b)$ such that if $\b=min(b)$ then setting $f'=f\cup\{(w, b)\}$,
\begin{enumerate}
\item $f'$ is a shift of an amenable function\footnote{This implies that $w$ is a sequential structure.},
\item $w$ is $(f', \phi)$-minimal with $\b^{\M, f', \phi}_w=\b$ (in particular, $\omega\b=\b$, see \rdef{important notation}),
\item $\omega\a=\b+\omega\gg^w$,\footnote{It follows from clause 5 of \rdef{important notation} that $\M'\models ``\cf(\gg)$ is not a measurable cardinal as witnessed by extenders in $A$".}
\item $\univ{\M'}=\J_{\b+\omega\gg^w}(\M'||\b)$ and $A\cap \univ{\M'}=A\cap \univ{\M'|\b}$.
\end{enumerate}
\end{enumerate}
For $w\in \dom(f')$, we say that $f'(w)$ is indexed at $\b_w+\omega\gg^w$ or that $\b_w+\omega\gg^w$ is the index of $f'(w)$.$\myqedhere$
\end{definition}

Suppose $\M$ is a hybrid $\mathcal{J}$-structure with an indexing scheme $\phi$. We will often say that ``$\M$ is indexed according to $\phi$" or that ``$\M$ is $\phi$-indexed". Notice that only the $f$ predicate is indexed according to $\phi$. In most situations that we will consider $A$ will be an extender sequence. Sometimes, however, we will need to consider cases where there are two or more $f$ predicates.

%

\begin{rem} Notice that it follows from  \rdef{hybrid j-structure} that the function $a\mapsto min(f(a))$ is injective on $\dom(f)$.  $\myqedhere$
\end{rem}

 Hod mice are a special blend of \textit{layered hybrid $\mathcal{J}$-structures} introduced below. 
 
\begin{definition}[Passive layered hybrid $\mathcal{J}$-structure]\label{passive layered hybrid j-structure}\index{passive layered hybrid $\mathcal{J}$-structures} We say $\M$
 is a \textbf{passive layered hybrid $\mathcal{J}$-structure}\index{passive layered hybrid $\mathcal{J}$-structure} over a self-well-ordered set $X$ with indexing scheme $\phi(x, y)$ if $\M=(\M', k)$ is an f.s. $\J$-structure such that
 \begin{enumerate}
 \item for some $\a$, $A\subseteq \univ{\M'}$ and $f\subseteq \univ{\M'}$, 
 \begin{center}
 $\M'=(\J_{\omega\a}^{A, f}(X), A, f, Y, X, \in)$,
 \end{center} 
 \item $Y\subseteq \J_{\omega\a}^{A, f}(X)$ and $dom(f)= Y \subseteq \{ \Q: \Q\inseg \M\}\cup X$,
 \item for all $\Q\in \dom(f)$, $f(\Q)=_{def}f_\Q$ is a shift of an amenable function\footnote{Below we will drop $(\M', f_\Q, \phi)$ from the $\b_w^{\M', f_\Q, \Q, \phi}$ notation.},
 \item for all $w\in \M'$ and $\Q\in Y$, $w\in \dom(f_\Q)$ if and only if  
 \begin{enumerate}
 \item  $w$ is $(f_\Q, \Q, \phi)$-minimal and $\b_{w}^{\Q}+\omega\gg^w<\omega\a$, and
 \item for all $\R\in Y$ such that $\R<_{\M'} \Q$ and for all $(f_\R, \R, \phi)$-minimal $u\in \M'|\omega\b^\Q_w$, $u\in \dom(f_\R\cap \univ{\M'|\b^\Q_w})$\footnote{Implying that $\b^\R_u<\b^\Q_w$.},
 \end{enumerate}
 
\item for all $\Q\in Y$ and for all $w\in \dom(f(\Q))$, 
\begin{enumerate}
\item $\b_w^\Q=\min(f(w))$,
\item $\univ{\M'|(\b_w^\Q+\omega\gg^w)}=\J_{\b_w^\Q+\omega\gg^w}(\M'||\b_w^\Q)$ and $A\cap \univ{\M'|(\b_w^\Q+\omega\gg^w)}=A\cap \univ{\M'|\omega\b_w^\Q}$\footnote{It also follows that $f\cap \univ{\M'|(\b_w^\Q+\gg^w)}=f\cap \univ{\M'|\omega\b^\Q_w}$.}.
\end{enumerate}
\end{enumerate}
$\myqedhere$
\end{definition}

\begin{definition}[Layered hybrid $\mathcal{J}$-structure]\label{layered hybrid j-structure}\index{layered hybrid $\mathcal{J}$-structures} We say $\M$
 is a \textbf{layered hybrid $\mathcal{J}$-structure}\index{passive layered hybrid $\mathcal{J}$-structure} over self-well-ordered set $X$ with indexing scheme $\phi(x, y)$ if $\M=(\M', k)$ is an f.s. $\J$-structure such that
 \begin{enumerate}
 \item for some $A\subseteq \univ{\M'}$ and $f\subseteq \univ{\M'}$, 
 \begin{center}
 $\M'=(\J_{\omega\a}^{A, f}, A, f, Y, B, F, X, \in)$,
 \end{center}
 \item $\M'=(\J_{\omega\a}^{A, f}, A, f, Y, X, \in)$ is a passive layered hybrid $\J$-structure over $X$,
 \item only one of $B$ and $F$ is non-empty,
\item if $F\not=\emptyset$ then $F$ is an ordered pair $(\Q, (w, b))$ such that $\Q\in Y$, $b\subseteq \ord(\M)$, and if $\b=min(b)$ then setting $f'=f\cup\{(\Q, (w, b))\}$,
\begin{enumerate}
\item $f'(\Q)$ is a shift of an amenable function,
\item $w$ is $(f',\Q, \phi)$-minimal with $\b^{\M, f', \Q, \phi}_w=\b$,
\item $\a=\b+\omega\gg^w$,
\item $\univ{\M'}=\J_{\b+\omega\gg^w}(\M'||\b)$ and $A\cap \univ{\M'}=A\cap \univ{\M'|\b}$,
\item  for all $\R\in Y$ such that $\R<_{\M'} \Q$ and for all $(f'_\R, \R, \phi)$-minimal $u\in \M'|\b$, $u\in \dom(f_\R\cap \univ{\M'|\b})$\footnote{Implying that $\b^\R_u<\b$.}.
\end{enumerate}
\end{enumerate}
$\myqedhere$
\end{definition}

Suppose $\M$ is a layered hybrid $\mathcal{J}$-structure with an indexing scheme $\phi$. We will often say that ``$\M$ is indexed according to $\phi$" or that ``$\M$ is $\phi$-indexed".

We will often omit $\phi$ when discussing a particular layered hybrid $\mathcal{J}$-structure. 
If $\M$ is a layered hybrid $\mathcal{J}$-structure then we let $f^\M$ and $Y^\M$ be as in \rdef{passive layered hybrid j-structure}. We again have that for each $\Q\in Y^\M$,  the function $a\mapsto min(f^\M(\Q)(a))$ is injective on $\dom(f(\Q))$. 

Notice that hybrid $\mathcal{J}$-structures can be viewed as a special case of layered hybrid $\mathcal{J}$-structures. Because of this, in the sequel we will only establish terminology for layered hybrid $\mathcal{J}$-structures though we might use the same terminology for hybrid $\mathcal{J}$-structures. 

Typically, when discussing hybrid $\mathcal{J}$-structures, $X$ will be an \textit{iterable} structure and $f$ will be the predicate coding its strategy.\footnote{\label{fn:cofMea} In this case, the $\gamma$ defined in Definition \ref{hybrid j-structure} is the length of a tree $\mathcal{T}$ according to $f$. The condition ``$\M \vDash$ cof$(\gamma)$ is not measurable" in Definition \ref{hybrid j-structure} ensures that fine structure is preserved under iterations.} 

As mentioned above, hod mice are a special type of layered hybrid $\mathcal{J}$-structures: the $f$ predicate of a hod mouse codes a strategy for its layers. When the $A$ predicate of a layered hybrid $\mathcal{J}$-structure is a coherent sequence of extenders then the resulting model is called a \textit{hybrid layered premouse}\index{hybrid layered premouse}. 

Results of this manuscript are independent of particular extender-indexing schemes, but for technical reasons we will use a mixture of Jensen indexing as developed in \cite[Definition 2.4]{ANS}, \cite{JensenBook}, \cite{NeemanSteelSub}, \cite{Zeman} and Mitchell-Steel indexing as developed in \cite{FSIT} and \cite[Definition 2.4]{OIMT}. Suppose $\M=\mathcal{J}^{\vec{E}, f}(X)$ is a $\phi$-indexed layered hybrid $\mathcal{J}$-structure over a self-well-ordered set $X$ and $\vec{E}$ is a sequence of extenders. We say $\eta$ is a cutpoint of $\M$ if there is no $\a\in \dom(\vec{E}^\M)$ such that $\cp(\vec{E}^\M(\a))<\eta\leq \a$. We say $\vec{E}$ is a \textit{mixed indexed extender sequence} if the following clauses hold:
\begin{enumerate}
\item (j-like indexing) If $\k=\cp(E)$ is a limit of Woodin cardinals of $\M$ and is a cutpoint of $\M$ then letting $\N=\pi^\M_E(\M|(\k^+)^\M)$, $E$ is indexed at $(\eta^+)^\N$ where $\eta=\sup\{ \a\in \dom(\vec{E}^{\N}): \cp(\vec{E}^\N(\a))=\k\}$.
\item (ms-indexing) All other extenders are indexed according to the ms-indexing. 
\end{enumerate}
The initial segment condition for $E$ is clause 3 of \cite[Definition 2.4]{OIMT}. 
There are many papers in the literature that connect the two indexing schemes. For example, the reader may consult \cite{Fuchs1}, \cite{Fuchs2}, and \cite{schlutzenberg2021background}. Our goal in this book is to present the theory of the minimal model of the Largest Suslin Axiom in as shorter space as we can, and because of this we will avoid fine structural issues that have been well-treated in the literature.  

\begin{definition}[Layered hybrid e-structure]\label{layered hybrid e-structure}\index{layered hybrid e-structure} Suppose $\M=\mathcal{J}^{\vec{E}, f}(X)$ is a $\phi$-indexed layered hybrid $\mathcal{J}$-structure over a self-well-ordered set $X$. $\M$ is called a \textbf{$\phi$-indexed layered hybrid potential e-structure} ($\sf{lhpes}$)\index{layered hybrid potential premouse} if $\vec{E}$ is a mixed indexed extender sequence. We write $\vec{E}^\M$ for $\vec{E}$ etc.

If $\M$ is an $\sf{lhpes}$ and $E=\vec{E}^\M(\gg)$ then we let $\ind^\M(E)=\gg$.

We say that $\M$ is a $\phi$-indexed layered hybrid e-structure ($\sf{lhes}$)\index{$\sf{lhes}$} if $\M=(\M', k)$ is an f.s. $\J$-structure such that $\M'$ is a $\phi$-indexed $\sf{lhpes}$ and for every $(\omega\b, m)<l(\M)$,  $\M||(\omega\b, m)$ is sound.$\myqedhere$
\end{definition}

Mixed indexing smoothens implementation of certain technical arguments. The most crucial property for us is the following. Suppose $E$ is an extender on the extender sequence of an $\sf{lhes}$ $\M$ such that $\M\models ``\cp(E)$ is a cutpoint and a limit of Woodin cardinals" (i.e. there is no $F\in \vec{E}^\M$ such that $\cp(F)<\cp(E)\leq \ind^\M(F)$). So $E$ has j-like indexing. Let $\gg=\sup\{ \xi: \xi<\ind^\M(E), \xi\in \dom(\vec{E}^\M)$ and $\cp(E_\xi)=\cp(E)\}$. Then $\gg=\sup\{\xi: \xi\in \dom(\vec{E}^{Ult(\M, E)})$ and $\cp(E_\xi)=\cp(E)\}$. The advantage of mixed indexing over other indexing schemes can be seen in \rdef{sec:pitb}.


For an $\sf{lhes}$ $\M$ with just one layer (that is, $|Y^\M|=1$), we say $\M$ is a \textit{hybrid e-structure} ($\sf{hes}$). Next we introduce $\sf{lhes}$ that are internally closed under sharps. We will use such a closure to introduce short tree strategy premice (see \rdef{strategy premouse} and \rdef{sts premouse}).

\begin{definition}[Closed $\sf{lhes}$]\label{closed under sharps} Suppose $\M$ is an $\sf{lhes}$ and $\a\leq \ord(\M)$. Then we say $\M$ is \textbf{closed}  below $\a$ if for all $\b<\a$ there is $\gg\in \dom(\vec{E}^\M)$ such that $\gg<\alpha$ and $\cp(E_\gg^\M)>\b$. We say $\M$ is  closed if $\M$ is closed below $\ord(\M)$.$\myqedhere$
\end{definition}

\section{Iteration trees and stacks}\label{sec: it trees}

 Below we review iteration trees. Our notation is mostly in line with most of the references we have quoted above in \rsec{fine structure: sec}. The only difference is that we incorporated the concepts of a stack of iteration trees into an iteration tree. 

Suppose $\M$ is an $\sf{lhes}$ (or ${\sf{hes}}$). Thus, $\M$ is an f.s. $\J$-structure that has a designated extender sequence $\vec{E}^\M$. For a limit ordinal $\eta$, we let $\eta-1=\eta$. 

\begin{definition}\label{putative it} We say $\T$ is a \textbf{putative iteration tree} on $\M$ if 
\begin{center}
$\T=((\M_\a)_{\a<\eta}, (E_\a)_{\a<\eta-1}, D, R, (\beta_\a, m_\a)_{\a\in R}, T)$
\end{center}
and the following conditions hold.
\begin{enumerate}
\item $T$ is a tree order on $\eta$. \\\\
Let $\T(\a+1)$ be the $T$-predecessor of $\a+1$ and $(\a, \b)_\T$ be the $T$-interval $(\a, \b)$\footnote{Similarly define all other combinations of $(\a, \b)_\T$, like $[\a, \b)_\T$ and etc.}. 
\item For all $\a$ such that $\a+1<\eta$, $\M_\a$ is a well-founded $\sf{lhes}$ (or ${\sf{hes}}$).
\item $R\subseteq \eta-1$, $0\in R$ and for all $\a\in R$ and for all $\b\geq \a$, $\T(\b+1)\geq \a$. 
\item For all $\a\in R$, $(\omega\b_\a, m_\a)\leq l(\M_\a)$.\\

Set $\M_\a'=\begin{cases} \M_\a &: \a\not \in R \vee (\a\in R \wedge \omega\b_\a= \ord(\M_\a))\\
\M_\a||(\omega\b_\a, m_\a)& : \a\in R \wedge \omega\b_\a<\ord(\M_\a)
\end{cases}$
\item $\M_0=\M$.
\item For all $\a+1<\eta$, $E_\a\in \vec{E}^{\M_\a'}$.
\item for all $\a+1<\eta$, setting $\b=\T(\a+1)$ and $\k_\a=\cp(E_\a)$, 
\begin{center}
$\M_\a|(\k_\a^+)^{\M_\a|\ind^{\M_\a}(E_\a)} \insegeq \M_\b'$.
\end{center}
\item for all $\a+1<\eta$, 
\begin{center}
$\M_{\a+1}=Ult(\M_\b'||(\omega\xi_\a, k_\a), E_\a)$ 
\end{center}
where 
\begin{enumerate}
\item $\b=\T(\a+1)$,
\item $\omega\xi_\a\leq \ord(\M'_\b)$ is the largest such that $(\k_\a^+)^{\M_\a|\ind^{\M_\a}(E_\a)}=(\k_\a^+)^{\M_\b'|\omega\xi_\a}$,
\item $k_\a$ is the largest such that $(\omega\xi_\a, k_\a)\leq l(\M_\b')$ and $\cp(E_\a)<\rho_{k_\a}(\M_\b'||(\omega\xi_\a, k_\a))$.
\end{enumerate}
\item $D=\{\a+1<\eta:$ letting $\b=\T(\a+1)$, $(\omega\xi_\a, k_\a)<l(\M_\b)\}$.

Let\begin{center} $\pi^{\T}_{\b, \a+1}=\pi_{E_\a}^{\M'_\b||(\omega\xi_\a, k_\a)}: \M'_\b||(\omega\xi_\a, k_\a)\rightarrow \M_{\a+1}$\end{center} be the ultrapower map and for $\a<_T\gg<\eta$ let $\pi^\T_{\a, \gg}:\M_\a\rightarrow \M_\gg$ be the embedding obtained by compositions.\footnote{Assuming these embeddings can be composed. $\pi^\T_{\a, \gg}$ is defined if and only if $D\cap (\a, \gg]_\T=\emptyset$.}

\item For limit $\l<\eta$, $D\cap (0, \l)_\T$ is finite and letting $\b\in [0, \l)_\T$ be the least such that $D\cap (\b, \l)_\T=\emptyset$, $\M_\l$ is the direct limit of the system $(\M_\gg, \pi_{\gg, \gg'}^\T: \gg<\gg', \gg, \gg'\in [\b, \l)_\T)$ and for $\gg\in [\b, \l)$, $\pi^\T_{\gg, \l}:\M_\gg\rightarrow \M_\l$ is the direct limit embedding.

More precisely, $j(\M_\l)$ is the direct limit of $(j(\M_\gg), \pi_{\gg, \gg'}^\T: \gg<\gg', \gg, \gg'\in [\b, \l)_\T)$ and $k(\M_\l)=k(\M_\b)$ (recall \rdef{fine structural j-structure} which defined $j(\M)$ and $k(\M)$).
\end{enumerate}
If $\a+1\in D$ then we say that there is a drop at $\a+1$. Suppose $\a+1\in D$ and $\b=\T(\a+1)$. If $\omega\xi_\a<\ord(\M_\b)$ then we say that there is a drop in model at $\a+1$ and otherwise we say there is a drop in degree. $\myqedhere$
\end{definition}

We set $\M_\a^\T=\M_\a$, $E_\a^\T=E_\a$, $\ind_\a^\T=\ind^{\M_\a}(E_\a)$, $\lh(\T)=\eta$, $\k_\a^\T=\cp(E_\a^\T)$ and $\nu_\a^\T=\nu(E_\a^\T)$\footnote{Here $\nu(E)$ is the natural length of $E$. See \cite[Definition 2.2]{OIMT}.}. We will drop superscript $\T$ when it is clear from context.
\begin{definition} We say that $\T$ is an \textbf{iteration tree} if it is a putative iteration tree such that for every $\a<\lh(\T)$, $\M_\a^\T$ is well-founded.$\myqedhere$
\end{definition}

\begin{definition}\label{normal iteration tree} Given a putative iteration $\T$ on $\M$ we say that $\T$ is \textbf{normal} if 
\begin{enumerate}
\item $R^\T=\{0\}$,
\item for all $\a<\b<\lh(\T)$, $\ind_\a^\T<\ind_\b^\T$, and
\item for all $\a$ such that $\a+1<\lh(\T)$, $\b=\T(\a+1)$ is the least $\b'$ such that $(\k_\a^+)^{\M_\a|\ind_\a}<\nu_{\b'}$.
\end{enumerate}
$\myqedhere$
\end{definition}

\begin{notation}\label{notation for iteration trees} Given a putative iteration tree 
\begin{center} $\T=((\M_\a)_{\a<\eta}, (E_\a)_{\a<\eta-1}, D, R, (\b_\a, m_\a)_{\a\in R}, T)$\end{center} and ordinals $\gg<\zeta\leq \eta$, we will use the following notations:
\begin{enumerate}
\item $\T_{\leq \gg}=((\M_\a)_{\a\leq \gg}, (E_\a)_{\a<\gg}, (\xi_\a)_{\a<\gg}, D\cap (0, \gg], R\cap [0, \gg), (\b_\a, m_\a)_{\a\in R\cap [0, \gg)}, T\cap (\gg+1)^2)$,
\item $\T_{[\gg, \zeta]}=((\M_\a)_{\a\in [\gg, \zeta]}, (E_\a)_{\a\in [\gg, \zeta)}, D\cap (\gg, \zeta], R\cap [\gg, \zeta), (\b_\a, m_\a)_{\a\in R\cap [\gg, \zeta)}, T\cap [\gg, \zeta]^2)$.
\item Other notations such as $\T_{<\gg}$, $\T_{\geq\gg}$, $\T_{(\gg, \xi)}$ and etc are defined in the obvious manner.
\item Given $\a\in R^\T$, 
\begin{center}
${\sf{next}}^\T(\a)=\begin{cases}
\min(\R^\T-(\a+1)) &:  \R^\T-(\a+1)\not =\emptyset\\
\lh(\T) &: \ \text{otherwise.}
\end{cases}$\end{center}
\item Given $\a\in R^\T$, ${\sf{nc}}^\T_\a=\T_{[\a, \a']}$ where $\a'={\sf{next}}^\T(\a)$\footnote{``nc" stands for ``normal component".}.
\item We say that $\U$ is a \textbf{normal component} of $\T$ if for some $\a\in R^\T$, $\U={\sf{nc}}^\T_\a$. 
\item If $\U$ is a normal component of $\T$ then we let $\a^\T(\U)=\a$ where $\a$ is as above. We say $\U$ is the last normal component of $\T$ if $\a^\T(\U)=\max(R^\T)$.
\end{enumerate}
$\myqedhere$
\end{notation}

\begin{definition}
We say that a putative iteration tree $\T$ is a \textbf{putative \textit{stack}} if 
\begin{enumerate}
\item $R^\T$ is closed,
\item for all $\a\in R^\T$, ${\sf{nc}}^\T_\a$, after obvious re-enumeration of its members, is a normal iteration tree on $\M_\a^\T||(\b_\a^\T, m_\a^\T)$.
\end{enumerate}
We say $\T$ is a \textbf{stack} if $\T$ is a putative stack that is an iteration tree.$\myqedhere$
\end{definition}

\begin{definition}\label{stacks on hod like lses} Suppose $\P$ is an ${\sf{lhes}}$. We say that $\T$ is a \textbf{semi-smooth}\footnote{``Smooth" is inspired by Jensen's terminology who uses ``smooth stack" for stacks that do not allow drops in model or degree at the begining of the rounds. See \cite{JensenBook}.} stack on $\P$ if for all $\a\in R^\T$, 
\begin{enumerate}
\item ${\sf{nc}}^\T_\a$ is normal\footnote{After trivial re-organization.}
\item $\M_\a^\T||(\omega\b_\a^\T, m_\a^\T)$ is a layer of $\M_\a^\T$.
\end{enumerate}
$\myqedhere$
\end{definition}

\begin{rem}[Semi-smooth convention]\label{semi-smooth convention}\index{semi-smooth convention} Because we will mostly work with semi-smooth stacks on ${\sf{lhes}}$, we make the convention that all stacks are semi-smooth. 

$\myqedhere$
\end{rem}

In the sequel, we will often say that $\R$ is a node of $\T$ to mean that $\R=\M_\a^\T$ for some $\a<\lh(\T)$.\\

\textbf{Branches of iterations.}\\
Suppose $\M$ is an $\sf{lhes}$ and $\T$ is an iteration tree on $\M$ such that $\lh(\T)$ is a limit ordinal. We say $b\subseteq \lh(\T)$ is a putative cofinal branch of $\T$ if $b$ is cofinal in $\lh(\T)$ and for every $\a\in b$, $[0, \a]_\T\subseteq b$. Given a putative cofinal branch $b$ of $\T$, we say $b$ is a cofinal branch of $\T$ if  $D^\T\cap b$ is finite. Given a cofinal branch $b$ of $\T$ we let $\M^\T_b$ be the direct limit of the directed system $((\M_\a)_{\a\in b'}, (\pi_{\a, \b}^\T)_{\a<\b, \a,\b\in b'})$ where $b'=b-\max(D^\T\cap b)$\footnote{Here the direct limit is defined analogously to clause 7 of \rdef{putative it}.}. We say $b$ is a well-founded cofinal branch if it is a cofinal branch such that $\M^\T_b$ is well-founded. If $b$ is a putative cofinal branch of $\T$ then we let $\T^\frown \{b\}$ be the unique putative iteration tree $\U$ such that 
\begin{enumerate}
\item $\lh(\U)=\lh(\T)+1$, 
\item $T^\T=T^\U\rest (\lh(\T)\times \lh(\T))$
\item for all $\a<\lh(\T)$, $\M_\a^\U=\M_\a^\T$ and $E_\a^\T=E_\a^\U$,
\item $D^\T=D^\U$ and $R^\T=R^\U$,
\item $\M_{\lh(\T)}^\U=\M^\T_b$,
\item $[0, \lh(\T))_\U= b$.
\end{enumerate}
We also let for $\a\in b$, $\pi^\T_{\a, b}=\pi^{\T^\frown\{b\}}_{\a, \lh(\T)}$ given the later embedding is defined, and if $\pi^\T_{0, b}$ is defined then we let $\pi^\T_b=\pi^\T_{0, b}$. \\

\textbf{Strategies.}\\
Given a stack $\T$ on $\M$ with last model $\N$ and a stack $\U$ on $\N$ we can form $\T$-followed-$\U$ stack $\T^\frown \U$. More formally, $\T^\frown \U$ is the unique stack $\W$ on $\M$ such that $\lh(\W)=\lh(\T)+\lh(\U)$,  $\W_{<\lh(\T)}=\T$, 
\begin{center}
$R^\W=R^\T\cup \{\lh(\T)+\a: \a\in R^\U\}$
\end{center}
 and $\W_{\geq \lh(\T)}=\U$\footnote{After obvious re-enumeration of its members.}. We may often say that $\T^\frown \U$ is a normal iteration tree if after straightforward re-enumeration of its members it becomes a normal iteration tree.\footnote{In this re-enumeration we must set $R^{\T^\frown \U}=\{0\}$.}

Suppose $\M$ is an $\sf{lhes}$ and $\Sigma$ is a function. We say $\Sigma$ is an iteration strategy for $\M$ if whenever $\T\in \dom(\Sigma)$ then $\T$ is a stack on $\M$, $\lh(\T)$ is a limit ordinal and $\Sigma(\T)$ is a cofinal well-founded branch of $\T$.

A putative iteration tree $\T$ on $\M$ is according to $\Sigma$ if for all limit ordinals $\a<\lh(\T)$, $\T_{<\a}\in \dom(\Sigma)$ and $\Sigma(\T_{<\a})=[0, \a)_\T$. We say $\Sigma$ is a $\k$-strategy for $\M$ if 
\begin{enumerate}
\item $\Sigma$ is an iteration strategy for $\M$ such that if $\T\in \dom(\Sigma)$ then $\T$ is a normal iteration tree on $\M$ of length $<\k$,
\item if $\T$ is a normal iteration tree on $\M$ such that $\lh(\T)<\k$, $\lh(\T)$ is a limit ordinal and $\T$ is according to $\Sigma$ then $\T\in \dom(\Sigma)$, 
\item if $\T\in \dom(\Sigma)$, $b=\Sigma(\T)$ and $\U$ is a normal finite putative iteration tree on $\M^\T_b$ such that $\lh(\T)+\lh(\U)<\k$ and $\T^\frown \U$ is a normal putative iteration tree on $\M$ (after obvious re-enumeration of it) then $\U$ is an iteration tree. 
\end{enumerate}
Alternatively, $\Sigma$ is a $\k$-strategy if it is a winning strategy for II in the iteration game $\mathcal{G}(\M, \k)$. This game is defined immediately after \cite[Definition 3.3]{OIMT}. The subscript $k$ that appears in this game is just $k(\M)$.  

Similarly we can define $(\k, \l)$-strategy for $\M$ that acts on stacks. Here $\k$ bounds the number of normal components and $\l$ bounds the length of the normal components of the stacks. The relevant iteration game, $\mathcal{G}(\M, \k, \l)$ appears soon after \cite[Remark 4.3]{OIMT}.\\

\textbf{The common part.}\\
Given a normal iteration tree $\T$ we let the common part of $\T$ be 
\begin{center}
$\m(\T)=\cup_{\a<\lh(\T)}\M^\T_\a|\ind_\a^\T$.
\end{center} Usually in literature, for example in \cite{OIMT}, $\m(\T)$ is denoted by $\M(\T)$. However, inner model theorists make $\M$ tired, and so in this book we give it less weight to carry than its usual heavy load. Following \cite{OIMT} we let $\d(\T)=\ord(\m(\T))=\sup\{\ind^\T_\a: \a<\lh(\T)\}$.\\

\textbf{Restrictions}

\begin{terminology}\label{above eta} Suppose $\M$ is an $\sf{lses}$ and $\T$ is a stack on $\M$. 
\begin{enumerate}
\item Given $\eta\leq \ord(\M)$, we say $\T$ is \textbf{above $\eta$} if for all $\a+1<\lh(\T)$, $\cp(E_\a^\T)\geq \eta$. 
\item We say that $\T$ is \textbf{below $\eta$} if for all $\a+1<\lh(\T)$, either $\pi^{\T}_{0, \a}$ is undefined or $\ind_\a^\T< \pi_{0, \a}^\T(\eta)$.
\item If $\N\insegeq \M$ then we say $\T$ is \textbf{based on} $\N$ if $\T$ is below $\ord(\N)$.  
\end{enumerate}
$\myqedhere$
\end{terminology}

\begin{definition}\label{restriction of a stack} Suppose $\M$ is an $\sf{lhes}$ and $\N\insegeq \M$ is such that $ord(\N)$ is a regular cardinal of $\M$ such that $\rho(\M)>\ord(\N)$. Suppose further that $\T$ is a stack on $\M$. We then let $\downarrow(\T, \N)$ be the portion of $\T$ that is based on $\N$. More precisely, if 
\begin{center}
$\T=((\M_\a)_{\a<\eta}, (E_\a)_{\a<\eta-1}, D, R, (\b_\a, m_\a)_{\a\in R}, T)$.
\end{center}
then 
\begin{center}
$\downarrow(\T, \N)=((\M'_\a)_{\a<\eta'}, (E'_\a)_{\a<\nu-1}, D', R', (\b_\a', m_\a')_{\a\in R'}, T')$.
\end{center}
is such that there is an order preserving map $\sigma:\eta'\rightarrow \eta$ such that
\begin{enumerate}
\item for all $\a<\a'<\eta'$, $(\a, \a')\in T'\iff (\sigma(\a), \sigma(\a'))\in T$,
\item for all $\a<\eta'$, $\a\in R'\iff \sigma(\a)\in R$,
\item for all $\a<\eta'$, $\M'_\a\insegeq \M_{\sigma(\a)}$ and $E'_\a=E_{\sigma(\a)}$,
\item for all $\a<\eta'$, $\b_\a'=\b_{\sigma(\a)}$ and $m'_\a=m_{\sigma(\a)}$,
\item for all $\a+1<\eta'$, $\sigma(\a)$ is the least $\b\in (\sup_{\gg<\a}\sigma(\gg), \eta)$ such that $\ind_\b^\T\leq \ord(\M'_{\a+1})$,
\item if $\a+1=\eta'$ and $\a$ is a limit ordinal then $\sigma(\a)=\sup_{\gg<\a}\sigma(\gg)$,
\item if $\a+1=\eta'$ and $\a=\b+1$ then $\sigma(\a)=\sigma(\b)+1$,
\item for all $\a<\eta'$, if there is $\b<\eta$ such that $\ind_\b^\T\leq \ord(\M'_\a)$ then $\a+1<\eta'$. 
\end{enumerate}
$\myqedhere$
\end{definition}

\begin{definition}\label{upward extension of a stack} Suppose $\M$ is an $\sf{lhes}$ and $\N\insegeq \M$ is such that ${\sf{ord}}(\N)$ is a regular cardinal of $\M$ such that $\rho(\M)>\ord(\N)$. Suppose further that $\T$ is a stack on $\N$. We then let $\uparrow(\T, \M)$ be the result of ``applying"  $\T$ to $\M$. More precisely, if 
\begin{center}
$\T=((\N_\a)_{\a<\eta}, (E_\a)_{\a<\eta-1}, D, R, (\b_\a, m_\a)_{\a\in R}, T)$.
\end{center}
then 
\begin{center}
$\uparrow(\T, \M)=((\M_\a)_{\a<\eta}, (E_\a)_{\a<\eta-1}, D, R, (\b_\a, m_\a)_{\a\in R}, T)$.
\end{center}
with $\M_0=\M$.$\myqedhere$
\end{definition}

Suppose for some $\eta\leq \ord(\Q)$, $\eta$ is a regular cardinal of $\P$, $\rho(\P)>\eta$ and $\T$ is below $\eta$. We then have that $\T\rest \Q$ is the unique stack $\U$ on $\Q$ such that the copy of $\U$ onto $\P$ via $id$ is $\T$. 
  
\section{Layered strategy e-structures}

In this manuscript, we are mostly concerned with $\sf{lhes}$ whose $f$ predicate codes a strategy.  The goal of this section is to introduce the language used to describe such structures.  

Suppose that $\M$ is an $\sf{lhes}$. We then say that a shifted amenable function $f$ \textit{codes a partial strategy function} for $\M$ if letting $g$ be the amenable component of $f$, the following conditions hold:
\begin{enumerate}
\item $\dom(f)\subseteq \{ (\mathcal{J}_\omega(\T), \T, \in): \T$ is a stack on $\M$ without a last model$\}$.
\item Whenever $\T$ is a stack on $\M$ such that $(\mathcal{J}_\omega(\T), \T, \in) \in \dom(f)$ and whenever $\U$ is an initial segment of $\T$ without a last model, $(\mathcal{J}_\omega(\U), \U, \in)\in \dom(f)$ and 
\begin{center}
$g((\mathcal{J}_\omega(\U), \U, \in))=[0, \lh(\U))_\T$.
\end{center}
\item For all $(\mathcal{J}_\omega(\T), \T, \in)\in \dom(f)$, $g((\mathcal{J}_\omega(\T), \T, \in))$ is a cofinal branch of $\T$.
\end{enumerate}
 Notice that we do not require that\\\\
 (a) $g((\mathcal{J}_\omega(\T), \T, \in))$ is a well-founded branch of $\T$,\\
 (b) if $(\mathcal{J}_\omega(\T), \T, \in)\in \dom(f)$ and $b=g((\mathcal{J}_\omega(\T), \T, \in))$ is a cofinal well-founded branch of $\T$ then any reasonable finite extension of $\T^\frown\{b\}$ has well-founded models.\\\\
 Conditions (a) and (b) will be part of a more restrictive notion. 
 
 When defining short tree strategy mice, we will encounter hybrid structures whose $f$ predicate doesn't necessarily code a strategy but a partial strategy. We make this notion more precise. First we make a useful definition.  
 
 
 \begin{definition}\label{semi-strategy}
 Suppose $\M$ is an $\sf{lhes}$. We then say that $\Sigma$ is a \textbf{semi-strategy} for $\M$ if the domain of $\Sigma$ consists of quadruples $(\M_0, \T_0, \M_1, \U)$ such that 
 \begin{enumerate}
 \item $\M_0=\M$, 
 \item $\T_0$ is a normal tree on $\M_0$, 
 \item $\M_1$ is either the last model of $\T_0$ or $\T_0$ doesn't have a last model and $\M_1=(\m(\T_0))^\#$\footnote{This is the true, $\omega_1$-iterable, sharp of $\m(\T_0)$.}, and
 \item $\U$ is a stack on $\M_1$ below $\d(\T_0)$\footnote{This means that all extenders used in $\U$ have lengths below the image of $\d(\T_0)$. I.e. for each $\a<\lh(\T_0)$ either $[0, \a)_\T\cap D^\T\not =\empty$ or $\ind^\T_\a<\pi_{0, \a}^\T(\d(\T_0))$.}.$\myqedhere$
 \end{enumerate}
 
  We can then consider amenable functions that code partial semi-iteration strategies. We will abuse our terminology and will treat semi-iteration strategies as if they were just strategies. 
 \end{definition}

Suppose then a shifted amenable function $f$ codes a partial strategy function for $\M$. We then let $\Sigma^f$ be the partial strategy function coded by $f$. More precisely, letting $g$ be the amenable component of $f$, 
\begin{enumerate}
\item $\dom(\Sigma^f)=\{ \T: (\mathcal{J}_\omega(\T), \T, \in)\in \dom(f)\}$ and
\item for all $\T\in \dom(\Sigma^f)$, $\Sigma^f(\T)=g((\mathcal{J}_\omega(\T), \T, \in))$.
\end{enumerate}
 We say $f$ \textit{codes a partial strategy} if $\Sigma^f$ chooses cofinal and well-founded branches. We say $f$ \textit{codes a total $A$-strategy} if $\Sigma^f(\T)$ is defined whenever $\T\in A$ is of limit length and is according to $\Sigma$. If $A$ is clear from context then we will drop it from our notation.
 
Following \cite{ATHM}, if $\M$ is an $\sf{lhes}$, $\N\insegeq \M$ and $\Sigma$ is an iteration strategy for $\M$ then $\Sigma_\N$ is the strategy of $\N$ we get by the copy construction. More precisely, $\Sigma_\N$ is the $id$-pullback of $\Sigma$.  Like in \cite{ATHM}, if a transitive structure $P$ has a distinguished sequence of extenders then when discussing iterability of $P$ we will always mean iterability with respect to that extender sequence.
 
 \begin{definition}[Strategic e-structure, $\sf{ses}$]\label{strategy premouse} \index{strategic e-structure, ses}Suppose $\P$ is a transitive structure, $X$ is a self-well-ordered set such that $\P\in X$ and $\M$ is a $\phi$-indexed ${\sf{hes}}$. We say $\M$ is a \textbf{$\phi$-indexed strategic e-structure} ($\sf{ses}$) over $X$ based on $\P$ if $f^\M$ codes a partial iteration strategy for $\P$ and for any $w\in \dom(f^\M)$ if $\b=\min(f^\M(w))$ then $\M|\b$ is closed\footnote{See \rdef{closed under sharps}. Also, recall that for such $\b$ we have $\omega\b=\b$}.
 
 We say $\M$ is based on $\P$ if $\M$ is over $\mathcal{J}_\omega[\P]$ and is based on $\P$.$\myqedhere$
\end{definition}

In \rsec{sec short tree strategy mice}, we will also need unindexed $\sf{ses}$\footnote{Notice that in \rdef{unindexed ses}, \textit{unindexed} simply means that no indexing is specified. It is possible that a given unindexed $\sf{ses}$ $\M$ is in fact $\phi$-indexed for some $\phi$.}.

 \begin{definition}[Unindexed $\sf{ses}$]\label{unindexed ses} Suppose $\P$ is a transitive structure, $X$ is a self-well-ordered set such that $\P\in X$ and $\M=\mathcal{J}^{\vec{E}, f}(X)$ is a hybrid $\mathcal{J}$-structure over $X$. We say $\M$ is an \textbf{unindexed strategic e-structure} (unindexed $\sf{ses}$) over $X$ based on $\P$ if the following clauses hold.
 \begin{enumerate}
 \item $f^\M$ codes a partial iteration strategy for $\P$ such that for any $w\in dom(f^\M)$ if $\b=min(f^\M(w))$ then $\M|\b$ is closed\footnote{See \rdef{closed under sharps}. Also, recall that for such $\b$ we have $\omega\b=\b$}.
 \item $\vec{E}$ is a mixed indexed extender sequence. 
 \item If $\M=(\M', k)$\footnote{See \rdef{fine structural j-structure}.} then for every $(\omega\b, m)<l(\M)$,  $\M||(\omega\b, m)$ is sound.
 \end{enumerate}
 
 We say $\M$ is based on $\P$ if $\M$ is over $\mathcal{J}_\omega[\P]$ and is based on $\P$.$\myqedhere$\\
\end{definition}

\begin{definition}[Layered strategic e-structure, $\sf{lses}$]\label{strategy lhp}\index{layered strategy e-structure, lsp} Suppose $\M$ is a $\phi$-indexed $\sf{lhes}$. We say $\M$ is a \textbf{$\phi$-indexed layered strategic e-structure} ($\sf{lses}$) if for all $\Q\in Y^\M$, in $\M$, 
\begin{enumerate}
\item  $f^\M(\Q)$ codes a partial iteration strategy  for $\Q$ such that for every $w\in \dom(f^\M(\Q))$, if $\b=\min(f^\M(\Q)(w))$ then $\M|\b$ is closed, and
\item if $\Q_0,\Q_1\in Y^\M-(X^\M\cup\{X^\M\})$\footnote{Recall $X^\M$ is the set or structure over which $\M$ is defined.} are such that $\Q_0\insegeq \Q_1$ then letting, for $i\in 2$, $\Sigma_i$ be the partial iteration strategy  coded by $f^{\M}(\Q_i)$ and $\Lambda$ be the $id$-pullback of $\Sigma_1$, then $\Lambda\subseteq \Sigma_0$\footnote{Here, we cannot demand equality as there maybe $\T\in \dom(\Sigma_0)$ such that if $\U$ is the $od$-copy of $\T$ on $\Q_1$, $\U\not \in \dom(\Lambda)$.}. 
\end{enumerate}
$\myqedhere$
\end{definition}


 If $\Q\in Y^\M$ then we let $\Sigma_\Q^\M$ be the partial strategy function coded by $f^\M(\Q)$ and let $\Sigma^\M$ be the function with domain $Y^\M$ such that $\Sigma^\M(\Q)=\Sigma^\M_\Q$. The next definition isolates the language of $\sf{lses}$ and $\sf{ses}$.
 
 \begin{definition}\label{language of lses} We let $\mathcal{L}_{\sf{ses}}$ be the language of $\sf{ses}$ intended for lightface $\sf{ses}$, where we say $\M$ is a lightface $\sf{ses}$ if for some $\P$, $\M$ is an $\sf{ses}$ over $\mathcal{J}_{\omega}[\P]$ based on $\P$. Thus, $\mathcal{L}_{\sf{ses}}$ augments the ordinary language for premice as introduced in \cite[Definition 2.10]{OIMT} by adding one constant symbol $\dot{\P}$ for $\P$ and a predicate symbol $\dot{f}$ for $f$. $\mathcal{L}_{\sf{ses}}$ can be further augmented by a constant symbol for $X$ (see \rdef{strategy premouse}), and this language can be used for boldface $\sf{ses}$. 
 
 We let $\mathcal{L}_{\sf{lses}}$ be the language of $\sf{lses}$ over $\emptyset$ (those are the $\sf{lses}$ whose $X$ predicate is the $\emptyset$). Thus, $\mathcal{L}_{\sf{lses}}$ is the language of premice augmented by symbols $\{\dot{B}, \dot{f}, \dot{Y}\}$. 
 
In some cases, it is convenient to use the symbol $\dot{\V}$ to  denote the universe of $\sf{lses}$ or $\sf{ses}$, and also the symbol $\dot{\Sigma}$ to indicate the strategy function coded by $\dot{f}$. Moreover, if $\Q\in \dom(\dot{f})$ then we will use $\dot{\Sigma}_\Q$ to denote the strategy function given by $\dot{f}(\Q)$. 

$\myqedhere$
 \end{definition} 

In most applications, $\sf{lses}$ have a very canonical indexing scheme which is originally due to Woodin. At each stage the stack whose branch is being indexed by $f$ is the least stack whose branch hasn't yet been indexed. We call this the \textit{standard indexing scheme} (see \rsec{hp indexing scheme:sec}). \index{standard indexing scheme}

\begin{remark} Unless indicated otherwise, we will always tacitly assume that the extenders used to witness the existence of large cardinals in $\sf{lses}$ belong to the extender sequence of the $\sf{lses}$. Thus, when we say ``$\k$ is a measurable cardinal in $\M$" we mean that there is an extender $E\in \vec{E}^\M$ such that $E$ witnesses that $\k$ is a measurable cardinal in $\M$. In \cite{FarmerThesis}, Schlutzenberg extensively studied the problem of whether in pure extender models all large cardinal properties are witnessed by extenders that are indexed on the extender sequence. In particular, he showed that measurability and Woodinness are witnessed by extenders that are on the extender sequence.$\myqedhere$
\end{remark}

\begin{remark}\label{the y predicate} Suppose $\M$ is an $\sf{lses}$ and $\b<\ord(\M)$. The notations $\M|\b$ and $\M||\b$ were introduced just before \rrem{omegaalpha}. In this remark, we would like to clarify the meaning of $Y^{\M|\b}$. It is not hard to re-formulate \rdef{layered hybrid j-structure} in a way that $\sf{lses}$ become hierarchical $\mathcal{J}$-structures (see \rdef{hierarchical}) with the property that $Y^{\M|\b}=X^\M\cup \{\Q\insegeq \M: \Q\in Y^\M\wedge \ord(\Q)<\omega\b\}$\footnote{One could for example index every $\M||\omega\b\in Y^\M$ at $\omega\b+\omega$.}.$\myqedhere$
\end{remark}


Suppose $\M$ is an $\sf{lses}$ and $\Sigma$ is a $(\k, \theta)$-iteration strategy for $\Q$ for some $\Q\in Y^\M$. Then it can be the case that $\Sigma^\M_\Q\subseteq \Sigma$. When this happens we get structures relative to $\Sigma$. 

\begin{definition}[$(\Sigma, \phi)$-premouse]\label{sigma-premouse}\index{$(\Sigma, \phi)$-premouse} Suppose $X$ is a transitive self-well-ordered structure and $\P\in X$ is an $\sf{ses}$ or $\sf{lses}$ or just a transitive self-well-ordered set. Suppose further that $\Sigma$ is a $(\k, \theta)$-iteration strategy for $\P$ and $\M$ is a $\phi$-indexed $\sf{ses}$ over $X$ based on $\P$. Then $\M$ is called a \textbf{$(\Sigma, \phi)$-premouse} over $X$ based on $\P$ if $\Sigma^\M\subseteq \Sigma\rest \M$.

Similarly, if $\M$ is a  $(\Sigma, \phi)$-premouse over $X$ based on $\P$, $X=(\mathcal{J}_\omega[\P], \P, \in)$ and $\Sigma$ is a $(\k, \theta)$-iteration strategy for $\P$ then $\M$ is called a  \textbf{$(\Sigma, \phi)$-premouse} over $X$.$\myqedhere$
\end{definition}

We then say $\M$ is a $(\Sigma, \phi)$-premouse if one of the cases in \rdef{sigma-premouse} holds.

\begin{definition}[$(\Sigma, \phi)$-mouse]\label{sigma-mouse}\index{$(\Sigma, \phi)$-mouse} Keeping the notation of \rdef{sigma-premouse},  
we say $\M$ is a \textbf{$(\Sigma, \phi)$-mouse} if $\M$ has an $\omega_1+1$-iteration strategy $\Lambda$ such that whenever $\N$ is a $\Lambda$-iterate of $\M$ then $\N$ is a $(\Sigma, \phi)$-premouse. $\myqedhere$
\end{definition}

We warn the reader that we will often omit $\phi$ from our notation and say $``\M$ is a $\Sigma$-mouse" instead of ``$\M$ is a $(\Sigma, \phi)$-mouse" if $\phi$ is clear from the context. 


\section{Iterations of $(\Sigma, \phi)$-mice}

Suppose $X$ is a transitive self-well-ordered structure such as $\sf{ses}$ or $\sf{lses}$ or just a transitive self-well-ordered set.   Suppose further that $\Sigma$ is an $(\omega_1, \omega_1)$-iteration strategy for some $\P\in X$ (which is also $\sf{ses}$ or $\sf{lses}$ or some transitive set) and $\phi$ is an indexing scheme. Given two $(\Sigma, \phi)$-mice, we can compare them using the usual comparison argument. 

\begin{theorem}[Theorem 3.11 of \cite{OIMT}] Suppose $\M$ and $\N$ are two countable $(\Sigma, \phi)$-mice with $(\omega_1+1)$-iteration strategies $\Lambda$ and $\Gamma$ respectively. Then there are iteration trees $\T$ and $\U$ on $\M$ and $\N$ respectively according to $\Lambda$ and $\Gamma$ respectively, having last models $\M^\T_\a$ and $\N_\eta^\N$ such that either
\begin{enumerate}
\item the iteration embedding $\pi^\T_{0, \a}$-exists and $\M_\a^\T$ is an initial segment of $\M_\eta^\U$, or
\item the iteration embedding $\pi^\U_{0, \eta}$-exists, and $\M_\eta^\U$ is an initial segment of $\M_\a^\T$.
\end{enumerate}
\end{theorem}

Comparison for $\sf{lses}$ is more involved and we do not know how to do it in general. Below we recall our primary method of identifying the good branches of iteration trees. Recall that the strategy for a sound mouse projecting to $\omega$ is determined by
{\em $\Q$-structures}\index{$\Q$-structures}. For $\T$ normal, let $\Phi(\T)$ be the phalanx of $\T$ (see Definition 6.6 of \cite{CMIP}). 

\begin{definition}\label{qstructures} Suppose $\M$ is an $\sf{lses}$ (or $\sf{ses}$). Let $\T$ be a normal tree of limit
length on $\M$  and let $b$ be a cofinal branch of $\T$. Then $\Q(b,\T)$ is the shortest
initial segment $\Q$ of $\M_b^{\T}$, if one exists, such that
$\Q$ projects strictly across $\d(\T)$ (i.e. $\rho(\Q)<\d(\T)$) or defines a function witnessing $\d(\T)$ is not a 
Woodin cardinal as witnessed by the extenders on the sequence of $\m(\T)$. Equivalently, $\Q(b, \T)=\M^\T_b||\omega\xi$ such that $\xi$ is the largest $\xi'$ with the property that $\M^\T_b||\omega\xi'\models ``\d(\T)$ is a Woodin cardinal".$\myqedhere$
\end{definition}

Next we would like to state a general result stating that branches identified by $\Q$-structures are unique. 
\begin{definition}\label{id pullback initial segment}
Suppose that $\M$ is an $\sf{lses}$ and $\Sigma$ is a strategy for $\M$.
If $\N$ is a $\Sigma$-iterate of $\M$ via $\T$ then we let $\Sigma_{\N, \T}$ be the strategy of $\N$ given by $\Sigma_{\N, \T}(\U)=\Sigma(\T^\frown \U)$. If then $\Q\insegeq \N$ then we let $\Sigma_{\Q, \T}$ be the $id$-pullback of $\Sigma_{\N, \T}$.$\myqedhere$
\end{definition}

\begin{definition} Suppose $\M$ is a $\phi$-indexed $\sf{lses}$ (perhaps over some set $X$ and based on some $\P\in X$) and $\Sigma$ is an iteration strategy for $\M$. We say $(\M, \Sigma)$ is a \textbf{layered strategy $\phi$-mouse} ($\phi$-$\sf{lsm}$) pair if $\Sigma$ has hull condensation (see Definition 1.30 of \cite{ATHM}) and whenever $\N$ is a $\Sigma$-iterate of $\M$ via $\T$ then $\N$ is a $\phi$-indexed $\sf{lses}$ and for any $\Q\in Y^\N-X$, $\Sigma^\N_\Q\subseteq \Sigma_{\Q, \T}$. We say $(\M, \Sigma)$ is sound if $\M$ is sound.

Similarly we can define $\phi$-${\sf{sm}}$. We will say that $\M$ is a $(\Sigma, \phi)$-$\sf{lsm}$ or $(\Sigma, \phi)$-$\sf{sm}$ if $(\M, \Sigma)$ is respectively a $\phi$-$\sf{lsm}$ or $\phi$-$\sf{sm}$. $\myqedhere$
\end{definition}

\begin{terminology}\label{above eta and other things} Suppose $\M$ is an $\sf{lses}$. 
\begin{enumerate}
\item We say $\gg$ is a \textbf{cutpoint} of $\M$ if there is no extender $E\in \vec{E}^\M$ such that $\cp(E)<\gg\leq \ind^\M(E)$.
\item We say $\gg$ is a \textbf{strong cutpoint} of $\M$ if there is  no extender $E\in \vec{E}^\M$ such that $\cp(E)\leq \gg\leq \ind^\M(E)$.
\item An extender $E\in \vec{E}^\M$ \textbf{overlaps} $\kappa$ if $\cp(E)< \kappa \leq \lh(E)$, and \textbf{weakly overlaps} $\k$ if $\cp(E)\leq \k\leq \lh(E)$.
\item $\ord(Y^\M)=\sup\{\ord(\Q): \Q\in Y^\M\}$. 
\end{enumerate}
$\myqedhere$
\end{terminology}

\begin{theorem}\label{strategyunique} Suppose $(\M, \Sigma)$ is a sound $\phi$-$\sf{lsm}$ pair, and suppose $\gg<\ord(\M)$ is a strong cutpoint of $\M$ such that 
\begin{center}
$\ord(Y^\M)\leq \gg$ and $\rho(\M)\leq \gg$.
\end{center}
 Then $\M$ has at most one $(\omega_1+1)$-iteration strategy $\Lambda$ that acts on iteration trees that are strictly above $\gg$ and whenever $\N$ is a $\Lambda$-iterate of $\M$ then $\N$ is a $\phi$-indexed $\sf{lses}$ and $\Sigma^\N\subseteq \Sigma\rest \N$.
 
Moreover, any such strategy $\Lambda$ is determined by: for countable length normal iteration trees $\T$, $\Lambda(\T)$ is the unique cofinal wellfounded $b$ such that the phalanx \begin{center}
 $\Phi(\T)^\frown ( \d(\T), \Q(b,\T))$\end{center} is $\omega_1+1$-iterable (as a $(\Sigma, \phi)$-phalanx, see \rdef{qstructures} for the meaning of $\Q(b, \T)$).\footnote{The meaning of this is left to the reader, but see \cite[Definition 6.7]{CMIP} or \cite[Definition 2.22]{ANS}.}
\end{theorem}  

In some cases, however, it is enough to assume that $\Q(b, \T)$ is countably iterable. This happens, for instance, when $\M$ has no local Woodin cardinals with extenders overlapping it. While the $\sf{lses}$ we will consider may have initial segments that have Woodin cardinals that are not cutpoints, no such cardinal will be Woodin in the entire model. This simplifies our situation somewhat, and below we describe exactly how this will be used. 

\begin{definition}[{Definition 2.1} of \cite{CMWMWC}]\label{qoft} Let $(\M, \Sigma)$ be a sound $\phi$-$\sf{lsm}$ pair and let $\gg<\ord(\M)$ be such that $\tau=\ord(Y^\M) \leq \gg$. Suppose $\T$ is a normal iteration tree on $\M$ that is above $\gg+1$;
then  $\Q(\T)$, if exists, is the unique $\Q$ that has the following properties.
\begin{enumerate}
\item $\Q$ is a $(\Sigma_{\M||\tau}, \phi)$-${\sf{sm}}$ over $\m(\T)$ based on $\M||\tau$ (in particular, $\d(\T)$ is a strong cutpoint of $\Q$).
\item $\mathcal{J}_\omega(\Q)\models ``\d(\T)$ is not a Woodin cardinal",
\item $k(\Q)$ is the least $k$ such that 
\begin{enumerate}
\item $\rho_k(\Q)< \d(\T)$ or 
\item $\rho_k(\Q)=\d(\T)$ and there is $r\Sigma_{k}^\Q$-definable function $f:\d(\T)\rightarrow \d(\T)$ witnessing that $\d(\T)$ is not a Woodin cardinal as witnessed by  the extenders of $\m(\T)$.
\end{enumerate}
\end{enumerate}
$\myqedhere$
\end{definition}

Countable iterability is usually enough to guarantee there is at most one $\sf{lses}$ with
the properties of $\Q(\T)$. If it exists, $\Q(\T)$ might identify the good branch of $\T$,
the one any sufficiently powerful iteration strategy must choose. This is the content of
the next lemma which can be proved by analyzing the proof of Theorem 6.12 of \cite{OIMT}. To state it we need to introduce fatal drops.

\begin{definition}[Fatal drop]\label{fatal drop}\index{fatal drop}
Suppose $\M$ is a $\phi$-indexed $\sf{lses}$ and $\T$ is an iteration tree on $\M$. We say $\T$ has a \textbf{fatal drop} if for some $\a<\lh(\T)$ and $\eta<\ord(\M_\a^\T)$,
\begin{enumerate}
\item $\eta$ is a cutpoint of $\M_\a^\T||\omega\xi_\a^\T$,
\item $\sup\{\ind^\T_\b: \b<\a\}\leq \eta$,
\item $\rho(\M_\a^\T||(\omega\xi_\a^\T, k^\T_\a))\leq \eta$\footnote{$\omega\xi_\a^\T$ and $k_\a^\T$ are defined in clause 8 of \rdef{putative it}.}, 
\item  $\T_{\geq\a}$ is a normal iteration tree on $\M_\a^\T||(\omega\xi_\a^\T, k_\a^\T)$ that is above $\eta$.
\end{enumerate}
 We then say $\T$ has a fatal drop at $(\a, \eta)$ if the pair is the lexicographically least satisfying the above condition. $\myqedhere$
\end{definition}

The following is the lemma mentioned above. 

\begin{lemma}\label{lem qstructures} Let $(\M, \Sigma)$ be a $\phi$-$\sf{lsm}$ pair such that $\ord(Y^\M)$ is a strong cutpoint of $\M$\footnote{The hod mice considered in the manuscript satisfy this condition.} and let $\gg<\ord(\M)$ be such that $\ord(Y^\M)\leq \gg$.  Suppose $\T$ is a normal iteration tree on $\M$ that is  above $\gg+1$ and has limit length.
\begin{enumerate}
\item  Suppose $\Q(\T)$ exists. Then there is at most one cofinal branch $b$ of $\T$ such that
either $\Q(\T) = \M_b^{\T}$ or $\Q(\T) = \M_b^{\T} ||\omega\xi$ for some $\xi$ in the wellfounded
part of $\M_b^\T$.
\item  Suppose further no measurable cardinal of $\M$ which is $\geq \gg$ is a limit of Woodin cardinals. Suppose further that $\T$ is according to $\Sigma$, $\T$ doesn't have a fatal drop and if $b=\Sigma(\T)$ then $\Q(b, \T)$-exists.  Then $\Q(b, \T)=\Q(\T)$. 
\end{enumerate}
\end{lemma} 

$\Q(\T)$ identifies $b$ because it determines a canonical cofinal subset
of $rng(\pi_{\a,b}^{\T} \cap \d(\T))$, for some $\a \in b$, to which we can apply Lemma 1.13 of \cite{ATHM} (which is an immediate consequence of the zipper argument from \cite{IT}). 

\begin{rem}Suppose $(\M, \Sigma)$ is a $\phi$-lsm pair and $\Q\in Y^\M-X^\M$. Let $\R=\M$ if $\Q$ is the largest initial segment of $\M$ in $Y^\M$ and otherwise, let $\R$ be the least member of $Y^\M$ properly extending $\Q$. Suppose $\T$ is a tree on $\M$ which is above $\ord(\Q)+1$ and is based on $\R$. Notice that in this case we can define $\Q(\T)$ just as in \rdef{qoft} by using $\R$ instead of $\M$. $\myqedhere$
\end{rem}

We end this section by introducing the $\mathcal{O}$-stack.
Suppose $\P$ is an $\sf{lses}$, $\a, \eta <\ord(\P)$ and $\Q\insegeq \P||\eta$. Let $e^\P_{\eta, \a}$ be the least ordinal $\b > \eta$, if it exists, such that $\b\in \dom(\vec{E}^\P)$, and letting $E=\vec{E}^\P(\b)$, $\cp(E)\in (\alpha, \eta)$. Thus, $e^\P_{\eta, \a}$ is the index of the first extender that overlaps $\eta+1$ and has a critical point $>\a$. Otherwise, if there is no such extender then set $e^\P_{\eta, \a}=\ord(\P)$.

Let $s^\P_{\eta, \Q}$ be the least ordinal $\b> \eta$, if it exists, such that for some $\R\in Y^\P-X^P$ with $\Q\inseg \R$ letting $F$ be the set indexed at $\b$ in $\P$, $F$ is a pair of the form $(\R, a)$. Thus, $s^\P_{\eta, \Q}$ is the first place above $\eta$ where a branch of some iteration tree $\T$ that is based on a strictly longer layer than $\Q$ is added.  If there is no such $\R$ then let $s^\P_{\eta, \Q}=\ord(\P)$. Let $\eta'=(\eta^+)^\P$ if $(\eta^+)^\P$ exists and otherwise let $\eta'=\ord(\P)$. Set $\a^{\P}_{\eta, \Q, \a}=\min\{ e^{\P}_{\eta, \a}, s^\P_{\eta, \Q}, \eta'\}$. 

Suppose $\M$ is f.s. $\mathcal{J}$-structure and $\eta<\ord(\M)$ is the largest cardinal of $\M$. We then let $\M|(\eta^+)^\M=\M|\ord(\M)$ and $\M||(\eta^+)^\M=\M$. 

\begin{definition}[$\mathcal{O}^\P$-stack]\label{the o stack}\index{$\mathcal{O}^\P$-stack}
Suppose $\P$ is an $\sf{lses}$, $\a, \eta <\ord(\P)$ and $\Q\insegeq \P||\eta$.

We now set 
\begin{center}
$\mathcal{O}^\P_{\eta, \Q, \a}=\P|(\eta^+)^{\P|\a^\P_{\eta, \Q, \a}}$.
\end{center}
Next we define the stack $(\mathcal{O}^{\P, \xi}_{\eta, \Q, \a} : \xi\leq \Omega^\P_{\eta, \Q, \a})$ according to the following recursion:
\begin{enumerate}
\item $\mathcal{O}^{\P, 0}_{\eta, \Q, \a}=\mathcal{O}^\P_{\eta, \Q, \a}$,
\item for $\xi+1\leq \Omega^{\P}_{\eta, \Q, \a}$, $\mathcal{O}^{\P, \xi+1}_{\eta, \Q, \a}=\mathcal{O}^\P_{\ord(\mathcal{O}^{\P, \xi}_{\eta, \Q, \a}), \Q, \a}$,
\item for limit $\l\leq \Omega^{\P}_{\eta, \Q, \a}$, $\mathcal{O}^{\P, \l}_{\eta, \Q, \a}=\bigcup_{\xi<\l}\mathcal{O}^{\P, \xi}_{\eta, \Q, \a}$, and
\item $\Omega^\P_{\eta, \Q, \a}$ is the least $\nu$ such that $\ord(\mathcal{O}^{\P, \nu}_{\eta, \Q, \a})=\a_{\eta, \Q, \a}^\P$.
\end{enumerate}
If $\Q = \P||\kappa$, then we write $\mathcal{O}^\P_{\eta,\kappa,\alpha}$ for $\mathcal{O}^\P_{\eta,\Q,\alpha}$; if $\alpha=0$, we also write $\mathcal{O}^\P_{\eta,\kappa}$ for $\mathcal{O}^\P_{\eta,\Q,\alpha}$. For $\xi\leq \Omega^\P_{\eta, \P||\eta, \a}$, we let $\mathcal{O}^{\P, \xi}_\eta=\mathcal{O}^{\P, \xi}_{\eta, \P||\eta, 0}$ with $\mathcal{O}^{\P}_\eta=\mathcal{O}^{\P}_{\eta, \P||\eta, 0}$.$\myqedhere$
\end{definition}

\section{Hod-like layered hybrid premice}\label{sec: hod-like layered hybrid premice}

 The difference between the $\sf{lses}$ considered here and those considered in \cite{ATHM} is that here we will have $\sf{lses}$ whose predicate codes the \textit{short tree strategy} of its initial segments. The hod mice we will consider in this paper are all \textit{layered}, and we start by introducing these objects. 

If $\M$ is an $\sf{lses}$ and $\k$ is an $\M$-cardinal then we set $E_\xi^\M=\vec{E}^\M(\xi)$ and
\begin{center}
$X_\k^\M=\{ \xi: E_\xi^\M\not =\emptyset$ and $\cp(E_\xi^\M)=\k\}$\index{$X^\M_\k$}.\index{$X_\k^\M$}
\end{center}
We also let 
\begin{center}
$o^\M(\k)=\max(\sup X^\M_\k , (\k^+)^\M)$\index{$o^\M(\k)$}.\index{$o^\M(\k)$}
\end{center}

Suppose $M$ is a transitive structure and $\eta$ is an ordinal. Then we let 
$(\eta^{+\a})^M$ be the $\a$th-cardinal successor of $\eta$ in $M$ if it exists and otherwise, we let it be $\ord(M)$. 

\begin{definition}[Pre-hod-like]\label{pre-hod-like}
Suppose $\P$ is an $\sf{lses}$. We say $\P$ is \textbf{pre-hod-like} if one of the following holds:
\begin{enumerate}
\item (Meek)\index{Meek} There is $\d$ such that 
\begin{enumerate}
\item $\P\models ``\d$ is a Woodin cardinal or a limit of Woodin cardinals",  
\item $\d$ is a cutpoint of $\P$\footnote{This condition follows from the other conditions, but we would like to isolate it.},
\item if $\k<\ord(\P)$ is a limit of Woodin cardinals of $\P$ then $o^\P(\k)<\d$,
\item $\P\models \sf{ZFC}-{\sf{Replacement}}$ and 
\item if $\d$ is a Woodin cardinal of $\P$ then $\P=\bigcup_{n<\omega}\P|(\d^{+n})^\P$, and if $\d$ is a limit of Woodin cardinals of $\P$ then $\d$ is the largest cardinal of $\P$.
\end{enumerate}
\item (Non-meek)\index{Non-meek} There is $\d\leq \ord(\P)$ such that
\begin{enumerate}
\item there is $\k<\d$ such that $\d\leq o^\P(\k)$,
\item if $\k$ is the least $\eta<\d$ such that $\d\leq o^\P(\eta)$ then $o^\P(\k)=\d$ and $\P\models ``\k$ is a limit of Woodin cardinals", 
\item letting $\k<\d$ be the least such that $o^\P(\k)=\d$, $\rho(\P)\in (\k, \d]$ or $\ord(\P)$ is a limit of ordinals $\xi$ such that $\rho(\P||(\xi, \omega))\in (\k, \d]$\footnote{Here, we implicitly assuming that $\xi=\omega\b$ for some $\b$. See \rrem{omegaalpha}.}. 
\item $\P$ is $\d$-sound,
\item if $\dom(\vec{E}^\P)\cap (\d^\P, \ord(\P)]=\emptyset$ then $\mathcal{J}_{\omega}[\P]\models``\d^\P$ is not a Woodin cardinal".
\end{enumerate}
\item (Gentle) $\d=_{def}\ord(\P)$ is a limit of Woodin cardinals of $\P$ and $\P\models \sf{ZFC}-{\sf{Replacement}}$.
\end{enumerate}
We let $\d^\P$\index{$\d^\P$} be the $\delta$ above. $\myqedhere$
\end{definition}
The next definition is somewhat technical. The meaning of it is that we will wait until we see the sharp of a layer before we will activate the strategy.

\begin{definition}[Properly non-meek]\label{proper type II}\index{Properly non-meek} Suppose $\P$ is a non-meek pre-hod-like $\sf{lses}$. We say $\P$ is \textbf{properly  non-meek} if there is $\xi\in \dom(\vec{E}^\P)$ ($\xi$ may be o$(\P)$) such that $\cp(E_\xi^\P)>\d^\P$ and $\P|\xi=\mathcal{J}_\xi[\P|\d^\P]$. $\myqedhere$
\end{definition}

The next definition isolates the type of hod premice that give rise to pointclasses satisfying the Largest Suslin Axiom.

\begin{figure}
\centering
\begin{tikzpicture}[]
\coordinate [label={left:$\kappa$}] (A) at (0, 0);
\coordinate [label={left:$\delta$}] (B) at (0,1.5);
\coordinate (C) at (0,2);
\coordinate [label = {left:$\xi$}](D) at (0,3);
\coordinate (E) at (0.3,0);
\coordinate (F) at (0.3, 1.5);
\coordinate (G) at (0,-1);
\coordinate [label={above:$\P$}] (H) at (0,4);
\coordinate (I) at (-3,4);
\coordinate (J) at (3,4);
\coordinate (K) at (0.3,2);
\coordinate (L) at (0.3,3);
\node [node distance=0.3cm, right of = L] {$E^\P_\xi$};
\draw (A) -- (E);
\draw[->] (E) -- (F);
\draw (H) -- (G);
\draw[dotted] (G) -- (I) -- (J) -- (G);
\draw (C) -- (K) -- (L) -- (D);

%
   \end{tikzpicture}
\caption{Lsa type $\sf{lses}$. Here, $\P$ is an lsa type $\sf{lses}$. $\kappa$ is a limit of Woodin cardinals in $\P$, $\delta=\delta^\P$ is Woodin in $\P$, and $o^\P(\kappa)=\delta$. $\P|\xi$ is the least active level of $\P$ above $\delta$.}
\label{fig:lsa_type}
\end{figure}

\begin{definition}[Lsa type, Figure \ref{fig:lsa_type}]\label{lsa type}\index{lsa type} Suppose $\P$ is a pre-hod-like $\sf{lses}$. We say $\P$ is of \textbf{lsa type} if 
\begin{enumerate}
\item $\P$ is properly non-meek, 
\item $\P\models ``\d^\P$ is a Woodin cardinal"
\end{enumerate}

Suppose $\P$ is a pre-hod-like $\sf{lses}$ of lsa type. We let $\P_{\sf{ex}}\insegeq \P$ be the longest initial segment $\P'$ of $\P$ such that $\P'$ is of lsa type, $\d^\P=\d^{\P'}$ and letting $k=k(\P')$, for every $\k<\d^\P$ there is no cofinal $f:\kappa\rightarrow \d^\P$ that is $r\Sigma_{k}^{\P'}$-definable over $\P'$\footnote{$\sf{ex}$ stands for ``exact".}. We then say that $\P$ is \textbf{exact} if $\P=\P_{\sf{ex}}$.

Continuing with $\P$, let $\a=\min(\dom(\vec{E}^\P)-\d^\P)$ and set $\P_{\#}=\P||\a$. We then say that $\P$ is of \textbf{$\#$-lsa type} if $\P_{\#}=\P$ and $\mathcal{J}_{\omega}[\P]\models ``\d^\P$ is a Woodin cardinal". 

If $\Sigma$ is a strategy of $\P$ then we let $\Sigma_{{\sf{ex}}}$ be the strategy of $\P_{\sf{ex}}$  with the property that $\Sigma_{{\sf{ex}}}= (id$-pullback of $\Sigma)$.$\myqedhere$
\end{definition}

In this paper we will consider hod mice that are lsa small. 

\begin{definition}[Lsa small]\label{lsa small} Suppose $\P$ is a pre-hod-like $\sf{lses}$. We say $\P$ is \textbf{lsa small} if for all $\P$-cardinals $\k$ such that $o^\P(\k)<\d^\P$ and $\P\models ``\k$ is a limit of Woodin cardinals", $\P\models ``o^\P(\k)$ is not a Woodin cardinal".$\myqedhere$
\end{definition}

 \begin{rem}\label{lsa small convention} From now on we tacitly assume that all $\sf{lses}$ considered in this paper are lsa-small. We will, from time to time, remind the reader of this. $\myqedhere$
 \end{rem}

We can now isolate the layers of pre-hod-like $\sf{lses}$. 
\begin{remark}\label{intuitive remarks} Before we give the definition we make the following intuitive remarks. Suppose $\P$ is a hod premouse, which are the objects that we eventually want to define (see \rdef{hod premouse}).
\begin{enumerate}
\item The philosophy behind ``layering" is the desire to make maximal complexity jumps in the Wadge hierarchy. Ordinary mice and premice are designed to reach large cardinals by using the least amount of information large cardinals give us, namely the extenders that induce those embeddings that we use to define the large cardinal in question. For example, to reach a measurable cardinal in a mouse we only use ultrafilters. However, measurability tells us much more than just that there is a nice ultrafilter on some cardinal. For example, if $\kappa$ is measurable then every $\bP^1_1$ set is $\kappa$-homogenously Suslin, and in trying to build mice with measurable cardinals we ignore this extra information. We justify our ignorance by claiming that our algorithms that produce mice (e.g. fully backgrounded constructions, $K^c$ constructions and etc) using extenders as oracles output structures that do inherit all the important properties of large cardinals. That this indeed happens has been verified by Neeman for large cardinals in the region of Woodin cardinal that is a limit of Woodin cardinals (see \cite{Neeman}). However, a priori, this dream-like solution may have been wrong, and more of the information given to us by large cardinals might have been required to reach them in canonical structures, and perhaps the fact that we cannot do significantly better than a Woodin cardinal that is a limit of Woodin cardinals  is a sign that only extenders won't do. 

Hod mice have an entirely different purpose. Instead of large cardinals the aim is to reach or rather ``capture" the Wadge hierarchy inside canonical structures. This is parallel to Shoenfield's Absoluteness, namely that $L$ is $\Sigma^1_2$-correct. Each layer of a hod mouse corresponds to a new level of the Wadge hierarchy. However, what the philosophy of layering claims to be possible is that we can reach all levels of the Wadge hierarchy by simply jumping to the most significant levels of it, and here the significant levels of the Wadge hierarchy are defined to be \textit{the Solovay pointclass}. 

\begin{definition}\label{solovay pointclass}
Assume $\sf{ZF+AD^+}$. We say $\Gamma$ is a \textbf{Solovay pointclass} if there is $\k$ such that $\k$ is a member of the Solovy sequence and $\Gamma=\{A\subseteq \bR: w(A)<\k\}$\footnote{See \rdef{solovay sequence}.}. $\myqedhere$
\end{definition} 
 The strategy of each layer of a hod mouse generates a Solovay pointclass in the sense that the named strategy has Wadge rank $\theta_\a$ for some $\a$. The dream of the ``layering" philosophy is that by only generating the Solovay pointclasses we will reach all levels of the Wadge hierarchy. Internalizing this idea would help the reader with a knowledge of ${\sf{AD^+}}$ theory to understand why layers are defined the way they are defined: every initial segment whose strategy corresponds to a Solovay pointclass is a layer.

Of course, at this stage the idea is vague. If $\P$ is our hod mouse, $\Sigma$ is an iteration strategy for $\P$, $\Q_0\inseg \Q_1\insegeq \P$ then the reader should expect that $\Sigma_{\Q_0}$ is not more complex then $\Sigma_{\Q_1}$, and one can easily build many situations where in fact $\Sigma_{\Q_0}$ is Wadge reducible to $\Sigma_{\Q_1}$ but not vice a versa. However, it may be the case that neither $\Q_0$ nor $\Q_1$ are layers of $\P$. To make the idea work we need to \textit{anticipate} the initial segments of $\P$ whose strategies generate the Solovay pointclasses. Below we spell out what initial segments should be layers in the minimal model of ${\sf{LSA}}$.

\item 
The basic phenomenon that guides us in our definition of layers is the following: Suppose $(\P, \Sigma)$ is a pair such that $\P$ is a hod mouse and $\Sigma$ is its iteration strategy\footnote{For explanatory reasons, we are being somewhat vague.}. Let $\M_\infty(\P, \Sigma)$ be the direct limit of all countable $\Sigma$-iterates of $\P$ and $\pi_{\P, \infty}:\P\rightarrow \M_{\infty}(\P, \Sigma)$ be the direct limit embedding. \\\\
\textbf{Key Phenomenon:} For $\d\leq \d^\P$, $\pi_{\P, \infty}(\d)$ is a member of the Solovay sequence if and only if $\d$ is either a cutpoint Woodin cardinal of $\P$ or a cutpoint limit of Woodin cardinals of $\P$.\\\\
The way we use the Key Phenomenon is as follows. Suppose we have declared $\Q$ a layer of $\P$.
\begin{enumerate}
\item If $\d^\Q$ is a cutpoint Woodin cardinal of $\P$ or a cutpoint limit of Woodin cardinals of $\P$ then $\Q$ is the unique layer $\Q'$ of $\P$ such that $\d^{\Q'}=\d^\Q$.
\item  If, however, $\d^\Q=o^\Q(\k)$ for some $\k$ and $\Q'$ is such that $\rho(\Q')\leq \d^\Q$ then $\pi_{\Q, \infty}(\k)\leq \pi_{\Q', \infty}(\k)$ and the strict inequality cannot be ruled out. Thus, we declare $\Q'$ a layer as it can generate a new Solovay pointclass.
\item Also, suppose $\k<\d^\P$ is a limit of cutpoint Woodin cardinals of $\P$ and suppose that $\a$ is such that either $\a=\ind^\P(E)$ for some $E\in \vec{E}^\P$ with $\cp(E)=\k$ or $\a$ is a limit of such points. Notice now that for every $\b< \a$, $\pi_{\P||(\omega\b, \omega), \infty}(\k)<\pi_{\P||(\omega\a, \omega), \infty}(\k)$. This is because we must have that $\pi_{\P||(\omega\b, \omega), \infty}=\pi_{Ult(\P, E)||(\omega\b, \omega), \infty}$. Therefore,  $\P||(\omega\a, \omega)$ must be a layer of $\P$. 
\end{enumerate}

\item Layers of $\P$  are those proper initial segments of $\P$ whose strategy is being indexed on the strategy predicate of $\P$, with the exception that $\P$ is also considered to be a layer of itself.
\item Meek layers of our hod mice were already studied in \cite{ATHM}. 
\item The key new ingredient of our hod mice is the way we treat the lsa type layers. Given an lsa type layer $\Q$ of $\P$, say $\Q$ is minimal in $\P$ if there is no layer $\Q'\inseg \Q$ such that $\d^\Q=\d^{\Q'}$. Given a minimal lsa type layer $\Q$, we start indexing the short-tree-strategy of $\Q$ into the strategy predicate, and this leads to our notion of short-tree-strategy (${\sf{sts}}$) premouse (see \rdef{sts premouse}). There are two possibilities here. Either (a) we reach a level $\Q'$ with the property that $\Q'$ is an ${\sf{sts}}$ premouse over $\Q$ such that $\d^\Q$ is not a Woodin cardinal definably over $\Q'$ or (b) $\d^\P=\d^\Q$, $\d^\P$ is Woodin in $\P$ and $\P$ above $\d^\Q$ is an ${\sf{sts}}$ premouse. If we do reach such a $\Q'$ then $\Q'$ becomes a layer and we start adding the strategy of $\Q'$. Notice that $\Q'$ itself is of lsa type.
\item The following conditions essentially characterize all proper layers of $\P$, but the conditions below do not spell out the actual definition and are given for explanatory purposes. 
\begin{enumerate}
\item (Woodin cardinals) If $\eta<\d^\P$ is a Woodin cardinal of $\P$ then there is a layer $\Q$ of $\P$ such that $\d^\Q=\eta$.
\item (Limit of Woodin cardinals) If $\k<\d^\P$ is a limit of Woodin cardinals of $\P$ then $\Q=_{def}\P|(\k^+)^\P$ is a layer of $\P$ such that $\d^\Q=\k$ and $\Q$ is the unique layer $\Q'$ of $\P$ with $\d^{\Q'}=\k$. Moreover, $\k$ is a strong cutpoint in $\Q$ (in particular, if $E\in \vec{E}^\P$ is such that $\cp(E)=\k$ then $E$ is total.)
\item (Active layers) If $\k<\d^\P$ is a limit of Woodin cardinals and $E\in \vec{E}^\P$ is such that $\cp(E)=\k$ then there is a layer $\Q$ of $\P$ such that $\d^\Q=\ind^\P(E)$.
\item (Limits of layers) If $\nu<\d^\P$ is a limit of ordinals of the form $\d^\Q$ where $\Q$ is a layer of $\P$ then there is a layer of $\P$ such that $\d^\Q=\nu$.
\item If $\Q$ is a layer of $\P$ with $\d^\Q<\d^\P$ and $\Q'\insegeq \mathcal{O}^\P_{\ord(\Q)}$ is such that $\rho(\Q')\leq \d^\Q$ then $\Q'$ is a layer of $\P$ with $\d^{\Q'}=\d^\Q$. 
\end{enumerate}
\item The layers of a hod-like $\sf{lses}$ are defined in a way that all non-meek layers are properly non-meek. There is no deep reason for doing this. The theory can be developed without this condition, but having more room above $\d^\P$ is a convenience. 
\end{enumerate}
$\myqedhere$
\end{remark}
 
\begin{definition}[Layers of $\sf{lses}$]\label{layers of hod-like lsp}
Suppose $\P$ is an lsa small pre-hod-like $\sf{lses}$. We define the \textbf{layers} $(\P_{\xi, \xi'}: \xi\leq \eta \wedge \xi'\leq \nu_\xi)$ of $\P$ as follows. As part of the definition, we will also define a sequence $(\d_\xi, \iota_{\xi, \xi'}: \xi\leq \eta \wedge \xi'\leq \nu_\xi)$. The sequences are subject to the following requirements:\\\\
\textbf{The Condition Defining the Sequence $(\d_\xi: \xi\leq \eta)$}\\\\
${\sf{R0:}}$ The sequence $(\d_\xi: \xi\leq \eta)$ enumerates in increasing order the set consisting of the following ordinals. \begin{enumerate}
\item Woodin cardinals of $\P$ that are $\leq \d^\P$\footnote{Examining \rdef{pre-hod-like}, one could see that clause 2c leaves open the possibility of $\P$ having Woodin cardinals $>\d^\P$.}.
\item The limits of Woodin cardinals of $\P$ that are $\leq \d^\P$.
\item Ordinals $\nu$ with the property that $\nu\in \dom(\vec{E}^\P)$ and $\cp(\vec{E}^\P(\nu))<\d^\P$ is a limit of Woodin cardinals of $\P$.
\item Ordinals  $\nu$ which are limits of ordinals as in clause 3 above. 
\end{enumerate}
\textbf{The Conditions Defining the Sequence $(\iota_{\xi, 0}: \xi< \eta)$}\\\\
${\sf{R1:}}$ Suppose $\xi<\eta$ and $\d_\xi$ is a Woodin cardinal of $\P$. Then
\begin{center}$\iota_{\xi, 0}=\ord(\mathcal{O}^{\P, \omega}_{\d_\xi, \P|\d_{\xi-1}})$\footnote{Since we do not have a Woodin limit of Woodin cardinals in our $\P$, $\xi-1$ makes sense. For $\xi=0$, we let $\d_{-1}=0$.}.\end{center}
${\sf{R2:}}$ Suppose $\xi<\eta$ and $\d_\xi$ is a limit of Woodin cardinals of $\P$. Then $\iota_{\xi, 0}=\d_\xi$.
$\sf{R3:}$ Suppose $\xi<\eta$ and $\d_\xi$ is neither a Woodin cardinal of $\P$ nor a limit of Woodin cardinals of $\P$. Then
\begin{center}$\iota_{\xi, 0}=\min(\dom(\vec{E}^\P)-(\d_{\xi}+1))$.\end{center}
\textbf{The Conditions Defining the Sequence $(\iota_{\xi, 1}: \xi< \eta)$ for $\xi$ as in $\sf{R3}$}\\\\
In $\sf{R4}$-$\sf{R5}$, suppose $\xi<\eta$ and $\d_\xi$ is neither a Woodin cardinal of $\P$ nor a limit of Woodin cardinals of $\P$. \\\\
$\sf{R4:}$ Suppose $rud(\P|\iota_{\xi, 0})\models ``\d_\xi$ is not a Woodin cardinal". Then $\iota_{\xi, 1}$ is the least  ordinal $\b>\iota_{\xi,0}$ such that $\rho(\P||(\b, \omega))\leq \d_\xi$.\\\\
$\sf{R5:}$ Suppose $rud(\P|\iota_{\xi, 0})\models ``\d_\xi$ is a Woodin cardinal". Then $\iota_{\xi, 1}$ is the largest ordinal $\b>\iota_{\xi,0}$ such that $\P|\b\models ``\d_\xi$ is a Woodin cardinal".\\\\
\textbf{The Conditions Defining the Sequence $(\iota_{\xi, \xi'}: \xi'\leq \nu_\xi)$ for $\xi<\eta$}\\\\
$\sf{R6:}$  Suppose $\d_\xi$ is a Woodin cardinal of $\P$ or a limit of Woodin cardinals of $\P$. Then $\iota_{\xi, 0}$ is defined as in $\sf{R1}$ and $\sf{R2}$. If $\d_\xi$ is a Woodin cardinal then $\nu_\xi=0$. If $\d_\xi$ is a limit of Woodin cardinals then $\nu_\xi=1$ and \begin{center} $\iota_{\xi, 1}=\ord(\mathcal{O}^{\P}_{\d_\xi+1, \P|\d_{\xi}+1})$.\end{center}
$\sf{R7:}$  Suppose $\d_\xi$ is neither a Woodin cardinal of $\P$ nor a limit of Woodin cardinals of $\P$. Then $\iota_{\xi, 0}$ is defined as in $\sf{R3}$, $\iota_{\xi, 1}$ is defined as in $\sf{R4}$-$\sf{R5}$, and the sequence $(\iota_{\xi, \xi'}: \xi'\in (1, \nu_\xi])$ enumerates in increasing order the closure of the set \begin{center}
$\{ \a<\d_{\xi+1}: \rho(\P||(\a, \omega))\leq \d_{\xi}\}$.
\end{center}
\textbf{When $\xi=\eta$}\\\\
$\sf{R8:}$ Suppose $\P$ is meek. If $\d^\P$ is a Woodin cardinal then $\nu_\eta=0$ and $\iota_{\eta, 0}=\ord(\P)$. If $\d^\P$ is a limit of Woodin cardinals then $\nu_\eta=1$, $\iota_{\eta, 0}=\d_\eta$ and $\iota_{\eta, 1}=\ord(\P)$.\\\\\
$\sf{R9:}$ Suppose $\P$ is non-meek and $\dom(\vec{E}^\P)-(\d_\eta+1)=\emptyset$\footnote{If $\ord(\P)=\d_\eta$ then this condition is satisfied.}. Then $\iota_{\eta, 0}=\ord(\P)$ and $\nu_\eta=0$ (in this case, we have that $\P=\mathcal{J}_{\iota_{\eta, 0}}[\P|\d_\eta]$).\\\\
$\sf{R10:}$ Suppose $\P$ is non-meek and $\dom(\vec{E}^\P)-(\d_\eta+1)\not=\emptyset$. Then $\iota_{\eta, 0}=\min(\dom(\vec{E}^\P)-(\d_\eta+1))$ and one of the following conditions holds:
\begin{enumerate}
\item If $rud(\P||\iota_{\eta, 0})\models ``\d_\eta$ is not a Woodin cardinal" then $(\iota_{\eta, \xi}: \xi\in [1, \nu_\eta])$ enumerates in increasing order the closure of the set 
\begin{center}
$\{\a\leq\ord(\P): \rho(\P||(\a, \omega))\leq \d_\eta\}$.
\end{center}
\item If $rud(\P||\iota_{\eta, 0})\models ``\d_\eta$ is a Woodin cardinal" but $\P\models ``\d_\eta$ is not a Woodin cardinal" then $\iota_{\xi, 1}$ is the largest ordinal $\b$ such that $\P|\b\models ``\d_\eta$ is a Woodin cardinal" and the sequence $(\iota_{\eta, \xi}: \xi\in [2, \nu_\eta])$ enumerates in increasing order the closure of the set
\begin{center}
$\{\a\leq \ord(\P): \rho(\P||(\a, \omega))\leq \d_\eta\}$.
\end{center}
\item If $\P|\models ``\d_\eta$ is a Woodin cardinal" then $\iota_{\eta, 1}=\ord(\P)$ and $\nu_\eta=1$.
\end{enumerate}
$\sf{R11:}$ Suppose $\P$ is gentle. Then $\nu_\eta=0$ and $\iota_{\eta, 0}=\d^\P$. \\\\
\textbf{The Definition of $(\P_{\xi, \xi'}: \xi\leq \eta, \xi'\leq \nu_\xi)$}\\\\
$\sf{R12:}$ $\P_{\xi,\xi'}=\P||\iota_{\xi, \xi'}$.\\

We say $\Q$ is a layer of $\P$ if for some $\xi\leq \eta$ and $\xi'\leq \nu_\xi$, 
\begin{center} $\Q=\P||\iota_{\xi, \xi'}$. \end{center}
We say $\Q$ is a proper layer of $\P$ if $\Q$ is a layer of $\P$ and $\Q\not=\P$. We write $\Q\insegeq_{hod} \P$ if and only if $\Q$ is a layer of $\P$, and we write  $\Q\inseg_{hod} \P$ if and only if $\Q$ is a proper layer of $\P$. $\myqedhere$
 \end{definition}
 
 \begin{remark}\label{single extender} Suppose $\P$ is an active $\sf{lses}$ such that
 \begin{enumerate}
 \item if $\a=\ord(\P)$ then $\P|\a$ is a hod-like $\sf{lses}$, 
 \item if $E=\vec{E}^\P(\a)$ then $o^{\P|\a}(\cp(E))=\d^{\P|\a}$ and
 \item $\rho(\P)>\cp(E)$.
 \end{enumerate}
  Then $\P$ itself is hod-like. It falls under clause 2c of \rdef{pre-hod-like}. Notice that $\a$ is enumerated in the $\d$-sequence of $\P$.$\myqedhere$
 \end{remark}

 
Next we introduce hod-like $\sf{lses}$. These will eventually turn into hod premice. To do this we need to impose conditions on the layers of $\sf{lses}$, which are just the members of $Y^\P$ where $\P$ is an ${\sf{lses}}$. 

  \begin{definition}[Hod-like $\sf{lses}$]\label{hod like lsp}
 Suppose $\P$ is a pre-hod-like $\sf{lses}$. We say $\P$ is \textbf{hod-like} if the following conditions hold.
 \begin{enumerate}
 \item  $\{\Q: \Q$ is a proper layer of $\P\}= (Y^\P- X^\P)$.
 \item For all layers $\Q$ of $\P$ such that $\d^\Q$ is a limit of Woodin cardinals of $\P$, $\ord(\Q)$ is a cardinal of $\P$.
 \end{enumerate}
 $\myqedhere$
 \end{definition}
 
  \begin{remark} Perhaps clause 2 of \rdef{hod like lsp} needs some more explanation. According to \rdef{layers of hod-like lsp} if $\xi$ is such that $\Q=\Q_{\xi, 0}$ then \begin{center}$\ord(\Q)=\ord(\mathcal{O}^{\P}_{\d_\xi+1, \P|\d_{\xi}})$\end{center} which is the longest initial segment of $\P$ whose strategy predicate codes a strategy for $\P|\d_{\xi}=\Q|\d^\Q$. A priori there is no reason for $\ord(\Q)$ to be a cardinal. Clause 3, following \cite{ATHM}, makes this demand. It is a fullness condition that we will have to verify every time we build a hod premouse. $\myqedhere$
\end{remark}
 \begin{remark}\label{no partial extenders with cp=bottom part} Continuing with the set up of clause 2 of \rdef{hod like lsp}, it follows that if $E\in \vec{E}^\P$ is such that $\cp(E)=\d^\Q$ then $E$ is total. For if $E$ is such that $\cp(E)=\d^\Q$ then 
 \begin{center}
 $\mathcal{O}^\P_{\d^\Q+1, \P|\d^\Q}\insegeq \P|\ind^\P(E)$.
 \end{center} But if $E$ is not total we must also have that $\P|\ind^\P(E)\insegeq \P|((\d^\Q)^+)^\P$.$\myqedhere$
\end{remark}

\begin{remark}\label{more on layers} Each $\sf{lses}$ comes with its own $Y$ predicate and the role of \rdef{layers of hod-like lsp} and  \rdef{hod like lsp} is to impose conditions that the $Y$ predicate of a hod like $\sf{lses}$ must have. One important point is that conditions like $\sf{R1}$ and $\sf{R6}$ depend on external factors. For example, in $\sf{R6}$ we demand that $\iota_{\xi, 1}=\ord(\mathcal{O}^{\P}_{\d_\xi+1, \P|\d_{\xi}+1})$ while there can be many ordinals $\a\in (\d_\xi, \iota_{\xi, 1})$ with the property that $\P|\a$ is hod-like, yet none of them determine a layer. When designing $\P$ via hod pair constructions (see \rdef{gamma-hod pair construction*}), we will need to choose $\iota_{\xi, 1}$, and its choice depends on the pointclass $\Gamma$ that we attempt to generate via the hod pair construction. Intuitively, $\iota_{\xi, 1}$ is defined to be the ordinal height of the stack of all sound $\sf{ses}$ over $\P|\d_\xi+1$ that are based on $\P|\d_{\xi}$, have a projectum $\leq \d_\xi$ and have an iteration strategy in $\Gamma$. $\myqedhere$
\end{remark}

\begin{notation}\label{l p} Suppose $\P$ is a hod-like $\sf{lses}$. Let 
\begin{center}
$L^\P=\{ \d: \exists \Q\in Y^\P-X^\P( \d^\Q=\d)\}\cup\{\d^\P\}$.
\end{center}
Let $\l^\P+1$ be the order type of $L^\P$.\index{$\l^\P$} We let $(\d_\a^\P: \a\leq\l^\P)$ be the increasing enumeration of $L^\P$.\index{$\d_\a^\P$} Also for $\xi\leq \l^\P$, set 
\begin{center}
$\P(\xi)=\cup\{\Q\in Y^\P-X^\P: \d_\xi=\d^\Q\}$.
\end{center}

We say $w=(\eta^w, \d^w)$ is a window of $\P$ if 
\begin{enumerate}
\item $\eta^w$ is the least $\eta$ such that $\d^w=o^\P(\eta^w)$ and
\item there is a layer $\Q$ of $\P$ such that $\d^w=\d^\Q$.
\end{enumerate}
 We say $w$ is the top window of $\P$ if $\d^w=\d^\P$\index{window}\index{top window}. Given a hod-like $\sf{lses}$ $\P$, we set $\ml(\P)=\cup(Y^\P)$\footnote{$\ml$ stands for \textit{maximal layer}.}.
 
 We say that $\Q$ is a \textbf{complete} layer of $\P$ if $\Q$ is a layer of $\P$ such that if $\Q$ is non-meek then there is no layer of $\R$ of $\P$ with the property that $\Q\inseg_{hod}\R$ and $\R^b=\Q^b$.
 
 If $\Q$ is a layer of $\P$ of successor type then letting $\xi$ be such that $\d^\Q=\d_{\xi+1}^\P$, $\Q^-=_{def}\Q(\xi)$ . Thus, $\Q^-$ is the longest complete layer that is in an initial segment of $\Q$. $\myqedhere$
\end{notation}

\begin{definition}[Germane $\sf{lses}$]\label{germane lses} Suppose $\M$ is an $\sf{lses}$. We say $\M$ is \textbf{germane} if letting $\a=\sup\{ \ord(\Q): \Q\in Y^\M-X^\M\}$, the following conditions hold:
\begin{enumerate}
\item If $\Q\in Y^\M-X^\M$ then $\Q$ is a hod-like $\sf{lses}$\footnote{Recall that if $\Q\in Y^\M-X^\M$ then $Y^\Q=\{\R\in Y^\M: \ord(\R)<\ord(\Q)\}$.}.
\item If $\Q\in Y^\M-X^\M$ and $\Q$ is meek then for all $\omega\b\in [\ord(\Q), \ord(\M))$, $\rho(\M||(\omega\b, \omega))>\d^\Q$.
\item $\a+1$ is a cutpoint of $\M$.
\item If $\M$ is pre-hod-like then it is hod-like.
\item One of the following conditions holds:
\begin{enumerate}
\item $\M$ is pre-hod-like.
\item $\M$ is not pre-hod like and one of the following holds:
\begin{enumerate}
\item $\a=\ord(\M)$ and $\a$ is a limit of Woodin cardinals of $\M$.
\item $\a<\ord(\M)$, $\M||\a\in Y^\M$ and $\a$ is a cardinal of $\M$ (see \rrem{remark on germane}). 
\end{enumerate}
\end{enumerate}
\end{enumerate}
If $\M$ is a germane $\sf{lses}$ then we let 
\begin{center}
${\sf{hl}}(\M)=\begin{cases}
\M&: \ \text{5.a holds}\\
\M||\a &: \text{otherwise}
\end{cases}$\end{center}
 where $\a=\sup\{ \ord(\Q): \Q\in Y^\M-X^\M\}$\footnote{$\sf{hl}$ stands for ``hod-like".}. $\myqedhere$
\end{definition}

\begin{remark}\label{remark on germane} Continuing with the set up of \rdef{germane lses}, suppose $\M$ is germane but not hod-like. Clause 5.b then says that either $\a$ is a limit of Woodin cardinals of $\M$ or $\M||\a$ is the longest hod-like initial segment of $\M$ and, moreover, it is declared to be a layer of $\M$\footnote{See $\sf{R9}$ and $\sf{R10}$ of \rdef{layers of hod-like lsp}. We demand that $\a$ be a cardinal of $\M$ because otherwise $\M$ would be pre-hod-like and hence, hod-like.}. $\myqedhere$
\end{remark}

It is not hard to create examples of germane $\sf{lses}$ that are not hod-like. For example, if $\P$ is hod-like and $\Sigma$ is its strategy then $\Sigma$-premie over $\P$ will be germane. This comment is not literally true as such premice can project in ways not allowed by  \rdef{germane lses}, but also such premice need to be re-organized into $\sf{lses}$. 

\begin{terminology}\label{types of lsa small premice} Suppose $\P$ is a hod-like $\sf{lses}$.
\begin{enumerate}
\item \textbf{(Successor type)}\index{successor type} We say $\P$ has a \textbf{successor type} if $\P$ has a top window $(\eta, \d)$ and $\eta$ is not a limit of Woodin cardinals of $\P$.
\item \textbf{(Limit type)}\index{limit type}  We say $\P$ has a \textbf{limit type} if either $\P$ doesn't have a top window or if $(\eta, \d)$ is the top window of $\P$ then $\eta$ is a limit of Woodin cardinals of $\P$.
\end{enumerate}  
If $\M$ is germane then we say $\M$ is of successor type if ${\sf{hl}}(\M)$ is of successor type hod-like $\sf{lses}$. Otherwise we say that $\M$ is of limit type. We say $\M$ is of \textbf{$b$-type}\footnote{$``b"$ stands for bottom, see below.} if $\M$ is of limit type and letting $\a=\sup\{ \ord(\Q): \Q\in Y^\M-X^\M\}$, $\a$ is not a limit of Woodin cardinals of $\M$\footnote{This means that $\M||\a$ is hod-like and is of limit type.}. $\myqedhere$
\end{terminology}

Next, we isolate the bottom part of $b$-type germane $\sf{lses}$. For non-meek hod-like $\sf{lses}$, this is essentially the part of $\P$ that is below the largest measurable limit of cutpoint Woodin cardinals.

 \begin{definition}[The bottom part of $\sf{lses}$]\label{the bottom part of lsp}\index{$\P^b$}
 Given a limit type hod-like $\sf{lses}$ $\P$ we let $\P^b=\P$ if $\P$ doesn't have a top window and otherwise, letting $(\eta, \d)$ be the top window of $\P$, we let \begin{center}
$\P^b=\P|(\eta^+)^\P$
\end{center}
where $``b"$ stands for ``bottom". We say that $\P^b$ is the \textbf{bottom} part of $\P$. It follows that $\P^b$ is a hod-like meek $\sf{lses}$ of limit type. 

Similarly, if $\M$ is germane of $b$-type then $\M^b=({\sf{hl}}(\M))^b$. $\myqedhere$
\end{definition}

\begin{definition}\label{projecting types} Suppose $\M$ is germane. We say $\M$ is \textbf{projecting well} if letting $k=k(\M)$\footnote{See \rsec{fine structure: sec}.} one of the following clauses holds:
\begin{enumerate}
\item $\M$ is of successor type and setting $\d=\d^{{\sf{hl}}(\M)}$, $\d$ is Woodin with respect to all $f:\d\rightarrow \d$ which are $r\Sigma_{k+1}^\M$-definable as witnessed by the extender sequence $\vec{E}^{\M|\d}$.
\item $\M$ is of $b$-type and $\rho_{k+1}(\M)> \d^{\M^b}$.
\item $\M$ is of limit type but not of $b$-type and $\rho_{k+1}(\M)>\ord({\sf{hl}}(\M))$.
\end{enumerate}
Otherwise we say that $\M$ \textbf{projects badly}. We say $\M$ \textbf{projects precisely} if $\M$ projects well and if there is $n$ such that $(\M, n)$ projects badly then letting $k=k(\M)$, $\M'=(\M, k+1)$ projects badly. $\myqedhere$
\end{definition}
Clearly if $\M$ projects badly then there is always an initial segment $\M'$ of $\M$ such that $\ord(\M')=\ord(\M)$ and $\M'$ projects precisely. 

\begin{remark}\label{projects precisely} We are interested in germane $\M$ that project precisely because we would like to apply stacks that are based on ${\sf{hl}}(\M)$ to $\M$ without changing the stack. 

For example, assume ${\sf{hl}}(\M)=_{def}\P$ and $\P$ is a meek hod-like $\sf{lses}$ of limit type. Suppose $\M$ projects badly. If now $E\in \vec{E}^{\P|\d^\P}$ then $Ult(\M, E)$ may have more layers than $Ult(\P, E)$, and $\P$'s strategy doesn't act on these new layers. On the other hand if $\M$ projects precisely then this is no longer the case as the functions used to compute $\pi^\M_E(\d^\P)$ and $\pi^\P_E(\d^\P)$ are the same, and they all are in $\P$.

We will use this sort of arguments later, when we need to show that if $\P$ is full, $\Sigma$ is its strategy, $\M$ is germane such that ${\sf{hl}}(\M)=\P$ and $\M$ is a $\Sigma$-mouse over $\P$ then $\M$ doesn't project badly.

Notice that our comment above concerns only to germane $\M$ which are not themselves hod-like. If $\M$ is of $b$-type and projects across $\ord({\sf{hl}}(\M))$ but it does not project badly then $\M$ itself is hod-like. 

Summarizing, if $\M$ projects precisely and $\T$ is a stack on ${\sf{hl}}(\M)$ then we define $\uparrow(\T, \M)$ just like we did in \rdef{upward extension of a stack}. $\myqedhere$
\end{remark}

\begin{definition}[Almost non-dropping stacks]\label{almost non-dropping stacks}\index{almost non-dropping stacks} Suppose $\M$ is germane of $b$-type and projects precisely. Suppose further that $\T$ is a stack on $\M$ that is based on ${\sf{hl}}(\M)$. 
We say that $\T$ is \textbf{almost non-dropping} if one of the following holds:
\begin{enumerate}
\item There is $\a\in R^\T$ such that $\pi^{\T_\leq \a}$ exists and $\T_{\geq\a}$ is above $\ord(\M_\a^b)$.
\item $\T$ has a last model and $\pi^{\T}$ exists. 
\end{enumerate}
If $\T$ is almost non-dropping and the first clause holds then let $\a(\T)$ witness it.  If $\T$ is almost non-dropping then we set
\begin{center}
$\pi^{\T, b}=\begin{cases}
\pi^\T\rest \M^b & : \pi^\T\ \text{exists}\\
\pi^{\T_{\leq \a(\T)}}\rest \M^b &:\ \text{otherwise}
\end{cases}$\end{center}
Suppose $\Sigma$ is an iteration strategy for $\M$\footnote{It is worth remembering that this entails that $\Sigma$-iterates of $\M$ have the same indexing scheme as $\M$.}. We then let 
\begin{center}
$I(\M, \Sigma)=\{(\T, \R):  \T$ is according to $\Sigma$, $\T$ is based on ${\sf{hl}(\M)}$, $\R$ is the last model of $\T$ and $\pi^{\T}$ is defined$\}$.\index{$I(\M, \Sigma)$}\\
$I^b(\M, \Sigma)=\{ (\T, \R): \T$ is according to $\Sigma$, $\T$ is based on ${\sf{hl}(\M)}$, $\R$ is the last model of $\T$ and $\pi^{\T, b}$ is defined$\}$.\index{$I^b(\M, \Sigma)$}
\end{center}
$\myqedhere$
\end{definition}

\begin{rem}
Notice that if $\T$ is almost non-dropping then it may only have drops in some image of the top window of $\P$. $\myqedhere$
\end{rem}

\begin{definition}\label{cutpoint of a stack} Suppose $\P$ is an $\sf{lses}$ and $\a<\lh(\T)$. We say $\a$ is a \textbf{cutpoint} of $\T$ if $\T_{\geq \a}$ is a stack (after trivial re-enumeration) on $\M_\a^\T$, or equivalently, if for every $\b+1\in (\a, \lh(\T)))$, $\T(\b+1)\geq \a$. $\myqedhere$
\end{definition}

%

The reader may benefit from reviewing \rnot{notation for iteration trees}.

\begin{definition}\label{proper stack} Suppose $\M$ is germane $\sf{lses}$ and 
\begin{center}
$\T=((\M_\a)_{\a<\eta}, (E_\a)_{\a<\eta-1}, D, R, (\b_\a, m_\a)_{\a\in R}, T)$
\end{center}  is a stack on $\M$ that is based on $\P=_{def}{\sf{hl}}(\M)$. We say $\T$ is a \textbf{proper stack} if the following conditions hold:
\begin{enumerate}
\item $\T$ is semi-smooth\footnote{This condition is already built into our definition of stack. See  \rrem{semi-smooth convention}. We are only making it explicit here.}.
\item $R=\{ \a: \a$ is a cutpoint of $\T\}$.
\item For all $\a\in R$ such that $\a\not=\max(R)$, if ${\sf{nc}}^\T_\a$ has a fatal drop then $\T_{\geq \a}$ is a normal stack.
\item If $\omega\b_\a<\ord(\M_\a)$ then $\M_\a||(\omega\b_\a, m_\a)$ is a non-meek hod-like $\sf{lses}$.
\item  For all $\a\in R$, if $\M_\a$ is of $b$-type\footnote{See \rnot{types of lsa small premice}.} and $\ind_\a<\d^{\M_\a^b}$ then letting $\gg$ be the least such that 
\begin{itemize}
\item $\ind_\a<\ord(\M_\a(\gg))$\footnote{See \rnot{l p} for the definition of $\P(\b)$. The definition obviously carries over to germane $\sf{lses}$.} and 
\item $\ord(\M_\a(\gg))$ is a cutpoint of $\M_\a$, 
\end{itemize}
 if $\M_\a(\gg)$ is of successor type and ${\sf{next}}^\T(\a)\in R$ then $\pi^{{\sf{nc}}^\T_\a}$ exists.
\end{enumerate}
$\myqedhere$
\end{definition}
The following lemma summarizes the properties of proper stacks.
\begin{lemma}\label{properties of proper stacks} Suppose 
\begin{center}
$\T=((\M_\a)_{\a<\eta}, (E_\a)_{\a<\eta-1}, D, R, (\b_\a, m_\a)_{\a\in R}, T)$
\end{center}
is a proper stack on a germane $\M$. Then the following conditions hold. 
\begin{enumerate}
\item  For all $\a\in R$, if $\M_\a$ is of $b$-type and $\ind_\a<\d^{\M_\a^b}$ then letting $\gg$ be the least such that 
\begin{itemize}
\item $\ind_\a<\ord(\M_\a(\gg))$ and 
\item $\ord(\M_\a(\gg))$ is a cutpoint of $\M_\a$, 
\end{itemize}
then the following conditions hold:
\begin{enumerate}
\item ${\sf{nc}}^\T_\a$ is based on $\M_\a(\gg)$.
\item If $\M_\a(\gg)$ is of successor type and ${\sf{next}}^\T(\a)\in R$ then $\pi^{{\sf{nc}}^\T_\a}$ exists and ${\sf{nc}}^\T_\a$ is above $\ord(\M_\a(\gg-1))$\footnote{This condition follows from the requirement that all cutpoints of $\T$ are in $R$. Similarly the last portion of the next clause.}.
\item If $\M_\a(\gg)$ is of limit type and ${\sf{next}}^\T(\a)\in R$ then $\pi^{{\sf{nc}}^\T_\a, b}$ exists and ${\sf{nc}}^\T_\a$ is above $\d^{\M_\a(\gg)^b}$\footnote{Here, $\pi^{{\sf{nc}}^\T_\a, b}$ is defined provided ${\sf{nc}}^\T_\a$ is on $\M_\a(\gg)$ which may not be the case. The meaning  of this here and in the sequel is that letting $\U=\downarrow({\sf{nc}}^\T_\a, \M_\a(\gg))$, $\pi^{\U, b}$ is defined.}.
\end{enumerate}
\item For all $\a\in R$, if $\M_\a$ is of $b$-type and $\ind_\a>\d^{\M_\a^b}$ then 
 \begin{enumerate}
 \item ${\sf{nc}}^\T_\a$ is above $\d^{\M_\a^b}$\footnote{This condition follows from clause 2 of \rdef{proper stack}.} and 
 \item if ${\sf{next}}^\T(\a)\in R$ then $\pi^{{\sf{nc}}^\T_\a, b}$ exists\footnote{This condition can be deduced from clause 3 of \rdef{proper stack}.}.
 \end{enumerate}
\end{enumerate}
\end{lemma}

\begin{notation}\label{sflayer notation} Suppose 
\begin{center}
$\T=((\M_\a)_{\a<\eta}, (E_\a)_{\a<\eta-1}, D, R, (\b_\a, m_\a)_{\a\in R}, T)$
\end{center}
is a proper stack on a germane $\M$. For $\a\in R$, we define ${\sf{layer}}^\T_\a$ to be the least complete\footnote{See \rnot{l p}.} layer $\N$ of $\M_\a$ such that $\ind_\a^\T\in \N$. We also let ${\sf{rnc}}^\T_\a=\downarrow({\sf{nc}}^\T_\a, {\sf{layer}}^\T_\a)$\footnote{See \rdef{restriction of a stack}.}. Often, we will represent $\T$ as
\begin{center}
$\T=((\M_\a)_{\a<\eta}, (E_\a)_{\a<\eta-1}, D, R, ({\sf{rnc}}_\a, {\sf{layer}}_\a)_{\a\in R}, (\b_\a, m_\a), T)$.
\end{center} 

Suppose $\a\in R$ is such that $\M_\a$ is of $b$-type. We say $\a$ is of \textbf{limit type}  if either $\ind_\a>\d^{\M_\a^b}$ or ${\sf{layer}}^\T_\a$ is of limit type. Otherwise we say that $\a$ is of \textbf{successor type}. We say $\a$ is of \textbf{bottom} type if $\ind_\a<\d^{\M_\a^b}$.$\myqedhere$
\end{notation}

\begin{remark}[Proper-stacks Convention]\label{proper stacks convention} In this book all stacks on hod-like $\sf{lses}$ are proper stacks. $\myqedhere$
\end{remark}

%

%

\section{The iteration embedding $\pi^{\T, b}$}\label{sec:pitb}

Recall our convention regarding stacks (see \rrem{proper stacks convention}). In this section, we define the embedding $\pi^{\T, b}$ via an inductive process. The reader may skip this section. The one important point that will come up later is the following. Suppose $\a\in R^\T$ and $\b<\lh(\T)$. Then if $\pi_{\a, \b}^{\T, b}$ is defined then its domain is $(\M_\a||(\omega\b_\a, m_\a))^b$ which in general may not be the same as $\M_\a^b$.

Assume that $\P$ is a limit type hod-like $\sf{lses}$ which isn't meek and suppose $\T$ is an iteration tree on $\P$. Again, we will not be concerned with the particular indexing scheme that $\P$ has. In some cases, regardless of whether $\T$ has a last model or not, it is possible to extract an embedding out of the iteration embeddings given by $\T$ that acts on $\P^b$. We describe this embedding below. First we define it by assuming that $\T$ is a normal iteration tree and then extend the definition to stacks. Recall that our $\sf{lses}$ are lsa-small (see \rdef{lsa small}).

\begin{definition}\label{embedding on the top for stacks}
Suppose $\P$ is a non-meek hod-like $\sf{lses}$\footnote{In particular, $\P\not=\P^b$.}. Suppose 
\begin{center}
$\T=((\M_\a)_{\a<\eta}, (E_\a)_{\a<\eta-1}, D, R, (\b_\a, m_\a), T)$
\end{center}
is a normal iteration tree on $\P$. We define $(\pi^{\T, b}_{\a, \a'}: \a<\a'<\eta \wedge (\a, \a')\in T)$ by induction maintaining that if $\pi^{\T, b}_{\a, \a'}\not =\emptyset$ then \\\\
(a) if $\a\not \in R$ then $\pi^{\T, b}_{\a, \a'}: \M_\a^b\rightarrow \M_{\a'}^b$ is an elementary embedding, and\\
(b) if $\a\in R$ then $\pi^{\T, b}_{\a, \a'}:(\M_\a||(\omega\b_\a, m_\a))^b\rightarrow \M_{\a'}^b$ is an elementary embedding. \\\\
\textbf{The successor case}\\\\
Suppose $\b+1<\eta$  and we have defined  $(\pi^{\T, b}_{\a, \a'}: \a<\a'\leq \b \wedge (\a, \a')\in T)$.  Let $\gg'=\T(\b+1)$. For $\gg<\b+1$ such that $(\gg, \b+1)\in T$, we define $\pi^{\T, b}_{\gg, \b+1}$ as follows. 
\begin{enumerate}
\item If $\pi^{\T, b}_{0, \gg'}=\emptyset$ then set $\pi^{\T, b}_{\gg, \b+1}=\emptyset$. \\\\
In the next three clauses we assume that $\pi^{\T, b}_{0, \gg'}\not =\emptyset$ and that $\gg'\not \in R$.
\item Suppose $\cp(E_\b)>\ord(\M_{\gg'}^b)$. Then set $\pi^{\T, b}_{\gg, \b+1}=\pi^{\T, b}_{\gg, \gg'}$.
\item Suppose $\cp(E_\b)\leq \ord(\M_{\gg'}^b)$ and $\b+1\in D$. Then $\pi^{\T, b}_{\gg, \b+1}=\emptyset$.
\item Suppose $\cp(E_\b)\leq \ord(\M_{\gg'}^b)$ and $\b+1\not \in D$. Then $\pi^{\T, b}_{\gg, \b+1}=(\pi^\T_{\gg', \b+1}\rest \M_{\gg'}^b)\circ \pi^{\T, b}_{\gg, \gg'}$.\\\\
In the next three clauses we assume that $\pi^{\T, b}_{0, \gg'}\not =\emptyset$ and that $\gg'\in R$.
\item Suppose $\cp(E_\b)>\ord((\M_{\gg'}||(\omega\b_{\gg'}, m_\gg'))^b)$. Then set $\pi^{\T, b}_{\gg, \b+1}=\pi^{\T, b}_{\gg, \gg'}$.
\item Suppose $\cp(E_\b)\leq \ord((\M_{\gg'}||(\omega\b_{\gg'}, m_\gg'))^b)$ and $\b+1\in D$. Then $\pi^{\T, b}_{\gg, \b+1}=\emptyset$.
\item Suppose $\cp(E_\b)\leq \ord((\M_{\gg'}||(\omega\b_{\gg'}, m_\gg'))^b)$ and $\b+1\not \in D$. Then $\pi^{\T, b}_{\gg, \b+1}=(\pi^\T_{\gg', \b+1}\rest \M_{\gg'}^b)\circ \pi^{\T, b}_{\gg, \gg'}$.\\
\end{enumerate}
\textbf{The limit case}\\\\
Suppose next that $\b<\eta$ is a limit ordinal and we have defined $(\pi^{\T, b}_{\a, \a'}: \a<\a'< \b \wedge (\a, \a')\in T)$. Then we define $\pi^{\T, b}_{\gg, \b}$ for $\gg\in [0, \b)_\T$ according to the following cases:
\begin{enumerate}
\item If $\gg\in [0, \b)_\T$ is such that there is $\gg'\in [0, \b)_\T$ with the property that $\pi^{\T, b}_{\gg, \gg'}=\emptyset$ then $\pi^{\T, b}_{\gg, \b}=\emptyset$.
\item Suppose $\gg\in [0, \b)_\T$ is such that for all $\gg'\in [0, \b)_\T$, $\pi^{\T, b}_{\gg, \gg'}$ is defined. Then letting $\nu+1\in b$ be such that $T(\nu+1)=\gg$, $\pi^{\T, b}_{\gg, \b}=\pi^{\T, b}_{\nu+1, \b}\circ \pi^\T_{\gg, \nu+1}$ where $\pi^{\T, b}_{\nu+1, \b}$ is the direct limit embedding given by the directed system $(\M_\xi^b, \pi_{\xi, \xi'}^{\T, b}: \nu+1\leq\xi<\xi' \wedge (\xi, \xi')\in c^2)$ where $c=[0, \b)_\T$.\\
\end{enumerate}
\textbf{The iteration embedding $\pi^{\T, b}$}\\\\
Continuing with the $\T$ above, we let $\pi^{\T, b}$ be defined according to the following clauses.
\begin{enumerate}
\item Suppose $\lh(\T)=\gg+1$. Then set $\pi^{\T, b}=\pi^{\T, b}_{0, \gg}$. \\\\
In all the clauses below we assume that $\lh(\T)$ is a limit ordinal.
\item Suppose there is $\gg<\lh(\T)$ such that $\pi^{\T, b}_{0, \gg}=\emptyset$ and $\T_{\geq \gg}$ is a normal stack on $\M_\gg$. Then set $\pi^{\T, b}=\emptyset$. \\\\
In all the clauses below we assume that if $\gg<\lh(\T)$ is such that $\T_{\geq \gg}$ is a normal stack on $\M_\gg$ then $\pi^{\T, b}_{0, \gg}\not =\emptyset$. 
\item Suppose there is $\gg<\lh(\T)$ such that $\T_{\geq \gg}$ is a normal stack on $\M_\gg$ based on $\M^b_\gg$. Then set $\pi^{\T, b}=\emptyset$. 
\item Suppose there is $\gg<\lh(\T)$ such that $\T_{\geq \gg}$ is a normal stack on $\M_\gg$ above $\ord(\M_\gg^b)$. Then set $\pi^{\T, b}=\pi^{\T, b}_{0, \gg}$. 
\item Suppose there is a cofinal $c\subseteq \lh(\T)$ such that $\{ \gg<\lh(\T): \exists \gg' \in c( (\gg, \gg')\in T)\}$ is a well-founded branch of $\T$ and for all $\gg<\gg'$ with $(\gg, \gg')\in c^2$, $\pi_{0, \gg}^{\T, b}\not =\pi_{0, \gg'}^{\T, b}$. Then set $\pi^{\T, b}=\pi^{\T}_c\rest \P^b$.
\end{enumerate}
Given $(\a, \a')\in T$, we say $\pi^{\T, b}_{\a, \a'}$ is \textbf{defined} or \textbf{exists} if $\pi^{\T, b}_{\a, \a'}\not=\emptyset$. Similarly we say $\pi^{\T, b}$ is \textbf{defined} or \textbf{exists} if $\pi^{\T, b}\not=\emptyset$. $\myqedhere$
\end{definition}

\begin{remark}\label{pitb is an iteration} Suppose $\P$ is a non-meek hod-like $\sf{lses}$ and $\T$ is a stack on $\P$. For all $\a<\a'$ such that $(\a, \a')\in T$, if $\pi^{\T, b}_{\a, \a'}$ exists then it is essentially the iteration embedding. However, given how $\pi_{\a, \a'}^\T$ is defined, it is possible that $\pi^{\T, b}_{\a, \a'}$ exists yet $\pi_{\a, \a'}^\T$ is undefined. 

In general, we have that $\pi^{\T, b}_{\a, \a'}$ is defined if and only if for all $\gg$ such that $\gg+1\in [\a, \a')_\T\cap D^\T$, $\cp(E_\gg)>\ord(\M_\gg^b)$. $\myqedhere$
\end{remark}

Notice that in \rdef{embedding on the top for stacks} we are not assuming that the stack has a last model. The fragment of the eventual iteration embedding $\pi^{\T}$ restricted to $\P^{b}$ can be seen without actually having the last branch. 


%
\section{Canonical singularizing sequence}

The following notion will be used throughout this paper. 

\begin{definition}[Canonical singularizing sequences]\label{canonical singularizing sequences} Suppose $\P$ is a germane $\sf{lses}$ of $b$-type that projects precisely and $\T$ is an almost non-dropping stack on $\P$. Let $\Q=\pi^{\T, b}(\P^{b})$. Then $\Q$ is a hod-like $\sf{lses}$. If $w=(\eta, \d)$ is a window of $\Q$ then we let 
\begin{center}
$s(\T, w)=\{ \a: \exists a\in \eta^{<\omega} \exists f\in \P^b( \a=\pi^{\T, b}(f)(a))\}\cap \d$
\end{center}
$\myqedhere$
\end{definition}

The following is an easy lemma, which is a consequence of our assumption that all hod-like ${\sf{lses}}$ are lsa small. It traces back to the fact that if $\P$ is a hod-like ${\sf{lses}}$ and $E\in \vec{E}^\P$ is an extender such that $\cp(E)=\d^\Q$ for some layer $\Q$ of $\P$ and $\nu\in [\cp(E), \ind^\P(E)]$ then $Ult(\P, E)\models ``\nu$ is not a Woodin cardinal". 


The following definition will be used in the next few lemmas.

\begin{definition}\label{t-critical}
Suppose  $\T$ is an almost non-dropping stack on (an lsa small) germane, $b$-type $\sf{lses}$ $\P$ that projects precisely, $\lh(\T)=\a+1$ and $\Q=\pi^{\T, b}(\P^b)$\footnote{We drop $\T$ from our notation.}. Suppose $\xi$ is a cardinal of $\Q$. Let  $\iota\leq \a$ be the least such that $\pi^{\T_{\leq \iota+1}, b}$ is defined, $\iota+1\in [0, \a]_\T$ and for some $\xi' \in \M_{\iota+1}^b$, $\pi_{\iota+1, \a}(\xi')=\xi$\footnote{This embedding may not be defined, but because both $\pi^{\T, b}$ and $\pi^{\T_{\leq \iota+1}, b}$ are defined, $\pi_{\iota+1, \a}\rest \M_{\iota+1}^b$ is meaningful.}. We say $\xi$ is $\T$-critical if $\xi'=\cp(E_{\iota}^\T)$. $\myqedhere$
\end{definition}

Given $A\subseteq X\times Y$ and $x\in X$, we set $A_x=\{y: (x, y)\in A\}$.

\begin{lemma}\label{t-critical sequence} Suppose  $\T$ is an almost non-dropping stack on (an lsa small) germane, $b$-type $\sf{lses}$ $\P$ that projects precisely, $\lh(\T)=\a+1$ and $\Q=\pi^{\T, b}(\P^b)$. Suppose $\xi$ is $\T$-critical. Then there is a finite sequence $(\gg_i, \gg_i', \xi_i: i\leq n+1)$ such that 
\begin{enumerate}
\item $(\gg_i :i\leq n)$ is increasing and for each $i\leq n$, $\gg_i\in [0, \a]_\T$,
\item $\gg_{n+1}=\gg_{n+1}'=\a$ and $\xi_{n+1}=\xi$,
\item for every $i\leq n$, $\gg_i=\T(\gg_i'+1)$,
\item $\xi_0$ is not $\T_{\leq \gg_0}$-critical,
\item for every $i\leq n$, $\pi^{\T_{\leq \gg_i}, b}$ is defined,
\item for every $i\leq n$, $\xi_i=\cp(E_{\gg_i'})$,
\item for every $i$ such that $i+1\leq n$, $\xi_{i+1}=\pi_{\gg_i'+1, \gg_{i+1}}(\xi_i)$,
\item $\xi=\pi_{\gg_n'+1, \a}(\xi_n)$,
\item for every $i\leq n$, 
\begin{center}
$\powerset(\xi_{i+1})\cap \M_{\gg_{i+1}}^\T=\{\pi_{\gg_i'+1, \gg_{i+1}}^\T(g)(t): g\in \M_{\gg'_i+1}^\T|(\xi_i^+)^{\M^\T_{\gg'_i+1}} \wedge t\in [\xi_{i+1}]^{<\omega}\}$.
\end{center}
\item for every $i\leq n+1$, for every $(m_0, m_1, ...,m_k)\in \mathbb{N}^{<\omega}$ and for every $A\in \M^\T_{\gg_i}$ such that
\begin{center}
$A\subseteq [\xi_i]^{m_0}\times [\xi_i]^{m_1}\times...\times[\xi_i]^{m_k}$
\end{center}
 there is $B\in \M_{\gg_0}^\T|(\xi_0^+)^{\M_{\gg_0}^\T}$ and $t\in [\xi_i]^{<\omega}$ such that 
\begin{center}
$A=\pi_{\gg_0, \gg_{i}}^\T(B)_t\cap ([\xi_i]^{m_0}\times [\xi_i]^{m_1}\times...\times[\xi_i]^{m_k})$.
\end{center}
\end{enumerate}
\end{lemma}
\begin{proof} We first get a finite sequence satisfying clauses 1-8, and then show that any such sequence also satisfies clauses 9 and 10. Because $\xi$ is $\T$-critical, we have some $(\iota, \xi')$ satisfying the clauses of \rdef{t-critical}. Let $\gg=\T(\iota+1)$. The claim now can be proven by induction. Assuming our claim is true for $\T_{\leq \gg}$ we have two cases. Suppose first that $\xi'$ is not $\T_{\leq \gg}$-critical. Set then $n=1$, $\gg_0=\gg$, $\gg'_0=\iota$ and $\xi_0=\xi'$. Otherwise let $(\gg_i, \gg_i', \xi_i: i\leq m)$ witness the claim for the pair $(\xi', \T_{\leq \gg'})$. Then set $n=m+2$, $\gg_{m+1}=\gg$, $\xi_{m+1}=\xi'$ and $\gg_{m+1}'=\iota$. This finishes the proof that there is a finite sequence satisfying clauses 1-8.

We now want to show that any sequence that satisfies clauses 1-8 also satisfies clauses 9 and 10. Let then $(\gg_i, \gg_i', \xi_i: i\leq n+1)$ be a sequence satisfying clauses 1-8. 9 is easy to show as the generators of $\pi_{\gg_i'+1, \gg_{i+1}}^\T\rest \M_{\gg'_i+1}^\T|(\xi_i^+)^{\M^\T_{\gg'_i+1}}$ are contained in $\xi_{i+1}=\pi_{\gg_i'+1, \gg_{i+1}}^\T(\xi_i)$. 

We now show clause 10. Fix $i+1\leq n+1$. Without loss of generality we can assume that clause 10 holds for all $j\leq i$. We then want to prove it for $i+1$.
Below we drop $\T$ from superscripts. Fix $(m_0, m_1, ...,m_k)\in \mathbb{N}^{<\omega}$ and $A\in \M_{\gg_{i+1}}$ such that
\begin{center}
$A\subseteq [\xi_{i+1}]^{m_0}\times [\xi_{i+1}]^{m_1}\times...\times[\xi_{i+1}]^{m_k}$.
\end{center}
 We want to find a $B\in \M_{\gg_0}|(\xi_0^+)^{\M_{\gg_0}}$ and $t\in [\xi_{i+1}]^{<\omega}$ such that 
\begin{center}
$A=\pi_{\gg_0, \gg_{i+1}}(B)_t\cap ([\xi_{i+1}]^{m_0}\times [\xi_{i+1}]^{m_1}\times...\times[\xi_{i+1}]^{m_k})$.
\end{center}
It follows from clause (9) that  for some $u\in [\xi_{i+1}]^{<\omega}$, $A=\pi_{\gg_i'+1, \gg_{i+1}}(g)(u)$. Let $p=\card{u}$. Notice that 
\begin{center}
$g:[\xi_i]^p\rightarrow \powerset([\xi_{i}]^{m_0}\times [\xi_{i}]^{m_1}\times...\times[\xi_{i}]^{m_k})$.
\end{center}
Let then $G\subseteq [\xi_i]^p\times [\xi_{i}]^{m_0}\times [\xi_{i}]^{m_1}\times...\times[\xi_{i}]^{m_k}$ be given by 
\begin{center}
$(x, y)\in G \iff y\in g(x)$. 
\end{center}
We thus have that\\\\
(1) for all $(x, y) \in  [\xi_i]^p\times [\xi_{i}]^{m_0}\times [\xi_{i}]^{m_1}\times...\times[\xi_{i}]^{m_k}$, 
\begin{center}
$y\in g(x) \iff (x, y)\in \pi_{\gg_i, \gg_i'+1}(G)$,\\
\end{center}
implying that\\\\
(2) for all $x \in  [\xi_i]^p$, 
\begin{center}
$g(x)=\pi_{\gg_i, \gg_i'+1}(G)_x\cap  [\xi_{i}]^{m_0}\times [\xi_{i}]^{m_1}\times...\times[\xi_{i}]^{m_k}$.\\
\end{center}
Because clause 10 holds for $i$, we get some $H\in \M_{\gg_0}|(\xi_0^+)^{\M_{\gg_0}}$ and $s\in [\xi_i]^{<\omega}$ such that
\begin{center}
$G=\pi_{\gg_0, \gg_i}(H)_s\cap [\xi_i]^p\times [\xi_{i}]^{m_0}\times [\xi_{i}]^{m_1}\times...\times[\xi_{i}]^{m_k}$.
\end{center}
Combining the above with (2) we get that \\\\
(3) for all $x \in  [\xi_i]^p$, 
\begin{center}
$g(x)=(\pi_{\gg_0, \gg_i'+1}(H)_s)_x\cap [\xi_{i}]^{m_0}\times [\xi_{i}]^{m_1}\times...\times[\xi_{i}]^{m_k}$.\\
\end{center}
Applying $\pi_{\gg_i'+1, \gg_{i+1}}$ to the equation above and recalling that $A=\pi_{\gg_i'+1, \gg_{i+1}}(g)(u)$ for some $u\in [\xi_{i+1}]^{<\omega}$, we get that
\begin{center}
$A=(\pi_{\gg_0, \gg_{i+1}}(H)_s)_u\cap ([\xi_{i+1}]^{m_0}\times [\xi_{i+1}]^{m_1}\times...\times[\xi_{i+1}]^{m_k})$
\end{center}
Because both $s, u\in [\xi_{i+1}]^{<\omega}$, we now can find some $B$ and $t$ such that 
\begin{center}
$A=\pi_{\gg_0, \gg_{i+1}}(B)_t\cap ( [\xi_{i+1}]^{m_0}\times [\xi_{i+1}]^{m_1}\times...\times[\xi_{i+1}]^{m_k})$
\end{center}
with $t\in [\xi_{i+1}]^{<\omega}$.

\end{proof}

\begin{lemma}\label{towards the main lemma}
 Suppose $\T$ is an almost non-dropping stack on (an lsa small) germane, $b$-type $\sf{lses}$ $\P$ that projects precisely and $\Q=\pi^{\T, b}(\P^b)$. Suppose $\xi$ is a limit of Woodin cardinals of $\Q$ and $A\in \Q$. Then the following holds.
 \begin{enumerate}
 \item Suppose $A\in \powerset(\xi)$ and $\xi$ is not $\T$-critical. Then $A=\pi^{\T, b}(B)_t$ where $B\in \P^b$ and $t\in [\xi]^{<\omega}$.
 \item Suppose for some $(m_0, m_1, ..., m_{k})\in \mathbb{N}^{<\omega}$, 
\begin{center}
$A\subseteq [\xi]^{m_0}\times [\xi]^{m_1}\times...\times[\xi]^{m_k}$.
\end{center}
 and $\xi$ is $\T$-critical. Then there is $B\in \P^b$ and $t\in [\xi]^{<\omega}$ such that 
\begin{center}
$A=\pi^{\T, b}(B)_t\cap ([\xi]^{m_0}\times [\xi]^{m_1}\times...\times[\xi]^{m_k})$.
\end{center}
 \end{enumerate}
 \end{lemma}
 \begin{proof} We drop $\T$ from superscripts.  Let $\a$ be the least such that $\Q=\M_{\a}^b$. Notice that because $\pi^{\T, b}$ is defined, $\pi_{0, \a}\rest \P^b$ makes sense and is equal to $\pi^{\T, b}$. Therefore, we will simply use $\pi_{\iota, \iota'}$ as if it is defined on all of $\M_\iota$. To prove our claim, we may just as well assume, without losing generality, that $\a+1=\lh(\T)$.
 
 Towards a contradiction suppose our claim is false. Without loss of generality we may assume that\\\\
 (*) for every $\iota$ such that $\iota+1<\lh(\T)$ and $\pi^{\T_{\leq \iota}, b}$ is defined, for every $\nu$ which is a limit of Woodin cardinals of $\Q'=_{def}\pi^{\T_{\leq \iota}, b}(\P^b)$, and for every $C\in \Q'$ the following holds: 
  \begin{enumerate}
 \item Suppose $C\in \powerset(\nu)$ and $\nu$ is a not $\T$-critical. Then $C=\pi^{\T_{\leq \iota}, b}(D)_t$ where $D\in \P^b$ and $t\in [\nu]^{<\omega}$.
 \item Suppose $(n_0, n_1, ..., n_{k})\in \mathbb{N}^{<\omega}$ is such that
\begin{center}
$C\subseteq [\nu]^{n_0}\times [\nu]^{n_1}\times...\times[\nu]^{n_k}$.
\end{center}
 and $\nu$ is $\T$-critical. Then there is $D\in \P^b$ and $t\in [\xi]^{<\omega}$ such that 
\begin{center}
$C=\pi^{\T_{\leq \iota}, b}(D)_t\cap ([\nu]^{n_0}\times [\nu]^{n_1}\times...\times[\nu]^{n_k})$.
\end{center}
 \end{enumerate}
 
 A direct limit argument then shows that  $\alpha=\b+1$. Let $\k=\cp(E_\b)$ and let $\gg=\T(\a)$.  Notice that if $\xi<\k$ then our claim follows from (*)\footnote{Notice that we must have that $\pi^{\T_{\leq \gg}, b}$ is defined as otherwise $\pi^{\T, b}$ cannot be defined.}. Thus, we assume that $\xi\geq \k$.  \\
 
 \textbf{Assume first that $\xi\leq \ind_\b$.}\\\\
  Because $\xi$ is a limit of Woodin cardinals of $\Q$ and there are no Woodin cardinals of $\Q$ in the interval $(\k, \ind_\b]$\footnote{This is consequence of the fact that $\P$ is lsa small.}, we have that in fact $\xi=\k$. Because $\k=\cp(E_\b)$, we have that $\xi$ is $\T$-critical.  Applying \rlem{t-critical sequence} to $(\xi, \T)$ we get a finite sequence $(\gg_i, \gg_i', \xi_i: i\leq n)$ satisfying the clauses of \rlem{t-critical sequence}. In particular, $\xi=\xi_{n+1}$, $\xi_n=\k$, $\gg_n=\gg$, $\gg_n'=\b$ and $\gg_{n+1}=\a$. Clause 10 of \rlem{t-critical sequence} implies that there is $B'\in \M_{\gg_0}|(\xi_0^+)^{\M_{\gg_0}}$ and $s\in [\xi]^{<\omega}$ such that\\\\
  (1)  $A=\pi_{\gg_0, \a}(B')_s\cap ([\xi]^{m_0}\times [\xi]^{m_1}\times...\times[\xi]^{m_k})$.\\\\
   Because $\xi_0$ is not $\T_{\leq \gg_0}$-critical and because $\xi_0$ is a limit of Woodin cardinals of $\M_{\gg_0}$, applying (*) to $(\xi_0, \T_{\leq \gg_0})$ we get some $B''\in \P^b$ and $s'\in [\xi_0]^{<\omega}$ such that\\\\
   (2) $\pi_{0, \gg_0}(B'')_{s'}=B'$.\\\\
    Putting (1) and (2) together and rearranging $B''$ we get some $B\in \P^b$ and $t\in [\xi]^{<\omega}$ such that
\begin{center}
$A=\pi^{\T, b}(B)_t\cap ([\xi]^{m_0}\times [\xi]^{m_1}\times...\times[\xi]^{m_k})$.
\end{center}
  
  \textbf{Assume now that $\xi>\ind_\b$.}\\\\ 
 Let $\l$ be the least such that $\pi_{\gg, \a}(\l)\geq \xi$. Because $\xi$ is a limit of Woodin cardinals of $\Q$, we have that \\\\
 (1) $\M_{\gg}\models ``\l$ is a limit of Woodin cardinals".\\\\
 We now have some $g\in \M_\gg$, $g:\kappa\rightarrow \l$ such that $\pi_{\gg, \a}(g)(s)=A$ where $s\in [\ind_\b]^{<\omega}\subset [\xi]^{<\omega}$. 
 
 Suppose first that $\l$ is not $\T$-critical. Since $\l$ is not $\T$-critical, applying (*) to $(\l, \T_{\leq \gg})$, we get some $f\in \P^b$ and some $t\in [\l]^{\omega}$ such that $g=\pi_{0, \gg}(f)(t)$. Therefore, $A=\pi_{0, \a}(f)(u)(s)$ where $u=\pi_{\gg, \a}(t)$. But because $u\in [\xi]^{<\omega}$, we can find some $f^*\in \P^b$ and some $u^*\in [\xi]^{<\omega}$ such that $A=\pi_{0, \a}(f^*)(u^*)$. Rearranging $f^*$ we get some $B\in \P^b$ and $t\in [\xi]^{<\omega}$ such that $A=\pi_{0, \a}(B)_t$.

Finally, suppose that $\l$ is $\T$-critical. In this case, $\l$ is a regular cardinal of $\M_\gg$ and so we have two cases. Either $\l=\k$ or $\pi_{\gg, \a}(\l)=\xi$. 
 In both cases, we have some $B'\in \M_\gg|(\l^+)^{\M_\gg}$ such that $A=\pi_{\gg, \a}(B')_t$ for $t\in [\xi]^{<\omega}$. 
 
We now have two cases. Suppose first that $\pi_{\gg, \a}(\l)=\xi$. Applying (*) to $(\l, \T_{\leq \gg})$ we get some $B''\in \P^b$ and some $t'\in [\l]^{<\omega}$ such that 
\begin{center}
$B'=\pi_{0, \gg}(B'')_{t'}\cap ([\l]^{\card{t}}\times [\l]^{m_0}\times [\l]^{m_1}\times...\times[\l]^{m_k})$.
\end{center}
and so rearranging $B''$ we get some $B\in \P^b$ and $s\in [\xi]^{<\omega}$ such that 
\begin{center}
 $A=\pi_{0, \a}(B)_s\cap ([\xi]^{m_0}\times [\xi]^{m_1}\times...\times[\xi]^{m_k})$.
\end{center}

Suppose next that $\pi_{\gg, \a}(\l)>\xi$. As $\l$ is a regular cardinal of $\M_\gg$ this is only possible if $\l=\k$. Since $\k$ is a $\T$-critical point we have that there is some $B''\in \P^b$ and some $t'\in [\k]^{<\omega}$
\begin{center}
$B'=\pi_{0, \gg}(B'')_{t'}\cap ([\k]^{\card{t}}\times [\k]^{m_0}\times [\k]^{m_1}\times...\times[\k]^{m_k})$.
\end{center}
Therefore, we get that 
\begin{center}
$A=(\pi_{0, \a}(B'')_{t'}\cap ([\pi_{\gg, \a}(\k)]^{\card{t}}\times [\pi_{\gg, \a}(\k)]^{m_0}\times ...\times [\pi_{\gg, \a}(\k)]^{m_k})_t$.
\end{center}
Therefore, for some $B'''\in \P^b$ and $t''\in [\xi]^{<\omega}$,
\begin{center}
$A=\pi_{0, \a}(B''')_{t''}\cap ([\xi]^{m_0}\times [\xi]^{m_1}\times...\times[\xi]^{m_k})$.
\end{center}
Since $\xi\in (\ind_\b, \pi_{\gg, \a}(\k))$, we have that $\xi=\pi_{\gg, \a}(g)(u)$ for some $u\in [\xi]^{<\omega}$. Therefore, rearranging $B'''$ we get some $B\in \P^b$ and $s\in [\xi]^{<\omega}$ such that 
$A=\pi_{0, \a}(B)_{s}$.

 \end{proof}

\begin{lemma}\label{the canonical singularizing sequence exists} Suppose $\P$ is a germane, $b$-type $\sf{lses}$ that projects precisely and $\T$ is an almost non-dropping stack on $\P$. Let $\Q=\pi^{\T, b}(\P^b)$. Then for any window $w=(\eta, \d)$  of $\Q$ (see \rnot{l p}) such that  $\Q\models ``\d$ is a Woodin cardinal",
\begin{center}
$\sup(s(\T, w))=\d$.
\end{center}
\end{lemma}
\begin{proof} We drop $\T$ from superscripts. Let $\a^*$ be the least such that $\Q=\M_{\a^*}^b$. To prove our claim, we may just as well assume, without losing generality, that $\a^*+1=\lh(\T)$. 

Suppose to the contrary that $w=(\eta, \d)$ is a window of $\Q$ such that $\Q\models ``\d$ is a Woodin cardinal" but $\sup(s(\T, w))<\d$. Let then $\a$ be the least $\a'$ such that $\d<\d^{(\M_{\a'})^b}$ and $\Q|\d=\M_{\a'}|\d$. As $\pi_{\a, \a^*}\rest \d+1=id$\footnote{Notice that while $\pi_{\a, \a^*}$ may not exist, $\pi_{\a, \a^*}\rest \M_\a^b$ must be defined, and so the use of $\pi_{\a, \a^*}$ is justified.}, we can now assume without loss of generality that $\a=\a^*$. 

We can also assume, without loss of generality, that\\\\
(*) for any $\b\in [0, \a)_\T$, $\d\not \in \rge(\pi_{\b, \a})$\footnote{Notice that while $\pi_{\b, \a}$ may not be defined, it nevertheless is defined on $\M_\b^b$, and so here and in the sequel we will ignore the fact that $\pi_{\b, \a}$ may not be defined.}.\\\\ This is because for any such $\b$, letting $\d'$ be such that $\pi_{\b, \a}(\d')=\d$, we must have that $\sup(\pi_{\b, \a}[\d'])=\d$\footnote{This is a consequence of the fact that because we are only considering lsa small hod-like ${\sf{lses}}$, $\d'$ is not a critical point of any $E\in \vec{E}^{\M_\b}$.}.

Because $\d$ has no pre-image in any $\M_\b$, it must be the case that $\a=\b+1$ for some $\b$. Let $\gg=\T(\b+1)$. We thus have that $\M_\a=Ult(\M_\gg, E_\b)$\footnote{Notice that $E_\b$ cannot cause a drop as we are assuming that $\pi^{\T, b}$ exists.} and that\\\\
(1) $\d\not \in \rge(\pi_{\gg, \a})$ and hence, $\d>\cp(E_\b)$\footnote{$\d=\cp(E_\b)$ is not possible because of lsa smallness.}.\\\\
Furthermore, notice that\\\\
(2) $\d>\ind_\b$\\\\ as otherwise in the case that $\d=\ind_\b$ we have that $\d$ is a successor cardinal in $\M_\a$ and hence not a Woodin cardinal, or in the case that $\d\in (\cp(E_\b), \ind_\b)$ we have that $\d$ is not a Woodin cardinal in $\M_\a$ as it is not a Woodin cardinal in $Ult(\M_\b, E_\b)$ and $\M_\a\cap \powerset(\d)=Ult(\M_\b, E_\b)\cap \powerset(\d)$.

Notice next that (2), and more relevantly the argument used to establish (2), also implies that\\\\
(3) $\eta>\ind_\b$.\\\\
It then follows that if $\xi$ is least such that $\pi_{\gg, \a}(\xi)>\eta$ then
\begin{center} $\d=\sup(\{\pi_{\gg, \a}(f)(s): f\in \M_\gg$, $f:\cp(E_\b)^{\card{s}}\rightarrow \xi$ and $s\in [\ind_\b]^{<\omega}\}\cap \d)$\end{center}
 and therefore, it follows from (3) that\\\\
(4) $\d=\sup(\{\pi_{\gg, \a}(f)(s): f\in \M_\gg$, $f:\cp(E_\b)^{\card{s}}\rightarrow \xi$ and $s\in [\eta]^{<\omega}\}\cap \d)$.\\\\
Notice that it follows from our choice of $\xi$ that $\M_\gg\models ``\xi$ is a cutpoint limit of Woodin cardinals". (4) and \rlem{towards the main lemma} now give a contradiction\footnote{We apply \rlem{towards the main lemma} to functions $f$ used in (4).}, as we get that $\sup(s(\T, w))=\d$.

Below we calculate the details in the case $\xi$ is $\T$-critical. In this case, if $f\in \M_\gg$ is such that $f:\cp(E_\b)^{\card{s}}\rightarrow \xi$ then letting $F$ be the graph of $f$, we can find $G\in \P^b$ and $t\in [\xi]^{<\omega}$ such that $F=\pi_{0, \gg}(G)_t\cap [\xi]^{\card{s}}\times [\xi]$. But then for every $u\in [\ind_\b]^{<\omega}$, if $\pi_{\gg, \a}(f)(u)$ is defined then it is the unique $x$ such that
\begin{center}
 $(s, x)\in \pi_{0, \a}(G)_t\cap [\pi_{\gg, \a}(\xi)]^{\card{s}}\times [\pi_{\gg, \a}(\xi)]$. 
 \end{center}
 Setting then $g(v)=G_v$, we get that for every $u\in [\ind_\b]^{<\omega}$, if $\pi_{\gg, \a}(f)(u)$ is defined then it is equal to $\pi_{0, \a}(g)(t)(u)$. It then follows that there is $h\in \P^b$ such that for every $u\in [\ind_\b]^{<\omega}$ there is $u'\in [\ind_\b]^{<\omega}$ such that if  $\pi_{\gg, \a}(f)(u)$ is defined then $\pi_{\gg, \a}(f)(u)=\pi_{0, \a}(h)(u')$. It then follows from (4) that $\sup(s(\T, w))=\d$.
\end{proof}

\section{The un-dropping game}\label{undropping game sec}

Recall our convention regarding proper stacks (see \rrem{proper stacks convention}).
Before we proceed, we explain the meaning of the un-dropping game. Suppose we are comparing the strategies of two lsa type hod-like $\sf{lses}$ $\P$ and $\Q$. Let $\Sigma$ be the strategy of $\P$ and $\Lambda$ be the strategy of $\Q$. Let us assume that the pointclasses generated by $(\P, \Sigma)$ and $(\Q, \Lambda)$ are the same. We are then searching for $\R$ which is an iterate of $\P$ and $\Q$ and $\Sigma_\R=\Lambda_\R$. In this comparison we might be forced to consider iteration trees $\T$ and $\U$ with last models $\M$ and $\N$ such that $\pi^\T$ and $\pi^\U$ don't exist and for some $\K\insegeq_{hod}\M$ and $\K\insegeq_{hod} \N$, $\Sigma_{\K}\not= \Lambda_{\K}$. We can continue the comparison by comparing $(\M, \Sigma_\M)$ and $(\N, \Lambda_\N)$ and producing $(\S, \Phi)$ which is a common tail of $(\M, \Sigma_\M)$ and $(\N, \Lambda_\N)$. However, $(\S, \Phi)$ cannot be thought of as a last model of a successful comparison of $(\P, \Sigma)$ and $(\Q, \Lambda)$ simply because $\pi^\T$ and $\pi^\U$ do not exist. What we need to do is to compare  $(\M, \Sigma_\M)$ and $(\N, \Lambda_\N)$ and then somehow get back to $\P$ and $\Q$. This is what the un-dropping game achieves.

\begin{figure}
\centering
\begin{tikzpicture}[node distance=2cm, auto]
  \node (A) {$\P=\M_0$};
  \node (B) [right of=A] {$\dots$};
  \node (C) [right of=B] {$\R_i$};
  \node (X) [node distance=0.5cm, below of = C] {$\triangledown$};
  \node (D) [right of=C] {};
  \node (E) [node distance=1cm, below of=C] {$\Q_i$};
  \node (G) [right of=E] {$\R_{i+1}$};
  \node (Y) [node distance=0.5cm, below of=G] {$\triangledown$};
  \node (H) [node distance = 1cm, below of=G] {$\Q_{i+1}$};
  \node (I) [right of = H] {$\dots$};
  \node (J) [node distance = 1cm, right of = I] {$\R_k$};
  \node (Z) [node distance =0.5cm, below of=J] {$\triangledown$};
  \node (K) [node distance=1cm, below of=J] {$\Q_k$};
  \node (L) [right of =K] {};
  \draw[->] (E) to node {$\T_i$}(G);
  \draw[->] (H) to node {$\T_{i+1}$} (I);
  \draw[->] (K) to node  {$\T_k$} (L);
   \end{tikzpicture}
\caption{A stack with neat drops.}
\label{fig:neat_drops}
\end{figure}

\begin{definition}[The main drops of a stack, Figure \ref{fig:neat_drops}]\label{main drops}\index{main drops} Suppose 
$\P$ is a  germane, $b$-type $\sf{lses}$ that projects precisely and 
\begin{center}
$\T=((\M_\a)_{\a<\eta}, (E_\a)_{\a<\eta-1}, D, R, ( \rnc_\a, \layer_\a)_{\a\in R}, (\b_\a, m_\a), T)$
\end{center}
is a (proper) stack on $\P$ based on ${\sf{hl}}(\P)$. 

We say that $\a\in R$ is a \textbf{main drop} if 
\begin{enumerate}
\item $\a$ is of limit type and of bottom type, 
\item $\lh(\T_\a)$ is a successor ordinal,
\item $\pi^{\T_\a}$ is undefined (see \rdef{proper stack}),
\item $\pi^{\T_\a, b}$ is defined.
\end{enumerate} 
We say $T$ has a \textbf{main drop} if there is $\a\in R$ which is a \textbf{main drop}.

Suppose $\T$ has main drops and let $(\a_i: i \in [1,  k])\subseteq R$ be the sequence of \textbf{main drops} of $\T$ enumerated in increasing order. We then set 
 \begin{enumerate}
 \item $\a_{0}=0$ and $\a_{k+1}=\lh(\T)-1$,
 \item for $i\leq k+1$, $\R_i=\M_{\a_i}$ and for $i\leq k$, $\Q_i=\layer_{\a_i}$\footnote{See \rnot{sflayer notation}.},
 \item for $i\leq k$, $\T_{i}=\T_{[\a_i, \a_{i+1}]}$,
 \item $\T_{k+1}$ and $\Q_{k+1}$ are undefined,
 \item $md^\T=(\a_i, \R_i, \T_i, \Q_i: i\leq k+1)$.
 \end{enumerate}
 We then say that $md^\T=(\a_i, \R_i, \T_i, \Q_i: i\leq k+1)$ is the \textbf{md}-sequence of $\T$. $\myqedhere$
\end{definition}

 Next we define the un-dropping extender of a stack. This is essentially the extender given by dovetailing the embeddings $\pi^{\T_{i}, b}$. The un-dropping extender allows us to get back to the original model, and hence it ``undrops" the main drops of $\T$. First notice that the following is true.
 
 \begin{definition}\label{one point extension} Suppose \begin{itemize}
\item $\P$ is a germane, $b$-type $\sf{lses}$ that projects precisely and
\item $\T$ is a stack on $\P$ such that $\T$ has a last model and it is based on ${\sf{hl}}(\P)$. 
\end{itemize}
 Let $\nu+1=\lh(\T)$. We say $\T$ has a \textbf{one point extension} if  letting
 \begin{center}
 $\T=((\M_\a)_{\a\leq \nu}, (E_\a)_{\a<\nu}, D, R, T)$,
 \end{center}
 and
  \begin{center}
 $\T^{ope}=((\M_\a)_{\a\leq \nu}, (E_\a)_{\a<\nu}, D, R\cup\{\nu\}, T)$,
 \end{center}
 is a proper stack (or according to our convention \rrem{proper stacks convention} just a stack)\footnote{See \rdef{proper stack}.}.$\myqedhere$
 \end{definition}
 The following can now be demonstrated by examining \rdef{proper stack} and the definition of $\pi^{\T, b}$. 
 \begin{lemma}\label{existence of main drops} Suppose
 \begin{itemize}
\item $\P$ is a germane, $b$-type $\sf{lses}$ that projects precisely and
\item $\T$ is a stack on $\P$ that has a one point extension and $\pi^{\T, b}$ is undefined.
\end{itemize} Then $\T$ has a main drop and letting $md^\T=(\a_i, \R_i, \T_i, \Q_i: i\leq k+1)$ be the $md$-sequence of $\T$, $\pi^{\T_k\rest \Q_k, b}$ exists.
 \end{lemma}
 
 We make the following convention.
 
 \begin{terminology}\label{extender convention} Suppose $j: M\rightarrow N$ is a map between two transitive sets or classes $M$ and $N$, and suppose $(\k, \l)$ is such that $j(\k)\geq \l$ and $j\rest \k=1\not id$. We then say that $E$ is the \textbf{$(\k, \l)$-extender} derived from $j$ if 
 \begin{center}
 $E=\{ (a, A): A\in \powerset([\k]^{\card{a}})\cap M \wedge a\in [\l]^{<\omega} \wedge a\in j(A)\}$.
 \end{center}
 We say $E$ is a \textbf{short} extender if $\cp(j)=\k$ and otherwise we say $E$ is \textbf{long}. All extenders used to build extender sequences that we consider in this book are short extenders. In particular, when discussing fully backgrounded constructions (e.g. \rdef{gamma-hod pair construction*}) we tacitly assume that all extenders are short. However, we may from time to time derive an extender from a given embedding and not specify whether it is short or long. For example, see the definition of $E^\T_\Q$ below. $\myqedhere$
 \end{terminology}

\begin{definition}[The un-dropping extender of a proper stack]\label{the un-dropping extender of a continuable stack}\index{un-dropping extender} Suppose 
\begin{itemize}
\item $\P$ is a germane, $b$-type $\sf{lses}$ that projects precisely and
\item $\T$ is a stack on $\P$ such that $\T$ is based on ${\sf{hl}}(\P)$ and $\T$ has a one point extension.
\end{itemize} 
\textbf{When $\pi^{\T, b}$ is undefined.}\\\\
Let $md^\T=(\a_i, \R_i, \T_i, \Q_i : i\leq k+1)$ be the $md$-sequence of $\T$. For $i\leq k+1$, set $\k_i=\d^{\R_i^b}$ and for $i\leq k$, let
\begin{center}
$\sigma^{\T}_{i}:(\powerset(\k_i))^{\R_i}\rightarrow (\powerset(\k_{i+1}))^{\R_{i+1}}$
\end{center}
be given by 
\begin{center}
$\sigma^{\T}_{i}(A)=\pi^{\T_{i}\rest \Q_i, b}(A)\cap \k_{i+1}$. 
\end{center}
Set $\sigma^{\T}=\sigma_k^\T\circ \sigma_{k-1}^{\T}\cdot \cdot\cdot \circ \sigma^{\T}_{0}$. 

Suppose $\Q\insegeq_{hod}\R_{k+1}^b$ is meek. We then let $E^{\T}_\Q$ be the $(\k_0, \d^\Q)$-extender derived from $\sigma^{\T}$. More precisely,
\begin{center}
$E^{\T}_\Q=\{ (a, A): a$ is a finite subset of $\d^\Q$, $A\in (\powerset([\k_0]^{\card{a}}))^{\P}$, and $a\in \sigma^{\T}(A)\}$\index{$E^{\T}_\Q$}. 
\end{center}
\textbf{When $\pi^{\T, b}$ is defined.}\\\\
Suppose $\Q\insegeq_{hod}\pi^{\T, b}(\P^b)$ is a complete layer\footnote{See \rnot{l p}.} of $\pi^{\T, b}(\P^b)$. We then let $E^{\T}_\Q$ be the $(\d^{\P^b}, \d^\Q)$-extender derived from $\sigma^{\T}$. More precisely,
\begin{center}
$E^{\T}_\Q=\{ (a, A): a$ is a finite subset of $\d^\Q$, $A\in (\powerset([\d^{\P^b}]^{\card{a}}))^{\P}$, and $a\in \pi^{\T, b}(A)\}$. 
\end{center}

We then say that $E^\T_\Q$ is the $\Q$-un-dropping extender of $\T$. We also say that $E$ is the \textbf{main un-dropping} extender of $\T$ if $E=E^\T_{\R_{k+1}^b}$ or if $E=E^\T_{\pi^{\T, b}(\P^b)}$. $\myqedhere$
\end{definition}

%

When comparing hod premice we need to consider iterations in which at certain stages $I$ is allowed to use the un-dropping extender of the resulting stack. 
The game producing such iterations is defined below. 

\begin{definition}[The un-dropping iteration game]\label{the un-dropping iteration game}\index{un-dropping iteration game} Suppose $\P$ is a germane, $b$-type $\sf{lses}$ that projects precisely. The un-dropping iteration game on $\P$, $\mathcal{G}^{u}(\P, \k, \l, \a)$, is an iteration game satisfying the following conditions:
\begin{enumerate}
\item In $\mathcal{G}^{u}(\P, \k, \l, \a)$, player I and II collaborate to produce a sequence 
\begin{center} $p=(\M_\b, \T_\b, \Q_\b, E_\b: \b< \gg)$\end{center} such that  \begin{enumerate}
\item $\gg\leq \k$,
\item $\M_0=\P$,
\item for all $\b<\gg$, $\T_\b$ is a stack on $\M_\b$ (and is produced via  the rules of $\mathcal{G}(\M_\b, \l, \a)$\footnote{This is the game defined in \cite[Chapter 4]{OIMT}.}),
\item for each $\b$ such that $\b+1<\gg$, the iteration embedding $\pi^p_{0, \b}:\M_0\rightarrow \M_\b$ is defined,
\item for each $\b$ such that $\b+1<\gg$, either
\begin{enumerate}
\item $E_\b$ is the $\Q_\b$- un-dropping extender of $\T_\b$ and $\M_{\b+1}=Ult(\M_\b, E_\b)$, \\
or
\item $E_\b=\Q_\b=\emptyset$, $\pi^{\T_\b}$ exists and $\M_{\b+1}$ is the last model of $\T_\b$,
\end{enumerate}
\item for a limit ordinal $\b<\gg$, $\M_\b$ is the direct limit of $(\M_\xi, \pi^p_{\xi, \zeta}: \xi<\zeta<\b)$ where $\pi^p_{\xi, \zeta}:\M_\xi\rightarrow \M_\zeta$ is the iteration embedding,
\item player I is the player that chooses extenders while playing $\mathcal{G}(\M_\b, \l, \a)$ to produce $\T_\b$,
\item  player I is the player that chooses to stop the run of $\mathcal{G}(\M_\b, \l, \a)$ by either playing the $\Q_\b$-un-dropping extender $E_\b$ or by letting $\M_{\b+1}$ be the last model of $\T_\b$ (in which case $\pi^{\T_\b}$ must be defined),
\item player II chooses branches while playing $\mathcal{G}(\M_\b, \l, \a)$.
\end{enumerate}
\item Player II loses a run $p$ of $\mathcal{G}^{u}(\P, \k, \l, \a)$ if one of the models appearing in $p$ is ill-founded.
\end{enumerate}
We say $\Sigma$ is a \textbf{$(\k, \l, \a)$-strategy} for $\P$ if it is a strategy for II in $\mathcal{G}^{u}(\P, \k, \l, \a)$ such that any run of $\mathcal{G}^{u}(\P, \k, \l, \a)$ in which player II plays according to $\Sigma$ is not a loss for II.

We say $p=(\M_\b, \T_\b, E_\b: \b< \gg)$ is a \textbf{generalized stack} on $\P$ if it is produced by a run of $\mathcal{G}^{u}(\P, \k, \l, \a)$ and $p$ is not a loss for $II$. Since $\lh(E_\b)=\d^{\Q_\b}$ there is no ambiguity in omitting $\Q_\b$s. $\myqedhere$
\end{definition}

\begin{remark}\label{some comments on gen strategies} Suppose $\P$ is a germane, $b$-type $\sf{lses}$ that projects precisely and $\Sigma$ is a $(\k, \l, \a)$-strategy-strategy. Suppose $\Q'$ is a $\Sigma$-iterate of $\P$ via $p=(\M_\b, \T_\b, E_\b: \b< \gg)$\footnote{Thus, $\Q$ is the last model of $p$.} such that $\pi^{p, b}$ is defined (see \rdef{almost non-dropping gen stacks}) and either $\Q=\Q'$ or $\Q\in Y^{\Q'}$ is such that $\Q^b=(\Q')^b$. Then $\Sigma_{\Q, p}$ is the $(\k', \l, \a)$-strategy of $\Q$ given by $\Sigma_{\Q, p}(q)=\Sigma(p^\frown q)$. Here, $\gg+\k'=\k$. 

Suppose next that $\Q'$ is a $\Sigma$-iterate of $\P$ via $p=(\M_\b, \T_\b, E_\b: \b< \gg)$ and $\Q\insegeq \Q'$ is such that at least one of the following holds:
\begin{enumerate}
\item $\pi^{p, b}$ doesn't exist\footnote{This means that if $\b+1=\gg$ then $\pi^{\T_\b}$ is undefined.}.
\item $\pi^{p, b}$ exists and $\Q\insegeq \pi^{p, b}(\P^b)$. 
\end{enumerate}
Then $\Sigma_{\Q, p}$ is defined like in the previous case but only for stacks produced by $\mathcal{G}(\Q, \l, \a)$. $\myqedhere$
\end{remark}

Just like with ordinary strategies, it also possible to pullback $(\k, \l, \a)$-strategies. The proof of the fallowing theorem is just like the proof of the same theorem for ordinary strategies. 
\begin{theorem}\label{pullback gen straetgies} Suppose $\P$ and $\Q$ are germane, $b$-type $\sf{lses}$ which project precisely, $\sigma:\Q\rightarrow \P$ is a  weak embedding\footnote{In the sense of \cite[Fact 2.13]{ANS}. See the paragraph after \cite[Fact 2.13]{ANS}.} and $\Sigma$ is a $(\k, \l, \a)$-strategy. Then $\Q$ has a $(\k, \l, \a)$-iteration strategy, $\Lambda$, with the following property. For all generalized stack $q=(\Q_\b, \U_\b, F_\b: \b< \gg)$ on $\Q$, $q$ is according to $\Lambda$ if and only if there is a generalized stack $p=(\P_\b, \T_\b, E_\b: \b< \gg)$ on $\P$ and sequences $(\sigma_\b: \b<\gg)$ and $(\tau_{\b, \iota}: \b<\gg \wedge \iota<\lh(\U_\b))$ such that the following clauses hold:
\begin{enumerate}
\item $\sigma_0=\sigma$ and for all $\b<\gg$, $\sigma_\b:\Q_\b\rightarrow\P_\b$ is a weak embedding.
\item For all $\b<\gg$, $\T_\b=\sigma_\b\U_\b$, i.e., $\T_\b$ is obtained from $\U_\b$ via the $\sigma_\b$-copying construction (see \cite[Chapter 4.1]{OIMT}). 
\item For all $\b<\gg$, $(\tau_{\b, \iota}: \iota<\lh(\U_\b))$ is the sequence of copy maps produced during the construction of $\T_\b$. 
\item For each $\b<\gg$, $F_\b$ is the undropping extender of $\U_\b$ if and only if $E_\b$ is the undropping extender of $\T_\b$.
\item For each $\b<\gg$, $F_\b=\emptyset$ and $\Q_\b$ is the last model of $\U_\b$ if and only if $E_\b=\emptyset$ and $\Q_\b$ is the last model of $\U_\b$.
\item For each $\b$ such that $\b+1<\gg$ and $F_\b$ is the undropping extender of $\U_\b$, letting $\nu=\lh(\U_\b)$, for all $a\in \lh(F_\b)^{<\omega}$ and $A\in \M_\b$, $(a, A)\in F_\b$ if and only if $(\tau_{\b, \nu}(a), \sigma_\b(A))\in E_\b$.
\item For each $\b$ such that $\b+1<\gg$ and $F_\b$ is the undropping extender of $\U_\b$, letting $\nu=\lh(\U_\b)$, $\sigma_{\b+1}:Ult(\Q_\b, F_\b)\rightarrow Ult(\P_\b, E_\b)$ is such that $\sigma_{\b+1}([a, f]_{F_\b})=[\sigma_{\b, \nu}(a), \sigma_\b(f)]_{E_\b}$. 
\item For each $\b$ such that $\b+1<\gg$ and $F_\b$ is the undropping extender of $\U_\b$, letting $\nu=\lh(\U_\b)$, $\sigma_{\b+1}=\sigma_{\b, \nu}$.
\end{enumerate}
\end{theorem}
\begin{notation}\label{notation for generalized stacks} Suppose $\T=(\M_\a, \T_\a, E_\a: \a<\gg)$ is a generalized stack. 
\begin{enumerate}
\item For $\a<\gg$ and $\a'<\lh(\T_\a)$, we let $\M_{\a, \a'}^\T=\M_{\a'}^{\T_\a}$. 
\item For $\a<\gg$ and $\iota_0\leq \iota_1< \lh(\T_\a)$ such that $\iota_0 \in [0, \iota_1)_\T$, $\pi^{\T, \a}_{\iota_0, \iota_1}:\M_{\a, \iota_0}^{\T}\rightarrow \M_{\a, \iota_1}^{\T}$ is the iteration embedding $\pi^{\T_\a}_{\iota_0, \iota_1}$ provided it is defined. 
\item Suppose next that $\a_0\leq \a_1<\gg$ and $\iota<\lh(\T_{\a_1})$. We then let $\pi^{\T}_{\a_0, \a_1}:\M_{\a_0}\rightarrow \M_{\a_1}$ be the iteration embedding and $\pi_{\a_0, (\a_1, \iota)}^\T=\pi^{\T_{\a_1}}_{\a_1, \iota}\circ \pi^{\T}_{\a_0, \a_1}$ given that $\pi^{\T_{\a_1}}_{\a_1, \iota}$ is defined.
\item We let $\T^{\sf{ue}}$ be the \textbf{un-dropping extension} of $\T$. More precisely, $\T^{\sf{ue}}$ is defined assuming $\eta=\b+1$ and $\T_\b$ has a one point extension\footnote{See \rdef{one point extension}.}, in which case $\T^{\sf{ue}}$ is obtained by letting $E_\b$ be the un-dropping extender of $\T_\b$, $\M_{\b+1}=Ult(\M_\b, E_\b)$ and $\T_{\b+1}=\emptyset$.
\item  We can also define $\T^{\sf{ue}}_\Q$ assuming $\Q\insegeq \R^b$ where $\R$ is the last model of $\T_\b$. Here, we let $E_\b$ be the $\Q$-un-dropping extender of $\T_\b$.
\item Again assuming $\lh(\T)=\b+1$ and $\T$ has a one point extension, letting $\R$ the last model of $\T_\b$ and $\Q\insegeq \R^b$ be a complete layer of $\R$, we can define the $\Q$-un-dropping extender of $\T$ by setting:
\begin{center}
$E^\T_\Q=\{ (a, A): a\in [\d^\Q]^{<\omega} \wedge A\in \powerset(\d^{\P^b})\cap \P \wedge a\in \sigma^{\T_\b}(\pi^\T_{0, \b}(A))\}$,
\end{center}
where $\sigma^\X$ is defined in \rdef{the un-dropping extender of a continuable stack}. We then set
\begin{center}
$\sigma^\T=\sigma^{\T_b}\circ \pi^{\T}_{0, \b}$.
\end{center} 
Alternatively, $E^\T_\Q$ is the $(\d^{\P^b}, \d^\Q)$-extender derived from $\pi^{\T^{\sf{ue}}}$. We say $E^\T$ is the un-dropping extender of $\T$ if $\T$ is the $\R^b$-un-dropping extender of $\T$.
\item As ordinary stacks are instances of generalized stacks, $\T^{\sf{ue}}$ and $\T^{\sf{ue}}_\Q$ can also be used for ordinary stacks.
\end{enumerate}
 Often, when $\T$ is clear from the context, we will omit it from our notation. $\myqedhere$
\end{notation}


%
The next definition introduces \textit{self-cohering} iteration strategies. The idea is as follows. Suppose $\P$ is a non-meek hod-like ${\sf{lses}}$ and suppose $\T=(\M_\a, \T_\a, F_\a: \a<\eta)$ is a generalized stack on $\P$ according to some iteration strategy $\Sigma$. Let $\R$ be the last model of $\T_0$. Then $\R^b\insegeq \M_1^b$. But it is not clear that $\Sigma_{\R^b, \T_0}=\Sigma_{\R^b, \T_0^\frown \{F_0\}}$. Self-cohering strategies have this property. We will use this property in our diamond comparison argument (see \rdef{sec:diamond comparison}).

\begin{definition}\label{initial segment of a stack given by a node} Suppose $\T$ is a stack and $\R=\M_\a^\T$ for some $\a<\lh(\T)$. We then say that $\R$ is a node of $\T$ and write $\T_{\leq \R}$ for $\T_{\leq \a}$. Similarly if $\R'=\M_\b^\T$ for $\b>\a$ then we can define $\T_{<\R}$, $\T_{\R, \R'}$ and $\T_{\geq \R}$. Similar notation can be introduced for generalized stacks in the obvious way. $\myqedhere$
\end{definition}

\begin{definition}\label{self-cohering}\index{self-cohering strategies} Suppose $(\P, \Sigma)$ is a a hod-like ${\sf{lses}}$ pair (see \rdef{hod-like lsp pair}). We say that $\Sigma$ is \textbf{self-cohering} if whenever 
\begin{itemize}
\item $\T=(\M_\a, \T_\a, F_\a: \a<\eta)$ is a generalized stack according to $\Sigma$, 
\item $\a_0, \a_1<\eta$,
\item $\xi_0< \lh(\T_{\a_0})$ and $\xi_1<\lh(\T_{\a_1})$,
\item $\R\inseg_{hod}\M^{\T_{\a_0}}_{\xi_0}=_{def}\S_0$ and $\R\inseg_{hod}\M^{\T_{\a_1}}_{\xi_1}=_{def}\S_1$,
\end{itemize}
\begin{center}
$\Sigma_{\R, \T_{\leq \S_0}}=\Sigma_{\R, \T_{\leq \S_1}}$\footnote{See \rdef{id pullback initial segment}.}.
\end{center}
where the equality is between the $(\omega_1, \omega_1)$ portions of both strategies. $\myqedhere$
\end{definition}

Self-cohering is a desired property and we will have to establish that our constructions produce strategies that are self-cohering. However, it is more convenient not to make it part of our definitions.

\begin{definition}[Hod-like $\sf{lses}$ pair]\label{hod-like lsp pair}\index{hod-like ${\sf{lses}}$ pair}
We say $(\P, \Sigma)$ is a \textbf{hod-like $\sf{lses}$ pair} (with an indexing scheme $\phi$) if 
\begin{enumerate}
\item $\P$ is a hod-like $\sf{lses}$ (with an indexing scheme $\phi$),
\item if $\P$ is non-meek then $\Sigma$ is a $(\k, \l, \nu)$-strategy,
\item if $\P$ is meek or gentle then $\Sigma$ is a $(\l, \nu)$-strategy,
\item if $\Q$ is a $\Sigma$-iterate of $\P$ via $\T$ and $\R\insegeq_{hod}\Q$ then $\Sigma^\R\subseteq \Sigma_{\R, \T}$\footnote{This clause is asserting that the internal strategy of $\R$ agrees with $\Sigma_{\R, \T}$.}.
\end{enumerate}
 We say $(\P, \Sigma)$ is a \textbf{simple hod-like $\sf{lses}$ pair} if $\P$ is a hod-like $\sf{lses}$ , $\Sigma$ is a $(\l, \nu)$-iteration strategy and clause 4 above holds. 
 
 In the context of $\sf{AD^+}$, unless otherwise specified, the strategy of a hod-like $\sf{lses}$ pair or a simple  hod-like $\sf{lses}$ pair is an $(\omega_1, \omega_1, \omega_1)$-strategy or an $(\omega_1, \omega_1)$-strategy. 
 
 $\myqedhere$
\end{definition}


Finally we finish this section by stating the version of \rlem{the canonical singularizing sequence exists} for generalized stacks. Its proof is just like the proof of \rlem{the canonical singularizing sequence exists}. First we generalize $\pi^{\T, b}$ and \rdef{canonical singularizing sequences} to generalized stacks.

 \begin{definition}[Almost non-dropping generalized stacks]\label{almost non-dropping gen stacks}\index{almost non-dropping generalized stacks} Suppose $\M$ is germane of $b$-type and projects precisely. Suppose further that 
 \begin{center}
 $\T=(\M_\a, \T_\a, E_\a: \a<\gg)$
 \end{center}
 is a generalized stack on $\M$ that is based on ${\sf{hl}}(\M)$. We say that $\T$ is \textbf{almost non-dropping} if either $\gg$ is a limit ordinal or $\gg=\a+1$ and $\pi^{\T_\a, b}$ exists. Assuming $\T$ is almost non-dropping we set
 \begin{center}
$\pi^{\T, b}=\begin{cases}
\pi^\T_c & :\ \gg\ \text{is a limit ordinal and $c$ is the unique branch of $\T$}\\
\pi^{\T_\a, b}\circ \pi_{0, \a}^\T\rest \M^b &:\ \text{otherwise}
\end{cases}$\end{center}
Suppose $\Sigma$ is a $(\k, \l, \a)$-iteration strategy for $\M$\footnote{It is worth remembering that this entails that $\Sigma$-iterates of $\M$ have the same indexing scheme as $\M$.}. We then let 
\begin{center}
$I(\M, \Sigma)=\{(\T, \R):  \T$ is according to $\Sigma$, $\T$ is based on ${\sf{hl}}(\M)$, $\R$ is the last model of $\T$ and $\pi^{\T}$ is defined$\}$.\index{$I(\M, \Sigma)$}\\
$I^b(\M, \Sigma)=\{ (\T, \R): \T$ is according to $\Sigma$, $\T$ is based on ${\sf{hl}}(\M)$, $\R$ is the last model of $\T$ and $\pi^{\T, b}$ is defined$\}$\index{$I^b(\M, \Sigma)$}\\
$I^{ope}(\M, \Sigma)=\{(\T, \R): \T$ is according to $\Sigma$, $\T$ is based on ${\sf{hl}}(\M)$, $\T$ has a one point extension\footnote{See \rdef{one point extension}. Here one point extension of a generalized stack $\T=(\M_\a, \T_\a, E_\a: \a<\eta)$ is $\T$ unless $\eta=\b+1$, in which case we let $\T^{ope}=(\M_\a, \T_\a, E_\a: \a<\b)^\frown (\M_\b, \T_\b^{ope})$.} and $\R$ is the last model of $\T\}$\index{$I^{ope}(\M, \Sigma)$}\\
$B^{ope}(\M, \Sigma)=\{(\T, \R):$ there is $(\T, \R')\in I^{ope}$ and $\R$ is a layer of $\R'\}$.\index{$B^{ope}(\M, \Sigma)$}\\
\end{center}
$\myqedhere$
\end{definition}

We remark that we will use $I^{ope}(\M, \Sigma)$ and $B^{ope}(\M, \Sigma)$ even when $\Sigma$ is an iteration strategy acting on stacks.

\begin{definition}[Canonical singularizing sequences]\label{canonical singularizing sequences gen} Suppose $\P$ is a germane $\sf{lses}$ of $b$-type that projects precisely and $\T$ is an almost non-dropping generalized stack on $\P$. Let $\Q=\pi^{\T, b}(\P^{b})$. Then $\Q$ is a hod-like $\sf{lses}$. If $w=(\eta, \d)$ is a window of $\Q$ then we let 
\begin{center}
$s(\T, w)=\{ \a: \exists a\in \eta^{<\omega} \exists f\in \P^b( \a=\pi^{\T, b}(f)(a))\}\cap \d$
\end{center}
$\myqedhere$
\end{definition}

\begin{lemma}\label{the canonical singularizing sequence for g-stacks} Suppose $\P$ is germane, $b$-type $\sf{lses}$ that projects precisely and
\begin{center}
 $\T=(\M_\a, \T_\a, E_\a: \a<\gg)$
 \end{center}
 is a generalized stack $\P$ such that $\pi^{\T, b}$ exists. Let $\Q=\pi^{\T, b}(\P^b)$. Then for any window $w=(\eta, \d)$  of $\Q$ (see \rnot{l p}) such that  $\Q\models ``\d$ is a Woodin cardinal",
\begin{center}
$\sup(s(\T, w))=\d$.
\end{center}
\end{lemma}

\chapter{Short tree strategy mice}\label{chap: shorttreestrategymice}

The main purpose of this chapter is to isolate the definition of short tree strategy mice. As was mentioned many times before, the main problem with defining this concept is the fact that it is possible that maximal iteration trees (which should not  have branches indexed in the strategy predicate) may \textit{core down} to short iteration trees (which must have branches indexed in the strategy predicate), thus causing indexing issues. To solve this issue we will design an authentication procedure which will carefully choose iteration trees and index their branches. Thus, if some iteration tree doesn't have a branch indexed in the strategy predicate then it is because the authentication procedure hasn't yet found an authenticated branch, and therefore, such iteration trees cannot core down to an iteration tree whose branch is authenticated.

The following is a rough roadmap of the chapter. \rsec{short tree strategy section} introduces the short tree component of an iteration strategy, while \rsec{short tree game and sts mice} introduces the short tree strategy as an abstract object. This is an important step as the strategy predicate of a short tree strategy mouse codes a short tree strategy in the sense of \rdef{short tree strategy}. The next important step is the isolation of two different kinds of iterations, those that are universally short (see \rdef{nus stacks}), i.e. short with respect to any strategy, and those that are ambiguous. As there is no ambiguity involved in determining whether a universally short iteration trees are short or not and moreover, since universal shortness is preserved under Mostowski collapses, we will simply add the branches of such iterations to the strategy predicate without authenticating them first. The branches of ambiguous iterations will  be authenticated before being added to the strategy predicate. 

A key tool in the authentication procedure is the fully backgrounded constructions that produce iterates of hod like ${\sf{lses}}$ (see \rdef{full short tree coherent background constructions}). Such constructions are used to find a branch of an iteration tree with the property that the branch model itself iterates to the same construction.  \rdef{sts indexing scheme} and \rdef{sts phi premouse} introduce the particular ways branches will be indexed. Our authentication procedure appears as \rdef{authentic lsp}, and \rdef{weak psi alpha indexing scheme a} defines indexing scheme. Then \rdef{sts premouse} introduces the short tree strategy mice and \rrem{how branches get indexed} explains exactly how branches get indexed. \rdef{hod premouse} finally introduces the concept of a hod premouse.

\begin{remark}\label{germane remark} All the notions introduced in this chapter can be routinely carried over to germane $\sf{lses}$ that project precisely. Thus, when discussing germane $\sf{lses}$, we will freely use the language developed in the sections of this chapter. $\myqedhere$
\end{remark}

\section{The short tree component of a strategy}\label{short tree strategy section}

Suppose $(\P, \Sigma)$ is a hod-like ${\sf{lses}}$ pair or a simple hod-like ${\sf{lses}}$ pair such that $\P$ is of lsa type (see \rdef{hod-like lsp pair}). 
Since the particular indexing scheme will not matter for what follows, we suppress the indexing scheme that the pair $(\P, \Sigma)$ has. The  next definition isolates the \textit{short tree component} of $\Sigma$ denoted by $\Sigma^{stc}$. Let $\k=\d^{\P^b}$ and $\d=\d^\P$. Recall that all our stacks are proper stacks (see \rrem{proper stacks convention}). The next few concepts will be introduced for generalized stacks, and as stacks are instances of generalized stacks, they can be used in connection with stacks. 

\begin{remark}\label{we need projects precisely} The short tree component of $\Sigma$ is a strategy that acts on $\P_\#$\footnote{See \rdef{lsa type}.}. Thus, the short tree component does not in general produce stacks that can be applied to $\P$ without dropping in degree. Such dropping can happen, for example, when $\rho_{k(\P)}(\P)<\d^\P$.$\myqedhere$
\end{remark}

\begin{definition}\label{sigma short} Suppose $(\P, \Sigma)$ is a hod-like ${\sf{lses}}$ pair or a simple hod-like ${\sf{lses}}$ pair such that $\P$ is of lsa type. Suppose $\T=(\M_\b, \T_\b, E_\b: \b<\gg)$ is a generalized stack on $\P_{\sf{ex}}$ according to $\Sigma_{\sf{ex}}$\footnote{See \rdef{lsa type}.}. We say $\T$ is \textbf{$\Sigma$-short} if $\T\in \dom(\Sigma_{\sf{ex}})$ and letting $b=\Sigma_{{\sf{ex}}}(\T)$ one of the following conditions holds: 
\begin{enumerate}
\item $\pi^{\T}_b$ is undefined.
\item $\gg$ is a limit ordinal.
\item $\gg$ is a successor ordinal and $\lh(\T_{\gg-1})$ is a limit ordinal.
\item $\gg$ is a successor ordinal,  $\lh(\T_{\gg-1})$ is a successor ordinal and letting 
\begin{center}
$\a=\max(R^{\T_{\gg-1}})$,
\end{center}
 $\pi^\T_b(\d)>\d(\T_{\geq \a})$.
\end{enumerate}

We then say that $\T$ is \textbf{$\Sigma$-maximal} if it is not $\Sigma$-short.  $\myqedhere$
\end{definition}

\begin{remark}\label{maximal alpha} Recall that according to our convention \rrem{proper stacks convention}, if $\T$ is a stack and $\a$ is a cutpoint of $\T$ then $\a\in R^\T$. Hence, if for some $\a<\lh(\T)$, $\T_{<\a}$ is $\Sigma$-maximal then $\a\in R^\T$. $\myqedhere$
\end{remark}
  
  Notice that if $\T$ is $\Sigma$-short then it does not follow that initial segments of $\T$ are also $\Sigma$-short. If $\T$ is a generalized stack or just a stack then we let $\T^-$\index{$\T^-$} be $\T$ without its last model if it exists and $\T$ otherwise. The next definitions describe exactly when a stack is according to the short tree strategy component of $\Sigma$. \rdef{the short tree component domain 1} introduces the domain of $\Sigma^{stc}$ restricted to ordinary indexable stacks and \rdef{the short tree component domain 2} introduces the domain of $\Sigma^{stc}$.

\begin{definition}\label{mtsharp} Suppose $\T$ is a normal iteration tree of limit length. We then let $\m^+(\T)=(\m(\T))^\#$. $\myqedhere$
\end{definition}

\rdef{the short tree component domain 1} needs a slight modification of the concept of a tree order. 

\begin{definition}\label{modified tree order} Suppose $(\iota_\tau: \tau<\nu)$ is an increasing sequence of ordinals and for all $\tau<\nu$ such that $\tau+1<\nu$, $I_\tau$ is either the interval $[\iota_\tau, \iota_{\tau+1})$ or the interval $[\iota_\tau, \iota_{\tau+1}]$. We say $I_\tau$ is right-open if $I_\tau=[\iota_\tau, \iota_{\tau+1})$ and otherwise we say $I_\tau$ is right-closed. Let $\iota=\sup\{\iota_\tau+1:\tau<\nu\}$.  We then say that $U$ is a \textbf{tree order} on $\prod_{\tau<\nu} I_\tau$ if $U\subseteq \iota^2$ such that the following clauses hold.
\begin{enumerate}
\item $U$ is a partial order preserving the usual order on ordinals.
\item If $(\a, \b)\in U$ then for some $\tau<\nu$, $(\a, \b)\in I_\tau\times I_\tau$. 
\item For all limit ordinals $\l<\iota$, either 
\begin{enumerate}
\item for some $\tau<\nu$, $\l=\iota_{\tau+1}$ and $I_\tau$ is right-open, or
\item $\{ \a<\l: (\a, \l)\in U\}$ is a closed unbounded subset of $\l$.
\end{enumerate}
\end{enumerate}
We will freely adopt the usual notation used for ordinary tree orders. For example, $<_U$ is the order given by $U$, $U(\a+1)$ is the $U$-predecessor of $\a+1$ and $[\a, \b]_U=\{\gg: \a\leq_U\gg\leq_U \b\}$.

Suppose $T$ is a tree order on $\iota$. We then let $T\rest \prod_{\tau<\nu}I_{\tau}$  be the unique tree order $U$ on $\prod_{\tau<\nu} I_\tau$ with the property that for all $\tau<\nu$ and for all $(\a, \b)\in I_\tau\times I_\tau$, $(\a, \b)\in U \iff (\a, \b)\in T$. $\myqedhere$
\end{definition}

\begin{definition}[The domain of the short tree component of a strategy I]\label{the short tree component domain 1} Suppose $(\P, \Sigma)$ is a hod-like ${\sf{lses}}$ pair or a simple hod-like ${\sf{lses}}$ pair such that $\P$ is of lsa type. We let
\begin{center}
$\U
=((\M_\a)_{\a<\eta}, (E_\a)_{\a<\eta-1}, D, R, (\beta_\a, m_\a)_{\a\in R}, \so, \ma, U )\in \dom(\Sigma^{stc})$
\end{center}
 if there is a stack 
 \begin{center}
 $\T=((\M'_\a)_{\a<\eta}, (E'_\a)_{\a<\eta-1}, D', R', (\beta'_\a, m_\a')_{\a\in R'}, T)\in \dom(\Sigma_{{\sf{ex}}})$\footnote{See \rdef{lsa type}.}
 \end{center} such that $\U$ is the same as $\T$ except it doesn't have the maximal branches of $\T$; more precisely, the following conditions hold.
\begin{enumerate}
\item $\M_0=\P_{\#}$ and $\T$ is below $\d^\P$\footnote{I.e., for every $\a<\eta$ if $\pi_{0, \a}^\T$ is defined then $\ind_\a^\T<\pi_{0, \a}^\T(\d^\P)$.}.
\item $D=D'$, $R'=R=\so\cup \ma$, $\so\cap \ma=\emptyset$ and $\ma$ is finite.
\item For all $\a<\eta$, $E_\a=E_{\a}'$, $\b_\a=\b_\a'$ and $m_\a=m_\a'$.
\item For all successor $\a<\eta$, $\M_\a=\M'_\a$.
\item For all limit $\a<\eta$ such that $\T_{<\a}$ is $\Sigma$-short, $\M_\a=\M_\a'$.
\item For all limit $\a<\eta$ such that $\T_{<\a}$ is $\Sigma$-maximal, letting $\mathcal{X}$ be the last normal component of $\T_{<\a}$, $\M_\a=\m^+(\mathcal{X})$\footnote{If $\T_{<\a}$ is $\Sigma$-maximal then because of clause 1 its last normal component cannot be normally continued implying that $\a\in R$.}.\\\\
 Let $\nu=o.t.(R)$ and $(\iota_\tau: \tau < \nu)$ be the increasing enumeration of $R$. If $\tau+1=\nu$ then set 
 \begin{center}$
 \iota_{\tau+1}=\begin{cases} \eta &:\eta\ \text{is a limit ordinal}\\
 \eta-1 &: \text{otherwise}
 \end{cases}$\end{center}
 We say $\tau+1$ is irrelevant if $\tau+1=\nu$ and $\iota_{\tau+1}=\eta$, and if $\tau+1$ is not irrelevant then we say $\tau$ is relevant. For $\tau<\nu$ such that $\iota_\tau\in \so$ let $I_\tau=[\iota_\tau, \iota_{\tau+1}]$ and otherwise set $I_{\tau}=[\iota_\tau, \iota_{\tau+1})$.
 \item $U=T\rest \prod_{\tau<\nu} I_\tau$. 
 \item For every $\tau\leq \nu$ if $\iota_\tau$ is defined\footnote{if $\nu$ is limit then $\iota_\tau$ is not defined.} then $\iota_\tau\in \ma$ if and only if $\tau$ is relevant and $\T_{<\iota_{\tau+1}}$ is $\Sigma$-maximal.

\end{enumerate}
If $\T$ and $\U$ are as above then we write $\U=\T^{sc}$\index{$\T^{sc}$} and say that $\U$ is the short component of $\T$. $\myqedhere$
\end{definition}
%
%
%
%
%
%
%

\begin{definition}[The short tree component of a strategy II]\label{the short tree component domain 2} Suppose $(\P, \Sigma)$ is a hod-like ${\sf{lses}}$ pair or a simple hod-like ${\sf{lses}}$ pair such that $\P$ is of lsa type.
We set
\begin{center}
$\U=(\N_\a, \U_\a, E_\a: \a<\eta)\in \dom(\Sigma^{stc})$
\end{center}
 if there is a generalized stack $\T=(\M_\a, \T_\a, F_\a : \a<\eta)\in \dom(\Sigma_{{\sf{ex}}})$ (see \rdef{the un-dropping iteration game}) such that $\U$ is the same as $\T$ except it doesn't have the maximal branches of $\T$; more precisely, 
\begin{enumerate}
\item $\M_0=\P||\a_0$,
\item for every $\a<\eta$, $\N_\a=\M_\a$ and $E_\a=F_\a$,
\item for every $\a<\eta$, $\U_\a=\T_\a^{sc}$,
\item there are at most finitely many $\a$ such that $\U_\a\not =\T_\a$, and
\item either $\eta$ is a limit ordinal or the last normal component of $\T_{\eta-1}$ has a limit length (this condition is redundant as $\T\in \dom(\Sigma)$).
\end{enumerate}
If $\T$ and $\U$ are as above then we write $\U=\T^{sc}$ and say that $\U$ is the short component of $\T$. $\myqedhere$
\end{definition}

Conditions (3-4) in \rdef{the short tree component domain 2}  ensure that if the relevant stacks are of limit length, we can take the direct limit. We will not be concerned with quasi-limits (cf. \cite{CMI}) here. The next definition defines the short tree component of $\Sigma$. 


\begin{definition}[The short tree component of a strategy]\label{the short tree component stacks}\index{the short tree component of a strategy} Suppose $(\P, \Sigma)$ is a hod-like ${\sf{lses}}$ pair or a simple hod-like ${\sf{lses}}$ pair such that $\P$ is of lsa type and is exact. Suppose 
\begin{center}
$\T=(\M_\a, \T_\a, F_\a: \a<\eta)$
\end{center} is a generalized stack on $\P$ such that $\T\in \dom(\Sigma_{\sf{ex}})$ and $\T^{sc}\in \dom(\Sigma^{stc})$. Let $b=\Sigma(\T)$. We then set $\Sigma^{stc}(\U)=x$ where $x$ is defined as follows:
\begin{enumerate}
\item If $\eta$ is a successor ordinal, $\T_{\eta-1}$ has a last normal component\footnote{See \rdef{notation for iteration trees}.} and letting $\X$ be the last normal component of $\T_{\eta-1}$, $\pi^{\T}_b$ is defined and $\pi^\T_b(\d)=\d(\mathcal{X})$ then $x=\m^+(\X)$.
\item Otherwise $x=b$\footnote{For reader's convenience we spell out the exact clauses of ``Otherwise". \begin{enumerate} \item $\eta$ is a limit ordinal.
\item $\eta$ is a successor ordinal and $\T_{\eta-1}$ doesn't have a last normal component or $\pi^{\T}_b$ is undefined.
\item $\eta$ is a successor ordinal, $\T_{\eta-1}$ has a last normal component $\mathcal{X}$, $\pi^{\T}_b$ is defined and $\pi^\T_b(\d)>\d(\mathcal{X})$.\end{enumerate}}. 
\end{enumerate}
$\myqedhere$
\end{definition}

Thus, $\Sigma^{stc}(\T)$ either returns the value of $\Sigma_{\sf{ex}}(\T)$ or $\m^+(\X)$ where $b=\Sigma_{{\sf{ex}}}(\T)$. From now on, we will use this notation even when $\Sigma$ is a partial iteration strategy. 

Notice the similarity with the short tree iterability for suitable mice in the context of core model induction or in the context of $\H$ analysis and $\Sigma^{stc}$. If $\P$ is a $\Sigma^2_1$-suitable premouse and $\Sigma$ is fullness preserving iteration strategy for $\P$,\footnote{Here $\Sigma^2_1$ and fullness preservation are relative to an $\sf{AD}^+$-model.} $\Sigma^{stc}$ is just the short tree iterability strategy of $\P$. 

\section{Short tree stacks and short tree strategies}\label{short tree game and sts mice}

In order to define the \textit{short tree strategy mice}, we will need to introduce the concept of \textit{short tree strategy} that is independent of a particular strategy. We start by defining short-tree-stacks, or just st-stacks. Recall our convention that all stacks are proper (see \rrem{proper stacks convention}). 
We will not take the usual route of first defining putative st-stacks and then defining st-stacks, and leave such matters to the reader. Our goal is to concentrate on the important new property that st-stacks have. 


\begin{definition}\label{st-stack} Suppose $\P$ is a hod-like $\#$-lsa type  ${\sf{lses}}$\footnote{See \rdef{lsa type}.}. Set $\d=\d^\P$. We say that $\T$ is an \textbf{st-stack}\index{st-stack}  on $\P$ if 
\begin{center}
$\T=((\M_\a)_{\a<\eta}, (E_\a)_{\a<\eta-1}, D, R, (\beta_\a, m_\a)_{\a\in R}, \so, \ma, T)$
\end{center}
and the following conditions hold.
\begin{enumerate}
\item $R$ is a closed subset of $\eta$ and $0\in R$.
\item $R=\so\cup \ma$, $\so\cap \ma=\emptyset$ and $\ma$ is finite.\\\\
 Let $\nu=o.t.(R)$ and $(\iota_\tau: \tau < \nu)$ be the increasing enumeration of $R$. For each $\a<\eta$, let $\tau_\a$ be the largest $\tau<\nu$ such that $\iota_{\tau}\leq \a$ and set $\iota^\a=\iota_{\tau_\a}$. If $\tau+1=\nu$ then set 
 \begin{center}$
 \iota_{\tau+1}=\begin{cases} \eta &:\eta\ \text{is a limit ordinal}\\
 \eta-1 &: \text{otherwise}
 \end{cases}$\end{center}
 We say $\tau+1$ is irrelevant if $\tau+1=\nu$ and $\iota_{\tau+1}=\eta$, and if $\tau+1$ is not irrelevant then we say $\tau$ is relevant. For $\tau<\nu$ such that $\iota_\tau\in \so$ let $I_\tau=[\iota_\tau, \iota_{\tau+1}]$ and otherwise set $I_{\tau}=[\iota_\tau, \iota_{\tau+1})$. We say that $(\iota_\tau: \tau\leq \nu)$ is the $\iota$-sequence of $\T$
 \item $T$ is a tree order on $\prod_{\tau<\nu} I_\tau$. 
\item For all $\a<\eta$, $\M_\a$ is a well-founded $\sf{lhes}$ (or ${\sf{hes}}$).
\item For all $\a\in R$, $(\omega\b_\a, m_\a)\leq l(\M_\a)$.\\\\
Set
\begin{center} $\M_\a'=\begin{cases} \M_\a &: \a\not \in R \vee (\a\in R \wedge \omega\b_\a= \ord(\M_\a))\\
\M_\a||(\omega\b_\a, \omega)& : \a\in R \wedge \omega\b_\a<\ord(\M_\a)
\end{cases}$
\end{center}
\item $\M_0=\P$.
\item For all $\a+1<\eta$, $E_\a\in \vec{E}^{\M_\a'}$.
\item Normality conditions hold. More precisely, the following conditions hold.
\begin{enumerate}
\item For all $\a+1<\eta$, letting $\b=T(\a+1)$ and $\k_\a=\cp(E_\a)$, then $\b$ is the least ordinal $\gg\geq \tau_\a$ such that 
\begin{center}
$(\k_\a^+)^{\M_\a|\ind^{\M_\a}(E_\a)}<\nu(E_\gg)$.
\end{center}
\item For all $\a<\b$ such that $\b+1<\eta$ and $\iota^\a=\iota^\b$, $\ind^{\M_\a}(E_\a)<\ind^{\M_\b}(E_\b)$. 
\end{enumerate}
\item For all $\a+1<\eta$,
\begin{center}
$\M_{\a+1}=Ult(\M_\b'||(\omega\xi_\a, k_\a), E_\a)$ 
\end{center}
where 
\begin{enumerate}
\item $\b=T(\a+1)$,
\item $\omega\xi_\a\leq \ord(\M'_\b)$ is the largest such that $(\k_\a^+)^{\M_\a|\ind^{\M_\a}(E_\a)}=(\k_\a^+)^{\M_\b'|\omega\xi_\a}$,
\item $k_\a$ is the largest such that $(\omega\xi_\a, k_\a)\leq l(\M_\b')$ and $\cp(E_\a)<\rho_{k_\a}(\M_\b'||(\omega\xi_\a, k_\a))$.
\end{enumerate}
\item $D=\{\a+1<\eta:$ letting $\b=T(\a+1)$, $(\omega\xi_\a, k_\a)<l(\M_\b)\}$.

Let\begin{center} $\pi^{\T}_{\b, \a+1}=\pi_{E_\a}^{\M'_\b||(\omega\xi_\a, k_\a)}: \M'_\b||(\omega\xi_\a, k_\a)\rightarrow \M_{\a+1}$\end{center} be the ultrapower map and for $\a<\gg<\eta$ such that $\tau_\a=\tau_\gg$ and $\a<_{T}\gg<\eta$ let $\pi^\T_{\a, \gg}:\M_\a\rightarrow \M_\gg$ be the embedding obtained by compositions.\footnote{Assuming these embeddings can be composed. $\pi^\T_{\a, \gg}$ is defined if and only if $D\cap (\a, \gg]_\T=\emptyset$.}

\item Suppose $\l<\eta$ is a limit ordinal. Then the following clauses hold.
\begin{enumerate}
\item Suppose $\l\not \in R$. Then $D\cap (\iota^\l, \l)_{T}$ is finite and letting $\b\in [\iota^\l, \l)_{T}$ be the least such that $D\cap (\b, \l)_{T}=\emptyset$, $\M_\l$ is the direct limit of the system $(\M_\gg, \pi_{\gg, \gg'}^\T: \gg<\gg', \gg, \gg'\in [\b, \l)_{T})$ and for $\gg\in [\b, \l)$, $\pi^\T_{\gg, \l}:\M_\gg\rightarrow \M_\l$ is the direct limit embedding.
\item Suppose $\l\in \so$. Then $\sup(\ma \cap \l)<\l$\footnote{This is a consequence of the fact that $\ma$ is finite.} and setting $\l_0=\sup(\ma \cap \l)$ and $\l_1=\sup(D\cap \l)$, $\M_\l$ is the direct limit of $(\M_\a, \pi_{\a, \b}: \max{\l_0, \l_1}<\a<\b, (\a, \b)\in \so^2\cap \l^2)$.\\
\end{enumerate}
For each $\tau<\nu$ such $\iota_\tau \in \so$, let $\T^\tau$ be the re-organization of $\T_{[\iota_\tau, \iota_{\tau+1}]}$\footnote{If $\iota_{\tau+1}=\eta$, we let $\T_{[\iota_\tau, \iota_{\tau+1}]}=\T_{[\iota_\tau, \iota_{\tau+1})}$. Also, see the discussion after \rdef{putative it}.}  as a normal iteration tree on $\M_{\iota_\tau}$ and for each $\tau<\nu$ such that $\iota_\tau\in \ma$, let $\T^\tau$ be the re-organization of $\T_{[\iota_\tau, \iota_{\tau+1})}$ as a normal iteration tree on $\M_{\iota_\tau}$.\\
\item For each $\tau<\nu$ such that $\iota_\tau\in \ma$, $\M_{\iota_{\tau+1}}=m^+(\T^\tau)$, $\d^{\M_{\iota_{\tau+1}}}=\d(\T^\tau)$ and 
\begin{center}
$\mathcal{J}_{\omega}[\M_{\iota_{\tau+1}}]\models ``\d(\T^\tau)$ is a Woodin cardinal".\\
\end{center}

\item For each $\tau<\nu$ such that $\iota_\tau \in \so$, $\tau+1$ is relevant and $\pi_{\iota_\tau, \iota_{\tau+1}}$ is defined, $\pi_{\iota_\tau, \iota_{\tau+1}}(\d^{\M_{\iota_\tau}})>\d(\T^\tau)$.
\item For all $\tau <\nu$ such that $\iota_\tau\in \ma$, $\pi^{\T^\tau, b}$ is defined\footnote{See \rdef{embedding on the top for stacks}.}.  
\item For every $\a<\eta$, if $\pi_{0, \a}$ is defined then $\ind^{\M_\a}(E_\a)<\pi_{0, \a}(\d)$. 
\item If $\tau < \nu$ is such that $\iota_\tau$ is the least member of $\ma$ then $\pi_{0, \iota_{\tau}}$ is defined. 
\item If $\tau_0<\tau_1$ are such that $\iota_{\tau_0}\in \ma$, $\iota_{\tau_1}\in \ma$ and $[\iota_{\tau_0+1}, \iota_{\tau_1})\cap \ma=\emptyset$ then (provided $\tau_0+1<\tau_1$) $\pi_{\iota_{\tau_0+1}, \iota_{\tau_1}}$ is defined and for every $\a\in [\iota_{\tau_0+1}, \iota_{\tau_1+1})$, if $\pi_{\iota_{\tau_0+1}, \a}$ is defined then $\ind^{\M_\a}(E_\a)<\pi_{\iota_{\tau_0+1}, \a}(\d^{\M_{\iota_{\tau_0+1}}})$.
\item If $\tau< \nu$ is such that $\iota_{\tau}$ is the largest member of $\ma$ then for every $\a\in [\iota_{\tau}, \eta)$ if $\pi_{\iota_{\tau}, \a}$ is defined then $\ind^{\M_\a}(E_\a)<\pi_{\iota_{\tau}, \a}(\d^{\M_{\iota_{\tau}}})$.
\end{enumerate}
We say $\T$ is (an ordinary) normal st-stack if $R^\T=\{0\}$ and $(\b_0^\T, m_0^\T)=l(\M_0^\T)$. 

We adopt our proper stacks convention, \rrem{proper stacks convention}, and in particular demand that all cutpoints of $\T$ are in $R^\T$. $\myqedhere$
\end{definition}

\begin{remark}\label{pitb for st-stack} $\pi^{\T, b}$ can also be defined for st-stacks. See \rdef{embedding on the top for stacks}. $\myqedhere$
\end{remark}

\begin{remark}[Proper st-stack convention]\label{proper st-stack convention} We again make the convention that st-stacks are proper stacks. Adopting the definition of proper stack to st-stacks is a straightforward matter which we leave to the reader. $\myqedhere$
\end{remark}

We will use superscript $\T$ to denote the objects introduced in \rdef{st-stack} (e.g. $\ma^\T$ or $\iota_\tau^\T$). Also, we write $\lh(\T)$ for the ordinal $\eta$. 

It is now straightforward to define the concept of \textit{generalized st-stacks} on $\P$ following the definition of \rdef{the un-dropping iteration game}. These have the form $(\M_\b, \T_\b, E_\b: \b< \gg)$ where $\T_\b$ is an st-stack on $\M_\b$ and $E_\b$ is the un-dropping extender. We leave the details of the definition to the reader. Next we define st-strategy and leave it to the reader to define \textit{generalized st-strategies}.

\begin{definition}[St-strategy]\label{short tree strategy} Suppose $\P$ is a hod-like $\#$-lsa type  ${\sf{lses}}$\footnote{See \rdef{lsa type}.}. We say that $\Lambda$ is an \textbf{st-strategy} for $\P$ if $\Lambda$ is a function with the following properties.
\begin{enumerate}
\item If $x\in \dom(\Lambda)$ then $x$ is an st-stack on $\P$ such that if 
\begin{center}$x=_{def}\T=((\M_\a)_{\a<\eta}, (E_\a)_{\a<\eta-1}, D, R, (\beta_\a)_{\a\in R}, \so, \ma, T )$\end{center} then $\eta$ is a limit ordinal.
\item If $\T\in \dom(\Lambda)$, 
\begin{center}$\T=((\M_\a)_{\a<\eta}, (E_\a)_{\a<\eta-1}, D, R, (\beta_\a)_{\a\in R}, \so, \ma, T )$\end{center} and $\Lambda(\T)=x$ then $\T^\frown \{x\}$ is an st-stack on $\P$. More precisely the following conditions hold.
\begin{enumerate}
\item If $o.t.(R)$ is a limit ordinal then letting $\a\in R$ be such that $\ma\cup D\subseteq \a$, $x$ is the direct limit of $(\M_\b, \pi_{\b, \gg}: \b<\gg, (\b, \gg)\in (R-\a)^2)$. 
\item If $\tau+1=o.t.(R)$ then $x$ is either a branch of $\T_{\iota_{\tau}}$ or $x=\m^+(\T_{\iota_\tau})$\footnote{See \rdef{mtsharp}.}.
\end{enumerate}
 \item If $\T\in \dom(\Lambda)$ then $\T$ is according to $\Lambda$, i.e., for every limit ordinal $\eta'<\eta$, $\T_{<\eta'}\in \dom(\Lambda)$ and $\T_{\leq \eta'}=\T^{\eta'\frown} \{x\}$ where $x=\Lambda(\T_{<\eta'})$. 
\end{enumerate}
$\myqedhere$
\end{definition}

We say that $\T$ is a $(\k, \l)$-st-stack on $\P$ if $\T$ is an st-stack on $\P$ such that $o.t.(R^\T)<\k$ and for every $\tau<o.t.(R^\T)$, $\lh(\T^\tau)< \l$. As we said above, we could define the concept of putative st-stack similarly to \rdef{putative it}. As doing this is straightforward, we leave it to the reader. Putative essentially means that all models in the stack except possibly the last one are well-founded.

\begin{definition}\label{st-strategy}  Suppose $\P$ is a hod-like lsa type $\phi$-indexed ${\sf{lses}}$. We say $\Lambda$ is a \textbf{$(\k, \l)$-st-strategy} for $\P$ if the following clauses hold.
 \begin{enumerate}
 \item $\Lambda$ is an st-strategy.
 \item If $\T$ is a putative $(\k, \l)$-st-stack that is according to $\Lambda$ then $\T$ is a $(\k, \l)$-stack. 
 \item If $\T$ is a $(\k, \l)$-st-stack that is according to $\Lambda$ such that 
 \begin{enumerate}
 \item $\lh(\T)$ is a limit ordinal and 
 \item if $o.t.(R^\T)=\tau+1$ then $\lh(\T^{\tau})+1<\l$,
 \end{enumerate}
  then $\T\in \dom(\Lambda)$.
\end{enumerate}
\end{definition}

As we said above, we can then define generalized $(\k, \l, \nu)$-st-strategy which acts on generalized st-stacks. The definition of this notion is rather straightforward. 


Suppose now $\P$ and $\Lambda$ are as in \rdef{st-strategy}. We let $b(\Lambda)$ be the set of all $\T\in \dom(\Lambda)$ such that $\T$ has a last normal component of limit length and $\Lambda(\T)$ is a cofinal wellfounded branch of $\T$. Let $m(\Lambda)=\dom(\Lambda)-b(\Lambda)$\index{$b(\Lambda)$}\index{$m(\Lambda)$}. We call $m(\Lambda)$ the \textit{model component of $\Lambda$}.\index{model component of a short tree strategy}
 Given $\U\in \dom(\Lambda)$ such that the last component of $\U$ has a limit length, we let
\begin{center}
$\M(\Lambda, \U)=\begin{cases}
\M^{\U}_b &: \Lambda(\U)=b\\
\Lambda(\U) &: \text{otherwise}.
\end{cases}$\index{$\M(\Lambda, \U)$}
\end{center} 
If $\Lambda$ is an st-strategy for $\P$ and $\T$ is a stack on $\P$ according to $\Lambda$ with last model $\N$ then we let $\Lambda_{\N, \T}$ be the short tree strategy of $\N$ induced by $\Lambda$, i.e., for every $\U$ on $\N$, $\Lambda_{\N, \T}(\U)=\Lambda(\T^\frown \U)$. Here $\T^\frown \U$ is an st-stack defined in a natural way so that the normal components of $\T$ and $\U$ are the normal components of $\T^\frown \U$. $\myqedhere$
\begin{rem}
In many situations, it is expected that finding $(\k, \l)$-st-strategies must be easy. For example, whenever $\T$ is normal iteration tree of length $\omega$ such that $\mathcal{J}_{\omega}(\m^+(\T))\models ``\d(\T)$ is a Woodin", we can set $\Lambda(\T)=\m^+(\T)$. Thus, instead of working hard to define the correct branch, we declare success by setting $\Lambda(\T)=\m^+(\T)$. However, we will be interested in st-strategies that have certain fullness preservation properties. For instance, suppose $\M$ is just a suitable mouse in the sense of $L(\mathbb{R})$. If we now demand that $\Lambda$ must have the property that whenever $\T\in \dom(\Lambda)$ is such that $\Q(\T)$ exists then $\Lambda(\T)$ must be a branch $b$ with the property that $\Q(b, \T)=\Q(\T)$ then $\Lambda$ would be a rather complex object. We will have that $\Lambda(\T)$ is a model only in the case when $\T$ is a maximal iteration tree. In this case, $\Lambda$ is in fact a ``short tree iterability strategy" in the sense of $L(\mathbb{R})$. Such strategies are difficult to construct, and in our current situation, we will be interested in a notion of fullness preservation with respect to a much more complicated pointclass than $(\Sigma^2_1)^{L(\mathbb{R})}$. $\myqedhere$
\end{rem}

\section{Hull and branch condensation for short tree strategy} 

The goal of this section is to introduce \textit{hull condensation} for st-strategies. Hull condensation for iteration strategies was introduced in Definition 1.31 of \cite{ATHM}. It is an important property that is used to show that when doing hod pair constructions no discrepancies arise due to the coring down process. Thus if $\T$ is according to a strategy with hull condensation and $\U$ is a \textit{hull} of $\T$ (cf. Definition \ref{hull of a normal tree}) then it is according to the strategy. 

The difference between strategies and st-strategies is that st-strategies have a model component, and this difference causes some complications when trying to outright generalize hull condensation: such a direct generalization leads to a very strong property. Our definition is based on our indexing scheme \rdef{sts indexing scheme}. In short tree strategy mice, we only index branches of a certain kinds of iterations, and we need to apply hull condensation to these types of iterations. We start by introducing these iterations. 

First we define the \textit{universally short normal trees} which are essentially those normal iteration trees that are short with respect to any iteration strategy. 

\begin{definition}\label{normal stack} We say that $\T$ is a \textbf{normal stack} on $M$ if letting
\begin{center}
$\T=((\M_\a)_{\a<\eta}, (E_\a)_{\a<\eta-1}, D, R, (\b_\a, m_\a)_{\a\in R}, T)$,
\end{center}
 for all $\a<\eta$, $\b_\a=\ord(\M_\a)$, $m_\a=k(\M_\a)$ and setting
 \begin{center}
$\U=((\M_\a)_{\a<\eta}, (E_\a)_{\a<\eta-1}, D, T)$,
\end{center}
$\U$ is a normal iteration tree\footnote{Recall our general convention that all cutpoints of a stack a $\W$ belong to $R^\W$.}.  Given an st-stack 
\begin{center}
$\T=((\M_\a)_{\a<\eta}, (E_\a)_{\a<\eta-1}, D, R, (\b_\a, m_\a)_{\a\in R}, \so, \ma, T)$,
\end{center}
we say $\T$ is a \textbf{normal stack} if
\begin{enumerate}
\item $\ma=\emptyset$ and letting 
 \begin{center}
$\U=((\M_\a)_{\a<\eta}, (E_\a)_{\a<\eta-1}, D, T)$,
\end{center}
$\U$ is a normal iteration tree, or
\item ${\card{\ma}}=1$ and letting $\a$ be the unique element of $\ma$, ${\sf{next}}^\T_\a=\lh(\T)$ and 
 \begin{center}
$\U=((\M_\a)_{\a<\eta-1}, (E_\a)_{\a<\eta-1}, D, T)$,
\end{center}
$\U$ is a normal iteration tree.
\end{enumerate}
$\myqedhere$
\end{definition}

\begin{definition}[Universally short stacks]\label{nus stacks} Suppose $\P$ is a hod-like $\#$-lsa type ${\sf{lses}}$ and  $\T$ is a normal stack on $\P$ (see \rdef{normal stack}) such that $\lh(\T)$ is a limit ordinal. We say 
\begin{center}
$\T=((\M_\a)_{\a<\eta}, (E_\a)_{\a<\eta-1}, D, R, (\b_\a, m_\a)_{\a\in R}, T)$,
\end{center}
is \textbf{universally short} \index{${\sf{uvs}}$} (${\sf{uvs}}$) if one of the following holds:
\begin{enumerate}
\item $\pi^{\T, b}$ is undefined.\\\\
Suppose next $\pi^{\T, b}$ is defined and let $\a<\lh(\T)$ be the least such that $\pi_{0, \a}$ is defined and $\M_\a^b=\pi^{\T, b}(\P^b)$. It then follows that $\T_{\geq \a}$ is a stack on $\M_\a$ that is above $\ord(\M_\a^b)$\footnote{The condition that $\pi_{0, \a}$ is defined follows from the equality $\S^b=\pi^{\T, b}(\P^b)$.}. 
\item $R$ is cofinal in $\lh(\T)$.  
\item $\T$ has a fatal drop (see \rdef{fatal drop}).
\item For some $\b\in R-(\a+1)$, $D\cap (\a, \b]_\T\not=\emptyset$. 
\item For some $\b\in R-(\a+1)$ and some $\eta\in (\d^{\S^b}, \d^{\M_\b})$, $\T_{\geq \b}$  is a normal stack on $\M_\b$ that is below $\eta$. 
\item There is $\Q\insegeq \m(\T)^\#$ such that $\Q\models ``\d(\T)$ is a Woodin cardinal"  and $\mathcal{J}_{\omega}(\Q)\models ``\d(\T)$ isn't a Woodin cardinal".  

\end{enumerate}
If $\T$ is not ${\sf{uvs}}$ then we say that $\T$ is \textbf{non-universally short} (${\sf{nuvs}}$). $\myqedhere$
\end{definition} 


\begin{definition}[Indexable stack]\label{indexable stack}\index{indexable stack} Suppose $\P$ is a hod-like $\#$-lsa type ${\sf{lses}}$\footnote{See \rdef{lsa type}.}. We say that an st-stack\footnote{See \rdef{st-stack}.} 
\begin{center}
$\T=((\M_\a)_{\a<\eta}, (E_\a)_{\a<\eta-1}, D, R, (\beta_\a, m_\a)_{\a\in R}, \so, \ma, T)$
\end{center}
 is an \textbf{indexable stack} on $\P$ if one of the following clauses hold:
\begin{enumerate}
\item ${\sf{max}}=\emptyset$ and there is $\a\in R^\T$ such that $\pi^{\T_{\leq \a}, b}$ is defined and $\T_{\geq \a}$ is based on $\pi^{\T_{\leq \a}, b}(\P^b)$.
\item $\card{{\sf{max}}}=1$, $\T$ is a normal stack\footnote{See \rdef{normal stack}.} and if $\a$ is the unique element of ${\sf{max}}$ then $\pi^\T_{0, \a}$ is defined and ${\sf{next}}^\T(\a)=\lh(\T)$\footnote{It follows that $\T_{\geq \a}$ is above $\pi^{\T}_{0, \a}(\d^{\P^b})$. See also \rnot{notation for iteration trees}.}.
\end{enumerate}

Below and elsewhere we will use the notation $\T=(\P_0, \T_0, \P_1, \T_1)$ to denote indexable stacks. Here $\T_0=\T_{\leq \a}$ where $\a$ is either as in clause 1 or 2  and $\T_1=\T_{\geq\a}$. We will say that the indexable stack is \textbf{ordinary} if ${\sf{max}}^\T=\emptyset$. $\myqedhere$
\end{definition}

The iterations that we will index in short tree strategy mice are finite st-stacks of length 2. We define hull condensation for such stacks. 

\begin{definition}[Hull of a stack]
\label{hull of a normal tree}
Suppose $\M$ and $\M'$ are hod-like ${\sf{lses}}$ and $\T$ and $\T'$ are stacks on $\M$ and $\M'$ respectively. Set 
\begin{center}
$\T=((\M_\a)_{\a<\eta}, (E_\a)_{\a<\eta-1}, D, R, (\beta_\a, m_\a)_{\a\in R},  T)$\\
$\T'=((\M'_\a)_{\a<\eta}, (E'_\a)_{\a<\eta-1}, D', R', (\beta'_\a, m_\a)_{\a\in R'}, T')$.
\end{center}
Let $(\iota_\beta :\b\leq o.t.(R))$ and $(\iota'_\gg: \gg\leq  o.t.(R'))$ be the $\iota$-sequences of $\T$ and $\T'$ respectively (see \rdef{st-stack}). Let $i_{\a, \b}=\pi_{\a, \b}^\T$ and $i'_{\a, \b}=\pi^{\T'}_{\a, \b}$ provided the aforementioned embeddings exist.

We say $(\M', \T')$ is a \textbf{hull} of $(\M,\T)$ if there is a tuple 
\begin{center}
$(\sigma,  (\tau_\a)_{\a<\lh(\T')})$
\end{center}
such that the following clauses hold.
\begin{enumerate}
\item $\sigma : \lh(\T') \rightarrow \lh(\T)$ is an injective map that preserves the tree order and is such that $\sigma[R']\subseteq R$ and $\sigma(0)=0$.
\item For all $\a, \b$ such that $\a+\b<\lh(\T')$, $\sigma(\a+\b)=\sigma(\a)+\sigma(\b)$.
\item For every $\b<o.t.(R')$, $\sigma(\iota'_{\b+1})=\iota_{\sigma(\b)+1}$.
\item For every $\a<\lh(\T')$, $\tau_\a:\M'_\a\rightarrow_{\Sigma_1} \M_{\sigma(\a)}$ and $E_{\sigma(\a)}=\tau_\a(E'_\a)$.
\item For all $\a<\b<\lh(\T')$, $[\a, \b]_{\T'}\cap D'=\emptyset \iff [\sigma(\a), \sigma(\b)]_\T\cap D=\emptyset$. 
\item For every $\a<\b<\lh(\T')$, if $\sup(R'\cap (\a+1))=\sup(R'\cap \b)$ then 
\begin{center}
$\tau_\a\rest \lh(E'_\a)+1=\tau_\b\rest \lh(E'_\a)+1$.
\end{center}
\item For every $\a<\b<\lh(\T')$ such that $\a\leq_{\T'} \b$ and $(\a, \b]_{\T'}\cap R'=\emptyset$, 
\begin{center}
$\tau_\b\circ i'_{\a, \b} =
i_{\sigma(\a), \sigma(\b)}\circ \tau_\a$.
\end{center}
\item For every $\a+1<\lh(\T')$, if $\b=\T'(\a+1)$ then $\sigma(\b)=\T(\sigma(\a)+1)$ and 
\begin{center}
$\tau_{\a+1}([a,f]_{E'_\a})=[\tau_\a(a), \tau_\b(f)]_{E_{\sigma(\a)}}$.
\end{center}
 \end{enumerate}
 We say $(\sigma, (\tau_\a)_{\a<\lh(\T')})$ \textbf{witnesses} that $(\M', \T')$ is a hull of $(\M, \T)$.
 
 If $\M=\M'$ then we say that $(\M, \T')$ is a \textbf{hull} of $(\M, \T)$ if there is a  tuple $(\sigma, (\tau_\a)_{\a<\lh(\T')})$ witnessing that $(\M, \T')$ is a hull of $(\M, \T)$ and such that $\tau_0=id$. 
 
Both in the case $\M=\M'$ and $\M\not =\M'$, it is not ambiguous to simply say that $\T'$ is a hull of $\T$ to mean that $(\M', \T')$ is a hull of $(\M, \T)$, and so we will use this terminology\footnote{Notice that in the case $\M=\M'$, we must have that $\tau_0=id$.}. $\myqedhere$
\end{definition}

\begin{figure}
\centering
\begin{tikzpicture}[node distance=2cm, auto]
  \node (A) {$\M$};
  \node (B) [below of=A] {$\M$};
  \draw[->] (B) to node {$id$}(A);
  \node (C) [right of=A] {$\M^\T_{\sigma^0(\alpha)}$};
  \node (D) [right of=B] {$\M^\U_\alpha$};
  \draw[->] (D) to node {$\tau^0_\alpha$}(C);
  \draw[->] (A) to node {$\pi^\T_\alpha$}(C);
   \draw[->] (B) to node {$\pi^\U_\alpha$}(D);
 \node (E) [right of=C] {$\M_2$};
  \node (F) [right of=D] {$\M_1$};
\begin{scope}[dashed]
\draw (C) -- (E);
\draw (D) -- (F);
\end{scope}
  \draw[->] (F) to node {$\tau_0^1$} (E);
  \node (G) [node distance=3cm,  right of =E] {$\M_{\sigma^1(\gg)}^\S$};
  \node (H) [node distance=3cm, right of =F] {$\M_{\gg}^\W$};
  \draw[->] (F) to node {} (H);
  \draw[->] (E) to node {} (G);
  \node (I) [right of =H] {};
  \node (K) [right of =G]{};
  \begin{scope}[dashed]
  \draw (H) -- (I);
  \draw (G) -- (K);
  \end{scope}
  \draw[->] (H) to node {$\tau^1_\gg$}(G);
  \draw[->,bend left=37] (A) to node {$\pi^{\T,b}$}(E);
  \draw[->, bend right=37] ( B) to node {$\pi^{\U,b}$} (F);
   \end{tikzpicture}
\caption{Hull of a stack of length $2$. $(\M, \U, \M_1, \W)$ is a hull of $(\M,\T, \M_2, \S)$.}
\label{fig:hull}
\end{figure}

\begin{definition}[Hull of an indexable stack]
\label{hull of a stack of length 2}
(See  Figure \ref{fig:hull}.) Suppose $\M$ is a hod-like $\#$-lsa type ${\sf{lses}}$ and 
\begin{center}
$u=(\M, \U, \M_1, \W)$\\
$t=(\M, \T, \M_2, \S)$
\end{center}
 are two indexable stacks. We say $(\M, u)$ is a \textbf{hull} of $(\M,t)$ if either
\begin{enumerate}
\item  both $u$ and $t$ are ordinary (see \rdef{indexable stack}) and $(\M, u)$ is a hull of $(\M, t)$ (in the sense of \rdef{hull of a normal tree}) or
\item both $u$ and $t$ are not ordinary, and there are two tuples $(\sigma^0, (\tau_\a^0)_{\a<\lh(\U)})$ and $(\sigma^1, (\tau_\a^1)_{\a<\lh(\W)})$ such that the following holds.
\begin{enumerate}
\item $(\sigma^0, (\tau_\a^0)_{\a<\lh(\U)})$ witnesses that $(\M, \U)$ is a hull of $(\M, \T)$.
\item $(\sigma^1, (\tau_\a^1)_{\a<\lh(\W)})$ witnesses that $(\M_1, \W)$ is a hull of $(\M_2, \S)$.
\item $\tau_0^1\rest (\M_1^b)\circ \pi^{\U, b}=\pi^{\T, b}$. 

\end{enumerate}
\end{enumerate}
$\myqedhere$
\end{definition}

To finally define hull condensation for short tree strategy, we need to introduce a few more definitions. Suppose $(\P, \Sigma)$ such that $\P$ is a hod-like $\#$-lsa type ${\sf{lses}}$ and $\Sigma$ is a st-strategy for $\P$. First we introduce two sorts of iterates of $(\P, \Sigma)$, $I^b(\P, \Sigma)$ and $I(\P, \Sigma)$. 

\begin{notation} Suppose $\P$ is a hod-like $\#$-lsa type ${\sf{lses}}$\footnote{See \rdef{lsa type}.} and $\Sigma$ is a st-strategy\footnote{See \rdef{st-strategy}.} for $\P$. We let
$\max(\P, \Sigma)$\index{$\max(\P, \Sigma)$} be the set of $\Sigma$-maximal iterations. More precisely, $\max(\P, \Sigma)$ consists of pairs $(\T, \Q)$ such that $\T\in m(\Sigma)$ and $\Q=\m^+(\T)$. $\myqedhere$
\end{notation}
 In the following definition, we recycle the notations used in Definition \ref{almost non-dropping stacks}. The difference here is that $\Sigma$ is the short-tree strategy.

\begin{definition}[$I^b(\P, \Sigma)$ and $I(\P, \Sigma)$]\label{short tree iterates} Suppose $(\P, \Sigma)$ is a pair such that $\P$ is a hod-like $\#$-lsa type ${\sf{lses}}$ and $\Sigma$ is an st-strategy for $\P$.  We then let 
\begin{center}
$I^b(\P, \Sigma)=\{(\T, \Q): \T$ is according to $\Sigma$, $\Q$ is the last model of $\T$ and $\pi^{\T, b}$ exists$\}$,
\end{center}
\begin{center}
$I(\P, \Sigma)=\{(\T, \Q):$ either $(\T, \Q)\in \max(\P, \Sigma)$ or for some $\b\in {\sf{max}}(\T)$, $\pi^{\T_{\geq \b}}$  exists$\}$.
\end{center}
\end{definition}

\begin{notation}\label{hc set up} We let ${\sf{HC}}$ be the set of all hereditarily countable sets. In \rdef{coding of countable objects}, we fix a coding of elements of $\sf{HC}$ by reals. This coding then induces a coding of elements of $\cup_{n\in \omega}\powerset({\sf{HC}}^n)$ by sets of reals. Let ${\sf{Code}}$ be the coding function introduced in \rdef{coding of countable objects}. Thus for $A\subseteq {\sf{HC}}^n$, ${\sf{Code}}(A)$ is the set of reals that codes $A$. $\myqedhere$
\end{notation}

\begin{definition}\label{gamma(p, sigma) and b(p, sigma) for sts} Suppose $(\P, \Sigma)$ is a pair such that $\P$ is a hod-like $\#$-lsa type ${\sf{lses}}$ and $\Sigma$ is an st-strategy for $\P$. We then let 
\begin{center}
$B(\P, \Sigma)=\{(\T, \Q) : \exists \R ((\T, \R)\in I^b(\P, \Sigma)\wedge \Q\insegeq_{hod}\R^b)\}$,
\end{center}
and 
\begin{center}
$\Gamma^b(\P, \Sigma)=\{ A\subseteq \mathbb{R} : \exists (\T, \Q) \in B(\P, \Sigma) (A\leq_w {\sf{Code}}(\Sigma_{\Q, \T}))\}$.
\end{center}
$\myqedhere$
\end{definition}

\begin{definition}[Hull condensation\index{hull condensation}]\label{hull condensation} Suppose $\P$ is a hod-like $\#$-lsa type ${\sf{lses}}$ and 
$\Sigma$ is a st-strategy for $\P$. We say $\Sigma$ has \textbf{hull condensation} if the following clauses hold.
\begin{enumerate}
\item For all $(\T, \Q)\in B(\P, \Sigma)$, $\Sigma_{\Q, \T}$ has hull condensation, and
\item Whenever $(\T, \Q)\in I(\P, \Sigma)$,  $u=(\Q, \U, \Q_1, \W)$ and $t=(\Q, \T', \Q_2, \W')$ are two indexable stacks on $\Q$ such that $t$ is according to $\Sigma_{\Q, \T}$ and $(\Q, u)$ is a hull of $(\Q, t)$ then $u$ is according to $\Sigma_{\Q, \T}$.
\end{enumerate} 
$\myqedhere$
\end{definition}

\begin{figure}
\centering
\begin{tikzpicture}[node distance=2cm, auto]
  \node (A) {$\P$};
  \node (B) [right of=A] {$\Q$};
  \node (C) [node distance=0.5cm, right of=B] {$\triangleright$};
  \node (D) [node distance=0.7 cm, right of = C] {$\Q(\beta)$};
  \node (X) [node distance=1.5cm, below of=B] {$\R$};
  \node (E) [node distance=0.5cm, below of=X] {$\triangledown$};
  \node (G) [node distance=0.5cm, below of=E] {$\R^b$};
  \node (H) [right of=G] {$\M^\S_c$};
  
  \draw[->] (A) to node {$\T$}(B);
  \draw[->] (G) to node {$\S, c$} (H);
  \draw[->] (H) to node  {$\tau$} (D);
  \draw[->] (A) to node {$\U$} (X);
   \end{tikzpicture}
\caption{Branch condensation for short tree strategies. Notations are as in Definition \ref{branch condensation for short tree strategies}. In the above, $\pi^{\T,b} = \pi \circ \pi^\S_c \circ \pi^{\U,b}$.}
\label{fig:branch_cond_sts}
\end{figure}

Next we introduce branch condensation for short tree strategies. We will need this notion in the definition of hod mice (see \rdef{hod premouse}). 

\begin{definition}[Branch condensation for st-strategies]\label{branch condensation for short tree strategies} (See Figure \ref{fig:branch_cond_sts}.) Suppose 
 $(\P, \Sigma)$ is such that $\P$ is a hod-like $\#$-lsa type ${\sf{lses}}$ and $\Sigma$ is a st-strategy for $\P$. We say $\Sigma$ has \textbf{branch condensation} if whenever $(\T, \Q, \U, \R, \tau, \S, c, \a, \b)$ is such that 
 \begin{enumerate}
 \item $(\T, \Q), (\U, \R)\in I^b(\P, \Sigma)$,
 \item $\a<\l^{\R^b}$ and $\d^{\R(\a+1)}$ is a Woodin cardinal of $\R$\footnote{See \rnot{l p} for the definition of $\R(\tau)$.},
 \item $\S$ is a normal iteration tree of limit length according to $\Sigma_{\R^b, \U}$ that is based on $\R(\a+1)$ and is above $\d_\a^\R$,
 \item $c$ is a branch of $\S$ such that $\pi^{\S}_c$ exists, and 
 \item $\tau:\M^\S_c\rightarrow \Q(\b)$ and $\pi^{\T, b}=\tau\circ \pi^{\S}_c \circ \pi^{\U, b}$
 \end{enumerate}
 then $c=\Sigma_{\R, \U}(\S)$. $\myqedhere$
\end{definition}

\section{St-type pairs}

Suppose $\P$ is a hod-like $\#$-lsa type ${\sf{lses}}$\footnote{See \rdef{lsa type}.} and suppose $\Lambda$ is an st-strategy for $\P$. We would like to introduce the notion of  short tree premice and in particular, $\Lambda$-premice. The main technical problem is that we do not have a reasonable notion of condensation for st-strategies.  In particular, if $\Lambda=\Sigma^{stc}$ for some strategy $\Sigma$, then it may well be that there is a $\Sigma$-maximal iteration tree $\T$ on $\P$ such  there is a $\Sigma$-short hull $\U$ of $\T$. 

The above scenario is the main difficulty with defining short tree strategy mice. We have to find a particular indexing of short tree strategies, or rather carefully skip over ``bad trees", in a way that when $\T$ above is ``cored down" to $\U$ above then our indexing is still preserved. In particular, the branch of $\T$ cannot be added too early. The idea is to wait until the branch of $\T$ or rather its correct $\Q$-structure is \textit{certified}. Before we define short tree hybrids, however, we have to make a few definitions that will be useful to us in the future.

We will only consider st-strategies $\Lambda$ with the property that whenever $\T\in \dom(\Lambda)$ is ${\sf{uvs}}$ then $\Lambda(\T)$ is a branch. 

\begin{definition}[Faithful short tree strategy]\label{faithful short tree strategy}\index{faithful strategy}
Suppose $\P$ is a hod-like $\#$-lsa type ${\sf{lses}}$  and $\Lambda$ is a $(\k, \l, \eta)$-st-strategy for $\P$. We say $\Lambda$ is a \textbf{\textit{faithful}\index{faithful} $(\k, \l, \eta)$-st-strategy} if whenever $\T\in \dom(\Lambda)$ is ${\sf{uvs}}$, $\T\in b(\Lambda)$. $\myqedhere$
%
\end{definition}

%
%

\begin{definition}[St-type pair]\label{lsa type pair} 
We say $(\P, \Lambda)$ is a \textbf{hod-like st-type pair} if 
\begin{enumerate}
\item $\P$ is a hod-like $\#$-lsa type ${\sf{lses}}$,
\item  $\Lambda$ is a faithful $(\omega_1, \omega_1, \omega_1)$-st-strategy,
\item if $\Q$ is a $\Lambda$-iterate of $\P$ via $\T$ and $\R\in Y^\Q$ then $\Sigma^\R\subseteq \Lambda_{\R, \T}$\footnote{This clause asserts that the internal strategy of $\R$ agrees with $\Lambda_{\R, \T}$.}.
\item $\Lambda$ has hull condensation\footnote{See \rdef{hull condensation}.}.
\end{enumerate}

Similarly we can define \textbf{simple} hod-like st type pairs by demanding that $\Lambda$ is a faithful $(\omega_1, \omega_1)$-strategy and that clause 3 above holds. $\myqedhere$

\end{definition}

\section{$(\P, \Sigma)$-hod pair construction}\label{(p, sigma)-hod pair constructions}

Suppose that $(\P, \Sigma)$ is a hod-like st-type pair. Below we describe a fully backgrounded construction that, if successful, constructs a $\Sigma$-iterate of $\P$. To learn more about such backgrounded constructions the reader may consult \cite[Chapter 11]{FSIT} and also various papers of Schlutzenberg and Schindler-Steel-Zeman that deal with certain fine structural issues present in \cite{FSIT} (for example, \cite[Chapter 2.2, Definition 2.4]{OIMT} and the discussion after it, and also \cite{FarmerRes} and \cite{SSZ}).  We say a $(\k, \l)$-extender $E$ coheres $\Sigma$ if $\P\in V_\k$, $V_\l\subseteq Ult(V, E)$ and $\pi_E(\Sigma)\rest V_\l=\Sigma\rest V_\l$. 

In this manuscript, our goal is to deal with novel issues arising  from the theory of short tree strategy mice, such as developing an indexing scheme for short tree strategies, proving a comparison theorem for lsa small hod pairs and obtain core model induction applications at the level of $\sf{LSA}$, to list a few. We don not have space to also carefully develop the theory of fully backgrounded constructions, but all issues that arise have been handled in literature. For example, to deal with issues arising from our mixing indexing we refer the reader to Schlutzenberg's \cite{schlutzenberg2021background} and to deal with issues regarding inheriting large cardinals we refer the reader to \cite[Chapters 9-12]{FSIT} and to \cite{FarmStrongs}.

Unlike in \cite{ANS} and \cite{FSIT}, and other similar places in literature where the convergence of the backgrounded constructions is established, here we will not be concerned with iterability issues of the backgrounded constructions and just simply assume that such constructions converge provided the background universe is iterable. Our assumption is justified by the results of \cite[Chapter 12]{FSIT}. The consequence of our assumption is that in clause (3) below we simply take the core rather than the dropdown sequence. See \rdef{the core} for the definition of $\C$. 

\begin{definition}[$(\P, \Sigma)$-coherent fully backgrounded constructions]\label{full short tree coherent background constructions}
Suppose $\k$ is an inaccessible cardinal and $(\P, \Sigma)$ is a hod-like st-type pair such that $\Sigma$ is a $(\k, \k, \k)$-st-strategy. Then for $\eta\leq \k$, we say $(( \M_\gg , \N_\gg : \gg\leq \eta), ( F_\gg: \gg<\eta), ( \T_\gg: \gg\leq\eta))$
is the output of the \textbf{$(\P, \Sigma)$-coherent fully backgrounded construction} of $V_\k$ if the following holds.
\begin{enumerate}
\item $\M_0=\emptyset$.
\item $\M_\gg$  is a hod-like ${\sf{lses}}$ such that there is a tree $\T_\gg$\footnote{Notice that if there is such a $\T_\gg$ then it is unique.} on $\P$ according to $\Sigma$ such that either 
\begin{enumerate}
\item $\T_\gg$ has a last model $\M$ such that if $\xi=\ord(\M_\gg)$ then $\M_\gg|\xi=\M|\xi$ or 
\item $\M_\gg=\m(\T_\gg)$. 
\end{enumerate}
\item Suppose $\gg\leq \eta$ is such that either $\T_\gg$ has a last model or $\T_\gg\in b(\Sigma)$. Let $\M$ be the last model of $\T_\gg$ if it exists and otherwise, setting $b=\Sigma(\T_\gg)$, let $\M=\M^{\T_\gg}_b$. Let $\xi=\ord(\M_\gg)$ and suppose $\M_\gg=\mathcal{J}^{\vec{E}, f}_\xi$. Then the following statements hold.
\begin{enumerate}
\item If $\M_\gg=\M$ then $\gg=\eta$.
\item If $\M_\gg$ is active and $\M_\gg\not=\M||\xi$ then $\gg=\eta$.
\item If $\M_\gg$ is active and $\M_\gg=\M||\xi$ then $\N_{\gg+1}=\mathcal{J}_\omega[\M_\gg]$ and $\M_{\gg+1}={\sf{core}}(\N_\gg)$.
\item Suppose $\M_\gg$ is passive and $\M_\gg\inseg \M$. Suppose there is no pair $(F^*, F)$ and an
ordinal $\zeta<\xi$ such that $F^*\in V_\k$ is an extender over $V$ that coheres $\Sigma$, $F$ is an extender over $\M_\gg$,    $V_{\zeta+\omega}\subseteq Ult(V, F^*)$ and 
    \begin{center}
    $F=F^*\cap ([\zeta]^{\omega}\times \mathcal{J}_\xi^{\vec{E}, f})$
    \end{center}
 such that $(\mathcal{J}_\xi^{\vec{E}, f}, \in, \vec{E}, f, \tilde{F})$ is a hod-like ${\sf{lses}}$\footnote{Here $\tilde{F}$ is the amenable code of $F$, see the discussion after \cite[Lemma 2.9]{OIMT}.}. Then $\N_{\gg}=\mathcal{J}_\omega(\M_\gg)$ and $\M_{\gg+1}={\sf{core}}(\N_\gg)$.
\item Again suppose $\M_\gg$ is passive and $\M_\gg\inseg \M$ but there is a pair $(F^*, F)$ and an ordinal $\zeta$ satisfying the above conditions. 
Then if $F\in \vec{E}^\M$ then we let 
    \begin{center}
    $\N_{\gg}=(\mathcal{J}_\xi^{\vec{E}, f}, \in, \vec{E}, f, \tilde{F})$ and $\M_{\gg+1}={\sf{core}}(\N_{\gg})$.
    \end{center}
\item Again suppose $\M_\gg$ is passive, $\M_\gg\inseg \M$ and that $\M||\xi$ is an active $\mathcal{J}$-structure such that its last predicate codes a set $A$ that is not an extender. Let then $\N_\gg=(\M_\gg, A, \in)$\footnote{Here we mean that $A$ is being indexed in the strategy predicate of $\N_\gg$.} and $\M_{\gg+1}={\sf{core}}(\N_\gg)$. 
\end{enumerate}
\item Suppose $\gg\leq \eta$ is such that $\T_\gg$ is of limit length and $\T_\gg\not \in b(\Sigma)$. Then $\gg=\eta$.
\end{enumerate}
$\myqedhere$
\end{definition} 

\begin{rem}\label{partial strategy remark} Notice that the constructions introduced in \rdef{full short tree coherent background constructions} can be carried out even when $\Sigma$ is a partial strategy. Thus, for example, we may say that $``(( \M_\gg , \N_\gg : \gg\leq \eta), ( F_\gg: \gg<\eta), ( \T_\gg: \gg\leq\eta))$
is the output of the $(\P, \Sigma)$-coherent fully backgrounded construction of $N"$ the meaning of which should be self-evident with one wrinkle. It may be that for some $\gg\leq \eta$, $\Sigma(\T_\gg)$ is undefined. In this case, we have that $\gg=\eta$ and we stop the construction. 

If the background universe has a distinguished extender sequence then we tacitly assume that the extenders appearing in the $(\P, \Sigma)$-coherent fully background construction come from this distinguished extender sequence.  $\myqedhere$
\end{rem}

\section{A short tree strategy indexing scheme}\label{stsis sec}

Our goal here is to introduce the notion of a \textit{short tree strategy premouse} (\textit{sts premouse}). As we mentioned in the previous section, the difficulty with doing this lies in the fact that maximal trees might ``core down" to short trees and thus, creating  indexing issues. The idea behind the solution presented here is to add a branch for a tree as soon as we see a certificate of shortness, which in our case will be a $\Q$-structure. As the $\Q$-structures that we will be looking for are themselves sts premice, this  inevitably leads to an induction.

Technically speaking $\M$ in \rdef{unambiguous hp} should not be ${\sf{ses}}$ (see \rdef{strategy premouse}) as $f^\N$ doesn't quite code an iteration strategy. Its domain consists of indexable stacks (see \rdef{indexable stack}). But recall the abuse of terminology proposed after \rdef{semi-strategy}. Also, recall the definition of $\m^+(\T)=\m(\T)^\#$ (see \rdef{mtsharp}).

The language of unindexed $\sf{ses}$\footnote{See \rdef{unindexed ses}.}  includes constant symbols for $\vec{E}$, $f$, $X$ and $\P$. We denote these symbols by $\dot{E}$, $\dot{f}$, $\dot{X}$ and $\dot{\P}$. Also, we let $\dot{<}$ be the symbol denoting the constructibility order and $\dot{\Sigma}$ be the partial strategy coded by $\dot{f}$.  $\dot{<}$ and $\dot{\Sigma}$ are not symbols in the language but they can be easily defined from the other symbols. 

\begin{definition}[$\phi^*$-formula]\label{phi*} We let $\phi^*(x)$ be the conjunction of the  following statements in the language of $\sf{ses}$. 
\begin{enumerate}
\item $x$ is a sequential structure of the form $(\mathcal{J}_{\omega}(t), t, \in)$ where $t=(\P_0, \T_0, \P_1, \T_1)$ is an indexable stack on $\dot{\P}$,
\item $t$ is according to $\dot{\Sigma}$ where $\dot{\Sigma}$ is the partial strategy coded by $\dot{f}$, and
\item $\cf(\lh(\T_0))$ and $\cf(\lh(\T_1))$ are not measurable cardinals.
\item there is $(\nu, \xi)$ such that letting  $(( \M_\gg , \N_\gg : \gg\leq \eta), ( F_\gg: \gg<\eta), ( \W_\gg: \gg<\eta))$ be the output of the $(\dot{\P}, \dot{\Sigma})$-coherent fully backgrounded construction of the universe \footnote{See \rrem{partial strategy remark}.} in which extenders used have critical points $>\nu$\footnote{See \rdef{full short tree coherent background constructions}.}, $\W_\xi=\T_0$.
\end{enumerate}
$\myqedhere$
\end{definition}

\begin{definition}[Unambiguous ${\sf{ses}}$]\label{unambiguous hp}\index{unambiguous ${\sf{ses}}$} Suppose $\M$ is an unindexed ${\sf{ses}}$\footnote{See \rdef{unindexed ses}.} over some self-well-ordered set $X$ based on a hod-like $\#$-lsa type ${\sf{lses}}$ $\P$. We say 
$\M$ is \textbf{unambiguous} 
if $\M$ is closed\footnote{See \rdef{closed under sharps}.} and  whenever $w$ is a sequential structure of the form $(\mathcal{J}_\omega(t), t, \in)$ where $t=(\P_0, \T_0, \P_1, \T_1)\in \M$ is an indexable stack according to $\Sigma^\M$ and such that 
\begin{enumerate}
\item $\M\models \phi^*[w]$ and
\item either
\begin{enumerate}
\item $\T_1=\emptyset$ and $\M\models ``\T_0$ is a ${\sf{uvs}}$\footnote{See \rdef{nus stacks}.}  of limit length" or
\item  $\T_1$ is a nonempty stack of limit length
\end{enumerate}
\end{enumerate}
then $t\in \dom(\Sigma^\M)$.
We say $\M$ is ambiguous if it is not unambiguous. $\myqedhere$
\end{definition}

Notice that ambiguity is a first order property of unindexed $\sf{ses}$. The next definition introduces an indexing scheme that we will use to define short tree premice. The indexing scheme only defines the strategy on certain carefully chosen stacks. It turns out that this much information is enough to extend the strategy to all stacks (see \rchap{lsa internal theory chapter}). 

\begin{remark}\label{intuition behind ambiguity} The reader may find the following remark helpful. \rdef{nus stacks} introduced the ${\sf{uvs}}$ stacks, which are stacks that are short with respect to all reasonable strategies. \rdef{unambiguous hp} introduces unambiguous ${\sf{ses}}$, which are the ${\sf{ses}}$ whose internal strategy predicate is total on all indexable ${\sf{uvs}}$ stacks that satisfy the formula $\phi^*$ (see \rdef{phi*}). Negating this, we have that if $\N$ is ambiguous ${\sf{ses}}$ then in $\N$ there is a ${\sf{uvs}}$ $\T$ of limit length that is according to the internal strategy of $\N$ yet no branch of $\T$ is indexed in the strategy predicate of $\N$. $\myqedhere$
\end{remark}

\begin{definition}[$\psi$-sts indexing scheme]\label{sts indexing scheme} Suppose $\phi(x)$ and $\psi(x, y)$ are two formulas in the language of ${\sf{ses}}$. We say $\phi$ is a \textbf{$\psi$-sts indexing scheme} if $\phi(w)$ is the conjunction of the following clauses: 
\begin{enumerate}
\item For all ordinals $\gg$ there is $\xi>\gg$ such that $\dot{E}(\xi)$ is defined\footnote{i.e. the universe is closed in the sense of \rdef{closed under sharps}.}. 
\item $\dot{\Sigma}$ is a partial faithful st-strategy such that $m(\dot{\Sigma})=\emptyset$\footnote{Notice that clause 4 below guarantees that $\dot{\Sigma}$ is really a partial strategy rather than an st-strategy. We emphasize the fact that $\dot{\Sigma}$ is an st-strategy to point out the fact that there is no iteration according to $\dot{\Sigma}$ that is $\dot{\Sigma}$-maximal.}.
\item $\phi^*(w)$.
\item Either
\begin{enumerate}
\item The universe is ambiguous and $w$ is the $\dot{<}$-least sequential structure $w'$ witnessing ambiguity of the universe.\\
 Or
\item The universe is unambiguous and $w$ is the $\dot{<}$-least sequential structure $w'$ of the form $(\mathcal{J}_{\omega}(t), t, \in)$ with the property that $t=(\P_0, \T_0, \P_1, \T_1)$ is an indexable stack on $\dot{\P}$, $\phi^*(w')$ holds, $\dot{\Sigma}(\T_0)$ is undefined and there is a unique cofinal well-founded branch $b$ of $\T_0$ such that  $\psi(\T_0, b)$ holds.
\end{enumerate}

\end{enumerate}
$\myqedhere$
\end{definition}

\begin{remark}\label{useful remark about indexing} The reader may find it useful to compare \rdef{sts indexing scheme} with \rdef{important notation}, \rdef{passive hybrid j-structure} and \rdef{hybrid j-structure}. The model over which we intend to evaluate $\phi$ in \rdef{sts indexing scheme} corresponds to $\M|\omega\b$ in \rdef{important notation}. More precisely, if $\phi$ is as in \rdef{sts indexing scheme} and $w$ is a sequential structure, then to decide whether we need to index a branch of $w$ or not we need to look for $\b$ such that $\M|\omega\b\models \phi[w]$.

The meaning of clause 4b is that $\psi$ is the certification of $b$ as the correct branch, but \rdef{sts indexing scheme} doesn't say anything about a particular certification procedure that we will use. The exact certification method is presented in \rdef{weak psi alpha indexing scheme a}. $\myqedhere$
\end{remark}


Notice that $\phi$ is uniquely determined by $\psi$.  The next definition uses ideas from \rdef{important notation} and \rdef{passive hybrid j-structure}, and it may be useful to review those definitions (in particular clause 4a of \rdef{passive hybrid j-structure}).

\begin{definition}[Sts $\psi$-premouse]\label{sts phi premouse1} Suppose $X$ is a self-well-ordered set, $\P\in X$ is a hod-like $\#$-lsa type ${\sf{lses}}$ and $\psi(x, y)$ is a formula in the language of unindexed $\sf{ses}$. Let $\phi$ be the $\psi$-sts indexing scheme. Then $\M$ is an \textbf{sts $\psi$-premouse} over $X$ based on $\P$ if $\M$ is a $\phi$-indexed  ${\sf{ses}}$ over $X$ based on $\P$ and if $w\in \dom(f^\M)$ is such that clause 4b of \rdef{sts indexing scheme} applies to $w=_{def}(\mathcal{J}_\omega(t), t, \in)$ where $t=_{def}(\P_0, \T_0, \P_1, \T_1)$ then letting $\b=\min(f^\M(w))$, 
\begin{center}
$f^\M(w)=\{\b+\omega\gg: \gg\in b\}$
\end{center}
 where $b\in \M|\b$ is the unique branch of $\T_0$ such that $\M|\b\models \psi[\T_0, b]$. $\myqedhere$
\end{definition}

If $\psi(x, y)=``0=1"$ then we say $\M$ has a trivial indexing scheme and also say that $\M$ is a trivial sts premouse. Notice that in a trivial sts premouse indexable ${\sf{nuvs}}$  stacks do not have branches indexed in the strategy predicate.\\\\
\textbf{$\Q$-structures are sts $\psi$-premouse}\\\\
Suppose $\P$ is an $\#$-lsa type hod like $\sf{lses}$ and $\T$ is a normal $\sf{nuvs}$ tree on $\P$. Suppose $b$ is a well-founded branch of $\T$ such that $\Q(b, \T)$ exists. Does it follow that $\Q(b, \T)$ is an sts $\psi$-premouse in some reasonable sense? The following lemma gives the answer we need.

\begin{definition}\label{uniformly sts} Suppose $\P$ is a hod-like $\#$-lsa type $\sf{lses}$ and $\psi(x, y)$ is a formula in the language of unindexed $\sf{ses}$. We say $\P$ is \textbf{uniformly $\psi$-organized} if for each $\#$-lsa type layer $\Q$ of $\P$ such that $\Q^b=\P^b$ and $\d^\Q<\d^\P$, if $\nu$ is the largest such that $\P||\nu\models ``\d^\Q$ is a Woodin cardinal" then $\P||\nu$ is an sts $\psi$-premouse over $\Q$. $\myqedhere$
\end{definition}

\begin{lemma}\label{qstructures are sts} Suppose $\P$ is a uniformly $\psi$-organized hod-like $\#$-lsa type ${\sf{lses}}$. Suppose $\T$ is a normal iteration tree on $\P$. Suppose $\a<\lh(\T)$ and $\R\insegeq_{hod} \M^\T_\a$ is such that letting $(\P_{\xi, \xi'}: \xi\leq \eta, \xi'\leq \nu_\xi)$ be the layers of $\M^\T_\a$, for some $\xi\leq \eta$, $\R=\P_{\xi, 0}$ and $\P_{\xi, 1}$ is defined either according to condition $\sf{R5}$ of \rdef{layers of hod-like lsp} or clause 2 of $\sf{R10}$ of  \rdef{layers of hod-like lsp}\footnote{This simply means that $\R$ is a $\#$-lsa type layer of $\M^\T_\a$.}. Then $\P_{\xi, 1}$ is an sts $\psi$-premouse over $\R$.
\end{lemma}
\begin{proof} We prove the claim by induction. Suppose first $\a=\b+1$ and the claim is true for $\b$ (i.e. the claim is true for all $\zeta\leq \b$ and hod initial segments of $\M^\T_\zeta$). In this case, we have that $\M^\T_{\a}=Ult(\N, E^\T_\b)$ where $\gg=\T(\a)$ and $\N\insegeq \M^\T_\gg$ is the appropriate initial segment of $\M^\T_\gg$. If now $\d^\R<\d^{\M^\T_\a}$ then the claim follows from elementarity of $\pi^\T_{\gg, \a}$ restricted to strict initial segments of $\N$. 

Assume then that $\d^\R=\d^{\M^\T_\a}$. If $\N\models ``\d^\N$ is not a Woodin cardinal" then once again elementarity implies  the claim (as $\P_{\xi, 1}\in \rge(\pi^\T_{\gg, \a})$).  Assume then $\N\models ``\d^\N$ is a Woodin cardinal". In this case, we have that 
$\P_{\xi, 1}=Ult(\N, E_\b^\T)$ and $\N$ is an sts $\psi$-premouse over $\N|\tau$ where $\tau=\min(\vec{E}^\N-\d^\N)$. 

Set now $\Q=_{def}\P_{\xi, 1}$, $E=E_\b^\T$ and $j=\pi^\T_{\gg, \a}$.  Let $f$ be the strategy predicate of $\Q$ and suppose that $\Q$ is not an sts $\psi$-premouse over $\R$. Notice that because $j$ is $\Sigma_1$-elementary, we must have that for every $\omega\zeta<\sf{ord}(\Q)$, $\Q||\omega \zeta$ is an sts $\psi$-premouse over $\R$. This is because otherwise $\N$ would satisfy that there is some $\omega\zeta'+\omega<\sf{ord}(\N)$ such that $\N|\omega\zeta'+\omega\models ``\N||\omega\zeta'$ is not an sts $\psi$-premouse". 

Thus, it is enough to show that if $\Q$ is active and its last predicate is a pair $(\T, b)\in f$ then $(\T, b)$ conforms the rules of sts $\psi$-premice. Set then $w=(\mathcal{J}_\omega(\T), \T, \in)$ and $\nu=\min(f(w))$. Let $\nu'=j^{-1}(\nu)$ and $\T'=j^{-1}(\T)$. Because $j\rest \N|\nu'$ is fully elementary, we have that\\\\
(1) $\T'$ is chosen in $\N|\nu'$ according to clause 4a of \rdef{sts indexing scheme} if and only if $\T$ is chosen in $\Q|\nu$ according to clause 4a of \rdef{sts indexing scheme}.\\
(2) $\T'$ is chosen in $\N|\nu'$ according to clause 4b of \rdef{sts indexing scheme} if and only if $\T$ is chosen in $\Q|\nu$ according to clause 4b of \rdef{sts indexing scheme}.\\\\
Thus, to finish, we need to verify that letting $b^*=j^{-1}[b]$ and $b'$ be the closure of $b^*$ in $\T'$ then \\\\
(*) $b'$ is as in clause 4b of \rdef{sts indexing scheme} if and only if $b$ is as in clause 4b of \rdef{sts indexing scheme}.\\\\
(*) is straightforward because if $b'$ is as in clause 4b then $b=j(b')$, and if $b$ is as in clause 4b then $b'=j^{-1}(b)$\footnote{The equivalence follows from the fact that because $\cf^\N(\lh(\T'))$ is not a measurable cardinal in $\N$, $j\rest \lh(\T')$ is cofinal in $\lh(\T)$.}. 

The case when $\a$ is a limit ordinal is very similar, and we leave it to the reader.
\end{proof}

\section{Authentic indexable stacks}\label{authentic iterations and finite stacks sec}

Suppose $(\P, \Sigma)$ is a hod-like limit type pair. Suppose $\T$ is a tree on $\P$ according to $\Sigma$ such that $\pi^{\T, b}$ exists and $\m^+(\T)\models ``\d(\T)$ is a Woodin cardinal" (see \rdef{mtsharp}). When defining short tree strategy mice, we will be faced with the following question. How can we guess the correct branches of iteration trees that are on $\m^+(\T)$ and are according to $\Sigma_{\m^+(\T), \T}$? In this section, we present an authentication process that allows us to guess the correct branches of such iterations. 

The main technical object used in our authentication process is $s(\T, w)$ introduced in \rdef{canonical singularizing sequences}. In the light of \cite{NormIter}, we could use \textit{strong hull condensation} instead of $s(\T, w)$, similar to the way optimal Suslin representations are obtained in \cite[Chapter 2]{Mousepairs}. We, however, do not know if the core model induction applications of this book could be done using the ideas of \cite{Mousepairs}.

We start by recalling $s(\T, w)$ (and slightly modifying it). Suppose $\P$ is a non-meek hod-like ${\sf{lses}}$ and $\T$ is a stack on $\P$ such that $\pi^{\T, b}$ exists. Let $\S=\pi^{\T, b}(\P^{b})$,  $w=(\eta, \d)$ be a window\footnote{See \rdef{l p}.} of $\S$ and $X\subseteq \P^b$. We then set 
\begin{center}
$s(\T, X, w)=\{ \a: \exists a\in \eta^{<\omega} \exists f\in X( \a=\pi^{\T, b}(f)(a))\}\cap \d$\index{$s(\T, X, \xi)$}
\end{center}
When $X=\P^b$ then we just write $s(\T, w)$. 

\begin{definition}\label{useful x} Suppose $\P$ is a hod-like limit type ${\sf{lses}}$  and $X\subseteq \P^b$. We then say that $X$ is \textbf{useful} if whenever $\T$ is a stack on $\P$ such that $\pi^{\T, b}$ is defined, $\d$ is a Woodin cardinal of $\S=_{def}\pi^{\T, b}(\P^{b})$ and $w$ is a window of $\S$ such that $\d^w=\d$ then $s(\T, X, w)$ is cofinal in $\d$. $\myqedhere$
\end{definition}

Recall that \rlem{the canonical singularizing sequence exists} shows that $X=\P^b$ is useful.

\begin{notation}\label{hull notation}
Here and elsewhere in the manuscript , given a collection of formulas $\Gamma$, by $cHull^\M_\Gamma(Y)$ we mean the transitive collapse of $X$ where $a\in X$ if and only if there is a formula $\phi\in \Gamma$ and $s\in Y^{<\omega}$ such that $a\in \M$ is the unique $b$ with the property that $\M\models \phi[b, s]$. If $\Gamma$ contains all formulas then we omit it from the notation. If $\Gamma$ is the set of all $\Sigma_n$ formulas then we just write $cHull_n^\M(Y)$. If $\Gamma$ is the set of all formulas then we just write $cHull^\M(Y)$. $\myqedhere$
\end{notation}

\begin{figure}
\centering
\begin{tikzpicture}[node distance=3cm, auto]
  \node (A) {$X\subset \P^b$};
  \node (B) [right of=A] {$\S=\pi^{\T,b}(\P^b)$};
  \node (X) [node distance=0.5cm, below of = B] {$\triangledown$};
  \node (D) [node distance=0.5cm, below of=X] {$\W$};
  \node (E) [node distance=2cm, below of=D] {$\R$};
  \draw[->] (A) to node {$\T$}(B);
  \draw[->] (E) to node {$\U$}(D);

   \end{tikzpicture}
\caption{$(\T,X)$ authenticates $\R$. The objects $\xi, \U$ etc. are as in \ref{authentic lsp}.}
\label{fig:authentic_lses}
\end{figure}

\begin{definition}[Authentic hod-like $\sf{lses}$]\label{authentic lsp}\index{authentic hod-like lsp} (see Figure \ref{fig:authentic_lses}) Suppose $(\P, \Sigma)$ is a hod-like st-type pair, $\T$ is a normal iteration tree on $\P$ according to $\Sigma$ such that $\pi^{\T, b}$ exists and $X\subseteq \P^b$ is useful. Let $\S=\pi^{\T, b}(\P^b)$ and suppose $\R$ is a hod-like ${\sf{lses}}$. We say $(\T, X)$ \textbf{authenticates} $\R$ if there is a normal iteration tree $\U$ on $\R$ such that the following clauses hold.
\begin{enumerate}
\item $\U$ has a last model $\W$ such that $\pi^\U$ is defined and $\W\insegeq_{hod}\S$.
\item If $\gg<\lh(\U)$ is a limit ordinal such that $\S\models ``\d(\U\rest \gg)$ is a Woodin cardinal"\footnote{This condition then implies that for some window $w=(\eta, \d)$, $\S\models ``\d$ is a Woodin cardinal" and $\m(\U\rest \gg)=\S|\d$. See \rdef{layers of hod-like lsp}.}, letting $w$ be the unique window of $\S$ such that $\d(\U\rest \gg)=\d^w$ and setting $b=[0, \gg)_\U$, for some $\tau\in b$, 
\begin{center}
$s(\T, X, w)\subseteq \rge(\pi^{\U}_{\tau, b})$.
\end{center}
\item If $\R$ is of limit type then 
\begin{center}$\W^b=cHull^{\S^b}(\pi^{\T, b}[X]\cup \d^{\W^b})$\end{center} and if $\sigma: \W^b\rightarrow \S^b$ is the uncollapse map then 
\begin{center}
$\sigma^{-1}[\pi^{\T, b}[X]]\subseteq \rge(\pi^{\U, b})$. 
\end{center}
\end{enumerate}
We say $\R$ is $(\P, \Sigma, X)$-authentic if there is $\T$ on $\P$ according to $\Sigma$ such that $(\T, X)$  authenticates $\R$. We also say that $\R$ is $(\P, \Sigma, X, \T)$-authentic. 

Notice that there is only one iteration tree $\U$ with the above properties. We then say that $\U$ is the $(\T, X)$-authentication tree on $\R$. When $X=\P^b$ we simply omit it from terminology. $\myqedhere$
\end{definition}

Clearly the tree $\U$ in \rdef{authentic lsp} is a tree built via a comparison process in which $\S$ doesn't move. A typical $\R$ that we would like to authenticate will be an iterate of $\P$. If $\Sigma$ has nice properties, such as \textit{strong branch condensation} (see \rdef{strong branch condensation}) then clauses 2 and 3 of \rdef{authentic lsp} hold for the iterates of $\P$. Next, we would like to define \textit{authentic iterations}. 

\begin{definition}[Authentic iterations]\label{authentic iterations}\index{authentic hod-like lsp} Suppose $(\P, \Sigma)$ is a hod-like  st-type pair, $\T$ is a normal tree on $\P$ according to $\Sigma$ such that $\pi^{\T, b}$ exists and $X\subseteq \P^b$ is useful. Let $\S=\pi^{\T, b}(\P^b)$. Suppose $\R$ is an ${\sf{lses}}$ and $\mathcal{X}$ is a stack on $\R$. We say $(\T, X)$ \textbf{authenticates} $(\R, \mathcal{X})$ if $(\T, X)$ authenticates $\R$ and, letting $\U$ be the $(\T, X)$-authentication tree on $\R$ and $\W$ be the last model of $\U$, $\mathcal{X}$ is according to $\pi^\U$-pullback of $\Sigma_{\W, \T}$. 

Again we omit $X$ when $X=\P^b$. We say $(\R, \mathcal{X})$ is a \textbf{$(\P, \Sigma, X)$-authenticated iteration} if there is a tree $\T$ on $\P$ according to $\Sigma$ such that $(\T, X)$ authenticates $(\R, \mathcal{X})$. We also say that $(\R, \mathcal{X})$ is \textbf{$(\P, \Sigma, X, \T)$-authentic}. When $X=\P^b$ we simply omit it from our terminology. $\myqedhere$
\end{definition}

Next we define \textit{authentic indexable stacks}. These are stacks that will be important in our definition of short tree strategy mice (see \rdef{weak psi alpha indexing scheme a}). It maybe helpful to review the notation introduced in \rnot{notation for iteration trees}.


 \begin{definition}[Authentic indexable stacks]\label{authentic stacks of length 2}\index{authentic iterations} Suppose $(\P, \Sigma)$ is a  hod-like st-type pair, $X\subseteq \P^b$ is useful and $\R$ is a hod-like $\#$-lsa type ${\sf{lses}}$. Suppose 
 \begin{center}
 $t=(\R_0, \U, \R_1, \W)$
 \end{center}
  is an indexable stack on $\R=\R_0$. We say $t$ is \textbf{$(\P, \Sigma, X)$-authenticated} if the following conditions hold.
\begin{enumerate}
\item Suppose $\a\in R^\U$ is such that $\pi^{\U_{\leq \a}, b}$ exists. Then for all 
\begin{center}
$\a'\in (R^\U-(\a+1))\cup \{\lh(\U)\}$
\end{center}
 such that $\K=_{def}\U_{[\a, \a']}$ is based on $\S=_{def}\M_\a^b$,  $(\S, \K)$ is $(\P, \Sigma, X)$-authenticated iteration.
\item Suppose $\a\in R^\U$ is such that $\pi^{\U_{\leq \a}, b}$ exists. Then for all 
\begin{center}
$\a'\in (R^\U-(\a+1))\cup \{\lh(\U)\}$
\end{center}
 such that $\K=_{def}\U_{[\a, \a']}$ is above $\ord(\S)$ where $\S=_{def}\M_\a^b$,  the following conditions hold.
\begin{enumerate}
\item Suppose $\K$ doesn't have any fatal drops\footnote{See \rdef{fatal drop}.}. Then for any limit $\a<\lh(\K)$, if $b$ is the branch of $\K\rest \a$ then $\Q(b, \K\rest \a)$ exists and is $(\P, \Sigma, X)$-authentic.  
\item Suppose $\K$ has a fatal drop at $(\a, \eta)$. Let $\Q=\M_\tau^\K||\omega\xi_\tau^\K$. Then $(\Q, \K_{\geq\Q})$ is a $(\P, \Sigma, X)$-authenticated iteration.  
\end{enumerate}
\item $(\R_1^b, \W)$ is a $(\P, \Sigma, X)$-authenticated iteration. 
\end{enumerate}
When $X=\P^b$ we simply omit it from our terminology. $\myqedhere$
\end{definition}

It is of course desirable that $(\P, \Sigma, X)$-authenticated stacks are according to $\Sigma$. In the next section, we will use our authentication idea to define certified stacks.

\section{Short-tree-strategy mice}\label{sec short tree strategy mice}

We now have developed enough terminology and tools to define sts premice. We use the following notation below. Suppose $\M$ is a transitive model of some fragment of set theory and $\l$ is a limit of Woodin cardinals. Let $g\subseteq Coll(\omega, <\l)$ be $\M$-generic. For $\a<\l$, let $g_\a=g\cap Coll(\omega, <\a)$. We let $D(\M, \l, g)$ stand for the derived model of $\M$ at $\l$ computed using $g$. More precisely, letting $\bR^*=\bigcup_{\a<\l}\bR^{\M[g_\a]}$, $D(\M, \l, g)$ is defined in $\M(\bR^*)$ by letting
\begin{enumerate}
\item $\Gamma$ be the set of $A$ such that for some $\a<\l$ and some $B\in \M[g_\a]$ such that $\M[g_\a]\models ``B$ is $<\l$-universally Baire", $A$ is the interpretation of $B$ in $\M(\bR^*)$\footnote{The meaning of this is the following. For each $\M$-cardinal $\b\in (\a, \l)$, let $(T_\b, S_\b)\in \M[g\cap Coll(\omega, <\a)]$ be $\b$-absolutely complementing trees such that $p[T_\b]=B$. We then have that $A=\cup_{\b<\l}(p[T_\b])^{\M(\bR^*)}$. It is customary to set $\Gamma=Hom^*$. See \cite{DMT}}.
\item $D(\M,\l, g)=L(\Gamma, \bR^*)$. 
\end{enumerate}
Woodin's derived model theorem says that $D(\M, \l, g)\models \sf{AD}^+$ (see \cite{DMT}). We will use this theorem throughout this book.

Before we introduce the notion of short tree strategy premouse, we take a moment to describe the intuition behind the definition. Suppose $\P$ is a hod-like $\#$-lsa type ${\sf{lses}}$ and $\T$ is a normal ${\sf{nuvs}}$  tree on $\P$. We would like to find the correct $\Q$-structure for $\T$. We first attempt to find this $\Q$-structure among ${\sf{ses}}$ that have the trivial indexing scheme $\psi_0$, i.e., no indexable ${\sf{nuvs}}$ stack has an indexed branch. However, there may never be such an ${\sf{ses}}$ that can be used as $\Q$-structure. Assume then that this is the case. We then immediately encounter two problems. 

The first problem has to do with determining the exact stage of the constructibility order where we must stop looking for a $\Q$-structure among the ${\sf{ses}}$ that have the trivial indexing scheme. We will do this as soon as we reach a sufficiently closed stage. To know that we have reached such a level, we need to address the second problem.

The second problem is to describe the next type of gadgets that can be used as $\Q$-structures. A natural choice is the collection of ${\sf{ses}}$ over $\m(\T)$ in which all ${\sf{nuvs}}$  trees have $\Q$-structures with the trivial indexing scheme. This is our second indexing scheme. Let us call it $\psi_1$. One wrinkle is that we need a certification method for the $\Q$-structures that are used in a $\psi_1$-sts premouse. This is done by using the ideas from \rdef{authentic stacks of length 2}. 

The way we put the two ideas together is as follows. We first search for a $\Q$-structure among ${\sf{ses}}$ with the trivial indexing scheme $\psi_0$. If we reach a level $\M_0$ that has a $\psi_0$-sts $\Q_1\in \M_0$ that can be used as a $\Q$-structure then we stop and see if $\M_0$ certifies $\Q_1$ (see Definition \ref{weak psi alpha indexing scheme a}). If yes, then we declare success. If no, then we continue with the trivial indexing. This naturally leads to an induction, in which we define more and more complex indexing schemes which themselves are indexed by ordinals. One issue is that the most straightforward approach to the problem of defining the indexing schemes involves extending the language of ${\sf{ses}}$ to have names for ordinals, and this creates several unpleasant issues. Instead, we will first introduce ${\sf{ses}}$ whose indexing scheme may not be first order definable, the externally-$\phi$-${\sf{ses}}$. Afterwards, it will be straightforward to verify that being a short-tree-strategy premouse is in fact first order. 

Another issue is to show that if there is a $\Q$-structure for some tree $\T$  then we will indeed reach this $\Q$-structure inside our short-tree-strategy mice. For this, we will use an appropriate notion of fullness. Finally, the reader may find \rrem{how branches get indexed} useful. What follows is parallel to \rsec{layered sec}. The reader may want to review \rdef{j-structures def}, \rdef{amenable function}, \rnot{smsphi notation}, \rdef{important notation},  \rdef{passive hybrid j-structure}, \rdef{hybrid j-structure}, \rdef{layered hybrid e-structure} and \rdef{unindexed ses}.\\\\ 
\textbf{Externally-indexed ${\sf{hes}}$}\\

The main difference between \rdef{smsphi notation a} and \rdef{smsphi notation} is clause 4 bellow. 
\begin{notation}\label{smsphi notation a}
 Suppose that $\M=\mathcal{J}_{\iota}^{A_0, ..., A_n, f}(X)$ is a $\J$-structure or an f.s. $\J$-structure, $P\in X$ and $\Phi$ is a set of triples $(x, y, z)$ such that if $(x, y, z)\in \Phi$ then $x$ is a sequential structure. Let $S^{\M}_{P, \Phi}$ be the set of pairs $(\b, w)$ such that
 \begin{enumerate}
 \item $\omega\beta+\omega\gg^w\leq \ord(\M)$, 
 \item $\M|\omega\b\models ``\cf(\gg^w)$ is not a measurable cardinal as witnessed by extenders in $A_0$"\footnote{See \rrem{measurable cofinality issue}.}, and 
 \item $\M|\omega\b\models {\sf{ZFC}}$, and 
 \item $(w, \M|\omega\b, P)\in \Phi$. 
 \end{enumerate}
 $\myqedhere$
 \end{notation}

 \begin{definition}\label{important notation a}  Suppose that $(\M, P, \Phi)$ are as in \rnot{smsphi notation a}. 
 Suppose further that $f$ is a shifted amenable function with amenable component $g$ such that $\dom(f)\subseteq \univ{M}$ and for all $w\in \dom(f)$, $\min(f(w))+\gg^w\leq \ord(\M)$\footnote{Recall our convention that $X^\M$ is self-well-ordered.}.  We say $w$ is \textbf{weakly $(f, P, \Phi)$-\textit{minimal}} if there is $\b$ such that
 \begin{enumerate}
 \item $(\b, w)\in S^\M_{P, \Phi}$ (in particular, because $\M|\omega\b\models {\sf{ZFC}}$, $\omega\b=\b$),
 \item $w\not \in \dom(f\cap \univ{\M|\b})$, 
 \item  $\{ u\in \univ{\M|\b}: u<_{\M|\b}w$ and there is $\xi< \b$ such that $(\xi, u)\in S^\M_{P, \Phi}\}\subseteq \dom(f\cap \univ{\M|\b})$.
 \end{enumerate}
 We say $w$ is \textbf{$(f, P, \Phi)$-\textit{minimal}} if there is $\b$ witnessing that $w$ is weakly $(f, P, \Phi)$-\textit{minimal} and such that $w$ is the $<_{\M|\b}$-minimal $w'$ which is weakly $(f, P, \Phi)$-\textit{minimal} as witnessed by $\b$.

 

If $w$ is $(f, P, \Phi)$-\textit{minimal} then we let $\b^{\M, f, P, \Phi}_w$ be the least $\b$ witnessing that $w$ is $(f, P, \Phi)$-\textit{minimal}. In many cases, $(\M, f, P, \Phi)$ will be clear from context and so we will drop it from our notation. $\myqedhere$
\end{definition}

 We are now in a position to introduce the  \textit{externally-$\Phi$-indexed passive hybrid $\mathcal{J}$-structures}, or just $e\Phi$-indexed passive hybrid  $\mathcal{J}$-structures. 

\begin{definition}[$e\Phi$-indexed Passive Hybrid $\mathcal{J}$-structures]\label{passive hybrid j-structure a} We say $\M$ is an \textbf{$e\Phi$-indexed passive hybrid $\mathcal{J}$-structure} over a self-well-ordered set $X$ based on $P$ if $\M=(\M', k)$ is an f.s. $\mathcal{J}$-structure such that the following conditions hold.
\begin{enumerate}
\item  For some $\a$, $A \subseteq \univ{\M'}$ and $f\subseteq \univ{\M'}$, \begin{center}
$\M'=(\mathcal{J}_{\omega\a}^{A, f}(X), A, f, X, \in)$\footnote{We would like to emphasize that $\M'$ has only the displayed predicates. Also, below $(\M', f, \phi)$ are omitted from $\b_w$ notation.},
\end{center}
\item $f$ is a shift of an amenable function.
\item For all $w\in \univ{\M'}$, $w\in \dom(f)$ if and only if  $w$ is $(f, P, \Phi)$-minimal.
\item For all $w\in \dom(f)$, 
\begin{enumerate}
\item $\b_w=\min(f(w))$ and $\b_w+\omega\gg^w<\ord(\M)$\footnote{Here $\b_w$ is defined in \rdef{important notation a}.},
\item $\univ{\M'|(\b_w+\omega\gg^w)}=\J_{\b_w+\omega\gg^w}(\M'||\omega\b_w)$ and $A\cap \univ{\M'|(\b_w+\omega\gg^w)}=A\cap \univ{\M'|\omega\b_w}$\footnote{It also follows that $f\cap \univ{\M'|(\b_w+\gg^w)}=f\cap \univ{\M'|\b_w}$.}.
\end{enumerate}
\end{enumerate}
$\myqedhere$
\end{definition}

\begin{definition}[$e\Phi$-indexed Hybrid $\mathcal{J}$-structures]\label{hybrid j-structure a} We say $\M$ is an \textbf{$e\Phi$-indexed hybrid $\mathcal{J}$-structure} over a self-well-ordered set $X$ based on $P$ if $\M=(\M', k)$ is an f.s.  $\mathcal{J}$-structure such that 
\begin{enumerate}
\item  for some $\a$,  $A \subseteq \univ{\M'}$ and $f\subseteq \univ{\M'}$, \begin{center}
$\M'=(\mathcal{J}_{\omega\a}^{A, f}(X), A, f, B, F, X, \in)$\footnote{Below $(\M', f, \P, \Phi)$ are omitted from $\b_w$ notation.},
\end{center}
\item $(\mathcal{J}_{\omega\a}^{A, f}(X), A, f, X, \in)$ is an $e\Phi$-indexed passive hybrid $\J$-structure, 
\item at most one of $B$ and $F$ is not empty,
\item if $F\not=\emptyset$ then $F$ is an ordered pair $(w, b)$ such that if $\b=min(b)$ then setting $f'=f\cup\{(w, b)\}$,
\begin{enumerate}
\item $f'$ is a shift of an amenable function\footnote{This implies that $w$ is a sequential structure.},
\item $w$ is $(f', P, \Phi)$-minimal with $\b^{\M, f', P, \Phi}_w=\b$ (in particular, $\omega\b=\b$, see \rdef{important notation a}),
\item $\omega\a=\b+\omega\gg^w$,\footnote{It follows from clause 5 of \rdef{important notation a} that $\M'\models ``\cf(\gg)$ is not a measurable cardinal as witnessed by extenders in $A$".}
\item $\univ{\M'}=\J_{\b+\omega\gg^w}(\M'||\b)$ and $A\cap \univ{\M'}=A\cap \univ{\M'|\b}$.
\end{enumerate}
\end{enumerate}
For $w\in \dom(f')$, we say that $f'(w)$ is indexed at $\b_w+\omega\gg^w$ or that $\b_w+\omega\gg^w$ is the index of $f'(w)$. $\myqedhere$\\
\end{definition}


 \begin{definition}[$e\Phi$-indexed Strategic e-structure, $e\Phi-\sf{ses}$]\label{strategy premouse a} Suppose $\P$ is a transitive structure, $X$ is a self-well-ordered set such that $\P\in X$ and $\M=\mathcal{J}^{\vec{E}, f}(X)$ is an $e\Phi$-indexed hybrid $\mathcal{J}$-structure over $X$ based on $P$. We say $\M$ is an \textbf{$e\Phi$-indexed strategic e-structure} ($e\Phi$-$\sf{ses}$) over $X$ based on $\P$ if the following clauses hold.
 \begin{enumerate}
 \item $f^\M$ codes a partial iteration strategy for $\P$ such that for any $w\in dom(f^\M)$ if $\b=min(f^\M(w))$ then $\M|\b$ is closed\footnote{See \rdef{closed under sharps}. Also, recall that for such $\b$ we have $\omega\b=\b$}.
 \item $\vec{E}$ is a mixed indexed extender sequence. 
 \item If $\M=(\M', k)$\footnote{See \rdef{fine structural j-structure}.} then for every $(\omega\b, m)<l(\M)$,  $\M||(\omega\b, m)$ is sound.
 \end{enumerate}
 
 We say $\M$ is based on $\P$ if $\M$ is over $\mathcal{J}_\omega[\P]$ and is based on $\P$. $\myqedhere$\\
\end{definition}
\textbf{External $\Psi$-sts indexing scheme}\\

The reader may want to review \rdef{unindexed ses}.

\begin{definition}[$e\Psi$-sts indexing scheme]\label{sts indexing scheme a} Suppose $\Phi$ and $\Psi$ are two sets. We say $\Phi$ is an \textbf{external $\Psi$-sts indexing scheme} (or just an $e\Psi$-sts indexing scheme) if for all triples $(w, \N, \P)$, $(w, \N, \P)\in \Phi$ if and only if the following clauses hold.
\begin{enumerate}
\item $X$ is a self-well-ordered set, $\P\in X$ is a hod-like $\#$-lsa type ${\sf{lses}}$ and $w \in \N$.
\item $\N$ is an unindexed $\sf{ses}$ over $X$ based on $\P$. 
\item $\N$ is closed\footnote{See \rdef{closed under sharps}.}.
\item $\N\models ``\Sigma^\N$ is a partial faithful st-strategy with $m(\Sigma^\N)=\emptyset$"\footnote{This comment was made before as well, but we remind the reader. Notice that clause 4 below guarantees that $\Sigma^\N$ is really a partial strategy rather than st-strategy. We emphasize the fact that $\Sigma^\N$ is an st-strategy to point out the fact that there is no iteration according to $\Sigma^\N$ that is $\Sigma^\N$-maximal.},
\item $\N\models {\sf{ZFC}}+\phi^*(w)$\footnote{See \rdef{phi*}.}.
\item Either
\begin{enumerate}
\item $\N$ is ambiguous and $w$ is the $<_\N$-least sequential structure witnessing the ambiguity of $\N$.\\
Or
\item $\N$ is unambiguous and $w$ is the $<_\N$-least sequential structure $w'\in \N$ of the form $w'=(\mathcal{J}_{\omega}(t), t, \in)$ where $t=(\P, \T_0, \P_1, \T_1)$ such that $\N\models \phi^*[w']$, $t$ is an indexable ${\sf{nuvs}}$ such that  $\dot{\Sigma}(\T_0)$ is undefined and there is a cofinal well-founded branch $b$ of $\T_0$ such that $b\in \N$ and $(\T_0, \N, b)\in \Psi$.
\end{enumerate}
\end{enumerate}
$\myqedhere$
\end{definition}

\begin{definition}[Sts $\Psi$-premouse]\label{sts phi premouse} Suppose $X$ is a self-well-ordered set, $\P\in X$ is a hod-like $\#$-lsa type ${\sf{lses}}$ and $\Psi$ is a set. Let $\Phi$ be the $e\Psi$-sts indexing scheme. Then $\M$ is an \textbf{sts $\Psi$-premouse} over $X$ based on $\P$ if $\M$ is a $e\Phi$-${\sf{ses}}$ over $X$ based on $\P$ and if $w\in \dom(f^\M)$ is such that clause 6.b of \rdef{sts indexing scheme a} applies to $w=_{def}(\mathcal{J}_\omega(t), t, \in)$ where $t=_{def}(\P_0, \T_0, \P_1, \T_1)$ then letting $\b=\min(f^\M(w))$, 
\begin{center}
$f^\M(w)=\{\b+\omega\gg: \gg\in b\}$
\end{center}
 where $b\in \M|\b$ is the unique branch of $\T_0$ such that $(\T_0, \M|\b, b)\in \Psi$. $\myqedhere$
\end{definition}

Notice that in \rdef{sts indexing scheme a}, $\Phi$ is uniquely determined by $\Psi$. We now by induction define a sequence of sets $(\Psi_\b: \b\in Ord)$ and  for $\b\in Ord$, we let $\Phi_\b$ be the e$\Psi_\b$-sts indexing scheme. To start we let $\Psi_0=\emptyset$. Thus, if $\M$ is an $e\Phi_0-{\sf{ses}}$ then $\M$ does not have branches for ${\sf{nuvs}}$  stacks. We will use the following concept.

\begin{definition}[Terminal tree]\label{terminal tree} Suppose $X$ is a self-well-ordered set, $\P\in X$ is a hod-like $\#$-lsa type ${\sf{lses}}$, $\N$ is an $\sf{ses}$ over $X$ based on $\P$. Given $\T\in \N$ on $\P$, we say $\T$ is \textbf{$\N$-terminal} if $\T$ is ${\sf{nuvs}}$, $\T$ is according to $\Sigma^\N$ and $\T\not \in \dom(\Sigma^\N)$. $\myqedhere$
\end{definition} 

\begin{definition}[Sts indexing scheme]\label{weak psi alpha indexing scheme a}\index{sts indexing scheme} Let $\Psi_0=\emptyset$ and suppose $( \Psi_\xi: \xi<\alpha)$ have been defined. For $\xi<\a$ let $\Phi_\xi$ be the e$\Psi_\xi$-sts indexing scheme (see \rdef{sts indexing scheme a}). We let $\Psi_\a$ be the set of triples $(\T, \M, b)$ such that  $\M$ is an  ${\sf{ses}}$ over $X$ based on $\P$ and $(\T, b)$ is the $\M$-lexicographically least\footnote{This is just the order defined by: first order the first coordinate by $<_\M$, the canonical well-order of $\M$, then order the second coordinate by $<_\M$.} pair such that $\T$ is a normal iteration tree on $\P$, $\T$ is $\M$-terminal, and $b$ is a cofinal branch through $\T$ such that for some pair $(\gg, \xi)\in \ord(\M)\times \a$ the following clauses hold:
\begin{enumerate}
\item $\M|\gg$ is unambiguous\footnote{See \rdef{unambiguous hp}.} and $\M|\gg\models \sf{ZFC} + ``$there are infinitely many Woodin cardinals $>\d(\T)$".
\item $\M|\gg\models ``\lh(\T)$ is not of measurable cofinality". 
\item $b\in \M|\gg$ and $\M|\gg\models ``b$ is a well-founded branch".
\item $\M|\gg\models ``\Q(b, \T)$ exists" and $\Q(b, \T)$ is an $e\Phi_\xi-\sf{ses}$ over $\m^+(\T)$\footnote{This last statement about $\Q(b, \T)$ may not be first order over $\M|\gg$.}.
\item  Letting $(\d_i: i<\omega)$ be the first $\omega$ Woodin cardinals of $\M|\gg$ that are strictly greater than $\d(\T)$, the following holds in $\M|\gg$: $\Q(b, \T)$  is $<Ord$-iterable above $\d(\T)$ via a strategy $\Lambda$ such that letting $\l=\sup_{i<\omega}\d_i$, for every generic $g\subseteq Coll(\omega, <\l)$, $\Lambda$ has an extension $\Lambda^+ \in D(\M|\gg, \l, g)$ such that 
\begin{enumerate}
\item $D(\M, \l, g)\models ``\Lambda^+$ is an $\omega_1$-iteration strategy" and 
\item whenever $\R\in D(\M|\gg, \l, g)$ is a $\Lambda^+$-iterate of $\Q(b, \T)$ above $\delta(\T)$ and $t\in \R$ is an indexable stack on $\m^+(\T)$ according to $\Sigma^\R$, 
\begin{center}
$\M|\gg[g]\models ``t$ is $(\P, \Sigma^{\M|\gg})$-authenticated"\footnote{The witness for $t$ being $(\P,\Sigma^{\M|\gg})$-authenticated is in $\M|\gg$}.
\end{center}
\end{enumerate}
\end{enumerate}
The lexicographically least pair $(\gg, \xi)$ satisfying the above conditions is called the least $(\M, \Psi_\a)$-shortness witness for $(\T, b)$. We also say that $(\gg, \xi, b)$ is an $\M$-minimal shortness witness for $\T$. We also say that $\T$ has an $\M$-shortness witness. 

We say $\M$ is a \textbf{potential} $\sf{sts}$ premous if $\M$ is an $\sf{sts}$ $\Psi_{\a}$-premouse for some $\a$. $\myqedhere$
\end{definition}

Notice that because we minimized $b$, $\M$ has at most one $\M$-shortness witness for $\T$. The next lemma can now be established via an induction on ordinals. 
\begin{lemma}\label{correct calculation of sd 0} Suppose $M$ is a transitive model of ${\sf{ZFC}}$ and $\R\in M$. Then
\begin{enumerate}
\item For every $\a<\ord(\R)$, $M\models ``\R$ is an $\sf{sts}$ $\Psi_{\a}$-premouse" if and only if  $\R$ is an $\sf{sts}$ $\Psi_{\a}$-premouse.   
\item For every $\a<\ord(\M)$, $M\models ``\R$ is an $e\Phi_\a-\sf{ses}$" if and only if $\R$ is an $e\Phi_\a-\sf{ses}$. 
\item For every $\a<\ord(\R)$, \rdef{sts indexing scheme a} and \rdef{weak psi alpha indexing scheme a} define the sequences $(\Psi_\b\cap M:\b\leq \a)$ and  $(\Phi_\b\cap M: \b\leq \a)$ in $M$.
\end{enumerate}
\end{lemma}
\begin{proof} The claim is obvious for $\a=0$. In this case, $\Psi_0=\emptyset$, and since clause 6b of \rdef{sts indexing scheme a} is not applicable, the statement $``(w, \N, \P)\in \Phi_0"$ is a first order (over  $\N$) property of $\N$. Thus,  the three clauses above follow.

Suppose now that for some $\a<\ord(M)$, the three clauses have been verified for all $\b<\a$. We want to verify it for $\a$. 

We start with clause 3 of \rlem{correct calculation of sd 0}. Notice that all clauses of \rdef{weak psi alpha indexing scheme a} except the second half of clause 4 are internal properties of $\M$, where $\M$ is as in \rdef{weak psi alpha indexing scheme a}. But our induction hypothesis implies that for every $\xi<\a$, 
being $e\Phi_\xi-\sf{ses}$ is absolute between $M$ and $V$, implying that the second half of clause 4 of \rdef{weak psi alpha indexing scheme a} is absolute between $M$ and $V$ (notice that in \rdef{weak psi alpha indexing scheme a} the branch $b$ is in $\M$). This means that $\Psi_\a^M=\Psi_\a\cap M$. 

Next, fix $(w, \N, P)\in M$. Notice that if $(w, \N, P)\in \Phi_\a^M$ then 
\begin{enumerate}
\item if clause 6a of \rdef{sts indexing scheme a} applies then $(w, \N, P)\in \Phi_\a$ and
\item if clause 6b of \rdef{sts indexing scheme a} applies then $(\T_0, \N, b)\in \Psi_\a^M$, where $(\T_0, b)$ are as in clause 6b (and in particular, $b\in \N$).
\end{enumerate}
The first statement above holds as all clauses of \rdef{sts indexing scheme a} except 6b are internal properties of $\N$ and as such are absolute between $M$ and $V$. Notice next that because we already have that $\Psi_\a^M=\Psi_\a\cap M$, the second statement above implies that $(\T_0, \N, b)\in \Psi_\a$ and hence, $(w, \N, P)\in \Phi_\a$. Thus, $\Phi_\a^M=\Phi_\a\cap M$.

The proof of the remaining clauses are very similar, and can be easily established by examining \rdef{sts phi premouse}.
\end{proof}
\begin{corollary}\label{correct calculation of sd 1} Suppose $M$ is a transitive model of ${\sf{ZFC}}$ and $\a<\ord(M)$. Then $M\models ``\M$ is an $e\Phi_\a-\sf{ses}$" if and only if $\M$ is an $e\Phi_\a-\sf{ses}$. Also, $M\models ``\M$ is an $\sf{sts}$ $\Psi_{\a}$-premouse" if and only if  $\M$ is an $\sf{sts}$ $\Psi_{\a}$-premouse,  $\a<\ord(M)$ and $\M\in M$. 
\end{corollary}

\begin{definition}\label{degree of shortness} Suppose $\M$ is a potential $\sf{sts}$ premouse. We say $\a$ is the \textbf{shortness degree} of $\M$ if $\a$ is the least for which $\M$ is a $e\Phi_{\a}-\sf{ses}$. We let $sd(\M)$ be the shortness degree of $\M$. $\myqedhere$
\end{definition}
\textbf{The shortness degree of a potential sts premouse}\\

Suppose $\M$ is a potential $\sf{sts}$. We now describe a well-founded tree $U(\M)$ whose rank bounds $sd(\M)$. The nodes in $U(\M)$ consist of finite sequences of the form $(x_0, x_1, ..., x_n)$ such that the following conditions hold:
\begin{enumerate}
\item For each $i\leq n$, $x_i=(t_i, b_i, \M_i)$ where $t_i=(\M_i, \T_i)$ is an indexable stack on $\M_i$, $\T_i$ is a normal tree on $\M_i$, $b_i$ is a branch of $\T_i$ and for $i+1\leq n$, $\M_{i+1}=\Q(b_i, \T_i)$.
\item  $\M_0=\M$.
\item For each $i\leq n$, $t_i\in \dom(\Sigma^{\M_i})$\footnote{In particular, $t_i\in \M_i$.} and $b_i=\Sigma^{\M_i}(\T_i)$.
\item For each $i\leq n$, $\mathcal{J}_{\omega}[\m^+(\T_i)]\models ``\d(\T_i)$ is a Woodin cardinal".
\end{enumerate}
Notice that if $(x_0, ..., x_n)\in U(\M)$ then for each $i< n$, $\M_{i+1}\in \M_i$\footnote{Notice that because $\d(\T_i)$ is a Woodin cardinal, $t_i$ is an $\sf{nuvs}$ and therefore, $b_i, \M_{i+1}\in \M_i$. See clause 3 of \rdef{weak psi alpha indexing scheme a}.}. Therefore, $U(\M)$ is well-founded. If $p\in U(\M)$ then we use superscript $p$ to denote the objects that appear in $p$. For example $\M_n^p$, $x_i^p$ or $\T_i^p$. Given a well-founded tree $S$ we let $rank(S)$ be its rank.

\begin{lemma}\label{bounding sd} $sd(\M)=rank(U(\M))$. 
\end{lemma}
\begin{proof} The proof is by induction on $sd(\M)$. Suppose $sd(\M)=0$. In this case, $U(\M)=\emptyset$ and hence, its rank is also $0$. On the other hand, if $U(\M)=\emptyset$ then clearly $\M$'s strategy predicate does not index any branch for an $\sf{nuvs}$ indexable stack, and therefore, $sd(\M)=0$. For $\b$ an ordinal, let $I(\b)$ be the conjunction of the following two statements.\\\\
(1) For all $\M'$, if $sd(\M')=\b$ then $rank(U(\M'))=\b$.\\
(2) For all $\M'$, if $rank(U(\M'))=\b$ then $sd(\M')=\b$.\\\\
We want to prove that for all $\b$, $I(\b)$ is true. Assume then that for some $\a\geq 1$, for all $\b<\a$, $I(\b)$ holds. We want to prove $I(\a)$, which amounts to proving that the following two statements hold.\\\\
(A) For any $\M$ such that $sd(\M)=\a$, $rank(U(\M))=\a$.\\
(B) For any $\M$ such that $rank(U(\M))=\a$, $sd(\M)=\a$.\\\\
We now prove (A). Fix $\M$ such that $sd(\M)=\a$. We want to see that $rank(U(\M))=\a$. Suppose first that $p=(x_0,..., x_n)\in U(\M)$. Then we have that\\\\
(*) $U(\M_{n+1}^p)=\{ q: p^\frown q\in U(\M)\}$.\\\\
We now show that the rank of $U(\M)$ is at least as big as $\a$. To see this, it is enough to show that for each $\b<sd(\M)$ there is a node $p=(x_0,..., x_n)\in U(\M)$ such that letting $\M_{n+1}=\Q(b^p_n, \T^p_n)$, $sd(\M_{n+1})\geq \b$. (*) then will imply that in fact $\b\leq rank(U(\M))$. But because $\b<sd(\M)$, we must have a pair $(\P, \T)\in \M$ such that $\T\in \dom(\Sigma^\M)$ and if $b=\Sigma^\M(\T)$ and $t=(\P, \T)$ then $(t, b, \M)\in U(\M)$ and $sd(\Q(b, \T))\geq \b$. It then follows that $p=(t, b, \M)$ is as dessired.

We now show that $rank(U(\M))\leq \a$. Indeed, let $p=(x_0)\in U(\M)$ and set $\M_1=\Q(b_0^p, \T_0^p)$. Because $sd(\M)=\a$, it follows that $sd(\M_1)<\a$. Therefore, it follows from (*) and $\forall \b<\a I(\b)$ that $rank(U(\M_1))<\a$. As this is true for any node of $U(\M)$ of length 1, we have that $rank(U(\M))\leq \a$.

The proof of (B) is very similar. Indeed, if $\M$ is such that $rank(U(\M))=\a$ then (A) implies that $sd(\M)\geq \a$. Suppose then $sd(\M)>\a$. We then claim that there is $p\in U(\M)$ such that $p$ has length $n+1$ and $sd(\M_{n+1}^p)=\a$. Suppose otherwise. Thus the following is true:\\\\
(**) whenever $p\in U(\M)$ is of length $n+1$, either $sd(\M_{n+1}^p)>\a$ or $sd(\M_{n+1}^p)<\a$.\\\\ 
We now inductively define $(p_i: i<\omega)$ such that for all $i<\omega$, 
\begin{enumerate}
\item $p_i\in U(\M)$,
\item $p_{i+1}$ extends $p_i$,
\item $p_i$ has length $i+1$,
\item $sd(\M_{i+1}^{p_i})>\a$.
\end{enumerate}
Let $p_0\in U(\M)$ be of length 1 and such that $sd(\M_1^{p_0})\geq \a$. There is indeed such a $p_0$ as  all $\Q$-structures used in $\M$ would have shortness degree $<\a$ implying that $sd(\M)\leq \a$. (**) now implies that in fact $sd(\M_1^{p_0})>\a$. Repeating this construction $\omega$ times produces our desired sequence. It now follows that $U(\M)$ is not well-founded, which is a contradiction and hence, (B) holds.
\end{proof}

\begin{lemma}\label{correct calculation of sd 2}  Suppose $M$ is a transitive model of ${\sf{ZFC}}$, $\M\in M$ and $\M$ is a potential $\sf{sts}$ premouse. Then $M\models ``\M$ is potential $\sf{sts}$ premouse".
\end{lemma}
\begin{proof} It follows from \rlem{correct calculation of sd 1} that it is enough to establish that $sd(\M)\in M$. But this follows from the fact that $U(\M)\in M$ (and hence, $rank(U(\M))\in M$) and $sd(\M)=rank(U(\M))$. 
\end{proof}

\textbf{Authenticated potential sts premouse}

\begin{definition}\label{authenticated potential sts} Suppose $\M$ is an unindexed $\sf{ses}$ over $X$ based on $\P$ and $\Q\in \M$ is a potential $\sf{sts}$ premouse over some $\#$-lsa type hod-like $\sf{lses}$ $\S$.  We say $\Q$ is $\M$-\textbf{authenticated} if the following clauses hold:
\begin{enumerate}
\item $\M$ has at least $\omega$ many Woodin cardinals $>\ord(\Q)$. 
\item  Letting $(\d_i: i<\omega)$ be the first $\omega$ Woodin cardinals of $\M$ that are strictly greater than $\ord(\Q)$, the following holds in $\M$: $\Q$  is $<Ord$-iterable above $\d^\S$ via a strategy $\Lambda$ such that letting $\l=\sup_{i<\omega}\d_i$, for every generic $g\subseteq Coll(\omega, <\l)$, $\Lambda$ has an extension $\Lambda^+ \in D(\M, \l, g)$ such that 
\begin{enumerate}
\item $D(\M, \l, g)\models ``\Lambda^+$ is an $\omega_1$-iteration strategy" and 
\item whenever $\R\in D(\M, \l, g)$ is a $\Lambda^+$-iterate of $\Q$ above $\d^\S$ and $t\in \R$ is an indexable stack on $\m^+(\T)$ according to $\Sigma^\S$, 
\begin{center}
$\M[g]\models ``t$ is $(\P, \Sigma^{\M})$-authenticated"\footnote{The witness for $t$ being $(\P,\Sigma^{\M})$-authenticated is in $\M$}.
\end{center}
\end{enumerate}
\end{enumerate}
Being $\M$-authenticated is a first order property of $\M$, and so we write $\M\models ``\Q$ is authenticated" for the statement $\Q$ is $\M$-authenticated. $\myqedhere$
\end{definition}

\textbf{The definition of short tree strategy premouse}\\

Let $U_0(x, y)$ be the formula in the language of unindexed ${\sf{ses}}$ expressing the statement that ``$x$ is an ordinal and $y$ is the universe up to $\omega x$". Thus, $U_0(x, y)$ defines the function $\gg\mapsto \M|\omega\gg$ over any ${\sf{ses}}$ $\M$. 

\begin{definition}\label{sts0} We let ${\sf{sts_0}}(x, y)$ be the formula in the language of ${\sf{ses}}$ expressing the following: there is an ordinal $\gg$ such that letting $M$ be such that $U_0(\gg, M)$, the following clauses hold:
\begin{enumerate}
\item   $M\models {\sf{ZFC}}$.
\item $x$ is a normal iteration tree of limit length and $y$ is a cofinal well-founded branch of $x$.
\item  $M\models ``\Q(y, x)$ exists and is an authenticated potential $\sf{sts}$ premouse".
\item For any well founded branch $y'$ of $x$ and an ordinal $\gg'$, letting  $M'$ be such that $U(\gg', M')$, if
\begin{enumerate}
\item $y\not =y'$,
\item $M'\models {\sf{ZFC}}$, and 
\item  $M'\models ``\Q(b', \T)$ exists and is an authenticated potential $\sf{sts}$ premouse",
\end{enumerate}
then letting $\Q=\Q(y, x)$ and $\Q'=\Q(y', x)$, either $(\gg, sd(\Q))<_{lex} (\gg', sd(\Q'))$ or $(\gg, sd(\Q))=(\gg', sd(\Q'))$ and $b<_{M}b'$. 
\end{enumerate} 
$\myqedhere$
\end{definition}

 \begin{definition}[Sts-indexed $\sf{ses}$, Sts mouse]\label{sts premouse}\index{sts mouse} Suppose $X$ is a self-well-ordered set and $\P\in X$ is a hod-like $\#$-lsa type ${\sf{lses}}$. Let $\sf{sts}$ be the $\sf{sts_0}$-sts indexing scheme\footnote{See \rdef{sts indexing scheme}.}. We say $\M$ is an \textbf{sts premouse} over $X$ based on $\P$ if $\M$ is an $\sf{sts}$-indexed $\sf{ses}$  over $X$ based on $\P$. If additionally $\M$ is $\omega_1+1$-iterable the we say that $\M$ is an  \textbf{sts mouse}.\footnote{Here implicit in this is the demand that iterates of $\P$ according to the strategy are sts premice.} $\myqedhere$
\end{definition}

The following is an easy lemma.

\begin{lemma}\label{abs of sts} Suppose $M$ is a transitive set and $\M\in M$. Then $\M$ is an $\sf{sts}$ premouse if and only if $M\models ``\M$ is an $\sf{sts}$ premosue".
\end{lemma}

The following is a corollary to \rlem{qstructures are sts}. It implies that the certified $\Q$-structures themselves are $\sf{sts}$ premice.

\begin{lemma}\label{qstructures are sts 1} Suppose $\P$ is a uniformly\footnote{See \rdef{uniformly sts}.} $\sf{sts}$-organized $\#$-lsa type hod like ${\sf{lses}}$. Suppose $\T$ is a normal $\sf{nuvs}$ tree on $\P$  and $b$ is well-founded branch of $\T$ such that $\Q(b, \T)$ exists. Then $\Q(b, \T)$ is an $\sf{sts}$ premouse based on $m^+(\T)$.
\end{lemma}

\begin{remark}[On how branches get indexed]\label{how branches get indexed} The first key point is that $\M|\gg$ in \rdef{weak psi alpha indexing scheme a} is not the analogue of $\M|\b$ in \rdef{important notation}. The analogue of $\M|\b$ in the sense of \rdef{important notation} is $\M$ itself. Recall that the indexing scheme is not $\Psi_\a$ but rather $\Phi_\a$, and so the relevant definitions for determining the analogue of $\M|\b$ in the sense of \rdef{important notation} are \rdef{sts0} and \rdef{sts premouse}. 

\rdef{hybrid j-structure} introduced layered hybrid $\mathcal{J}$-structures, and a key aspect of that definition is the indexing of branches. The indexing scheme $\phi$ (in the sense of \rdef{hybrid j-structure}) is only picking the iteration trees that we would like to index, where the branches are indexed is then uniquely determined by the procedure described in \rdef{hybrid j-structure}. \rdef{sts0} and \rdef{sts premouse} are relevant definitions, and explain what the $\phi$ in \rdef{hybrid j-structure} should be. 

The reader may wonder why we have concentrated so much on ${\sf{nuvs}}$   iterations. The point is that clause 4b of \rdef{sts indexing scheme} requires that we add the branches of ${\sf{uvs}}$   iterations, and these branches are not branches that we intend to certify. These branches are told to the model by consulting an outside strategy. It is only the branches of ${\sf{nuvs}}$   iterations, the ones that appear in clause 4b of \rdef{sts indexing scheme}, are being certified. The schemes introduced in \rdef{weak psi alpha indexing scheme a} determine our certification procedures. $\myqedhere$
\end{remark}

\begin{definition}[$\Lambda$-sts premouse]\label{lambda sts premouse}\index{$\Lambda$-sts premouse} Suppose $X$ is a self-well-ordered set, $\P\in X$ is a hod-like $\#$-lsa type ${\sf{lses}}$, $\Lambda$ is an st-strategy for $\P$ and $\M$ is an sts premouse over $X$ based on $\P$. Then we say $\M$ is a \textbf{$\Lambda$-sts premouse} over $X$ based on $\P$ if $\Sigma^M\subseteq \Lambda\rest \M$. $\myqedhere$
\end{definition}

\begin{definition}[$\Lambda$-sts mouse]\label{lambda sts mouse}\index{$\Lambda$-sts mouse} Suppose $X$ is a self-well-ordered set, $\P\in X$ is a hod-like $\#$-lsa type ${\sf{lses}}$, $\Lambda$ is an st-strategy for $\P$ and $\M$ is a $\Lambda$-sts premouse over $X$ based on $\P$.  Then we say $\M$ is a \textbf{$\Lambda$-sts mouse} over $X$ based on $\P$ if $\M$ has an $\omega_1+1$-iteration strategy $\Sigma$ such that whenever $\N$ is a $\Sigma$-iterate of $\M$ via $\Sigma$, $\N$ is a $\Lambda$-sts premouse over $X$ based on $\P$. 

We say $\M$ is a $\Lambda$-sts (pre)mouse over $\P$ if $\M$ is a $\Lambda$-sts (pre)mouse over $\mathcal{J}_{\omega}[\P]$ based on $\P$. $\myqedhere$
\end{definition}

\section{The hod premouse indexing scheme}\label{hp indexing scheme:sec}

The goal of this short section is to introduce \textit{the hod premouse indexing scheme} ($\sf{hp}$ indexing scheme). This scheme combines the standard indexing scheme with the $\sf{sts}$ indexing scheme The standard indexing scheme which is used in \cite{ATHM} is due to Woodin. According to this scheme we must pick the least iteration whose branch has not yet been indexed in the strategy predicate and index the branch of this iteration in the strategy predicate. Below we give a formal definition of the $\sf{hp}$ indexing scheme.


The reader may find it helpful to review \rdef{important notation} and \rdef{language of lses}. In particular, the reader should keep in mind that the intended universes where indexing schemes are evaluated are the models of the form $\M|\omega\b$ of \rnot{smsphi notation}. Thus, these universes themselves are not hod like $\sf{lses}$ (see \rdef{hod like lsp}). But each such $\M|\omega\b$ has a its own predicate $Y^{\M|\omega\b}$ which is what we will use below to describe the $\sf{hp}$-indexing scheme. Perhaps reviewing \rrem{the y predicate} may clarify some of the questions that the reader might have. 

\begin{definition}\label{index-ready} We say that $\sf{lses}$ $\M$ is \textbf{strategy-ready} if letting $\iota=\ord(Y^\M)$, $\omega\iota+\omega^2<\ord(\M)$.  $\myqedhere$
\end{definition} 


\begin{definition}[Hod premouse indexing scheme, $\sf{hp}$ indexing scheme]\label{sis}\index{standard indexing scheme}\index{hp-indexed} 
We say $\phi(x, y)$ is the \textbf{hod premouse indexing scheme} ($\sf{hp}$ indexing scheme) if $\phi$ is the conjunction of the following clauses.
\begin{enumerate}
\item The universe is closed (see \rdef{closed under sharps}). 
\item The universe is strategy-ready (see \rdef{index-ready}).
\item $x=\cup Y^{\dot{\V}}$.
\item If $x$ is lsa like then 
\begin{enumerate}
\item $\dot{\V}$ is an ${\sf{sts}}$ premouse over $\dot{\V}|\iota+\omega$ based on $x$ (see \rdef{sts premouse}),
\item ${\sf{sts}}[y]$.
\end{enumerate} 
\item If $x$ is not lsa like then $y$ is the $<_{\dot{\V}}$-least sequential structure of the form $(\mathcal{J}_\omega(\T^y), \T^y, \in)$ where $\T^y$ is a stack on $x$ that is according to $\dot{\Sigma}_x$ and doesn't have a last model. 
\end{enumerate}
We let $\sf{hp}$ denote the $\sf{hp}$-indexing scheme. $\myqedhere$
\end{definition}

\begin{remark} The determination of the $Y$ predicate of the models appearing in the hod pair constructions (see \rdef{gamma-hod pair construction*}) is an important step in such constructions. $\myqedhere$
\end{remark}

The next definition isolates the standard indexing scheme. It is defined in the language of $\sf{ses}$ which has a constant symbol for a structure whose strategy is indexed on the sequence of $\sf{ses}$. We let $\dot{\P}$ be this constant. 

\begin{definition}\label{standard indexing scheme} We say $\phi(y)$ is the \textbf{standard indexing scheme} ($\sf{sis}$-indexing scheme) if $\phi$ expresses the following statement: $y$ is then $<_{\dot{\V}}$-least sequential structure of the form $(\mathcal{J}_\omega(\T^y), \T^y, \in)$ where $\T^y$ is a stack on $\dot{\P}$ that is according to $\dot{\Sigma}$ and $\lh(\T^y)$ is a limit ordinal. We let $\sf{sis}$ denote the $\sf{sis}$ indexing scheme. $\myqedhere$
\end{definition}

%

\begin{remark}
Woodin's method of feeding the branch information into the model (as described in clause 4 of \rdef{sis}) is easy to comprehend and allows us to develop the basic theory of hod mice in this manuscript; however, it does not seem to allow for the proof of $\square$ to generalize easily. An alternative method to feeding in branch information that does allow for the $\square$ proof to generalize is to use the $\mathfrak{B}$-operator (see \cite{trang2013}). This method is summarized in Section \ref{sec:hod_mice}; we also describe where the $\square$ proof seems to break down if Woodin's method was used. Nevertheless, a hod mouse constructed using Woodin's method constructs the same sets as the one using the $\mathfrak{B}$-operator (given that everything else is the same). Woodin's method is used from now on to the end of Chapter 10 because of its simplicity. $\myqedhere$
\end{remark}

\section{Hod mice}\label{hod mice sec}

The main goal of this section is to introduce \textit{lsa small hod premice}. The reader might find it helpful to review \rsec{sec: hod-like layered hybrid premice}. In particular, we will use \rdef{the o stack}, \rdef{pre-hod-like}, \rdef{proper type II}, \rdef{lsa type}, \rdef{layers of hod-like lsp}, \rdef{hod like lsp}, \rnot{l p}, and \rter{types of lsa small premice}. Also recall our convention introduced in \rrem{lsa small convention}. According to this convention all our hod-like ${\sf{lses}}$ are lsa small. 

We start by isolating the types of points in $Y^\P$ where $\P$ is hod-like ${\sf{lses}}$. 

\begin{notation}[Meek and lsa points]\label{points} Suppose $\P$ is a hod-like ${\sf{lses}}$. 
\begin{enumerate}
\item $meek(\P)=\{\Q\in Y^\P: \Q$ is meek\footnote{See \rdef{pre-hod-like}.}$\}$. 
\item $lsa(\P)=\{ \Q\in Y^\P : \Q$ is of $\#$-lsa type\footnote{See \rdef{lsa type}.}$\}$. 
\item ${\sf{\ml}}(\P)=\bigcup Y^\P$\footnote{This object was introduced in \rdef{l p}.}.
\end{enumerate}
$\myqedhere$
\end{notation}

%

\rdef{pre-hod-like} and \rdef{hod like lsp} do most of the job that we need to do to define hod premice. Essentially what is missing from \rdef{hod like lsp} is the exact nature of premice at lsa layers. In the next definition, we will not repeat what has already been introduced in \rdef{pre-hod-like} and \rdef{hod like lsp}.

\begin{definition}[Hod premouse]\label{hod premouse}\index{hod premouse} Suppose $\P$ is a (lsa small) hod-like ${\sf{lses}}$\footnote{See \rdef{pre-hod-like} and \rdef{hod like lsp}.}. Let $(\P_{\xi, \xi'}: \xi\leq \eta \wedge \xi'\leq \nu_\xi)$ be the sequence of layers of $\P$ and $(\d_\xi, \iota_{\xi, \xi'}: \xi\leq \eta \wedge \xi'\leq \nu_\xi)$ be the sequence of ordinal parameters associated with it (see \rdef{layers of hod-like lsp}).
 We say $\P$ is an \textbf{lsa small hod premouse} or just a \textbf{hod premouse} if $\P$ is $\sf{hp}$-indexed\footnote{See \rdef{sis}.} hod-like $\sf{lses}$ that has the following properties: 
 \begin{enumerate}
 \item Suppose $\nu$ is a cutpoint of $\P$. Then the following holds.
 \begin{enumerate}
 \item If $\P$ is meek and $\nu<\d^\P$ then $\P\models ``\mathcal{O}^\P_{\nu, \nu}$ has an $Ord$-strategy (sts strategy respectively) acting on iteration trees that are above\footnote{See \rter{above eta}.} $\nu"$.
 \item If $\P$ is non-meek and $\nu<\d^\P$ then $\P|\d^\P\models ``\mathcal{O}^\P_{\nu, \nu}$ has a $\d^\P$-strategy acting on trees that are above $\nu$".

 \end{enumerate}
 \item If $\P$ is  of successor type\footnote{See \rter{types of lsa small premice}}, $\xi+1=\eta$ and $\Q=\P_{\xi, \iota_\xi}$ then for any $\eta\in (\d^\Q, \d^\P)$, $\P\models ``\P|\eta^+$ is $(Ord, Ord)$-iterable for stacks that are above $\ord(\Q)$".
 \item If $\P$ is of lsa type and $\eta \in (\ord(\P^b), \d^\P)$ then $\P|\d^\P\models ``\P|\eta^+$ is $(Ord, Ord)$-iterable for stacks that are above $\ord(\P^b)$"
 
%


 \end{enumerate}
 $\myqedhere$
\end{definition}

Next we define hod pairs.

\begin{definition}[Hod pairs]\label{hod pairs}\index{hod pairs} We say $(\P, \Sigma)$ is a (simple) hod pair if $(\P, \Sigma)$ is a (simple) hod-like $\sf{lses}$ pair\footnote{See \rdef{hod-like lsp pair}.}, $\P$ is a hod premouse and $\Sigma$ has hull condensation. $\myqedhere$
%
\end{definition}

Next we introduce the collection of sets generated by hod pairs.

\begin{definition}[$\Gamma(\P, \Sigma)$ and $B(\P, \Sigma)$]\label{gamma(p, sigma) and b(p, sigma)}\index{$\Gamma(\P, \Sigma)$}\index{$B(\P, \Sigma)$} Suppose $(\P, \Sigma)$ is a hod pair of limit type. We then let 
\begin{center}
$B(\P, \Sigma)=\{(\T, \Q) : \exists \R ((\T, \R)\in I(\P, \Sigma)\wedge \Q\unlhd_{hod}\R^b)\}$, and\\
$\Gamma(\P, \Sigma)=\{ A\subseteq \mathbb{R} : \exists (\T, \Q) \in B(\P, \Sigma) (A\leq_w {\sf{Code}}(\Sigma_{\Q, \T})\}$.
\end{center}
$\myqedhere$
\end{definition}

\begin{definition}[Pre-sts hod pairs]\label{pre sts hod pairs} We say $(\P, \Sigma)$ is a \textbf{pre-sts hod pair} if $(\P, \Sigma)$ is a hod-like st-type pair\footnote{See \rdef{lsa type pair}.} and $\Sigma$ is a $(\k, \l, \nu)$-st-strategy for $\P$ with hull condensation.
%

We say $(\P, \Sigma)$ is a \textbf{simple pre-sts hod pair} if $(\P, \Sigma)$ is a hod-like st-type pair and $\Sigma$ is a $(\l, \nu)$-st-strategy for $\P$ with hull condensation. $\myqedhere$
\end{definition}

To define sts hod pairs, we will make use of the notation introduced in \rdef{gamma(p, sigma) and b(p, sigma) for sts}.
Recall that in \rdef{gamma(p, sigma) and b(p, sigma) for sts}, we introduced $\Gamma^b(\P, \Sigma)$ but not $\Gamma(\P, \Sigma)$. We will define $\Gamma(\P, \Sigma)$ for sts hod pairs in \rsec{gamma(p, sigma) in the successor case sec}. 

Suppose now that $X$ is a self-well-ordered set, $(\P, \Sigma)$ is a pre-sts pair such that $\P\in X$ and $\Q$ is a $\Sigma$-sts mouse over $X$ based on $\P$. Let $\Lambda$ be the strategy of $\Q$. We then let $\Gamma(\Q, \Lambda)$ be the collection of all sets of reals $A$ such that for some $\Lambda$-iterate $\R$ of $\Q$, there is $(\T, \S)\in B(\P, \Sigma^\R)$ such that $A\leq_w \Sigma_{\S, \T}$.

\begin{definition}[Sts hod pairs]\label{sts hod pairs}\index{sts hod pairs}  We say $(\P, \Sigma)$ is an \textbf{sts hod pair} if $(\P, \Sigma)$ is a pre-sts pair such that whenever $(\T, \R, \tau)$ is such that letting $(\R_{\xi, \xi'}: \xi\leq \eta \wedge \xi'\leq \nu_\xi)$ be the sequence of layers of $\R$ and $(\d_\xi, \iota_{\xi, \xi'}: \xi\leq \eta \wedge \xi'\leq \nu_\xi)$ be the sequence of ordinal parameters associated with it (see \rdef{layers of hod-like lsp}),
\begin{enumerate}
\item $(\T, \R)\in I(\P, \Sigma)$\footnote{See \rdef{short tree iterates}.} and
\item  $\R_{\tau, 0}\in lsa(\R)$ and $\d_\tau<\d^\R$,
\end{enumerate}
 then $\R_{\tau, 1}$ has an iteration strategy $\Phi\in \Gamma^b(\P, \Sigma)$ witnessing that $\R_{\tau, 1}$ is a $\Sigma_{\R_{\tau, 0}, \T}$-sts mouse based on $\R_{\tau, 0}$\footnote{Thus, all the iterates of $\R_{\tau, 1}$ via $\Phi$ are above $\ord(\R_{\tau, 0})=\iota_{\tau, 0}$.} and such that $\Gamma(\R_{\tau, 1}, \Phi)\subset \Gamma^b(\P, \Sigma)$.
 
 Similarly we can define \textbf{simple} sts hod pairs. $\myqedhere$
\end{definition}

\begin{definition}\label{allowable pair} We say $(\P, \Sigma)$ is an \textit{allowable pair} if it is one of the hod pairs introduced above. More precisely, one of the following holds:
\begin{enumerate}
\item $(\P, \Sigma)$ is a hod pair. 
\item $(\P, \Sigma)$ is a simple hod pair.
\item $(\P, \Sigma)$ is an sts hod pair.
\item $(\P, \Sigma)$ is a simple sts hod pair.
\end{enumerate}
In the context of $\sf{AD}^+$, unless otherwise specified, the strategy component of any of the above pairs will always be $(\omega_1, \omega_1, \omega_1)$ or $(\omega_1, \omega_1)$ strategy or st-strategies. $\myqedhere$
\end{definition}

\rdef{sts hod pairs} imposes conditions on sts hod pairs that may seem unnatural. However, these conditions are needed to prove that sts hod pairs behave nicely. These clauses will be used in Chapter \ref{lsa internal theory chapter}.

\chapter{A comparison theory of hod mice}\label{chap:comparison}

This section is devoted to proving a comparison theorem for hod pairs. We will have two comparison theorems, \rcor{comparison holds 1} and \rcor{diamond comparison}. \rcor{comparison holds 1} is useful in determinacy context while \rcor{diamond comparison}  is useful in Core Model Induction applications. The following is a key hypothesis used in many of the theorems of this chapter.

\begin{definition}\label{no mouse with a superstrong} We let $\sf{NsesS}$ stand for the statement ``there is no $\omega_1$-iterable $\sf{ses}$ with a superstrong cardinal".  $\myqedhere$
\end{definition} 

\section{Backgrounds and Suslin capturing}

The goal of this section is to introduce \textit{backgrounds} and the concept of \textit{Suslin, co-Suslin capturing}. We will use these notions to build hod pairs with desired properties, such as \textit{fullness preservation} and \textit{branch condensation}. Before we do this, we fix a coding of hereditarily countable sets by reals. We will use this coding throughout this book.
\begin{definition}\label{coding of countable objects}
Given a real $x\in \bR$, we let $E_x=\{ (m, n): x(2^m3^n)=1\}$. We let $x\in \sf{Code}$ if $m_x=_{def}(\omega, E_x)$ is a well-founded model satisfying the Axiom of Extensionality. If $x\in \sf{Code}$ then we let $\pi_x: m_x\rightarrow  M_x$ be the transitive collapse of $m_x$ and let $c_x=\pi_x(0)=\{ \pi_x(m): x(2^m)=1\}$. We then say that $x$ \textbf{codes} $c_x$.   $\myqedhere$
\end{definition}

Recall that $\sf{HC}$ is the set of hereditarily countable sets (see \rdef{hc set up}). Given $n\in \omega$ and  $A\subseteq {\sf{HC}}^n$ we let ${\sf{Code}}(A)=\{ x\in \sfc : c_x\in A\}$. Notice that \begin{center}$\sfc: \cup_{n\in \omega}\powerset({\sf{HC}}^n)\rightarrow \powerset(\bR)$\end{center} is an injective function. 

\begin{definition}\label{pairing function}
Let $(p_i: i<\omega)$ be the sequence of prime numbers. Let 
\begin{center}
${\sf{merge}}:\bR^{\leq \omega}\rightarrow \bR$
\end{center}
 be given by ${\sf{merge}}(q)=y$ if letting $q=(y_i)_{i<n}$, 
\begin{center}
$y(j)=
\begin{cases}
y_{j_0}(j_1)  &: j=p_{j_0}^{j_1}\\
0 &: \text{otherwise}.
\end{cases}$
\end{center}
$\myqedhere$
\end{definition}

\begin{notation}\label{rank notation}
If $M$ is a transitive set and $\a\leq \ord(M)$ then we let $M|\a=V_\a^M$. 

$\myqedhere$
\end{notation}

\begin{definition}[Background]\label{background}\index{background}
We say\begin{center} $\mathbb{M}=(M, \d, \vec{G})$\end{center} is a \textbf{background} if 
\begin{enumerate}
\item $M\models {\sf{ZFC}}+``\d$ is a Woodin cardinal",
\item $\vec{G}:\d\rightarrow V_\d^M$ is a partial function such that for each $\a\in \dom(\vec{G})$, for some $(\k, \l)$, $M\models ``\vec{G}(\a)$ is a $(\k, \l)$-extender such that $M|\l\subseteq Ult(M, \vec{G}(\a))$ and $\l$ is inaccessible",
\item $\M\models ``\vec{G}$ witnesses that $\d$ is a Woodin cardinal"\footnote{I.e., $M\models ``$ for every $A\subseteq \d$ there is $\k$ such that for every $\l<\d$ there is a $\a$ with the property that letting $E=\vec{G}(\a)$, $\cp(E)=\k$, $\lh(E)\geq \l$ and $A\cap \lh(E)= \pi_E(A)\cap \lh(E)$.},
\end{enumerate}
We say 
\begin{center}
$\mathbb{M}=(M, \d, \vec{G}, \Sigma)$\end{center} is an \textbf{internally iterable background} if in addition to the three clauses above, the following clauses hold:
\begin{enumerate}
 \item $\Sigma\in M$ and $M\models ``\Sigma$ is a winning strategy for $II$ in the version of the iteration game $\mathcal{G}(M, \d, \d+1)$ in which player $I$ is required to choose extenders whose (natural) lengths are inaccessible cardinals in the model they are chosen from and are also below the image of $\d$".
 \item $\Sigma$ has hull condensation,
 \item $\dom(\Sigma)\subseteq \mathcal{J}_\omega(V^M_\d)$. 
 \end{enumerate}

We say $\mathbb{M}=(M, \d,  \vec{G}, \Sigma)$ is an \textbf{externally iterable  background} if $(M, \d,  \vec{G})$ is a background and $\Sigma$ is a winning strategy for $II$ in the version of $\mathcal{G}(M, \omega_1, \omega_1)$ mention in clause 1 above. We say  $\mathbb{M}=(M, \d,  \vec{G}, \Sigma)$ is an \textbf{iterable background} if it is either internally or externally iterable background.

We say that an externally iterable background $(M, \d,  \vec{G}, \Sigma)$ is \textbf{self-knowledgable} if $(M, \d,  \vec{G}, \Sigma\rest \mathcal{J}_{\omega}(M|\d))$ is an internally iterable background. 

Suppose $(M, \d,  \vec{G}, \Sigma)$ is a background and $N$ is a $\Sigma$-iterate of $M$. Let $i: M\rightarrow N$ be the iteration embedding. We set 
\begin{center} $\mathbb{M}_N=(N, i(\d), i(\vec{G}), \Sigma_N)$.\end{center} 
In most cases considered in this book, $\Sigma_N$ won't depend on the iteration producing $N$. $\myqedhere$
\end{definition}

Suppose $\mathbb{M}=(M, \d,  \vec{G}, \Sigma)$ is an externally iterable background and $A\subseteq \bR$. We review the standard capturing notions (for example see \cite{OIMT} or \cite{DMATM} and references presented in those papers). We say $\mathbb{M}$ Suslin captures $A$\index{Suslin capturing} at an $M$-cardinal $\eta$ if there is a tree $T\in M$ such that whenever $N$ is a $\Sigma$-iterate of $M$ with $i:M\rightarrow N$ the iteration embedding and whenever $g$ is $<i(\eta)$-generic over $N$, $(p[i(T)])^{N[g]}=A\cap N[g]$. We say $\mathbb{M}$ Suslin, co-Suslin captures $A$ at $\eta$ if it Suslin captures both $A$ and $A^c$ at $\eta$. We say $\mathbb{M}$ Suslin captures $A$ if $\mathbb{M}$ Suslin captures $A$ at $(\d^+)^M$, and similarly $\mathbb{M}$ Suslin, co-Suslin captures $A$ if $\mathbb{M}$ Suslin, co-Suslin captures $A$ at $(\d^+)^M$. 

Finally we recall the notion of \textit{self-capturing background} (Definition 2.24 of \cite{ATHM}). 

\begin{definition}\label{self-capturing bt}
Suppose $\mathbb{M}=(M, \d,  \vec{G}, \Sigma)$ is a self-knowledgable, externally iterable background. We say $\mathbb{M}$ is \textbf{self-capturing}\index{self-capturing background} if $\Sigma$ is positional\footnote{I.e. whenever $N$ is a $\Sigma$-iterate of $M$ via $\X$, $\Sigma_{N, \X}$ is independent of $\X$. See \cite[Definition 2.35]{ATHM}.} and for every $M$-inaccessible cardinal $\l<\d$ there is a name $\dot{X}\in M^{Coll(\omega, M|\l)}$ such that for any $M$-generic $g\subseteq Coll(\omega, M|\l)$, $(M[g], \d, \vec{G}, \Sigma)$\footnote{Here we abuse the notation a bit. In reality we should use $\Sigma'$ which is the portion of $\Sigma$ that acts on stacks above $\l+1$.} Suslin, co-Suslin captures ${\sf{Code}}(\Sigma_{M|\l})$ at $(\d^+)^M$ as witnessed by $\dot{X}_g=(T, S)$.  $\myqedhere$
\end{definition}

\rthm{n*x} is the main method for producing self-capturing backgrounds.

\subsection{Capturing pointclasses}\label{subsec: capturing pointclasses}

We recall the definition of a good pointclass (see \cite[Definition 9.12]{DMATM}). Unlike  \cite[Definition 9.12]{DMATM} we include scale property into the definition of good pointclass. 

\begin{definition}\label{good pointclass} We say $\Gamma$ is a good pointclass if $\Gamma$ is closed under recursive substitutions, is closed under quantification over $\omega$, is closed under existential quantification over $\bR$, is $\omega$-parametrized\footnote{This means that there is $U\subseteq \omega\times \bR$ such that $U\in \Gamma$ and $\{ A\subseteq \bR: A\in \Gamma\}=\{ U_e: e\in \omega\}$.} and has the scale property. $\myqedhere$
\end{definition}

Suppose $\Gamma$ is a good pointclass. For $x\in \bR$, we let $C_\Gamma(x)$ be the largest countable $\Gamma(x)$-set of reals. For transitive $a\in {\sf{HC}}$\footnote{${\sf{HC}}$ is the set of hereditarily countable sets.} and surjection $g:\omega\rightarrow a$, we let $a_g$ be the real coding $(a, \in)$ via $g$. More precisely, 
\begin{center}
$a_g(k)=
\begin{cases} 
1 :& k=2^m3^n\ \text{and}\ g(m)\in g(n)\\
0 :& \text{otherwise}.
\end{cases}$
\end{center}
Clearly $M_{a_g}=(a, \in)$. If $b\subseteq a$, then we let $b_g=\{ m: g(m)\in b\}$. We then  let $C_\Gamma(a)=\{ b\subseteq a:$ for comeager many $g:\omega\rightarrow a$, $b_g\in C_\Gamma(a_g)\}$.

Continuing with $\Gamma$, we say $P$ is a $\Gamma$-Woodin if there is a $P$-cardinal $\d_P$ such that
\begin{enumerate}
\item $P$ is countable, 
\item  $P=C_\Gamma(C_\Gamma(V_{\d_P}^P))$,
\item $P\models ``\d_P$ is the only Woodin cardinal" and
\item for every $\eta<\d_P$, $C_\Gamma(V_\eta)\models ``\eta$ is not a Woodin cardinal".
\end{enumerate}
We say $(P, \Psi)$ is a \textit{$\Gamma$-Woodin pair}\index{$\Gamma$-Woodin pair} if
\begin{enumerate}
\item $\Psi$ is an $\omega_1$-iteration strategy for $P$ and
\item for every $\Psi$-iterate $Q$ of $P$, $Q$ is a $\Gamma$-Woodin\footnote{$P$ is a coarse structure, there is no notion of dropping for iterations of $P$, so $P$-to-$Q$ embedding always exists.}. 
\end{enumerate}
Woodin, assuming  $\sf{AD}^+$,  showed that if $\Gamma$ is a good pointclass not closed under $\forall^\mathbb{R}$ then there are $\Gamma$-Woodin pairs (see \cite[Theorem 10.3]{DMATM}). To learn more on Woodin's work one may consult \cite{CoarseAD}.

Suppose $\Gamma$ is a good pointclass and $(P, \Psi)$ is a $\Gamma$-Woodin pair. Let $\mathcal{L}_{\Psi}$ be the extension of the language of set theory obtained by adding one predicate symbol $\dot{\Psi}$ and one constant symbol $e$. The intended interpretation of $\dot{\Psi}$ is ${\sf{Code}}(\Psi)$. $e$ wlll denote a real number. Given $u\in \bR$, we define $T'_n(\Psi, u)$ to be the set of $(\phi, \vec{x})$ such that $\phi$ is a $\Sigma_n$-formula in $\mathcal{L}_\Psi$, $\vec{x}\in \bR^m$ where $m$ is the number of free variables of $\phi$ and 
\begin{center}
$({\sf{HC}}, {\sf{Code}}(\Psi), u, \in)\models \phi[\vec{x}]$.
\end{center}
We let $T'_n(\Psi)=T'_n(\Psi, 0)$. 

Next we code $T'_n(\Psi, u)$ by a set of reals as follows. First let $G_\Psi$ be the set of natural numbers that are G\"odel numbers for $\mathcal{L}_\Psi$-formulae. We say $y\in \bR$ is $\Psi$-appropriate if $y(0)$ is a G\"odel number of an $\mathcal{L}_{\Psi}$ formula. If $y$ is $\Psi$-appropriate then we let $\phi_y$ be the formula that $y(0)$ codes and $l_y$ be the number of free variables of $\phi_y$. Let $(p_i: i<\omega)$ be the sequence of prime numbers in increasing order. For $i\leq l_y$, let $y_i\in \bR$ be such that for all $k\in \omega$, $y_i(k)=y(p_i^{k+1})$. If $y$ is $\Psi$-appropriate then we say $y$ is neat if for all $k'$ such that $k'\not=0$ and $k'\not \in \{ p_i^k: i< l_y \wedge k\in \omega\}$, $y(k')=0$. Let then $T_n(\Psi, u)$ be the set of $\Psi$-appropriate neat $y\in \bR$ such that 
\begin{center}
$(\phi_y, {\sf{merge}}(y_i: i< l_y))\in T_n'(\Psi, u)$. 
\end{center}
Again, set $T_n(\Psi)=T_n(\Psi, 0)$.

Suppose $z\in \bR$, $\phi$ is an $\mathcal{L}_\Psi$-formula with $l+1$ free variables and $(x_i: 2\leq i\leq l)\in \bR^m$. Let $y_0\in \bR$ be such that $y_0(0)$ is the G\"odel number of $\phi$ and for $i>0$, $y_0(i)=0$. Let $y_1=z$ and  for $2\leq i\leq l$, $y_i=x_i$. Set $a(\phi, z, \vec{x})={\sf{merge}}((y_i: i\leq l))$. Notice that $(\phi, z, \vec{x})$ is uniquely determined by $a(\phi, z, \vec{x})$. In fact, the function $(\phi, z, \vec{x})\mapsto a(\phi, z, \vec{x})$ is a $\Pi^0_1$ injection. 

 Assuming ${\sf{AD}}$, if $A\subseteq \bR$ then $w(A)$ is its Wadge rank, and if $\Gamma$ is a pointclass then $w(\Gamma)=\sup\{ w(A): A\in \Gamma\}$.
 \begin{notation} \label{lub gamma}
Suppose $\Gamma$ is a pointclass closed under continous preimages and $A\subseteq \bR$. We say $A$ is a least upper bound for $\Gamma$ if $\Gamma=\{ B\subseteq \bR: w(B)< w(A)\}$. Set then $lub(\Gamma)=\{ A\subseteq \bR: A$ is a least upper bound for $\Gamma\}$. $\myqedhere$
\end{notation}

\begin{definition}\label{capturing gamma}
Suppose $\Gamma$ is any pointclass closed under the continuous preimages. We say that the tuple $(\mathbb{M}, (\P, \Psi), \Gamma^*, A)$ Suslin, co-Suslin captures $\Gamma$ if the following conditions hold:
\begin{enumerate}
\item $A\in lub(\Gamma)$,
\item $\Gamma^*$ is the least good pointclass such that $\Gamma\subseteq \utilde{\Delta}_{\Gamma^*}$. 
\item $(P, \Psi)$ is a $\Gamma^*$-Woodin pair.
\item  $(P, \d_P, \Psi)$ Suslin, co-Suslin captures  $A$.
\item $\mathbb{M}$ is a self-capturing background. 
\item  $\mathbb{M}$ Suslin, co-Suslin captures the sequence $(T_n(\Psi): n<\omega)$.
\end{enumerate}
$\myqedhere$
\end{definition}

\begin{notation}\label{cn for capturing triples} Suppose  $\Gamma$ is a pointclass closed under the continuous preimages, ${\sf{C}}=(\mathbb{M}, (\P, \Psi), \Gamma^*, A)$ Suslin, co-Suslin captures $\Gamma$ and $\mathbb{M}=(M, \d,  \vec{G}, \Sigma)$. If $N$ is a $\Sigma$-iterate of $M$ then we set 
${\sf{C}}_N=(\mathbb{M}_N, (\P, \Psi), \Gamma^*, A)$. $\myqedhere$
\end{notation}

The following is an important yet straightforward lemma that we will use throughout this book.

\begin{terminology}\label{leq generic}
Below and throughout this book we say that $``g$ is $<\eta$-generic" to mean that the poset for which $g$ is generic has size $<\eta$. Similarly we say that $``g$ is $\leq\eta$ generic" to mean that the  poset for which $g$ is generic has size $\leq\eta$. $\myqedhere$
\end{terminology}

\begin{lemma}[Correctness of backgrounds]\label{correctness of backgrounds} Suppose $(\mathbb{M}, (\P, \Psi), \Gamma^*, A)$ Suslin, co-Suslin captures $\Gamma$ and set $\mathbb{M}=(M, \d,  \vec{G}, \Sigma)$. Suppose $x\in \bR\cap M$. Let $(S_n, U_n: n<\omega)\in M$ be the sequence of trees on $\omega\times (\d^+)^M$ such that $(S_n, U_n)$ Suslin, co-Suslin captures $T_n(\Psi)$. Let $g$ be $<\d$-generic over $M$. Then for any real $u\in M[g]$,
\begin{center}
$({\sf{HC}}^{M[g]}, {\sf{Code}}(\Psi)\cap M[g], u, \in) \prec  ({\sf{HC}}, {\sf{Code}}(\Psi), u, \in)$. 
\end{center}
\end{lemma}
\begin{proof} It is enough to verify that if $\phi$ is a formula, $m+1$ is the number of its free variables and $\vec{x}\in \bR^m\cap M[g]$ then if $({\sf{HC}}, {\sf{Code}}(\Psi), \in)\models \exists v \phi[v, \vec{x}]$ then there is $v\in M[g]\cap \bR$ such that $({\sf{HC}}, {\sf{Code}}(\Psi), \in)\models \phi[v, \vec{x}]$.
Let $n$ be such that $\phi$ is $\Sigma_n$. Then there is $v$ such that $a(\phi, v, \vec{x})\in T_n(\Psi)$. 

Working in $M[g]$, let $S'=\{ (s, h)\in S_n: s(0)$ is the G\"odel number of $\phi\}$, and let $S$ be the tree on $\omega\times \d$ whose branches are pairs $(y', f)\in \bR\times \d^\omega$ with the property that if $y=a(\phi, y', \vec{x})$ then $(y, f)\in [S']$.

We now have that whenever $h$ is any $Coll(\omega, \d)$-generic extension of $M[g]$, in $M[g][h]$, $p[S]$ is the set of $y'\in \bR^{M[g][h]}$, such that if $y=a(\phi, y', \vec{x})$ then $y\in p[S']$. Because $S_n$ Suslin captures $T_n(\Psi)$ we have that $(p[S])^{M[g]}\not =\emptyset$\footnote{This follows from genericity iterations. One can iterate $M[g]$ via $\Sigma$ to obtain $i: M[g]\rightarrow N$ such that $(p[i(S)])^{N}\not =\emptyset$. For example if $v$ is generic over $N$ for the extender algebra at $i(\d)$ then $v\in (p[i(S)])^{N}$.}. Notice next that if $v\in p[S]\cap M[g]$ then $({\sf{HC}}, {\sf{Code}}(\Psi), \in)\models \phi[v, \vec{x}]$.
\end{proof}

Self-capturing backgrounds are very useful for building hod pairs and proving comparison. The following theorem of Woodin shows that under $\sf{AD}^+$, self-capturing backgrounds are abundant. \cite{CoarseAD} has an outline of  the proof of \rthm{n*x}.

\begin{theorem}[Woodin, Theorem 10.3 of \cite{DMATM}]\label{n*x} Assume $ \sf{AD}^+$. Suppose $\Gamma$ is a good pointclass and there is a good pointclass $\gG^*$ such that $\Gamma\subseteq \Delta_{\gG^*}$. Suppose $(N, \Psi)$ is $\gG^*$-Woodin which Suslin, co-Suslin captures some $A\in lub(\Gamma)$. There is then a function $F$ defined on $\mathbb{R}$ such that for a Turing cone of $x$, $F(x)=(\N^*_x, \M_x, \d_x, \Sigma_x)$ is such that
\begin{enumerate}
\item $N\in L_1[x]$,
\item $\N^*_x|\d_x=\M_x| \d_x$,
\item $\M_x$ is a $\Psi$-mouse over $x$: in fact, $\M_x=\M_1^{\Psi, \#}(x)|\k_x$ where $\k_x$ is the least inaccessible cardinal of $\M_1^{\Psi, \#}(x)$ that is $>\d_x$,
\item $\N^*_x\models ``\d_x$ is the only Woodin cardinal",
\item $\Sigma_x$ is the unique iteration strategy of $\M_x$,
\item $\N^*_x= L(\M_x, \Lambda)$ where $\Lambda=\Sigma_x\rest \dom(\Lambda)$ and 
\begin{center}
$\dom(\Lambda)=\{\T\in \M_x: \T$ is a normal iteration tree on $\M_x$, $\lh(\T)$ is a limit ordinal and $\T$ is below $\d_x\}$,
\end{center}
\item setting $\vec{G}=\{(\a, \vec{E}^{\N^*_x}(\a)): \N^*_x\models  ``\lh(\vec{E}^{\N^*_x}(\a))$ is an inaccessible cardinal $<\d_x"\}$ and $\mathbb{M}_x=(\N^*_x, \d_x, \vec{G}, \Sigma_x)$, $(\mathbb{M}_x, (N, \Psi), \Gamma^*, A)$ Suslin, co-Suslin captures $\Gamma$\footnote{Hence, $(\N^*_x, \d_x, \vec{G}, \Sigma_x)$ is a self-capturing background.}. 
\end{enumerate}
\end{theorem}

\subsection{The meaning of ${\sf{Lp}}^\Gamma$, ${\sf{HP}}^\Gamma$ and ${\sf{Mice}}^\Gamma$}

The reader may find it helpful to review \rdef{standard indexing scheme} and \rdef{allowable pair}. Recall that we say $X$ is \textit{self-well-ordered} if there is a wellordering of $\univ X$ in $\mathcal{J}_1(X)$ definable over $\J_0(X)$. 

 \begin{definition}[The ${\sf{Lp}}$ function]\label{the lp function} Suppose $\Gamma$ is a pointclass and $(\P, \Sigma)$ is an allowable pair\footnote{See \rdef{allowable pair}.} such that ${\sf{{\sf{Code}}}}(\Sigma)\in \Gamma$. Suppose $X$ is a self-well-ordered set such that $\P\in X$.
 \begin{enumerate}
 \item If $\Sigma$ is an iteration strategy then ${\sf{Lp}}^{\Gamma, \Sigma}(X)$ is the stack of all sound $(\Sigma, \sf{sis})$-mice $\M$ over $X$ based on $\P$\footnote{See \rdef{strategy premouse}, \rdef{sigma-premouse} and \rdef{standard indexing scheme}.} such that $\rho(\M)=\ord(trc(X))$\footnote{Our fine structural notation was introduced in \rdef{the core}.} and $\M$ has a strategy in $\Gamma$.
 \item If $\Sigma$ is a st-strategy ${\sf{Lp}}^{\Gamma, \Sigma}(X)$ is the stack of all sound $\Sigma$-sts mice $\M$ over $X$ based on $\P$ such that $\rho(\M)=\ord(trc(X))$ and $\M$ has a strategy in $\Gamma$\footnote{From here on, ``Lp" means ``g-organized Lp" as defined in \cite{trang2013} unless explicitly stated otherwise. We will occasionally remind the reader of this convention. The reason we need to use g-organization is so that $S$-constructions go through.}.
 \end{enumerate} 
 We set ${\sf{Lp}}^{\Gamma, \Sigma}(\P)={\sf{Lp}}^{\Gamma, \Sigma}(\mathcal{J}_\omega[\P])$. $\myqedhere$
 \end{definition}
 
  Below if $\Psi$ is an iteration strategy or an st-strategy then we let $M_\Psi$ be the structure that $\Psi$ is iterating.

\begin{notation}\label{mice relative to gamma} Suppose $\Gamma$ is a pointclass. Following Section 2.5 of \cite{ATHM} we let 
\begin{center}
${\sf{Hp}}^\Gamma=\{ (\P, \Sigma): (\P, \Sigma)$ is an allowable pair such that ${\sf{Code}}(\Sigma)\in \Gamma\}$\\

${\sf{Mice}}^\Gamma=\{(a, \Sigma, \M): a\in {\sf{HC}}\wedge a \text{\ is a self-well-ordered} \wedge (M_\Sigma, \Sigma)\in {\sf{Hp}}^\Gamma\wedge \M_\Sigma\in a\wedge \M\insegeq {\sf{Lp}}^{\Gamma, \Sigma}(a) \wedge \rho(\M)=\ord(trc(a))\}$
\end{center}
and given $(\P, \Sigma)\in {\sf{Hp}}^{\Gamma}$, 
\begin{center}
${\sf{Mice}}^\Gamma_{\Sigma}=\{(a, \M): a\in {\sf{HC}}\wedge a \text{\ is a self-well-ordered} \wedge \P \in a\wedge \M\insegeq {\sf{Lp}}^{\Gamma, \Sigma}(a) \wedge \rho(\M)=\ord(trc(a))\}$
\end{center}
When $\Gamma=\powerset(\bR)$, we omit it from our notation. 

Suppose $A\subseteq \bR$ with $w(\Gamma)\leq w(A)$. We say $\sigma\in \bR$ is an $A$-code if $\sigma(0)$ is a G\"odel number for some formula $\phi$, and if $B$ is the set of reals definable over $({\sf{HC}}, A, \sigma, \in)$ via $\phi$\footnote{I.e., $u\in B\iff  ({\sf{HC}}, A, \sigma, \in)\models \phi[u]$.} then $B\in \rge(\sfc)$. We then let $C_{\sigma}=\sfc^{-1}(B)$ and $A{\sf{Code}}$ be the set of $A$-codes.

Given a set $A\subseteq \bR$ with $w(\Gamma)\leq w(A)$, we let ${\sf{Code}}({\sf{Hp}}^\Gamma, A)$ be the set of $\sigma\in A{\sf{Code}}$ such that $C_{\sigma}\in {\sf{Hp}}^\Gamma$. If $\sigma\in{\sf{Code}}({\sf{Hp}}^\Gamma, A)$ then we let $(\P_\sigma, \Sigma_\sigma)$ be the pair determined by $\sigma$.   

Given a set $A\subseteq \bR$ with $w(\Gamma)\leq w(A)$, we let ${\sf{Code}}({\sf{Mice}}^\Gamma, A)$ be the set of $(\sigma_0, \sigma_1)$ such that $\sigma_0\in A{\sf{Code}}$, $\sigma_1\in A{\sf{Code}}$ and $C_{\sigma_1}={\sf{Mice}}^\Gamma_{\Sigma_{\sigma_0}}$.  

Given a set $A\subseteq \bR$ with $w(\Gamma)\leq w(A)$, we let $A_\Gamma$ be the set of triples $(\sigma_0, \sigma_1,  \sigma_2)$ such that 
\begin{enumerate}
\item For each $i<3$, $\sigma_i\in A\sf{Code}$,
\item $\sigma_0\in {\sf{HP}}^\Gamma$,
\item $C_{\sigma_1}=(a, \Sigma_{\sigma_0}, \M)\in {\sf{Mice}}^\Gamma$,
\item $C_{\sigma_2}$ is the unique $\omega_1$-iteration strategy of $\M$. 
\end{enumerate}
$\myqedhere$
\end{notation}

The following is an easy consequence of Lemma \ref{correctness of backgrounds}. It follows from the fact that each of
\begin{center}
 ${\sf{Code}}(A)_\Gamma$,
${\sf{Code}}({\sf{Hp}}^\Gamma, A)$ and ${\sf{Code}}({\sf{Mice}}^\Gamma, A)$.
\end{center}
is definable over $({\sf{HC}}, \sfc(\Psi), \in)$, where $(P, \Psi)$ is as below.
\begin{corollary}\label{capturing particular sets} Suppose $\mathbb{M}=(M, \d,  \vec{G}, \Sigma)$ and $(\mathbb{M}, (P, \Psi), \Gamma^*, A)$ Suslin, co-Suslin captures $\Gamma$. Then $\mathbb{M}$ Suslin, co-Suslin captures
\begin{center}
 ${\sf{Code}}(A)_\Gamma$,
${\sf{Code}}({\sf{Hp}}^\Gamma, A)$ and ${\sf{Code}}({\sf{Mice}}^\Gamma, A)$.
\end{center}
\end{corollary}


\subsection{Internalizing ${\sf{HP}}^\Gamma$}\label{internalizing gamma sets}
Suppose next that $\Gamma$ is a pointclass and $(\mathbb{M}, (P, \Psi), \Gamma^*, A)$ Suslin, co-Suslin captures $\Gamma$. In \rsec{sec: hod pair constructions}, we will describe the $\Gamma$-hod pair construction of $\mathbb{M}$ that produces hod pairs. When describing this construction, we will use the following concepts and simple observations. 

\begin{definition}\label{generically absolute pair} Suppose $M$ is a transitive model of $\sf{ZFC}$, $X\in M$ and $\phi$ is a formula. We say $(X, \phi)$ is $(M, \eta)$-\textbf{generically absolute} if for some $\theta\geq \eta$ such that $X\in V_\theta^M$, for all $Y\prec (V^M_\theta, X, \eta \in)$ such that $Y\in M$ and $M\models ``Y$ is countable", letting $N_Y$ be the transitive collapse of $Y$ and $\pi_Y: Y\rightarrow N_Y$ be the collapse map, whenever $g\in M$ is $\leq\pi(\eta)$-generic over $N_Y$ and $x\in N_Y[g]\cap \bR$,
\begin{center}
$N_Y[g]\models \phi[\pi_Y(X), x] \iff M\models \phi[X, x]$.
\end{center}
$\myqedhere$
\end{definition}

\begin{definition}\label{generically absolute pair+} Suppose $M$ is a transitive model of $\sf{ZFC}$, $X\in M$ and $\phi$ is a formula. We say $(X, \phi)$ is $(M, \eta, \a)$-\textbf{generically absolute} if $\a<\eta$ and whenever $g\subseteq Coll(\omega, M|\a)$ is $M$-generic, $((X,g), \phi)$ is $(M[g], \eta)$-\textbf{generically absolute}. $\myqedhere$
\end{definition}

The following theorem can be proven by using the \textit{Tree Production Lemma} (see \cite[Lemma 4.1]{DMT}).

\begin{lemma}\label{getting trees for m-ub} Suppose $M$ is a transitive model of $\sf{ZFC}$ and $(X, \phi)$ is $(M, \eta)$-\textbf{generically absolute}. 
There is then a pair $(T, S)\in OD^{M}_{X}$ such that $(T, S)$ is $\leq\eta$-absolutely complementing and whenever $g$ is $\leq\eta$-generic
\begin{center}
$(p[T])^{M[g]}=\{ x: M[g]\models \phi[X, x]\}$.
\end{center}
\end{lemma}

\begin{definition}\label{satisfaction rel for self-capturing triples} Suppose $\Gamma$ is a pointclass,  $(\mathbb{M}, (P, \Psi), \Gamma^*, A)$ Suslin, co-Suslin captures $\Gamma$, $\mathbb{M}=(M, \d,  \vec{G}, \Sigma)$ and $(X, \phi)$ is $(M, \d, \a)$-generically absolute. 


We then write
\begin{center}
$M\models (X, \phi)\in {\sf{Hp}}^\Gamma$
\end{center}
to mean that the following holds. 

Whenever $g\subseteq Coll(\omega, M|\a)$ is $M$-generic,  there is a real $\sigma\in M[g]\cap \sfc({\sf{Hp}}^\Gamma, A)$
 such that letting $\tau$ be the formula coded by $\sigma(0)$, whenever $h$ is $\leq \d$-generic over $M[g]$, in $M[g][h]$, 
\begin{center}
$\{x\in \bR: \phi[(X, g), x]\}=\{x: ({\sf{HC}}^{M[g][h]}, A\cap M[g][h], \sigma, \in )\models \tau[c_x]\}$.
\end{center}
where $c_x$ is the set coded by $x$. $\myqedhere$
\end{definition}

The following lemma is a straightforward consequence of genericity iterations.

\begin{lemma}\label{uniqueness of sigma code} Suppose $\Gamma$ is a pointclass,  $(\mathbb{M}, (P, \Psi), \Gamma^*, A)$ Suslin, co-Suslin captures $\Gamma$, $\mathbb{M}=(M, \d,  \vec{G}, \Sigma)$, $(X, \phi)$ is $(M, \d, \a)$-generically absolute and $M\models (X, \phi)\in {\sf{Hp}}^\Gamma$. Suppose $g\subseteq Coll(\omega, M|\a)$ is $M$-generic and $\sigma_0, \sigma_1\in M[g]\cap \sfc({\sf{Hp}}^\Gamma, A)$ are two reals witnessing that $M\models (X, \phi)\in {\sf{Hp}}^\Gamma$. Then $(\P_{\sigma_0}, \Sigma_{\sigma_0})=(\P_{\sigma_1}, \Sigma_{\sigma_1})$.
\end{lemma}

The following now is not hard to show. It follows from \rlem{correctness of backgrounds}, which implies that
\begin{center}
$({\sf{HC}}^{M[g][h]}, A\cap M[g][h], \sigma, \in)\prec ({\sf{HC}}, A, \sigma, \in)$,
\end{center}
and also from \rlem{getting trees for m-ub}.

\begin{lemma}\label{computing hpgamma sets correctly} Suppose $\Gamma$ is a pointclass,  $(\mathbb{M}, (P, \Psi), \Gamma^*, A)$ Suslin, co-Suslin captures $\Gamma$, $\mathbb{M}=(M, \d,  \vec{G}, \Sigma)$ and $(X, \phi)$ is $(M, \d, \a)$-generically absolute. Suppose further that $g\subseteq Coll(\omega, M|\a)$ is $M$-generic, $\sigma\in M[g]\cap \sfc({\sf{Hp}}^\Gamma, A)$ witnesses that $M\models (X, \phi)\in {\sf{Hp}}^\Gamma$ and $\tau$ is the formula coded by $\sigma(0)$. Set $u\in C$ if and only if
\begin{center}
 there is an iteration $i: M\rightarrow N$ according to $\Sigma$ such that $\cp(i)>\a$ and for some $N[g]$-generic $h\subseteq Coll(\omega, i(\d))$, $u\in N[g][h]$ and $N[g][h]\models \tau[c_x]$.\end{center}
 
  Then $\sfc^{-1}(C)=(\P_\sigma, \Sigma_\sigma)\in {\sf{HP}}^\Gamma$ and $C$ is Suslin, co-Suslin captured by $(M[g], \d, \Sigma)$.
\end{lemma}

\section{Fully backgrounded constructions relative to short tree strategy}

Suppose $(M, \d,  \vec{G}, \Sigma)$ is an iterable background and $\P\in V_\d^M$ is a $\#$-lsa type hod premouse ( see \rdef{lsa type}). Suppose $\Lambda\in M$ is a $(\d, \d, \d)$ st-strategy for $\P$ and $X\in V_\d^M$ is a transitive self-well-ordered set such that $\P\in X$. We can then define the model $\mathcal{J}^{\vec{E}, \Lambda}(X)$ exactly like in the case $\Lambda$ is an iteration strategy. The construction will ensure that the model $\mathcal{J}^{\vec{E}, \Lambda}(X)$ is an sts premouse over $X$ based on $\P$. Here is the precise definition. 

Recall that if $(\M_\a: \a<\xi)$ is a sequence of $\mathcal{J}$-structures and $\xi$ is a limit ordinal then $\M=lim_{\a\rightarrow \xi}\M_\a$ is the $\mathcal{J}$-structure with the property that for each $\b$ such that $\mathcal{J}_\b^\M$ is defined, there is $\gg<\xi$ such that for all $\a\in (\gg, \xi)$, $\J_\b^{\M_\a}=\J_\b^\M$.

Suppose $(M, \d,  \vec{G}, \Sigma)$ is an internally or externally iterable background, $A\subseteq V_\d^M$ and $E\in V_\d^M$ is an extender. Then we say $E$ \textit{coheres} or \textit{reflects} $A$ if $\nu(E)$ is an inaccessible cardinal of $M$, $V_{\nu(E)}^M\subseteq Ult(\M, E)$ and $A\cap V_{\nu(E)}^M=\pi_E(A)\cap V_{\nu(E)}^M$. Recall that an ${\sf{lses}}$ $\M$ is reliable if for all $k$, ${\sf{core}}_k(\M)$ exists and $({\sf{core}}_k(\M), k)$ is $\omega_1+1$-iterable (see \rdef{the core} and \cite[Chapter 11]{FSIT}). Finally recall our notation $\univ{M}$ denoting the universe of $M$. This notation was introduced in \rsec{sec: jstructures}. Finally recall that sts premice are $\sf{sts}$-indexed (see \rdef{sts0} and \rdef{sts premouse}). 

As was stated many times, in this book we are mostly concerned with new issues that arise from dealing with sts mice. Reproving  all the well-established facts will add 1000s of more pages to this book without adding any new ideas. In particular, our exposition of the fully backgrounded constructions heavily relies on \cite{FSIT} and \cite[Chapter 5]{FarmDef}. The later proves the uniqueness of the next extender in full generality. 

\begin{definition}\label{fully backgrounded sts construction} Suppose $(M, \d,  \vec{G})$ is a background  and $(\P, \Lambda)\in M$ is an sts hod pair\footnote{In particular, $\Lambda$ has hull condensation. An easy Skolem hull and a realizability argument implies that if $E$ is a countably complete total extender in $M$ then $\pi_E(\Lambda)=\Lambda\rest Ult(M, E)$.}, $\Lambda$ is $(\d, \d, \d)$ st-strategy for $\P$ and $X\in V_\d^M$ is a transitive self-well-ordered set such that $\P\in X$. Suppose further that $\Lambda$ has hull condensation. Then 
\begin{center}
${\sf{Le}}((\P, \Lambda), X)=(\M_\gg , \N_\gg, F^+_\gg, F_\gg, b_\gg: \gg\leq \delta)$
\end{center}
is the output of \textbf{the fully backgrounded construction of $(M, \d,  \vec{G})$ relative to
$\Lambda$ done over $X$ using the coherence condition} if the following conditions hold. 
\begin{enumerate}
\item $\M_0=\mathcal{J}_\omega(X)$, and for all $\gg<\d$, each of $\M_\gg$ and $\N_\gg$ is either undefined or is a $\Lambda$-sts premouse\footnote{See \rdef{lambda sts premouse}.}.
\item If for some $\xi\leq \eta$, $\N_{\xi}$ is defined but is not a reliable sts premouse over $X$ based on $\P$  then all other objects with index $\geq \xi$ are undefined.
\item Suppose for some $\xi<\d$, for all $\gg\leq \xi$, both $\M_\gg$ and $\N_\gg$ are defined. Then $\M_{\xi+1}$, $\N_{\xi+1}$, $F^+_\xi$, $F_\xi$ and $b_{\xi}$ are determined as follows.
\begin{enumerate}
\item Suppose $\M_\xi=(\mathcal{J}^{\vec{E}, f}_{\omega\a}, \in, \vec{E}, f)$ is a passive ${\sf{ses}}$\footnote{I.e., with
no last predicate} and there is an extender $F^*\in \vec{G}$,
 an extender $F$ over $\M_\xi$, and an
ordinal $\nu<\omega\a$ such that 
\begin{enumerate}
\item $\nu<\nu(F^*)$,
\item $F=F^*\cap ([\nu]^{\omega}\times \univ{\M_\xi})$, and
 \item setting
    \begin{center}
    $\N_{\xi+1}=(\mathcal{J}_{\omega\a}^{\vec{E}, f}, \in, \vec{E}, f, \tilde{F})$
    \end{center}
     where $\tilde{F}$ is the amenable code of $F$\footnote{For the definition of the ``amenable code" see the last paragraph on page 14 of \cite{OIMT}.}, clause 2 fails for $\xi+1$.
     \end{enumerate} 
     Then
$\M_{\xi+1}=\C(\N_{\xi+1})$\footnote{Recall that $\C(\M)$ is the core of $\M$.}, $F_\xi^+=\vec{G}(\xi)$ where $\xi$ is the least such that $F^*=\vec{G}(\xi)$ has the above properties, $F_\xi=F^+\cap ([\nu]^{\omega}\times \univ{\M_\xi})$ where $\nu$ is chosen so that the above clauses hold and $b_\xi=\emptyset$.
\item Suppose $\M_\xi=(\mathcal{J}^{\vec{E}, f}_{\omega\a}, \in, \vec{E}, f)$ is a passive ${\sf{ses}}$, $\a=\b+\gg$ and there is  $t=(\P_0, \T, \P_1, \U)\in \univ{\M_\xi|\omega\b} \cap \dom(\Lambda)$ such that setting $w=(\mathcal{J}_{\omega}(t), t, \in)$, $w$ is $(f, {\sf{sts}})$-minimal as witnessed by $\b$\footnote{See \rdef{important notation}. In particular, this means that we have to index the branch of $t$ at $\omega\a$.} and $\gg=\lh(t)$. 
Set $b=\Lambda(t)$ and
     \begin{center}
     $\N_{\xi+1}=(\mathcal{J}^{\vec{E}, f^+}_{\omega\b+\omega\gg}, \in, \vec{E}, f, \tilde{b})$
     \end{center}
     where $\tilde{b}\subseteq \omega\b+\omega\gg$ is defined by $\omega\b+\omega\nu\in \tilde{b}\iff \nu \in b$. Assuming clause 2 fails for $\xi+1$, 
$\M_{\xi+1}={\sf{core}}(\N_{\xi+1})$, $F^+_\xi=F_\xi=\emptyset$ and $b_\xi=\tilde{b}$.\\\\
\textbf{Important Anomaly:} Suppose $t$ is $\sf{nuvs}$ and suppose $e\in \M_\xi|\omega\b$ is such that $\M_\xi|\omega\b\models {\sf{sts}_0}(t, e)$\footnote{See \rdef{sts0}. This means that $e$ is the branch of $t$ we must choose.}.
 If $e\not =b$ then $\N_{\xi+1}$ is not an sts premouse over $X$ based on $\P$, and so clause 2 holds. \\
\item If $\M_\xi$ doesn't satisfy clause 2a or 2b then set $\N_{\xi+1}=\mathcal{J}_{\omega}[\M_\xi]$. Assuming clause 2 fails for $\xi+1$, 
$\M_{\xi+1}=\C(\N_{\xi+1})$, $F_\xi^+ =F_\xi=b_\xi=\emptyset$. 
\end{enumerate}


\item Suppose $\xi\leq \d$ is a limit ordinal and for all $\gg<\xi$, both $\M_\gg$ and $\N_\gg$ are defined. Then $\M_{\xi}$ and $\N_{\xi}$ are determined as follows\footnote{$F_\xi, b_\xi$ will be defined at the next stage of the induction as in clause 2.}. Set $\N_\xi=lim_{\a\rightarrow \xi}\M_\a$. Assuming clause 2 fails for $\xi$, $\M_\xi=\C(\N_\xi)$. 
\item $\M_\d=\N_\d$ and $F^+_\d=F_\d=b_\d=\emptyset$.
\end{enumerate}
We say that ${\sf{Le}}((\P, \Lambda), X)$ is \textbf{successful} if for all $\xi<\d$ clause 2 above fails. Given $\k<\d$, we can also define ${\sf{Le}}((\P, \Lambda), X)_{\geq \k}$ by requiring that in clause 3.a, $\cp(F)\geq \k$. 

We will use the following terminology. We say $\Q$ is an $\N$-\textbf{model} of ${\sf{Le}}((\P, \Lambda), X)_{\geq \k}$ if for some $\gg\leq \d$, $\Q=\N_\gg$. Similarly we define $\M$-model and other such expressions. We say $\Q$ is \textbf{the last model} of ${\sf{Le}}((\P, \Lambda), X)_{\geq \k}$ if $\Q=\N_\d$. $\myqedhere$
\end{definition}

The fully backgrounded constructions of both \cite{FSIT} and \cite[Chapter 5]{FarmDef} do not use the coherence condition. In most cases considered in this book, we also do not need the coherence condition. The following theorem is essentially a corollary to \cite[Chapter 12]{FSIT}.

\begin{theorem}\label{successful le constructions} Suppose $(M, \d,  \vec{G}, \Sigma)$ is an iterable background  and $(\P, \Lambda)\in M$ is an sts hod pair, $\Lambda$ is $(\d, \d, \d)$ st-strategy for $\P$ and $X\in V_\d^M$ is a transitive self-well-ordered set such that $\P\in X$. Then for any $\k<\d$, ${\sf{Le}}((\P, \Lambda), X)_{\geq \k}$ is not successful if and only if for some $\xi<\d$, the Anomaly stated in clause 3.b of \rdef{fully backgrounded sts construction} holds.
\end{theorem}

\begin{remark}\label{general fully backgrounded constructions} Assuming that $(\P, \Lambda)$ is a pair with the property that $\P$ is an $\sf{lses}$ and $\Lambda$ is an iteration strategy for $\P$ with hull condensation, we could define ${\sf{Le}}((\P, \Lambda), X)_{\geq \k}$  just like above except in clause 3.b we require $t$ to be a stack on $\P$ according to $\Lambda$, $b=\Lambda(t)$ and $w$ be $(f, \sf{sis})$-minimal (see \rdef{standard indexing scheme}). ${\sf{Le}}((\P, \Lambda), X)_{\geq \k}$ can be defined for various types of strategies; in particular, it can be defined for $\omega_1$-strategies and for $(\omega_1, \omega_1)$-strategies. We also let ${\sf{Le}}$ be the construction relative the $\emptyset$. Thus, the models of $\sf{Le}$ are simply ordinary premice. We leave the details of the above mentioned constructions to the reader who may want to consult \cite[Definition 2.3]{ATHM}. $\myqedhere$
\end{remark}

The important comment in clause 3.b is a non-trivial matter. Recall that according to our $\sf{sts}$ indexing scheme (see \rdef{weak psi alpha indexing scheme a}), the branch we have to index at stage $\xi$ in clause 3.b is $e$ not $b$. However, if $e\not =b$ then the resulting structure cannot be a $\Lambda$-sts mouse. Thus, if $e\not=b$ then we have to halt the construction. When $\Lambda$ has nice properties such as \textit{strong branch condensation} (see \rdef{strong branch condensation}) then such anomaly will never arise. See \rrem{on hod pair constructions} for an in-depth  discussion of this issue.

\section{Hod pair constructions}\label{sec: hod pair constructions}

In  this section we introduce the $\Gamma$-hod pair constructions. The goal of such a construction is to produce a hod pair $(\P, \Sigma)$ such that $w({\sf{Code}}(\Sigma))\geq w(\Gamma)$ but for any hod initial segment $\Q\inseg \P$, $w({\sf{Code}}(\Sigma_\Q))<w(\Gamma)$ (or equivalently $(\Q, \Sigma_\Q)\in {\sf{Hp}}^\Gamma$). 

%

The reader may benefit from reviewing the concept of fully backgrounded constructions as presented in \cite[Chapter 11 and Chapter 12]{FSIT}. Such constructions inherit a strategy from the background model\footnote{The model where the construction is being done.} via the procedure described in \cite[Chapter 12]{FSIT}. Other forms of such constructions also have appeared in \cite{schlutzenberg2021background} and \cite{FarmStrongs}.

Suppose $\Gamma$ is a good pointclass and $\mathbb{M}=(M, \d,  \vec{G}, \Sigma)$ is a background Suslin, co-Suslin capturing $\Gamma$ (see \rdef{capturing gamma}). We will work with $\mathbb{M}$ and $\Gamma$, but we will omit both from our notations. 

All concepts introduced here depend on $\mathbb{M}$. For instance, $\sf{E}$ below should really be $\sf{E}^\mathbb{M}$. Also, all fully backgrounded constructions that we will use are fully backgrounded constructions in the sense of $V_\d^M$, and if $M$ is equipped with a distinguished extender sequence then we tacitly assume that all the backgounded constructions use extenders from this particular extender sequence. 
 
 The reader may find it helpful to review \rdef{layers of hod-like lsp}, \rdef{l p}, \rter{types of lsa small premice}, \rdef{the bottom part of lsp} and \rdef{allowable pair}. We start by introducing those hod premice that can be used as layers in the $\Gamma$-hod pair construction.
 \begin{definition}\label{cbl}
($\sf{CBL}$:) We say that an allowable pair $(\R, \Lambda)$ \textit{can be $\Gamma$-layered} or is just $\Gamma-\sf{cbl}$ if one of the following conditions hold:
\begin{enumerate}
\item $\R$ is a hod premouse of successor type and $\R={\sf{Lp}}_{\omega}^{\Gamma, \Lambda_{\R^-}}(\R|\d^\R)$. 
\item $\R$ is a properly non-meek\footnote{See \rdef{proper type II}.} hod premouse of limit type and letting $\k=\d^{\R^b}$,
\begin{center}
$\R^b={\sf{Lp}}^{\Gamma, \Lambda_{\R|\k}}(\R^b|\k)$.
\end{center}
\item $\R$ is a gentle hod premouse such that if $\Q\in Y^\R$ then $(\Q, \Lambda_\Q)$ is $\Gamma-\sf{cbl}$.
\end{enumerate}
$\myqedhere$
\end{definition}
Recall that an $\sf{lses}$ $\M$ over $\emptyset$ set has a predicate, $Y^\M$, whose members are the layers of the $\sf{lses}$ (see \rdef{layered hybrid j-structure}). Thus, below, when describing $\M_\gg$, we must also simultaneously define $Y^{\M_\gg}$. 

$\Phi_\gg$ below will be the iteration strategy induced by $\Sigma$ essentially via the resurrection procedure describe in \cite[Chapter 12]{FSIT}. The procedure described in \cite[Chapter 12]{FSIT} only induces $(\omega_1, \omega_1)$-iteration strategies, but  it is not hard to modify it to obtain an $(\omega_1, \omega_1, \omega_1)$-strategy (see \rdef{the un-dropping iteration game}). We will give an outline of how to do this after \rdef{gamma-hod pair construction*}.

\begin{terminology}\label{below strategy} Suppose $\M$ is an $\sf{lses}$ and $\a\leq \ord(\M)$. We say that a stack $\T$ on $\M$ is \textbf{below} $\a$ if for every $\gg<\lh(\T)$ such that $[0, \gg)_\T\cap D^\T=\emptyset$, $\ind_\gg^{\T}<\pi^\T_{0, \gg}(\a)$.  Similarly we define the meaning of ``below $\a$" for generalized stacks.

Suppose $\M$ is an $\sf{lses}$ and $\Q\insegeq \M$. We say $\Sigma$ is a strategy of $\M$ \textbf{based} on $\Q$ if whenever $\T$ is according to $\Sigma$,  $\T$ is below $\ord(\Q)$. $\myqedhere$
\end{terminology}

As was mentioned before, our exposition of the fully backgrounded constructions heavily relies on \cite{FSIT} and \cite[Chapter 5]{FarmDef} . As was mentioned before, the later reference proves the uniqueness of the next extender in full generality. 

\begin{definition}\label{gamma-hod pair construction*} Suppose $\Gamma$ is a pointclass,  ${\sf{C}}=(\mathbb{M}, (P, \Psi), \Gamma^*, A)$ Suslin, co-Suslin captures $\Gamma$ and $\mathbb{M}=(M, \d,  \vec{G}, \Sigma)$. Then 
\begin{center}${\sf{hpc}}=(\M_\gg , \N_\gg, Y_\gg, \Phi_\gg, F^+_\gg, F_\gg, b_\gg : \gg\leq \delta)$
\end{center}
 is the output of the $\Gamma$-\textbf{hod pair construction} ($\Gamma-\sf{hpc}$) of $\mathbb{M}$ if the following conditions hold (the construction is over $\emptyset$).
 \begin{enumerate}
 \item $\M_0=\mathcal{J}_\omega$, and for all $\gg\leq \d$, each of $\M_\gg$ and $\N_\gg$ is either undefined or is an $\sf{hp}$-indexed $\sf{lses}$ (see \rdef{sis}).
 \item For all $\gg\leq \d$, if $\M_\gg$ is defined then $Y_\gg=Y^{\M_\gg}$ (see \rdef{layered hybrid j-structure}).
 \item For all $\gg\leq \d$, if $\M_\gg$ is defined then $\Phi_\gg$ is the strategy defined in \rdef{gammastrategy}\footnote{This strategy is induced by $\Sigma$ essentially via the resurrection procedure of \cite[Chapter 12]{FSIT}.}.
 \item For all $\gg\leq \d$, if $\N_\gg$ is defined and either
 \begin{enumerate}
 \item $\N_\gg$ is not a reliable $\sf{hp}$-indexed $\sf{lses}$\footnote{Recall clause 2 of \rdef{strategy lhp}. To verify that $\N_\gg$ is $\sf{lses}$, we need to verify that clause 2 of \rdef{strategy lhp} holds.} or
 \item $\N_\gg$ is a reliable $\sf{hp}$-indexed $\sf{lses}$ but for some $\Q\in Y^{\N_\gg}$ such that $\Q$ is meek or gentle\footnote{See \rdef{pre-hod-like}.} and for some $n<\omega$, $\rho_n(\N_\gg)\leq \d^\Q$,
 \end{enumerate}
  then all remaining objects with index $\geq \gg$ are undefined. \\\\
  For all $\gg\leq \eta$ for which clause 4 (the above statement) fails, $\pi_\gg:\C(\N_\gg)\rightarrow \N_\gg$ is the uncollapse map. 
 \item Suppose for some $\xi<\d$, for all $\gg\leq \xi$, both $\M_\gg , \N_\gg$ are defined. Then $\M_{\xi+1}$, $\N_{\xi+1}$, $Y_{\xi+1}$, $\Phi_{\xi+1}$, $F_{\xi}^+$, $F_\xi$ and $b_{\xi}$ are deteremined as follows.
 \begin{enumerate}
\item Suppose $\M_\xi=(\mathcal{J}^{\vec{E}, f}_{\omega\a}, \in, \vec{E}, f, Y_\xi, \in)$ is a passive ${\sf{hp}}$-indexed $\sf{lses}$\footnote{I.e., with
no last predicate.}, there is an extender $H^*\in \vec{G}$ 
an extender $H$ over $\M_\xi$, and an
ordinal $\nu<\omega\a$ such that $\nu<\lh(H^*)$ and setting
    \begin{center}
    $H=H^*\cap ([\nu]^{\omega}\times \univ{\M_\xi})$, and
    $\N_{\xi+1}=(\mathcal{J}_{\omega\a}^{\vec{E}, f}, \in, \vec{E}, f, Y_\xi, \tilde{H}, \in)$
    \end{center}
    where $\tilde{H}$ is the amenable code of $H$, clause 4.a fails for $\xi+1$. 
    Then letting $\iota\in \dom(\vec{G})$ be the least such that $H^*=_{def}\vec{G}(\iota)$ has the above properties, 
    \begin{center}
    $\N_{\xi+1}=(\mathcal{J}_{\omega\a}^{\vec{E}, f}, \in, \vec{E}, f, Y_\xi, \tilde{H}, \in)$
    \end{center}
    where $\tilde{H}$ is the amenable code of $H$\footnote{Here $H$ is what is determined by $H^*$. For the definition of the ``amenable code" see the last paragraph on page 14 of \cite{OIMT}.}. Assuming clause 4 fails for $\xi+1$, the remaining objects 
     are defined as follows.
    \begin{enumerate} 
    \item $\M_{\xi+1}=\C(\N_{\xi+1})$\footnote{Recall that $\C(\M)$ is the core of $\M$.}, 
    \item $F^+_\xi= H^*$ and $F_\xi=H$, 
    \item $b_\xi=\emptyset$ and
    \item $Y_{\xi+1}=\pi^{-1}_{\xi+1}(Y_\xi)$.
\end{enumerate}
\item Suppose $\M_\xi=(\mathcal{J}^{\vec{E}, f}_{\omega\a}, \in, \vec{E}, f, Y_\xi, \in)$ is a passive ${\sf{hp}}$-indexed $\sf{lses}$, $\M_\xi$ is strategy-ready\footnote{See \rdef{index-ready}.}, $\a=\b+\gg$ and there is  $t\in \univ{\M_\xi|\omega\b}$ such that setting $w=(\mathcal{J}_{\omega}(t), t, \in)$, $w$ is $(f, {\sf{hp}})$-minimal as witnessed by $\b$\footnote{See \rdef{important notation}. In particular, this means that we have to index the branch of $t$ at $\omega\a$.} and $\gg=\lh(t)$. 
Set $b=\Phi_\xi(t)$ and
     \begin{center}
     $\N_{\xi+1}=(\mathcal{J}^{\vec{E}, f^+}_{\omega\b+\omega\gg}, \in, \vec{E}, f, Y_\xi, \tilde{b}, \in)$
     \end{center}
     where $\tilde{b}\subseteq \omega\b+\omega\gg$ is defined by $\omega\b+\omega\nu\in \tilde{b}\iff \nu \in b$. Assuming clause 4 fails for $\xi+1$, the remaining objects 
     are defined as follows. 
         \begin{enumerate} 
    \item $\M_{\xi+1}=\C(\N_{\xi+1})$, 
    \item $F_\xi=F^+_\xi=\emptyset$, 
    \item $b_\xi=\tilde{b}$ and
    \item $Y_{\xi+1}=\pi^{-1}_{\xi+1}(Y_\xi)$. \\
    \end{enumerate}

\textbf{Important Anomaly:} Suppose $\cup Y_\xi$ is $\#$-lsa type\footnote{See \rdef{lsa type}.} and $t$ is $\sf{nuvs}$. Suppose $e\in \M_\xi|\omega\b$ is such that $\M_\xi|\omega\b\models {\sf{sts}_0}(t, e)$\footnote{See \rdef{sts0}. This means that $e$ is the branch of $t$ we must choose.}.
 If $e\not =b$ then $\N_{\xi+1}$ is not an sts premouse over $\mathcal{J}_{\omega}(\cup Y_\xi)$ based on $\cup Y_{\xi}$, and so the construction must stop.\\
\item If $\M_\xi$ doesn't satisfy clause 2a or 2b then set $\N_{\xi+1}=\mathcal{J}_{\omega}[\M_\xi]$ (this presupposes that $Y^{\N_{\xi+1}}=Y_\xi$). Assuming clause 4 fails for $\xi+1$, the remaining objects 
     are defined as follows.
  \begin{enumerate} 
    \item $\M_{\xi+1}=\C(\N_{\xi+1})$\footnote{Recall that $\C(\M)$ is the core of $\M$.}, 
    \item $F_\xi=F^+_\xi=\emptyset$,
    \item $b_\xi=\emptyset$,
    \end{enumerate}
    and $Y_{\xi+1}$ is defined as follows.
\begin{enumerate}
\item If $M\models ``(\M_\xi, \Phi_\xi)$ is $\Gamma$-{\sf{cbl}}\ $+(X_\xi, \phi_\xi)\in {\sf{Hp}}^\Gamma"$\footnote{$(X_\xi, \phi_\xi)$ is defined in \rdef{xgamma}. The meaning of  $M\models (X_\xi, \phi_\xi)\in {\sf{Hp}}^\Gamma$ is essentially that $M\models (\M_\xi, \Phi_\xi)\in {\sf{Hp}}^\Gamma$.} then $Y_{\xi+1}=\pi^{-1}_{\xi+1}(Y_\xi)\cup\{\pi^{-1}_{\xi+1}(\M_\xi)\}$.
\item If $M\models ``(\M_\xi, \Phi_\xi)$ is $\Gamma$-{\sf{cbl}}\ $+(X_\xi, \phi_\xi)\not \in {\sf{Hp}}^\Gamma"$ then all remaining objects with index $\geq \xi$ are undefined. 
\item If both 5.c.A and 5.c.B fail then $Y_{\xi+1}=\pi^{-1}_{\xi+1}(Y_\xi)$. 
\end{enumerate}
\end{enumerate}
\item Suppose $\xi\leq \d$ is a limit ordinal and for all $\gg<\xi$, both $\M_\gg$ and $\N_\gg$ are defined. Then $\M_{\xi}$ and $\N_{\xi}$ are determined as follows\footnote{The rest of the objects will be defined at the next stage of the induction as in clause 4.}. Set $\N_\xi=lim_{\a\rightarrow \xi}\M_\a$. Assuming clause 4 fails for $\xi+1$, the remaining objects 
     are defined as follows.
     \begin{enumerate}
     \item $\M_\xi=\C(\N_\xi)$ and
     \item $Y_\xi=\pi^{-1}_{\xi}(Y^{\N_\xi})$\footnote{$F_\xi$ and $b_\xi$ are defined at step $\xi+1$.}.
     \end{enumerate}
     \item $\M_\d=\N_\d$ and $Y_\d, \Phi_\d, F^+_\d, F_\d$, and $b_\d$ are undefined.
 \end{enumerate}
Let  \begin{center}${\sf{hpc}}=(\M_\gg , \N_\gg, Y_\gg, \Phi_\gg, F^+_\gg, F_\gg, b_\gg: \gg\leq \delta)$\end{center} be the output of the  $\Gamma-{\sf{hpc}}$ of $\mathbb{M}$. We say that the $\Gamma-{\sf{hpc}}$ of $\mathbb{M}$ is \textbf{successful} if  clause 4 fails for all $\gg<\d$. We say that the $\Gamma-{\sf{hpc}}$ of $\mathbb{M}$ \textbf{reaches its goal} if the $\Gamma-{\sf{hpc}}$ of $\mathbb{M}$ is successful and for some $\xi<\d$, clause 5.c.ii holds.

For each $\gg\leq \d$, we let $\Phi^+_\gg$ be the extension of $\Phi_\gg$ defined in \rsec{the induced strategy for the undropping game}. We then set
\begin{center}${\sf{hpc}}^+=(\M_\gg , \N_\gg, Y_\gg, \Phi^+_\gg, F^+_\gg, F_\gg, b_\gg: \gg\leq \delta)$.\end{center}
Notice that ${\sf{hpc}}\in M$ while ${\sf{hpc}}^+\not \in M$.

Also, given $\xi\leq \d$ and $\a\leq \d$, we set 
\begin{center}
${\sf{hpc}}\rest (\xi, \a)=(\M_\gg , \N_\gg, Y_\gg, \Phi_\gg\rest M|\a, F^+_\gg, F_\gg, b_\gg: \gg\leq \xi)$
\end{center}
and finally we let 
\begin{center}
${\sf{hpc}}^-=(\M_\gg , \N_\gg, Y_\gg, F^+_\gg, F_\gg, b_\gg: \gg\leq \delta)$
\end{center}
We will often use $({\sf{C}}, \Gamma)$ as a subscript to emphasize the dependence on $({\sf{C}}, \Gamma)$. Thus, we will write ${\sf{hpc}}_{{\sf{C}}, \Gamma}$ and etc. Also, to emphasize the dependence on ${\sf{C}}$, we may also say that $\Gamma-{\sf{hpc}}$ of ${\sf{C}}$ is successful or reaches its goal.

We say that $\Q$ is an $\N$-\textbf{model} of $\sf{hpc}$ if for some $\gg\leq \d$, $\Q=\N_\d$. We define other such expressions (e.g. $\M$-model and etc) in a similar fashion. We say $\W$ is \textbf{the last model} of  $\sf{hpc}$ if $\Q=\N_\gg$, the last defined $\N$-model of $\sf{hpc}$.$\myqedhere$
\end{definition}

\begin{remark}\label{gamma cbl} \rsubsec{internalizing gamma sets} defines the meaning of $M\models (X_\xi, \phi_\xi)\in {\sf{Hp}}^\Gamma$, which in reality formalizes $M\models (\M_\xi, \Phi_\xi)\in {\sf{Hp}}^\Gamma$ . Using very similar ideas, one can also easily formalize the meaning of $M\models ``(\M_\xi, \Phi_\xi)$ is $\Gamma$-{\sf{cbl}}". Such a formalism will refer to some set $Z_\xi$ and a formula $\psi_\xi$. The definition of these will be similar to the definitions of $X_\xi$ and $\phi_\xi$ (see \rdef{xgamma}). We leave the details to the reader. 

$\myqedhere$
\end{remark}

\begin{remark}\label{germane models} Notice that each $\M_\gg$ and $\N_\gg$ are germane (see \rdef{germane lses}), and so we can use the concepts introduced in \rsec{undropping game sec}. $\myqedhere$
\end{remark}

%

\subsection{The construction of $\Phi_\gg^+$}\label{the induced strategy for the undropping game}

We are continuing with the objects defined in \rsec{sec: hod pair constructions} and in particular, in \rdef{gamma-hod pair construction*}. Recall that $\Phi_\gg^+$ must be an $(\omega_1, \omega_1, \omega_1)$-iteration strategy. Its $(\omega_1, \omega_1)$ component can be defined using the procedure of \cite[Chapter 12]{FSIT}. Also $\Phi^+_\gg$ is the strategy of $\M_\gg$ that is based on $\cup Y^{\M_\gg}$. Thus, to define $\Phi_\gg^+$ we may just as well assume that $\cup Y^{\M_\gg}$ is a limit type hod premouse, as this is when an $(\omega_1, \omega_1, \omega_1)$-iteration strategy is used. As the process is a straightforward adaptation of  \cite[Chapter 12]{FSIT}, we will only give a short outline.

The procedure of \cite[Chapter 12]{FSIT} gives an $(\omega_1, \omega_1)$-strategy $\Phi_\gg$ for $\M_\gg$. Set $\P^+=\M_\gg$ and $\P=\cup Y^{\M_\gg}$. Suppose $N$ is a $\Sigma$-iterate of $M$ via $\X$ and $i:M\rightarrow N$ is the iteration embedding. Recall that we had
\begin{center}${\sf{hpc}}_{{\sf{C}}, \Gamma}=(\M_\gg , \N_\gg, Y_\gg, \Phi_\gg, F^+_\gg, F_\gg, b_\gg : \gg\leq \delta)$
\end{center}
and our background is $\mathbb{M}=(M, \d, \vec{G}, \Sigma)$. Suppose $\a<\lh(\X)$ and $\b\leq \pi^{\X}_{0, \a}(\d)$. We then let $\R^\X_{\a, \b}$ be the $\b$-th $\N$-model of $\pi^{\X}_{0, \a}({\sf{hpc}}_{{\sf{C}}, \Gamma})$, and also we let $\Phi^N_{\a, \b}$ be the $(\omega_1, \omega_1)$-iteration strategy of $\R_{\a, \b}^\X$ induced by $\Sigma_N$. We let $\Phi^N$ be the strategy of $i(\P^+)$ induced by $\Sigma_N$.

Given $N$ as above and a stack $\T$ on $i(\P^+)$ that is based on $i(\P)$ and is according to $\Phi^N$, we let $r\T$ be the \textit{resurrection of $\T$}. The reader may wish to review properties H1-H7 on page 113-115 of \cite{FSIT}, which outline the construction of $r\T$. Below we outline the description of  $\Phi^M$ and leave $\Phi^N$ to the reader. Assuming 
\begin{center}
$\T=((\M_\a)_{\a<\eta}, (E_\a)_{\a<\eta-1}, D, R, (\beta_\a, m_\a)_{\a\in R}, T)$\footnote{Recall that our stacks are proper, see \rdef{proper stacks convention}.}\\ and\\
$r\T=((r\M_\a)_{\a<\eta}, (rE_\a)_{\a<\eta-1}, rD, rR, (r\beta_\a, rm_\a)_{\a\in rR}, rT)$\footnote{Here we only use $(\omega_1, \omega_1)$-portion of $\Lambda_N$.} 
\end{center}
there are sequences $\vec{\sigma}=(\sigma_\a: \a<\eta)$ and $\vec{\nu}=(\nu_\a: \a<\eta)$ satisfying the following conditions:
\begin{enumerate}
\item $rD=\emptyset$, $T=rT$ and $R=rR$.
\item For each $\a<\eta$, $\nu_\a\leq \pi^{r\T}_{0, \a}(\gg)$ and $\sigma_\a:\M_\a\rightarrow \R^{r\T}_{\a, \nu_\a}$ is a weak embedding\footnote{For example, see the discussion after Fact 2.13 of \cite{ANS}.}.
\item If $[0, \a)_\T\cap D=\emptyset$ then $\nu_\a=\pi^{r\T}_{0, \a}(\gg)$ and $\R^{r\T}_{\a, \nu_\a}=\pi^{r\T}_{0, \a}(\P)$.
\item For each $\a<\a'$ such that $(\a, \a')\cap R=\emptyset$, $\sigma_\a\rest \ind^\T_\a=\sigma_{\a'}\rest \ind^\T_{\a'}$.
\item For each $\a<\a'$ such that $\a\T\a'$ and $\pi_{\a, \a'}^\T$ is defined, $\pi^{r\T}_{\a, \a'}\circ \sigma_\a =\sigma_{\a'}\circ \pi_{\a, \a'}^\T$.
\item Moreover, $(\sigma_\a: \a<\eta)$ and $(\nu_\a: \a<\eta)$ are uniquely determined via the procedure described on pages 113-115 of \cite{FSIT}. 
\end{enumerate}
We then say that $(\vec{\sigma}, \vec{\nu})$ are the $r\T$-sequences. 

\begin{definition}\label{res and emb}
 Suppose now that $p=(\P_\b, \T_\b, E_\b: \b< \gg)$ is a generalized stack on $\P^+$ that is based on $\P$. We say $p$ is \textbf{correct} if there is a stack $q=((\Q_\a)_{\a<\eta}, (F_\a)_{\a<\eta-1}, D, R, (\beta_\a, m_\a)_{\a\in R}, Q)$ according to $\Sigma$ and a sequence of embeddings $(\sigma_\b: \b<\gg)$ such that the following conditions hold:
 \begin{enumerate}
 \item $\eta=\Sigma_{\b<\gg}\lh(\T_\b)$ and $\eta_\b=_{def}\Sigma_{\b'<\b}\lh(\T_{\b'})$.
 \item For all $\b<\gg$, $\sigma_\b:\P_\b\rightarrow \pi^q_{0, \eta_\b}(\P^+)$ is a weak embedding.
 \item  $\sigma_0=id$.
 \item For all $\b<\gg$, $\eta_\b\in R$.
 \item For all $\b<\gg$, $q_{[\eta_\b, \eta_{\b+1})}=r(\sigma_\b\T_\b)$.
 \item For all $\b<\gg$ such that $\b+1<\gg$ and $E_\b$ is an un-dropping extender\footnote{The case when $\pi^{\T_\b, b}$ is defined is easier and very similar, and we leave it to the reader.}, letting 
 \begin{enumerate}
 \item $md^{\T_\b}=(\a_i, \R_i, \W_i, \S_i : i\leq k+1)$ be the main drops of $\T_\b$,
 \item for $i\leq k+1$, $\k_i=\d^{\R_i^b}$\footnote{See \rdef{the un-dropping extender of a continuable stack}.},
 \item $\xi+1=\lh(\T_\b)$,
 \item $m:\M_{\xi}^{\T_\b}\rightarrow \M_{\xi}^{\sigma_\b\T_\b}$ is the map obtained via the copying process,
 \item $(\vec{n}, \vec{\nu})$ are the $r(\sigma_\b\T_\b)$-sequences,
 \item $k=\sigma^\T$,
 \end{enumerate}
the embedding 
  \begin{center} $\sigma_{\b+1}: \P_{\b+1}\rightarrow \pi^q_{0, \eta_{\b+1}}(\P^+)$\end{center} is given by 
  \begin{center}
 $\sigma_{\b+1}(\pi_{E_\b}(f)(a))=\pi_{\eta_\b, \eta_{\b+1}}^q(\sigma_\b(f))(n_\xi(m(a)))$
 \end{center}
 \end{enumerate}
The definition of $\sigma_{\b+1}$ works because we have that  $(a, A)\in E_\b$ if and only if $n_\xi\circ m(a)\in \pi^q_{\eta_\b, \eta_{\b+1}}(\sigma_\b(A))$.

  Notice that both $q$ and the embeddings $\vec{\sigma}=(\sigma_\b: \b<\gg)$ are uniquely determined. We then set $q={\sf{res}}(p)$ and $\vec{\sigma}={\sf{emb}}(p)$.  $\myqedhere$
 \end{definition}
 The following is an easy lemma. It uses the objects introduced above.
 \begin{lemma}\label{sigmabeta} $\sigma_{\b+1}\rest (\P_{\b+1}|\lh(E_\b))=n_\xi\circ m \rest (\P_{\b+1}|\lh(E_\b))$.
 \end{lemma}
 
 It is now straightforward to show, using the resurrection process of \cite[Chapter 12]{FSIT}, that if $p$ is a correct generalized stack on $\P^+$ based on $\P$ of limit length then there is a unique branch $b$ of $p$ such that $p^\frown \{b\}$ is also correct. We then let $\Phi_\gg^+$ be the unique $(\omega_1, \omega_1, \omega_1)$-strategy of $\P^+$ with the property that $p$ is according to $\Phi^+_\gg$ if and only if $p$ is a correct generalized stack on $\P^+$ based on $\P$. Notice finally that the definition of $\Phi_\gg^+$ can be done locally inside $M$.

 \begin{definition}\label{gammastrategy} If $\cup Y^{\M_\gg}$ is not of $\#$-lsa type then $\Phi_\gg=\Phi_\gg^+\rest M|(\d^+)^M$. If $\cup Y^{\M_\gg}$ is of $\#$-lsa type then $\Phi_\gg=(\Phi_\gg^+)^{stc}\rest M|(\d^+)^M$. $\myqedhere$
 \end{definition}
 The following lemma summarizes \rdef{res and emb}. 
 \begin{lemma}\label{summary res and emb} Suppose $\Gamma$ is a pointclass,  ${\sf{C}}=(\mathbb{M}, (P, \Psi), \Gamma^*, A)$ Suslin, co-Suslin captures $\Gamma$ and $\mathbb{M}=(M, \d,  \vec{G}, \Sigma)$. Set
\begin{center}${\sf{hpc}}^+_{{\sf{C}}, \Gamma}=(\M_\gg , \N_\gg, Y_\gg, \Phi^+_\gg, F^+_\gg, F_\gg, b_\gg : \gg\leq \delta)$.
\end{center}
 Suppose $\gg\leq \d$ is such that $Y_\gg$ is defined. Set $\P=\cup Y_\gg$ and suppose $\M_\gg$ is of $b$-type. Suppose $\T$ is a generalized stack according to $\Phi^+_\gg$ with last model $\Q$. There is then a $\Sigma$-iterate $N$ of $M$ such that letting $i: M\rightarrow N$ be the iteration embedding and 
 \begin{center}${\sf{hpc}}^+_{{\sf{C_N}}, \Gamma}=(\R_\gg , \S_\gg, Z_\gg, \Psi^+_\gg, E^+_\gg, E_\gg, c_\gg : \gg\leq i(\delta))$,
\end{center} 
there is $\nu\leq i(\gg)$ and a weak embedding $\sigma:\Q\rightarrow  \S_\nu$ such that the following holds.
 \begin{enumerate}
 \item If $\pi^\T$ is defined then $\nu=i(\gg)$ and $\sigma\circ \pi^\T=i\rest \P$.
 \item If $\pi^{\T, b}$ is defined then $i(\P^b)=\S_\nu^b$ and $\sigma\circ \pi^{\T, b}=i\rest \P^b$.
 \item $\Phi^+_{\Q, \T}$ is the $\sigma$-pullback of $\Psi_{\nu}^+$.\\
 \end{enumerate}
 \end{lemma}
 
 We remark that a similar result holds for all $\gg$. We now have the following lemma connecting different strategies to each other.
 
 \begin{lemma}\label{different strategies agree with each other} Suppose $\Gamma$ is a pointclass,  ${\sf{C}}=(\mathbb{M}, (P, \Psi), \Gamma^*, A)$ Suslin, co-Suslin captures $\Gamma$ and $\mathbb{M}=(M, \d,  \vec{G}, \Sigma)$. Set
\begin{center}${\sf{hpc}}^+_{{\sf{C}}, \Gamma}=(\M_\gg , \N_\gg, Y_\gg, \Phi^+_\gg, F^+_\gg, F_\gg, b_\gg : \gg\leq \delta)$.
\end{center}
 Suppose $\a<\b\leq \d$ are such that $\N_{\a}$ and $\N_\b$ are defined. Let $\Q\in \Y_\a$ be a meek hod premouse. Set $\Psi^0=\Phi^+_{\a+1}$, $\P_1=\M_{\a+1}$ and define $\Psi^1$ as follows:
 \begin{itemize}
 \item If $\rho(\N_\b)\leq \d^\Q$ then let $n$ be the largest such that for every $\k<\d^\Q$, any $r\Sigma_n^{\N_\b}$-definable $f:\k\rightarrow \d^\Q$ is in $\Q$ and let $\Psi^1$ be the strategy of $\P_1=_{def}{\sf{core}}_n(\N_\b)$ defined via the resurrection procedure described above.
 \item If $\rho(\N_\b)>\d^\Q$ then let $\Psi^1$ be the strategy of $\P_1=_{def}{\sf{core}}(\N_\b)$  defined via the resurrection procedure described above.
 \end{itemize}
 Then $\Psi^0_\Q=\Psi^1_\Q$.
 \end{lemma}
 The proof of the lemma is straightforward. Let $\gg$ be such that $\M_\gg=\Q$ and let $\zeta=\sup\{\lh(F_\iota^+): \iota<\gg\}$. Observe now that because we assume that $\Q$ is meek, if $\T$ is a stack on $\Q$ then the $id$-copy of $\T$ onto $\P_0$ and onto $\P_1$ is simply $\T_0=_{def}\uparrow(\T, \P_0)$ and $\T_1=_{def}\uparrow(\T, \P_1)$ respectively, and these stacks use exactly the same extenders as $\T$. Therefore the resurrection procedure resurrects both $\T_0$ and $\T_1$ to stacks based on $M|\zeta$. Hence both $\Psi^0_\Q$ and $\Psi^1_\Q$ are determined by $\Sigma_{M|\zeta}$. 
 
 Lastly we state the following consequence of \rlem{sigmabeta}.
 \begin{definition}\label{weakly selfcohering} Suppose $(\P, \Sigma)$ is a hod pair or an sts pair. We say $\Sigma$ is \textbf{weakly self-cohering} if the following clauses hold:
 \begin{enumerate}
 \item Whenever $\T$ is a generalized stack on $\P$ according to $\Sigma$ with last model $\S$ and $\Q$ is a complete layer of $\S^b$ such that $\T^{\sf{ue}}_\Q$ is defined\footnote{See \rnot{notation for generalized stacks}.}, $\Sigma_{\Q, \T}=\Sigma_{\Q, \T^{{\sf{ue}}}}$. 
 \item Whenever $\T$ is a generalized stack on $\P$ according to $\Sigma$ with last model $\S$ and such that $\T$ has a one point extension\footnote{See \rdef{one point extension}.}, $\Q\insegeq_{hod} \S^b$ is of limit type, and $\U$ is a stack on $\Q$ according to $\Sigma_{\Q, \T}$ such that $\U$ has a one point extension then letting $E$ be the un-dropping extender of $\U$, $Ult(\S, E)$ is well-founded.
\end{enumerate}

Suppose next that $(\P, \Sigma)$ is a simple hod pair or an sts hod pair. Then we say that $\Sigma$ is \textbf{weakly self cohering} if the following clauses hold:
 \begin{enumerate}
 \item Whenever $\T$ is a stack on $\P$ according to $\Sigma$ with last model $\S$ and $\Q$ is a complete layer of $\S^b$ such that $\T^{\sf{ue}}_\Q$ is defined, the last model of $\T^{\sf{ue}}_\Q$ is well-founded. 
 \item Whenever $\T$ is a stack on $\P$ according to $\Sigma$ with last model $\S$ and such that $\T$ has a one point extension, $\Q\insegeq_{hod} \S^b$ is of limit type, and $\U$ is a stack on $\Q$ according to $\Sigma_{\Q, \T}$ such that $\U$ has a one point extension then letting $E$ be the un-dropping extender of $\U$, $Ult(\S, E)$ is well-founded.
\end{enumerate}
$\myqedhere$
 \end{definition}
 
The following now is an easy consequence of \rlem{sigmabeta}.
\begin{lemma}\label{weakly self cohering} Suppose $\Gamma$ is a pointclass,  ${\sf{C}}=(\mathbb{M}, (P, \Psi), \Gamma^*, A)$ Suslin, co-Suslin captures $\Gamma$ and $\mathbb{M}=(M, \d,  \vec{G}, \Sigma)$. Set
\begin{center}${\sf{hpc}}^+_{{\sf{C}}, \Gamma}=(\M_\gg , \N_\gg, Y_\gg, \Phi^+_\gg, F^+_\gg, F_\gg, b_\gg : \gg\leq \delta)$.
\end{center}
 Suppose $\gg\leq \d$ is such that $Y_\gg$ is defined. Set $\P=\cup Y_\gg$ and suppose $\M_\gg$ is of $b$-type. Then $(\Phi_\gg)_{\P}$ is weakly self-cohering.
\end{lemma}
 
\subsection{The definition of $(X_\gg, \phi_\gg)$.}
 
 Notice that the map $\T\mapsto {\sf{res}}(\T)$ can be defined without any reference to any strategy for $\P$ or $M$. In this view, ${\sf{res}}(\T)$ may not have well-founded models. Moreover, the construction of ${\sf{res}}(\T)$ only depends on ${\sf{hpc}}\rest \gg+1$. 
 \begin{definition}\label{xgamma}
 Let $\xi_\gg<\d$ be the least inaccessible cardinal of $M$ such that
 \begin{center}
  ${\sf{hpc}}^-\rest \gg+1\in M|\xi_\gg$. 
  \end{center}
  Let $\dot{X}_\gg\in M^{Coll(\omega, M|\xi_\gg)}$ witness that $M$ is self-capturing for $M|\xi_\gg$ (see \rdef{self-capturing bt}) and set 
 \begin{center}
 $X_\gg=(\dot{X}_\gg,  \M_\gg, {\sf{hp}c}\rest (\gg+1, \xi_\gg), M|\xi_\gg)$. 
 \end{center}
 Let $\psi(x, y, z, w)$ be a formula such that
 \begin{center}
  $\psi[\M_\gg, {\sf{hpc}}\rest (\gg+1, \xi_\gg), M|\xi_\gg]$ 
  \end{center}
  expresses all the clauses of \rdef{gamma-hod pair construction*} except the portion of clause 5.c that defines $Y_{\xi+1}$. Let $\phi_{\gg}(u, v, w)$ be the conjunction of the following formulas.
 \begin{enumerate}
 \item $u=(Y, g)$ such that  $Y=(Z, Q, h, f, N)$, $Z$ is  $Coll(\omega, N)$-name and $g\subseteq Coll(\omega, N)$ is a filter,
 \item $\psi(Q, h, f, N)$,
 \item $w$ is a stack on $Q$, and
 \item letting $Z_g=(U, W)$, ${\sf{res}}(w)\in p[U]$. 
 \end{enumerate} 
 $\myqedhere$
 \end{definition}

\section{On backgrounded constructions}

The following sequence of lemmas will be used in the proof of \rthm{existence of thick sets}.

\begin{definition}\label{property *}
We say $(M, \d, \vec{G}, \Sigma, \P)$ has the property $(*)$ if
\begin{itemize}
\item $(M, \d, \vec{G})$ is a background\footnote{See \rdef{background}.}, 
\item $\Sigma\in M$ is a $(\d, \d)$-iteration or $\d$-iteration strategy for $\P$ with hull condensation\footnote{The exact nature of $\P$ is irrelevant.}, or $(\P, \Sigma)$ is an sts hod pair and $\Sigma$ is a $(\d, \d, \d)$ st-strategy. 
\item ${\sf{Le}}(\Sigma, \mathcal{J}_{\omega}[\P])^{(M, \d, \vec{G})}$ is successful. 
\end{itemize}
We say that $(M, \d, \vec{G}, \Sigma, \P)$ has the property $(*+)$ if in addition to the above clauses $M\models ``\Sigma$ is a $(\d^+, \d^+)$-iteration strategy, $\d^+$-iteration strategy or $(\d^+, \d^+, \d^+)$ st-strategy".
If $Q\subseteq M$ then we let $\Sigma^Q=\Sigma\rest (Q|\d)$. $\myqedhere$
\end{definition}

\begin{lemma}\label{le coherence} Assume $(M, \d, \vec{G}, \Sigma, \P)$ has the property $(*)$. Set $P= \mathcal{J}_{\omega}[\P]$. Suppose $\l<\d$ is such that $P\in M|\l$ and $\N$ is the last model of ${\sf{Le}}(\Sigma^M, P)^{(M, \d, \vec{G})}$. Suppose $F^*\in \vec{G}$ is such that
\begin{enumerate}
\item $\lh(F^*)=\eta$ is an inaccessible cardinal of $M$, 
\item $\pi_{F^*}(\N)|\eta=\N|\eta$.
\end{enumerate}
Set $\kappa=\cp(F^*)$ and let $F'$ be the $(\k, \eta)$ extender derived from\begin{center} $\pi_{F^*}\rest \N:\N\rightarrow \pi_{F^*}(\N)$.\end{center} Then for any $\rho\in [(\k^+)^\N, \eta)$ such that $\rho$ is the natural length of $F'\rest \rho$, letting $F$ be the trivial completion of $F'\rest \rho$, one of the following conditions hold:
\begin{enumerate}
\item $\lh(F)\in \dom(\vec{E}^\N)$ and $F=\vec{E}^\N(\lh(F))$ or
\item $\lh(F)\not\in \dom(\vec{E}^\N)$, $\rho$ is a limit ordinal $>(\k^+)^\N$, $\rho$ is a generator of $F$, $\rho\in \dom(\vec{E}^\N)$ and letting $E=\vec{E}^{\N}(\rho)$, $F=\pi_E^{\N|\rho}(\vec{E}^{\N|\rho})(\lh(F))$.
\end{enumerate}
\end{lemma}

Suppose $(M, \d, \vec{G})$ is a background. We write $(M, \vec{G})\models ``\k$ reflects $A$" to mean that $\k$ reflects $A$ using extenders in $\vec{G}$\footnote{I.e., the set of $\k$ such that for every $\l<\d$ there is $F\in \vec{G}$ such that $\pi_F(\vec{G})\rest \l=\vec{G}\rest \l$.}. Working in $M$, let $(A^M_i : i<\omega)$ be defined by the following induction:
\begin{enumerate}
\item $A^M_0\subseteq \d$ is the set of $<\d$-strong cardinals $\k$ such that $(M, \vec{G})\models ``\k$ reflects $\vec{G}$. 
\item $A^M_{i+1}\subseteq \d$ is the set of $<\d$-strong cardinals $\k$ such that  $(M, \vec{G})\models ``\k$ reflects $A^M_i"$.
\end{enumerate}
Clearly, $A_i^M$ depends on both $\d$ and $\vec{G}$, but in all of the lemmas below $(\d, \vec{G})$ is clear from the context. Sometimes, when $(\d, \vec{G})$ is not clear from the context, we will write $A_i^{M, \d, \vec{G}}$.  We have the following straightforward lemma.

\begin{lemma}\label{a sets} Suppose $(M, \d, \vec{G})$ is a background. Then the following holds in $M$. 
\begin{enumerate}
\item Suppose $A\subseteq \d$ and $\k_0<\k_1<\d$ are such that $(M, \vec{G})\models ``\k_1$ reflects $(A, \vec{G})"$ and $(M|\k_1, \vec{G}\rest \k)\models ``\k_0$ reflects $(A\cap \k_1, \vec{G}\rest \k_1)"$. Then $(M, \vec{G})\models ``\k_0$ reflects $(A, \vec{G})"$.
\item For each $i<\omega$, if $\l\in A_{i+1}^M$ or is a limit point of $A_{i+1}^M$ then $\l$ is a limit point of $A_i^M$.
\item For all $i<\omega$ and for every $\l$, which is a member or a limit point of $A_{i}^M$, $A_i^M\cap \l=A_i^{M|\l,\vec{G}\rest \l}$. 
\item For all $i\in \omega$, $A_{i+1}^M\subseteq A_i^M$.
\item $\k\in \cap_{i\in \omega} A_i^M$ if and only if for each $i\in \omega$, $(M, \vec{G})\models ``\k$ reflects $A^M_i"$. Hence, $\cap_{i\in \omega} A^M_i\not =\emptyset$.
\item If $(M, \vec{G})\models ``\k<\d$ reflects $(A_i: i\in \omega)"$ then $\k\in \cap_{i<\omega} A_i^M$.
\end{enumerate}
\end{lemma}

\begin{lemma}\label{extenders go up and down} Assume $(M, \d, \vec{G}, \Sigma, \P)$ has the property $(*)$. Set $P= \mathcal{J}_{\omega}[\P]$. Suppose $\l<\d$ is such that $P\in M|\l$, $\N'$ is the last model of ${\sf{Le}}(\Sigma^M, P)^{(M, \d, \vec{G})}$ and $\N=L_{\ord(M)}[\N']$\footnote{$\d$ is a Woodin cardinal of $\N$ and all bounded subsets of $\d$ in $\N$ are in $\N'$. The first claim can be shown by the results of \cite[Chapter 11]{FSIT}, and the second follows from the fact that $\d$ is a regular cardinal, which allows us to take Skolem hulls of M that are transitive below $\d$.}. Let $\vec{H}=\{E\in\vec{E}^\N: \N\models ``\nu(E)$ is inaccessible$"\}$. Then \begin{center}
 $\cap_{i<\omega} A_i^M=\cap_{i<\omega} A_i^{\N, \vec{H}}$.
\end{center}
\end{lemma}
\begin{proof} We will use $A_n^\N$ for $A_n^{\N, \vec{H}}$. It is enough to show that $i<\omega$, $A_{i+1}^M\subseteq A_i^\N\subseteq A_i^M$. Notice first that\\\\
(1) in $M$, if $\k\in A_0^M-(\l+1)$ and $\Q$ is an $\N$-model of 
${\sf{Le}}(\Sigma^M, P)^{(M, \d, \vec{G})}$ such that $\Q\in M|\k$ then $\Q$ is an $\N$-model of 
${\sf{Le}}(\Sigma^M, P)^{(M, \k, \vec{G}\rest \k)}$.\\\\
(1) then easily implies that\\\\
(2) if $\k\in A_0^M-(\l+1)$ then $\N|\k$ is the last model of ${\sf{Le}}(\Sigma^M, P)^{(M, \k, \vec{G}\rest \k)}$\footnote{This is a mild abuse of our notation as $\k$ may not be a Woodin cardinal of $M$. But $\sf{Le}$ construction do not depend on the Woodinness of $\d$.} and the $\k$th $\N$-model of ${\sf{Le}}(\Sigma^M, P)^{(M, \d, \vec{G})}$.\\\\
 $A_0^\N\subseteq A_0^M$ follows from the fact that the backgrounding extenders used in 
 \begin{center}${\sf{Le}}(\Sigma^M, P)^{(M, \d, \vec{G})}$\end{center} are all total $M$-extenders.

Suppose now that $\k\in A_1^M$. We want to see that $\k\in A_0^\N$. Let $\eta_0<\eta_1$ be two members of $A_0^M$ such that $\k<\eta_0$. Let $F^*\in \vec{G}$ be an extender such that $\pi_{F^*}(A_0^M)\cap \eta_1+1=A_0^M\cap \eta_1+1$ and $\pi_{F^*}(\vec{G})\rest \eta_1=\vec{G}\rest \eta_1$. (2) then implies that\\\\
(3) $\pi_{F^*}(\N)|\eta_1=\N|\eta_1$.\\\\
Indeed, it follows from (2) that $\pi_{F^*}(\N)|\eta_1$ is the last model of 
\begin{center}
$({\sf{Le}}(\Sigma^{Ult(M, F^*)}, P)_{>\l})^{(Ult(M, F^*)|\eta_1, \eta_1, \pi_{F^*}(\vec{G})\rest \eta_1)}$,
\end{center}
 and since $Ult(M, F^*)|\eta_1=M|\eta_1$, we have that $\pi_{F^*}(\N)|\eta_1$ is the last model of ${\sf{Le}}(\Sigma^M, P)^{(M|\eta_1, \eta_1, \vec{G}\rest \eta_1)}$, which according to (2) is just $\N|\eta_1$. 
 
Let now $F$ be the $(\k, \eta_1)$ extender derived from $\pi_{F^*}\rest \N$. Since $\eta_0$ is a regular cardinal of $\N$ and hence, $\eta_0\not \in \dom(\vec{E}^\N)$, it follows from \rlem{le coherence} that the trivial completion of $F\rest \eta_0$ is on $ \vec{E}^\N$. As $\eta_0$ was arbitrary, we have that $\d=\sup\{ \lh(E): E\in \vec{E}^\N \wedge \cp(E)=\k\}$, implying that $\k\in A_0^\N$.

Assume now that $A_{n+1}^M\subseteq A_n^\N\subseteq A_n^M$. We want to see that \\\\
(a) $A_{n+1}^\N\subseteq A_{n+1}^M$\\
(b) $A_{n+2}^M\subseteq A_{n+1}^\N$.\\\\
First suppose $\k\in A_{n+1}^\N$. To show that $\k\in A_{n+1}^M$, we need to show that in $M$, $\k$ reflects $A_n^M$.
Let $\eta\in A_0^\N$ be a limit point of $A_n^\N$. Let $E\in \vec{H}$ be a $(\k, \eta)$-extender that reflects $A_n^\N$\footnote{By this we mean an extender whose natural length is $\eta$. As $\eta$ is a regular cardinal of $\N$, there are no $(\k, \eta)$-extenders on the sequence of $\N$.}. Thus, $\pi_E(A_n^\N)\cap \eta=A_n^\N\cap \eta$. Let $E^*\in \vec{G}$ be the resurrection of $E$. We then have that $E=E^*\cap (\eta^{<\omega}\times \N)$ and an embedding $\sigma: Ult(\N, E)\rightarrow \pi_{E^*}(\N)$ such that $\cp(\sigma)\geq \eta$. Because $A_n^\N\subseteq A_n^M$ we have that\\\\
(5) $A_n^{\pi_{E^*}(\N)}\subseteq A_n^{Ult(M, E^*)}$,\\
(6) $\sigma(A_n^\N\cap \eta)=\sigma(A_n^{Ult(\N, E)}\cap \eta)=A_n^{\pi_{E^*}(\N)}\cap \sigma(\eta)\subseteq A_n^{Ult(M, E^*)}$. \\\\
Since $\sigma[A_n^\N\cap \eta]=A_n^\N\cap \eta$, it follows that \\\\
(7) $\pi_{E^*}(A_n^M)\cap \eta$ is cofinal in $\eta$.\\\\
It then follows that $\pi_E^*(A_n^M)\cap \eta=A_n\cap \eta$ (see \rlem{a sets}). Thus, $\k\in A_{n+1}^M$.

Finally suppose that $\k\in A_{n+2}^M$. We want to see that $\k\in A_{n+1}^\N$. Let $\eta_0<\eta_1$ be two members of $A_{n+1}^M$ such that $\eta_0$ is a limit of $A_{n+1}^M$ and $\k<\eta_0$. Let $F^*$ be such that $\pi_{F^*}(A_{n+1}^M)\cap \eta_1+1=A_{n+1}^M\cap \eta_1+1$. Like in the $n=0$ case, we have that if $F'$ is the $(\k, \eta_1)$-extender derived from $\pi_{F^*}\rest \N$ and $F$ is the trivial completion of $F'\rest \eta_0$ then $F\in \vec{E}^\N$. Let now $\sigma:Ult(\N, F)\rightarrow \pi_{F^*}(\N)$ be the canonical factor map. We have that $\cp(\sigma)\geq \eta_0$. We also have that $A_{n+1}^M\cap \eta_1\subseteq A_n^\N\cap \eta_1$. Arguing as above, we get that, in $\N$, $F$ reflects $A_n^\N$.  
\end{proof}

\begin{lemma}\label{extenders go up and down1} Assume $(M, \d, \vec{G}, \Sigma, \P)$ has the property $(*)$. Set $P= \mathcal{J}_{\omega}[\P]$. Suppose $\l<\d$ is such that $P\in M|\l$, $\N'$ is the last model of ${\sf{Le}}(\Sigma^M, P)^{(M, \d, \vec{G})}$ and $\N=L_{\ord(M)}[\N']$. Suppose $F^*\in \vec{G}$ is such that
\begin{enumerate}
\item $\lh(F^*)=\eta\in A_1^M$, 
\item $\pi_{F^*}(\N)|\eta=\N|\eta$.
\end{enumerate}
Set $\kappa=\cp(F^*)$ and let $F'$ be the $(\k, \eta)$-extender derived from $\pi_{F^*}\rest \N:\N\rightarrow \pi_{F^*}(\N)$. Let $F$ be the trivial completion of $F'\rest \eta$. Then $F\in \vec{E}^\N$.
\end{lemma}
\begin{proof}  Let $\eta'\in A_1^M-(\eta+1)$ and let $H\in \vec{G}$ be an extender such that $\cp(H)=\eta$, $\lh(H)>\eta'$, and $\pi_H(A_0^M)\cap (\eta'+1)=A_0^M\cap \eta'+1$. We have that\\\\
(1) $\N|\eta'$ is the last model of both 
\begin{center}
$({\sf{Le}}(\Sigma^M, P)_{>\l})^{(M|\eta', \eta', \vec{G}\rest \eta')}$ and $({\sf{Le}}(\Sigma^{Ult(M, H)}, P)_{>\l})^{(Ult(M, H)|\eta', \eta', \vec{G}\rest \eta')}$,
\end{center}
(2) $\pi_H(\N)|\eta'=\N|\eta'$. \\\\
It follows from \rlem{le coherence} that all initial segments of $F$ are on the sequence of $\N$ or an ultrapower away. Thus, in $Ult(M, H)$, we have that all initial segments of $\pi_H(F)$ are on the extender sequence of $\pi_H(\N)$ or an ultrapower away. As $F\rest \eta=\pi_H(F)\rest \eta$, we have that in $Ult(M, H)$, the trivial completion of $F\rest \eta$ is on the sequence of $\pi_H(\N)$. But the trivial completion of $F\rest \eta$ both in $M$ and in $Ult(M, H)$ is $F$, as it only depends on $\pi_{F|\eta}(\N|(\k^+)^\N)$ which is computed the same way in both models. Thus, $F$ is on the extender sequence of $\pi_H(\N)$. Since $\pi_H(\N)|(\eta^+)^{\pi_H(\N)}=\N|(\eta^+)^\N$, we have that $F$ is on the extender sequence of $\N$.
\end{proof}

\begin{definition}\label{core le sequence} Suppose $(M, \d, \vec{G})$ is a background. Let $\vec{K}^M$ consist of all extenders $E\in \vec{G}$ such that 
\begin{itemize}
\item $\nu(E)\in \cap_{i<\omega} A_i^M$ and is a limit point of $\cap_{i<\omega} A_i^M$,
\item $E$ reflects $(A_i^M: i<\omega)$.
\end{itemize}
We say $\S$ is the \textbf{fully backgrounded $\l$-core} of $(M, \d, \vec{G})$ if $\S$ is the last model of ${\sf{Le}}^{(M, \d, \vec{K}^M)}_{>\l}$. We let ${\sf{LeCore}}^{(M, \d, \vec{G})}_{>\l}$ be the fully backgrounded $\l$-core of $(M, \d, \vec{G})$. $\myqedhere$
\end{definition}
Clearly $(M, \d, \vec{K}^M)$ is a background.

\begin{lemma}\label{core l e} Assume $(M, \d, \vec{G}, \Sigma, \P)$ has the property $(*)$. Set $P= \mathcal{J}_{\omega}[\P]$. Suppose $\l<\d$ is such that $P\in M|\l$, $\R'$ is the last model of ${\sf{Le}}(\Sigma^M, P)_{>\l}^{(M, \d, \vec{G})}$ and $\R=L_{\ord(M)}[\R']$. Then ${\sf{LeCore}}^{(M, \d, \vec{G})}$ is the last model of ${\sf{Le}}^{(\R, \d, \vec{K}^\R)}_{>\l}$ where $\vec{K}^\R$ is computed relative to $\vec{H}^\R=\{ E\in \vec{E}^\R: \R\models ``\nu(E)$ is an inaccessible cardinal$"\}$.
\end{lemma}
\begin{proof} It is enough to show that if $\Q$ is an $\M$-model of both ${\sf{Le}}^{(M, \d, \vec{K}^M)}_{>\l}$ and ${\sf{Le}}^{(\R, \d, \vec{K}^\R)}_{>\l}$ then the $\N$-models of ${\sf{Le}}^{(M, \d, \vec{K}^M)}_{>\l}$ and ${\sf{Le}}^{(\R, \d, \vec{K}^\R)}_{>\l}$ constructed immediately after $\Q$ coincide. Assume then the $\N$-model of ${\sf{Le}}^{(M, \d, \vec{K}^M)}_{>\l}$ constructed immediately after $\Q$ is $\Q'$. The only non-trivial case is when $\Q'$ is obtained by adding an extender to $\Q$. Thus, assume $\Q'=(\Q, F)$. We need to see that $(\Q, F)$ is the $\N$-model of ${\sf{Le}}^{(\N, \d, \vec{K}^\R)}_{>\l}$ constructed immediately after $\Q$. Let $F^*$ be the background extender of $F$. It follows that 
\begin{itemize}
\item $\nu(F)< \nu(F^*)$ and $F^*\in \vec{K}^M$,
\item $\nu(F^*)\in \cap A_i^M$ and is a limit point of $\cap_{i<\omega} A_i^M$,
\item $F^*$ reflects $(A_i^M: i<\omega)$.
\end{itemize}
Set $\eta=\lh(F^*)$ and let $F'$ be the $(\k, \eta)$-extender derived from $\pi_{F^*}\rest \R$. Let $E$ be the trivial completion of $F'|\eta$. It follows from \rlem{extenders go up and down1} that in fact $E\in \vec{E}^\R$ and it also follows from \rlem{extenders go up and down} that $E\in \vec{K}^\R$. Since $F=F^*\cap (\nu(F)^{<\omega}\times \univ{\Q})$, we have that $\Q\subseteq \R$ and  in $\R$, $E$ is a background certificate of $F$. It then follows from the uniqueness of the next extender (see \cite[Chapter 9]{FSIT} and \cite[Theorem 5.1]{FarmDef}\footnote{This reference contains the proof of non-existence of mixed bicephali, completing \cite[Chapter 9]{FSIT}.}) that in fact that $\Q'$ is the $\N$-model of ${\sf{Le}}^{(\R, \d, \vec{K}^\R)}_{>\l}$ constructed immediately after $\Q$. 

Conversely, suppose  the $\N$-model of ${\sf{Le}}^{(\R, \d, \vec{K}^\R)}_{>\l}$ constructed immediately after $\Q$ is $(\Q, F)$ and let $F^*\in \vec{K}^\R$ be the background extender of $F$. We then have that 
\begin{itemize}
\item $\nu(F)< \nu(F^*)$,
\item $\nu(F^*)\in \cap_{i<\omega} A_i^\R$ and is a limit point of $\cap_{i<\omega} A_i^\R$,
\item $F^*$ reflects $(A_i^\R: i<\omega)$.
\end{itemize}
It then follows from \rlem{extenders go up and down} that letting $F^{**}$ be the background extender of $F^{*}$ and $E=F^{**}|\lh(F^*)$, $E\in \vec{K}^M$ and $E$ backgrounds $F$. It then follows from the uniqueness of the next extender (see the above references) that $(\Q, F)$ is indeed the $\N$-model of ${\sf{Le}}(\Sigma^M, P)_{>\l}^{(M, \d, \vec{G})}$ constructed immediately after $\Q$.
\end{proof}

\begin{corollary}\label{the core is a class} Suppose $(M, \d, \vec{G})$ is a background and $\l<\d$. Then for any $(\P, \Sigma)$ such that $(M, \d, \vec{G}, \Sigma, \P)$ has the property (*) and $\mathcal{J}_{\omega}[\P]\in M|\l$, letting $\R$ be the last model of ${\sf{Le}}(\Sigma^M, P)_{>\l}^{(M, \d, \vec{G})}$, ${\sf{LeCore}}^{(M, \d, \vec{G})}_{>\l}$ is a definable class of $\R$.
\end{corollary}

\section{On the existence of thick sets}\label{sec: thick hulls}

\cite[Chapter 5.1]{ATHM} develops a methodology for proving branch condensation and various uniqueness results for iteration strategies. The basic idea, due to Jensen\footnote{Jensen developed similar ideas for the $K^c$ constructions, see \cite{JSSS}.} and Steel\footnote{The first author learnt about the main idea behind \cite[Chapter 5.1]{ATHM} from Steel sometime between 2004-2006. To the author's best knowledge  \cite[Chapter 5.1]{ATHM} is the first written account of this material.}, is that the \textit{stack} over a fully backgrounded construction has covering properties. However, both \cite{ATHM} and our current exposition needs, in addition, that \textit{thick sets} exist. While \cite{ATHM} uses their existence, it seems that \cite{ATHM} does not establish their existence. In this section, we take a moment to fill this gap.

First we import one important definition from \cite[Chapter 5.1]{ATHM}. Recall that if $M$ is a transitive set then we let $M|\a$ be $V_\a^M$.

\begin{definition}\label{the stack in scb} Suppose $\k$ is a regular cardinal, $\Sigma$ is a $\kappa^+$-iteration strategy\footnote{The nature of the structure that $\Sigma$ is a strategy of is not important.} and $\M$ is a $\Sigma$-premouse (possibly over some set $X$) such that $\M\subseteq H_\kappa$. We let ${\sf{stack}}(\M, \Sigma)$ be the union of all sound countably iterable $\Sigma$-premice $\S$ such that $\M\insegeq \S$ and $\rho(\S)=\k$. 

If the stack is computed inside an inner model $M$ then to emphasize the dependence on $M$, we will write ${\sf{stack}}^M(\M, \Sigma)$. $\myqedhere$
\end{definition}

\begin{definition}\label{thick set} Suppose $\k$ is a regular cardinal, $\Sigma$ is a $\kappa^+$-iteration strategy and $\M$ is a $\Sigma$-premouse such that $\M\subseteq H_\kappa$. We say $\M$ is $\kappa$-\textbf{fat} if $\kappa=\ord(\M)$ and letting $\M'={\sf{stack}}(\M, \Sigma)$, $\cf(\ord(\M'))\geq \k$. To emphasize the dependence on $\Sigma$, we say that $\M$ is $(\kappa, \Sigma)$-fat. 

We say $\M$ has \textbf{thick sets} (or $\kappa$-\textbf{thick sets} or $(\kappa, \Sigma)$-\textbf{thick sets}) if $\M$ is $\kappa$-fat and $\M'=_{def}{\sf{stack}}(\M, \Sigma)$ has a $(\k, \k+1)$-iteration strategy $\Lambda$ (as a $\Sigma$-premouse) such that whenever $\X$ is a stack on $\M'$ according to $\Lambda$ such that $\X$ is below $\k$, $\pi^\X$ exists and $\pi^\X(\k)=\k$, 
\begin{enumerate}
\item $\pi^\X(\M')={\sf{stack}}(\pi^\X(\M), \Sigma)$, and
\item for some club $C\subseteq \k$, whenever $\tau\in C$ is a non-measurable inaccessible cardinal, $\pi^\X[\ord(\M')]$ contains a $\tau$-club. 
\end{enumerate}
If $\Lambda$ is as above then we say that $(\M, \Lambda)$ has \textbf{thick sets}. $\myqedhere$
\end{definition}

The following lemma is due to Steel. Its proof can be found in \cite[Lemma 5.2]{ATHM}. Below $H_\l$ is the set of hereditarily size $<\l$ sets. 

\begin{lemma}\label{cofinality of the stack} Assume $\sf{NsesS}$ and suppose $(M, \d, \vec{G}, \Sigma, \P)$ has the property $(*+)$\footnote{See \rdef{property *}.}. Let $\l<\d$ and $\M$ be the last model of $({\sf{Le}}(\Sigma^M, \mathcal{J}_{\omega}[\P])_{>\l})^{M|\d}$. Then $M\models ``\M$ is $\d$-fat". 
\end{lemma}

\begin{definition}\label{kappa universal} Suppose $\M$ is a $\Sigma$-premouse and $\k$ is a cardinal such that $\M\subseteq H_\kappa$. We say $\M$ is $(\kappa, \Sigma)$-\textbf{universal} if $\M$ has a $(\kappa, \kappa+1)$-iteration strategy $\Lambda$ (as a $\Sigma$-mouse) such that for all  $(\N, \Phi, \Q, \T)$ with the property that
\begin{itemize}
\item $\N\subseteq H_\kappa$ is a $\Sigma$-premouse,
\item $\Phi$ is a $\kappa+1$-iteration strategy for $\N$ (as a $\Sigma$-premouse),
\item $\X$ is an iteration of $\M$ according to $\Lambda$ such that $\pi^\X$ is defined and $\pi^\X(\kappa)=\kappa$,
\item $\Q$ is the last model of $\X$,
\end{itemize}
$(\Q, \Lambda_{\Q, \X})$ wins the coiteration with  $(\N, \Phi)$. More precisely, if $(\T, \U)$ are the normal stacks on $\Q$ and $\N$ respectively that are produced according to the ordinary comparison procedure by using $\Lambda_{\Q, \X}$ on the $\Q$ side and $\Phi$ on the $\N$ side, then letting $\Q'$ and $\N'$ be the last models of $\T$ and $\U$ respectively,  $\N'\insegeq\Q$. 

If $\Lambda$ is as above then we say that $(\M, \Lambda)$ is $(\kappa, \Sigma)$-\textbf{universal}. $\myqedhere$
\end{definition}

The following simple lemma will be used in the proof of \rthm{existence of thick sets}.
\begin{lemma}\label{measurables have high cofinality} Suppose $\l<\d$ are cardinals, $\d$ is a regular cardinal, $M\subseteq \d$, $E\in V_\l$ is an (possible long) $M$-extender and $N=Ult(M, E)$\footnote{This is the ultrapower that is constructed using functions in $M$.}. Then the following holds:
\begin{enumerate}
\item Suppose $\k\in (\l, \d)$ and $\cf^M(\k)\geq \l$. Then $\sup(\pi_E[\k])=\pi_E(\k)$. 
\item If $\k\in (\l, \d)$ is an inaccessible cardinal then $\pi_E(\k)=\k$.  
\item Suppose $\k>\l$ is an inaccessible cardinal and $\eta\in (\k, \d) $ is a measurable cardinal of $N$ such that $\cf(\eta)<\eta$. Then there is $\eta'>\k$ such that $M\models ``\eta'$ is a measurable cardinal" and $\cf(\eta')<\eta'$.
\end{enumerate}
\end{lemma}
\begin{proof} As clause 1 and 2 are straightforward, we only prove clause 3.  Let $f \in M$ be such that for some $a\in \lh(E)^{<\omega}$, $\eta=\pi_E(f)(a)$. Let $(f_i: i<\k)\subseteq M$ and $(a_i: i<\k)\subseteq \lh(E)^{<\omega}$ be such that $(\pi_E(f_i)(a_i): i<\kappa)$ is increasing and cofinal in $\pi_E(f)(a)$. Let $\tau\leq \l$ be the least such that $\pi_E(\tau)\geq \lh(E)$ and set $h_i(s)=\sup\{ f_i(t): t\in \tau^{<\omega} \wedge f_i(t)<f(s)\}$. We then have that $(\pi_E(h_i)(a): i<\k)$ is cofinal in $\pi_E(f)(a)$. It follows that for $E_a$ measure one many $s$, $(h_i(s): i<\k)$ is cofinal in $f(s)$, as otherwise if $h(s)=\sup\{ h_i(s)+1: i<\k\}$ then we would have $\pi_E(h)(a)<\pi_E(f)(a)$ and for each $i$, $\pi_E(h_i)(a)<\pi_E(h)(a)$. Fix one such $s$ with the property that $f(s)>\k$  and $f(s)$ is a measurable cardinal of $M$ ($E_a$-measure one many $s$ have this property).  Because $(h_i(s): i<\k)$ is cofinal in $f(s)$, we have that $\cf(f(s))<f(s)$. Hence, $\eta'=f(s)$ is as desired. 
\end{proof}

\rthm{existence of thick sets} is the main theorem on thick sets that we will use throughout this book. 
\begin{theorem}\label{existence of thick sets} Assume $\sf{NsesS}$. Suppose
\begin{center}
$(M, \d, \vec{G}, \P, \Sigma, \P', \P^+, \Q, \Lambda, \Lambda', E, \R, \Phi)$
\end{center} 
has the following properties:
\begin{enumerate}
\item $(M, \d,  \vec{G})$ is a background.
\item $\l<\d$, $\mathbb{P}\in M|\l$ is a poset and $g\subseteq \mathbb{P}$ is $M$-generic.
\item  $(\P, \Sigma)$ and $(\Q, \Lambda)$ are allowable pairs with the property that 
\begin{enumerate}
\item $\P\in M|\l$ and $(M, \d, \vec{G}, \Sigma, \P)$ has the property $(*+)$,
\item $\Q\in M|\l[g]$ is a successor type and $M[g]\models ``\Lambda$ is a $(\d^+, \d^+)$-strategy".
\end{enumerate}
\item $\P'$ is the last model of ${\sf{Le}}((\P, \Sigma), \mathcal{J}_\omega[\P])_{> \l}^{(M, \d, \vec{G})}$ and $\P^+=L_{\ord(M)}[\P']$.
\item $E\in M|\l[g]$ is a $\P^+$-extender such that 
\begin{enumerate}
\item $\P^+_E=_{def}Ult(\P^+, E)$ is well-founded, 
\item $\Lambda'=_{def}\Lambda\rest \P^+_E\in \P^+_E$,
\item letting $\vec{H}=\{ E'\in \vec{E}^{\P^+_E}: \cp(E')>\pi_E(\ord(\P))$ and $\P^+_E\models ``\nu(E')$ is an inaccessible cardinal$"\}$, $(\P^+_E, \d, \vec{H}, \Lambda', \Q)$ has the property (*),
\item letting ${\sf{Le}}((\Q^{-}, \Lambda_{\Q^-}), \mathcal{J}_{\omega}[\Q])_{>\l}^{(\P^+_E, \d, \vec{H})}=(\Q_\gg , \Q'_\gg, F^+_\gg, F_\gg, b'_\gg: \gg<\delta)$, 
\begin{enumerate}
\item for all $\gg<\delta$, $\rho(\Q'_\gg)>\d^{\Q}$, 
\item $\R=\Q'_\d$.
\end{enumerate}
\end{enumerate}
\item $\Phi\in M[g]$ is a $(\d, \d+1)$-iteration strategy of $\R$.
\end{enumerate}
Suppose that $\Phi_{\Q^-}=\Lambda_{\Q^-}$. Then for every stack $\X$ according to $\Phi$ such that $\lh(\X)<\d$ and $\pi^\X$ exists, letting $\R_1$ be the last model of $\X$,
\begin{enumerate}
\item $(\R_1, \Phi_{\R_1, \X})$ has $(\d, \Lambda_{\pi^\X(\Q^-), \X})$-thick sets, and
\item (consequently) $(\R_1, \Phi_{\R_1, \X})$ is $(\d, \Lambda_{\pi^\X(\Q^-), \X})$-universal\footnote{Universality follows from the existence of thick sets, for example see the proof of  \cite[Lemma 5.4]{ATHM}.}.
\end{enumerate}

Furthermore, in $M[g]$, $\Phi_{\Q}$ is the unique $(\d, \d+1)$-strategy $\Psi_0$ of $\Q$ such that for some $\S$ and a $(\d, \d+1)$-iteration strategy $\Psi$ of $\S$,
\begin{enumerate}
\item $\Q=\S|\d^\Q$ and $\d^\Q$ is a regular cardinal of $\S$,
\item $\Psi_\Q=\Psi_0$,
\item $\Psi_{\Q^-}=\Lambda_{\Q^-}$, 
\item for every stack $\X$ according to $\Psi$ such that $\lh(\X)<\d$ and $\pi^\X$ exists, letting $\S_1$ be the last model of $\X$, $(\S_1, \Psi_{\S_1, \X})$ has $(\d, \Lambda_{\pi^\X(\Q^-), \X})$-thick sets.
\end{enumerate}
\end{theorem}
\begin{proof} The proofs of all of the claims made above are essentially contained in \cite{ATHM}. We first prove the statements made before the ``furthermore" clause. The following is the first important step.  Set $\Lambda_\X=\Lambda_{\pi^\X(\Q^-), \X}$.

\begin{sublemma}\label{1st important sublemma} In $M[g]$, $\R_1$ is $(\d, \Lambda_\X)$-fat.
\end{sublemma}
\begin{proof} Towards a contradiction assume not. Let $(Z_\a: \a<\d)$ be a continuous chain of submodels of $H_{\d^+}[g]$ of size $<\d$ such that for a club of $\a$, letting $N_\a$ be the transitive collapse of $Z_\a$ and $\tau_\a:N_\a\rightarrow Z_\a$ be the inverse of the collapse, $\a=\cp(\tau_\a)$ and $\powerset(\a)^{\R_1}\subseteq N_\a$. Such a sequence can be constructed following the construction given in the proof of \cite[Lemma 5.2]{ATHM}. 

Let $\W={\sf{LeCore}}_{>\l}^{(\P^+, \d, \vec{H})}$ where $\vec{H}$ consists of those extenders of $\vec{E}^{\P^+}$ whose natural length is an inaccessible cardinal of $\P^+$. It follows from \rlem{core l e} that $\pi_E(\W)$ is a class of $\R$ and therefore, $\W_\X=_{def}\pi^\X(\pi_E(\W))$ is a class of $\R_1$. Hence, for a club of $\a<\d$ the following conditions are true:\\\\
(1) $E\in Z_\a$, $\pi^\X\circ \pi_E(\a)=\a$, $\cp(\tau_\a)=\a$ and $\powerset(\a)^{\W_\X}\subseteq N_\a$. \\\\
If $\a$ is as in (1) then we in fact have that $\powerset(\a)^\W\subseteq N_\a$. However, as in the proof of \cite[Lemma 5.2]{ATHM}, we can find an extender $F^*\in \vec{G}$ such that for some $\nu$, the trivial completion of $F^*\cap (\nu^{<\omega}\times \univ{\W})$ is on $\vec{E}^\W$ and witnesses that $\cp(F)$ is a superstrong cardinal in $\W$, contradicting $\sf{NsesS}$.
\end{proof}
Set $\R^+_1={\sf{stack}}(\R_1, \Lambda_\X)$ and let $\X^+=\uparrow(\X, \R_1^+)$. Applying the proof of \cite[Lemma 5.3]{ATHM} we get the following.
\begin{sublemma}\label{2nd imp sublemma} Suppose $\Y$ is an iteration of $\R$ according to $\Phi$ such that $\pi^\Y$ is defined and $\pi^\Y(\d)=\d$. Then all models of $\Y^+=_{def}\uparrow(\Y, \R^+)$ are well-founded and if $\S$ is the last model of $\Y^+$ then $\S={\sf{stack}}(\S|\d, \Phi_{\pi^{\X^\frown \Y}(\Q^-), \X^\frown \Y})$.
\end{sublemma}
We now have that $(\R_1, \Phi_{\R_1, \X})$ is $(\d, \Lambda_{\Q^-})$-universal (e.g. see the proof of \cite[Lemma 5.4]{ATHM}). Next we show that $(\R_1, \Phi_{\R_1, \X})$ has $(\d, \Lambda_{\X})$-thick sets. Let $\U$ be a stack on $\R_1$ according to $\Phi_{\R_1, \X}$ and let $\U^+=\uparrow(\U, \R^+)$. We are assuming that $\pi^{\U^+}(\d)=\d$\footnote{We in fact should also assume that $\U$ is above $\pi^\X(\ord(\Q^-))$ but this is irrelevant to the proof.}  and want to show that\\\\
(a) for some club $C\subseteq \d$, whenever $\k\in C$ is a non-measurable inaccessible cardinal, $\pi^{\U^+}[\ord(\R^+_1)]$ contains a $\k$-club. \\\\
Set  $\sigma=(\pi^{\U^+}\rest \d+1)\circ (\pi^{\X^+}\rest \d+1) \circ (\pi_E\rest \d+1)$ and $\nu=\ord(\R^+_1)$. Notice that since $\sigma(\d)=\d$, we have a club $C\subseteq \d$ such that for each $\a\in C$, $\sigma[\a]\subseteq \a$. Let $\l'<\d$ be such that $\max(\l, \ord(\P))<\l'$ and $\X\in M|\l'[g]$. We want to show that $C-(\l'+1)$ witnesses (a). Suppose then $\k\in (\l', \d)$ is an inaccessible cardinal of $M$ which is not measurable in $M$ and $\k\in C$. It then follows that $\sigma(\k)=\k$. Indeed, because $E, \X\in M|\l'[g]$ we have that $\pi^\X(\pi_E(\k))=\k$. Notice now that because $\k$ is not measurable in $M$, $\k$ is not measurable in $\P^+$ and therefore, in $\P^+_E$ and consequently in $\R$ and $\R_1$. Hence, it follows from $\pi^\U[\k]\subseteq \k$ that $\pi^\U(\k)=\k$. 

Suppose now that $\a\in [\d, \nu)$ and $\cf^M(\a)=\k$. We claim that $\sup(\pi^\U[\a])=\pi^\U(\a)$. The claim is clear if $\cf^{\R^+_1}(\a)=\k$. Suppose then that  $\eta=_{def}\cf^{\R^+_1}(\a)>\k$. Notice that we have that $\cf^M(\eta)=\k$. We claim that \\\\
(b) $\eta$ is not a measurable cardinal of $\R_1$.\\\\
Assume $\eta$ is measurable in $\R_1$. Then it follows from \rlem{measurables have high cofinality} that there is $\eta'>\k$ such that $\eta'$ is a measurable cardinal of $\R$ and $\cf^M(\eta')<\eta'$. Because $\eta'$ is measurable in $\R$, $\eta'$ is a measurable cardinal of $\P^+_E$. Since $\cf^M(\eta')<\eta'$, \rlem{measurables have high cofinality} implies that there is a measurable cardinal $\eta''$ of $\P^+$ such that $\eta''>\k$ and $\cf^M(\eta'')<\eta''$. But each measurable cardinal of $\P^+$ that is $>\k$ is a measurable cardinal of $M[g]$, contradiction! Thus, (b) holds.

Since $\eta$ is not  a measurable cardinal of $\R_1$ we get that $\sup(\pi^{\U^+}[\eta])=\eta$. Hence, $\sup(\pi^{\U^+}[\a])=\pi^{\U^+}(\a)$. It then follows that $\pi^{\U^+}[\nu]$ is a $\k$-club.

Next, we prove that in $M[g]$, $\Phi_{\Q}$ is the unique $(\d, \d+1)$-strategy $\Psi_0$ of $\Q$ such that for some $\S$ and a $(\d, \d+1)$-iteration strategy $\Psi$ of $\S$,
\begin{enumerate}
\item $\Q=\S|\d^\Q$ and $\d^\Q$ is a regular cardinal of $\S$,
\item $\Psi_\Q=\Psi_0$,
\item $\Psi_{\Q^-}=\Lambda_{\Q^-}$, 
\item for every stack $\X$ according to $\Psi$ such that $\lh(\X)<\d$ and $\pi^\X$ exists, letting $\S_1$ be the last model of $\X$, $(\S_1, \Psi_{\S_1, \X})$ has $(\d, \Lambda_{\pi^\X(\Q^-), \X})$-thick sets.
\end{enumerate}
Fix then $(\S, \Psi)$ that satisfies clause 1, 3, and 4 above. We have that $(\R, \Phi)$ also satisfies those clauses. It is then enough to show that $\Psi_\Q=\Phi_\Q$. Assume not. Let $\U$ be a stack on $\Q$ such that $\Phi_\Q(\U)\not =\Psi_\Q(\U)$. Let $\R^+={\sf{stack}}(\R, \Lambda_{\Q^-})$, $\S^+={\sf{stack}}(\S, \Lambda_{\Q^-})$, $\U_0=\uparrow(\U, \R^+)$ and $\U_1=\uparrow(\U, \S^+)$. Because $\Phi_{\Q^-}=\Psi_{\Q^-}$, we have some $\a<\lh(\U)$ such that $\pi^\U_{0, \a}$ is defined and $\U_{\geq \a}$ is a normal stack on $\M_\a^\U$ and is above $\ord(\pi^\U_{0, \a}(\Q^-))$. Let then $\X_0=(\U_0)_{\leq \a}$, $\X_1=(\U_1)_{\leq \a}$ and $\Y=\U_{\geq \a}$. Set $\R_1=\M_\a^{\X_0}$, $\S_1=\M_\a^{\X_1}$, $\Y_0=\uparrow(\Y, \R_1)$ and $\Y_1=\uparrow(\Y, \S_1)$. Finally set $\Phi_\Q(\U)=b_0$ and $\Psi_\Q(\U)=b_1$. 

We claim that $\Q(b_0, \U)$ doesn't exist. Towards a contradiction assume it does exist. Assume first that $\Q(b_1, \U)$ doesn't exist. It follows that $\d(\Y)$ is not a limit of Woodin cardinals of $\m(\Y)$, and therefore, $\Q(b_0, \U)$ is a $\Lambda_{\pi^\U_{0, \a}(\Q^-), \U_{\leq \a}}$-mouse over $\m(\Y)$, and since $\M_{b_1}^{\Y_1}$ is universal, $\Q(b_0, \U)\insegeq \M_{b_1}^{\Y_1}$\footnote{More precisely, setting $\S'=\M_{b_1}^{\Y_1}$, $(\S', \Psi_{\S', \U_1^\frown\{b\}})$ is $(\d, \Lambda_{\pi^\U_{0, \a}(\Q^-), \U_{\leq \a}})$-universal.}. Thus, we must have that both $\Q(b_0, \U)$ and $\Q(b_1, \U)$ exist. A similar argument shows that $\U$ cannot have a fatal drop, implying that $\Q(b_0, \U)$ and $\Q(b_1, \U)$ are $\Lambda_{\pi^\U_{0, \a}(\Q^-), \U_{\leq \a}}$-mice over $\m(\Y)$. Hence, $\Q(b_0, \U)=\Q(b_1, \U)$ implying that $b_0=b_1$, contradiction. Hence, $\Q(b_0, \U)$ doesn't exist.  A symmetric argument shows that $\Q(b_1, \U)$ also does not exist.

We thus have that for $i\in 2$, $\pi^{\Y_i}_{b_i}$ is defined. Let $\R_2=\M^{\Y_0}_{b_0}$ and $\S_2=\M^{\Y_1}_{b_1}$. Both $\R_2$ and $\S_2$ are $\Lambda_{\pi^\U_{0, \a}(\Q^-), \U_{\leq \a}}$-mice. We can then find $\W$ such that \\\\
(1) $\W$ is a $\Phi_{\R_2, \U_0^\frown \{b_0\}}$-iterate of $\R_2$ and the iteration embedding $j_0: \R_2\rightarrow \W$ exists and has the property that $j_0(\d)=\d$, and\\
(2) $\W$ is a $\Psi_{\S_2, \U_1^\frown \{b_1\}}$-iterate of $\S_2$ and the iteration embedding $j_1: \S_2\rightarrow \W$ exists and has the property that $j_1(\d)=\d$.\\\\
Because of our assumption on thick sets, we have a club $C_0\subseteq \d$ and a club $C_1\subseteq \d$ such that for every $\k\in C_0\cap C_1$ that is an inaccessible cardinal of $M$ but not a measurable cardinal of $M$, \\\\
(3) $j_0 \circ \pi^{\Y_0}_{b_0}[\ord(\R_1)]$ and $j_1 \circ \pi^{\Y_1}_{b_1}[\ord(\S_1)]$ contain a $\k$-club.\\\\
(3) then implies that\\\\
(4) $(j_0 \circ \pi^{\Y_0}_{b_0}[\ord(\R_1)])\cap (j_1 \circ \pi^{\Y_1}_{b_1}[\ord(\S_1)])$ contains a $\k$-club.\\\\
Let then $D\subseteq ((j_0 \circ \pi^{\Y_0}_{b_0}[\ord(\R_1)])\cap (j_1 \circ \pi^{\Y_1}_{b_1}[\ord(\S_1)])$ be a $\k$-club and set $D_0=(j_0 \circ \pi^{\Y_0}_{b_0})^{-1}[D]$ and $D_1=(j_1 \circ \pi^{\Y_1}_{b_1})^{-1}[D]$. Let now $\Q_0=\pi^\U_{0, \a}(\Q)$. Notice that $\Y$ is a normal stack on $\Q_0$ that is above $\ord(\Q^-_0)$ and below $\d^{\Q_0}$. We now have that\\\\
(5) $\d^{\Q_0}=\sup(Hull^{\R_1}(D_0 \cup \Q_0^-)\cap \d^{\Q_0})$,\\
(6) $\d^{\Q_0}=\sup(Hull^{\S_1}(D_1 \cup \Q_0^-)\cap \d^{\Q_0})$,\\
(7) $\d(\Y)=\sup(Hull^{\R_2}(\pi^{\Y_0}_{b_0}[D_0] \cup \Q_0^-)\cap \d(\Y))$,\\
(8) $\d(\Y)=\sup(Hull^{\S_2}(\pi^{\Y_1}_{b_1}[D_1] \cup \Q_0^-)\cap \d(\Y))$.\\\\
(5)-(8) are consequences of universality. For example, (5) can be shown as follows. Suppose $\d^{\Q_0}>\sup(Hull^{\R_1}(D_0 \cup \Q_0^-)\cap \d^{\Q_0})$ and set $\gg=\sup(Hull^{\R_1}(D_0 \cup \Q_0^-)\cap \d^{\Q_0})$. Let $\R'=cHull^{\R_1}(D_0 \cup \gg)$ and let $\tau:\R'\rightarrow \R_1$ be the inverse of the transitive collapse. Then because $\tau(\Q^-_0)=\Q_0^-$, we have that $\R'$ is a $\Lambda_{\Q_0^-, \U_{\leq \a}}$-mouse as witnessed by $\Phi'=(\tau$-pullback of $\Phi_{\R_1, \X_0})$. Moreover, it follows from \cite[Lemma 5.4]{ATHM} that $\R'={\sf{stack}}(\R'|\d, \Lambda_{\Q_0^-, \U_{\leq \a}})$ and $(\R', \Phi')$ is $(\d, \Lambda_{\Q_0^-, \U_{\leq \a}})$-universal. But because $\R'\models ``\gg$ is a Woodin cardinal$"$ and $\R_1\models ``\gg$ is not a Woodin cardinal$"$, we have a contradiction.

(5)-(8) easily imply that $\rge(\pi^{\Y_0}_{b_0})\cap \rge(\pi^{\Y_1}_{b_1})$ is cofinal in $\d(\Y)$. Hence, because $\pi^{\Y_0}_{b_0}\rest \d^{\Q_0}=\pi^\Y_{b_0}\rest \d^{\Q_0}$ and $\pi^{\Y_1}_{b_1}\rest \d^{\Q_0}=\pi^\Y_{b_1}\rest \d^{\Q_0}$, we have that $b_0=b_1$. 
\end{proof}

The next few chapters are essentially applications of \rthm{existence of thick sets}.

\section{Fullness preservation}\label{fullness preservation}

Throughout this section we assume $\sf{AD}^+$. Below, we use $\R^*$ to denote the \textasteriskcentered-translation of $\R$ (cf. \cite{Selfiter} or \cite[Remark 12.7]{DMATM}.). Suppose $\eta$ is a cutpoint cardinal of a hod premouse $\R$. The \textasteriskcentered-translation is used to translate $\R|(\eta^+)^\R$ into an $\sf{lses}$ over $\R|\eta$. More precisely, $\univ{(\R|(\eta^+)^\R)^*}=\univ{\R|(\eta^+)^\R}$ but $\eta$ is a strong cutpoint of $(\R|(\eta^+)^\R)^*$. If in fact $\eta$ is already a strong cutpoint of $\R$ then $(\R|(\eta^+)^\R)^*=\R|(\eta^+)^\R$. Thus, as far as grasping the main ideas are concerned, the reader will lose little by treating all cutpoint cardinals as strong cutpoint cardinals.

\begin{definition}\label{projectively closed pointclass} We say $\Gamma$ is \textbf{projectively closed} if whenever $A$ is a set of reals such that for some $B\in \Gamma$, $A$ is first order definable over $({\sf{HC}}, B, \in)$ (with parameters), $A\in \Gamma$. $\myqedhere$
\end{definition}

\begin{definition}[$\Gamma$-Fullness preservation]\label{gamma fullness preservation}\index{$\Gamma$-fullness preservation} Suppose $(\P, \Sigma)$ is a hod pair or an sts hod pair\footnote{Recall that if $(\P, \Sigma)$ is an sts hod pair then $\P=(\P|\d^\P)^\#$. See \rdef{pre sts hod pairs}.} such that $\P\in {\sf{HC}}$ and $\Gamma$ is a projectively closed pointclass. We say $\Sigma$ is \textbf{$\Gamma$-fullness preserving} if the following holds for all $(\T, \Q)\in I(\P, \Sigma)$. 
\begin{enumerate}
\item For all meek\footnote{See \rdef{pre-hod-like}.} layers $\R$ of $\Q$ such that $\R$ is of successor type\footnote{See \rdef{types of lsa small premice}.}, letting $\S=\R^-$\footnote{This is the longest proper layer of $\R$. See \rnot{l p}.}, for all $\eta\in (\ord(\S), \ord(\R))$ if $\eta$ is a cutpoint cardinal of $\R$ then
\begin{center}
$(\R|(\eta^+)^\R)^*={\sf{Lp}}^{\Sigma_{\S, \T}}(\R|\d)$.
\end{center} 
\item  For all meek\footnote{See \rdef{pre-hod-like}.} layers $\R$ of $\Q$ such that $\R$ is of limit type, 
\begin{center}
$\R={\sf{Lp}}^{\Sigma_{\R|\d^\R, \T}}(\R|\d^\R)$.
\end{center}
\item If $\P$ is of $\#$-lsa type then ${\sf{Lp}}^{\Gamma, \Sigma^{stc}_{\Q, \T}}(\Q)\models ``\d^\Q$ is a Woodin cardinal"\footnote{Here, if $\Sigma$ is a short tree strategy then $\Sigma^{stc}=\Sigma$.}. 
\end{enumerate}
If only conditions 1 and 2 hold then we say that $\Sigma$ is \textbf{almost $\Gamma$-fullness preserving}. We say that $\Sigma$ is \textbf{lower-level $\Gamma$-fullness preserving} if the above clauses hold for $\R\inseg_{hod}\Q$\footnote{We will use this version of fullness preservation when studying anomalous hod pairs (see \rsec{sec: anomalous hod premice}). For now, the reader may ignore it. The concept will became important in \rthm{the generation of mouse full pointclasses ii}.}. 

Suppose $(\P, \Sigma)$ is a hod pair such that $\P$ is gentle. Then we say that $\Sigma$ is $\Gamma$-\textbf{fullness preserving} if for every $\Q\in Y^\P$, $\Sigma_\Q$ is $\Gamma$-fullness preserving.

If $\Gamma$ is a Solovay pointclass then we will omit it from the terminology. $\myqedhere$
\end{definition}

\begin{theorem}[Fullness preservation of induced strategies]\label{fullness preservation of background constructions}\label{thm:fullness_preserving} Assume $\sf{AD}^+$.\\ Suppose $\Gamma$ is a pointclass such that for some $\a$ with $\theta_\a<\Theta$, $\Gamma=\{ A\subseteq \bR: w(A)<\theta_\a\}$, ${\sf{C}}=(\mathbb{M}, (P, \Psi), \Gamma^*, A)$ Suslin, co-Suslin captures $\Gamma$ and $\mathbb{M}=(M, \d,  \vec{G}, \Sigma)$. Set
\begin{center} 
${\sf{hpc}}_{{\sf{C}}, \Gamma}^+=(\M_\gg , \N_\gg, Y_\gg, \Phi^+_\gg, F^+_\gg,  F_\gg, b_\gg: \gg\leq \d)$.
\end{center} 
Suppose $\b<\d$, $\P\in Y_\b$ and $M\models ``(\P, (\Phi_\b)_\P)\in {\sf{HP}}^\Gamma$". Then $(\Phi_\b^+)_\P$ is almost $\Gamma$-fullness preserving. 

Moreover, assuming that 
\begin{itemize}
\item $\P$ is of $\#$-lsa type,
\item letting $\Psi=(\Phi^+_\b)^{stc}_\P$, ${\sf{Lp}}^{\Gamma, \Psi}(\P)\models ``\d^\P$ is a Woodin cardinal" and
\item ${\sf{Le}}((\P, (\Phi_\b^{stc})_\P), \mathcal{J}_{\omega}[\P])$ does not break down because of the anomaly stated in clause 3.b of \rdef{fully backgrounded sts construction},
\end{itemize}
$\Psi$ is $\Gamma$-fullness preserving. Also, the above clauses hold for $\xi$ as in the definition of $Y_{\xi+1}$ that appears in clause 5.c.ii of \rdef{gamma-hod pair construction*}.
\end{theorem}
\begin{proof} Below we will use the universality clause of \rthm{existence of thick sets}. Towards a contradiction, assume $\Lambda=_{def}(\Phi_\b^+)_\P$ is not $\Gamma$-fullness preserving. We have that $\Lambda$ is the $id$-pullback of $\Phi_\b^+$\footnote{See \rdef{id pullback initial segment}.}. It follows by absoluteness\footnote{See \rlem{correctness of backgrounds} and \rcor{capturing particular sets}. Here we use the fact that $M\models (\P, (\Phi_\b)_\P)\in {\sf{HP}}^\Gamma$.} that there is a counterexample in $M[g]$ where $g\subseteq Coll(\omega, \nu)$ is $M$-generic and $\nu<\d$.  All the clauses of $\Gamma$-fullness preservation are very similar and follow from the universality of background constructions. Below we derive a contradiction from the failure of  clause 2 of \rdef{gamma fullness preservation} and leave the rest to the reader. We also leave the ``moreover" clause to the reader as it is very similar to the other cases. We can then further assume that $\P$ is of limit type as otherwise we would just be re-proving \cite[Lemma 5.7]{ATHM}. 

Fix $(\T, \Q)\in I(\P, \Lambda)$\footnote{It is irrelevant whether $\T$ is an ordinary stack or a generalized stack.} and fix $\R$ which is as in clause 2 of \rdef{gamma fullness preservation} (so $\R$ is a meek layer of $\Q$ of limit type, implying that $\R=\R^b$). Let $\k=\d^\R$. We need to see that
\begin{center}
$\R={\sf{Lp}}^{\Gamma, \Lambda_{\R|\k, \T}}(\R|\k)$.
\end{center}
Using \rlem{summary res and emb}, we can find a $\Sigma$-iterate $N$ of $M$ such that letting $i: M\rightarrow N$ be the iteration embedding and 
 \begin{center}${\sf{hpc}}^+_{{\sf{C_N}}, \Gamma}=(\S'_\gg , \S_\gg, Z_\gg, \Psi^+_\gg, E^+_\gg, E_\gg, c_\gg : \gg\leq i(\delta))$,
\end{center} 
there is a $\gg< i(\d)$ and a weak embedding $\sigma:\Q\rightarrow  i(\P)$ such that\\\\
(1) $\sigma\circ \pi^\T=i\rest \P$,\\
(2) $\Lambda_{\Q, \T}$ is the $\sigma$-pullback of $(\Psi_{\gg}^+)_{i(\P)}$,\\
(3) $i(\P)\in Z_\gg$ .\\\\
Suppose first that $\M\insegeq \R$ is such that $\R|\k\inseg \M$. We need to see that \\\\
(a) $\M$ as a $\Lambda_{\R|\k, \T}$-premouse has an $\omega_1$-iteration strategy in $\Gamma$.\\\\
Clearly (2) and (3) above easily imply (a)\footnote{In fact this also follows from our choice of $\Gamma$ as since ${\sf{Code}}(\Lambda_{\R|\k, \T})\in \Gamma$, any $\Lambda_{\R|\k, \T}$-mouse $\M$ such that $\rho(\M)=\k$ has an iteration strategy in $\Gamma$.}. 

Fix now $\M\insegeq {\sf{Lp}}^{\Gamma, \Lambda_{\R|\k, \T}}(\R|\k)$ such that $\rho(\M)=\k$. We want to see that\\\\
(b) $\M\insegeq \R$. \\\\
We let $\pi=\pi^\T$, $\tau=\d^{\P^b}$, $\zeta=\sup\{\lh(F^+_\gg): \gg<\b\}$, $\vec{G}'=\{F\in \vec{G}: \cp(F)>\max(\zeta, \nu)\}$ and $\N$ be the last model of 
\begin{center}
 $({\sf{Le}}((\P|\tau, \Lambda_{\P|\tau}), \P^b)_{>\zeta})^{(M[g], \d, \vec{G}')}$. 
 \end{center} 
Notice that because of our choice of $\Gamma$ (see the footnote above), the fact that $(\P, \Lambda)$ is a $\Gamma-{\sf{cbl}}$ and the  $(\d, \Lambda_{\P|\tau})$-universality of $\N$,\\\\
 (4) $\P=\N|(\tau^+)^\N$. \\\\
Notice next that if $E$ is the $(\tau, \d^{\Q^b})$-extender derived from $\pi^\T$  then \\\\
 (5) $M[g]\models ``Ult(\N, E)$ is $\d$-iterable".\\\\
 This is because $\sigma:\Q\rightarrow i(\P)$ can be extended to  $\sigma^+:Ult(\N, E)\rightarrow i(\N)$. 
 
 Let then $\pi^+=\pi^\N_{E}$, $\vec{H}=\{E\in \vec{E}^{Ult(\N, E)}: \cp(E)>\d^\Q$ and $\nu(E)$ is an inaccessible cardinal of $Ult(\N, E)\}$, and $\N^*$ be the last model of
 \begin{center}
 $({\sf{Le}}((\R|\k, \Lambda_{\R|\k, \T}), \R|\k))^{Ult(\N, E), \d, \vec{H}}$.
 \end{center} 

It then follows from $(\d, \Lambda_{\R|\k, \T})$-universality of $\N^*$ that $\M\insegeq \N^*$. Therefore, $\M\in Ult(\N, E)$, and since $\R=Ult(\N, E)|(\k^+)^{Ult(\N, E)}$, $\M\in \R$. Since $\M$ is $\omega_1$-iterable, it follows that  $\M\insegeq\R$.
\end{proof}

The proof actually gives more. 

\begin{definition}[Strongly $\Gamma$-fullness preserving]\label{strongly fullness preserving}\index{strong $\Gamma$-fullness preservation} Suppose $(\P, \Sigma)$ is a hod pair or an sts hod pair and $\Gamma$ is a pointclass. We say $\Sigma$ is strongly $\Gamma$-fullness preserving if $\Sigma$ is $\Gamma$-fullness preserving and whenever
\begin{enumerate}
\item  $\T$ is a stack according to $\Sigma$ with last model $\S$ such that if $\P$ is of limit type then $\pi^{\T, b}$ exists and otherwise $\pi^{\T}$ exists, and
\item $\R$ is such that there are elementary embedding $(\sigma, \tau)$ with the property that 
\begin{enumerate}
\item if $\P$ is of limit type then $\sigma:\P^b\rightarrow \R$, $\tau:\R\rightarrow \S^b$ and $\pi^{\T, b}=\tau\circ \sigma$, and
\item if $\P$ is of successor type then $\sigma:\P\rightarrow \R$, $\tau:\R\rightarrow \S$ and $\pi^{\T}=\tau\circ \sigma$,
\end{enumerate}
\end{enumerate}
then the $\tau$-pullback strategy of $\Sigma_{\S^b, \T}$ if 2(a) holds and of $\Sigma_{\S,\T}$ if 2(b) holds is $\Gamma$-fullness preserving. Following \rdef{gamma fullness preservation} we can also define the meaning of \textbf{strongly almost $\Gamma$-fullness preserving} as well as the meaning of  \textbf{strongly low-level $\Gamma$-fullness preserving}.$\myqedhere$
\end{definition}

The following is then a corollary to the proof of \rthm{fullness preservation of background constructions} and we leave it to the reader.

\begin{theorem}[Strong fullness preservation of induced strategies]\label{strong fullness preservation}\index{strong fullness preservation} Assume $\sf{AD}^+$ and suppose $\Gamma$ is a pointclass such that for some $\a$ with $\theta_\a<\Theta$, $\Gamma=\{ A\subseteq \bR: w(A)<\theta_\a\}$, ${\sf{C}}=(\mathbb{M}, (P, \Psi), \Gamma^*, A)$ Suslin, co-Suslin captures $\Gamma$ and $\mathbb{M}=(M, \d,  \vec{G}, \Sigma)$. Set
\begin{center} 
${\sf{hpc}}_{{\sf{C}}, \Gamma}^+=(\M_\gg , \N_\gg, Y_\gg, \Phi^+_\gg, F^+_\gg,  F_\gg, b_\gg: \gg\leq \d)$.
\end{center} 
Suppose $\b<\d$ and $\P\in Y_\b$. Then $(\Phi_\b^+)_\P$ is almost $\Gamma$-fullness preserving. 

Moreover, assuming that 
\begin{itemize}
\item $\P$ is of $\#$-lsa type,
\item letting $\Psi=(\Phi^+_\b)^{stc}_\P$, ${\sf{Lp}}^{\Gamma, \Psi}(\P)\models ``\d^\P$ is a Woodin cardinal" and
\item for every $\zeta<\d$ the ${\sf{Le}}((\P, (\Phi_\b^{stc})_\P), \mathcal{J}_{\omega}[\P])_{>\zeta}$ does not break down because of the anomaly stated in clause 3.b of \rdef{fully backgrounded sts construction}\footnote{We only need this condition for $\zeta=\sup\{\lh(F^+)_\gg: \gg<\b\}$.},
\end{itemize}
$\Psi$ is $\Gamma$-fullness preserving. 
\end{theorem}

The following is an easy yet useful consequence of strong fullness preservation. 

\begin{lemma}\label{easy consequence of strong fullness preservation}
Assume $\sf{AD}^+$ and suppose $\Gamma$ is a pointclass. Suppose further that $(\P, \Sigma)$ is a hod pair or an sts hod pair such that $\Sigma$ is strongly $\Gamma$-fullness preserving. Let $\T$ be a stack on $\P$ according to $\Sigma$ with last model $\S$ such that if $\P$ is of limit type then $\pi^{\T, b}$ exists and otherwise $\pi^{\T}$ exists. Suppose $(\R, \sigma, \tau)$ is such that 
\begin{enumerate}
\item if $\P$ is of limit type then $\sigma:\P^b\rightarrow \R$, $\tau:\R\rightarrow \S^b$ and $\pi^{\T, b}=\tau\circ \sigma$, and
\item if $\P$ is of successor type then $\sigma:\P\rightarrow \R$, $\tau:\R\rightarrow \S$ and $\pi^{\T}=\tau\circ \sigma$.
\end{enumerate}
Let $E$ be such that
\begin{enumerate}
\item if $\P$ is of limit type then $E$ is the $(\d^{\P^b}, \d^\R)$-extender derived from $\sigma$, and
\item if $\P$ is of successor type then $E$ is the $(\d^\P, \d^\R)$-extender derived from $\sigma$
\end{enumerate} 
Then $\R=Ult(\P, E)$. In particular, $\R=\{\pi_E(f)(a):f\in \P$ and $a\in (\d^\R)^{<\omega}\}$.
 \end{lemma}
 \begin{proof}
 Let $k:Ult(\P, E)\rightarrow \R$ be the factor map, i.e., $k(\pi(f)(a))=\sigma(f)(a)$. Then if $\P$ is of limit type then
 $\pi^{\T, b}=\tau\circ k\circ \pi_E$ and if $\P$ is of successor type then $\pi^{\T}=\tau\circ k\circ \pi_E$. Notice that $\cp(k)>\d^\R$. It now follows from strong $\Gamma$-fullness preservation of $\Sigma$ that $\Sigma_{\S, \T}^{\tau\circ k}$, the $\tau\circ k$-pullback of $\Sigma_{\S, \T}$, is $\Gamma$-fullness preserving. But because $k\rest \d^\R=id$, we have that for every $\R'\in Y^{\R}$, 
 \begin{center}
 $(\Sigma^{\tau\circ k}_{\S, \T})_{\R'}=(\Sigma^{\tau}_{\S, \T})_{\R'}$.
 \end{center}
 It then follows that $\R=Ult(\P, E)$.
 \end{proof}
 
 \section{Tracking disagreements}

Here we introduce terminology that we will use to track the disagreements between strategies. The reader may wish to review \rnot{l p}, \rdef{allowable pair} and  \rter{types of lsa small premice}.


\begin{definition}[Low level disagreement between strategies]\label{low level disagreement between strategies}\index{low level disagreement} Suppose $(\P, \Sigma)$ and $(\P, \Lambda)$ are two allowable pairs. 
We say that there is a \textbf{low level disagreement} between $\Sigma$ and $\Lambda$ if one of the following conditions holds:
\begin{enumerate}
\item $\P$ is of successor type and $\Sigma_{\P^-}\not =\Lambda_{\P^-}$.
\item $\P$ is gentle and for some complete proper layer $\Q$ of $\P$, $\Sigma_\Q\not=\Lambda_\Q$.
\item $\P$ is of limit type, $\P$ is meek and there is $(\T, \Q)\in B(\P, \Sigma)\cap B(\P, \Lambda)$ such that $\Sigma_{\Q, \T}\not = \Lambda_{\Q, \T}$.
\item $\P$ is of limit type, $(\P, \Sigma)$ and $(\P, \Lambda)$ are hod pairs or sts hod pairs and there is $(\T_1, \P_1)\in I^b(\P, \Sigma)$ and $(\T_2, \P_2)\in I^b(\P, \Lambda)$ such that 
\begin{enumerate}
\item $\Q=_{def}\P_1^b=\P_2^b$, 
\item $\pi^{\T_1, b}=\pi^{\T_2, b}$, and
\item $\Sigma_{\Q, \T_1}\not =\Lambda_{\Q, \T_2}$.
\end{enumerate}
\item $\P$ is of limit type, $(\P, \Sigma)$ and $(\P, \Lambda)$ are simple hod pairs or simple sts hod pairs and there is $(\T_1, \P_1)$ and $(\T_2, \P_2)$ such that 
\begin{enumerate}
\item $(\T_1, \P_1)\in I^{ope}(\P, \Sigma)$\footnote{See \rdef{one point extension} and \rdef{almost non-dropping gen stacks}. We mainly use this to conclude that the un-dropping extender exists.},
\item $(\T_2, \P_2)\in I^{ope}(\P, \Lambda)$,
\item $\Q=_{def}\P_1^b=\P_2^b$, 
\item the $\Q$-un-dropping extenders of $\T_1$ and $\T_2$ are the same,
\item $\Sigma_{\Q, \T_1}\not =\Lambda_{\Q, \T_2}$.
\end{enumerate}
\end{enumerate} 
If clause 4 or 5 holds then we say that $(\T_1, \P_1, \T_2, \P_2)$ is a low level disagreement between $\Sigma$ and $\Lambda$. Suppose next that $\P$ is of limit type. We say $(\T_1, \P_1, \T_2, \P_2, \Q)$ is a \textbf{minimal low level disagreement} \index{minimal low level disagreement} if, 
\begin{enumerate}
\item $(\T_1, \P_1, \T_2, \P_2)$ is a low level disagreement between $\Sigma$ and $\Lambda$,
\item $\Q$ is of successor type and $\Q\insegeq \P_1^b=\P_2^b$,
\item $\Sigma_{\Q^-, \T_1} = \Lambda_{\Q^-, \T_2}$, 
\item $\Sigma_{\Q, \T_1}\not =\Lambda_{\Q, \T_2}$.
\end{enumerate}
$\myqedhere$
\end{definition}

Next we show that the existence of a disagreement translates into the existence of a minimal low level disagreement.  The reader may wish to review \rdef{initial segment of a stack given by a node}, \rdef{main drops}, \rnot{notation for generalized stacks},  \rrem{some comments on gen strategies} and \rdef{allowable pair}.

\begin{lemma}[Disagreement implies low level disagreement]\label{disagreement implies low level disagreement} Suppose $\Gamma$ is a projectively closed pointclass, and $(\P, \Sigma)$ and $(\P, \Lambda)$ are allowable pairs such that  both $\Sigma$ and $\Lambda$ are almost $\Gamma$-fullness preserving. Suppose that one of the following conditions holds:
\begin{enumerate}
\item $\P$ is of limit type but not of lsa type, and $\Sigma\not =\Lambda$.
\item $(\P, \Sigma)$ and $(\P, \Lambda)$ are sts pairs or simple sts pairs, $\Sigma\not=\Lambda$ and both $\Sigma$ and $\Lambda$ are fullness preserving. 
\item $(\P, \Sigma)$ is a (simple) sts pair, $(\P, \Lambda)$ is a (simple) hod pair, $\Sigma\not=\Lambda$ and both $\Sigma$ and $\Lambda$ are fullness preserving. 
\item $\P$ is of lsa type, $\mathcal{J}_{\omega}[\P]\models ``\d^\P$ is not a Woodin cardinal" and  $\Sigma\not=\Lambda$.
\end{enumerate}
Then there is a low level disagreement between $\Sigma$ and $\Lambda$. 
\end{lemma}
\begin{proof}
We give the proof from clause 2, which is the hardest, and leave the rest to the reader. The proof from clause 1 is easier and is similar to \cite[Proposition 2.41]{ATHM}). We also assume that $(\P, \Sigma)$ and $(\P, \Lambda)$ are sts pairs (as apposed to simple sts pairs).

Thus, we assume that $(\P, \Sigma)$ and $(\P, \Lambda)$ are sts hod pairs and $\Sigma\not=\Lambda$. We then have that $\P=(\P|\d^\P)^\#$. Assume there is no low level disagreement between $\Sigma$ and $\Lambda$ and let $\T=(\S_\a, \Y_\a, E_\a: \a<\eta)$ be any disagreement between $\Sigma$ and $\Lambda$. Because $\Sigma(\T)\not =\Lambda(\T)$ we must have that\\\\
(1) $\eta=\gg+1$, $\Y_\gg\not=\emptyset$, $E_\gg=\emptyset$ and $\lh(\Y_\gg)$ is a limit ordinal. \\\\
Set $\U=\Y_\gg$. For $\xi<\lh(\U)$ we let $\M_\xi=\M_\xi^\U$.  
\begin{sublemma}\label{sublemma about u} The following holds. 
\begin{enumerate}
\item If $\a\in R^{\U}$ is such that $\pi_{0, \a}^{\U}$ is defined then $\Sigma_{\M_{\a}^b}=\Lambda_{\M_{\a}^b}$. 
\item $\U$ does not have a main drop\footnote{See \rdef{main drops}.}.
\item $R^{\U}$ has a largest element and if  $\a=\max(\R^{\U})$ then $\U_{\geq \a}$ is above $\ord(\M_{\a}^b)$. 
\end{enumerate}
\end{sublemma}
\begin{proof} Clause 1 is an immediate consequence of our assumption that there are no low level disagreements between $\Sigma$ and $\Lambda$. To see that $\U$ does not have a main drop, suppose that it does and let  
\begin{center}
 $md^{\U}=(\a_i, \R_i, \W_i, \R'_i: i\leq k+1)$
 \end{center}
 be the $md$-sequence of $\U$. It follows that $(\U)_{\geq \a_1}$ is based on $\R'_1\insegeq \R_1^b$ and therefore, $\Sigma_{\R'_1, \T_{\leq \R'_1}}\not =\Lambda_{\R'_1, \T_{\leq \R'_1}}$. 
 Let $\T'=(\S'_\xi, \Y'_\xi, E'_\xi: \xi< \eta)$ be such that
 \begin{enumerate}
 \item for $\xi\leq \gg$, $\S'_\xi=\S_\xi$ and $E_\xi'=E_\xi$,
 \item for $\xi<\gg$, $\Y'_\xi=\Y_\xi$,
 \item $\Y'_\gg=\U_{\leq \b}$ where $\b$ is the least such that $\M_\b^b=\R_1^b$\footnote{Notice that it follows that $\pi^{\U_{\leq \b}}$ is defined.}.
 \end{enumerate} Then $(\T', \R_1^b)$  is a low level disagreement between $\Sigma$ and $\Lambda$. 
 
 To see that clause 3 holds, notice that if $R^{\U}$ doesn't have a maximal element then $\U$ has a unique branch which must be chosen by both $\Sigma$ and $\Lambda$. Suppose now that $\a=\max(R^\U)$. Recall our convention on proper stacks (see \rrem{proper stacks convention}). Thus, every cutpoint of $\U$ belongs to $R^{\U}$. Therefore, as $\a=\max(R^\U)$ and as every cutpoint of $\U$ belongs to $R^\U$, we have the following four possibilities.
 \begin{enumerate}
 \item $\U_{\geq \a}$ is above $\ord(\M^b_\a)$.
 \item $\U_{\geq \a}$ is above $\d^{\M^b_\a}$ but below $\ord(\M_\a^b)$.
 \item $\U_{\geq\a}$ is below $\d^{\M^b_\a}$.
 \item $\lh(\U_{\geq\a})=2$ and $\cp(E_0^{\U_{\geq\a}})=\d^{\M_\a^b}$.
 \end{enumerate}
 If 1 holds then there is nothing to prove. Clearly 4 fails as $\U$ has a limit length. 
 
We now show that neither 2 nor 3 can hold. Assume 2 holds. Because both $\Sigma$ and $\Lambda$ are $\Gamma$-fullness preserving, $\Sigma(\T)=\Lambda(\T)$. 

Assume now that 3 holds. Let $\T'=(\S'_\xi, \Y'_\xi, E'_\xi: \xi< \eta)$ be such that
 \begin{enumerate}
 \item for $\xi\leq \gg$, $\S'_\xi=\S_\xi$ and $E_\xi'=E_\xi$,
 \item for $\xi<\gg$, $\Y'_\xi=\Y_\xi$,
 \item $\Y'_\gg=\U_{\leq \b}$ where $\b$ is the least such that $\M_\b^b=\M_\a^b$\footnote{Notice that it follows that $\pi^{\U_{\leq \b}}$ is defined.}.
 \end{enumerate}
 Then $(\T', \M_\b^b)$ constitutes a low level disagreement between $\Sigma$ and $\Lambda$.
 \end{proof}
Let $\a_0=\max(R^\U)$ and $\mathcal{X}=\U_{\geq \a_0}$. Set $\P_1=\m^+(\X)$.

\begin{sublemma}\label{t is a branch} There are ordinary stacks\footnote{As apposed to generalized stacks.} $\T_1$ and $\T_2$ on $\M_{\a_0}^\U$ such that
\begin{enumerate}
\item $(\T_{\leq \M_\a^\U})^\frown \T_1$ is according to $\Sigma$ and $(\T_{\leq \M_\a^\U})^\frown \T_2$ is according to $\Lambda$, 
\item $\T_1$ and $\T_2$ use the same extenders,
\item both $\pi^{\T_1, b}$ and $\pi^{\T_2, b}$ exist and $\pi^{\T_1, b}=\pi^{\T_2, b}$,
\item $(\T_{\leq \M_\a^\U})^\frown \T_1\in b(\Sigma)$ and $(\T_{\leq \M_\a^\U})^\frown \T_2\in b(\Lambda)$, i.e., $\Sigma((\T_{\leq \M_\a^\U}) ^\frown \T_1)$ and $\Lambda((\T_{\leq \M_\a^\U})^\frown \T_2)$ are branches rather than models,
\item $\T_1$ and $\T_2$ have last normal components $\X_1$ and $\X_2$, 
\item letting $b=\Sigma((\T_{\leq \M_\a^\U})^\frown \T_1)$ and $c=\Lambda((\T_{\leq \M_\a^\U})^\frown \T_2)$, $\Q(b, \X_1)\not=\Q(c, \X_2)$\footnote{Because $\T_1$ and $\T_2$ use the same extenders, we have that $\m^+(\X_1)=\m^+(\X_2)$.}.
\end{enumerate}
\end{sublemma}
\begin{proof} Towards a contradiction, suppose not. As $\T$ is a disagreement between $\Sigma$ and $\Lambda$, we have that $\T\not \in b(\Sigma)\cap b(\Lambda)$ as otherwise we could just take $\T_1=\X=\T_2$. Notice that since $\T$ is a disagreement between $\Sigma$ and $\Lambda$, $\T\not \in m(\Sigma)\cap m(\Lambda)$, as otherwise $\Sigma(\T)=\P_1=\Lambda(\T)$. Assume without loss of generality that $\T\in m(\Sigma)$ and $\T\in b(\Lambda)$. Then letting $c=\Lambda(\T)$, $\Q(c, \X)$ exists. Let $\Sigma_1$ be the $(\omega_1, \omega_1)$-portion of $\Sigma_{\P_1, \T}$ and $\Lambda_1$ be the $(\omega_1, \omega_1)$-portion of $(\Lambda_{\Q(c, \X), \T})_{\sf{ex}}$\footnote{Here and below, if $\Psi$ is an st-strategy then $\Psi^{stc}=\Psi$. Also, if for example $\T\in m(\Sigma)$ then $\Sigma_{\P_1, \T}$ is an $(\omega_1, \omega_1, \omega_1)$-st-strategy. The definition of $\Psi_{\sf{ex}}$ appeared in \rdef{lsa type}.}. 

It follows from $\Gamma$ fullness preservation that $\Sigma_1\not=\Lambda_1^{stc}$. Indeed, if $\Sigma_1=\Lambda_1^{stc}$ then  $\Q(c, \T)$  is a $\Sigma_1$-sts mouse over $\P_1$ with an iteration strategy in $\Gamma$\footnote{Recall that our sts indexing scheme indexes branches of $(\omega_1, \omega_1)$-iterations and not generalized stacks.}. Hence, $\T\in b(\Sigma)$ and $\Sigma(\T)=b$. 

Notice now that there is no low level disagreement between $\Sigma_1$ and $\Lambda_1^{stc}$ since if $(\T_1, \P_1, \T_2, \P_2)$ is a low level disagreement between $\Sigma_1$ and $\Lambda_1^{stc}$ then letting $E$ be the $\P_1^b=\P_2^b$-un-dropping extender of $\T_1$ and $\T_2$, $\T^\frown \T_1^\frown \{Ult(\S_\gg, E), E \}$ and $\T^\frown \T_2^\frown \{Ult(\S_\gg, E), E \}$ induce a low level disagreement between $\Sigma$ and $\Lambda$. 

Let $\U_1$ be a disagreement between $\Sigma_1$ and $\Lambda_1^{stc}$. Arguing as we have argued for $\T$, we get that $\pi^{\U_1, b}$ is defined. Let $\R_1$ be the least node of $\U_1$ such that $\R_1^b=\pi^{\U_1, b}(\P_1^b)$. It then follows that $(\U_1)_{\geq \R_1}$ is a stack on $\R_1$ that is above $\ord(\R_1^b)$. Let $\X_1$ be the last normal component of $(\U_1)_{\geq \R_1}$. It follows from $\Gamma$-fullness preservation that $\X_1$ doesn't have a fatal drop and $\m^+(\X_1)\models ``\d(\X_1)$ is a Woodin cardinal".  Set $\P_2=\m^+(\X_1)$. We now claim that \\

\textit{Claim.} $\U_1\in b(\Lambda_1^{stc})$. \\\\
\begin{proof} Assume that $\U_1\in m(\Lambda_1^{stc})$.  Because $\U_1$ is a disagreement between $\Sigma_1$ and $\Lambda^{stc}_1$, we must have that $\U_1\in b(\Sigma_1)$. Let $\U^*\in \dom(\Lambda_1)$ be such that $(\U^*)^{sc}=\U_1$\footnote{See \rdef{the short tree component domain 1}.}. It then follows that both $\Sigma(\T^\frown \U_1)$ and $\Lambda(\T^\frown \{c\}^\frown \U^*)$ are branches, and therefore, letting $b_1=\Sigma_1(\U_1)$ and $c_1=\Lambda_1(\U^*)$\footnote{Recall that $\Lambda_1$ is a strategy.}, we must have that $\Q(b_1, \X_1)=\Q(c_1, \X_1)$. Indeed, if $\Q(b_1, \X_1)\not =\Q(c_1, \X_1)$ then letting $\T_1=\X^\frown \U_1$ and $\T_2=\X^\frown \{c\} ^\frown \U^*$\footnote{In this iteration, player $I$ starts a new round of the iteration after player $II$ plays $c$. At the begining of this round, player $I$ drops to $\Q(c, \X)_{{\sf{ex}}}$.}, $\T_1$ and $\T_2$ are as desired. It follows that\\\\
(2) $b_1=c_1$. \\\\
We now have that because $\U_1\in m(\Lambda_1^{stc})$, letting $\b$ be the largest member of ${\sf{max}}^{\U_1}$,\\\\
(3) $\pi^{ \U^*}_{\b, c_1}$ is defined and $\pi^{ \U^*}_{\b, c_1}(\d^{\M_{\b}^{\U_1}})=\d^{\P_2}$.\\\\
 Because $\U_1\in b(\Sigma_1)$, we must have that\\\\
 (4) either $\pi^{\U_1}_{\b, b_1}$ is undefined or $\pi^{\U_1}_{\b, b_1}(\d^{\M_{\b}^{\U_1}})>\d^{\P_2}$. \\\\
 But because of (2)\\\\
 (5) $\pi^{\U_1}_{\b, b_1}$ is undefined  if and only if $\pi^{ \U^*}_{\b, c_1}$ is undefined, and if $\pi^{\U_1}_{\b, b_1}$ is defined then $\pi^{ \U^*}_{\b, c_1}(\d^{\M_{\b}^{\U_1}})=\pi^{\U_1}_{\b, b_1}(\d^{\M_{\b}^{\U_1}})$,\\\\
  as the calculation of both depends on the functions in $\M_\b^{\U_1}|\d^{\M_{\b}^{\U_1}}$. Clearly (2), (3), (4) and (5) contradict each other.
\end{proof}

Since $\U_1$ is a disagreement, we have that $\U_1\in m(\Sigma_1)$. Let then $c_1=\Lambda_1^{stc}(\U_1)$. Notice that\\\\
(6) $\Q(c_1, \X_1)$ exists and if $\U^*\in \dom(\Lambda_1)$ is such that $(\U^*)^{sc}=\U_1$ then either $\pi^{\U^*}_{c_1}$ is undefined or $\pi^{\U^*}_{c_1}(\d^{\P_1})>\d^{\P_2}$.\\\\
We now continue in the above manner by letting 
\begin{center}
$\Sigma_2=(\Sigma_1)_{\P_2, \U_1}$ and $\Lambda_2=((\Lambda_1)_{\Q(c_1, \X_1), \U^{*\frown} \{c_1\}})_{\sf{ex}}$.\end{center}
 Notice that $\Gamma$-fullness preservation once again implies that $\Sigma_2\not =\Lambda_2^{stc}$. By repeating in the above manner we obtain sequences $(\U^*_i: i\in [1, \omega))$, $(\Lambda_i: i\in [1, \omega))$ and $(c_i: i\in [1, \omega))$ such that the following conditions are satisfied:
\begin{enumerate}
\item $\U^*_1=\U^*$ where $\U^*$ is as in (6).
\item  For each $i<\omega$, $\U^*_{i}$ is according to $\Lambda_i$ and $c_i=\Lambda_i(\U_i^*)$.
\item For each $i<\omega$, $\U^*_i$ has a last normal component $\X_i$ and $\Q(c_i, \X_i)$ exists.
\item For each $i<\omega$, $\Lambda_{i+1}=((\Lambda_{i})_{\Q(c_i, \X_i), \U_i^{*\frown}\{c_i\}})_{\sf{ex}}$,
\item For each $i\in \omega$, letting $\d_{i+1}=\d(\X_i)$, either $\pi^{\U^*_i}_{c_i}$ is undefined or $\pi^{\U^*_i}_{c_i}(\d_i)>\d_{i+1}$. 
\end{enumerate}
Concatenating the $\U^*_i$s we get $\U$ according to $\Lambda_1$ without a well-founded branch.
\end{proof}

Let $\T_1$ and $\T_2$ be as in \rsublem{t is a branch}. Set $\U_1=(\T_{\leq \M_{\a_0}^\U})^\frown \T_1$, $\U_2=(\T_{\leq \M_{\a_0}^\U})^\frown \T_2$ $b=\Sigma(\U_1)$ and $c=\Lambda(\U_2)$. Let $\X_1$ and $\X_2$ be the last normal components of $\T_1$ and $\T_2$. It follows that $\m^+(\X_1)=\m^+(\X_2)$, both $\Q(b, \mathcal{X}_1)$ and $\Q(c, \X_2)$ exist and $\Q(b, \X_1)\not =\Q(c, \X_2)$. 

Let $\P_2=\m^+(\X)$. Notice that it follows from our smallness assumption on hod mice, namely that hod mice do not have lsa hod initial segments, that $\d^{\P_2}$ is a strong cutpoint of both $\Q(b, \X_1)$ and $\Q(c, \X_2)$. We then have that $\Q(b, \X_1)$ is a $\Sigma^{stc}_{\P_2, \U_1}$-sts mouse over $\P_2$, $\Q(c, \X_2)$ is a $\Lambda_{\P_2, \U_2}^{stc}$-sts mouse over $\P_2$, and the comparison of $\Q(b, \X_1)$ and $\Q(c, \X_2)$ does not halt (as otherwise we would have $\Q(b, \X_1)=\Q(c, \X_2)$).  Set $\nu=\d^{\P_2}$, $\M_0=\Q(b, \X_1)$ and $\M_1=\Q(c, \X_2)$. We now have that\\\\
(7) $\M_0\not \insegeq \M_1$, $\M_1\not\insegeq \M_0$, $\M_0||\nu=\M_1||\nu$, $\M_0$ and $\M_1$ are  $\nu$-sound and project to $\nu$, and $\nu$ is a strong cutpoint of both $\M_0$ and $\M_1$.\\
(8) $\M_0$ is a  $\Sigma^{stc}_{\P_2, \U_1}$-sts mouse over $\P_2$ and $\M_1$ is a $\Lambda_{\P_2, \U_2}^{stc}$-sts mouse over $\P_2$.\\
(9) The comparison of $\M_0$ and $\M_1$ cannot halt.\\\\
(9) holds as otherwise its failure implies that either $\M_0\insegeq \M_1$ or $\M_1\insegeq \M_0$, both of which are impossible (because of (7)). 

It follows that the comparison of $\M_0$ and $\M_1$ encounters disagreements involving strategies, as otherwise the usual comparison argument would imply that the comparison halts. Let $\Psi_0$ and $\Psi_1$ be the canonical strategies of $\M_0$ and $\M_1$ respectively. Thus, $\Psi_0$ witnesses that $\M_0$ is a $\Sigma^{stc}_{\P_2, \U_1}$-sts mouse, and $\Psi_1$ witnesses that $\M_1$ is a  $\Lambda_{\P_2, \U_2}^{stc}$-sts mouse. 

We can then find $\Psi_0$-iterate $\K_0$ of $\M_0$ and $\Psi_1$-iterate $\K_1$ of $\M_1$ such that $\K_0$ and $\K_1$ are produced via the usual extender comparison procedure (this implies that both iterations are above $\nu$) and for some $\a$, \\\\
(10) $\K_0|\a=\K_1|\a$, $\K_0||\a\not =\K_1||\a$, $\a\not \in \dom(\vec{E}^{\K_0})\cup  \dom(\vec{E}^{\K_1})$.\\\\
Notice that it follows from our indexing scheme (see \rdef{sts indexing scheme}) that there must be a branch indexed at $\a$ in both $\K_0$ and $\K_1$. Let then $t=(\P_2, \W, \P_3, \W') \in \K_0|\a$ be such that its branch is indexed at $\a$ in both $\K_0$ and $\K_1$. 

 We now have to analyze exactly what kind of stack $t$ is. Recall that our indexing scheme is so that we add branches for two kinds of stacks that we now list.\\\\
\textbf{Case 1.} $\W$ is a $\K_0|\a$-terminal tree\footnote{See \rdef{terminal tree}.} and $\W'$ is undefined.\\
\textbf{Case 2.} $\W'$ is defined and is a stack on $(\P_3)^b$. \\\\
We can immediately rule out case 1 above: $\K|\a=\N|\a$ and the branch of $\W$ just depends on $\K|\a$\footnote{See \rdef{weak psi alpha indexing scheme a}.}. On the other hand, case 2, just like in the proof of \rlem{sublemma about u},  leads to a low level disagreement between $\Sigma$ and $\Lambda$, which is contrary to our assumption. This contradiction implies that the comparison of $\M_0$ and $\M_1$ does not encounter strategy disagreement implying that (7) is false.  This contradiction also completes our proof of \rlem{disagreement implies low level disagreement}.
\end{proof}

\begin{lemma}[Minimal low level disagreement]\label{minimal low level disagreement exist} Suppose $\Gamma$ is a pointclass projectively closed pointclass, and $(\P, \Sigma)$ and $(\P, \Lambda)$ are allowable pairs such that  both $\Sigma$ and $\Lambda$ are almost $\Gamma$-fullness preserving. Suppose that one of the following conditions holds:
\begin{enumerate}
\item $\P$ is of limit type but not of lsa type, and $\Sigma\not =\Lambda$.
\item $(\P, \Sigma)$ and $(\P, \Lambda)$ are sts pairs or simple sts pairs, $\Sigma\not=\Lambda$ and both $\Sigma$ and $\Lambda$ are fullness preserving and are weakly self-cohering.  
\item $(\P, \Sigma)$ is a (simple) sts pair, $(\P, \Lambda)$ is a (simple) hod pair, $\Sigma\not=\Lambda$ and both $\Sigma$ and $\Lambda$ are fullness preserving and are weakly self-cohering. 
\item $\P$ is of lsa type, $\mathcal{J}_{\omega}[\P]\models ``\d^\P$ is not a Woodin cardinal" and  $\Sigma\not=\Lambda$.
\end{enumerate}
Then there is a minimal low level disagreement between $\Sigma$ and $\Lambda$. 
\end{lemma}
\begin{proof} Again we give the proof from clause 2. Assume there is no minimal low level disagreement between $\Sigma$ and $\Lambda$. It follows from \rlem{low level disagreement between strategies} that there is a low level disagreement between $\Sigma$ and $\Lambda$. Let $(\T_1, \P_1)\in I^b(\P, \Sigma)$ and $(\U_1, \R_1)\in I^b(\P, \Lambda)$ be a low level disagreement. Set $\Q=\P_1^b(=\R_1^b)$. We thus have that $\Sigma_{\Q, \T_1}\not =\Lambda_{\Q, \U_1}$. Notice that if $\Sigma_{\Q|\d^\Q, \T_1}=\Lambda_{\Q|\d^\Q, \U_1}$ then $\Gamma$-fullness preservation implies that $\Sigma_{\Q, \T_1}=\Lambda_{\Q, \U_1}$. Thus, there is $\b<\l^\Q$ such that\footnote{See \rnot{l p}.}
\begin{itemize}
\item $\Sigma_{\Q(\b), \T_1}\not=\Lambda_{\Q(\b), \U_1}$, 
\item $\ord(\Q(\b))$ is a cutpoint of $\Q$.
\end{itemize}
Let $\b_1$ be the least such ordinal and set $\Q_1=\Q(\b_1)$. If $\Q_1$ is of successor type then by minimality of $\b_1$ we get that $(\T_1, \U_1, \Q_1)$ is a minimal low level disagreement. Thus, we have that $\Q_1$ is limit type. The minimality of $\b_1$ then implies that\\\\
(1) $\Q_1$ is non-meek,\\
(2) $\Sigma_{\Q_1^b, \T_1}=\Lambda_{\Q_1^b, \U_1}$\footnote{Notice that these are ordinary strategies not generalized strategies. The reason is that $\Sigma_{\P'}$ for $\P'\inseg \P^b$ is an ordinary strategy}.\\\\
Applying \rlem{low level disagreement between strategies}, we get $(\T_2, \P_2, \U_2, \R_2)$ that constitute a low level disagreement between $\Sigma_{\Q_1, \T_1}$ and $\Lambda_{\Q_1, \U_1}$. Let $\b_2$ be the least $\b$ such that 
\begin{itemize}
\item $\P_2(\b)=\R_2(\b)$,
\item $\Sigma_{\P_2(\b), \T_1^\frown \T_2}\not=\Lambda_{\R_2(\b), \U_1^\frown \U_2}$, 
\item $\ord(\P_2(\b))$ is a cardinal of both $\P_2$ and $\R_2$.
\end{itemize}
Thus,\\\\
(3) $\P_2(\b_2)\inseg  \P_2^b$.\\\\
Set $\Q_2=\P_2(\b)$. We claim that\\

\textit{Claim.} $\Q_2$ is not of successor type. \\
\begin{proof} To see this, suppose that $\Q_2$ is of successor type. Let $\T_1=(\M_\a, \X_\a, F_\a: \a\leq \eta)$ and $\U_1=(\N_\a, \Y_\a, G_\a: \a\leq \eta)$. We have that $\M_\eta, F_\eta$, $\N_\eta$ and $G_\eta$ are undefined. Let $\X=\X_\eta^\frown \T_2$ and $\Y=\Y_\eta^\frown \U_2$. In forming $\X$, we let player $I$ start a new round on $\P_1$ by dropping to $\Q_1$. The same happens in $\Y$ as well. Let then $F_\eta$ be the $\Q_2$-un-dropping extender of $\X$ and $G_\eta$ be the $\Q_2$-un-dropping extender of $\Y$ and set $\X'=\X^\frown\{Ult(\M_\eta, F_\eta), F_\eta\}$ and $\Y'=\Y^\frown \{ Ult(\N_\eta, G_\eta), G_\eta\}$. Notice that $F_\eta=G_\eta$ as $\pi^{\T_1, b}=\pi^{\U_1, b}$ and the $\Q_2$-un-dropping extenders of $\T_2$ and $\U_2$ are the same. Because $\Sigma$ and $\Lambda$ are weakly self-cohering, we have that $\Sigma_{\Q_2, \X'}=\Sigma_{\Q_2, \X}$ and $\Lambda_{\Q_2, \Y'}=\Lambda_{\Q_2, \Y}$. Thus, $\Sigma_{\Q_2, \X'}\not =\Lambda_{\Q_2, \Y'}$ and hence, $(\X, \Y, \Q_2)$ is a minimal low level disagreement. 
\end{proof}

Continuing in this fashion we can now produce a sequence $(\P_i, \T_i, \Q_i: i\in [2, \omega))$ such that the following conditions hold.
\begin{enumerate}
\item For all $i\in [2, \omega)$, $\Q_i\inseg_{hod} \P^b_i$.
\item $\Q_i$ is non-meek.
\item $\T_i$ is a stack on $\Q_i$ such that $\pi^{\T_i, b}$ exists and $\P_{i+1}$ is the last model of $\T_i$. 
\end{enumerate}
Clearly, the concatenation of $\T_i$'s is an iteration according to $\Sigma_{\P_2, \X}$ without a well-founded branch.
\end{proof}

Next we introduce several definitions that will be useful in the sequel. 

\begin{definition}[Comparison stack]\label{comparison stack}\index{comparison stack} Suppose $(\P, \Sigma)$ and $(\Q, \Lambda)$ are two hod pairs or sts hod pairs. Then we say $(\T, \R, \U, \S)$ are comparison stacks for 
\begin{center}
$((\P, \Sigma), (\Q, \Lambda))$ 
\end{center}
with last models $(\R, \S)$
if $(\T, \R)\in I(\P, \Sigma)$, $(\U, \S)\in I(\Q, \Lambda)$, and either
\begin{enumerate}
\item $\S\in Y^\R$ and $\Sigma_{\S, \T}=\Lambda_{\S, \U}$.
\item $\R\in Y^\S$ and $\Sigma_{\R, \T}=\Lambda_{\R, \U}$.
\end{enumerate}
$\myqedhere$
\end{definition}

\begin{definition}[Agreement up to the top]\label{agreement up to the top}\index{agreement up to the top} Suppose $\P$ and $\Q$ are two hod premice of limit type. Then we say $\P$ and $\Q$ agree up to the top if $\P^b=\Q^b$. Suppose further that $\Sigma$ and $\Lambda$ are such that $(\P, \Sigma)$ and $(\Q, \Lambda)$ are two hod pairs or sts hod pairs. Then we say $(\P, \Sigma)$ and $(\Q, \Lambda)$ agree up to the top if $\P$ and $\Q$ agree up to the top and $\Sigma_{\P^b}=\Lambda_{\Q^b}$. $\myqedhere$
\end{definition}

\begin{definition}[Extender and strategy disagreement]\label{extender and strategy disagreement}\index{extender disagreement}\index{strategy disagreement}
Given two hod premice $\P$ and $\Q$ such that $\P\not =\Q$, we let $\b(\P, \Q)$ be the least ordinal $\gg$ such that $\P|\gg =\Q|\gg$ but $\P||\gg\not=\Q||\gg$. We say $\P$ and $\Q$ have an extender disagreement if $\b(\P, \Q) \in \dom(\vec{E}^\R)\cup\dom(\vec{E}^\Q)$. We say $\P$ and $\Q$ have a strategy disagreement if $\b(\P, \Q)\not \in \dom(\vec{E}^\R)\cup \dom(\vec{E}^\Q)$. In this case, we let 
\begin{center}
$\R_{\P, \Q}=\cup Y^{\P|\b(\P, \Q)}(=\cup Y^{\Q|\b(\P, \Q)})$
\end{center} 
Thus, both $\P$ and $\Q$ have a branch indexed at $\b(\P, \Q)$ for some $\T$ on $\R_{\P,\Q}$. We say $\R_{\P,\Q}$ is the disagreement layer of $\P$ and $\Q$. $\myqedhere$
\end{definition}

\begin{definition}[Extender comparison]\label{extender comparison}\index{extender comparison} Suppose that $(\P, \Sigma)$ and $(\Q, \Lambda)$ are two allowable pairs which agree up to the top. Then we say $(\T, \R, \U, \S)$ are the trees of the extender comparison of $(\P, \Sigma)$ and $(\Q, \Lambda)$ if 
\begin{enumerate}
\item $\T$ is according to $\Sigma$ and $\R$ is its last model,
\item  $\U$ is according to $\Lambda$ and $\S$ is its last model, and
\item $\T$ and $\U$ are obtained by using the usual extender comparison process (i.e., by removing the least extender disagreements) for comparing the top windows of $\P$ and $\Q$ until a strategy disagreement appears. 
\end{enumerate}
$\myqedhere$
\end{definition}

It follows that if in \rdef{extender comparison}, $\R\not =\S$ then $\R$ and $\S$ have a strategy disagreement.

\section{Self-cohering}

Here our goal is to show that the strategies appearing in hod pair constructions are self-cohering\footnote{See \rdef{self-cohering}.}.

\begin{theorem}\label{hpc strategy is self-cohering} Assume $\sf{AD}^+$. Suppose $\Gamma$ is a pointclass such that for some $\a$ with $\theta_\a<\Theta$, $\Gamma=\{ A\subseteq \bR: w(A)<\theta_\a\}$, ${\sf{C}}=(\mathbb{M}, (P, \Psi), \Gamma^*, A)$ Suslin, co-Suslin captures $\Gamma$ and $\mathbb{M}=(M, \d,  \vec{G}, \Sigma^*)$. Set
\begin{center} 
${\sf{hpc}}_{{\sf{C}}, \Gamma}^+=(\M_\gg , \N_\gg, Y_\gg, \Phi^+_\gg, F^+_\gg,  F_\gg, b_\gg: \gg\leq \d)$.
\end{center} 
Suppose $\b<\d$, $\P\in Y_\b$ and and $M\models ``(\P, (\Phi_\b)_\P)\in {\sf{HP}}^\Gamma$. Then $(\Phi_\b^+)_\P$ is self-cohering. 
\end{theorem}
 \begin{proof} Set $\Sigma=(\Phi_\b^+)_\P$. The hardest case is when $\P$ is non-meek and $\Sigma$ is generalized strategy. Suppose
\begin{itemize}
\item $\T=(\M_\a, \T_\a, F_\a: \a<\eta)$ is a generalized stack according to $\Sigma$, 
\item $\a_0, \a_1<\eta$,
\item $\xi_0< \lh(\T_{\a_0})$ and $\xi_1<\lh(\T_{\a_1})$, and
\item $\R\inseg_{hod}\M^{\T_{\a_0}}_{\xi_0}=_{def}\S_0$ and $\R\inseg_{hod}\M^{\T_{\a_1}}_{\xi_1}=_{def}\S_1$.
\end{itemize}
By absoluteness, we can find such a $\T$ in $M[g]$ where $g\subseteq Coll(\omega, \zeta_0)$ is $M$-generic and $\zeta_0<\d$. We want to see that
\begin{center}
$\Sigma_{\R, \T_{\leq \S_0}}=\Sigma_{\R, \T_{\leq \S_1}}$\footnote{See \rdef{id pullback initial segment}.}.
\end{center}
Assume then that $\Sigma_{\R, \T_{\leq \S_0}}\not =\Sigma_{\R, \T_{\leq \S_1}}$. Again, the hard case is when $\R$ is of limit type, and so we assume this.

 It follows from \rlem{minimal low level disagreement exist} that there is a minimal low level disagreement $(\U_0, \R_0, \U_1, \R_1, \Q)$ between $\Sigma_{\R, \T_{\leq \S_0}}$ and $\Sigma_{\R, \T_{\leq \S_1}}$. We thus have that
\begin{itemize}
\item $\Q$ is of successor type, 
\item $\Sigma_{\Q, (\T_{\leq \S_0})^\frown \U_0}\not =\Sigma_{\Q, (\T_{\leq \S_1})^\frown \U_1}$ and
\item $\Sigma_{\Q^-, (\T_{\leq \S_0})^\frown \U_0} =\Sigma_{\Q^-, (\T_{\leq \S_1})^\frown \U_1}$.
\end{itemize}
Let $\X_0=((\T_{\a_0})_{\leq \S_0})^\frown \U_0$ and $\X_1=((\T_{\a_1})_{\leq \S_1})^\frown \U_1$. Let $H_0$ and $H_1$ be the
 $\Q$-un-dropping extenders of $\X_0$ and $\X_1$\footnote{See \rdef{the un-dropping extender of a continuable stack}.}. Finally, let for $i\in 2$, $\Y_i=(\M^i_\xi, \T^i_\xi, F^i_\xi: \xi\leq \a_i+1)$ be the generalized stack that has the following properties:
 \begin{itemize}
 \item For $\a\leq \a_i$, $\M^i_\xi=\M_\xi$.
 \item For $\a< \a_i$, $\T_\a=\T^i_\a$ and $F^i_\xi=F_\xi$. .
 \item $\T_{\a_i}^i=\X_i$, $\M^i_{\a_i+1}=Ult(\M^i_{\a_i}, H_i)$ and $F_{\a_i}^i=H_i$. 
 \end{itemize}
 Let for $i\in 2$, $E_i$ be the $(\d^{\P^b}, \pi^{Y_i, b}(\d^{\P^b}))$-extender derived from $\pi^{\Y_i, b}$. Because $\Sigma$ is a weakly self-cohering\footnote{See \rlem{weakly self cohering}.}, we have that 
 \begin{itemize}
\item $\Q$ is of successor type, 
\item $\Sigma_{\Q, \Y_0}\not =\Sigma_{\Q, \Y_1}$ and
\item $\Sigma_{\Q^-, \Y_0} =\Sigma_{\Q^-, \Y_1}$.
\end{itemize}
 
Set $\k=\d^{\P^b}$ and $\zeta_1=\sup\{\lh(F_\gg^+):\gg\leq \b\}$. Set $\zeta=\max((\zeta_0^+)^M, (\zeta^+_1)^M)$. Let now $\N$ be the last model of 
\begin{center}
$({\sf{Le}}((\P|\k, \Sigma_{\P|\k}), \mathcal{J}_{\omega}[\P^b])_{>\zeta})^{(M, \d, \vec{G})}$.
\end{center} We have that $\N|(\k^+)^\N=\P^b$\footnote{See \rthm{fullness preservation of background constructions}, which implies that $\P$ is full.}. We then set for $i\in 2$, $\N_i=Ult(\N, E_i)$. Because for $i\in 2$, $\Y_i$ is an iteration of $\P$ according to the strategy induced by $\Sigma^*$, we have a $\Sigma^*$-iterate $M_i$ of $M$ such that letting $j_i: M\rightarrow M_i$ be the iteration embedding the following clauses hold:\\\\
(1) For each $i\in 2$, there is an $\M$-model $\M_i$ of ${\sf{hpc}}_{\Gamma}^{M_i}$ and an elementary embedding $\sigma_i: \Q\rightarrow \M_i$.\\
(2) For each $i\in 2$, $\M_i\inseg_{hod} j_i(\P)$\footnote{This follows from the fact that we are not allowed to project across $\d^{\M_i}$.}.\\
(3) For each $i\in 2$, there is an $\M$-model $\M_i'$ of ${\sf{hpc}}_{\Gamma}^{M_i}$ with index $\leq j_i(\b)$ such that $\M_i\in Y^{\M_i'}$ and $\Sigma_{\Q, \Y_i}$ is the $\sigma_i$-pullback of $\Phi_i$, where $\Phi_i$ is the strategy of $\M_i'$ induced by $\Sigma^*_{M_i}$.\\
(4) For each $i\in 2$, $\sigma_i$ extends to $\sigma_i^+: \N_i\rightarrow j_i(\N)$ and $j_i=\sigma_i^+\circ \pi_{E_i}$.\\\\
For each $i\in 2$, let $\Lambda_i^*$ be the strategy of $j_i(\N)$ induced by $\Sigma^*_{M_i}$ and let $\Lambda_i$ be the $\sigma_i^+$-pullback of $\Lambda_i^*$.  
 \begin{lemma}\label{strategy agreement footnote} For each $i\in 2$, $(\Lambda_i)_\Q=\Sigma_{\Q, \Y_i}$.
 \end{lemma}
 \begin{proof}
It is enough to show that $(\Phi_i)_{\M_i}=(\Lambda^*)_{\M_i}$. This follows easily from the fact that $\d^{\M_i}$\footnote{Recall that because $(\U, \Q)$ is a minimal disagreement, $\Q$ is of successor type. Thus, in fact $\d^{\M_i}$ is a Woodin cardinal of $j_i(\N_i)$ and $\M_i'$.} is a regular cardinal both in $\M_i'$ and in $j_i(\N)$. Because of this, both $(\Phi_i)_{\M_i}$ and $(\Lambda^*)_{\M_i}$ are the strategy of $\M_i$ induced by $\Sigma^*_{M_i|\tau_i}$ where $\tau_i$ is the least such that $\M_i$ is constructed inside $M_i|\tau_i$.
 \end{proof}
 We now let for $i\in 2$, $\W_i$ be the last model of 
 \begin{center}
 $({\sf{Le}}((\Q^-, (\Lambda_i)_{\Q^-}), \mathcal{J}_{\omega}[\Q])_{>\pi_{E_i}(\zeta)})^{(\N'_i, \d, \vec{K_i})}$
 \end{center}
 where $\N'_i=L_{\ord(M)}[\N_i]$ and $\vec{K}_i=\{ K\in \vec{E}^{\N_i}: \nu(K)$ is an inaccessible cardinal of $\N_i\}$.
 Let $\Omega_i$ be the strategy of $\W_i$ induced by $\Lambda_i$. We once again have that $(\Omega_i)_{\Q}=(\Lambda_i)_\Q$. Applying the ``furthermore" clause of \rthm{existence of thick sets} to $((\Omega_0)_{\Q}, \Q)$ and $((\Omega_1)_{\Q}, \Q)$, we get that $(\Omega_0)_{\Q}=(\Omega_1)_{\Q}$. However, since $(\Lambda_i)_\Q=\Sigma_{\Q, \Y_i}$ and $(\Omega_i)_{\Q}=(\Lambda_i)_\Q$, we have that 
 $(\Omega_0)_{\Q}\not=(\Omega_1)_{\Q}$. This contradiction completes the proof of \rthm{hpc strategy is self-cohering}.
 \end{proof}
 
 \section{Branch condensation}\label{proving branch condensation sec}

In this subsection we prove that the hod pair constructions produce strategies with branch condensation and in fact more. In order, however, to prove that hod pair constructions converge, we will need to establish the solidity and universality of the standard parameter of the models appearing in such constructions. Establishing such fine structural facts wasn't an issue in \cite{ATHM} as the fine structure for hod mice considered in that paper was a routine generalization of the fine structure theory developed in \cite{FSIT}. Here the matter is somewhat more complicated as the fine structure of non-meek hod mice cannot be viewed as a routine generalization of the fine structure of \cite{FSIT}. Nevertheless, the matter isn't too complicated as a simple generalization of branch condensation, \textit{strong branch condensation}, allows us to reduce our case to the one in \cite{FSIT}. In this subsection, we will establish that hod pair constructions produce strategies with strong branch condensation. The next definition will use concepts from \rnot{l p}, \rdef{one point extension}, \rdef{almost non-dropping gen stacks} and \rdef{allowable pair}.


\begin{definition}\label{induces a commuting diagram} Suppose $\P$ is a non-gentle hod premouse.

Suppose next that  either
\begin{itemize}
\item $\P$ is of successor type or 
\item $\P$ is of lsa type and $\mathcal{J}_{\omega}[\P]\models ``\d^\P$ is a Woodin cardinal".
\end{itemize}
Suppose $\sigma: \R\rightarrow \Q$ is an elementary embedding. We say that there is a \textbf{total} $(\Q, \R, \sigma)$-\textbf{b-condensation diagram} on $\P$ if there is $(\pi, \tau)$ such that
\begin{itemize}
\item $\pi:\P\rightarrow \Q$ is an elementary embedding,
\item $\tau:\P\rightarrow \R$ is an elementary embedding,
\item $\pi=\sigma\circ \tau$,
\end{itemize}
We then say that $(\pi, \tau)$ \textbf{supports} a total $(\Q, \R, \sigma)$-b-condensation diagram on $\P$.

Suppose next that $\P$ is of limit type and $\sigma: \R\rightarrow \Q$ is an elementary embedding. We say 
there is a \textbf{bottom-type} $(\Q, \R, \sigma)$-\textbf{b-condensation diagram} on $\P$ if there is $((\pi, \Q'), (\tau, \R'), \sigma')$ such that
\begin{itemize}
\item $\pi:\P^b\rightarrow \Q'$ is an elementary embedding,
\item $\tau:\P^b\rightarrow \R'$ is an elementary embedding,
\item $\sigma':\R'\rightarrow \Q'$ is an elementary embedding,
\item $\pi=\sigma'\circ \tau$,\\\\
and either
\item $\Q$ and $\R$ are of successor type, $\Q\insegeq_{hod} \Q'$, $\R\insegeq_{hod}\R'$ and $\sigma'\rest \R=\sigma$, or
\item $\Q$ and $\R$ are of of limit type, $\Q^b\insegeq_{hod} \Q'$, $\R^b\insegeq_{hod}\R'$ and $\sigma'\rest \R^b=\sigma\rest \R^b$.
\end{itemize}
We then say that $((\pi, \Q'), (\tau, \R'), \sigma')$ \textbf{supports} a bottom-type $(\Q, \R, \sigma)$-b-condensation diagram on $\P$. We say that $((\pi, \Q'), (\tau, \R'), \sigma')$ supports a \textbf{strict} bottom-type $(\Q, \R, \sigma)$-b-condensation diagram on $\P$ if in clause 6, $\Q^b\inseg_{hod}\Q'$.

Suppose now that $\P$ is as above and $(\P, \Sigma)$ is an allowable pair. We then say that there is a $(\P, \Sigma)$-supported $(\Q, \R, \sigma)$-\textbf{b-condensation diagram} on $\P$ if there is $(\T, \Q^*)\in I^{ope}(\P, \Sigma)$ such that one of the following clauses holds:
\begin{enumerate}
\item  \begin{itemize} 
\item $(\P, \Sigma)$ is a hod pair, 
\item $\P$ is either of successor type or of lsa type and such that $\mathcal{J}_{\omega}[\P]\models ``\d^\P$ is a Woodin cardinal", 
\item $\Q^*=\Q$,
\item  $(\T, \Q)\in I(\P, \Sigma)$, and 
\item there is $\tau:\P\rightarrow \R$ such that $(\pi^\T, \tau)$-supports a total $(\Q, \R, \sigma)$-b-condensation diagram on $\P$.
\end{itemize}
\item \begin{itemize}
\item $\P$ is of limit type\footnote{This clause also works for simple hod pairs and simple sts hod pairs.}
\item $\Q$ is a complete layer of $\Q^*$ and
\item letting 
\begin{center}
$E=\begin{cases}
E^\T_{\Q^b} &:\ \text{$\Q$ is of limit type}\\
E^\T_\Q &:\ \text{$\Q$ is of successor type},
\end{cases}$\end{center}
there are $\tau:\P^b\rightarrow \R'$ and $\sigma':\R'\rightarrow \Q'=_{def}\pi_E(\P^b)$ such that $((\pi_E\rest \P^b, \Q'), (\tau, \R'), \sigma')$ supports a bottom type $(\Q, \R, \sigma)$-b-condensation diagram on $\P$,
\item (\textbf{The sts conditions})\footnote{We need this conditions in order to make sense of $\sigma$-pullback of $\Sigma_{\Q, \T}$.} if $(\P, \Sigma)$ is an sts hod pair or a simple sts hod pair then $\Q\not =\Q^*$ provided one of the following holds:
\begin{itemize}
\item $\pi^\T$ exists.
\item $\pi^\T$ doesn't exist but letting $\T=(\P_\a, \X_\a, G_\a: \a\leq \b)$ and $\gg$ be the largest element of ${\sf{max}}^\X_\b$, $\pi^{(\X_\b)_{\geq \gg}}$ exists.
\end{itemize}
\end{itemize}
\end{enumerate}
$\myqedhere$
\end{definition}

\begin{definition}[Strong branch condensation]\label{strong branch condensation}\index{strong branch condensation} Suppose $(\P, \Sigma)$ is an allowable pair and $\P$ is not gentle. We say $\Sigma$ has \textbf{strong branch condensation with low-level-agreements} if 
\begin{enumerate}
\item $\Sigma$ has branch condensation\footnote{See \cite[Definition 2.14]{ATHM}. If $\Sigma$ is an st-strategy then we apply \cite[Definition 2.14]{ATHM} to stacks $\T$ and $\U$ such that ${\sf{max}}^\T={\sf{max}}^\U=\emptyset$.},
\item whenever 
\begin{itemize}
\item $(\T, \Q), (\U, \R)\in I^{ope}(\P, \Sigma)$,
\item $\pi:\R^b\rightarrow \Q^b$ is such that $\pi^{\T, b}=\pi\circ \pi^{\U, b}$,
\item $\X$ is a stack on $\R^b$ according to $\Sigma_{\R^b, \U}$,
\item $c$ is a branch of $\X$ such that $\pi^\X_c$ is defined and there is $\sigma: \M_c^\X\rightarrow \Q^b$ with the property that $\pi=\sigma\circ \pi^\X_c$,
\end{itemize}
 $c=\Sigma(\U^\frown \X)$. 
\item whenever $(\Q, \R, \sigma)$, $(\T, \Q^*)\in I^{ope}(\P, \Sigma)$ and $(\W, \R)\in B^{ope}(\P, \Sigma)\cup I^{ope}(\P, \Sigma)$ are such that
\begin{itemize}
\item there is a $(\P, \Sigma)$-supported  $(\Q, \R, \sigma)$-b-condensation diagram on $\P$ as witnessed by $(\T, \Q^*)$  and
\item letting $\Lambda$ be the $\sigma$-pullback of $\Sigma_{\Q, \T}$, there is no low level disagreement between $\Sigma_{\R, \W}$ and $\Lambda$,
\end{itemize}
then one of the following holds:
\begin{enumerate}
\item If \begin{itemize}
\item $\P$ is of lsa type, $(\P, \Sigma)$ is a hod pair and $\mathcal{J}_{\omega}[\P]\models ``\d^\P$ is a Woodin cardinal",
\item $(\T, \Q^*)$ supports a bottom-type  $(\Q, \R, \sigma)$-b-condensation diagram \begin{center}$((\pi, \Q'), (\tau, \R'), \sigma')$\end{center} on $\P$ and
\item  $\pi^{\T}$ is defined,
\end{itemize} then $\R$ is of lsa type and $(\Sigma_{\R, \W})^{stc}=\Lambda^{stc}$.
\item In all other cases, $\Sigma_{\R, \W}=\Lambda$.
\end{enumerate}
\end{enumerate}

We say $\Sigma$ has \textbf{strong branch condensation} if  $\Lambda=\Sigma_{\R, \W}$ holds without the requirement that  there is no low level disagreement between $\Lambda$ and $\Sigma_{\R, \W}$. $\myqedhere$
\end{definition}

\begin{remark}
The proof of \rthm{strong condensation for backgrounded strategies} only establishes clause 3 of strong branch condensation, but the proof can be easily modified to show clause 1 and 2 as well. $\myqedhere$
\end{remark}

The following is an easily provable lemma, which establishes the equivalence between strong branch condensation and strong branch condensation with low-level-agreements.  The reader may wish to review \rdef{weakly selfcohering}.

\begin{lemma}\label{no low level disagreement is not necessary} Suppose $(\P, \Sigma)$ is a hod pair or an sts hod pair, $\Sigma$ is weakly self-cohering and $\Gamma$ is a projectively closed pointclass. Suppose that 
\begin{itemize}
\item $\Sigma$ has strong branch condensation with low-level-agreements,
\item $\Sigma$ is $\Gamma$-strongly fullness preserving,
\item if $\P$ is of successor type then $\Sigma_{\P^-}$ has strong branch condensation. 
\end{itemize} 
Then $(\P, \Sigma)$ has strong branch condensation. 
\end{lemma}
\begin{proof} Suppose that $(\Q, \R, \sigma)$ is such that there is a $(\P, \Sigma)$-supported  $(\Q, \R, \sigma)$-b-condensation diagram on $\P$ as witnessed by $(\T, \Q^*)$. Let $\Lambda$ be the $\sigma$-pullback of $\Sigma_{\Q, \T}$. Fix a pair $(\W, \R)\in B^{ope}(\P, \Sigma)\cup I^{ope}(\P, \Sigma)$. Our goal is to argue that $\Sigma_{\R, \W}=\Lambda$. Towards a contradiction assume that $\Sigma_{\R, \W}\not =\Lambda$. Thus, we must have that there is a lower level disagreement between $\Sigma_{\R, \W}$ and $\Lambda$. 

Suppose first that $\P$ is of successor type. Because there is a lower level disagreement between $\Sigma_{\R, \W}$ and $\Lambda$, we must have that $\Sigma_{\R^-, \W}\not =\Lambda_{\R^-}$. However, it is not hard to see that there is a $(\P^-, \Sigma_{\P^-})$-supported  $(\Q^-, \R^-, \sigma\rest \R^-)$-b-condensation diagram on $\P^-$ as witnessed by $(\downarrow(\T, \P^-), \Q^-)$. Because $\sigma\rest \R^-$-pullback of $\Sigma_{\Q^-, \downarrow(\T, \P^-)}$ is just $\Lambda_{\R^-}$ and because $\Sigma_{\P^-}$ has strong branch condensation, we have that $\Sigma_{\R^-, \W} =\Lambda_{\R^-}$

We now assume that $\P$ is of limit type. Since all the cases are very similar, we will examine two representative cases, namely:\\\\
(A) $\P$ is of lsa type, $\mathcal{J}_{\omega}[\P]\models ``\d^\P$ is a Woodin cardinal", $(\P, \Sigma)$ is a hod pair and there is a total $(\P, \Sigma)$-supported  $(\Q, \R, \sigma)$-b-condensation diagram on $\P$  as witnessed by $(\T, \Q^*)$.\\
(B) $(\P, \Sigma)$ is an sts hod pair and $\Q$ is of limit type.\\\\
We start with (A). In this case, $\Q^*=\Q$ and $\pi^\T$ exists. Let $\tau:\P\rightarrow \R$ be such that $\pi^\T=\sigma\circ \tau$. Let then $(\W_1, \R_1, \W_1', \R_1', \R_2)$ be a minimal low level disagreement between $(\R, \Lambda)$ and $(\R, \Sigma_{\R, \W})$. Thus,
\begin{itemize}
\item (since $(\R, \Lambda)$ and $(\R, \Sigma_{\R, \W})$ are hod pairs), we have that $\W_1$ and $\W_1'$ are generalized stacks.
\item $\R_2$ is of successor type,
\item $\Lambda_{\R_2^-, \W_1}=\Sigma_{\R_2^-, \W^\frown \W_1'}$ and
\item $\Lambda_{\R_2, \W_1}\not =\Sigma_{\R_2, \W^\frown \W_1'}$.
\end{itemize}
 Let $\T_1=\sigma\W_1$ and let $\Q_1$ be the last model of $\T_1$. Let $\sigma_1:\R_1\rightarrow \Q_1$ be the copy map and set $\Q_2=\sigma_1(\R_2)$. Notice that both $\R_2$ and $\Q_2$ are of successor type. 
%
 
 Let  $\T=(\P_\a, \X_\a, G_\a: \a<\eta)$ and let $\U=\T^\frown \T_1$ where we construct this by setting $\P_\eta=\Q$ and $\X_{\eta+1}=\T_1$. Thus, $\U\rest \eta=\T$. Combining $\T$ and $\T_1$ this way is a legal way of producing a generalized stack because $\pi^\T$ is defined\footnote{If $(\P, \Sigma)$ was a simple hod pair then at this step we would let $\T^\frown \T_1$ be a stack.}.  Let then
 \begin{itemize}
  \item $E$ be the $\Q_2$-un-dropping extender $\U$ and $\Q'=\pi_E(\P^b)$,
 \item  $F'$ be the $\R_2$-un-dropping extender of $
 \W_1$,
 \item $F=\{ (a, A): (a, \tau(A))\in F'\}$,
 \item $\sigma_2: \R'\rightarrow \Q'$ be the map given by $\sigma_2([a, f]_F)=[\sigma_1(a), f]_E$.
 \end{itemize}
It follows that \\\\
(A1) $(a, X)\in F\iff (\sigma_1(a), X)\in E$, and hence $\sigma_2$ is an elementary embedding, and\\
(A2) $\sigma_2\rest \R_2=\sigma_1\rest \R_2$ and $\pi_E\rest \P^b=\sigma_2\circ \pi_F\rest \P^b$. \\\\
Therefore, $(\T^\frown \T_1, \Q_1)$ and $((\pi_E\rest \P^b, (\Q_2')^b), ((\pi_F\rest \P^b, (\R_2')^b), \sigma_2)$ support a bottom-type $(\Q_2, \R_2, \sigma_1\rest \R_2)$-b-condensation diagram on $\P$. Therefore, since $\Lambda_{\R_2^-, \W_1}=\Sigma_{\R_2^-, \W^\frown \W_1'}$ (i.e. there is no low level disagreement between $\Lambda_{\R_2, \W_1}=\Sigma_{\R_2, \W^\frown \W_1'}$) and $\Lambda_{\R_2, \W_1}$ is the $\sigma_2$-pullback of $\Sigma_{\Q_2, \T^\frown \T_1}$, $\Lambda_{\R_2, \W_1}=\Sigma_{\R_2, \W^\frown \W_1'}$, contradiction!\\
 
 We now work assuming (B). Most of what we say below is very similar to the above with only minor differences. In this case we have that clause 2 of \rdef{induces a commuting diagram} holds. Thus,
 \begin{itemize}
\item $\Q$ is a complete layer of $\Q^*$,
\item $\Q$ is of limit type and 
\item letting $E=E^\T_{\Q^b}$, there are $\tau:\P^b\rightarrow \R'$ and $\sigma':\R'\rightarrow \Q'=_{def}\pi_E(\P^b)$ such that $((\pi_E\rest \P^b, \Q'), (\tau, \R'), \sigma')$ supports a bottom type $(\Q, \R, \sigma)$-b-condensation diagram on $\P$.
\end{itemize}
Set $\pi=\pi_E\rest \P^b$. Let then $(\W_1, \R_1, \W_1', \R_1', \R_2)$ be a minimal low level disagreement between $(\R, \Lambda)$ and $(\R, \Sigma_{\R, \W})$. Let $\T_1=\sigma\W_1$ and let $\Q_1$ be the last model of $\T_1$. Let $\sigma_1:\R_1\rightarrow \Q_1$ be the copy map and let $\Q_2=\sigma_1(\R_2)$. Notice that both $\R_2$ and $\Q_2$ are of successor type.

We now define $\U$ as follows. If $\pi^\T$ is defined then we let $\U=\T^\frown \T_1$ be as in (A). Assume then $\pi^\T$ is not defined. In this case, $\T=(\P_\a, \X_\a, G_\a: \a\leq \b)$, $\Q^*$ is the last model of $\X_\b$ and $\pi^{\X_\b}$ is not defined. Let then $\U$ be the same as $\T$ except that the $\b$th stack used in $\U$ is $\X_\b^\frown \T_1$\footnote{Thus, we have that for some $\gg\in R^{\X_\b^\frown \T_1}$, $\Q^*=\M_\gg^{\X_\b^\frown \T_1}$ and $(\omega\b_\gg^{\X_\b^\frown \T_1}, m_\gg^{\X_\b^\frown \T_1})=(\ord(\Q), \omega)$.}. 

Just like in case (A), we have that if $E'$ is the $\Q_2$-un-dropping extender of $\U$, $F'$ is the $\R_2$-un-dropping extender of $\W_1$, $F=\{ (a, A): (a, \tau(A))\in F'\}$ and  $\sigma_2:  \pi_{F}(\P^b)\rightarrow \pi_{E'}(\P^b)$ is the map given by $\sigma_2([a, f]_F)=[\sigma_1(a), f]_E$ then\\\\
(B1) $(a, X)\in F\iff (\sigma_1(a), X)\in E'$, and hence $\sigma_2$ is an elementary embedding, and\\
(B2) $\sigma_2\rest \R_2=\sigma_1\rest \R_2$ and $\pi_{E'}\rest \P^b=\sigma_2\circ \pi_F\rest \P^b$.  \\\\
Here the situation may seem somewhat more complicated as $\W_1$ is on $\R$ and not on $\R'$. But since $\R^b\insegeq_{hod}\R'$, $F'$ is an $\R'$-extender. Moreover, since $\T_1=\sigma\W_1$ and $\sigma'\rest \R^b=\sigma\rest \R^b$, we have that for each $A\in \powerset(\d^{\P^b})\cap \P$, 
\begin{center}
$\sigma_2(\sigma^{\W_1}(\tau(A)))=\sigma^{\T_1}(\pi(A))$\footnote{Here, $\sigma^\X$ is defined in \rdef{the un-dropping extender of a continuable stack} and \rnot{notation for generalized stacks}.}.
\end{center}
We then once again, just like in (A), have that 
\begin{center}
$(\T^\frown \T_1, \Q_1)$ and $((\pi_{E'}\rest \P^b, \pi_{E'}(\P^b)), ((\pi_F\rest \P^b, \pi_F(\P^b)), \sigma_2)$ 
\end{center}
support a bottom-type $(\Q_2, \R_2, \sigma_1\rest \R_2)$-b-condensation diagram on $\P$, and which implies, just like in (A), that $\Lambda_{\R_2, \W_1}=\Sigma_{\R_2, \W^\frown \W_1'}$.
\end{proof}

\begin{theorem}\label{strong condensation for backgrounded strategies} Assume $\sf{AD}^++{\sf{NsesS}}$. Suppose
\begin{itemize}
\item  for some $\a_0$ such that $\theta_{\a_0}<\Theta$, $\Gamma=\{ A\subseteq \bR: w(A)<\theta_{\a_0}\}$,
\item  ${\sf{C}}=(\mathbb{M}, (P, \Psi), \Gamma^*, A)$ Suslin, co-Suslin captures $\Gamma$,
\item $\mathbb{M}=(M, \d,  \vec{G}, \Sigma^*)$, 
\item ${\sf{hpc}}=(\M_\gg , \N_\gg, Y_\gg, \Phi_\gg, F^+_\gg, F_\gg, b_\gg : \gg\leq \delta)$ is the output of the $\Gamma-{\sf{hpc}}$ of $\mathbb{M}$,
\item $\xi<\d$ is such that $(\M_\xi, \Phi_\xi^+)$ is a hod pair, $\M_\xi$ is not gentle and $M\models (\M_\xi, \Phi_\xi)\in {\sf{Hp}}^{\Gamma}$.
\end{itemize}
 Then $\Phi_\xi^+$ has strong branch condensation. 
 
 Also if $\xi<\d$ is such that $(\M_\xi, (\Phi_\xi^+)^{stc})$ is an sts pair and $M\models (\M_\xi, \Phi^{stc}_\xi)\in {\sf{Hp}}^{\Gamma}$ then $(\Phi_\xi^+)^{stc}$ has strong branch condensation.
\end{theorem}
\begin{proof} The proof of the second half of the theorem is similar to the first and so we will prove the first and leave the second to the reader. The proof of the branch condensation is very similar to the proof of the second half of strong branch condensation, and so we give the proof of the second half of strong branch condensation.

 Towards a contradiction, suppose that for some $\xi'$, $\M_{\xi'}$ is a hod premouse and $\Phi^+_{\xi'}$ doesn't have strong branch condensation, and  let $\xi$ be the least such $\xi'$. Because of \rlem{no low level disagreement is not necessary}, it is enough to show that $\Phi_\xi^+$ has strong branch condensation with low-level-agreements. 

 Just like in the proof of fullness preservation (see \rthm{fullness preservation of background constructions}), if $\Phi^+_\xi$ does not have strong branch condensation then for some $\zeta_0<\d$ the witness can be found in some $M[g]$ where $g\subseteq Coll(\omega, \zeta_0)$ is $M$-generic. Let $\zeta_1=\{ \sup(F_\gg^+): \gg<\xi\}$ and set $\zeta=\max((\zeta_0^+)^M, (\zeta_1^+)^M)$. 

Let $\P=\M_\xi$ and $\Sigma=\Phi^+_\xi$. The difficult case is when $\P$ is non-meek, and so we assume this. We start working in $M[g]$. What we need to show is that whenever 
\begin{itemize}
\item $(\Q, \R, \sigma)$ is such that there is a $(\P, \Sigma)$-supported  $(\Q, \R, \sigma)$-b-condensation diagram on $\P$ as witnessed by $(\T, \Q^*)$, and 
\item $(\W, \R)\in B^{ope}(\P, \Sigma)\cup I^{ope}(\P, \Sigma)$ is such that letting $\Lambda$ be the $\sigma$-pullback of $\Sigma_{\Q, \T}$, there is no low level disagreement between $\Sigma_{\R, \W}$ and $\Lambda$,
\end{itemize}
then $\Sigma_{\R, \W}=\Lambda$.

Fix then $(\Q, \R, \sigma)\in M|\zeta[g]$ such that there is a $(\P, \Sigma)$-supported  $(\Q, \R, \sigma)$-b-condensation diagram on $\P$ as witnessed by $(\T, \Q^*)\in M|\zeta[g]$ and let $(\W, \R)\in B^{ope}(\P, \Sigma)\cup I^{ope}(\P, \Sigma)$ be  such that $(\W, \R)\in M|\zeta[g]$ and letting $\Lambda$ be the $\sigma$-pullback of $\Sigma_{\Q, \T}$,\\\\
(1) there is no low level disagreement between $\Sigma_{\R, \W}$ and $\Lambda$ but $\Sigma_{\R, \W}\not =\Lambda$. \\\\
It follows from \rlem{minimal low level disagreement exist} that\\\\
(2) either $\R$ is of successor type or $\R$ is of lsa type and $\mathcal{J}_{\omega}[\R]\models ``\d^\R$ is a Woodin cardinal".\\\\
%
%
\textbf{Case 1:}  $\pi^\T$ is defined, and for some $\tau$, $(\pi^\T, \tau)$ supports a total $(\Q, \R, \sigma)$-b-condensation diagram on $\P$.\\\\
We thus have that $\tau:\P\rightarrow \R$ and $\pi^\T=\sigma\circ \tau$. Let then
\begin{center}
$\Sigma'=\begin{cases}
\Sigma_{\P^-} &: \text{$\P$ is of successor type}\\
\Sigma^{stc} &: \text{otherwise}
\end{cases}$
\end{center}
\begin{center}
$\P'_0=\begin{cases}
\P &: \text{$\P$ is of successor type}\\
(\P|\d^\P)^\#&: \text{otherwise}
\end{cases}$
\end{center}
\begin{center}
$\P'_1=\begin{cases}
\P^- &: \text{$\P$ is of successor type}\\
(\P|\d^\P)^\#&: \text{otherwise}
\end{cases}$
\end{center}
\begin{center}
$\Lambda'=\begin{cases}
\Lambda_{\R^-} &: \text{$\P$ is of successor type}\\
\Lambda^{stc} &: \text{otherwise}.\\
\end{cases}$
\end{center}
Let $\P^+$ be the last model of 
\begin{center}
$({\sf{Le}}((\P'_1, \Sigma'), \mathcal{J}_\omega[\P'_0])_{>\zeta})^{(M[g], \d, \vec{G})}$\footnote{See \rdef{the stack in scb}. Here we are assuming that if $\P$ is of lsa type then the above construction doesn't break down because of the anomaly stated in clause 3.b of \rdef{fully backgrounded sts construction}. In the sequel, we will prove that such constructions indeed converge. See \rthm{sts fb constructions converge}.}. 
\end{center}
Define $\Q'_1$ and $\R'_1$ the same way $\P_1'$ is defined. We then let $\Sigma'_{\Q'_1}$ and $\Sigma'_{\R'_1}$ be defined the same way $\Sigma'$ is defined but relative to $\Sigma_{\Q, \T}$ and $\Sigma_{\R, \W}$. It follows from (2) that  $\Lambda'=\Sigma'_{\R'_1}$.

 Let 
 \begin{itemize}
 \item $E$ be the $(\d^\P, \d^\R)$-extender derived from $\tau$, 
 \item $F$ be the $(\d^\P, \d^\Q)$-extender derived from $\pi^\T$ and  
 \item $H$ be the $(\d^\P, \d^{\R})$-extender derived from $\pi^{\W}$\footnote{In the case that $(\W, \R)\in B^{ope}$, we let $H$ be the $\R$-un-dropping extender of $W$ and continue as below.}.
 \end{itemize}
  We let 
 \begin{center}
 $\R^+=Ult(\P^+, E)$, $\Q^+=Ult(\P^+, F)$ and $\S^+=Ult(\P^+, H)$. 
 \end{center}
 We also have $\sigma^+:\R^+\rightarrow \Q^+$ such that 
 \begin{center}
 $\pi^{\P^+}_F=\sigma^+\circ \pi^{\P^+}_E$ and $\sigma^+\rest \R=\sigma$. 
 \end{center}
 More precisely, $\sigma^+(x)=\pi^{\P^+}_F(f)(\sigma(a))$ where $f\in \P^+$, $a\in (\R)^{<\omega}$ and $x=\pi^{\P^+}_E(f)(a)$.
 
Notice now that both $\Q^+$ and $\S^+$ have strategies induced by $\Sigma^*$ via the resurrection procedure of \cite[Chapter 12]{FSIT} that we have outlined in \rlem{summary res and emb}. Let $\Psi^*$ and $\Phi$ be these strategies. We then have that $\Psi^*_\Q=\Sigma_{\Q, \T}$ and $\Phi_\R=\Sigma_{\R, \W}$.  Let now $\Psi$ be the $\sigma^+$-pullback of $\Psi^*$. Applying the ``furthermore" clause of \rthm{existence of thick sets} to $(\R^+, \Psi)$ and $(\S^+, \Phi)$,  we conclude that $\Psi_\R=\Phi_\R$. \\\\
\textbf{Case 2:} There is $\tau:\P\rightarrow \R'$ such that $(\tau, \R')\in M|\zeta[g]$ and letting
\begin{itemize}
\item $F$ be the $\Q^b$-un-dropping extender of $\T$ if $\Q$ is of limit type and 
\item $F$ be the $\Q$-un-dropping extender otherwise,
\end{itemize}
 there is a $\sigma':\R'\rightarrow\Q'=_{def}\pi_{F}(\P^b)$ such that $((\pi_{F}\rest \P^b, \Q'), (\tau, \R'), \sigma')$ supports a $(\Q, \R, \sigma)$-b-condensation diagram on $\P$.\\\\
This case is very similar. Notice that (2) implies that $\Sigma_{\R, \W}^{stc}=\Lambda^{stc}$ assuming the hypothesis of clause 2 of \rdef{strong branch condensation} holds (as there are no low level disagreements between $\Sigma_{\R, \W}$ and $\Lambda$). Thus, we assume that the hypothesis of clause 2 is not applicable. Notice now that the sts conditions\footnote{See \rdef{induces a commuting diagram}.} and the fact that $\Q$ has to be a complete layer of $\Q^*$ imply that $\Q$ and $\R$ are not of lsa type. (2) then implies that $\Q$ and $\R$ must be of successor type\footnote{See \rlem{disagreement implies low level disagreement}.}. It then follows that $\Q\inseg_{hod}\pi_{F}(\P^b)$.

 Let now $\P'$ be the last model of
\begin{center}
$({\sf{Le}}((\P^b, \Sigma_{\P^b}), \mathcal{J}_\omega[\P^b])_{>\zeta})^{(M[g], \d, \vec{G})}$.
\end{center}
Let  \begin{itemize}
 \item $E$ be the $(\d^{\P^b}, \d^{\R'})$-extender derived from $\tau$,  and  
 \item $H$ be the $\R$-un-dropping extender of $\W$.
 \end{itemize}
\begin{itemize}
\item $\R_1=Ult(\P', E)$ and $\pi_0=\pi_E^{\P'}$,
\item $\Q_1=Ult(\P', F)$ and $\sigma_1:\R_1\rightarrow \Q_1$ is the canonical factor map, and 
\item $\S_1=Ult(\P', H)$ and $k=\pi_H^{\P'}$.
\end{itemize}

We then define $\R^+$ and $\S^+$ as follows. Let $\R^+$ be the last model of  
\begin{center}
$({\sf{Le}}((\R^-, \Lambda_{\R^-}), \mathcal{J}_\omega[\R])_{>\zeta})^{(L_{\ord(M)}[\R_1], \d, \vec{H}')}$
\end{center}
where $\vec{H}'=\{ K\in \vec{E}^{\R_1}: \nu(K)$ is inaccessible in $\R_1\}$. 

Let $\S^+$ be the last model of  
\begin{center}
$({\sf{Le}}((\R^-, \Sigma_{\R^-, \W}), \mathcal{J}_\omega[\R])_{>\zeta})^{(L_{\ord(M)}[\S_1], \d, \vec{H}'')}$
\end{center}
where $\vec{H}''=\{ K\in \vec{E}^{\S_1}: \nu(K)$ is inaccessible in $\S_1\}$. 

Notice that because $\Sigma_{\R^-, \W}=\Lambda_{\R^-}$ we have that both $\R^+$ and $\S^+$ are $\Sigma_{\R^-, \W}$-mice over $\R$. Once again, in $M[g]$, both $\R^+$ and $\S^+$ have $(\d, \d+1)$-iteration strategies $\Phi$ and $\Psi$ such that $\Lambda_\R=\Phi_\R$ and $\Sigma_{\R, \W}=\Psi_\R$. It then again follows from the ``furthermore" clause of \rthm{existence of thick sets} that $\Phi_\R=\Psi_\R$. 
\end{proof}

\section{Positional and commuting}\label{sec positional and commuting}

In this section, our goal is to show that strong branch condensation implies commuting. Recall \cite[Definition 2.35]{ATHM}: if $M$ is a transitive model of a fragment of $\sf{{\sf{ZFC}}}$ and $\Sigma$ is an iteration strategy for $M$ then we say $\Sigma$ is \textit{positional} if whenever $Q$ is a $\Sigma$-iterate of $M$ via $\W$ and $(\T, R), (\U, R) \in I(Q, \Sigma_{Q, \W})$, $\Sigma_{R, \W^\frown \T}=\Sigma_{R, \W^\frown \U}$. Recall that commuting means that in the above scenario, $\pi^{\T}=\pi^{\W}$. If the above only holds for $Q=M$, then we say that $\Sigma$ is \textit{weakly positional} (and \textit{weakly commuting} respectively). Using the usual proof of the Dodd-Jensen lemma, we get that (weakly) positional implies (weakly) commuting.  

\begin{remark}\label{gentle branch cond} In the previous section, we only studied branch condensation for non-gentle hod premice. This is because if $(\P, \Sigma)$ is a hod pair and $\P$ is gentle then $\Sigma$ is essentially $\oplus_{\Q\inseg_{hod}\P}\Sigma_{\Q}$. Thus, we can say that $\Sigma$ has strong branch condensation if for every complete layer $\Q\inseg_{hod}\P$, $\Sigma_\Q$ has strong branch condensation. We will state our theorems for hod pairs or sts hod pairs, but the proofs will be given for pairs $(\P, \Sigma)$ such that $\P$ is non-gentle. $\myqedhere$
\end{remark}

\begin{proposition}\label{positional} Suppose $(\P, \Sigma)$ is an allowable pair, $\Gamma$ is a projectively closed pointclass and $\Sigma$ has strong branch condensation and is strongly $\Gamma$-fullness preserving. Then $\Sigma$ is positional.  Moreover, if $\Sigma$ is an iteration strategy then it is also commuting.
\end{proposition}
\begin{proof} We just prove weak positionality and hence weak commuting. The proof of the general case is only notationally more complicated. 

Suppose $(\T, \Q), (\U, \Q)\in I(\P, \Sigma)$. We want to see that $\Sigma_{\Q, \T}=\Sigma_{\Q, \U}$. Towards a contradiction, suppose not. Suppose first that $\P$ is of limit type and if it is of the lsa type then $\Sigma_{\Q, \T}^{stc}\not=\Sigma_{\Q, \U}^{stc}$. Let then $((\T_1, \R_1), (\T_2, \R_2), \R_3)$ a  a minimal lower level disagreement\footnote{See \rlem{minimal low level disagreement exist}.} between $\Sigma_{\Q, \T}$ and $\Sigma_{\Q, \U}$. We can then apply strong branch condensation to $(\R_3, \R_3, id)$. Notice that $(\T^\frown \T_1, \Q)$ supports a $(\R_3, \R_3, id)$-b-condensation diagram on $\P$ as witnessed by $((\pi, \R), (\pi, \R), id)$ where letting $E$ be the $\R_3$-un-dropping extender of $\T^\frown \T_1$, $\R=\pi_E(\P^b)$ and $\pi=\pi_E\rest \P^b$. 

Next suppose that $\P$ is of successor type or of lsa type but $\Sigma_{\Q, \T}^{stc}=\Sigma_{\Q, \U}^{stc}$. It then follows that $(\pi^\T, \pi^\T)$ supports a total $(\Q, \Q, id)$-b-condensation-diagram on $\P$. It then again follows that $\Sigma_{\Q, \T}=\Sigma_{\Q, \U}$.
\end{proof}

The proof actually gives more.

\begin{proposition}\label{positional1} Suppose $(\P, \Sigma)$ is an allowable pair, $\Gamma$ is a projectively closed pointclass and $\Sigma$ has strong branch condensation and is strongly $\Gamma$-fullness preserving. Suppose that $(\T, \Q)\in B^{ope}(\P, \Sigma)\cup I^{ope}(\P, \Sigma)$ and $(\U, \Q)\in B^{ope}(\P, \Sigma)\cup I^{ope}(\P, \Sigma)$. Then $\Sigma_{\Q, \T}=\Sigma_{\Q, \U}$.
\end{proposition}
%
%

\begin{definition}\label{sigmaq} Suppose $(\P, \Sigma)$ is an allowable pair and $\Gamma$ is a projectively closed pointclass. Suppose further that $\Sigma$ has strong branch condensation and is strongly $\Gamma$-fullness preserving. 
Given $\Q\in p[I^{ope}(\P, \Sigma)]\cup p[B^{ope}(\P, \Sigma)]$\footnote{$p[A]$ is the projection of $A$. In general, the coordinate onto which we are projecting will be clear from the context.}, we let $\Sigma_\Q=\Sigma_{\Q, \T}$ where $\T$ is such that $(\T, \Q)\in I^{ope}(\P, \Sigma)\cup B^{ope}(\P, \Sigma)$. $\myqedhere$
\end{definition}

We need commuting not only for iteration strategies but also for short tree strategies. 

\begin{definition}\label{positional and commuting for sts pairs} Suppose $(\P, \Sigma)$ is an sts hod pair. We say $\Sigma$ is \textbf{weakly commuting} if whenever $(\T, \Q)\in I^b(\P, \Sigma)$ and $(\U, \R)\in I^b(\P, \Sigma)$ are such that $\R^b\insegeq_{hod} \Q^b$ and $\R^b=cHull^{\Q^b}(\pi^{\T, b}[\P^b]\cup \d^{\R^b})$, then letting
\begin{itemize}
\item $k':Hull^{\Q^b}(\pi^{\T, b}[\P^b]\cup \d^{\R^b})\rightarrow \R^b$ be the transitive collapse and
\item $k=_{def}k'\circ \pi^{\T, b}:\P^b\rightarrow \R^b$,
\end{itemize}
 $k=\pi^{\U, b}$.

In the above situation we say that $k$ is the \textbf{collapse} of $\pi^{\T, b}[\P^b]$. We say $(\P, \Sigma)$ is \textbf{commuting} if whenever $(\T, \Q)\in I^{ope}(\P, \Sigma)$, $\Sigma_{\Q}$ is weakly commuting. $\myqedhere$
\end{definition}

It is not known to us if strong branch condensation and $\Gamma$-fullness preservation for sts pairs implies commuting. Nevertheless, hod pair constructions produce sts pairs that are commuting (also see \rprop{commuting from strong branch condensation}).
\begin{theorem}\label{commuting for backgrounded st strategies} Assume $\sf{AD}^++{\sf{NsesS}}$. Suppose
\begin{itemize}
\item  for some $\a_0$ such that $\theta_{\a_0}<\Theta$, $\Gamma=\{ A\subseteq \bR: w(A)<\theta_{\a_0}\}$,
\item  ${\sf{C}}=(\mathbb{M}, (P, \Psi), \Gamma^*, A)$ Suslin, co-Suslin captures $\Gamma$,
\item $\mathbb{M}=(M, \d,  \vec{G}, \Sigma^*)$, 
\item ${\sf{hpc}}=(\M_\gg , \N_\gg, Y_\gg, \Phi_\gg, F^+_\gg, F_\gg, b_\gg : \gg\leq \delta)$ is the output of the $\Gamma-{\sf{hpc}}$ of $\mathbb{M}$,
\item $\xi<\d$ is such that $(\M_\xi, \Phi_\xi^+)$ is a hod pair, $\M_\xi$ is not gentle and $M\models (\M_\xi, \Phi_\xi)\in {\sf{Hp}}^{\Gamma}$.
\end{itemize}
Suppose $\xi<\d$ is such that $(\M_\xi, (\Phi_\xi^+)^{stc})$ is an sts pair and $M\models (\M_\xi, \Phi^{stc}_\xi)\in {\sf{Hp}}^{\Gamma}$. Then $(\Phi_\xi^+)^{stc}$ is commuting.
\end{theorem}
\begin{proof} The proof is very similar to the proofs we have already given. Set $(\P, \Sigma)=(\M_\xi, (\Phi_\xi^+)^{stc})$. We prove that $\Sigma$ is weakly commuting as the general case is only notationally more complicated. Towards a contradiction assume $\Sigma$ is not weakly commuting. There is then some $\zeta_0$ such that for some $g\subseteq Coll(\omega, \zeta_0)$, $M[g]$ has a witness to the fact that $\Sigma$ is not weakly commuting. Let $\zeta_1=\sup\{ \lh(F_\iota): \iota<\xi\}$. Let $\zeta=\max((\zeta_0^+)^M, (\zeta_1^+)^M)$ and $\vec{G}'=\{ H\in \vec{G}: \cp(H
)>\zeta\}$. 

Fix then $(\T, \Q)\in I^b(\P, \Sigma)$ and $(\U, \R)\in I^b(\P, \Sigma)$ such that 
\begin{itemize}
\item  $\R^b\insegeq_{hod} \Q^b$,
\item  $\R^b=cHull^{\Q^b}(\pi^{\T, b}[\P^b]\cup \d^{\R^b})$ and
\item $(\T, \Q, \U, \R)\in M[g]$.
\end{itemize} 
Let $k:\P^b\rightarrow \R^b$ be the collapse of $\pi^{\T, b}[\P^b]$. We want to see that $k=\pi^{\U, b}$.

Let $\P'$ be the last model of $({\sf{Le}}((\P^b, \Sigma_{\P^b}), \mathcal{J}_{\omega})_{>\zeta})^{(M[g], \d, \vec{G}')}$ and $\P^+={\sf{stack}}(\P', \Sigma_{\P^b})$. Next let $E_0$ be the $(\d^{\P^b}, \d^{\R^b})$-extender derived from $k$ and $E_1$ be the $(\d^{\P^b}, \d^{\R^b})$-extender derived from $\pi^{\U, b}$. We need to show that\\\\
(a) $\pi_{E_0}\rest \P^b=\pi_{E_1}\rest \P^b$. \\\\
Notice that we have that\\\\
(1) $\pi_{E_0}\rest \d^{\P^b}=\pi_{E_1}\rest \d^{\P^b}$.\\\\
 This is because $k\rest \d^{\P^b}=\pi^{\T, b}\rest \d^{\P^b}$, so if $\pi_{E_0}\rest \d^{\P^b}\not = \pi_{E_1}\rest \d^{\P^b}$, we have that for some $\S\inseg_{hod}\P$, $\S$ is a complete layer of $\P$ and $\Sigma_\S$ is not commuting. But since $\Sigma_\S$ has strong branch condensation and is fullness preserving\footnote{See \rthm{strong condensation for backgrounded strategies} and \rthm{fullness preservation of background constructions}.}, \rprop{positional} implies that $\Sigma_\S$ is commuting. 

Let then $\N_0=Ult(\P^+, E_0)$ and $\N_1=Ult(\P^+, E_1)$. Notice that it follows from \rthm{strong condensation for backgrounded strategies} that\\\\
(2) both $\N_0$ and $\N_1$ are $\Sigma_{\R^b}$-mice over $\R^b$,\\
(3) for $i\in 2$, $\N_i={\sf{stack}}(\N_i|\d, \Sigma_{\R^b})$,\\
(4) both $\N_0$ and $\N_1$ have $(\d, \d+1)$-strategies as $\Sigma_{\R^b}$-mice\footnote{These strategies act on iterations below $\d$.}.\\\\
 Let then $\M$ be a common iterate of $\N_0$ and $\N_1$ (via these strategies). Let $j_0:\N_0\rightarrow \M$ and $j_1:\N_1\rightarrow \M$. We then have some $\k$ such that 
\begin{itemize}
\item $j_0(\pi_{E_0}(\k))=j_1(\pi_{E_1}(\k))$,
\item for $i\in 2$, $\cp(j_i)>\d^{\R^b}$, and 
\item $\rge(j_0\circ \pi_{E_0})\cap \rge(j_1\circ \pi_{E_1})$ contains a $\k$-club $C$\footnote{See \rthm{existence of thick sets}.}.
\end{itemize}
Let $D_0=(j_0\circ \pi_{E_0})^{-1}[C]$ and $D_1=(j_1\circ \pi_{E_1})^{-1}[C]$ and set $D=D_0\cap D_1$. It follows from universality of $\P^+$, $\N_0$ and $\N_1$ that\\\\
(5) $\P^b\subseteq Hull^{\P^+}(D\cup \d^{\P^b})$, and for $i\in 2$, $\R^b\subseteq Hull^{\N_i}(\pi_{E_i}[D])$.\\\\
Suppose now that $A\in \P^b\cap \powerset(\d^{\P^b})$ and fix $s\in D^{<\omega}$, $t\in [\d^{\P^b}]^{<\omega}$ and a term $\phi$ such that $A=\phi^{\P^+}[s, t]$. We then have that $\pi_{E_0}(A)=\phi^{\N_0}(\pi_{E_0}(s), \pi_{E_0}(t))$. It follows that $j_0(\pi_{E_0}(A))\in \rge(j_1\circ \pi_{E_1})$. Let $B\in \P^b\cap \powerset(\d^{\P^b})$ be such that $j_0(\pi_{E_0}(A))=j_1(\pi_{E_1}(B))$.  Because $\cp(j_i)>\d^{\R^b}$, we have that $\pi_{E_0}(A)=\pi_{E_1}(B)$. But (1) now implies that $A=B$. Therefore, $\pi_{E_0}(A)=\pi_{E_1}(A)$.
\end{proof}

It follows from \rprop{positional} that iterates of $(\P, \Sigma)$ can be successfully compared with one another. To prove it we simply compare $(\Q, \Sigma_\Q)$ with $(\R, \Sigma_\R)$ by using least-extender-disagreement comparison.

\begin{corollary}\label{self comparison} Suppose $(\P, \Sigma)$ is an allowable pair, $\Gamma$ is a  projectively closed pointclass and $\Sigma$ has strong branch condensation and is strongly $\Gamma$-fullness preserving. Suppose $(\T, \Q)\in I^{ope}(\P, \Sigma)$ and $(\U, \R)\in I^{ope}(\P, \Sigma)$. Then there is $(\T_1, \Q_1)\in I^{ope}(\Q, \Sigma_\Q)$ and $(\U_1, \R_1)\in I^{ope}(\R, \Sigma_\R)$ such that $\T_1$ and $\U_1$ are normal stacks and the following holds: 
\begin{enumerate}
\item Suppose $(\P, \Sigma)$ is a hod pair or a simple hod pair. Then one of the following holds:
\begin{enumerate}
\item $\Q_1\insegeq_{hod}\R_1$, $\pi^{\T_1}$ exists and $(\Sigma_{\R_1})_{\Q_1}=\Sigma_{\Q_1}$.
\item $\R_1\insegeq_{hod}\Q_1$, $\pi^{\U_1}$ exists and $(\Sigma_{\Q_1})_{\R_1}=\Sigma_{\R_1}$.
\end{enumerate}
Moreover, if in addition $(\T, \Q)\in I(\P, \Sigma)$ and $(\U, \R)\in I(\P, \Sigma)$, then $\Q_1=\R_1$ and both $\pi^{\T_1}$ and $\pi^{\U_1}$ are defined. 
\item Suppose $(\P, \Sigma)$ is an sts hod pair or a simple sts hod pair. Then one of the following holds:
\begin{enumerate}
\item $\Q_1\insegeq_{hod}\R_1$ and $(\Sigma_{\R_1})_{\Q_1}=\Sigma_{\Q_1}$. 
\item $\R_1\insegeq_{hod}\Q_1$ and $(\Sigma_{\Q_1})_{\R_1}=\Sigma_{\R_1}$.
\end{enumerate}
Moreover, if in addition $(\T, \Q)\in I(\P, \Sigma)$ and $(\U, \R)\in I(\P, \Sigma)$, then $\Q_1=\R_1$. Consequently, if $\Sigma$ is commuting then $\pi^{\T^\frown \T_1, b}=\pi^{\U^\frown \U_1, b}$.
\end{enumerate}
\end{corollary}
In clause 2 of \rcor{self comparison}, the conclusion $\Q_1=\R_1$ is a consequence of $\Gamma$-fullness preservation and our minimality assumption. If, for example, $\Q_1\inseg_{hod}\R_1$ then there is a $\Sigma_{\Q_1}$-sts $\W\inseg \R_1$ such that 
\begin{itemize}
\item $\W\models ``\d^{\Q_1}$ is a Woodin cardinal" ,
\item $\mathcal{J}_{\omega}[\W]\models ``\d^{\Q_1}$ is not a Woodin cardinal", and
\item $\W$ has a strategy in $\Gamma$.
\end{itemize}
It then follows that $\W\insegeq \Q_1$, contradiction.

The following is a corollary of \rcor{comparison holds 1} and \rthm{commuting for backgrounded st strategies}. 

\begin{proposition}\label{commuting} Suppose $(\P, \Sigma)$ is a hod pair or an sts hod pair, $\Gamma$ is a projectively closed pointclass and $\Sigma$ has strong branch condensation and is strongly $\Gamma$-fullness preserving. Then for some $(\T, \Q)\in I(\P, \Sigma)$, $\Sigma_{\Q, \T}$ is commuting.  
\end{proposition}

The next lemma will be used in the proof of \rthm{main theorem on gen int}.

\begin{lemma}\label{lemma used in gen int} Suppose 
\begin{itemize}
\item $(\P, \Sigma)$ is an allowable pair and $\P$ is non-meek, 
\item $\Gamma$ is a  projectively closed pointclass, 
\item $\Sigma$ has strong branch condensation and is strongly $\Gamma$-fullness preserving,
\item  if $(\P, \Sigma)$ is an sts hod pair or a simple sts hod pair then $\Sigma$ is commuting,
\item $(\T, \Q)\in I^{ope}(\P, \Sigma)$, $(\U, \R)\in I^{ope}(\P, \Sigma)$ and $(\W, \S)\in I(\R, \Sigma_{\R})$ are such that $\W$ is based on $\R^b$ and $\S\insegeq_{hod}\Q$.
\end{itemize}
 Then $\S=cHull^{\Q^b}(\pi^{\T, b}[\P^b]\cup \d^{\S})$.
\end{lemma}
\begin{proof} The proof is an easy corollary of commuting. Let $\W^+=\uparrow(\W, \R)$ and let $\S^+$ be the last model of $\W^+$. Notice that $(\S^+)^b=\S$. Let $\M$ be a common iterate of $(\Q, \Sigma_\Q)$ and $(\S^+, \Sigma_{\S^+})$ via respectively $\X$ and $\Y$, which we can find because of \rcor{self comparison}. We have that \\\\
(1) $\S=cHull^{\M^b}(\pi^{\Y, b}[\S])$ and $\Q^b=cHull^{\M^b}(\pi^{\X, b}[\Q^b])$.\\\\
It follows from \rprop{positional1} and commuting that\\\\
(2) $\pi^{\T^\frown \X, b}=\pi^{\U^\frown (\W^+)^\frown \Y, b}$ and $\pi^{\X, b}\rest \d^\S=\pi^{\Y, b}\rest \d^\S$\footnote{This equality holds as we have $(\Sigma_\Q)_{\S}=(\Sigma_{\S^+})_{\S}$.}.\\\\
It follows from (1), (2) and the fact that $\Sigma$ is strongly $\Gamma$-fullness preserving that\\\\
(3) $\S=cHull^{\M^b}(\pi^{\T^\frown \X, b}[\P^b]\cup \pi^{\X, b}[\d^\S])$.\\\\
Therefore, $\S=cHull^{\Q^b}(\pi^{\T, b}[\P^b] \cup \d^{\S})$.
\end{proof}

\section{Solidity and condensation}

The main contributions of this section are \rthm{solidity} and \rthm{condensation} that can be used to show that fully backgrounded hod pair constructions are successful, which amounts to showing that  clause 4 of \rdef{gamma-hod pair construction*} never occurs. We start with the following version of \rlem{general theorem on dodd-jensen} for phalanxes that is used in the proof of solidity and universality. We omit the actual proofs of  \rthm{solidity} and \rthm{condensation}  as, in the light of \rlem{dodd-jensen for certified phalanxes}, the proofs of solidity and universality are trivial generalizations of the usual proofs of these facts (see \cite[Chapter 5]{OIMT}).

\begin{remark} This section is devoted to showing that hod pair constructions of a background $(M, \d,\vec{G})$ converge. We thus think of the hod pairs that appear in the statement of lemmas and propositions of this section as hod pairs constructed by hod pair constructions, and since we would like to show that hod pair constructions are successful, which amounts to showing that clause 4 of \rdef{gamma-hod pair construction*} never occurs, the pairs we consider here are models appearing in the intermediate stages of hod pair constructions. This, in particular, means that first of all, we must deal with hod pairs (as opposed to sts hod pairs) and also the hod premouse is non-meek and not of lsa type, which is a consequence of our minimality assumption. The reason that it is enough to consider non-meek hod premice is that clause 4 of \rdef{gamma-hod pair construction*} for meek or gentle type hod premice is not new, and the usual proofs of solidity and condensation can be used. $\myqedhere$
\end{remark}

\begin{definition}[Certified phalanxes]\label{certified phalanx}\index{certified phalanx} Suppose $(\P, \Sigma)$ is a hod pair such that $\P$ is non-meek and $\R$ is a hod premouse. Let $\pi,\zeta$ be such that $\pi:\R\rightarrow \P$ is a $\Sigma_1$-embedding, and $\zeta\leq \cp(\pi)$. We say $(\P, \R, \zeta)$ is a $(\pi, \P, \Sigma)$-certified phalanx if $\zeta>o(\P^b)$. We also say $(\P,\R,\zeta)$ is a $(\P,\Sigma)$-certified phalanx witnessed by $\pi$. $\myqedhere$
\end{definition}
Continuing with the set up of \rdef{certified phalanx}, we let $\pi^+: (\P, \R, \zeta)\rightarrow (\P, \P, \zeta)$ be given by $(id, \pi)$, and also, we let $\Sigma^{\pi^+}$ be the $\pi^+$-pullback of $\Sigma$.

\begin{lemma}[No strategy disagreement]\label{no strategy disagreement} Suppose $(\P, \Sigma)$ is a hod pair such that  $\P$ is non-meek, $\Sigma$ has strong branch condensation and $\Sigma$ is strongly $\Gamma$-fullness preserving for some pointclass $\Gamma$ that is projectively closed. Suppose $(\P, \R, \zeta)$ is a $(\P, \Sigma)$ certified phalanx as witnessed by $\pi:\R\rightarrow \P$. Let $\Lambda=\Sigma^{\pi^+}$. Then no strategy disagreement appears in the comparison of $\P$ and $(\P, \R, \zeta)$ where $\Sigma$ is used on the $\P$ side and $\Lambda$ is used on the $(\P, \R, \zeta)$ side.
\end{lemma}
\begin{proof}
Towards a contradiction suppose not. It follows from the proof of \rlem{disagreement implies low level disagreement} that we can find a minimal low level disagreement $((\T, \Q), (\U, \S), \W)$ between $\Sigma$ and $\Lambda$. Let then $E=E^{\U}_\W$, be the $\W$-un-dropping extender of $\U$. We have that $\W\inseg_{hod} Ult(\P, E)$. Let now $\X=\pi^+\U$, $\P_1$ be the last model of $\X$, $\sigma: \S\rightarrow \P_1$ be the copy map and $F$ be the $\sigma(\W)$-un-dropping extender of $\X$. Let $\sigma': Ult(\P, E)\rightarrow Ult(\P, F)$ be given by $\sigma'([a, f]_E)=[\sigma(a), f]_E$. 

We now have that $(\X, \P_1)$ and $((\pi_F\rest \P^b, Ult(\P, F)^b), (\pi_E\rest \P^b, Ult(\P, E)^b), \sigma')$ support a $(\sigma(\W), \W, \sigma\rest \W)$-b-condensation diagram on $\P$. Because $\sigma\rest \W$ pullback of $\Sigma_{\sigma(\W), \X}$ is $\Lambda_{\W, \U}$, it follows from strong branch condensation that $\Sigma_{\W, \T}=\Lambda_{\W, \U}$.
\end{proof}

\begin{definition}[Certified pairs]\label{certified pairs}\index{certified pair} Suppose $(\P, \Sigma)$ is a hod pair and $\R$ is a hod premouse such that both $\P$ and $\R$ are of limit type. Suppose that there is $\pi$ such that $\pi:\P^b\rightarrow \R^b$ is elementary. We say the pair $(\pi, \R)$ is $(\P, \Sigma)$-certified by $(\sigma, \T, \Q, \Q')$ if
 \begin{enumerate}
 \item  $(\T, \Q)\in I(\P, \Sigma)$, $\Q'\insegeq_{hod}\Q$ and $\sigma:\R\rightarrow \Q'$ is $\Sigma_1$-elemnetary,
 \item $(\Q')^b=Hull^{\Q}(\pi^{\T}[\P^b]\cup \d^{(\Q')^b})$, and
 \item  letting $k:\P^b\rightarrow (\Q')^b$ be the collapse of $\pi^{\T}[\P^b]$, $k=(\sigma\rest \R^b)\circ \pi$. 
 \end{enumerate}
We say $(\R, \Lambda)$ is a $(\P, \Sigma)$-certified hod pair if for every $(\U, \S)\in I(\R, \Lambda)$, there is $\pi$ and $(\sigma,\T,\Q,\Q')$ such that $(\pi, \S)$ is $(\P,\Sigma)$-certified by $(\sigma, \T, \Q,\Q')$ and 
\begin{center}$\Lambda_{\S^b, \U}=(\sigma$-pullback of $\Sigma_{\Q^b, \T})$. \end{center}
$\myqedhere$
\end{definition}

%

\begin{lemma}\label{general theorem on dodd-jensen} Suppose $(\P, \Sigma)$ is a hod pair such that $\P$ is non-lsa type non-meek hod premouse, $\Gamma$ is a projectively closed pointclass and $\Sigma$ has strong branch condensation and is strongly $\Gamma$-fullness preserving. Suppose $(\T, \R)\in I^b(\P, \Sigma)$\footnote{Thus, $\pi^{\T, b}$ exists, see \rdef{almost non-dropping stacks}.} is such that for some $\Lambda$, $(\R, \Lambda)$ is $(\P, \Sigma)$-certified and there is a $\Sigma_0$-elementary embedding $\pi:\P\rightarrow \R$ such that $(\pi\rest \P^b,\R)$ is $(\P,\Sigma)$-certified by $(\sigma, \U, \Q,\Q')$. Then $\pi^{\T}$ exists and $\pi^{\T}\leq \pi$.
\end{lemma}
\begin{proof} To implement the usual proof of the Dodd-Jensen property (see \cite[Chapter 4.2]{OIMT}), we need to know that\\\\
(a) $\Sigma$ is the $\pi$-pullback of $\Lambda$.\\\\
 Because $(\R, \Lambda)$ is  $(\P, \Sigma)$-certified, (a) easily follows from strong branch condensation of $\Sigma$. 
\end{proof}

\begin{lemma}[Dodd-Jensen for certified phalanxes]\label{dodd-jensen for certified phalanxes} 
Suppose $\Gamma$ is a projectively closed pointclass and $(\P, \Sigma)$ is a hod pair such that $\Sigma$ has strong branch condensation and is strongly $\Gamma$-fullness preserving. Suppose that $(\P, \R, \zeta)$ is a $(\P, \Sigma)$-certified phalanx as witnessed by $\pi:\R\rightarrow \P$. Suppose that
\begin{itemize}
\item $\T$ is a stack on $(\P, \R, \zeta)$ according to $\Sigma^{\pi^+}$ with last model $\Q$,
\item $\U$ is a stack on $\P$ according to $\Sigma$ with last model $\S$, and
\item  the last branch of $\T$ is on $\P$ and either
\begin{enumerate}
\item $\Q\insegeq_{hod}\S$ and $\pi^\T$ exists or 
\item $\S\insegeq_{hod}\Q$ and $\pi^\U$ exists. 
\end{enumerate}
\end{itemize}
Then $\Q=\S$ and $\pi^{\T}=\pi^{\U}$. 
\end{lemma}
\begin{proof} Let $\T^*=\pi^+\T$. Let $\Q^*$ be the last model of $\T^*$ and let $\sigma:\Q\rightarrow \Q^*$ come from the copying construction. Suppose first that $\Q\insegeq_{hod}\S$ and $\pi^\T$ exists.  Notice next that $\pi^{\T}$-pullback of $(\Sigma_{\Q^*})^{\sigma}$ is $\Sigma$. Hence, applying the ordinary proof of the Dodd-Jensen property we get that $\S=\Q$, $\pi^\U$ exists and $\pi^\U\leq \pi^\T$. 

Suppose now $\S\insegeq_{hod}\Q$ and $\pi^\U$ exists. Notice that $\pi^\U$ induces an embedding $\pi^*:(\P, \R, \zeta)\rightarrow (\S, \S, \pi^\U(\zeta))$ such that $\pi^*\rest \P=\pi^\U$ and $\pi^*\rest \R=\pi^\U\circ \pi$. Notice that\\\\
(1) $\Sigma^{\pi^+}=(\pi^*$-pullback of $\Sigma_{\S})$.\\\\
  Applying \rlem{general theorem on dodd-jensen} to the embedding $(\sigma\rest \S)\circ \pi^\U$ and $(\T^*, \Q^*)$, we get \\\\
  (2) $\sigma(\S)=\Q^*$, $\pi^{\T^*}$ exists and $\pi^{\T^*}\leq \sigma\circ \pi^\U$.\\\\
  (2) now implies that $\S=\Q$.  Since $\pi^{\T^*}=\sigma\circ \pi^\T$, we have that $\pi^\T$ exists and $\pi^\T\leq \pi^\U$. Putting the two arguments together we see that $\pi^\U=\pi^\T$.
\end{proof}

It is clear that it follows from \rlem{dodd-jensen for certified phalanxes} and from \rlem{no strategy disagreement} that the usual proofs of condensation, universality and solidity go through for hod mice. We state the results without proofs (see \cite[Chapter 5]{OIMT} for the usual proofs of these results.)

\begin{theorem}[Solidity and universality]\label{solidity} Suppose $\Gamma$ is a projectively closed pointclass, $k<\omega$ and $(\P, \Sigma)$ is a hod pair such that 
\begin{enumerate}
\item $\P$ is $k$-sound non-meek hod premouse,
\item $\P$ is not of lsa type and $\rho(\P)>o(\P^b)$, and 
\item $\Sigma$ is strongly $\Gamma$-fullness preserving and has  strong branch condensation.
\end{enumerate} Let $r$ be the $k+1$st standard parameter of $(\P, u_k(\P))$; then $r$ is $k+1$-solid and $k+1$-universal over $(\P, u_k(\P))$. 
\end{theorem}

\begin{theorem}[Condensation]\label{condensation}  Suppose $\Gamma$ is a projectively closed pointclass and $(\P, \Sigma)$ is a hod pair such that 
\begin{enumerate}
\item $\P$ is non-meek hod premouse,
\item $\P$ is not of lsa type and $\rho(\P)>o(\P^b)$, and 
\item $\Sigma$ is strongly $\Gamma$-fullness preserving and has  strong branch condensation.
\end{enumerate} 
 Suppose $(\P, \R, \zeta)$ is a $(\P, \Sigma)$ certified phalanx as witnessed by $\pi:\R\rightarrow \P$ such that $\zeta=\cp(\pi)=\rho_\omega^\R$. Then either
\begin{enumerate}
\item $\R\insegeq_{hod}\P$ or
\item there is an extender $E$ on the sequence of $\P$ such that $lh(E)=\rho_\omega^\R$ and $\R\insegeq_{hod}Ult(\P, E)$.
\end{enumerate}
\end{theorem}

 \section{Backgrounded constructions relative to st-strategies}\label{sec: background constructions st-strategy}

In this section, we show that if $(\P, \Sigma)$ is an sts pair constructed by a hod pair construction then if the fully backgrounded construction relative to $\Sigma$ breaks down then it does so because it reaches an sts mouse that destroys the Woodiness of $\d^\P$. The reader may find it helpful to review \rdef{authentic lsp}, \rdef{authentic iterations}, \rdef{authentic stacks of length 2}, clause 5 of \rdef{weak psi alpha indexing scheme a} and \rdef{sts0}.

\begin{theorem}\label{sts fb constructions converge} Assume $\sf{AD}^+$. Suppose
\begin{itemize}
\item  for some $\a$ such that $\theta_\a<\Theta$, $\Gamma=\{ A\subseteq \bR: w(A)<\theta_\a\}$, 
\item  ${\sf{C}}=(\mathbb{M}, (P, \Psi), \Gamma^*, A)$ Suslin, co-Suslin captures $\Gamma$ and 
\item $\mathbb{M}=(M, \d,  \vec{G}, \Sigma^*)$, 
\item ${\sf{hpc}}=(\M_\gg , \N_\gg, Y_\gg, \Phi_\gg, F^+_\gg, F_\gg, b_\gg : \gg\leq \delta)$ is the output of the $\Gamma-{\sf{hpc}}$ of $M$,
\item $\xi<\d$ is such that $(\M_\xi, \Phi_\xi)$ is an sts pair and $M\models (\M_\xi, \Phi_\xi)\in {\sf{Hp}}^{\Gamma}$.
\end{itemize}
Set $(\M_\xi, \Phi_\xi^+)=(\P, \Sigma)$ and let $\zeta\geq \sup\{\lh(F^+_\nu): \nu<\xi\}$. Let  
\begin{center}
$({\sf{Le}}((\P, \Sigma^{stc}), \mathcal{J}_\omega[\P])_{>\zeta})^M=(\P_\gg , \P'_\gg, X^+_\gg, X_\gg, b_\gg: \gg<\eta)$.
\end{center}
Suppose there is $\xi_0$ such that the anomaly stated in clause 3.b of \rdef{fully backgrounded sts construction} occurs at $\xi_0$.
Then $\P_{\xi_0}\models ``\d^\P$ is not a Woodin cardinal". 

Consequently, if ${\sf{Lp}}^{\Gamma, \Sigma^{stc}}(\P)\models ``\d^{\P}$ is a Woodin cardinal" then the anomaly stated in clause 3.b of \rdef{fully backgrounded sts construction} does not occur.
\end{theorem}
\begin{proof} Suppose that $\P_{\xi_0}\models ``\d^\P$ is a Woodin cardinal". Set $\S=(\P_{\xi_0})_{\sf{ex}}$\footnote{See \rdef{lsa type}.}. We assume that $\xi$ is the least such that our claim fails for $(\M_\xi, \Phi_\xi^+)$. Suppose $\S=(\mathcal{J}^{\vec{E}, f}_{\omega\a_0}, \in, \vec{E}, f)$, $\a_0=\b_0+\gg_0$ and $t=(\P, \T)\in \univ{\S|\omega\b} \cap \dom(\Sigma)$ is such that setting $w=(\mathcal{J}_{\omega}(t), t, \in)$, $w$ is $(f, {\sf{sts}})$-minimal as witnessed by $\b_0$\footnote{See \rdef{important notation}. In particular, this means that we have to index the branch of $t$ at $\omega\a_0$.}, $\T$ is $\P|\omega\b_0$-terminal\footnote{See \rdef{terminal tree}.} and $\gg_0=\lh(t)$. 
Set $b=\Sigma(\T)$, $p=\uparrow(\T, \S)$, $\P_1=\m^+(\T)$ and
\begin{center}
$\S_1=\begin{cases}
\M^p_b &: \Q(b, p)\ \text{doesn't exist}\\
\Q(b, p) &:\ \text{otherwise}.
\end{cases}$\end{center} 
Notice that if $\mathcal{J}_{\omega}[\S]\models ``\d^\P$ is not a Woodin cardinal" then $\Q(b, p)$ is defined. 

We then have that
     \begin{center}
     $\P'_{\xi_0+1}=(\mathcal{J}^{\vec{E}, f^+}_{\omega\b_0+\omega\gg_0}, \in, \vec{E}, f, \tilde{b})$
     \end{center}
     where $\tilde{b}\subseteq \omega\b_0+\omega\gg_0$ is defined by $\omega\b_0+\omega\nu\in \tilde{b}\iff \nu \in b$. Since we are assuming that the anomaly occurs, we have that there is $e\in \S|\omega\b_0$ such that $\S|\omega\b_0\models {\sf{sts}_0}(t, e)$\footnote{See \rdef{sts0}. This means that $e$ is the branch of $t$ we must choose.} and $e\not =b$. Let $\Phi$ be the strategy of $\S$ induced by $\Sigma^*$\footnote{Notice that $\P$ is constructed in $M|\zeta$ and $\S$ above $\d^\P$ is constructed using extenders with critical points $>\zeta$. It follows that $\Sigma^*$ indeed induces a strategy $\Phi$ for $\S$ via the ordinary resurrection procedure of \cite[Chapter 12]{FSIT}. See also \rsubsec{the induced strategy for the undropping game}.}. Notice that because $\d^\P$ is a cutpoint in $\S$, $\Phi$ extends $\Sigma$. 
     
\begin{sublemma}\label{iterates are correct sts premice} Whenever $\S'$ is a $\Phi$-iterate of $\S$ via a stack that is above $\d^\P$, $\S'$ is a $\Sigma^{stc}$-sts premouse over $\P$. Thus, any two $\Phi$-iterates of $\S$ can be compared to each other.

Similarly if $\U$ is a generalized stack on $\S$ according to $\Phi$ such that $\U$ is based on $\P$ and has a last normal component\footnote{See \rnot{notation for iteration trees}.} $\U'$, $d=\Sigma(\U)$, $\S'=\M^\U_d$ and $\m^+(\U')\in Y^{\S'}$ is $\#$-lsa like then whenever $\S''$ is a $\Phi_{\S', \U^\frown \{d\}}$-iterate of $\S'$ via a stack that is above $\d(\U')$, $\S''$ is a $\Sigma^{stc}_{\m^+(\U'), \U}$-sts premouse over $\m^+(\U')$.
\end{sublemma}
\begin{proof} This follows from hull condensation of $\Sigma$. We do the proof for stacks, and the more general proof is only notationally more complicated.  Suppose $\U$ is a stack on $\S$ according to $\Phi$ with last model $\S'$. Let $\U^*$ be the resurrection of $\U$ onto $M$ and let $M^*$ be the last model of $\U^*$\footnote{See \cite[Chapter 12]{FSIT} and also \rsubsec{the induced strategy for the undropping game}.}. Set
\begin{center}
$\pi^{\U^*}(({\sf{Le}}((\P, \Sigma^{stc}), \mathcal{J}_\omega[\P])_{>\zeta})^M)=(\R_\tau , \R'_\tau, Z^+_\tau, Z_\tau, c_\tau: \tau<\eta')$.
\end{center}

We then have some $\iota\leq \pi^{\U^*}(\xi_0)$ and $\sigma:\S'\rightarrow \R_\iota$. Let $s=(\P, \X_0, \P_1, \X_1)$ be an indexable stack such that $s\in \dom(\Sigma^{\S'})\cap \dom(\Sigma^{stc})$. We want to see that
\begin{center}
$\Sigma^{stc}(s)=\Sigma^{\S'}(s)$.
\end{center} 
Let $d=\Sigma^{\S'}(s)$. Notice next that $s^\frown\{ d\}$ is hull of $\sigma(s)^\frown \{\sigma(d)\}$ and that $\sigma(s)^\frown \{\sigma(d)\}$ is according to $\Sigma$\footnote{Set $s'=\sigma(s)^\frown \{\sigma(d)\}$. We have that $\P|\d^\P$ is constructed inside $M|\zeta$ while the construction producing $\S$ and $\S''$ uses extender with critical points $>\zeta$. It follows that $\U^*$ is above $\zeta+1$ while $\Sigma$ is determined by the pair $(M|\zeta, \Sigma^*_{M|\zeta})=(M^*|\zeta, (\Sigma^*_{M^*})_{M|\zeta})$. Thus, $(\sigma^*(\Phi_\xi))^+=\Sigma$. Notice also that it follows from the elementarity of $\pi^{\U^*}$ that $\R_\iota$ is indeed a $\Sigma^{stc}$-sts as the first time this breaks down is at $\pi^{\U^*}(\xi_0)$. Thus, any indexable stack that has been indexed in $\R_\iota$ is according to $\Sigma^{stc}$.}. But $\Sigma$ has hull condensation (see \cite[Lemma 2.9]{ATHM}), implying that $\Sigma(s)=d$.

The second part of the claim is very similar. This time, letting $\U^*$ be the resurrection of $\U$, we have $\sigma: \S'\rightarrow \Q_\iota$ where 
 \begin{center}
 $\pi^{\U^*}({\sf{hpc}})=(\Q_\iota, \Q_\iota', Z_\iota, \Omega_\iota, H^+_\iota, H_\iota, e_\iota : \iota\leq \pi^{\U^*}(\delta))$
 \end{center}
  is the output of the $\Gamma-{\sf{hpc}}$ of the last model of $\U^*$ and $\iota\leq \pi^{\U^*}(\xi)$. But now we repeat the same argument as before noting that $\Sigma^{stc}_{\m^+(\U'), \U}$ is the short-tree component of the $\sigma$-pullback of $\Omega^+_\iota$.
\end{proof}

We now have that $\Q(e, \T)$ exists and it is an sts premouse over $\m^+(\T)$\footnote{See \rdef{mtsharp}.}. Set $\Q=\Q(e, \T)$. Notice that because $w$ is $(f, {\sf{sts}})$-minimal, we must have that for some $\tau_0<\a_0$, $\rho(\S||\tau_0)\leq \d(\T)$ and $\S||\tau_0\models {\sf{sts}_0}(t, e)$. Let $\tau_1$ witnesses that $\S|\tau_0\models {\sf{sts}_0}(t, e)$. Set $\S'=\S||\tau_0$. 

Using \rlem{n*x}, we can find a self-capturing background $(M_0, \d_0, \vec{G}_0, \Sigma_0^*)$ which Suslin, co-Suslin captures ${\sf{Code}}(\Sigma^*)$ and 
\begin{center}
$({\sf{HC}}^{M_0}, {\sf{Code}}(\Sigma^*), \in)\prec^{\bR^{M_0}} ({\sf{HC}}, {\sf{Code}}(\Sigma^*), \in)$
\end{center}
where $\prec^Z$ means elementarity with respect to parameters in $Z$. 

Let $\l$ be the supremum of the first $\omega$-Woodin cardinals of $\S'|\omega\tau_1$.  Let $h\subseteq Coll(\omega, \bR^{M_0})$ be generic and let $\S''$ be an $\mathbb{R}^{M_0}$-genericity iterate of $\S'$ via $\Phi_{\S'}$. Thus, we have a stack $\U\in M_0[h]$ on $\S'$ such that the following holds in $M_0[h]$:
\begin{enumerate}
\item $\lh(\U)=\omega_1^{M_0}+1$,
\item $\S''$ is the last model of $\U$,
\item for every $\a<\omega_1^{M_0}$, $\U_{\leq \a}\in M_0$,
\item $D^\U=\emptyset$,
\item $\pi^\U(\l)=\omega_1^{M_0}$,
\item for some $\S''$-generic $k\subseteq Coll(\omega, <\omega_1^{M_0})$, $\bR^{M_0}$ is the set of symmetric reals of $\S''[k]$.
\end{enumerate}
We let $N=D(\S'', \omega_1^{M_0}, k)$\footnote{This is the derived model of $\S''$ as computed by $k$. See \rsec{sec short tree strategy mice}. We need to work inside $M_0$ to guarantee that $\S''\in V$.}. We now have a strategy $\Lambda\in N$ for $\Q$ and some $\nu$ with the property that
\begin{enumerate}
\item $N\models ``\Lambda$ is an $\omega_1$-iteration strategy" and 
\item whenever $\R\in N$ is a $\Lambda$-iterate of $\Q$ above $\delta(\T)$ and $s\in \R$ is an indexable stack on $\P_1=\m^+(\T)$ according to $\Sigma^\R$, 
\begin{center}
$\S''[k]\models ``s$ is $(\P, \Sigma^{\S''|\nu})$-authenticated"\footnote{See \rdef{weak psi alpha indexing scheme a}.}.
\end{center}
\end{enumerate}
 Notice that $N\in M_0$\footnote{This follows from the fact that $\rho(\S')\leq \d(\T)$ and from \rsublem{iterates are correct sts premice}. See \cite[Proposition 3.0.1]{DMATM}.} implying that $\Lambda\in M_0$. Moreover, because ${\sf{Code}}(\Lambda)$ is projective in ${\sf{Code}}(\Sigma^*)$, we have that ${\sf{Code}}(\Lambda)$ is $\d$-universally  Baire in $M_0$. Thus, $\Lambda$ also acts on length $\omega_1^{M_0}$ iterations. 

\begin{sublemma}\label{an incorrect iterate} There is a $\Lambda$-iterate $\Q^*$ of $\Q$ and a $\Phi_{\S_1, p}$-iterate $\S^*$ of $\S_1$ such that $\S^*=\Q^*|\ord(\S^*)$ and $\Q^*||\ord(\S^*)$ is not a $\Sigma_{\m^+(\T), \T}^{stc}$-sts premouse over $\m^+(\T)$ (implying that $\Q^*|\ord(\S^*)\not = \Q^*||\ord(\S^*)$).
\end{sublemma}
\begin{proof} 
Our goal now is to compare $\S_1$ with $\Q$. We use  $\Phi_{\S_1, p}$ for $\S_1$ and $\Lambda$ for $\Q$. Assume for a moment that the comparison is successful. If $\S_1=\Q(b, p)$ then we in fact have that $\Q(b, p)=\Q$ and since $\Q=\Q(e, \T)$, we get that $b=e$, which is a contradiction. Hence, $\S_1=\M^p_b$ and $\Q(b, p)$ doesn't exits (and hence $\pi^{p}_b$ is defined). In this case, we must have that $\S_1$ loses the comparison and if $(\S^*, \Q^*)$ are the last models of the comparison then $\S^*\inseg \Q^*$\footnote{Equality is not possible because $\S^*$ is not a $\Q$-structure for $\T$.}. Because the $\S_1$-side loses we have that the iteration embedding $j_0:\S_1\rightarrow \S^*$ is defined. Let then $j=j_0\circ \pi^p_b:\S\rightarrow \S^*$.

We now argue that $\Q^*$ is as desired. To show this we only need to show that $\Q^*||\ord(\S^*)$ is not a $\Sigma_{\m^+(\T), \T}^{stc}$-sts premouse over $\m^+(\T)$. Notice that our $\sf{sts}$-indexing scheme is so that if $\W$ is an sts premouse properly extending $\S^*$ and $\a^*=\ord(\S^*)$ then $(j(\T), j(e))$ has to be indexed in $\W$ at $\a^*$. Thus, $(j(\T), j(e))$ must be indexed at $\a^*$ in $\Q^*$. Assume then that $(j(\T), j(e))$ is according to $\Sigma^{stc}_{\m^+(\T), \T}$. Because $\T^\frown\{e\}$ is a hull of $\T^\frown\{b\}^\frown j(\T)^\frown\{j(e)\}$ as witnessed by $j$, it follows from hull condensation\footnote{See \cite[Lemma 2.9]{ATHM}.} of $\Sigma$ that $e=\Sigma(\T)$, contradiction.

Thus, we must have that the comparison between $\S_1$ and $\Q$ does not terminate. But then \rsublem{iterates are correct sts premice} implies that this can only happen because the comparison of $\S_1$ and $\Q$ produces an iterate $\Q^*$ of $\Q$ which is not a $\Sigma^{stc}_{\P_1, \T}$-sts over $\P_1$. Let then $\S^*$ be the corresponding iterate of $\S_1$. We thus have some $\iota<\ord(\S^*)$ such that \\\\
(1) $\iota\not \in \dom(\vec{E}^{\S^*})\cap \dom(\vec{E}^{\Q^*})$ and $\S^*|\iota=\Q^*|\iota$ but $\S^*||\iota\not = \Q^*||\iota$. \\\\
Let $t_1=(\P_1, \T_1, \P_1', \T_1')\in \dom(\Sigma^{\Q^*})\cap \dom(\Sigma^{\S^*})$ be such that $\Sigma^{\S^*}(t_1)\not=\Sigma^{\Q^*}(t_1)$. Assume first that $\T_1'$ is defined. In this case we have that $\pi^{\T_1, b}$ is defined and $\T_1'$ is based on $(\P_1')^b$. We then have that $((\P_1')^b, \T_1')$ is  $(\P, \Sigma^{\S''|\nu})$-authenticated iteration\footnote{In fact, an authentication exists in $\S''[k]$. See \rdef{authentic lsp}.}. Suppose $\Y$ authenticates $((\P_1')^b, \T_1')$. Set $\R=(\P_1')^b$ and let $\W$ be the last model of $\Y$. Notice that $\T^\frown \T_1$ is according to $\Sigma^{stc}$ and moreover, $\R=\pi^{\T^\frown \T_1, b}(\P^b)$. 

It follows that there is $\pi:\P^b\rightarrow \R$ and $\sigma:\R\rightarrow \W^b$ such that $\pi^{\Y, b}=\sigma\circ \pi$ and $\sigma\T_1'$ is according to $(\Sigma^{\S''|\nu})_{\W^b, \Y}$\footnote{Recall that $\Sigma^\X$ is the strategy predicate of $\X$.}. It follows from \rsublem{iterates are correct sts premice} that $(\Sigma^{\S''|\nu})_{\W^b, \Y}\subseteq \Sigma_{\W^b, \Y}$. Applying \rthm{strong condensation for backgrounded strategies} to $(\Y, \W^b, \R, \R, \sigma)$ and $\T^\frown \T_1$, we get that $\T_1'$ is according to $\Sigma_{\R, \T^\frown \T_1}$. Hence, we must have that $\T_1'=\emptyset$.

We thus have that $t_1=(\P_1, \T_1)$. If $\T_1$ is ${\sf{uvs}}$\footnote{See \rdef{nus stacks}.} then by arguing as above we once again prove that $\T_1$ is according to $\Sigma_{\P_1, \T}$. Assume then that $\T_1$ is $\sf{nuvs}$. It follows that $\m^+(\T_1)$ is a $\#$-lsa type hod premouse. In this case, our sts scheme guarantees that there are branches $c_1\in \S^*$ and $c_2\in \Q^*$ such that $(\T_1, c_1)$ is indexed at $\iota$ in $\S^*$ and $(\T_1, c_2)$ is indexed at $\iota$ in $\Q^*$. But because $\S^*|\iota=\Q^*|\iota$, we must have that $c_1=c_2$ as what branch is indexed at $\iota$ in either of the models depends solely on $\S^*|\iota=\Q^*|\iota$ and not on any external factors. We thus have that $\S^*||\iota=\Q^*||\iota$ contradicting (1). This contradiction implies that in fact the comparison between $\S_1$ and $\Q$ is successful. 
\end{proof}

Let then $\Q^*$ and $\S^*$ be as in the sublemma above. Because $\Q^*$ wins the coiteration we have that the iteration embedding $j_0:\S_1\rightarrow \S^*$ exists. $j_0$ is according to $\Phi_{\S_1, p}$. As we have argued in the proof of the sublemma, we must also have that $\pi^p_b$ exists. Set then $j=j_0\circ \pi^p_b$. As in the proof of the sublemma, the pair $(j(\T), j(e))$ must be indexed in $\Q^*$ at $\ord(\S^*)$.

 It then follows that setting $\T_1=j(\T)$, $e_1=j(e)$ and $\Q_1=_{def}(j(\Q(e, \T)))_{\sf{ex}}=(\Q(e_1, \T_1))_{\sf{ex}}$\footnote{See \rdef{lsa type}.}, $\Q_1$ is $(\P, \Sigma^{\S''|\nu})$-authenticated. This means that we can find a normal stack $\Y$ on $\P$ with last model $\W$ and an ordinal $\iota$ such that for some normal stack $\X$ on $\Q_1$,\\\\
 (2) $\X$ is based on $\P_2=_{def}\m^+(j(\T))$,\\
 (3) $\W||\iota$ is the last model of $\X$ and $\pi^\X$ is defined\footnote{See \rdef{authentic lsp}. Clause 1 applies to our current situation.},\\
 (4) $\W||\iota$ is a $\Sigma^{stc}_{\W_0, \Y}$-sts over $\W_0$ where $\W_0=_{def}(\W||\iota)_\#$.\\\\\ 
 (4) follows from the fact that $\S''$ is a $\Sigma^{stc}$-sts premouse over $\P$ (see \rsublem{iterates are correct sts premice}). We claim that\\\\
 (5) $\X$ is according to $\Sigma^{stc}_{\P_2, \T^\frown \T_1}$. \\\\
 (5) follows from \rsublem{next lemma}. Assuming (5) we finish the argument. Let 
 \begin{itemize}
 \item $\X_0=\downarrow(\X, \P_2)$
 \item $p_1=\uparrow(\T_1, \S_1)$,
 \item $b_1=\Sigma(p^\frown\{b\}^\frown p_1)$,
 \item $\S_2=\M_{b_1}^{p^\frown\{b\}^\frown p_1}$,
 \item $\X'=\uparrow(\X_0, \S_2)$.
 \end{itemize}
 Arguing just like for $(b, \S_1)$ we have that $\Q(b_1, p^\frown\{b\}^\frown p_1)$ does not exist and $\pi^{p^\frown\{b\}^\frown p_1}_{b_1}$ is defined.
 It follows from (5) that $\X'$ is according to $\Phi_{\S_2, q}$ where $q=p^\frown\{b\}^\frown p_1^\frown\{b_1\}$. Let $\S_3$ be the last model of $\X'$. Because $\pi^{\X'}$ is defined and because $\mathcal{J}_{\omega}[\W||\iota]\models ``\d(\X_0)$ is not a Woodin cardinal", we have that there is a normal $\Phi_{\S_3, q^\frown \X'}$-iterate $\S_4$ of $\S_3$ and a normal $\Sigma_{\W||\iota, \Y}$-iterate $\W'$ of $\W||\iota$ such that $\S_4\inseg \W'$ and the iteration embedding $k:\S\rightarrow \S_4$ exists. Because $\W'$ is a $\Sigma$-iterate of $\P$, we have that $(k(\T), k(e))$, which according to our sts indexing scheme must be indexed in $\W'$, is according to $\Sigma$.  It then follows that $\T^\frown \{e\}$ is a hull of $\T^\frown \{b\}^\frown \T_1^\frown\{b_1\}^\frown \X_0^\frown k(\T)^\frown \{k(e)\}$ implying that $b=e$\footnote{Notice that $\W||\iota$-to-$\W'$ iteration is above $\d(\X_0)$.}. 
 
 This completes the proof of the theorem assuming (5). \rsublem{next lemma} implies (5). To simplify the matter, the symbols used in the statement of \rsublem{next lemma}, with the exception of $(\P, \Sigma)$, do not have the same meaning as the same symbols in the proof given above. 

\begin{sublemma}\label{next lemma} Suppose 
\begin{itemize}
\item $\X$ is a generalized stack on $\P$ according to $\Sigma^{stc}$,
\item $\X$ has a last normal component $\X'$ with the property that $\R=\m^+(\X')$ is a $\#$-like lsa type hod premouse,
\item $\T$ is a stack on $\P$ according to $\Sigma^{stc}$ that authenticates $\R$.
\end{itemize}
 Let $\S$ be the last model of $\T$ and let $\U$ be a normal stack on $\R$ witnessing that $\T$ authenticates $\R$\footnote{See \rdef{authentic lsp}.}. Then $\U$ is according to $\Sigma^{stc}_{\R, \X}$.  
\end{sublemma}
\begin{proof}
The proof uses ideas from \rlem{the canonical singularizing sequence for g-stacks}, \rthm{strong fullness preservation}, \rthm{strong condensation for backgrounded strategies}. We will use  \rlem{the canonical singularizing sequence for g-stacks} and \rthm{strong fullness preservation} to conclude that $\U$ picks the branches according to $\Sigma^{stc}_{\R, \X}$ in successor windows. Because the proofs are very similar to the proofs already given in the above mentioned theorems, we will sketch the arguments.

Set $\eta+1=\lh(\U)$. Suppose $\a\leq \eta$ is a limit ordinal and $\U_{<\a}$ is according to $\Sigma^{stc}_{\R, \X}$. We want to see that \\\\
(a) $[0, \a)_\U=\Sigma^{stc}_{\R, \X}(\U_{<\a})$. \\\\
Let $\S'\insegeq_{hod}\S$ be the longest such that $\S'\insegeq \m(\U_{<\a})$ and either $\S'$ is a layer of $\S$ or a limit of layers of $\S$. There are two essential cases.\\\\
\textbf{Case 1:} $\S'\inseg_{hod} \S^b$ is a complete layer\footnote{See \rdef{l p}.} of $\S^b$.\\\\
 Let $\S''$ be the least layer of $\S^b$ such that $\S'\inseg \S''$ and $\d^{\S''}$ is a Woodin cardinal of $\S^b$. If now $\d(\U_{<\a})<\d^{\S''}$ then (a) follows from fullness preservation of $\Sigma$\footnote{Notice that \rthm{strong condensation for backgrounded strategies} implies that $\Sigma_{\S', \T}=\Sigma_{\S', \X^\frown \U}$. This is because we can apply \rthm{strong condensation for backgrounded strategies} to $\Q=_{def}\S'$, $E=($the $(\d^{\P^b}, \d^\Q)$-extender derived from $\pi^{\T, b})$, $\R=_{def}\S'$ and $\sigma=_{def}id$.}. Suppose then that $\d(\U_{<\a})=\d^{\S''}$ then letting $w$ be the window of $\S$ such that $\d^w=\d^{\S''}$, we need to see that\\\\
(b) if $b=\Sigma(\U_{<\a})$ then $s(\T, w)\subseteq \rge(\pi^{\U_{<\a}}_b)$.\\\\
(b) is a consequence of the second clause of strong branch condensation\footnote{See \rdef{strong branch condensation}.}. Let $\b$ be the least such that $\S'\insegeq \M_\b^\U$. Set $\N=\M^\U_\b$. Notice that both $\b$ and $\alpha$ are on the main branch of $\U$ implying that both $\pi^\U_{\b, \eta}$ and $\pi^\U_{\a, \eta}$ are defined. Let $\pi=\pi^\U_{\b, \eta}\rest \N^b$ and $\sigma=\pi^\U_{\a, \eta}\rest (\M_\a^\U)^b$. Let $c=[0, \a)_\U$ and set $\Y=\downarrow(\U_{[\b, \a)}, \N^b)$. We can now apply clause 2 of strong branch condensation to $(\pi, \sigma, \U_{\leq \b}, \N, \T, \S, \Y, c)$.\\\\
\textbf{Case 2:} $\S'\inseg_{hod} \S^b$ is not a complete layer\footnote{See \rdef{l p}.} of $\S^b$.\\\\
Let $c=[0, \a)_\U$. In this case, we have that $\Q(c, \U_{<\a})$ exists and $\Q(c, \U_{<\a})\insegeq \S$. The dificult case is when $\Q(c, \U_{<\a})$ is an sts premouse over $\m^+(\U_{<\a})$, and so we assume it. We then have that $\Q(c, \U_{<\a})$ is a $\Sigma^{stc}_{\m^+(\U_{<\a}), \T}$-sts premouse over $\m^+(\U_{<\a})$. It then follows from \rprop{positional} that in fact $\Sigma_{\m^+(\U_{<\a}), \T}=\Sigma_{\m^+(\U_{<\a}), \X^\frown \U_{<\a}}$\footnote{Notice that we used \rthm{sts fb constructions converge} in the proof of \rthm{strong condensation for backgrounded strategies}, namely in the proof of Case 1. However, we use \rprop{positional} for low level strategies or for the short-tree-component of our strategy, while Case 1 of \rthm{strong condensation for backgrounded strategies} deals with the full lsa type hod premice.}. 

The last remaining case is when $\S^b\insegeq \S'$ and this case is very similar to Case 2 above.
\end{proof}
This finishes the proof of \rthm{sts fb constructions converge}.
\end{proof}

We finish this section by recording some consequences of the proof given above. Suppose $(\P, \Sigma)$ is an sts hod pair. There is one potential problem with our definition of short tree strategy indexing scheme\footnote{See \rdef{sts indexing scheme}.}. Suppose $\M$ is an unambiguous\footnote{See \rdef{unambiguous hp}.} $\Sigma$-sts premouse and $\T$ is an ${\sf{nuvs}}$\footnote{See \rdef{nus stacks}.} stack on $\P$. Suppose $(\gg, \xi, b)$ is an $\M$-minimal shortness witness for $\T$ and let $\Q=\Q(b, \T)$. It is not clear that $\Q$ is a $\Sigma_{\m^+(\T), \T}$-sts premouse. More precisely, it is not clear that $\Sigma^{\Q}\subseteq \Sigma_{\m^+(\T), \T}\rest \Q$. However, the proof of \rthm{sts fb constructions converge} shows that in many situations it is indeed the case that\\\\
(A) $\Q$ a $\Sigma_{\m^+(\T), \T}$-sts premouse, and\\
(B)  $\Sigma(\T)=b$.\\\\
As the proofs are very similar to the proofs already given in the proof of \rthm{sts fb constructions converge}, we will simply state our results. 

\begin{proposition}\label{summary of fb convergence} Suppose $(\P, \Sigma)$ is an sts hod pair and $\Gamma$ is a projectively closed pointclass. Suppose that $\Sigma$ has strong branch condensation and is strongly almost $\Gamma$-fullness preserving. Then the following holds:
\begin{enumerate}
\item Suppose $t=(\P, \T, \P_1, \T_1)$ is $(\P, \Sigma)$-authenticated indexable stack\footnote{See \rdef{authentic stacks of length 2}.}. Then $t$ is according to $\Sigma$.
\item  Suppose $\M$ is a $\Sigma$-sts premouse over some set $X$ and based on $\P$, $\T\in \M$ is an ${\sf{nuvs}}$\footnote{See \rdef{nus stacks}.} stack on $\P$, $(\gg, \xi, b)$ is an $\M$-shortness witness for $\T$\footnote{See \rdef{weak psi alpha indexing scheme a}.} and $\Q=\Q(b, \T)$. Then $\Q$ is a $\Sigma_{\m^+(\T), \T}$-sts premouse over $\m^+(\T)$.
\item Suppose $\M$ is an $\sf{hp}$-indexed germane $\sf{lses}$ such that ${\sf{hl}}(\M)=\P$\footnote{In particular, $\M$ can be viewed as a $\Sigma$-sts premouse over $\P$.} and $\mathcal{J}_{\omega}[\M]\models ``\d^\P$ is a Woodin cardinal". Suppose further that $\Lambda$ is an $\omega_1+1$-iteration strategy for $\M$ such that $\Lambda_\P=\Sigma$ and suppose $(\T, b)\in \M$ is  such that
\begin{itemize}
\item $\T$ is an ${\sf{nuvs}}$,
\item for some $\b$ and $\gg$ such that $\omega\b+\omega\gg\leq \ord(\M)$, setting $t=(\P, \T)$ and $w=(\mathcal{J}_{\omega}(t), t, \in)$, $w$ is $(f, {\sf{sts}})$-minimal as witnessed by $\b$\footnote{See \rdef{important notation}. In particular, this means that we have to index the branch of $t$ at $\omega(\b+\gg)$.},
\item $\gg=\lh(\T)$,
\item $b\in \M|\omega\b$ and $\M|\omega\b\models {\sf{sts}_0}(\T, b)$\footnote{See \rdef{sts0}. This means that $e$ is the branch of $t$ we must choose.} .
\end{itemize}
Then $\Sigma(\T)=b$.
\item Suppose $\Sigma$ is strongly $\Gamma$-fullness preserving and $\M$ is an $\sf{hp}$-indexed germane $\sf{lses}$ such that ${\sf{hl}}(\M)=\P$\footnote{In particular, $\M$ can be viewed as a $\Sigma$-sts premouse over $\P$.} and $\mathcal{J}_{\omega}[\M]\models ``\d^\P$ is a Woodin cardinal". Suppose further that $\Lambda$ is an $\omega_1$-iteration strategy for $\M$ that acts on iterations above $\d^\P$ and suppose $(\T, b)\in \M$ is  such that
\begin{itemize}
\item if $\Lambda^*$ is the $\omega_1$ fragment of $\Lambda$ then ${\sf{Code}}(\Lambda^*)\in \Gamma$,
\item $\T$ is an ${\sf{nuvs}}$,
\item for some $\b$ and $\gg$ such that $\omega\b+\omega\gg\leq \ord(\M)$, setting $t=(\P, \T)$ and $w=(\mathcal{J}_{\omega}(t), t, \in)$, $w$ is $(f, {\sf{sts}})$-minimal as witnessed by $\b$\footnote{See \rdef{important notation}. In particular, this means that we have to index the branch of $t$ at $\omega(\b+\gg)$.},
\item $\gg=\lh(\T)$,
\item $b\in \M|\omega\b$ and $\M|\omega\b\models {\sf{sts}_0}(\T, b)$\footnote{See \rdef{sts0}. This means that $e$ is the branch of $t$ we must choose.} .
\end{itemize}
Then $\Sigma(\T)=b$.
\end{enumerate}
\end{proposition}

The proof of clause 1 of \rprop{summary of fb convergence} is contained in the proof of \rsublem{next lemma}. Clause 2 easily follows from Clause 1 and the relevant definitions. The hypothesis of Clause 3 is exactly what we have at the begining of the proof of \rthm{sts fb constructions converge} (e.g. see \rsublem{iterates are correct sts premice}). Clause 4 follows from the fact that $\Sigma$ is strongly $\Gamma$-fullness preserving. As in the proof of \rthm{sts fb constructions converge} we have that $\Lambda^{*}$ induces a strategy for $\Q(b, \T)$. Thus, if $\Phi$ is this strategy then ${\sf{Code}}(\Phi)\in \Gamma$. Therefore, by strong $\Gamma$-fullness preservation, $\Sigma(\T)=\b$. 

\begin{remark}[On hod pair constructions]\label{on hod pair constructions} Suppose $(\P, \Sigma)$ is an sts hod pair and $\Gamma$ is a projectively closed pointclass. Suppose that $\Sigma$ has strong branch condensation and is strongly almost $\Gamma$-fullness preserving. Recall \rdef{fully backgrounded sts construction}, which introduces fully backgrounded constructions relative to $\Sigma$. In particular, recall the Important Anomaly in clause 3.b of \rdef{fully backgrounded sts construction}. It follows from the clause 4 of \rprop{summary of fb convergence} that, in the terminology of clause 3.b of \rdef{fully backgrounded sts construction}, as long as $\M_\xi$ has an $\omega_1$-iteration strategy (as a $\Sigma$-sts premouse over $\P$) the Important Anomaly cannot occur. $\myqedhere$
\end{remark}


%

\section{The normal-tree comparison theory}\label{sec:normal_comparison}

As in Theorem 2.2.2 of \cite{ATHM}, under $\sf{AD}^+$ and in several other contexts, we can prove a comparison theorem where comparison is achieved via normal trees. In this section we state a comparison theorem for hod pairs that can be applied inside models of $\sf{AD}^+$ and also, inside models satisfying sufficiently rich extensions of $\sf{{\sf{ZFC}}}$, like hod mice themselves. 
Such comparison arguments, among other things, are useful in core model induction arguments and in the analysis of $\H$ of models of $\sf{AD}^+$.

 We start with some general definitions and facts. One warning is that our exposition differs from the one in \cite{ATHM} mainly because we would like to set up our arguments here in a more general setting than the ones stated in \cite{ATHM}. The notation $\insegeq_{hod}$ was introduced in \rdef{layers of hod-like lsp}.

\begin{definition}[Comparison]\label{comparison def} Suppose $(\P, \Sigma)$ and $(\Q, \Lambda)$ are two hod pairs. Then we say that \textbf{comparison holds} for $(\P, \Sigma)$ and $(\Q, \Lambda)$ if there are $(\T, \R)$  and $(\U, \S)$ such that
\begin{enumerate}
\item $\T$ is a stack on $\P$ according to $\Sigma$ with last model $\R$,
\item $\U$ is a stack on $\Q$ according to $\Lambda$ with last model $\S$,

and one of the following holds:
\item $(\Q, \Lambda)$ wins: More precisely the following clauses hold:
\begin{enumerate}
\item $\R\insegeq_{hod}\S$,  
\item $\Lambda_{\R, \U}=\Sigma_{\R, \T}$,
\item $\pi^\T$ is defined,
\item If $\P$ is meek or gentle then $\pi^\U$ is defined,
\item If $\P$ is non meek then letting $\a<\lh(\U)$ be the least such that $\P^b\insegeq \M_\a^\U$, $\pi^\U_{0, \a}$ is defined.
\end{enumerate}
\item $(\P, \Sigma)$ wins: More precisely the following clauses hold:
\begin{enumerate}
\item $\S\insegeq_{hod}\R$,  
\item $\Lambda_{\S, \U}=\Sigma_{\S, \T}$,
\item $\pi^\U$ is defined,
\item If $\Q$ is meek or gentle then $\pi^\T$ is defined,
\item If $\Q$ is non meek then letting $\a<\lh(\T)$ be the least such that $\Q^b\insegeq \M_\a^\T$, $\pi^\T_{0, \a}$ is defined.
\end{enumerate}
\end{enumerate}
If clause 1 holds then we say that $(\Q, \Lambda)$ \textbf{wins the comparison}, and otherwise we say that $(\P, \Sigma)$ wins.
We say normal comparison for $(\P, \Sigma)$ and $(\Q, \Lambda)$ holds if we can take $\T$ and $\U$ to be normal. 

Similarly we define the meaning of ``\textbf{comparison holds} for $(\P, \Sigma)$ and $(\Q, \Lambda)$" in the case $(\P, \Sigma)$ or $(\Q, \Lambda)$ are allowable pairs. For example, if $(\P, \Sigma)$ is a hod pair and $(\Q, \Lambda)$ is an sts hod pair then we say that \textbf{comparison holds} for $(\P, \Sigma)$ and $(\Q, \Lambda)$ if there are $(\T, \R)$ and $(\U, \S)$ such that in the case $(\P, \Sigma)$ wins, $\Sigma_{\S, \T}^{stc}=\Lambda_{\S, \U}$. $\myqedhere$
\end{definition}

%


As in \cite{ATHM}, we can prove comparison for pairs whose corresponding strategies are fullness preserving. Here we show that the fully backgrounded constructions are universal in the sense that they win the comparison with hod pairs that they capture. To establish this fact, we will use \textit{the strategy absorption} argument. The strategy absorption argument was first presented in \cite{ATHM} (see the proof of Theorem 2.28 of \cite{ATHM}) and it builds  on unpublished ideas of Steel. Because we will use the strategy absorption argument several times in this paper and in the next proof, it is important to understand how it works. The general form of the argument is as follows. We have a hod pair $(\P, \Lambda)$ captured by some \textit{background} $(M, \d,  \vec{G}, \Sigma)$. There is also an iteration tree $\T$ on $\P$ according to $\Lambda$ with last model $\Q$ and $\R\insegeq_{hod}\Q$ such that $\R$ is constructed via some hod pair construction of $M$. It is additionally required that the background extenders used to build $\R$ cohere $\Lambda$\footnote{However, the fact that ${\sf{Code}}(\Lambda)$ is Suslin, co-Suslin captured by $(M, \d,  \vec{G}, \Sigma)$ implies that all extenders in $\vec{G}$ cohere $\Lambda$.}. The goal of the argument is to show that the strategy $\R$ inherits from the background universe is the same as $\Lambda_{\R, \T}$. In many cases, this can be done by appealing to branch condensation and the existence of minimal disagreements. Here is how a typical argument works. 

Let $\Phi$ be the iteration strategy of $\R$ induced by the background strategy. Fix $\U$ on $\R$ that is according to both $\Lambda_{\R, \T}$ and $\Phi$ but $\Lambda_{\R, \T}(\U)\not =\Phi(\U)$. Let $\U^*$ be the stack on $M$ obtained by resurrection  process. Thus, $\U^*=r\U$ (see \rsubsec{the induced strategy for the undropping game}) . Let $b=\Phi(\U^*)$. We then have that $\pi^{\U^*}_b(\T)$ is according to $\Lambda$ (this is where we use coherence). Then branch condensation is applied to the equality \begin{center}
$\pi^{\pi^{\U^*}_b(\T)}=\sigma\circ \pi^{\U}_b\circ \pi^{\T}$ 
\end{center}
where $\sigma:\M^{\U}_b\rightarrow \pi^{\U^*}_b(\R)$ is the canonical factor map that the resurrection process gives (in particular, $\pi^{\U^*}_b\rest \R=\sigma\circ \pi^{\U}_b$). The reader may wish to review \rsubsec{the induced strategy for the undropping game}.  Recall that strong branch condensation and $\Gamma$-fullness preservation implies positional (see \rprop{positional}).

\begin{theorem}[Universality of backgrounded construction]\label{universality of background construction} Assume $\sf{AD}^+$.\\ Suppose that 
\begin{itemize}
\item $\Gamma$ is a pointclass, 
\item $(\P, \Lambda)$ is an allowable pair, 
\item $k(\P)={\sf{ep}}(\P)$\footnote{Recall that $\P$ is f.s $\mathcal{J}$-structure. To define ${\sf{ep}}(\P)$ we ignore its fine-structural component $k(\P)$ and treat $\P$ as just a $\mathcal{J}$-structure. See \rdef{the core}.},
\item $\Lambda$ is $\Gamma$-fullness preserving and has strong branch condensation\footnote{We could instead assume just the first two clauses of strong branch condensation and also that $\Lambda$ is self-cohering. However, our proof will use self-cohering in an indirect way. Strategies with strong branch condensation are positional and therefore, self-cohering. The reader may wish to review \rdef{self-cohering}, \rdef{strong branch condensation}, \rthm{strong condensation for backgrounded strategies} and \rprop{positional}.}, 
\item ${\sf{C}}=(\mathbb{M}, (P, \Psi), \Gamma^*, A)$ Suslin, co-Suslin captures both $\Gamma$ and ${\sf{Code}}(\Lambda)$, and \item $\mathbb{M}=(M, \d,  \vec{G}, \Sigma)$.
\item $\sf{NsesS}$.
\end{itemize} Let
\begin{center}${\sf{hpc}}_{{\sf{C}}, \Gamma}^+=(\M_\gg , \N_\gg, Y_\gg, \Phi^+_\gg, F^+_\gg, F_\gg, b_\gg : \gg\leq \delta)$
\end{center} 
be the output of $\Gamma-\sf{hpc}$ of $\mathbb{M}$ with the property that each $F^+_\gg$ coheres $\Lambda\rest M$\footnote{This actually follows from the fact that ${\sf{Code}}(\Lambda)$ is Suslin, co-Suslin captured.}.
Then there is $\gg\leq \d$ such that the following holds.
\begin{enumerate}
\item If $(\P, \Lambda)$ is not an sts hod pair then $\gg<\d$ and there is a normal stack $\X$ such that $(\X, \M_\gg)\in I(\P, \Lambda)$ and $\Phi^+_\gg=\Lambda_{\M_\gg}$.
\item  If $(\P, \Lambda)$ is an sts hod pair then there is a normal stack $\X$ such that letting
\begin{center}
$\N=\begin{cases}
\M_\gg&: \gg<\d\\
\M_\gg^{\#} &: \gg=\d,
\end{cases}$
\end{center}
$(\X, \N)\in I(\P, \Lambda)$ and $\Phi^+_\gg=\Lambda_{\N}$.
\item For every $\b\leq \gg$, there is a $\Lambda$-iterate $\R$ of $\P$ via a normal stack $\T$ such that $\M_\b\insegeq \R$ and if $\S\in Y_{\b}$ then $(\Phi^+_\b)_\S=\Lambda_{\S}$.
\item For every $\b\leq \gg$, there is a $\Lambda$-iterate $\R$ of $\P$ via a normal stack $\T$ such that $\N_\xi\insegeq \R$.
\end{enumerate}
\end{theorem}
\begin{proof} As in the proof of Lemma 2.10 of \cite{ATHM}, in the comparison of $\P$ with the models of ${\sf{hpc}}_{{\sf{C}}, \Gamma}$ no extender disagreement appears on ${\sf{hpc}}_{{\sf{C}}, \Gamma}$ side. Many of the details of the argument have appeared in \cite[Lemma 3.21]{trang2013}, and because of this we only concentrate on the new aspects of the proof. We then assume that $\P$ is non-meek. 

We first show that clause 3 holds and then show that clause 1 and 2 hold. Clause 4 is similar to clause 3. To prove clause 3, we only verify that\\\\
(1) for every $\b<\d$, if $(\T, \R, \S)$ are such that $\T$ is a normal stack on $\P$ according to $\Sigma$, $\R$ is the last model of $\T$, $\M_\b\insegeq \R$ and $\S\in Y_\b$ then $(\Phi^+_\b)_\S=\Lambda_{\S}$.\\\\
The proof of (1) is the portion of the proof that goes beyond \cite[Theorem 2.28]{ATHM} and \cite[Lemma 3.21]{trang2013}, and so we prove (1). 
 
 Because $\Lambda$ is self-cohering (see \rdef{self-cohering}) we can in fact assume that\\\\
 (2) for every $\a+1<\lh(\T)$, $\S\not \insegeq \M_\a^\T$.\\\\
  For simplicity, we prove $(\Phi^+_\b)_\S=\Lambda_{\S}$ for ordinary stacks as opposed to generalized stacks. The more general proof is only notationally more complex. The reader may wish to review \rsec{the induced strategy for the undropping game} and \rlem{summary res and emb}.

Towards a contradiction, we assume that (1) fails. Let $(\b, \S, \T, \R)$ witness the failure of (1) such that $\b$ is the least possible and (2) holds. We assume that $\S\in Y^\R\cap Y_{\b}$ is the least layer for which (1) fails. Let $\Phi=\Phi^+_\b$ and $\Q=\M_\b$. \\\\
\textbf{Case 1:} $\S$ is of successor type. \\\\
Then we get a contradiction using branch condensation of $\Lambda$. Let $\U$ be a stack on $\S$ such that it is according to both $\Phi_\S$ and $\Lambda_{\S}$ but $\Phi(\U)\not =\Lambda_{\S}(\U)$. Let $b=\Phi(\U)$ and $c=\Lambda_{\S}(\U)$. Let $\U^*=r\U$\footnote{This is the resurrection of $\U$. See \rsec{the induced strategy for the undropping game}.}. Then because extenders used to construct $\Q$ cohere $\Lambda$, we have that  $\pi^{\U^*}(\T)$ is according to $\Lambda$. Let $N$ be the last model of $\U^*$.\\ 

\textit{Claim.} $\pi^{\U}_b$ exists.\\
\begin{proof}
The claim is a consequence of $\Gamma$-fullness preservation and the fact that $\Phi_{\S^-}=\Lambda_{\S^-}$. Towards a contradiction assume that $\pi^\U_b$ is undefined. Because $\Phi_{\S^-}=\Lambda_{\S^-}$ and because $\Phi(\U)\not =\Lambda_{\S}(\U)$, we must have some $\iota\in R^\U$ such that $\pi^\U_{0, \iota}$ is defined and $\U_{\geq \iota}$ is above $\M_\iota^-$. Moreover, because $\Phi(\U)\not =\Lambda_{\S}(\U)$, $R^\U$ must have a maximal element. Let then $\iota=\max(R^\U)$ and set $\X=\U_{\geq \iota}$.

Suppose now that $\pi^\U_c$ is not defined. $\Gamma$-fullness preservation implies that $\X$ does not have fatal drop, and so  $\d(\X)$ is a strong cutpoint in both $\Q(b,\U)$ and $\Q(c,\U)$. Hence $\Gamma$-fullness preservation implies that $b=c$. Contradiction. 

Thus,  $\pi^\U_c$ is defined. It then follows from $\Gamma$-fullness preservation that $\Q(b, \X) \insegeq \M^{\U}_c$. Therefore, $b = c$, which is a contradiction.
\end{proof}

Let then $\U^+= \uparrow (\R, \U)$\footnote{See \rdef{upward extension of a stack}.} be the unique stack on $\R$ whose tree structure and extenders are exactly those of $\U$. Let $\R^*$ be the last model of $\U^+$. We then have $\sigma:\R^*\rightarrow \pi^{\U^*}_b(\R)$ such that, assuming $\pi^\T$ is defined,
\begin{center}
$\pi^{\pi^{\U^*}_b(\T)}=\sigma\circ \pi^{\U^+}_b\circ \pi^{\T}$.
\end{center} 
Notice next that it follows from (2) and the fact that $\S$ is of successor type that $\pi^\T$ is defined. Branch condensation of $\Lambda$ and the displayed equality implies that $b=c$\\\\
\textbf{Case 2.} $\S$ is of limit type.\\\\
 Then by appealing to \rlem{minimal low level disagreement exist}, we can fix some $((\U_1, \S_1), (\U_2, \S_2), \S_3)$ that constitutes a minimal low level disagreement between $\Phi$ and $\Lambda_{\S}$. Let $\U_1^*=r\U_1$, $N$ be the last model of $\U_1^*$ and $\sigma: \S_1\rightarrow \Q$ be elementary such that $\Q$ is an $\M$-model appearing in the $\Gamma-{\sf{hpc}}$ of $N$. We then have that letting $\T^*=\pi^{\U_1^*}(\T)$, for some $\xi$, $\Q\insegeq \M_\xi^{\T^*}$. Let $\Q'=\sigma(\S_3)$. Notice next that\\\\
 (3) $(\T^*, \Q)$ supports a bottom type $(\Q', \S_3, \sigma\rest \S_3)$-b-condensation diagram on $\P$. \\\\
 Let $\Phi^*$ be the strategy of $\Q'$ induced by $\Sigma_N$. We then have that\\\\
 (4) $\Phi_{\S_3, \U_1}=(\sigma$-pullback of $\Phi^*)$ and $\Lambda_{\Q'}=\Phi^*$.\\\\
 $\Lambda_{\Q'}=\Phi^*$ follows from the fact that $\Q'\inseg \Q$ and also from the fact that $\b$ is the least satisfying (1).
 Thus, the strong branch condensation of $\Lambda$ implies that $\Lambda_{\S_3}=\Phi_{\S_3, \U_1}$. \\

 Next, we need to  verify that clause 1 and 2 hold for some $\gg<\d$. Set $\N=\M_\d$. Assume first that $(\P, \Lambda)$ is not an sts hod pair. This means that $\Lambda$ is an iteration strategy. Assume then clause 1 fails. It follows that we have a normal stack $\X$ on $\P$ such that $\lh(\X)=\d$ and $\m(\X)=\N$. Let $b=\Lambda(\X)$. Let $\a<\lh(\X)$ be least such that $\d\in \rge(\pi^\X_{\a, b})$. Because the entire construction takes place in $M$ and because $\d$ is regular, we have that letting $\eta$ be such that $\d=\pi^{\X}_{\a, b}(\eta)$, $\eta$ must be a measurable cardinal of $\M_\a^\X$. 
 
 Notice that $\M_\a^\X$ is germane\footnote{See \rdef{germane lses}.} and  because $\X$ may drop in model, $\M_\a^\X$ may not be hod-like. Let $\R\insegeq \N$ be the longest such that
 \begin{itemize}
 \item $\R\in Y^{\N}$, 
 \item $\R$ is meek or gentle, 
 \item $\R\insegeq \M_\a^\X$ and 
 \item $\d^\R\leq \eta$.
 \end{itemize}
 We now break into cases. Let $\a'\leq \a$ be the least such that $\R\insegeq \M_{\a'}^\X$. 
 
 Suppose first that $\R$ is of successor type. We must then have $\R'\in Y^{\M^\X_{\a'}}$ such that $(\R')^-=\R$. But now $\X_{\geq \a'}$ is based on $\R'$ and is above $\d^\R$. Because in this case $\R'$ out-iterates $\N$, this contradicts our assumption that $\sf{NsesS}$ \footnote{This is a consequence of the ordinary universality of the background constructions. If a mouse outiterates a fully backgrounded construction then it generates a mouse with a superstrong. See \cite[Lemma 2.13]{ATHM}.}. 
 
 Suppose then $\R$ is gentle. In this case, we must have $\R'$ such that $\R'$ is meek of limit type, $\d^{\R'}=\d^\R$ and $\R'\in Y^{\M^\X_{\a'}}$ or $\R'=\M^\X_{\a'}$. If $\d^\R<\eta$ then we get that $\X_{\geq \a'}$ is based on $\R'$ and is above $\d^\R+1$, and once again this leads to a contradiction. 
 
 Suppose now that $\d^\R=\eta$. Let now $\k> \eta$ be such that it reflects 
 \begin{itemize}
 \item $\X$, and
 \item ${\sf{hpc}}_{{\sf{C}}, \Gamma}=(\M_\gg , \N_\gg, Y_\gg, \Phi_\gg, F^+_\gg, F_\gg, b_\gg : \gg\leq \delta)$.
 \end{itemize}
Let $\xi=o^{\N}(\k)$\footnote{Notice that $\xi<\d$ as otherwise $\d$ would be a Woodin cardinal of $\M^\X_b$ and since it is also a measurable cardinal, we would get a contradiction to our minimality assumption on hod mice.} and $\zeta+1\in b$ be such that $\X(\zeta+1)=\k$.  Let $E\in \vec{G}$ be an extender such that $\cp(E)=\k$, $\lh(E)>\xi$ and it reflects the above sets.  It follows that\\\\
 (5) $\k\in b$ and $\cp(\pi^\X_{\k, b})=\k$,\\
 (6)  $E_{\zeta}^\X$ agrees with $E$,\\
 (7) $\pi^\X_{\a', \k}$ is defined and $\pi^{\X_{\a', \k}}(\d^\R)=\k$,\\
 (8) $\M_\k^\X|\xi=\N|\xi$\\
 (9) $\M^\X_{\zeta+1}|\ind_\zeta^\X=\N|\ind_\zeta^\X$.\\\\
It follows from (6), (8) and (9) that $E_\zeta^\X\in \vec{E}^{\N}$ as $E$ can serve as a background certificate for it. Clearly this is a contradiction.
 
 Finally suppose $\R$ is of limit type. In this case we have that $\d^\R<\eta$. We also have $\R'\in Y^{\M^\X_{\a'}}$ such that either $\R'$ is of successor type and $(\R')^-=\R$ or $\R'$ is of limit type and $(\R')^b=\d^\R$. The first case leads to a contradiction via a similar argument as the one given above. Let then $\R'$ be a complete layer of $\M^\X_{\a'}$ (see \rnot{l p}) such that $(\R')^b=\R$. It follows that $\X_{\geq \a'}$ is based on $\R'$, $\X_{\geq \a'}$ is above $\d^\R$ and also, that $\N^b=\R^b$. This case once again leads to a contradiction because assuming $\sf{NsesS}$ universality implies that $\R'$ cannot out-iterate $\N$. 

The case that $(\P, \Lambda)$ is an sts hod pair is very similar. In this case, we note that $\X$ must be $\Lambda$-maximal as otherwise $\Lambda(\X)$ is a branch and all of the above arguments can be repeated. If $\X$ is $\Lambda$-maximal then $\Lambda(\X)=\N^\#$, which is one of the possibilities in clause 2.
\end{proof}

As a corollary to \rthm{universality of background construction} we get that comparison holds.

\begin{corollary}\label{comparison holds 1} Assume $\sf{AD}^+$ and suppose 
\begin{itemize}
\item $\Gamma$ is a pointclass, 
\item $(\P, \Sigma)$ and $(\Q, \Lambda)$ are two allowable pairs such that both $\Sigma$ and $\Lambda$ are $\Gamma$-fullness preserving and have strong branch condensation, 
\item $k(\P)={\sf{ep}}(\P)$ and $k(\P)={\sf{ep}}(\Q)$,
\item there is a good pointclass $\Gamma'$ such that $\Gamma\cup \{{\sf{Code}}(\Lambda), {\sf{Code}}(\Sigma)\} \subseteq {\mathbf{\Delta}}_{\Gamma'}$.
\end{itemize} Then the normal comparison holds for $(\P, \Sigma)$ and $(\Q, \Lambda)$.
\end{corollary}
\begin{proof} Using \rthm{n*x} we can find ${\sf{C}}=(\mathbb{M}, (P, \Psi), \Gamma^*, A)$ which Suslin, co-Suslin captures $\Gamma$, ${\sf{Code}}(\Lambda)$ and ${\sf{Code}}(\Sigma)$. 
Let
\begin{center}${\sf{hpc}}_{{\sf{C}}, \Gamma}^+=(\M_\gg , \N_\gg, Y_\gg, \Phi^+_\gg, F^+_\gg, F_\gg, b_\gg : \gg\leq \delta)$
\end{center} 
be the output of $\Gamma-\sf{hpc}$ of $\mathbb{M}$ with the property that each $F^+_\gg$ coheres both $\Sigma\rest M$ and $\Lambda\rest M$.

It follows from \rthm{universality of background construction} that there are $\b, \gg\leq \d$ and normal stacks $\T$ and $\U$ such that
\begin{enumerate}
\item $(\T, \M_\b)\in I(\P, \Sigma)$ and $\Phi^+_\b=\Sigma_{\M_\b}$ and
\item  $(\U, \M_\gg)\in I(\Q, \Lambda)$,  $\Phi^+_{\gg}=\Lambda_{\M_\gg}$.
\end{enumerate}
If $\b=\gg$ then clearly the normal comparison for $(\P, \Sigma)$ and $(\Q, \Lambda)$ holds. Suppose then $\b<\gg$. Then there is $(\U', \Q')$ such that 
\begin{itemize}
\item $\U'$ is a normal stack on $\Q$ according to $\Lambda$,
\item $\Q'$ is the last model of $\U'$, and
\item $\M_\b\insegeq_{hod}\Q'$ and $\Phi_\b=\Lambda_{\M_\b}$. 
\end{itemize}
Let $\a<\lh(\U')$ be the least such that $\M_\b\insegeq \M_\a^{\U'}$. Set $\X=\U'_{\leq \a}$ and $\S=\M_\a^\X$ and $\R=\M_\b^\T$. In order to show that $(\T, \R)$ and $(\X, \S)$ witnesses that comparison holds for $(\P, \Sigma)$ and $(\Q, \Lambda)$, we need to show that $[0, \a)_\X\cap D^\X=\emptyset$. However, to obtain this condition we may need to change $\X$.

First observe that if $\P$ is meek or gentle then indeed $[0, \a)_\X\cap D^\X=\emptyset$. We give the argument in the case $\P$ is of successor type and as the rest is similar, we leave the rest to the reader. Since $\P$ is of successor type, we have that $\d^\R$ is a cardinal of $M$. Notice that $\a$ is the least $\a'$ such that $\M_{\a'}^\X|\d^\R=\R|\d^\R$. This follows from $\Gamma$-fullness preservation, which implies that if $\M_{\a'}^\X|\d^\R=\S|\d^\R$ then $\R\insegeq \M_{\a'}^\X$. Thus, $\a$ must be a limit ordinal. Suppose then $[0, \a)_\X\cap D^\X\not =\emptyset$. It follows that $\rho(\S)<\d^\R$. But hod premice do not project across layers of successor type (or rather meek or gentle type)\footnote{See \rdef{pre-hod-like}.}.

Suppose then that $\P$ is non-meek. Let $\iota$ be the least such that $\M_\iota^\X|\ord(\R^b)=\R^b$. It follows from the argument above that $[0, \iota)_{\X}\cap D^\X=\emptyset$. Moreover, $\X_{\geq \iota}$ is above $\ord(\R^b)$. Set $\k=\d^{\R^b}$.

 Suppose now that there is $E\in \vec{E}^{\S}$ such that $\cp(E)=\k$ and $\ind^\S(E)>\ord(\R)$. Let $\X'$ be the continuation of $\X$ obtained by using $E$ at stage $\a$. Notice that $E$ must be applied to $\M_{\iota}^\X$. As $\X_{\leq \iota}$ doesn't have drop on its main branch, we have that $\X'$ also doesn't have a drop on its main branch and moreover, $(\T, \R)$ and $(\X', Ult(\M_\b^\X, E))$ witness that comparison holds for $(\P, \Sigma)$ and $(\Q, \Lambda)$.
\end{proof}

Using reflection, we can eliminate the extra assumptions on $\Gamma$ and the two strategies. 

\begin{corollary}[Comparison]\label{comparison holds} Assume $\sf{AD}^+$ and suppose $\Gamma$ is a pointclass. Suppose further that $(\P, \Sigma)$ and $(\Q, \Lambda)$ are two hod pairs such that 
\begin{itemize}
\item both $\Sigma$ and $\Lambda$ are $\Gamma$-fullness preserving and have branch condensation, 
\item  $k(\P)={\sf{ep}}(\P)$ and $k(\P)={\sf{ep}}(\Q)$,
\end{itemize}
Then the normal comparison hold for $(\P, \Sigma)$ and $(\Q, \Lambda)$.
\end{corollary}
\begin{proof} Suppose not. Applying $\Sigma^2_1$-reflection, we can find $\Gamma^*$ and two hod pairs $(\P_1, \Sigma_1)$ and $(\Q_1, \Lambda_1)$ such that $\Gamma^*\cup\{ {\sf{Code}}(\Sigma_1), {\sf{Code}}(\Lambda_1)\}\subseteq {\mathbf{\Delta}}^2_1$ and the claim of the corollary fails for $(\Gamma^*, (\P_1, \Sigma_1), (\Q_1, \Lambda_1))$. We then apply \rcor{comparison holds 1}. We use \rthm{n*x} to get a ${\sf{C}}=(\mathbb{M}, (P, \Psi), \Gamma^*, A)$ that Suslin, co-Suslin captures $\Gamma$, ${\sf{Code}}(\Lambda)$ and ${\sf{Code}}(\Sigma)$.
\end{proof}

\begin{remark}\label{ep remark} In most situations, our allowable pairs $(\P, \Sigma)$ will have the property that $k(\P)={\sf{ep}}(\P)$. Thus, we make a convention that unless otherwise specified, all allowable pairs have the property that $k(\P)={\sf{ep}}(\P)$. When it is necessary we will remind the reader of this.

However, \rthm{universality of background construction} can also be proven in the case that $k(\P)<{\sf{ep}}(\P)$. In this case, what we get is that letting $k=k(\P)$, $(\X, (\C_k(\N_\gg), k))\in I(\P, \Sigma)$. Similar results can also be proven for germane $\sf{lses}$. $\myqedhere$
\end{remark}

\section{Diamond comparison}\label{sec:diamond comparison}

Our goal here is to provide another comparison argument, \textit{diamond comparison}, that doesn't rely on branch condensation as heavily as our other argument (see \rcor{comparison holds}). The new comparison argument follows the same line of thought as the proof of a similar comparison argument from \cite{ATHM} (see Theorem 2.47 of \cite{ATHM}). 


As in \cite{ATHM}, the diamond comparison argument can be used to show that $\sf{AD}^+ + \sf{LSA}$ is consistent relative to a Woodin cardinal that is a limit of Woodin cardinals. This will appear as \rthm{lst from wlw}. In \cite{ATHM}, a similar argument gave the consistency of ${\sf{AD}}_\bR+``\Theta$ is a regular cardinal" relative to a Woodin cardinal that is a limit of Woodin cardinals. 

Following the proof of Theorem 2.47 of \cite{ATHM}, we first define a \textit{bad block} and a \textit{bad sequence} and show that there cannot be such a bad sequence of length $\omega_1$. We then show that the failure of comparison produces such bad sequences of length $\omega_1$.

\subsection{Bad sequences}

For the purposes of this subsection, we make a definition of a bad block and a bad sequence. In later subsections, we will redefine these names for different objects. Below and elsewhere, if $\T$ is a stack of successor length then we let $\T^-$ be $\T_{<\a}$ where $\a+1=\lh(\T)$. 

\begin{definition}[Bad block]\label{bad block} \index{bad blocks}Suppose $(\P, \Sigma)$ and $(\Q, \Lambda)$ are two hod pairs such that both $\P$ and $\Q$ are of limit type and are not gentle. Then 
\begin{center}
$B=(( (\P_i, \Q_i, \Sigma_i, \Lambda_i): i<4), ( \T_i , \U_i: i<3), (c, d) )$
\end{center} 
is a bad block on $((\P, \Sigma), (\Q, \Lambda))$ if the following holds:
\begin{enumerate}
\item $(\P_0, \Sigma_0)=(\P, \Sigma)$ and $(\Q_0, \Lambda_0)=(\Q, \Lambda)$.
\item $\T_0$ is a stack according to $\Sigma_0$ on $\P$.
\item  $\U_0$ is a stack according to $\Lambda_0$ on $\Q$.
\item Let $\T_0=( \M_\b, \T^*_\b, E_\b: \b\leq \nu)$ and $\U_0=( \N_\b, \U^*_\b, F_\b: \b\leq \nu)$. Then $\T^*_\nu$ and $\U^*_\nu$ are undefined, $\P_1=\M_{\nu}$ and $\Q_1=\N_{\nu}$.
\item There is $\K$ such that $\K\inseg_{hod}\P_1$, $\K\inseg_{hod} \Q_1$, $\K$ is of successor type, $\Sigma_{\K, \T_0}\not=\Lambda_{\K, \U_0}$ and 
$\Sigma_{\K, \T_0}=\Sigma_{\K, \U_0}$.

\item $\T_1$ and $\U_1$ are stacks on $\P_1$ and $\Q_1$ respectively with last models $\P_2$ and $\Q_2$ such that  $\pi^{\T_1}$ and $\pi^{\U_1}$ are defined, $\pi^{\T_1}(\K)=\pi^{\U_1}(\K)$ and setting $\K'=\pi^{\T_1}(\K)$, $\Sigma_{\K', \T_0^\frown \T_1}=\Lambda_{\K', \U_0^\frown \U_1}$\footnote{Because of \rthm{comparison holds} we can take $\T_1$ and $\U_1$ to be normal trees. We will always use the diamond comparison argument in situations where \rthm{comparison holds} applies to low level strategies.}.
\item $\T_1$ and $\U_1$ have a last normal component of successor length whose predecessor is a limit ordinal\footnote{Recall that in \rdef{comparison stack}, we required that comparison stacks have a last model.} and $\T_1^-=\U_1^-$.
\item $c=\Sigma_{\P_1, \T_0}(\T_1^-)$, $d=\Lambda_{\Q_1, \U_0}(\U_1^-)$\footnote{Thus, $\P_2=\M^{\T_1^-}_{c}$ and $\Q_2=\M^{\T_1^-}_{d}$.}.
\item $\Sigma_1=\Sigma_{\P_1, \T_0}$, $\Sigma_2=\Sigma_{\P_2, \T_0^{\frown}\T_1}$, $\Lambda_1=\Sigma_{\Q_1, \U_0}$, and $\Lambda_2=\Sigma_{\Q_2, \U_0^{\frown}\U_1}$,
\item  $(\T_2, \P_3)\in I(\P_2, \Sigma_2)\cap I^{ope}(\P_2, \Sigma_2)$ and $(\U_2, \Q_3)\in I(\Q_2, \Lambda_2)\cap I^{ope}(\Q_2, \Lambda_2)$,
\item $\Sigma_3=(\Sigma_2)_{\P_3, \T_2}$ and $\Lambda_3=(\Lambda_2)_{\Q_3, \U_2}$. 
\item $\P_3^b=\Q_3^b$ and $(\Sigma_3)_{\P_3^b}=(\Lambda_3)_{\Q_3^b}$.
\end{enumerate}
We set $\T^{B}=\T_0^{\frown}\T_1^{\frown}\T_2$ and $\U^{B}=\U_0^{\frown}\U_1^{\frown}\U_2$. We say $\T^{B}$ is the stack on the top of $B$ and $\U^{B}$ is the stack in the bottom of $B$. $\myqedhere$
\end{definition}

Next we show that there cannot be a bad sequence of length $\omega_1$. 

\begin{lemma}[No bad sequences, $\sf{ZF+DC}$]\label{no bad sequences} Suppose $(\P, \Sigma)$ and $(\Q, \Lambda)$ are two hod pairs of limit type such that $\P$ and $\Q$ are countable, and both $\Sigma$ and $\Lambda$ are $(\omega_1, \omega_1, \omega_1)$-strategies. There is then no bad sequence, i.e., a sequence $( B_\b: \b<\omega_1)$ satisfying the following conditions:
\begin{enumerate}
\item For all $\b<\omega_1$, $B_\b=(( (\P_{\b, i}, \Q_{\b,i}, \Sigma_{\b, i}, \Lambda_{\b, i}): i<4), ( \T_{\b, i} , \U_{\b, i}: i<3), (c_\b, d_\b) )$.
\item For all $\b<\omega_1$, $B_\b$ is a bad block on $((\P_{\b, 0}, \Sigma_{\b, 0}), (\Q_{\b, 0}, \Lambda_{\b, 0}))$.
\item For all $\b<\omega_1$, $\P_{\b+1, 0}=\P_{\b, 3}$ and $\Q_{\b+1, 0}=\Q_{\b, 3}$.
\item For $\b<\a<\omega_1$, let $\pi_{\b, \a}:\P_{\b, 0}\rightarrow \P_{\a, 0}$ be the composition of the embeddings on the ``top" and $\sigma_{\b, \a}:\Q_{\b, 0}\rightarrow \Q_{\a, 0}$ be the composition of the embeddings on the ``bottom". Then for all limit $\l<\omega_1$, $\P_{\l, 0}$ is the direct limit of $(\P_\a, \pi_{\a, \b} : \a<\b<\l)$. Similarly, for all limit $\l<\omega_1$, $\Q_{\l, 0}$ is the direct limit of $(\Q_\a, \sigma_{\a, \b}: \a<\b<\l)$ under the maps $\sigma_{\b, \a}$.
\item For all limit ordinals $\l<\omega_1$, $\P_{\l, 0}^b=\Q_{\l, 0}^b$. 
\item For all $\b<\omega_1$, $\Sigma_{\b, 0}=\Sigma_{\P_{\b, 0}, \oplus_{\gg<\b}\T^{B_\gg}}$ and $\Lambda_{\b, 0}=\Sigma_{\Q_{\b, 0}, \oplus_{\gg<\b}\U^{ B_\gg}}$.
\end{enumerate}
\end{lemma}
\begin{proof}
Towards a contradiction, suppose $\vec{B}=( B_\b: \b<\omega_1)$ is a bad sequence. Let $\P_{\omega_1}$ be the direct limit of $( \P_{\a, 0}, \pi_{\a, \b} : \a<\b<\omega_1)$ and $\Q_{\omega_1}$ be the direct limit of $( \Q_{\a, 0}, \sigma_{\a, \b} : \a<\b<\omega_1)$. Let $N=L((\P, \Sigma), (\Q, \Lambda), \vec{B}, \bR)$, $\zeta=\Theta^N$ and $X$ be a countable submodel of $N|(\zeta^+)^N$ such that letting $\tau:M\rightarrow N|(\zeta^+)^N$ be the uncollapse map, $\vec{B}\in \rge(\tau)$. Let $\k=\omega_1^M$ and notice that for every $\b<\k$, 
\begin{center}
$B_\b^-=_{def}(( (\P_{\b, i}, \Q_{\b,i}): i<4), ( \T_{\b, i} , \U_{\b, i}: i<3), (c_\b, d_\b) )\in M$
\end{center}
 and $B_\b^-$ is countable in $M$. It then follows that $\tau^{-1}(\P_{\omega_1})=\P_{\k, 0}$ and $\tau^{-1}(\Q_{\omega_1})=\Q_{\k, 0}$. Let 
 \begin{center}
 $\pi_{\b}:\P_{\b, 0}\rightarrow \P_{\omega_1}$ and $\sigma_{\b}:\Q_{\b, 0}\rightarrow \Q_{\omega_1}$
 \end{center}
  be the direct limit embeddings. 

Standard arguments show that for all $x\in\P_{\k, 0}\cap \Q_{\k, 0}$, 
\begin{center}
$\pi_{\k}(x)=\tau(x)=\sigma_\k(x)$.
\end{center}
Notice that $\P_{\k, 0}^b=\Q_{\k, 0}^b$ (see clause 5 of our hypothesis). Set $\d=\d^{\P^b_{\k, 0}}$ and let $\phi=\pi^{\T_{\k, 0}}$ and $\psi=\pi^{\U_{\k, 0}}$. We now have that \\\\
(1) $\P_{\k, 0}^b=\Q_{\k, 0}^b$ and $\pi_\k\rest \P_{\k, 0}^b=\sigma_\k\rest \Q_{\k, 0}^b$ .\\\\
 Let 
 \begin{itemize}
 \item $\K$ witness clause 5 of \rdef{bad block} for $B_\k$,
 \item  $p=\pi^{\T_{\k, 1}^-}_{c_{\k}}$ and $q=\pi^{\T_{\k, 1}^-}_{d_{\k}}$, 
 \item $i:\P_{\k, 2}\rightarrow \P_{\omega_1}$ and $j:\Q_{\k, 2}\rightarrow \Q_{\omega_1}$ be the iteration embeddings along the top and bottom of $\vec{B}$. 
 \end{itemize}
  Notice that because 
 \begin{center}
 $(\Sigma_{\k, 1})_{\K^-}=(\Lambda_{\k, 1})_{\K^-}$,
\end{center}
we have that\\\\
(2) $i \circ p\rest \K^-= j\circ q \rest \K^-$.\\\\
Next it follows from \rlem{the canonical singularizing sequence for g-stacks} that\\\\
(3) $\d^{\K}=\sup\{\phi(f)(a) : f\in \P_{\k, 0}\wedge f:\d\rightarrow \d \wedge a\in (\K^-)^{<\omega}\}$\\
(4) $\d^{\K}=\sup\{\psi(f)(a) : f\in \Q_{\k, 0}\wedge f:\d\rightarrow \d \wedge a\in (\K^-)^{<\omega}\}$\\\\
Because 
\begin{center}
$\K'=_{def}p(\K)=q(\K)$ and $(\Sigma_{\k, 2})_{\K'}=(\Lambda_{\k, 2})_{\K'}$,
\end{center}
 we have that\\\\
(5) $i\rest \K'= j\rest \K'$.\\\\
Let then 
\begin{center}
$s=\{\phi(f)(a) : f\in \P_{\k, 0}\wedge f:\d\rightarrow \d \wedge a\in (\K^-)^{<\omega}\}$\\
 $t=\{\psi(f)(a) : f\in \Q_{\k, 0}\wedge f:\d\rightarrow \d \wedge a\in (\K^-)^{<\omega}\}$.
 \end{center}
 (1) and (2) then imply that\\\\
(6) $i\circ p [s]=j\circ q[t]$.\\\\
(5) and (6) now imply that\\\\
(7) $p[s]=q[t]$.\\\\
It follows from (3), (4) and (7) that\\\\
(8) $p[s] \cap q[t]$ is cofinal in $\d^{\K'}$.\\\\
It then follows that $c_{\k}=d_{\k}$, contradiction.
\end{proof}

\subsection{The comparison argument}

In this subsection we prove the following comparison theorem under the hypothesis that the \textit{lower level comparison} holds.
Suppose $(\P, \Sigma)$ and $(\Q, \Lambda)$ are two hod pairs of limit type such that $\Gamma(\P, \Sigma)=\Gamma(\Q, \Lambda)=_{def}\Gamma$, both $\Sigma$ and $\Lambda$ are $\Gamma$-fullness preserving.

\begin{definition}[Lower Level Comparison]\label{low level comparison}\index{low level comparison}
We say low level comparison holds for hod pairs or sts hod pairs $(\P, \Sigma)$ and $(\Q, \Lambda)$ if 
\begin{enumerate}
\item for every $(\T, \P_1)\in B(\P, \Sigma)$ and $(\U, \Q_1)\in B(\Q, \Lambda)$, comparison holds for $(\P_1, \Sigma_{\P_1, \T})$ and $(\Q_1, \Lambda_{\Q_1, \U})$, and
\item whenever $(\T, \P_1)\in I(\P, \Sigma)$, $(\U, \Q_1)\in I(\Q, \Lambda)$ and $\K$ are such that
\begin{itemize}
\item $\K\insegeq_{hod} \P_1$ and $\K\insegeq_{hod}\Q_1$,
\item $\K$ is of successor type and,
\item $\Sigma_{\K^-, \T}=\Lambda_{\K^-, \U}$, 
\end{itemize}
there is a normal stack $\S$ of limit length according to both $\Sigma_{\P_1, \T}$ and $\Lambda_{\Q_1, \U}$ that is based on $\K$ and is such that letting $b=\Sigma_{\P_1, \T}(\S)$ and $c=\Lambda_{\Q_1, \U}(\S)$,
\begin{enumerate}
\item $\pi^{\S}_b$ and $\pi^\S_{c}$ exist,
\item $\pi^\S_b(\K)=\pi^\S_c(\K)$, and
\item letting $\K'=\pi^\S_b(\K)$, $\Sigma_{\K', \T^\frown \S^\frown\{b\}}=\Lambda_{\K', \U^\frown \S^\frown\{c\}}$.
\end{enumerate}
\end{enumerate}
$\myqedhere$
\end{definition}

The following is the comparison theorem we will prove in this subsection. The theorem uses concepts defined in \rdef{gamma(p, sigma) and b(p, sigma) for sts} and \rdef{gamma(p, sigma) and b(p, sigma)}.

\begin{theorem}[Diamond comparison]\label{diamond comparison} Suppose $(\P, \Sigma)$ and $(\Q, \Lambda)$ are two hod pairs such that $\Gamma(\P, \Sigma)=\Gamma(\Q, \Lambda)=_{def}\Gamma$, both $\Sigma$ and $\Lambda$ are $\Gamma$-fullness preserving ($\omega_1, \omega_1, \omega_1)$-strategies, $\P$ and $\Q$ are countable and are of limit type, and lower level comparison holds between $(\P, \Sigma)$ and $(\Q, \Lambda)$. Then there are $(\T, \R)\in I(\P, \Sigma)$ and $(\U, \R)\in I(\Q, \Lambda)$ such that either
\begin{enumerate}
\item $\P$ and $\Q$ are of lsa type and $\Sigma_{\R, \T}^{stc}=\Lambda^{stc}_{\R, \U}$ or
\item $\P$ and $\Q$ are not of lsa type and $\Sigma_{\R, \T}=\Lambda_{\R, \U}$.
\end{enumerate}
\end{theorem}

There are several other variations of the above theorem that works for sts hod pairs and also for a hod pair and an sts hod pair. We will state these theorems after the proof of \rthm{diamond comparison}. We prove \rthm{diamond comparison} by showing that the failure of its conclusion produces a bad sequence of length $\omega_1$. Towards showing this, we prove two useful lemmas. 

 We say that \textit{weak comparison}\index{weak comparison} holds between $(\P, \Sigma)$ and $(\Q, \Lambda)$ if there is $(\T, \U, \R, \S)$ such that 
\begin{enumerate}
\item $(\T, \R)\in I(\P, \Sigma)$,
\item $(\U, \S)\in I(\Q, \Lambda)$,
\item $\R^b=\S^b$ and $\Sigma_{\R^b, \T}=\Lambda_{\S^b, \U}$.
\end{enumerate}
Our first lemma says that lower level comparison implies that weak comparison holds.

\begin{lemma}\label{weak comparison} Suppose $(\P, \Sigma)$ and $(\Q, \Lambda)$ are two hod pairs such that $\Gamma(\P, \Sigma)=\Gamma(\Q, \Lambda)=_{def}\Gamma$\footnote{See \rdef{gamma(p, sigma) and b(p, sigma)}.}, both $\Sigma$ and $\Lambda$ are $\Gamma$-fullness preserving, $\P$ and $\Q$ are of limit type, and that lower level comparison holds between $(\P, \Sigma)$ and $(\Q, \Lambda)$. Then weak comparison holds between $(\P, \Sigma)$ and $(\Q, \Lambda)$.
\end{lemma}
\begin{proof} We inductively construct $(\P_i, \T_i: i<\omega)$ and $(\Q_i, \U_i: i<\omega)$ such that the following conditions hold.
\begin{enumerate}
\item $\P_0=\P$ and $\Q_0=\Q$.
\item Suppose $i=2n$. Then the following holds.
\begin{enumerate}
\item $\T_i$ is a stack on $\P_i$ based on $\P_i^b$ and according to $\Sigma_{\P_i, \oplus_{k<i}\T_k}$ with last model $\P_{i+1}$. 
\item $\U_i$ is a stack on $\Q_i$ according to $\Lambda_{\Q_i, \oplus_{k<i}\U_i}$ with last model $\Q_{i+1}$. 
\item Letting $\P'_i$ be the least non gentle layer of $\P_{i+1}$ such that $\pi^{\T_i}[\P_i^b]\subseteq \P_i'$, $\P_i'\insegeq_{hod}\Q_{i+1}^b$ and $\Lambda_{\P_i', \oplus_{k\leq i}\U_{k}}=\Sigma_{\P_i', \oplus_{k\leq i}\T_k}$.
\end{enumerate}
\item Suppose $i=2n+1$. Then the following holds. 
\begin{enumerate}
\item $\T_i$ is a stack on $\P_i$ according to $\Sigma_{\P_i, \oplus_{k<i}\T_k}$ with last model $\P_{i+1}$. 
\item $\U_i$ is a stack on $\Q_i$ based on $\Q_i^b$ and according to $\Lambda_{\Q_i, \oplus_{k<i}\U_i}$ with last model $\Q_{i+1}$. 
\item Letting $\Q_i'$ be the least non gentle layer of $\Q_{i+1}$ such that $\pi^{\U_i}[\Q_i^b]\subseteq \Q_i'$, $\Q_i'\insegeq_{hod}\P_{i+1}^b$ and $\Lambda_{\Q', \oplus_{k\leq i}\U_{k}}=\Sigma_{\Q', \oplus_{k\leq i}\T_k}$.
\end{enumerate}
\end{enumerate}

We show how to carry out the inductive step. Suppose we have constructed $(\P_i, \Q_i: i\leq 2n)$ and $(\T_i, \U_i: i< 2n)$. We now consider two cases.\\

\textbf{Case 1.} $\cf^{\P_{2n}}(\d^{\P^b_{2n}})$ is not a measurable cardinal in $\P_{2n}$. \\

Notice that in this case, we have that $\P_1=\Q_1$ and $\Sigma_{\P_1, \T_0}=\Lambda_{\Q_1, \U_0}$. Thus, weak comparison holds for $(\P, \Sigma)$ and $(\Q, \Lambda)$ provided we can take care of $n=0$ case. Notice also that in this case $\P_0=\P_0^b$.

Let $(\N_i: i<\omega)$ be a sequence of layers of $\P(=\P_0)$ such that 
\begin{itemize}
\item for all $i<\omega$, $\d^{\N_i}$ is a Woodin cardinal of $\P$, 
\item for all $i<\omega$, $\N_i\inseg_{hod}\N_{i+1}$ and
\item $\P|\d^\P=\cup_{i<\omega}\N_i$.
\end{itemize}
 By induction we construct a sequence $(\T^*_k, \W_k, \S_k, \R_k, \S^*_k, \R^*_k, \R^{**}_k: k<\omega)$ such that the following hold.
\begin{enumerate}
\item $(\S^*_0, \R^*_0)\in I(\Q, \Lambda_\Q)$, $\R^{**}_0\insegeq_{hod}\R^*_0$ and 
\begin{center}
$\Gamma(\N_0, \Sigma_{\N_0})=\Gamma(\R^{**}_0, \Lambda_{\R_0^{**}, \S^*_0})$.
\end{center}
Also, $(\T^*_0, \W_0)\in I(\P, \Sigma)$, $(\S_0, \R_0)\in I(\R^*_0, \Lambda_{\R_0^*, \S^*_0})$ and the following conditions hold:
\begin{enumerate}
\item $\T_0^*$ is based on $\N_0$ and $\S_0$ is based on $\R^{**}_0$. 
\item $\pi^{\T_0^*}(\N_0)=\pi^{\S_0}(\R_0^{**})$ and
\item letting $\K=\pi^{\T_0^*}(\N_0)$,
\begin{center}
$\Sigma_{\K, \T_0^*}= \Lambda_{\K, \S_0^{*\frown} \S_0}$.
\end{center}
\end{enumerate}
\item For $k+1<\omega$, 
$(\S^*_{k+1}, \R^*_{k+1})\in I(\R_k, \Lambda_{\R_k, \oplus_{m\leq k}(\S_m^{*\frown}\S_m)})$, $\R^{**}_{k+1}\inseg_{hod}\R^*_{k+1}$ and 
\begin{center}
$\Gamma(\N^*_{k+1}, \Sigma_{\N^*_{k+1}, \oplus_{m\leq k}\T^*_m})=\Gamma(\R^{**}_{k+1}, \Lambda_{\R_{k+1}^{**},  \oplus_{m\leq k}(\S_m^{*\frown}\S_m)})$.
\end{center}
where $\N_{k+1}^*=\pi^{\oplus_{m\leq k}\T_m^*} (\N_{k+1})$.
Also, 
\begin{center}
$(\T_{k+1}^*, \W_{k+1})\in I(\W_{k}, \Sigma_{\W_k, \oplus_{m\leq k}\T^*_m})$,\\
 $(\S_{k+1}, \R_{k+1})\in I(\R^*_{k+1}, \Lambda_{\R_{k+1}^*,\oplus_{m\leq k}(\S_m^{*\frown}\S_m))^\frown \S_{k+1}^*})$
 \end{center} 
 and the following conditions hold:
\begin{enumerate}
\item $\T_{k+1}^*$ is based on $\N^*_{k+1}$ and $\S_{k+1}$ is based on $\R^{**}_{k+1}$. 
\item $\pi^{\T_{k+1}^*}(\N_{k+1}^*)=\pi^{\S_{k+1}}(\R_{k+1}^{**})$ and
\item letting $\K=\pi^{\T_{k+1}^*}(\N_{k+1}^*)$,
\begin{center}
$\Sigma_{\K, \oplus_{m\leq k+1}\T_m^*}= \Lambda_{\K, \oplus_{m\leq k+1}(\S_m^{*\frown}\S_m)}$.
\end{center}
\end{enumerate}
\end{enumerate}
We then let $\T_{0}=\oplus_{k<\omega}\T^*_k$ and $\U_{0}=\oplus_{m<\omega}\S_k^{*\frown}\S$. Also, we let $\P_{1}$ be the last model of $\T_{0}$ and $\Q_{1}$ be the last model of $\U_{0}$. \\

\textbf{Case 2.} $\cf^{\P_{2n}}(\d^{\P^b_{2n}})$ is a measurable cardinal in $\P$. \\

The difference between this case and the previous case is that here we cannot start by fixing $(\N_i: i<\omega)$ as above. Here is the outline of the construction of $(\T_{2n}, \U_{2n}, \P_{2n+1}, \Q_{2n+1})$. 

Because $\Gamma(\P_{2n}, \Sigma_{\P_{2n}, \oplus_{i<2n}\T_i})=\Gamma(\Q_{2n}, \Lambda_{\Q_{2n}, \oplus_{i<2n}\U_i})$, we can find
\begin{center}
 $(\S_0, \R_0)\in I(\Q_{2n}, \Lambda_{\Q_{2n}, \oplus_{i<2n}\U_i})$
 \end{center}
 and $\R_0^*\inseg_{hod}\R_0$ such that letting $E\in \vec{E}^{\P_{2n}}$ be the extender of Mitchel order 0 with $\cp(E)=\cf^{\P_{2n}}(\d^{\P_{2n}})$, 
\begin{center}
$\Gamma(\P_{2n}^b, \Sigma_{\P_{2n}, (\oplus_{i<2n}\T_i)^\frown\{Ult(\P_{2n}, E)\} })=\Gamma(\R_0^*, \Lambda_{\R_0^*, (\oplus_{i<2n}\U_i)^\frown\{ \S_0\}})$
\end{center}
Appealing to low level comparison, we can find 
\begin{center}
$(\T^*_{2n}, \P_{2n+1})\in I(Ult(\P_{2n}, E), \Sigma_{\P_{2n}, (\oplus_{i<2n}\T_i)^\frown\{Ult(\P_{2n}, E)\} })$ and\\ $(\S_1, \R_1)\in  I(\R_0, \Lambda_{\R_0, (\oplus_{i<2n}\U_i)^\frown\S_0})$
\end{center}
such that 
\begin{enumerate}
\item $\T^*_{2n}$ is based on $\P_{2n}^b$, 
\item $\S_1$ is based on $\R_0^*$, 
\item $\pi^{\T^*_{2n}}(\P_{2n}^b)=\pi^{\S_1}(\R_0^*)=_{def}\K$, and 
\item $\Sigma_{\K, (\oplus_{i<2n}\T_i)^\frown\{E\}^\frown \T_{2n}^*}=\Lambda_{\K, (\oplus_{i<2n}\U_i)^\frown\S_0^\frown\S_1 }$
\end{enumerate}
Let then $\T_{2n}= \{E\}^\frown\T_{2n}^*$, $\U_{2n}=\S_0^\frown \S_1$ and $\Q_{2n+1}=\R_1$. 

The two cases above finish the construction of $(\T_{2n}, \U_{2n}, \P_{2n+1}, \Q_{2n+1})$. The construction of $(\T_{2n+1}, \U_{2n+1}, \P_{2n+2}, \Q_{2n+2})$ is very similar and we leave it to the reader.

Notice now that if $\T=\oplus_{i<\omega}\T_i$, $\U=\oplus_{i<\omega}\U_i$, $\R$ is the last model of $\T$ and $\S$ is the last model of $\U$ then $(\T, \R)$ and $(\U, \S)$ witness that weak comparison holds for $(\P, \Sigma)$ and $(\Q, \Lambda)$.
\end{proof}

\begin{lemma}\label{getting a bad block} Suppose $(\P, \Sigma)$ and $(\Q, \Lambda)$ are two hod pairs such that $\Gamma(\P, \Sigma)=\Gamma(\Q, \Lambda)=_{def}\Gamma$, both $\Sigma$ and $\Lambda$ have strong branch condensation and are strongly $\Gamma$-fullness preserving, both $\P$ and $\Q$ are of limit type and low level comparison holds for $(\P, \Sigma)$ and $(\Q, \Lambda)$.  Suppose further that
$\P^b=\Q^b$ and $\Sigma_{\P^b}=\Lambda_{\Q^b}$. Let $(\T, \R, \U, \S)$ be the trees of the extender comparison of $\P$ and $\Q$\footnote{Thus, $\T$ is on $\P$ with last model $\R$ and $\U$ is on $\Q$ with last model $\S$. See \rdef{extender comparison}.}. Suppose that either
\begin{enumerate}
\item $\R\not =\S$ or 
\item $\R=\S$ and $\Sigma_{\R, \T}\not =\Lambda_{\S, \U}$.
\end{enumerate}
Then there is a bad block on $((\P, \Sigma) , (\Q, \Lambda))$.
\end{lemma}
\begin{proof}
It follows from \rlem{disagreement implies low level disagreement} that we can find minimal low level disagreement $((\T^*, \P^*), (\U^*, \Q^*), \K)$  between $(\R, \Sigma_{\R, \T})$ and $(\S, \Lambda_{\S, \U})$. Let $E$ be the $\W$-un-dropping extender of $\T^\frown \T^*$ and $F$ be the $\W$-un-dropping extender of $\U^\frown \U^*$, and let $\T_0$ be the extension of $\T^\frown \T^*$ obtained by applying $E$ and $\U_0$ be the extension of $\U^\frown \U^*$ obtained by applying $F$. We then let $\P_1$ and $\Q_1$ be the last models of $\T_0$ and $\U_0$. 

Let $\X$ be a normal stack as in clause 2 of \rdef{low level comparison}. Let $b=\Sigma(\T_0^\frown\S)$, $c=\Lambda(\U_0^\frown\X)$, $\P_2=\M^{\T_1}_b$ and $\Q_2=\M^{\X}_c$. Set $\T_1=\X^\frown \{b\}$ and $\U_1=\X^\frown \{c\}$. We thus have that $\pi^{\T_1}$ and $\pi^{\U_1}$ exist, $\pi^{\T_1}(\K)=\pi^{\U_1}(\K)$ and 
\begin{center}
$\Sigma_{\pi^{\T_1}(\K), \T_0^\frown\T_1}=\Lambda_{\pi^{\U_1}(\K), \U_0^\frown\U_1}$
\end{center}

Next (appealing to \rlem{weak comparison}) let $(\T_2, \P_3)$ and $(\U_2, \Q_3)$ witness that the weak comparison holds for 
\begin{center}
$(\P_2, \Sigma_{\P_2, \T_0^\frown\T_1})$, and$ (\Q_2, \Lambda_{\Q_2, \U_0^\frown\U_1})$. 
\end{center}

Next let $\P_0=\P$, $\Q_0=\Q$, $\Sigma_0=\Sigma$, $\Lambda_0=\Lambda$, and for $i\in \{1, 2, 3\}$ let $\Sigma_i=\Sigma_{\P_i, \oplus_{k<i} \T_k}$ and $\Lambda_i=\Lambda_{\Q_i, \oplus_{k<i} \U_k}$. 
It is then easy to see that 
\begin{center}
$(( (\P_i, \Q_i, \Sigma_i, \Lambda_i): i<4), ( \T_i , \U_i: i<3), (b, c) )$
\end{center}
is a bad block on $((\P, \Sigma) , (\Q, \Lambda))$. 
\end{proof}

\textbf{The proof of  \rthm{diamond comparison}}\\\\
 Suppose that the conclusion of \rthm{diamond comparison} fails. This means that\\\\
(1) whenever $(\T, \R)\in I(\P, \Sigma)$ and $(\U, \R)\in I(\Q, \Lambda)$,
\begin{enumerate}
\item if $\P$ and $\Q$ are of lsa type then $\Sigma_{\R, \T}^{stc}\not=\Lambda^{stc}_{\R, \U}$ or
\item if $\P$ and $\Q$ are not of lsa type then $\Sigma_{\R, \T}\not=\Lambda_{\R, \U}$.
\end{enumerate}

It follows from \rlem{weak comparison} that, without loss of generality, we can assume that $\P^b=\Q^b$ and  $\Sigma_{\P^b}=\Lambda_{\Q^b}$. We now by induction construct a bad sequence $(B_\a: \a<\omega_1)$ on $((\P, \Sigma), (\Q, \Lambda))$. 

It follows from \rlem{getting a bad block} that there is a bad block on $((\P, \Sigma), (\Q, \Lambda))$. Let $B_0$ be any bad block on $((\P, \Sigma), (\Q, \Lambda))$. Suppose next that we have constructed $(B_\b :\b<\l)$ for $\l$ a limit. Let $\P_\l$ and $\Q_\l$ be the direct limit of respectively $(\P_\b:\b<\l)$ and $(\Q_\b:\b<\l)$ under the corresponding iteration embeddings. Then letting $\Sigma_{\l, 0}$ and $\Lambda_{\l, 0}$ be the corresponding tails of $\Sigma$ and $\Lambda$, we have that $(\P_\l, \Sigma_\l)$ and $(\Q_\l, \Lambda_\l)$ satisfy the hypothesis of \rlem{getting a bad block}. Let then $B_\l$ be a bad block on $((\P_\l, \Sigma_{\l}), (\Q_\l, \Lambda_{\l}))$. 

Next suppose that we have constructed $(B_\b :\b<\l+1)$. Let $\P_{\l+1}=\P_{\l, 3}$, $\Q_{\l+1}=\Q_{\l, 3}$ and let $\T$ and $\U$ be the stacks respectively on the top of $(B_\b :\b<\l+1)$ and in the bottom of $(B_\b :\b<\l+1)$. We then again can find, using \rlem{getting a bad block}, a bad block $B_{\l+1}$ on $((\P_{\l+1}, \Sigma_{\P_{\l+1}, \T}), (\Q_{\l+1}, \Lambda_{\Q_{\l+1}, \U}))$. It then follows that the resulting sequence $( B_\b: \b<\omega_1)$ is a bad sequence on $((\P, \Sigma), (\Q, \Lambda))$. This is a contradiction to Lemma \ref{no bad sequences}, and this contradiction completes the proof of \rthm{diamond comparison}.

\section{Some concluding remarks}\label{sec: concluding remarks}

The proof of \rthm{diamond comparison} can be used to show that fullness preserving strategies that have strong branch condensation become commuting on a tail. We end this section by a an outline of this useful fact.

\begin{proposition}\label{commuting from strong branch condensation} Suppose $(\P, \Sigma)$ is an sts hod pair and $\Gamma$ is a projectively closed pointclass. Suppose that $\Sigma$ has strong branch condensation and is $\Gamma$-fullness preserving. Then for some $(\T, \Q)\in I^{ope}(\P, \Sigma)$, $\Sigma_\Q$ is commuting\footnote{See \rdef{positional and commuting for sts pairs}.}.
\end{proposition}
\begin{proof} Towards a contradiction assume not. Then we can find a sequence
\begin{center}
 $c=(\P_\a, \T_\a, \Q_\a, \U_\a, \R_\a, k'_\a, k_\a: \a\leq \omega_1)$
 \end{center}
  such that the following conditions hold:
\begin{enumerate}
\item For each $\a<\omega_1$, $\P_{\a+1}$ is the result of comparing the pairs $(\Q_\a, \Sigma_{\Q_\a})$ and $(\R_\a, \Sigma_{\R_\a})$\footnote{The comparison is possible because of \rcor{self comparison}.}
\item For each $\a<\omega_1$, $(\P_\a, (\T_\a, \Q_\a), (\U_\a, \R_\a), k'_\a, k_\a)$ witnesses that $\Sigma_{\P_\a}$ is not commuting.
\item For each limit ordinal $\a\leq \omega_1$\footnote{The rest of the objects are undefined for $\a=\omega_1$.}, $\P_\a$ is the direct limit of $(\P_\b, \pi_{\b, \gg}^t: \b<\gg<\a)$\footnote{Notice that in \rdef{positional and commuting for sts pairs} we can assume that $\pi^{\T}$ is defined, possibly by applying undropping extenders. This is because  commuting for sts hod pairs is a principle about the bottom parts not the entire model.} where $\pi_{\b, \gg}^t:\P_\b\rightarrow \P_\gg$ is the embedding given by concatenating the $\P_\b$-to-$\Q_\b$-to-$\P_{\b+1}$ stacks. 
\end{enumerate}
It follows from \rprop{positional1} that in clause 3 above we could define $\P_\a$ as the direct limit of $(\P_\b, \pi_{\b, \gg}^b: \b<\gg<\a)$ where $\pi_{\b, \gg}^b:\P_\b\rightarrow \P_\gg$ is the embedding given by concatenating the $\P_\b$-to-$\R_\b$-to-$\P_{\b+1}$ stacks.

Suppose now that $\sigma: H\rightarrow H_{\omega_2}$ is such that $H$ is countable and transitive, and $c\in \rge(\sigma)$. Let $\k=\omega_1^H$. It follows that 
\begin{center}
$\pi^t_{\k, \omega_1}=\sigma\rest \P_\k=\pi_{\k, \omega_1}^b$.
\end{center}
Let now $j:\Q_\k\rightarrow \P_{\omega_1}$ and $i:\R_\k\rightarrow \P_{\omega_1}$ be the two iteration embeddings. It follows from strong branch condensation and in particular from \rprop{positional1} that letting $\d=\d^{\R_\k^b}$,\\\\
(1) $i\rest \R_\k|\d=j\rest \Q_\k|\d$\footnote{Notice that $k_\k\rest \d=id$.}.\\\\
Hence, we have that $\pi^{\U_\k}\rest \P_\k|\d^{\P_\k^b}=\pi^{\T_\k}\rest \P_\k|\d^{\P_\k^b}$. It remains to show that for $A\in \powerset(\d^{\P_\k^b})\cap \P_\k$, $\pi^{\U_\k}(A)=k_\k(A)$. But we have that \\\\
(2) $i(\pi^{\U_k}(A))=j(\pi^{\T_\k}(A))$ and $k_\k(A)=\pi^{\T_\k}(A)\cap \d$.\\\\
It then follows from (1) and (2) that $\pi^{\U_\k}(A)=k_\k(A)$. 
\end{proof}

The following proposition implies that in many situations we can construct authenticating iterations as described in \rdef{authentic lsp}. We will use it in the proof of \rthm{main theorem on gen int}.

\begin{proposition}\label{authenticated iteration construction} Suppose $(\P, \Sigma)$ is an sts pair and $\Gamma$ is a projectively closed pointclass. Suppose $\Sigma$ is 
\begin{itemize}
\item strongly $\Gamma$-fullness preserving,
\item has strong branch condensation and
\item is commuting\footnote{See \rdef{positional and commuting for sts pairs}.}.
\end{itemize}
Suppose $(\T, \Q)\in I^{ope}(\P, \Sigma)$, $(\U, \R)\in I^{ope}(\P, \Sigma)$ and $(\W, \S)\in I^{ope}(\R, \Sigma_{\R})$ are such that for some $\d<\d^{\Q^b}$ the following conditions hold:
\begin{enumerate}
\item $\Q\models ``\d$ is a Woodin cardinal",
\item $\W$ is a normal stack, and
\item $\S|\d=\Q|\d$.
\end{enumerate}
Let 
\begin{itemize}
\item $\K\insegeq_{hod}\Q$ be such that $\d^\K=\d$,
\item $\a<\lh(\W)$ be the least such that $\K^-\insegeq \M_\a^\W$,
\item $\b< \lh(\W)$ be such that $\m(\W_{<\b})=\K|\d$,
\item $w$ is the window of $\Q$ such that $\d^w=\d$\footnote{See \rnot{l p}.}, and
\item $b=[0, \b)_\W$.
\end{itemize}
Then $s(\T, w)\subseteq \pi^{\W_{[\a, \b)}}_b$\footnote{See \rdef{canonical singularizing sequences}.}
\end{proposition}

The \rprop{authenticated iteration construction} can be proven by simply comparing $(\S, \Sigma_\S)$ and $(\Q, \Sigma_\Q)$ and then using commuting and \rcor{self comparison}.

\chapter{Hod mice revisited}

In this section we generalize the result of \cite[Chapter 3]{ATHM} to our current context. As in \cite{ATHM}, these results lead towards showing that given a hod pair $(\P, \Sigma)$, $\Gamma(\P, \Sigma)$ is an $OD$-full pointclass (see Definition 3.16 of \cite{ATHM}). 

Recall the effect of \rprop{positional}; if $(\P, \Sigma)$ is a hod pair such that $\Sigma$ has strong branch condensation and if $\Q\in pI(\P, \Sigma)$, then the strategy of $\Q$ induced by $\Sigma$ is independent of the particular iteration producing $\Q$. In \rsec{sec positional and commuting}, this strategy was denoted by $\Sigma_\Q$. In this chapter, whenever the strategy of a hod mouse has a strong branch condensation, we will make use of the aforementioned notation without giving any further explanation. 

\section{The uniqueness of the internal strategy}

The first theorem, \rthm{uniqueness of internal strategies}, is just a direct generalization of \cite[Theorem 3.3]{ATHM}. It says that the internal strategies are unique. First we prove a useful lemma.

\begin{lemma}\label{bigger cofinality} Suppose $\P$ is a hod premouse, $\Q\insegeq_{hod}\P$, $\U\in \P$ is a stack on $\Q$ with last model $\R$ such that $\U$ has a one point extension\footnote{See \rdef{one point extension}.}, and $\R'\insegeq_{hod}\R$ is such that $\R\models ``\d^{\R'}$ is a Woodin cardinal". Suppose further that if $\pi^\U$ is undefined then letting $E$ be the $\R'$-un-dropping extender of $\U$, $Ult(\P, E)$ is well-founded. Then $\cf^\P(\d^{\R'})>\omega$.
\end{lemma}
\begin{proof}
Towards a contradiction, assume not. We give the proof assuming that $\pi^\U$ is defined. If not, then one could instead work with $Ult(\P, E)$ instead of $\P$ and $\pi_E$ instead of $\pi^\U$, where $E$ is the $\R'$-un-dropping extender of $\U$. 

Notice that it cannot be the case that $\R'\in \rge(\pi^\U)$ as $\pi^\U$ is continuous at the Woodin cardinals of $\P$. Therefore, by minimizing $\Q$, we can assume that $\Q$ is of limit type. We now apply \rlem{the canonical singularizing sequence exists} to $(\U, w)$ where $w$ is the window of $\R$ such that $\d^w=\d^{\R'}$. Let $\eta=\eta^w$. We thus have that there is a sequence $(h_i: i<\omega)\in \Q^b$  and a sequence $(a_i: i<\omega) \in (\eta^{<\omega})^{\omega}$ such that 
\begin{center}
$\d^{\R'}=\sup\{ \pi^\U(h_i)(a_i): i<\omega\}$.
\end{center}
Notice now that $(\pi^\U(h_i): i<\omega)\in \R$. Therefore, 
\begin{center}
$\d^{\R'}=\sup\{ \pi^\U(h_i)(a): a\in [\eta]^{<\omega} \wedge \pi^\U(h_i)(a)<\d^{\R'}\}$.
\end{center}
It then follows that $\R\models \cf(\d^{\R'})\leq \eta$, which is a contradiction as $\d^{\R'}$ is a Woodin cardinal of $\R$.
\end{proof}

\begin{theorem}[Uniqueness of internal strategies]\label{uniqueness of internal strategies} Suppose $\P$ is a hod premouse such that $\P\models \sf{ZFC-Powerset}$, $\d^\P$ is a regular cardinal of $\P$ and  $\W\inseg_{hod}\P$ is such that $\P\models ``\Sigma^\P_\W$ is a $((\d^\P)^+, (\d^\P)^+)$-strategy"\footnote{If $(\d^\P)$ is the largest cardinal then we assume that $((\d^\P)^+)^\P=\ord(\P)$.}. Then 
$\P\models$ ``$\W$ has a unique iteration strategy ".
\end{theorem}
\begin{proof} Working in $\P$, suppose $\Lambda\not =\Sigma_\W^\P$ is another iteration strategy for $\W$. Let $\Sigma=\Sigma_\W^\P$. Notice that \rlem{bigger cofinality} implies that if 
\begin{itemize}
\item $\U$ is a stack on $\W$ according to both $\Lambda$ and $\Sigma$,
\item $\lh(\U)$ is a limit ordinal, and
\item $b=\Sigma(\U)$ and $c=\Lambda(\U)$
\end{itemize}
then \\\\
(*) either\\\\
(A) both $\Q(b, \U)$ and $\Q(c, \U)$ exist, or\\
(B) $b=c$.\\\\
This is because if $b\not =c$ then $\cf^\P(\d(\U))=\omega$ and hence, we have that
\begin{enumerate}
\item either $\pi^\U_b$ is undefined or $\d(\U)$ is not a Woodin cardinal of $\M^\U_b$, and
\item either $\pi^\U_c$ is undefined or $\d(\U)$ is not a Woodin cardinal of $\M^\U_c$.
\end{enumerate}
The above clauses imply that $\d(\U)$ is not a Woodin cardinal neither in $\M^\U_b$ nor in $\M^\U_c$. Therefore, both $\Q(b, \U)$  and $\Q(c, \U)$ exist. 

 It then follows from the proof of \rlem{disagreement implies low level disagreement}\footnote{The use of $\Gamma$-fullness preservation can be substituted by (*).} that we can find a minimal low-level disagreement  $(\T_1, \W_1, \T_2, \W_2, \Q)$ between $(\W, \Sigma)$ and $(\W, \Lambda)$. Moreover, we can assume that $\lh(\T_1)<\d^\P$ and $\lh(\T_2)<\d^\P$\footnote{If not, then we can reflect below $\d^\P$. Recall that $\W\inseg_{hod}\P$, so the desired Skolem hull of $\P$ can be required to contain $\W$.}. Let $\S\in \P$ be a stack on $\Q$ according to both $\Sigma_{\Q, \T_1}$ and $\Lambda_{\Q, \T_2}$ and such that $\Sigma_{\Q, \T_1}(\S)\not =\Lambda_{\Q, \T_2}(\S)$.
 It then follows from (*) that letting $b=\Sigma_{\Q, \T_1}(\S)$ and $c=\Sigma_{\Q, \T_2}(\S)$, both $\Q(b, \S)$ and $\Q(c, \S)$ exist. However, since $\Sigma_{\Q^-, \T_1}=\Lambda_{\Q^-, \T_2}$ and also both $\Q(b, \S)$ and $\Q(c, \S)$ are $\d^\P+1$-iterable in $\P$, we have that $\Q(b, \S)=\Q(c, \S)$
%
%
%
\end{proof}

%

\section{Generic interpretability}\label{generic interpretability sec}

We now move to generic interpretability. We start by recalling and generalizing the definition of a pre-hod pair (see \cite[Definition 3.7]{ATHM}). 

\begin{definition}[Prehod pair]\index{prehod pair}\label{prehod pair}
$(\P, \Sigma)$ is a prehod pair if
\begin{enumerate}
\item $\P$ is a countable hod premouse of successor type,
\item If $\P^-$ is not of limit type then $\Sigma$ is an $(\omega_1, \omega_1)$-strategy for $\P$ acting on stacks based on $\P^-$.
\item If $\P^-$ is of limit type then $\Sigma$ is an $(\omega_1, \omega_1, \omega_1)$-strategy for $\P$ acting on stacks based on $\P^-$.
\item If $i:\P\rightarrow \Q$ comes from an
iteration according to $\Sigma$, $\Sigma_{\Q^-}^\Q=\Sigma_{\Q^-}\rest \Q$\footnote{Thus, $\P$ is a $\Sigma$-mouse over $\P^-$.},
\item For any $\P$-cardinal $\eta\in (\d^{\P^-}, \d^\P)$, considering $\P|\eta$ as a $\Sigma$-mouse over $\P^-$, there is an $\omega_1$-strategy $\Lambda$ for $\P|\eta$\footnote{Thus, $\Lambda$ acts on stacks above $\d^{\P^-}$.}.
\end{enumerate}
$\myqedhere$
\end{definition}

Notice that there must be a unique strategy $\Lambda$ as in clause 5 of \rdef{prehod pair}.\footnote{$\Lambda$ is the $\Q$-structure guided strategy.} Also, recall the definition of Generic Interpretability, \cite[Definition 3.8]{ATHM}. In our current context it takes the following form.

\begin{definition}[Generic Interpretability]\label{generic interpretability} Suppose 
$(\P, \Sigma)$ is a pre-hod pair, a meek hod pair of limit type or an sts hod pair. We say \textit{generic interpretability} holds for $(\P, \Sigma)$ if there is a function $F$ such that
\begin{enumerate}
\item $F$ is definable over $\P$ with no parameters,
\item $dom(F)$ consists of pairs $(\Q, \k)$ such that $\Q\in Y^\P$, $\Q\insegeq \P|\d^\P$ and $\k\in (\d^\Q, \d^\P)$ is a $\P$-cardinal,
\item for $(\Q, \k)\in dom(F)$, $F(\Q, \k)=( \dot{T}, \dot{S})$ such that ,
    \begin{enumerate}
    \item $\dot{T}, \dot{S}\in \P^{Coll(\omega, \ord\Q))}$,
     \item $\P\models ``\forces_{Coll(\omega, \ord\Q))} \dot{T}$ and $\dot{S}$ are $\k$-complementing",
     \item for any $\nu\in( \ord\Q),\k)$ and any $\P$-generic $g\subseteq Coll(\omega, \nu)$,
\begin{center}
$\P[g]\models ``p[\dot{T}_{g}]$ is an $(\omega_1, \omega_1, \omega_1)$-iteration strategy for $\Q$ which extends $\Sigma_{\Q}^\P$"
\end{center}
and
\begin{center}
$(p[\dot{T}_g])^{\P[g]}=\Sigma_{\Q}\rest HC^{\P[g]}$.
\end{center}
 \end{enumerate}
\end{enumerate}
$\myqedhere$
\end{definition}

The proof that the generic interpretability holds is just like the proof of \cite[Theorem 3.10]{ATHM} using \rthm{universality of background construction} and \rthm{uniqueness of internal strategies} instead of  \cite[Lemma 2.15]{ATHM}  and \cite[Theorem 3.3]{ATHM}.  First the proof of \cite[Lemma 3.9]{ATHM} can be used with no changes to establish the following useful lemma. 

\begin{lemma}\label{strengthening of clause 6}
Suppose $(\P, \Sigma)$ is a prehod pair. Let $\k\in (\d^{\P^-}, \d^\P)$ be a $\P$-cardinal and $\Lambda^*$ be the iteration strategy of $\P|\k$ as in 5 of \rdef{prehod pair}. Let $\Lambda$ be the fragment of $\Lambda^*$ that acts on non-dropping stacks.
Let $g\subseteq Coll(\omega, \k)$ be $\P$-generic. Then $\P[g]$ locally Suslin, co-Suslin captures
${\sf{Code}}(\Lambda)$ and its complement at any cardinal of $\P$ greater than $\k$\footnote{Recall that this means that for every $\P$-cardinal $\nu>\k$, there are $\nu$-complementing trees $U, V\in \P[g]$ such that for any $<\nu$-generic $h$, ${\sf{Code}}(\Lambda)\cap P[g][h]=(p[U])^{\P[g][h]} = (\mathbb{R}^{\P[g][h]}-p[V])^{\P[g][h]}.$}.
\end{lemma}

Fix now a prehod pair $(\P, \Sigma)$ and let $\Q\in Y^\P$.  Let $\k<\d^\P$ be a $\P$-cardinal such that 
\begin{itemize}
\item $\kappa> \ord(\Q)$ and
\item in the case $\Q$ is of limit type, $\P$ has no extenders with critical point $\d^{\Q^b}$ and index greater than $\k$.
\end{itemize}
Let $\vec{G}=\{ E\in \vec{E}^{\P|\d^\P}: \nu(E)$ is an inaccessible cardinal of $\P$ and $\cp(E)>\k\}$. Notice that $(\P, \d^\P, \Sigma, \vec{G})$ is a self-capturing background triple. 
Let 
\begin{center}${\sf{hpc}}^+=(\M_\gg , \N_\gg, Y_\gg, \Phi^+_\gg, F^+_\gg, F_\gg, b_\gg : \gg\leq \delta)$
\end{center}
be the output of ${\sf{hpc}}$ of $(\P, \d^\P, \Sigma, \vec{G})$\footnote{See \rdef{gamma-hod pair construction*}. The aforementioned definition requires a pointclass $\Gamma$ but one can simply ignore all the clauses of \rdef{gamma-hod pair construction*} that mention $\Gamma$.}. 

Here we abuse the notation and write $\Phi_\b$ both for the strategy of $\M_\b$ that is internal to $\P$ and also for the external strategy. It follows from \rthm{universality of background construction} , \rlem{uniqueness of internal strategies} and \rlem{strengthening of clause 6} that for some $\b$, $(\M_\b, \Phi^+_{\b})$ is a tail of $(\Q, \Sigma_\Q)$. We then set 
\begin{center}
$\N_{\k, \Q}^\P=\M_\b$ and $\Lambda_{\k, \Q}=\Phi^+_{\b}$. 
\end{center}
In what follows, we will omit superscript $\P$, but ask the reader to keep in mind that certain notions depend on $\P$. 
Also let $\pi_{\k, \Q}:\Q\rightarrow \N_{\k, \Q}$ be the iteration embedding according to $\Sigma_{\Q}$ and let $\T_{\k, \Q}$ be the normal stack on $\Q$ with last model $\N_{\k, \Q}$. The following is a consequence of \rlem{strengthening of clause 6}, hull condensation of $\Sigma$ and the proof of \rthm{universality of background construction}.

\begin{corollary}\label{towards getting generic interpretability} Suppose $(\P, \Sigma)$ is a pre-hod pair such that for some projectively closed pointclass $\Gamma$, $\Sigma$ has branch condensation and is $\Gamma$-fullness preserving. Suppose $\Q\in Y^\P$ and $\k>\ord(\Q)$ are such that
\begin{itemize}
\item $\kappa> \ord(\Q)$ and
\item in the case $\Q$ is of limit type, $\P$ has no extenders with critical point $\d^{\Q^b}$ and index greater than $\k$.
\end{itemize}
Let $\eta\in (\ord(\N_{\k, \Q}), \d^\P)$ and $n<\omega$. Then there are names $(\dot{T}, \dot{S})\in \P^{Coll(\omega, \eta)}$ such that  
 \begin{enumerate}

     \item $\P\models ``\forces_{Coll(\omega, \eta)} \dot{T}$ and $\dot{S}$ are $(\d^\P)^{+n}$-complementing",
     \item for any $\l\in (\eta, ((\d^\P)^{+n})^\P)$ and any $\P$-generic $g\subseteq Coll(\omega, \l)$,
     \begin{center}
$\P[g]\models ``p[\dot{T}_{g}]$ is an $(\omega_1, \omega_1)$-iteration strategy for $\N_{\k, \Q}$"
\end{center}
and letting $\Phi$ be the $\pi^\P_{\k, \Q}$-pullback of the strategy given by $(p[\dot{T}_g])^{\P[g]}$ then 
\begin{center}
$\Phi=\Sigma_{\Q}\rest HC^{\P[g]}$.
\end{center}
\end{enumerate}
\end{corollary} 

Our generic interpretability result can now be proved using the tree production lemma (\cite[Theorem 3.3.15 ]{Stationarytower}) and \rcor{towards getting generic interpretability}. We leave the details to the reader.

\begin{theorem}[The generic interpretability]\label{generic interpretability holds} Suppose $(\P, \Sigma)$ is a prehod pair or is a non-gentle hod pair of limit type or is an sts hod pair. Also, suppose that for some projectively closed pointclass $\Gamma$, $\Sigma$ is $\Gamma$-fullness preserving. Assume that for every $\Q\in Y^\P$, $\Sigma_{\Q}$ has strong branch condensation. Then generic interpretability holds for $(\P, \Sigma)$.
\end{theorem}

Next, we present our result on internal fullness preservation. The proof follows the same line of thought as the proof of \cite[Theorem 3.12 ]{ATHM}. Below $\S^*(\R)$ is the $*$-transform of $\S$ into a hybrid mouse over $\R$ and it is defined  when $\ord(\R)$ is a cutpoint of $\S$ (see \cite[Remark 12.7]{DMATM} and \cite{Selfiter}).

\begin{definition}\label{internal fullness preservation} Suppose $\P$ is a hod premouse and $\Q\in Y^\P$. We say $\Lambda=\Sigma_{\Q}^\P$ is \textbf{internally fullness preserving} if the following holds for $(\T, \R)\in I(\Q, \Lambda)$\footnote{Thus, $(\T, \R)\in \P$.} such that $\P\models ``|{\T}|^+$ exists".
\begin{enumerate}
\item For all limit type $\S\in Y^\R$, if $\M\in \P$ is a sound $\max(\d^\P+1, (|{\T}|^+)^\P)$-iterable $\Lambda_{\S|\d^{\S^b}, \T}$-mouse over $\S|\d^{\S^b}$ then $\M\insegeq \S$.
\item Suppose $\W\inseg_{hod}\S$ is of lsa type and $\W=((\W|\d^\W)^\#)^\S$. Suppose $\M\in \P$ is a sound $\max(\d^\P+1, (\card{{\T}}^+)^\P)$-iterable $\Lambda_{\W, \T}$-sts mouse over $\W$. Then $\M\insegeq \S$.
\item Suppose $\R_1\inseg_{hod} \R$ is of successor type  and $\eta\in (\ord(\R_1^-), \d^{\R_1}]$ is a cutpoint cardinal of $\R$. Suppose $\M\in \P$ is a sound $\max(\d^\P+1, (|{\T}|^+)^\P)$-iterable $\Lambda_{\R^-_1, \T}$-mouse over $\R|\eta$. Then $\M\insegeq (\R|(\eta^+)^\R)^*(\R|\eta)$.
\end{enumerate}
$\myqedhere$
\end{definition}

\begin{theorem}[Internal fullness preservation]\label{internal fullness preservation1} Suppose $\P$ is a hod premouse and $\Q\in Y^\P$ is such that $(\ord(\Q)^+)^\P$ exists. Then $\Sigma_\Q^\P$ is internally fullness preserving. 
\end{theorem}

\section{The derived models of hod mice}\label{sec: derived models of hod mice}

In this section, we state, without a proof, a version of \cite[Theorem 3.19]{ATHM}. Suppose $(\P, \Sigma)$ is an allowable pair\footnote{See \rdef{allowable pair}.} such that $\Sigma$ has strong branch condensation and is fullness preserving\footnote{See \rdef{gamma fullness preservation}.}. Suppose $\Q\insegeq_{hod} \P$ is such that $\Q$ is meek and is of limit type. Thus, $\d^\Q$ is a limit of Woodin cardinals of $\P$. Suppose further that $\cf^\P(\d^\Q)$ is not a measurable cardinal in $\P$. We then let $D^*(\P, \Sigma, \Q)$ be the set of all $A\subseteq \mathbb{R}$ such that for some strong cutpoint $\tau<\d^\Q$ of $\Q$ and $g\subseteq Coll(\omega, \tau)$-generic over $\P$ there are  trees $T,U\in \P[g]$ such that 
\begin{enumerate}
\item $\P[g]\models ``(T, U)$ is $\d^\Q$-complementing" and
\item $x\in A$ if and only if there is $(\S, \R)\in I(\Q, \Sigma_{\Q})$ and a Woodin cardinal $\d$ of $\R$ such that
\begin{itemize}
 \item $\pi^\S$ is above $\tau$,
 \item $x$ is generic for the extender algebra of $\R[g]$ at $\d$ and
 \item $\R[g, x]\models x\in p[\pi^{\S}(T)]$.
\end{itemize}
\end{enumerate}
It follows from \rcor{comparison holds} and \rthm{positional} that for $x\in \bR$, the right hand side of the above equivalence is independent of the choice of $(\S, \R)$.

We let $D(\P, \Sigma, \Q)$\index{$D(\P, \Sigma, \Q)$} be the derived model of $\Q$ as computed by $\Sigma_\Q$, i.e., for $A\subseteq \mathbb{R}$, $A\in D(\P, \Sigma, \Q)$ if there is $(\S, \R)\in I(\P, \Sigma)$ such that $\S$ is based on $\Q$ and $A\in D^*(\R, \Sigma_{\R}, \pi^{\S}(\Q))$. 

Next recall \cite[Definition 3.18]{ATHM}. Essentially a pointclass is completely mouse-full if the next model of determinacy has the same mice relative to common iteration strategies. We introduce this notion more carefully. 

Given a set of reals $A\subseteq \bR$, we let $W_A=\{ B\subseteq \bR: w(B)< w(A)\}$. Next following Definition 3.13 of \cite{ATHM}, we say $A\subseteq\mathbb{R}$ is a new set\index{new set} if
\begin{enumerate}
  \item $L(A, \mathbb{R})\models {\sf{AD^+}}$,
  \item $\powerset(\mathbb{R}) \cap L(W_A, \mathbb{R})=W_A$,
  \item $\Theta^{L(W_A, \mathbb{R})}$ is a Suslin cardinal of $L(A, \mathbb{R})$.
\end{enumerate}

The following is \cite[Definition 3.17]{ATHM}. 

\begin{definition}\label{completely mouse full} Given a pointclass $\Gamma$, we say $\Gamma$ is \textbf{completely mouse full}\index{completely mouse full} if either $\Gamma=\powerset(\mathbb{R})$ or there is a new set $A$ such that
  \begin{enumerate}
  \item $\Gamma=W_A$,
  \item if $(\P, \Sigma)$ is allowable such that ${\sf{Code}}(\Sigma)\in \Gamma$ and $L(A, \mathbb{R})\models ``\Sigma$ has strong branch condensation and is $\Gamma$-fullness preserving" then for every $a\in HC$, 
  \begin{center}
  $Lp^{\Gamma, \Sigma}(a)=(Lp^\Sigma(a))^{L(A, \bR)}$. 
  \end{center}
  \end{enumerate}
  $\myqedhere$
  \end{definition}

 Given two pointclasses $\Gamma_1$ and $\Gamma_2$, we write $\Gamma_1\insegeq_{mouse} \Gamma_2$ if $\Gamma_1\subseteq \Gamma_2$ and $\Gamma_2$ has the same mice as $\Gamma_1$ relative to common iteration strategies. More precisely, if $(\P, \Sigma)\in \Gamma_1$ is an allowable pair such that $L(\Gamma_2, \bR)\models ``\Sigma$ has strong branch condensation and is $\Gamma_1$-fullness preserving" then for any $a\in HC$,
\begin{center}
$Lp^{\Gamma_1, \Sigma}(a)=Lp^{\Gamma_2, \Sigma}(a)$.
\end{center}
 Finally, following \cite[Definition 3.18]{ATHM}, 
 \begin{definition}\label{mouse full pointclass}\label{mouse full}
 $\Gamma$ is \textbf{mouse full}\index{mouse full} if either it is completely mouse full or is a union of completely mouse full pointclasses $(\Gamma_\a: \a<\Omega^\Gamma)$ such that for all $\a$, $\Gamma_\a\insegeq_{mouse}\Gamma_{\a+1}$ and for all limit $\a$, $\Gamma_\a=\bigcup_{\b<\a}\Gamma_\b$. $\myqedhere$
 \end{definition} We can now state our generalization of \cite[Theorem 3.19]{ATHM}.

\begin{theorem}\label{the derived model theorem for hod pairs} Suppose $(\P, \Sigma)$ is an allowable pair and $\Gamma$ is a pointclass closed under continuous preimages.\footnote{We define the Solovay sequence $(\theta^\Gamma_\alpha : \alpha\leq \Omega)$ relative to $\Gamma$ as the Solovay sequence defined in the model $L(\Gamma,\mathbb{R})$ if $\Gamma$ is constructibly closed (i.e., $\powerset(\bR)\cap L(\Gamma, \bR)=\Gamma$). We can aslo make sense of the Solovay sequence relative to $\Gamma$ in the case $\Gamma$ is a limit of constructibly closed pointclasses; here for $A\in \Gamma$, we say a set $B$ is $OD^\Gamma(A)$ if $B$ is $OD(A)^{L(\Lambda,\mathbb{R})}$ for some constructibly closed $\Lambda \lhd \Gamma$. From here on, when we talk about the Solovay sequence relative to a pointclass $\Gamma$, $\Gamma$ is assumed to have one of the two properties above. Notice that if $\Gamma$ is a constructibly closed pointclass which is a union of constructibly closed pointclasses strictly contained in it, then the two ways of computing the Solovay sequence relative to $\Gamma$ are equivalent.} Suppose further that $\P$ is non-gentle and of limit type, and that $\Sigma$ has strong branch condensation and is $\Gamma$-fullness preserving. Then
\begin{enumerate}
\item  $\Gamma(\P, \Sigma)=\bigcup_{\Q\in pI(\P, \Sigma), \Q'\insegeq \Q}D(\Q, \Sigma_{\Q}, \Q')$.
\item For any $\Q\in pI(\P, \Sigma)$, if $\Q'\inseg_{hod}\Q$ is non-gentle and is of limit type then $D(\Q, \Sigma_{\Q}, \Q')$ is completely mouse full.
\item For any $\Q\in pI(\P, \Sigma)$, if $\Q'\inseg_{hod}\Q''\insegeq_{hod}\Q$ are such that
\begin{itemize}
\item $\Q'$ and $\Q''$ are non-gentle and are of limit type,
\item $\Q''\insegeq_{hod}\Q$ is the least non-gentle layer of $\Q$ that has $\omega$ more Woodin cardinals than $\Q'$,
\end{itemize}
then letting $\Gamma'=D(\Q, \Sigma_{\Q}, \Q'')$, if $\xi$ is such that $\theta_{{\sf{Code}}\Sigma_{\Q'})}^{\Gamma'}=\theta_\xi^{\Gamma'}$ then for every $n$, letting $\Q'_n\insegeq_{hod} \Q''$ be the layer of $\Q''$ that has exactly $n$ Woodin cardinals above $\ord(\Q')$,
    \begin{center}
    $\theta_{{\sf{Code}}\Sigma_{\Q_n'})}^{\Gamma'}=\theta_{\xi+n}^{\Gamma'}$ and $\Omega^{\Gamma'}=\xi+\omega$.
    \end{center}
\item $\Gamma(\P, \Sigma)$ is a mouse full pointclass.
\end{enumerate}
\end{theorem}

%

We finish with a theorem generalizing \cite[Theorem 3.20]{ATHM}. It shows that $\Gamma(\P, \Sigma)$ satisfies mouse capturing for any $\Sigma_\Q$ where $\Q\in pI(\P, \Sigma)$. Recall from \cite{ATHM} (the first page of the introduction of \cite{ATHM}) that \textsf{MC} stands for mouse capturing, i.e., for the statement that for $x, y\in \bR$, $x\in OD_y$ if and only if there is an $\omega_1$-iterable $y$-mouse $\M$ such that $x\in \M$. Given an allowable pair $(\P, \Sigma)$ such that $\Sigma$ has strong branch condensation and is $\Gamma^*$-fullness preserving for some projectively closed pointclass $\Gamma^*$, we say \textsf{MC} holds for $\Sigma$\footnote{The statement ``\textsf{MC} holds for $\Sigma$" can be made precise for an arbitrary strategy with hull condensation. Our definition also includes st-strategies.} if for $x, y\in \bR$, $x\in OD_{y, \Sigma}$ if and only if there is an $\omega_1$-iterable $\Sigma$-mouse $\M$ over $y$ such that $x\in \M$. Given a mouse full pointaclass $\Gamma$ and a allowable pair $(\P, \Sigma)\in \Gamma$ such that $\Sigma$ is $\Gamma$-fullness preserving and has strong branch condensation, we write
\begin{center}
$\Gamma\models ``\sf{MC}$ for $\Sigma$"
\end{center}
if one of the following holds:
\begin{enumerate}
\item $\Gamma$ is completely mouse full and whenever $A$ is a new set such that $\Gamma=W_A$ then $L(A, \bR)\models ``\sf{MC}$ for $\Sigma$".
\item $\Gamma$ is not completely mouse full and if $(\Gamma_\a: \a<\Omega)$ are the completely mouse full pointclasses witnessing that $\Gamma$ is mouse full then for some $\a<\Omega$, $L(\Gamma_\a, \bR)\models ``\sf{MC}$ for $\Sigma$".
\end{enumerate}

\begin{theorem}\label{smc in the derived model} Suppose $(\P, \Sigma)$ is an allowable pair of limit type and $\Sigma$ has strong branch condensation and is $\Gamma^*$-fullness preserving for some projectively closed pointclass $\Gamma^*$. Suppose further that there is a good pointclass $\Gamma$ such that ${\sf{Code}}(\Sigma)\in \Delta_{\utilde{\Gamma}}$. Then for every $\Q\in pB(\P, \Sigma)$, 
\begin{center}
$\Gamma(\P, \Sigma)\models ``\sf{MC}$ for $\Sigma_\Q$".
\end{center}
\end{theorem}

\section{Anomalous hod premice}\label{sec: anomalous hod premice}

In this paper, we use anomalous hod premice the same way we used them in \cite{ATHM}, to generate pointclasses that are mouse full but not completely mouse full. The reader may wish to review \rdef{germane lses} and \rdef{sis}.

\begin{definition}[Anomalous hod premouse of type I]\label{anomalous hod premouse of type I}
$\P$ is an \textbf{anomalous hod premouse of type I} if $\P$ is a germane ${\sf{hp}}-\sf{lses}$ such that letting $\Q={\sf{hl(\P)}}$, $\Q$ is of successor type, $\P\models ``\d^\Q$ is Woodin" and either $\rho(\P)<\d^\Q$ or $\mathcal{J}_\omega[\P]\models ``\d^\Q$ is not a Woodin cardinal".  $\myqedhere$
\end{definition}


\begin{definition}[Anomalous hod premouse of type II]\label{anomalous hod premouse of type II}
$\P$ is an \textbf{anomalous hod premouse of type II} if $\P$ is a germane ${\sf{hp}}-\sf{lses}$ such that letting $\Q={\sf{hl(\P)}}$, $\Q$ is a gentle hod premouse, $\rho(\P)< \d^\Q$ but for every $\xi\in (\d, \ord\P))$, $\rho(\P||\omega\xi)\geq \d^\Q$. $\myqedhere$
\end{definition}


\begin{definition}[Anomalous hod premouse of type III]\label{anomalous hod premouse of type III}
$\P$ is \textbf{an anomalous hod premouse of type III} if $\P$ is a germane ${\sf{hp}}-\sf{lses}$ such that letting $\Q={\sf{hl(\P)}}$, $\Q$ is non-gentle limit type hod premouse, $\rho(\P)<\d^{\Q^b}$ but for every $\omega\xi<\ord(\P)$, $\rho(\P||\omega\xi)>\d^{\Q^b}$\footnote{It follows from the arguments on page 142 of \cite{ATHM} that $\rho(\P||\omega\xi)=\d^{\Q^b}$ is not possible in situations that will arise in this book.}. $\myqedhere$
\end{definition}

Thus, in the language of \rdef{projecting types}, if $\P$ is an anomalous hod premouse then $\P$ is not projecting well but all of its initial segments do project well. We say $\P$ is an anomalous hod premouse if it is an anomalous hod premouse of type $i$ where $i\in \{ I, II, III\}$. 

\begin{definition}[Anomalous hod pair]\label{anomalous hod pair}
$(\P, \Sigma)$ is \textbf{an anomalous hod pair} if one of the following conditions holds:
\begin{enumerate}
\item $\P$ is an anomalous hod premouse of type I or II, $\Sigma$ is an $(\omega_1, \omega_1)$-iteration strategy with hull condensation and whenever $\Q$ is a $\Sigma$ iterate of $\P$, $\Sigma^\Q\subseteq\Sigma\rest \Q$\footnote{Recall that $\Sigma^\Q$ is the internal strategy of $\Q$.}.
\item $\P$ is an anomalous hod premouse of type III, $\Sigma$ is a $(\omega_1, \omega_1, \omega_1)$-iteration strategy\footnote{See \rdef{the un-dropping iteration game}.} with hull condensation and whenever $\Q$ is a $\Sigma$ iterate of $\P$, $\Sigma^\Q\subseteq\Sigma\rest \Q$.
\end{enumerate}
We then say that $(\P, \Sigma)$ is a \textbf{simple anomalous hod pair} if either 
\begin{itemize}
\item it is an anomalous hod pair and $\P$ is of type I or II, or
\item $\P$ is an anomalous hod premouse of type III, $\Sigma$ is a $(\omega_1, \omega_1)$-iteration strategy with hull condensation and whenever $\Q$ is a $\Sigma$ iterate of $\P$, $\Sigma^\Q\subseteq\Sigma\rest \Q$.
\end{itemize}
$\myqedhere$
\end{definition}

The following lemma is due to Mitchell and Steel. It appears as Claim 5 in the proof of Theorem 6.2 of \cite{FSIT}. In the current work, the lemma is used to show that certain hod pair constructions converge, which leads to showing that generation of pointclasses holds (see \rthm{the generation of mouse full pointclasses}). It was used in \cite{ATHM} in a similar fashion (see \cite[Lemma 3.25]{ATHM}).

\begin{lemma}\label{resurrection of strong uniqueness} Suppose $(\P, \Sigma)$ is an anomalous hod pair or a simple hod pair such that for $n<k(\P)$, $(\P, n)$ is not anomalous. Let $k=k(\P)$, $\P'=(\P, k-1)$ and $\Sigma'=\Sigma_{\P'}$, and suppose  $(\T, \Q)\in I(\P', \Sigma')$. Then 
\begin{itemize}
\item if $\P$ is of type I or II then $\rho_k(\Q)<\d^\Q$ and 
\item if $\P$ is of type III then $\rho_k(\Q)<\d^{\Q^b}$.
\end{itemize}
\end{lemma}

The next theorem is the adaptation of \cite[Theorem 3.27]{ATHM} to our current setting. It generalizes our results from previous sections to anomalous hod pairs.

\begin{theorem}\label{derived model of gamma full anomalous hod pairs} Suppose $(\P, \Sigma)$ is an anomalous hod pair of type II or III. Suppose that there is a projectively closed pointclass $\Gamma$ such that for any $(\T, \Q)\in B(\P, \Sigma)$ there is a hod pair $(\R, \Lambda)$ such that $\Lambda$ has (strong) branch condensation and is low-level $\Gamma$-fullness preserving\footnote{See \rdef{gamma fullness preservation}.}, and there is $\pi: \Q\rightarrow \R$ such that $\Lambda^\pi=\Sigma_{\Q, \T}$. Then
\begin{enumerate}
\item For every $(\T, \Q)\in B(\P, \Sigma)$, $\Sigma_{\Q, \T}$ has (strong) branch condensation, is positional and is commuting.
\item $\Sigma$ is strongly low-level $\Gamma(\P, \Sigma)$-fullness preserving and $\Gamma(\P, \Sigma)$ is a mouse full pointclass.
\end{enumerate}
\end{theorem}

We omit the proof of \rthm{derived model of gamma full anomalous hod pairs} as it is only notationally more complicated than the proof of \cite[Theorem 3.10]{ATHM}.  We remind the reader that the proof of \cite[Theorem 3.27]{ATHM} depended on the generic interpretability result, which appeared as \cite[Theorem 3.10]{ATHM}. In our current context we need to use \rthm{generic interpretability holds}. The general idea is that we can translate the properties of $\Sigma$ into the derived model of $\P$ as computed via $\Sigma$. This fact then just gets preserved under pull-back embeddings. 

It is also possible to prove a version of \rthm{derived model of gamma full anomalous hod pairs} for sts hod pairs. To prove it, we again need to use \rthm{generic interpretability holds}. We state it without a proof.

\begin{theorem}\label{gamma fullness preservation sts hod pairs} Suppose $(\P, \Sigma)$ is an sts hod pair and $\Gamma$ is a projectively closed pointclass. Suppose that for any $(\T, \Q)\in B(\P, \Sigma)$ there is a hod pair $(\R, \Lambda)$ such that $\Lambda$ has strong branch condensation and is (strongly) $\Gamma$-fullness fullness preserving, and there is $\pi: \Q\rightarrow \R$ such that $\Lambda^\pi=\Sigma_{\Q, \T}$. Then
\begin{enumerate}
\item For every $(\T, \Q)\in B(\P, \Sigma)$, $\Sigma_{\Q, \T}$ has (strong) branch condensation, is positional and is commuting.
\item $\Sigma$ is strongly $\Gamma^b(\P, \Sigma)$-fullness preserving and $\Gamma^b(\P, \Sigma)$ is a mouse full pointclass.
\end{enumerate}
\end{theorem}

The following is an easy corollary of \rthm{derived model of gamma full anomalous hod pairs}.

\begin{corollary}[Branch condensation pulls back]\label{branch condensation pulls back} Suppose $(\P, \Sigma)$ is a hod pair of limit type and $\Sigma$ has (strong) branch condensation. Suppose $\pi:\Q\rightarrow \P$ is elementary. Then for every $\R\inseg_{hod} \Q$ such that $\d^\R$ is a cutpoint of $\Q$, $(\Sigma^\pi)_{\R}$ has (strong) branch condensation.
\end{corollary}

\section{Branch condensation on a tail}

The main theorem of this section, \rthm{branch condensation on a tail}, will be used in several places (e.g. the proof of \rthm{the generation of mouse full pointclasses ii}) in this book as well as in core model induction applications. First we need to introduce a new concept, which fortunately for us, Farmer Schlutzenberg has developed independently and much more generally.

\begin{definition}\label{supported anomalous pair} Suppose $(\P, \Sigma)$ is an anomalous pair of type $III$. We say $(\P, \Sigma)$ has a \textbf{supporting bicephalous} if there is a bicephalous $B=(\rho, \M, \P)$ in the sense of \cite[Definition]{FarmStrongs} such that
\begin{enumerate}
\item $\rho=\d^{\P^b}$,
\item $\M$ is germane\footnote{See \rdef{germane lses}.} and such that $\rho(\M)<\rho$, ${\sf{hl}}(\M)=\P|\rho$ and $\M\inseg {\sf{Lp}}^{\Sigma_{\P|\rho}}(\P|\rho)$,
\item for every $n<\omega$, $\rho_n(\M)\not =\rho$,
\item $k(\M)$ is the least $n$ such that $\rho_{n+1}(\M)<\rho$,
\item for every $\gg\in [{\sf{ord}}(\P^b), {\sf{ord}}(\M))$, $\rho(\M||\gg)>\rho$, 
\item $B$ has an $\omega_1$-iteration strategy $\Sigma^+$ which extends $\Sigma$. 
\end{enumerate}
$\myqedhere$
\end{definition}

\begin{remark} The reader unfamiliar with \cite{FarmStrongs} may treat $B$ in \rdef{supported anomalous pair} as a pair constructed by some $\Gamma$-hod pair construction. After reaching $\P|\rho$ the hod pair construction aims to reach the next $\Gamma-{\sf{cbl}}$\footnote{See \rdef{cbl}.}. The construction proceeds as a fully backgrounded construction relative to $\Sigma_{\P|\rho}$. Once $\P^b$ is reached it is declared to be a layer and a new strategy appears, the strategy of $\P^b$. To reach $\M$ we just simply need to continue the construction relative to $\Sigma_{\P|\rho}$. This will all be relevant in the proof of  \rthm{the generation of mouse full pointclasses ii}. Also notice that clause 4 implies that we iterate $\M$ using one fine structural level lower than one would normally do. $\myqedhere$
\end{remark}

\begin{theorem}[Branch condensation on a tail]\label{branch condensation on a tail} Suppose $(\P, \Sigma)$ is an anomalous hod pair of type II or III. Suppose that for every $(\T, \Q)\in B(\P, \Sigma)$, $\Sigma_{\Q, \T}$ has strong branch condensation. Moreover, assume either\\\\
(1) ${\sf{ZFC}}$ holds and $\P$ is of type $II$, or\\
(2) ${\sf{AD^+}}$ holds and if $\P$ is of type $III$ then it has a supporting bicephalous.\\\\
Then if $\P$ is of type $II$ then there is $(\T, \Q)\in I(\P, \Sigma)$ such that $\Sigma_{\Q, \T}$ has strong branch condensation, and if $\P$ is of type $III$ then $\Sigma$ has strong branch condensation.
\end{theorem}
\begin{proof}
The case when $\P$ is of type $II$ is very similar to the proof of \cite[Theorem 3.28]{ATHM}. The case when $\P$ is of type $III$ is similar to the proof of \rthm{strong condensation for backgrounded strategies}. In order not to repeat the entire \rsec{proving branch condensation sec}, we outline the proof of branch condensation and leave the rest to reader. Let $B=(\rho, \M, \P)$ and $\Sigma^+$ witness that $(\P, \Sigma)$ has a supporting bicephalous.

Suppose $(\T, \Q, \U, c, \sigma)$ is such that 
\begin{enumerate}
\item $(\T, \Q)\in I(\P, \Sigma)$,
\item $\U$ is a stack according to $\Sigma$ and $\lh(\U)$ is a limit ordinal,
\item $c$ is a cofinal well-founded branch of $\U$,
\item $\sigma:\M^\U_b\rightarrow \Q$ is elementary\footnote{$\sigma$ may not be fully elementary, just at the right fine structural level.} and such that $\pi^\T=\sigma\circ \pi^{\U}_b$.
\end{enumerate}
We would like to show that $c=\Sigma(\U)$. Let $d=\Sigma(\U)$. The most dificult case, which also represents the difficulties involved in other cases that are left to the reader, is the case when $\pi^{\U, b}$ exists, $R^\U$ has a maximal element, and if $\a=\max(R^\U)$ then $\U_{\geq \a}$ is above $(\M_\a^\U)^b$ and $\Q(c, \U_{\geq \a})$-exists. Set then $\W=\m^+(\U)$ and let $\Phi_0'$ be the $\sigma$-pullback of $\Sigma_{\sigma(\W), \T}$ and $\Phi_1'=\Sigma_{\W, \U^\frown\{d\}}$. Finally, set $\Phi_0=(\Phi_0')^{stc}$ and $\Phi_1=(\Phi_1')^{stc}$. We then have that if $\Phi_0=\Phi_1$ then in fact, as $\Q(d, \U)$ exists, $\Q(d, \U)=\Q(c, \U)$ and therefore, $c=d$. Assume then that $\Phi_0\not=\Phi_1$.

Let then $(\X_0, \W_0, \X_1, \W_1, \R)$ be a minimal low level disagreement\footnote{See \rdef{low level disagreement between strategies}.} between $\Phi_0$ and $\Phi_1$. Recall the notation $\X^{\sf{ue}}$ introduced in clause 4 and 5 of \rnot{notation for generalized stacks}. Let $\Y_0=\U^\frown \{c\}^\frown (\X_0)_\R^{\sf{ue}}$ and $\Y_1=\U^\frown \{c\}^\frown (\X_1)_\R^{\sf{ue}}$. Let $\Y_0^+$ and $\Y_1^+$ be the stacks on $B$ obtained by applying $\Y_0$ and $\Y_1$ to $B$. $\Y_1^+$ is according to $\Sigma^+$ while $\Y_0^+$ is not according to a strategy but it has well-founded models because of $\sigma$. Let $B_0=(\nu_0, \M_0, \P_0)$ and $B_1=(\nu_1, \M_1, \P_1)$ be the last models of $\Y_0^+$ and $\Y_1^+$. Let $\Lambda_0'$ be the strategy of $B_0$\footnote{This strategy comes from the copying procedure; $B_0$ embeds into a $\Sigma^+$-iterate of $B$ that starts out by applying $\T$ to $B$ and then copies $(\X_0)_\R^{\sf{ue}}$.} and let $\Lambda_1'=\Sigma^+_{B_1, \Y_1^+}$. 

We now have that $\pi^{\Y_0, b}=\pi^{\Y_1, b}$, which implies that  $\nu_0=\nu_1$, $\M_0=\M_1$ and $\P_0^b=\P_1^b$. In fact, letting $F$ be the $(\rho, \pi^{\Y_0, b}(\rho))$-extender\footnote{$\rho=\d^{\P^b}$.} derived from $\pi^{\Y_0, b}$ then for $i\in 2$, $\M_i=Ult(\M, F)$ and $\P_i=Ult(\P^b, F)$. Let then $\nu=\nu_0$, $\N=\M_0$, $\S=\P_0^b$, $\Lambda^0$  be the strategy of $ \N$ induced by $\Lambda_0'$ and $\Lambda^1$ be the strategy of $\N$ induced by $\Lambda_1'$. We have that $\Lambda^0_\R\not= \Lambda^1_\R$.

Notice next that because $\Sigma_{\P|\rho}$ has branch condensation, $\rho(\N)<\d^\R$\footnote{See \rlem{resurrection of strong uniqueness}. This follows from the proof of Claim 5 in the proof of Theorem 6.2 of \cite{FSIT}.} and moreover, letting $n=k(\N)$\footnote{We abuse our notation and think of $\N$ as both fine structural and non-fine structural.}, and $\N'$ be the $n$th reduct of $\N$ then \\\\
(1) $\sup(Hull^{\N'}_1(p_1(\N')\cup \R^-)\cap \d^\R)=\d^\R$. \\\\
Let then $\K'=cHull^{\N'}_1(p_1(\N')\cup \R)$ and $i':\K'\rightarrow \N'$ be the uncollapsing embedding. Let $\K$ be decoding of $\K'$ and $i:\K\rightarrow \N$ be the canonical uncoring embedding. (1) then implies that\\\\
(2) $\K\models ``\d^\R$ is a Woodin cardinal" and $\mathcal{J}_{\omega}[\K]\models ``\d^\R$ is not a Woodin cardinal". \\\\
Let now $\Lambda^2=(i$-pullback of $\Lambda^0$) and $\Lambda^3=(i$-pullback of $\Lambda^0$). Let $\mathcal{Z}$ be a normal stack\footnote{We can choose $\mathcal{Z}$ to be normal because of \rthm{universality of background construction}.} on $\K$ based on $\R$ such that $\mathcal{Z}^d=_{def}\downarrow(\mathcal{Z}, \R)$ is according to both $\Lambda^2$ and $\Lambda^3$ and such that setting $e=\Lambda^2(\mathcal{Z})$ and $f=\Lambda^3(\mathcal{Z})$ then \\\\
(3) $e\not = f$.\\\\
It follows from \rthm{derived model of gamma full anomalous hod pairs} that $\Q(e, \mathcal{Z}^d)$ exists if and only if $\Q(f, \mathcal{Z}^d)$ exists, and therefore, neither exists. Let $\K_e=\M_e^{\mathcal{Z}}$ and $\K_f=\M_f^{\mathcal{Z}}$.

 Let $\Lambda^e=\Lambda^2_{\K_e, \mathcal{Z}^\frown \{e\}}$ and $\Lambda^f=\Lambda^3_{\K_f, \mathcal{Z}^\frown \{f\}}$. Letting $\tau=\d(\mathcal{Z})$, we now have that\\\\
(4) $\K_e, \K_f\models ``\d^\R$ is a Woodin cardinal" and $\mathcal{J}_{\omega}[\K_e], \mathcal{J}_{\omega}[\K_f]\models ``\d^\R$ is not a Woodin cardinal",\\
(5) $\K_e$ and $\K_f$ are $\tau$-sound,\\
(6) $\pi^{\mathcal{Z}}_e(\R)=\pi^{\mathcal{Z}}_f(\R)=_{def}\R_1$ and $\Lambda^e_{\R_1}=\Lambda^f_{\R_1}=_{def}\Lambda_1$\footnote{This is a consequence of \rthm{universality of background construction}. $\mathcal{Z}$ is produced by iterating $\R$ into a universal model.}.\\\\
Thus, if we argue that $\K_e=\K_f$ then we would be done as it would show that $e=f$, contradicting (3). Set $\Gamma_e=\Gamma(\K_e, \Lambda^e)$ and $\Gamma_f=\Gamma(\K_f, \Lambda^f)$. Suppose first that $\Gamma_e=\Gamma_f$. Then because $\Lambda^e$ and $\Lambda^f$ both are $\Gamma_e$-fullness preserving and therefore, $(\K_e, \Lambda^e)$ and $(\K_f, \Lambda^f)$ can be compared as in \rthm{diamond comparison}, (5) implies that $\K_e=\K_f$.

We now assume that $\Gamma_e\not =\Gamma_f$. Without losing generality, lets suppose that $\Gamma_e\subset \Gamma_f$. It follows from the above argument that $\K_e$ is ordinal definable in $\Gamma_f$ from $\Lambda_1$.
Indeed, let $A\in \Gamma_f$ be such that every set in $\Gamma_e$ has Wadge rank $<w(A)$. Then in $L(A, \bR)$, $\K_e$ is the unique anomalous hod premouse $\V$ that has an $\omega_1$-iteration strategy $\Pi$ such that
\begin{enumerate}
\item $\Gamma(\V, \Pi)=\{ C\subseteq \bR: w(C)<w(\Gamma_e)\}$,
\item $\R_1\insegeq^c_{hod}\V$\footnote{See \rdef{complete layer notation}.},
\item $\Pi_{\R_1}=\Lambda_1$,
\item $\V$ is $\tau$-sound. 
\end{enumerate}
It then follows from \rthm{derived model of gamma full anomalous hod pairs} that $\K_e\in \K_f$, which contradicts (4).
\end{proof}

\chapter{The internal theory of lsa hod mice}\label{lsa internal theory chapter}

A major shortcoming of our treatment of short-tree-strategy mice is that we did not add branches to all trees. Suppose $(\P, \Sigma)$ is an sts hod pair, $X$ is a self-well-ordered set such that $\P\in X$ and $\M$ is a $\Sigma$-sts premouse over $X$ based on $\P$. Recall the short tree strategy indexing scheme \rdef{weak psi alpha indexing scheme a}. Recall that our strategy for indexing branches was to consider two kinds of iterations, $\sf{uvs}$ and $\sf{nuvs}$\footnote{See \rdef{nus stacks}.}. We outright index the branches of $\sf{uvs}$ iterations. However, we only consider a subclass of $\sf{nuvs}$ iterations. If for some $\b<o(\M)$, $\T\in \dom(\Sigma^{\M|\b})$ is an $\M|\b$-ambiguous tree then (i) $\T$ is a result of comparing $\P$ with a certain background construction of $\M|\b$ and (ii) we index the branch of $\T$ after we find a certain certificate of shortness (recall \rdef{weak psi alpha indexing scheme a}). It is then not clear from our definition that $\Sigma\rest \M\subseteq \M$. The main goal of this chapter is to show that, provided $\M$ is sufficiently closed, $\Sigma\rest \M$ is a definable class of  $\M$. Below we make our goal more precise. \\\\
\textbf{Motivational Question.} Suppose $(\P, \Sigma)$ is a hod pair or an sts hod pair, $X$ is a self-well-ordered set such that $\P\in X$ and $\M$ is a $\Sigma$ or $\Sigma$-sts mouse over $X$ (see \rdef{lambda sts premouse}). Is $\Sigma\rest \N$ definable over $\N$? Is $\Sigma\rest \N[g]$ definable over $\N[g]$ where $g$ is $\N$-generic? \\

In \rsec{generic interpretability sec} we gave an answer to Motivational Question in the case $\M$ is $\P$ itself (see \rthm{generic interpretability holds}). Another answer was given by \cite[Lemma 3.35]{ATHM}, where it was shown that $\Sigma\rest \N[g]$ is definable over $\N[g]$ provided $\P$ doesn't have non-meek layers. Here, we are mainly concerned with proving a version of  \cite[Lemma 3.35]{ATHM} in the case of a non-meek hod premice. Because of this we will state many of our definitions and theorems for hod pairs or sts hod pairs $(\P, \Sigma)$ such that $\P$ is non-meek\footnote{See \rdef{pre-hod-like}.}. To simplify our terminology, we will say $(\P, \Sigma)$ is a non-meek hod pair\index{non-meek hod pair} if $\P$ is a non-meek hod premouse and $\Sigma$ is either an iteration strategy or an sts-strategy (this is only allowed in the case $\P$ is of lsa type). 

While a positive answer to the Motivational Questions is desirable, it is naive to hope that one exists for all such $\N$. A positive answer depends on how closed $\N$ is. If for instance the branch of $\T$ is given via a $\Q$-structure that is beyond the $\#$-operator while our $\N$ is only closed under the $\#$-operator then, in most cases, identifying the correct branch of $\T$ inside $\N$ via a procedure that is uniform in $\T$ will be impossible. In this chapter, we give a positive answer to the Motivational Question provided our $\N$ is sufficiently closed. We make this notion more precise. 

Suppose $(\P, \Sigma)$ is a non-meek hod pair and  $\N$ is a $\Sigma$-mouse such that $\N\models \sf{ZFC}$$-Replacement$. We say $\N$ is $\Sigma$-closed if $\Sigma\rest \N\subseteq \N$. 
We say $\N$ is generically $\Sigma$-closed if $\N$ is $\Sigma$-closed and whenever $g$ is $\N$-generic, $\Sigma\rest \N[g]$ is definable over $(\N[g], \in)$ (in the language of $\Sigma$-premice). It is worth remarking that the structure $(\N[g], \in)$ is a structure in the language of $\Sigma$-premice and in particular, there are names for $\vec{E}^\N$ and $\Sigma^\N$.

\begin{definition}\label{generically sigma-closed} 
We say $\N$ is uniformly generically $\Sigma$-closed if $\N$ is generically $\Sigma$-closed and there are formulas $\phi$ and $\psi$ (in the language of $\Sigma$-premice) such that for any $\N$-generic $g$, any stack $\T\in \N[g]$ on $\P$ and any $b\in \N[g]$,
\begin{center}
$\T\in \dom(\Sigma)\iff (\N[g], \in)\models \phi[\T]$\\
$\Sigma(\T)=b\iff (\N[g], \in)\models \psi[\T, b]$
\end{center} 
$\myqedhere$
\end{definition}


%

The main theorem of this chapter is \rthm{main theorem on gen int}. It gives a positive answer to our Motivational Question in the case $\N$ is $\Sigma$-closed and has fullness preserving iteration strategy (see  \rdef{sigma closed mouse} and  \rdef{fullness preserving strategy for sigma mice}). The main idea behind the proof of \rthm{main theorem on gen int} is that the branch of an iteration tree $\T$ on $\P$ can be identified by the authentication process introduced in \rdef{authentic iterations}.

Given a transitive set $X$, we let $X^\#$ be the least sound active mouse over $X$. Also recall that if $X$ is any set and $A\subseteq X^2$ then $p[A]$ is the projection of $A$ onto one of the coordinates of $A$. The specific coordinate onto which we project will always be clear from the context.

\section{Internally $\Sigma$-closed mice}\label{sigma-closed mice sec}

In this section we introduce a kind of closure property of hybrid mice for which we can give a positive answer to our motivational question. The first such closure property is \text{internal closure}, which postulates that our mouse has enough of the strategy.

\begin{definition}[Internally $\Sigma$-closed mouse]\label{sigma closed mouse} \index{internal $\Sigma$ closure} Suppose 
\begin{itemize}
\item $(\P, \Sigma)$ is an allowable pair\footnote{See \rdef{allowable pair}.}, 
\item $\P$ is a non-meek hod premouse,
\item if $\Sigma$ is an iteration strategy then $\N$ is a $\Sigma$-premouse over some $X$ based on $\P$, and
\item if $\Sigma$ is an sts premouse then $\N$ is a $\Sigma$-sts premouse over some $X$ based on $\P$.  
\end{itemize}
\begin{enumerate}
\item We say $\N$ is \textbf{internally $\Sigma$-closed premouse} if  for every $\k<\ord(\N)$ there is $\M\insegeq \N$ such that
\begin{enumerate}
\item $\M\models \sf{ZFC}$, 
\item $\N|\k\insegeq \M$,
\item for every $(\T, \S)\in B(\P, \Sigma^\M)$\footnote{Thus, $(\T, \S)\in \M$.}, $\Sigma^\M_{\S}$ is total in $\M$\footnote{Thus, $\Sigma^\M_{\S}=\Sigma_{\S, \T}\rest \M$.}, 
\item $\M$ has at least three Woodin cardinals that are greater than $\k$,
\item letting $\d_0<\d_1<\d_2$ be the first three Woodin cardinals of $\M$ that are greater than $\k$, for every $i\in 3$ and $\eta\in [\k, \d_i)$, letting $(( \M_\gg , \N_\gg : \gg\leq \nu), ( F_\gg: \gg<\nu), ( \T_\gg: \gg\leq\nu))$
be  the output of the $(\P, \Sigma)$-coherent fully backgrounded construction of $\M|\d_i$\footnote{See Definition \ref{full short tree coherent background constructions}.} in which extenders used have critical points $>\eta$,  the following conditions hold:
\begin{enumerate}
\item If $\Sigma$ is an iteration strategy then $\pi^{\T_\nu}$-exists and $\M_\nu$ is the last model of $\T_\nu$.
\item If $\Sigma$ is an st-strategy then $\pi^{\T_\nu, b}$ exists and $\T_\nu$ is $\M$-terminal\footnote{See \rdef{terminal tree}.}.
\item If $(\T, \S)\in B(\P, \Sigma^{\N|\eta})$\footnote{See \rdef{gamma(p, sigma) and b(p, sigma) for sts}.} then for some $\b$, $\M_\nu|\b$ is a $\Sigma^{\N}_{\S}$-iterate of $\S$ via a normal stack. 
\end{enumerate}
\end{enumerate}
\item If $\M, \N$ and $\k$ are as above then we say $\M$ \textbf{witnesses} the internal $\Sigma$-closure of $\N$ at $\k$. 
\item We say $\N$ is an \textbf{internally $\Sigma$-closed mouse} (sts mouse) if it is an internally $\Sigma$-closed premouse and has an $\omega_1$-iteration strategy $\Lambda$ witnessing that $\N$ is a $\Sigma$-mouse (sts mouse). 
\end{enumerate}$\myqedhere$
\end{definition}

Two remarks are in order. First notice that internal $\Sigma$-closure is a first order property of $\N$, and in clause 3 above we do not need to require that $\Lambda$-iterates of $\N$ are internally $\Sigma$-closed as this is just a consequence of elementarity. 

Secondly, we cannot in general hope to prove that generic interpretability holds for internally $\Sigma$-closed mice. The reason is that there might be $\Q\in B(\P, \Sigma)$ such that $\Sigma_\Q$ is beyond the iteration strategy of $\N$ (in the sense that $\Lambda <_w \Sigma_\Q$), and if such a $\Q$ is generic over $\N$ then it is not wise to hope that $\Sigma_\Q\rest \N$ would be definable over $\N[\Q]$. In order to prove generic interpretability result for internally $\Sigma$-closed premice we need to find a \textit{fullness condition} that would let us take care of examples as above. In particular, we seem to need to require that any $\Sigma_\Q$ as above is strictly below the strategy of $\N$. The next couple of paragraphs make this intuitive notion more precise. 

Suppose $\N$ is an internally $\Sigma$-closed mouse, $\k$ is an $\N$-cardinal and $\M$ is as in \rdef{sigma closed mouse}.  Let $\d_0<\d_1<\d_2$ be the first three Woodin cardinals of $\M$ that are greater than $\k$, and let $\eta\in [\k, \d_2)$ and $i\in 3$ be the least such that $\eta<\d_i$. We then let $\S^\M_\eta$ be the $\Sigma^\M$-iterate of $\P$ constructed via the $(\P, \Sigma^\M)$-coherent fully backgrounded construction of $\M|\d_i$ where critical points of extenders used are $>\eta$. We let $\U^\M_\eta$ be the normal tree on $\P$ with last model $\S^\M_\eta$ and 
\begin{center}
$\pi^\M_\eta=\begin{cases} 
\pi^{\U^\M_\eta, b}&: \P\ \text{is of lsa type}\\
\pi^{\U^\M_\eta} &: \text{otherwise.}
\end{cases}$
\end{center}
Notice that $\pi^{\M}_\eta\in \N$.


Keeping the notation and terminology of \rdef{sigma closed mouse}, suppose $\Lambda$ is an iteration strategy for $\N$ (witnessing that $\N$ is an internally $\Sigma$-closed mouse). Suppose $\xi<\ord(\N)$ and $\Lambda^\xi$ is the fragment of $\Lambda$ that acts on stacks above $\xi$.
We then let $\Gamma(\N, \Lambda^\xi)$ be the collection of all sets $A\subseteq \bR$ such that for some $(\T, \R)\in I(\N, \Lambda)$, $\k<\ord(\R)$ and $\M\insegeq \R$ witnessing that $\R$ is internally $\Sigma$-closed at $\k$ the following holds: letting $\d_0<\d_1<\d_2$ be the first three Woodin cardinals of $\M$ that are greater than $\k$, whenever $\eta\in[\k, \d_2)$, there is $\Q\inseg_{hod}(\S^\M_\eta)^b$ 
 such that 
 \begin{enumerate}
 \item $\S^\M_\eta\models ``\d^\Q$ is a Woodin cardinal" and
\item $w(A)\leq w({\sf{Code}}(\Sigma_{\Q, \U^\M_\eta}))$.
\end{enumerate}

\begin{remark}\label{convenient notation} For convenience, we will use the notation $\Gamma(\P, \Sigma)$ for both sts pairs and hod pairs. In the case of sts hod pairs, it is just $\Gamma^b(\P, \Sigma)$. $\myqedhere$
\end{remark}

\begin{definition}\label{fullness preserving strategy for sigma mice} Suppose $\N$ is as in \rdef{sigma closed mouse} and $\Lambda$ is an $\omega_1$-iteration strategy for $\N$ (witnessing that $\N$ is an internally $\Sigma$-closed mouse). We then say that $\Lambda$ is a \textbf{fullness preserving} iteration strategy for $\N$ if  for every $\xi<\ord(\N)$, letting $\Lambda^\xi$ be the fragment of $\Lambda$ that acts on stacks above $\xi$,  $\Gamma(\N, \Lambda^\xi)=\Gamma(\P, \Sigma)$. $\myqedhere$
\end{definition}

The following is our generic interpretability result for internally $\Sigma$-closed mice $\N$ that have a fullness preserving iteration strategy. 

\begin{theorem}\label{main theorem on gen int} Assume $\sf{NsesN}$\footnote{See \rdef{no mouse with a superstrong}.} Suppose $(\P, \Sigma)$ is an allowable pair, $\Gamma$ is a projectively closed pointclass and $\N$ is an internally $\Sigma$-closed premouse (possibly over some set $X$). Suppose $\Sigma$ is 
\begin{itemize}
\item strongly $\Gamma$-fullness preserving,
\item has strong branch condensation and
\item is commuting\footnote{See \rdef{positional and commuting for sts pairs}.}.
\end{itemize}
 Then the following conditions hold.
\begin{enumerate}
\item If $(\P, \Sigma)$ is a hod pair and $\N$ is a $\Sigma$-mouse then for any $\N$-generic $g$, $\N[g]$ is $\Sigma$-closed and $\Sigma\rest \N[g]$ is uniformly in $g$ definable over $\N[g]$.
\item If $(\P, \Sigma)$ is an sts hod pair and $\N$ is a $\Sigma$-sts mouse with a fullness preserving iteration strategy then for any $\N$-generic $g$, $\N[g]$ is $\Sigma$-closed and $\Sigma\rest \N[g]$ is uniformly in $g$ definable over $\N[g]$.
\end{enumerate}
\end{theorem}

In the next few sections, we will develop the terminology we need to prove \rthm{main theorem on gen int}. We will not give the proof of clause 1 of \rthm{main theorem on gen int}. It is much easier than the proof of clause 2 of \rthm{main theorem on gen int} and it is very much like the proof of \cite[Theorem 3.10]{ATHM}. Thus, we only concentrate on sts hod pairs.

\section{Authentication procedure revisited}\label{sec: authentication revisited}

Suppose $(\P, \Sigma)$ is an sts hod pair,  $\N$ is an internally $\Sigma$-closed premouse, $g$ is $\N$-generic and  $\T\in \dom(\Sigma)\cap \N[g]$ is a normal stack on $\P$ above $\P^b$ such that $\T$ doesn't have fatal drops. Suppose first that $\T\in b(\Sigma)$. In this case, we would like to identify $\Q(b, \T)$ in $\N[g]$ via a procedure that is uniform in $\T$. Here $b=\Sigma(\T)$. Clearly if $\Q(b, \T)\insegeq \m^+(\T)$ then we can easily identify $\Q(b, \T)$. Suppose then $\m^+(\T)\inseg \Q(b, \T)$. We now face two problems. 

The first problem is showing that $\Q(b, \T)\in \N[g]$ and the second is showing that $\Q(b,\T)$ can be identified in $\N$ in a uniform manner. Both of these require more of $\N$ than just internal $\Sigma$-closure. To prove both of these facts, we will need that $\N$ has a fullness preserving iteration strategy. Our strategy for finding $\Q(b, \T)$ in $\N$ is that if $\N$ is sufficiently rich then some backgrounded construction will reach $\Q(b, \T)$. To execute this plan, we first need to describe the sort of backgrounded constructions that we will consider. In what follows, we borrow ideas from \rsec{authentic iterations and finite stacks sec}. In particular, it will be helpful to recall \rdef{authentic stacks of length 2} and other definitions from that section.

\begin{definition}[$(\N, X)$-authenticated iteration strategy]\label{n certified iteration strategy} Suppose $(\P, \Sigma)$ is an sts hod pair, $X\subseteq \P^b$ and $\N$ is a $\Sigma$-sts  premouse such that $X\in \N$. Suppose that $g\subseteq\mathbb{P}$ is $\N$-generic for some poset $\mathbb{P}\in \N$ and $\R\in \N[g]$ is an lsa type hod premouse. We define a partial short tree strategy $\Phi^{\N, X, g}_\R$ without a model component for $\R$ as follows. $\Phi^{\N, X, g}_\R$ acts on indexable stacks\footnote{See \rdef{indexable stack}.}.
\begin{enumerate}
\item $t=(\R, \T, \R_1, \T_1)\in \dom(\Phi^{\N, X, g}_\R)\cap  \N[g]$ if and only if $t$ is $(\P, \Sigma^\N, X)$-authenticated\footnote{See \rdef{authentic lsp}.}.
\item Given $t=(\R, \T_0, \R_1, \T_1)\in \dom(\Phi^{\N, X, g}_\R)\cap  \N[g]$ with $\T_1\not =\emptyset$, $\Phi^{\N, X,  g}_\R(t)=b$ if and only if $t^\frown\{b\}$ is $(\P, \Sigma^\N, X)$-authenticated. 

\end{enumerate}
When $X=\P^b$ we simply omit it from our terminology. $\myqedhere$
\end{definition}

Continuing with the $\R, \N$ of \rdef{n certified iteration strategy}, we next define an $\N$-authenticated backgrounded construction over $\R$. This is essentially a fully backgrounded construction relative to $\Phi^{\N, g}_\R$ (see \rdef{fully backgrounded sts construction}).

\begin{definition}\label{certified backgrounded constructions}
Suppose $(\P, \Sigma)$ is an sts hod pair, $X\subseteq\P^b\cap \N$ and $\N$ is a $\Sigma$-sts premouse such that $X\in \N$. Suppose that $g\subseteq\mathbb{P}$ is $\N$-generic for some poset $\mathbb{P}\in \N$ and $Y, \R\in \N[g]$ are such that $Y$ is a self-well-ordered set and $\R\in Y$ is an lsa type hod premouse. Suppose further that $\k$ is an $\N$-cardinal such that  $\{ \mathbb{P}, \R, Y\}\in \N|\k[g]$.

We then say that 
$(( \M_\gg , \N_\gg : \gg\leq \nu), ( F_\gg: \gg<\nu), ( \T_\gg: \gg\leq\nu))$
is the output of the $(\N, \k, X)$-authenticated fully backgrounded construction over $Y$ based on $\R$ in which extenders used have critical points $>\k$ if $(( \M_\gg , \N_\gg : \gg\leq \nu), ( F_\gg: \gg<\nu), ( \T_\gg: \gg\leq\nu))$
is the output of $(\R, \Phi^{\N, X, g}_\R)$-coherent fully backgrounded construction\footnote{See \rdef{full short tree coherent background constructions}.} of $\N$ done over $Y$ using extenders with critical points $>\k$\footnote{If the construction reaches a stack $\T$ such that $\Phi^{\N, X, g}_\R(\T)$ is undefined we stop the construction.}. 

Finally, we say $\Q$ is an $(\N, X)$-authenticated sts mouse over $Y$ based on $\R$ if $\Q\in \N$ and for some $\nu$, 
\begin{itemize}
\item $\{ \mathbb{P}, \R, Y, \Q\}\in \N|\nu[g]$ and 
\item $\Q$ appears as a model in the $(\N, \nu, X)$-authenticated fully backgrounded construction over $Y$ based on $\R$.
\end{itemize} When $X=\P^b$ we simply omit it from our terminology. $\myqedhere$

\end{definition}

Suppose now that $(\P, \Sigma)$ is an sts hod pair, $X\subseteq \P^b$ and $\N$ is an internally $\Sigma$-closed mouse with a fullness preserving iteration strategy $\Lambda$ such that $X\in \N$. Suppose $\mathbb{P}\in \N$ is a poset and $g\subseteq \mathbb{P}$ is $\N$-generic. Suppose further that $Y\in \N[g]$. We then let
\begin{center}
${\sf{Lp}}^{\N, g, X, sts, +}(Y, \R)=\bigcup\{\K\in \N[g]:$ there is an $\N$-cardinal $\k$ such that $\{ \mathbb{P}, \R, Y, \K\}\in \N|\k[g]$ such that $\K$ is an $(\N, \k, X)$-authenticated sound sts mouse over $Y$ based on $\R$ such that $\rho(\K)=\ord(Y)\}$
\end{center}

Again, if $X=\P^b$, then we omit it from the notation.

Notice that we do not know that ${\sf{Lp}}^{\N, g, X, sts, +}(Y, \R)$ is a meaningful object, since we do not know that if $\Q_0$ and $\Q_1$ are authenticated by $\M_0$ and $\M_1$ respectively then they are compatible. This, however, is true when $\R$ is an iterate of $\P$ and $\Sigma$ has strong branch condensation and is strongly $\Gamma$-fullness preserving for some $\Gamma$. This fact will also be verified in the next section.

We can then define $({\sf{Lp}}^{\N, g, X, sts, +}_\a(Y, \R): \a<\ord(\N))$ by induction as usual. More precisely, the sequence is defined via the following recursion.
\begin{enumerate}
\item ${\sf{Lp}}^{\N, g, X, sts, +}_0(Y, \R)=trc(Y,\R)$.
\item ${\sf{Lp}}^{\N, g, X, sts, +}_1(Y, \R)={\sf{Lp}}^{\N, g, X, sts, +}(Y, \R)$.
\item ${\sf{Lp}}^{\N, g, X, sts, +}_{\a+1}(Y, \R)={\sf{Lp}}^{\N, g, X, sts, +}({\sf{Lp}}^{\N, g, X, sts, +}_\a(Y, \R))$.
\item ${\sf{Lp}}^{\N, g, X, sts, +}_\l(Y, \R)=\bigcup_{\a<\l}{\sf{Lp}}^{\N, g, X, sts, +}_\a(Y, \R)$.
\end{enumerate}
When $Y=\mathcal{J}_{\omega}[\R]$ or $X=\P^b$, we omit them from the above notation. 

The ``+" version of the ${\sf{Lp}}$ operator defined above may stack more sts mice than we need. To get the proper operator we need to only consider $\K\in \N[g]$ which have an iteration strategy in $\Gamma(\P, \Sigma)$. This fact can be expressed in a first order manner over $\N[g]$. 

\begin{definition}\label{simple q} Suppose $(\P, \Sigma)$, $\N$ and $(\mathbb{P}, g, X, Y, \R)$ are as above. Let $\K\insegeq {\sf{Lp}}^{\N, g, X, sts, +}(Y, \R)$. We say $\K$ is \textit{simple} (in $\N$) if there is $\k$ and $\M\insegeq \N$ such that
\begin{itemize}
\item $(\mathbb{P}, X, Y, \R, \K)\in \N|\k[g]$,
\item $\M$ witnesses the internal $\Sigma$-closure of $\N$ at $\k$,
\item letting $\d_0<\d_1<\d_2$ be the first three Woodin cardinals of $\M$ that are above $\k$, there is some $\Q\inseg_{hod} (\S^\M_\k)^b$ such that if $\eta\in (\ord(\Q), \d_0)$ is the least such that ${\sf{Lp}}^{\Gamma(\P, \Sigma), \Sigma_\Q}(\M|\eta)\models ``\eta$ is a Woodin cardinal"\footnote{Let $\Q^+$ be the least hod initial segment of $\S^\M_\k$ such that $\Q\inseg \Q^+$ and $\d^{\Q^+}$ is a Woodin cardinal of $\S^\M_\k$. Since $\M$ is closed under $\Sigma_{\Q^+}$, the condition  ${\sf{Lp}}^{\Gamma(\P, \Sigma), \Sigma_\Q}(\M|\eta)\models ``\eta$ is a Woodin cardinal" is first order over $\N[g]$.} then $\K$ is an $(\M|\eta, \ord(\Q), X)$-authenticated sound sts mouse over $Y$ based on $\R$ such that $\rho(\K)=\ord(\R)$.
\end{itemize}
$\myqedhere$
\end{definition}
We then let
\begin{center}
${\sf{Lp}}^{\N, g, X, sts}(Y, \R)=\bigcup\{\K\in \N[g]:\K$ is simple and $\K\inseg {\sf{Lp}}^{\N, g, X, sts, +}(Y, \R)\}$.
\end{center}
We will omit $g$ and $X$ when they are clear from the context. The effect of clause 3 of \rdef{simple q} is that since $\K$ is built by a fully backgrounded construction of $\M|\eta$ using extenders with critical points $>\ord(\Q)$, the strategy $\K$ acquired from the strategy of $\M|\eta$ via the resurrection procedure of \cite[Chapter 12]{FSIT} is in $\Gamma(\P, \Sigma)$. This is because the strategy of $\M|\eta$ that acts on stacks above $\ord(\Q)$ is in $\Gamma(\P, \Sigma)$. Thus, the following claim is true.

\begin{proposition}\label{matching lp operators} Suppose $(\P, \Sigma)$, $\N$ and $(\mathbb{P}, g, X, Y, \R)$ are as in \rdef{simple q}. Suppose $\R=\m^+(\T)$ where $\T$ is a normal stack according to $\Sigma$. Suppose further that $\Sigma$ has strong branch condensation and is $\Gamma$-fullness preserving for some projectively closed pointclass $\Gamma$. Then ${\sf{Lp}}^{\N, g, X, sts}(\R)\insegeq {\sf{Lp}}^{\Gamma(\P, \Sigma), \Sigma_{\R}}(\R)$.
\end{proposition}
\begin{proof} We have already explained that every $\K\insegeq {\sf{Lp}}^{\N, g, X, sts}(\R)$ has an iteration strategy in $\Gamma(\P, \Sigma)$. Moreover, because $\Sigma$ has strong branch condensation and is $\Gamma$-fullness preserving, the strategy $\K$ acquired from the strategy of $\M|\eta$ witnesses that $\K$ is a $\Sigma_{\R}$-sts\footnote{While this is non-trivial, most of the proof is contained in the proofs of  \rthm{sts fb constructions converge} and \rprop{summary of fb convergence}.}.
\end{proof}

We can now describe the $\N$-authenticated iterations of $\P$. 

\begin{definition}[$\N$-authenticated iteration]\label{certified iterations without fatal drops ii} Suppose $(\P, \Sigma)$ is an sts pair, $\Gamma$ is a projectively closed pointclass and $\N$ is an internally $\Sigma$-closed mouse with a fullness preserving iteration strategy $\Lambda$. Suppose further that $\Sigma$ has strong branch condensation and is strongly $\Gamma$-fullness preserving. Also suppose that $g\subseteq\mathbb{P}$ is $\N$-generic for some poset $\mathbb{P}\in\N$ and $\T\in \N[g]$ is a stack on $\P$. We say $\T$ is \textbf{$\N$-authenticated} if the following conditions holds.
\begin{enumerate}
\item For every $\a\in {\sf{max}}^\T$\footnote{See \rdef{the short tree component domain 1}.}, 
\begin{center}
${\sf{Lp}}^{\N, sts}(\M_\a^\T)\models ``\d^{\M_\a^\T}$ is a Woodin cardinal".
\end{center} 
\item For every $\a\in {\sf{max}}^\T$, $\pi^{\T_{<\a}, b}$ exists.
\item For all $\a\in R^\T$ such that $\pi^{\T_{\leq \a}, b}$ exists, letting $\W={\sf{nc}}^\T_\a$\footnote{See \rnot{notation for iteration trees}.}, if $\W$ is above $\ord((\M_\a^\T)^b)$ then for all limit ordinals $\gg<lh(\W)$ such that $\W\rest \gg$ is ${\sf{nuvs}}$, 
\begin{enumerate}
\item ${\sf{Lp}}^{\N, sts}(\m^+(\W\rest \gg))\models ``\d(\W\rest \gg)$ is not a Woodin cardinal", and
\item letting $b=[0, \gg)_\T$, $\Q(b, \W\rest \gg)$ exists and 
\begin{center}$\Q(b, \W\rest \gg)\insegeq {\sf{Lp}}^{\N, sts}(\m^+(\W\rest \gg))$.\end{center}
\end{enumerate}
\item For every $\a\in R^\T$ such that $\pi^{\T_{\leq \a}, b}$ exists, if ${\sf{nc}}^\T_\a$ is based on $\S=_{def}(\M_\a^\T)^b$ then $(\S, {\sf{nc}}^\T_\a)$ is a $(\P,\Sigma^\N)$-authenticated iteration\footnote{See \rdef{authentic iterations}.}.
\item For every $\a\in R^\T$ such that $\pi^{\T_{\leq \a}, b}$ exists, letting $\U={\sf{nc}}^\T_\a$ and $\S=_{def}\M_\a^\T$, if $\U$ is above $\d^{\S^b}$ and is such that for some $\eta\in (\d^{\S^b}, \d^\S)$, $\U$ is based on $\mathcal{O}^\S_{\eta, \S|\eta, \eta}$\footnote{See \rdef{the o stack}.} and is above $\eta$, then $(\mathcal{O}^\S_{\eta, \S|\eta, \eta}, \U)$ is a $(\P,\Sigma^\N)$-authenticated iteration.
\item For every $\a\in R^\T$ such that $\pi^{\T_{\leq \a}, b}$ exists, letting $\U={\sf{nc}}^\T_\a$ and $\S=_{def}\M_\a^\T$, if $\U$ is a normal tree on $\S^b$ above $\d^{\S^b}$, then $(\S^b, \T_{\geq\S})$ is a $(\P,\Sigma^\N)$-authenticated iteration\footnote{This clause follows from the one above it.}.
\end{enumerate}
$\myqedhere$
\end{definition}

%

\section{Generic interpretability in internally $\Sigma$-closed premice}\label{gen int sec}

In this section, we prove our main theorem, \rthm{main theorem on gen int}. As we said before, we will only prove clause 2. We start by fixing an sts hod pair $(\P, \Sigma)$ such that $\Sigma$ has strong branch condensation, a projectively closed pointclass $\Gamma$ such that $\Sigma$ is strongly $\Gamma$-fullness preserving and an internally $\Sigma$-closed premouse $\N$ such that $\N$ has a fullness preserving iteration strategy $\Lambda$\footnote{See \rdef{fullness preserving strategy for sigma mice}.}. We want to show that $\N$ is uniformly generically $\Sigma$-closed. 

Fix a poset $\mathbb{P}\in \N$ and an $\N$-generic $g\subseteq \mathbb{P}$. We start by defining a short tree  strategy $\Phi$ for $\P$. $\Phi$ is defined over $\N[g]$ in a uniform manner. Its domain consists of $\N$-authenticated iterations (see \rdef{certified iterations without fatal drops ii}). Given an $\N$-authenticated iteration $\T$ of limit length, we set
$\Phi(\T)=x$ if and only if one of the following conditions holds.
\begin{enumerate}
\item $\T$ is $\sf{nuvs}$ and letting $\a=\max(R^\T)$, ${\sf{Lp}}^{\N, sts}(\m^+(\T_{\geq \a}))\models ``\d(\T_{\geq \a})$ is a Woodin cardinal" and $x=\m^+(\T)$. 
\item Clause 1 above fails, $x\in \N[g]$ is a branch of $\T$ such that $\N[g]\models ``x$ is a cofinal well-founded branch of $\T$" and  $\T^\frown \{x\}$ is $\N$-authenticated.
\end{enumerate}
To complete the proof of \rthm{main theorem on gen int} we need to show that\\\\ 
(a) whenever $\T\in \dom(\Phi)\cap \dom(\Sigma)$, $\Phi(\T)$ is defined and is equal to $\Sigma(\T)$. \\\\
Fix then $\T'\in \N[g]$ such that $\T'\in \dom(\Phi)\cap \dom(\Sigma)$. Suppose first that \\\\
(*1) $\T'$ is $\sf{nuvs}$.\\\\
 Let $\a=\max(R^{\T'})$ and set $\T=\T'_{\geq \a}$ and $\S=_{def}\M^{\T'}_\a$. We thus have that $\mathcal{J}_{\omega}[\m^+(\T)]\models ``\d(\T)$ is a Woodin cardinal". We want to conclude that $\Phi(\T)$ is defined and $\Phi(\T)=\Sigma(\T)$.
 
 Suppose first that $\Sigma(\T)=\m^+(\T)$. Then because $\Sigma$ is $\Gamma$-fullness preserving, there is no $\Sigma_{\m^+(\T)}$-sts $\Q$ over $\m^+(\T)$ such that 
\begin{itemize}
\item $\Q$ is sound, 
\item $\Q$ has a strategy $\Lambda\in \Gamma(\P, \Sigma)$ witnessing that $\Q$ is a $\Sigma_{\m^+(\T)}$-sts over $\m^+(\T)$,
\item $\Q\models ``\d(\T)$ is a Woodin cardinal" but $\mathcal{J}_{\omega}[\Q]\models ``\d(\T)$ is not a Woodin cardinal".
\end{itemize} 
It then follows from \rprop{matching lp operators} that $\Phi(\T)=\m^+(\T)$.

Suppose next that $\Sigma(\T)=b$ where $b$ is a cofinal well-founded branch of $\T$. We thus have that $\Q(b, \T)=_{def}\W$ exists and want to conclude that $\Phi(\T)=b$. Notice that if 
\begin{center}
${\sf{Lp}}^{\N, sts}(\m^+(\T'_{\geq \a}))\models ``\d(\T'_{\geq \a})$ is not a Woodin cardinal"
\end{center}
then $\W\insegeq {\sf{Lp}}^{\N, sts}(\m^+(\T'_{\geq \a}))$ and, therefore, $\Phi(\T)=b$. Assume then that \\\\
(1) ${\sf{Lp}}^{\N, sts}(\m^+(\T'_{\geq \a}))\models ``\d(\T'_{\geq \a})$ is a Woodin cardinal".\\\\
 Set $\T=\T'_{\geq \S}$. The following claim finishes the proof of (a) assuming (*1).\\

\textit{Claim.} $\W\insegeq {\sf{Lp}}^{\N, sts}(\m^+(\T))$.\\\\
\begin{proof}
Recall from \rdef{sts hod pairs} that $\W$ has a strategy in $\Psi\in \Gamma^b(\P, \Sigma)$ witnessing that $\W$ is a $\Sigma_{\m^+(\T)}$-sts mouse over $\m^+(\T)$. Let $\k$ be an $\N$-cardinal such that $\{ \mathbb{P}, \T\}\in \N|\k[g]$. Using fullness preservation of $\Lambda$, fix an iteration tree $\U_0$ on $\N$ above $\k$ and according to $\Lambda$ with last model $\N_1$ such that $\pi^{\U_0}$ exists and there is an $\M\insegeq \N_1$ such that
\begin{enumerate}
\item $\M$ witnesses the internal $\Sigma$-closure of $\N_1$ at $\k$ and 
\item for some $\Q\inseg (\S^\M_\k)^b$, $w({\sf{Code}}(\Psi))\leq w({\sf{Code}}(\Sigma_{\Q}))$.
\end{enumerate}
Fix a real $x$ witnessing  $w({\sf{Code}}(\Psi))\leq w({\sf{Code}}(\Sigma_{\Q}))$. 

Let $\d_0<\d_1<\d_2$ be the first three Woodin cardinals of $\M$ that are greater than $\k$. Let $\U_1$ be an iteration tree on $\M$ based on $\M|\d_1$ according to $\Lambda_{\M, \U_0}$ and above $\d_0$ that is constructed according to the rules of $x$-genericity iteration. Let $\pi=\pi^{\U_1}$ and let $\M_2$ be the last model of $\U_1$. We then have that $x$ is generic for the extender algebra of $\M_2$ at $\pi(\d_1)$. It follows that \\\\
(2) $\Psi\rest (\M_2|\pi(\d_2))[g][x]\in \M_2[g][x]$\footnote{To make this conclusion, we use the fact that $\M_2[g][x]$ is closed under $\Sigma_\Q$. This can be established using the generic interpretability results of \cite{ATHM}.}.\\\\
Let $(( \M_\gg , \N_\gg : \gg\leq \nu), ( F_\gg: \gg<\nu), ( \T_\gg: \gg\leq\nu))$ be the output of the $(\M_2|\pi(\d_2), \pi(\delta_1))$-authenticated fully backgrounded construction over $\mathcal{J}_{\omega}[\m^+(\T)]$ based on $\m^+(\T)$\footnote{See \rdef{certified backgrounded constructions}.}. Next we have that.\\\\
(3) For some $\gg\leq \nu$, $\W=\M_\gg$.\\\\
(3) is a consequence of the following facts:\\\\
(3A) Extenders used in the construction of $\S$ have critical points $>\pi(\delta_1)$ (so the construction side doesn't move).\\
(3B) For each $\gg$, $\M_\gg$ is a $\Sigma_{\m^+(\T)}$-sts mouse over $\m^+(\T)$\footnote{See \rthm{sts fb constructions converge} and \rprop{summary of fb convergence}.}.\\
(3C) $\W$ side loses (because (2) implies $\W$ is $\pi(\d_2)$-iterable in  $\M_2|[g][x]$ and the construction is universal\footnote{See the universality clause of \rthm{existence of thick sets}.}).\\\\ 
Clearly (3) finishes the proof of the claim.
\end{proof}
We now assume the following:\\\\
(*2)  $\T'$ is $\sf{uvs}$. \\\\
Again, our goal is to show that $\Sigma(\T')=\Phi(\T')$. As many components of the proof are similar to the proofs of \rthm{sts fb constructions converge} and \rprop{next lemma}, we will not give the full proof. It is in fact enough to show the following:\\\\
(b) Suppose $\a\in R^{\T'}$ is such that $\pi^{\T', b}$ is defined and for all $\b\in R^{\T'}\cap \a$, $(\M^{\T'}_\b)^b \not =(\M^{\T'}_\a)^b$.  Let $\T=\T'_{\leq \a}$ and $\S=\M^{\T'}_\a$. There is then $\k<\ord(\N)$ and $\M\insegeq \N$ witnessing that $\N$ is internally $\Sigma$-closed at $\k$ and such that
\begin{enumerate}
\item $\{\mathbb{P}, \T\}\in \N|\k[g]$,
\item there is a normal stack $\U$ on $\S$ according to $\Sigma_\S$ such that either
\begin{enumerate}
\item $\U$ has a last model $\K\inseg \S_\k^\M$ or
\item $\U$ is of limit length and $\Sigma(\U)=\m^+(\U)=\S_\k^\M$.
\end{enumerate}
\end{enumerate} 
The tree $\U$ is built using the authentication procedure described in \rdef{authentic lsp}. \rprop{authenticated iteration construction} guarantees that the prescription for finding the branches of $\U$ as in clause 2 of \rdef{authentic lsp} produces a branch which is according to $\Sigma_\S$. The fact that $\N$ has a fullness preserving iteration strategy implies that $\S$ cannot outiterate $\S_\k^\M$.

\section{$S$-constructions}

Our definition of sts mice makes heavy use of the fact that the set $X$ is a self-well-ordered set. In particular, our definition cannot be used to define sts mice over $\mathbb{R}$. Another shortcoming of our definition is that it does not explain how to do $S$-constructions. In this short section, motivated by Section 3.38 of \cite{ATHM}, we indicate how to use \rthm{main theorem on gen int} to redefine hod mice in a way that one can define sts mice over $\mathbb{R}$ and perform $S$-constructions. 

Recall the difficulty with defining hybrid mice over $\bR$. In our definition, we always choose the least stack of some sort for which the branch has not been added and index a branch. Since $\bR$ may not be self-well-ordered, we do not have the luxury of choosing the least such stack. 

The problem with $S$-constructions is very similar. Suppose $(\P, \Sigma)$ is a hod pair or an sts hod pair and $N$ and $M$ are two transitive models of some fragment of set theory such that $\mathcal{J}_\omega(M)\subseteq \mathcal{J}_\omega(N)$ and for some poset $\mathbb{P}\in \mathcal{J}_\omega(M)$ and some $M$-generic $G\subseteq \mathbb{P}$, $\mathcal{J}_{\omega}(N)=\mathcal{J}_\omega(M)[G]$. Suppose further that both $M$ and $N$ are $\Sigma$-closed and $\P\in N\cap M$. For us, $S$-constructions are constructions that translate $\Sigma$-mice over $N$ to $\Sigma$-mice over $M$. For more details consult Section 3.38 of \cite{ATHM}.\footnote{In \cite{Selfiter}, this process is called $P$-constructions.} 

The difficulty in performing $S$-constructions is the following. Suppose $\N$ is a $\Sigma$-mouse over $N$, and we want to translate it onto a $\Sigma$-mouse over $M$. Suppose our translation has produced a $\Sigma$-mouse $\M$ over $M$, and our indexing scheme demands that a branch of some stack $\T\in \N$ be indexed in the very next step in the translation procedure. The problem is that $\T$ may not be a stack in $\M$ nor may it be the stack whose branch is indexed in $\M$. 

To solve this problem, we change the definition of hybrid premouse in a way that the iterations whose branches are indexed do not depend on generic extensions. In particular, instead of indexing iterations according to $\Sigma$, we considered generic genericity iterations on $\M_1^{\#, \Sigma}$. Such iterations make levels of the model generically generic and do not depend on generic extensions. This move solves both problems. In the first case what is important is that the indexed iterations do not depend on the well-ordering of the model, and in the second case what is important is that the indexed iterations do not depend on generic extensions. For more details consult Definition 3.37 of \cite{ATHM} or \cite{trang2013} for a similar construction.

Here our solution is similar. Suppose $(\P, \Sigma)$ is an sts hod pair such that $\Sigma$ has strong branch condensation and is $\Gamma$-fullness preserving for some projectively closed pointclass $\Gamma$ and $\M$ is an $\Sigma$-sts mouse over some set $X$ such that $\P\in X$. Then the iterations of $\P$ that are indexed in $\M$ are of the form $t=(\P, \T, \Q, \U)$, where $t$ is an indexable stack\footnote{See \rdef{indexable stack}.}. $\T$ is always the result of comparing $\P$ with a certain backgrounded construction. Notice that this neither depends on the well-ordering of $\M$ nor on small generic extensions. $\U$ is a stack on $\Q^b$ and, in \rdef{weak psi alpha indexing scheme a}, we chose the least such stack. Thus the choice of $\U$ depends on both the well-ordering of $\M$ and small generic extensions\footnote{Small in the sense that the generic is smaller than the critical point of the first background extender used in the construction.}. To solve the issue, we will start considering stacks $s=(\P, \T, \Q, \U)$ where $\T$ is as before but now $\U$ is a generic genericity iteration on $\M_2^{\#, \Sigma_{\Q^b}}$ to make a level of the model generically generic. We only consider such generic genericity iterations of $\M_2^{\#, \Sigma_{\Q^b}}$ that are based on the first Woodin of $\M_2^{\#, \Sigma_{\Q^b}}$. 


The reason we choose $\M_2^{\#, \Sigma_{\Q^b}}$ is that we want to use clause 1 of \rthm{main theorem on gen int}. It is not hard to see that if $\d_0<\d_1$ are the first two Woodin cardinals of $\M_2^{\#, \Sigma_{\Q^b}}$ and $g\subseteq Coll(\omega, \d_0)$ is $\M_2^{\#, \Sigma_{\Q^b}}$-generic then $\M_2^{\#, \Sigma_{\Q^b}}|\d_1[g]$ is internally $\Sigma_{\Q^b}$-closed. Clause 1 of \rthm{main theorem on gen int} is a weaker result than \cite[Lemma 3.35]{ATHM}, which is what is used to reorganize hod mice in \cite{ATHM}. We could prove an equivalent of \cite[Lemma 3.35]{ATHM}, but doing this is much harder than proving clause 1 of \rthm{main theorem on gen int}.

To show that the resulting structure $\M$ is closed under $\Sigma$, we will first show that we can find branches of indexable stacks. Given such a stack $t=(\P, \T, \Q, \U)$ let $\W$ be an iteration of $\M_2^{\#, \Sigma_{\Q^b}}$ such that $(\P, \T, \Q, \W)$ is indexed in $\M$ and if $\S$ is the last model of $\W$ then $\U$ is generic over $\S$ for $\mathbb{B}^\S_{\d}$ where $\d$ is the least Woodin of $\S$ and $\mathbb{B}^\S_{\d}$ is the extender algebra of $\S$ at $\d$. It then follows from \rthm{main theorem on gen int} that $\Sigma_{\Q^b}\rest \S|\eta[\U]\in \S[\U]$ where $\eta$ is the second Woodin cardinal of $\S$. The rest of the proof is just repeating the proof of \rthm{main theorem on gen int}. 

Instead of re-developing the entire theory of sts mice, we will simply give the definition of indexable stack. The rest of the definitions, those developed in \rsec{stsis sec}, \rsec{authentic iterations and finite stacks sec}, \rsec{sec short tree strategy mice}, \rsec{hp indexing scheme:sec} and \rsec{hod mice sec}, stay more or less the same. It is important to note that the theory of sts mice as well as hod mice does not in general depend on particular indexing schemes. 

\begin{definition}[Revised Indexable stack]\label{revised indexable stack}\index{revised indexable stack} Suppose $\P$ is a hod-like $\#$-lsa type ${\sf{lses}}$\footnote{See \rdef{lsa type}.}. We say that an st-stack\footnote{See \rdef{st-stack}.} 
\begin{center}
$\T=((\M_\a)_{\a<\eta}, (E_\a)_{\a<\eta-1}, D, R, (\beta_\a, m_\a)_{\a\in R}, \so, \ma, T)$
\end{center}
 is a \textbf{revised indexable stack} on $\P$ if one of the following clauses hold:
\begin{enumerate}
\item ${\sf{max}}=\emptyset$ and there is $\a\in R^\T$ such that $\pi^{\T_{\leq \a}, b}$ is defined and letting\footnote{$\M_2^{\#, \Sigma_{(\M_\a^\T)^b}})^{\M_\a^\T}$ is $(\M_2^{\#, \Sigma_{(\M_\a^\T)^b}}$ in the sense of $\M_\a^\T$.}   $\M=(\M_2^{\#, \Sigma_{(\M_\a^\T)^b}})^{\M_\a^\T}$, $\T_{\geq \a}$ is a normal stack on $\M$ that is above $\ord((\M_\a^\T)^b)$ and is based on $\M|\d$ where $\d$ is the least Woodin cardinal of $\M$.
\item $\card{{\sf{max}}}=1$, $\T$ is a normal stack\footnote{See \rdef{normal stack}.} and if $\a$ is the unique element of ${\sf{max}}$ then $\pi^\T_{0, \a}$ is defined and ${\sf{next}}^\T(\a)=\lh(\T)$\footnote{It follows that $\T_{\geq \a}$ is above $\pi^{\T}_{0, \a}(\d^{\P^b})$. See also \rnot{notation for iteration trees}.}.
\end{enumerate} 
$\myqedhere$
\end{definition}

 We say $\P$ is a revised hod premouse if it is indexed according to our revised indexing scheme, which will only index \textit{authentic revised indexable stacks}\footnote{See \rsec{authentic iterations and finite stacks sec}.}. We say $(\P, \Sigma)$ is revised hod pair if $\P$ is revised hod premouse and $\Sigma$ is an iteration strategy for $\P$.

\begin{theorem}\label{closure of revised hod premouse} Suppose $(\P, \Sigma)$ is a revised hod premouse such that $\Sigma$ is strongly $\Gamma$-fullness preserving for some projectively closed pointclass $\Gamma$ and $\Sigma$ has strong branch condensation. Then for any $\Q\in Y^\P$ and $\P$-generic $g$, 
\begin{enumerate}
\item if $\Q$ is not of lsa type then $\Sigma_\Q\rest \P[g]$ is uniformly in $\Q$ definable over $\P[g]$, and
\item if $\Q$ is of lsa type then the fragment of $\Sigma_\Q^{stc}\rest \P[g]$ that acts on revised indexable stacks is uniformly in $\Q$ definable over $\P[g]$.
\end{enumerate}
\end{theorem}

We now just import our lemmas on $S$-construction from Section 3.8 of \cite{ATHM} to our current context. Let $(\P, \Sigma)$ be a hod pair or an sts pair such that $\Sigma$ has the strong branch condensation and is strongly $\Gamma$-fullness preserving for some pointclass $\Gamma$. Suppose $\M$ is a sound $\Sigma$-mouse and $\d$ is a cutpoint cardinal of $\M$. Suppose further that $\N\in \M|\d+1$ is such that $\d\subseteq \N\subseteq H_\d^\M$, $\N$ models a sufficiently strong fragment of ${\sf{ZF}}-{\sf{Replacement}}$, $\N$ is a $\Sigma$-mouse or a $\Sigma$-sts mouse and there is a partial ordering $\mathbb{P}\in L_{\omega}[\N]$ such that $\M|\d$ is $\mathbb{P}$-generic over $L_{\omega}[\N]$. We would like to define $S$-construction of $\M$ over $\N$ relative to $\Sigma$.

\begin{definition}\label{s-construction}
An $S$-construction of $\M$ over $\N$ relative to $\Sigma$ is a sequence $( \S_\a, \bar{\S}_\a: \a\leq \eta)$ of $\Sigma$-mice over $\N$ such that
\begin{enumerate}
\item $\S_0=L_{\omega}[\N]$,
\item if $\M|\d$ is generic over $\bar{\S}_\a$ for a forcing in $L_{\omega}[\N]$ then
\begin{enumerate}
\item if $\M||(\omega \cdot \a)$ is active and has a last branch $b$ then $\S_\a$ is the expansion of $\bar{\S}_\a$ by $b$ and  $\bar{\S}_{\a+1}=rud(\S_\a)$.
\item if $\M||(\omega \cdot \a)$ is active and has a last extender $E$ then $\S_\a$ is the expansion of $\bar{\S}_\a$ by $E$ and $\bar{\S}_{\a+1}=rud(\S_\a)$,
\item if $\M||(\omega\times \a)$ is passive then $\S_\a=\bar{\S}_\a$ and $\bar{\S}_{\a+1}=rud(\S_\a)$,
\end{enumerate}
\item if $\l$ is limit then $\bar{\S}_\l=\bigcup_{\a<\l}\S_\a$.
\end{enumerate}
$\myqedhere$
\end{definition}

The following is the restatement of Lemma 3.42 of \cite{ATHM}.

\begin{lemma}\label{s-construction lemma}
Suppose $(\P, \Sigma)$, $\M, \N$  are as above and $\d$ is a strong cutpoint cardinal of $\M$.  Suppose further that $\N\in \M|\d+1$ is such that $\d\subseteq \N\subseteq H_\d^\M$ and there is a partial ordering $\mathbb{P}\in L_{\omega}[\N]$ such that whenever $\Q$ is a $\Sigma$-mouse over $\N$ such that $H_\d^\Q=\N$ then $\M|\d$ is $\mathbb{P}$-generic over $\Q$. Then there is a $\Sigma$-mouse $\S$ over $\N$ such that $\M|\d$ is generic over $\S$ and $\S[\M|\d]=\M$.
\end{lemma}

The following is just the restatement of Lemma 3.43 of \cite{ATHM}.

\begin{lemma}\label{s-constructions give hybrid q-structures} Suppose $(\P, \Sigma)$, $\M$ and $\N$ are as above. Suppose further that $\M\models {\sf{ZFC}}-{\sf{Replacement}}$ is a $\Sigma$-mouse and $\eta$ is a strong cutpoint non-Woodin cardinal of $\M$. Suppose $\gg>\eta$ is a cardinal of $\M$ and $\N=(\mathcal{J}^{\vec{E}, \Sigma})^{\M|\gg}$. Suppose $\mathcal{J}_{\omega}(\N|\eta)\models ``\eta$ is Woodin". Let $( \S_\a, \bar{\S}_\a :\a<\nu)$ be the $\S$-construction of $\M|(\eta^+)^\M$ over $\N|\eta$ relative to $\Sigma$. Then for some $\a<\nu$, $\S_\a\models ``\eta$ isn't Woodin".
\end{lemma}

\chapter{Analysis of $\H$}

In this chapter we analyze $V_\Theta^\H$ of the minimal model of the Largest Suslin Axiom. The analysis is very much like the analysis of $V_\Theta^\H$ in the minimal model of $AD^++\theta_1=\Theta$, which appeared in \cite[Chapter 4]{ATHM}. Just like in \cite[Chapter 4]{ATHM}, we need to introduce the notion of suitable pair, $B$-iterable pair and etc. The proof of \rthm{computation of hod} is very much like the proof of \cite[Theorem 4.24]{ATHM}.

\section{$B$-iterability}

In this section, we import $B$-iterability technology to our current context. Most of what we will need was laid out in  \cite[Section 4.1 and Section 4.2]{ATHM}. Here we will only sketch the necessary arguments.

\begin{definition}[Suitable pair] \label{suitable pair}
 $(\P, \Sigma)$ is a \emph{suitable pair} if the following clauses hold:
\begin{enumerate}
 \item Either $\P$ is a hod premouse of successor type or $\P$ is a $\#$-like lsa type hod premouse\footnote{See \rdef{lsa type}. Thus, in both cases $\P\models ``\d^\P$ is a Woodin cardinal". Also, if $\P$ is not of lsa type then $\P$ is a $\Sigma_{\P^-}$-mouse above $\P^-$.}.
 \item If $\P$ is not of lsa type then 
 \begin{itemize}
 \item $(\P, \Sigma)$ is a pre-hod pair\footnote{See \rdef{prehod pair}.},
 \item  $\Sigma$ has strong branch condensation and is strongly fullness preserving\footnote{Thus, $\powerset(\bR)$-fullness preserving.},
 \item For any $\P$-cardinal $\eta>\d_{\l-1}^\P$, if $\eta$ is a strong cutpoint
     then $\P|(\eta^+)^\P={\sf{Lp}}^{\Sigma}(\P|\eta)$.
     \end{itemize}
 \item If $\P$ is of lsa type then $(\P, \Sigma)$ is an sts hod pair such that $\Sigma$ has strong branch condensation and is strongly fullness preserving\footnote{In this fullness preservation implies that ${\sf{Lp}}^{\Sigma}(\P)\models ``\d^\P$ is a Woodin cardinal".}.
\end{enumerate}
$\myqedhere$
\end{definition}

For convenience, we extend the notation $\P^-$ to lsa type (see \rnot{l p}).

\begin{notation}\label{p- notation} Suppose $\P$\footnote{As was decided many pages before, we develop the theory of hod mice over $\emptyset$. Thus, $X^\P=\emptyset$.} is either of lsa type or of successor type. We then let 
\begin{center}
$\P^-=\begin{cases}
\P&:\text{$\P$ is of lsa type}\\
\bigcup_{\Q\inseg_{hod}\P}\Q&:\text{otherwise}.
\end{cases}$\end{center}
Also, if $(\P, \Sigma)$ is a suitable sts pair then we let ${\sf{lp}}(\P, \Sigma)={\sf{Lp}}_\omega^{\Sigma}(\P)$. $\myqedhere$
\end{notation}

Suppose $(\P, \Sigma)$ and $(\Q, \Lambda)$ are hod pairs or sts hod pairs such that $\Sigma$ and $\Lambda$ have strong branch condensation and are strongly fullness preserving.  We then let 
\begin{center}
$(\P, \Sigma)\leq_{DJ} (\Q, \Lambda)$
\end{center}
 if and only if $(\P, \Sigma)$ loses the coiteration with $(\Q, \Lambda)$. Notice that $\leq_{DJ}$ is a well-founded relation. We then let $\a(\P, \Sigma)=\card{(\P, \Sigma)}_{\leq_{DJ}}$\index{$\a(\P, \Sigma)$}, and we let $[\P, \Sigma]$\index{$[\P, \Sigma]$} be the $=_{DJ}$\index{$=_{DJ}$} equivalence class of $(\P, \Sigma)$, i.e.,
\begin{center}
$(\Q, \Lambda)\in [\P, \Sigma]$ iff $(\Q, \Lambda)$ is a hod pair such that $\Lambda$ has branch condensation and is strongly fullness preserving and $\a(\Q, \Lambda)=\a(\P, \Sigma)$.
\end{center}
Notice that $[\P, \Sigma]$ is independent of $(\P, \Sigma)$. We let
\begin{center}
$\mathbb{B}(\P, \Sigma)=\{B\subseteq [\P, \Sigma] \times \mathbb{R} : B$ is $OD\}$.
\end{center}
Note that $\mathbb{B}(\P, \Sigma)$ is defined for hod pairs or sts hod pairs, but not for suitable pairs that are not sts hod pairs\footnote{This is because if $(\P, \Sigma)$ is a non-sts suitable pair then $\Sigma$ does not act on $\P$ but only on $\P^-$. For such pairs, comparison is somewhat meaningless.}. 

 The following standard lemma features prominently in our computations of $\H$. The proof is very much like the proof of Lemma 4.16 of \cite{ATHM}. Below $\sf{SMC}$ stands for \textit{Strong Mouse Capturing}\index{SMC}\index{strong mouse capturing}. More precisely, $\sf{SMC}$ states that for any hod pair or sts hod pair $(\P, \Sigma)$ such that $\Sigma$ is strongly fullness preserving and has strong branch condensation then for any $x, y\in \bR$, $x\in OD_{y, \Sigma}$ if and only if $x\in {\sf{Lp}}^{\Sigma}(y)$.

\begin{lemma}\label{term capturing} Assume $\sf{SMC}$ and suppose $(\P, \Sigma)$ is a suitable pair and $B\in \mathbb{B}(\P^-, \Sigma)$. Set
\begin{center}
$\P_+=\begin{cases}
\P &: \P\ \text{is of successor type}\\
{\sf{lp}}(\P, \Sigma) &:\ \text{otherwise}.
\end{cases}$
\end{center}
 Suppose $\k$ is a $\P_+$-cardinal such that if $\P$ is of lsa type then for some $n>0$, $\k=((\d^\P)^{+n})^{\P_+}$ and otherwise $\k>\d^{\P^-}$. Then there is $\tau\in \P_+^{Coll(\omega, \k)}$ such that $(\P_+, \tau)$ locally term captures $B_{(\P, \Sigma)}$ at $\k$ for a comeager set of $\P_+$-genetics $g\subseteq Coll(\omega,\kappa)$.
\end{lemma}

If $B$ is locally term captured for comeager many set generics over a suitable pair $(\P, \Sigma)$ then we let $\tau_{B, \k}^{\P, \Sigma}$ be the invariant term in $\P_+$ locally term capturing $B$ at $\k$ for comeager many set generics. One way to get term capturing for all generics is to show that a suitable pair can be extended to a structure that has one more Woodin.

\begin{definition}[$n$-Suitable pair] \label{n-suitable pair}
 $(\P, \Sigma)$ is an \emph{$n$-suitable pair} if there is $\d$ such that 
 \begin{center}
 $\P\models ``\d$ is a Woodin cardinal"
 \end{center}
 and the following clauses hold:
\begin{enumerate}
\item Either $(\P|(\d^{+\omega})^\P, \Sigma)$ is a suitable pair or letting $\a=\min(\dom(\vec{E}^\P)-\d)$, $(\P|\a, \Sigma)$ is suitable\footnote{Because $\P|(\d^{+\omega})^\P$ has infinitely many cardinals above $\d$, in the first case $\P|(\d^{+\omega})^\P$ is of successor type and in the second, $\P|\a$ is of an lsa type.}. Set 
\begin{center}
$\P_0=\begin{cases}
\P|(\d^{+\omega})^\P &: (\P|(\d^{+\omega})^\P, \Sigma)\ \text{is a suitable pair}\\
\P|\a&:\ \text{otherwise}
\end{cases}$\end{center}
 \item If $\P_0$ is of lsa type then $\P$ is a $\Sigma$-sts premouse over $\P_0$ and otherwise $\P$ is a $\Sigma$-premouse over $\P_0$,
 \item $\P\models {\sf{ZFC}}-{\sf{Replacement}}+ ``$there are exactly $n$ Woodin cardinals, $\eta_0<\eta_1<...<\eta_{n-1}$ that are  strictly greater than $\d"$,
 \item $\ord(\P)=\sup_{i<\omega}(\eta_{n-1}^{+i})^\P$ (here we set $\eta_{-1}=\delta$),
 \item For any $\P$-cardinal $\eta\geq \d$, if $\eta$ is a strong cutpoint
     then $\P|(\eta^+)^\P={\sf{Lp}}^{\Sigma}(\P|\eta')$ where $\eta'=\min(\dom(\vec{E}^\P)-\eta)$.
 \end{enumerate}
 We say $\P$ is of lsa type if $\P_0$ is of lsa type. Otherwise we say that $\P$ is of successor type. $\myqedhere$
 
\end{definition}
%

If $(\P, \Sigma)$ is $n$-suitable then we let $\d^\P$ be the $\d$ of \rdef{n-suitable pair} and $\P_0$ be as in \rdef{n-suitable pair}.
 Clearly $0$-suitable pair is just a suitable pair. The following are easy consequences of \rlem{term capturing}. 

\begin{lemma}\label{more term capturing} Assume $\sf{SMC}$. Suppose $(\P, \Sigma)$ is an $n$-suitable pair and $B\in \mathbb{B}(\P^-, \Sigma)$. Suppose $\k$ is a $\P$-cardinal such that 
\begin{itemize}
\item if $\P$ is of lsa type then $\k>((\d^\P)^+)^\P$ and 
\item otherwise $\k>\d^{\P_0^-}$.
\end{itemize}
 Then there is $\tau\in \P^{Coll(\omega, \k)}$ such that $(\P, \tau)$ locally term captures $B_{(\P, \Sigma)}$ at $\k$ for comeager set of $\P$-generics $g\subseteq Coll(\omega,\kappa)$.
\end{lemma}

\begin{corollary}\label{term capturing for all generics}
Assume $\sf{SMC}$. Suppose $(\P, \Sigma)$ is an $n$-suitable pair and $B\in \mathbb{B}(\P^-, \Sigma)$. Let $\nu=((\d^\P)^{+\omega})^\P$. 
Suppose $\k$ is a $\P$-cardinal such that 
\begin{itemize}
\item if $\P$ is of lsa type then $\k\in(((\d^\P)^+)^\P, \nu)$ and 
\item otherwise $\k\in (\d^{\P_0^-}, \nu)$.
\end{itemize}
 Then $(\P|\nu, \tau_{B, \k}^{\P, \Sigma})$ locally term captures $B_{(\P, \Sigma)}$ at $\k$ for all $\P$-generic $g\subseteq Coll(\omega,\kappa)$.
\end{corollary}

\rcor{term capturing for all generics} is our main method of showing that various $B$ are term captured over the hod mice that we will construct.  Suppose now that $(\P, \Sigma)$ is a hod pair. It is now a trivial matter to import the terminology of \cite[Section 4.1]{ATHM} to our current context. We will have that $S(\Sigma)$ consists of those $\Q$ such that $\Q_0\in pI(\P, \Sigma)$ and $(\Q, \Sigma_{\Q_0})$ is a suitable pair. Given $\Q\in S(\Sigma)$, we let  $f_B(\Q)=\oplus_{\k<\ord(\Q)}\tau_{B, \k}^{\Q, \Sigma_{\Q_0}}$\footnote{Here $\Q_0$ is defined in \rdef{n-suitable pair}.}. Then the rest of the notions are defined for $F=\{ f_B: B\in \mathbb{B}(\P, \Sigma)\}$. Therefore, in the sequel, we will freely use the terminology of  \cite[Section 4.1]{ATHM}. 

\section{The computation of $\H$}\label{computation of hod sec}

Throughout this section we assume $\sf{AD}^++\sf{SMC}$ and let $\la \theta_\a: \a\leq \Omega\ra$ be the Solovay sequence. Our goal is to compute $V_{\theta_\b}^{\H}$ for $\b\leq\Omega$. We will do it under some additional hypothesis described below. In the next few chapters, we will prove that our additional hypothesis follows from $\sf{AD}^++``$No initial segment of the Solovay sequence satisfies $\sf{LSA}$". 

Suppose $(\P, \Sigma)$ is a hod pair or an sts pair such that $\Sigma$ has strong branch condensation and is strongly fullness preserving. We will continue using the notation $\a(\P, \Sigma)$, $\P_0$ and $\P^-$ from the previous section\footnote{See \rdef{p- notation}.}. 

Suppose first that $\b+1=\Omega$. We then let $\mathcal{I}=\{ (\Q, \Lambda, \vec{B}) :$
\begin{enumerate}
\item $(\Q, \Lambda)$ is suitable, $\Lambda$ is strongly fullness preserving and has strong branch condensation, and $\a(\Q^-, \Lambda)=\b$,
\item for some integer $n$, $\vec{B}=(B_0, ..., B_n)$ and for every $i<n$, $B_i\in \mathbb{B}(\Q^-, \Lambda)$, and
\item $(\Q, \Lambda)$ is strongly $\vec{B}$-iterable $\}$.
\end{enumerate}
 \rthm{getting strong b-iterability} and the results of \rsec{generation of mouse full pointclass sec} show that $\mathcal{I}\not=\emptyset$. Define $\preceq$ on $\mathcal{I}$ by
\begin{center}
$(\P, \Sigma, \vec{B}) \preceq (\Q, \Lambda, \vec{C}) \iff \vec{B}\subseteq \vec{C}$ and $(\Q, \Lambda, \vec{B})$ is a $\vec{B}$-tail of $(\P, \Sigma, \vec{B})$.
\end{center}
When $(\R, \Psi, \vec{B}) \preceq (\Q, \Lambda, \vec{C})$, there is a canonical map 
\begin{center}
$\pi: H_{\vec{B}}^{\R, \Psi} \rightarrow H_{\vec{B}}^{\Q, \Lambda}$,
\end{center}
 which is independent of $\vec{B}$-iterable branches. We let $\pi_{(\R, \Psi, \vec{B}), (\Q, \Lambda, \vec{B})}$ be this map. We then have that $(\mathcal{I}, \preceq)$ is directed. Let
\begin{center}
$\mathcal{F}=\{ H_{\vec{B}}^{\Q, \Lambda}: (\Q, \Lambda, \vec{B})\in \mathcal{I}\}$.
\end{center}
and also let $\M_{\infty}$ be the direct limit of $\mathcal{F}$ under the iteration maps $\pi_{(\R, \Psi, \vec{B}), (\Q, \Lambda, \vec{B})}$. Let $\d_\infty=\d^{\M_\infty}$. For $(\Q, \Lambda, B)\in \mathcal{I}$, we let $\pi_{(\Q, \Lambda, B), \infty}: H^{\Q, \Lambda}_{B} \rightarrow \M_\infty$. Standard arguments show that $\M_\infty$ is well-founded.

 Following \cite[Section 4.4]{ATHM}, we let $\phi$ be the following sentence: for every $\b+1<\Omega$ there is a hod pair $(\P, \Sigma)$ such that 
   \begin{enumerate}
   \item $\P$ is of successor type,
   \item $\a(\P^-, \Sigma_{\P^-})=\b$,
   \item $\Sigma$ is strongly fullness preserving and has strong branch condensation,
     \item for any $\Q\in pI(\P, \Sigma)\cup pB(\P, \Sigma)$, if $\Q$ is of successor type then
   \begin{enumerate}
   \item there is a sequence $\la B_i: i<\omega\ra\in \mathbb{B}(\Q^-, \Sigma_{\Q^-})^\omega$ which guides $\Sigma_\Q$
   and
    \item for any $B\in\mathbb{B}(\Q^-, \Sigma_{\Q^-})$ there is $\R\in pI(\Q, \Sigma_{\Q})$ such that $\Sigma_{\R}$ respects $B$.
      \end{enumerate}
  \end{enumerate}

Our additional hypothesis, $\psi$, is the conjunction of $\phi$ with the following statement:
 If $\Omega=\b+1$ then there is a suitable $\emptyset$-iterable $(\P, \Sigma)$ such that 
    \begin{enumerate}
    \item $\a(\P^-, \Sigma_{\P^-})=\b$ and $\Sigma_{\P^-}$ is strongly fullness preserving and has strong branch condensation,
    \item for any $B\in \mathbb{B}(\P^-, \Sigma_{\P^-})$ there is an $\emptyset$-iterate $(\Q, \Phi)$ of $(\P, \Sigma)$ such that $(\Q, \Phi)$ is strongly $B$-iterable.
    \item $\M_\infty$ is well-founded and $\d_\infty=\Theta=\theta_{\b+1}$.
\end{enumerate}

We will use the following lemma to establish $\psi$.  It can be proved exactly the same way as \cite[Lemma 4.23]{ATHM}. 
\begin{lemma}\label{proving psi} Suppose $\Gamma\subseteq \powerset(\mathbb{R})$ is such that 
\begin{center}
$L(\Gamma, \mathbb{R})\models {\sf{AD}^+}+{\sf{SMC}}+\Omega=\b+1$ and $\Gamma=\powerset(\mathbb{R})\cap L(\Gamma, \mathbb{R})$.
\end{center} Suppose $\Gamma^*\subseteq \powerset(\mathbb{R})$ is such that $\Gamma\subseteq \Gamma^*$, $L(\Gamma^*, \mathbb{R})\models \sf{AD}^+$ and there is a hod a pair $(\P, \Sigma)\in \Gamma^*$ such that the following hold.
\begin{enumerate}
\item $\Sigma$ has strong branch condensation and is strongly $\Gamma$-fullness preserving.
\item Either $\P$ is of successor type or of lsa type. 
\item ${\sf{Code}}(\Sigma_{\P^-})\in \Gamma$ and
\begin{enumerate}
\item if $\P$ is of successor type then $L(\Gamma, \mathbb{R})\models ``(\P, \Sigma_{\P^-})$ is a suitable pair such that $\a(\P^-, \Sigma_{\P^-})=\b$" and 
\item  if $\P$ is of lsa type then $L(\Gamma, \mathbb{R})\models ``(\P, \Sigma_{\P}^{stc})$ is a suitable pair such that $\a(\P^-, \Sigma_{\P^-})=\b$".
\end{enumerate}
Letting 
\begin{center}
$\Lambda=\begin{cases}
\Sigma_{\P^-} &: \P\ \text{is of successor type}\\
\Sigma^{stc}&:\ \text{otherwise}
\end{cases}$
\end{center}
the following clauses hold:
\item There is a sequence $\la B_i: i<\omega\ra\in [(\mathbb{B}(\P^-, \Lambda))^{L(\Gamma, \mathbb{R})}]^{\omega}$ guiding $\Sigma$.
\item For any $B\in(\mathbb{B}(\P^-, \Lambda))^{L(\Gamma, \mathbb{R})}$ there is $\R\in pI(\P, \Sigma)$ such that $\Sigma_{\R}$ respects $B$.
\end{enumerate}
Then $L(\Gamma, \mathbb{\R})\models \psi$ and $\M_\infty^{L(\Gamma, \bR)}=\M^+_\infty(\P, \Sigma)$.\footnote{Recall that $\M^+_\infty(\P, \Sigma)$ is the direct limit of all $\Sigma$-iterates of $\P$}
\end{lemma}

The next theorem is the adaptation of \cite[Theorem 2.24]{ATHM} to our current context. It can be proved via exactly the same proof. Because of this, we omit the proof.

\begin{theorem}[Computation of $\H$]\label{computation of hod} Assume $\sf{AD}^+$. Suppose $\Gamma\subseteq \powerset(\mathbb{R})$ is such that $\Gamma=\powerset(\mathbb{R})\cap L(\Gamma, \mathbb{R})$. Set $W=L(\Gamma, \mathbb{R})$ and let $(\theta_\b: \b\leq \Omega)$ be the Solovay sequence of $W$. Then the following holds:
\begin{enumerate}
\item Suppose $W\models \phi$ and $\b+1<\Omega$. Let $(\P, \Sigma)$ witness $\phi$ for $\b$. Then letting $\M=\M_\infty^+(\P, \Sigma)$, $\vec{E}=\vec{E}^\M$ and $\Lambda=\Sigma^\M$, for every $\a\leq \b$ there $\Q\inseg_{hod} \M$ such that 
\begin{enumerate}
\item $\d^\Q=\theta_\a$,
\item $\d^\Q$ is either a Woodin cardinal of $\M$ or a limit of Woodin cardinals of $\M$, and
\item $\M|\theta_\a=(V_{\theta_\a}^{\H^W}, \vec{E}\rest \theta_\a, \Lambda\rest V_{\theta_\a}^{\H^W}, \in)$.
\end{enumerate}
\item If $W\models \psi$ then letting $\M=\M_\infty^{W}$
$\vec{E}=\vec{E}^\M$ and $\Lambda=\Sigma^\M$, for every $\a\leq \Omega$ there is $\Q\insegeq_{hod} \M$ such that 
\begin{enumerate}
\item $\d_\Q=\theta_\a$,
\item $\d^\Q$ is either a Woodin cardinal of $\M$ or a limit of Woodin cardinals of $\M$, and
\item $\M|\theta_\a=(V_{\theta_\a}^{\H^W}, \vec{E}\rest \theta_\a, \Lambda\rest V_{\theta_\a}^{\H^W}, \in)$.
\end{enumerate}
\end{enumerate}
\end{theorem}

Thus, working in a model of $\sf{AD}^+$, if $\a<\Omega$ then to compute $\H|\theta_\a$ we only need to produce a hod pair $(\P, \Sigma)$ satisfying \rlem{proving psi}. In the next chapter, in particular in \rthm{getting strong b-iterability} and \rsec{generation of mouse full pointclass sec}, we will show that this is indeed true in the minimal model of the Largest Suslin Axiom.


\chapter{Models of ${\sf{LSA}}$ as derived models}\label{chap:LSA_DM}

In this chapter, we show that certain derived models satisfy the ${\sf{LSA}}$. We also prove results that are important elsewhere. The results of \rsec{generation of mouse full pointclass sec} and \rthm{getting strong b-iterability} are needed to carry out the computation of $\H$ (see \rthm{computation of hod}). We start with introducing the pointclass $\Gamma(\P, \Sigma)$ where $(\P, \Sigma)$ is an sts hod pair. 

\section{$\Gamma(\P, \Sigma)$ revisited}\label{gamma(p, sigma) in the successor case sec}

In this section, we translate the results of \cite[Section 5.6]{ATHM} to our current context. Suppose $(\P, \Sigma)$ is a hod pair such that either $\P$ is of  successor type or of $\#$-lsa type\footnote{See the discussion after \rdef{lsa type}.} and $\Sigma$ is strongly fullness preserving and has strong branch condensation. Recall the notation $\P^-$.

Suppose first that $\P$ is of successor type. We now generalize the results of \cite[Section 5.6]{ATHM}. Recall the notation ${\sf{Mice}}_{\Sigma}$ (see \rnot{mice relative to gamma}). Because $\P$ is not of lsa type, it follows that ${\sf{Code}}(\Sigma)$ is Suslin, co-Suslin (this can be proved using the proof of \cite[Lemma 5.9]{ATHM}). It follows that there is a scaled pointclass closed under continuous images and pre-images and under $\exists^\bR$, and also contains ${\sf{Mice}}_{\Sigma_{\P^-}}$. We then let $\Gamma_\Sigma^*$ be the least such pointclass.  Also, let 
\begin{center}
$\Gamma_\Sigma=(\Sigma^2_1({\sf{Code}}(\Sigma_{\P^-})))^{L({\sf{Mice}}_{\Sigma_{\P^-}}, \mathbb{R})}$.
\end{center}
Notice that $\Gamma_\Sigma$ is a lightface good pointclass, and so we set
\begin{center}
$\utilde{\Gamma}_{\Sigma}=\cup_{x\in \bR}(\Sigma^2_1({\sf{Code}}(\Sigma_{\P^-}), x))^{L({\sf{Mice}}_{\Sigma_{\P^-}}, \mathbb{R})}$.
\end{center}
 Also ${\sf{Mice}}_{\Sigma_{\P^-}}$ belongs to $\Gamma_\Sigma$ and is a universal $\Gamma_\Sigma$ set. We let
\begin{center}
$\Gamma(\P, \Sigma)=\{ A:$ for cone of $x\in \mathbb{R}$, $A\cap C_{\Gamma_\Sigma}(x)\in C_{\Gamma_\Sigma}(C_{\Gamma_\Sigma}(x))\}=Env(\Gamma_\Sigma)$\footnote{Here, $C_\Gamma(x)$ is the largest countable $\Gamma(x)$ set. It is defined to be the set of $y\in \bR$ such that for some set $A\in \Gamma\cap \powerset(\bR^3)$ and some ordinal $\a<\omega_1$, for every $z\in \bR$ coding $\a$, $y$ is the unique real such that $(y, x, z)\in A$.}.
\end{center}
Notice that if $(\Q, \Lambda)$ is a tail of $(\P, \Sigma)$ then $\Gamma(\Q, \Lambda)=\Gamma(\P, \Sigma)$. The next theorem is essentialy the conjunction  of \cite[Lemma 5.13-5.16]{ATHM}. 
\begin{theorem}\label{gamma(p, sigma) in the case p is not lsa} Assume ${\sf{AD^+}}$ and suppose $(\P, \Sigma)$ is a hod pair of successor type and $\Sigma$ is strongly fullness preserving and has strong branch condensation. Then the following holds.
\begin{enumerate}
\item There is a tail $(\Q, \Lambda)$ of $(\P, \Sigma)$ such that $\Gamma^*_\Lambda=\utilde{\Gamma}_\Lambda$.
\item Suppose $\Gamma^*_\Sigma=\utilde{\Gamma}_\Sigma$. Then for any real $x$ coding $\P^-$, 
\begin{center}
$C_{\Gamma_\Sigma}(x)={\sf{Lp}}^{\Gamma, \Sigma_{\P^-}}(x)$. 
\end{center}
\item Suppose $\Gamma^*_\Sigma=\utilde{\Gamma}_\Sigma$. Then ${\sf{Code}}(\Sigma)\not \in \Gamma(\P, \Sigma)$.
\item Suppose $\Gamma^*_\Sigma=\utilde{\Gamma}_\Sigma$. Then there is a tail $(\Q, \Lambda)$ of $(\P, \Sigma)$ such that 
\begin{center}
$\Gamma(\Q, \Lambda)=\powerset(\bR)\cap L(\Gamma(\Q, \Lambda), \bR)$.
\end{center}
Because $\Gamma(\Q, \Lambda)=\Gamma(\P, \Sigma)$, it follows that $\Gamma(\P, \Sigma)=\powerset(\bR)\cap L(\Gamma(\P, \Sigma), \bR)$.
\end{enumerate}
\end{theorem}

We spend the rest of this section defining $\Gamma(\P, \Sigma)$ in the case $\P$ is of $\#$-lsa type. The reader may wish to review \rdef{lsa type pair}, \rdef{pre sts hod pairs} and \rdef{sts hod pairs}. The difficulty with representing the $\sf{LSA}$ pointclasses as $\Gamma(\P, \Sigma)$ is the following. Suppose $\Gamma$ is an $\sf{LSA}$ pointclass, i.e., $\Gamma=\powerset(\bR)\cap L(\Gamma, \bR)$ and $L(\Gamma, \bR)\models \sf{AD}^++\sf{LSA}$. Let $\a$ be such that $\a+1=\Omega^\Gamma$ and set $\Gamma^{b}=\{ A\subseteq \bR: w(A)<\theta_\a\}$\footnote{The superscript ``b" stands for bottom.}. ${\sf{LSA}}$ pointclasses are peculiar because the pair that generates $\Gamma^b$ is essentially the same as the pair that generates $\Gamma$. More precisely, if $(\P, \Sigma)$ generates $\Gamma$ then $(\P, \Sigma^{stc})$ generates $\Gamma^b$. 

\begin{definition}
Suppose $(\P, \Sigma)$ is a hod pair such that $\P$ is of $\#$-lsa type and $\Sigma$ has strong branch condensation and is strongly fullness preserving\footnote{See \rdef{strong branch condensation} and \rdef{strongly fullness preserving}.}. We then let 
\begin{center}
$\Gamma(\P, \Sigma)=\{ A:$ for cone of $x\in \mathbb{R}$, $A\cap {\sf{Lp}}^{\Sigma^{stc}}(x)\in {\sf{Lp}}_2^{\Sigma^{stc}}(x)\}$.
\end{center}
$\myqedhere$
\end{definition}
Notice that the definition of $\Gamma(\P, \Sigma)$ depends on $\Sigma^{stc}$ and hence, can also be defined for sts pairs. It is not immediately clear that $L(\Gamma(\P, \Sigma))\cap \powerset(\bR)=\Gamma(\P, \Sigma)$. \rthm{gamma(p, sigma) in the lsa case} shows that it is indeed true. Before we prove it, we prove some useful lemmas. The first lemma shows that various $\Sigma$-sts mice are internally $\Sigma$-closed.

\begin{lemma}\label{general lemma on internal closure} Assume ${\sf{AD^+}}+{\sf{NsesS}}$\footnote{See \rdef{no mouse with a superstrong}.}. Suppose $(\R, \Phi)$ is an sts hod pair such that $\Phi$ has strong branch condensation and is strongly fullness preserving\footnote{See \rprop{positional}.} and $\M$ is a $\Phi$-sts mouse over $\R$. Suppose $\eta$ is a Woodin cardinal of $\M$ and $(\eta^+)^\M$ exists. Suppose further that whenever $\Q\in B(\R, \Sigma^{\M|\eta})$ and $\Q$ is of successor type, then $\Sigma^\M_\Q=\Phi_\Q\rest \M$. Given $\nu<\eta$, let $\S^\M_\nu$ be the last model of $(\R, \Sigma^\M)$-coherent fully backgrounded construction of $\M|\eta$ that uses extenders with critical points $>\nu$\footnote{See \rdef{full short tree coherent background constructions}. See also \rsec{sec: background constructions st-strategy}.} and let $\T_\nu$ on $\R$ be the normal tree leading to $\S^\M_\nu$.  Then for all $\nu<\eta$, $\pi^{\T_{\nu}, b}$ exists and $\pi^{\T_{\nu}, b}(\d^{\R^b})=\d^{(\S^\M_\nu)^b}$.
\end{lemma}
\begin{proof}
Towards a contradiction assume that for some $\nu$ our claim fails. Suppose first that $\pi^{\T_{\nu}, b}$ is undefined. We omit $\nu$ and $\M$ from subscripts and superscripts. Let $B$ be the set of layers $\P$ of $\S^\M$ such that $\pi^{\T}_{\leq \d^\P, b}$ exists. We then have that $\cup_{\P\in B}\P\not =\S$, and so letting $\a=\sup\{\d^\P: \P\in B\}$, $\a<\eta$ and $\pi^{\T_{\leq \a}, b}$ is defined. 

Suppose first that $\d^{(\M^\T_\a)^b}>\a$. Let $\Q\inseg_{hod}\M^\T_\a$ be the least complete layer\footnote{See \rnot{l p}.} of $\M^\T_\a$ such that $\cup B\inseg\Q$. It follows that $\T_{\geq \a}$ is a normal tree based on $\Q$. But since $\Sigma^M_\Q=\Phi\rest \M$, it follows from universality\footnote{See \rthm{existence of thick sets} and \rthm{universality of background construction}.} that $\lh(\T_{\geq \a})<\eta$ and $\pi^{\T_{\geq \a}}$ is defined. This is a contradiction, as it implies that there is $\Q'\in B$ such that $\d^{\Q'}>\a$.

 Assume now that $\d^{(\M^\T_\a)^b}=\a$. Since $\pi^{\T, b}$ does not exist, $\T_{\geq \a}$ must be based on $(\M^\T_\a)^b$ and be above $\d^{(\M^\T_\a)^b}$. However, it follows from our assumption that $\Sigma^\M_\Q=\Phi_\Q\rest \M$, and once again we get a counterexample to the universality of $\S$. 
 
 The proof that $\pi^{\T_{\nu}, b}(\d^{\R^b})=\d^{(\S^\M_\nu)^b}$ is very similar and we leave it to the reader.

\end{proof}
The following set up will be used in \rlem{fullness preservation of phi for n}, \rcor{fullness pres of phi}, \rcor{cor to full pres of phi}, \rlem{iterability of sts mice}, \rcor{cor to iterability of sts mice}, \rlem{gamma p sigma as a derived model}, \rcor{cor to gamma p sigma as a derived model} and \rlem{wadge rank computation}.\\\\
Assume ${\sf{AD}^++NsesS}$. Suppose $(\P,  \Sigma)$ is a hod pair such that $\P$ is of lsa type, $\P=(\P|\d^\P)^\#$ and $\Sigma$ has strong branch condensation and is strongly fullness preserving. Suppose ${\sf{Code}}(\Sigma)$ is Suslin, co-Suslin. Let $\Gamma$ be any good pointclasses such that ${\sf{Code}}(\Sigma)\in \utilde{\Delta}_{\Gamma}$. Let $\mathbb{M}=(M, \d, \vec{G}, \Omega)$ and let ${\sf{C}}=(\mathbb{M}, (P_0, \Psi_0), \Gamma^*, A)$ Suslin, co-Suslin capture both $\Gamma$\footnote{See \rdef{self-capturing bt}, \rdef{capturing gamma} and \rlem{n*x}.} and ${\sf{Code}}(\Sigma)$. We then have that the fully backgrounded hod pair construction of $\mathbb{M}$ reaches a tail of $(\P, \Sigma)$ (see \rthm{comparison holds}). Let $(\Q, \Lambda)$ be this tail. Let $\N$ be the last model of\footnote{See \rdef{fully backgrounded sts construction}.} 
\begin{center}
$({\sf{Le}}((\Q, \Lambda^{stc}), \mathcal{J}_{\omega}[\Q]))^{(M, \d, \vec{G})}$.
\end{center}

Because $\Sigma$ is fullness preserving we have that $\N\models ``\d^\Q$ is a Woodin cardinal".  Let $\Phi$ be the strategy of $\N$ induced by $\Psi$.  Notice that $\Phi$ is fullness preserving in the sense of ${\sf{Lp}}$ operator, i.e., whenever $\M$ is a $\Phi$-iterate of $\N$ and $\eta$ is a strong cutpoint of $\M$ then $\M|(\eta^+)^\M={\sf{Lp}}^{\Lambda^{stc}}(\M|\eta)$. This can be shown using the proof of \rthm{fullness preservation of background constructions}. We now prove several lemmas about $(\N, \Phi)$ leading up to showing that $\Gamma(\Q, \Lambda^{stc})$ can be realized as a derived model of $\N$. Let $\kappa$ be the least strong cardinal of $\N$. The first lemma is quite standard.

\begin{lemma}\label{least strong limit of woodins} $\N\models ``\k$ is a limit of Woodin cardinals".
\end{lemma}
\begin{proof} It is enough to show that $\d$ is a limit of $M$-cardinals $\eta$ such that ${\sf{Lp}}^{\Lambda^{stc}}(M|\eta)\models ``\eta$ is a Woodin cardinal". Fix $\nu<\d$. Because ${\sf{Code}}(\Sigma)\in \utilde{\Delta}_{\Gamma}$, we have that for cone of $z$, ${\sf{Lp}}^{\Sigma^{stc}}(z)\in C_{\Gamma}(z)$. We can assume, using absoluteness\footnote{See \rlem{correctness of backgrounds}.}, that the base of this cone is in $M$. Let $T, S\in M$ be $\d$-complementing trees witnessing that $A$ is Suslin, co-Suslin captured by $(M, \d, \vec{G}, \Omega)$. Let $\pi: R\rightarrow H_{(\d^+)^{M}}$ be a Skolem hull such that $\cp(\pi)>\nu$ is an $M$-cardinal and $\{ T, S\}\in \rge(\pi)$. Let $\eta=\cp(\pi)$. Then it follows that $C_{\Gamma}(M|\eta)\in M$ and hence, $C_{\Gamma}(M|\eta)\models ``\eta$ is a Woodin cardinal". It follows that ${\sf{Lp}}^{\Lambda^{stc}}(M|\eta)\models ``\eta$ is a Woodin cardinal". 
\end{proof}

The next lemma uses language introduced in \rdef{fullness preserving strategy for sigma mice}. 
\begin{lemma}\label{fullness preservation of phi for n}
$\Phi$ is fullness preserving, i.e., $\Phi$ witnesses that $\Gamma(\N|\k, \Phi)=\Gamma^b(\Q, \Lambda^{stc})$.
 \end{lemma}
\begin{proof}
Clearly, because $\Phi$ witnesses that $\N$ is a $\Lambda^{stc}$-sts mouse, $\Gamma(\N|\kappa, \Phi)\subseteq \Gamma^b(\Q, \Lambda^{stc})$. Fix then $(\T, \R)\in B(\Q, \Lambda^{stc})$. We want to see that\\\\
(1) there is a $\Phi$-iterate $\N_1$ of $\N|\k$ such that for some $t=(\Q, \T, \S, \U) \in \N_1$, $t$ is according to $\Sigma^{\N_1}$ and $\Lambda_\R\leq_w \Lambda_{\S^b}$.\\\\
Suppose (1) fails. We can then assume, without loss of generality, that for some $\nu<\d$ and some $g\subseteq Coll(\omega, \nu)$, $(\T, \R)\in M[g]$. Again without losing generality we can assume that $\R$ is of successor type. Let now $\S$ be the output of the $(\Q, \Sigma^\N)$-coherent fully backgrounded construction of $\N$ that uses extenders with critical points $>\nu$. Let $\U$ be a normal tree on $\Q$ with last model $\S$. We claim that\\\\
(2) $\pi^{\U, b}$ exists, $\pi^{\U, b}(\d^{\Q^b})=\d^{\S^b}$ and $\S'\inseg_{hod}{\S^b}$ is a $\Lambda_\R$-iterate of $\R$.\\\\
The first two clauses of (2) are consequences of \rlem{general lemma on internal closure}. The third is a straightforward consequence of the fact that $\Lambda$ is both positional and fullness preserving and of the fact that $\S$ side never moves in the comparison with $\R$\footnote{See \rprop{positional} and \cite[Lemma 2.6]{ATHM}.}. This finishes the proof of \rlem{fullness preservation of phi for n}.

\end{proof}

Before we proceed, we record some lemmas that can now be proved. Since these lemmas are standard, we will state these results without proofs and instead will give references. The next lemma can be proved following the proof of clause 2 of \rthm{main theorem on gen int} and also standard arguments like the proofs of Corollary 1.2 and Proposition 1.4, 1.5 of \cite{negres} and \cite[Chapter 3.1]{ATHM}. 

\begin{lemma}\label{fullness pres of phi} Suppose $\pi: \N|(\k^{+})^\N\rightarrow \M$ is an iteration via $\Phi_{\N|\k}$ and $g$ is $\M$-generic. Then letting $F$ be the function $F(X)={\sf{Lp}}^{\Lambda^{stc}}(X)$, $F\rest \M[g]$ is  definable over $\M[g]$ uniformly in $(\M, g)$\footnote{I.e., the definition works for any such $\M$ and $g$.}.
\end{lemma}
Below ${\sf{HC}}$ stands for the set of all hereditarily countable sets.
\begin{corollary}\label{cor to full pres of phi} Suppose $\pi: \N|(\k^+)^\N\rightarrow \M$ is an iteration via $\Phi_{\N|(\k^+)^\N}$ and $F$ is as in \rlem{fullness pres of phi}. Then if $h\subseteq Coll(\omega,<\pi(\k))$ is $\M$-generic then $F\rest {\sf{HC}}^{\M[h]}\in \M[\bR^{\M[h]}]$\footnote{Because $\k$ is a regular cardinal in $\N$, we have that $\bR^{M[h]}=(\bR^*)^{M[h]}$.}.
\end{corollary}

\rlem{fullness pres of phi} can be used to prove the following lemma. See also the proof of clause 2 of \rthm{main theorem on gen int}, \rprop{matching lp operators}, \rdef{certified iterations without fatal drops ii} and \cite[Proposition 1.5]{negres}. 

\begin{lemma}\label{iterability of sts mice} Suppose $\pi: \N|(\k^+)^\N\rightarrow \M$ is an iteration via $\Phi_{\N|(\k^+)^\N}$ and $\d$ is a cutpoint Woodin cardinal of $\M$. Let $\xi$ be a cutpoint cardinal of $\M$ such that $\M$ has no Woodin cardinals in the interval $(\xi, \d)$. Let $\eta\in (\xi, \d)$ be an $\M$-cardinal and let $\Psi$ be the fragment of $\Phi$ that acts on normal non-dropping trees based on $\M|(\eta^+)^\M$ that are above $\xi$. Then letting $h\subseteq Coll(\omega, (\eta^+)^\M)$ be $\M$-generic, $\Phi\rest \M|\pi(\k)[h]\in \M$ and is $\pi(\k)$-universally Baire in $\M[h]$. 
\end{lemma}

\begin{corollary}\label{cor to iterability of sts mice} Suppose $\pi: \N|(\k^+)^\N\rightarrow \M$ is an iteration via $\Phi_{\N|(\k^+)^\N}$. Suppose $g$ is $\M|\pi(\k)$-generic, $X\in (\M|\pi(\k))[g]$ and $\R\in {\sf{Lp}}^{\Lambda^{stc}}(X)$ is such that $\rho(\R)={\sf{ord}}(X)$. Let $h\subseteq Coll(\omega, \card{X})$ be $(\M|\pi(\k))[g]$-generic. Then $\R\in \M[g][h]$ and $\M[g][h]\models ``\R$ has a $\pi(\k)$-universally Baire iteration strategy $\Psi$ witnessing that $\R$ is a $\Lambda^{stc}$-sts mouse over $X$ based on $\Q"$. 

Moreover, if $\R\in (\M|\pi(\k))[g]$ is a sound $\Lambda^{stc}$-sts premouse over $X$ such that $\rho(\R)={\sf{ord}}(X)$ and for some $(\M|\pi(\k))[g]$-generic $h\subseteq Coll(\omega, \card{X})$, $\M[g][h]\models ``\R$ has a $\pi(\k)$-iteration strategy" then $\R\insegeq {\sf{Lp}}^{\Lambda^{stc}}(X)$\footnote{Notice that $\R$ has a unique $\pi(\k)$-iteration strategy in $\M[g][h]$.}. 
\end{corollary}

The next lemma shows that $\Gamma(\Q, \Lambda^{stc})$ can be realized as the set of reals of \textit{a derived model} of a $\Phi$-iterate of $\N$. We introduced the notation $D(\M, \l, h)$ in \rsec{sec short tree strategy mice}. The derived model theorem says that $D(\M, \l, h)\models {\sf{AD}^+}$, but we need a stronger version of this theorem. 

Suppose $\V$ is a transitive inner model of ${\sf{ZFC}-{Powerset}}$, $\lambda$ is a limit of Woodin cardinals of $\V$, $(\l^{++})^\V$ exists and $h\subseteq Coll(\omega, <\l)$ is $\V$-generic. Let 
\begin{center}
$\bR^*=\cup_{\a<\l}\bR^{\V[g\cap Coll(\omega, <\a)]}$
\end{center}
 and $\Gamma=\{A\in \V(\bR^*): \V(\bR^*)\models ``L(A, \bR^*)\models {\sf{AD^+}}"\}$. Set $D^+(\V, \l, h)=_{def}(L(\Gamma, \bR^*))^{\V(\bR^*)}$. Then the stronger version of Woodin's derived model theorem says that $D^+(\V, \l, h)\models {\sf{AD}^+}$. Sometimes $D^+(\V, \l, h)$ is called the \textit{new} derived model.
 
 Suppose now that in addition to the above, $\V$ is countable and $\Phi$ is an $\omega_1+1$-iteration strategy for $\V$. Let $g:\omega\rightarrow \bR$ be generic for $Coll(\omega, \bR)$ and let $(x_i: i<\omega)$ be the enumeration of $\bR$ given by $x_i=g(i)$. We can now perform an $\bR$-genericity iteration of $\V$ via $\Phi$ much like it is done in \cite[Chapter 7.2 and Corollary 7.17]{OIMT}. Let $\V'$ be this iterate of $\V$ and let $h\subseteq Coll(\omega, <\omega_1^V)$ be $\V'$-generic such that $(\bR^*)^{\V'[g]}=\bR^V$. We then let $D^+(\V, \Phi, \l, g)=D^+(\V', \omega_1^V, h)$.

\begin{lemma}\label{gamma p sigma as a derived model} The new derived model of $\N|(\k^+)^\N$ as computed via $\Psi=_{def}\Phi_{\N|(\k^+)^\N}$ is $L(\Gamma(\Q, \Lambda^{stc}))$. More precisely, for any $g\subseteq Coll(\omega, \bR)$, $D^+(\N|(\k^+)^\N, \Psi, \k, g)=L(\Gamma(\Q, \Lambda^{stc}))$ and $\powerset(\bR)\cap D^+(\V, \Psi, \l, g)=\Gamma(\Q, \Lambda^{stc})$. In particular, $\Gamma(\Q, \Lambda^{stc})=\powerset(\bR)\cap L(\Gamma(\Q, \Lambda^{stc}))$.
\end{lemma}
\begin{proof}
We will use clause 2 of \rthm{main theorem on gen int}. First we verify that clause 2 of \rthm{main theorem on gen int} applies. For this we need to verify that \\\\
(1) $\N$ is internally $\Lambda^{stc}$-closed, and\\
(2) $\Phi$ is a fullness preserving strategy for $\N$. \\\\
Notice that (1) is a consequence of \rlem{general lemma on internal closure} and (2) is just \rlem{fullness preservation of phi for n}.
We thus have that clause 2 of \rthm{main theorem on gen int} applies. 

To prove \rlem{gamma p sigma as a derived model} we need to show that given an $\bR$-genericity iteration $\pi:\N|(\k^+)^\N\rightarrow \W$ according to $\Phi_{\N|(\k^+)^\N}$,\\\\
(3) if $A\in \Gamma(\Q, \Lambda^{stc})$ then $A\in \W(\bR)$, and\\
(4) if $A\in \W(\bR)$ is such that $L(A, \bR)\models \sf{AD}^+$ then $A\in \Gamma(\Q, \Lambda^{stc})$. \\\\
We start with (3). Towards a contradiction, assume not and let $A\in\Gamma(\Q, \Lambda^{stc})$ witness this. We have that for cone of $z\in \bR$, $A\cap {\sf{Lp}}^{\Lambda^{stc}}(z)\in {\sf{Lp}}_2^{\Lambda^{stc}}(z)$. Let $z$ be some base of the aforementioned cone. Let $\xi>\Theta$ be such that $L_\xi(\powerset(\bR))\models {\sf{ZF}-Replacement}$ and $\sigma: M\rightarrow L_\xi(\powerset(\bR))$ be a countable hull such that $\N, z\in {\sf{HC}}^M$ and $\{\Phi, A\}\in \rge(\sigma)$. Let $A^M=\sigma^{-1}(A)$.

Let $g\in L(\powerset(\bR))$ be $M$-generic for $Coll(\omega, \bR^M)$.
Let $(y_i: i<\omega)$ be the generic sequence enumerating $\bR^M$ and let $(\d_i: i<\omega)$ be a sequence of cutpoint Woodin cardinals of $\N|(\k^+)$ with sup $\k$. Let $(\N_i, \T_i: i<\omega)$ be the $\bR^M$-genericity iteration. Thus, $\N_0=\N|(\k^+)^\N$, $\T_i$ is a tree on $\N_i$ that is based on $\N_i|\pi^{\oplus_{j<i}\T_j}(\d_i)$ and is above $\pi^{\oplus_{j < i}\T_j}(\d_{i-1})$\footnote{Let $\d_{-1}=0$.} and $\T_i$ is built according to the rules of $y_i$-genericity iteration. Let $\pi_{i, k}:\N_i\rightarrow \N_k$ be the composition of the iteration embeddings. Let $\N_\omega$ be the direct limit of $\N_i$ under $\pi_{i, k}$. 

Because $z\in \bR^M$, we have that $A\cap (\N_\omega|\omega_1^M)(\bR^M)\in {\sf{Lp}}^{\Lambda^{stc}}((\N_\omega|\omega_1^M)(\bR^M)))$. Notice that it follows from \rlem{fullness pres of phi} that if $\N_\omega^+$ is the last model of $\uparrow(\oplus_{i<\omega}\T_i, \N)$\footnote{See \rdef{upward extension of a stack}.} then 
\begin{center}
${\sf{Lp}}^{\Lambda^{stc}}((\N_\omega|\omega_1^M)(\bR^M))\in \N_\omega^+(\bR^M)$.
\end{center}
 It follows that $A^M \in D(\N_\omega, \omega_1^M, h)$ where $h\subseteq Coll(\omega, <\omega_1^M)$ is an $\N_\omega$-generic such that $\bR^{\N_\omega[h]}=\bR^M$. This finishes the proof of (3).

We keep the notation used to prove (3) and start proving (4). To prove (4), we need to show that if $A$ is as in (4) and $M,\sigma$ etc were defined as before relative to $A$ then\\\\
(5) $\sigma^{-1}(A)\in (\Gamma(\Q, \Lambda^{stc}))^M$.\\\\
Suppose that (5) fails. We then have that there is $B\in \N_\omega(\bR^M)$ such that $L(B, \bR^M)\models \sf{AD}^+$ and $B\not \in (\Gamma(\Q, \Lambda^{stc}))^M$. We first claim that\\

\textit{Claim.} in $L(B, \bR^M)$, for cone of $y$, $B\cap {\sf{Lp}}^{\Lambda^{stc}}(y)\in {\sf{Lp}}_2^{\Lambda^{stc}}(y)$.\\\\
\begin{proof}
Suppose not. Working in $L(B, \bR^M)$, fix $y\in \bR^M$ such that for any $y^*\in \bR^M$ Turing above $y$, $B\cap {\sf{Lp}}^{\Lambda^{stc}}(y)\not \in {\sf{Lp}}_2^{\Lambda^{stc}}(y)$. Fix $i<\omega$ such that $y\in \N_\omega[h\cap Coll(\omega, \d_i)]$. Notice that\\\\
(6) for every $y\in \bR^M$, $({\sf{Lp}}^{\Lambda^{stc}}(y))^{L(B, \bR^M)}={\sf{Lp}}^{\Lambda^{stc}}(y)$.\\\\
(6) is a consequence of \rcor{cor to iterability of sts mice}. This is because if $\R\insegeq ({\sf{Lp}}^{\Lambda^{stc}}(y))^{L(B, \bR^M)}$ is such that $\rho(\R)=\omega$ then $\R$ has an iteration strategy in $\N_\omega[y]$ as the iteration strategy of $\R$ is ordinal definable from $\Lambda^{stc}, y$ in the derived model of $\N_\omega$.  

Let $k<\omega$ be such that there is a symmetric name $\tau$ for $B$ in $\N_\omega[h\cap Coll(\omega, \d_k)]$. Let $j=\max(i, k)+1$. We then have that\\\\
(7) in $L(B, \bR^M)$, $B\cap (\N_\omega|\d_j)[h\cap Coll(\omega, \d_j)]\not \in {\sf{Lp}}^{\Lambda^{stc}}((\N_\omega|\d_j)[h\cap Coll(\omega, \d_j)])$.\\\\
However, it follows from \rlem{s-construction lemma} that \\\\
(8) ${\sf{Lp}}^{\Lambda^{stc}}((\N_\omega|\d_j)[h\cap Coll(\omega, \d_j)])=\N_\omega|(\d_j^+)^{\N_\omega}[h\cap Coll(\omega, \d_j)]$.\\\\
(8) and (7) contradict (6) (as $\tau_{h\cap Coll(\omega, \d_j)}=B\cap (\N_\omega|\d_j)[h\cap Coll(\omega, \d_j)]$).
\end{proof}

We will now make use of \cite[Theorem 0.1]{lpR}. It follows from the proof of the aforementioned theorem (applied to all sets of reals in $L(B, \bR^M)$) that \\\\
(9) $\powerset(\bR)^{L(B, \bR^M)}\subseteq ({\sf{Lp}}^{\Lambda^{stc}}(\bR^M))^{L(B, \bR^M)}$ and\\
(10) if $\K\insegeq({\sf{Lp}}^{\Lambda^{stc}}(\bR^M))^{L(B, \bR^M)}$ is such that $\rho(\K)=\bR$ and $k: \K'\rightarrow \K$ is such that $\K'$ is countable in $L(B, \bR^M)$ then $\K'\insegeq {\sf{Lp}}^{\Lambda^{stc}}(k^{-1}(\bR^M))$.\\\\
A Skolem hull argument done inside $\N_\omega^+$ shows that (10) implies that, \\\\
(11) $({\sf{Lp}}^{\Lambda^{stc}}(\bR^M))^{L(B, \bR^M)}\insegeq {\sf{Lp}}^{\Lambda^{stc}}(\bR^M)$.\\\\
Suppose now that \\\\
(a) ${\sf{Lp}}^{\Lambda^{stc}}(\bR^M)\in M$.\\\\
Then clearly (11) implies that \\\\
(12) $B\in M$.\\\\
(12) and the Claim imply (5). Thus, it is enough to prove that (a) holds. (a) easily follows from the fact that ${\sf{Lp}}^{\Lambda^{stc}}(\bR^M)\in \N_\omega^+(\bR)$ implying that ${\sf{Lp}}^{\Lambda^{stc}}(\bR^M)\in M[g]$. But since $g$ is generic for a homogenous poset, it follows that ${\sf{Lp}}^{\Lambda^{stc}}(\bR^M)\in M$.
\end{proof}

The following is a simple corollary of the proof of \rlem{gamma p sigma as a derived model}.

\begin{corollary}\label{cor to gamma p sigma as a derived model} Suppose $(\eta_i: i<\omega)$ is a sequence of consecutive Woodin cardinals of $\N|\k$ and $\l=\sup_{i<\omega}\eta_i$. The derived model of $\R=_{def}\N|(\l^+)^\N$ as computed via $\Phi_\R$ is $L(\Gamma(\Q, \Lambda^{stc}))$. In particular, $\Gamma(\Q, \Lambda^{stc})=\powerset(\bR)\cap L(\Gamma(\Q, \Lambda^{stc}))$.
\end{corollary}

 Let $\Q_\infty$ be the direct limit of all $\Lambda$-iterates of $\Q$ and let $\pi:\Q\rightarrow \Q_\infty$ be the iteration embedding. Let $\Omega$ be the $(\omega_1, \omega_1)$ fragment of $\Lambda_{\Q_\infty^b}$\footnote{See \cite{CoarseAD} where it is shown that $\Lambda$ is $<\Theta$-uB.}. Notice that $\pi\rest \Q^b$ depends only on $\Lambda^{stc}$\footnote{If $\T$ is the  $\Q$-to-$\Q_\infty$ stack then $\pi\rest \Q=\pi^{\T, b}$.} and hence (by the coding lemma), it is in $L(\Gamma(\Q, \Lambda^{stc}))$. Also, because $\Lambda^{stc}$ is fullness preserving, it follows that $\pi[\Q^b]$ can be coded as a subset of $w(\Gamma^b(\Q, \Lambda))$. This is because $\Q^b_\infty|\d^{\Q^b_\infty}=\bigcup\{\M_\infty(\R, \Lambda_\R): \R\in pB(\Q, \Lambda)\}$ and $\d^{\Q^b}=w(\Gamma^b(\Q, \Lambda))$. 

\begin{lemma}\label{wadge rank computation}  $\Lambda^{stc}\in \mathcal{J}_\omega(\pi[\Q^b], \Q^b_\infty, \Gamma^b(\Q, \Lambda))$.
\end{lemma}
\begin{proof} Set $\Psi=\Lambda^{stc}$. Notice that if $(\T, \S)\in I(\Q, \Psi)$ and $\W$ is a tree on $\S$ of limit length according to $\Psi_\S$ such that $\W$ is above $\d^{\S^b}$ and $\W\in b(\Psi_\S)$ then letting $b=\Psi_\S(\W)$\footnote{Thus, $b$ is a branch.}, $\Q(b, \W)$ exists and has an iteration strategy in $\Gamma^b(\Q, \Lambda)$. This is simply because there is an extender $E\in \vec{E}^{\M^\W_b}$ with critical point $\d^{\S^b}$ such that $\Q(b, \W)\inseg (Ult(\M^\W_b, E))^b$. We can then define $\Psi$ in $\mathcal{J}_\omega(\pi[\Q^b], \Q^b_\infty, \Gamma^b(\Q, \Lambda))$ with the following procedure.  We work in $\mathcal{J}_\omega(\pi[\Q^b], \Q^b_\infty, \Gamma^b(\Q, \Lambda))$.

Suppose first $X$ is a transitive set and $\R\in X$ is an lsa type hod mouse. Suppose that there is an embedding $\tau:\Q^b\rightarrow \R^b$. Suppose further that $\M$ is an sts mouse over $X$ based on $\R$. We say $\M$ is \textit{good} if it has an iteration strategy  $\Delta\in \Gamma^b(\Q, \Lambda)$ such that if $\S$ is a $\Delta$-iterate of $\M$, $t=(\R, \T, \R_1^*, \U)\in \S$ is according to $\Sigma^\S$, and $\R_1=\pi^{\T, b}(\R^b)$ then letting $\Delta_1=\Delta_{\R_1}$,
\begin{enumerate}
\item $(\R_1, \Delta_1)$ is a hod pair such that $\Delta_1$ has strong branch condensation and is strongly fullness preserving,
\item $\R_1=Hull^{\R_1}(\pi^{\T, b}\circ \tau[\Q^b]\cup \d^{\R_1})$,
\item letting $\sigma:\R_1\rightarrow \Q_\infty^b$ be given by
\begin{center}
$\sigma(x)=\pi(f)(\pi^{\Delta_1}_{\R_1, \infty}(a))$,
\end{center}
where $f\in \Q^b$ and $a\in (\d^{\R_1^b})^{<\omega}$ are such that $x=\pi^{\T, b}\circ \tau(f)(a)$,
\begin{center}
$\pi\rest \Q^b=\sigma\circ \pi^{\T, b}\circ \tau$.
\end{center}
\item $\U$ is according to $\Delta_1$.
\end{enumerate}
We can now define ${\sf{Lp}}^{good, sts, \tau}(X)$ which is the stack of good sts mice over $X$ that are based on $\R$. Then we can define ${\sf{Lp}}_{\omega}^{good, sts, \tau}(X)$. 

Suppose next that $\R$ is an lsa type hod premouse and $\tau:\Q^b\rightarrow \R^b$ is an embedding. Suppose $\U$ is a stack on $\R^b|\d^{\R^b}$. We say $(\R^b, \U)$ is a $\tau$-\textit{good} iteration if there is $k:\R^b\rightarrow \Q^b_\infty$ such that $\pi\rest \Q^b=k\circ \tau$ and for some $(\S, \Delta)\in \Gamma^b(\Q, \Lambda)$ such that $\Delta$ has strong branch condensation and is strongly fullness preserving, $k\rest (\R^b|\d^{\R^b})\subseteq \pi^\Delta_{\S, \infty}[\S]$ and if $\sigma:\R|\d^{\R^b}\rightarrow \S$ is given by
\begin{center}
$\sigma(x)=(\pi^\Delta_{\S, \infty})^{-1}(k(x))$
\end{center}
then $\U$ is according to $\sigma$-pullback of $\Delta$.

We can similarly define $\tau$-\textit{good} iterations when $\U$ is above $\d^{\R^b}$. In this case, we simply demand that $\U$ be according to the unique strategy of $\Psi'$ of $\R^b$ which acts on stacks that are above $\d^{\R^b}$ and letting $\Delta'=(\sigma$-pullback of $\Delta)$\footnote{Here, $\sigma$ is as above.}, $\Psi'$ witnesses that $\R^b$ is a $\Delta'$-premouse above $\d^{\R^b}$.

Suppose now that 
\begin{center}
$\T=((\M_\a)_{\a<\eta}, (E_\a)_{\a<\eta-1}, D, R, (\beta_\a, m_\a)_{\a\in R}, \so, \ma, T)$
\end{center}
 is an st-stack\footnote{See \rdef{st-stack}.} on $\Q$ of countable length. Recall \rrem{proper st-stack convention} and \rdef{proper stack}. These conventions stipulate that $R$ consists of cutpoints of $\T$. Also recall \rnot{notation for iteration trees}. We say $\T$ is $\pi$-\textit{b-realizable} if there is a sequence $(\sigma_\a: \a\in R)$ such that the following clauses hold\footnote{For the definition of $\pi^{\T, b}_{\a,\a'}$, see \rsec{sec:pitb}.}:
 \begin{enumerate}
 \item $\T$ doesn't have a fatal drop\footnote{See \rdef{fatal drop}.}, 
 \item $\sigma_\a:(\M_\a)^b\rightarrow \Q_\infty^b$ is an elementary embedding.
 \item For all $\a, \a'\in R$ with $\a<\a'$, $\sigma_\a=\sigma_{\a'}\circ \pi^{\T, b}_{\a,\a'}$.
 \item For all $\a\in R$, letting $\Lambda_\a=(\sigma_\a \rest \M_\a|\d^{\M_\a^b}$-pullback of $\Omega)$, for each complete layer $\R\inseg \M_\a^b$, $\sigma_\a\rest \R=\pi^{\Lambda_\a}_{\R, \infty}$ where $\pi^{\Lambda_\a}_{\R, \infty}: \R\rightarrow \M_\infty(\R, (\Lambda_\a)_\R)$ is the iteration map according to $(\Lambda_\a)_\R$.
 \item For all $\a\in R$ such that $\a\not =\max(R)$, letting $\a'=\min(R-(\a+1))$, if $\T_{\a, \a'}$ is based on $\M_\a^b|\d^{\M_\a^b}$ then $\T_{\a, \a'}$ is according to $\Lambda_\a$. 
 \item For all $\a\in R$ such that $\a\not =\max(R)$, letting $\a'=\min(R-(\a+1))$, if $\T_{\a, \a'}$ is based on $\M_\a^b$ and is above $\d^{\M_\a^b}$ then $\T_{\a, \a'}$ is according to the unique strategy of $\M_\a^b$ that acts on stacks above $\d^{\M_\a^b}$ and witnesses that $\M_\a^b$ is a $(\Lambda_\a)_{\M_\a|\d^{\M_\a^b}}$-mouse over $\M_\a|\d^{\M_\a^b}$.
 \item For every $\a\in R$ such that $\a+1<\lh(\T)$, letting $\W={\sf{nc}}^\T_\a$\footnote{See \rnot{notation for iteration trees}.}, for all limit ordinal $\gg<lh(\W)$ such that $\W\rest \gg$ is $\sf{nuvs}$, letting $\tau=\pi^{\T, b}_{0, \a}$,
\begin{enumerate}
\item if ${\sf{Lp}}^{good, sts, \tau}(\m^+(\W\rest \gg))\models ``\d(\W\rest \gg)$ is a Woodin cardinal" then $\lh(\W)=\gg+1$ and $\gg\in R$ and
\item if ${\sf{Lp}}^{good, sts, \tau}(\m^+(\W\rest \gg))\models ``\d(\W\rest \gg)$ is not a Woodin cardinal" then setting $b=[0, \gg)_\W$,  $b$ is a cofinal branch  for $\W\rest \gg$ such that $\Q(b, \W\rest \gg)$ exists and $\Q(b, \W\rest \gg)\insegeq {\sf{Lp}}^{good, sts, \tau}(\m^+(\W\rest\gamma))$.
\end{enumerate}
 \end{enumerate}
 
Let then $\Delta$ be an iteration strategy for $\Q$ such that its domain consists of st-stacks $\T$ which are $\pi$-b-realizable and for $\T\in  \dom(\Delta)$, $\Delta(\T)=b$ if and only if $\T^\frown \{b\}$ is a $\pi$-b-realizable st-stack. It can now be shown that $\Delta$ is the fragment of $\Psi$ that acts on st-stacks that do not have a fatal drop. The proof is very much like the proof of clause 2 of \rthm{main theorem on gen int} and it also uses \rdef{sts hod pairs}\footnote{In particular, see the conclusion of \rdef{sts hod pairs}.}. We leave it to the reader.

To compute $\Lambda^{stc}$, notice that $\Lambda^{stc}$ is the unique short-tree strategy $\Lambda'$ of $\Q$ such that $\Lambda'$ is fullness preserving and $\Delta$ as defined above is the fragment of $\Lambda'$ that acts on st-stacks without fatal drops. This easily follows from \rlem{disagreement implies low level disagreement}.  
\end{proof}

We are now in a position to state the main theorem of this section.

\begin{theorem}\label{gamma(p, sigma) in the lsa case} Assume ${\sf{AD}^++NsesS}$. Suppose $(\P,  \Sigma)$ is a hod pair such that $\P$ is of $\#$-lsa type\footnote{See \rdef{lsa type}.} and $\Sigma$ has strong branch condensation and is strongly fullness preserving. Suppose ${\sf{Code}}(\Sigma)$ is Suslin, co-Suslin. Then for some $\Q\in pI(\P, \Sigma)$,
\begin{enumerate}
\item $L(\Gamma(\Q, \Sigma_\Q))\cap \powerset(\bR)=\Gamma(\Q, \Sigma_\Q)$,
\item the set $\{ (x, y): x\in \bR$ and $y\not \in {\sf{Lp}}^{\Sigma^{stc}_\Q}(x)\}$ cannot be uniformized in $L(\Gamma(\Q, \Sigma_\Q))$, and 
\item $L(\Gamma(\Q, \Sigma_\Q))\models \sf{LSA}$.
\end{enumerate}
\end{theorem}
\begin{proof} Assume that one of 1-3 above is false. Let $\Gamma_0=\{A\subseteq \bR: A$ is ordinal definable from $\Sigma$ and a real$\}$. Then one of 1-3 is false inside $L(\Gamma_0, \bR)$, which means we can assume that $V=L(\Gamma_0, \bR)$. Let $(\a_0, \b_0)$ be lexicofraphically least such that letting $\Gamma^*_0=\{A\subseteq \bR: w(A)<\a_0\}$ the following holds:
\begin{enumerate}
\item $W=_{def}L_{\b_0}(\Gamma_0^*, \bR)\models {\sf{ZF}-Powerset}+``\Theta$ exists", 
\item $\Sigma\in \Gamma_0^*$ and $\a_0=\Theta^{L_{\b_0}(\Gamma_0^*, \bR)}$, and 
\item one of clauses 1-3 fails in $L_{\b_0}(\Gamma_0^*, \bR)$.
\end{enumerate}
Let $\Gamma_0=(\Sigma^2_1({\sf{Code}}(\Sigma)))^W$ and let $\Gamma$ be any good pointclass such that $\Gamma_0\subseteq \Delta_{\Gamma}$. Using \rthm{n*x} we can find $\mathbb{M}=(M, \d, \vec{G}, \Omega)$ and ${\sf{C}}=(\mathbb{M}, (P_0, \Psi_0), \Gamma^*, A)$ such that ${\sf{C}}$ Suslin, co-Suslin capture both $\Gamma$\footnote{See \rdef{self-capturing bt}, \rdef{capturing gamma} and \rlem{n*x}.}, ${\sf{Code}}(\Sigma)$ and the set $D$ consisting of triples $(u, v, w)\in \bR^3$ where $u$ codes $\Q\in pI(\P, \Sigma)$, $v$ codes a self-well-ordered $X\in {\sf{HC}}$ with $\Q\in X$ and $w$ codes ${\sf{Lp}}^{\Sigma^{stc}}(X)$. 

We then have that the fully backgrounded hod pair construction of $\mathbb{M}$ reaches a tail of $(\P, \Sigma)$ (see \rthm{comparison holds}).
Let $(\Q, \Lambda)$ be this tail (so $\Lambda=\Sigma_\Q$). Let $\N$ be the last model of 
\begin{center}
$({\sf{Le}}((\Q, \Lambda^{stc}), \mathcal{J}_{\omega}[\Q]))^{(M, \d, \vec{G})}$.
\end{center}
Because $\Sigma$ is fullness preserving we have that $\N\models ``\d^\Q$ is a Woodin cardinal".  Let $\Phi$ be the strategy of $\N$ induced by $\Omega$. We now start proving that $(\Q, \Lambda)$ is as desired. 

Clause 1 is just \rlem{gamma p sigma as a derived model}. We prove clause 2 of \rthm{gamma(p, sigma) in the lsa case}, which amounts to showing that the set $B=\{ (x, y):  x\in \bR \wedge y\not \in {\sf{Lp}}^{\Lambda^{stc}}(x)\}$ as computed in $W$ cannot be uniformized in $L(\Gamma(\Q, \Lambda^{stc}))$. Towards a contradiction assume that $B$ can be uniformized in $L(\Gamma(\Q, \Lambda^{stc}))$. It follows that we can find a set of reals $A\in \Gamma(\Q, \Lambda^{stc})$ such that $A$ codes a sjs $(A_i: i<\omega)$ with the property that $A_0=B$.

 Let $\pi:\N|(\k^+)\rightarrow \M$ be an $\bR$-genericity iteration. We then have that $A$ is in the (new) derived model of $\M$. Fix then a $<\pi(\k)$-generic $g$ over $\M$ such that there is a term relation $\tau\in \M[g]$ realizing $A$. Let $\d<\pi(\k)$ be a cutpoint Woodin cardinal of $\M$ which is not a limit of Woodin cardinals of $\M$ and such that $g$ is a $<\d$-generic. Let $\xi<\d$ be a cutpoint $\M$-cardinal such that $\M$ has no Woodin cardinals in the interval $(\xi, \d)$. Let $\M^*\inseg \M$ be such that $\tau\in \M^*$, $\M^*\models {\sf{ZFC}-Powerset}$ and $\M|\pi(\k)\insegeq \M^*$. Let now $\sigma:\S\rightarrow \M^*$ be such that $\cp(\sigma)\in (\xi, \d)$, $\sigma(\cp(\sigma))=\d$, $\cp(\sigma)$ is an $\M$-cardinal and $\tau\in \rge(\sigma)$. It follows that ${\sf{Lp}}^{\Lambda^{stc}}(\M|\cp(\sigma))\in \S$ and ${\sf{Lp}}^{\Lambda^{stc}}(\M|\cp(\sigma))\models ``\cp(\sigma)$ is a Woodin cardinal", contradiction! This finishes the proof of clause 2 of  \rthm{gamma(p, sigma) in the lsa case}.

To finish the proof of \rthm{gamma(p, sigma) in the lsa case} we need to show that $L(\Gamma(\Q, \Lambda))\models \sf{LSA}$. Suppose first that\\\\
(a) for every transitive $X\in {\sf{HC}}$ such that $\Q\in X$ and for every $\R\insegeq {\sf{Lp}}^{\Lambda^{stc}}(X)$ such that $\rho(\R)={\sf{ord}}(X)$, $\R$ has an iteration strategy $\Phi^*\in \Gamma^b(\Q, \Lambda)$ such that $\Phi^*$ witnesses that $\R$ is a $\Lambda^{stc}$-sts mouse over $X$ based on $\Q$. \\\\
We claim that (a) implies $L(\Gamma(\Q, \Lambda))\models \sf{LSA}$. Towards a contradiction assume not and set $B=\{ (x, y):  x\in \bR \wedge y\not \in {\sf{Lp}}^{\Lambda^{stc}}(x)\}$. We claim that\\\\
(1) $B$ is Suslin, co-Suslin in $L(\Gamma(\Q, \Lambda))$.\\\\
Clearly (2) contradicts clause 2 of \rthm{gamma(p, sigma) in the lsa case}. Set $\Psi=\Lambda^{stc}$. It follows from (a) that\\\\
(2) $B$ is projective in $\Psi$.\\\\
Let $\Q_\infty$ be the direct limit of all $\Lambda$-iterates of $\Q$ and let $\pi:\Q\rightarrow \Q_\infty$ be the iteration embedding. Notice that $\pi\rest \Q^b$ depends only on $\Psi$ and hence, because of \rlem{wadge rank computation}, it is in $L(\Gamma(\Q, \Lambda))$. Also, because $\Psi$ is fullness preserving, it follows that $\pi[\Q^b]$ can be coded as a subset of $w(\Gamma^b(\Q, \Lambda))$. This is because $\Q^b_\infty|\d^{\Q^b_\infty}=\bigcup\{\M_\infty(\R, \Lambda_\R): \R\in pB(\Q, \Lambda)\}$ and $\d^{\Q^b}=w(\Gamma^b(\Q, \Lambda))$. 

It follows from (2) and \rlem{wadge rank computation} that $B\in \mathcal{J}_\omega(\pi[\Q^b], \Q^b, \Gamma^b(\Q, \Lambda))$. Since we are assuming $L(\Gamma(\Q, \Lambda))\models \neg \sf{LSA}$ and since, in $L(\Gamma(\Q, \Lambda))$, $\d^{\Q^b_\infty}$ is both $<\Theta$ and is a limit of Suslin cardinals, $B$ must be Suslin, co-Suslin in $L(\Gamma(\Q, \Lambda))$, implying (1). Thus, it is enough to prove (a). 

Suppose (a) fails. We can then assume that the witness is in some $Coll(\omega, \Q)$-generic extension of $M$. Let $g\subseteq Coll(\omega, \Q)$ be $M$-generic and let $X\in {\sf{HC}}^{M[g]}$ be a counterexample to (a). We then have that $X$ is $<\k$-generic over $\N$. In fact, if $\eta\in ({\sf{ord}}(\Q), \k)$ is any Woodin cardinal of $\N$, then $X$ can be added to $\N$ by the extender algebra of $\N$ at $\eta$. Let then $\R\insegeq {\sf{Lp}}^{\Lambda^{stc}}(X)$ be the least such that $\rho(\R)=o(X)$ yet if $\Delta$ is the strategy of $\R$ witnessing that $\R$ is a $\Lambda^{stc}$-sts mouse over $X$ based on $\Q$ then $\Delta\not \in \Gamma^b(\Q, \Lambda)$. Notice that we have that\\\\
(3) ${\sf{Code}}(\Delta)$ is Suslin, co-Suslin in $L(\Gamma(\Q, \Lambda))$ (this follows from \rcor{cor to iterability of sts mice}). \\\\
It is then enough to show that the Suslin, co-Suslin sets of $L(\Gamma(\Q, \Lambda))$ are exactly those of $\Gamma^b(\Q, \Lambda)$. Assume otherwise. Let $\Q_\infty=\M_\infty(\Q, \Lambda)$. Because every set in $\Gamma^b(\Q, \Lambda)$ is $\d^{\Q_\infty^b}$-Suslin, co-Suslin we have that there is some $\eta<\k$ such that if $h\subseteq Coll(\omega, \eta)$ is $\N|\kappa$-generic then there is  
\begin{center}
$(\T, \S)\in I(\Q,  \Lambda^{stc})\cap {\sf{HC}}^{\N|\k[h]}$
\end{center}
 such that $\Lambda_{\S}\in L(\Gamma(\Q, \Lambda))$\footnote{This can be shown using \rthm{universality of background construction} and the fact that $\Lambda^{stc}$ is Suslin, co-Suslin in $L(\Gamma(\Q, \Lambda))$, which follows from our assumption and \rlem{wadge rank computation}.}. It then follows that 
 \begin{center}
 $\Lambda_\S\rest \N|\k[h]\in \N[h]$\footnote{$\Lambda_\S\rest \N|\k[h]\in \N[h]$ because of \rlem{fullness pres of phi}.}.
 \end{center} 

Let now $\nu>o(\S)$ be a cutpoint Woodin cardinal of $\N|\k$. Let $\S_1$ be an iterate of $\S$ above $\d^{\S}$ that is built according to the rules of $\N|\nu$-genericity iteration\footnote{This iteration starts by iterating the least measurable of $\S$ that is $>\d^{\S^b}$ $\nu+1$ times.}. For this genericity iteration we use the extender algebra at $\d^\S$ that uses extenders with critical points $>\d^{\S^b}$. Thus, the $\S$-to-$\S_1$ iteration is above $\d^{\S^b}$. We have that $\S_1\in \N[h]|(\nu^+)^\N$. Let $\N_1$ be the output of $({\sf{Le}}((\Q, \Lambda^{stc}), \mathcal{J}_{\omega}[\Q])^{(L[\N], \d, \vec{G'})}$ where $\vec{G'}$  consists of those extenders of $\N$ that have an inaccessible length (in $\N$) and a critical point $>\nu^{+}$. It follows from fullness preservation that $\N_1\models ``\d^\S$ is a Woodin cardinal".  

Let $\N_2$ be the $(\N_1, \nu, \pi^{\T, b}[\Q^b])$-authenticated backgrounded construction over $\N|\nu$ based on $\Q$\footnote{This makes sense as $\N|\nu$ is generic over $\S_1$ and $\pi^{\T, b}\in \N|\d$, see \rdef{certified backgrounded constructions}.}. Then it follows from universality of $\N_2$ that $\N|(\nu^+)^\N\subseteq \N_2\subseteq \N_1[\N|\nu]$. However, $\d^{\S_1}$ is not a cardinal of $\N$ yet it is a cardinal of $\N_1[\N|\nu]$, contradiction! This finishes the proof of (a) and hence, the proof of \rthm{gamma(p, sigma) in the lsa case}. 
\end{proof}

The next theorem can now be proved using \rcor{cor to gamma p sigma as a derived model} and the proof of Theorem 5.20 of \cite{ATHM}. 

\begin{theorem}\label{getting strong b-iterability}
Assume ${\sf{AD}^++NsesS}$. Suppose $(\P, \Sigma)$ is a hod pair such that $\P$ is either of successor type or of $\#$-lsa type  and $\Sigma$ has branch condensation and is fullness preserving. Suppose $B\in \mathbb{B}(\P^-, \Sigma_{\P^-})$. 
There is then $\Q\in pI(\P, \Sigma)$ and 
 $\vec{B}=\la B_i: i<\omega\ra\subseteq \mathbb{B}(\P, \Sigma_{\P^-})$
such that $\vec{B}$ strongly guides $\Sigma_{\Q}$ and $B_0=B$.
\end{theorem}

\section{A hybrid upper bound for LSA}

The main theorem of this section, \rthm{lsa from min omega woodins over lsa}, is a corollary to the proofs given in the previous section. It can be used in core model induction applications to show that certain hypotheses imply that there is a model of $\sf{LSA}$. We give a fairly detailed proof of \rthm{lsa from min omega woodins over lsa}.


\begin{definition}\label{omega woodins over lsa1} Suppose $(\P, \Sigma)$ is an sts hod pair\footnote{See \rdef{sts hod pairs}.}. We let $\N_{\omega.2, lsa}^{\#}(\P, \Sigma)$ be the minimal active $\Sigma$-sts mouse $\M$ over $\P$ such that $\M$ has $\omega.2$ many Woodin cardinals greater than $\d^\P$. $\myqedhere$
\end{definition}

Recall \rdef{closed under sharps} and \rdef{strategy premouse}.  Suppose now that $\M=\N_{\omega.2, lsa}^{\#}(\P, \Sigma)$. Let $\l$ be the supremum of the Woodin cardinals of $\M$. Because the only total extender of $\M$ whose  critical point is $>\l$ is the last extender of $\M$, the strategy predicate above $\l$ is empty. Thus, $\M=(\M|\l)^{\#}$. We use $\omega.2$ many Woodin cardinals because we need to produce proper initial segments of $\M$ that are unambiguous and satisfy the properties listed in clause 5 of \rdef{weak psi alpha indexing scheme a}. Notice that the way we stated clause 5 of \rdef{weak psi alpha indexing scheme a} implies that the strategy predicate of $\M|\gg$ cannot be empty above $\gamma$. We remark that we strongly believe that one could re-organize the manuscript in a way that we could prove all the lemmas in this section for $\N_{\omega, lsa}^\#(\P, \Sigma)$ which is the minimal active $\Sigma$-sts premouse over $\P$ that has $\omega$ Woodin cardinals above $\d^\P$. 

\begin{definition}\label{omega woodins over lsa} We say $\N$ is \textbf{an active $\omega.2$ Woodin lsa mouse} if it has an iteration strategy $\Sigma$ such that 
\begin{enumerate}
\item $\N$ has a Woodin cardinal $\d$ such that letting $\P=((\N|\d)^\#)^\N$, $(\P, \Sigma_\P^{stc})$ is an sts hod pair such that $\Sigma^{stc}_\P$ has strong branch condensation and is strongly $\Gamma^b(\P, \Sigma^{stc}_\P)$-fullness preserving,
\item $\N=\N_{\omega.2, lsa}^{\#}(\P, \Sigma^{stc}_\P)$,
\item for every $\P'\inseg_{hod} \P$ such that $\P'$ is of $\#$-lsa type\footnote{See \rdef{lsa type}. This means that $\P'=((\P'|\d^{\P'})^\#)^\P$.} layer of $\P$, $\N_{\omega.2, lsa}^\#(\P', \Sigma_{\P'}^{stc})\insegeq \P$ and 
\begin{center}
$\N_{\omega.2, lsa}^\#(\P', \Sigma_{\P'}^{stc})\models ``\xi$ is not a Woodin cardinal".  
\end{center}
\end{enumerate}
We say $\P$ is the lsa part of $\N$. We say $(\N, \Sigma)$ is \textbf{an active $\omega.2$ Woodin lsa pair}. 

$\myqedhere$
\end{definition}
Notice that if $(\N, \Sigma)$ is an active $\omega.2$ Woodin lsa pair then $\rho(\N)\leq (\k^+)^{\N}$ where, letting $\P$ be the lsa part of $\N$, $\k$ is the least $<\d^\P$-strong cardinal of $\P$\footnote{The fact that $\rho(\N)\leq (\k^+)^{\N}$ can be proved as follows. Suppose that $\rho(\N)> (\k^+)^{\N}$. Let $\M=Hull^\N((\k^+)^\N)$. Clearly $\M$ is also an active $\omega$ Woodin lsa mouse. We would be done if we had $\M\insegeq \N$. To show this, we use the proof of \rthm{condensation}, and compare $(\N, \M, (\k^+)^\N)$ with $\N$. We need to verify that a version of \rlem{dodd-jensen for certified phalanxes} holds for $(\N, \M, (\k^+)^\N)$. However, this can be done via exactly the same proof. We leave the details to the reader.}. 


In what follows, we let the statement \textit{there is an active $\omega.2$ Woodin lsa pair} be shortening for the statement that there is a pair $(\N, \Sigma)$ such that $\N$ is an active $\omega$ Woodin lsa mouse and $\Sigma$ witnesses the clauses of \rdef{omega woodins over lsa}.

Notice that it follows from \rthm{diamond comparison} that if $(\N, \Sigma)$ and $(\M, \Lambda)$ are two active $\omega$ Woodin lsa pairs with common lsa part $\P$ such that $\Sigma^{stc}=\Lambda^{stc}$ then $\N=\M$ and $\Sigma=\Lambda$. Let $I=\omega.2-\{\omega\}$.

\begin{lemma}\label{always authenticated} Suppose  $(\bar{\N}, \Sigma)$ is an active $\omega.2$ Woodin lsa pair and $\P$ is the lsa part of $\bar{\N}$. Let $\N$ be the result of iterating the last extender of $\bar{\N}$ through the ordinals. Let $(\d_i: i\in I)$ be the Woodin cardinals of $\N$ above $\d^\P$ and let $\l$ be their supremum. Let $\pi: \N\rightarrow \M$ be an iteration via $\Sigma$ that is above $\d^\P$. Suppose $g$ is $<\pi(\l)$-generic over $\M$ and $\W\in (\M|\l[g])\cap pB(\P, \Sigma^{stc})$\footnote{See \rdef{gamma(p, sigma) and b(p, sigma) for sts}. Recall that $pT$ is the projection of $T$.}. Let $k\in \omega$ be such that $g$ is generic for a poset in $\M|\pi(\d_k)$ and let $\S^\M_k$ be the last model of the $(\P, \Sigma^\M)$-coherent fully backgrounded construction of $\M|\pi(\d_{k+1})$ using critical points $>\d_k$\footnote{See \rdef{full short tree coherent background constructions}.}. Then the following holds:
\begin{enumerate}
\item Suppose $\T=_{def}\T^\M_k$ is the normal $\P$-to-$\S^\M_k$ stack. Then
\begin{enumerate}
\item $\lh(\T)$ is a limit ordinal,
\item $\T$ is $\sf{nuvs}$\footnote{See \rdef{nus stacks}.},
\item $\pi^{\T, b}$ exists,
\item  $\pi^{\T b}(\P^b)=(\S^\M_k)^b$, and
\item $\N_{\omega.2, lsa}^\#(\m^+(\T), \Sigma^{stc}_{\m^+(\T)})\models ``\d(\T)$ is a Woodin cardinal$"$.
\end{enumerate}
\item There is $\U\in \M[h]$ such that $\U$ is according to $\Sigma_\W$ and the last model of $\U$ is a layer of $(\S^\M_k)^b$. 
\end{enumerate}
\end{lemma}
\begin{proof} To make the proof notationally more pleasant, we ignore $\pi$ and assume $\N=\M$. The general case is very similar.

Clause 2 above follows from clause 1 and from the fact that $\Sigma$ is positional\footnote{See \rsec{sec positional and commuting}.} and that $(\S^\N_k)^b$-side doesn't move in the comparison of $\W$ and $(\S^\N_k)^b$. As proofs like this have appeared in the manuscript many times before we omit most of it. The exact procedure used to recover $\U\in \N[h]$ is the authentication process used to define sts mice\footnote{See \rsec{authentic iterations and finite stacks sec} and also the proof of \rsublem{next lemma}.}. 

 
 Clause 1.a and clause 1.b follows from standard arguments. Clause 1.b is a consequence of the fact that assuming $\T$ is a ${\sf{uvs}}$, $(\P, \T)$ is an indexable stack and since  $\N$ has more than $\d_{\k+1}$ many inaccessible cardinals, $\T\in \dom(\Sigma^\N)$ and hence, the construction producing $\S^\N_k$ can go further\footnote{See \rdef{indexable stack}.}. Clause 1.c and 1.d are straightforward consequences of clause 1.b\footnote{See \rlem{properties of proper stacks}.}. We verify clause 1.e.
 
 Let $\P_1=\m^+(\T)$. Notice that $\P_1|\d^{\P_1}=\S^\N_k$ and also $\P_1$ is a $\#$-lsa type. We want to see that  $\N_{\omega.2, lsa}^\#(\P_1, \Sigma^{stc}_{\P_1})\models ``\d^{\P_1}$ is a Woodin cardinal$"$. Towards a contradiction suppose\\\\
 (*) $\N_{\omega.2, lsa}^\#(\P_1, \Sigma^{stc}_{\P_1})\models ``\d^{\P_1}$ is not a Woodin cardinal$"$.\\\\
  Let $b=\Sigma(\T)$. (*) then implies that $\Q(b, \T)$ exists and is a $\Sigma_{\P_1}^{stc}$-sts mouse over $\P_1$. We now work towards showing that $\N$ has a branch indexed for $\T$, which is a contradiction as then the construction of $\S^\N_k$ can go further. 
  
 Working in $\N$, let $\Sigma_1$ be the $\N$-authenticated st-iteration strategy\footnote{See \rdef{n certified iteration strategy}.} of $\P_1$ and let $\K'$ be the output of the fully backgrounded construction of $\N|\l$ relative to $\Sigma_1$ done over $\mathcal{J}_{\omega}[\P_1]$ using extenders with critical point $>\d_k$\footnote{See \rdef{fully backgrounded sts construction}.} and let $\K=\mathcal{J}[\K']$. Notice that $\Sigma_1=\Sigma_{\P_1}\rest \N|\l$\footnote{See \rthm{main theorem on gen int}.}. \\

\textit{Claim 1.}  $\K$ has $\omega.2$ Woodin cardinals. In fact, for every $k'>k$, $\d_{k'}$ is a Woodin cardinal of $\K$.\\\\
\begin{proof} Suppose not. This means that the construction producing $\K$ doesn't reach $\lambda$. As iterability cannot be an issue (recall that $\N$ is iterable), the construction fails to reach $\l$ because the construction reaches  a model $\K^*$ such that there is an indexable stack $t=(\m^+(\T), \T_1, \P_2, \U)\in \K^*$  whose branch must be indexed but $t\not \in \dom(\Sigma_1)$. 
Notice now that $t$ cannot be ${\sf{nuvs}}$ as branches of such iterations are determined internally in $\K^*$\footnote{Recall that there can be an issue here. It can be the case that the branch determined by $\K^*$ does not agree with the branch determined by $\Sigma_1$. To show this, we use an argument like the one used in the proof of \rthm{sts fb constructions converge}.}. Thus, $t$ must be ${\sf{uvs}}$. Notice, however, that because $\P_2^b\in pB(\P_1, \Sigma_{\P_1})$, we have that $\P_2$ is $(\P,\Sigma^\N)$-authenticated and so, we must have that $(\P_2^b, \U)$ is an $(\P,\Sigma^\N)$-authenticated iteration. 
\end{proof}

Our goal now is to compare the construction producing $\K$ and $\Q(b, \T)$. Let $\Psi$ be the strategy of $\Q(b, \T)$ witnessing that $\Q(b, \T)$ is a $\Sigma^{stc}_{\P_1}$-sts mouse. Notice that we do not know that $\Q(b, \T)\in \N[h]$. The comparison that we use is the one used in \cite{SchSteel14}.\\

\textit{Claim 2.} The comparison of the construction producing $\K$ and $\Q(b, \T)$ is successful.\\\\
\begin{proof} Towards a contradiction assume not. We can then find a normal tree $\T_1$ on $\Q(b, \T)$ with last model $\Q_1$ and a normal tree $\U_1$ on $\N$ with last model $\N_1$ such that 
\begin{itemize}
\item $\T_1$ is according to $\Psi$, 
\item $\U_1$ is according to $\Sigma$ and has no drops,
\item for some $\b\not \in \dom(\vec{E}^{\Q_1})$, letting $\K_1=\pi^{\U_1}(\K)$, 
\begin{itemize}
\item $\Q_1|\b=\K_1|\b$,
\item $\b\not \in \dom(\vec{E}^{\K_1})$ and 
\item $\Q_1||\b\not =\K_1||\b$.
\end{itemize}
\end{itemize}
 Let then $t=(\P_1, \W, \R, \W_1)\in \Q_1|\b$ be an indexable stack  whose branch is indexed at $\b$ (either in $\Q_1$ or $\K_1$). As our indexing schema is local, it follows that a branch of $t$ must be indexed at $\b$ in both $\K_1$ and $\Q_1$. Since both $\Q(b, \T)$ and $\K_1$ are $\Sigma^{stc}_{\P_1}$-sts mice over $\P_1$, we have that $\K_1$ and $\Q_1$ cannot disagree on the branch of $t$.
\end{proof}

Because $\K$ has $\omega.2$ Woodin cardinals and is  a proper class model, it follows from Claim 2 and clause 3 of \rdef{omega woodins over lsa} that $\Q(b, \T)\insegeq \K$\footnote{This is because $\Q(b, \T)$ is $\omega.2$ small.}. We thus have that $\Q(b, \T)\in \N$. It follows that to show that $\N$ has a branch indexed for $\T$, it is enough to show that clause 5 of  \rdef{weak psi alpha indexing scheme a} holds for $\W=_{def}\Q(b, \T)$. To do this, we need to show that\\\\
(a) there is $\M\insegeq \N$ and a pair $(\b, \gg)$ such that,
\begin{enumerate}
\item $\b<o(\M)$ and $b\in \M|\b$,
\item $\M|\b$ is unambiguous (see \rdef{unambiguous hp}) and $\M|\b\models \sf{ZFC}$$+``$there are infinitely many Woodin cardinals $>\d(\T)$",
\item  letting $(\eta_i: i<\omega)$ be the first $\omega$ Woodin cardinals $>\d(\T)$ of $\M|\b$,  $\M|\b \models ``\W$  is $<Ord$-iterable above $\d(\T)$ via a strategy $\Phi$ such that letting $\nu=\sup_{i<\omega}\eta_i$, for every generic $g\subseteq Coll(\omega, <\nu)$, $\Phi$ has an extension $\Phi^+ \in D(\M|\b, \nu, g)$ such that $D(\M, \nu, g)\models ``\Phi^+$ is an $\omega_1$-iteration strategy" and whenever $\R\in D(\M|\b, \nu, g)$ is a $\Phi^+$-iterate of $\W$ and $t\in \R$ is an indexable stack on $\P_1$ then $t$ is $(\P, \Sigma^\M)$-authenticated.
\end{enumerate}

To show the existence of such an $\M$, it is enough to show that $\N|\d_{\omega+1}$ satisfies clauses 1-3 and first two clauses are straightforward. Let $(\eta_i: i\in \omega)$ enumerate $(\d_{i}: i\in (k+2, \omega))$ in increasing order.  We show that $(\eta_i: i\in \omega)$ witnesses clause 3 holds. \\

\textit{Claim 3.} if $\K_1$ is the $(\N, \d_{k+1})$-authenticated\footnote{See \rdef{certified backgrounded constructions}.} construction of $\N|\d_{k+2}$ done over $\mathcal{J}_{\omega}[\P_1]$ based on $\P_1$ then ${\sf{ord}}(\K_1)=\d_{k+2}$ and $\K_1\insegeq \K$. \\\\
\begin{proof}
Suppose not. It follows from the proof of Claim 2 that $\K_1$ has height $\d_{k+2}$. If $\K_1\not \insegeq \K$ then there is some model $\Q$ appearing in the construction producing $\K$ such that $\rho(\Q)<\d_{k+2}$. Let $p$ be the standard parameter of $\Q$. Let $X\prec \Q$ be such that $\rho(\Q)< X\cap \delta_{k+2}\in\delta_{k+2}$ and  $X\cap \delta_{k+2}$ is a cardinal in $\N$\footnote{This is possible because $\delta_{k+2}$ is strongly inaccessible in $\N$.} and $\bar{\Q}$ be the transitive collapse of $X$. By condensation (using the fact that $X$ contains solidity witnesses for $p$), $\bar{\Q}\inseg \Q$. Since $\bar{\Q}$ is sound and $\rho(\bar{\Q})=\rho(\Q) < X\cap \delta$, $X\cap \delta$ is not a cardinal in $\N$. Contradiction. 
\end{proof}

It follows from Claim 3 that $\W\insegeq \K_1$. To complete the proof of Clause 3 of (a), it is now enough to show the following claim.\\

\textit{Claim 4.} Suppose $\eta\in (\d_{k+1}, \d_{k+2})$ is an $\N$-cardinal and $g\subseteq Coll(\omega, (\eta^+)^\N)$. Let $\Phi$ be the fragment of $\Sigma$ that acts on non-dropping trees that are based on $\N|(\eta^+)^\N$ and are above $\d_{k+1}$. Then $\Phi\rest \N|\l[g]\in \N|\l[g]$ and if $\Lambda=\Phi\rest {\sf{HC}}^{\N|\l[g]}$ then in $\N[g]$, $\Lambda$ is a $<\l$-universally Baire iteration strategy such that for any poset $\mathbb{P}\in \N|\l[g]$, if $k\subseteq \mathbb{P}$ is $\N[g]$-generic and $\Lambda^k$ is the canonical extension of $\Lambda$ to ${\sf{HC}}^{\N|\l[g*k]}$ then $\Lambda^k=\Phi\rest {\sf{HC}}^{\N|\l[g*k]}$.\\\\
\begin{proof}
We only prove that $\Phi\rest \N|\l[g]\in \N|\l[g]$ and leave the rest to the reader. Let $\Q=\N|(\eta^+)^\N$ and let $\W_1\in \N[g]$ be a tree on $\Q$ of limit length and according to $\Phi$. Let $e=\Phi(\W_1)$. We want to show that $e\in \N[g]$ and $\N[g]$ has uniform way of identifying $e$. Notice that $\Q(e, \W_1)$ exists. Let $\K_2$ be the $\N$-authenticated background construction over $\M(\W_1)$. The proof of Claim 1 and Claim 2 show that $\Q(e, \W_1)\insegeq \K_2$. It is now easy to find the uniform definition of $e$. The reader may wish to consult the proof of \cite[Proposition 1.4]{negres}.
\end{proof}

Claim 4 finishes the proof of \rlem{always authenticated}.
\end{proof}

\begin{corollary}\label{bottom part fullness preservation} Suppose  $(\bar{\N}, \Sigma)$ is an active $\omega.2$ Woodin lsa pair and $\P$ is the lsa part of $\bar{\N}$. Let $\N$ be the result of iterating the last extender of $\bar{\N}$ through the ordinals. Let $\Phi$ be the fragment of $\Sigma$ that acts on stacks above $\d^\P$. Then $\Phi$ is $\Gamma^b(\P, \Sigma^{stc}_\P)$-fullness preserving\footnote{See \rdef{fullness preserving strategy for sigma mice}.}. 
\end{corollary}
\begin{proof} Given $\S\in pB(\P, \Sigma^{stc}_\P)$, let $\pi:\N\rightarrow \M$ be a $\Sigma$-iterate of $\N$ above $\d^\P$ such that $\S$ is generic over $\M$ for the extender algebra at the first Woodin of $\M$ that is larger than $\d^\P$. It follows from \rlem{always authenticated} that $\S$ is $\M$-authenticated\footnote{See \rdef{n certified iteration strategy}.}. 
\end{proof}

\begin{lemma}\label{reconstructing} Suppose  $(\bar{\N}, \Sigma)$ is an active $\omega.2$ Woodin lsa pair and $\P$ is the lsa part of $\bar{\N}$. Let $\N$ be the result of iterating the last extender of $\bar{\N}$ through the ordinals. Let $\d<\eta$ be two consecutive Woodin cardinals of $\N$ such that $\d>\d^\P$. Let $\N^*$ be the output of $(\N, \d)$-authenticated\footnote{See \rdef{certified backgrounded constructions}.} background construction of $\N|\eta$ done over $\mathcal{J}_{\omega}[\P]$ based on $\P$. Then
\begin{enumerate}
\item $\N^*$ has height $\eta$ and
\item if $\N_1$ is the result of translating $\N$ onto a structure over $\N^*$ via $S$-constructions\footnote{See \rdef{s-construction}.} then $\N_1$ is a normal iterate of $\N$ via a tree that is based on $\N|\d_0$ where $\d_0$ is the least Woodin cardinal of $\N$ above $\delta^\P$.
\end{enumerate}
\end{lemma}
\begin{proof}
We start by verifying clause 1. Suppose $\N^*$ fails to reach height $\eta$. This can only happen if at some stage of the construction we reach a model $\M$ such that there is some indexable stack $t=(\P, \T, \P_1, \VU)\in \M$ that is according to $\Sigma^\M$, it is required by the rules of the sts indexing scheme that we add a branch of $t$ to $\M$ but $t$ does not have an $\N$-authenticated branch. Notice that this can only happen when $t$ is $\sf{uvs}$ but \rlem{always authenticated} implies that any $\sf{uvs}$ has an $\N$-authenticated branch.

We now verify clause 2. Notice that $\N_1[\N|\eta]=\N$. Thus $\N_1$ is $\eta$-sound $\omega.2$ Woodin mouse. It is then enough to show that there is a tree $\U\in \N$ on $\N|\delta_0$ such that $\m(\U)=\N^*$. 

Suppose not. Let $\U\in \N$ be the normal stack on $\N|\d_0$ that is built by comparing $\N|\d_0$ with the construction producing $\N^*$. Since the aforementioned comparison fails, we must have that $\Sigma(\U)\not \in \N$. Let $b=\Sigma(\U)$. It follows from \rlem{s-construction lemma} that $\Q(b, \U)\in \N$\footnote{$\Q(b, \U)$ can be obtained via an $S$-construction, translating $\N$ to an sts mouse over $\m(\U)$.}. Hence, $b\in \N$, contraditicon. 

We must have that $\Q(b, \U)$ exists and $\Q(b, \U)\not \insegeq \N^*$. It follows that $\N^*\models ``\d(\U)$ is a Woodin cardinal". Thus, in the further comparison of $\Q(b, \U)$ and the construction producing $\N^*$, $\N^*$ side does not move. 
\end{proof}

\begin{theorem}\label{lsa from min omega woodins over lsa} Suppose  $(\bar{\N}, \Sigma)$ is an active $\omega.2$ Woodin lsa pair and $\P$ is the lsa part of $\bar{\N}$. Let $\N$ be the result of iterating the last extender of $\bar{\N}$ through the ordinals and let $\Sigma^>$ be the strategy of $\N$ that acts on iterations above $\d^\P$.  Let $\l$ be the supremum of the Woodin cardinals of $\N$ and let $\l'$ be the supremum of the first $\omega$ Woodin cardinals of $\N$. Then whenever $g\subseteq Coll(\omega, \bR)$ is generic, $D^+(\N, \Sigma^{>}, \l', g)\models \sf{LSA}$\footnote{See \rdef{gamma p sigma as a derived model}.}.
\end{theorem}
\begin{proof}
Let $(\d_i: i<\omega.2)$ be the Woodin cardinals of $\N$ and their limits that are greater than $\d^\P$. It follows from \rlem{general lemma on internal closure} that $\N|\l$ is internally $\Sigma^{stc}_\P$-closed. It follows from \rcor{bottom part fullness preservation} that $\Sigma^{>}$ is $\Gamma^b(\P, \Sigma^{stc}_\P)$-fullness preserving. 

Suppose $X$ is a transitive countable set such that $\P\in X$. Let for $i\in 2$, $\pi_i:\N\rightarrow \M_i$ be an iteration according to $\Sigma$ such that $\cp(\pi_i)>\d^\P$ and $X$ is $<\pi(\l')$-generic over $\M_i$.\\

\textit{Claim 1.} ${\sf{Lp}}^{\M_0, sts}(X, \P)={\sf{Lp}}^{\M_1, sts}(X, \P)$\footnote{See \rdef{simple q} and the discussion following it.}.\\\\
\begin{proof}
Let $\K_0$ be the $\M_0$-authenticated background construction over $X$ based on $\P$ and $\K_1$ be the $\M_1$-authenticated background construction over $X$ based on $\P$. We compare the construction producing $\K_0$ with the one producing $\K_1$. Notice that it follows from the proof of Claim 1 of \rlem{always authenticated} that both constructions reach proper class models. It then follows from the proof of Claim 2 of \rlem{always authenticated} that the aforementioned comparison produces $\sigma_0:\M_0\rightarrow \M_2$ and $\sigma_1:\M_1\rightarrow \M_3$ such that $\cp(\sigma_i)>{\sf{ord}}(X)$ and $\sigma_0(\K_0)$ and $\sigma_1(\K_1)$ are lined up (i.e. one is an initial segment of the other). Because they both have exactly $\omega.2$ Woodin cardinals it follows from our minimality assumption on $\N$ that $\sigma_0(\K_0)=\sigma_1(\K_1)$. The claim now follows. 
\end{proof}

Given a transitive $X\in {\sf{HC}}$, we let $\W(X)={\sf{Lp}}^{\M, sts}(X, \P)$ where $\M$ is such that there is an iteration $\pi:\N\rightarrow \M$ according to $\Sigma$ such that $\cp(\pi)>\d^\P$ and $X$ is $<\pi(\l')$-generic over $\M$. Suppose now that $\S'\in pI(\P, \Sigma)$ and $\S$ is a $\#$-lsa type\footnote{See \rdef{lsa type}.} proper layer of $\S'$. Let $\eta=\d^\S$. We then claim that \\\\
\textit{Claim 2.} $\W(\S)\models ``\eta$ is not a Woodin cardinal". \\\\
\begin{proof}
Suppose otherwise. Notice that $\S'\models ``\eta$ is not a Woodin cardinal". Let $\Q\insegeq \S'$ be the longest initial  segment $\Q^*$ of $\S'$ such that $\Q^*\models ``\eta$ is a Woodin cardinal". Then $\Q$ is a $\Sigma^{stc}_{\S}$-sts mouse. Let now $\pi:\N\rightarrow \M$ be an iteration according to $\Sigma^{>}$ such that $\S'$ is $<\pi(\l')$-generic over $\M$. Let $\K$ be the $\M$-authenticated background construction done over $\mathcal{J}_{\omega}[\S]$ based on $\S$. Because we are assuming that the claim fails, we must have that $\K\models ``\eta$ is a Woodin cardinal". 

We now compare $\Q$ with the construction of $\M$ producing $\K$. Notice that this comparison halts (this follows from the proof of Claim 2 that appears in the proof of \rlem{always authenticated}). Now, $\Q$ has to win this comparison. Since $\K$ is proper class and has $\omega.2$ Woodin cardinals, the fact that $\Q$ wins contradicts the minimality assumption on $\N$ (more precisely, contradicts clause 3 of Definition \ref{omega woodins over lsa}).
\end{proof}

Suppose next that $\S\in pI(\P, \Sigma)$ and $\eta=\d^\S$. 
We then have that\\\\
\textit{Claim 3.} $\W(\S)\models ``\eta$ is a Woodin cardinal$"$. \\\\
\begin{proof}
Let $\sigma:\N\rightarrow \S^+$ be the result of applying the iteration producing $\S$ to the entire model $\N$. Thus $\S$ is the lsa part of $\S^+$. Let now $\pi:\N\rightarrow \M$ be an iteration according to $\Sigma$ above $\d^\P$ such that $\S$ is $<\pi(\l)$-generic over $\M$. Let $\K$ be the $\M$-authenticated background construction done over $\mathcal{J}_{\omega}[\S]$ based on $\S$. We now compare the construction producing $\K$ with $\S^+$. As before this construction has to halt. It then follows from our minimality condition on $\N$ that $\W(\S)\models ``\eta$ is a Woodin cardinal$"$.
\end{proof}

The next claim computes the powerset of the Woodin cardinals of $\N$. The proof is very similar to the proof of Claim 3 and we omit it. \\\\
\textit{Claim 4.} Let $\pi:\N\rightarrow \M$ be an iteration according to $\Sigma$ above $\d^\P$. Then for any $k<\omega$, $\M|(\d_k^+)^\M=\W(\M|\d_k)$.\\\\
The next claim can be proved using the proof of Claim 3 and the proof of \rlem{cor to iterability of sts mice}. Also see the proof of Claim 4 of \rlem{always authenticated}.\\\\
\textit{Claim 5.} Suppose $X\in HC$ is a transitive set and $\R\insegeq \W(X)$ is such that $\rho(\R)=o(X)$. Let $\pi:\N\rightarrow \M$ be an iteration according to $\Sigma$ above $\d^\P$ such that $X$ is $<\pi(\l')$-generic over $\M$. Let $k<\omega$ be such that for some $g\subseteq Coll(\omega, <\pi(\d_k))$, $X\in {\sf{HC}}^{\M|\pi(\d_k)[g]}$. Then $\R$ has a $<\pi(\l')$-universally Baire iteration strategy in $\M[g]$.\\\\
Suppose $g\subseteq Coll(\omega, \bR)$ generic. Let $(x_i: i<\omega)$ be an enumeration of $\bR$ in $V[g]$. Let $\pi:\N\rightarrow \M$ be $\bR$-genericity iteration according to $\Sigma$ that is below $\l'$ and is guided by $(x_i: i<\omega)$. The next claim is a corollary to Claim 5 and clause 2 of \rthm{main theorem on gen int}.\\\\
\textit{Claim 6.} Set $B=\{ (x, y)\in \bR^2: y\not \in \W(x)\}$. Then $B\in \M(\bR)$ and $\Sigma^{stc}_\P\in \M(\bR)$.\\\\
 Let $\P_\infty$ be the direct limit of all $\Psi=_{def}\Sigma_\P$-iterates of $\P$ and let $\pi:\P\rightarrow \P_\infty$ be the iteration embedding. Notice that $\pi\rest \P^b$ depends only on $\Psi$. Also, because $\Psi$ is strongly $\Gamma^b(\P, \Psi)$-fullness preserving, it follows that $\pi[\P^b]$ can be coded as a subset of $w(\Gamma^b(\P, \Sigma^{stc}))$\footnote{$\d^{\P_\infty^b}$ is the largest cardinal of $\P_\infty^b$ and $\pi[\P^b]$ is cofinal in $\P_\infty^b$. Thus, we have $A\subseteq \d^{\P_\infty^b}$ which codes the pair $(\P^b, \pi[\P^b])$.}. This is because $\P^b_\infty|\d^{\P^b_\infty}=\bigcup\{\M_\infty(\R, \Lambda_\R): \R\in pB(\P, \Sigma^{stc})\}$ and $\d^{\P^b}=w(\Gamma^b(\P, \Sigma^{stc}))$. It follows from \rlem{wadge rank computation} that\\\\
\textit{Claim 7.} $\Psi\in \mathcal{J}_\omega(\P_\infty^b, \pi[\P^b], \Gamma^b(\P, \Psi))$.\\\\
Next we establish a crucial claim.\\\\
\textit{Claim 8.} $\mathcal{J}(\P_\infty^b, \pi[\P^b], \Gamma^b(\P, \Psi))\models \sf{AD}^+$.\\\\
\begin{proof} Suppose not. Let $A\in \mathcal{J}(\P_\infty^b, \pi[\P^b], \Gamma^b(\P, \Psi))$ be a set of reals that is not determined.  Let $X=\pi[\P^b]$. Fix $x\in \bR$ and $\Q\in pB(\P, \Psi)$ such that $A$ is definable from $(X, x, (\Q, \Psi_\Q), \P_\infty^b)$ and a finite sequence of ordinals over $\mathcal{J}(\P_\infty^b, \pi[\P^b], \Gamma^b(\P, \Psi))$. By minimizing the sequence of ordinals we can suppose that $A$ is definable without ordinal parameters.

Let $(\M_i, \T_i: i<\omega)$ be the $\bR$-genericity iteration of $\N$ relative to a generic enumeration $(x_i: i<\omega)$ of $\bR$ (this iteration is according to $\Sigma^{>}$, is below $\l'$ (the sup of the first $\omega$-Woodins of $\N$ and is above $\d^\P$). For $i<\omega$ let $\pi_i=\pi^{\oplus_{j\leq i}\T_j}$ and for $i<j\leq \omega$ let $\pi_{i, j}:\M_i\rightarrow \M_j$ be the composition of iteration embeddings. We then have that $A\in \M(\bR)$, where $\M$ is the direct limit of $\M_i$'s under the embeddings $\pi_{i, j}$. 

 Let $i$ be large enough so that 
$x, \Q\in {\sf{HC}}^{\M_i[(x_j: j\leq i)]}$ and $\Sigma_\Q\rest {\sf{HC}}^{\M_i[(x_j: j\leq i)]}$ is $<\pi_i(\l')$-universally Baire. Let $\tau\in \M_i[(x_j:j\leq i)]$ be a name such that $\pi_{i, \omega}(\tau)$ is a term relation for $A$. Let $\eta=\pi_i(\d_{i+1})$.  We claim that if 
\begin{center}
$\R=(\M_i|(\eta^+)^\M))[(x_j: j\leq i)]$
\end{center} then letting $\Phi$ be the fragment of $\Sigma_{\M_i}$ that acts on trees based on $\R$ that are above $\pi_i(\d_i)$, $(\R, \Phi, \tau_\R)$ term captures $A$ where $(p, u)\in \tau_\R$ if and only if $p\in Coll(\omega, \eta)$, $u\in \R^{Coll(\omega, \eta)}$ and $p\forces ``\forces_{Coll(\omega, <\pi_i(\l'))} u\in \tau"$. It then follows from a result of Neeman that $A$ is determined (see \cite{OPD}). 

Let then $\T$ be an iteration tree on $\M_i$ based on $\R$ according to $\Phi$. Let $\S$ be the last model of $\T$. We want to see that\\\\
(a) if $h\subseteq Coll(\omega, \pi^\T(\eta))$ is $\S$-generic then $(\pi^\T(\tau_\R))_h=A\cap \S[h]$.\\\\
Let $k>i$ be large enough that $\S\in \M_k[(x_j: j\leq k)]$. Let $\S^*$ be the output of $\M_k|\pi_k(\d_{k+1})$-authenticated backgrounded construction done over $\S|\pi^\T(\eta)$ based on $\P$. We then have that $\S^*$ is an iterate of $\S|\pi^\T(\pi_i(\d_{i+2}))$\footnote{See \rlem{reconstructing}. More precisely, there is a normal stack $\U$ on $\S|\pi^\T(\pi_i(\d_{i+2}))$ that is above $\pi^\T(\eta)$, $\lh(\U)$ is a limit ordinal and $\m(\U)=\S^*$.}. Let $\S^{**}=\pi_{k, k+1}(\S^*)$. Finally, let $\S_1$ be the result of translating $\M_{k+1}$ over $\S^{**}$ via $S$-constructions. We then have that\\\\
(1) $\S_1[\M_{k+1}|\pi_{k+1}(\d_{k+1})]=\M_{k+1}$\\
(2) $\S_1$ is an iterate of $\S$ such that if $\nu:\S\rightarrow \S_1$ is the iteration embedding then $\cp(\nu)>\pi^\T(\eta)$.\\\\
(2) is a consequence of the fact that $\S^{**}=\m(\U')$ where $\U'=\pi_{k, k+1}(\U)$ and $\U$ is as in the footnote above. It then follows that if $b$ is the branch of $\U'$ given by $\Sigma_{\M_i}$ then $\M^\U_b$ is $\pi_k(\d_{k+1})$-sound. 

It follows that we can think of $p=(\T_j: j\in (k+1, \omega))$ as an $\bR$-genericity iteration on $\S_1$ guided by $(x_j: j\in (k+1, \omega))$. Let then $\S_2$ be the last model of this genericity iteration and let $m: \N\rightarrow \S_2$ be the iteration embedding.
More precisely, $m=\pi^p\circ \nu \circ \pi^\T\circ \pi_i$. 

Because $\M|\pi_{k+1}(\d_{k+1})=\M_{k+1}|\pi_{k+1}(\d_{k+1})$, we  have that $\S_2[\M|\pi_{k+1}(\d_{k+1})]=\M$. Let $\sigma: \M_i\rightarrow \S_2$ be the iteration embedding. It then follows that in $\S_2[(x_j: j\leq i)]$, $\sigma(\tau)$ is the term relation that is forced by $Coll(\omega, <m(\l'))$ to be the least set in $\mathcal{J}(\P_\infty^b, \pi[\P^b], \Gamma^b(\P, \Psi))$ which is not determined and is definable from $(X, x, (\Q, \Psi_\Q), \P_\infty^b)$. It then follows that\\\\
(3) $\sigma(\tau)$ is realized as $A$.\\\\
 (a) now follows from (2) and the fact that $\cp(\pi^p)>\pi^\T(\eta)$.
\end{proof}

The proof of the next claim is exactly like the proof of (a) that appeared in the proof of \rthm{gamma(p, sigma) in the lsa case} and \rlem{always authenticated}. We leave it to the reader.\\\\
\textit{Claim 9.} For any transitive $X\in \sf{HC}$ such that $\P\in X$ and for any $\R\insegeq \W(X)$ such that $\rho(\R)=o(X)$, $\R$ has an iteration strategy in $\Gamma^b(\P, \Psi)$.\\\\

It follows from Claim 9 that the set $B=\{ (x, y)\in \bR^2: y\not \in \W(x)\}$ is projective in $\Psi$ and hence, $B\in \mathcal{J}(\P_\infty^b, \pi[\P^b], \Gamma^b(\P, \Sigma))$. It follows from Claim 9 that $\mathcal{J}(B)\models \sf{AD}^+$. We now have the following:\\\\
\textit{Claim 10.}  In $\M(\bR)$, let $\Gamma=\{ A\subseteq \bR: L(A, \bR)\models \sf{AD}^+\}$. Then $\Psi, B\in L(\Gamma, \bR)$.\\\\

It follows from the proof of clause 2 of \rthm{gamma(p, sigma) in the lsa case} that $B$ cannot be uniformized in $L(\Gamma, \bR)$. Hence, $L(\Gamma, \bR)\models \sf{LSA}$.
\end{proof}

\chapter{Condensing sets}\label{chap:condensing_sets}

The goal of this chapter is to introduce the theory of condensing sets. Such sets were first considered in \cite[Section 10, 11.1]{CuBF}, where they were presented in the form of a condensation property for elementary embeddings (see \cite[Definition 11.14]{CuBF}).  The current presentation dates back to an unpublished note by the first author.

Prior to this work, condensing sets have been used in the context of the core model induction. As a convenience to the reader, we recap some of the basic machinery used in the core model induction. We model our presentation on \cite{CuBF} but we will also use the set up of \cite{Trang2015PFA}. A typical situation is as follows. We have an embedding $j:M\rightarrow N$ with critical point $\kappa$ and such that $H_{\kappa^+}^M=H_{\kappa+}^N$. In $M$, we consider the maximal model of determinacy that has been built via core model induction. While the exact definition of the maximal model is somewhat case specific, it can be essentially described as follows.

Let $g\subseteq Coll(\omega, <j(\kappa))$ be $N$-generic. For $\nu<j(\k)$ let $g_\nu=g\cap Coll(\omega, <\nu)$. We then can extend $j$ to act on $M[g_\k]$. We denote this extension by $j$ again and we have that $j:M[g_\k]\rightarrow N[g]$. 

Consider the set of hod pairs $(\Q, \Lambda)$ such that 
\begin{enumerate}
\item $\Q\in {\sf{HC}}^{M[g_\kappa]}$,
\item for some $\nu<\k$ such that $\Q\in M[g_\nu]$, letting $\Psi=\Lambda\rest {\sf{HC}}^{M[g_\nu]}$, $\Psi\in M[g_\nu]$ and $M[g_\nu]\models ``{\sf{Code}}(\Psi)$ is $\k$-uB" and
\item if $T, S\in M[g_\nu]$ witness that ${\sf{Code}}(\Psi)$ is $\k$-uB then ${\sf{Code}}(\Lambda)=p[T]^{M[g_\k]}$. 
\end{enumerate}
Let $\Gamma$ be the set of such pairs $(\Q, \Lambda)$. An additional requirement is that $\Lambda$ is fullness preserving and has branch condensation. While the branch condensation is the same as before, fullness preservation is not the same as the definition given in earlier chapters. We refer the interested reader to \cite{CuBF} for more details on how to define $\Gamma$. It is in fact somewhat more involved. 

The goal of a core model induction is to show that $\Gamma$ is rich. This is done as follows. First a target theory is fixed. The theory used in \cite{CuBF} is $``\sf{AD}_{\bR}+``\Theta$ is regular". In \rchap{chap:lsa_from_pfa}, our target is $\sf{LSA}$. Suppose then there is no lsa type hod pair $(\Q, \Lambda)\in \Gamma$. Preliminary arguments, such as those used in \cite[Theorem 4.1]{UBH}, show that $\Gamma$ is of limit type, i.e., for any $(\Q, \Lambda)\in \Gamma$ there is $(\R, \Psi)\in \Gamma$ such that $\Gamma(\Q, \Lambda)\subset \Gamma(\R, \Psi)$. 

Next we let $\P^-=\bigcup_{(\Q, \Lambda)\in \Gamma}\M_\infty(\Q, \Lambda)$. Fixing a complete layer\footnote{See \rdef{l p}.} $\R$ of $\P^-$ and $(\Q, \Lambda)\in \Gamma$ such that $\R=\M_\infty(\Q, \Lambda)$, we let $\Sigma_\R=\Lambda_{\R}$. It follows from comparison that $\Sigma_R$ is independent of $(\Q, \Lambda)$. Let $\Sigma=\oplus_{\R}\Sigma_\R$ where the joint ranges over the complete layers of $\P^-$. 

We now define $\P$ as follows. Suppose next that there is $\M\insegeq {\sf{Lp}}^{\Sigma}(\P^-)$ such that $\rho(\M)<{\sf{ord}}(\P^-)$. We then let $\P$ be the least such $\M$. Otherwise we let $\P={\sf{Lp}}^{\Sigma}(\P^-)$. 

The next major step is to build an iteration strategy for $\P$ that extends $\Sigma$. We let $\Sigma^+$ be this new strategy.  $\Sigma^+$ is constructed as follows.

\begin{definition}[The construction of the strategy] \label{the construction of the strategy} 
Suppose  $\T\in {\sf{HC}}^{N[g]}$ is a stack on $\P$ where
\begin{center}
$\T=((\M_\a)_{\a<\eta}, (E_\a)_{\a<\eta-1}, D, R, (\beta_\a, m_\a)_{\a\in R}, T)$.
\end{center}
Recall \rnot{notation for iteration trees}. Suppose $j\rest \P\in N[g]$. Working in $N[g]$, we say $\T$ is \textbf{$j$-realizable} if there is a sequence $(\sigma_\a: \a\in R)$ such that the following clauses hold\footnote{For the definition of $\pi^{\T, b}_{\a,\a'}$, see \rsec{sec:pitb}.}:
 \begin{enumerate}
 \item $\T$ doesn't have a fatal drop\footnote{See \rdef{fatal drop}.}, 
 \item $\sigma_\a:\M_\a\rightarrow j(\P)$ is an elementary embedding.
 \item For all $\a, \a'\in R$ with $\a<\a'$, $\sigma_\a=\sigma_{\a'}\circ \pi^{\T}_{\a,\a'}$.
 \item For all $\a\in R$, letting $\Lambda_\a=(\sigma_\a \rest \M_\a|\d^{\M_\a}$-pullback of $j(\Sigma))$, for each complete layer $\R\inseg \M_\a$, $\sigma_\a\rest \R=\pi^{\Lambda_\a}_{\R, \infty}$ where $\pi^{\Lambda_\a}_{\R, \infty}: \R\rightarrow \M_\infty(\R, (\Lambda_\a)_\R)$ is the iteration map according to $(\Lambda_\a)_\R$\footnote{This condition assumes that $\Lambda_\a$ is fullness preserving.}.
 \item For all $\a\in R$ such that $\a\not =\max(R)$, letting $\a'=\min(R-(\a+1))$, $\T_{\a, \a'}$ is according to $\Lambda_\a$\footnote{Notice that because we are assuming $\T$ does not have fatal drops, $\T_{\a, \a'}$ is based on $\M_\a|\d^{\M_\a}$.}. 
 \end{enumerate}

Given  a stack  $\T\in {\sf{HC}}^{N[g]}$, we set $\T\in \dom(\Sigma^+)$ if $\T$ is $j$-realizable.  For $\T\in \dom(\Sigma)$, we set $\Sigma^+(\T)=b$ if $\T^\frown \{\M^\T_b\}$ is $j$-realizable. $\myqedhere$
\end{definition}

It is not hard to extend $\Sigma^+$ to act on all stacks, not just those without fatal drops. $\Sigma^+$ may not be a total strategy simply because we may not be able to satisfy clauses 4 and 5 of \rdef{the construction of the strategy}. Moreover, it may also depend on the realization maps. However, the proof of \cite[Lemma 11.6]{CuBF} gives the following.

\begin{theorem}\label{existence of sigma+} Suppose $\card{\P}<(\k^+)^M$. Then $j\rest \P\in N[g]$ and  $\Sigma^+$ is a total $(\omega_1, \omega_1)$-strategy in $N[g]$.
\end{theorem}

Then there are two arguments that we run as part of the proof of Theorem \ref{existence of sigma+}. First we show that $\P={\sf{Lp}}_\omega^{\Sigma}(\P^-)$. The reader can see, for example \cite[Lemma 3.78]{Trang2015PFA}, for an argument. 
Roughly, if not, suppose $n$ is such that $\rho_{n+1}(\P)<\delta^\P\leq\rho_n(\P)$, then in $j(\Gamma)$, we can find a complete layer $\R$ of $\P$ and an $OD^{j(\Gamma)}_{\Sigma_\R}$ set $A\subseteq \delta^\R$ such that $A\notin \P$. By fullness of $\P$ and $\sf{SMC}$ in $j(\Gamma)$, $A\in \P$. Contradiction.

The next argument attempts to show that $\P\models ``\d^\P$ is regular". Showing this finishes the proof of the main theorem of \cite{CuBF}. In this book we present an argument for obtaining a model of $\sf{LSA}$ from $\sf{PFA}$ (see \rthm{thm:square_lsa}). To prove Theorem \ref{thm:square_lsa}, we need to do more in order to finish the argument. It is in this step that the theory of condensing sets is used. A reader interested in more details may consult \cite[Section 10, 11.1]{CuBF} and \cite[Lemma 3.81]{Trang2015PFA}.

\section{Condensing sets}\label{condensing sets sec}

We introduce the notion of condensing set in the most general setting. Suppose $\phi$ is a formula in the language of set theory and $A$ is a set. We let $\mathcal{F}_{\phi, A}$ be a collection of hod pairs $(\Q, \Lambda)$ such that $\Q$ is countable, $\Lambda$ is an $(\omega_2, \omega_2, \omega_2)$-iteration strategy having strong branch condensation and such that $\phi[A, (\Q, \Lambda)]$ holds. 

\begin{terminology}\label{closure properties of f}
\begin{enumerate}
\item We say $(\phi, A)$ is \textbf{bottom part closed} if whenever $(\Q, \Lambda)\in \mathcal{F}_{\phi, A}$ and $\R\in pB(\Q, \Lambda)$ then $(\R, \Lambda_\R)\in \mathcal{F}_{\phi,A}$. 
\item We say $(\phi, A)$ is of \textbf{limit type} if for every $(\Q, \Lambda)\in \mathcal{F}_{\phi, A}$, there is $(\R, \Psi)\in \mathcal{F}_{\phi, A}$ such that $\R$ is of limit type and ${\sf{Code}}(\Lambda)\in \Gamma^b(\R, \Psi)$.  
\item Let $\Gamma_{\phi, A}=\bigcup\{ \Gamma(\R, \Psi): (\R, \Psi)\in  \mathcal{F}_{\phi, A}\wedge \R$ is of limit type$\}$. We say $(\phi, A)$ is \textbf{stable} if whenever $(\R, \Psi)\in \mathcal{F}_{\phi, A}$, $\Psi$ is strongly $\Gamma_{\phi, A}$-fullness preserving. 
\item We say $(\phi, A)$ is \textbf{directed} if whenever $(\Q, \Lambda), (\P, \Sigma)\in \mathcal{F}_{\phi, A}$, there are $\R\in pI(\Q, \Lambda)$ and $\S\in pI(\P, \Sigma)$ such that either
\begin{enumerate}
\item $\R\insegeq_{hod}\S$ and $\Sigma_\R=\Lambda_\R$ or
\item $\S\insegeq_{hod}\R$ and $\Lambda_\S=\Sigma_\S$.
\end{enumerate} 
\end{enumerate}
$\myqedhere$
\end{terminology}
\begin{notation}\label{complete layer notation} Given a hod premouse $\P$, we write $\R\insegeq^c_{hod}\P$ if $\R$ is a complete\footnote{See \rnot{l p}.} layer of $\P$.$\myqedhere$
\end{notation}

\begin{notation}\label{direct limit}
Suppose $(\phi, A)$ is bottom part closed, is of limit type, is stable and is directed. 
\begin{enumerate}
\item Let $\P^-_{\phi, A}=\bigcup_{(\Q, \Lambda)\in \mathcal{F}_{\phi, A}}\M_\infty(\Q, \Lambda)$. 
\item Fix $\R\inseg^c_{hod} {\P^-_{\phi, A}}$ and $(\Q, \Lambda)\in \mathcal{F}_{\phi, A}$ such that $\R=\M_\infty(\Q, \Lambda)$. Let $\Sigma_{\R, \phi, A}=\Lambda_{\R}$ and let $\Sigma_{\phi, A}=\oplus_{\R\inseg^c_{hod}\P^-_{\phi, A}}\Sigma_{\R,\phi,A}$.
\item Suppose there is $\M\insegeq {\sf{Lp}}^{\Gamma_{\phi, A}, \Sigma_{\phi,A}}(\P^-_{\phi, A})$ such that $\rho(\M)<{\sf{ord}}(\P^-_{\phi, A})$. Then let $\P_{\phi, A}$ be the least such $\M$. Otherwise let $\P_{\phi, A}={\sf{Lp}}^{\Gamma_{\phi, A}, \Sigma_{\phi, A}}(\P^-_{\phi, A})$.
\end{enumerate}
In clause 3 above $\M\insegeq {\sf{Lp}}^{\Gamma_{\phi, A}, \Sigma_{\phi,A}}(\P^-_{\phi, A})$ if and only if whenever $\pi:\M'\rightarrow \M$ is an elementary embedding and $\M'$ is countable, $\M'\insegeq{\sf{Lp}}^{\Gamma_{\phi, A}, \Sigma^\pi_{\phi,A}}(\pi^{-1}(\P^-_{\phi, A}))$.$\myqedhere$
\end{notation}

\begin{definition}\label{full collection of hod pairs} Suppose $(\phi, A)$ is bottom-part closed, is of limit type, is stable and is directed. We say $(\phi, A)$ is \textbf{full} if $\P_{\phi, A}={\sf{Lp}}^{\Gamma_{\phi, A}, \Sigma_{\phi, A}}(\P^-_{\phi, A})$.$\myqedhere$
\end{definition}

\begin{definition}\label{f covering} We say \textbf{lower part $(\phi, A)$-covering} holds if $(\phi, A)$ is full and $cf({\sf{ord}}(\P_{\phi, A}))\geq \omega_1$. $\myqedhere$
\end{definition}

\begin{notation}\label{qx notation}
Suppose now that lower part $(\phi, A)$-covering fails. Given $X\in \powerset_{\omega_1}(\P)$, we let
\begin{itemize}
\item $\P_X$ be the transitive collapse of $Hull^{\P_{\phi, A}}(X)$,
\item $\tau_X:\P_X\rightarrow \P_{\phi, A}$ be the inverse of the transitive collapse,
\item $\Sigma_X$ be the $\tau_X$-pullback of $\Sigma_{\phi, A}$,
\item $\d_X=\d^{\P_X}$.
\end{itemize}
$\myqedhere$
\end{notation}
\begin{remark}\label{important rem about sigma x} Thus, $\Sigma_X$ is a strategy that acts on stacks that are based on $\P_X|\d_X$. It follows that if $\P_X\models ``\d_X$ is a regular cardinal" then $\Sigma_X$ is (essentially\footnote{As defined, $\Sigma_X$ still does not act on iterations that are above $\d_X$.}) a strategy for $\P_X$.$\myqedhere$
\end{remark}

\begin{definition}[Weakly condensing set]\label{weakly condensing set}\index{weakly condensing set}
Suppose $(\phi, A)$-covering fails and set $\Gamma=\Gamma_{\phi, A}$, 
 $\P=\P_{\phi, A}$ and $\Sigma=\Sigma_{\phi, A}$. We say that $X\in \powerset_{\omega_1}(\P)$ is a \textbf{$(\phi, A)$-weakly condensing set} if $\P = Hull^\P(X\cup \delta^\P)$ and whenever $X\subseteq Y\in \powerset_{\omega_1}(\P)$,  $\Sigma_Y$ is a strongly $\Gamma$-fullness preserving iteration strategy with strong branch condensation. $\myqedhere$
\end{definition}

\begin{notation}\label{sigmaxy notation}
Suppose $(\phi, A)$-covering fails and set $\Gamma=\Gamma_{\phi, A}$, 
 $\P=\P_{\phi, A}$ and $\Sigma=\Sigma_{\phi, A}$.
Suppose $X\subseteq Y\in \powerset_{\omega_1}(\P)$.  Let $\tau_{X,Y}:\P_X\rightarrow \P_Y$ be $\tau_Y^{-1}\circ \tau_X$. Let \begin{itemize}
\item $\sigma_Y^{X,-} = \bigcup_{\R\inseg^c_{hod}<\P_Y} \pi^{\Sigma_Y}_{\R,\infty}$,
\item  $\sigma^X_Y: \P_Y\rightarrow \P$ be given by: for any $f\in \P_X$ and any $a\in (\P_Y|\delta_Y)^{<\omega}$, and $x=\tau_{X,Y}(f)(a)$,
\begin{center}
$\sigma^X_Y(x) = \tau_X(f)(\sigma^{X,-}_Y(a))$.
\end{center}
\end{itemize}
$\myqedhere$
\end{notation}
\begin{definition}\label{def: extension of} Suppose $(\phi, A)$-covering fails and set $\Gamma=\Gamma_{\phi, A}$, 
 $\P=\P_{\phi, A}$ and $\Sigma=\Sigma_{\phi, A}$.
Let $X\subseteq Y\in \powerset_{\omega_1}(\P)$. We say that $Y$ \textbf{extends} $X$ or $Y$ is an \textbf{extension of} $X$\index{extension} if 
\begin{enumerate}
\item $\tau_{X,Y}\rest (\P_X|\delta^{\P_X})$ is the iteration map via $\Sigma_X$,
\item letting $\nu=\sup\tau_{X, Y}[\d^{\P_X}]$, $\tau_Y\rest \P_Y|\nu$ is the iteration embedding according to $(\Sigma_X)_{\P_Y|\nu}$, and 
\item $\P_Y=Hull_1^{\P_Y}(\delta^{\P_Y}\cup\tau_{X,Y}[\P_X])$.
\end{enumerate}
$\myqedhere$
\end{definition}

\begin{definition}\label{def:honest_extension}\index{honest extension}
Suppose $(\phi, A)$-covering fails and set $\Gamma=\Gamma_{\phi, A}$, 
 $\P=\P_{\phi, A}$ and $\Sigma=\Sigma_{\phi, A}$. Suppose $Y$ is an extension of a weakly condensing set $X$. Let $\delta_Y=\delta^{\P_Y}$. We say that $Y$ is \textbf{an honest extension} of $X$ if 
\begin{enumerate}[(a)]
\item $\P_Y|\sup(\tau_{X, Y}[\d_X])$ is a $\Sigma_X$-iterate of $\P_X|\d_X$,
\item $\tau_{X, Y}\rest (\P_X|\delta_X) = \pi^{\Sigma_X}_{\P_X|\d_X, \P_Y|\sup(\tau_{X, Y}[\d_X])}$ and
\item $\sigma^X_Y$ is an elementary embedding\footnote{We clearly have that $\tau_X = \sigma^X_Y\circ \tau_{X,Y}$.}.
\end{enumerate}
$\myqedhere$
\end{definition}
\begin{remark}
$X$ is obviously an honest extension of itself, but there are other (non-trivial) honest extensions of $X$. For example, if $X = X'\cap \P$ where $X'\prec H^V_{\lambda}$ for some regular $\lambda$ (this will be the case for our intended $X$) and $Y = Y'\cap \P$ for some $X'\prec Y'$, then $Y$ is an honest extension of $X$. $\myqedhere$
\end{remark}
\begin{definition}[Condensing set]\label{def:condensing_set}\index{condensing set} Suppose $(\phi, A)$-covering fails and set $\Gamma=\Gamma_{\phi, A}$, 
 $\P=\P_{\phi, A}$ and $\Sigma=\Sigma_{\phi, A}$. 
Suppose $X\in \powerset_{\omega_1}(\P)$ is a $(\phi, A)$-weakly condensing set. We say that $X$ is a \textbf{$(\phi, A)$-condensing set} if whenever $Y$ extends $X$, $Y$ is an honest extension of $X$. 

We say that $X$ is a \textbf{strongly $(\phi, A)$-condensing set} if whenever $Y$ extends $X$, $Y$ is a $(\phi, A)$-condensing set.$\myqedhere$
\end{definition}


We expect that under many hypothesis such as $\sf{PFA}$  lower part $(\phi, A)$-covering fails. We also expect that under many hypothesis, failure of lower part $(\phi, A)$-covering implies the existence of $(\phi, A)$-condensing sets. In the next few chapters, we explore some specific situations where we know how to prove the existence of $(\phi, A)$-condensing sets. 

We finish by remarking that $(\phi, A)$ depends on the specific situation we are in. For instance, in \cite{CuBF}, $\phi$ isolates those hod pairs that have certain extendability and self-determining properties (see \cite[Definition 3.1, 3.5, 3.8]{CuBF}).

We finish here by showing that below $\sf{LSA}$, pullback strategies are unique.

\begin{lemma}[Uniqueness of strategies]\label{lem:uniqueness of strategies} Suppose $(\phi, A, X)$ is such that $\phi$ is a formula in the language of set theory,  $(\phi, A)$ is full, lower part $(\phi, A)$-covering fails and $X$ is a $(\phi, A)$-condensing set. Suppose further that whenever $(\Q, \Lambda)\in \Gamma_{\phi, A}$, $\Q$ is not of lsa type. Then whenever $Y_1$ and $Y_2$ are two honest extensions of $X$ such that $\P_{Y_1}=\P_{Y_2}$, then $\Sigma_{Y_1}=\Sigma_{Y_2}$. 
\end{lemma}
\begin{proof} Suppose that $\Sigma_{Y_1}\not =\Sigma_{Y_2}$. Let $\P_1=\P_{Y_1}$, $\P_2=\P_{Y_2}$, $\Phi_1=\Sigma_{Y_1}$ and $\Phi_2=\Sigma_{Y_2}$. Because we can trace disagreement of strategies to minimal disagreements (using our smallness assumption on hod mice)\footnote{See \rlem{disagreement implies low level disagreement}.}, we can find a minimal low level disagreement\footnote{See \rdef{low level disagreement between strategies}.} $(\T_1, \Q'_1, \T_2, \Q'_2, \R)$ between $\Phi_1$ and $\Phi_2$\footnote{See \rrem{important rem about sigma x}. It follows that $\T_1$ and $\T_2$ are based on a proper initial segment of $\P_{1}=\P_2$.}. Let $E$ be the $\R$-un-dropping extender of $\T_1$ and $\T_2$\footnote{See clause 5d of \rdef{low level disagreement between strategies}.}, and set for $i=1, 2$, $\W_i=Ult(\P_{i}, E)$. We thus have that\\\\
(1) $\W_1=\W_2$, $\R$ is of successor type and $(\Phi_1)_{\R^-}=(\Phi_2)_{\R^-}$.\\\\
Because both $Y_1$ and $Y_2$ are extensions of $X$, we have that both $\tau_{X, Y}\rest (\P_X|\d_X)$ and $\tau_{X, Z}\rest (\P_X|\d_X)$ are the iteration embedding according to $\Sigma_X$. Because $\Sigma_X$ has strong branch condensation and  is strongly $\Gamma_{\phi, A}$-fullness preserving, we have that $\tau_{X, Y}\rest (\P_X|\d_X)=\tau_{X, Z}\rest (\P_X|\d_X)$\footnote{See \rprop{positional}.}. Let then $\tau=_{def}\tau_{X, Y}\rest (\P_X|\d_X)=\tau_{X, Z}\rest (\P_X|\d_X)$. 


Next, because of the smallness assumption on hod pairs in $\Gamma_{\phi,A}$, it follows from $(\phi,A)$-condensation of $X$ that\\\\
(2) for $i\in \{1,2\}$, $\sup (Hull^{\W_i}(\pi_{E}\circ \tau[\P_X])\cap \delta^\R)=\d^\R$\footnote{It is easier to first establish that $\sup (Hull^{\W_i}(\pi_{E}\circ \tau[\P_X]\cup \d^{\R^-})\cap \delta^\R)=\d^\R$. See \rlem{the canonical singularizing sequence exists}.}.\\\\
Set for $i=1, 2$, $\X_i=\T_i^\frown \{E\}$. We can now find, using \rthm{comparison holds}, a normal stack $\U_1$ on $\W_1$ according to $(\Phi_1)_{\W_1, \X_1}$ and a normal stack $\U_2$ on $\W_2$ according to $(\Phi_2)_{\W_2, \X_2}$ such that setting $b_1=(\Phi_1)_{\W_1, \X_1}(\U_1)$, $b_2=(\Phi_2)_{\W_2, \X_2}(\U_2)$, $\R_1=\M^{\U_1}_{b_1}$ and $\R_2=\M^{\U_2}_{b_2}$ then setting $\Psi_1=(\Phi_1)_{\R_1, \X_1^\frown \U_1^\frown\{b\}}$ and $\Psi_2=(\Phi_2)_{\R_2, \X_2^\frown \U_2^\frown\{b_2\}}$,\\\\
(3) for $i\in \{1, 2\}$, $\U_i$ is based on $\R$, $\downarrow(\U_1, \R)=\downarrow(\U_2, \R)$, $b_1\not =b_2$ and $\pi^{\U_1}_{b_1}(\R)=\pi^{\U_2}_{b_2}(\R)$\\
(4) letting $\S=\pi^{\U_1}_{b_1}(\R)$, $(\Psi_1)_{\S}=(\Psi_2)_{\S}$.\\\\
Notice now that we can find $k_1:\R_1\rightarrow \P$ and $k_2:\R_2\rightarrow \P$ such that letting for $i=1,2$, $\tau_{Y_i}=\tau_i$,\\\\
(5) for $i=1,2$, $\tau_i=k_i\circ (\pi^{\U_i}_{b_i}\circ \pi_E)$,\\
(6) for $i=1, 2$, $k_i\rest \S=\pi^{(\Psi_i)_{\S}}_{\S, \infty}$\footnote{(5) and (6) easily follows from the fact that $\T_1$ and $\T_2$ are based on proper initial segment of $\P_1|\d^{\P_1}$.}, and\\
(7) for $i=1, 2$, $\pi^{\U_i}_{b_i}\circ \pi_E\circ \tau$ is according to $\Sigma_X$.\\\\
It follows that letting for $i=1,2$,  $Z_i=Y_i\cup \rge(\pi^{\Psi_i}_{\S, \infty})$, $Z_i$ extends $X$, and moreover, because $X$ is a condensing set, for $i=1, 2$, $k_i=\sigma^X_{Z_i}$\footnote{Again, this easily traces back to the fact that $\T_i$ is based on a proper initial segment of $\P_i|\d^{\P_i}$.}. 

Notice that it follows from (4) that $k_1\rest \S=k_2\rest \S$. Also, notice that\\\\
(8) $k_1\rest (Hull^{\R_1}(\S^-\cup \pi^{\U_1}_{b_1}\circ \pi_E\circ \tau [\P_X]))=k_2\rest (Hull^{\R_2}(\S^-\cup \pi^{\U_2}_{b_2}\circ \pi_E\circ \tau [\P_X]))$.\\\\
Combining (2) and (8) we get that (using (3))\\\\
(9) $\rge(\pi^{\U_1}_{b_1})\cap \rge(\pi^{\U_2}_{b_2})$ is cofinal in $\d^\S$. \\\\
Clearly (9) and parts of (3)\footnote{That for $i\in \{1, 2\}$, $\U_i$ is based on $\R$, $\downarrow(\U_1, \R)=\downarrow(\U_2, \R)$.} imply that $b_1=b_2$, while other parts of (3) state that $b_1\not=b_2$.
\end{proof}

The following is a useful corollary of the definition of a condensing set. We will apply this corollary in many applications later.
\begin{corollary}\label{cor:strategy_condensation}
Suppose $Y\prec Z$ are extensions of a $(\phi,A)$-condensing set $X$ and $Z$ is an extension of $Y$. Suppose $B\in \powerset(\delta^\P)\cap \P$ and $B\in Y$. Let $a\in  (\delta^{\Q_Y})^{<\omega}$. Then $\pi^{\Sigma_Y}_{\Q_Y,\infty}(a)\in B$ if and only if $\pi^{\Sigma_Z}_{\Q_Z,\infty}(\tau_{Y,Z}(a))\in B$.
\end{corollary}

\section{Condensing sets from elementary embeddings}\label{sec:condensing_zfc}

The following two theorems can be proved using the proof of \cite[Lemma 11.15]{CuBF}. First we introduce some terminology.

\begin{terminology}\label{homogenous formula} Suppose $\kappa$ is an inaccessible cardinal and $G\subseteq Col(\omega, <\kappa)$ is $V$-generic. Suppose $(\phi, A)$ is such that $V[G]\models ``(\phi, A)$ is full and lower part $(\phi, A)$-covering fails".   We say $(\phi, A)$ is \textbf{homogenous} if $\P_{\phi, A}\in V$, $\Sigma_{\phi, A}\rest V\in V$ and for any $(\Q, \Lambda)\in \mathcal{F}_{\phi, A}$, there is $(\R, \Psi)\in \mathcal{F}_{\phi, A}$ such that $\R\in V$, $\Psi\rest H_{\kappa}^V\in V$ and $V[G]\models \Gamma(\Q, \Lambda)\subseteq \Gamma(\R, \Psi)$. $\myqedhere$
\end{terminology}

\begin{theorem}\label{existence of condensing sets in n} Suppose $N\subseteq M$ are transitive models of set theory and $j:M\rightarrow N$ is an elementary embedding with critical point $\k$ such that $j$ is amenable to $M$, i.e., for every $X\in M$, $j(X)\in M$. Suppose $g\subseteq Coll(\omega, <j(\k))$ is $N$-generic. Let $j^+:M[g_\kappa]\rightarrow N[g]$ be the extension of $j$ where for $\a<j(\k)$, $g_\a=g\cap Coll(\omega, <\a)$. Suppose $\phi$ is a formula in the language of set theory and $A\in M[g]$. Suppose further that $M[g_\kappa]\models ``(\phi, A)$ is full, $(\phi, A)$ is homogenous and lower part $(\phi, A)$-covering fails". Then $j[\P_{\phi, A}]$ is a strongly $(\phi, j(A))$-condensing set in $N[g]$. Hence, 
 $M[g]\models ``$there is a strongly $(\phi, A)$-condensing set".
\end{theorem}

\begin{terminology}
We say $(\phi, A)$ is \textbf{maximal} if there is no hod pair or an sts hod pair $(\Q, \Lambda)$ such that $\Q$ is of limit type, $\Lambda$ has strong branch condensation and is strongly $\Gamma_{\phi, A}$-fullness preserving and $\Gamma(\Q, \Lambda)=\Gamma_{\phi, A}$.$\myqedhere$
\end{terminology} 

\begin{theorem}\label{regularity from condensing set} Assume $\sf{ZF}+DC$ and suppose $(\phi, A)$ is maximal and full, lower part $(\phi, A)$-covering fails and $X$ is a $(\phi, A)$-condensing set. Then $\P_{\phi, A}\models ``\d^{\P_{\phi, A}}$ is regular". 
\end{theorem}

We will not prove \rthm{regularity from condensing set} but will give a fairly complete proof of \rthm{existence of condensing sets in n}. \\\\
\textbf{The proof of \rthm{existence of condensing sets in n}.}\\\\
We fix $(M, N, j, \kappa, g, \phi, A)$ as in the statement of \rthm{existence of condensing sets in n}. The proof follows the proof of   \cite[Theorem 10.3]{CuBF}. Throughout this section we will use the following notation:
\begin{notation}\label{not: notation for this section}
Working in $M[g_\kappa]$, let 
\begin{itemize}
\item $\P^-=\P^-_{\phi, A}$, 
\item $\P=\P_{\phi, A}$, 
\item $\Sigma=\Sigma_{\phi, A}$, 
\item $\mathcal{F}=\mathcal{F}_{\phi, A}$,
\item $\Gamma=\Gamma_{\phi, A}$.
\end{itemize}
$\myqedhere$
\end{notation}

\begin{theorem}\label{thm:weakly_condensing_set}
$N[g]\models ``j[\P]$ is a weakly condensing set".
\end{theorem}
\begin{proof}
Notice that that $j[{\sf{ord}}(\P)]$ is cofinal in ${\sf{ord}}(j(\P))$. Below, we often confuse strategies with their interpretations in relevant generic extensions or in relevant inner models. However, in some cases, the distinction between the two strategies is important, and in those situations we will either separate the two strategies or point out that the distinction is important. Also, below if $Y\in \powerset_{\omega_1}(j(\P))$ then we let $\P_Y=j(\P)_Y$ and $\Sigma_Y=j(\Sigma)_Y$.

We want to show that \\\\
(a) if $Y\in (\powerset_{\omega_1}(j(\P)))^{N[g]}$ is such that $j[\P]\subseteq Y$ then $L(j^+(\Gamma))\models ``\P_Y$ is $\Sigma_Y$-full"\footnote{I.e., $\P_Y={\sf{Lp}}^{j^+(\Gamma), \Sigma_Y}(\P_Y|\d_Y)$.}.\\\\
Towards a contradiction assume that (a) is false. Notice that if $Y$ witnesses that (a) is false then $\P_Y$ may not be in $M[g_\kappa]$. Fix one such $Y$ that is a counterexample to (a), and let $\M$ be a sound $\Sigma_Y$-mouse over $\P_Y|\d_Y$ that has an iteration strategy in $j^+(\Gamma)$ but such that $\M\not \insegeq \P_Y$ and $\rho(\M)=\d_Y$. Let $\dot{Y}\in N^{Coll(\omega, <(\kappa, \l))}$ be a name for $Y$ and $\dot{\M}$ be a name for $\M$.  We can then find some $\Sigma_{j[\P]}$-hod pair $(\P^+, \Pi)\in N$ and a hod pair $(\S, \Phi)\in N$ such that
\begin{enumerate}
\item $\P^+\in H_{j(\kappa)}^N$,
\item $\Pi$ has strong branch condensation, 
\item $\P^+$ is meek and of limit type, 
\item $\cf^{\P^+}(\d^{\P^+})=\omega$,
\item $(Y\cap j(\P|\d^\P))\subseteq \rge(\pi^{\Phi}_{\S, \infty})$ and no proper complete layer of $\S$ has this property\footnote{I.e., if $\S'\inseg^c_{hod}\S$ then $(Y\cap j(\P|\d^\P))\not \subseteq \rge(\pi^{\Phi}_{\S', \infty})$. See \rnot{complete layer notation}.},
\item $\Pi\in N$ is a $(j(\kappa), j(\kappa))$-strategy for $\P^+$ that can be uniquely extended to a strategy $\Pi^{g}\in j^+(\Gamma)$\footnote{$\Pi$ can be obtained by computing the direct limit of all those hod pairs $(\Q, \Lambda)\in j^+(\mathcal{F}_{\phi, A})$ with the property that $\Q\in N[g_\nu]$ where $\nu\in (\kappa, j(\kappa))$ is chosen in a way that $\Sigma_Y$ and the strategy of $\M$ appear in $N[g_\nu]$. We might then have to take some initial segment of this direct limit to satisfy clauses 3-5 above.}, and moreover, $\Pi$ witnesses that $\P^+$ is a $\Sigma_{j[\P]}$-hod mouse,
\item $N[g_\kappa]\models$ it is forced by $Coll(\omega, <(\kappa, j(\k))$ that 
\begin{enumerate}
\item $\dot{\M}$ is a sound $\Sigma_{\dot{Y}}$-mouse over $\P_{\dot{Y}}|\d_{\dot{Y}}$ that projects to $\d_{\dot{Y}}$.
\item $\dot{\M}$ has an iteration strategy in the derived model of $(\P^+, \Pi)$\footnote{We confuse the extension of $\Pi$ to this extension with $\Pi$-itself.} as computed by any $\bR$-genericity iteration,
\item $\Phi$ is in the derived model of $(\P^+, \Pi)$ as computed by any $\bR$-genericity iteration,
\item $\dot{\M}$ is not an initial segment of $\P_{\dot{Y}}$. 
\end{enumerate}
\end{enumerate}
Because $\P^+$ might have cardinality $>\kappa$, when we form $\P_Y^+=_{def}Ult(\P^+, E)$, where $E$ is the $(\cp(\tau_{j[\P], Y}), \d_Y)$-extender derived from $\tau_{j[\P], Y}$, we cannot conclude that $\P_Y^+$ is iterable in $N[g]$. This is because we do not know that $j\rest \P^+\in N$. To resolve this issue we take a hull of size $\kappa$. Let $\kappa_1=(\kappa^+)^M$.

We work in $N[g_\kappa]$. We can now find $\pi: W[g_\kappa]\rightarrow (H_{j(\kappa_1)})^{N[g_\kappa]}$ (in $N[g_\kappa]$) such that
\begin{itemize}
\item $W\in M$ is transitive and $\kappa+1\subseteq W$, 
\item $(j(\P), j\rest \P, \dot{Y}, (\P^+, \Pi), (\S, \Phi))\in \rge(\pi)$.
\end{itemize}
Let $\dot{Z}=\pi^{-1}(\dot{Y})$, $\dot{\N}=\pi^{-1}(\dot{\M})$, $\R=\pi^{-1}(j(\P))$ and $k: \P\rightarrow \R$ be $\pi^{-1}(j\rest \P)$. Working in $N[g_{\kappa_1}]$, let $h\subseteq Coll(\omega, <(\kappa, k(\k))$ be $W$-generic, and set 
\begin{itemize}
\item $Z=\dot{Z}_h$, $\dot{\N}_h=\N$, $\Q=(\P_{Z})^{W[h]}$, 
\item $\sigma=(\tau_{k[\P], Z})^{W[h]}$ and $\tau=(\tau_Z)^{W[h]}$,
\item $\overline{\P^+}=\pi^{-1}(\P^+)$ and $\overline{\Pi}=\pi^{-1}(\Pi)$,
\item $(\overline{\S}, \overline{\Phi})=\pi^{-1}(\S, \Phi)$.
\end{itemize}
Thus, we have that\\\\
(A) $k=\tau\circ \sigma$, $\sigma:\P\rightarrow \Q$ and $\tau:\Q\rightarrow \R$,\\\
(B) in $W[g_\kappa*h]$, 
\begin{enumerate}
\item $\N$ is a sound $\Sigma_{Z}$-mouse over $\Q|\d^\Q$ that projects to $\d^\Q$.
\item  in any derived model of $(\overline{\P^+}, \overline{\Pi})$ as computed by an $\bR$-genericity iteration, $\N$ has an $\omega_1$-iteration strategy witnessing that it is a $\Sigma_Z$-mouse,
\item $\N$ is not an initial segment of $\Q$. 
\item $\overline{\Phi}$ is in the derived model of $(\overline{\P^+}, \overline{\Pi})$ as computed by any $\bR$-genericity iteration,
\item letting $\xi:\Q|\d^\Q\rightarrow \overline{\S}|\d^{\overline{\S}}$ be such that $\xi=(\pi^{\overline{\Phi}}_{\overline{\S}, \infty})^{-1} \circ \tau$, $\Sigma_Z=(\xi$-pullback of $\overline{\Phi}_{\overline{\S}|\d^{\overline{\S}}}$. 
\end{enumerate}
Let now $F$ be the $(\cp(\sigma), \d^\Q)$-extender derived from $\sigma$, and set $\Q^+=Ult(\overline{\P^+}, F)$. Let $\sigma^+=\pi_F^{\P^+}$. Notice that because $\pi\circ k=j\rest \P$, we have $\phi^+:\Q^+\rightarrow j(\overline{\P^+})$ such that\\\\ (C) $j\rest \overline{\P^+}=\phi^+\circ \sigma^+$.\\\\
Let $\overline{\Pi}^+$ be the $\pi\rest \overline{\P^+}$-pullback of $\Pi$\footnote{We confuse $\Pi$ with its extension to $N[g]$. Similarly, we think of $\overline{\Pi}^+$ as a strategy in $N$ as well as in $N[g]$. Same comment applies below to $\overline{\Pi}$ and $\overline{\Phi}$.} and let $\overline{\Phi}^+$ be the $\pi$-pullback of $\Phi$. Notice that\\\\
(D1) $\overline{\Pi}^+\rest {\sf{HC}}^{W[g_\kappa*h]}= \overline{\Pi}$\footnote{See proof of Claim 2 in the proof of \cite[Lemma 10.4]{CuBF}. The same equation for $W[g]$ follows easily from hull condensation of $\overline{\Pi}^+$, but this equation for $W[g_\kappa*h]$ needs more work.},\\
(D2) $\overline{\Pi}^+$ witnesses that $\overline{\P^+}$ is a $\Sigma$-hod mouse\footnote{This follows from the fact that $\Pi$ witnesses that $\P^+$ is a $\Sigma$-hod mouse and $\pi\rest \P=id$.},\\
(D3) $\overline{\Phi}^+\rest {\sf{HC}}^{W[g_\kappa*h]}= \overline{\Phi}$.\\\\
Notice now that we have \\\\
(F) in $N[g]$, $j^+(\overline{\Pi}^+\rest (H_{\kappa_1}^{M[g_\k]}))$ is a $(j(\kappa), j(\kappa))$-iteration strategy witnessing that $j(\overline{\P^+})$ is a $j(\Sigma)$-hod mouse, and moreover, $j\rest \overline{\P^+}\in N[g]$\footnote{Because $\card{\overline{\P^+}}^M=\k$.}.\\\\
We let $\Psi=(\Sigma_Z)^{W[g_\kappa*h]}$. Notice that in $W[g_\kappa*h]$, $\Psi$ is the $\tau$-pullback of $\pi^{-1}(j(\Sigma))$. Let $\Psi^+$ be the $\phi^+\rest (\Q|\d^\Q)=\pi\circ \tau\rest (\Q|\d^\Q)$-pullback of $j(\Sigma)$. It follows that\\\\
(G) $\Psi^+$ is the $\pi\circ \xi$-pullback of $\Phi$, and it is also $\xi$-pullback of $\overline{\Phi}^+$.\\\\
We now claim that\\\\
(b) in $N[g]$, in any derived model of $(\overline{\P^+}, \overline{\Pi}^+)$ as computed by an $\bR$-genericity iteration, $\N$ has an $\omega_1$-iteration strategy witnessesing that $\N$ is a $\Psi^+$-mouse.\\\\
The proof of (b) is like the proof of Claim 1 of \cite[Lemma 10.4]{CuBF}. We outline it below.
Working in $W[g_\kappa]$, let $\W=\M_\omega^{\#, \overline{\Pi}, \overline{\Phi}}$ and let $\Lambda$ be the unique iteration strategy of $\W$. Because $\pi(\W)=\M_\omega^{\#, \Pi, \Phi}$, we have that letting $\Lambda^+$ be the $\pi$-pullback of $\pi(\Lambda)$,\\\\
(H) $\W=\M_\omega^{\#, \overline{\Pi}^+, \overline{\Phi}^+}$, and $\Lambda^+$ witnesses that $\W$ is a $\overline{\Pi}^+\oplus \overline{\Phi}^+$-mouse. \\\\
Working in $W[g_\kappa*h]$ and using (B2), we can find a $\Lambda$-iterate $\W_1$ of $\W$, a Woodin cardinal $\eta$ of $\W_1$ and $\W_1$-generic $m\subseteq Coll(\omega, \eta)$ such that letting $\l$ be the sup of the Woodin cardinals of $\W_1$,\\\\
(I1) $\N, \xi\in \W_1[m]$,\\
(I2) $\W_1[m]\models ``$the derived model at $\l$ satisfies that any derived model of $(\overline{\P^+}, \overline{\Pi})$ as computed by an $\bR$-genericity iteration has an $\omega_1$-iteration strategy for $\N$ witnessing that $\N$ is a mouse relative to the $\xi$-pullback of $\overline{\Phi}$".\\\\
Let then $\W_\infty$ be a $\Lambda^+$-iterate of $\W_1$ which is obtained via some $\bR^{\N[g]}$-genericity iteration in such a way that letting $i:\W_1\rightarrow \W_\infty$ be the iteration embedding, $\cp(i)>\eta$.  It then follows from (I1), (I2) and (H) that \\\\
(J) $\W_\infty[\N, \xi]\models ``$the derived model at $\l$ satisfies that any derived model of $(\overline{\P^+}, \overline{\Pi}^+)$ as computed by an $\bR$-genericity iteration has an $\omega_1$-iteration strategy for $\N$ witnessing that $\N$ is a mouse relative to the $\xi$-pullback of $\overline{\Phi}^+$".\\\\
(b) now easily follows from (J), (H) and (G).

To finish the proof of \rthm{thm:weakly_condensing_set}, it remains to implement the last portion of the proof of \cite[Theorem 10.3]{CuBF}. Let $\Delta_0$ be $\phi^+$-pullback of $j^+(\overline{\Pi}^+\rest (H_{\kappa_1}^{M[g_\k]}))$. Notice that it follows from (F) that $\Delta_0$ witnesses that $\Q^+$ is a $\Psi^+$-hod mouse. It then follows from (b) that\\\\
(K) in $N[g]$, in any derived model of $(\Q^+, \Delta_0)$ as computed by an $\bR$-genericity iteration, $\N$ has an $\omega_1$-iteration strategy $\Delta$ witnessing that $\N$ is a $\Psi^+$-mouse.\\\\
(K) gives contradiction, as it implies that\\\\
(L) $\Q^+\models ``{\sf{ord}}(\Q)$ is not a cardinal"\footnote{This is because (K) implies that $\N$ is ordinal definable in $\Q^+$ and therefore, $\N\in \Q$.},\\\\
while clearly $\overline{\P^+}\models ``{\sf{ord}}(\P)$ is a cardinal", contradicting the elementarity of $\phi^+$
\end{proof}

The main theorem of this chapter is.
\begin{theorem}\label{thm:condensing_set}
In $N[g]$, $j[\P]$ is a strongly condensing set.
\end{theorem}
\begin{proof}
We will show that $j[\P]$ is a condensing set. A very similar proof, which is only notationally more complicated, shows that $j[\P]$ is strongly condensing. To prove the theorem, we need the following definition, due to the first author (cf. \cite{CuBF} or \cite{Trang2015PFA}). The proof is based on \cite[Lemma 11.15]{CuBF}. For completeness, we give a fairly detailed argument here. The reader may wish to recall \rnot{complete layer notation}.  Below if $Y\in \powerset_{\omega_1}(j(\P))$ then we let $\P_Y=j(\P)_Y$ and $\Sigma_Y=j(\Sigma)_Y$.

We work in $N[g]$. Suppose $X\in \powerset_{\omega_1}(j(\P))$ is a weakly condensing set and $B\in X \cap \powerset(\delta^{j(\P)})$. We say that $X$ has \textit{$B$-condensation} if whenever $Y\in \powerset_{\omega_1}(j(\P))$ is such that $X\prec Y$, $\tau_{X, Y}(T_{X,B}) = T_{Y,B}$, where  for $Z\in \powerset_{\omega_1}(j(\P))$,
\begin{center}
$T_{Z,B} = \{(\varphi,s) \ | \ s\in [\delta^\R]^{<\omega} \textrm{ for some } \R\inseg^c_{hod}\P_Z \wedge j(\P) \vDash \varphi[\pi^{\Sigma_{Z}}_{\R,\infty}(s), B]\}$.
\end{center}
We say $X$ has \textit{term condensation} if it has $B$-condensation for every $B\in X \cap \powerset(\delta^{j(\P)})$.

To prove that a weakly condensing set $X$ is condensing, it is enough to prove that $\tau_X$ has term condensation. It is not hard to show that if for every $A\in j(\P)$ there is $X$ with $A$-condensation then $j[\P]$ has term condensation\footnote{See the proof of \cite[Lemma 11.15]{CuBF}. This essentially follows from the elementarity of $j$.}. We say, working in $N[g]$, $X\in \powerset_{\omega_1}(j(\P))$ is \textit{good} if for every $\R\inseg^c_{hod}\P_X$, $\tau_X\rest \R=\pi^{\Sigma_X}_{\R, \infty}$. It follows from \cite[Lemma 11.9]{CuBF} that the set of good $X$ is a club. Notice that $j[\P]$ is good.

Towards a contradiction, assume that (in $N[g]$) there is a set $A\in \P$ such that no $X\in \powerset_{\omega_1}(j(\P))$ with the property $j[\P]\subseteq X$, has $A$-condensation. We now fix such a set $A$. We say (in $N[g]$) that a tuple $\{\langle \P_i,\Q_i, X_i,  Y_i, \xi_i,\pi_i, \phi_i \ | \ i<\omega \rangle, B, \M \}$ is a \textit{bad tuple} (relative to $A$) if
\begin{enumerate}
\item $X_0=Y_0=j[\P]$,
\item for all $i<\omega$, $X_i\in \powerset_{\omega_1}(j(\P))$ is good,
\item  for all $i<\omega$, $\P_i = \P_{X_i}$  and $\Q_i = \P_{Y_i}$;
\item for all $i < j<\omega$, $X_i \prec Y_i \prec X_j$;
\item for all $i<\omega$, $\xi_i=\tau_{X_i, Y_i}$, $\pi_i=\tau_{Y_i, X_{i+1}}$ and $\phi_{i} =\tau_{X_i,X_{i+1}}$\footnote{Thus, $\xi_i:\P_i\rightarrow \Q_i$, $\pi_i: \Q_i \rightarrow \P_{i+1}$ and $\phi_i= \pi_i \circ \xi_i$.};
\item $\M\in j(\P|\d^\P)$ and letting $\eta=\sup_{i<\omega}(X_i\cap \d^{j(\P)})$, $\M|\eta=j(\P)|\eta$;
\item letting $\eta$ be as above, for every formula $\phi$ and for every $s\in \eta^{<\omega}$, $\M\models \phi[B, s]$ if and only if $j(\P)\models \phi[A, s]$; 
\item for all $i\in [1, \omega)$, $\xi_i(T_{X_i, A}) \neq T_{Y_i,A}$. 
\end{enumerate}

\begin{claim}\label{claim:bad_tuple} There is a bad tuple.
\end{claim}
\begin{proof}
It is easy to construct a bad tuple $\{\langle \P_i,\Q_i, X_i,  Y_i, \xi_i,\pi_i, \phi_i \ | \ i<\omega \rangle, A, j(\P) \}$ with $j(\P)$ playing the role of $\M$ and $A$ playing the role of $B$. Once this is done, letting $\eta=\sup(X_i\cap \d^{j(\P)})$, we set $\M=cHull^{j(\P)}(\{A\}, \eta)$. $B$ then is the transitive collapse of $A$.


\end{proof}
Fix a bad tuple $\mathcal{A}^*=\{\langle \P_i,\Q_i, X_i,  Y_i, \xi_i,\pi_i, \phi_i \ | \ i<\omega \rangle, B, \M\}$. 
Let 
\begin{itemize}
\item $\eta=\sup(X_i\cap \d^{j(\P)})$, 
\item $Z_i=X_i\cap \eta$, $W_i=Y_i\cap \eta$, 
\item $\Phi_i=\Sigma_{X_i}$ and $\Psi_i=\Sigma_{Y_i}$,
\item $T_i=\{ (\phi, s): s\in [\d^{\P_i}]^{<\omega} \wedge \M\models \phi[B, \pi^{\Phi_i}_{\P_i|\d^{\P_i}, \infty}(s)]\}$,
\item $S_i=\{ (\phi, s): s\in [\d^{\Q_i}]^{<\omega} \wedge \M\models \phi[B, \pi^{\Psi_i}_{\Q_i|\d^{\Q_i}, \infty}(s)]\}$.
\end{itemize}
and  set
\begin{center}
$\mathcal{A}=\{\langle \P_i, \Q_i, \Phi_i, \xi_i, \pi_i, \phi_i, T_i, S_i\ | \ i<\omega \rangle, B, \M \}$
\end{center}
Notice that it follows that for all $i<\omega$, $T_i=T_{X_i, A}$ and $S_i=T_{Y_i, A}$. 
Let $C\in j^+(\Gamma)$ be such that $\eta<\Theta^{L(C, \bR)}$ and $\M\subseteq \H^{L(C, \bR)}$. Then because $\Phi_i$ and $\Psi_i$ can be recovered from $j(\Sigma)_{\M|\eta}$ and respectively $Z_i$ and $W_i$, $\mathcal{A}\in L(C, \bR)$ and $L(C, \bR)\models *(\mathcal{A})$ where $*(\mathcal{A})$ is the conjunction of the following clauses:
\begin{enumerate}
\item for all $i<\omega$, $\xi_i:\P_i\rightarrow \Q_i$, $\pi_i: \Q_i \rightarrow \P_{i+1}$ and $\phi_i= \pi_i \circ \xi_i$;
\item for all $i<\omega$, $T_i\in \P_i$ and $\xi_i(T_i) \neq S_i$;
\item for all $i<\omega$, letting $\Psi_i$ be the $\pi_i$-pullback of $\Phi_{i+1}$, $S_i=\{ (\phi, s): s\in [\d^{\Q_i}]^{<\omega} \wedge \M\models \phi[B, \pi^{\Psi_i}_{\Q_i|\d^{\Q_i}, \infty}(s)]\}$;
\item for all $i<\omega$, $T_i=\{(\phi, s): s\in [\d^{\P_i}]^{<\omega}$ and $\M\models \phi[B, \pi^{\Phi_i}_{\P_i|\d^{\P_i}, \infty}(s)]\}$;
\item $\eta<\Theta$ and $\M\subseteq \H^{L(C, \bR)}$.
\end{enumerate}
Notice that $\P_0=\P$ and $\Phi_0=\Sigma$. Let now $(\P_0^+,\Pi_0)$ be a $\Phi_0$-hod pair such that 
\begin{center}
$L(\Gamma(\P_0^+,\Pi_0), \bR) \vDash *(\mathcal{A})$.
\end{center}
We may also assume $(\P_0^+,\Pi_0\rest N)\in N$, $\P_0^+$ is of limit type and $\cf^{\P_0^+}(\d^{\P_0^+})$ is not a measurable cardinal of $\P_0^+$. This type of reflection is possible because we replaced $j(\P)$ by $\M$. Let $u=(\phi_i, T_i: i<\omega)$ and set \begin{center}
$\W = \M_{\omega}^{\sharp, \Pi_0, \oplus_{i<\omega}\Phi_{i}}(u)$.\end{center} 
Let $\Lambda$ be the unique strategy of $\W$ witnessing that $\W$ is a $\Pi_0\oplus(\oplus_{i<\omega}\Phi_{i})$-mouse over $u$. We now have that\\\\
(A) in $N[g]$, whenever $D$ is obtained as a derived model of  $(\W, \Lambda)$ via some $\mathbb{R}$-genericity iteration,
\begin{center}
$D \vDash *(\mathcal{A})$.\\
\end{center} 
We remark that the following objects are in $N$:
\begin{itemize}
\item $\M$, $\W$ and $\Lambda\rest N$. 
\item $\langle \P_i, \Phi_i\rest N, \phi_i, T_i\ | \ i<\omega \rangle$.
\end{itemize}
However, $(\Q_i, \xi_i, \pi_i, S_i)$ are not in $V$. Notice also that the objects listed above are in $D(\mathcal{Z}, \omega_1^{N[g]}, h)$. We set $\mathcal{B}= \{\langle \P_i, \Phi_i, \phi_i, T_i\ | \ i<\omega \rangle, B, \M\}$ and given 
$t=(\N_i, \psi_i, \sigma_i, U_i)$ we write $\mathcal{B}^t$ for the set  $\{\langle \P_i, \N_i, \Phi_i, \psi_i, \sigma_i, \phi_i, T_i, U_i\ | \ i<\omega \rangle, B, \M\}$.
Thus, we have the following:\\\\
(B) in $N[g]$,  whenever $D$ is obtained as a derived model of  $(\W, \Lambda)$ via some $\mathbb{R}$-genericity iteration,
\begin{center}
$D \vDash$ ``there is $t=(\N_i,  \psi_i, \sigma_i, U_i)$ such that $*(\mathcal{B}^t)$".\\
\end{center} 

Let now $\pi: W[g_\kappa]\rightarrow (H_{j(\kappa_1)}^{N[g_\k]})$ be such that all the relevant objects are in the range of $\pi$, $W\in N$, $\card{W}^N=\kappa$, $j\rest \P\in \rge(\pi)$ and $\cp(\pi)>\kappa$. By``all relevant objects" we mean those objects that are in $N$, and in particular those listed above. For $a\in \rge(\pi)$, let $\overline{a}=\pi^{-1}(a)$. For $i<\omega$, let $\overline{\Phi_i}^+$ be the $\pi$-pullback of $\Phi_i$, and also Let $\overline{\Lambda}^+$ be the $\pi$-pullback of $\Lambda$. 

We thus have that\\\\
(C) $\overline{\W}=\M_{\omega}^{\sharp, \overline{\Pi_0}^+, \oplus_{i<\omega}\overline{\Phi_{i}}^+}$, and $\overline{\Lambda}^+$ witnesses that $\overline{\W}$ is a $\overline{\Pi_0}^+ \oplus (\oplus_{i<\omega}\overline{\Phi_{i}}^+)$-mouse,\\
(D) in $N[g]$, whenever $D$ is obtained as a derived model of $(\overline{\W}, \overline{\Lambda}^+)$ via some $\mathbb{R}$-genericity iteration, there is $\S\in \H^D$ and $F\in \S$ such that letting $\mathcal{B}_0=\{\langle \overline{\P_i}, \overline{\Phi_i}^+, \overline{\phi_i}, \overline{T_i}\ | \ i<\omega \rangle, F, \S\}$
\begin{center}
$D \vDash$ ``there is $t=(\N_i, \psi_i, \sigma_i, U_i)$ such that $*(\mathcal{B}_0^t)$",
\end{center}  
(E) in $N[g]$, whenever $D$ is obtained as a derived model of $(\overline{\P_0^+}, \overline{\Pi_0}^+)$ via some $\mathbb{R}$-genericity iteration, there is $\S\in \H^D$ and $F\in \S$ such that letting $\mathcal{B}_0=\{\langle \overline{\P_i}, \overline{\Phi_i}^+, \overline{\phi_i}, \overline{T_i}\ | \ i<\omega \rangle, F, \S\}$
\begin{center}
$D \vDash$ ``there is $t=(\N_i, \psi_i, \sigma_i, U_i)$ such that $*(\mathcal{B}_0^t)$".\\
\end{center}  

The proof of (E) is like the proof of (b) in the proof of \rthm{thm:weakly_condensing_set}. Notice that in $N[g]$, the derived models of $(\overline{\P_0^+}, \overline{\Pi_0}^+)$ obtained via $\bR$-genericity iteration have the form $D=_{def}L(\Gamma(\overline{\P_0^+}, \overline{\Pi_0}^+))$.  Fix then some $(F, \S)\in D$ and $t=(\N_i,  \psi_i, \sigma_i, U_i) \in D$  such that letting $\mathcal{B}_0=\{\langle \overline{\P_i}, \overline{\Phi_i}^+, \overline{\phi_i}, \overline{T_i}\ | \ i<\omega \rangle, F, \S\}$,
$D \vDash *(\mathcal{B}_0^t)$.

We thus have that the following clauses hold:
\begin{enumerate}
\item for all $i<\omega$, $\psi_i:\overline{\P_i}\rightarrow \N_i$, $\sigma_i: \N_i \rightarrow \overline{\P_{i+1}}$ and $\overline{\phi_i}= \sigma_i \circ \psi_i$;
\item for all $i<\omega$, $\overline{T_i}\in \overline{\P_i}$ and for all $i\in [1, \omega)$, $\psi_i(\overline{T_i}) \neq U_i$;
\item for all $i<\omega$, $U_i=\{ (\phi, s): s\in [\d^{\N_i}]^{<\omega} \wedge \S\models \phi[F, \pi^{\Sigma_i}_{\N_i|\d^{\N_i}, \infty}(s)]\}$ where $\Sigma_i$ is the $\sigma_i$-pullback of $\overline{\Phi_i}^+$;
\item for all $i<\omega$, $\overline{T_i}=\{\phi, s)\in : s\in [\d^{\overline{\P_i}}]^{<\omega}$ and $\S\models \phi[F, \pi^{\overline{\Phi_i}^+}_{\overline{\P_i}|\d^{\overline{\P_i}}, \infty}(s)]\}$;
\item $\S\in \H^{L(\Gamma(\overline{\P_0^+}, \overline{\Pi_0}^+))}$.
\end{enumerate}

Now we define by induction $\psi_i^+: \overline{\P_i^+} \rightarrow \N_i^+$, $\sigma_i^+: \N_i^+ \rightarrow \overline{\P_{i+1}^+}$, $\overline{\phi_{i}^+}: \overline{\P_i^+}\rightarrow \overline{\P_{i+1}^+}$ as follows. $\overline{\phi_{0}^+}: \overline{\P_0^+}\rightarrow \overline{\P_{1}^+}$ is the ultrapower map by the $(\cp(\overline{\phi_0}),\d^{\P_1})$-extender derived from $\overline{\phi_0}$. Note that $\overline{\phi_{0}^+}$ extends $\overline{\phi_{0}}$. Let $\psi_0^+: \overline{\P_0^+} \rightarrow \N_0^+$  be the ultrapower map by the $(\cp(\psi_0),\delta^{\N_0})$-extender derived from $\psi_0$. Again $\psi_0^+$ extends $\psi_0$. Finally let $\sigma_0^+ = (\overline{\phi^+_{0}})^{-1}\circ \psi_0^+$. The maps $\psi_i^+, \sigma_i^+, \overline{\phi_{i}^+}$ are defined similarly. Let $\overline{\P_\omega^+}$ be the direct limit of the linear system $(\overline{\P_i^+}, \overline{\phi^+_{i, k}}: i<k<\omega)$ where $\overline{\phi^+_{i, k}}$ is the composition of $(\overline{\phi^+_m}: m\in [i, k))$. Let $\overline{\phi_{i, \omega}^+}:\overline{\P_i^+}\rightarrow \overline{\P_\omega^+}$ and $\sigma^+_{i, \omega}:\N_i^+\rightarrow \overline{\P_\omega^+}$ be the direct limit embeddings.

Let now ${\sf{Hypo}}$ be the following statement:\\\\
$\sf{Hypo:}$ There is an $(\omega_1, \omega_1+1)$-iteration strategy for $\overline{\Pi_\omega}$ for $\overline{\P_\omega^+}$ such that the following clauses hold: \\  
$\sf{Hypo1:}$ $\overline{\Pi_\omega}$ acts on stacks that are above $\d^{\overline{\P_\omega}}$ where $\overline{\P_\omega}=\overline{\phi_{0, \omega}}(\overline{\P_0})$.\\
$\sf{Hypo2:}$  For every $i<\omega$\footnote{Including $i=0$.}, letting $\overline{\Pi_i}$ be the $\overline{\phi_{i, \omega}^+}$-pullback of $\overline{\Pi_\omega}$,  $\overline{\Pi_i}$ witnesses that $\overline{\P_i^+}$ is a $\overline{\Phi_i}^+$-hod mouse over $\overline{\P_i}$.\\\\
We have that $N[g]\models {\sf{Hypo}}$. Indeed, let for $i<\omega$, $m_i=\tau_{X_i}\circ (\pi\rest \overline{\P_i})$ and let $m_\omega:\overline{\P_\omega}\rightarrow j(\P)$ be the canonical embedding built via the direct limit construction. We thus have that for each $i$, $m_i=m_\omega\circ \overline{\phi_{i, \omega}}$. Just like $\overline{\phi_i}$, we can extend $m_i$ (for $i\leq \omega$) to $m_i^+:\overline{\P_i^+}\rightarrow j(\overline{\P_0^+})$. The desired strategy $\overline{\Pi_\omega}$ is $m_\omega^+$-pullback of $j(\overline{\Pi_0}^+\rest (H_{\kappa^+}^M))$. Because $m_i^+$ extends $\tau_{X_i}\circ (\pi\rest \overline{\P_i})$, we have that $\sf{Hypo}$ holds.

We now show how to finish the proof assuming $\sf{Hypo}$. By a similar argument as in \cite[Theorem 3.1.25]{trangThesis2013} or as in \cite[Page 663, just before (8) in the proof of Lemma 11.15]{CuBF}, we can use the strategies $\overline{\Pi_i}^+$'s to simultaneously execute a $\mathbb{R}^{V[G]}$-genericity iterations. The process yields a sequence of models $\langle\overline{\P^+_{\omega, i}},\N_{\omega, i}^+ \ | \ i\leq \omega\rangle$ and maps $\psi^+_{\omega, i}:\overline{\P^+_{\omega, i}}\rightarrow \N^+_{\omega, i}$, $\sigma^+_{\omega, i}:\N^+_{\omega, i}\rightarrow \overline{\P^+_{\omega, i+1}}$, and $\overline{\phi^+_{\omega, i}} = \psi^+_{\omega, i}\circ \sigma^+_{\omega, i}$\footnote{This embedding should not be confused with $\overline{\phi^+_{ i, \omega}}$.}. The iteration described above uses $\sigma_i^+$-pullback of $\overline{\Pi_i}$ to iterate $\N_i^+$. We denote this strategy by $\Sigma_i^+$.

Because the genericity iterations are above ${\sf{ord}}(\overline{\P_i})$ and ${\sf{ord}}(\N_i)$ for all $i\leq \omega$ and by \cite[Theorem 3.26]{ATHM}, the interpretation of the strategy of $\overline{\P_i}$ ($\N_i$ respectively) in the derived model of $\overline{\P_{\omega, i}^+}$ ($\N_{\omega, i}^+$, respectively) is $\overline{\Phi_i}^+$ ($\Sigma_i$, respectively). Let $C_i$ be the derived model of $\overline{\P^+_{\omega, i}}$ and $D_i$ be the derived model of $\N^+_{\omega, i}$ (at the sup of the Woodin cardinals of each model). Then $\mathbb{R}^{V[G]} = \mathbb{R}^{C_i} = \mathbb{R}^{D_i}$. Furthermore, $C_i\cap \powerset(\mathbb{R})\subseteq D_i\cap \powerset(\mathbb{R})\subseteq C_{i+1}\cap \powerset(\mathbb{R})$ for all $i$. 

Notice that we in fact have that $C_i=L(\Gamma(\overline{\P_i^+}, \overline{\Pi_i}))$ and $D_i=L(\Gamma(\N_i^+, \Sigma^+_i))$.  
Therefore, it follows from our choice of $\Pi_0$, $(F, \S)\in \cap_{i\leq \omega}(C_i\cap D_i)$, and since $\S$ is ordinal definable in each of $C_i$ and $D_i$, $(F, \S)\in \cap_{i\leq \omega}(\overline{\P_{\omega, i}^+}\cap \N_{\omega, i}^+)$. It follows that for each $i<\omega$, $\overline{T_i}$ and $U_i$ are definable respectively in $\overline{\P_{\omega, i}^+}$ and $\N_{\omega, i}^+$ from $(F, \S)$. Indeed, we have that\\\\
(F) for every $i<\omega$, $(\phi, s)\in \overline{T_i}$ if and only if $s\in [\d^{\overline{\P_i}}]^{<\omega}$ and in the derived model of $\overline{\P_{\omega, i}^+}$ at $\d^{\overline{\P_{\omega, i}^+}}(=\omega_1^{N[g]})$, $\S\models \phi[F, \pi^{\overline{\Phi_i}^+}_{\overline{\P_i}|\d^{\overline{\P_i}}, \infty}(s)]\}$\footnote{Here we abuse notation and use $\overline{\Phi_i}^+$ for the extension of $\overline{\Phi_i}$ to the derived model of $\overline{\P_{\omega, i}^+}$.},\\
(G) for every $i<\omega$, $(\phi, s)\in U_i$ if and only if $s\in [\d^{\N_i}]^{<\omega}$ and in the derived model of $\N_{\omega, i}^+$ at $\d^{\N_{\omega, i}^+}(=\omega_1^{N[g]})$, $\S\models \phi[F, \pi^{\Sigma_i}_{\N_i|\d^{\N_i}, \infty}(s)]\}$\footnote{Similar comments like above apply here as well.}.\\\\
Notice next that\\\\
(H) for every $i<\omega$, $\psi_{\omega, i}^+(\overline{\P_i})=\N_i$ and $\psi_{\omega, i}^+(\overline{\Phi_i})=\Sigma_i$\footnote{By this equation we mean that the internal strategy of $\overline{\P_i^+}$ is mapped to the internal strategy of $\N_i$.},\\
(I) there is $i_0<\omega$ such that for every $i\geq i_0$, $\phi_{\omega, i}^+(\S, F)=(\S, F)$ and $\psi_{\omega, i}^+(\S, F)=(\S, F)$\footnote{Because $\S$ and $F$ are ordinal definable in the respective derived models.}. \\\\
Thus, if $i_0$ is as in (I) then $\psi_{i_0}(\overline{T_{i_0}})=U_{i_0}$, contradiction! 
\end{proof}
It is now easy to derive \rthm{existence of condensing sets in n} from \rthm{thm:condensing_set} and \rthm{thm:weakly_condensing_set}.

Theorem \ref{existence of condensing sets in n} is typically applied in core model induction applications where there exists a mild large cardinal (e.g. a measurable cardinal) that gives rise to the embedding $j$ as in the hypothesis of the theorem. Below, we outline the proof of the following theorem, which gives the existence of condensing sets in some situations where large cardinals may not exist (e.g. under $\sf{PFA}$). One applies \rthm{existence of condensing sets in n} in applications where the core model induction is carried out in $V^{Coll(\omega,<\kappa)}$ and applies \rthm{thm:condensing_set_no_j} applications where the core model induction is carried out in $V^{Coll(\omega,\kappa)}$. 

In the following, we use the notations as in \ref{not: notation for this section} and in the previous section (in particular, $\P = \P_{(\phi,A)}, \Sigma=\Sigma_{(\phi,A)}$ etc. In the case $A\in V$, we define $\P^+ = \P^+_{\phi,A} = \sf{Lp}$$^{\Sigma,\Gamma,c}(\P^-)$ to be the union of sound $\Sigma$-premice $\M$ such that $\rho_\omega(\M)\leq \sf{ord}$$(\P^-)$ such that whenever $\pi: \M^*\rightarrow \M$ is elementary and $\M^*\in V$ is countable, then there is a unique iteration strategy $\Lambda\in \Gamma$ witnessing $\M^*$ is a $\Sigma^\pi$-mouse. We note that $\P,\P^+\in V$ and  $$\P\unlhd \P^+$$ though in general, equality may not hold.

\begin{theorem}\label{thm:condensing_set_no_j}
Suppose $\kappa$ is a  cardinal such that $\kappa^\omega = \kappa$, and $\kappa \geq 2^{2^{\aleph_0}}$. Let $g\subseteq Coll(\omega,\kappa)$ be $V$-generic. Suppose $\phi$ is a formula in the language of set theory and $A\in V$ such that $V[g] \models ``(\phi,A)$ is full, homogeneous and lower part $(\phi,A)$-covering fails". Furthermore, suppose  cof$(\sf{ord}$$(\P^+))\leq\kappa$. Then $V[g] \models$ ``there is a strongly $(\phi,A)$-condensing set".
\end{theorem}

\begin{remark}
The assumption ``cof$(\sf{ord}$$(\P^+))\leq\kappa$" in the above theorem holds in many situations, e.g. $\sf{PFA}$. If $\P = \P^+$ then this clause is superfluous as it is implied by the failure of lower part $(\phi,A)$-covering. $\myqedhere$
\end{remark}

The rest of the section is dedicated to outlining the proof of the theorem. We assume the hypothesis of the theorem from now to the end of this section. Let $\lambda >> \kappa$ and $X\prec (H_\lambda, \epsilon)$. We say that $X$ is \textit{good} if $|X|=\kappa$, $\kappa\subset X$, $X^\omega\subseteq X$, $\{\P,\P^+, \Sigma\rest V\} \in X$, and $X\cap \P^+$ is cofinal in $\P^+$. 

Let $X$ be good. Let $\pi_X: M_X\rightarrow (H_\lambda, \epsilon)$ be the uncollapse map. $\pi_X$ extends uniquely to a map $\pi_X^+: M_X[g]\rightarrow H_\lambda[g]$. We let $\gamma_X$ be the critical point of $\pi_X$ and $\pi^+_X(\P_X^-,\P_X,\P_X^+,\Sigma_X,\Gamma_X) = (\P^-,\P,\P^+,\Sigma,\Gamma)$; in general, if $a\in H_\lambda[g]$ is in the range of $\pi_X^+$, then we let $a_X = \pi_X^{+,-1}(a)$. We say that a good $X$ is \textit{$c\Gamma$-full} if $\P_X^+ = \sf{Lp}$$^{\Sigma^{\pi_X},\Gamma,c}(\P_X^-)$ and is \textit{$\Gamma$-full} if $\P_X = \sf{Lp}$$^{\Sigma^{\pi_X},\Gamma}(\P_X^-)$. It is clear that 
\begin{center}
$ \sf{Lp}$$^{\Sigma^{\pi_X},\Gamma}(\P_X^-)\unlhd \sf{Lp}$$^{\Sigma^{\pi_X},\Gamma,c}(\P_X^-)$.
\end{center}

\begin{lemma}\label{lem:full_X}
The set $S$ of $c\Gamma$-full $X$ is stationary. Furthermore, there is a stationary $T\subseteq S$ such that for each $X\in T$, $X$ is $\Gamma$-full.
\end{lemma}
\begin{proof}
We first show the first clause implies the second. Suppose the set $S$ of  $c\Gamma$-full $X$ is stationary and for contradiction, suppose that there is a club $C$ such that for all $X\in C\cap S$, $X$ is not $\Gamma$-full. Let $(X_\alpha, \M_\alpha: \alpha<\kappa^+)$ be such that 
\begin{itemize}
\item $(X_\alpha: \alpha<\kappa^+)$ is an increasing and continuous sequence in $C$ such that for all successor $\alpha$, $X_\alpha\in S$.
\item For each $\alpha$, $\M_\alpha\in \sf{Lp}$$^{\Sigma^{\pi_{X_\alpha}},\Gamma}(\P_{X_\alpha}^-)\backslash \P_{X_\alpha}$.
\end{itemize}
Letting $\P_{\alpha} = \P_{X_\alpha}, \P_{\alpha}^+ = \P^+_{X_\alpha}, \pi_{X_\alpha}=\pi_\alpha$ etc., we note that for any successor $\alpha$, $$\M_\alpha \lhd \P^+_{\alpha}.$$ This easily implies that letting $\pi_{\alpha,\beta} = \pi_{\beta}^{-1}\circ \pi_\alpha$ for $\alpha < \beta$, then for all successor ordinals $\alpha < \beta$, $$\pi_{\alpha,\beta}(\M_\alpha) = \M_\beta.$$
Now let $\M$ be the direct limit of $(\pi_{\alpha,\beta}, \M_\alpha: \alpha < \beta \wedge X_\alpha, X_\beta\in S)$, then $$ \P\unlhd \M \unlhd \P^+.$$ In particular, $\M$ is not an initial segment of $\P$. Now we show $\M\lhd \P$, which is a contradiction. Let $\pi:\M^*\rightarrow \M$ be elementary with $\M^*\in V[g]$ countable, transitive. Then there is $X_\alpha \in S$ and $\tau: \M^*\rightarrow \M_\alpha$ such that $\pi_\alpha\circ \tau = \pi$. This implies:
\begin{itemize}
\item $\Sigma^\pi = \Sigma_\alpha^\tau$, and
\item by the definition of $\M_\alpha$, $\M^*$ is a $\Sigma^\pi$-mouse with unique strategy in $\Gamma$.
\end{itemize}
This shows $\M\lhd \P$. Contradiction.

Now we prove the first clause; the idea of the proof is basically that of \cite[Theorem 3.4]{JSSS}. Suppose the set $W$ of good $X$ such that $X$ is not $c\Gamma$-full contains a $\kappa$-club. Let $\eta = cof^V(\sf{ord}$$(\P^+))$ and $(\M_i : i < \eta)$ be an enumeration of a cofinal sequence of sound $\M$ such that $\rho_1(\M)\leq \sf{ord}$$(\P^-)$ and $\P^- \lhd \M \lhd \P^+$. Let $(X_\alpha: \alpha < \kappa^+)$ enumerate an increasing, continuous sequence such that for successor ordinal $\alpha$ or limit $\alpha$ of cofinality $\geq\omega_1$, $X_\alpha\in W$. We use the notation as above, writing for example $\P_\alpha = \P_{X_\alpha}$. We also write $\Theta$ for $\sf{ord}$$(\P^-)$, $\theta_\alpha$ for $\pi_\alpha^{-1}(\Theta)$, and $\gamma_{\alpha}$ for $\gamma_{X_\alpha}$. We assume $(\M_i : i < \eta)\in X_0$ and let $(\M_i^\alpha : i<\eta) = \pi_\alpha^{-1}((\M_i : i<\eta)$. 

For each $\alpha$, let $\N_\alpha$ be the least sound $\N$ such that $\P_\alpha^+\lhd \N \lhd \sf{Lp}$$^{\Sigma^{\pi_X},\Gamma,c}(\P_X^-)$ and $\rho_\omega(\N_\alpha)\leq \sf{ord}$$(\P_\alpha^-)$. Let $n_\alpha$ be the least $n$ such that $\rho_{n+1}(\N_\alpha)\leq \sf{ord}$$(\P_\alpha^-)$. Let $\pi^*_{\alpha}: \N_\alpha \rightarrow \Q_\alpha$ be the corresponding $r\Sigma_{n_\alpha}$ ultrapower map given by the extender of length $\Theta$ derived from $\pi_\alpha$; similarly, we define $\pi^*_{\alpha,\beta}: \N_\alpha\rightarrow \Q^\beta_\alpha$ from $\pi_{\alpha,\beta}$. Note that the objects $\Q_\alpha, \Q^\beta_\alpha$ are all well-founded and hence we identify them with their transitive isomorph. By the assumption that $\eta < \kappa$ and $X_\alpha$ is good, the map $\pi_\alpha$ is cofinal in $\sf{ord}$$(\P^+)$ and therefore, $\neg (\Q_\alpha \lhd \P^+$). 

So there is an $n$ such that $$C = \{\alpha < \kappa^+ : \gamma_\alpha = \alpha \wedge n_\alpha = n\wedge cof(\alpha)\geq\omega_1\}.$$ contains a $\omega_1$-club. Fix an $\alpha \in C$ for now and let $(Y_\beta : \beta < \kappa^+)$ be a continuous, increasing sequence of $Y\prec H_{\lambda}$ such that $Y\cap \kappa^+\in \kappa^+$, $Y^\omega\subset Y$ and $(\N_\alpha,\Q_\alpha,\pi^*_\alpha)\in Y$.  For each $\beta$, let $\sigma_\beta: H_\beta \rightarrow Y_\beta$ be the uncollapse map and $\kappa_\beta = \textrm{crt}(\sigma_\beta)$. As in the proof of \cite[Theorem 3.4]{JSSS}, we get a club $$C_\alpha = \{\beta < \kappa^+ : \kappa_\beta = \beta \wedge \textrm{rng}(\pi_\beta)\cap \P^+  = Y_\beta \cap \P^+\}$$ and furthermore, using the agreement between $\sigma_\beta$ with $\pi_\beta$ on points in $\P^+$, we also get for $\beta\in C_\alpha$, $$\Q^\beta_\alpha = \sigma_\beta^{-1}(\Q_\alpha).\footnote{As in the proof of \cite[Theorem 3.4]{JSSS}, the map $\varphi: \Q^\beta_\alpha \rightarrow \sigma_\beta^{-1}(\Q_\alpha)$ defined as: $\varphi(\pi^*_{\alpha,\beta}(f)(a)) = \sigma_\beta^{-1}\circ \pi^*_{\alpha}(f)(a)$ for $a\in [\Theta_\beta]^{<\omega}$ and letting $\delta$ be such that $\pi_{\alpha,\beta}(\delta)>\textrm{max}(a)$, $f: [\delta]^{|a|}\rightarrow \N_\alpha$ come from the level $n$ Skolem term over $\N_\alpha$, is a well-defined elementary map and surjective. Therefore, it must be the identity.}$$

Fix $\beta \in \triangle_{\alpha<\kappa^+} C_\alpha \cap C$ such that $cof(\beta) \neq \eta$. This is possible because $\eta \leq \kappa$ and $\kappa \geq \omega_2$. For simplicity, let us assume $\rho_1(\N_\beta) = \Theta_\beta$\footnote{The general case where $n>1$ is the least such that $\rho_n(\N_\beta) = \Theta_\beta$ is handled just as in \cite[Theorem 3.4]{JSSS} by working with the $n$-reduct.}. So we have $\rho_0(\N_\beta) = \sf{ord}$$(\N_\beta)$ and 
\begin{center}
$\eta = cof(\sf{ord}$$(\P^+_\beta)) = cof(\sf{ord}$$(\N_\beta)) \leq \kappa.\footnote{The second equality follows from \cite[Lemma 1.2]{JSSS}.}$
\end{center}
Now we let $(\delta_i : i < \eta)$ be cofinal in $\sf{ord}$$(\N_\beta)$ and for each $i<\eta$, let $\overline{\sigma_i}: \overline{\N_i}\rightarrow Hull_1^{\N_\beta|\delta_i}(\Theta_\beta \cup p_1(\N_\beta))$ be the uncollapse map. By condensation, $\overline{\N_i}\lhd \P_\beta^+$ for each $i$.

Since $cof(\beta)\neq \eta$ and $\Theta_\beta = \bigcup_{\alpha<\beta} \pi_{\alpha,\beta}''\Theta_\alpha$, so $\P_\beta^+$ is the direct limit of the $\P_\alpha^+$ under the maps $\pi_{\alpha,\beta}$, there is an $\alpha < \beta$ and there are cofinal sets $T,T'\subset \eta$ such that 
\begin{center}
$i\in T \Rightarrow \M^\beta_i, p(\M^\beta_i)\in Hull_1^{\N_\beta}(\pi_{\alpha,\beta}''\Theta_\alpha\cup p_1(\N_\beta))$
\end{center}
and
\begin{center}
$i\in T' \Rightarrow \overline{\N_i}, \overline{\sigma_i}^{-1}(p_1(\N_\beta))\in \textrm{rng}(\pi_{\alpha, \beta})$.
\end{center}
Now we claim that for the $\alpha$ above,
\begin{equation}\label{eqn:capture_pi}
Hull_1^{\N_\beta}(\pi_{\alpha,\beta}''\Theta_\alpha\cup p_1(\N_\beta))\cap \textsf{ord}(\P_\beta^+) = \textrm{rng}(\pi_{\alpha,\beta})\cap \textsf{ord}(\P_\beta^+).
\end{equation}

Suppose $\xi\in \textrm{rng}(\pi_{\alpha,\beta})\cap \textsf{ord}(\P_\beta^+)$. Let $\pi_{\alpha,\beta}(\overline{\xi})=\xi$ for some $\overline{\xi} < \sf{ord}$$(\P_\alpha^+)$. There is some $i\in T$ such that $$\overline{\xi}\in Hull_1^{\M^\alpha_i}(\Theta_\alpha\cup p_1(\M_\alpha^1))$$. This implies $$\xi\in Hull_1^{\M^\beta_i}(\pi_{\alpha,\beta}''\Theta_\alpha  \cup p_1(\M_\beta^i))\subset Hull_1^{\N^\beta}(\pi_{\alpha,\beta}''\Theta_\alpha \cup p_1(\N_\beta)).$$

Conversely, let $\xi\in Hull_1^{\N_\beta}(\pi_{\alpha,\beta}''\Theta_\alpha\cup p_1(\N_\beta))\cap \textsf{ord}(\P_\beta^+) $. So there is $i\in T'$, a Skolem term $\tau$, and parameter $\vec{\epsilon}\in [\pi_{\alpha,\beta}''\Theta_\alpha]^{<\omega}$ such that $\xi = \tau^{\N_\beta|\delta_i}[\vec{\epsilon},p_1(\N_\beta)]$. We may also assume $\Theta_\beta \in Hull_1^{\N^\beta|\delta_i}(\pi_{\alpha,\beta}''\Theta_\alpha \cup p_1(\N_\beta))$; this is possible by the ``$\supseteq$" direction of \ref{eqn:capture_pi}, which we just proved. We easily get that $\xi\in   Hull_1^{\N_\beta|\delta_i}(\Theta_\beta \cup p_1(\N_\beta))$, hence $\xi \in Hull_1^{\overline{\N_i}}(\Theta_\beta\cup \overline{\sigma_i}^{-1}(p_1(\N_\beta)))$. So in fact, $\xi = \tau^{\overline{\N_i}}[\vec{\epsilon},\overline{\sigma_i}^{-1}(p_1(\N_\beta))]$. This implies $$\xi\in Hull_1^{\overline{\N_i}}(\pi_{\alpha,\beta}''\Theta_\alpha \cup \overline{\sigma_i}^{-1}(p_1(\N_\beta)))\subset \textrm{rng}(\pi_{\alpha,\beta})$$
as desired.

Now we finish the proof of the theorem. Let $\overline{\sigma}:\overline{\N}\rightarrow \N_\beta$ be the uncollapse map. By \ref{eqn:capture_pi}, $\overline{\sigma}$ and $\pi_{\alpha,\beta}$ agree on rng$(\pi_{\alpha,\beta})\cap \sf{ord}$$(\P_\beta^+)$. Therefore, $\overline{\N}|\sf{ord}$$(\P_\alpha^+) = \N_\alpha | \sf{ord}$$(\P_\alpha^+)$ and $\Q^\beta_\alpha|\sf{ord}$$(\P_\beta^+) = \N_\beta|\sf{ord}$$(\P_\beta^+)$. Now let $\pi: M\rightarrow H_\lambda$ be elementary with $M$ countable transitive and rng$(\pi)$ containing all relevant objects. We let $\pi(\overline{\M})=\overline{\N}$ and $\pi(\overline{\M}^*)=\N_\alpha$. Then note that $\overline{\M}$ is a $\Sigma_\beta^{\pi\circ \pi_{\alpha,\beta}}$-mouse and $\overline{\M}^*$ is a $\Sigma_\alpha^\pi$-mouse. But $\Sigma_\alpha = \Sigma_\beta^{\pi_{\alpha,\beta}}$ so in fact, $\overline{\M}$ is a $\Sigma_\alpha^\pi$-mouse. This easily implies $\overline{\M} = \overline{\M}^*$. By elementarity, $\overline{\N} = \N_\alpha$. Finally, using  $\overline{\N} = \N_\alpha$ and the agreement between $\overline{\sigma}$ and $\pi_{\alpha,\beta}$, we have
\begin{center}
$\Q^\beta_\alpha = \N_\beta$.
\end{center}
By pressing down, there is an $\alpha$ and a stationary set $Y$ of $\beta$ such that for all $\beta\in Y$, $\Q^\beta_\alpha=\N_\beta$ is the ultrapower of $\N_\alpha$ using the extender of length $\Theta_\beta$ derived from $\pi_{\alpha,\beta}$. Let $\N$ be the direct limit of such the $\N_\beta$ under these ultrapower maps. Then we easily get $\N = \Q_\alpha$ and hence $$\Q_\alpha \lhd \P^+.$$ Contradiction.

\end{proof}

\begin{remark}
The proof of the above lemma just requires a bit less of $\kappa$ than the hypothesis of Theorem \ref{thm:condensing_set_no_j}, namely we just need $\kappa \geq \omega_2$ and $\kappa^\omega=\kappa$.$\myqedhere$
\end{remark}


The following theorems are the corresponding versions of \rthm{thm:weakly_condensing_set} and \rthm{thm:condensing_set} and immediately imply Theorem \ref{thm:condensing_set_no_j}. Before proving \rthm{thm:weak_condensing_X}, we note that the  set of good $X$ contains an $\omega_1$-club and if $X$ is good, then $X$ contains $\powerset(\mathbb{R})$ because $\kappa \geq 2^{2^\omega}$. The following lemma will be used in the proof of both theorems.

\begin{lemma}\label{lem:key_realizing_lemma}
Suppose $X$ is good and $\Gamma$-full. Suppose $(\P_X^*,\Pi)$ is a $\Sigma_X$-hod pair in $\Gamma$ such that $(\P_X^*,\Pi\rest V)\in V$. Let $\P^* = \textrm{Ult}(\P_X^*,E)$ be the ultrapower of $\P_X^*$ via the extender of length $\Theta$ derived from $\pi_X$ and $\pi^*_X$ be the ultrapower map. Let $k: \R^*\rightarrow \P^*$ be elementary and $\R^*$ is countable, transitive in $V$. Then there is a $\Sigma_0$-elementary map $\overline{\pi}:\overline{\R^*}\rightarrow \P_X^+$ such that letting $\R = k^{-1}(\P)$, then $\Sigma^{k\rest \R} = \Sigma_X^{\overline{\pi}\rest \R}$.\footnote{D. Adolf has observed that this lemma holds and can be used to prove Lemma \ref{lem:full_X}. However, Lemma \ref{lem:key_realizing_lemma} uses essentially that $\kappa\geq 2^{2^\omega}$, while Lemma \ref{lem:full_X} holds with less required of $\kappa$.}
\end{lemma}
\begin{proof}
First we note that $E$ is a total extender over $\P_X^*$ because Lemma \ref{lem:full_X} implies that $\powerset^{\P^*_X}(\Theta_X)\subset M_X$. So the definition of $\P^*$ makes sense.

The proof of this lemma is essentially that of \cite[Lemma 8.12]{SZ} but with an additional detail. We will use the notations as introduced in \cite[Section 8]{SZ} regarding extenders. First let $W = \{(\P_\alpha,\Sigma_\alpha): \alpha < 2^{2^\omega}\}$ enumerate all countable hod pairs in $V$ such that $\Sigma_\alpha\in \Gamma$\footnote{We confuse $\Sigma_\alpha$ with its canonical extension in $V[g]$.}. Since $X$ is good, $W\subseteq X$; this is where we use $\kappa \geq 2^{2^\omega}$ in an essential way. Let $\alpha$ be such that $\Sigma_\alpha = \Sigma^{\pi^*\rest \R}$ and $\R = \Q_\alpha$.

Let $U = \textrm{rng}(\overline{\pi})$ and $((a_i, A_i) : i<\omega)$ enumerate all pairs $(c,A)$ such that there is a $\Sigma_0$-formula $\psi$ and $[a^1,f_1]_E, \dots, [a^k,f_k]_E\in U$ such that 
\begin{center}
$A = \{u \in \sf{ord}$$(\P_X)^{|c|} : P_X^* \models \psi[f_1^{a^1,c}(u), \dots, f_k^{a_k,c}(u)]\}\in E_c$.
\end{center}

Let $a\subset \omega$ be the set of $n$ such that $[a_n, f_n]_E$ represents some element of $\P$. Let $\{\tau_n : n<\omega\}$ enumerate all the Skolem functions of $\P^*_X$ and $b = \{i : \exists n\in a \  f_n = \tau_{n_i}\}$. So $\{\pi_X^*(f_n)(a_n) : n\in a\}$ is an elementary substructure of $\P$. In $H_\lambda$, the following first order statement with parameters $(\Q_\alpha, \Sigma_\alpha), (\P,\Sigma)$ holds: `` there is a sequence $(a_n : n<\omega)$ of finite sets of ordinals such that for each $n$, $a_n \in \pi_X(A_n)$ and $\Q_\alpha = \{\pi_X(\tau_{n_i})(a_n) : i\in b \wedge n\in a\}\prec_{\Sigma_0} \P$ and $\Sigma_\alpha$ is the pullback of $\Sigma$ under the uncollapse map". So by elementarity, the corresponding statement holds in $M_X$: `` there is a sequence $(a_n : n<\omega)$ of finite sets of ordinals such that for each $n$, $a_n \in A_n$ and $\Q_\alpha = \{\tau_{n_i}(a_n) : i\in b \wedge n\in a\}\prec_{\Sigma_0} \P_X$ and $\Sigma_\alpha$ is the pullback of $\Sigma_X$ under the uncollapse map". Let $(\bar{a_n} : n<\omega)$ witness the above statement.  

The embedding $\overline{\pi}$ is defined by: $\overline{\pi}([a_n,f_n]_E) = f_n(\bar{a_n})$ is the desired embedding with the property that 
\begin{center}
$\Sigma_X^{\overline{\pi}\rest \R} = \Sigma^{k\rest \R}$.
\end{center}
\end{proof}

\begin{theorem}\label{thm:weak_condensing_X}
There is a stationary set $S'\subset S$ such that whenever $X\in S'$, $X\cap \P$ is a weakly $(\phi,A)$-condensing set.
\end{theorem}

\begin{proof}
Suppose not. Fix a good $X$ such that $X$ is $\Gamma$-full but $X\cap \P$ is not a weakly condensing set. Note that $\pi_X\rest\P_X$ is cofinal in $\P$. Let $Y$ be an extension of $X\cap \P$ such that $(\Q_Y,\Sigma_Y)$ has the following properties:
\begin{enumerate}[(i)]
\item letting $k = \tau_Y$, $\Sigma_Y = \Sigma^k$;
\item $\Q_Y$ is not $\Gamma$-full, so letting $\R = \Q_Y$, there is a sound $\Sigma_Y$-mouse $\M$ such that $\neg (\M\unlhd \R)$ and $\rho_\omega(\M) = \delta^\R$.
\end{enumerate} 

By definition, $\tau_X = \tau_Y \circ \tau_{X,Y}$ ($\tau_X = \pi_X\rest \P_X$ here). Let $(\P^*_X,\Lambda_{X})\in V$ be a $\Sigma_X$-hod pair such that 
\begin{itemize}
\item $\Gamma(\P^*_X,\Lambda_X)\vDash \R$ is not full as witnessed by $\M$. \footnote{For brevity, we suppress mentioning the pair $(\S,\Phi)$ as in the proof of Theorem \ref{thm:weakly_condensing_set} and instead focus on the main points of the proof.}
\item $\Lambda_X\in \Gamma$ is $\Gamma$-fullness preserving and has strong branch condensation.
\item $\P^*_X$ is meek, is of limit type, and cof$^{\P^*_X}(\delta^{\P^*_X})=\omega$.
\end{itemize}
Such a pair $(\P^*_X,\Lambda_X)$ exists by boolean comparisons. In particular, $\P_X^*$ is a $\Sigma_X$-hod premouse over $\P_X$.

By arguments similar to before or that used in \cite[Lemma 3.78]{Trang2015PFA}, no $\M \lhd \P^*_X$ is such that  $\rho_\omega(\M) < \sf{ord}$$(\P^-_X)$ and in fact, $\sf{ord}$$(\P_X)$ is a cardinal of  $\P^*_X$.

By the above argument, $\P^*_X$ thinks $\P_X$ is full. Let 
\begin{center}
$\pi_X^*:\P_X^*\rightarrow \P^*$
\end{center}
be the ultrapower map by the extender $E$ of length $\Theta$ induced by $\pi_X$. Note that $\pi_X^*$ extends $\pi_X\rest \P_X$ (since $\pi_X$ is cofinal in $\P$) and $\P^*$ is wellfounded since $X$ is closed under $\omega$-sequences. Let 
\begin{center}
$i^*: \P^*_X\rightarrow \R^+$
\end{center}
be the ultrapower map by the extender of length $\delta^\R$ induced by $i =_{def} \tau_{X,Y}$. Note that $\R\lhd \R^+$ and $\R^+$ is wellfounded since there is a natural map 
\begin{center}
$k^*: \R^+ \rightarrow \P^*$
\end{center} 
extending $k$ such that $\tau_X^* = k^*\circ i^*$. Without loss of generality, we may assume $\M$'s unique strategy $\Sigma_\M \leq_w \Lambda_X$. Also, let $(\dot{\R},\dot{\M})$ be the canonical $Col(\omega,\kappa)$-names for $(\R,\M)$. Let $K$ be the transitive closure of $H^V_{\kappa}\cup (\dot{\R},\dot{\M})$.

Let $\W = \M_\omega^{\Lambda_X,\sharp}$ and $\Lambda$ be the unique strategy of $\W$. Let $\W^*$ be a $\Lambda$-iterate of $\W$ below its first Woodin cardinal that makes $K$-generically generic. Then in $\W^*[K]$, the derived model $D(\W^*[K])$ satisfies
\begin{center}
$L(\Gamma(\P^*_X,\Lambda_X),\mathbb{R}) \vDash \dot{\R} \textrm{ is not full as witnessed by } \M$.\footnote{This is because we can continue iterating $\W^*$ above the first Woodin cardinal to $\W^{**}$ such that letting $\lambda$ be the sup of the Woodin cardinals of $\W^{**}$, then there is a $Col(\omega,<\lambda)$-generic $h$ such that $\mathbb{R}^{V[G]}$ is the symmetric reals for $\W^{**}[h]$. And in $\W^{**}(\mathbb{R}^{V[G]})$, the derived model satisfies that $L(\Gamma(\P^+_X,\Lambda_X)) \vDash \R$ is not full.}
\end{center}
So the above fact is forced over $\W^*[K]$ for $\dot{\R}$.

Let $H\prec H_{\lambda}$ be countable (in $V$) such that all relevant objects are in $H$. Let $\pi: M\rightarrow H$ invert the transitive collapse and for all $a\in H$, let $\overline{a}=\pi^{-1}(a)$. By Lemma \ref{lem:key_realizing_lemma}, there is a map $\overline{\pi}:\overline{\R^+}\rightarrow \P_X^*$\footnote{We abuse notation a bit here. Technically, $\R$ is not in $V$. $\overline{\R}$ is the interpretation the name $\overline{\dot{\R}}$ over a $M[h]$ where $h\in V$ is $Coll(\omega,\overline{\kappa})$ generic over $M$. A similar comment applies to the maps $\overline{i^*}, \overline{k^*}$.} such that letting $\Lambda_0$ be the $\pi$-pullback of $\Lambda_X$ and $\Lambda_1$ be the $\overline{\pi}$-pullback of $\Lambda_X$, then
\begin{center}
$\Lambda_1\rest \overline{\R} = \Sigma^{\pi\rest \overline{\P}\circ \overline{k}}$, \footnote{This fact was missing from the proof of \cite[Lemma 3.80]{Trang2015PFA}. We need this to know that $(\overline{\R^+},\Lambda_1)$ is a $\Sigma^{\pi\rest \overline{\P}\circ \overline{k}}$-hod pair.}
\end{center}
and furthermore since $\pi\rest \overline{\P_X^*} = \overline{\pi}\circ \overline{i^*}$, \footnote{This follows from the definition of $\overline{\pi}$ and the fact that $\pi^*_X\circ \pi\rest \overline{\P^*_X} =  \pi\rest \bar{\P^+}\circ\overline{k^*}\circ \overline{i^*}$. }
\begin{center}
$\Lambda_0 = \Lambda_1^{\overline{i}}$.
\end{center}
In particular, $\Lambda_0\leq_w \Lambda_1$ and letting $\Sigma_{\overline{\R}} = \Lambda_1\rest \overline{\R} = \Sigma^{\pi\rest \overline{\P}\circ \overline{k}}$, $(\overline{\R^+},\Lambda_1)$ is a $\Sigma_{\overline{\R}}$-hod pair and that $\sf{ord}$$(\overline{\R})$ is a cardinal in $\overline{\R^+}$. 

We also confuse $\overline{\Lambda}$ with the $\pi$-pullback of $\Lambda$. Hence $\Gamma(\overline{\P_X^*},\Lambda_0)$ witnesses that $\overline{\R}$ is not full and this fact is forced over $\bar{\W^*}[\bar{K}]$ for the name $\bar{\dot{\R}}$. This means if we further iterate $\overline{\W^*}$ via $\overline{\Lambda}$ to $\mathcal{Y}$ such that $\mathbb{R}^{V[G]}$ can be realized as the symmetric reals over $\mathcal{Y}$ then in the derived model $D(\mathcal{Y})$, 
\begin{equation}\label{notFull}
L(\Gamma(\overline{\P^*_X},\Lambda_0)) \vDash \overline{\R} \textrm{ is not full}. 
\end{equation}
In the above, we have used the fact that the interpretation of the UB-code of the strategy for $\overline{\P_X^*}$ in $\mathcal{Y}$ to its derived model is $\Lambda_0\rest\mathbb{R}^{V[G]}$; this key fact is proved in \cite[Theorem 3.26]{ATHM} and Chapter 6.

Now we iterate $\overline{\R^*}$ to $\S$ via $\Lambda_1$ to realize $\mathbb{R}^{V[G]}$ as the symmetric reals for the collapse $Col(\omega,<\delta^\S)$, where $\delta^\S$ is the sup of $\S$'s Woodin cardinals. By  the fact that $\Lambda_0\leq_w \Lambda_1$ and $(\overline{\R^*},\Lambda_1)$ is a $\Sigma_{\overline{\R}}$-hod pair,  we get that in the derived model $D(\S)$,
\begin{center}
$\overline{\R}$ is not full as witnessed by $\bar{\M}$.
\end{center}
So $\Sigma_{\bar{\M}}$ is OD$_{\Sigma_{\overline{\R}}}$ in $D(\S)$ and hence $\overline{\M}\in \overline{\R^*}$.\footnote{We note that it is crucial here  that both $\overline{\M}$ and $\overline{\R^*}$ are $\Sigma_{\overline{\R}}$-mice.} This contradicts internal fullness of $\overline{\R}$ in $\overline{\R^*}$ since $\overline{\M}$ collapses $\sf{ord}$$(\overline{\R})$ in $\overline{\R^*}$ but $\sf{ord}$$(\overline{\R})$ is a cardinal in $\overline{\R^*}$. 

\end{proof}

For a good $X$, using the embedding $\pi_X$ we can define a $\pi_X$-realizable strategy $\Sigma^+_X$ for $\P_X$ using the construction of \rdef{the construction of the strategy}. We have that $\Sigma_X^+$ is such that 
\begin{itemize}
\item $\Sigma_X^+$ extends $\Sigma_X$;
\item for any $\Sigma_X^+$ iterate $\Q$ of $\P_X$ via stack $\vec{\T}$ such that the iteration embedding $\pi^{\vec{\T}}$ exists, there is an embedding $\sigma:\Q\rightarrow \P$ such that $\pi_X = \sigma\circ \pi^{\vec{\T}}$. Furthermore, letting $\Psi=(\Sigma_X^+)_{\VT, \Q}$, for all $\S\lhd^c_{hod}\Q$, $\Psi_{\S}$ has branch condensation.
\item $\Sigma_X^+$ is $\Gamma(\P_X,\Sigma_X^+)$-fullness preserving.
\end{itemize}
Theorem \ref{thm:weak_condensing_X} then implies that $\Sigma_X^+$ is $\Gamma$-fullness preserving.


\begin{theorem}\label{thm:condensing_X}
There is a stationary set $S'\subset S$ such that whenever $X\in S'$, $X\cap \P$ is a strongly $(\phi,A)$-condensing set.
\end{theorem}
\begin{proof}
The proof of this theorem is an adaptation of the proof of Theorem \ref{thm:weak_condensing_X} in a similar way one adapts the proof of \rthm{thm:weakly_condensing_set} to prove  \rthm{thm:condensing_set}. For completeness, we give a fairly detailed argument here. We will omit $(\phi,A)$ from our notations.

Suppose $X$ is a weakly condensing set and $B\in \P_X \cap \powerset(\Theta_X)$.\footnote{For the rest of this proof, whenever $X$ is weakly condensing, we automatically assume that $X=X'\cap \P$ for some good $X'$.} We say that $\tau_X$ has \textbf{$B$-condensation} if whenever $\Q=\Q_Y$ (where $Y$ is an extension of $X$) is such that there are elementary embeddings $\upsilon: \P_X \rightarrow \Q$, $\tau:\Q\rightarrow \P$ such that $\Q$ is countable in $V[g]$ and $\tau_X = \tau\circ \upsilon$, then $\upsilon(T_{\P_X,B}) = T_{\Q,\tau,B}$, where  
\begin{center}
$T_{\P_X,B} = \{(\psi,s) \ | \ s\in [\Theta_X]^{<\omega}\wedge \P_X \vDash \psi[s,B]\}$,
\end{center}
and
\begin{center}
$T_{\Q,\tau,B} = \{(\psi,s) \ | \ s\in [\delta_\alpha^\Q]^{<\omega} \textrm{ for some } \alpha<\lambda_\Q \wedge \P \vDash \phi[\pi^{\Sigma^{\tau}_{\Q}}_{\Q(\alpha),\infty}(s),\tau_X(B)]\}$,
\end{center}
where $\Sigma^{\tau}_\Q$ is the $\tau$-pullback strategy of $\Sigma$. We say $\tau_X$ has \textbf{condensation} if it has $B$-condensation for every $B\in \P_X \cap \powerset(\delta_X)$.

As before, we just prove the condensing part . To prove that a weakly condensing set $X$ is condensing, it is enough to prove that $\tau_X$ has condensation. Suppose for contradiction that the set $T$ of $X'\in S$ such that $X=X'\cap \P$ is cofinal in $\P$ and is not a condensing set is stationary. For each $X'\in T$, let $X=X'\cap \P$ (we will use this type of notations throughout this proof without mentioning again) and $A_X$ be the $\lesssim_X$-least such that $\tau_X$ fails to have $A_X$-condensation, where $\lesssim_X$ is the canonical well-ordering of $\P_X$. We say that a tuple $\{\langle \P_i,\Q_i,\tau_i,\xi_i,\pi_i,\sigma_i \ | \ i<\omega \rangle, \M_{\infty,Y}\}$ is a \textbf{bad tuple} if
\begin{enumerate}
\item $Y\in T$;
\item $\P_i = \P_{X_i}$ for all $i$, where $X'_i\in T$ and $\Q_i = \Q_{Y_i}$ for $Y_i$ an extension of $X_i$;
\item for all $i < j$, $X_i \prec Y_i \prec X_j \prec Y$;
\item $\M_{\infty,Y}$ be the direct limit of iterates $(\Q,\Lambda)$ of $(\P_Y,\Sigma^+_Y)$ such that $\Lambda$ has branch condensation; 
\item for all $i$, $\xi_i:\P_i\rightarrow \Q_i$, $\sigma_i:\Q_i \rightarrow \M_{\infty,Y}$,  $\tau_i: \P_{i+1}\rightarrow \M_{\infty,Y}$, and $\pi_i: \Q_i \rightarrow \P_{i+1}$;
\item for all $i$, $\tau_i = \sigma_i\circ \xi_i$, $\sigma_i = \tau_{i+1} \circ \pi_i$, and $\tau_{X_i,X_{i+1}}\rest \P_i=_{\textrm{def}} \phi_{i,i+1} = \pi_i \circ \xi_i$;
\item $\phi_{i,i+1}(A_{X_i}) = A_{X_{i+1}}$;
\item for all $i$, $\xi_i(T_{\P_i,A_{X_i}}) \neq T_{\Q_i,\sigma_i,A_{X_i}}$. 
\end{enumerate}
In (8), $T_{\Q_i,\sigma_i,A_{X_i}}$ is computed relative to $\M_{\infty,Y}$, that is
\begin{center}
$T_{\Q_i,\sigma_i,A_{X_i}} = \{(\phi,s) \ | \ s\in [\delta_\alpha^{\Q_i}]^{<\omega} \textrm{ for some } \alpha<\lambda^{\Q_i} \wedge \M_{\infty,Y} \vDash \phi[\pi^{\Sigma^{\sigma_i}_{\Q_i}}_{\Q_i(\alpha),\infty}(s),\tau_i(A_{X_i})]\}$
\end{center}
\begin{claim}\label{claim:bad_tuple} There is a bad tuple.
\end{claim}
\begin{proof}
For brevity, we first construct a bad tuple $\{\langle \P_i,\Q_i,\tau_i,\xi_i,\pi_i,\sigma_i \ | \ i<\omega \rangle, \P\}$ with $\P$ playing the role of $\M_{\infty,Y}$. We then simply choose a sufficiently large, good $Y$ and let $i_Y:\P_Y\rightarrow \M_{\infty,Y}$ be the direct limit map, $m_Y: \M_{\infty,Y}\rightarrow \P$ be the natural factor map, i.e. $m_Y\circ i_Y = \pi_Y$. It's easy to see that for all sufficiently large $Y$, the tuple $\{\langle \P_i,\Q_i,m_Y^{-1}\circ \tau_i,m_Y^{-1}\circ \xi_i,m_Y^{-1}\circ \pi_i,m_Y^{-1}\circ \sigma_i \ | \ i<\omega \rangle, \M_{\infty,Y}\}$ is a bad tuple.

The key point is (6). Let $A_X^* = \tau_X(A_X)$ for all $X\in T$. By Fodor's lemma, there is an $A^*$ such that $\exists^* X\in T \ A^*_X = A^*$.\footnote{``$\exists^* X\in T$" means ``stationarily many $X\in T$".} So there is an increasing and cofinal sequence $\{X_\alpha \ | \ \alpha < \kappa^+ \} \subseteq T$ such that for $\alpha < \beta$, $\tau_{X_\alpha,X_\beta}(A_{X_\alpha}) = A_{X_\beta} = \tau_{X_\beta}^{-1}(A)$. This easily implies the existence of such a tuple $\{\langle \P_i,\Q_i,\tau_i,\xi_i,\pi_i,\sigma_i \ | \ i<\omega \rangle, \P\}$.
\end{proof}
Fix a bad tuple $\mathcal{A} = \{\langle \P_i,\Q_i,\tau_i,\xi_i,\pi_i,\sigma_i \ | \ i<\omega \rangle, \M_{\infty,Y}\}$. Let $(\P_0^+,\Pi)$ be a (g-organized) $\Sigma_{\P_0}$-hod pair (cf. \cite{trang2013}) such that 
\begin{center}
$\Gamma(\P_0^+,\Pi) \vDash \mathcal{A}$ is a bad tuple.
\end{center}
We may also assume $(\P_0^+,\Pi\rest V)\in V$, $\delta^{\P_0^+}$ is limit of Woodin cardinals and is of nonmeasurable cofinality in $\P_0^+$ and there is some $\alpha<\lambda^{\P_0^+}$ such that $\Sigma_{Y} \leq_w \Pi_{\P_0^+(\alpha)}$. This type of reflection is possible because we replace $\P$ by $\M_{\infty,Y}$. Let $\W = \M_{\omega}^{\sharp, \Sigma_Y, \Pi, \oplus_{n<\omega}\Sigma_{X_n}}$ and $\Lambda$ be the unique strategy of $\W$. If $\mathcal{Z}$ is the result of iterating $\W$ via $\Lambda$ to make $\mathbb{R}^{V[G]}$ generic, then letting $h$ be $\mathcal{Z}$-generic for the Levy collapse of the sup of $\mathcal{Z}$'s Woodin cardinals to $\omega$ such that $\mathbb{R}^{V[G]}$ is the symmetric reals of $\mathcal{Z}[h]$, then in $\mathcal{Z}(\mathbb{R}^{V[G]})$,

\begin{center}
$\Gamma(\P_0^+,\Pi) \vDash \mathcal{A}$ is a bad tuple.
\end{center}

Now we define by induction $\xi_i^+: \P_i^+ \rightarrow \Q_i^+$, $\pi_i^+: \Q_i^+ \rightarrow \P_{i+1}^+$, $\phi_{i,i+1}^+: \P_i^+\rightarrow \P_{i+1}^+$ as follows. $\phi_{0,1}^+: \P_0^+\rightarrow \P_{1}^+$ is the ultrapower map by the extender derived from $\pi_{X_0,X_1}$ of length $\Theta_{X_1}$. Note that $\phi_{0,1}^+$ extends $\phi_{0,1}$. Let $\xi_0^+: \P_0^+ \rightarrow \Q_0^+$ extend $\xi_0$ be the ultrapower map by the extender derived from $\xi_0$ of length $\delta^{\Q_0}$. Finally let $\pi_0^+ = (\phi^+_{0,1})^{-1}\circ \xi_0^+$. The maps $\xi_i^+, \pi_i^+, \phi_{i,i+1}^+$ are defined similarly. Let also $\M_Y = \textrm{Ult}(\P_0^+,E)$, where $E$ is the extender derived from $\pi_{X,Y}$ of length $\Theta_Y$. There are maps $\epsilon_{2i}: \P_i^+ \rightarrow \M_Y$, $\epsilon_{2i+1}:\Q_i^+\rightarrow \M_Y$ for all $i$ such that $\epsilon_{2i} = \epsilon_{2i+1}\circ \xi^+_i$ and $\epsilon_{2i+1} = \epsilon_{2i+2}\circ \pi_i^+$. When $i = 0$, $\epsilon_0$ is simply $i_E$. Letting $\Sigma_i = \Sigma_{\P_i}$ and $\Psi=\Sigma_{\Q_i}$, $A_i = A_{X_i}$, there is a finite sequence of ordinals $t$ and a formula $\theta(u,v)$ such that in $\Gamma(\P_0^+,\Pi)$
\begin{enumerate}
\setcounter{enumi}{8}
\item for every $i<\omega$, $(\phi,s)\in T_{\P_i,A_i} \Leftrightarrow \theta[\pi^{\Sigma_i}_{\P_i(\alpha),\infty},t]$, where $\alpha$ is least such that $s\in [\delta_\alpha^{\P_i}]^{<\omega}$;
\item for every $i$, there is $(\phi_i,s_i)\in T_{\Q_i,\xi_i(A_i)}$ such that $\neg \theta[\pi^{\Psi_i}_{\Q_i(\alpha)}(s_i),t]$ where $\alpha$ is least such that $s_i\in [\delta_\alpha^{\Q_i}]^{<\omega}$.
\end{enumerate}
The pair $(\theta,t)$ essentially defines a Wadge-initial segment of $\Gamma(\P_0^+,\Pi)$ that can define the pair $(\M_{\infty,Y}, A^*)$, where $\tau_i(A_i)=A^*$ for some (any) $i$. 

Now let $X\prec H_{\lambda}$ be countable that contains all relevant objects and $\pi: M\rightarrow X$ invert the transitive collapse. For $a\in X$, let $\overline{a}=\pi^{-1}(a)$. By countable completeness of the extender $E$ and Lemma \ref{lem:key_realizing_lemma}, there is a map $\pi^*: \overline{\M_Y} \rightarrow \P^+_0$ with the property specified in Lemma \ref{lem:key_realizing_lemma}. Let $\overline{\Pi_i}$ be the $\pi^*\circ \overline{\epsilon_i}$-pullback of $\Pi$, so in $V[g]$,  
\begin{center}
$(\overline{\M_Y}, \Pi^{\pi^*})$ is a $\Sigma_{\M_{\infty,Y}}^\pi$-hod pair,\footnote{\cite[Lemma 3.82]{Trang2015PFA} concludes this by claiming $\pi\rest \overline{\M_Y} = \epsilon_0\circ \pi^*$, which is not true. One needs Lemma \ref{lem:key_realizing_lemma} to conclude this.}
\end{center}
\begin{center}
$\forall i < \omega$, $(\overline{\P_i^+},\overline{\Pi_i})$ is a $\Sigma^\pi_{\P_i}$-hod pair,
\end{center}
and $$\overline{\Sigma_Y} \leq_w \overline{\Pi_0} \leq_w \overline{\Pi_1} \dots \leq_w \Pi^{\pi^*}.$$ 

Let $\dot{\mathcal{A}}\in (H_{\bar{\kappa}})^M$ be the canonical name for $\bar{\mathcal{A}}$. It's easy to see (using the assumption on $\W$) that if $\W^*$ is a result of iterating $\bar{\W}$ via $\bar{\Lambda}$ (we confuse $\bar{\Lambda}$ with the $\pi$-pullback of $\Lambda$; they coincide on $M$) in $M$ below the first Woodin of $\bar{\W}$ to make $H$-generically generic, where $H$ is the transitive closure of $H_{\omega_2}^M\cup \dot{A}$, then in $\W^*[H]$, the derived model of $\W^*[H]$ at the sup of $\W^*$'s Woodin cardinals satisfies:
\begin{center}
$L(\bar{\P}_0,\mathbb{R}) \vDash \dot{\mathcal{A}}$ is a bad tuple.
\end{center} 

Now we stretch this fact out to $V[G]$ by iterating $\W^*$ to $\W^{**}$ to make $\mathbb{R}^{V[G]}$-generic. In $\W^{**}(\mathbb{R}^{V[G]})$, letting $i: \W^* \rightarrow \W^{**}$ be the iteration map then
\begin{center}
$\Gamma(\bar{\P_0}^+,\bar{\Pi}) \vDash i(\bar{\mathcal{A}})$\footnote{We abuse the notation slightly here. Technically, $\bar{\mathcal{A}}$ is not in $\W^*$ but $\W^*$ has a canonical name $\dot{{\mathcal{A}}}$ for $\bar{\mathcal{A}}$. Hence by $i(\bar{\mathcal{A}})$, we mean the interpretation of $i(\dot{{\mathcal{A}}})$.} is a bad tuple.
\end{center} 

By a similar argument as in \cite[Theorem 3.1.25]{trangThesis2013}, we can use the strategies $\overline{\Pi_i}^+$'s to simultanously execute a $\mathbb{R}^{V[G]}$-genericity iterations. The last branch of the iteration tree is wellfounded. The process yields a sequence of models $\langle\overline{\P^+_{i,\omega}},\overline{\Q_{i,\omega}^+} \ | \ i<\omega\rangle$ and maps $\overline{\xi^+_{i,\omega}}:\overline{\P^+_{i,\omega}}\rightarrow \overline{\Q^+_{i,\omega}}$, $\overline{\pi^+_{i,\omega}}:\overline{\Q^+_{i,\omega}}\rightarrow \overline{\P^+_{i+1,\omega}}$, and $\overline{\phi^+_{i,i+1,\omega}} = \overline{\pi^+_{i,\omega}}\circ \overline{\pi^+_{i,\omega}}$. Furthermore, each $\overline{\P^+_{i,\omega}}, \overline{\Q^+_{i,\omega}}$ embeds into a $\Pi^{\pi^*}$-iterate of $\overline{\M_Y}$ and hence the direct limit $\P_\infty$ of $(\overline{\P^+_{i,\omega}}, \overline{\Q^+_{j,\omega}} \ | \ i,j<\omega)$ under maps $\overline{\pi^+_{i,\omega}}$'s and $\overline{\xi^+_{i,\omega}}$'s is wellfounded. As mentioned above, $\overline{\P^+_{i,\omega}}$ is a (g-organized) $\Sigma^\pi_i$-premouse and $\overline{\Q^+_{i,\omega}}$ is a $^{g}\Psi^\pi_i$-premouse. Let $C_i$ be the derived model of $\overline{\P^+_{i,\omega}}$, $D_i$ be the derived model of $\overline{\Q^+_{i,\omega}}$ (at the sup of the Woodin cardinals of each model), then $\mathbb{R}^{V[G]} = \mathbb{R}^{C_i} = \mathbb{R}^{D_i}$. Furthermore, $C_i\cap \powerset(\mathbb{R})\subseteq D_i\cap \powerset(\mathbb{R})\subseteq C_{i+1}\cap \powerset(\mathbb{R})$ for all $i$.

(9), (10) and the construction above give us that there is a $t\in [\textrm{OR}]^{<\omega}$, a formula $\theta(u,v)$ such that
\begin{enumerate}
\setcounter{enumi}{10}
\item for each $i$, in $C_i$, for every $(\phi,s)$ such that $s\in \delta^{\overline{\P_i}}$, $(\phi,s)\in T_{\overline{\P_i},\overline{A_i}}\Leftrightarrow \theta[\pi^{\overline{\Sigma_i}}_{\overline{\P_i}(\alpha),\infty}(s),t]$ where $\alpha$ is least such that $s\in [\delta_\alpha^{\overline{\P_i}}]^{<\omega}$.
\end{enumerate}
Let $n$ be such that for all $i\geq n$, $\overline{\xi^+_{i,\omega}}(t) = t$. Such an $n$ exists because the direct limit $\P_\infty$ is wellfounded as we can arrange that $\P_\infty$ is embeddable into a $\Pi^{\pi^*}$-iterate of $\overline{\M_Y}$. By elementarity of $\overline{\xi^+_{i,\omega}}$ and the fact that $\overline{\xi^+_{i,\omega}}\rest \P_i = \overline{\xi_i}$,
\begin{enumerate}
\setcounter{enumi}{11}
\item for all $i\geq n$, in $D_i$, for every $(\phi,s)$ such that $s\in \delta^{\overline{\Q_i}}$, $(\phi,s)\in T_{\overline{\Q_i},\overline{\xi_i}(\overline{A_i})}\Leftrightarrow \theta[\pi^{\overline{\Psi_i}}_{\overline{\Q_i}(\alpha),\infty}(s),t]$ where $\alpha$ is least such that $s\in [\delta_\alpha^{\overline{\Q_i}}]^{<\omega}$.
\end{enumerate}
However, using (10), we get
\begin{enumerate}
\setcounter{enumi}{12}
\item for every $i$, in $D_i$, there is a formula $\phi_i$ and some $s_i\in [\delta^{\overline{\Q_i}}]^{<\omega}$ such that $(\phi_i,s_i)\in T^{\overline{\Q_i},\overline{\xi_i}(\overline{A_i})}$ but $\neg \phi[\pi^{\overline{\Psi_i}}_{\overline{\Q_i}(\alpha),\infty}(s_i),t]$ where $\alpha$ is least such that $s\in [\delta_\alpha^{\overline{\Q_i}}]^{<\omega}$.
\end{enumerate}
Clearly (12) and (13) give us a contradiction. This completes the proof of the lemma.

\end{proof}
\section{Condensing sets in models of $\sf{AD}^+$}\label{condensing sets under ad+ sec}

Thus far we have built condensing sets while working in models of $\sf{ZFC}$. In this section, we prove their existence in models of $\sf{AD}^+$. The material presented in this section will be used in the proof of generation of pointclasses (see \rthm{the generation of mouse full pointclasses}). Throughout this section we assume $\sf{AD}^++V=L(\powerset(\bR))$. Recall the notation $\Gamma_1\insegeq_{mouse}\Gamma_2$ (see \cite[Page 82]{ATHM} or \rsec{sec: derived models of hod mice}). 

Suppose $\Gamma$ is a mouse full pointclass (\rdef{mouse full}) such that:\\\\
\noindent (\textasteriskcentered)$_\Gamma$ \ there is a good pointclass $\Gamma^*$ containing $\Gamma$ and there is a sequence $(\Gamma_\a: \a<\Omega)$ with the property that
\begin{enumerate}
\item $\Omega$ is a limit ordinal,
\item $\Gamma_\a\inseg_{mouse}\Gamma$,
\item for $\a<\b<\Omega$, $\Gamma_{\a}\inseg_{mouse}\Gamma_\b$,
\item $\forall -1\leq \alpha < \Omega$, $\Gamma_{\alpha+1}$ is completely mousefull\footnote{Set $\Gamma_{-1}=\emptyset$.},
\item there is no completely mouse-full pointclass $\Psi\inseg_{mouse}\Gamma$ such that for some $\a$, $\Gamma_{\a}\inseg_{mouse}\Psi\inseg_{mouse}\Gamma_{\a+1}$,
\item if $\a<\Omega$ is a limit ordinal then $\Gamma_\a=\bigcup_{\b<\a}\Gamma_\b$,
\item $\Gamma=\bigcup_{\a<\Omega}\Gamma_\a$.
\end{enumerate}

Recall the definitions of ${\sf{HP}}^\Gamma$ and ${\sf{Mice}}^\Gamma$ (see \rnot{mice relative to gamma}). Let $\mathcal{F}=\{(\P, \Sigma)\in {\sf{HP}}^\Gamma: \Sigma$ is strongly $\Gamma$-fullness preserving and has strong branch condensation$\}$. We then let $\M^-=\bigcup_{(\P, \Sigma)\in \mathcal{F}}\M_\infty(\P, \Sigma)$. It follows from $\sf{AD}^+$ theory that if $(\P, \Sigma)\in \mathcal{F}$ then $\Sigma$ can be extended to a $(\Theta, \Theta, \Theta)$-iteration strategy\footnote{For example, see \cite{CoarseAD}.}. In what follows, we assume that if $(\P, \Sigma)\in \mathcal{F}$ then $\Sigma$ is a $(\Theta, \Theta, \Theta)$-iteration strategy.

Recall \rnot{complete layer notation}. Given $\R\inseg^c_{hod}\M^-$, we let $\Sigma_\R$ be the strategy of $\R$ such that whenever $(\P,\Lambda)\in \mathcal{F}$ is such that $\M_\infty(\P, \Sigma)=\R$ then $\Lambda_{\R}=\Sigma_\R$. Next we let ${\sf{Lp}}^{\Gamma, \oplus_{\R\inseg^c_{hod}\M^-}\Sigma_\R}(\M^-)$ be the stack of all sound $\oplus_{\R\inseg^c_{hod}\M^-}\Sigma_\R$-premice $\N$ over $\M^-$ such that $\rho(\N)\leq {\sf{ord}}(\M^-)$ and whenever $\pi:\S\rightarrow \N$ is elementary and $\S$ is countable then $\S$, as a $\oplus_{\R\inseg^c_{hod}\pi^{-1}(\M^-)}\Sigma^\pi_{\R}$-mouse, has an $\omega_1$-iteration strategy in $\Gamma$. Finally, if there is $\N\insegeq {\sf{Lp}}^{\Gamma, \oplus_{\R\inseg^c_{hod}\M^-}\Sigma_\R}(\M^-)$ such that $\rho(\N)<{\sf{ord}}(\M^-)$ then let $\M$ be the least such $\N$ and otherwise let $\M={\sf{Lp}}^{\Gamma, \oplus_{\R\inseg^c_{hod}\M^-}\Sigma_\R}(\M^-)$.

We let $\phi(u, v)$ be the formula that expresses the fact that $u$ is a mouse full pointclass such that (\textasteriskcentered)$_u$ holds and $v$ is a hod pair $(\Q, \Lambda)$ such that ${\sf{Code}}(\Lambda)\in u$ and $\Lambda$ has strong branch condensation and is strongly $u$-fullness preserving. 

\begin{remark} We have developed the concept of a hod mouse below $\sf{LSA}$. In the next theorem, hod pairs are all lsa small. However, the proof is general enough and uses this hypothesis only because we have not set up a general theory of hod mice. Because of this, we omit the extra hypothesis that we are in the minimal model of $\sf{LSA}$.$\myqedhere$
\end{remark}

\begin{theorem}\label{condensing set for gamma} Assume ${\sf{ZF}}+{\sf{AD}^+}$\footnote{Also, see the above remark.}. Suppose $\Gamma$ is a mouse full pointclass such that $(*)_\Gamma$ holds. Let 
\begin{itemize}
\item $\mathcal{F}=\mathcal{F}_{\phi, \Gamma}=\{(\P, \Sigma)\in {\sf{HP}}^\Gamma: \Sigma$ is strongly $\Gamma$-fullness preserving and has strong branch condensation$\}$,
\item $\M^-=\bigcup_{(\P, \Sigma)\in \mathcal{F}}\M_\infty(\P, \Sigma)$, and
\item let $\M$ be defined as follows: if there is $\N\insegeq {\sf{Lp}}^{\Gamma, \oplus_{\R\inseg^c_{hod}\M^-}\Sigma_\R}(\M^-)$ such that $\rho(\N)<{\sf{ord}}(\M^-)$ then let $\M$ be the least such $\N$ and otherwise let $\M={\sf{Lp}}^{\Gamma, \oplus_{\R\inseg^c_{hod}\M^-}\Sigma_\R}(\M^-)$.
\end{itemize}
 Then one of the following holds\footnote{What follows is not intended as an ``either or" conclusion.}. 
\begin{enumerate}
\item There is a hod pair or an anomalous hod pair $(\P, \Sigma)$ such that $\Sigma$ has strong branch condensation and is strongly $\Gamma$-fullness preserving, and $\Gamma(\P, \Sigma)=\Gamma$ (i.e., $(\phi, \Gamma)$ is not maximal).
\item $\M={\sf{Lp}}^{\Gamma, \oplus_{\R\inseg^c_{hod}\M^-}\Sigma_\R}(\M^-)$, lower part $(\phi, \Gamma)$-covering fails and there is a strongly $(\phi, \Gamma)$-condensing set $X\in \powerset_{\omega_1}(\M)$.
\item For some $(\Q, \Lambda)\in {\sf{HP}}^\Gamma$ such that $\Lambda$ has strong branch condensation and is strongly $\Gamma$-fullness preserving and for some $x\in \bR$, ${\sf{Lp}}^{\Lambda}(x)\not ={\sf{Lp}}^{\Gamma, \Lambda}(x)$.
\end{enumerate}
\end{theorem}
\begin{proof}
Towards a contradiction assume that all three clauses are false. We drop $(\phi, \Gamma)$ from our terminology. We will abuse our terminology and will say ``$\Gamma$-hod pair construction of $M$". Whenever we do this we mean the $\Gamma$-hod pair construction of $\mathbb{M}$ as defined in \rdef{gamma-hod pair construction*}. Here, $\mathbb{M}$ is a background whose universe is $M$, and it will always be clear exactly what $\mathbb{M}$ should be. 

 Let $A_0\subseteq \mathbb{R}$ be such that $A_0\in lub(\Gamma)$. Let $\Gamma_0, \Gamma_0^*, (N_0, \Phi_0), A_0^*,  \Gamma_1$ be such that
\begin{itemize}
\item $\Gamma_0$, $\Gamma^*_0$ and $\Gamma_1$ are good pointclasses,
\item $\Gamma\subseteq \utilde{\Delta}_{\Gamma_0}$,
\item $\Gamma_0\subseteq \utilde{\Delta}_{\Gamma^*_0}$,
\item $A_0^*\in lub(\Gamma_0^*)$,
\item $(N_0, \Phi_0)$ is a $\Gamma^*_0$-Woodin Suslin, co-Suslin capturing the sequence $(T_n(A_0): n\in \omega)$\footnote{See \rsubsec{subsec: capturing pointclasses}. There $T_n(X)$ is defined for $X$ a strategy but the same definition can be applied to any set of reals.}.
\end{itemize}
Let $F_0$ be as in  \rthm{n*x} for $(\Gamma_0, \Gamma^*_0, (N_0, \Phi_0), A_0^*)$, and fixing some $(N_1, \Phi_1), \Gamma_1^*, A_1^*$ let  $F_1$ be as in \rthm{n*x} for $(\Gamma_1, \Gamma^*_1, (N_1, \Phi_1), A_1^*)$.
\begin{itemize}
\item Let $x\in \dom(F_0)$ be such that if $F_0(x)=(\N', \M', \d', \Psi')$ then letting $\vec{G}$ be as in clause 7 of \rthm{n*x} and setting $\mathbb{M}_0=(\N, \d, \vec{G}, \Psi')$, $(\mathbb{M}_0, (N_0, \Phi_0), \Gamma_0, A_0^*)$ Suslin, co-Suslin captures $\Gamma$ and $A_0$. 
\item Let $y\in \dom(F_1)$ be such that if $F_1(y)=(\N^*_y, \M_y, \d_y, \Psi_y)$ then letting $\vec{G}_y$ be as in clause 7 of \rthm{n*x} and setting $\mathbb{M}_y=(\N^*_y, \d_y, \vec{G}_y, \Sigma_y)$,\begin{center} $(\mathbb{M}_y, (N_1, \Phi_1), \Gamma_1^*, A_1^*)$\end{center} Suslin, co-Suslin captures $\Gamma$ and $(N_1, \Phi_1)$ Suslin, co-Suslin captures ${\sf{Code}}(\Psi^*)$ where $\Psi^*$ is the $\omega_1$-strategy of $\M_2^{\#, \Phi_0}$. 
\end{itemize}
We record the following fact, which is a consequence of the proof of \rlem{correctness of backgrounds}\footnote{The lemma follows because letting $\d_1$ be the second Woodin cardinal of $\W'$, $\Psi$ allows us to define a $\d_1$-uB representation for $T_n(\Phi_0)$ (see \rlem{correctness of backgrounds}).}.
\begin{lemma}\label{consequence of correctness of backgrounds} Suppose $u$ is a set, $\W=\M_1^{\#, \Phi_0}(u)$ and $\Lambda$ is the unique strategy of $\W$ witnessing that $\W$ is a $\Psi$-mouse. Let $\d$ be the least Woodin cardinal of $\W$ and let $\W'$ be a $\Lambda$-iterate of $\W$ such that the iteration embedding $j:\W\rightarrow \W'$ exists. Let $h\subseteq Coll(\omega, j(\d))$ be $\W'$-generic. Then for any real $\tau\in \W'[h]$,
\begin{center}
$({\sf{HC}}^{\W'[h]}, A_0\cap \W'[h], \tau, \in)\prec ({\sf{HC}}, A_0, \tau, \in)$.
\end{center}
\end{lemma}

Let $\kappa$ be the least $<\d_y$-strong cardinal of $\N^*_y$. Let $g\subseteq Coll(\omega, <\k)$ be $\N^*_y$-generic. Let $\mathcal{F}_0\in \N^*_y[g]$ be the set of $(\Q, \Lambda)\in \N^*_y[g]$ such that 
\begin{itemize}
\item $\Q\in {\sf{HC}}^{\N^*_y[g]}$,
\item $\N^*_y[g]\models (\Q, \Lambda)\in {\sf{HP}}^\Gamma$,
\item $\N^*_y[g]\models ``\Lambda$ is $\Gamma$-fullness preserving and has strong branch condensation".
\end{itemize}
We use the methodology of \rsubsec{internalizing gamma sets} to obtain $(D, \psi)$ such that $\mathcal{F}_0=(\mathcal{F}_{\psi, D})^{\N^*_y[g]}$. Notice that $(\Q, \Lambda)\in \mathcal{F}_0$ if and only if there is a real $\sigma\in \N^*_y[g]$ such that $\sigma(0)$ is  a G\"odel number for some formula $\zeta$ and (in $\N^*_y[g]$) letting $A_0^{y, g}=A_0\cap \N^*_y[g]$, \\\\ 
(A) ${\sf{Code}}(\Lambda)$ is definable over $({\sf{HC}}, A_0^{y, g}, \sigma, \in)$ via $\zeta$ without parameters and\\
(B) $({\sf{HC}}, A_0^{y, g}, \sigma, {\sf{Code}}(\Lambda), \in)\models ``\Lambda$ is $\Gamma$-fullness preserving and has strong branch condensation",\\
(C) (A) and (B) hold in all further generic extensions of $\N^*_y[g]$. \\\\
We have that $(\psi, D)$ is lower part closed and stable. The next claim shows that it is directed. \\

%

\textit{Claim 1.} $\N^*_y[g]\models ``(\psi, D)$ is directed". \\\\
\begin{proof}
Fix $(\Q_0, \Lambda_0), (\Q_1, \Lambda_1)\in \mathcal{F}_0$. We now compare $(\Q_i, \Lambda_i)$ with the hod pair construction of $\N^*_y$. It follows from \rthm{comparison holds} that for each $i<2$, $\Q_i$ iterates, via $\Lambda_i$, to some model $\Q_i^+$ in the aforementioned hod pair construction such that $(\Lambda_i)_{\Q_i^+}$ is the strategy $\Q_i$ inherits from the background construction. Let $\nu_i<\k$ be such that $\Q_i\in\N^*_y[g\cap Coll(\omega, \nu_i)]$, and let $g_i=g\cap  Coll(\omega, \nu_i)$. To complete the proof it is enough to show that\footnote{This is because then by a Skolem hull argument we can obtain common iterates of $(\Q_0, \Lambda_0), (\Q_1, \Lambda_1)$ that are in ${\sf{HC}}^{\N^*_y[g]}$, and apply \rlem{correctness of backgrounds}.}.  \\\\
(a) for each $i$, $(\Q^+_i, (\Lambda_i)_{\Q_i^+})$ appears in the $\Gamma$-hod pair construction of $\N^*_y|\k[g_i]$ in which all extenders used have critical point $>\max(\nu_0, \nu_1)$. \\\\
 Let $\eta\in (\k, \d_y)$ be such that $(\Q^+_i, \Lambda_i)$ appears in the $\Gamma$-hod pair construction of $\N^*_y|\eta[g_i]$. Let then $E\in \vec{E}^{\N^*_y}$ be such that $\cp(E)=\k$ and $\nu_E>\nu_i$. It follows that in $Ult(\N^*_y, E)[g_i]$, $(\Q^+_i, \Lambda_i)$ appears in the $\Gamma$-hod pair construction of $(Ult(\N^*_y, E)|\pi_E(\k))[g_i]$. (a) now follows from elementarity. 
\end{proof}

Working in $\N^*_y$, let $\P^-=\P^-_{\psi, D}$. For $\a<\k$, let $g_\a=g\cap Coll(\omega, <\a)$. Our next claim implies that $(\psi, D)$ is of limit type.\\

\textit{Claim 2.} $\P^-$ is of a limit type. \\\\
\begin{proof}
Suppose not. It follows that there is $(\Q, \Lambda)\in \mathcal{F}_0$ such that $\P^-=\M_\infty(\Q, \Lambda)$. Let $\nu<\k$ be a cutpoint cardinal of $\N^*_y$ such that $\Q\in {\sf{HC}}^{\N^*_y[g_\nu]}$\footnote{It follows that $\Lambda\rest {\sf{HC}}^{\N^*_y[g_\nu]}\in \N^*_y[g_\nu]$ and $\Lambda=(\Lambda\rest {\sf{HC}}^{\N^*_y[g_\nu]}\in \N^*_y[g_\nu])^g$. We leave the details of such calculations to the reader. The methodology behind such calculations is presented in \rsubsec{internalizing gamma sets}.}. It follows from the proof of (a) in Claim 1 above that the $\Gamma$-hod pair construction of $\N^*_y|\k$ in which extenders used have critical point $>\nu$  reaches a pair $(\R, \Phi)$ such that $\R$ is a $\Lambda$-iterate of $\Q$ and $\Phi=\Lambda_\R$.

 Because of our condition on $\Gamma$ (namely that $\Omega$ is a limit ordinal) there is $\a+1<\Omega$ such that $\Gamma_\a=\Gamma(\Q, \Lambda)$. It follows that the $\Gamma$-hod pair construction of $\N^*_y$ using extenders with critical point $>\nu$ reaches $(\S, \Delta)\in \mathcal{F}$ such that $\Gamma(\S, \Delta)=\Gamma_{\a+1}$. It follows from the proof of (a) in Claim 1 above that the $\Gamma$-hod pair construction of $\N^*_y|\k$ in which extenders used have critical point $>\nu$ reaches such a pair $(\S, \Delta)$. It is then enough to show that $\N^*_y[g]\models (\S, \Delta)\in {\sf{HP}}^\Gamma$\footnote{Here we confuse $\Delta$ with its extension to $\N^*_y[g]$. Fullness preservation and branch condensation follow from \rthm{fullness preservation of background constructions} and \rthm{strong condensation for backgrounded strategies}. Recall that we are assuming that clause 3 of \rthm{condensing set for gamma} is false.}. Let $\nu_1\in (\nu, \k)$ be an $\N^*_y$-cutpoint cardinal such that $\S\in \N^*_y|\nu_1$.
 
 Let $\eta\in (\nu_1, \k)$ be the least $\N^*_y$-cardinal such that $\M_1^{\#, \Phi_0}(\N^*_y|\eta)\models ``\eta$ is a Woodin cardinal". Let $\N_1$ be the output of the fully backgrounded construction of $\N^*_y|\eta$ relative to $\Phi_1$ using extenders with critical points $>\nu_1$\footnote{See \rrem{general fully backgrounded constructions}}. We now compare $(\S, \Delta)$ with the $\Gamma$-hod pair construction of $\N_1$. Notice that all extenders of $\N_1$ have critical points $>\nu_1$. Let $\S_1$ be the output of the aforementioned $\Gamma$-hod pair construction. We claim that \\\\
 (b) some proper initial segment of $\S_1$ is a $\Delta$-iterate of $\S$.\\\\
  Suppose not. Let $z\in \dom(F_1)$ be such that $y<_T z$ and letting
  \begin{itemize}
  \item $F_1(z)=(\N^*_z, \M_z, \d_z, \Psi_z)$, 
  \item $\vec{G}_z$ be as in clause 7 of \rthm{n*x} and 
  \item $\mathbb{M}_z=(\M_z, \d_z, \vec{G}_z, \Psi_z)$,
  \end{itemize}
   then $(\mathbb{M}_z, (N_1, \Phi_1), \Gamma_1^*, A_1^*)$ Suslin, co-Suslin captures ${\sf{Code}}(\Delta)$ and $\N^*_y \in {\sf{HC}}^{\N^*_z}$.
 
Working in $\N^*_z$, let $\eta_1$ be the least $\N^*_z$-cardinal such that $\M_1^{\#, \Phi_0}(\N^*_z|\eta_1)\models ``\eta_1$ is a Woodin cardinal".  Let $\N^*$ be the output of the fully backgrounded construction of $\N^*_z|\eta_1$ relative to $\Phi_1$ done over $\N^*_y|\nu_1$. Comparing $\N^*_y$ with the construction producing $\N^*$ we get a normal stack $\T$ on $\N^*_y$ according to $\Psi_y$ such that $\T$ is based on $\N^*_y|\eta$ and if $\T^-$ is $\T$ without its last branch then $\m(\T^-)=\N^*|\eta_1$. 
 
 We now have that $\M_1^{\#, \Phi_0}(\N^*|\eta_1)\models ``\eta_1$ is a Woodin cardinal" (this can be shown by considering $S$-constructions). Yet, by elementarity $(\S, \Delta)$ wins the comparison against the $\Gamma$-hod pair construction of $\N^*|\eta_1$, contradicting universality of the latter. This contradiction implies that some initial segment of $\S_1$ is a $\Delta$-iterate of $\S$. Let $\S_2$ be this initial segment. This finishes the proof of (b).
 
 We now want to show that there is a real $q\in \mathbb{R}^{\N^*_y[g]}$ such that ${\sf{Code}}(\Delta_{\S_2})$ is definable over $({\sf{HC}}, A_0, q, \in)$ without parameters. Fix $r\in \bR$ such that ${\sf{Code}}(\Delta_{\S_2})$ is definable over $({\sf{HC}}, A_0, r, \in)$ without parameters, and let $\zeta$ be the formula defining ${\sf{Code}}(\Delta_{\S_2})$. Let $\xi$ be a cutpoint of $\N_1$ such that $\S_2\in \N_1|\xi$. Let $\N_1^+=\M_1^{\#, \Phi_0}(\N_1|\eta)$ and let $\Psi^+$ be the strategy of $\N_1^+$. Let $\pi:\N_1^+\rightarrow \N_2$ be an iteration of $\N_1^+$ via $\Psi^+$ such that $r$ is generic over $\N_1^+$ for the extender algebra at $\pi(\eta)$. We now have that\\\\
 (1) $\N_2[r]\models ``{\sf{Code}}(\Delta_{\S_2})$ is definable over $({\sf{HC}}, A_0\cap \N_2[r], r, \in)$\footnote{Here and below, we confuse $A_0$ with its interpretations in  relevant models.} via formula $\zeta"$\footnote{See \rlem{consequence of correctness of backgrounds}.}. \\\\
 It follows from elementarity of $\pi$ that\\\\
 (2) $\N_1^+\models ``$it is forced by $Coll(\omega, \eta)$ that there is a real $s$ such that ${\sf{Code}}(\Delta_{\S_2})$ is definable via $\zeta$ over 
 $({\sf{HC}}, A_0, s \in)$".\\\\
 Because $\N_1^+$ is countable in $\N^*_y[g]$, we can fix $q\in \bR^{\N^*_y[g]}$ such that \\\\
 (3) $q$ is in some $\leq \eta$-generic extension of $\N_1^+$  and 
 $\N_1^+[q]\models ``{\sf{Code}}(\Delta_{\S_2})$ is definable via $\zeta$ over $({\sf{HC}}, A_0, q, \in)$".\\\\
 Now $\delta$ is a Woodin cardinal in $\N_1^+[q]$, and so using genericity iterations we can show that 
 ${\sf{Code}}(\Delta_{\S_2})$ is definable over $({\sf{HC}}, A_0, q, \in)$ via $\zeta$. This finishes the proof of Claim 2.\footnote{The proof is a bit more involved. Notice that $\N_1^+[q]$ captures Suslin, co-Suslin captures ${\sf{Code}}(\Delta_{\S_2})$. This is because for some $\N_1^+$-successor cardinal $\nu'\in [\nu_1, \eta)$, ${\sf{Code}}(\Delta_{\S_2})$ is determined by the fragment of $\Psi^+_{\N_1^+|\nu'}$ that acts on iteration that are above $\nu_1$. It follows that  $\N_1^+[q]$ has a way of determining ${\sf{Code}}(\Delta_{\S_2})$ in its generic extensions.  \rlem{consequence of correctness of backgrounds} then gives what we want.}
 \end{proof}
 
 Our discussion before Claim 1, Claim 1 and Claim 2 show that $(\psi, D)$ is  lower part closed, is of limit type, is stable and is directed. We now work in $\N^*_y[g]$. 
 \begin{notation}\label{notation for this lemma: sigma} Let
 \begin{enumerate}
 \item $\Sigma=\Sigma_{\psi, D}$ (see clause 2 of \rnot{direct limit}) and
 \item $\P=\P_{\psi, D}$. 
 \end{enumerate}
 Notice that if $h$ is $Coll(\omega, \bR^{\N^*_y[g]})$-generic over $\N^*_y[g]$ then there is a real $z\in \N^*_y[g]$ such that $z(0)$ is a G\"odel number for a formula $\zeta$ such that $\Sigma$ is definable over $({\sf{HC}}^{\N^*_y[g*h]}, A_0\cap {\sf{HC}}^{\N^*_y[g*h]}, z, \in)$ via $\zeta$ without parameters. Notice that if $\Sigma^+$ is the strategy for $\P|\d^\P$ definable over $({\sf{HC}}, A_0, z, \in)$ via $\zeta$ without parameters then $\Sigma^+\rest (\N^*_y|\d_y)[g]=\Sigma$. We will confuse $\Sigma^+$ with $\Sigma$.$\myqedhere$
 \end{notation}
 
 \textit{Claim 3.} ${\sf{Code}}(\Sigma)\in \Gamma$.\\\\
 \begin{proof}
  Towards a contradiction assume ${\sf{Code}}(\Sigma)\not \in \Gamma$. It then follows that $\Gamma(\P|\d^\P, \Sigma)=\Gamma$, and hence clause 1 of \rthm{condensing set for gamma} holds.
  \end{proof}
  
Since $\Gamma_1^* \not =\powerset(\bR)$, there is a $C\subseteq \bR$ such that $\Gamma_1^*, F_1\in L(C, \bR)$. We then have that $L(C, \bR)\models \sf{DC}$. Work in $W=L(C, \bR)$ and let $G\subseteq Coll(\omega_1, \bR)$ be  $W$-generic. Notice that $W[G]\models \sf{ZFC}$. Recall $\mathcal{F}$ from the statement of \rthm{condensing set for gamma}. Let $((\Q_\a, \Lambda_\a):\a<\omega_1)\in W[G]$ be an enumeration of $\mathcal{F}$ and $(z_\a:\a<\omega_1)\in W[G]$ be an enumeration of $\bR$. In $W[G]$, choose a sequence $(y_\a: \a<\omega_1)$ of reals such that 
  \begin{enumerate}
  \item $y_0=y$ and $g\in \N^*_{y_1}$,
  \item for all $\a<\omega_1$, letting $F_1(y_\a)=(\N^*_{y_\a}, \M_{y_\a}, \d_{y_\a}, \Psi_{y_\a})$, $(z_\b: \b\leq \a)\in \N^*_{y_\a}$ and $\oplus_{\b\leq \a}\Lambda_\a$ is Suslin, co-Suslin captured by $(\N^*_{y_\a}, \d_{y_\a}, \Psi_{y_\a})$, and
  \item for $\b<\omega_1$, $(\N^*_{y_\gamma}: \gamma<\b)\in {\sf{HC}}^{\N^*_{y_\b}}$.
  \end{enumerate}
  
 We now construct a sequence of $\Phi_1$-mice $(\M_\a, \N_\a: \N_\a\inseg \M_\a \wedge \a<\omega_1)$ and a sequence of commuting embeddings $\pi_{\a, \b}:\M_\a\rightarrow \M_\b$ such that $\pi_{\a, \b}(\N_\a)=\N_\b$ and if $\k_\a=\cp(\pi_{\a, \b})$ then $\N_\a=\M_\a|\k_\a$. For $\a>0$ we will have that $\M_\a$ is the output of a fully backgrounded  construction of $\N^*_{y_\a}$ relative to $\Phi_1$ and also that $\N_\a\insegeq \M_\a$, and $\M_\a$ will be a $\Psi_y$-iterate of $\N^*_{y}$. Below we describe the construction.
 \begin{itemize}
 \item Set $\M_0=\N^*_{y_0}$ and $\N_0=\M_0|\k$.
 \item For $\a<\omega_1$, let $\vec{G}_\a$ consist of those $E\in \vec{E}^{\N^*_{y_\a}}$ such that $\cp(E)>{\sf{ord}}(\N_\a)$ and $\nu(E)$ is an inaccessible cardinal of $\N^*_{y_\a}$. 
 \item Given $\M_\a$ and $\N_\a$, let $\M_{\a+1}=({\sf{Le}}((N_1, \Phi_1), \N_\a))^{(\N^*_{y_{\alpha}}, \d_\a, \vec{G}_\a)}$. 
 \item Let $\pi_{\a, \a+1}:\M_{\a}\rightarrow \M_{\a+1}$ be the iteration embedding according to $(\Psi_y)_{\M_\a}$\footnote{Notice that because $\N^*_{y_\a}\in {\sf{HC}}^{\N^*_{y_{\a+1}}}$, $\M_{\a+1}$ is a $(\Psi_y)_{\M_\a}$-iterate of $\M_\a$.}. 
 \item Let $\k_{\a+1}$ be the least $\d_{y_{\a+1}}$-strong cardinal of $\M_{\a+1}$ and let 
 \begin{center}
 $\N_{\a+1}=\M_{\a+1}|\k_{\a+1}$.
 \end{center}
 It follows that $\N_{\a+1}=\pi_{\a, \a+1}(\N_\a)$\footnote{Notice that if $E\in \vec{E}^{\M_{\a}}$ is the extender with the least index on the extender sequence of $\M_\a$ such that $\cp(E)=\k_\a$ then $E$ is the first extender used in the $\M_\a$-to-$\M_{\a+1}$ iteration. Here we assume that all extenders with $\cp(E)$ are total. Otherwise we can translate them away as is done in \cite[Remark 12.7]{DMATM}.}. 
 \item Suppose now that $\l<\omega_1$ is a limit ordinal and we have constructed a sequence $(\M_\a, \N_\a: \N_\a\inseg \M_\a \wedge \a<\l)$ and a sequence of commuting embedding $\pi_{\a, \b}:\M_\a\rightarrow \M_\b$ for $\a<\b<\l$. Let $\M^*_\l$ be the direct limit of $\M_\a$ under $\pi_{\a, \b}$. Let $\pi^*_{\a, \l}:\M_\a\rightarrow \M_\l^*$ be the embedding given by the direct limit construction. Let then $\N_{\l}=\pi_{0, \l}(\N_0)$ and let $\M_\l=({\sf{Le}}((N_1, \Phi_1), \N_\l))^{(\N^*_{y_{\l}}, \d_\l, \vec{G}_\l)}$. Letting $k:\M_\l^*\rightarrow \M_\l$ be the iteration embedding according to $(\Psi_{y})_{\M_\l^*}$, we set $\pi_{\a, \l}=k\circ \pi_{\a, \l}^*$.
 \end{itemize}
 

Finally let $\M_{\omega_1}$ be the direct limit of the system $(\M_\a, \pi_{\a, \b}: \a<\b<\omega_1)$ and let $\pi_{\a, \omega_1}:\M_\a\rightarrow \M_{\omega_1}$ be the direct limit embedding. Let $\P_{\omega_1}=\pi_{0, \omega_1}(\P)$ and $\P^-_{\omega_1}=\pi_{0, \omega_1}(\P^-)$. \\

\textit{Claim 4.} Fix $\a<\omega_1$ and let $h\subseteq Coll(\omega, <\k_\a)$ be $\N_{y_\a}^*$-generic. Then $\pi_{0, \a}(\P)=(\P_{\psi, D})^{\N^*_{y_\a}[h]}$ and $\pi_{0, \a}(\Sigma)=(\Sigma_{\psi, D})^{\N^*_{y_\a}[h]}$.\\\\
We leave the proof of Claim 4 to the reader as it is very similar to the proofs of Claim 1 and Claim 2. For $\a\leq \omega_1$, we let $\P_\a=\pi_{0, \a}(\P_\a)$, $\P^-_\a=\pi_{0, \a}(\P^-)$ and $\Sigma^\a=\pi_{0, \a}(\Sigma)$.\\

\textit{Claim 5.} $\P_{\omega_1}=\M$.\\\\
\begin{proof}
 Notice that\\\\
 (1) for $\a<\b<\omega_1$ and for $\R\inseg^c_{hod}\P_\a$, $\pi_{\a, \b}\rest \R$ is the iteration embedding according to $(\Sigma^\a)_{\R}$, and\\
 (2) if $\a<\omega_1$, $\R\inseg^c_{hod}\P_\a$ and $\Q$ is a $(\Sigma^\a)_{\R}$-iterate of $\R$ then there is $\b<\omega_1$ such that some $\pi_{\a, \b}(\R)$ is a $(\Sigma^\a)_\Q$-iterate of $\Q$.\\
 (3) for all $\a<\omega_1$ there is $\R\inseg^c_{hod}\P_\a$ such that $\R$ is a $\Lambda_\a$-iterate of $\Q_\a$.\\\\
  To see (2), let $\b$ be large enough such that $(\Q_\b, \Lambda_\b)=(\Q, (\Sigma^\a)_\Q)$. It then follows that $\pi_{\a, \b}(\R)$ is a $(\Sigma^\a)_\Q$-iterate of $\Q$. It follows from (1) and (2) that $\P_{\omega_1}|\d^\M=\M|\d^\M$.
  
  If $\rho(\P_{\omega_1})<{\sf{ord}}(\P^-_{\omega_1})$ then we must have that $\P_{\omega_1}=\M$. Suppose then $\rho(\P_{\omega_1})>{\sf{ord}}(\P^-_{\omega_1})$. Clearly $\P_{\omega_1}\insegeq {\sf{Lp}}^{\Gamma, \Sigma}(\M^-)$. Suppose then $\P_{{\omega_1}}\inseg {\sf{Lp}}^{\Gamma, \Sigma}(\M^-)$. By a standard Skolem hull argument, it follows that for some $\a<\omega_1$, $\P_\a\inseg {\sf{Lp}}^{\Gamma, \Sigma^\a}(\pi_{0, \a}(\P^-))$. However, because $\rho(\P_{{\omega_1}})>{\sf{ord}}(\P_{\omega_1}^-)$, $\N^*_{y_\a}\models ``\P_\a= {\sf{Lp}}^{\Gamma, \Sigma^\a}(\pi_{0, \a}(\P^-))"$, contradiction.
\end{proof}

\textit{Claim 6.} $\rho(\M)>{\sf{ord}}(\M^-)$.\\\\
\begin{proof} Assume $\rho(\M)<{\sf{ord}}(\M^-)$ (it follows from the definition of $\M$ that equality is impossible). We now have that $\rho(\P)<\d^\P$. The argument now takes place in $\N^*_y[g]$. Let $N=\N^*_y$ and let $U\in N$ be the Mitchell order 0 ultrafilter on $\k$. Let $j: N\rightarrow Ult(N, U)$ be the ultrapower embedding and $j^+: N[g]\rightarrow Ult(N, U)[g']$ be its lift up to $N[g]$. Notice that $j^+(\Gamma)$ makes sense. As in core model induction applications $\Sigma$ can be extended to a strategy $\Sigma'$ for $\P$\footnote{For example, see \rdef{the construction of the strategy}, \rsec{sec: e-certified strategies}, \cite[Definition 6.14]{ATHM} and also \cite{CuBF} and \cite{Trang2015PFA}. Recall that $\Sigma$ is a strategy for $\P^-$ which in this case is just $\P|\d^\P$.}. It follows from clause 2 of \rthm{branch condensation on a tail} that there is a tail $(\Q, \Lambda)\in Ult(N, U)[g]$ of $(\P, \Sigma')$ such that $\Lambda$ has strong branch condensation. Because we are assuming that (in $Ult(N, U)[g]$) clause 1 of \rthm{condensing set for gamma} fails and because $\Sigma'$ is a $j$-realizable strategy, $\Gamma(\Q, \Lambda)\subset j^+(\Gamma)$ and ${\sf{Code}}(\Lambda)\in \Gamma$. Notice next that $\rho(\Q)<\d^\Q$. We can then finish by using the argument given on page 143 of \cite{ATHM}\footnote{This is a standard argument in core model induction. The reader can also consult \cite{CuBF} and  \cite{Trang2015PFA}.}. 
\end{proof}
  
 We thus have that $\P={\sf{Lp}}^{\Gamma, \Sigma}(\P^-)$. \\
  
  \textit{Claim 7.} $\N^*_y\models \card{\P}=\k$.\\\\
  \begin{proof}
  Recall the real $z$ introduced before the statement of Claim 3. We have that $z\in \N^*_y[g][h]$  where $h$ is $Coll(\omega, \bR^{\N^*_y[g]})$-generic. It then follows that $\P$ is definable over $({\sf{HC}}^{\N^*_y[g][h]}, A_0\cap{\sf{HC}}^{\N^*_y[g][h]}, z, \in)$ and hence, $\P\in {\sf{HC}}^{\N^*_y[g][h]}$. Thus $\card{\P}^{\N^*_y[g]}=\kappa$.
  \end{proof}
  
 Notice that Claim 7 implies that lower part $(\phi, \Gamma)$-covering fails as it implies that $\cf({\sf{ord}}(\P_{\omega_1}))=\omega$\footnote{This is because $\pi_{0, \omega_1}$ is continuous at ${\sf{ord}}(\P)$.}. 
 It follows from \rthm{thm:condensing_set} that $X=_{def}\pi_{0, 1}[\P]\in \powerset(\P_1)\cap \M_1$ is such that\\\\
(A) for any $\M_1$-generic $h\subseteq Coll(\omega, <\k_1)$, $\M_1[h]\models ``X$ is countable and is a $(\psi, D)$-condensing set". \\\\
 It follows from Claim 7 that\\\\
  (B) for every $\a\in [1, \omega_1)$ and for every $\M_{y_\a}$-generic $h\subseteq Coll(\omega, <\k_\a)$, $\M_{y_\a}[h]\models ``\pi_{1, \a}[X]$ is a $(\psi, D)$-condensing set". \\

   \textit{Claim 8.} For every $\a\in [1, \omega_1)$ and for every $\N^*_{y_\a}$-generic $h\subseteq Coll(\omega, <\k_\a)$, $\N^*_{y_\a}[h]\models ``\pi_{1, \a}[X]$ is a weakly $(\psi, D)$-condensing set". \\\\
   \begin{proof}
   We give the proof for $\a=1$ and leave the rest to the reader. Let $h\subseteq Coll(\omega, <\k_1)$ be $\N^*_{y_1}$-generic and let $Y\in (\powerset_{\omega_1}(\P_1))^{\N^*_{y_1}[h]}$ be an extension of $X$. In what follows we will use the notation introduced in  \rsec{condensing sets sec} relative to $\N^*_{y_1}[h]$. Thus, $\Sigma_Y\in \N^*_{y_1}[h]$ is the $\tau_Y:\P_Y\rightarrow \P_1$-pullback of $\pi_{0, 1}(\Sigma)$. However, we will also confuse $\Sigma_Y$ and $\pi_{0, 1}(\Sigma)$ with their canonical extensions that act on all stacks. 

The proof of the claim follows the steps of \rthm{thm:weakly_condensing_set}. Recall that in that proof the key step is to find a universal model extending $\P$ such that $\pi_{0, 1}$ acts on it. Here, we describe how to find this universal model and leave the rest, which is just like the proof of \rthm{thm:weakly_condensing_set}, to the reader.  To simplify, we only show that if $\Q_Y^-=\tau_Y^{-1}(\P^-_1)$ then $\Q_Y={\sf{Lp}}^{\Gamma, \Sigma_Y}(\Q_Y^-)$. The rest of the proof is very similar.
  
  Suppose then that $\Q_Y\inseg {\sf{Lp}}^{\Gamma, \Sigma_Y}(\Q_Y^-)$ and let $\S\inseg {\sf{Lp}}^{\Gamma, \Sigma_Y}(\Q_Y^-)$ be the least such that $\rho(\S)\leq {\sf{ord}}(\Q_Y^-)$ and $\S\not \insegeq \Q_Y$. Let $(\R, \Lambda)\in {\sf{HP}}^\Gamma$ be such that 
  \begin{itemize}
  \item $\R$ is meek and of limit type,
  \item $(\R, \Lambda)$ be a $\Sigma$-hod pair,
  \item $L(\Gamma(\R, \Lambda), \bR)\models ``\S$, as a $\Sigma_Y$-mouse, has an $\omega_1$-iteration strategy."
  \end{itemize}
  Let $\a<\omega_1$ be such that ${\sf{Code}}(\Lambda)$ is Suslin, co-Suslin captured by $(\N^*_{y_\a}, \d_{y_\a}, \Psi_{y_\a})$. Recalling \rdef{fully backgrounded sts construction} and \rrem{general fully backgrounded constructions}, let
  \begin{itemize}
  \item $\W^*$ be the output of $({\sf{Le}}((N_1, \Phi_1)\oplus (\P, \Sigma), \mathcal{J}_{\omega}(N_1, \P)))^{(\N_y^*, \d_y, \vec{G}_y)}_{>\k}$ and 
  \item $\W^{**}$ be the output of $({\sf{Le}}((N_1, \Phi_1)\oplus (\P, \Sigma), \mathcal{J}_{\omega}(N_1, \P)))^{(\N^*_{y_\a}, \d_{y_\a}, \vec{G}_{y_\a})}_{>\k}$.
  \end{itemize}
   Notice that it follows that ${\sf{ord}}(\W^*)=\d_y$ and ${\sf{ord}}(\W^{**})=\d_{y_\a}$. We now compare the construction producing $\W^*$ and the construction producing $\W^{**}$. The  comparison produces a tree $\T$ on $\N^*_y$ of limit length such that\\\\
   (T1) $\T\in \N^*_{y_\a}$,\\
   (T2) setting $b=\Psi_y(\T)$, $\pi^\T_b(\W^*)=\W^{**}$,\\
   (T3) $\T$ is above $\k$.\\\\
  Let $\W$ be the $\Gamma$-hod pair construction of $(\W^{**}, \d_{y_\a}, \vec{G}^*, \Sigma^*)$ done over $\P$ and relative to $\Sigma$\footnote{See \rdef{gamma-hod pair construction*}.} where 
  \begin{itemize}
  \item $\vec{G}^*$ is the set of those extenders from $\vec{E}^{\W^{**}}$ whose critical point is $>{\sf{ord}}(\P)$ and $\nu(E)$ is an inaccessible cardinal and 
  \item $\Sigma^*$ is the strategy of $\W^{**}$ induced by $\Psi_{y_\a}$.
  \end{itemize} It follows from \rthm{universality of background construction} that there is $\K\inseg_{hod}\W$ which is a  $\Lambda$-iterate of $\R$, and hence, \\\\ 
(1) $L(\Gamma(\K, \Lambda_\K))\models ``\S$, as a $\Sigma_Y$-mouse, has an $\omega_1$-iteration strategy".\\\\
 $\K$ is our universal model but we cannot yet apply $\pi_{0, 1}$ to it. To do this, let $\U=\pi_{0, 1}\T$. The copying construction produces $\sigma:\M^\T_b\rightarrow \M^\U_b$ such that $\pi^\U_b\circ \pi_{0,1 }=\sigma\circ \pi^\T_b$. Moreover, because of (T3) above, $\cp(\sigma)=\k$, $\sigma(\P)=\P_1$ and $\sigma\rest \P=\pi_{0, 1}\rest \P$. It then follows that\\\\
  (2) $\sigma(\K)$ is a $\pi_{0, 1}(\Sigma)$-hod premouse over $\P_1$, and\\
  (3) $\Lambda_\K$ is the $\sigma$-pullback of the strategy of  $\sigma(\K)$ induced by $(\Psi_{y})_{\M^\U_b}$.\\\\
  The reason (3) holds is the following. First notice that $\Lambda_\K$ is the strategy of $\K$ induced by $\Sigma^*$. But for some $\nu<\d_{y_\a}$, we build $\K$ via $\Gamma$-hod pair construction of $\W^{**}|\nu$, and hence $\Lambda_\K$ is the strategy of $\K$ induced by $\Sigma^*_{\W^{**}|\nu}$. $\W^{**}|\nu$ has a unique $\omega_1$-strategy as a $\Phi_1\oplus \Sigma$-mouse, and therefore $\Sigma^*_{\W^{**}|\nu}$ is the strategy induced by $(\Psi_y)_{\M^\T_b}$, and this strategy is the $\sigma$-pullback of the strategy of $\sigma(\W^{**}|\nu)$ which is induced by $(\Psi_y)_{\M^\U_b}$.
  
It now follows that we can lift $\pi_{X, Y}$ to $\K$ and obtain $\pi^{+}_{X, Y}:\K\rightarrow \Y$ and $\tau^+_Y:\Y\rightarrow \sigma(\K)$ such that\\\\
(4) $\sigma \rest \K=\tau^+_Y\circ \pi^{+}_{X, Y}$,\\
(5) $\sigma\rest \P=\pi_{0, 1}\rest \P$.\\\\
 The rest of the proof follows very closely to the proof of \rlem{thm:weakly_condensing_set} and uses (1). This finishes our outline of the proof of Claim 8.
 \end{proof}
 
 Working in $V$, given $A\in \P_{\omega_1}$ and $X\in \powerset_{\omega_1}(\P_{\omega_1})$, we let $T_{X, A}$ be the set of $(\phi, s)$ such that $s\in [\d_X]^{<\omega}$ and $\P_{\omega_1}\models \phi[A, \pi^{\Sigma_X}_{\P_X|\d_X, \infty}(s)]$. We then say that $X$ has \textit{$A$-condensation} if for every $Y\in \powerset_{\omega_1}(\P_{\omega_1})$, $\tau_{X, Y}(T_{X, A})=T_{Y, A}$. To show that there is a strongly $(\phi, \Gamma)$-condensing set, it is enough to show that for each $A$ there is an $X\in \powerset_{\omega_1}(\P_{\omega_1})$ with $A$-condensation. Assuming this, it is not hard to show that for some $\a<\omega_1$, $\pi_{\a, \omega_1}[\P_\a]$ is a condensing set.\\
 
  \textit{Claim 9.} Suppose $A\in \P_{\omega_1}$. There is $\a_0<\omega_1$ such that $A\in \rge(\pi_{\a_0, \omega_1})$ and for every $\a\in (\a_0, \omega_1)$ and for every $\N_{y_\a}^*$-generic $h\subseteq Coll(\omega, <\k_\a)$, $\N_{y_\a}^*[h]\models ``\pi_{\a_0, \a}[\P_{\a_0}]$ is an $A_\a$-condensing set" where $A_\a=\pi_{\a, \omega_1}^{-1}(A)$. \\\\
  \begin{proof}
 Towards a contradiction assume otherwise. Let $(\a_i: i<\omega)$ be such that
 \begin{itemize}
 \item for all $i<\omega$, $\a_i<\omega_1$,
 \item for all $i< \omega$ and for every $\N^*_{y_{\a_{i+1}}}$-generic $h\subseteq Coll(\omega, <\k_{\a_{i+1}})$, $\N_{y_{\a_{i+1}}}[h]\models ``\pi_{\a_i, \a_{i+1}}[\P_{\a_i}]$ is not a $A_{\a_{i+1}}$-condensing set".
 \end{itemize}
 Let $\a=\sup_{i<\omega}\a_i$. Let $\nu_i<\k_{\a_{i+1}}$ be such that for some $\N^*_{y_{\a_{i+1}}}$-generic $h\subseteq Coll(\omega, \nu_i)$ there is $W\in (\powerset_{\omega_1}(\P_{\a_{i+1}}))^{\N^*_{y_{\a_{i+1}}}[h]}$ such that $\N^*_{y_{\a_{i+1}}}[h]\models ``$it is forced by $Coll(\omega, <\k_{\a_{i+1}})$ that $W$ witness that $\pi_{\a_i, \a_{i+1}}[\P_{\a_i}]$ is not a $A_{\a_{i+1}}$-condensing set". Fix then $(h_i, W_i)$ that play the role of $(h, W)$ and such that $h_i\in \N^*_{y_\a}$. Set $Z_i=\pi_{\a_i, \a_{i+1}}[\P_{\a_i}] $. We thus have that\\\\
 (1) in $\N^*_{y_{\a_{i+1}}}[h_i]$, $\tau_{Z_i, W_i}(T^i_{Z_i, A_{\a_{i+1}}})\not =T^i_{W_i, A_{\a_{i+1}}}$ where $T^i_{U, A_{\a_{i+1}}}$ is defined like $T_{U, A}$ above only inside $\N^*_{y_{\a_{i+1}}}[h_i]$\footnote{More precisely, $(s, \phi)\in T^i_{U, A_{\a_{i+1}}}$ if and only if $s\in [\d_U]^{<\omega}$ and $\P_{\a_{i+1}}\models \phi[A_{\a_{i+1}}, \pi^{\Sigma_U}_{\P_U|\d_U, \infty}(s)]$. All of the relevant objects are computed in $\N^*_{y_{\a_{i+1}}}[h_i]$}.\\\\
 Let $X'_i=\pi_{\a_i, \a}[\P_{\a_i}]$ and $Y'_i=\pi_{\a_{i+1}, \a}[W_i]$. It is not hard to verify that\\\\
 (2) in $\N^*_{y_\a}$, $\tau_{X'_i, Y'_i}(T^\a_{X'_i, A_{\a}})\not =T^\a_{Y'_i, A_{\a}}$ where $T^\a_{U, A_\a}$ is defined like $T_{U, A}$ above only inside $\N^*_{y_\a}$.\\\\
 Like in the proof of \rthm{thm:condensing_set}, we can find some $\Y\in \P_\a|\d^{\P_\a}$ and $B\in \Y$ with the property that letting $\sup_{i<\omega}X'_i=_{def}\eta<\d^{\P_\a}$, $\Y|\eta=\P_\a|\eta$ and for every $s\in [\eta]^{<\omega}$ and every $\phi$, $\Y\models \phi[B, s]$ if and only if $\P_\a\models \phi[A_\a, s]$. Let now (in $\N^*_{y_\a}$) $\P_i=\P_{X'_i}$, $\Q_i=\P_{Y'_i}$, $\xi_i=\tau_{X'_i, Y'_i}: \P_i\rightarrow \Q_i$, $\pi_i=\tau_{Y'_i, X'_{i+1}}:\Q_i\rightarrow \P_{i+1}$ and $\phi_i=\tau_{X'_i, X'_{i+1}}$. Finally, set $(C, \X)=\pi_{\a, \omega_1}(B, \Y)$, $X_i=\pi_{\a, \omega_1}(X_i')$ and $Y_i=\pi_{\a, \omega_1}(Y_i')$. It is not hard to verify that\\\\
 (3) in $V$, $\mathcal{A}=\{ (\P_i, \Q_i, X_i, Y_i, \xi_i, \pi_i, \phi_i), C, \X\}$ is a bad tuple relative to $A$ in the sense that
 \begin{enumerate}
\item for all $i<\omega$, $X_i\in \powerset_{\omega_1}(\P_{\omega_1})$ is such that $\tau_{X_i}\rest \P_{X_i}=\pi^{\Sigma_{X_i}}_{\P_{X_i}|\d_{X_i}, \infty}$,
\item  for all $i<\omega$, $\P_i = \P_{X_i}$  and $\Q_i = \P_{Y_i}$;
\item for all $i < j<\omega$, $X_i \prec Y_i \prec X_j$;
\item for all $i<\omega$, $\xi_i=\tau_{X_i, Y_i}$, $\pi_i=\tau_{Y_i, X_{i+1}}$ and $\phi_{i} =\tau_{X_i,X_{i+1}}$\footnote{Thus, $\xi_i:\P_i\rightarrow \Q_i$, $\pi_i: \Q_i \rightarrow \P_{i+1}$ and $\phi_i= \pi_i \circ \xi_i$.};
\item $\X\in \P_{\omega_1}$ and letting $\nu=\sup_{i<\omega}(X_i\cap \d^{\P_{\omega_1}})$, $\X|\nu=\P_{\omega_1}|\nu$;
\item letting $\nu$ be as above, for every formula $\phi$ and for every $s\in \nu^{<\omega}$, $\X\models \phi[C, s]$ if and only if $j(\P)\models \phi[A, s]$; 
\item for all $i\in [1, \omega)$, $\xi_i(T_{X_i, A}) \neq T_{Y_i,A}$. 
\end{enumerate}
 
 As in the proof of Claim 8 we can find some normal stack $\T$ on $\N^*_{y}$ with last model $\W$ such that\\\\
 (3) $\T$ is above $\k$ and $\pi^\T$ is defined,\\
 (4) the $\Gamma$-hod pair construction of 
 \begin{center}
 $(\W, \pi^\T(\d_y), \pi^\T(\vec{G}_y))$
 \end{center}
  done over $\P$ and relative to $\Sigma$ using extenders with critical point $>\k$ reaches a $\Sigma$-hod pair $(\K, \Lambda)$ such that 
 \begin{center}
 $L(\Gamma(\K, \Lambda))\models ``$the sequence $\mathcal{A}=\{ (\P_i, \Q_i, X_i, Y_i, \xi_i, \pi_i, \phi_i), C, \X\}$ is a bad tuple".
 \end{center}
 Let $\U=\pi_{0, \a}\T$, $\W'$ be the last model of $\U$ and let $\sigma:\W\rightarrow \W'$ be the copy map. We have that $\sigma$-pullback of $\sigma(\Lambda)$ is $\Lambda$\footnote{Here we confuse the local strategies with their interpretations in $V$.}. We now finish the proof by performing the following steps:\\\\
Step 1:  lift $\K$ to each $\P_i$ via $\pi_{0, \a_i}$ and obtain $\P_i^+$,\\
Step 2: lift $\K$ to each $\Q_i$ via $\xi_i\circ \pi_{0, \a_i}$ and obtain $\Q_i^+$,\\
Step 3: extend $(\xi_i, \pi_i, \phi_i)$ to $\xi_i^+:\P_i^+\rightarrow \Q_i^+$, $\pi_i^+:\Q_i^+\rightarrow \P_{i+1}^+$ and $\phi_i^+:\P_{i}^+\rightarrow \P_{i+1}^+$,\\
Step 4: find maps $p_i:\P_i^+\rightarrow \sigma(\K)$ and $q_i: \Q_i^+\rightarrow \sigma(\K)$ such that $p_i=q_i\circ \xi_i^+=p_{i+1}\circ \pi_i^+$,\\
Step 5: let $\P^+_\omega$ be the direct limit of $(\P_i^+, \phi_{i, k}: i<k<\omega)$ where $\phi_{i, k}:\P^+_i\rightarrow \P^+_k$ be the composition of $(p_n: n\in [i, k))$, \\
Step 6: let $\phi_{i, \omega}:\P^+_i\rightarrow \P^+_\omega$ and $\psi_{i, \omega}: \Q_i^+\rightarrow \P_\omega^+$ be the direct limit embeddings,\\
Step 7: let $p_\omega^+:\P_\omega^+\rightarrow \sigma(\K)$ be constructed via the direct limit construction,\\
Step 8: set for $i\leq \omega$, $\Pi_i^p=(p_i$-pullback of $\sigma(\Lambda))$ and $\Pi^q_i=(q_i$-pullback of $\sigma(\Lambda))$,\\
Step 9: apply the three dimensional argument from the last portion of the proof of \rthm{thm:condensing_set} to derive a contradiction. \\\\
The above steps finish the proof of Claim 9.
 \end{proof}
 The discussion before Claim 9 implies that there is a $(\phi, \Gamma)$-condensing set, which is clause 3 of \rthm{condensing set for gamma}. Since we were assuming all three clauses of \rthm{condensing set for gamma} are false, this is clearly a contradiction and finishes the proof of \rthm{condensing set for gamma}.
\end{proof}

\chapter{Applications}\label{applications}

\section{The generation of the mouse full pointclasses}\label{generation of mouse full pointclass sec}

In this section, our goal is to show that under $\sf{Strong\ Mouse\ Capturing}$ ($\sf{SMC}$) if $\Gamma$ is a mouse full pointclass (see \rdef{mouse full}) such that $\Gamma\not =\powerset(\bR)$ and there is a good pointclass $\Gamma^*$ with the property that $\Gamma\subset \Gamma^*$ then there is a hod pair or an sts pair $(\P, \Sigma)$ such that $\Gamma(\P, \Sigma)=\Gamma$. Recall that $\sf{SMC}$ states that for any hod pair or sts hod pair $(\P, \Sigma)$ such that $\Sigma$ is strongly fullness preserving and has strong branch condensation then for any $x, y\in \bR$, $x\in OD_{y, \Sigma}$ if and only if $x\in Lp^{\Sigma}(y)$.
We work under the following two minimality assumptions.
\begin{definition}\label{lsa nwlw} 
${\sf{\#_{lsa}}}$ is the following statement: There is a pointclass $\Gamma\subset \powerset(\bR)$ such that there is a Suslin cardinal bigger than $w(\Gamma)$ and $L(\Gamma, \bR)\models \sf{LSA}$.

${\sf{NWLW}}$ is the following statement: There is no iteration strategy for an active mouse with a Woodin cardinal that is a limit of Woodin cardinals. $\myqedhere$
\end{definition}

As in \cite[Section 6.1]{ATHM}, we will construct $(\P, \Sigma)$ as above via a hod pair construction of some sufficiently strong background universe. 
Here is our theorem on generation of pointclasses. 

\begin{theorem}[The generation of the mouse full pointclasses I]\label{the generation of mouse full pointclasses} Assume 
\begin{center}
${\sf{AD^+}}+\neg {\sf{\#_{lsa}}}+{\sf{NWLW}}$\footnote{See \rthm{lst from wlw}, which removes the hypothesis that ${\sf{NWLW}}$ holds.}. 
\end{center}
Suppose $\Gamma\not =\powerset(\mathbb{R})$ is a mouse full pointclass such that $\Gamma\models \sf{SMC}$. Then one of the following holds:
 \begin{enumerate}
\item For some $(\Q, \Lambda)\in {\sf{{\sf{HP}}}}^\Gamma$ such that $\Lambda$ has strong branch condensation and is strongly $\Gamma$-fullness preserving and for some $x\in \bR$, ${\sf{Lp}}^{\Lambda}(x)\not ={\sf{Lp}}^{\Gamma, \Lambda}(x)$.
 \item $\Gamma$ is completely mouse full and letting $A\subseteq \mathbb{R}$ witness the fact that $\Gamma$ is completely mouse full, the following holds in $L(A, \bR)$:
 \begin{enumerate}
 \item $\neg \sf{LSA}$ and there is a hod pair $(\P, \Sigma)$ such that $\Sigma$ has strong branch condensation and is strongly fullness preserving and $\Gamma(\P, \Sigma)=\Gamma$.
 \item $\sf{LSA}$ and there is an sts hod pair $(\P, \Sigma)$ such that $\Sigma$ has branch condensation and is fullness preserving, $\P$ is of $\#$-lsa type\footnote{See \rdef{lsa type}.}  and $\Gamma^b(\P, \Sigma)=\Gamma$. 
 \end{enumerate}
 \end{enumerate}
 Additionally, assuming (i) clause 1 fails, (ii) if $A$ is as in clause 2 then $L(A, \bR)\models {\sf{LSA}}$, and (iii) there is a good pointclass $\Gamma^*$ such that $\Gamma\subset \utilde{\Delta}_{\Gamma^*}$, then there is a hod pair $(\P, \Sigma)$ such that $\P$ is of $\#$-lsa type, $(\P, \Sigma^{stc})\in L(A, \bR)$ and $(\P, \Sigma^{stc})$ satisfies the conditions in clause 2.b. 
\end{theorem}
\begin{proof}
Our proof has the same structure as the proof of \cite[Theorem 6.1]{ATHM}. However, unlike that proof, we will make an important use of \rthm{condensing set for gamma}. The proof is again by induction. Suppose $\Gamma\not =\powerset(\bR)$ is a mouse full pointclass such that whenever $\Gamma^*$ is properly contained in $\Gamma$ and is a mouse full pointclass then there is a hod pair $(\P, \Sigma)$ as in 1 or 2. We want to show that the claim holds for $\Gamma$. Towards a contradiction assume the conclusion of  \rthm{the generation of mouse full pointclasses} is false. We examine several cases.\\

\textbf{Case 1.} There is a sequence of mouse full pointclass $(\Gamma_\a: \a<\Omega)$ such that $\Gamma_\a\subseteq \Gamma$, $\Gamma=\bigcup_{\a<\Omega}\Gamma_\a$ and for $\a<\b<\Omega$, $\Gamma_\a\insegeq_{mouse} \Gamma_\b$.\\

We will use the terminology of \rsec{condensing sets under ad+ sec}. Let $\phi(u, v)$ be the formula that expresses the fact that $u$ is a mouse full pointclass having the properties that $\Gamma$ has and $v$ is a hod pair $(\Q, \Lambda)$ such that ${\sf{Code}}(\Lambda)\in u$ and $\Lambda$ has strong branch condensation and is strongly $u$-fullness preserving. 

Let $\M^-=\P^-_{\phi, \Gamma}$ and $\M=\P_{\phi, \Gamma}$. Because we are assuming that $\Gamma$ is not generated by a hod pair, it follows from clause 2 of \rthm{condensing set for gamma} that $\rho(\M)>o(\M^-)$ and that there is a condensing set $X\in \powerset_{\omega_1}(\M)$. In what follows we will use the notation introduced in \rsec{condensing sets sec}. In particular, recall the definition of $\tau_Y$ and $\sigma^X_Y$. 

Following the proof of \rthm{condensing set for gamma} let $\Gamma_0, \Gamma_0^*,\Gamma_1, \Gamma_1^*, A_0, A_1^*, (N_0, \Phi_0), (N_1, \Phi_1), F_0,$ and $F_1$ be as in that proof. We introduce two more kinds of set of reals that we need to be captured. 

Let $(\a_i: i<\omega)$ be an enumeration of $X$ and let $x_i=(\a_k: k\leq i)$. Let $(\phi_i: i<\omega)$ be an enumeration of formulas in the language of hod mice. Let $B_{i, k}$ be the set of pairs $((\Q, \Lambda, \b), (\R, \Psi, \gg))$ such that $(\Q, \Lambda), (\R, \Psi)\in {\sf{HP}}^\Gamma$, $\b<\d^\Q$, $\gg<\d^\R$ and $\pi^\Psi_{\R, \infty}(\gg)$ is the unique ordinal $\xi$ such that $\M\models \phi_k[x_i, \pi^{\Lambda}_{\Q, \infty}(\b), \xi]$. We then let $A_{i, k}$ be the set of reals $\sigma$ such that $\sigma(0)$ is a G\"odel number of some formula $\zeta$ such that $B_{i, k}$ is definable over $({\sf{HC}}, A_0, \sigma, \in)$ via $\zeta$ without parameters.  

Next, let $B'$ be the set of $(\Q, \Lambda)\in {\sf{HP}}^\Gamma$ such that $X\cap \d^\M\subseteq \pi^\Lambda_{\Q, \infty}[\Q|\d^\Q]$ and the transitive collapse of $Hull^\M(X\cup \pi^\Lambda_{\Q, \infty}[\Q|\d^\Q])$ is $\Q$. Given $(\Q, \Lambda)\in B'$, let $Y_{\Q, \Lambda}=\pi^\Lambda_{\Q, \infty}[\Q|\d^\Q]$. Let $B''$ be the set of $((\Q, \Lambda), X_{\Q, \Lambda})$ such that $(\Q, \Lambda)\in B'$ and $X_{\Q, \Lambda}=\tau_{Y_{\Q, \Lambda}}^{-1}[X]$. Let $B$ be the set of reals $\sigma$ such that $\sigma(0)$ is a G\"odel number of some formula $\zeta$ such that $B$ is definable over $({\sf{HC}}, A_0, \sigma, \in)$ via $\zeta$ without parameters.

We now define our final set $C$. Given $x\in \bR$, let $A_x=\{ u\in \bR: \{u\}$ is $OD^\Gamma_{x, X}\}$. We let $C=\{ (x, y)\in \bR^2: y$ codes $A_x\footnote{$A_x$ is countable. Here we just mean that $y$ lists the members of $A_x$ via the coding introduced in \rdef{pairing function}.}\}$. Let now $x\in \dom(F_0)$ be such that if $F_0(x)=(\N, \M, \d, \Psi)$ then letting $\vec{G}$ be as in clause 7 of \rthm{n*x}, $((\N, \d, \vec{G}, \Psi), (N_0, \Phi_0), \Gamma_0^*, A_0)$ Suslin, co-Suslin captures $\Gamma, A_0, B$ and $C$. Let $\Psi^*$ be the iteration strategy of $\M_3^{\#, \Phi_0}$ and let $y\in \dom(F_1)$ be such that if $F_1(y)=(\N^*_y, \M_y, \d_y, \Sigma_y)$ and letting $\mathbb{M}_y$ be as in clause 3 of \rthm{n*x}, then $(\mathbb{M}_y, (\N_1, \Phi_1), \Gamma_1^*, A_1^*)$ captures ${\sf{Code}}(\Psi^*)$.

We claim that some hod pair appearing on the $\Gamma$-hod pair construction of $\N^*_y|\d_y$ generates $\Gamma$. Here the proof is somewhat different than the proof of Theorem 6.1 of \cite{ATHM}. There the contradictory assumption that such constructions do not reach $\Gamma$ led to a construction of a hod pair $(\P, \Sigma)$ such that $\l^\P=\d^\P$ and $\P\models ``\d^\P$ is regular". This meant that a pointclass satisfying $\sf{AD}_{\bR}+``\Theta$ is a regular cardinal" had been reached giving the desired contradiction. In our current situation, if the constructions never stops then we will reach an lsa type hod premouse $\P$ of height $\d_y$. We need techniques to argue that this cannot happen.  
 
 We proceed by assuming that the $\Gamma$-hod pair construction of $\N^*_y|\d_y$ does not reach a pair generating $\Gamma$. Let $\P^{*}$ be the final model of the $\Gamma$-hod pair construction of $\N^*_y|\d_y$.  By this we mean that either (i) ${\sf{ord}}(\P^*)=\d_y$ and $\P^{*})^\#$ is a hod premouse or (ii) ${\sf{ord}}(\P^*)<\d_y$ and the $\Gamma$-hod pair construction of $\N^*_y|\d_y$ after reaching $\P^*$ does not produce a hod premouse $\Q$ such that $\d^\Q=\d^{\P^*}$. If (i) is true then set $\P=(\P^{*})^\#$ and otherwise set $\P=(\P^*|\d^{\P^*})^\#$.  Let $\Sigma$ be the strategy of $\P$ induced by $\Sigma_y$. Note that $\delta_y$ is not a limit of Woodin cardinals in $\P$ as otherwise $\P\models ``\d_y$ is a Woodin cardinal that is limit of Woodin cardinals", contradicting our smallness assumption.\\

 \textit{Claim 1.} ${\sf{ord}}(\P^*)=\d_y$. \\\\
 \begin{proof}
Suppose not. It follows from \rthm{strong condensation for backgrounded strategies}, \rthm{fullness preservation of background constructions} and \rthm{sts fb constructions converge} that the only way our construction could break down before reaching $\d_y$ is if 
\begin{itemize}
\item $\d^\P<\d_y$ and $\rho(\P)={\sf{ord}}(\P^*)$,
\item $\P$ is of $\#$-lsa type, and 
\item letting $\Lambda=\Sigma_{\P}$, $\Lambda$ has strong branch condensation and is strongly $\Gamma$-fullness preserving,
\item ${\sf{Lp}}^{\Lambda, sts}(\P)\models ``\d^\P$ is a Woodin cardinal". 
\end{itemize}
Because $\Gamma^b(\P, \Lambda^{stc})\subseteq \Gamma$ and $\Gamma^b(\P, \Lambda^{stc})\not =\Gamma$, we can fix $(\R, \Phi)\in {\sf{HP}}^\Gamma$ such that $\R$ is meek and of limit type and $(\P, \Lambda)\in L(\Gamma(\R, \Phi))$. We have that in $L(\Gamma(\R, \Phi))$, $\Lambda$ has strong branch condensation and is strongly fullness preserving. It now follows from \rthm{gamma(p, sigma) in the lsa case} applied in $L(\Gamma(\R, \Phi))$ that for some $\S\in pI(\P, \Lambda)$, $L(\Gamma(\S, \Lambda_\S))\models \sf{LSA}$, contradicting our assumption that $\neg \#_{lsa}$ holds. 
 \end{proof}
 
We thus have that $\d^\P=\d_y$. Let $\k$ be the least $<\d^\P$-strong cardinal of $\P$. For $\a<\k$, let $g_\a=g\cap Coll(\omega, <\a)$.\\
 
 \textit{Claim 2.} $(\P^b, \Sigma_{\P^b})\in B$. \\\\
 \begin{proof}
 Let $g\subseteq Coll(\omega, <\k)$ be $\N^*_y$-generic. We let $(\psi(u, v), D)$ be as in the proof of \rthm{condensing set for gamma}. Following the notation used in the proof of \rthm{condensing set for gamma}, let $\S=\P_{\psi, D}$\footnote{Recall that this is defined in $\N^*_y[g]$.} and $\S^-=\P^-_{\psi, D}$. It follows from the proof of \rthm{condensing set for gamma} that $\rho(\S)>o(\S)$. 
 
 We claim that $\S$ is an iterate of $\P^b$. Let, in $\N^*_y[g]$, $\M'_\infty(\P^b, \Sigma_{\P^b})$ be the direct limit of $\Sigma_{\P^b}$-iterates $\Q$ of $\P^b$ such that  $\P^b$-to-$\Q$ iteration has countable length. Notice that $\M'_\infty(\P^b, \Sigma_{\P^b})\insegeq \S$. This is simply because for every $\Q\insegeq_{hod}\P^b$, $(\Q, \Sigma_\Q)\in {\sf{HP}}^\Gamma$. Suppose then that $\M_\infty(\P^b, \Sigma_{\P^b})\inseg \S$. Let $(\R, \Pi)\in {\sf{HP}}^\Gamma\cap \N^*_y[g]$\footnote{Here we are abusing the notation and use $\Pi$ for both the strategy in $\N^*_y[g]$ as well as its extension in $V$.} be such that $\M_\infty(\P^b, \Sigma_{\P^b})\inseg \M_\infty(\R, \Pi)$. Let $\eta<\k$ be such that $\R\in \N^*_y[g_\eta]$ and there is a $\sigma\in \bR\cap \N^*_y[g_\eta]$ such that $\sigma(0)$ is a G\"odel number for a formula $\zeta$ with the property that ${\sf{Code}}(\Pi)$ is definable over $({\sf{HC}}, A_0, \sigma, \in)$ via $\zeta$ without parameters.
 
 Let $\Q$ be the output of the hod pair construction of $\P$\footnote{Following \rsec{sec: authentication revisited} it can be shown that for each $\nu\in (\k, \d_y)$ which is a successor cardinal of $\P$, the fragment of $\Sigma$ that acts on non-droping iterations based on $\P|\nu$ and are above $\k$ is in $\P$. This allows us to make sense of hod pair constructions. See also the results of \cite[Section 1]{negres}. The above outline uses \rthm{existence of thick sets}.} in which extenders used have critical points $>\kappa$. It follows from \rthm{universality of background construction} that for some $\Q'\inseg_{hod}\l^\Q$, $\Q'$ is a $\Pi$-iterate of $\R$. Let $E\in \vec{E}^\P$ be an extender with critical point $\k$ such that $\nu_E$ is an inaccessible cardinal of $\P$ and $\Q'$ is constructed by the hod pair construction of $\P|\nu_E$. Let $E^*\in \vec{E}^{\N^*_y}$ be the resurrection of $E$. It follows that in $Ult(\N^*_y, E^*)$, some hod pair appearing on the hod pair construction of $\pi(\P^b)$ in which extenders used are bigger than $\kappa$ is a $\Pi$-iterate of $\R$. It then follows that some hod pair appearing on some hod pair construction of $\P^b$ is a $\Pi$-iterate of $\R$. It then follows that ${\sf{Code}}(\Pi)<_w{\sf{Code}}(\Sigma_{\P^b})$ implying that $\M'_\infty(\R, \Pi)\inseg  \M'_\infty(\P^b, \Sigma_{\P^b})$, contradiction (here $\M'_\infty(\R, \Pi)$ is defined similarly to $\M'_\infty(\P^b, \Sigma_{\P^b}$). This contradiction proves the claim.
 \end{proof}

  It is not hard to see, by a simple Skolem hull argument using the fact that $\P\in \N^*_y$, that\\\\
  (1) for a club of $\eta<\d_y$, $(\P|\eta)^\#\models ``\eta$ is a Woodin cardinal". \\\\
Let $C$ be the club in (1). For $\eta\in C$, let $\R_\eta=(\P|\eta)^\#$, $\Sigma_\eta=\Sigma^{stc}_{\R_\eta}$ and $\Q_\eta\insegeq \P$ be the longest $\Sigma_\eta$-sts mouse such that $\Q_\eta\models ``\eta$ is a Woodin cardinal". Using \rlem{s-construction lemma}, we can translate $\Q_\eta$ into $\Sigma_\eta$-sts mouse $\bar{\Q_\eta}$ over $(\N^*_y|\eta)^\#$. Notice that\\\\
(2) for every $\eta$, $\bar{\Q_\eta}$ has an iteration strategy $\Delta$ witnessing that $\bar{\Q_\eta}$ is a $\Sigma_\eta$-sts mouse over $\mathcal{J}_\omega[(\N^*_y|\eta)^\#]$ based on $\R_\eta$. \\\\
(2) is a consequence of the fact that $\Q_\eta$ appears on a $\Gamma$-hod pair construction of $\N^*_y$. Moreover,\\\\
(3) for every $\eta$ and for every real $x$ coding $\N^*_y|\eta$, $\bar{\Q}_\eta$ is $OD^\Gamma_{x, X}$.\\\\
(3) follows from proofs that have already appeared in the book. For instance, see the notion of goodness that appeared in the proof of \rlem{wadge rank computation}. We now claim that\\\
 
\textit{Claim 3.} for a club of $\eta\in C-(\k+1)$, $\Q_\eta \in \mathcal{J}^{\Psi^*}_{\nu_\eta}(\N^*_y|\eta)$ where $\nu_\eta$ is the least ordinal such that $\mathcal{J}^{\Psi^*}_{\nu_\eta}(\N^*_y|\eta)\models \sf{ZFC}$.\\\\
\begin{proof}
Let $\l$ be least such that $\mathcal{J}_\l^{\Psi^*}(\N^*_y|\d_y)\models \sf{ZFC}$. Let $\eta\in C$ be such that there is a map $\pi: \mathcal{J}^{\Psi^*}_{\nu_\eta}(\N^*_y|\eta)\rightarrow \mathcal{J}_\l^{\Psi^*}(\N^*_y|\d_y)$. Thus\\\\
(4) $\mathcal{J}^{\Psi^*}_{\nu_\eta}(\N^*_y|\eta)\models ``\eta$ is a Woodin cardinal". \\

Using genericity iterations, we can find $\N\in \mathcal{J}^{\Psi^*}_{\nu_\eta}(\N^*_y|\eta)$ such that $\N$ is a $\Psi^*$-iterate of $\M_3^{\#, \Phi_0}$ such that $(\N^*_y|\eta)^\#$ is generic over the extender algebra $\mathbb{B}_{\d_0}^\N$ where $\d_0$ is the least Woodin cardinal of $\N$ that is $>\eta$. Let $h\subseteq Coll(\omega, \eta)$ be $\N^*_y$-generic. Fix a real $x\in \N[\N^*_y|\eta][h]$ coding $\N^*_y|\eta$. It follows that there is $y\in \bR$ such that $(x, y)\in C \cap \N[\N^*_y|\eta][h]$. Therefore $\Q_\eta\in \N[\N^*_y|\eta][h]$. As $x$ is arbitrary, we have that $\Q_\eta\in \mathcal{J}^{\Psi^*}_{\nu_\eta}(\N^*_y|\eta)$. It follows that $\mathcal{J}^{\Psi^*}_{\nu_\eta}(\N^*_y|\eta)\models ``\eta$ is not a Woodin cardinal", contradicting (4). 
\end{proof}

The rest of the proof is easy. It follows from Claim 3 that we can find an $\eta$ such that $\Q_\eta \in \mathcal{J}^\Psi_{\nu_\eta}(\N^*_y|\eta)$ and there is an elementary embedding $\pi: \mathcal{J}^{\Psi^*}_{\nu_\eta}(\N^*_y|\eta)\rightarrow \mathcal{J}_\l^{\Psi^*}(\N^*_y|\d_y)$ where $\l$ is the least such that $\mathcal{J}_\l^{\Psi^*}(\N^*_y|\d_y)\models \sf{ZFC}$. Because $\Q_\eta \in \mathcal{J}^\Psi_{\nu_\eta}(\N^*_y|\eta)$, we have that $\mathcal{J}^{\Psi^*}_{\nu_\eta}(\N^*_y|\eta)\models ``\eta$ is not a Woodin cardinal", and because $\pi: \mathcal{J}^{\Psi^*}_{\nu_\eta}(\N^*_y|\eta)\rightarrow \mathcal{J}_\l^{\Psi^*}(\N^*_y|\d_y)$, we have that $\mathcal{J}^{\Psi^*}_{\lambda}(\N^*_y|\delta_y)\models ``\delta_y$ is not a Woodin cardinal". This is an obvious contradiction! Thus, we must have that the $\Gamma$ hod pair construction of $\N^*_y$ reaches a generator for $\Gamma$. We now move to case 2.\\

\textbf{Case 2.} $\Gamma$ is a completely mouse full pointclass such that for some $\a$, $L(\Gamma, \bR)\models \theta_{\a+1}=\Theta$.\\

Because we are assuming $\neg\#_{lsa}$, we must have that $L(\Gamma, \bR)\models \neg \sf{LSA}$. The rest of the proof is very much like the proof of \cite[Theorem 6.1]{ATHM}. To complete it, we need to use \rthm{computation of hod} instead of \cite[Theorem 4.24]{ATHM}. We leave the details to the reader. The proof of ``additionally" clause is similar to Case 1 and uses \rthm{universality of background construction} and \rlem{wadge rank computation}.
\end{proof}

\rthm{the generation of mouse full pointclasses} has one shortcoming. It cannot be used to compute $\H$ of the minimal model of $\sf{LSA}$ as it only generates pointclasses whose Wadge ordinal is strictly smaller than the largest Suslin cardinal. To compute $\H$ of the minimal model of $\sf{LSA}$ we will need the following theorem.

\begin{theorem}\label{generating min model of lsa} Assume $\sf{AD}^++\sf{LSA}$ and suppose $\neg\#_{lsa}+{\sf{NWLW}}\footnote{As in the previous theorem, \rthm{lst from wlw} removes the extra assumption that $\sf{NWLW}$ holds.}$. Let $\a$ be such that $\theta_{\a+1}=\Theta$, and suppose that there is a hod pair or an sts hod pair $(\P, \Sigma)$ such that 
\begin{itemize}
\item $\Sigma$ is strongly fullness preserving and has strong branch condensation and
\item $\Gamma^b(\P, \Sigma)=\{ A\subseteq \bR: w(A)<\theta_\a\}$. 
\end{itemize}
Then $(\P, \Sigma)$ is an sts hod pair and for any $B\in \mathbb{B}[\P, \Sigma]$ there is $\Q\in pI(\P, \Sigma)$ such that $(\Q, \Sigma_\Q)$ is strongly $B$-iterable. 
\end{theorem}
\begin{proof}
Towards a contradiction, assume not. We reflect the failure of our claim to $\utilde{\Delta}^2_1$. Let $(\b, \gg)$ be lexicographically least such that letting $\Gamma=\{A\subseteq \bR: w(A)<\gg\}$,
\begin{enumerate}
\item $\Gamma=\powerset(\bR)\cap \mathcal{J}_\b(\Gamma, \bR)$ and $L_\b(\Gamma, \bR)\models \sf{LSA+ZF-Powerset}$,
\item letting $\a$ be such that $L_\b(\Gamma, \bR)\models ``\theta_{\a+1}=\Theta"$, $L_\b(\Gamma, \bR)\models ``$there is a hod pair or an sts hod pair $(\P, \Sigma)$ such that $\Sigma$ is strongly fullness preserving and has strong branch condensation and $\Gamma^b(\P, \Sigma)=\{ A\subseteq \bR: w(A)<\theta_\a\}$ but either
\begin{enumerate}
\item $(\P, \Sigma)$ is not an sts hod pair or
\item there is a $B\in \mathbb{B}[\P, \Sigma]$ such that whenever $\Q\in pI(\P, \Sigma)$, $(\Q, \Sigma_\Q)$ is not strongly $B$-iterable". 
\end{enumerate}
\end{enumerate}
Because $(\b, \gg)$ is minimized, we have that $\Gamma\subset \utilde{\Delta}^2_1$. Fix $(\P, \Sigma)$ as above. First we claim that \\

\textit{Claim.} $\Sigma$ is not an iteration strategy. \\\\
\begin{proof}
Towards a contradiction suppose not. Let $A_0\in lub(\Gamma)$, $\Gamma^*$ be a good pointclass beyond $\Gamma$ and $(N_0, \Psi_0)$ be a $\Gamma^*$-Woodin which Suslin, co-Suslin captures $(T_n(A_0): n<\omega)$\footnote{Here, $T_n(A_0)$ is defined the way $T_n(\Psi)$ is defined in \rsubsec{subsec: capturing pointclasses}.}. Let $F$ be as in \rthm{n*x} for $\Gamma^*$, and let $x\in \dom(F)$ be such that letting $F(x)=(\N^*_x, \M_x, \d_x, \Sigma_x)$ and $\mathbb{M}_x$ be as in clause 7 of \rthm{n*x}, $(\mathbb{M}_x, (N_0, \Psi_0), \Gamma^*, A_0)$ Suslin, co-Suslin captures ${\sf{Code}}(\Sigma)$ and $\Gamma$. It follows that ${\sf{Le}}((\P, \Sigma), \mathcal{J}_{\omega}[\P])^{\N^*_x|\d_x}$ reaches $\M_2^{\#, \Sigma}$. Let $\Psi$ be the iteration strategy of $\M_2^{\#, \Sigma}$. Notice that\\\\
(1) $\Psi\in L_\b(\Gamma, \bR)$.\\\\
 Because $\Sigma$ is an iteration strategy, it follows from clause 1 of \rthm{main theorem on gen int} that there are trees $(T, S)\in \M_2^{\#, \Sigma}$ such that letting $\d_0<\d_1$ be the Woodin cardinals of $\M_2^{\#, \Sigma}$
\begin{enumerate}
\item $\M_2^{\#,\Sigma}\models ``(T, S)$ are $\d_1$-complementing",
\item whenever $\pi: \M_2^{\#, \Sigma}\rightarrow \N$ is an iteration according to $\Psi$ and $g\subseteq Coll(\omega, \pi(\d_0))$ is $\N$-generic then ${\sf{Code}}(\Sigma)\cap \bR^{\N|\d_1[g]}=p[\pi(T)]$ and $({\sf{Code}}(\Sigma))^c\cap \bR^{\N|\d_1[g]}=p[\pi(S)]$.
\end{enumerate}
Let $\M_\infty$ be the direct limit of all $\Psi$-iterates of $\M_2^{\#,\Sigma}$ and let $\pi:\M_2^{\#, \Sigma}\rightarrow \M_\infty$ be the direct limit embedding. It then follows that ${\sf{Code}}(\Sigma)=p[\pi(T)]$ and $({\sf{Code}}(\Sigma))^c=p[\pi[S]]$. It follows from (1)  that $\pi(T), \pi(S)\in L(\Gamma, \bR)$, implying that $L(\Gamma, \bR)\models ``{\sf{Code}}(\Sigma)$ is Suslin, co-Suslin". It follows that ${\sf{Code}}(\Sigma)\in \Gamma(\P, \Sigma)$, contradiction!
\end{proof}

It follows from Claim 1 that $(\P, \Sigma)$ is an sts hod pair. Hence, we must have that\\\\
(2) there is $B\in \mathbb{B}[\P, \Sigma]$ such that whenever $\Q\in pI(\P, \Sigma)$, $(\Q, \Sigma_\Q)$ is not strongly $B$-iterable.\\\\ 
We can now finish by following the proof of \rthm{getting strong b-iterability}. The only issue is that in  \rthm{getting strong b-iterability} we require that $\Sigma$ be a strategy, but this is only needed externally. In our current context, we need a strategy $\Sigma^*$ that extends $\Sigma$ and is $L_\b(\Gamma, \bR)$-fullness preserving and has branch condensation. Obtaining such a $\Sigma^*$ might require passing to a $\Sigma$-iterate of $\P$. We can obtain such a strategy by further iterating $\P$ via $\Sigma$ to a hod pair construction of a sufficiently strong background triple using the theory developed in \rsec{sec:normal_comparison}. Notice that (2) holds even for this new pair, and so without loss of generality we may just as well assume that $\Sigma^*$ exists. The rest is just like in the proof of  \rthm{getting strong b-iterability}.
 \end{proof}
 
 \begin{theorem}[The generation of the mouse full pointclasses II]\label{the generation of mouse full pointclasses ii} Assume 
\begin{center}
${\sf{AD^+}}+\neg {\sf{\#_{lsa}}}+{\sf{NWLW}}$\footnote{See \rthm{lst from wlw}, which removes the hypothesis that ${\sf{NWLW}}$ holds.}. 
\end{center}
Suppose that
\begin{itemize}
\item $\Gamma\not =\powerset(\mathbb{R})$ is a mouse full pointclass such that $\Gamma\models \sf{SMC}$ and
\item for some $(\Q, \Lambda)\in {\sf{{\sf{HP}}}}^\Gamma$ such that $\Lambda$ has strong branch condensation and is strongly $\Gamma$-fullness preserving and for some $x\in \bR$, ${\sf{Lp}}^{\Lambda}(x)\not ={\sf{Lp}}^{\Gamma, \Lambda}(x)$.
\end{itemize}
Then there is an anomalous pair\footnote{See \rdef{anomalous hod pair}.} $(\P, \Sigma)$ such that 
\begin{itemize}
\item $\Sigma$ has strong branch condensation and is strongly $\Gamma$-fullness preserving and
\item $\Gamma(\P, \Sigma)=\Gamma$.
\end{itemize}
\end{theorem}
\begin{proof}
The proof is similar to the proof of \rthm{the generation of mouse full pointclasses} but here we need to revise \rthm{fullness preservation of background constructions} and \rthm{strong condensation for backgrounded strategies}. There, the strong fullness preservation and strong branch condensation are proved using the method of thick hulls developed in \rsec{sec: thick hulls}. Here, we need to use the fact that\\\\
(*) for some $(\Q, \Lambda)\in {\sf{{\sf{HP}}}}^\Gamma$ such that $\Lambda$ has strong branch condensation and is strongly $\Gamma$-fullness preserving and for some $x\in \bR$, ${\sf{Lp}}^{\Lambda}(x)\not ={\sf{Lp}}^{\Gamma, \Lambda}(x)$.\\\\
 The following is the main way (*) affects the proof of \rthm{fullness preservation of background constructions}: For example, we can no longer assume that if $\tau$ and $\N$ are as in that proof (just after (b)) then (4) of that proof holds.  Below we outline our method of dealing with the aforementioned issue. 
 
 Towards contradiction assume not and suppose $\Gamma$ is the least pointclass satisfying our hypothesis which is not generated as stated above. The following are the two lemmas that we need to prove.
 
\begin{lemma}\label{fp with projecting mouse} Suppose ${\sf{C}}=(\mathbb{M}, (P, \Psi), \Gamma^*, A)$ Suslin, co-Suslin captures $\Gamma$ and $\mathbb{M}=(M, \d,  \vec{G}, \Sigma)$. Set
\begin{center} 
${\sf{hpc}}_{{\sf{C}}, \Gamma}^+=(\M_\gg , \N_\gg, Y_\gg, \Phi^+_\gg, F^+_\gg,  F_\gg, b_\gg: \gg\leq \d)$.
\end{center} 
Suppose $\b<\d$, $\P\in Y_\b$ and $M\models ``(\P, (\Phi_\b)_\P)\in {\sf{HP}}^\Gamma"$. Then $(\Phi_\b^+)_\P$ is almost low-level strongly $\Gamma$-fullness preserving.\footnote{In this context, it may not be the case that $\P$ is $\Gamma$-full at the top. Meaning, if $\P$ is meek of limit type then it may not be the case that $\P={\sf{Lp}}^{\Gamma, (\Phi_\b)_{\P|\d^\P}}(\P|\d^\P)$. For example, this may happen if $\P=(\P|\d^\P)^\#$ and $\Gamma=\Gamma(\P, (\Phi_\b)_\P)$. In this context, fullness preservation is meant to be for lower level strategies. See \rdef{gamma fullness preservation}.}.
\end{lemma}

\begin{lemma}\label{bc with projecting mouse} Suppose ${\sf{C}}=(\mathbb{M}, (P, \Psi), \Gamma^*, A)$ Suslin, co-Suslin captures $\Gamma$ and $\mathbb{M}=(M, \d,  \vec{G}, \Sigma)$. Set
\begin{center} 
${\sf{hpc}}_{{\sf{C}}, \Gamma}^+=(\M_\gg , \N_\gg, Y_\gg, \Phi^+_\gg, F^+_\gg,  F_\gg, b_\gg: \gg\leq \d)$.
\end{center} 
Suppose $\b<\d$, $\P\in Y_\b$ and $M\models ``(\P, (\Phi_\b)_\P)\in {\sf{HP}}^\Gamma"$. Then $(\Phi_\b^+)_\P$ has branch condensation. 
\end{lemma}

In both of those cases, the hard case is when $\P$ is of limit type and $\Gamma$ is of limit type. We assume this and set $\Lambda=(\Phi_\b)_\P$. We will only outline the proof of fullness preservation. Strong fullness preservation can be established by a very similar argument. We then have that there is a witness to non fullness preservation or non branch condensation added by a collapse of some $\nu<\d$. Let then $g\subseteq Coll(\omega, \nu)$ be this collapse. We now outline how to proceed assuming that $\P$ is not gentle (there is not much to prove in this case). This in particular implies that $\P={\sf{Lp}}^{\Gamma, \Lambda}(\P|\d^{\P^b})$. 

In the case of fullness preservation this witness is a tuple $(\T, \M)$ such that $\T$ is according to $\Lambda$ and $\M$ is a mouse witnessing the failure of one of the clauses of fullness preservation. The key is that $\M$ has an iteration strategy coded by a set in $\Gamma$. Let $\Gamma'$ be a good pointclass contained in $\Gamma$ and such that the iteration strategy of $\M$ is in $\Delta_{\Gamma'}$. Let $\eta$ be the least $\Gamma'$-Woodin cardinal of $M$ above $\iota=_{def}\max(\nu, \zeta)$ where  $\zeta=\sup\{\lh(F^+_\gg): \gg<\b\}$. We now repeat the proof of  \rthm{fullness preservation of background constructions} while working inside $M'=C_{\Gamma'}(M|\eta)$. Let $\tau=\d^{\P^b}$ and $\vec{G}'=\{F\in \vec{G}: \cp(F)>\iota\}$. The key point is that when in that proof we let $\N$ be the last model of 
\begin{center}
 $({\sf{Le}}((\P|\tau, \Lambda_{\P|\tau}), \P^b)_{>\zeta})^{(M'[g], \eta, \vec{G}')}$,
 \end{center} 
we have that no level of $\N$ projects across $\P^b$. This is because if $\Sigma'$ is the fragment of $\Sigma_{M'}$ that acts on stacks that are above $\iota$ then ${\sf{Code}}(\Sigma')\in \Gamma$.

The issue with branch condensation can be resolved similarly. In the case of branch condensation, the witness is $(\T, \U, b, \sigma)\in M[g]$ such that\\\\
(i) $\T$ is according to $\Lambda$, $\pi^\T$ exists and $\T$ has a last model $\R$,\\
(ii)  $\U$ is according to $\Lambda$ and is of limit length,\\
(iii) $b$ is a cofinal branch of $\U$ and $\Lambda(\U)\not =b$,\\
(iv) $\sigma: \M^\U_b\rightarrow \R$ is such that $\pi^\T=\sigma\circ \pi^\U_b$.\\\\
The dificult case is when $\P$ is non-meek, and so we assume this. We assume $\Lambda$ is an iteration strategy as the other case is very similar. Let $\zeta$ and $\tau$ be as above.

 The most dificult case is when $\pi^{\U, b}$ is defined, $\Q(b, \U)$ exists and is an sts premouse over $\m^+(\U)$. Other cases follow the same pattern but this one is the most involved. Here we need to show that $\Q(b, \U)$ is in fact $(\Lambda_{\m^+(\U), \U})^{stc}$-sts mouse over $\m^+(\U)$. Let $\eta$ be such that $\sigma\rest \Q(b, \U):\M^\U_b\rightarrow \R||\eta$. Because ${\sf{Code}}((\Lambda_{\R|\eta, \T})^{stc})\in \Gamma$, we have that if $\Phi$ is the $\sigma\rest \Q(b, \U)$-pullback of $(\Lambda_{\R|\eta, \T})^{stc}$ then ${\sf{Code}}(\Phi)\in \Gamma$. Thus, it is enough to show that, setting $\Psi=(\Lambda_{\m^+(\U), \U})^{stc}$\\\\
 (a) $\Phi=\Psi$.\\\\
Let $\W=\m^+(\U)$. Suppose $\Phi\not =\Psi$, and let $(\X_0, \W_0, \X_1, \W_1, \Y)$ be a minimal disagreement\footnote{See \rdef{low level disagreement between strategies}.} between $\Phi$ and $\Psi$. We have that ${\sf{Code}}(\Phi_{\Y, \U^\frown \X_0})$ and ${\sf{Code}}(\Psi_{\Y, \U^\frown \X_1})$ are in $\Gamma$. Let then $\Gamma'$ be a good pointclass contained in $\Gamma$ such that 
\begin{center}
$\{{\sf{Code}}(\Phi_{\Y, \U^\frown \X_0}), {\sf{Code}}(\Psi_{\Y, \U^\frown \X_1})\}\subseteq \Delta_{\Gamma'}$.
\end{center}
Let now $\eta, M', \vec{G}', \Sigma'$ be as above defined relative to the new meaning of $\Gamma'$. Again the proof now simply follows the proof of \rthm{strong condensation for backgrounded strategies}.

The next major issue to deal with is when we pass from gentle stage to the next hod premouse. Just like in the proof of \rthm{the generation of mouse full pointclasses}, we build our desired generator for $\Gamma$ via a hod pair construction of some background triple. We find some ${\sf{C}}=(\mathbb{M}, (P, \Psi), \Gamma^*, A)$ that Suslin, co-Suslin captures $\Gamma$ where $\mathbb{M}=(M, \d,  \vec{G}, \Sigma)$. Fix a hod pair $(\Q, \Phi)\in {\sf{HP}}^\Gamma$ such that $\Phi$ has strong branch condensation and is $\Gamma$-fullness preserving and for some $x\in \bR$, ${\sf{Lp}}^{\Phi}(x)\not ={\sf{Lp}}^{\Gamma, \Phi}(x)$. We now get that\\\\
(1) whenever $i:\Q\rightarrow \Q'$ is an iteration according to $\Phi$ and whenever $y\in \bR$ is a real Turing above $x$, 
\begin{center}
${\sf{Lp}}^{\Phi_{\Q'}}(y)\not ={\sf{Lp}}^{\Gamma, \Phi_{\Q'}}(y)$.
\end{center}

(1) can be established via more or less standard arguments. For example, see \cite[Lemma 6.21 ]{ATHM}. The key ingredient of the proof is that if (1) fails for $\Q'$ and $y$ then any universal $\Phi_\Q$-mouse over which $i$ and $\Q$ are set generic is also $\Phi$-universal. To find such a universal $\Phi_\Q$-mouse, we can choose some good pointclass $\Gamma'$ such that $\Phi, \Phi'\in \Delta_{\Gamma'}$ where $\Phi'$ is the iteration strategy of ${\sf{Lp}}^{\Phi}(x)$. Let then $F$ be as in \rthm{n*x} for $\Gamma'$ and let $z\in \dom(F)$ be such that $(\mathbb{M}_z, (N, \Psi), \Gamma^*, A)$ Suslin, co-Suslin captures some set of reals coding $(\Q, \Phi_{\Q'}, y, \Phi', i, x)$, and where the rest of the objects are defined as in clause 7 of \rthm{n*x}. The desired universal $\Phi_{\Q'}$-mouse is 
\begin{center}
 $({\sf{Le}}((\Q', \Phi_{\Q'}), y)^{(M, \d, \vec{G})}$.
 \end{center} 
 Letting $\N$ be that model, we have that $\N$ has a Woodin cardinal and $(i, \Q, x)$ is set generic over $\N$. It then follows from the universality of $\N$ that ${\sf{Lp}}^{\Phi}(x)\in \N[(i, \Q, x)]$ and if $\K\insegeq{\sf{Lp}}^{\Phi}(x)$ then $\K$ appears in some fully backgrounded construction of $\N[(i, \Q, x)]$.  We leave the details to the reader.
 
For each $\Q'\in pI(\Q, \Phi)$ and for each $y\in \bR$ Turing above $x$ let $\M_{\Q', y}\insegeq {\sf{Lp}}^{\Phi_{\Q'}}(y)$ be the least such that $\M_{\Q', y}$ does not have an iteration strategy in $\Gamma$ (as a $\Phi_{\Q'}$-mouse). Let $\Psi_{\Q', y}$ be the unique iteration strategy of $\M_{\Q', y}$. In addition to the requirements mentioned above, we demand that $A$ that appeared in ${\sf{C}}$ codes the set of all triple that have the form $(\Q', y, \Psi_{\Q', y'})$. In particular, $\Q, x\in M$ and $\Phi_{\Q, x}$ is Suslin, co-Suslin captured by ${\sf{C}}$. 

Set now
\begin{center} 
${\sf{hpc}}_{{\sf{C}}, \Gamma}^+=(\M_\gg , \N_\gg, Y_\gg, \Phi^+_\gg, F^+_\gg,  F_\gg, b_\gg: \gg\leq \d)$.
\end{center} 
Because of our set up, there is $\gg<\d$ such that $\M_\gg\in pI(\Q, \Phi)$. This implies that the construction cannot last $\d$ steps. If it did, then because $\Phi^+_\gg=\Phi_{\M_\gg}$\footnote{See \rthm{universality of background construction}.}, letting $\N$ be the last model of
\begin{center}
 $({\sf{Le}}((\M_\gg, \Phi_{\M_\gg}))^{(\N_\d, \d, \vec{G}')}$,
 \end{center}  where $\vec{G}'=\{ F\in \vec{E}^{\N_\d}: \forall \gg'<\gg (\nu(F^+_{\gg'})<\cp(F))$ and $\nu(F)$ is an inaccessible cardinal of $\N_\d\}$, some fully backgrounded construction of some set generic extension of $\N$ would reach $\M_{\M_\gg, x}$. This would imply that ${\sf{Code}}(\Phi_{\M_\gg, x})\in \Gamma$, contradiction. We thus have that the construction has to stop. 

Because clause 4a of \rdef{gamma-hod pair construction*} never occurs, we must have that clause 4b occurs. Let then $\xi$ be such that $\N_\xi$ has the property described in clause 4b of \rdef{gamma-hod pair construction*}. We have that $\N_\xi$ is germane\footnote{See \rdef{germane lses}.}. Let $\P=\N_\xi$. Let $\Psi=\Phi^+_\xi$. The following is our main claim.\\
 
  \textit{Claim.} $\Gamma(\P, \Psi)=\Gamma$.\\\\
  \begin{proof} Because $\Psi$-iterates of $\P$, via the resurrection process, embed into hod mice whose iteration strategies are in $\Gamma$, we have that  $\Gamma(\P, \Psi)\subseteq \Gamma$. It remains to show that $\Gamma\subseteq \Gamma(\P, \Psi)$. Assume then that $\Gamma(\P, \Psi)\subset \Gamma$. 
  
 Using \rthm{branch condensation on a tail}\footnote{\rthm{branch condensation on a tail} is applicable because $(\P, \Psi)$ is an anomalous hod pair and if it is of type III then we can produce a supporting bicephalous via fully backgrounded construction done over $\P^b$ relative to $\Psi_{\P|\d^{\P^b}}$. The proof of (1) above shows that this construction will reach an $\M$ with the property that $\rho(\M)\leq \d^{\P^b}$. The arguments presented on page 142 of \cite{ATHM} then show that if $\M$ is the least such level of the aformentioned backgrounded construction then in fact $\rho(\M)<\d^{\P^b}$.}, we can find some tail $(\P', \Psi')$ of $(\P, \Psi)$ such that $\P'\in pI(\P, \Psi')$, $\Psi'$ has branch condensation and $\Gamma(\P', \Psi')=\Gamma(\P, \Psi)$. It then follows that ${\sf{Code}}(\Psi)\in \Gamma$, which can be shown by using the proof of \rlem{wadge rank computation}. It now follows that ${\sf{Code}}(\Psi)\in \Gamma$. 
 
 Notice next that by induction we can assume that if $\S\inseg^c_{hod} \P$ then $\Psi_\S$ is $\Gamma$-fullness preserving and has branch condensation. This means that we can now apply the proof of 2a and 2b on page 142 of \cite{ATHM} to conclude that $\rho(\P)\geq {\sf{ord}}({\sf{hl}}(\P))$, which is a contradiction\footnote{${\sf{hl}}$ is defined in \rdef{germane lses}.}. 
  \end{proof}
  The desired hod pair generating $\Gamma$ is the tail of $(\P, \Psi)$ provided by \rthm{branch condensation on a tail}.
\end{proof}

\section{A proof of the Mouse Set Conjecture below LSA}

Throughout we will assume $\sf{AD}^{++}=_{def}\sf{AD}^++V=L(\powerset(\bR))$.\index{$\sf{AD}^{++}$} Recall the definition of ${\sf{\#_{lsa}}}$ and ${\sf{NWLW}}$ defined in the previous section\footnote{See \rdef{lsa nwlw}.}. Recall that $\sf{Strong\ Mouse\ Capturing}$ ($\sf{SMC}$) is the statement that for any hod pair or an sts hod pair $(\P, \Sigma)$ such that $\Sigma$ has strong branch condensation and is strongly fullness preserving, and for any reals $x, y$, $x$ is ordinal definable from $\Sigma$ and $y$ if and only if $x$ is in some $\Sigma$-mouse over $y$. The following is the main theorem of this section.

\begin{theorem}\label{msc} Assume $\sf{AD}^{++}+\neg \#_{lsa}+NWLW$. Then the $\sf{Strong\ Mouse\ Capturing}$ holds. 
\end{theorem}

The rest of this section is devoted to the proof of \rthm{msc}.  We assume familiarity with the proof of \cite[Theorem 6.19]{ATHM} and build directly on it.  We start by stating the main steps of \cite[Theorem 6.19]{ATHM}. We will follow these steps and provide proofs only for the new cases. 

Towards a contradiction assume that $\sf{SMC}$ is false. Our first step is to locate the minimal level of the Wadge hierarchy over which $\sf{SMC}$ becomes false. For simplicity we assume that the ${\sf{Mouse\ Capturing}}$, instead of the ${\sf{Strong\ Mouse\ Capturing}}$, is false. ${\sf{Mouse\ Capturing}}$ is the same as $\sf{SMC}$ when the pair $(\P, \Sigma)=\emptyset$. The general case is only different in one aspect, it needs to be relativized to some strategy or a short tree strategy $\Sigma$. 

\begin{notation}\label{the gamma} Throughout this section, we let $\Gamma$ be the least Wadge initial segment such that for some $\a$
\begin{enumerate}
\item $\Gamma=\powerset(\bR)\cap L_\a(\Gamma, \bR)$,
\item $L_\a(\Gamma, \bR)\models \sf{SMC}$,
\item there are reals $x$ and $y$ such that $L_{\a+1}(\Gamma,\mathbb{R})\models ``y$ is $OD(x)$" yet no $x$-mouse has $y$ as a member. 
\end{enumerate}
$\myqedhere$
\end{notation}
For the purposes of this section we make the following definition.
\begin{definition}\label{perfect hod pair} Suppose $(\P, \Sigma)$ is a hod pair and $\Gamma^*$ is a projectively closed pointclass\footnote{See \rdef{projectively closed pointclass}.}. We say $(\P, \Sigma)$ is \textbf{$\Gamma^*$-perfect} if the following conditions are met.
\begin{enumerate}
\item $\Sigma$ is $\Gamma^*$-strongly fullness preserving and has strong branch condensation.
\item For every $\Q\in pI(\P, \Sigma)\cup pB(\P, \Sigma)$ such that $\Q$ is of successor type, there is $\vec{B}=(B_i: i\leq \omega)\subseteq \mathbb{B}[\Q^-, \Sigma_{\Q^-})$ such that $\vec{B}$ strongly guides $\Sigma_{\Q}$. 
\end{enumerate}
If $\Gamma^*=\powerset(\bR)$ then we omit $\Gamma^*$ from our notation.  $\myqedhere$
\end{definition}

The following theorem was heavily used in \cite{ATHM}. It is essentially due to Steel and Woodin (see \cite{TWMS}).

\begin{theorem}\label{mc in l(sigma, r)} Assume $\sf{AD}^+$ and suppose $(\P, \Sigma)$ is a hod pair or an sts hod pair such that $L(\Sigma, \bR)\models ``(\P, \Sigma)$ is perfect". Then $L(\Sigma, \bR)\models \sf{MC}$$(\Sigma)$. 
\end{theorem}

A key theorem used in the proof of \rthm{msc} is the following capturing theorem. Its precursor is stated as  \cite[Theorem 6.5]{ATHM}.

\begin{theorem}\label{capturing of hod pairs} Suppose $(\P, \Sigma)$ is a perfect hod pair and $\Gamma_1$ is a good pointclass such that ${\sf{Code}}(\Sigma)\in \utilde{\Delta}_{\Gamma_1}$. Suppose $F$ is as in \rthm{n*x} for $\Gamma_1$ and $z\in \dom(F)$ is such that if $F(z)=(\N^*_z, \M_z, \d_z, \Sigma_z)$ then $(\N^*_z, \d_z, \Sigma_z)$ Suslin, co-Suslin captures ${\sf{Code}}(\Sigma)$\footnote{We abuse the terminology and omit the other object used to express this type of capturing. In the sequel, if the nature of these other objects, like the pair $(N, \Psi)$, is not important we will omit them from the discussions.}. Let $\N=({\sf{Le}}(\emptyset))^{\N^*_z|\d_z}$\footnote{This is just the ordinary fully backgrounded construction. See \rdef{fully backgrounded sts construction}.}. Then there is $\Q\in pI(\P, \Sigma)\cap \N$ such that $\Sigma_\Q\rest \N\in \mathcal{J}[\N]$. 
\end{theorem} 

The next key lemma that is used in the proof of \rthm{msc} is the following generation lemma that can be traced to \cite[Lemma 6.23]{ATHM}. Below $\Gamma$ is as in \rnot{the gamma}.

\begin{lemma}\label{generating gamma} There is a perfect pair $(\P, \Sigma)$ such that 
\begin{center}
$\Gamma(\P,\Sigma)\subseteq \Gamma \subseteq L(\Sigma, \bR)$.
\end{center}
\end{lemma}

Our goal now is to give an outline of the way \rthm{capturing of hod pairs} and \rlem{generating gamma} are used to prove \rthm{msc}.

\subsection{The structure of the proof of the Mouse Set Conjecture}
  
First we outline the proof of the following general theorem. 
 
 \begin{theorem}\label{mc in min models} Suppose $(\P, \Sigma)$ is a perfect pair. Then $L(\Sigma, \bR)\models ``$for every $\R\inseg^c_{hod} \P$\footnote{See \rnot{complete layer notation}.}, ${\sf{Mouse\ Capturing}}$ holds for $\Sigma_{\R}"$.
 \end{theorem}
 \begin{proof}
 We only outline the proof as the full proof is presented in \cite[Section 6.4]{ATHM}. For simplicity we outline the proof for the least complete layer of $\P$. Let $\R\inseg^c_{hod}\P$ be the least layer of $\P$. We want to show that\\\\
(1) $L(\Sigma, \bR)\models ``{\sf{Mouse\ Capturing}}$ holds for $\Sigma_{\R}"$.\\\\
 The general case is only notationally more complex. Suppose $x, y\in \bR$ are such that $L(\Sigma, \bR)\models ``y\in OD_{\Sigma_\R}(x)"$. It follows from \rthm{mc in l(sigma, r)} that there is a $\Sigma$-mouse $\M$ over $(\P, x)$ containing $y$ such that $\M$ has an iteration strategy in $L(\Sigma, \bR)$. In fact, it follows from \rthm{mc in l(sigma, r)} that\\\\
(2) for every $\Q\in pI(\P, \Sigma)$ there is a $\Sigma_\Q$-mouse $\M$ over $(\Q, x)$ such that $y\in \M$ and $\M$ has an iteration strategy in $L(\Sigma, \bR)$.\footnote{This is because $L(\Sigma_\Q,\mathbb{R})=L(\Sigma,\mathbb{R})$ and $L(\Sigma_\Q,\mathbb{R})\vDash \sf{MC}$$(\Sigma_\Q)$.}\\\\
Let $\M_\Q$ be the least $\Sigma_\Q$-mouse over $(\Q, x)$ such that $y$ is definable over $\M_\Q$. Let $\Lambda_\Q$ be the iteration strategy of $\M_\Q$ (witnessing that $\M_\Q$ is a $\Sigma_\Q$-mouse). Let $\Gamma^*\in L(\Sigma, \bR)$ be a good pointclass such that the set
\begin{center}
 $A=\{(z, u)\in \bR^2: z$ codes some $\M_\Q$ and $u$ is an iteration according to $\Lambda_\Q\}$
 \end{center}
  is in $\utilde{\Delta}_{\Gamma^*}$. Let $F$ be as in \rthm{n*x} for $\Gamma^*$ and let $z\in \dom(F)$ be such that if $F(z)=(\N_z^*, \M_z, \d_z, \Sigma_z)$ then $(\N^*_z, \d_z, \Sigma_z)$ Suslin, co-Suslin captures $\Sigma$ and the set $A$. Let $\N=({\sf{Le}}(\emptyset, x))^{\N^*_z|\delta_z}$. It follows from \rthm{capturing of hod pairs} that\\\\
 (3) there is a $\Q\in \N$ such that $\Sigma_\Q\rest \N\in \mathcal{J}[\N]$.\\\\
 It follows from the universality of $\N$ that $\M_\Q\in \N$ (this is because $({\sf{Le}}((\Q, \Sigma_\Q))^{\N}$ is universal in $\N^*_z$ and the strategy $\Lambda_\Q$ of $\M_\Q$ is captured by $\N^*_z$ (via $A$)). It then follows that $y\in \N$. As $\N$ is an $x$-mouse, this completes the proof. 
 \end{proof}
 
Suppose now that $(\P, \Sigma)$ is a $\Gamma$-perfect pair such that $\Gamma(\P, \Sigma)\subseteq \Gamma\subseteq L(\Sigma, \bR)$. Such a pair is given to us by \rlem{generating gamma}. 
%
%

We now apply \rthm{mc in l(sigma, r)}. For each $\Q\in pI(\P, \Sigma)$ there is a $\Sigma_\Q$-mouse $\M_\Q$ over $(\Q, x)$ such that $y$ is definable over $\M_\Q$. We then again can find an $x$-mouse $\N$ such that for some $\Q\in\N\cap pI(\P, \Sigma)$, $\M_\Q\in \N$. It follows that $y\in \N$. Thus, to finish the proof of \rthm{msc}, it is enough to establish \rthm{capturing of hod pairs} and \rlem{generating gamma}. 

\subsection{Review of basic notions}
 
In this subsection we review basic notions introduced in \cite[Theorem 6.5]{ATHM} for proving a version of \rthm{capturing of hod pairs}. 

\begin{terminology}\label{objects for the proof}
We are in fact working towards the proof of \rthm{capturing of hod pairs}, and the notation and the terminology of this subsection will be used in the later subsections. 
Fix $(\P, \Sigma)$, $\Gamma_1$, $F$ and $z$ as in the statement of \rthm{capturing of hod pairs}. Let $\N=({\sf{Le}}(\emptyset))^{\N^*_z}$. $\myqedhere$
\end{terminology}
\textbf{Goal:} We are looking for $\Q\in pI(\P, \Sigma)\cap \N$ such that $\Sigma_\Q\rest \N\in \mathcal{J}[\N]$.

 We start working in $\N^*_z$. Without loss of generality we can assume that\\\\
 (1) whenever $\R\in pB(\P, \Sigma)\cap (\N^*_z|\d_z)$ there is $\S\in pI(\R, \Sigma_\R)\cap \N$ such that $\Sigma_\S\rest \N\in \mathcal{J}[\N]$.\\\\
 As in \cite{ATHM}, there are several cases. 
 \begin{enumerate}
 \item $\P$ is of successor type.
 \item $\P$ is of limit type  and $\P$ is meek. 
 \item $\P$ is non-meek but $\P$ is not of $\#$-lsa type.
 \item $(\P, \Sigma)$ is an sts hod pair. 
 \end{enumerate}
 The first two cases are just like the cases considered in \cite[Theorem 6.5]{ATHM}, we leave those to the reader. Here we analyze the remaining two cases. To start, we need to import some notions from \cite[Section 6.3]{ATHM}.
 
 \begin{definition}\label{hod pair at kappa}
 Suppose for a moment that we are working in some model of $\sf{ZFC}$. Suppose $\k$ is an inaccessible cardinal. We say that $(\Q, \Lambda)$ is a \textbf{hod pair at $\k$} if
\begin{enumerate}
\item $(\Q, \Lambda)$ is a hod pair,
\item $\Q\in {\sf{HC}}$\footnote{We will later apply this definition to $\Q$ which are not countable. The reason we make this assumption is so that we can have clause 4 below. It follows that the current definition makes sense in a variety of situations, and in particular when clause 4 holds after collapsing $\Q$ to be countable.}
\item $\Lambda$ is a $(\kappa, \kappa)$-iteration strategy,
\item ${\sf{Code}}(\Lambda)$ is a $\k$-universally Baire set of reals.
\end{enumerate}
$\myqedhere$
\end{definition}
Suppose $(\Q, \Lambda)$ is a hod pair at $\k$.  Then we let
\begin{center}
$Lp^{\Lambda, \kappa}(a)=\bigcup \{\M: \M$ is a sound $\Lambda$-mouse over $a$ such that $\rho_\omega(\M)={\sf{ord}}(a)$ and $\M\trianglelefteq ({\sf{Le}}((\Q, \Lambda), a)^{V_\k}\}$.
\end{center}
As is customary, we let ${\sf{Lp}}^{\Lambda, \kappa}_\a(a)$ be the $\a$th iterate of ${\sf{Lp}}^{\Lambda, \kappa}(a)$. Below $\S^*(\R)$ is the $*$-transform of $\S$ into a hybrid mouse over $\R$, it is defined when $\R$ is a cutpoint of $\S$ (cf. \cite{Selfiter}).

\begin{definition}[Fullness preservation in models of $\sf{ZFC}$]\label{fullness press in models of zfc} Suppose now that $(\P, \Sigma)$ is a hod pair at $\kappa$. We then say $\Sigma$ is \textit{$\kappa$-fullness preserving} if the following holds for all $(\T, \Q)\in I(\P, \Sigma)\cap V_\k$. 
\begin{enumerate}
\item For all meek\footnote{See \rdef{pre-hod-like}.} layers $\R$ of $\Q$ such that $\R$ is of successor type\footnote{See \rdef{types of lsa small premice}.}, letting $\S=\R^-$\footnote{This is the longest proper layer of $\R$. See \rnot{l p}.}, for all $\eta\in (\ord(\S), \ord(\R))$ if $\eta$ is a cutpoint cardinal of $\R$ then
\begin{center}
$(\R|(\eta^+)^\R)^*={\sf{Lp}}^{\Sigma_{\S, \T}, \k}(\R|\d)$.
\end{center} 
\item  For all meek\footnote{See \rdef{pre-hod-like}.} layers $\R$ of $\Q$ such that $\R$ is of limit type, 
\begin{center}
$\R={\sf{Lp}}^{\Sigma_{\R|\d^\R, \T}, \k}(\R|\d^\R)$.
\end{center}
\item If $\P$ is of $\#$-lsa type then ${\sf{Lp}}^{\Sigma^{stc}_{\Q, \T}, \k}(\Q)\models ``\d^\Q$ is a Woodin cardinal"\footnote{Here, if $\Sigma$ is a short tree strategy then $\Sigma^{stc}=\Sigma$.}. 
\end{enumerate}
$\myqedhere$
\end{definition}

We continuing our work inside some model of $\sf{ZFC}$.  
\begin{definition}[Universal tail]\label{universal tail}
Suppose $(\Q, \Lambda)$ is a hod pair at $\k$ such that $\Lambda$ has branch condensation and is $\k$-fullness preserving. 
Suppose $\l<\k$ is an inaccessible cardinal. Then we say $(\Q^*, \Lambda^*)$ is a \textbf{$\l$-universal tail} of $(\Q, \Lambda)$ if there is a (possibly generalized) stack
\begin{center}
$\T=(\M_\b, \T_\b, E_\b: \b< \l)$  
\end{center} 
on $\Q$ according to $\Lambda$ with last model $\Q^*$ such that $\lh(\T)=\l$ and 
 for any $(\S, \R)\in I(\Q, \Lambda)\cap V_\l$ there is a stack $\U$ on $\R$ according to $\Lambda_{\R, \S}$ such that for some $\a<\l$, $\M_\a$ is the last model of $\U$. 

If $\T$ is as above then we say $\T$ is a \textbf{$\l$-universal stack }on $\Q$ according to $\Lambda$. $\myqedhere$
\end{definition}
 
We now resume the proof of \rthm{capturing of hod pairs}, and continue with the objects introduced in \rter{objects for the proof}.
and start working in $\N^*_z$. Observe that because of our assumption on $(\P, \Sigma)$, whenever $\Q, \R \in pI(\P, \Sigma)$, $(\Q, \Sigma_\Q)$ and $(\R, \Sigma_\R)$ have a common tail in $\N^*_z|\d_z$. In fact more is true. Suppose $\kappa$ is a strong cardinal of $\N^*_z$. Then it follows from \rcor{comparison holds} that if   $\Q, \R \in pI(\P, \Sigma)\cap \N^*_z|\k$ then $(\Q, \Sigma_\Q)$ and $(\R, \Sigma_\R)$ have a common tail in $\N^*_z|\k$. This means that whenever $\kappa<\d_z$ is a cardinal of $\N^*_z$ and $\Q\in (pI(\P, \Sigma)\cup pB(\P, \Sigma))\cap \N^*_z|\k$, we can form the direct limit of all $\Sigma_\Q$ iterates of $\Q$ that are in $\N^*_z|\k$. Let $\R^{\Q, \Sigma_\Q}_\kappa$ be this direct limit. The next lemma shows that the universal tails are unique. It appeared as \cite[Lemma 6.8]{ATHM}.

In what follows, we will often abuse the terminology introduced in \rdef{universal tail}. Usually when, working inside $\N^*_z$, we talk about $\k$-universal tail of some $(\Q, \Lambda)$ with $\Q\in \N^*_z|\k$ then we mean that we are working in $\N^*_z[h]$ where $h\subseteq Coll(\omega, \Q)$ is $\N^*_z$-generic. Then take our current $\k$ be $\l$ of \rdef{universal tail} and $\k$ of \rdef{universal tail} be $\d_z$. 

\begin{lemma}[Uniqueness of universal tails]\label{uniqueness of universal tails} Suppose $\Q\in pI(\P, \Sigma)\cap \N^*_z|\d_z$. Then for each $\S\inseg^c_{hod}\Q$\footnote{See \rdef{complete layer notation}.} and $\N$-strong $\k<\d_z$ such that $\S\in\N^*_z|\k$, there is a unique $\k$-universal tail of $(\S, \Sigma_{\S})$. In fact, letting $\R=\R^{\S, \Sigma_{\S}}_{\kappa}$, $(\R, \Sigma_{\R})$ is the unique $\k$-universal tail of $(\S, \Sigma_{\S})$
\end{lemma}

\begin{definition}
Suppose $\Q\in (pI(\P, \Sigma)\cup pB(\P, \Sigma))\cap \N^*_z|\d_z$ and $\k$ is an $\N$-strong cardinal such that $\Q\in \N^*_z|\k$. 
Then we say $\N$ \textbf{captures} a tail of $(\Q, \Sigma_\Q)$ below $\k$ if there is a hod pair $(\R, \Lambda)\in \N$ such that $\Lambda$ is a $(\k, \k)$-iteration strategy and there is a term relation $\tau\in \N^{Coll(\omega, <\k)}$ such that whenever $g\subseteq Coll(\omega, \card{\R}^+)$ is $\N$-generic, \begin{enumerate}
\item $\N[g]\models ``(\R, \tau_g)$ is a hod pair at $\k$ such that $\tau_g$ is $\k$-fullness preserving" and $\tau_g\rest \N=\Lambda$,
\item for some $\l<\k$, $\R=\R^{\Q, \Lambda}_\l$ and letting $T, U\in \N[g]$ witness that $\tau_g$ is $\k$-uB, whenever $h\subseteq Coll(\omega, <\k)$ is $\N[g]$-generic, $(p[T])^{\N[g][h]}={\sf{Code}}(\Sigma_\R)\cap \N[g][h]$.
\end{enumerate}
We say $\N$ captures $B(\Q, \Sigma_\Q)$ below $\k$ if whenever $\R\in pB(\Q, \Sigma_\Q)\cap \N^*_z|\k$, $\N$ captures $(\R, \Sigma_{\R})$ below $\k$. $\myqedhere$
\end{definition}

Towards a contradiction, we assume that $\N$ does not capture a tail of $(\P, \Sigma)$ and that either
\begin{enumerate}
\item  $\P$ is non-meek but $\P$ is not of $\#$-lsa type,
\item $(\P, \Sigma)$ is an sts hod pair. 
\end{enumerate}

\begin{notation}\label{basic objects used in capturing}
For each $\Q\in pB(\P, \Sigma)$, we let $\l_{\Q}$ be the least $\N$-strong cardinal $\nu$ such that $\N$ captures the $\nu$-universal tail of $(\Q, \Sigma_{\Q})$. We let $(\R^{\Q, \Sigma}, \Phi^{\Q, \Sigma})$ be the $\l_{\Q}$-universal tail of $(\Q, \Sigma_{\Q})$. For each inaccessible cardinal $\nu$ such that $ \Q\in \N|\nu$, we let $(\R^{\Q, \Sigma}_\nu, \Phi^{\Q, \Sigma}_\nu)$ be the $\nu$-universal tail of $(\Q, \Sigma_{\Q})$. If $\l\geq \l_\Q$ then $\pi^{\Sigma_\Q}_{\Q, \R_\l^{\P, \Sigma}}$ is the iteration map $\pi^{\Sigma_\Q}_{\Q, \R^{\Q, \Sigma}}$. $\myqedhere$
\end{notation}
 
\subsection{The ideas behind the proof}

The notation and the terminology introduced in this subsection will be used in the next few subsections. We are continuing with the set up of \rter{objects for the proof}. 

\begin{notation} Suppose now that $\kappa_0$ is an $\N$-strong cardinal that reflects the set of $\N$-strong cardinals. Let 
\begin{center}
$\mathcal{E}_0=\{ E\in \vec{E}^\N: \N\models ``\nu(E)$ is inaccessible" and for all $\eta\in (\k_0, \nu(E))$, $\N\models ``\eta$ is a strong cardinal" if and only if $Ult(\N, E)\models ``\eta$ is a strong cardinal"$\}$.  
\end{center}
$\myqedhere$
\end{notation}

\begin{notation}\label{definition of f in n}
Working in $\N$, let 
\begin{center}
$\mathcal{F}=\{ (\Q, \Lambda) : \Q\in \N|\d_z\wedge \mathcal{J}[\N]\models ``(\Q, \Lambda)$ is a hod pair at $\d_z$ and $\Lambda$ has branch condensation and is $\d_z$-fullness preserving$"\}$.
\end{center}
We have that $\mathcal{F}$ is a directed system. Let for $\l\leq \d_z$, 
\begin{center}
$\mathcal{F}\rest \l=\{(\Q, \Lambda)\in \mathcal{F}: \Q\in \N|\l\}$.
\end{center}
 We let $\R^*$ be the direct limit of $\F\rest \k_0$ under the iteration maps. $\myqedhere$
\end{notation}
Recall that we fixed an $\N$-strong cardinal $\k_0$ that reflects the set of strong cardinals of $\N$. The equality below is computed inside $\N^*_z$. 
 \begin{definition}\label{definition of r}  Let $\R_0=(\R^{\P, \Sigma}_{\k_0})^b$.
$\myqedhere$
\end{definition}

 The next lemma summarizes what was proved in \cite{ATHM}. 
\begin{lemma}\label{summary of sec 6.3} The following holds.
\begin{enumerate}
\item 
Suppose $\Q\in pB(\P, \Sigma)\cap \N^*_z|\k_0$. Then $\R^{\Q, \Sigma_\Q}\in \N|\k_0$.
\item Suppose $\Q\in pB(\P, \Sigma)$, $\l>\k_0$ is a strong cardinal of $\N$ such that $\Q\in \N|\l$, and $E\in \mathcal{E}_0$ is an extender with critical point $\k_0$ such that $\nu(E)>(\l^+)^{\N^*_z}$. Then $\Phi^{\Q, \Sigma}\rest Ult(\N, E)\in Ult(\N, E)$.
\item Let $\R^*$ be as in \rnot{definition of f in n}. Then either $\R_0\trianglelefteq_{hod} \R^*$ or $\R_0|\d^{\R_0}=\R^*$. Moreover, $\R_0\in \N$ and $\Sigma_{\R_0}\rest \N\in \mathcal{J}[\N]$.
\end{enumerate}
\end{lemma}
Clause 1 is just \cite[Lemma 6.11]{ATHM}, clause 2 is \cite[Lemma 6.12]{ATHM} and clause 3 is \cite[Lemma 6.13]{ATHM}. 

Below we will develop a technology for recovering the full iterate of $\P$. Let $\R_0^+=\R^{\P, \Sigma}_{\k_0}$ be the iterate of $\P$ extending $\R_0$ and let $i: \P\rightarrow \R_0^+$ be the iteration embedding. We will recover an iterate of $\R_0^+$ inside $\N$ as an output of a backgrounded construction that is done over $\R_0$. Such constructions are called mixed hod pair constructions. The details of this construction appear in \rsec{sec: mixed hod pair constructions}. 

There are two kinds of extenders that we will use in this construction. The extenders with critical point $>\d^{\R_0}$ will have traditional background certificates. We will use the total extenders on the sequence of $\N$ to certify such extenders. The extenders with critical point $\d^{\R_0}$ will come from a different source. The following key lemma illustrates the idea. 

\begin{lemma}\label{prescription for extenders} Let $\d=\d^{\R_0}$. Suppose $\S \in pI(\R_0^+, \Sigma_{\R_0^+})$ is a normal iterate of $\R_0^+$ that is obtained by iterating entirely above $\d^{\R_0}$. Suppose that $\a\in \dom(\vec{E}^\S)$ is such that letting $E=_{def} \vec{E}^\S(\a)$, $\cp(E)=\d$, $\S|\a \in \N$ and $\Sigma_{\S|\a}\rest \N\in \mathcal{J}[\N]$. Then $E\in \N$. Moreover, $(a, A)\in E$ if and only if $a\in \nu_{E}^{<\omega}$, $A\in [\d]^{\card{a}}$ and whenever $F\in \mathcal{E}_0$ is such that $\cp(F)=\k_0$ and 
\begin{center}
$\N\models ``$there is a strong cardinal $\l$ in the interval $(\k_0, \nu_F)$ such that $\S\in \N|\l"$,
\end{center}
$\pi^{\Sigma_{\S|\a}}_{\S|\a, \pi_F(\R_0)}(a)\in \pi_F(A)$\footnote{The embedding $\pi^{\Sigma_{\S|\a}}_{\S|\a, \pi_F(\R_0)}$ is just $\pi^{\Sigma_{\S|\a}}_{\S|\a, \R_{\pi_F(\k_0)}^{\S|\a, \Sigma_{\S_\a}}}$. We will often abuse our notation this way.}.
\end{lemma}
\begin{proof} 
Set $\M^+=Ult(\R^+_0, E)$ and $\M=Ult(\R_0, E)$. Let $F^*$ be the resurrection of $F$ and let $\sigma: Ult(\N, F)\rightarrow \pi_{F^*}(\N)$ be the canonical factor map. We have that $\sigma\rest \nu_F=id$. Thus, $\pi_{F^*}\rest \N=\sigma\circ \pi_F$. It follows that $\pi_{F^*}\rest \R^+$ is the iteration embedding implying\\\\
(1) $\pi_{F^*}\rest \R_0^+=\pi^{\Sigma_{\M^+}}_{\M^+, \pi_{F^*}(\R_0^+)}\circ \pi^{\R_0^+}_{E}$. \\\\
We now have that
\begin{eqnarray*}
(a, A)\in E & \iff & a \in \pi^{\R_0^+}_{E}(A) \\
&\iff & \pi^{\Sigma_{\M^+}}_{\M^+, \pi_{F^*}(\R_0^+)}(a) \in  \pi^{\Sigma_{\M^+}}_{\M^+, \pi_{F^*}(\R_0^+)}(\pi^{\R_0^+}_{E}(A))\\
&\iff&  \sigma(\pi^{\Sigma_{\M}}_{\M, \pi_F(\R_0)}(a)) \in \pi_{F^*}(A)  \\
&\iff& \sigma(\pi^{\Sigma_{\M}}_{\M, \pi_F(\R_0)}(a)) \in \sigma(\pi_F(A))\\
&\iff& \pi^{\Sigma_{\M}}_{\M, \pi_F(\R_0)}(a) \in \pi_F(A)\\
\end{eqnarray*}
Therefore, 
\begin{center}
$(a, A)\in E\iff \pi^{\Sigma_{\M}}_{\M, \pi_F(\R_0)}(a) \in \pi_F(A)$.
\end{center}
By our assumption, the right hand side of the equivalence can be computed in $\N$. Hence $E\in \N$.
\end{proof}

Thus, the extenders with critical point $\d^{\R_0}$ that we will use in our mixed hod pair construction have the following property. If $\Q$ is the current level of the construction and $\Lambda$ is its strategy then let $E$ be the set of pairs $(a, A)$ such that $a\in (\d^{\R_0})^{<\omega}$ and for every $F\in \mathcal{E}_0$ such that $\cp(F)=\k_0$ and 
\begin{center}
$\N\models ``$there is a strong cardinal $\l$ in the interval $(\k_0, \nu_F)$ such that $\Q\in \N|\l"$,
\end{center}
$\pi^{\Lambda}_{\Q, \pi_F(\R_0)}(a)\in \pi_F(A)$.

There is one problem with this approach. We need to know the strategy $\Lambda$ of $\Q$ before we can find the relevant extender. To resolve this problem, we will first define the strategy $\Lambda$. Essentially $\Lambda$ will pick branches that, for some $\eta$, are $\pi_E$-realizable for all $E\in \mathcal{E}_0$ such that $\nu_E>\eta$. We will call such strategies $\mathcal{E}_0$-certified.

In the sequel, we will first introduce the $\mathcal{E}_0$-certified strategies. Then we will prove basic facts about them. Then we will introduce the mixed hod pair constructions and show that some model appearing on this construction is an iterate of $\R^+_0$ via an iteration that is entirely above $\d^{\R_0}$. 

\subsection{A condensing set in $\N$.}

The following is the main lemma of this section. Our goal is still to prove \rthm{capturing of hod pairs} and our set up is as in \rter{objects for the proof} and \rdef{definition of r}. 
\begin{lemma}\label{main lemma on condensing sets in n}  Suppose $\eta>\k_0$ is such that $\N\models ``\eta$ is a strong cardinal that reflects the set of strong cardinals". Set $\S^+=\R_{\eta}^{\P, \Sigma}$, $i^+=\pi^{\Sigma_{\R_0^+}}_{\R_0^+, \S^+}$, $\S=(\S^+)^b$ and $i=i^+\rest \R_0$. Then $i\in \N$ and $\N\models \card{\S}<(\eta^+)^\N$.
Moreover, the following holds:
\begin{enumerate}
\item If $E\in \vec{E}^{\N^*_z}$ is such that $\cp(E)=\eta$ and $E$ is total over $\N^*_z$ then $\pi_E\rest \S$ is a strongly condensing set in $Ult(\N^*_z, E)[g]$ where $g\subseteq Coll(\omega, \pi_E(\eta))$ is any $Ult(\N^*_z, E)$-generic\footnote{See \rdef{def:condensing_set}. This definition uses a formula $\phi$ and a parameter $A$. In the current case, $(\phi, A)$ is chosen in a way that the resulting directed system is $\mathcal{F}\rest \eta$ where $\mathcal{F}$ is as in \rnot{definition of f in n}.}.
\item  If $E\in \vec{E}^{\N}$ is such that $\cp(E)=\eta$ and $\nu(E)$ is an inaccessible cardinal of $\N$ then $\pi_E\rest \S$ is a strongly condensing set in $Ult(\N, E)[g]$ where $g\subseteq Coll(\omega, \pi_E(\eta))$ is any $Ult(\N^*_z, E)$-generic.
\end{enumerate}
\end{lemma}
\begin{proof} The fact that $\S\in \N$ follows from \cite{ATHM} and it is the same argument that shows that $\R_0\in \N$ (see \rlem{summary of sec 6.3}). Suppose now that $i\in \N$. It then follows that $\N\models \card{\S}<\eta^+$. Hence, \rthm{existence of condensing sets in n} implies both clause 1 and 2 above. Thus, it is enough to show that $i\in \N$.

Let $F\in \vec{E}^\N$ be any extender such that $\cp(F)=\k_0$ and $Ult(\N, F)\models ``\eta$ is a strong cardinal". Let $F^*$ be the background certificate of $F$ and let $k: Ult(\N, F)\rightarrow \pi_{F^*}(\N)$ be the canonical factor map. We now have that $\pi_{F^*}\rest \R_0^+=\pi^{\Sigma_{\R_0^+}}_{\R^+, \pi_{F^*}(\R_0^+)}$. We thus have that\\\\
(1) $\pi_{F^*}\rest \R_0^+=\pi^{\Sigma_{\S^+}}_{\S^+, \pi_{F^*}(\R_0^+)}\circ \pi^{\Sigma_{\R_0^+}}_{\R_0^+, \S^+}$.\\\\
Let $m=\pi^{\Sigma_{\S^+}}_{\S^+, \pi_{F^*}(\R_0^+)}\rest \S|\d^\S$. We have that\\\\
(2) $m=\pi^{\Sigma_{\S|\d^\S}}_{\S|\d^\S, \pi_{F^*}(\R_0)|\xi}$ where $\xi=\sup(m[\d^\S])$.\\\\
Because $\Sigma_{\S|\d^\S}\rest \N \in \N$, we have that $k(\Sigma_{\S|\d^\S}\rest \N)=\Sigma_{\S|\d^\S}\rest \pi_{F^*}(\N)$ and therefore, $m\in \pi_{F^*}(\N)$ and $m\in \rge(k)$. Let $n=k^{-1}(m)$. Thus,\\\\
(3) $n=\pi^{\Sigma_{\S|\d^\S}}_{\S|\d^\S, \pi_{F}(\R_0)|k^{-1}(\xi)}$.\\\\
Notice now that for $x\in \R_0$,\\\\
(4) $\pi_{F^*}(x)=\pi^{\Sigma_{\S^+}}_{\S^+, \pi_{F^*}(\R_0^+)}(i(x))$\\\\
implying that\\\\
(5) $\S$ is the transitive collapse of $\{\pi_{F^*}(f)(m(a)): f\in \R_0 \wedge a\in (\d^\S)^{<\omega}\}$ and $\pi^{\Sigma_{\S^+}}_{\S^+, \pi_{F^*}(\R_0^+)}\rest \S$ is the inverse of the transitive collapse.\\\\
(5) now implies that\\\\
(6) $\S$ is the transitive collapse of $\{\pi_{F}(f)(n(a)): f\in \R_0 \wedge a\in (\d^\S)^{<\omega}\}$.\\\\
Since $\{\pi_{F}(f)(n(a)): f\in \R_0 \wedge a\in (\d^\S)^{<\omega}\}\in \N$, we have that if $\sigma:\S\rightarrow \pi_F(\R_0)$ is the inverse of the transitive collapse then $\sigma\in \N$. Moreover,\\\\
(7) $\pi_{F^*}\rest \R_0=k\circ \sigma\circ i$ and $k\circ \sigma=\pi^{\Sigma_{\S^+}}_{\S^+, \pi_{F^*}(\R_0^+)}\rest \S$. \\\\
It now follows that\\\\
(8) $i(x)=\sigma^{-1}(\pi_F(x))$. \\\\
Since both $\sigma$ and $\pi_F$ are in $\N$, we get that $i\in \N$.
\end{proof}

\begin{notation}\label{notation for kappa} Suppose now that $\kappa$ is an $\N$-strong cardinal that reflects the set of $\N$-strong cardinals such that $\k>\k_0$. Let 
\begin{center}
$\mathcal{E}=\{ E\in \vec{E}^\N: \N\models ``\nu(E)$ is inaccessible" and for all $\eta\in (\k, \nu(E))$, $\N\models ``\eta$ is a strong cardinal" if and only if $Ult(\N, E)\models ``\eta$ is a strong cardinal"$\}$.  
\end{center}
Set $\R^+=\R^{\P, \Sigma}_\k$ and let $\R=(\R^+)^b$. It follows from \rlem{summary of sec 6.3} that $\R\in \N$. Let $\Phi^+=(\Phi^{\P, \Sigma}_\k)_{\R|\d^\R}$ and $\Phi=\Phi^+_\R$\footnote{See \rnot{basic objects used in capturing}.}. $\myqedhere$
\end{notation}

Notice that $\Phi\rest \N\in L[\N]$. We will abuse our notation and write $\Phi$ for both $\Phi\rest \N^*_z$ and $\Phi\rest \N$, but we encourage the reader to keep this subtle difference between the three versions of $\Phi$ in mind. 
\begin{definition}\label{e-realizable system} Working in $\N$, we say $(\sigma, \Q)$ is \textbf{$\mathcal{E}$-realizable} if 
\begin{itemize}
\item $\sigma:\R\rightarrow \Q$ is an elementary embedding,
\item for some $\N$-strong cardinal $\xi\in (\k, \d_z)$, $\Q\in \N|\xi$ and for all $E\in \mathcal{E}$ such that $\xi<\nu(E)$ and for all $\N$-generic $g\subseteq Coll(\omega, \Q)$, there is $j:\Q\rightarrow \pi_E(\R)$ such that $j\in Ult(\N, E)[g]$ and $\pi_E\rest \R=\sigma\circ j$.\footnote{Notice that $\pi_E\rest \R\in Ult(\N, E)$, see \rlem{main lemma on condensing sets in n}.}
\end{itemize}
We say that $j$ is $(\pi_E, \sigma)$-realizable. Continuing our work in $\N$, let $\mathcal{F}'_{\mathcal{E}}$ be the set of $\pi_E$-realizable pairs $(\sigma, \Q)$. Given $(\sigma, \Q)\in \mathcal{F}'_{\mathcal{E}}$, let $\xi(\sigma, \Q)<\d_z$ witness that clause 2 above holds for $(\sigma, \Q)$. Given $E\in \mathcal{E}$ such that $\xi(\sigma, \Q)<\nu(E)$, letting $j:\Q\rightarrow \pi_E(\R)$ be any $(\pi_E, \sigma)$-realizable embedding, set $\Psi_{\sigma, \Q, E, j}=(j$-pullback of $\pi_E(\Phi))$\footnote{$\Psi_{\sigma, \Q, E, j}$ is defined in $Ult(\N, E)$.}. $\myqedhere$
\end{definition}
The following is an easy consequence of \rlem{lem:uniqueness of strategies}.
\begin{lemma}\label{definability of the system} Suppose $(\sigma, \Q, E)$ are as in \rdef{e-realizable system} and that $\sigma\rest \R|\d^\R$ is an iteration embedding according to $\Phi_{\R|\d^\R}$. Let $j_0:\Q\rightarrow \pi_E(\R)$ and $j_1:\Q\rightarrow \pi_E(\R)$ be two $(\pi_E, \sigma)$-realizable embeddings. Then $\Psi_{\sigma, \Q, E, j_0}=\Psi_{\sigma, \Q, E, j_1}$. 
\end{lemma}
Given $(\sigma, \Q, E)$ as in \rlem{definability of the system}, we let $\Psi_{\sigma, \Q, E}$ be the common value of all $\Psi_{\sigma, \Q, E, j}$ where $j$ is any $(\pi_E, \sigma)$-realizable embedding. The next definition integrates $\Psi_{\sigma, \Q, E}$ with respect to $E$.

\begin{definition}\label{neatly e-realizable system} Working in $\N$, we say $(\sigma, \Q)$ is \textbf{neatly $\mathcal{E}$-realizable} if $(\sigma, \Q)$ is $\mathcal{E}$-realizable and for all $E_0, E_1\in \mathcal{E}$ with $\nu(E_0)\leq \nu(E_1)$, 
\begin{center}$\Psi_{\sigma, \Q, E_0}\rest \N|\nu(E_0)=\Psi_{\sigma, \Q, E_1}\rest \N|\nu(E_0)$.\end{center}
Let $\mathcal{F}_{\mathcal{E}}$ be the set of neatly $\mathcal{E}$-realizable pairs, and for $(\sigma, \Q)\in \mathcal{F}_{\mathcal{E}}$, let
\begin{center}
 $\Psi_{\sigma, \Q}=\cup\{\Psi_{\sigma, \Q, E}: E\in \mathcal{E} \wedge \xi(\sigma, \Q)<\nu(E)\}$\footnote{$\Psi_{\sigma, \Q}$ is defined in $\mathcal{J}[\N]$ and $\Psi_{\sigma, \Q}\rest \N$ is total.}. 
 \end{center}
 $\myqedhere$
\end{definition}
The following is the key lemma of this section.

\begin{lemma}\label{iterates are neatly e-realizable} Suppose $\S$ is a $\Phi^+$-iterate of $\R^+$ via $\T$ such that $\pi^{\T, b}$ is defined and $\S^b\in \N$. Then $\pi^{\T, b}\in \N$, $(\pi^{\T, b}, \S^b)$ is neatly $\mathcal{E}$-realizable and 
\begin{center}
$\Psi_{\pi^{\T, b}, \S^b}=\Phi^+_{\S^b}\rest \N$.
\end{center}
\end{lemma}
\begin{proof} The proof of $\pi^{\T, b}\in \N$ is exactly the proof of \rlem{main lemma on condensing sets in n}. The proof of the fact that $(\pi^{\T, b}, \S^b)$ is neatly $\mathcal{E}$-realizable is via a simple absoluteness argument. Let $E\in \mathcal{E}$ be such that $\S^b\in \N|\nu(E)$ and let $E^*$ be the background certificate of $E$. Let $k:Ult(\N, E)\rightarrow \pi_{E^*}(\N)$ be the canonical factor map. We have that $\cp(k)\geq \nu(E)$. Set $\sigma=\pi^{\T, b}$. Notice that\\\\
(1) in $\pi_{E^*}(\N)$, it is forced by $Coll(\omega, \S^b)$ that there is a $(\pi_{E^*}, \sigma)$-realizable $j:\S^b\rightarrow \pi_{E^*}(\R)$, and\\
(2) if $g\subseteq Coll(\omega, \S^b)$ is $\pi_{E^*}(\N)$-generic and $j:\S^b\rightarrow \pi_{E^*}(\R)$ is any $(\pi_{E^*}, \sigma)$-realizable embedding, then the $j$-pullback of $\pi_{E^*}(\Phi)$ is $\Phi^+_{\S^b}$. \\\\
It follows that\\\\
(3) in $\N$, it is forced by $Coll(\omega, \S^b)$ that there is a $(\pi_{E^*}, \sigma)$-realizable $j:\S^b\rightarrow \pi_{E^*}(\R)$, and\\
(4) if $g\subseteq Coll(\omega, \S^b)$ is $\N$-generic and $j:\S^b\rightarrow \pi_{E^*}(\R)$ is any $(\pi_{E^*}, \sigma)$-realizable embedding, then the $j$-pullback of $\pi_{E}(\Phi)$ is independent of $j$. \\\\
Let $\Pi$ in $Ult(\N, E)$ be the strategy of $\S^b$ such that it is forced by $Coll(\omega, \S^b)$, that for some $(\pi_E, \sigma)$-realizable $j$, $\Pi$ is the $j$-pullback of $\pi_E(\Phi)$. Let $\tau:\S^b\rightarrow \pi_E(\R)$ be defined by setting $\tau(x)=\pi_E(f)(\pi^{\Pi}_{\S|\d^{\S^b}, \pi_E(\R)}(a))$ where $x=\sigma(f)(a)$, $f\in \R$ and $a\in (\d^{\S^b})^{<\omega}$. It follows from clause 2 of \rlem{main lemma on condensing sets in n} that $\tau$ is a $(\pi_E, \sigma)$-realization and $\tau\in Ult(\N, E)$. It then follows from (2) that $k(\Pi)=\Phi^+_{\S^b}$ and therefore, $\Pi=\Phi^+_{\S^b}\rest Ult(\N, E)$.
\end{proof}

\begin{definition}\label{canonical e-realization} Suppose $(\sigma, \Q)\in \mathcal{F}_{\mathcal{E}}$ and $E\in \mathcal{E}$ is such that $\xi(\sigma, \Q)<\nu(E)$. We say that $\tau=\tau^E_{\sigma, \Q}$ is the \textbf{canonical $E$-realization} of $(\sigma, \Q)$ if $\tau:\Q\rightarrow\pi_E(\R)$ and $\tau(x)=\pi_E(f)(\pi^{\Psi_{\sigma, \Q, E}}_{\Q|\d^\Q, \R(\sigma, \Q)}(a))$ where $\R(\sigma, \Q)\inseg_{hod}\R$ is the $\Psi_{\sigma, \Q, E}$-iterate of $\Q|\d^\Q$, $f\in \R$, $a\in (\d^\Q)^{<\omega}$ and $x=\sigma(f)(a)$. $\myqedhere$
\end{definition}

\subsection{$\mathcal{E}$-certified iteration strategies}\label{sec: e-certified strategies}

Our goal is still to prove \rthm{capturing of hod pairs} and our set up is as in \rter{objects for the proof}, \rdef{definition of r} and \rnot{notation for kappa}. The following is a modification of \cite[Definition 6.14]{ATHM}.

\begin{definition}[$\pi_E$-realizable iterations]\label{semi j-realizable iterations}  Suppose 
\begin{enumerate}
\item $\V\in \N$ is a hod premouse extending $\R$ such that $\R=\V^b$,
\item $\T\in \N$ is either a stack on $\V$ or an st-stack\footnote{See \rdef{st-stack}.} on $\V$\footnote{If $\T$ is an st-stack then $\V$ must be of $\#$-lsa type.},
\item $E\in \mathcal{E}$.
\end{enumerate}
Recall \rdef{proper stack}, \rrem{proper stacks convention} and \rnot{notation for iteration trees}. Suppose that 
\begin{center}
$\T=((\M_\a)_{\a<\eta}, (E_\a)_{\a<\eta-1}, D, R, (\beta_\a, m_\a)_{\a\in R}, T)$
\end{center} is a stack. Set $R^b=\{\a\in R: \pi^{\T, b}_{0, \a}$ is defined$\}$. We say $\T$ is \textbf{$\pi_E$-realizable} if 
 the following holds:
\begin{enumerate}
\item $\N\models ``\l$ is a strong cardinal".
\item $\T\in \N|\lh(E)$.
\item For all $\a\in R^b$, $(\pi^{\T_{\leq \a}, b}, \M_\a^b)\in \mathcal{F}_{\mathcal{E}}$\footnote{See \rdef{neatly e-realizable system}.}.
\item For all $\a<\b$ such that $\a, \b\in R^b$, setting $\tau_\a=\tau_{\sigma, \Q}^E$, $\tau_{\a}=\tau_\b \circ \pi^{\T, b}_{\a, \b}$\footnote{See \rsec{sec:pitb}.}. 
\item For all $\a\in R^b$, letting $\Psi_\a=\Psi_{\sigma_\a, \M_\a^b}$,
\begin{enumerate}
\item if $\a\not =\max(R^b)$ and ${\sf{nc}}^\T_\a$ is based on $\M_\a^b|\d^{\M_\a^b}$ then ${\sf{nc}}^\T_\a$ is according to $\Psi_\a$, 
\item if $\a=\max(R^b)$ and $\U=\downarrow(\T_{\geq \a}, \M_\a^b)$\footnote{See \rdef{restriction of a stack}.} then 
\begin{enumerate}
\item if $\U$ is based on $\M_\a^b$ and is above $\d^{\M_\a^b}$ then it is according to the unique strategy $\Pi$ of $\M_\a^b$ witnessing that $\M_\a^b$ is a $\Psi_\a$-mouse over $\M_\a^b|\d^{\M_\a^b}$, and
\item if $\U$ is based on $\M_\a^b|\d^{\M_\a^b}$ then $\U$ is according to $\Psi_\a$.
\end{enumerate}
\end{enumerate}
\end{enumerate} 
We say that $(\sigma_\a : \a\in R^b)$ are the \textbf{$\pi_E$-realizable embeddings} of $\T$ and $(\Psi_\a: \a\in R^b)$ are the \textbf{$\pi_E$-realizable strategies} of $\T$. We say $\T$ is \textbf{$\mathcal{E}$-realizable} if for some $\eta$, $\T$ is $\pi_E$-realizable for every $E\in \mathcal{E}$ with the property that $\lh(E)>\eta$. 

The definition of the above concepts for st-stacks is very similar. The embeddings $\sigma_\a$ are once again defined for $\a\in R^b$ which once again consists of those $\a<\lh(\T)$ with the property that $\pi^{\T, b}_{0, \a}$ is defined. We leave the details to the reader.$\myqedhere$
\end{definition}

We now introduce a kind of backgrounded constructions reminiscent to the backgrounded construction introduced in  \rdef{fully backgrounded sts construction}. We will use them to find the $\Q$-structures of various iterations.  

\begin{definition}[$\mathcal{E}$-realizable backgrounded constructions]\label{e-certified background constructions} 

Suppose $\V, \T,$\footnote{We will omit $\T$ from most superscripts. Thus, $\M_\a$ is just $\M_\a^\T$.}  $\tau, \Q, \eta$ are such that
\begin{enumerate}
\item $\V\in \N$ is a hod premouse extending $\R$ such that $\V^b=\R$,
\item $\T$ is a $\mathcal{E}$-realizable stack or an st-stack on $\V$ such that $\pi^{\T, b}$ exists,
\item $\tau\in (R^b)^\T$\footnote{$R^b$ was introduced in \rdef{semi j-realizable iterations}.} and $\U=_{def}{\sf{nc}}^\T_\tau$ is based on $\M_\tau$ and is above $\d^{\M_\tau^b}$,
\item if $\U$ is of limit length then $\Q=\m(\U)$ and otherwise for some $\xi<\lh(\U)$, $\Q=\M_\xi^\U$,
\item $(\Q|\eta)^{\#}\models ``\eta$ is a Woodin cardinal".
\end{enumerate}
Then 
\begin{center}
${\sf{Le}}^{\mathcal{E}, c}((\Q|\eta)^\#)(\X_\gg , \Y_\gg, F^+_\gg, F_\gg, b_\gg: \gg\leq \delta_z)$
\end{center}
is the \textbf{fully backgrounded $\mathcal{E}$-realizable construction over $(\Q|\eta)^\#$ done in $\N$} if the following is true. 
\begin{enumerate}
\item $\X_0=\mathcal{J}_\omega((\Q|\eta)^\#)$, and for all $\xi<\delta_z$, $\X_\xi$ and $\Y_\xi$ are sts premice such that if $\W$ is a stack indexed either in $\X_\xi$ or $\Y_\xi$ then $\T_{\leq \Q}\ ^\frown \W$ is $\mathcal{E}$-realized.
\item If for some $\xi\leq \delta_z$, $\Y_{\xi}$ is defined but is not a reliable sts premouse over $(\Q|\eta)^\#$  then all other objects with index $\geq \xi$ are undefined.
\item Suppose for some $\xi<\delta_z$, for all $\gg\leq \xi$, both $\X_\gg$ and $\Y_\gg$ are defined. Then $\X_{\xi+1}$, $\Y_{\xi+1}$, $F^+_\xi$, $F_\xi$ and $b_{\xi}$ are determined as follows.
\begin{enumerate}
\item Suppose $\X_\xi=(\mathcal{J}^{\vec{E}, f}_{\omega\a}, \in, \vec{E}, f)$ is a passive ${\sf{ses}}$\footnote{I.e., with
no last predicate} and there is an extender $F^*\in \vec{G}$,
 an extender $F$ over $\X_\xi$, and an
ordinal $\nu<\omega\a$ such that 
\begin{enumerate}
\item $\nu<\nu(F^*)$,
\item $F=F^*\cap ([\nu]^{\omega}\times \univ{\X_\xi})$, and
 \item setting
    \begin{center}
    $\Y_{\xi+1}=(\mathcal{J}_{\omega\a}^{\vec{E}, f}, \in, \vec{E}, f, \tilde{F})$
    \end{center}
     where $\tilde{F}$ is the amenable code of $F$\footnote{For the definition of the ``amenable code" see the last paragraph on page 14 of \cite{OIMT}.}, clause 2 fails for $\xi+1$.
     \end{enumerate} 
     Then
$\X_{\xi+1}=\C(\Y_{\xi+1})$\footnote{Recall that $\C(\M)$ is the core of $\M$.}, $F_\xi^+=\vec{G}(\xi)$ where $\xi$ is the least such that $F^*=\vec{G}(\xi)$ has the above properties, $F_\xi=F^+\cap ([\nu]^{\omega}\times \univ{\X_\xi})$ where $\nu$ is chosen so that the above clauses hold and $b_\xi=\emptyset$.
\item Suppose $\X_\xi=(\mathcal{J}^{\vec{E}, f}_{\omega\a}, \in, \vec{E}, f)$ is a passive ${\sf{ses}}$, $\a=\b+\gg$ and there is  $t=(\P_0, \T, \P_1, \U)\in \univ{\X_\xi|\omega\b} \cap \dom(\Lambda)$ such that setting $w=(\mathcal{J}_{\omega}(t), t, \in)$, $w$ is $(f, {\sf{sts}})$-minimal as witnessed by $\b$\footnote{See \rdef{important notation}. In particular, this means that we have to index the branch of $t$ at $\omega\a$.} and $\gg=\lh(t)$. 

Suppose there is a branch $b$ of $t$ (in $\N$) such that $(\T_{\leq \Q})^\frown t^\frown\{b\}$\footnote{When re-organizing $(\T_{\leq \Q})^\frown t^\frown\{b\}$ as a stack, there maybe a drop at $\Q$, as we have to drop to $(\Q|\eta)^\#$.} is $\mathcal{E}$-certified. Then set
     \begin{center}
     $\Y_{\xi+1}=(\mathcal{J}^{\vec{E}, f^+}_{\omega\b+\omega\gg}, \in, \vec{E}, f, \tilde{b})$
     \end{center}
     where $\tilde{b}\subseteq \omega\b+\omega\gg$ is defined by $\omega\b+\omega\nu\in \tilde{b}\iff \nu \in b$. Assuming clause 2 fails for $\xi+1$, 
$\X_{\xi+1}={\sf{core}}(\Y_{\xi+1})$, $F^+_\xi=F_\xi=\emptyset$ and $b_\xi=\tilde{b}$.\\\\
\textbf{Important Anomaly:} Suppose $t$ is $\sf{nuvs}$ and suppose $e\in \X_\xi|\omega\b$ is such that $\X_\xi|\omega\b\models {\sf{sts}_0}(t, e)$\footnote{See \rdef{sts0}. This means that $e$ is the branch of $t$ we must choose.}.
 If $e\not =b$ then $\Y_{\xi+1}$ is not an sts premouse over $X$ based on $\V$, and so clause 2 holds. \\
\item If $\X_\xi$ doesn't satisfy clause 2a or 2b then set $\Y_{\xi+1}=\mathcal{J}_{\omega}[\X_\xi]$. Assuming clause 2 fails for $\xi+1$, 
$\X_{\xi+1}=\C(\Y_{\xi+1})$, $F_\xi^+ =F_\xi=b_\xi=\emptyset$. 
\end{enumerate}
\item Suppose $\xi\leq \delta_z$ is a limit ordinal and for all $\gg<\xi$, both $\X_\gg$ and $\Y_\gg$ are defined. Then $\X_{\xi}$ and $\Y_{\xi}$ are determined as follows\footnote{$F_\xi, b_\xi$ will be defined at the next stage of the induction as in clause 2.}. Set $\Y_\xi=lim_{\a\rightarrow \xi}\X_\a$. Assuming clause 2 fails for $\xi$, $\X_\xi=\C(\Y_\xi)$. 
\item $\X_{\delta_z}=\Y_{\delta_z}$ and $F^+_{\d_z}=F_{\d_z}=b_{\d_z}=\emptyset$.
\end{enumerate}
We say that ${\sf{Le}}^{\mathcal{E}, c}((\Q|\eta)^\#)$ is \textbf{successful} if either for all $\xi<\d_z$ clause 2 above fails or letting $\xi_0$ be the least for which clause 2 holds, there is $\xi<\xi_0$ such that $\mathcal{J}_{\omega}[\X_\xi]\models``\eta$ is not a Woodin cardinal". The $``c"$ in ${\sf{Le}}^{\mathcal{E}, c}$ stands for ``certified". Given $\k<\d_z$, we can also define ${\sf{Le}}^{\mathcal{E}, c}((\Q|\eta)^\#)_{\geq \k}$ by requiring that in clause 3.a, $\cp(F)\geq \k$. 

We will use the following terminology. We say $\K$ is a $\Y$-\textbf{model} of ${\sf{Le}}^{\mathcal{E}, c}((\Q|\eta)^\#)_{\geq \k}$ if for some $\gg\leq \d_z$, $\Q=\Y_\gg$. Similarly we define $\X$-model and other such expressions. We say $\K$ is \textbf{the last model} of ${\sf{Le}}^{\mathcal{E}, c}((\Q|\eta)^\#)_{\geq \k}$ if $\K=\Y_{\d_z}$. $\myqedhere$
\end{definition}

We can now define the $\mathcal{E}$-certified iterations. 

\begin{definition}\label{e-certified iterations} Suppose $\V\in \N$ is a hod premouse extending $\R$ such that $\R=\V^b$. Suppose $\T\in \N$ is a stack or an st-stack on $\V$ and $E\in \mathcal{E}$. We say $\T$ is \textbf{$E$-certified} if the following conditions are satisfied.
\begin{enumerate}
\item $\T$ is $\pi_E$-realizable.
\item Suppose $\tau \in (R^b)^\T$ is such that letting $\U=_{def}{\sf{nc}}^\T_\tau$,  $\U$ is above $\M_\tau^b$. Let $\a<\lh(\U)$ be a limit ordinal and let $c=[0, \a)_\U$. Then the following conditions hold.
\begin{enumerate}
\item If $\m^+(\U\rest \a)\models ``\d(\U\rest \a)$ is not a Woodin cardinal"\footnote{See \rdef{mtsharp}.} then $\Q(c, \U\rest \a)$ exists and $\Q(c, \U\rest \a)\insegeq \m^+(\U\rest \a)$.
\item If $\m^+(\U\rest \a)\models ``\d(\U\rest \a)$ is a Woodin cardinal" and there is $\W$ such that 
\begin{enumerate}
\item $\W$ appears on the ${\sf{Le}}^{\mathcal{E}, c}(\m^+(\U|\rest \a))$ construction of $\N$ and
\item $\W\models ``\d(\U\rest \a)$ is a Woodin cardinal" but $\mathcal{J_{\omega}}[\W]\models ``\d(\U\rest \a)$ is not a Woodin cardinal", 
\end{enumerate}
then $\Q(c, \U\rest \a)$ exists and $\Q(c, \U\rest \a)=\W$.
\item The above two clauses fail. Then $\T$ is an st-stack, $\a+1=\lh(\U)$ and $\tau+\a\in R^\T\cap {\sf{max}}^\T$.  
\end{enumerate}
\end{enumerate}
We say that $\T$ is \textbf{$\mathcal{E}$-certified} if for some $\l$, $\T$ is $E$-certified for every $E\in \mathcal{E}$ such that $\lh(E)>\l$. $\myqedhere$
\end{definition}

And finally we define $\mathcal{E}$-certified strategies. 

\begin{definition}\label{e-certified strategies} Suppose $\V\in \N$ is a hod premouse extending $\R$ such that $\R=\V^b$. We let $\Lambda_\V$ be the partial strategy of $\V$ with the property that 
\begin{enumerate}
\item $\dom(\Lambda_\V)$ consists of $\mathcal{E}$-certified stacks $\T$ of limit length, and
\item for all $\T\in \dom(\Lambda_\V)$, $\Lambda_\V(\T)=b$ if $b$ is the unique $x$ such that $\T^\frown\{x\}$ is $\mathcal{E}$-certified.
\end{enumerate}
We say $\Lambda_\V$ is the $\mathcal{E}$-certified strategy of $\V$. $\myqedhere$
\end{definition}

\begin{remark} According to our definition, the $\mathcal{E}$-certified strategy of $\V$ is an ordinary strategy acting on smooth iterations. However, it is straightforward to generalize our definition to obtain a strategy acting on generalized stacks\footnote{See \rsec{undropping game sec}.}. Below we assume that the strategy described in \rdef{e-certified strategies} acts on generalized stacks, and is a partial $(\d_z, \d_z, \d_z)$-iteration strategy. However, there is one subtle point that we address.

$\Lambda_\V$, the $\mathcal{E}$-certified strategy of $\V$, must be self-cohering\footnote{See \rdef{self-cohering}.}. To achieve this, we use the following idea. Suppose $\T$ is a $\pi_E$-realizable generalized stack and $F$ is the un-dropping extender of $\T$. Let $\Q$ be the last model of $\T$. We then have $\sigma:\Q^b\rightarrow \R'\insegeq \pi_E(\R)$ coming from the definition of $\pi_E$-realizability. To ensure that $\Lambda_\V$ will be self-cohering it is enough to use $\sigma':Ult(\V, F)^b\rightarrow \pi_E(\R)$ given by $\sigma'(x)=\pi_E(f)(\sigma(a))$ where $f\in \V$ and $a\in [\d^{\Q^b}]^{<\omega}$ is such that $x=\pi_F(f)(a)$. $\myqedhere$
\end{remark}

Suppose now that
\begin{center}
$\T=((\M_\a)_{\a<\eta}, (E_\a)_{\a<\eta-1}, D, R, (\beta_\a, m_\a)_{\a\in R}, T)$
\end{center}
is $\pi_E$-realizable as witnessed by $(\sigma_\a: \a\in R^b)$ and $(\Psi_\a:\a\in R^b)$. Using the language of \rsec{condensing sets sec} applied in $Ult(\N, E)$ to $\pi_E(\mathcal{F}\rest \k)$, it is not hard to see that for $\a\in (R^b)^\T$, $\M_\a^b=\P_{Y_\a}$ where $Y_\a=\sigma_\a[\M_\a^b]$ and $\Psi_\a=\Sigma_{Y_\a}$. It now follows from \rlem{lem:uniqueness of strategies} that there are unique $\mathcal{E}$-certified strategies.

\begin{lemma}\label{uniqueness of e-certified} Suppose $\V\in \N$ is a hod premouse extending $\R$ such that $\R=\V^b$. Suppose $\Lambda$ and $\Psi$ are two $\mathcal{E}$-certified strategies for $\V$. Then $\Lambda=\Psi$. 
\end{lemma}

It is also not hard to show that $\mathcal{E}$-certified iterations are according to $\Phi^+$, which is the topic of the next subsection.

\subsection{Correctness of certified strategies}

Our goal is still to prove \rthm{capturing of hod pairs} and our set up is as in \rter{objects for the proof} and \rdef{definition of r}.

\begin{lemma}\label{correctness of certified strategies: msc} Suppose $\S^*$ is a $\Phi^+$-iterate of $\R^+$ via an iteration that is entirely above $\d^{\R}$. Suppose further that $\S\insegeq_{hod} \S^*$ is such that $\S^b=\R$ and $\S\in \N$. Let $\T\in \N$ be a stack on $\S$\footnote{We assume that $\T$ is a stack, but the proof works for generalized stacks as well.}. Suppose $\T$ is $\mathcal{E}$-certified. Then $\T$ is according to $\Phi^+_\S$.  Thus, $\Lambda_\S=\Phi^+_\S\cap \N^2$.\footnote{This equation does not imply that $\Lambda_\S=\Phi^+_\S\rest \N$, simply because it does not imply that if $x\in \dom(\Phi^+_\S)\cap \N$ then $\Phi^+_\S(x)\in \N$. To get the aforementioned equality, we need to show that $\Lambda_\S$ is total.}
\end{lemma}
\begin{proof} Suppose 
\begin{center}
$\T=((\M_\a)_{\a<\eta}, (E_\a)_{\a<\eta-1}, D, R, (\beta_\a, m_\a)_{\a\in R}, T)$, 
\end{center}
and suppose $\a\in R^b$ is such that $\T_{\leq \a}$ is according to $\Phi^+_\S$. We want to show that $\U={\sf{nc}}^\T_\a$ is according to $\Phi^+_{\M_\a}$. 

Suppose first that $\U$ is based on $\M_\a^b$\footnote{There is yet another case: namely, $\a=\max{R^b}$ and $\U=\T_{\geq \a}$. But this case is very similar to our two cases.}. Let $E\in \mathcal{E}$ be such that $\T$ is $\pi_E$-realizable
 as witnessed by $(\sigma_\a: \a\in R^b)$ and $(\Psi_\a:\a\in R^b)$. Let $E^*$ be the background certificate of $E$ and let $k: Ult(\N, E)\rightarrow \pi_{E^*}(\N)$ be the canonical factor map. Notice that for $\a\in R^b$,\\\\
(1) $k\rest \N|\xi=id$ where $\xi$ is the least such that $\T\in \N|\xi$.\\
(2) In $\pi_{E^*}(\N)$, $k(\sigma_\a):\M_\a^b\rightarrow \pi_{E^*}(\N)$ and $k(\Psi_\a)$ is the $k(\sigma_\a)$-pullback of $\pi_{E^*}(\Phi)$.\\
(3) $k(\sigma_\a)\rest \M_\a|\d^{\M_\a^b}$ is the iteration embedding according to $k(\Psi_\a)$.\\\\ 

Let $F$ be the un-dropping extender of $\T_{\leq \a}$ and set $\K^+=Ult(\R^+, F)$ and $j=\pi^{\Phi^+_{\K^+}}_{\K^+, \pi_{E^*}(\R^+)}\rest \M_\a|\d^{\M_\a^b}$. Notice now that\\\\
(4) $\Phi_{\M_\a|\d^{\M_\a^b}}$ is the $j$-pullback of $\pi_{E^*}(\Phi)$ and $j$ is the iteration embedding according to $\Phi_{\M_\a|\d^{\M_\a^b}}$.\\\\
 As the pairs $(k(\sigma_\a), k(\Psi_\a))$ and $(j, \Phi_{\M_\a|\d^{\M_\a^b}})$ have the same property, it follows from \rlem{lem:uniqueness of strategies} that $k(\sigma_\a)=j$ and $k(\Psi_\a)=\Phi_{\M_\a|\d^{\M_\a^b}}\rest \pi_{E^*}(\N)$. Since $k(\U)=\U$, in the case $\U$ is based on $\M_\a|\d^{\M_\a^b}$, we have that $\U$ is according to $k(\Psi_\a)$ and therefore, $\U$ is according to $\Phi_{\M_\a|\d^{\M_\a^b}}$, and in the case $\U$ is above $\d^{\M_\a^b}$, we have that $\U$ is according to the unique strategy of $\M_\a^b$ that witnesses the fact that $\M_\a^b$ is a $\Phi_{\M_\a|\d^{\M_\a^b}}$-mouse over $\M_\a|\d^{\M_\a^b}$. 
 
 Suppose now that $\U$ is above $\sf{ord}(\M_\a^b)$. Here, we need to see that\\\\
 (a) if $\b<\lh(\U)$ is a limit ordinal then letting $b=[0, \b)_\U$, either $\Q(b, \U)\insegeq \m^+(\U)$ or else $\Q(b, \U)\insegeq {\sf{Lp}}^{\Gamma, (\Phi^+)^{sts}_{\m^+(\U)}}(\m^+(\U))$.\\\\
The following lemma establishes (a). For convenience, we will ignore the objects introduced above and treat next lemma in a general context. Thus $\T$ in the next lemma is not the $\T$ fixed above. 
\begin{lemma}\label{q-structure correctness lemma} Suppose $\T$ is an $\mathcal{E}$-certified iteration of $\S$, $\a\in R^b$ and $\U={\sf{nc}}^\T_\a$ is above $\sf{ord}(\M_\a^b)$. Suppose further that $\b<\lh(\U)$ is a limit ordinal and $\U_{<\b}$ is according to $\Phi^+_{\M_\a}$. Let $\Q=\M_\b^\U$ and $\eta>\d^{\Q^b}$ be such that $\mathcal{J}_{\omega}[(\Q|\eta)^\#]\models ``\eta$ is a Woodin cardinal" and let $\W\insegeq \Q$ be an sts mouse over $(\Q|\eta)^\#$.  Then $\W$ is a $(\Phi^+)_{(\Q|\eta)^\#}^{stc}$-sts mouse. 
\end{lemma}
\begin{proof}
Towards a contradiction assume that $\W$ is not a $(\Phi^+)_{(\Q|\eta)^\#}^{stc}$-sts mouse. It follows that $b=[0, \b)_\U$ is not the branch chosen by $\Phi^+_{\S}$. For convenience, we change our notation and let $\U$ be $\U\rest \b$  and $\Q=\m^+(\U)$. It follows from \rdef{e-certified iterations} that\\\\
(1) $\W$ is a model appearing in the fully backgrounded $\mathcal{E}$-realizable construction over $(\Q|\eta)^\#$ done in $\N$. \\\\
What we need to see is that $\W$ is a $(\Phi^+)_{\Q}^{stc}$-sts mouse over $\Q$. To show this it is enough to show that every stack indexed in $\W$ is according to $(\Phi^+)_{\Q}^{stc}$. To show this later fact, it is enough to show that\\\\
(b) if $t=(\Q, \U_0, \Q_1, \U_1)$ is an indexable stack\footnote{See \rdef{authentic stacks of length 2}.} on $\Q$ appearing in the fully backgrounded $\mathcal{E}$-realizable construction over $\Q$ (done in $\N$) and $c$ is the branch of $t$ indexed in this construction then $t^\frown \{c\}$ is according to $(\Phi^+)_{\Q}^{stc}$.\\\\
(b) is indeed enough. To see this, notice that if $s=(\Q, \U'_0, \Q_1', \U'_1)$ is indexed in $\W$ and $c'$ is the branch of $s$ indexed in $\W$ then for some stack $t=(\Q, \U_0, \Q_1, \U)$ as in (b) if $e$ is the branch of $t$ then $s^\frown \{c\}$ is a hull of $t^\frown \{e\}$. If $t$ is according to $(\Phi^+)_{\Q}^{stc}$ then it follows from hull condensation of $(\Phi^+)_{\Q}^{stc}$ that $s$ is also according to $(\Phi^+)_{\Q}^{stc}$. We now work towards showing that $t$ is according to $(\Phi^+)_{\Q}^{stc}$. 

Suppose first that $\U_0$ is according to $(\Phi^+)_{\Q}^{stc}$. We then have that $\U_1$ is a stack based on $\Q_1^b$. Because $(\T_{\leq \a})^\frown t$ is $\mathcal{E}$-certified, we can fix an extender $E\in \mathcal{E}$ such that $(\T_{\leq \a})^\frown t$ is $\pi_E$-realizable. We then have $\sigma:\Q_1^b\rightarrow \pi_{E}(\R)$ such that $\pi_E\rest \R=\sigma\circ \pi^{\U_0, b}\circ \pi^{\T_{\leq \Q}, b}$. We also have that $\U_1^\frown \{c\}$ is according to the $\sigma$-pullback of $\pi_E(\Phi_\R)$. Therefore, $t$ is according to $(\Phi^+)_{\Q}^{stc}$. 

It remains to show that $\U_0$ is according to $(\Phi^+)_{\Q}^{stc}$. Without loss of generality, we assume that
\begin{itemize}
\item $\lh(\U_0)=\gg+1$, 
\item $\gg$ is a limit ordinal, 
\item $\U_0\rest \gg$ is according to  $(\Phi^+)_{\Q}^{stc}$,
\item $[0, \gg)_{\U_0}\not =(\Phi^+)_{\Q}^{stc}(\U_0)$,
\item there is $\zeta\in R^{\U_0\rest \gg}$ such that $(\U_0)_{\geq \zeta}={\sf{nc}}^{\U_0}_{\zeta}$ and $\pi^{\U_0, b}$ exists,
\item $\mathcal{J}_{\omega}[\m^+(\U_0)]\models ``\d(\U_0)$ is a Woodin cardinal".
\end{itemize}
The last two clauses can be shown by examining the proof given for $\U_1$. Set $c_0=[0, \gg)_{\U_0}$, $\Q_0=_{def}\Q$, $\Q_2=\m^+(\U_0)$, $\W_0=_{def}\W$ and $\W_2=\Q(c_0, \U_0)$ . We then have that \\\\
(2) $\W_2$ appears in the fully backgrounded $\mathcal{E}$-realizable construction over $\Q_2$ (done in $\N$).\\\\ 
Clearly (2) leads to an infinite descend. 
\end{proof}
\end{proof}

\begin{rem}\label{remark on the important anomaly in e-realizable constructions} It is important to note that \rlem{q-structure correctness lemma} does not resolve the Important Anomaly stated in clause 3b of \rdef{e-certified background constructions}.
\end{rem}

The following straightforward yet important lemmas are steps towards showing that the Important Anomaly stated in clause 3b of \rdef{e-certified background constructions} does not occur. The reader may wish to review \rdef{authentic lsp} and \rdef{nus stacks}.

\begin{lemma}\label{authentic iterations are realizable1} Suppose 
\begin{enumerate}
\item $\V\in \N$ is a hod premouse extending $\R$ such that $\R=\V^b$,
\item $\T\in \N$ is either a stack on $\V$ or an st-stack\footnote{See \rdef{st-stack}.} on $\V$\footnote{If $\T$ is an st-stack then $\M$ must be of lsa type.},
\item $\pi^{\T, b}$ is defined, $\T$ has a last model and $\mathcal{E}$-realizable.
\end{enumerate}
Let $\S$ be the last model of $\T$ and suppose $\Q$ is authenticated\footnote{See \rdef{authentic lsp}.} by $\T$ and is meek and of limit type\footnote{Thus, clause 3 of \rdef{authentic lsp} holds.}. Then $\W, \U, \sigma$ be as in \rdef{authentic lsp} and letting $k:\R\rightarrow \Q$ be given by $k(x)=y$ if and only $\sigma^{-1}(\pi^{\T, b}(x))=\pi^{\U}(y)$, $(k, \Q)$ is $\mathcal{E}$-realizable.
\end{lemma}

The reader may wish to review \rnot{notation for iteration trees} and \rdef{e-certified strategies}. The lemma below implies that  the Important Anomaly stated in clause 3b of \rdef{e-certified background constructions} does not occur. 

\begin{lemma}\label{authentic iterations are realizable} Suppose 
\begin{enumerate}
\item $\V\in \N$ is a hod premouse extending $\R$ such that $\R=\V^b$,
\item $\T\in \N$ is either a stack on $\V$ or an st-stack\footnote{See \rdef{st-stack}.} on $\V$\footnote{If $\T$ is an st-stack then $\V$ must be of $\#$-lsa type.},
\item $\pi^{\T, b}$ is defined, $\T$ has a last model and $\mathcal{E}$-realizable.
\end{enumerate}
Let $\S'$ be the last model of $\T$ and suppose $\eta<{\sf{ord}}(\S')$ is such that $\mathcal{J}_{\omega}[(\S'|\eta)^\#]\models``\eta$ is a Woodin cardinal". Suppose $\S=_{def}(\S'|\eta)^\#$ is such that $\S\insegeq_{hod}\S'$ and $\U\in \N$ is an $\sf{nuvs}$ stack according to $(\Lambda_\V)_\S$ such that $\pi^{\U, b}$ is defined. Let $\Q=\m^+(\U)$ and suppose $t\in \N$ be an indexable stack on $\Q$ which is $(\S, (\Lambda_\V)_\S)$-authenticated\footnote{Notice that we, at this point, do not know that $\Lambda_\V$ is a total strategy in $\mathcal{J}[\N]$.}. Then $\T^\frown \U^\frown t$ is according to $\Lambda_\V$. 
\end{lemma}
\begin{proof}
Suppose $t=(\Q_0, \X_0, \Q_1, \X_1)$. Assume first that $\X_0$ is according to $(\Lambda_\V)_\Q$. Set $p=\T^\frown \U^\frown \X_0$ and let $\sigma=\pi^{p, b}$. It follows from \rlem{correctness of certified strategies: msc} that $p$ is according to $\Lambda_\V$, and the previous lemma implies that $(\sigma, \Q_1^b)\in \mathcal{F}_{\mathcal{E}}$. Because $(\Q_1^b, \X_1)$ is a $(\S, (\Lambda_\V)_\S)$-authenticated iteration, it follows from \rlem{correctness of certified strategies: msc} that $\X_1$ is according to $\Psi_{\sigma, \Q_1^b}$, and therefore, $\T^\frown \U^\frown t$ is according to $\Lambda_\V$. 

Thus, it is enough to show that $\X_0$ is according to $(\Lambda_\V)_\Q$. The argument given above implies that it is enough to show that for every $\a\in R^{\X_0}$ such that $\pi^{\X_0, b}_{0, \a}$ is defined and ${\sf{nc}}^{\X_0}_\a$ is a stack on $\M_\a^{\X_0}$ above ${\sf{ord}}((\M_\a^{\X_0})^b)$ then ${\sf{nc}}^{\X_0}_\a$ is according to $(\Lambda_\V)_{\M_\a^{\X_0}}$.

Assume then $\a$ is as above and  $(\X_0)_{\leq \a}$ is according to $(\Lambda_\V)_{\Q}$. Set $\M=\M_\a^{\X_0}$, $\X=(\X_0)_{\leq \a}$ and let $\Y={\sf{nc}}^{\X_0}_\a$. We want to see that $\Y$ is according to $(\Lambda_\V)_{\M}$. Let $\b<\lh(\Y)$ be a limit ordinal such that $\Y_{<\b}$ is according to $(\Lambda_\V)_{\M}$. We want to see that if $b=[0, \b)_\Y$ then $b=(\Lambda_\V)_\M(\Y_{<\b})$. The dificult case is when $\Q(b, \Y_{<\b})$ exists and is an sts mouse over $\m^+(\Y_{<\b})$. In this case, we want to see that $\Q(b, \Y_{<\b})$ is a model appearing in the fully backgrounded $\mathcal{E}$-realizable construction over $\m^+(\U_0)$ (done in $\N$). This would follows from the proof of the previous lemma. Our strategy for showing this is by showing (a) and (b) where these are the following statements:\\\\
(a) $\Q(b, \Y_{<\b})$ is a $(\Phi^+_{\m^+(\Y_{<\b})})^{stc}$-mouse over $\m^+(\Y_{<\b})$.\\\\
(b) If $\W$ is a $(\Phi^+_{\m^+(\Y_{<\b})})^{stc}$-mouse over $\m^+(\Y_{<\b})$ then $\W$ appears in the fully backgrounded $\mathcal{E}$-realizable construction over $\m^+(\U_0)$ (done in $\N$). More precisely, letting 
\begin{center}
${\sf{Le}}^{\mathcal{E}, c}(\m^+(\Y_{<\b}))=(\mathcal{Z}_\gg , \K_\gg, F^+_\gg, F_\gg, b_\gg: \gg\leq \delta_z)$
\end{center}
be the fully backgrounded $\mathcal{E}$-realizable construction over $\m^+(\Y_{<\b})$ done in $\N$ then for some  $\gg<\d_z$, $\mathcal{Z}_\gg=\W$.\\\\
(a) is a consequence of strong branch condensation of $\Phi^+$ and can be shown using the proof of \rsublem{next lemma} and \rlem{correctness of certified strategies: msc}. 

(b) is a consequence of the fact that $\Lambda_\V$ is total, and hence $\Phi^+_\V\rest \mathcal{J}[\N]=\Lambda_\V$ (see \rlem{correctness of certified strategies: msc}). Assuming that $\Lambda_\V$ is total, (b) can be proven by simply comparing $\W$ with the ${\sf{Le}}^{\mathcal{E}, c}(\m^+(\Y_{<\b}))$ construction. The stationarity of ${\sf{Le}}^{\mathcal{E}, c}(\m^+(\Y_{<\b}))$ implies that the construction side doesn't move, and the fact that $\Lambda_\V$ is total implies that the construction doesn't break down because in clause 3b of \rdef{e-certified background constructions} we are unable to find the desired branch. The Important Anomaly stated in clause 3b of \rdef{e-certified background constructions} does not occur (at least doesn't occur before reaching $\W$) because these type of branches are chosen internally and both the construction side and the $\W$-side must be choosing the same branch. But on the $\W$-side, the branch is according to $\Phi^+_\V$ and therefore, according to $\Lambda_\V$.  In the next subsection, we will prove that $\Lambda_\V$ is total, and more details will be given. 
\end{proof}

\subsection{$\Lambda_\V$ is total}

The goal of this subsection is to show that $\Lambda_\V$, the $\mathcal{E}$-certified strategy of $\V$, is total.

\begin{lemma}\label{lambdav is total} $\V\in \N$ is a hod premouse extending $\R$ such that $\R=\V^b$. Then $\Lambda_\V$ is total, and hence, $\Phi^+_{\V}\rest \N=\Lambda_\V$.
\end{lemma}
\begin{proof} The equality $\Phi^+_{\V}\rest \N=\Lambda_\V$ follows from \rlem{correctness of certified strategies: msc}. Suppose $\T\in \N$ is a stack according to $\Lambda_\V$ and of limit length. We want to show that $\Lambda_\V(\T)$ is defined. Let $b=\Phi^+_\V(\T)$. It is enough to show that $b\in \N$ and $\T^\frown \{b\}$ is $\mathcal{E}$-certified. When discussing objects in $\T$, we omit $\T$ from superscripts. The two non-straightforward cases are the following:\\\\
Case 1: There is $\a\in R$ such that $\pi^{\T, b}_{0, \a}$ is defined and $\T_{\geq \a}$ is based on $\M_\a^b$.\\
Case 2: There is $\a\in R$ such that $\a=\max(R)$, $\pi^{\T, b}_{0, \a}$ is defined, $\T_{\geq \a}$ is above ${\sf{ord}}(\M_\a^b)$, $\T_{\geq \a}$ doesn't have a fatal drop and $\mathcal{J}_{\omega}[\m^+(\T_{\geq \a}))]\models ``\d(\T_{\geq \a})$ is a Woodin cardinal". \\\\
 \rlem{iterates are neatly e-realizable} implies that under Case 1, $b\in \N$ and $\T$ is $\mathcal{E}$-certified. Assume then Case 2. Let $\W=\Q(b, \T)$. We have that $\W$ is a $(\Phi^+_{\m^+(\T)})^{stc}$-mouse over $\m^+(\T)$. It is now enough to show that $\W$ appears in the  ${\sf{Le}}^{\mathcal{E}, c}(\m^+(\T))$ construction of $\N$. Let 
 \begin{center}
${\sf{Le}}^{\mathcal{E}, c}(\m^+(\T))=(\mathcal{X}_\gg , \Y_\gg, F^+_\gg, F_\gg, b_\gg: \gg\leq \delta_z)$
\end{center}
be the ${\sf{Le}}^{\mathcal{E}, c}(\m^+(\T))$ construction of $\N$. For each $\gg\leq \d_z$, we let $\U_\gg$ be the normal stack on $\W$ above $\d(\T)$ that is constructed by comparing $\W$ with $\Y_\gg$. What we need to show is that\\\\
(a) the construction side never moves, i.e., for any $\gg<\d_z$, letting $\W_\gg$ be the last model of $\U_\gg$, if for some $\tau$, $\W_\gg|\tau=\Y_\gg|\tau$ and $\W_\gg||\tau\not= \Y_\gg||\tau$ then $\tau\not \in \dom(\vec{E})^{\Y_\gg}$ and $\tau$ is not an index of a branch in $\W_\gg$. \\\\
The fact that $\tau\not \in \dom(\vec{E})^{\Y_\gg}$ follows from the stationarity of the fully backgrounded constructions\footnote{See \rthm{existence of thick sets} and the references given there.}. It is then enough to show that $\tau$ is not an index of a branch in $\W_\gg$. Suppose to the contrary, and let \\\\
(1) $t\in \W_\gg$ be an indexable stack whose branch is indexed in $\W_\gg$ at $\tau$.\\\\
 Because $\W_\gg|\tau=\Y_\gg|\tau$, it follows that all initial segments of $t$ are according to $(\Lambda_\V)_{\m^+(\T)}$. We need to show that\\\\
 (b) there is a branch indexed at $\tau$ in $\Y_\gg$ and that branch is according to $(\Phi^+_{\m^+(\T)})^{stc}$.\\\\
 Again the two dificult cases are the cases stated under Case 1 and Case 2, and Case 1 can be analyzed as above, and so we only state Case 2. We thus have that\\\\
 (2) setting $\K=\m^+(\T)$, $t=(\K, \mathcal{Z})$ and $\mathcal{J}_{\omega}[\m^+(\mathcal{Z}_0)]\models ``\d(\mathcal{Z}_0)$ is a Woodin cardinal".\\\\
 Notice that $t\in \W_\gg$ and the branch chosen for $t$ both in $\W_\gg$ and in $\Y_\gg$ depends on $\W_\gg|\tau=\Y_\gg|\tau$. Hence, both $\W_\gg$ and $\Y_\gg$ must have the same branch of $t$ indexed in their strategy predicates. 
\end{proof}

\subsection{Mixed hod pair constructions}\label{sec: mixed hod pair constructions} 

We devote this entire subsection to the definition of a construction producing the iterate of $\R^+$. In this construction, we use $\mathcal{E}$-certification method to acquire extenders with critical point $\d^{\R}$, and we use the total extenders on the sequence of $\N$ to generate extenders with critical points $>\d^\R$. First we define $\mathcal{E}$-certified extenders. The reader may wish to review \rdef{canonical e-realization}.

\begin{definition}\label{e-certified extenders} Suppose $\Q\in \N$ is a hod premouse such that $\Lambda_\Q$ (see \rdef{e-certified strategies}) is total and $\Q^b=\R$. Suppose $F$ is an extender such that $(\Q, \tilde{F})$ is a reliable ${\sf{lses}}$ where $\tilde{F}$ is the amenable code of $F$. We say $F$ is \textbf{$\mathcal{E}$-certified} if
\begin{itemize}
\item $(\pi_F\rest \R, \pi_F(\R))\in \mathcal{F}_{\mathcal{E}}$ and
\item  for some $\N$-strong cardinal $\l$, for any $E\in \mathcal{E}$ such that $\lh(E)>\l$, setting $\tau=\tau_{\pi_F\rest \R, \pi_F(\R)}$\footnote{See \rdef{canonical e-realization}.},
\begin{center}
$(a, A)\in F\iff \tau(a)\in \pi_E(A)$.
\end{center}
\end{itemize}
We say that $\tau$ is the \textbf{$E$-realizability map} of $F$. $\myqedhere$
\end{definition}

The next lemma shows that $\mathcal{E}$-certified extenders are on the sequence of $\R^+$ and its iterates.

\begin{lemma}\label{ecertified extenders are on the sequence} Suppose $\S^*\in pI(\R^+, \Phi^+)$ and $\S\inseg_{hod}\S^*$ is such that $\S\in \N$ and $\S^b=\R$. Suppose $F$ is such that $(\S, \tilde{F})$ is a reliable $\sf{lses}$ where $\tilde{F}$ is the amenable code of $F$ and $F$ is $\mathcal{E}$-certified. Then $F\in \vec{E}^{\S^*}$.
\end{lemma}
\begin{proof} Let $\gg={\sf{ord}}(\S)$ and suppose $F^*\in \vec{E}^{\S^*}(\gg)$. Then $F^*$ has exactly the same property as $F$ and therefore, $F=F^*$. Thus, it is enough to show that $\gg\in \dom(\vec{E}^{\S^*})$. Suppose first that there is $\gg'\in \dom(\vec{E}^{\S^*})$ such that $\S\inseg_{hod}\S^*|\gg'$ and if $G'=\vec{E}^{\S^*}(\gg')$ then $\cp(G')=\d^{\R}$. Let $\gg^*$ be the least such $\gg'$ and set $G=\vec{E}^{\S^*}(\gg^*)$. As $F$ and $G$ both have the property described in \rdef{e-certified extenders}, $F$ is an initial segment of $G$, and therefore, $\gg=\gg^*$ and $\gg\in \dom(\vec{E}^{\S^*})$. Suppose then that\\\\
(1) there is no $\gg'\in \dom(\vec{E}^{\S^*})$ such that $\cp(\vec{E}^{\S^*}(\gg'))=\d^{\R}$.\\\\
 Because $F$ is $\mathcal{E}$-certified, we have that for some $\N$-strong cardinal $\l$, whenever $E\in \mathcal{E}$ is such that $\lh(E)>\l$, some proper initial segment of $\pi_E(\R)$ is a $\Phi^+_\S$-iterate of $\S$. Therefore, $(\S, \Phi^+_\S)$ is in ${\sf{HP}}^\Gamma$, and hence, $(\S^*, \Phi^+_{\S^*})\in {\sf{HP}}^\Gamma$\footnote{$\Gamma$ was introduced in \rnot{the gamma}.}. This is because (1) implies that $\S^*\inseg {\sf{Lp}}^{\Gamma, \Phi^+_\S}(\S)$ or $\S^*\inseg {\sf{Lp}}^{\Gamma, (\Phi^+_\S)^{stc}}(\S)$.
\end{proof}

Next we introduce the mixed hod pair constructions.

\begin{definition}\label{mixed hod pair construction} We say that 
\begin{center}${\sf{mhpc}}=(\M_\gg , \N_\gg, Y_\gg, \Phi_\gg, F^+_\gg, F_\gg, b_\gg : \gg\leq \delta)$
\end{center}
 is the output of the \textbf{mixed hod pair construction} of $\mathbb{N}$ over $\R$ if the following conditions hold.
 \begin{enumerate}
 \item $\M_0=\mathcal{J}_\omega[\R]$, and for all $\gg\leq \d$, each of $\M_\gg$ and $\N_\gg$ is either undefined or is an $\sf{hp}$-indexed $\sf{lses}$ (see \rdef{sis}).
 \item For all $\gg\leq \d$, if $\M_\gg$ is defined then $Y_\gg=Y^{\M_\gg}$ (see \rdef{layered hybrid j-structure}).
 \item For all $\gg\leq \d$, if $\M_\gg$ is defined then $\Phi_\gg=\Phi_{\M_\gg}$ is the $\mathcal{E}$-certified strategy of $\M_\gg$\footnote{See \rdef{e-certified strategies}.}.
 \item For all $\gg\leq \d$, if $\N_\gg$ is defined and either
 \begin{enumerate}
 \item $\N_\gg$ is not a reliable $\sf{hp}$-indexed $\sf{lses}$\footnote{Recall clause 2 of \rdef{strategy lhp}. To verify that $\N_\gg$ is $\sf{lses}$, we need to verify that clause 2 of \rdef{strategy lhp} holds.} or
 \item $\N_\gg$ is a reliable $\sf{hp}$-indexed $\sf{lses}$ but for some $\Q\in Y^{\N_\gg}$ such that $\Q$ is meek or gentle\footnote{See \rdef{pre-hod-like}.} and for some $n<\omega$, $\rho_n(\N_\gg)\leq \d^\Q$, or
 \item $\Phi_\gg$ is not total,
 \end{enumerate}
  then all remaining objects with index $\geq \gg$ are undefined. \\\\
  For all $\gg\leq \eta$ for which clause 4 (the above statement) fails, $\pi_\gg:\C(\N_\gg)\rightarrow \N_\gg$ is the uncollapse map. 
 \item Suppose for some $\xi<\d$, for all $\gg\leq \xi$, both $\M_\gg , \N_\gg$ are defined. Then $\M_{\xi+1}$, $\N_{\xi+1}$, $Y_{\xi+1}$, $\Phi_{\xi+1}$, $F_{\xi}^+$, $F_\xi$ and $b_{\xi}$ are deteremined as follows.
 \begin{enumerate}
\item Suppose $\M_\xi=(\mathcal{J}^{\vec{E}, f}_{\omega\a}, \in, \vec{E}, f, Y_\xi, \in)$ is a passive ${\sf{hp}}$-indexed $\sf{lses}$\footnote{I.e., with
no last predicate. See \rdef{sis}.}, there is an extender $H^*\in \mathcal{E}$ 
an extender $H$ over $\M_\xi$, and an
ordinal $\nu<\omega\a$ such that $\nu<\lh(H^*)$ and setting
    \begin{center}
    $H=H^*\cap ([\nu]^{\omega}\times \univ{\M_\xi})$, and
    $\N_{\xi+1}=(\mathcal{J}_{\omega\a}^{\vec{E}, f}, \in, \vec{E}, f, Y_\xi, \tilde{H}, \in)$
    \end{center}
    where $\tilde{H}$ is the amenable code of $H$, clause 4.a fails for $\xi+1$. 
    Then letting $\iota\in \dom(\vec{E}^\N)$ be the least such that $H^*=_{def}\vec{E}^\N(\iota)\in \mathcal{E}$ has the above properties, 
    \begin{center}
    $\N_{\xi+1}=(\mathcal{J}_{\omega\a}^{\vec{E}, f}, \in, \vec{E}, f, Y_\xi, \tilde{H}, \in)$
    \end{center}
    where $\tilde{H}$ is the amenable code of $H$\footnote{Here $H$ is what is determined by $H^*$. For the definition of the ``amenable code" see the last paragraph on page 14 of \cite{OIMT}.}. Assuming clause 4 fails for $\xi+1$, the remaining objects 
     are defined as follows.
    \begin{enumerate} 
    \item $\M_{\xi+1}=\C(\N_{\xi+1})$\footnote{Recall that $\C(\M)$ is the core of $\M$.}, 
    \item $F^+_\xi= H^*$ and $F_\xi=H$, 
    \item $b_\xi=\emptyset$ and
    \item $Y_{\xi+1}=\pi^{-1}_{\xi+1}(Y_\xi)$.
\end{enumerate}
\item Suppose $\M_\xi=(\mathcal{J}^{\vec{E}, f}_{\omega\a}, \in, \vec{E}, f, Y_\xi, \in)$ is a passive ${\sf{hp}}$-indexed $\sf{lses}$\footnote{I.e., with
no last predicate.} and there is an extender 
 $H$ over $\M_\xi$ such that setting
    \begin{center}
    $\N_{\xi+1}=(\mathcal{J}_{\omega\a}^{\vec{E}, f}, \in, \vec{E}, f, Y_\xi, \tilde{H}, \in)$
    \end{center}
    where $\tilde{H}$ is the amenable code of $H$, clause 4.a fails for $\xi+1$ and $H$ is $\mathcal{E}$-certified as defined in \rdef{e-certified extenders}. Assuming clause 4 fails for $\xi+1$, the remaining objects 
     are defined as follows.
    \begin{enumerate} 
    \item $\M_{\xi+1}=\C(\N_{\xi+1})$\footnote{Recall that $\C(\M)$ is the core of $\M$.}, 
    \item $F^+_\xi= H^*$ and $F_\xi=H$, 
    \item $b_\xi=\emptyset$ and
    \item $Y_{\xi+1}=\pi^{-1}_{\xi+1}(Y_\xi)$.
\end{enumerate}
\item Suppose $\M_\xi=(\mathcal{J}^{\vec{E}, f}_{\omega\a}, \in, \vec{E}, f, Y_\xi, \in)$ is a passive ${\sf{hp}}$-indexed $\sf{lses}$, $\M_\xi$ is strategy-ready\footnote{See \rdef{index-ready}.}, $\a=\b+\gg$ and there is  $t\in \univ{\M_\xi|\omega\b}$ such that setting $w=(\mathcal{J}_{\omega}(t), t, \in)$, $w$ is $(f, {\sf{hp}})$-minimal as witnessed by $\b$\footnote{See \rdef{important notation}. In particular, this means that we have to index the branch of $t$ at $\omega\a$.} and $\gg=\lh(t)$. 
Set $b=\Phi_\xi(t)$ and
     \begin{center}
     $\N_{\xi+1}=(\mathcal{J}^{\vec{E}, f^+}_{\omega\b+\omega\gg}, \in, \vec{E}, f, Y_\xi, \tilde{b}, \in)$
     \end{center}
     where $\tilde{b}\subseteq \omega\b+\omega\gg$ is defined by $\omega\b+\omega\nu\in \tilde{b}\iff \nu \in b$. Assuming clause 4 fails for $\xi+1$, the remaining objects 
     are defined as follows. 
         \begin{enumerate} 
    \item $\M_{\xi+1}=\C(\N_{\xi+1})$, 
    \item $F_\xi=F^+_\xi=\emptyset$, 
    \item $b_\xi=\tilde{b}$ and
    \item $Y_{\xi+1}=\pi^{-1}_{\xi+1}(Y_\xi)$. \\
    \end{enumerate}

\textbf{Important Anomaly:} Suppose $\cup Y_\xi$ is $\#$-lsa type\footnote{See \rdef{lsa type}.} and $t$ is $\sf{nuvs}$. Suppose $e\in \M_\xi|\omega\b$ is such that $\M_\xi|\omega\b\models {\sf{sts}_0}(t, e)$\footnote{See \rdef{sts0}. This means that $e$ is the branch of $t$ we must choose.}.
 If $e\not =b$ then $\N_{\xi+1}$ is not an sts premouse over $\mathcal{J}_{\omega}(\cup Y_\xi)$ based on $\cup Y_{\xi}$, and so the construction must stop.\\
\item If $\M_\xi$ doesn't satisfy clause 2a, 2b or 2c then set $\N_{\xi+1}=\mathcal{J}_{\omega}[\M_\xi]$ (this presupposes that $Y^{\N_{\xi+1}}=Y_\xi$). Assuming clause 4 fails for $\xi+1$, the remaining objects 
     are defined as follows.
  \begin{enumerate} 
    \item $\M_{\xi+1}=\C(\N_{\xi+1})$\footnote{Recall that $\C(\M)$ is the core of $\M$.}, 
    \item $F_\xi=F^+_\xi=\emptyset$,
    \item $b_\xi=\emptyset$,
    \end{enumerate}
    and $Y_{\xi+1}=\pi^{-1}_{\xi+1}(Y_\xi)\cup\{\pi^{-1}_{\xi+1}(\M_\xi)$ in the case $\M_{\xi+1}$ is a hod premouse and otherwise, $Y_{\xi+1}=\pi^{-1}_{\xi+1}(Y_\xi)$.
\end{enumerate}
\item Suppose $\xi\leq \d$ is a limit ordinal and for all $\gg<\xi$, both $\M_\gg$ and $\N_\gg$ are defined. Then $\M_{\xi}$ and $\N_{\xi}$ are determined as follows\footnote{The rest of the objects will be defined at the next stage of the induction as in clause 4.}. Set $\N_\xi=lim_{\a\rightarrow \xi}\M_\a$. Assuming clause 4 fails for $\xi+1$, the remaining objects 
     are defined as follows.
     \begin{enumerate}
     \item $\M_\xi=\C(\N_\xi)$ and
     \item $Y_\xi=\pi^{-1}_{\xi}(Y^{\N_\xi})$\footnote{$F_\xi$ and $b_\xi$ are defined at step $\xi+1$.}.
     \end{enumerate}
     \item $\M_\d=\N_\d$ and $Y_\d, \Phi_\d, F^+_\d, F_\d$, and $b_\d$ are undefined.
 \end{enumerate}
We say that the ${\sf{mhpc}}$ is \textbf{successful} if for some $\gg$, $\M_\gg$ is a $\Phi^+$-iterate of $\R^+$. $\myqedhere$
\end{definition}

The following is the main fact we need, which is a corollary to several lemmas established before.

\begin{lemma}\label{mhpc is successful} ${\sf{mhpc}}$ is successful. 
\end{lemma}
\begin{proof}
The lemma follows easily from \rlem{ecertified extenders are on the sequence}, \rlem{lambdav is total} and (b) that appears in the proof of \rlem{authentic iterations are realizable} (which was also established in the proof of \rlem{lambdav is total}).

To prove the lemma, we simply compare $\R^+$ with ${\sf{mhpc}}$-construction of $\N$ and argue that ${\sf{mhpc}}$ side reaches an iterate of $\R^+$. As all extender used in ${\sf{mhpc}}$ with critical point $>{\sf{ord}}(\R)$ have background certificates, the usual stationarity argument shows that such extenders cannot be part of a disagreement in the resulting comparison process. \rlem{ecertified extenders are on the sequence} shows that extenders with critical point $\d^\R$ also cannot be part of a disagreement, while \rlem{lambdav is total} shows that there cannot be a strategy disagreement. Therefore, $\R^+$ iterates to some model appearing on the ${\sf{mhpc}}$-construction.
\end{proof}

\rlem{mhpc is successful} and \rlem{lambdav is total} now imply \rthm{capturing of hod pairs}, and this finishes our proof of \rthm{capturing of hod pairs}.

\subsection{A proof of \rlem{generating gamma}}

In this subsection we outline the proof of \rlem{generating gamma}. The proof is very similar to the proof of \cite[Lemma 6.23]{ATHM}. Suppose that there is no hod pair or an sts hod pair $(\P, \Sigma)$ such that 
\begin{enumerate}
\item $\Sigma$ has strong branch condensation and is strongly fullness preserving,
\item $\Gamma(\P, \Sigma)\subseteq \Gamma\subseteq L(\Sigma, \bR)$
\end{enumerate}   
Just like in the proof of \cite[Lemma 6.23]{ATHM}, it follows from \rthm{the generation of mouse full pointclasses} that $\Gamma$ is not a mouse full pointclass (as we are assuming that $L_\a(\Gamma, \bR)\models \sf{SMC}$). Following the proof of \cite[Lemma 6.23]{ATHM}, we let $A$ be the set of hod pairs or sts hod pairs $(\P, \Sigma)$ such that ${\sf{Code}}(\Sigma)\in \Gamma$ and $\Sigma$ has strong branch condensation and is strongly fullness preserving. It follows from Claim 1 on page 158 of \cite{ATHM} that $A\not = \emptyset$. It follows from Claim 2 on the same page of \cite{ATHM} that if
\begin{center}
$\Gamma_1=\bigcup_{(\P, \Sigma)\in A}\Gamma(\P, \Sigma)$
\end{center}
then\\\\
(1) $\Gamma_1$ is a mouse full pointclass such that for some limit ordinal $\a$ there is a sequence of mouse full pointclasses $(\Gamma_\b:\b<\a)$ such that for $\b<\gg<\a$, $\Gamma_\b\insegeq_{mouse}\Gamma_\gg$ and $\Gamma_1=\bigcup_{\b<\a}\Gamma_\b$.\\\\
It follows from \rthm{the generation of mouse full pointclasses} that there is a possibly anomalous hod pair $(\P, \Sigma)$ such that either 
\begin{enumerate}
\item $\P$ is of lsa type and $\Gamma^b(\P, \Sigma)=\Gamma_1$ or
\item $\P$ is not of lsa type and $\Gamma(\P, \Sigma)=\Gamma_1$. 
\end{enumerate}
Because $\Gamma\models \sf{SMC}$ and because $\Gamma_1\inseg_{mouse}\powerset(\bR)$, we must have that $\Sigma$ is strongly fullness preserving (for instance see \cite[Lemma 6.21]{ATHM}). Notice that even if clause 1.b of \rthm{the generation of mouse full pointclasses} applies, we still get a hod pair as opposed to an sts pair. This is because we have good pointclasses beyond $\Gamma$. 

Notice also that ${\sf{Code}}(\Sigma)\not \in \Gamma$, as otherwise it follows from Claim 2 on page 158 of \cite{ATHM} that $(\P, \Sigma)\in A$. Thus, it must be the case that $\P$ is an anomalous hod premouse. We now get a contradiction as in page 159 of \cite{ATHM}, where it is argued that the computation of $\H^{L(\Sigma, \bR)}$ gives a contradiction.

\section{A proof of \sf{LSA} from large cardinals}\label{sec: lsa from lc}

In this section, we generalize \cite[Theorem 6.26]{ATHM}.
 
 \begin{theorem}\label{lst from wlw} The theory  $\sf{AD}^++\sf{LSA}+V=L(\powerset(\bR))$ is consistent relative to a Woodin cardinal that is a limit of Woodin cardinals. 
 \end{theorem}
 \begin{proof} Woodin showed that it is consistent relative to a Woodin cardinal that is a limit of Woodin cardinals that there are  divergent models of $\sf{AD}^+$, i.e., there are sets of reals $A, B\subseteq \mathbb{R}$ such that $L(A, \bR)\models \sf{AD}^+$, $L(B, \bR)\models \sf{AD}^+$, $A\not \in L(B, \bR)$ and $B\not \in L(A, \bR)$. Moreover, his construction shows that we can assume that both $L(A, \bR)$ and $L(B, \bR)$ satisfy $\sf{MC}$$+\Theta=\theta_0+{\sf{NWLW}}$\footnote{See \rdef{lsa nwlw}. The proof of Woodin's theorem appeared in \cite{EA} as Theorem 6.1 and it is obtained as a forcing extension of the minimal active mouse with a Woodin cardinal that is a limit of Woodin cardinal. If we also would like to have ${\sf{NWLW}}$ then we simply take the $\in$-minimal model in which divergent models exists.}. Thus, we assume that such a pair of models exists. 
 
Suppose towards a contradiction that there is no inner model satisfying $\sf{AD}^++\sf{LSA}+V=L(\powerset(\bR))$. Let $\Gamma=L(A, \bR)\cap L(B, \bR)\cap \powerset(\bR)$.  It is an unpublished theorem of Woodin that $L(\Gamma, \bR)\models \sf{AD}_\bR$ but see \cite[Theorem 8.1]{Wenvelope}. We also have that $\Gamma=\powerset(\bR)\cap L(\Gamma, \bR)$. Applying \rlem{the generation of mouse full pointclasses} in $L(A, \bR)$ and in $L(B, \bR)$ we get two hod pairs or sts hod pairs $(\P, \Sigma)\in L(A, \bR)$ and $(\Q, \Lambda)\in L(B, \bR)$ such that both $\P$ and $\Q$ are of limit type and $\Gamma=\Gamma(\P, \Sigma)=\Gamma(\Q, \Lambda)$. 

Working in $L(A, \bR)$, let $\M^*=\bigcup_{(\S, \Psi)\in B(\P, \Sigma)}\M_\infty(\S, \Psi)$ and for $\S'\inseg^c_{hod}\M^*$\footnote{See \rdef{complete layer notation}.} let $\Psi_{\S'}$ be the iteration strategy of $\S'$ obtained from any $(\S, \Psi)$ such that $\S'=\M_\infty(\S, \Psi)$. Notice that $\M^*$ and the strategies $(\Psi_{\S'}: \S'\inseg^c_{hod}\M^*)$ are independent of $(\P, \Sigma)$ i.e. working in $L(B, \bR)$ and using $(\Q, \Lambda)$ instead of $(\P, \Sigma)$ yields the same model $\M^*$ and the same strategies $(\Psi_{\S'}: \S'\inseg^c_{hod}\M^*)$. Let\footnote{$\W\insegeq {\sf{Lp}}^{\Gamma, \oplus_{\S\inseg^c_{hod}\M^*}\Psi_\S}(\M^*)$ if and only if  $\W$ is a sound $\oplus_{\S\inseg^c_{hod}\M^*}\Psi_\S$-premouse over $\M^*$ such that $\rho(\W)\leq {\sf{ord}}(\M^*)$ and whenever $\pi:\W'\rightarrow \W$ is elementary and $\W'$ is countable, $\W'\insegeq {\sf{Lp}}^{\Gamma, \Psi'}(\pi^{-1}(\M^*))$ where $\Psi'$ is the $\pi$-pullback of $\oplus_{\S\inseg^c_{hod}\M^*}\Psi_\S$.}
\begin{center}
$\M_A=({\sf{Lp}}^{\Gamma, \oplus_{\S\inseg^c_{hod}\M^*}\Psi_\S}(\M^*))^{L(A, \bR)}$  and $\M_B=({\sf{Lp}}^{\Gamma,  \oplus_{\S\inseg^c_{hod}\M^*}\Psi_\S}(\M^*))^{L(B, \bR)}$.
\end{center}
We then have that either $\M_A\insegeq \M_B$ or $\M_B\insegeq \M_A$. Without loss of generality we assume that $\M_A\insegeq \M_B$. 
 
Let $\pi:\P^b\rightarrow \M_A$ be the iteration embedding given by $\Sigma$. It follows from the proof of Claim 7 appearing in the proof of \rthm{lsa from min omega woodins over lsa} that $\Sigma\in L(\pi[\P], \M_A, \Gamma)$. By our assumption, $\M_A\in L(B,\mathbb{R})$. Because $\pi[\P]$ is a countable set we have that $\pi[\P]\in L(B, \bR)$. It follows that $\Sigma\in L(B, \bR)$. Therefore, ${\sf{Code}}(\Sigma)\in \Gamma$ implying that $\Gamma(\P, \Sigma)\subset \Gamma$, contradiction! 
 \end{proof}
 
We remark that just like in the proof of \cite[Theorem 6.26]{ATHM}, we could have used \rthm{diamond comparison} instead of \rthm{lsa from min omega woodins over lsa}. 

\chapter{A proof of square in lsa-small hod mice}
\label{chap:square}

\begin{definition}\label{def:square}\index{$\square_\kappa$}
For a cardinal $\kappa$ and a cardinal $\gamma\leq \kappa$, the principle $\square_{\kappa,\gamma}$ states that there is a sequence $\langle \mathcal{C}_\alpha : \alpha<\kappa^+\rangle$ such that for each $\alpha<\kappa^+$
\begin{enumerate}
\item $\mathcal{C}_\alpha \neq \emptyset$ and for each $C\in \mathcal{C}_\alpha$, $C$ is a closed unbounded subset of $\alpha$ of order type at most $\kappa$,
\item $|\mathcal{C}_\alpha|\leq \gamma$,
\item for each $C\in\mathcal{C}_\alpha$, for each $\beta \in \lim(C_\alpha)$, $ C\cap \beta\in \mathcal{C}_\beta$.
\end{enumerate}
If $\gamma=1$, then the principle $\square_{\kappa,\gamma}$ is simply $\square_\kappa$.$\myqedhere$
\end{definition}

Pure extender models are models constructed from a canonical sequence of extenders. Jensen (cf. \cite{Jensen}) initiated the program of understanding square principles in pure extender models by proving $L\vDash \forall \kappa \ \square_\kappa$. Building on works of several people, Schimmerling and Zeman (cf. \cite{schimmerling2004characterization}) give the most optimal characterization of $\square$ in (short) extender models, namely they prove that in an iterable, short extender model, $\square_\kappa$ holds if and only if $\kappa$ is not subcompact. Results on squares in extender models are important in understanding structure theory of such models and have found many applications in set theory. The reader can see, for instance, \cite{schimmerling2007coherent} and \cite{JSSS}, for some of the applications of square in extender models in computing lower-bound consistency strength of theories like $\sf{PFA}$. Recent advances in the core model induction methods have indicated that to improve the lower-bounds of combinatorial principles like $\sf{PFA}$, failure of square at a singular cardinal, the existence of guessing models etc., one way is to prove square holds in the hod mice that are currently being studied and constructed.

All known square proofs in extender models rely heavily on the fine-structure of such models, in particular, they make essential use of condensation properties of these models (cf. \cite[Lemma 1.6]{schimmerling2004characterization}). Unfortunately, the full condensation lemma,\cite[Lemma 1.6]{schimmerling2004characterization}, does not hold in hod mice. However, it is possible to overcome this shortcoming. We present here a proof of $\square_{\kappa,2}$ in an \textit{lsa-small hod mouse} $\P$ for all cardinals $\kappa$ of $\P$. In this chapter by lsa-small hod mouse, we mean that $\P$ does not contain an active $\omega$ Woodin lsa mouse as defined in Definition \ref{omega woodins over lsa}.

We first set up some terminology.  Our hod premice $\P$ are \textit{lsa-small} and hence for no $\alpha<\lambda^\P$, $\P(\alpha)$ is an lsa hod premouse, though $\P$ can be of lsa type. Throughout this paper, if $\Q$ is an initial segment of $\P$, we let $\Sigma_\Q$ denote the restriction of $\Sigma$ to $\Q$. If $\P$ is of limit type and has a top window $[\delta_\alpha^\P,\delta_{\alpha+1}^\P)$, then we let $\P^b = \P|({{\delta^\P_\alpha}^+})^\P$.  See Section \ref{sec:hod_mice} for a more detailed discussion of hod mice along with the definitions used in statements of this section. In the definitions below, we adapt the $\Sigma^*$-language (see \cite{schimmerling2004characterization}) to hod mice in the obvious way. Let $\rho^n_\Q$\index{$\rho^n_\Q$} be the $n^{th}$-projectum of $\Q$, and $p^n_\Q$\index{$p^n_\Q$} be the $n^{th}$-standard parameter of $\Q$.\footnote{Other notations for the $n^{th}$-projectum and $n^{th}$-standard parameter of $\Q$ used elsewhere in this book are $\rho_n(\Q)$ and $p_n(\Q)$ respectively. For this chapter, we stick to the more compact notations $\rho^n_\Q$ and $p^n_\Q$; this notation is compatible with the notation used in \cite{schimmerling2004characterization}.} Semantically, suppose $\Q$ is an initial segment of $\P$, a relation $A\subset |\Q|$ is $\Sigma^{(n)}_l(\Q)$ from $p$, or $\mathbf{\Sigma}^{(n)}_l(\Q)$, if it is $\Sigma_l$ from $p$ (or $\mathbf{\Sigma}_l$) over the $n^{th}$-reduct $\langle H^n_\Q,A^{n}_\Q \rangle$ of $\Q$, where $H^n_\Q = |\Q|\rho^n_{\Q}|$ and $A^n_\Q$ is the $n^{th}$ standard master code (with respect to $p^n_\Q$) of $\Q$. 
\begin{definition}\label{def:strong_fullness_preserving}
Suppose $\Sigma$ is an iteration strategy for a hod premouse $\P$. Suppose $\Gamma$ is an inductive-like pointclass. We say that $\Sigma$ is \textit{locally strongly $\Gamma$-fullness preserving}\index{locally strongly $\Gamma$-fullness preserving} if $\Sigma$ is $\Gamma$-fullness preserving and if $\P$ is of limit type with a top window and whenever $(\vec{\T},\S)\in I(\P,\Sigma)$, and 
\begin{center}
$\pi^{\vec{\T},b}: \P^b \rightarrow \S^b$ exists, 
\end{center}
then letting $\pi = \pi^{\vec{\T},b}$, whenever $\S^b \lhd \W \unlhd \S$ is such that for some $n$ 
\begin{center}
$o(\S^b)\leq \omega\rho^{n+1}_\W <  \omega\rho^n_\W$,
\end{center}
$\W$ is $n$-sound, and $\tau:\R\rightarrow \W$ is cardinal preserving and $\Sigma_0^{(n)}$ and $\omega\rho^n_{\R} > $ cr$(\tau) \geq \omega\rho^{n+1}_\R=\omega\rho^{n+1}_\W$, then the $\tau$-pullback of the strategy $\Sigma_{\W,\vec{\T}}$ is $\Gamma$-fullness preserving in the following sense: whenever $\vec{\U}$ is according to $\Sigma_{\W,\vec{\T}}^\tau$, then letting $\R^*$ be the last model of $\vec{\U}$, $(\R^*)^b$ is $\Gamma$-full.$\myqedhere$
\end{definition}

\begin{definition}\label{def:strong_branch_condensation}
Suppose $\Sigma$ is an iteration strategy for a hod premouse $\P$. We say that $\Sigma$ has \textit{locally strong branch condensation}\index{locally strong branch condensation} if $\Sigma$ has branch condensation and if $\P$ is of limit type with a top window and $\Q$ is such that $\P^b \lhd \Q \unlhd \P$, and $n$ is such that $\omega\rho^{n+1}_\Q < \omega\rho^n_\Q$, $\Q$ is $n$-sound, $\omega\rho^{n+1}_\Q$ is a cardinal of $\P$, and $\S$ is a $\Sigma_\Q$-iterate along a stack $\vec{\T}$ such that $\pi^{\vec{\T},b}$ exists, and  $\tau:\Q\rightarrow \R$ is a cardinal preserving, $\Sigma^{(n)}_0$-embedding such that $\R \unlhd \S$ and $(\Q^*)^b= \R^b$ for some non-dropping $\Sigma_\Q$-iterate $\Q^*$ of $\Q$. Suppose also that letting $j:\Q\rightarrow \Q^*$ be the iteration map, then $j\rest \Q^b = \tau\rest \Q^b$. Then $\Sigma_{\R,\vec{\T}}^{\tau} = \Sigma_\Q$.$\myqedhere$
\end{definition}

The proof of Theorem \ref{thm:fullness_preserving} can be modified to get hod mice $(\P,\Sigma)$ with $\Sigma$ being locally strongly $\Gamma$-fullness preserving. Similarly the proof of strong branch condensation can be modified to obtain strategies with locally strong branch condensation. We leave the easy proofs to the reader. The term ``locally" refers to the fact that these forms of fullness preservation and branch condensation apply to initial segments of $\P$ (like $\W$ and $\Q$ in these definitions). In the proof of Theorem \ref{thm:square}, initial segments like $\Q\lhd \P$ are typically collapsing structures of some $\tau\in (\kappa,(\kappa^+)^\P)$, where $\kappa = \rho_Q^{n+1}$ is a cardinal of $\P$. We seem to need to modify the usual notions of fullness preservation and branch condensation (as in Definitions \ref{def:strong_fullness_preserving} and \ref{def:strong_branch_condensation}) to ensure that various phalanx comparison arguments involving initial segments like $\Q$ (which is not $\P$ in most cases of interest) go through in the proof of Theorem \ref{thm:square}. In most (but not all) applications, the map $\pi$ in Definitions \ref{def:strong_fullness_preserving} and \ref{def:strong_branch_condensation}, is the identity and $\tau$ is the uncollapse map associated to a sufficiently elementary hull. The main theorem is the following.
\begin{theorem}\label{thm:square}
Suppose $(\P,\Sigma)$ is an lsa-small hod pair such that $\Sigma$ has locally strong branch condensation and is locally strongly $\Gamma$-fullness preserving for some inductive-like pointclass $\Gamma$ that satisfies ``$\sf{AD}^+ + \sf{SMC}$". Then $\P \vDash \forall \kappa \ \square_{\kappa,2}$.\footnote{The assumption that $\P$ is lsa-small implies that there are no subcompact cardinals in $\P$ and all extenders on the $\P$-sequence are short.}
\end{theorem} 	

Many techniques in the proof of \ref{thm:square} come from the Schimmerling-Zeman's proof in \cite{schimmerling2004characterization}. In Section \ref{sec:hod_mice}, we import some results from the theory of hod mice we need. In Section \ref{sec:S_Z}, we will import some terminology, results from \cite{schimmerling2004characterization} that we need here. We also explain in this section why a straightforward adaptation of \cite{schimmerling2004characterization} fails in the context of hod mice. In Section \ref{sec:proof}, we give the actual proof of Theorem \ref{thm:square}.

Finally, we remark that hod pairs constructed in practice (those constructed in sufficiently strong $\sf{AD}^+$ models or in the core model induction settings) do have the properties in the hypothesis of Theorem \ref{thm:square}. The main application of Theorem \ref{thm:square} in this book is to improve the lower-bound consistency strength of various theories such as $\sf{PFA}$ to that of $\sf{LSA}$ (see Chapter \ref{chap:lsa_from_pfa}).

\section{Ingredients from hod mice theory}\label{sec:hod_mice}

We summarize some definitions and results of the hod mice theory developed above that we need to prove Theorem \ref{thm:square}. The language for hod premice's is $\mathcal{L}_1$ with symbols
 $\in, \dot{Y}, \dot{E},\dot{F},\dot{\Sigma},\dot{B}, \dot{\gamma}$. 
 
Suppose $(\P,\Sigma)$ is an lsa-small hod pair. $\P$ is constructible from a sequence of extenders and a sequence of strategies of its own hod initial segments. As before $\dot{Y}^\P$ codes the layers of $\P$.  In the following, all extenders in $\dot{E}^\P$ and $\{\dot{F}^\P\}$ will be indexed on the $\P$-sequence  according to the $\lambda$-indexing scheme  described in \cite{Zeman}. $\dot{\gamma}^\P$, just as in \cite{Zeman}, codes where $\dot{F}^\P$ restricted to the largest cutpoint of $\dot{F}^\P$ would be indexed. There are two ways in which an initial segment $\Q$ of $\P$ can be active: $B$-active and $E$-active. $\Q$ is \textit{$B$-active} the top predicate $\dot{B}^\Q$ for $\Q$ (amenably) codes a branch for some tree on an initial segment of $\Q$. $\Q$ is \textit{$E$-active} if the top predicate $\dot{F}^\Q$ of $\Q$ codes an extender. Otherwise, we say that $\Q$ is \textit{passive}. $B$-active levels and passive levels are more or less treated the same way in the proof of Theorem \ref{thm:square}. In our situation, we note that since all initial segments of $\P$ and $\P$ itself are lsa-small, all extenders on the $\P$-sequence are of type $A$, i.e. they have no cutpoints. So $\dot{\gamma}^\P = \emptyset$. This aspect somewhat simplifies our proof, compared to \cite{schimmerling2004characterization}.

A few words about how the $B$-predicate codes up branches for an iteration tree $\T$ in $\P$ are in order. Suppose $\lambda = \lh(\T)$ is limit and $\P|\gamma$ is $B$-active such that $B^{\P|\gamma}$ codes a cofinal branch $b$ of $\T$. The traditional way that $B$ codes $b$ is that letting $\gamma^*+\lambda=\gamma$, $B^{\P|\gamma}=\{\gamma^* + \alpha \ | \ \alpha\in b\}$. While this approach is sufficient for developing the basic theory of strategic premice and certainly is sufficient for the theory of hod mice we have developed so far, it seems to create significant obstructions in the proof of $\square$ in this chapter. So instead, we use the coding method developed in \cite{trang2013}. Using \cite[Definition 2.26]{trang2013}, we let $\P|\gamma=\mathfrak{B}(\P|\gamma^*,\T,b)$. The reader is advised to consult \cite{trang2013} for the precise definition of $\mathfrak{B}(\P|\gamma^*,\T,b)$. Roughly, for every $0 < \alpha < \lambda$, $\P|(\gamma^*+\omega\alpha)$ is $B$-active and $B^{\P|(\gamma^*+\omega\alpha)}$ codes the branch $[0,\alpha)_\T$ and $B^{\P|\gamma}$ codes $b$ in the manner described above. We make this a bit more precisely here. 

$\dot{\Sigma}^\P$ and $\dot{B}^\P$ are used to record information about an iteration strategy
 $\Omega$ of $\P$. $\dot{\Sigma}^\P$ codes the strategy information added at
  earlier stages; $\dot{\Sigma}^\P$ acts on stacks in $\P$ based on a hod layer $\Q\in \dot{Y}^\P$. $\dot{\Sigma}^\P(s,b)$ implies 
  that $s= \mathcal{T} $,
   where $\mathcal{T}$ is a stack on
    $\Q$ in $\P$ of limit length and $\mathcal{T}^\smallfrown b$ is according to the strategy.
 We say that $s$ is an $\P$-tree, and write $s =\mathcal{T}(s)$. We 
  write $\dot{\Sigma}^\P_{\nu,k}$ for the partial iteration strategy for $\P|(\nu,k)$
   determined by $\dot{\Sigma}$ (here $\Q = \P|(\nu,k)$ for some $\nu, k$). We write $\Sigma^\P(s)=b$ when $\dot{\Sigma}^\P(s,b)$, and
    we say that $s$ is according to  $\Sigma^\P$ if $\T(s)$ is according to $\dot{\Sigma}^\P_{\Q}$. 

We say $\P$ is \textit{branch-active} (or just $B$-active) iff 
\begin{itemize}
\item[(a)] there is a largest $\eta < o(\P)  $ such that $\P|\eta \models {\sf KP}$,
and letting $N=\P|\eta$,
\item[(b)] there is a
 $<_N$-least $N$-tree $s$ such that $s$ is by $\Sigma^N$,
		$\T(s)$ has limit length, and $\Sigma^N(s)$ is undefined.
		\item[(c)] for $N$ and $s$ as above, $o(\P) \le
			o(N) + lh(\T(s))$.
\end{itemize}

Note that being branch-active can be expressed by a $\Sigma_2$ sentence in
$\mathcal{L}_1 -\lbrace \dot{B} \rbrace$. This contrasts with being extender-active,
which is not a property of the premouse with its top extender removed. In contrast
with extenders, we know when branches must be added before we do so.

\begin{definition}\label{branchactivelpm}
Suppose that $\P$ is branch-active. We set
\begin{align*}
	\eta^P &= \text{the largest $\eta$ such that $\P|\eta \models {\sf KP}$,}\\
	\nu^\P &= \text{unique $\nu$ such that $\eta^\P + \nu = o(\P)$},\\
	s^\P &= \text{least $\P|\eta^M$-tree such that $\dot{\Sigma}^{\P|\eta^\P}$ is undefined, and}\\ 
	b^\P &= \text{$\lbrace \alpha \mid \eta + \alpha \in \dot{B}^\P \rbrace.$}	
\end{align*}
	Moreover,
	\begin{itemize}
	\item[(1)]
$\P$ is a {\em potential hod $B$-active premouse} iff $b^M$ is a cofinal
	branch of $\T(s) \rest \nu^\P$.
	\item[(2)] $\P$ is {\em honest}
	iff $\nu^\P = \lh(\T(s))$, or $\nu^\P < lh(\T(s))$ and
	$b^\P = [0,\nu^\P)_{T(s)}$.
	\item[(3)] $\P$ is a hod premice iff $\P$ is an honest potential lpm.
	\item[(4)] $\P$ is {\em strategy active} iff $\nu^\P = \textrm{lh}(\T(s))$.
	\end{itemize}$\myqedhere$
	\end{definition}

Note that $\eta^\P$ is a
$\Sigma_0^\P$ singleton, because it is the least ordinal in $\dot{B}^\P$ (because
0 is in every branch of every iteration tree), and thus $s^M$ is also a
$\Sigma_0^\P$ singleton. We have separated honesty from the other conditions
because it is not expressible by a $Q$-sentence, whereas the rest is.
Honesty is expressible by a Boolean combination of $\Sigma_2$ sentences.

The definition of $B$-active hod premice defined in previous chapters 
required that when $o(\P) < \eta^\P + lh(\T(s))$, $\dot{B}^\P$ is empty,
 whereas here we require that it code 
$[0,o(\P))_{T(s)}$, in the same way that $\dot{B}^\P$ will have to code a new
branch when $o(\P) = \eta^M + lh(\T(s))$. Of course, $[0,\nu^\P)_{T(s)} \in \P$ 
when $o(\P) < \eta^\P + lh(\T(s))$ and $\P$ is honest, so the current
$\dot{B}^\P$ seems equivalent to the original $\dot{B}^\P = \emptyset$.
However, $\dot{B}^\P=\emptyset$ leads to $\Sigma_1^\P$ being too weak,
with the consequence that a $\Sigma_1$ hull of $\P$ might collapse to
something that is not a hod premouse.\footnote{The hull could satisfy
$o(H) = \eta^H + lh(\T(s^H))$, even though
$o(\P) < \eta^\P + lh(\T(s^\P))$. But then being a hod premouse
requires $\dot{B}^H \neq \emptyset$.} Our current choice for
$\dot{B}^\P$ solves that problem. Furthermore, the indexing of branches for $\P$-trees $\T$ using the $\mathfrak{B}$-operator is done for all eligible $\T$ regardless of whether cof$(\textrm{lh}(\T))$ is measurable in $\T$ as we required using the original definition. This makes the definition more uniform.

\begin{remark} Suppose $N$ is an hod premouse, and $N \models {\sf KP}$.
It is very easy to see that $\dot{\Sigma}^N$ is defined on all $N$-trees $s$
 that are by
$\dot{\Sigma}^N$ iff there are arbitrarily large
$\xi < o(N)$ such that $N|\xi \models {\sf KP}$.
Thus if $M$ is branch-active, then $\eta^M$ is a successor
admissible; moreover, we do add branch information, related to
exactly one tree, at each successor
admissible. Waiting until the next admissible
to add branch information is just a convenient way to make 
sure we are done coding in the branch information for a
given tree before we move on to the next one. One could go faster.$\myqedhere$
\end{remark}

As mentioned above and discussed in more details in \cite{trang2013}, we can prove stronger forms of condensation for hod mice where the $\mathfrak{B}$-operator is used to index branches. For instance, we have the following condensation lemma for hod pairs.

\begin{lemma}\label{lem:preservehodpairs} Let $(M,\Omega)$ be a hod pair\footnote{Recall $\Omega$ has hull condensation.} with $k(M) = k$\footnote{$k(M)$ is the largest $k$ such that $M$ is $k$-sound.}, and let $\pi \colon H \to M$ be $\Sigma^{(k)}_0$ and cardinal preserving.\footnote{This is the corresponding version of what is called ``weak embeddings" in \cite{trang2013}.} Suppose that one of the following holds:
\begin{itemize}
\item[(a)] $M$ is passive or branch-active, or
\item[(b)] $H$ is a hod premouse.
\end{itemize}
Then $(H,\Omega^\pi)$ is a hod pair.
\end{lemma}

\begin{proof} We show first that $H$ is a hod premouse. If (b) holds,
this is rather easy (in fact, a tautology). If $M$ is passive, then so is $H$, noting that $Q$ sentences are preserved downwards
under $\pi$ and being a passive hod premouse can be expressed by a $Q$ sentence. So let us assume that
$M$ is branch-active.

Note that $H$ is a potential branch active
hod premouse by the same reasoning as above: the property of being a potential branch active hod premouse is expressible by a $Q$ sentence and this is preserved downwards by $\pi$. So we just need to see that $H$ is honest.
Let $\nu = \nu^H$, $b = b^H$, and $\T = \T(s^H)$. If $\nu = lh(\T)$,
there is nothing to show,
so assume $\nu < lh(\T)$. We must show
that $b = [0,\nu)_{T}$. We have by induction that for $N=H|\eta^H$,
$(N,\Omega^\pi_N )$ is a hod pair. Thus $\T$ is by $\Omega^\pi$, and so we just need to see
that for $\U = \T \rest \nu$ and $\mathcal{W} =\U ^\frown b$,
$\mathcal{W}$ is by $\Omega^\pi$, or equivalently,
that $\pi\mathcal{W}$ is by $\Omega$. But it is easy to see that
$\pi\mathcal{W}$ is a hull of $\pi(U)^\frown b^M$, and $\Omega$
has hull condensation, so we are done.

\end{proof}
The use of the $\mathfrak{B}$-operator in coding branches of iteration trees in the $\square$-construction will be explained in Section \ref{sec:proof} (see, for example Lemma \ref{lem:0_unbounded}, where Lemma \ref{lem:preservehodpairs} is used).

We briefly discuss indexing schemes for extenders on the $\P$-sequence. We recall the mixed indexing scheme for hod premice being used the the previous chapters. Suppose $\kappa$ is a cardinal limit of cutpoint Woodin cardinals of $\P$, and if $E$ is an extender on the $\P$'s sequence such that cr$(E)=\kappa$, then the index of $E$ is $\gamma$ where $\gamma$ is the successor cardinal of the least cutpoint above $\kappa$ in Ult$(\P|\xi,E)$ (we call this \textit{cutpoint indexing scheme}), where $\xi\leq o(\P)$ is the largest such that $E$ measures all sets in $\P|\xi$. It turns out that such extenders are all total over $\P$ (see Chapters 2, 3 for more discussions). Suppose $E$ is an extender with critical point $\xi$ and $E$ is indexed according to the cutpoint indexing scheme. Then according to \cite{steel2013lst}, for all $\gamma < \rm{lh}$$(E)$, $E\rest\gamma$ is not on the $\P$-sequence, though $E\rest\gamma \in \P$ (for $\gamma$ below the sup of the generators of $E$) and the trivial completion of $E\rest\gamma$ is on the $\P$-sequence for various $\gamma$ (this is similar to the initial segment condition for Jensen indexing). Also, the set of indices of extenders with a fixed critical point $\xi$ indexed according to the cutpoint indexing scheme is nowhere continuous. For other extenders on $\P$'s sequence, we use \textit{the Mitchell-Steel indexing scheme} (ms-indexing). 

In the following, all extenders on a hod premouse's extender sequence will be indexed according to the \textit{Jensen indexing scheme} ($\lambda$-indexing). By results of \cite{Fuchs1,Fuchs2}, one can translate a hod pair $(\P,\Sigma)$ in the mixed indexing scheme described above to a hod pair $(\P^*,\Sigma^*)$ in the $\lambda$-indexing scheme. The hod pair $(\P^*,\Sigma^*)$ is obtained by applying the $\Lambda$-function to $(\P,\Sigma)$. See for example, \cite[Lemmata 3.4, 6.3, 6.4]{Fuchs2}.\footnote{Technically speaking, Fuchs shows the intertranslatability for ms-mice and $\lambda$-mice. But the same proof techniques can be used without virtually any change to translated mice with mixed-indexing to $\lambda$-mice.} As in \cite{Fuchs2}, one can show that the universe of $\P$ and $\P^*$ are the same and much more. $\Sigma^*$ is essentially the same as $\Sigma$, so it has all the properties $\Sigma$ has. The change in indexing schemes is for convenience of importing terminology and results from \cite{schimmerling2004characterization}. The result we prove here for hod mice with the Jensen indexing scheme will also hold for hod mice with the Mitchell-Steel indexing scheme by \cite{Fuchs2} and the above discussion, i.e. if $\kappa$ is a cardinal of $\P^*$ (equivalently of $\P$) then $$\P^*\models \square_{\kappa,2} \textrm{ implies } \P\models \square_{\kappa,2}.$$ Suppose $E$ is an extender with critical point $\xi$ on the sequence of $\P$ and $E$ is indexed by the Jensen indexing scheme, that is the index of $E$ in $\P$ is the successor cardinal of $i_E(\xi)$ in Ult$(\P,E)$. For a summary of the fine structure, see \cite[Section 1]{schimmerling2004characterization}. A couple of remarks regarding the adaptation of \cite[Section 1]{schimmerling2004characterization} into our situation are in order. First, we still demand extenders indexed according to the Jensen indexing to satisfy the \textit{initial segment condition} (ISC) in the sense of \cite[Section 1.4]{schimmerling2004characterization}; that is for all $\gamma < \rm{lh}$$(E)$, if $\gamma$ is a cutpoint of $E$, then $E\rest \gamma \in |\P|\rm{lh}$$(E)|$. Secondly, under this initial segment condition, using the assumption that our hod premice are lsa-small, it's easy to see that these extenders $E$ are all of type $A$, that is the set of cutpoints is empty; this is because there are no superstrong cardinals in lsa-small hod mice. The initial segment condition is needed to prove comparisons terminate.

We recall some concepts related to layers discussed in previous chapters. If $\P$ is a hod premouse, we let $\lambda^\P$ denote the order type of the layer Woodin cardinals of $\P$ and $(\delta^\P_\alpha : \alpha < \lambda^\P)$ enumerate the closure of the set of Woodin cardinals and indices of extenders whose critical point is a limit of layer Woodin cardinals of $\P$. Recall that $\delta$ is a layer Woodin cardinal of $\P$ if there is some $\Q\in Y^\P$ (i.e. $\Q$ is a layer of $\P$) such that $\delta = \delta^\Q$. Intuitively, $Y^\P$ is the set of initial segments $\Q$ of $\P$ such that the strategy (or sts strategy) of $\Q$ is activated in $\P$.  See \ref{layers of hod-like lsp}. If $\P$ has a largest Woodin cardinal, we denote that $\delta^\P$. Recall we use $\P^b$ to denote the ``bottom part" of $\P$ in the case that $\P$ has a top window $[\kappa=\delta_\alpha^\P,\delta^\P)$, where recall that by definition, $\kappa$ is either a Woodin or limit of Woodins in $\P$. In this case, $\P^b = \rm{Lp}^{\Sigma^\P_\kappa,\P}(\P|\kappa)$, where $\Sigma^\P_\kappa= \oplus_{\beta<\alpha} \Sigma^\P_{\P(\beta)}$. In the case $\alpha$ is a limit ordinal, $\P^b = \P|((\kappa)^+)^\P$. In this case, if $\kappa$ happens to be measurable in $\P$, then all extenders $E$ on the $\P$ sequence with critical point $\kappa$ are indexed according to the cutpoint indexing scheme. Notice that since $\P$ is lsa-small, $\kappa$ is a cutpoint (but not a strong cutpoint) in $\P$, though $\kappa$ is a strong cutpoint in $\P^b = \P|(\kappa^+)^\P$. Let $o(\kappa)$ be the supremum of the indices of extenders on the $\P$ sequence with critical point $\kappa$. If $\P$ is of \textit{lsa type} (i.e. $o(\kappa)=\delta^\P$) then there may be local large cardinals in the interval $(\kappa,o(\kappa))$, e.g. there may be a $\gamma\in (\kappa,o(\kappa))$ which is Woodin in some initial segment $\Q$ of $\P$; such large cardinals are witnessed by the extender sequence and \textit{the short tree strategy} of initial segments of $\Q$, but not the full strategy. This point is crucial in many arguments given below (see Lemma \ref{lem:no_strat_disagreement}).

Suppose $(\P,\Sigma)$ is a hod pair such that $\Sigma$ is $\Gamma$-fullness preserving for some inductive-like pointclass $\Gamma$ and has branch condensation. Suppose $\R\lhd \P$ is an initial segment of $\P$, then we let $\Sigma_\R$ denote the restriction of $\Sigma$ to trees based on $\R$. Let $I(\P,\Sigma)$ denote the set of $(\vec{\T},\R)$ where $\vec{\T}$ is a stack according to $\Sigma$ with last model $\R$. In this case, the ``$\vec{\T}$-tail" of $\Sigma$, denoted $\Sigma_{\vec{\T},\R}$, is a strategy for $\R$. We let $B(\P,\Sigma)$ denote the set of $(\vec{\T},\R)$ where $\vec{\T}$ is according to $\Sigma$ and $\R$ is a strict hod-initial segment of $\N^{\vec{\T}}$, the last model of $\vec{\T}$. We let $\Gamma(\P,\Sigma)$ be the set of $A\subseteq \mathbb{R}$ such that $A <_w \Sigma_{\vec{\T},\R}$ for some $(\vec{\T},\R)\in B(\P,\Sigma)$. Note that $\Gamma(\P,\Sigma)$ is a Wadge initial segment of $\Gamma$. We say that \textit{$\P$ generates $\Omega$} if $\Gamma(\P,\Sigma) = \Omega$.

The following fact will be used in many places throughout this chapter, and whose proof is essentially that of \ref{no strategy disagreement}. 
\begin{lemma}[No strategy disagreement]\label{lem:no_strat_disagreement}
Suppose $(\P,\Sigma)$ is an lsa-small hod pair such that $\P$ has a top window $[\delta^\P_\alpha,\delta^\P)$ and $\delta^\P_\alpha$ is not a strong cutpoint of $\P$, $\Sigma$ has locally strong branch condensation and is locally strongly $\Gamma$-fullness preserving for some constructibly closed pointclass $\Gamma\vDash ``\sf{AD}^+ + \sf{SMC}$". Suppose $\pi:\P'\rightarrow \P^*$ for some cardinal preserving, $\Sigma^{(n)}_0$ map $\pi$ such that $\P^b\lhd\P^*\unlhd \P$, and $\omega\rho^n_{\P^*} >$ cr$(\pi) =_{\rm{def}} \gamma > \omega\rho^{n+1}_{\P'} = \omega\rho^{n+1}_{\P^*}\geq o(\P^b)$ and $\rho^{n+1}_{\P^*}$ is a cardinal of $\P$. Then letting $\Lambda = \Sigma_{\P^*}^\pi$, the comparison of the phalanx $(\P^*,\P',\gamma)$ (using $\Lambda$) versus $\P^*$ (using $\Sigma_{\P^*}$) does not involve disagreements of strategies.
\end{lemma} 

Lemma \ref{lem:no_strat_disagreement} is useful since it reduces such comparisons to ordinary extender comparisons. Such phalanx comparisons will appear in many places in the proof of Theorem \ref{thm:square}. A corollary of this is the following version of the Condensation Lemma for hod mice (cf. \cite[Lemma 9.3.2]{Zeman}). For notations used in the statement of the lemma, see \cite[Section 1.3]{schimmerling2004characterization}.
\begin{theorem}\label{thm:condensation_lemma}
Suppose $(\P,\Sigma)$ is an lsa-small hod pair such that $\P$ has a top window $[\delta^\P_\alpha,\delta^\P)$ and $\delta^\P_\alpha$ is not a strong cutpoint of $\P$, $\Sigma$ has locally strong branch condensation and is locally strongly $\Gamma$-fullness preserving for some constructibly closed pointclass $\Gamma\vDash ``\sf{AD}^+ + \sf{SMC}$". Suppose $\P^b \lhd \M \unlhd \P$, $\tilde{\M}$ is a hod premouse, and $\sigma:\tilde{\M}\rightarrow \M$ is a cardinal preserving and $\Sigma_0^{(n)}$ embedding such that $\sigma\rest \omega\rho^{n+1}_{\bar{\M}} = \rm{id}$, where $\omega\rho^{n+1}_{\tilde{\M}} = \omega\rho^{n+1}_\M \geq o(\P^b)$ is a cardinal of $\P$.\footnote{In the Mitchell-Steel language, one requires $\sigma$ to be a weak $n$-embedding such that $\sigma''T_n^{\bar{\M}} \subseteq T_n^\M$.}\footnote{If $\omega\rho^{n+1}_{\tilde{\M}} = \omega\rho^{n+1}_\M = o(\P^b)$, then since $o(\P^b)$ is a cardinal of $\P$, cr$(\sigma)>o(\P^b)$. Equality can happen in other cases.} Then $\tilde{\M}$ is solid and $p_\M$ is $k$-universal for all $k\in\omega$. Furthermore, if $\tilde{\M}$ is sound above $\nu=\rm{cr}(\sigma)$ then one of the following holds:
\begin{enumerate}[(a)]
\item $\tilde{\M}$ is $\nu$-core of $\M$ and $\sigma$ is the uncollapse map.
\item $\tilde{\M}$ is an initial segment of $\M$.
\item $\tilde{\M} = \rm{Ult}^*(\M||\eta,E^\M_\alpha)$ where $\nu\leq \eta < o(\M)$, $\alpha\leq\omega\eta$ and $\nu=(\kappa^+)^{\M||\eta}$ where $\kappa=\rm{cr}$$(E^\M_\alpha)$; moreover, $E^\M_\alpha$ has a single generator $\kappa$.

\end{enumerate}
\end{theorem}

\begin{remark}\label{rmk:crit_pt}
If $\delta_\alpha$ is a strong cutpoint of $\P$, then it follows simply from the definition of hod premice that for all $\kappa\in [\delta_\alpha,\delta)$, $\square_\kappa$ holds in $\P$; this is because $\P$ is a $\Sigma_\alpha^\P$-premouse ($\Sigma_\alpha$ is the strategy for $\P(\alpha)$) and the $\square$ proof of \cite{schimmerling2004characterization} adapts straightforwardly. On the other hand, if $\delta_\alpha$ is not a strong cutpoint of $\P$, then Theorem \ref{thm:condensation_lemma} is false if for example one required that $\M=\P$, the embedding $\sigma$ have critical point $\delta_\alpha$, $\delta_\alpha$ is a limit of Woodin cardinals, $\tilde{\M}\in \M$ is sound, and the cardinality $\tilde{\M}$ in $\M$ is less than $\delta_\alpha^+$. In this case, $\rho^{\tilde{\M}}_\omega=\delta_\alpha = \nu$. Case (a) is immediately ruled out because $\tilde{\M}\in \M$. Case (c) cannot happen because of the fact that no extenders on the $\M$-sequence can be indexed at $\delta_\alpha$, which is a cardinal of $\P$. Case (b) also fails because otherwise, $\tilde{\M} \lhd \M|\delta_\alpha^+$. $\tilde{\M}$ and hence $\M$ has extenders indexed in the interval $[\delta_\alpha, o(\tilde{\M}))\subset [\delta_\alpha,\delta_\alpha^+)$. But by the definition of hod premice, no extenders on the $\M$-sequence can be indexed in the interval $ [\delta_\alpha,\delta_\alpha^+)$.

Since we are below superstrong and the extenders on the model are indexed according to the Jensen indexing, the possibility that 
\begin{center}
$\tilde{\M}$ is a proper initial segment of $\rm{Ult}$$(\M,E^\M_{\rm{cr}(\sigma)})$
\end{center}
cannot happen. In \cite[Lemma 9.3.2]{Zeman}, the aforementioned case can occur; in that case, $E^\M_{\rm{cr}(\sigma)}$ is a superstrong extender.$\myqedhere$
\end{remark}

The proof of the theorem is essentially that of \cite[Theorem 9.3.2]{Zeman}. The idea is one compares the phalanx $(\M,\tilde{\M},\rm{cr}(\sigma))$ against $\M$. Depending on how the comparison terminates, one gets one of the four possibilities in the statement of the theorem. Using locally strong $\Gamma$-fullness preservation and the fact that cr$(\sigma)>\omega\rho^{n+1}_\M$, Lemma \ref{lem:no_strat_disagreement} shows that the comparison is an extender comparison (no strategy disagreements are encountered). This puts us the in the situation to apply the proof of \cite[Theorem 9.3.2]{Zeman} (the Dodd-Jensen-like property we assume as part of locally strong branch condensation is enough to carry out the proof of \cite[Theorem 9.3.2]{Zeman}). To illustrate the main ideas, we present a proof of a special case, which often shows up in the $\square$-constructions.

\begin{proof}
We assume $\tilde{\M}$ is sound. Let $\tilde{\tau} = \rm{cr}(\sigma)$ and let $\tau = \sigma(\tilde{\tau})$. We further assume that: letting $\kappa = \omega\rho^{n+1}_\M$, $\tau = (\kappa^+)^\M$, and hence $\tilde{\tau} = (\kappa^+)^{\tilde{\M}}$. In this case, we prove that $\tilde{\M} \lhd \M$. The reader can see \cite[Theorem 9.3.2]{Zeman} for the full argument.

\begin{claim}\label{claim:0_phalanx_comp}
Let $\Lambda = \Sigma_{\M}^{\sigma}$. Then the comparison of the phalanx $(\M,\tilde{\M},\tilde{\tau})$ and $\M$ using $\Lambda$ and $\Sigma_{\M}$ respectively is successful. Furthermore, the main branch on the phalanx side doesn't drop (in model or degree) and is above $\tilde{\M}$, and the $\M$ side doesn't move.
\end{claim}
\begin{proof}
Using locally strong fullness preservation of $\Sigma$, $\Lambda$ is fullness preserving; so the comparison can be carried out (see Chapter 4). By Lemma \ref{lem:no_strat_disagreement}, the comparison is an extender comparison (no strategy disagreements show up in the comparison). Now we use locally strong branch condensation to prove the claim. The proof is a fairly standard argument. Let $\T$ and $\U$ be the trees on $(\M,\tilde{\M},\tilde{\tau})$ and $\M$ respectively that are generated by the comparison (via $\Lambda$ and $\Sigma_{\M}$ respectively). The comparison terminates successfully with $\Q$ being the last model of $\T$ and $\S$ being the last model of $\U$. 

Let $\sigma\T$ be the copy tree and $\sigma^*:\Q\rightarrow \Q^*$ be the copy map, where $\Q^*$ is the last model of $\sigma\T$. Note then that $\sigma\T$ is via $\Sigma_{\M}$. 

Suppose $\Q$ is above $\M$. We prove this case is impossible. Suppose $\Q\lhd \S$ and hence the branch embedding $\pi^\T$ exists. Note that $(\Q^*)^b\lhd \Q$ and $\Q^*$ is a non-dropping $\Sigma_\M$ iterate. Hence by strong branch condensation, 
\begin{center}
$\Sigma_{\Q,\T}^{\pi^{\T}} = \Sigma_{\M}$.
\end{center}
The usual Dodd-Jensen argument yields a contradiction. The main point is that the tree $\pi^\T \U$ is via $\Sigma_{\Q,\T}$. 

Suppose now $\S\lhd \Q$ and hence the branch embedding $\pi^\U$ exists. Note then that $\sigma^*(\S)\lhd \Q^*$. Again, by strong branch condensation, 
\begin{center}
$\Sigma_{\sigma^*(\S),\sigma\T}^{\sigma^*\rest\S \circ \pi^\U} = \Sigma_{\M}$.\footnote{Note that in this case, $\S^b = \sigma^*(\S)^b$. The last equality follows from the fact that $\sigma^*$ has critical point $>o(\Q^b)\geq o(\S^b)$.}
\end{center}
The usual Dodd-Jensen argument then yields a contradiction. The main point is that $(\sigma^*\rest\S\circ\pi^\U)\sigma\T$ is according to $\Sigma_{\Q^*,\sigma\T}$.

The above arguments easily give us that: $\Q = \S$ and $\pi^\T$, $\pi^\U$ both exist and they are equal. We can then find a pair of extenders $(E, F)$ used on $\T$ and $\U$ respectively such that $E$ and $F$ are compatible. By a standard argument using the $\sf{ISC}$, this is not possible.

Hence $\Q$ is on the main branch above $\tilde{\M}$. Note then that if $\pi^\T$ exists, then cr$(\pi^\T)>\tilde{\tau}$. Say $b$ is the main branch of $\T$. Then $b$ cannot drop (in model or degree) as otherwise, we have $\S\lhd \Q$ and $\pi^\U$ exists. As before, $\sigma^*(\S)\lhd \Q^*$ and $\Sigma_{\M} = \Sigma_{\Q^*,\sigma\T}^{\sigma^*\rest\S\circ \pi^\U}$. We get a contradiction as before.

So $b$ doesn't drop. Since $\tilde{\M}$ is $\kappa$-sound, $\rho^\omega_{\tilde{\M}} = \kappa < \tilde{\tau}$ and the branch $b$ is above $\tilde{\tau}$ and does not drop in model or degree, we get that $b=\{0\}$. And hence $\Q = \tilde{\M}$. Now it's not the case that $\S$ is a strict segment of $\Q=\tilde{\M}$; otherwise, $\pi^\U$ exists and 
\begin{center}
$\sigma^*\circ \pi^\U: \M \rightarrow \sigma^*(\S)\lhd \Q^*$.
\end{center}
We get a contradiction as before.

If $\S=\Q=\tilde{\M}$, then $\U$'s main branch doesn't drop. This is because $\tilde{\M}$ is sound. Note also that $\U\neq \emptyset$ since otherwise, $\M = \tilde{\M}$ which is impossible (after all, $\tau = (\kappa^+)^{\M} > \tilde{\tau} = (\kappa^+)^{\tilde{\M}}$). Now $\rho^\omega_\M = \rho^\omega_{\tilde{\M}} = \kappa$ and if there is an extender $E$ used along the main branch of $\U$ such that $\nu(E) > \kappa$ \footnote{$\nu(E)$ is the sup of the generators of $E$.} then either $\S$ is not $\kappa$-sound or $\rho(\S) > \kappa$. Contradiction. 

So for all $E$ used along the main branch of $\U$, $\nu(E)\leq \kappa$. If for all such $E$, $\nu(E)<\kappa$, then since $\M|\kappa = \tilde{\M}|\kappa = \Q|\kappa$, $\S|\kappa\neq \Q|\kappa$. Contradiction. If there is some such $E$ such that $\nu(E)=\kappa$, then $\U$ must drop since otherwise, $\rho^\omega_\S>\kappa$. Contradiction.

So $\Q\lhd \S$. We claim that $\Q=\tilde{\M} \lhd \M$. It suffices to show $\U=\emptyset$. Otherwise, let $E = E^\U_0$. Then 
\begin{equation}\label{eqn:one1}
\textrm{lh}(E)\geq \tilde{\tau} \textrm{ and lh}(E) < o(\Q) = o(\tilde{\M}).
\end{equation}
Note that lh$(E)$ is a cardinal of $\S$ strictly larger than $\kappa$ and $|o(\Q)|^\S = \rho^\omega_{\tilde{\M}} = \kappa$. This contradicts \ref{eqn:one1}. This completes the proof of Claim \ref{claim:0_phalanx_comp}.
\end{proof}

Using the claim, it is easy to see that $\tilde{\M}\lhd \M$ (that is, case (b) holds). This is because the branch embedding on the phalanx side must have critical point $> \kappa$ and $\tilde{\M}$ is $\kappa$-sound, so the branch is trivial with end model $\tilde{\M}$.
\end{proof}

\section{Ingredients from the Schimmerling-Zeman construction}\label{sec:S_Z}

In this section, we briefly remind the reader of the $\square$-construction in \cite{schimmerling2004characterization}.  First, the reader should recall from \cite{schimmerling2004characterization} the notions of a \textit{protomouse} and a \textit{pluripotent level} of $L[E]$ (we give definitions of these notions in the context of hod premice in Section \ref{sec:set_up}). See the beginning of \cite[Section 2]{schimmerling2004characterization} for a fairly detailed discussion on how protomice appear in interpolation arguments. Basically, protomice arise in interpolation arguments where the target structure is a pluripotent level. The reader should see the definition of \textit{divisor}, \cite[Section 2.1]{schimmerling2004characterization}, and \textit{strong divisor}, \cite[Section 2.4]{schimmerling2004characterization} (these notions are also defined in Section \ref{sec:set_up} for hod premice). Divisors identify protomice in interpolation arguments and (canonical) strong divisors in some sense are those (amongst many possible divisors of a given collapsing structure) that one uses in the course of the construction.  

We proceed to briefly outline the proof of $\square_\kappa$ in $L[E]$ as done in \cite{schimmerling2004characterization}. To get the main ideas across in a reasonable amount of space, we will be imprecise at various places. The reader can see \cite[Section 3]{schimmerling2004characterization} for a precise construction of the $\square_\kappa$-sequence $(C_\tau : \tau<\kappa^+)$. The proof starts by choosing the collapsing structure $\N_\tau$ for $\kappa < \tau < (\kappa^+)^{L[E]}$: $\N_\tau$ is the first level of $L[E]$ that satisfies ``$\tau = \kappa^+$" and $\rho^\omega_{\N_\tau}=\kappa$. There is a club $\mathcal{S}\subset\kappa^+$ of such $\tau$ in $L[E]$. We further require that for each $\tau\in \S$, $\J^E_\tau \prec \J_{\kappa^+}^E$. For each $\tau\in \mathcal{S}$, let $\mathcal{S}_1\subseteq \mathcal{S}$ be the set of $\tau$ for which the strong divisors of $\N_\tau$ exists (and let $(\mu(\N_\tau),q(\N_\tau))$ be the canonical strong divisor and $\N_\tau(\mu(\N_\tau),q(\N_\tau))$ be the unique associated protomouse as defined at the end of \cite[Section 2]{schimmerling2004characterization}). Let $\mathcal{S}_0=\mathcal{S} - \mathcal{S}_1$. 

For $\tau \in \S_0$, the associated club $C_\tau\subset \tau$ can be constructed by Jensen's method of constructing $\square$-sequences in $L$. In this case, $C_\tau$ is obtained canonically from the set $B_\tau$ of all $\bar{\tau}\in \S_0\cap \tau$ such that:
\begin{itemize}
\item $\N_{\bar{\tau}}$ is a premouse of the same type as $\N_\tau$ and $n_{\bar{\tau}}=n_\tau$, where for a $\sigma\in \S$, $n_\sigma$ is the least $n$ such that $\omega\rho^{n+1}_{\N_\sigma}\leq \kappa < \omega\rho^n_{\N_\sigma}$.

\item There is a map $\sigma_{\bar{\tau},\tau}:\N_{\bar{\tau}}\rightarrow \N_\tau$ that is $\Sigma_0^{(n_\tau)}$-preserving with respect to the language of premice and such that: $\bar{\tau} = \rm{cr}(\sigma_{\bar{\tau},\tau})$, $\sigma_{\bar{\tau},\tau}(\bar{\tau})=\tau$, $\sigma_{\bar{\tau},\tau}(p(\N_{\bar{\tau}}))=p(\N_\tau)$, and each $\alpha\in p(\N_\tau)$ has a generalized witness with respect to $(\N_\tau,p(\N_\tau))$ in the range of $\sigma_{\bar{\tau},\tau}$. Here, and later, $p(\N_\tau)$ is the $n_\tau^{th}$-standard parameter of $\N_\tau$.
\end{itemize} 

For $\tau\in \S_1$, the set $C_\tau$ is obtained canonically from the set $B_\tau$ of $\bar{\tau}\in \S_1\cap \tau$ that satisfies: 
\begin{itemize}
\item $(\mu(\N_{\bar{\tau}}),|q(\N_{\bar{\tau}})|) = (\mu(\N_{\tau}),|q(\N_{\tau})|)$; here by definition of divisors, $q(\N_\tau)$ is a bottom initial segment of $d(\N_\tau)$, the Dodd-parameter of $\N_\tau$.
\item There is a map $\sigma_{\bar{\tau},\tau}:\N_{\bar{\tau}}(\mu(\N_{\bar{\tau}}),q(\N_{\bar{\tau}})) \rightarrow\N_\tau(\mu(\N_{\tau}),q(\N_{\tau}))$ that is $\Sigma_0$-preserving with respect to the language for coherent structures such that: $\bar{\tau} = \rm{cr}(\sigma_{\bar{\tau},\tau})$, $\sigma_{\bar{\tau},\tau}(\bar{\tau})=\tau$, $\sigma_{\bar{\tau},\tau}(q(\N_{\bar{\tau}}))=q(\N_\tau)$, and each $\alpha\in q(\N_\tau)$ has a generalized witness (with respect to $(\N_\tau(\mu(\N_{\tau}),q(\N_{\tau})),q(\N_\tau))$ in the range of $\sigma_{\bar{\tau},\tau}$.
\end{itemize}

See \cite[Section 3]{schimmerling2004characterization} for how $C_\tau$ is defined from $B_\tau$. Now we focus on the key point: the proof that 
\begin{center}
$B_\tau$ is unbounded in $\tau$ if $\tau\in \S^1$ and cof$(\tau)>\omega$ in $L[E]$. 
\end{center}
Fix such a $\tau$ and let $\kappa<\gamma < \tau$ be arbitrary. We want to find a $\gamma<\bar{\tau}<\tau$ in $B_\tau$. Working in $L[E]$, fix some $\theta>>\kappa$ and let $X\prec H_\theta$ be countable such that all relevant objects are in $X$, in particular $\{\kappa,\tau,\gamma\}\in X$. Let $\sigma: \bar{\M}\rightarrow \M$ be the uncollapse map of $X\cap \M$, where $\M = \N(\mu(\N_\tau),q(\N_\tau))$. We write $\sigma^{-1}(x) =\bar{x}$ for each $x$ in the range of $\sigma$. Let $\tilde{\tau}=\sup(\sigma''\bar{\tau})$. Let $\tilde{\sigma}:\bar{\M}\rightarrow\tilde{\M}$ come from the $(\rm{cr}(\sigma),\tilde{\tau})$-extender derived from $\sigma$. Also, let $\sigma':\tilde{\M}\rightarrow \M$ be given by the interpolation lemma \cite[Lemma 1.2]{schimmerling2004characterization}. In this case, $\tilde{\M}=(\N,\tilde{F})$ is a protomouse (even if $\N(\mu(\N_\tau),q(\N_\tau)) = \N_\tau$ since in this case, $\N_\tau$ is a pluripotent level of $L[E]$ and the map $\sigma'$ is not cofinal). The way one shows $\tilde{\tau}\in B_\tau$ is as follows. Let $\M^*$ be the largest segment of $\N$ such that $\tilde{F}$ measures all sets in $\M^*$. One then shows that Ult$(\M^*,\tilde{F})$ is $\N_{\tilde{\tau}}$. Say $\M=(\M^-,F)$. This is accomplished by applying the condensation lemma \cite[Lemma 1.6]{schimmerling2004characterization} to $\phi: \rm{Ult}$$(\M^*,\tilde{F}) \rightarrow \pi_F(\M^*)$ such that
\begin{center}
$\phi(\pi_{\tilde{F}}(f)(a)) = \sigma'(f)(\pi_F(a))$
\end{center}
where $\pi_F$ is the $F$-ultrapower embedding applied to the largest initial segment of $\M^-$ that makes sense.

The key for the proof above is that we can always compare two iterable pure extender models; in this case, we compare the phalanx $(\pi_F(\M^*),\rm{Ult}$$(\M^*,\tilde{F}),\tilde{\tau})$ against $\pi_F(\M^*)$. If one adapted this argument to hod mice, it fails because the hod mice $\pi_F(\M^*)$ and $\rm{Ult}$$(\M^*,\tilde{F})$ generally belong to two different pointclasses, and hence cannot be directly compared. This traces back to the fact that if $\P$ is a hod premouse of limit type with a top window, then $(\delta)^{\P^b}$ is a strong cutpoint of $\P^b$. The fix for this, as done in the next section, is to sometimes allow for the collapsing structure of $\tau$, $\N_\tau$, to \textit{not be an initial segment of} the hod mouse and incorporate this kind of collapsing structures into the construction. It is this aspect that forces the construction to yield a weaker result, i.e. $\square_{\kappa,2}$, rather than $\square_\kappa$.

One other new situation in the hod mouse case that does not come up in the $L[E]$ case is the following. Suppose in the above, $\M = \N_\tau$ is $B$-active. Then the way branches are coded into the model (using the $\mathfrak{B}$-operator as discussed in the previous section) allows us to show that $\tilde{M}$ is a $B$-active hod premouse. If one used the traditional coding of branches, then $\tilde{M}$ may fail to be a hod premouse due to the lack of condensation one can prove with the traditional coding device; this is the reason we switch to the coding of branches via the $\mathfrak{B}$-operator. We will discuss this in more details in the next section.

\section{The proof}\label{sec:proof}

We give a proof of Theorem \ref{thm:square}, making use of the notions, notations, and proofs in \cite{schimmerling2004characterization} whenever applicable. We only focus on the details that are new in our situation and direct the reader to constructions in \cite{schimmerling2004characterization} that are obviously generalizable to our situation. To make the situation more concrete, we make the following assumption about $\P$:
\begin{equation}\label{eqn:lsatype}
\P \textrm{ is of lsa type.}
\end{equation}
So under (\ref{eqn:lsatype}), $\P$ has the top window $(\mu,\delta^\P)$, $\mu$ is strong to $\delta^\P$. This is the hardest case and we will focus on the $\square$-construction under this assumption. The other cases where $\P$ is not of lsa type are easier and the proof is a simpler version of what we are about to give.
\subsection{Some set-up}\label{sec:set_up}

We will use the fine-structure terminology and notations from \cite[Section 1]{schimmerling2004characterization}, generalized to our context in an obvious way. For example, notions in \cite{schimmerling2004characterization} that are defined using the language of premice are defined here using the language of hod premice; when we talk about a \textit{coherent structure}\index{coherent structure} in this paper, we mean a structure $M$ of the form $(\Q,F)$ where $\Q$ is an amenable structure in the language of hod premice and $F$ is a whole extender at $(\kappa,\lambda)$ (in the language of \cite[Section 1]{schimmerling2004characterization} ) and $\Q = \rm{Ult}$$(\Q|\bar{\alpha},F)$, where $\bar{\alpha}$ is the largest $\tau \leq \kappa^{+\Q}$ such that dom$(F)=\powerset(\kappa)\cap \Q|\tau$. We say $N$ is \textit{the hod premouse associated with $M$}\index{hod premouse associated to a coherent structure}. The notion of \textit{a generalized witness} for some ordinal $\alpha$ with respect to a pair $(M,s)$ where $M$ is a coherent structure, $s$ is a finite set of ordinals (or a generalized witness for an ordinal $\alpha$ with respect to a hod premouse $N$ associated with $M$ and some finite set of ordinals $r\cup s$) in \cite{schimmerling2004characterization} is generalized in an obvious way to our context.\footnote{Let $M, N, \kappa, \lambda$ be as above and $s\subset \lambda$ is finite. The \textit{standard witness}\index{standard witness} $W^{\alpha,s}_M$\index{$W^{\alpha,s}_M$} for $\alpha$ with respect to $M$ and $s$ to be the transitive collapse of $h_M(\alpha\cup\{s\})$, where $h_M$ is the canonical $\Sigma_1$-Skolem function of the coherent structure $M$. Similarly, $W^{\alpha,r\cup s}_N$ denotes the standard witness for $\alpha$ with respect to $N$ and $r\cup s$ and is the transitive collapse of $\tilde{H}^{n+1}_N(\alpha\cup \{r\cup s\})$, where $\tilde{h}^{n+1}_N$\index{$\tilde{h}^{n+1}_N$} is the canonical $\Sigma^{(n)}_1$-Skolem function of the hod premouse $N$. A \textit{generalized witness}\index{generalized witness} for $\alpha$ with respect to $M$ and $s$ is a pair $(Q,t)$, where $t\subset Q$ is a finite set of ordinals and such that for any $\xi_1,\dots,\xi_l<\alpha$, if $M\vDash \Phi(i,\xi_1,\dots,\xi_l,s)$ then $Q\vDash \Phi(i,\xi_1,\dots,\xi_l,t)$, where $\Phi$ is the universal $\Sigma_1$-formula. A generalized witness for $\alpha$ with respect to $N$ and $r\cup s$ is a pair $(Q,t)$, where $t\subset Q$ is a finite set of ordinals such that given any $\xi_1,\dots,\xi_l < \alpha$, if $N\vDash \Phi(i,\xi_1,\dots,\xi_l,r\cup s)$ then $Q\vDash \Phi(i,\xi_1,\dots,\xi_l,t)$, where $\Phi$ is the universal $\Sigma^{(n)}_1$-formula.}  A \textit{protomouse} $\P=(\Q,F)$ \index{protomouse}is a coherent structure where $F$ is an extender with critical $\kappa$ such that $F$ does not measure $\powerset(\kappa)^\Q$. A \textit{pluripotent level}\index{pluripotent level} of a hod premouse $\P$ is an $E$-active initial segment $\Q$ of $\P$ such that cr$(E^\Q_{\rm{top}}) < \kappa$ and $\omega\rho^1_\Q=\kappa$, where $\kappa$ is a cardinal of $\P$. The language used for treating protomice and pluripotent levels is the language of coherent structures, namely $\mathcal{L}_1 - \{\dot{\gamma}\}$.\footnote{We note again that for type $A$ hod premice, which are all the hod premice that we encounter in this book, $\dot{\gamma} = \emptyset$, so there is essentially no distinction between the language of hod premice and coherent structures.}

Fix $(\P,\Sigma)$ as in the hypothesis of Theorem \ref{thm:square}. Fix $\kappa\geq \delta^{\P^b}$, a cardinal of $\P$.  Working in $\P$, let $\mu = \delta^{\P^b}$ and $\S\subset \kappa^{+}$ be the club of $\kappa<\tau<\kappa^+$ such that for some $\tau<\tau' < \kappa^+$, $\P|\tau' \prec \P|\kappa^{++}$ and $\tau = (\kappa^+)^{\P|\tau'}$. We note that since we are below superstrong, the set of indices below $\kappa^+$ is non-stationary in $\kappa^+$; therefore, we can assume the club $\S$ consists of $\tau$ such that $\P|\tau$ is passive.

Let $\N^*_\tau\lhd \P$\index{$\N^*_\tau$} be the collapsing level for $\tau$, that is $\N^*_\tau$ the least initial segment $\N$ of $\P$ such that $\N\vDash \tau = \kappa^+$ and $\rho^\omega_\N=\kappa$. Let $\gamma_\tau$ be the supremum of  indexes $\gamma$ of extenders $E\in \dot{E}^{\N^*_\tau}$ such that cr$(E)=\mu$. Note that $\gamma_\tau < o(\N^*_\tau)$ and if $\tau < \sigma$ in $\S$, then $\N^*_\tau \lhd \N^*_\sigma$ and therefore, $\gamma_\tau \leq \gamma_\sigma$. Without loss of generality, we may assume throughout this chapter that 
\begin{center}
$\kappa \geq o(\P^b)$ and sup$_{\tau\in\S}(\gamma_\tau) \geq \kappa^+$.\footnote{If $\kappa = \delta^\P_{\alpha}$, where $\mu  = \delta_{\alpha}^\P$, then since $\P^b=\P|(\kappa^+)^\P = \rm{Lp}$$^{\oplus_{\beta<\alpha}\Sigma^\P_{\P(\beta)}}(\P|\delta^\P_\alpha)$, then $\P\vDash \square_\kappa$ since $\kappa$ is a strong cutpoint cardinal of $\P^b$. If $\kappa > \delta_\alpha^\P$ and sup$_{\tau\in\S}(\gamma_\tau)<\kappa^+$, then the proof is significantly easier. One constructs the $\square_\kappa$-sequence using points $\tau\in\S$ above sup$_{\tau\in\S}(\gamma_\tau)$ mimicking essentially the Schimmerling-Zeman construction and use Theorem \ref{thm:condensation_lemma}. If $\kappa < \mu$, we will then be showing $\P(\alpha) \models \square_{\kappa,2}$, where $\P(\alpha)$ plays the role of $\P$.}
\end{center}
The following statements can be easily verified using the definitions and our assumption.
\begin{proposition}\label{prop:basic_facts}
For a club of $\tau \in \S$:
\begin{enumerate}
\item $o(\N^*_\tau) > \tau$,
\item $\gamma_\tau \geq \tau$.
\end{enumerate}
\end{proposition}
\begin{proof}
By the construction of $\S$, for any $\tau\in \S, \P|\tau$ is passive and is a $\sf{ZFC}^-$-model, and therefore $\rho_\omega^{\P|\tau} = \tau$. This means $o(\P|\tau) < o(\N^*_\tau)$. This proves (1).

For (2), suppose not. Then there is a stationary set of $\tau\in \S$ such that $\gamma_\tau < \tau$. By pressing down, there is a stationary set $T\subset \S$ and a $\gamma < \kappa^+$ such that for all $\tau\in T$, $$\gamma_\tau = \gamma < \tau.$$ But then since the sequence $(\gamma_\tau : \tau \in \S)$ is monotonically non-decreasing, this means that $$\textrm{sup}_{\tau \in \S}(\gamma_\tau) = \gamma < \kappa.$$ This contradicts our assumption.
\end{proof}

Extenders $E$ with cr$(E)=\mu$ play a special role in this construction. Note that $\mu$ is a strong cutpoint of $\P^b$, that is, there are no partial extenders with critical point $\mu$ on the sequence of $\P$. This is the main difference between our situation and the $L[E]$-situation. It is this situation that forces us to consider collapsing structures that are not initial segments of our hod mouse $\P$.

Some discussions regarding protomice and divisors are in order. Following \cite{schimmerling2004characterization}, for a hod premouse $\N$ such that $\omega\rho^{n+1}_\N\leq\kappa<\omega\rho^n_\N$, we say that $(\nu,q)$ is a \textit{divisor} of $\N$ \index{divisor}if and only if there is an ordinal $\lambda=\lambda_\N(\nu,q)$\index{$\lambda_\N(\mu,q)$} such that letting $p_\N$ be the ($n+1$)-st standard parameter of $\N$, setting $r=p_\N-q$, the following hold:
\begin{enumerate}[(a)]
\item $\nu\leq\kappa < \lambda <\omega\rho^n_\N$;
\item $q = p_\N\cap \lambda$;
\item $\tilde{h}^{n+1}_\N(\nu\cup\{r\})\cap \omega\rho^n_\N$ is cofinal in $\omega\rho^n_\N$;
\item $\lambda = \rm{min}(\rm{OR}$$\cap \tilde{h}^{n+1}_\N(\nu\cup\{r\})-\nu)$.
\end{enumerate}
As in \cite{schimmerling2004characterization}, both $\nu$ and $\lambda$ are (inaccessible) cardinals in $\N$. Let $\N^*(\nu,q)$\index{$\N^*(\nu,q)$} be the transitive collapse of $\tilde{h}^{n+1}_\N(\nu\cup\{r\})$.

The notion of strong divisors in \cite{schimmerling2004characterization} generalize in an obvious way to our context. We recall it now. A divisor $(\nu,q)$ of $\N$ is \textit{strong}\index{strong divisor} if and only if for every $\xi<\nu$ and every $x$ of the form $\tilde{h}^{n+1}_\N(\xi,p_\N)$ we have $x\cap \nu\in \N^*(\nu,q)$. If $\N$ is pluripotent, we define the notion of strong divisor in the same way, but with $h^*_\N$ (the $\Sigma_1$-Skolem function of $\N$ computed in the language of coherent structures) and $d_\N$ (the Dodd-parameter of $\N$) in place of $\tilde{h}^{n+1}_\N$ and $p_\N$, respectively. It turns out, see \cite[Section 2]{schimmerling2004characterization}, that the notion of divisor/strong divisor does not depend on which language one uses to compute it.

As in \cite{schimmerling2004characterization}, if $\N$ has strong divisors, \textit{the canonical strong divisor}\index{canonical strong divisor} $(\mu_\N,q_\N)$\index{$(\mu_\N,q_\N)$} of $\N$ is chosen as follows: $q_\N$ is the shortest initial segment of $p_\N$ such that for some $\nu^*$, $(\nu^*,q_\N)$ is a strong divisor of $\N$ and $\mu_\N$ is the largest $\nu^*$ such that $(\nu^*,q_\N)$ is a strong divisor of $\N$. Now we define our collapsing structure $\N_\tau$ for $\tau\in\S$.
\begin{definition}\label{def:collapsing_structures}
Fix $\tau\in \S$. Suppose there is a pointclass $\Omega\subsetneq\Gamma$ such that there is a hod pair $(\R,\Sigma_\R)$ such that 
\begin{enumerate}[(i)]
\item $\N^*_\tau|\gamma_\tau\lhd \R$, 
\item $\rho^\omega_\R=\kappa$, 
\item $\gamma_\tau$ is a cutpoint of $\R$ and $\Sigma_{\R|{\gamma_\tau}}=\Sigma_{\P|{\gamma_\tau}}$, 
\item $\R$ is $\gamma_\tau$-sound, 
\item the order type of $\R$'s layer Woodin cardinals above $\gamma_\tau$ is a limit ordinal,
\item $\R$ has a strong divisor of the form $(\mu,q)$ where $p_\R=q\cup r$ for $r$ above the supremum $\lambda$ of the layer Woodin cardinals of $\R$ and max$(q)$ is below $(\gamma_\tau^+)^\R$ and $(\gamma_\tau^+)^\R<\lambda$, 
\item $\Sigma_\R$ has branch condensation, is $\Omega$-fullness preserving, and $(\R,\Sigma_\R)$ generates $\Omega$; that is $\Gamma(\R,\Sigma_\R)=\Omega$. 
\end{enumerate}
We call $(\R,\Sigma_\R)$ with the above properties the pointclass generator of $\Omega$\index{pointclass generator}. Let $\Gamma_\tau$ be the Wadge-minimal such pointclass and $\N_\tau$ be the pointclass generator of $\Gamma_\tau$, $(\mu_\tau,q_\tau,\lambda_\tau)$ be the $(\mu,q,\lambda)$ associated with $\N_\tau$ as above (note that $\N_\tau$ must be distinct from $\N^*_\tau$ in this case). If $(\Gamma_\tau,\R, \mu,q,\lambda)$ doesn't exist, we let $\N_\tau =\N^*_\tau$.$\myqedhere$
\end{definition}
The properties of pointclass generators seem technical; these properties are abstracted from various situations in interpolation arguments. It seems hard to do with much less. Here is a very rough explanation for why we would consider such objects before going into details: suppose $(\Gamma_\tau,\R, \mu,q,\lambda)$ exists, then letting $\pi: \N^*_\tau(\mu,q)\rightarrow \R$ be the uncollapse map and $F$ be the $(\mu,\gamma_\tau)$-extender derived from $\pi$. Letting $\R^*\lhd \P|\mu^+$ be the largest such that $\rho_\omega^\S=\mu$ and $\powerset(\mu)\cap \R^* = \powerset(\mu)\cap \N^*_\tau(\mu,q)$., then we can show that $\R = \textrm{Ult}(\R^*, F)$. Since $\powerset(\mu)\cap \R^*\subsetneq \powerset(\mu)\cap \P$, $\R$ belongs to $\Gamma$ and will generate a pointclass ($\Gamma_\tau$) strictly smaller than $\Gamma$. It turns out that there is a unique such $(\R,\Sigma_\R)$ that generates $\Gamma_\tau$. This canonicity is important in the $\square$-construction.

The following proposition justifies the uniqueness of pointclass generators.
\begin{proposition}\label{prop:unique_ptclass_gen}
Let $\P,\tau, \Omega$ be as in Definition \ref{def:collapsing_structures}. Let $(\R_0,\Sigma_0)$ and $(\R_1,\Sigma_1)$ be pointclass generators of $\Omega$ with the properties described in \ref{def:collapsing_structures}. Then $(\R_0,\Sigma_0) = (\R_1,\Sigma_1)$.
\end{proposition}
\begin{proof}
We compare the pair $(\R_0,\Sigma_0)$ against $(\R_1,\Sigma_1)$, lining up the models and the strategies (as done in Section \ref{sec:normal_comparison}). The comparison is possible by the assumption, namely $\Gamma(\R_0,\Sigma_0)=\Gamma(\R_1,\Sigma_1)=\Omega$, and is above $\gamma_\tau$ (this is because $\gamma_\tau$ is a cutpoint, in fact a strong cutpoint, of both models and $\R_0|\gamma_\tau = \R_1|\gamma_\tau$). The end model is, say, $\S$ and the tail strategies of $\Sigma_0$ and $\Sigma_1$ on $\S$ are the same. The usual proof using the fact that $\R_0$ and $\R_1$ are $\gamma_\tau$-sound and the comparison is above $\gamma_\tau$ shows that $\S = \R_0 = \R_1$ (the comparison is trivial) and $\Sigma_0 = \Sigma_1$. The equality of models (i.e. $\R_0 = \R_1$) follows from the fact that $(\R_0,\Sigma_0)$ and $(\R_1,\Sigma_1)$ generate the same pointclass. Otherwise, say $\R_0 = \R_1(\alpha) \lhd \R_1$ for some $\alpha$ and $\Sigma_0 = (\Sigma_1)_{\R_0}$. It is easy to see that from (ii) and (v), there is a layer Woodin cardinal $\beta$ of $\R_1$ such that $\R_0 \lhd \R_1(\beta)$. But this means that $\Omega = \Gamma(\R_0, \Sigma_0) \subsetneq \Gamma(\R_1(\beta), (\Sigma_1)_{\R_1(\beta)})\subseteq \Gamma(\R_1,\Sigma_1) = \Omega$. Contradiction.
\end{proof}

We simply use the notations from \cite[page 49]{schimmerling2004characterization} in the definition of our square sequence below. For instance, $(\mu_\tau,q_\tau)$ denotes the canonical strong divisor of $\N_\tau$ (if exists) in the case $\N_\tau=\N_\tau^*$ and denotes the $(\mu_\tau,q_\tau)$ in Definition \ref{def:collapsing_structures} in the case $\N_\tau\neq \N_\tau^*$ (note that $(\mu_\tau,q_\tau)$ is the unique strong divisor of $\N_\tau$ with the properties as in Definition \ref{def:collapsing_structures}). If $\N_\tau=\N^*_\tau$ is a pluripotent level that has no strong divisors, then $(\mu_\tau,q_\tau)$ denotes $(\rm{cr}$$(E^{\rm{top}}_{\N_\tau}),d(\N_\tau))$, where $d(\N_\tau)$ is the Dodd-parameter of $\N_\tau$ . 

Suppose $(\nu,q)$ is a divisor of $\N_\tau$; let $r, \lambda, n$ be as in the definition of divisor. Let $\pi:\N_\tau^*(\nu,q)\rightarrow \tilde{h}^{n+1}_\N(\nu\cup\{r\})$ be the uncollapse map. We let the associated protomouse $\N_\tau(\nu,q)$ be the coherent structure $(\N_\tau|\xi,F)$ where $\xi = \pi((\nu^+)^{M^*})$ and $F = E_\pi\rest (\powerset(\nu)\cap \N_\tau^*(\nu,q))$. We denote the $\lambda$ associated to $\nu,q$ in the definition of divisor $\lambda_{\N_\tau}(\nu,q)$. 
We let $\M_\tau=\N_\tau(\mu_\tau,q_\tau)$ be the protomouse associated with $(\mu_\tau,q_\tau)$. If $\N_\tau$ is pluripotent, we let $\M_\tau = \N_\tau$.


The following proposition is easy to see and justifies that the structure $\N_\tau(\nu,q)$ are protomice (and not hod premice). See \cite[Section 2.1]{schimmerling2004characterization} for a detailed discussion and proof.
\begin{proposition}\label{prop:proper_protomice}
Suppose $(\nu,q)$ be a divisor of $\N_\tau$ and $\pi:\N_\tau^*(\nu,q)\rightarrow \tilde{h}^{n+1}_{\N_\tau}(\nu\cup\{r\})$ be the uncollapse map (and in the case $\N_\tau\neq \N^*_\tau$, assume $\nu = \mu$). Then $\powerset(\nu)\cap \N_\tau^*(\nu,q) \subsetneq \powerset(\nu)\cap \N_\tau$. Furthermore, $\nu$ is an (inaccessible) cardinal of $\N_\tau^*(\nu,q)$ and a limit cardinal of $\N_\tau$, and $\lambda_{\N_\tau}(\nu,q)$ is an  (inaccessible) cardinal of $\N_\tau$.
\end{proposition}

\begin{definition}\label{def:S_0_S_1}
Let $\S^1\subset \S$ be the set of $\tau$ such that $(\mu_\tau,q_\tau)$ is defined and $\S^0=\S - \S^1$. $\myqedhere$
\end{definition}

Suppose $\N_\tau=\N^*_\tau$, then no divisors of $\N_\tau$ are of the form $(\mu,q)$. This is because otherwise, $\lambda=\lambda_{\N_\tau}(\mu,q)$ is a limit of Woodin cardinals because it is the image of $\mu$ under the uncollapse map and $\mu$ itself is a limit of Woodin cardinals. Let $\gamma_0<\gamma_1$ be consecutive Woodin cardinals in the interval $(\mu, \lambda)$; then by definition of $\P$, $\P|\gamma_1$ is a $\Lambda^{sts}$-mouse where $\Lambda$ is the strategy of $M^+(\P|\gamma_0)$. On the other hand, by elementarity, $\P|\gamma_1$ is a $\Lambda$-mouse. Contradiction. \footnote{Another argument is as follows. Note that each Woodin cardinal in the interval $(\mu,\lambda)$ is $> (\mu^+)^\P$, and	 hence $\mu$ is strong to $\lambda$ (in $\P$) by the initial segment condition. This contradicts the definition and smallness assumption on $\P$ since one can easily find an active $\omega$ Woodin lsa mouse in $\P$ (as defined in Definition \ref{omega woodins over lsa}).}

A similar argument applies to show that no divisors for $\N_\tau$ are of the form $(\xi,q)$ for $\xi < \mu$; though we don't need this fact in our construction as no divisors $(\nu,q)$ in this proof will have the property that $\nu<\mu$. So if $(\nu,q)$ is a divisor of $\N_\tau$, then $\nu>\mu$. This allows us to simply quote results of \cite[Section 2]{schimmerling2004characterization} in this case (in light of Theorem \ref{thm:condensation_lemma}). In the case that $\mu_{\tau} = \mu$ (so $\N_\tau\neq \N_\tau^*$), more care needs to be taken since it's not obvious that all results in \cite[Section 2.4]{schimmerling2004characterization} can be generalized to this case.

Using the remarks above, it is easy to see that if $\N_\tau\neq \N^*_\tau$, then $\tau\in\S^1$ and in fact $\N_\tau$ is not an initial segment of $\P$ (though $\N_\tau\in \P$ by Proposition \ref{prop:N_tau_def}); also, if $\N_\tau=\N^*_\tau$ is pluripotent, then $\tau\in \S^1$. For $\tau\in \S^0$, $\N_\tau=\N^*_\tau$ is not pluripotent and does not admit a strong divisor.

The following lemma allows us to define our $\square_{\kappa,2}$-sequence in a uniform manner.
\begin{proposition}\label{prop:N_tau_def}
Suppose $\N_\tau\neq \N_\tau^*$. Then $\N_\tau$ is definable over $\P$ (in fact, over any $\N^*_\xi$ or $\N_\xi$ with $\gamma_\xi\geq \xi$ for $\xi>\tau$ in $\S$) unformly from $\{\tau,\gamma_\tau\}$. 
\end{proposition}
\begin{proof}
Fix $\xi > \tau$ in $\S$ with $\gamma_\xi\geq \xi$. We first claim that $\gamma_\xi > \gamma_\tau$. To see this, note that $\tau\leq \gamma_\tau < o(\N^*_\tau) < \xi$. This is because $\xi$ is a cardinal (successor of $\kappa$) in $\N_\xi$ while there is a surjection from $\kappa$ onto $\gamma_\tau$ in $\N_\xi$. Since $\xi \leq \gamma_\xi$, the claim follows.

Now let $E$ be the extender on the $\N_\xi$-sequence such that cr$(E)=\mu$, lh$(E) > \gamma_\tau$, and is the least such.\footnote{\label{fnt:indices_discontinuous}We note that the set of indices for extenders with critical point $\mu$ is nowhere continuous.} Let $\S = \rm{Ult}$$(\N_\xi,E)$ (this is a $\Sigma_0$-ultrapower). Let $i:\S\rightarrow \S_\infty$ be an $\mathbb{R}$-genericity iteration (above $\gamma_\tau$). Now it is easy to see that in the derived model of $\S_\infty$ (at the sup of its Woodin cardinals), the pointclass $\Omega$ in the definition of $\N_\tau$ is a strict Wadge initial segment of $\powerset(\mathbb{R})$ and is definable there from $\{\tau,\gamma_\tau\}$ (by Lemma \ref{prop:unique_ptclass_gen}). Then $\N_\tau \in \S_\infty$ and in fact is definable there from parameters $\{\tau,\gamma_\tau\}$. The same holds in $\S$ by elementarity and the fact that cr$(i)>\gamma_\tau$. Finally, $\N_\tau \in \N_\xi$ and is definable there from parameters $\{\tau,\gamma_\tau, E\}$. But $E$ is definable in $\N_\xi$ from $\{\tau,\gamma_\tau\}$. So $\N_\tau$ is definable in $\N_\xi$ from $\{\tau,\gamma_\tau\}$.
\end{proof}

\begin{remark}\label{rem:active_levels_nonstationary}
By our smallness assumption on $\P$, the set $\mathfrak{A} = \{ \xi \ | \ \kappa < \xi < \kappa^+ \wedge \P|\xi \textrm{ is } E\textrm{-active}\}$ is non-stationary in $\P$. The reason is $\mathfrak{A}=\mathfrak{A}_0\cup\mathfrak{A}_1$. Here $\mathfrak{A}_0$ consists of $\xi$'s such that the top extender of $\P|\xi$ has critical point $\mu$ and $\mathfrak{A}_1=\mathfrak{A}-\mathfrak{A}_0$. $\mathfrak{A}_0$ in nonstationary by Footnote \ref{fnt:indices_discontinuous}. $\mathfrak{A}_1$ is nonstationary because otherwise, $\kappa$ is subcompact by \cite{schimmerling2004characterization}. As in \cite{schimmerling2004characterization}, the fact that $\mathfrak{A}$ is nonstationary is crucial in our construction. We use this fact in various arguments to follow. $\myqedhere$      
\end{remark}

\subsection{Approximation of a $\square_{\kappa,2}$ sequence}
We use the notation established in the previous section. Below, as in \cite{schimmerling2004characterization}, $n_\tau$ is the unique $n$ such that $\rho^{n+1}_{\N_\tau}=\kappa < \rho^n_{\N_\tau}$ and $p_\tau$ is the standard parameter of $\N_\tau$. Let also $p^*_\tau$ be the standard parameter of $\N^*_\tau$.

\begin{definition}\label{def:S0}\index{$\S^0$}
Suppose $\tau\in \S^0$, let $\vec{B}_\tau=\{B^0_\tau\}$ be the set of $\bar{\tau}\in\S\cap\tau$ satisfying:
\begin{itemize}
\item $\N_{\bar{\tau}}$ is a hod premouse of the same type as $\N_\tau$. \footnote{In this case, it simply means: $\N_\tau$ is $E$ ($B$)-active if and only if $\N_{\bar{\tau}}$ is $E$ ($B$)-active. If $\N_\tau$ is $E$-active (equivalently, $\N_{\bar{\tau}}$ is $E$-active), then $E_{\N_\tau}^{\rm{top}}$ is indexed according to the cutpoint (Jensen) indexing scheme if and only if $E_{\N_{\bar{\tau}}}^{\rm{top}}$ is indexed according to the cutpoint (Jensen, respectively) indexing scheme. Recall that all $E$-active hod mice, where $E$ is indexed according to the Jensen indexing scheme, in this proof will be of type $A$, i.e. the set of cutpoints is empty.}
\item $n_\tau = n_{\bar{\tau}}$.
\item There is a map $\sigma^0_{\bar{\tau}\tau}:\N^*_{\bar{\tau}}\rightarrow \N_\tau$ that is $\Sigma_0^{(n_\tau)}$-preserving with respect to the language of hod premice such that
\begin{enumerate}[(a)]
\item $\bar{\tau}=\rm{cr}(\sigma^0_{\bar{\tau}\tau})$ and $\sigma^0_{\bar{\tau}\tau}(\bar{\tau})=\tau$.
\item $\sigma^0_{\bar{\tau}\tau}(p^*_{\bar{\tau}})=p_\tau$.
\item for each $\alpha\in p_\tau$, there is a generalized witness for $\alpha$ with respect to $\N_\tau$ and $p_\tau$ in the range of $\sigma_{\bar{\tau}\tau}$.
\end{enumerate}
\end{itemize}$\myqedhere$
\end{definition}

Note that if $\tau\in \S^0$, then $\N^*_\tau=\N_\tau$ and either crt$(E^{top}_{\N_\tau})\geq\kappa$ or $\rho_1^{\N_\tau}>\kappa$. Recall the definition of $(\mu_\tau,q_\tau)$, $p_\tau, d_\tau, \M_\tau$ for $\tau\in \S^1$ in Section \ref{sec:set_up}. Below, $m_\tau$ is $|q_\tau|$. We also let $r_\tau = d_\tau - q_\tau$ be the top part of $d_\tau$.
\begin{definition}\label{def:S1}\index{$\S^1$}
Suppose $\tau\in \S^1$. Let $B^1_\tau$ be the set of $\bar{\tau}\in\S^1\cap\tau$ satisfying:
\begin{itemize}
\item $(\mu_{\bar{\tau}},m_{\bar{\tau}}) = (\mu_\tau,m_{\tau})$.
\item There is a map $\sigma^1_{\bar{\tau}\tau}:\M_{\bar{\tau}}\rightarrow \M_\tau$ that is $\Sigma_0$-preserving with respect to the language of coherent structures such that
\begin{enumerate}[(a)]
\item $\bar{\tau}=\rm{cr}(\sigma^1_{\bar{\tau}\tau})$ and $\sigma^1_{\bar{\tau}\tau}(\bar{\tau})=\tau$.
\item $\sigma^1_{\bar{\tau}\tau}(q_{\bar{\tau}})=q_\tau$.
\item for each $\alpha\in q_\tau$, there is a generalized witness for $\alpha$ with respect to $\M_\tau$ and $q_\tau$ in the range of $\sigma^1_{\bar{\tau}\tau}$.
\end{enumerate}
\end{itemize}

Suppose in addition that either crt$(E^{top}_{\N^*_\tau})\geq\kappa$ or $\rho_1^{\N^*_\tau}>\kappa$, let $B^0_\tau$ be the set of $\bar{\tau}\in \S\cap \tau$ satisfying:

\begin{itemize}
\item $\N^*_{\bar{\tau}}$ is a hod premouse of the same type as $\N^*_\tau$. 
\item $n_\tau = n_{\bar{\tau}}$.
\item There is a map $\sigma^0_{\bar{\tau}\tau}:\N^*_{\bar{\tau}}\rightarrow \N^*_\tau$ that is $\Sigma_0^{(n_\tau)}$-preserving with respect to the language of hod premice such that
\begin{enumerate}[(a)]
\item $\bar{\tau}=\rm{cr}(\sigma^0_{\bar{\tau}\tau})$ and $\sigma^0_{\bar{\tau}\tau}(\bar{\tau})=\tau$.
\item $\sigma^0_{\bar{\tau}\tau}(p^*_{\bar{\tau}})=p_\tau$.
\item for each $\alpha\in p_\tau$, there is a generalized witness for $\alpha$ with respect to $\N^*_\tau$ and $p_\tau$ in the range of $\sigma^0_{\bar{\tau}\tau}$.
\end{enumerate}
\end{itemize}

Finally, if $B^0_\tau$ exists, let $\vec{B}_\tau=\{B^0_\tau,B^1_\tau\}$. Otherwise, let $\vec{B}_\tau=\{B^1_\tau\}$.$\myqedhere$

\end{definition}

As in \cite{schimmerling2004characterization}, it is easy to see that in both cases $\sigma_{\bar{\tau}\tau}, \sigma^0_{\bar{\tau}\tau}, \sigma^1_{\bar{\tau}\tau}$ (if exist) are uniquely determined, $\Sigma_0$ (and not $\Sigma_1$), and non-cofinal.  By \cite[Lemma 3.3]{schimmerling2004characterization}, for each $\tau\in \S$ such that $B^0_\tau$ is defined, and $\bar{\tau}\in B^0_\tau$, 
\begin{equation}\label{eqn:tail_agreement_0}
B^0_\tau\cap \bar{\tau} = B^0_{\bar{\tau}}-\textrm{min}B^0_{\tau}.
\end{equation}
And similarly, if $B^1_\tau$ is defined, then for all $\bar{\tau}\in B^1_\tau$,
\begin{equation}\label{eqn:tail_agreement_1}
B^1_\tau\cap \bar{\tau} = B^1_{\bar{\tau}}-\textrm{min}B^1_{\tau}.
\end{equation}
The following is the key lemma (cf. \cite[Lemma 3.5]{schimmerling2004characterization}).
\begin{lemma}\label{lem:key_lemma}
For each $\tau\in \S$ of uncountable cofinality, for $i\in\{0,1\}$, if $B^i_\tau$ is defined, then $B^i_\tau$ is a club subset of $\tau$ on a tail end. That is, there is a $\bar{\tau}<\tau$ such that $B^i_\tau - \bar{\tau}$ is closed and unbounded in $\tau$. If $i=0$, we can take $\bar{\tau}=0$.
\end{lemma}
Using the lemma and \ref{eqn:tail_agreement_0}, \ref{eqn:tail_agreement_1}, by the argument on \cite[pg 52-55]{schimmerling2004characterization} , we can construct a $\square'_{\kappa,2}$-sequence on $\S$. We summarize the construction next. First for $\tau\in\S$, for $i$ such that $B^i_\tau$ is defined, let
\begin{itemize}
\item $\tau^i(0)=\tau$;
\item $\tau^i(j+1)=\rm{min}$$(B^i_{\tau(j+1)})$;
\item $l^i_\tau = $ the least $j$ such that $B^i_{\tau(j)}=\emptyset$.
\end{itemize}
Now let
\begin{itemize}
\item $B^{i,*}=B^i_{\tau^i(0)}\cup \dots \cup B^i_{\tau^i(l^i_\tau-1)}$;
\item $\sigma^{i,*}_{\bar{\tau}\tau}=\sigma^i_{\tau^i(1)\tau^i(0)}\circ \dots \circ \sigma^i_{\tau^i(j)\tau^i(j-1)}\circ \sigma^i_{\bar{\tau}\tau^i(j)}$ whenever $\bar{\tau}\in B^{i,*}_\tau$ and $j$ is such that $\bar{\tau}\in B^i_{\tau(j)}$.
\end{itemize}

By the exact same proof as in \cite[Lemma 3.4]{schimmerling2004characterization}, we get the coherency of the $B^{i,*}_\tau$ sets.
\begin{lemma}\label{lem:coherence}
For $\tau\in \S$, for $i$ such that $B^i_\tau$ is defined, suppose $\bar{\tau}\in B^{i,*}_\tau$. Then $B^i_{\bar{\tau}}$ is defined and $B^{i,*}_{\bar{\tau}} = B^{i,*}_\tau\cap \bar{\tau}$.
\end{lemma}

For each $\tau\in\S$, for $i$ such that $B^i_\tau$ is defined, let $\beta^i_\tau$ be the least $\beta$ in $B^{i,*}_\tau\cup\{\tau\}$ such that $B^{i,*}_\tau-\beta$ is closed in $\tau$. Using Lemmata \ref{lem:key_lemma} and \ref{lem:coherence}, we easily get that letting 
\begin{equation}\label{eqn:C_star}
C^{i,*}_\tau = B^{i,*}_\tau - \beta^i_\tau,
\end{equation}
then for $\bar{\tau}\in \beta^{i,*}_\tau$, $\bar{\tau}\geq \beta_\tau$,
\begin{equation}
\beta^i_\tau=\beta^i_{\bar{\tau}} \textrm{ and } C^{i,*}_\tau\cap \bar{\tau} = C^{i,*}_{\bar{\tau}}. 
\end{equation}

Now note that if $C^{0,*}_\tau$ is defined, then o.t.$(C^{0,*}_\tau)$ may not be $\leq \kappa$, while if $C^{1,*}_\tau$ is defined then o.t.$(C^{0,*}_\tau)\leq \kappa$. As in \cite[pg 54-55]{schimmerling2004characterization}, we can shrink $C^{0,*}_\tau$ to a set $C^{0,'}_\tau\subseteq C^{0,*}_\tau$ such that
\begin{itemize}
\item o.t.$(C^{0,'}_\tau) \leq \kappa$;
\item $C^{0,'}_\tau$ is a closed subset of $\S\cap\tau$ and if cof$(\alpha)>\omega$, then $C^{0,'}_\tau$ is also unbounded in $\tau$;
\item $C^{0,'}_\tau\cap \bar{\tau} =C^{0,'}_{\bar{\tau}}$. 
\end{itemize}

So letting $\vec{C}'_\tau=\{C^{i,'}_\tau \ | \ i\in\{0,1\} \wedge C^{i,'}_\tau \textrm{ is defined}\}$, we get that the sequence $\langle \vec{C}'_\tau \ | \ \tau < \kappa^+ \rangle$ is a $\square'_{\kappa,2}$-sequence on $\S$. Since $\S$ is club subset of $\kappa^+$, by a standard combinatorial argument (cf. \cite{devlin1984constructibility}), the $\square'_{\kappa,2}$-sequence on $\S$ can be turned into a $\square_{\kappa,2}$-sequence. Our main task is to prove Lemma \ref{lem:key_lemma}. This will take up the rest of the section.

\begin{remark}
It's clear from \cite[pg 54-55]{schimmerling2004characterization}, Definitions \ref{def:S0} and \ref{def:S1} and Proposition \ref{prop:N_tau_def} that the square sequence $\square_{\kappa,2}$ is definable from $\kappa$ in $\P$ and the definition is uniform in $\kappa$. $\myqedhere$
\end{remark}

\subsection{When $\tau\in \S^0$}
Fix $\tau\in \S^0$. Assume $\tau$ is a limit point of $\S$ uncountable cofinality. Recall $B^0_\tau$ is defined to be the set of $\bar{\tau}\in \S$ such that
\begin{itemize}
\item $n_\tau = n_{\bar{\tau}}$.
\item $\N^*_{\bar{\tau}}$ is a hod premouse of the same type as $\N_\tau$.
\item There is an embedding $\sigma^0_{\bar{\tau}\tau}:\N^*_{\bar{\tau}}\rightarrow \N_\tau$ such that $\sigma^0_{\bar{\tau}\tau}$ is $\Sigma_0^{(n_\tau)}$-preserving (in the language of hod premice) and
\begin{enumerate}[(a)]
\item $\bar{\tau}=\rm{cr}(\sigma^0_{\bar{\tau}\tau})$ and $\sigma^0_{\bar{\tau}\tau}(\bar{\tau})=\tau$.
\item $\sigma^0_{\bar{\tau}\tau}(p^*_{\bar{\tau}})=p_\tau$, where recall $p^*_{\bar{\tau}}$ is the standard parameter of $\N^*_{\bar{\tau}}$.
\item for each $\alpha\in p_\tau$, there is a generalized witness for $\alpha$ with respect to $\N_\tau$ and $p_\tau$ in the range of $\sigma_{\bar{\tau}\tau}$.
\end{enumerate}
\end{itemize}
To simplify the notation, let $D$ denote $B^0_\tau$ and $\sigma_{\bar{\tau},\tau}$ denote $\sigma^0_{\bar{\tau},\tau}$. 

\begin{lemma}\label{lem:0_unbounded}
$D$ is unbounded in $\tau$.
\end{lemma} 
\begin{proof}
Given $\tau' < \tau$, we find $\tilde{\tau}\geq \tau'$ in $D$. In $\P$, form an countable elementary hull of $\{\N_\tau,\tau',\S\}$ in $H_{(\kappa^{++})}$ (in which everything relevant is present). Let $H$ be the transitive collapse of the hull and $\sigma_0:H \rightarrow H_{\kappa^{++}}$ be the uncollapse map. Set:
\begin{itemize}
\item $\bar{x}=\sigma_0^{-1}(x)$ for any $x$ in range of $\sigma_0$,
\item $\sigma=\sigma_0\rest \bar{\N_\tau} : \bar{\N_\tau}\rightarrow \N_\tau$,
\item $\tilde{\tau} = \rm{sup}$$(\sigma''\bar{\tau})$.
\end{itemize}
Note that since $\tau\in \S^0$, $\N_\tau =\N^*_\tau$ and either cr$(E^{\rm{top}}_{\N_\tau})\geq \kappa$ or $\omega\rho^1_{\N_\tau}>\kappa$. Set $n=n_\tau$. Using the interpolation lemma (Lemma \cite[Lemma 1.2]{schimmerling2004characterization}), we can find a map $\tilde{\sigma}: \bar{\N_\tau}\rightarrow \tilde{\N}$ which is $\Sigma_0^{(n)}$-preserving and cofinal (the map $\tilde{\sigma}$ is the ultrapower map via the $(\rm{cr}(\sigma),\tilde{\tau})$-extender derived from $\sigma$). Note that $\tilde{\tau} = (\kappa^+)^{\tilde{\N}}$. Also, by the interpolation lemma, there is a map $\sigma': \tilde{\N}\rightarrow \N_\tau$ satisfying $\sigma'\rest \tilde{\tau} = \rm{id}$, $\sigma'(\tilde{\tau}) =\tau$, and $\sigma' \circ \tilde{\sigma} = \sigma$.

We have that
\begin{itemize}
\item $\tilde{\N}$ is a hod premouse of the same type as $\N_\tau$.
\item $\tilde{\N}$ is sound.
\item $\omega\rho^\omega_{\tilde{\N}} = \omega\rho^{n+1}_{\tilde{\N}}\leq \kappa$.
\end{itemize}

The above follow from the proof of \cite[Lemma 3.7]{schimmerling2004characterization} for the most part, except for the first item in the case when $\N_\tau$ is $B$-active. In this case, the first item follows from Lemma \ref{lem:preservehodpairs} (or see \cite[Lemma 2.36]{trang2013}) and hull condensation of $\Sigma$.\footnote{Indexing branches using the $\mathfrak{B}$-operator allows the proof of Lemma \ref{lem:preservehodpairs} and \cite[Lemma 2.36]{trang2013} to go through in this situation. The traditional approach to indexing branches does not seem to imply that $\tilde{N}$ is a hod premouse.}

It remains to see that $\tilde{\N}$ is indeed $\N^*_{\tilde{\tau}}$. We apply Theorem \ref{thm:condensation_lemma}. (a) cannot hold since the map $\sigma'$ is a $\Sigma^{(n)}_0$, non-cofinal map and hence cannot be a core map. (c) cannot hold because $\tilde{\N}$ is sound. So (b) holds; in fact, $\tilde{\N}$ is a strict segment of $\N_\tau$ because $\tilde{\tau} = \rm{cr}(\sigma') = (\kappa^+)^{\tilde{\N}} < \tau = (\kappa^+)^{\N_\tau}$. This easily implies $\tilde{\N}  = \N^*_{\tilde{\tau}}$.
\end{proof}
\begin{lemma}\label{lem:0_closed}
$D$ is a closed subset of $\tau$.
\end{lemma}
\begin{proof}
Let $\tilde{\tau}$ be a limit point of $D$. We show that $\tilde{\tau}\in D$. Form the direct limit $\langle \tilde{\N}, \sigma_{\bar{\tau}\tilde{\tau}} \ | \ \bar{\tau}\in D\cap \tilde{\tau} \rangle$ of the system $\langle \N^*_{\bar{\tau}}, \sigma_{\tau^*\bar{\tau}} \ | \ \tau^* \leq \bar{\tau} \wedge \tau^*, \bar{\tau} \in D\cap \tilde{\tau} \rangle$. The direct limit is well-founded and there is a $\Sigma_0$ embedding $\sigma: \tilde{\N}\rightarrow \N_\tau$ (defined by $\sigma(\sigma_{\bar{\tau}\tilde{\tau}}(x))=\sigma_{\bar{\tau},\tau}(x)$). It is easy to check that:
\begin{enumerate}[(a)]
\item $\sigma\circ \sigma_{\bar{\tau}\tilde{\tau}} = \sigma_{\bar{\tau}\tau}$.
\item $\tilde{\tau}=\sigma_{\bar{\tau}\tilde{\tau}}(\bar{\tau})$, $\sigma_{\tilde{\tau}\tau}(\tilde{\tau})=\tau$, and $\tilde{\tau} = \rm{cr}(\sigma)$.
\item $\sigma$ is $\Sigma_0^{(n)}$ preserving where $n = n_\tau$ (with respect to the language of coherent structures).
\end{enumerate} 

We need to see that $\tilde{\N} = \N^*_{\tilde{\tau}}$. First, we show that $\tilde{\N}$ is a hod premouse of the same type as $\N_\tau$. Note that $\Pi_2$-properties which hold on a tail end are upward preserved under direct limit maps (cf. \cite[pg 8-9]{schimmerling2004characterization}). Furthermore, $\N^*_{\bar{\tau}}$ is of the same type as $\N_\tau$ for each $\bar{\tau}\in D\cap \tau$. So $\tilde{\N}$ is of the same type as $\N_\tau$ (as either a passive hod premouse, or a $B$-active hod premouse, or an $E$-active hod premouse with cr$(E^{\rm{top}}_{\tilde{\N}}) > \mu$, in which case $\tilde{\N}$ is of type $A$, or else an $E$-active hod premouse with cr$(E^{\rm{top}}_{\tilde{\N}}) = \mu$, in which case $\omega\rho^1_{\tilde{\N}}>\kappa$; these statements can be expressed in a $\Pi_2$-fashion).

Recall that for $\bar{\tau}\in D$, we use $\tilde{h}_{\bar{\tau}}$ to denote $\tilde{h}^{(n_{\bar{\tau}}+1)}_{\N^*_{\bar{\tau}}}$, the $\Sigma_1^{(n_{\bar{\tau}})}$-Skolem function of $\N^*_{\bar{\tau}}$. Here note that $n_{\bar{\tau}} = n_\tau = n$. Let $\tilde{p} = \sigma_{\bar{\tau}\tilde{\tau}}(p^*_{\bar{\tau}})$ for $\bar{\tau}\in D\cap \tilde{\tau}$. Given any $x\in \tilde{\N}$, there is $\bar{\tau}\in D\cap \tilde{\tau}$ and $\bar{x}\in \N^*_{\bar{\tau}}$ such that $x = \sigma_{\bar{\tau}\tau}(\bar{x})$. There is $\xi<\kappa$ such that $\bar{x} = \tilde{h}_{\bar{\tau}}(\xi,p_{\bar{\tau}})$. This $\Sigma^{(n)}_1$-statement is preserved by $\sigma_{\bar{\tau}\tilde{\tau}}$, so $x = \tilde{h}^{n+1}_{\tilde{\N}}(\xi,\tilde{p})$. So $\tilde{\N} = \tilde{h}^{n+1}_{\tilde{\N}}(\kappa\cup\{\tilde{p}\})$.

This gives $\omega\rho^{n+1}_{\tilde{\N}} =\omega\rho^{\omega}_{\tilde{\N}} \leq\kappa$. But $\kappa$ is a cardinal in $\P$, so we indeed have equality. For each $\alpha \in p_\tau$, there is a generalized witness for $\alpha$ with respect to $(\N_\tau,p_\tau)$ in range of $\sigma$. This is because rng$(\sigma)$ contains rng$(\sigma_{\bar{\tau},\tau})$ for any $\bar{\tau}\in D\cap \tilde{\tau}$ and rng$(\sigma_{\bar{\tau},\tau})$ contains such a witness. This takes care of (c) in the definition of $D$. This easily implies that $\tilde{\N}$ is sound and $\tilde{p}$ is the standard paramter of $\tilde{\N}$. We can now apply Theorem \ref{thm:condensation_lemma} as in the proof of Lemma \ref{lem:0_unbounded} to conclude that $\tilde{\N}=\N^*_{\tilde{\tau}}$.
\end{proof}

Lemmata \ref{lem:0_unbounded}, \ref{lem:0_closed} together complete the proof of Lemma \ref{lem:key_lemma} in the case $\tau\in\S^0$.

\subsection{When $\tau\in \S^1$}
Fix $\tau\in \S^1$ a limit point of $\S$ of uncountable cofinality. If $B^0_\tau$ is defined, then as in the previous section, using the fact that crt$(E^{top}_{\N^*_\tau})\geq\kappa$ or $\rho_1^{\N^*_\tau}>\kappa$, we can show that $B^0_\tau$ is closed and unbounded in $\tau$. So let us now focus on the case $B^1_\tau$ is defined. Define $D\subset\tau$ to be the set of $\bar{\tau}\in \S$ such that
\begin{itemize}
\item $(\mu_\tau,q^*_{\bar{\tau}})$ is a strong divisor of $\N_{\bar{\tau}}$ where $q^*_{\bar{\tau}}$ is the bottom segment of $p_{\bar{\tau}}$ of length $m_\tau$ (recall $m_\tau$ is the length of $q_\tau$).
\item Letting $\M^*_{\bar{\tau}}$ be the protomouse of $\N_{\bar{\tau}}$ associated with $(\mu_\tau,q^*_{\bar{\tau}})$, there is a map $\sigma_{\bar{\tau}\tau}:\M^*_{\bar{\tau}}\rightarrow \M_\tau$ that is $\Sigma_0$-preserving (with respect to the language of coherent structures) such that
\begin{enumerate}[(a)]
\item $\bar{\tau}=\rm{cr}(\sigma_{\bar{\tau}\tau})$ and $\sigma_{\bar{\tau}\tau}(\bar{\tau})=\tau$.
\item $\sigma_{\bar{\tau}\tau}(q^*_{\bar{\tau}})=q_\tau$.
\item for each $\alpha\in q_\tau$, there is a generalized witness for $\alpha$ with respect to $\M_\tau$ and $q_\tau$ in the range of $\sigma_{\bar{\tau}\tau}$ (in the language of coherent structures).
\end{enumerate}
\end{itemize}
We will show that there is some $\bar{\tau}<\tau$ such that $B^1_\tau-\bar{\tau}=D-\bar{\tau}$. Part of this is to show that for all sufficiently large $\bar{\tau} \in D$, $(\mu_\tau,q^*_{\bar{\tau}}) = (\mu_{\bar{\tau}},q_{\bar{\tau}})$.
\begin{lemma}\label{lem:1_unbounded}
$D$ is unbounded in $\tau$.
\end{lemma} 
\begin{proof}
Let $\tau' < \tau$. As before, we find $\tilde{\tau}\in D$ above $\tau'$. Since protomice are present, we carry out the argument in the language of coherent structures. 

We let $\sigma_0, H$ be defined as in Lemma \ref{lem:0_unbounded}. Again, we denote $\bar{x}$ for $\sigma_0^{-1}(x)$. We let $\sigma:\bar{\M_\tau}\rightarrow \M_\tau$ and $\tilde{\tau} = \rm{sup}\ \sigma''\bar{\tau}$. As before, $\tau' \leq \tilde{\tau} < \tau$. Let $\tilde{\sigma}:\bar{\M_\tau}\rightarrow \tilde{\M}$ be the $(\rm{cr}(\sigma),\tilde{\tau})$-ultrapower map derived from $\sigma$ and $\sigma':\tilde{\M}\rightarrow \M_\tau$ be the factor map. As in \cite[Lemma 3.10]{schimmerling2004characterization}, we have:
\begin{itemize}
\item $\tilde{\sigma}(\bar{\kappa},\bar{\tau}) = (\kappa,\tilde{\tau})$.
\item cr$(\sigma')=\tilde{\tau}$ and $\sigma(\tilde{\tau})=\tau$.
\item $h_{\tilde{\M}}(\kappa\cup\{\tilde{q}\}) = \tilde{\M}$ where $\tilde{q}=\tilde{\sigma}(\bar{q_\tau})$; in other words, $\tilde{\M}$ is $\Sigma_1$-generated by $\kappa\cup\{\tilde{q}\}$.
\item $\omega\rho^\omega_{\tilde{\M}}=\omega\rho^1_{\tilde{\M}} = \kappa$ and $\tilde{q}\in R_{\tilde{\M}}$, the set of very good parameters for $\tilde{\M}$.
\item The range of $\tilde{\sigma}$ contains a generalized solidity witness for $\alpha$ with respect to $(\M_\tau,q_\tau)$ for each $\alpha\in q_\tau$.
\item $\tilde{q}=p_{\tilde{\M}}$ and $\tilde{\M}$ is solid and sound.
\end{itemize}

Note that as in Lemma \ref{lem:0_unbounded}, $\tilde{\sigma}$ is $\Sigma_0$ (but not $\Sigma_1$) and is not cofinal. This implies that $\tilde{\M}$ is a protomouse, even if $\M_\tau$ is a hod premouse (in which case, $\M_\tau = \N_\tau = \N_\tau^*$ is pluripotent).

We show $\tilde{\M} = \N_{\tilde{\tau}}(\mu_\tau,\tilde{q})$. Let $\R_0, \R_1$ be the hod premice associated with $\M_\tau, \tilde{\M}$, respectively. We have that $\R_0 = \rm{Ult}$$_n(\N_0^*,F)$, where $F$ is the top extender (fragment) of $\M_\tau$ and $\N_0^*$ is largest (strict) segment of $\M_\tau$ such that $\omega\rho^\omega_{\N_0^*}\leq\rm{cr}$$(F)=\mu$ and $F$ measures all sets in $\N_0^*$ if exists, otherwise, $\N_0^*= \M_\tau$;\footnote{The first case is the case $\M_\tau$ is a protomouse and the second case is when $\M_\tau$ is a pluripotent level.} in the other case, $\R_1 = \rm{Ult}$$(\N_1^*,\tilde{F})$, where $\tilde{F}$ is the top extender (fragment) of $\tilde{\M}$ and $\N_1^*$ is the largest (strict) segment of $\tilde{\M}$ (equivalently, of $\M_\tau$) such that $\omega\rho^{\omega}_{\N_1^*}\leq\rm{cr}$$(\tilde{F})=\mu$ and $\tilde{F}$ measures all sets in $\N_1^*$. Let $\pi_i:\N_i^*\rightarrow \R_i$ be the ultrapower maps and $\pi_2: \R_1\rightarrow \pi_0(\N_1^*)$ be the factor map
\begin{center}
$\pi_2(\pi_1(f)(a)) = \pi_0(f)(\tilde{\sigma}(a))$.
\end{center}
Note that $\pi_2\rest\lambda^{\R_1} = \tilde{\sigma}\rest\lambda^{\R_1}$ and therefore crt$(\pi_2) = \tilde{\tau}$. Furthermore, if $\omega\rho^\omega_{\N^*_i}=\mu$, then $\N^*_i\lhd \P^b$ and therefore if $\N^*_i$ is $E$-active, $E^{\textrm{top}}_{\N^*_i}$ has critical point $>\mu$ because $\mu$ is a strong cutpoint of $\P^b$. This easily gives that $\R_i$ is of the same type as $\N_i^*$ (as potential premice).

Note that $p_{\R_1}=\pi_1(p_{\N_1^*})\cup p_{\tilde{\M}}$ (cf. \cite[Lemma 2.16, 2.19]{schimmerling2004characterization}). In the case $\N^*_\tau\neq\N_\tau$, and hence $\mu_\tau=\mu$, $\pi_1(p_{\N_1^*})$ is the part of $p_{\R_1}$ above $\pi_1(\mu)$, the supremum of $\R_1$'s layer Woodin cardinals, and $p_{\tilde{\M}}$ is the part below $\pi_1(\mu)$. 

The argument in \cite[Lemma 3.10]{schimmerling2004characterization} then shows that $(\mu_\tau,\tilde{q})$ is a strong divisor of $\R_1$.\footnote{The proof of this fact does not depend on whether $\mu_\tau > \mu$.} To show $\tilde{\M}=\N_{\tilde{\tau}}(\mu_\tau,\tilde{q})$, we show $\R_1 = \N_{\tilde{\tau}}$. This then will show $\tilde{\tau}\in D$ as desired. There are two cases to consider.
\\
\\
\noindent \textbf{Case 1.} $\N_\tau = \N^ *_\tau$.

If cr$(F) = \rm{cr}$$(\tilde{F}) > \mu$, then it is easy to see that $\P^b \lhd \R_0, \R_1$. Note that in this case, $\R_0 = \N_\tau = \N^*_\tau$ (see \cite[Section 2]{schimmerling2004characterization}). So we can apply Theorem \ref{thm:condensation_lemma} as in the proof of Lemma \ref{lem:0_unbounded} and conclude that $\R_1 = \N_{\tilde{\tau}} = \N_{\tilde{\tau}}^*$. Now suppose cr$(F)=\mu$ (so $\mu_\tau = \mu$). Recall from the discussion above that we know $(\mu_\tau,\tilde{q})$ is a strong divisor of $\R_1$ and $\tilde{q}$ is the bottom part of the standard parameter of $\R_1$ below $\pi_1(\rm{cr}$$(\tilde{F}))$. We show that $\R_1 = \N_{\tilde{\tau}}\neq \N_{\tilde{\tau}}^*$ by the following claims. We also will get then that $(\mu_\tau,\tilde{q})=(\mu_{\tilde{\tau}},q_{\tilde{\tau}})$ in this case.

Let $\gamma_\tau$ be defined as in the paragraphs before \ref{prop:basic_facts} for $\N_\tau$; let $\gamma_{\tilde{\tau}},\tilde{\gamma}$ be defined similarly for $\N^*_{\tilde{\tau}}, \R_1$, respectively. Let $\Lambda$ be $\R_0$'s iteration strategy.
\begin{claim}\label{claim:gamma_same}
$\tilde{\gamma} = \gamma_{\tilde{\tau}}$.
\end{claim}
\begin{proof}
Suppose not. Assume $\tilde{\gamma} < \gamma_{\tilde{\tau}}$ (the other case is similar). Let $E$ be least on the extender sequence of $\N_{\tilde{\tau}}$ (equivalently of $\N^*_{\tilde{\tau}}$) such that
\begin{itemize}
\item cr$(E) = \mu$,
\item lh$(E)\geq \tilde{\gamma}$.
\end{itemize}
Let $\S = \rm{Ult}$$(\R_0,E)$. Note that $\tilde{\gamma}$ is a cutpoint of $\S$ and $i_E(\mu)$ is a limit of $\Gamma$-full Woodin cardinals above $\tilde{\gamma}$. By $\sf{SMC}$ in $\Gamma$, we can conclude that $\R'\in \S$, where $\R'$ is a sound hod mouse extending $\R_1|\tilde{\gamma}$, $\tilde{\tau}=(\kappa^+)^{\R'}$, $\tilde{\gamma}$ as a cutpoint of $\R'$, and $\R'$ projects to $\kappa$. \footnote{\label{fnt:in} By genericity iterations, without loss of generality, we may assume that a real witnessing the Wadge reduction of $\Lambda^{\pi_2}$ to $\Lambda$ is generic over $\S$. In $\S$'s derived model at $i_E(\mu)$, we can find $\R_1$. This means, in the derived model of $\S$, there is some hod mouse $\R$ extending $\R_1|\tilde{\gamma}$, having $\tilde{\tau}=\kappa^+$, $\tilde{\gamma}$ as a cutpoint, and projects to $\kappa$; furthermore, we can demand that $(\mu_\tau,\tilde{q})$ is a strong divisor of $\R$ and $\tilde{q}$ is the bottom part of the standard parameter of $\R$ below the supremum of $\R$'s layer Woodin cardinals. Let $\Omega$ be the Wadge-minimal pointclass that has a pointclass generator with these properties. Note that this determines the unique pointclass generator $\S_\Omega$ for $\Omega$. This implies that $\S_\Omega\in \S$.}

Fix $\R'\in \S$ as above. $\R'$ defines a surjection $f$ from $\kappa$ onto $\tilde{\tau}$. Since $\R'\in \S$, $f\in \S$. This contradicts the fact that $\S \vDash \tilde{\tau} = \kappa^+$.
\end{proof}
\begin{claim}\label{claim:model_same}
There is a pointclass $\Omega$ with pointclass generator a sound hod mouse that projects to $\kappa$, extends $\P|\tilde{\gamma}$, having $\tilde{\tau}=\kappa^+$, $\tilde{\gamma}$ as a cutpoint, and the set of layer Woodin cardinals above $\tilde{\gamma}$ has limit order type. $\R_1$ is the generator for the Wadge minimal such pointclass.
\end{claim}
\begin{proof}
Clearly, such $\Omega$ exists since the pointclass generated by $\R_1$ is such. Let $\Omega_0$ be the pointclass $\R_1$ generates and $\Omega_1$ be the minimal pointclass satisfying the hypothesis of the claim. Let $(\N,\Psi)$ generate $\Omega_1$ with the properties in the statement of the claim. Note that at this point, we know $\R_1$ and $\N$ are sound, projects to $\kappa$, extends $\P|\gamma_{\tilde{\tau}}$, satisfies $\kappa^+ = \tilde{\tau}$, have $\gamma_{\tilde{\tau}}$ as cutpoint, and the set of layer Woodin cardinals above  $\gamma_{\tilde{\tau}}$ of both models is of limit order type. 

We claim that $\Omega_0 = \Omega_1$. Suppose for contradiction that $\Omega_0 \subsetneq \Omega_1$ (the other case is similar). Then, using $\mathbb{R}$-genericity iteration above  $\gamma_{\tilde{\tau}}$ and elementarity, in the derived model of $\N$ (at the supremum of its Woodin cardinals) there is a pointclass with a generator $\S$ that is sound, projects to $\kappa$, extends $\P|\gamma_{\tilde{\tau}}$, satisfies $\kappa^+ = \tilde{\tau}$, and have $\gamma_{\tilde{\tau}}$ as cutpoint. Some such $\S$ is in $\N$ by a similar argument as in Footnote \ref{fnt:in}. This implies as in Claim \ref{claim:gamma_same} that $\tilde{\tau}$ is not a cardinal in $\N$. Contradiction. We have shown $\Omega_0 = \Omega_1$.

Now we can compare $\R_1$ against $\N$. Since $\Omega_0 = \Omega_1$, both models are $\kappa$-sound, projects to $\kappa$, and $\kappa<\gamma_{\tilde{\tau}}$, just as Lemma \ref{prop:unique_ptclass_gen}, we conclude that $(\N,\Psi) = (\R_1, \Lambda)$.
\end{proof}

Using the claims and  the fact that $(\mu_\tau,\tilde{q})$ is a strong divisor of $\R_1$ (note that max$(\tilde{q}) < (\gamma_{\tilde{\tau}}^+)^{\R_1}$ and $\tilde{q}$ is the bottom part below $\pi_1(\rm{cr}$$(\tilde{F}))$ of the standard parameter of $\R_1$) we easily verify that $\R_1=\N_{\tilde{\tau}}$ and hence $\tilde{\M} = \N_{\tilde{\tau}}(\mu_\tau,\tilde{q})$. Hence $\tilde{\tau}\in D$ as desired.
\\
\\
\noindent \textbf{Case 2.} $\N_\tau \neq \N^*_\tau$.

In this case, $\mu_\tau = \mu$. As above, $\R_0 = \N_\tau$ and $(\mu_\tau,q_\tau)$ is a strong divisor of $\N_\tau$. We aim to show that $\R_1 = \N_{\tilde{\tau}}$. As above, $(\mu,\tilde{q})$ is a strong divisor of $\R_1$ by the proof of \cite[Lemma 3.10]{schimmerling2004characterization}. Furthermore, max$(\tilde{q}) < (\tilde{\gamma}^+)^{\R_1}=(\gamma_{\tilde{\tau}}^+)^{\R_1}$ and $\tilde{q}$ is the bottom part of the standard parameter of $\R_1$ below $\R_1$'s limit of layer Woodin cardinals $\pi_1(\rm{cr}$$(\tilde{F}))$. This easily implies, using Claim \ref{claim:model_same}, that $\N_{\tilde{\tau}}\neq \N_{\tilde{\tau}}^*$, $\R_1 = \N_{\tilde{\tau}}$, $\mu=\mu_{\tilde{\tau}}$, $\tilde{q} = q_{\tilde{\tau}}$ and hence $\tilde{\M} = \N_{\tilde{\tau}}(\mu_{\tilde{\tau}},q_{\tilde{\tau}})$. So $\tilde{\tau}\in D$ as desired.
\end{proof}

\begin{lemma}\label{lem:1_closed}
$D$ is a closed subset of $\tau$.
\end{lemma}
\begin{proof}
Let $\tilde{\tau}$ be a limit point of $D$. We show that $\tilde{\tau}\in D$. As in Lemma \ref{lem:0_closed}, form the direct limit $\langle \tilde{\M}, \sigma^1_{\bar{\tau}\tilde{\tau}} \ | \ \bar{\tau}\in D\cap \tilde{\tau} \rangle$ of the system $\langle \M^*_{\bar{\tau}}, \sigma^1_{\tau^*\bar{\tau}} \ | \ \tau^* \leq \bar{\tau} \ \wedge \ \{\tau^*, \bar{\tau}\} \subset D\cap \bar{\tau} \rangle$. The direct limit is well-founded (so we identify $\tilde{\M}$ with its transitive collapse) and there is a $\Sigma_0$ embedding $\sigma: \tilde{\M}\rightarrow \M_\tau$ (defined by $\sigma(\sigma^1_{\bar{\tau}\tilde{\tau}}(x))=\sigma^1_{\bar{\tau},\tau}(x)$). It is easy to check that (cf. \cite[Lemma 3.11]{schimmerling2004characterization}):
\begin{enumerate}[(a)]
\item $\tilde{\M}$ is a coherent structure.
\item $\sigma\circ \sigma^1_{\bar{\tau}\tilde{\tau}} = \sigma^1_{\bar{\tau}\tau}$.
\item $\tilde{\tau}=\sigma^1_{\bar{\tau}\tilde{\tau}}(\bar{\tau})$, $\sigma^1_{\tilde{\tau}\tau}(\tilde{\tau})=\tau$, and $\tilde{\tau} = \rm{cr}(\sigma)$.
\item $h_{\tilde{\M}}(\kappa\cup\{\tilde{q}\})=\tilde{\M}$ where $\tilde{q}=\sigma^1_{\bar{\tau}\tilde{\tau}}(q^*_{\bar{\tau}})$, so $\omega\rho^\omega_{\tilde{\M}}=\omega\rho^1_{\tilde{\M}}=\kappa$ and $\tilde{q}\in R_{\tilde{\M}}$.
\item For every $\alpha \in q_\tau$, there is a generalized witness for $\alpha$ with respect to $(\M_\tau,q_\tau)$ in the range of $\sigma$. Hence $\tilde{q} = p_{\tilde{\M}}=\sigma^{-1}(q_\tau)$ and $\tilde{\M}$ is sound and solid.
\end{enumerate} 
The first four clauses are clear. The last follows from the fact that rng$(\sigma)$ contains rng$(\sigma^1_{\bar{\tau},\tau})$ for sufficiently large $\bar{\tau}<\tau$ and rng$(\sigma^1_{\bar{\tau},\tau})$ has all relevant generalized witnesses. 


Note that $\tilde{\M}$ is always a protomouse (this is because $\sigma$ is not cofinal). If $\mu_\tau > \mu$ (or equivalently $\N_\tau = \N_\tau^*$), we can appeal to the proof of \cite[Lemma 3.11]{schimmerling2004characterization} to get that $\tilde{\M} = \N_{\tilde{\tau}}(\mu_\tau,\tilde{q})$ and $(\mu_\tau,\tilde{q})$ is a strong divisor of $\tilde{\M}$. Otherwise, the same conclusion follows from the proof of Claim \ref{claim:model_same}.

The previous paragraph gives $\tilde{\tau}\in D$ as desired.
\end{proof}
\begin{lemma}\label{lem:1_tail}
There is a $\bar{\tau}<\tau$ such that for all $\tau'\in D-\bar{\tau}$, $(\mu_\tau,q^*_{\tau'})=(\mu_{\tau'},q_{\tau'})$. Consequently, $B^1_\tau-\bar{\tau}=D-\bar{\tau}$.
\end{lemma}
\begin{proof}
We need to prove that there is $\bar{\tau}<\tau$ such that for every $\tau'\in D-\bar{\tau}$, $(\mu_\tau,q^*_{\bar{\tau}}) = (\mu_{\tau'},q_{\tau'})$. Assume for contradiction that there is a sequence $\langle \tau_i \ | \ i<\delta\rangle$ that is increasing, cofinal in $\tau$ such that $(\mu_{\tau_i},q_{\tau_i}) \neq (\mu_\tau,q^*_{\tau_i})$. We may assume without loss of generality that the sequence $\langle \mu_{\tau_i} \ | \ i<\delta\rangle$ is monotonic and all $q_{\tau_i}$'s have the same length, say $m$.

If $\mu_\tau = \mu$, then we claim that for each $i<\delta$, $(\mu_{\tau_i},q_{\tau_i}) = (\mu_\tau,q^*_{\tau_i})$. This follows from the proof of Lemma \ref{lem:1_unbounded}, where we prove that in this case,  for each $i<\delta$, $\N_{\tau_i}\neq \N^*_{\tau_i}$ and $\mu_{\tau_i} =\mu = \mu_\tau$ and $q_{\tau_i} = q^*_{\tau_i}$. This contradicts the assumption that $(\mu_{\tau_i},q_{\tau_i}) \neq (\mu_\tau,q^*_{\tau_i})$. So we must have that $\mu_\tau > \mu$, so $\N_\tau = \N^*_\tau$. This implies that for each $i<\delta$, $\N_{\tau_i} = \N^*_{\tau_i}$ (again, by remarks in Section \ref{sec:set_up} and the argument in Lemma \ref{lem:1_unbounded}). So it must be the case then that $\mu_{\tau_i}>\mu$ (note that since $\N^*_{\tau_i}=\N_{\tau_i}$, $\N^*_{\tau_i}$ cannot have strong divisors of the form $(\mu,q)$ for some $q$) and so $(\mu_{\tau_i},q_{\tau_i})$, by definition, is the canonical strong divisor of $\N_{\tau_i}$.

By the definition of $(\mu_{\tau_i},q_{\tau_i})$, each $q_{\tau_i}$ is a bottom part of $q^*_{\tau_i}$, say $q^*_{\tau_i}= q_{\tau_i}\cup s_{\tau_i}$ ($s_{\tau_i}$ may be empty). Recall we have shown $\mu_\tau,\mu_{\tau_i}>\mu$ (so we can freely quote results of \cite[Section 2.4 and Lemma 3.12]{schimmerling2004characterization} in the arguments that follow). Now we observe that $\mu_{\tau_i}>\mu_\tau$ for all $i<\delta$.  This is because the argument in \cite[Lemma 3.12]{schimmerling2004characterization} shows: if $q_{\tau_i}=q^*_{\tau_i}$, then $\mu_{\tau_i}$ must be $> \mu_\tau$ by maximality of $\mu_{\tau_i}$ for $\N_{\tau_i}$ and the assumption that $(\mu_\tau,q^*_{\tau_i})\neq (\mu_{\tau_i},q_{\tau_i})$; otherwise, $q_{\tau_i}$ is a strict bottom segment of $q^*_{\tau_i}$, and hence \cite[Lemma 2.26]{schimmerling2004characterization} shows that no $\nu\leq \mu_\tau$ is such that $(\nu,q_{\tau_i})$ is a strong divisor of $\N_{\tau_i}$. 

Set for some (equivalently for all sufficiently large) $i<\delta$, $q = \sigma_{\tau_i\tau}(q_{\tau_i})$, $s = \sigma_{\tau_i\tau}(s_{\tau_i})$, $r = r_\tau$, $\nu = \rm{sup}$$_{i<\delta} \mu_{\tau_i}$. Now $(\nu,q)$ is a divisor of $\N_\tau$ (see \cite[Lemma 3.12]{schimmerling2004characterization}). Since $\nu > \mu_\tau>\mu$, $(\nu,q)$ cannot be a strong divisor of $\N_\tau$. Then a calculation as in \cite[Lemma 3.12]{schimmerling2004characterization} shows that for some $i<\delta$, $(\mu_{\tau_i},q_{\tau_i})$ is not a strong divisor of $\N_{\tau_i}$. Contradiction.
\end{proof}
Lemmata \ref{lem:1_unbounded}, \ref{lem:1_closed}, \ref{lem:1_tail} together complete the proof of Lemma \ref{lem:key_lemma} in the case $\tau\in \S^1$. This finishes the construction of our $\square_{\kappa,2}$ sequence.

\chapter{$\sf{LSA}$ from $\sf{PFA}$}\label{chap:lsa_from_pfa}

For a cardinal $\kappa$, let $\powerset_0(\kappa)=\kappa$; $\powerset_{n+1}(\kappa) = 2^{\powerset_n(\kappa)}$ for all $n<\omega$. 
\begin{definition}\label{def:coherent_seq}
A sequence $\langle \vec{C}_\alpha \ | \ \alpha \in \lambda \rangle$ is a $\square(\kappa,\lambda)$ sequence if it satisfies the following properties.
\begin{enumerate}[(i)]
\item $0<|\vec{C}_\alpha| < \kappa$ for all limit $\alpha\in \lambda$.
\item $C\subseteq \alpha$ is club in $\alpha$ for all limit $\alpha\in \lambda$ and $C\in \vec{C}_\alpha$.
\item $C\cap \beta\in \vec{C}_\beta$ for all limit $\alpha \in \lambda$, $C\in \vec{C}_\alpha$ and $\beta\in\rm{Lim}$$(C)$.
\item There is no club $D\subseteq \lambda$ such that $D\cap \alpha\in \vec{C}_\alpha$ for all $\alpha\in \rm{Lim}$$(D)$.
\end{enumerate}
We say that $\square(\kappa,\lambda)$ holds if a $\square(\kappa,\lambda)$-sequence exists. $\myqedhere$
\end{definition}
Clearly, $\square_{\lambda,<\kappa}$ implies $\square(\kappa,\lambda^+)$ and if $\kappa\leq\kappa'$, then $\square(\kappa,\lambda)$ implies $\square(\kappa',\lambda)$. $\square(2,\lambda)$ is $\square(\lambda)$. The following is the main result of this chapter.
\begin{theorem}\label{thm:square_lsa}
Suppose $\kappa$ is a  cardinal such that $\kappa^\omega = \kappa$ and $2^{2^{\aleph_0}}\leq \kappa$. Suppose there is a regular cardinal $ \gamma\in [\omega_2 , \kappa)$ such that for all $\alpha \in [\gamma, (\powerset_4(\kappa))^{+}]$, $\neg \square(3,\alpha)$. Then there is a model $M$ containing $\rm{OR}\cup\mathbb{R}$ such that $M\vDash \sf{LSA}$.
\end{theorem}
We immediately have the following corollary.
\begin{corollary}\label{cor:pfa_lsa}
Assume one of the following theories.
\begin{enumerate}
\item $\sf{PFA}$. 
\item There is a strongly compact cardinal.
\end{enumerate}
Then there is a model $M$ containing $\rm{OR}\cup\mathbb{R}$ such that $M\vDash \sf{LSA}$.
\end{corollary}
\begin{proof}
It is well-known that (1) implies the hypothesis of Theorem \ref{thm:square_lsa} (cf. \cite{VW}); this is because $\sf{PFA}$ implies $\neg \square(3,\gamma)$ for all $\gamma\geq \omega_2$. For (2), let $\kappa$ be a cardinal above a strongly compact cardinal $\gamma$ such that $\kappa^{\omega}=\kappa$. The hypothesis for Theorem \ref{thm:square_lsa} holds at $\kappa$ by the construction in \cite{solovay1974strongly}.
\end{proof}


In the previous chapter, we show  $\sf{LSA}$ is consistent assuming the existence of a Woodin cardinal which is a limit of Woodin cardinals by analyzing the HOD of the minimal model of $\sf{LSA}$, here we use the core model induction method to construct some model of determinacy that satisfies $\sf{LSA}$. The proof of Theorem \ref{thm:square_lsa} is built on that of \cite{Trang2015PFA}, which in turns is inspired by \cite{PFA} and \cite{sargsyan2014nontame}.

The rest of the chapter is dedicated to proving Theorem \ref{thm:square_lsa}. We assume the hypothesis of Theorem \ref{thm:square_lsa} along with the following simplifying assumption:
\begin{equation}\label{eqn:simplifying}
\kappa \textrm{ is measurable and $\forall \xi\in [\kappa,\kappa^{++}] \ 2^\xi = \xi^+$.}
\end{equation}
From the theorem's assumption, we have that 
\begin{center}
$\forall \alpha\in [\gamma,\kappa^{+4}]$, $\neg \square(3,\alpha)$. 
\end{center}
Later, we show how to get rid of assumption \ref{eqn:simplifying}. Our smallness assumption throughout this section is:

\begin{adjustwidth}{1cm}{1cm}
\begin{enumerate}
\item[$(\dag) \ \ \ $] in $V[G]$, where $G\subseteq Col(\omega, <\kappa)$ is $V$-generic, there is no model $M$ containing $\mathbb{R}\cup\rm{OR}$ such that $M \vDash ``\sf{ZF} + \sf{AD}^+$$ + \Theta=\theta_{\alpha+1} + \theta_{\alpha}$ is the largest Suslin cardinal below $\theta_{\alpha+1}$ ".
\end{enumerate}
\end{adjustwidth}

Before plunging in the the details, we give a very rough outline of the proof of Theorem \ref{thm:square_lsa}. Fix $\kappa$ as in the hypothesis of Theorem \ref{thm:square_lsa}. We operate under assumptions $(\dag)$ and \ref{eqn:simplifying}. Let $\mathbb{P} = Col(\omega, < \kappa)$. In $V^{\mathbb{P}}$, let $\Omega$ be the ``maximal pointclass of determinacy" (to be defined in the next section). Let $\P^-$ be the direct limit of hod pairs $(\M,\Sigma)$ such that $\Sigma\rest \rm{HC}\in \Omega$ and $\Sigma$ is $\Omega$-fullness preserving and has branch condensation. Let $\P$ be the appropriate ``Lp"-closure of $\P^-$ (defined in Section \ref{sec:cmi_background}). So $\P^-$ is an initial segment of $\P$. \cite{Trang2015PFA} and the results of Chapter \ref{chap:condensing_sets} show that $\P \vDash ``o(\P^-)$ is a regular limit of Woodin cardinals." In $V^{\mathbb{P}}$, we carry out a mixture of validated sts constructions and hybrid $K^c$-constructions over some transitive set $W$ containing $\P$ (to be explained in Section \ref{sec:stacking_mice}). Either the constructions stop prematurely (before stage $\kappa^{+++}$ for various reasons to be specified in Section \ref{sec:stacking_mice}), in which case we show that we can obtain a model of $\sf{LSA}$; otherwise, we reach a model $\P^+$ (extending $\P$) of height $\kappa^{+++}$. Then we consider the stack $\S$ of (appropriately defined) hod mice over $\P^+$. Using the proof of \cite[Theorem 3.4]{JSSS}, we show that cof$(o(\S))\geq \kappa^{+++}$. Using the fact that $\S\in V$ and our hypothesis $\forall \alpha\in [\gamma,\kappa^{+4}]$, $\neg \square(3,\alpha)$, we show that cof$(o(\S)) < \kappa^{+++}$. This contradiction shows that the second case cannot occur; therefore, we must reach a model of $\sf{LSA}$.

\section{Some core model induction backgrounds}\label{sec:cmi_background}
We continue to assume $(\dag)$ and \ref{eqn:simplifying} in this section. We recall some notions and results from \cite{Trang2015PFA}. In $V[G]$, where $\mathbb{P} = Col(\omega, < \kappa)$ and $G\subseteq \mathbb{P}$ is $V$-generic, let

\begin{center}$\Omega = \bigcup \{\powerset(\mathbb{R})\cap M \ | \  \mathbb{R}\cup \textrm{OR} \subset M \wedge M \vDash \textsf{AD}^+\}$.\end{center}
\cite{Trang2015PFA} shows that, under $(\dag)$,\footnote{\cite{Trang2015PFA} uses a stronger assumption, namely no models of ``$\sf{AD}_\mathbb{R}$$ +\Theta$ is regular" exist. But the proof there combined with the results in Section \ref{chap:condensing_sets}, particularly Theorem \ref{existence of condensing sets in n}, work using $(\dagger)$; the main point is that the HOD analysis now can be carried out up to models of $\sf{LSA}$.} the Solovay sequence $\langle \theta_\alpha^\Omega\ | \ \alpha \leq \gamma^*\rangle$ of $\Omega$ is of limit length. Furthermore, we get that if $A\in \Omega$, then there is a hod pair (or sts hod pair) $(\P,\Sigma)\in \Omega$ such that $A\in \Gamma^b(\P,\Sigma)$.

Let $\P^-$ be the direct limit of all hod pairs $(\M,\Sigma)$ such that $\M$ is countable in $V^{\mathbb{P}}$ and $\Sigma$ is an $(\omega_1,\omega_1+1)$-strategy of $\M$ that is strongly $\Omega$-fullness preserving, has strong branch condensation, and $\Sigma\rest \rm{HC}\in\Omega$. We will say that a pair $(\M,\Sigma)$ with these properties is \textit{nice} and let $\mathcal{F}$ be the direct limit system of all nice hod pairs. \cite{Trang2015PFA} shows that if $(\M,\Sigma\rest V)\in V$, then $\Sigma$ can be uniquely extended to a $(\kappa^{+4},\kappa^{+4})$-strategy $\Sigma^+$ (and hence $\Sigma^+\rest V\in V)$ (see Remark \ref{rem:extension} and the discussion before it for a somewhat more general argument). Say $\M$ iterates (via $\Sigma^+$) to a complete layer $\P^-(\alpha)$ of $\P^-$, where $\alpha<\gamma^*$, we let $\Sigma_\alpha$ be the $\Sigma^+$-tail of $\Sigma^+$.\footnote{The complete layers of $\P^-$ are $(\P^-(\beta) : \beta < \gamma^*)$ and for each $beta$, $\theta^\Omega_\beta$ is either the largest Woodin cardinal or the limit of Woodin cardinals in $\P^-(\alpha)$.} $\Sigma_\alpha$ only depends on $\alpha$ and does not depend on any particular choice of $(\M,\Sigma^+)$ as long as $\Sigma^+$ is nice. Let 
\begin{center}
$\Sigma = \oplus_{\alpha<\gamma^*} \Sigma_\alpha$,
\end{center}
and
\begin{center}
$\P = \rm{Lp}$$^{\Omega,\Sigma}(\P^-),$
\end{center}
be defined as in Section 9.1. So for every countable (in $V[G]$), transitive $\M^*$ embeddable into a level $\M\lhd \P$ via $\pi$, $\M^*$ is $(\omega_1,\omega_1+1)$-iterable as an (anomalous) hod mouse with strategy $\Lambda$ such that $\Lambda$ has a unique iteration strategy $\Lambda$ such that $\Lambda\rest \rm{HC} \in \Omega$; furthermore, if $\pi\in V$, then $\Lambda\rest V\in V$ and $\Lambda$ can be uniquely extended to a $(\kappa^{+4},\kappa^{+4})$-strategy in $V[G]$.

\begin{remark}\label{rem:phi_A}
As in Chapter \ref{chap:condensing_sets}, we let $\phi(U, W)$ be the formula that expresses the fact that $U$ is a mousefull pointclass with all the properties that $\Omega$ has and $W$ is a hod pair $(\Q, \Lambda)$ such that $Code(\Lambda)\in U$ and $\Lambda$ is strongly $U$-fullness preserving and has strong branch condensation. Then the $\F$ above is $\F_{\phi,\Omega}$ etc. From this point on, we will often suppress $\phi,\Omega$ from our notations; this should not be confusing since all the notations that come into the following constructions are relative to $(\phi,\Omega)$. 

The hypothesis of Theorem \ref{thm:square_lsa} implies that ``$(\phi, \Omega)$ is full, maximal, homogenous and lower part $(\phi, \Omega)$-covering fails". In particular, the conclusions of Theorem \ref{existence of condensing sets in n} hold for $V[G]$. $\myqedhere$
\end{remark}

\begin{lemma}\label{lem:small}
Let $\lambda$ be the ordinal height of $\Omega$, so $\lambda = \sf{ord}$$(\P^-) = \delta^\P$.
\begin{enumerate}
\item No levels $\M\lhd \P$ is such that $\rho_\omega(\M)< \lambda$. Hence $\rho_{\omega}(\P) = \sf{ord}$$(\P)$ and $\P\vDash \sf{ZFC}^-$.
\item $\P\vDash \delta^\P$ is a regular limit of Woodin cardinals
\item $\lambda <  \kappa^{+}$.
\item In $V$, $\sf{ord}$$(\P) < \kappa^+$ and cof$(\sf{ord}$$(\P)) < \gamma$.
\end{enumerate}
\end{lemma}
\begin{proof}
(1) follows from \cite[Theorem 3.78]{Trang2015PFA}. (2) follows by adapting the core model induction argument in \cite{Trang2015PFA} to get that $\gamma^*$ is a limit ordinal and using Theorem \ref{existence of condensing sets in n} to show the existence of condensing sets in $V[G]$, which in turns will give us that $\gamma^*$ is a regular cardinal in $\P$ (see for example the paragraph above \cite[Remark 3.86]{Trang2015PFA}). For (3) first note that $2^{<\kappa}=\kappa$ and $\omega_1 = \kappa$ in $V[G]$. For any $\Q\lhd^c_{hod} \P$, note that there is a hod pair $(\R,\Psi)\in \Omega$ such that $\R\in H_\kappa[G\rest \alpha]$, $\Psi\rest V[G\rest \alpha]\in V[G\rest\alpha]$ for some $\alpha<\kappa$. This easily implies (3).

For (4), first note that $\P \in V$. Let $\vec{C}$ be the canonical $\square_{\lambda}$-sequence built in $\P$ (using the construction in \cite{schimmerling2004characterization}), where $\lambda=\sf{ord}$$(\P^-)$ is the ordinal height of $\Omega$ as defined above. $\vec{C}$ is not threadable (by the maximality of $\P$). So if $\sf{ord}$$(\P) = \kappa^+$ or cof$(\sf{ord}$$(\P))\geq \gamma$, then using our hypothesis $\forall \alpha\in[\gamma,\kappa^{+4}], \neg\square(3,\alpha)$, we can find a thread for $\vec{C}$ by standard arguments. Contradiction.
\end{proof}

\section{Condensing sets}\label{sec:condensing_sets}

In $V[G]$, as done in Chapter \ref{chap:condensing_sets}, for each $X\in \powerset_{\omega_1}(\P)$, we let $\Q_X$ be the transitive collapse of $Hull_1^\P(X)$, $\delta_X = \delta^{\Q_X}$, and $\tau_X:\Q_X\rightarrow \P$ be the uncollapse map. Let $\Sigma_X$ be the $\tau_X$-pullback strategy for $\Q_X$.\footnote{Typically, $X = X^*\cap \P$ for some countable $X^*\prec H_{\kappa^{+4}}$. And $\Sigma_X$ is the $\tau_{X^*}$-realization map, where $\tau_{X^*}$ is the uncollapse map of $X^*$.} For $X\subseteq Y\in \powerset_{\omega_1}(\P)$, let $\tau_{X,Y}= \tau^{-1}_Y\circ \tau_X$. 

If $Y \in \powerset_{\omega_1}(\P^-)$ and $X\cap \P^- \subseteq Y$, we let $\Q^X_Y$ be the transitive collapse of $Hull_1^\P(X\cup Y)$, $\Sigma^X_Y$ be the $\tau_{X\oplus Y}$-pullback of $\Sigma$, and $\sigma^X_{Y}:\Q^X_Y\rightarrow \P$ be given by
\begin{center}
$\sigma^X_{Y}(q) = \tau_X(f)(\pi^{\Sigma^X_Y}_{\Q^X_Y,\infty}(a))$
\end{center}
where $a\in (\Q_Y|\delta^{\Q_Y})^{<\omega}$ and $q = \tau_{X,X\oplus Y}(f)(a)$. We also write $\tau^X_Y$ for $\tau_{X\oplus Y}$ and $\pi^X_Y$ for $\tau_{X,X\oplus Y}$.

Recall the notions of ($(\phi,\Omega)$)-condensing sets, extensions, and honest extensions discussed in Chapter \ref{chap:condensing_sets}. Suppose $Y \subset Z$ are extensions of $X$ and $X\in Cnd(\P)$. We write $\pi^X_{Y,Z}$ for the natural, uncollapse map from $\Q^X_Y$ to $\Q^X_Z$. By Lemma \ref{lem:condensing_set}, $\pi^X_{Y,Z}\rest \delta^{\Q^X_Y}$ agrees with the iteration map $\pi^{\Sigma_Y}_{\Q^X_Y,\Q^X_Z}$. We note that $\Q^X_X = \Q_X$ and if $Y$ is an extension of $X$, then $\pi^X_{X,Y}$ is just $\tau_{X,Y}$. When $X$ is a fixed condensing set and $Y$ extends $X$, we sometimes write $\Q_Y$ for $\Q^X_Y$ when no confusion arises.

Let $\mathfrak{S}\in V$ be the set of $X\prec H_{\kappa^{+4}}$ such that 
\begin{itemize}
\item $\kappa\cap X \in \kappa$, 
\item $\gamma < |X|< \kappa$, 
\item $\P\in X$, $X\cap \P$ is cofinal in $o(\P)$, and
\item $X^\xi \subset X$ for any $\xi < |X|$. 
\end{itemize}
Note that $\mathfrak{S}$ is stationary. We say that $\mathfrak{S}$ is the collection of \textit{good} hulls. For $X\in \mathfrak{S}$, we let $\pi_X: M_X\rightarrow H_{\kappa^{+4}}$ be the uncollapse map. Note that letting $\kappa_X$ be the critical point of $\pi_X$, $\pi_X$ extends to an elementary map from $M_X[G\rest \kappa_X]\rightarrow H_{\kappa^{+4}}[G]$, where $G\rest \kappa_X = G\cap Coll(\omega,<\kappa_X)$. We also call this map $\pi_X$.

The following facts follow easily from \cite{Trang2015PFA} and Chapter \ref{chap:condensing_sets}. The point is $j[\P]$ is a (strongly) $(\phi,\Omega)$-condensing set in $M[H]$ where $H\subseteq Coll(\omega,<j(\kappa))$ is $V$-generic and $G = H\rest \kappa$. Furthermore, because $\P\in V$ and $\sf{ord}$$(\P)<\kappa^+$ by Lemma \ref{lem:small}, $j[\P]\in M$.

\begin{lemma}\label{lem:condensing_set}
\begin{enumerate}[(i)]
\item $(\phi, \Omega)$ is full, maximal, homogenous and lower part $(\phi, \Omega)$-covering fails.
\item $\forall^* X'\in \mathfrak{S}_{\phi,\Omega}$, $X = X'\cap \P$ is a (strongly) condensing set.
\item Suppose $Y$ is an honest extension of a (strongly) condensing set $X$ and there are elementary maps $i:\Q_Y\rightarrow \R$ and $\sigma:\R\rightarrow \P$ such that $\sigma\circ i = \tau_Y\rest \Q_Y$ and every $x\in \R$ has the form $i(f)(a)$ for $f\in \Q_Y$ and $a\in [\delta^\R]^{<\omega}$. Then letting $\Lambda$ be the $\sigma$-pullback strategy of $\R$, and $\tau(i(f)(a)) = \tau_Y(f)(\pi^{\Lambda}_{\R|\delta^\R,\infty}(a))$, then $\tau$ is well-defined, elementary, and $\tau\rest \R|\delta^\R = \pi^{\Lambda}_{\R|\delta^\R,\infty}\rest \R|\delta^\R$.
\item Suppose $X$ is (strongly) condensing and $Y,Z$ are honest extensions of $X$ such that $\Q^X_Y = \Q^X_Z$, then $\Sigma^X_Y=\Sigma^X_Z$.
\end{enumerate}
\end{lemma}

\begin{remark}\label{rem:honest_extension}
Let $X$ be as in (ii) of the lemma. Then it is easy to se that any $Y = Y^*\cap \P$ where $Y^*\prec H_{\kappa^{+4}}$ is such that $Y^*$ is countable (in $V[G]$) is an honest extension of $X$. 

From now on, by ``condensing set", we mean ``strongly condensing set" and we will omit $(\phi,\Omega)$ from our terminology.

We let $Cnd(\P)$ be the collection of condensing sets $X\in \powerset_{\omega_1}(\P)$ in $V[G]$. If $X'$ is a good hull, and $X = X'\cap \P$ is a condensing set, then we say that $X'$ is a $X$-good hull. Similarly, any good hull $Y$ such that $X' \subset Y$ and $\{\P, X\}\in Y$ is called an \textit{$X$-good} hull. $\myqedhere$
\end{remark}

We also get the following easy consequences.

\begin{proposition}\label{prop:unique}
Suppose $X$ is a condensing set and $\Q$ is such that for some extension $Y$ of $X$, $\Q = \Q^X_Y$. Then there is a unique honest extension $W$ of $X$ such that $\Q = \Q^X_W$.
\end{proposition}

\begin{proposition}\label{prop:pullback}
Suppose $X$ is a condensing set. Suppose $Y$ and $W$ are extensions of $X$ such that there is a $\Sigma_1$-embedding $i: \Q^X_Y\rightarrow \Q^X_W$ and $\tau^X_W\circ i[\Q^X_Y]$ is an extension of $X$. Then $(\Sigma^X_W)^i = \Sigma^X_Y$.
\end{proposition}
\begin{proof}
Letting $Y^* = \tau^X_W\circ i[\Q^X_Y]$, then $\Q^X_Y = \Q^X_{Y^*}$. Moreover, $\Sigma^X_Y = \Sigma^X_{Y^*}$ by Lemma \ref{lem:condensing_set}, and hence $\Sigma^X_Y = (\Sigma^X_W)^i$.
\end{proof}

An easy modification of the proof of Theorem \ref{existence of condensing sets in n} also gives us the following useful fact. A proof of this can also be found in \cite[Lemma 2]{remarkSteel2015}. Note also that Corollary \ref{cor:strategy_condensation} is a consequence of Lemma \ref{lem:derived_model}. A standard argument, using the Vopenka algebra and the fact that $\P \models ``\delta^\P$ is a regular limit of Woodin cardinals", gives us that $L(\Omega,\P) \models ``\sf{AD}_\mathbb{R} + $$\Theta$ is regular." See \cite{remarkSteel2015} for a proof.

\begin{lemma}\label{lem:derived_model}
Suppose $X$ is a condensing set. Suppose $\Q, \R$ are countable, transitive in $V[G]$ with the property that there are elementary maps $i: \P_X\rightarrow \Q$, $k: \Q\rightarrow \R$, $\tau: \Q\rightarrow \P$, $\sigma: \R\rightarrow \P$ such that $\tau\circ i = \pi_X\rest \P_X$, and $\tau = \sigma \circ k$. Letting $\Sigma_\Q = \Sigma^\tau$ and $\Sigma_\R = \Sigma^\sigma$, then for any $A\subseteq \delta_X$ in $\P_X$ and any formula $\varphi$,
\begin{center}
$L(\Omega,\P) \models ``\forall s\in \Q ( \varphi[\Q,s,\Sigma_\Q,\P,\pi_X(A)] \Leftrightarrow  \varphi[\R,k(s),\Sigma_\R,\P,\pi_X(A)])"$.
\end{center}
\end{lemma}

\section{$X$-suitable hod premouse and $X$-validated iterations}\label{sec:X_suitable}

We summarize some key notions developed in \cite{sargsyan2019exact}. The reader can read \cite{sargsyan2019exact} for more details. We continue using the notations from the previous section. 

\begin{definition}\label{def_niceext}
Let $X$ be a condensing set. $\Q$ \textbf{nicely extends} $\Q_X$ if $\Q$ is non-meek and $\Q^b = \Q_X$. We also say that $\Q$ is a \textbf{nice extension} of $\Q_X$. Similarly, we can define what it means for $\Q$ to nicely extend $\Q^X_Y$ for an extension $Y$ of $X$ and $\Q$ to nicely extend $\P$. $\myqedhere$
\end{definition} 
 

Suppose $X$ is a condensing set and $Y$ is an extension of $X$. Suppose further that $\Q$ nicely extends $\Q^X_Y$. A stack (of normal iteration trees) $\T$ on $\Q$ has the form
\begin{center}
$\VT=((\M_\a)_{\a<\eta}, (E_\a)_{\a<\eta-1}, D, R, (\beta_\a, m_\a)_{\a\in R}, T)$,
\end{center}
where the displayed objects are introduced in \ref{putative it}. The above notation is quite standard. $D$ is the set of drops, $R$ is the set of stages where player $I$ starts a new round of the iteration game, $(\beta_\a, m_\a)$ is the place player $I$ drops at the beginning of the $\a$th round, and $T$ is the tree order. Sometimes to make clear these objects are associated with $\vec{\T}$, we write $D^{\vec{\T}}$ etc. Recall the convention in Chapter 2, we assume that all our stacks are proper (see Remark \ref{proper stacks convention}). One of the key aspects of being proper is that if $\b<lh(\VT)$ is such that $\VT_{\geq \b}$ is a stack on $\M^\VT_\b$ then $\b\in R$\footnote{Thus, no normal component of $\VT$ can be split into two normal components.}. We will also use the notation introduced in \ref{notation for iteration trees}. In particular, for $\a\in R^\VT$, ${\sf{next}}^{\VT}(\a)=min(R^{\VT}-(\a+1))$ if this minimum exists and otherwise ${\sf{next}}^{\VT}(\a)=lh(\VT)$. For $\a\in R^\VT$, we also set ${\sf{nc}}^{\VT}_\a=\VT_{[\a, \a']}$ where $\a'={\sf{next}}^{\VT}(\a)$.


\begin{definition} \label{dfn:realizable_ext} Suppose $X\in Cnd(\P)$ and $Y$ is an extension of $X$. Suppose further that $\Q$ nicely extends $\Q^X_Y$. Given $E\in \vec{E}^\Q$ such that $\cp(E)=\d^{\Q^X_Y}$, we say $E$ is \textbf{$(X, Y)$-realizable} if there is $W$, an extension of  $X\oplus Y$ such that $E=E^X_{Y, W}$, where $E^X_{Y,W}$ is the extender defined by\:
\begin{equation}\label{dfn:extender}
(a,A)\in E^X_{Y,W} \Leftrightarrow \tau^X_{W}(a) = \pi^{\Sigma^X_W}_{\Q^X_W,\infty}(a)\in \tau^X_Y(A),
\end{equation}
for any $a\in [lh(E)]^{<\omega}$ and $A\in \powerset(\cp(E))^{|a|}\cap \Q$. 

We are continuing with the notation of \rdef{dfn:realizable_ext}. Suppose $\VT$ is a stack on $\Q$. We say $\VT$ is a \textbf{$(X, Y)$-realizable} iteration if there is a sequence $(W_\a: \a \in R^{\VT})$ such that
\begin{enumerate}
\item $W_0=Y$,
\item if $\a, \b\in R^{\VT}$ and $\a<\b$ then $W_\b$ is an extension of $X\oplus W_\a$, 
\item if $\a, \b\in R^{\VT}$, $\a<\b$ and $\pi^{\VT, b}_{\a, \b}$ is defined then $\pi^{\VT, b}_{\a, \b}=\pi^X_{W_\a, W_\b}$\footnote{The embedding $\pi^{\VT}_{\a, \b}$ is defined similarly to $\pi^{\VT, b}$, it is essentially the embedding $\pi^{\VT}_{\a, \b}\rest \M_\a^b$. See Section \ref{sec:pitb}.}, and
\item if $\a\in R^\VT$ and $\VU$ is the largest fragment of $\VT_{\geq \a}$ that is based on $\M_\a^b$ then $\VU$ is according to $\Psi^X_{W_\a}$.
\end{enumerate}

We say $\VT$ is \textbf{$X$-realizable} if $Y$ is an honest extension of $X$ and $\VT$ is  $(X, Y)$-realizable.  $\myqedhere$


%

\end{definition}

The following lemma shows that $(X,Y)$-realizability is equivalent to $X$-realizability. The proof can be found in \cite[Section 7]{sargsyan2019exact}.

\begin{lemma}\label{lem:XY-realizability}
Suppose $Y$ is an extension of $X$ and $\Q$ nicely extends $\Q^X_Y$. Suppose $\VT$ is a $(X, Y)$-realizable iteration as witnessed by $(W'_\a: \a \in R^{\VT})$. For $\a\in R^{\VT}$ let $W_\a$ be the unique honest extension of $X$ with the property that $(\M_\a^{\VT})^b=\Q^X_{W_\a}$. Then $(W_\a: \a\in R^{\VT})$ witnesses that $\Q$ is $X$-realizable. 
\end{lemma}

We fix a condensing set $X$ throughout this section. Suppose $\Q$ nicely extends $\Q^X_Y$ and $\VT$ is a $X$-realizable iteration of $\Q$. We cannot in general prove that $\VT$ picks unique branches mainly because we say nothing about $\Q$-structures that appear in $\VT$ when we iterate above $\d^{\S^b}$ for some $\S = \M^\T_\beta$ and $\beta\in R^{\vec{\T}}$. The next definition introduces a notion of a premouse that resolves this issue. 

\begin{definition}\label{weakly suitable} We say $\R$ is \textbf{weakly $X$-suitable} if $\R$ is a hod premouse of lsa type such that $\R=(\R|\d^\R)^{\#}$, $\R$ has no Woodin cardinals in the interval $(\d^{\R^b}, \d^\R)$ and for some extension $Y$ of $X$, $\R$ nicely extends $\R^b=\Q^X_Y$. We say $\R$ is \textbf{weakly suitable} if $\R$  is a hod premouse of lsa type such that $\R=(\R|\d^\R)^{\#}$, $\R$ has no Woodin cardinals in the interval $(\d^{\R^b}, \d^\R)$, and $\R^b = \P$.$\myqedhere$
\end{definition}

We now define the notion of \textit{$X$-approved sts premouse of depth $n$} by induction on $n$. The induction ranges over all weakly $X$-suitable hod premice. Suppose $\R$ is weakly $X$-suitable hod premouse and $Y$ is an extension of $X$ such that $\R$ nicely extends $\R^b = \Q^X_Y$. \\\\
(1)  We say that $\M$ is a \textbf{$X$-approved sts premouse over $\R$ of depth $0$} if $\M$ is an sts premouse over $\R$\footnote{This in particular means that the strategy indexed on the sequence of $\M$ is a strategy for $\R$.} such that if $\T\in \M$ is according to $S^\M$ then $\T$ is $(X,Y)$-realizable.\\
(2) We say that $\M$ is a \textbf{$X$-approved sts premouse over $\R$  of depth $n+1$} if $\M$ is a $X$-approved sts premouse over $\R$ of depth $n$ such that if $\T\in \M$ is a $\sf{nuvs}$ and $S^\M(\T)$ is defined then letting $b=S^\M(\T)$, $\Q(b, \T)$ is a $X$-approved sts premouse over $\textrm{m}^+(\T) = (\M(\T))^\#$ of depth $n$.

\begin{definition}\label{z-validated sts mice} We say $\M$ is a \textbf{$X$-approved sts premouse} over $\R$ if for each $n<\omega$, $\M$ is an $X$-approved sts premouse over $\R$ of depth $n$. We say $\M$ as above is a \textbf{$X$-approved sts mouse} (over $\R$) if $\M$ has a $\mu$-strategy $\Sigma$ such that whenever $\N$ is a $\Sigma$-iterate of $\M$, $\N$ is a $X$-approved sts premouse over $\R$. 

We say $\M$ is an \textbf{$X$-approved hod premouse} if whenever $\T\in \M$ is according to $S^\M$, then $\T$ is $X$-realizable. We say $(\M,\Sigma)$ is an \textbf{$X$-approved hod mouse} if whenever $\mathcal{U}$ is according to $\Sigma$ with last model $\N$, then $\N$ is an $X$-approved hod premouse and $\Sigma_{\mathcal{U},\N}\rest \N = S^\N$.  \index{$X$-approved hod premouse}

$\myqedhere$
\end{definition}

We let $Lp^{Xa, sts}(\R)$ be the union of all $X$-approved sound sts mice over $\R$ that project to $ \leq \sf{ord}$$(\R)$. Finally, we can define the \textit{correctly guided $X$-realizable iterations}.

\begin{definition}\label{correctly guided iterations} Suppose $\R$ is a weakly $X$-suitable hod premouse and $\VT$ is a $X$-realizable iteration of $\R$. We say $\VT$ is \textbf{correctly guided} if whenever $\alpha \in R^{\vec{\T}}$, $\U = \textsf{nc}$$_\alpha^{\vec{\T}}$ is above $\delta^{\M^b_\alpha}$, and $\a<lh(\U)$ is a limit ordinal such that $\textrm{m}^+(\M(\U\rest \a))\models ``\d(\U\rest\alpha)$ is a Woodin cardinal", then letting $b=[0, \a]_\U$, $\Q(b, \U\rest \alpha)$ is an $X$-approved sts mouse over $\textrm{m}^+(\M(\U\rest \a))$. $\myqedhere$
\end{definition} 

The following facts follow straightforwardly from the definitions above (see \cite[Section 7]{sargsyan2019exact} for proofs).
\begin{proposition}\label{preservation of zv under embeddings}
\begin{enumerate}[(i)]
\item \label{hull of realizable iterations} Suppose $\R$ and $\S$ are weakly $X$-suitable hod premice. Suppose further that $\VT$ is a $X$-realizable iteration of $\S$ and $\VU$ is an iteration of $\R$ such that $(\R, \VU)$ is a hull of $(\S, \VT)$. Then $\VU$ is also $X$-realizable. 

\item \label{weakly suitable condenses} Suppose $\R$ and $\S$ are weakly $X$-suitable, $\N$ is an sts premouse over $\R$ and $\M$ is a $X$-approved premouse (mouse) over $\S$. Suppose $\pi:\N\rightarrow_{\Sigma_1} \M$. Then $\N$ is also a $X$-approved premouse (mouse). 

\item \label{hulls of corectly guided iterations} Suppose $\R$ and $\S$ are weakly $X$-suitable hod premice. Suppose further that $\VT$ is a correctly guided $X$-realizable iteration of $\S$ and $\VU$ is an iteration of $\R$ such that $(\R, \VU)$ is a hull of $(\S, \VT)$. Then $\VU$ is also correctly guided $X$-realizable iteration. 

\end{enumerate}
\end{proposition}
Our uniqueness theorem applies to $\R$ that are not \textit{infinitely descending}.

\begin{definition}\label{inf desc} We say that a weakly $X$-suitable hod premouse $\R$ is \textbf{infinitely descending} if there is a sequence $(p_i, \R_i, Y_i: i<\omega)$ such that 
\begin{enumerate}
\item $\R_0=\R$,
\item for every $i<\omega$, $\R_i$ is weakly $X$-suitable and nicely extends $\R_i^b = \Q^X_{Y_i}$,
\item for every $i<\omega$, $p_i$ is a correctly guided $X$-realizable iteration of $\R_i$,
\item for every $i<\omega$, $p_i$ has a last normal component $\T_i$ of successor length such that $\a_i=_{def}lh(\T_i)-1$ is a limit ordinal and 
$\R_{i+1} = \textrm{m}^+(\M(\T_i\rest \a_i))$,
\item for every $i<\omega$, setting $b_i=_{def}[0, \a_i)_{\T_i}$, $b_i$ is a cofinal branch of $\T_i$ such that $\Q(b_i, \T_i)$ exists and is $X$-approved.
\end{enumerate}
$\myqedhere$
\end{definition} 

Note that in the above definition, for each $i$, $\R_{i+1}$ is a strict initial segment of $\Q(b_i, \T_i)$. The following is the uniqueness result we need. 
\begin{proposition}\label{branch uniqueness} Suppose $\R$ is a weakly $X$-suitable hod premouse that is not infinitely descending and $\VT$ is a correctly guided $X$-realizable iteration of limit length on $\R$. There is then a unique branch $b$ of $\VT$ such that $\VT^\frown \{b\}$ is correctly guided and $X$-realizable.  \index{weakly $X$-suitable hod premouse}
\end{proposition}
\begin{proof} The proof easily follows from results of Section 4, particularly Lemma \ref{disagreement implies low level disagreement}, and Section 9.2. Suppose first that
\begin{enumerate}[(a)]
\item if $\VT$ doesn't have a last component or
\item if there is $\a\in R^\VT$ such that $\VT_{\geq\a}$ is based on $\M_\a^b$.
\end{enumerate}
In case (a), there is nothing to prove. In case (b), let $\S = \M^{\vec{\T}}_\alpha$, and let $W_\S$ be as in \rdef{dfn:realizable_ext}, $\Psi^X_{W_S}$ only depends on $\S^b$ (by Lemma \ref{lem:condensing_set})\footnote{Notice that in this case there is a branch $b$ such that $\VT^\frown \{b\}$ is correctly guided and $X$-realizable.}. Suppose we are not in case (a) or (b). Let now $\T={\sf{nc}}^\VT_\a$ be the last normal component of $\VT$. If $b, c$ are two different branches of $\T$ such that $\VT^\frown\{b\}$ and $\VT^\frown \{c\}$ are correctly guided $X$-realizable iterations then $\Q(b, \T)\not =\Q(c, \T)$ and both are $X$-approved sts mice over $\m^+(\T)$. It now follows from Lemma  \ref{disagreement implies low level disagreement} and the fact that $\R$ is not infinitely descending that we can reduce the disagreement of $\Q(b, \T)$ and $\Q(c, \T)$ to a disagreement between $\Psi^X_Y$ and $\Psi^X_Z$ for some extensions $Y,Z$ of $X$ with $\Q^X_Y=\Q^X_Z$. However, this cannot happen by Lemma \ref{lem:condensing_set}. 
\end{proof}


\begin{definition}\label{z-validated sts mice above mu} Suppose $X$ is a condensing set. Suppose $\R_0$ extends $\P$, $p$ is an iteration of $\R_0$ such that if $p$ is $\sf{nuvs}$, then setting $\R = \textrm{m}^+(p)$, $\M$ is an sts premouse over $\R$. Suppose $(\R, \M, p)\in H_{\kappa^{+4}}$. 
\begin{enumerate}
\item We say $\R$ is not infinitely descending if whenever  $U$ is an $X$-good hull such that $\R\in U$, $\pi_U^{-1}(\R)$ is not infinitely descending.
\item We say \textbf{$p$ is $X$-validated} if whenever $U$ is an $X$-good hull such that $\{\R, p\}\subseteq U$, $\pi_U^{-1}(p)$ is a correctly guided $X$-realizable iteration of $\pi_U^{-1}(\R)$.
\item Suppose $\R$ is weakly suitable. We say \textbf{$\M$ is a $X$-validated sts premouse} over $\R$ if  for every $X$-good hull $U$  such that $\{\R, \M\}\subseteq U$, letting $\N=\pi_U^{-1}(\M)$, $\N$ is a $X$-approved sts premouse over $\pi^{-1}_U(\R)$. 
\item Suppose $\R$ is weakly suitable. We say \textbf{$\M$  is a $X$-validated sts mouse} over $\R$ if whenever $U$ is an $X$-good hull such that $\{\R, \M\}\subseteq U$, letting $\N=\pi_U^{-1}(\M)$, $\N$ is a $X$-approved sts mouse over $\pi^{-1}_U(\R)$.
\item Suppose $\M$ is a $X$-validated sts mouse over $\R$ and $\xi$ is an ordinal. We say $\M$ has an $X$-validated $\xi$-strategy $\Sigma$ if whenever $\N$ is an iterate of $\M$ via $\Sigma$, $\N$ is a $X$-validated sts mouse over $\R$. \index{$X$-validated sts premouse}
\item Suppose $q$ is an iteration of $\R$. We say $q$ is \textbf{$X$-validated} if $p^\smallfrown q$ is $X$-validated.
\item We say a hod premouse $\M$ such that $\mathcal{P}\unlhd \M$ is an \textbf{$X$-validated} hod premouse (mouse) if  for every $X$-good hull $U$  such that $\{\mathcal{P}, \M\}\subseteq U$, letting $\N=\pi_U^{-1}(\M)$, $\N$ is a $X$-approved hod premouse (mouse, respectively).
\end{enumerate}
$\myqedhere$
\end{definition}

The following proposition is very useful and is an immediate consequence of \rprop{preservation of zv under embeddings}. When $U$ is a good hull we will use it as a subscript to denote the $\pi_U$-preimages of objects that are in $U$.

\begin{proposition}\label{one hull witness for premice} Suppose $\R_0,p,\R,\M$ are as in Definition \ref{z-validated sts mice above mu}. Suppose $U$ is an $X$-good hull such that $\{\R, \M\}\subseteq U$ and  $\M_U$ is not $X$-approved. Then whenever $U^*$ is an $X$-good hull such that $U\cup\{U\}\subseteq U^*$, $\M_{U^*}$ is not $X$-approved. Hence, $\M$ is not $X$-validated.
\end{proposition}

A similar result holds for iterations.

\begin{proposition}\label{one hull witness for iterations} Suppose $\R_0$ is as in \rprop{one hull witness for premice}. Suppose $p\in H_{\kappa^{+4}}$ is an iteration of $\R_0$. Suppose $U$ is an $X$-good hull such that $\{\R_0, p\}\subseteq U$ and  $p_U$ is not $X$-realizable. Then whenever $U^*$ is an $X$-good hull such that $U\cup\{U\}\subseteq U^*$, $p_{U^*}$ is not $X$-realizable. Hence, $p$ is not $X$-validated.
\end{proposition}

\begin{definition}\label{z-full} We say $\R$ is \textbf{$X$-suitable} if it is weakly $X$-suitable and whenever $\M$ is an $X$-approved sts mouse over $\R$ then $\M\models ``\d^\R$ is a Woodin cardinal". We say $\R$ is \textbf{suitable} if it is weakly suitable and whenever $\M$ is an $X$-validated sts mouse over $\R$, then $\M\models ``\d^\R$ is a Woodin cardinal". $\myqedhere$ \index{$X$-suitable hod premouse}
\end{definition}

 We let $Lp^{Xv, sts}(\R)$ be the union of all $X$-validated sound sts mice over $\R$ that project to $\leq \sf{ord}$$(\R)$. The following proposition is a consequence of \rprop{branch uniqueness}.

\begin{proposition}\label{unique branches above mu} Suppose $(\R_0,p,\R)$ are as in Definition \ref{z-validated sts mice above mu} and $\R$ is not infinitely descending. Suppose $\VT$ is an $X$-validated iteration of $\R$ of limit length. Then there is at most one branch $b$ of $\VT$ such that $\VT^\frown\{b\}$ is $X$-validated.
\end{proposition}

In the next two sections, we describe two kinds of constructions: \textit{the hybrid $K^c$-construction} over $\P$ or some $\P'$ extending $\P$,\footnote{This is a variation of the mixed hod pair construction in Definition \ref{mixed hod pair construction}.} and the \textit{$X$-validated sts constructions} over some weakly suitable $\R$. We use the notations and definitions from the previous section. We fix a condensing set $X\in V$ ($X$ exists by the previous section); and we assume that $X = X'\cap \P$ where $X'\prec H_{\kappa^{+4}}$ is of size $\kappa$ in $V$. The hybrid $K^c$-construction proceeds more or less according to the usual procedure for building hod pairs (as described many times in this book) except that extenders $F$ we put on the sequence of the models are correctly backgrounded (described below) instead of being fully backgrounded. This construction will reach a weakly suitable stage $\R$. We then continue with the $X$-validated sts construction over $\R$. If this construction produces a $\Q$-structure $\Q$ extending $\R$, then we attempt to construct an $X$-validated strategy for $\Q$. If this is successful, we then continue with the hybrid $K^c$-construction over $\Q$. If not, we show that an honest suitable $\R'$ must exist (see below). Producing such an $\R'$ means that the $X$-validated sts construction over $\R'$ will no longer produce $\Q$-structures for $\R'$ and will either reach a model that produces $\N_{\omega.2,lsa}^{\sharp}$ and hence a model of $\sf{LSA}$ or go on for $\kappa^{+++}$ many steps. We will rule out the latter by an argument using the technique of stacking mice developed in \cite{JSSS} (see the next section). This argument also rules out the case that the two aforementioned constructions alternate for $\kappa^{+++}$ many times.

The two constructions described above will produce a sequence of models $(\M_\xi, \N_\xi: \xi\leq \Upsilon)$. Before defining the sequence, we discuss the kind of background extenders being used in this construction. Suppose $\M_\xi$ has been constructed as part of one of such constructions and is in $V$, is passive, $\powerset(\delta^\P)^{\M_\xi}=\powerset(\delta^\P)^\P$. Suppose $F$ is a $(\cp(F),\sf{ord}$$(\N_\xi))$-extender that coheres the sequence of $\N_\xi = \mathcal{C}(\M_\xi)$.\footnote{$\N_\xi$ is the appropriate fine-structural core of $\M_\xi$, as dictated by the construction. See the next section.} Suppose $Y\prec H_{\kappa^{+4}}$ (or sometimes, we'll let $Y\prec H_\gamma$ for $\gamma\geq \kappa^{+4}$) is in V and is a good $X$-hull and $Y$ contains all relevant objects. Let $\pi_Y$ be the corresponding uncollapse map.  Suppose $\N^Y_\xi =_{\textrm{def}} \pi_Y^{-1}(\N_\xi)$ has a unique $X$-realizable strategy $\Sigma_\xi^Y$ such that $\Sigma_\xi^Y\rest\rm{HC}\in \Omega$ (these properties will be maintained during the course of our construction). Alternatively, we sometimes write $\Sigma_{Y,\xi}$ for $\Sigma_\xi^Y$.

\begin{definition}\label{def:correctly_backgrounded_ext}\index{correctly backgrounded extender}
We say that an extender $F$ is \textbf{correctly backgrounded}\index{correctly backgrounded extender} if one of the following holds:
\begin{itemize}
\item if $\cp(F) = \delta^\P$ and the least cutpoint above $\delta^\P$ is the largest cardinal in $\M_\xi$, then $(a,A)\in F$ if and only if $\forall^*Y$, $Y$ is $X$-good, letting $a_Y = \pi_Y^{-1}(a)$, $\pi^{\Sigma^Y_\xi}_{\N^Y_\xi,\infty}(a_Y)\in A$. We say that $F$ is \textbf{$X$-certified.}

\item if $\cp(F) > \delta^\P$, then say, $\lambda=F(\cp$$(F))$, $F$ is \textbf{certified by a collapse} in the sense of \cite{JSSS}, that is, there is $Z\prec H^V_{\kappa^{+4}}$ (in $V$) such that $|Z| < \kappa^{+++}$, where $o(\P)+1\subset Z$, $Z\cap \kappa^{+++}\in \kappa^{+++}$,  $Z^{<\kappa}\subseteq Z$,\footnote{Note that $|o(\P)| = |\Omega| =  \kappa$. Futhermore, $Z[G]^{<\omega_1} \subseteq Z[G]$ in $V[G]$.} and letting $\pi_Z:M_Z\rightarrow Z$ be the uncollapse, we have: $\M_\xi|(\cp$$(F))^{+,\M_\xi}\in M_Z$, $\cp(\pi_Z) = \cp(F)$, and
\begin{center}
$F$ is the Jensen completion of $(\pi\rest \powerset(\cp$$(\pi_Z))\cap \N_\xi)\rest \lambda$.
\end{center}
\end{itemize}
In either case, we will sometimes say ``$F$ is certified".$\myqedhere$
\end{definition}

Suppose $F$ is $X$-certified and $Y$ is $X$-good. Let $F_Y = \pi_Y^{-1}(F)$. Then it is easy to check that $(a,A)\in F_Y$ iff $\pi^{\Sigma^Y_\xi}_{\N^Y_\xi,\infty}(a)\in \pi_Y(A)$. We say that $F_Y$ is \textit{$\pi_Y$-certified over $(\N^Y_\xi,\Sigma^Y_\xi)$}.

\section{The $X$-validated sts constructions}\label{sec:sts_constructions}\index{$X$-validated sts construction}

We first describe the $X$-validated sts construction for a fixed condensing set $X$. Suppose $\R$ is weakly suitable and $Y\in H_{\kappa^{+++}}$ is transitive such that either $\R= Y$ or $\R\in Y$. Recall the conventions for hod  premice introduced earlier in the book. A (hod) premouse has the form $(\M,k)$, where $\M$ is a $k$-sound, acceptable $J$-structure. $k(M) = k$ is the degree of soundness of $\M$. We write the core $\mathcal{C}(\M)$ for the ($k(\M)+1$-)core of $\M$ (if this makes sense, i.e. when $\M$ is $k(\M)+1$-solid, or just solid; the same abbreviation will be applied when we say $\M$ is universal, meaning $\M$ is $k(\M)+1$-universal). Similarly, we write $\rho(\M)$ for the $k(\M)+1$-projectum and $p(\M)$ for the $k(\M)+1$-standard parameters of $\M$. When $\mathcal{C}(\M)$ exists, $k(\mathcal{C}(\M)) = k(\M)+1$. $\M$ is sound iff $\M = \mathcal{C}(\M)$. For brevity, we suppress the degree of soundness of the models constructed below. For instance, if $k(\M_\xi) = k$, then we write $\M_\xi$ for $(\M_\xi, k)$. Before, describing the next construction, we advise the reader to consult \cite[Chapter 10]{sargsyan2019exact} for a similar construction.

\begin{definition}[$X$-validated sts construction]\label{dfn:Kc_constr} We say $(\M_\xi, \N_\xi: \xi\leq \Upsilon)$ are the models of the \textit{$X$-validated sts construction over $Y$} if the following conditions hold:
\begin{enumerate}
\item $\Upsilon \leq \kappa^{+++}$, and for all $\xi<\Upsilon$ if $\M_\xi, \N_\xi$ are defined then $\M_\xi, \N_\xi\in H_{\kappa^{+++}}$.
\item For every $\xi\leq \Upsilon$, $\M_\xi$ and $\N_\xi$ are $X$-validated sts hod premice over $\R$ or are $X$-validated sts hod premice over $Y$ based on $\R$.
\item Suppose the sequence $(\M_\xi, \N_\xi: \xi<\eta)$ has been constructed. Suppose further that
there is a total $(\k, \nu)$-extender $F$ such that $F$ is certified by a collapse and letting $G$ be the Jensen completion of $\N_{\eta-1}\cap F$, $(\N_{\eta-1}, G)$ is a $X$-validated sts hod premouse over $\R$ or over $Y$ based on $\R$. Let then $\M_\eta=(\N_{\eta-1}, F)$ and 
$\N_{\eta}=\mathcal{C}(\M_\eta)$.

\item Suppose the sequence $(\M_\xi, \N_\xi: \xi<\eta)$ has been constructed, and $\T\in \N_{\eta-1}$ is the $<_{\N_{\eta-1}}$-least $\sf{uvs}$ without an indexed branch. Suppose further that there is a branch $b$ of $\T$ such that $(\N_{\eta-1}, B_b)$\footnote{$B_b$ is a code for $b$ as done in \cite[Section 2]{trang2013} and outlined in Chapter \ref{chap:square}. We only note that this amenable coding ensures condensation under very weak hull embeddings, cf. \cite[Lemma 3.10]{trang2013} and this fact is in turns used to show that $\square_{\kappa,2}$ holds in hod mice. From now on, we may confuse the structure $(\N_{\eta-1}, B_b)$ with $(\N_{\eta-1}, b)$.}  is a $X$-validated sts hod premouse\footnote{This in particular implies that $b\in \N_{\eta-1}$. For brevity, we suppress the other predicates that are part of the hod premouse, like $\epsilon, \vec{E}$ etc.} over $\R$ or over $Y$ based on $\R$. Let then $\M_\eta=(\N_{\eta-1}, B_b)$ and $\N_{\eta}=\mathcal{C}(\M_\eta)$. 
\item Suppose the sequence $(\M_\xi, \N_\xi: \xi<\eta)$, and for some $\sf{nuvs}$ tree $\T\in \N_{\eta-1}$ there is a branch $b$ such that $(\N_{\eta-1}, B_b)$ is an $X$-validated sts hod premouse over $\R$ or over $Y$ based on $\R$. Let $\T$ be the $\N_{\eta-1}$-least such tree and $b$ be such a branch for $\T$. Then $\M_\eta=(\N_{\eta-1}, B_b)$ and $\N_{\eta}=\mathcal{C}(\M_\eta)$.
\item Suppose the sequence $(\M_\xi, \N_\xi: \xi<\eta)$ has been constructed and all of the above cases fail. In this case we let $\M_\eta=\mathcal{J}_1(\N_{\eta-1})$ and provided $\M_\eta$ is a $X$-validated sts hod premouse over $\R$ or over $Y$ based on $\R$, $\N_{\eta}=\mathcal{C}(\M_\eta)$.
\item Suppose the sequence $(\M_\xi, \N_\xi: \xi<\eta)$ has been constructed and $\eta$ is a limit ordinal. Then $\M_\eta=liminf_{\xi\rightarrow \eta}\M_\xi$.
\end{enumerate}$\myqedhere$
\end{definition}

The construction fails at $\eta$ if one of the following holds.

\noindent $(i)_\eta$  $\M_\eta$ is not solid or universal.\\
$(ii)_\eta$ $\M_\eta$ is not $X$-validated.\\
$(iii)_\eta$ There is a $\sf{uvs}$  $\T\in \N_{\eta-1}$ such that the indexing scheme demands that a branch of $\T$ must be indexed yet $\T$ has no branch $b$ such that $(\N_{\eta-1}, B_b)$ is a $X$-validated sts premouse over $\R$ or over $Y$ based on $\R$. \\
$(iv)_\eta$ There is a $\sf{nuvs}$ $\T\in \N_{\eta-1}$ such that the indexing scheme demands that a branch of $\T$ must be indexed yet $\T$ has no branch $b$ such that $(\N_{\eta-1}, B_b)$ is a $X$-validated sts premouse over $\R$ or over $Y$ based on $\R$.\\
$(v)_\eta$ $\rho(\M_\eta)\leq \d^{\mH}$.\\

We note that in (3) above, $\cp(F) > \delta^\P$. If the construction holds at $\eta$ (i.e. $(1)_\eta-(5)_\eta$ all fail) then $\M_\eta, \N_\eta$ are $\sf{hp}$-indexed $\sf{lses}$. Indeed, we inductively maintain that the models in our construction are $\sf{hp}$-indexed $\sf{lses}$. Verifying $(iii)_\eta, (iv)_\eta$ fail for all $\eta$ in the $X$-validated sts construction roughly corresponds to showing that the important anomaly doesn't occur in the construction in \ref{mixed hod pair construction}.\footnote{We note that $\T$ as in $(iii)_\eta,(iv)_\eta$ are of the form $(\R,\T_0,\S,\T_1)$ but we will suppress this notation for brevity.}

One shows $(i)_\eta$ fails by the usual arguments, namely showing that countable substructures of $\M_\eta$ are iterable. That $(v)_\eta$ fails will be shown when $Y$ is specified. The proof that $(v)_\eta$ fails for certain precisely defined $Y$ is to show that $\M_\eta$ is iterable via a fullness preserving strategy. The proof that $(v)_\eta$ fails will subsume the proof that $(i)_\eta$ fails and will be given later (see Lemma \ref{cor:solidity_universality}). \cite[Section 9]{sargsyan2019exact} shows that $(ii)_\eta, (iii)_\eta, (iv)_\eta$ cannot fail if $\R$ is honest as witnessed by $(\vec{\mathcal{V}},p)\in H_{\kappa^{+4}}$, where $\vec{\mathcal{V}} = (\mathcal{V}_\alpha : \alpha \leq \xi)$ is an array with the $X$-realizability property (defined in \cite[Definition 9.1, 9.2]{sargsyan2019exact}) and either $\R = \mathcal{V}_\xi$ and $p=\emptyset$ or $p$ is an $X$-validated iteration of $\mathcal{V}_\xi$ of limit length such that $\pi^{p,b}$ exists and $\R = \textrm{m}^+(p)$. We summarize some main points of the arguments in \cite{sargsyan2019exact} below and outline the proof that $(iii)_\nu$ fails and $(iv)_\eta$ fails. First we define the objects $(\vec{\mathcal{V}},p)$ more precisely.

\begin{definition}\label{array} We say $\vec{\V}=(\V_\a:\a\leq \eta)$ is an \textbf{array} of length $\eta$ if the following conditions hold. 
\begin{enumerate}
\item For every $\a<\eta$, $(\V_\b:\b\leq \a)$ is an array of length $\a$ at $\mu$.
\item $\V_\eta$ nicely extends $\P$ and is a $X$-validated hod premouse.
\item For all $\a<\eta$, if $\V_\a$ is weakly $X$-suitable then there is $\b\leq \eta$ such that $\V_\b$ is a $X$-validated sts mouse over $\V_\a$ and $\J_1(\V_\b)\models ``$there are no Woodin cardinals $>\d^{\P}$".
\item For all $\a<\eta$, if $\J_1(\V_\a)\models ``$there are no Woodin cardinals $>\d^{\P}$" then $\V_\a$ has a $X$-validated iteration strategy. 
\end{enumerate}
We say $\vec{\V}$ is \textbf{small} if $rud(\V_\eta)\models ``$there are no Woodin cardinals $>\d^{\V_\eta^b}$". We let $\eta=lh(\vec{\V})$ and for $\a\leq \eta$, we let $\vec{\V}\rest \a=(\V_\b:\b\leq \a)$. $\myqedhere$
\end{definition}

Recall the notions of $k$-maximal iteration trees in \cite[Definition 3.4]{OIMT}, weak $k$-embeddings \cite[Definition 4.1]{OIMT}. For an iteration tree $\T$ on $\M$, letting $\M^\T_\alpha$ be the $\alpha$-th model in the tree; for $\alpha+1<lh(\T)$, recall the notion of degree $deg^\T(\alpha+1)$ \cite[Definition 3.7]{OIMT}. Recall the definition of $D^\T$:  if $\alpha+1\in D$, then the extender $E^\T_{\alpha+1}$ is applied to a strict initial segment of $\M^\T_\beta$ where $\beta = T-pred(\alpha+1)$. For $\lambda$ limit, $deg^\T(\lambda)$ is the eventual values of $deg^\T(\alpha+1)$ for $\alpha+1\in [0,\lambda]_\T$. For a cofinal branch $b$ of $\T$, $deg^\T(b)$ is defined to be the eventual value of $deg^\T(\alpha+1)$ for $\alpha+1\in b$. We write $\mathcal{C}_k(\M)$ for the $k$-th core of $\M$. Sometimes, we confuse $\mathcal{C}_0(\M)$ with $\M$ itself.

\begin{definition}\label{realizability property}\index{$X$-realizability property} Suppose $\vec{\V}$ is an array. We say $\vec{\V}$ has the $X$-\textbf{realizability property} if for all $\a<lh(\V)$, $\vec{\V}\rest \a$ has the $X$-realizability property and whenever $g\subseteq Coll(\omega, <\kappa)$ is generic, in $V[g]$, whenever $\pi:\W\rightarrow \mathcal{C}_k(\V_\eta)$ is a weak $k$-embedding and $\T$ are such that 
\begin{enumerate}
\item $X\subseteq rng(\pi)$
\item $\W, \T\in HC$,
\item $\T$ is $X$-approved normal, $k$-maximal iteration of $\W$ that is above $\d^{\W^b}$,
\end{enumerate}
then one of the following holds (in $V[g]$).
\begin{enumerate}
\item $\T$ is of limit length and there is a cofinal well-founded branch $c$ such that $c$ has no drops in model (i.e. $D^\T\cap b = \emptyset$); letting $l=deg^\T(b)$, there is a weak $l$-embedding $\tau:\M^\T_c\rightarrow \mathcal{C}_l(\V_\eta)$ such that $\pi\rest \W=\tau\circ \pi^\T_c$.
\item $\T$ is of limit length and there is a cofinal well-founded branch $c$ such that $c$ has a drop in model, and there is $\b<\eta$ and a weak $l$-embedding $\tau:\M^\T_c\rightarrow \mathcal{C}_l(\V_\b)$ such that $\tau\rest (\M^\T_c)^b=\pi\rest (\M^\T_c)^b$, where $l=deg^\T(c)$.
\item $\T$ has a last model and letting $\gg=lh(\T)-1$, $[0, \gg]_\T\cap \mathcal{D}^\T=\emptyset$ and there is a weak $l$-embedding $\tau:\M^\T_\gg\rightarrow \mathcal{C}_l(\V_\eta)$ such that $\pi\rest \W=\tau\circ \pi^\T$, where $l = deg^\T(\gamma)$.
\item  $\T$ has a last model and letting $\gg=lh(\T)-1$, $[0, \gg]_\T\cap \mathcal{D}^\T\not =\emptyset$ and for some $\b<\eta$ there is a weak $l$-embedding $\tau:\M^\T_\gg\rightarrow \mathcal{C}_l(\V_\b)$ such that $\tau\rest (\M^\T_\gg)^b=\pi\rest (\M^\T_\gg)^b$, where $l = deg^\T(\gamma)$.
\end{enumerate}
When the above 4 clauses hold we say that $\T$ is $(\pi, \vec{\V})$-realizable. In cases where $\vec{\V}$ is clear from the context, we omit it from our notation.$\myqedhere$
\end{definition}

\begin{definition}\label{honest wzc} Suppose $\R$ is a weakly $X$-suitable hod premouse. We say $\R$ is \textbf{honest} if there is an array $\vec{\V}=(\V_\a:\a\leq \eta)$ at $\mu$ with the $X$-realizability property such that $\R, \vec{\V}\in H_{\kappa^{+4}}$, the following conditions hold. 
\begin{enumerate}
\item Either $\V_\eta=\R$ or there is a $X$-validated iteration $p$ of $\V_\eta$ of limit length such that $\pi^{p, b}$ exists and $\R=\textrm{m}^+(p)$.
\item $\vec{\V}$ is small if and only if  $\V_\eta \not=\R$.
\end{enumerate}
If $\R$ is honest and $\vec{\V}$ is as above then we say that $\vec{\V}$ is an honesty certificate for $\R$.$\myqedhere$\index{honest certificate} 
\end{definition}

We fix $\R, \vec{\mathcal{V}}, p$ as in Definition \ref{honest wzc} and the models $(\M_\eta, \N_\eta : \eta \leq \Upsilon)$ are constructed by the above $X$-validated sts construction over $\R$. We proceed to verify the failure of $(ii)_\eta, (iii)_\eta, (iv)_\eta$ for $\eta\leq \Upsilon$.

%

\subsubsection{Failure of $(ii)_\eta$} \label{failure of ii}

Let us first verify $(ii)_\eta$ fails. Towards contradiction assume that there is some model $\W^*$ that appears in the $X$-validated sts construction such that $\W^*$ is not $X$-validated. Let $\W$ be the least such model. 

Suppose first that $\W$ is $\W=\M_\a$ for some $\a$. Suppose first $\a$ is a limit ordinal. Let $U$ be an $X$-good hull such that $ \{\R,  \W\}\subseteq U$ and $(\M_\b,\N_\beta:\b<\a)\in U$. We have that $\M_\b,\M_\b$ are $X$-validated for every $\b<\a$. Let $(\K_\xi: \xi\leq \a_U)=\pi^{-1}_U(\M_\b:\b<\a)$. Fix $\T\in \K_{\a_U}$ according to $S^{\K_{\a_U}}$.  We need to see that $\T$ is $X$-approved. Fix $\xi<\a_U$ such that $\T\in \K_\xi$ and is according to $S^{\K_\xi}$. Then $\pi_U(\T)\in \M_{\pi_U(\xi)}$ and is according to $S^{\M_{\pi_U(\xi)}}$. Therefore, $\T$ is $X$-approved. In this case, it is also easy to see that $\N_\alpha = \mathcal{C}(\M_\alpha)$ is $X$-validated.

Suppose next that $\a=\b+1$. Suppose the least model that is not $X$-validated is $\N_\a$. We must have that $\M_\a$ is $X$-validated and that all models $(\M_\xi,\N_\xi : \xi \leq \beta)$ are $X$-validated. Let now $U$ be an $X$-good hull such that $\{\R, \W\}\subseteq U$. But then $\pi_U^{-1}(\N_\a)=\mathcal{C}(\pi^{-1}_U(\M_\a))$. It then follows that $\pi_U^{-1}(\N_\a)$ is $X$-approved (see \rprop{preservation of zv under embeddings}).

We now assume that $\W=\M_\a$ for some $\alpha$. Suppose first that $\M_\a=(\N_{\a-1}, b)$ where $b$ is a branch. Then it follows from the definition of $X$-validated sts constructions that $\M_\a$ is $X$-validated. The case that $\M_\a=(\N_{\a-1}, E)$ for some extender $E$ is trivial as no new iterations of $\R$ have been introduced. 

Finally suppose $\M_\a=\J_1(\N_{\a-1})$. If $\M_\a$ is not $X$-validated then it is because there is $\T\in \M_\a-\N_{\a-1}$ such that $\T$ is according to $S^{\M_\a}$ yet $\T$ is not $X$-validated. Let $\xi=\sup\{\zeta: \T\rest \zeta\in \N_{\a-1}\}$. Then all  proper initial segments of $\T\rest \xi$ is in $\N_{\a-1}$ and hence, all of the proper initial segments of $\T\rest \xi$ are $X$-validated. Because $\T$ is not $X$-validated, $\xi+1\leq lh(\T)$. 

The following is the key point. There is no limit ordinal $\b\in (\xi, lh(\T))$. This is because to define $[0, \b]_\T$ we need to ``leave behind" a level that at the minimum is a model of $\sf{ZFC}$, while there is no such level between $\M_\a$ and $\N_{\a-1}$. Thus, it must be that $lh(\T)=\xi+n$ for some $n\in [1, \omega)$. Let then $m$ be least such that $\T\rest \xi+m$ is $X$-validated but $\T\rest \xi+m+1$ is not. We then have three cases. 

The first case is the following. 
\begin{enumerate}
\item $\pi^{\T\rest \xi, b}$ is defined, 
\item $\T_{\geq \xi}$ is above $(\M^\T_\xi)^b$, and
\item  letting $\Q=\M^\T_{\xi+m}$, $\cp(E_{\xi+m}^\T)=\d^{\Q^b}$.
\end{enumerate}
Let now $U$ be an $X$-good hull such that $\{\R, \W\}\subseteq U$. Let $E=\pi^{-1}(E_{\xi+m}^\T)$ and let $\Y$ be the least node of $\T_U$ to which $E$ must be applied. Using \cite[Proposition 9.11]{sargsyan2019exact}, fix $\b$, an $l<\omega$, and a weak $l$-embedding $k:\Y\rightarrow \mathcal{C}_l(\V_\b)$ such that $X\subseteq rng(k)$. Using the fact that $(\T_U)_{\geq\Y}$ is $(k, \vec{\V}\rest \b)$-realizable, we can find $\gg\leq \b$ and a weak $n$-embedding $\tau:\Q_U\rightarrow \mathcal{C}_n(\V_\gg)$ such that $X\subseteq rng(\tau)$. Therefore, as $\tau(E)$ is $X$-validated, letting $Z=k[\Y^b]\cap \mathcal{P}^-$, there is $Y$ an extension of $Z\oplus X$ such that $Ult(\Y^b, E)=\Q^X_Y$. Hence, $\T$ is $X$-validated.  

The next possibility is when $\T_{\geq \xi}$ is either a tree of finite length based on $\pi^{\T_{\geq \xi}, b}(\R^b)$ or it only uses extenders with critical points $>\pi^{\T_{\geq \xi}, b}(\d^b)$. The second case is easy because the tree $\T_{\geq \xi}$ is finite. The first case follows from an argument similar to the one given above. 

Finally we could have that $\xi+1=lh(\T)$, where $\xi$ is a limit ordinal, and for cofinal set of $\b<\xi$, letting $\gg_\b=pred_T(\b+1)$, $\cp(E_\b)=\d^{\pi^\T_{0, \gg_\b}(\d^{\R^b})}$. This case, however, easily follows from the direct limit construction. This completes our argument that $(ii)_\eta$ fails.

\subsubsection{Failure of $(iii)_\eta$} \label{failure of iii}

We now consider $(iii)_\eta$. Suppose $\R$ is honest as witnessed by $(\vec{\V}, p)$. Then we say $\T$ is a $X$-validated iteration of $\R$ if $p^\frown \T$ is a $Z$-validated iteration of $\V_\eta$ where $\eta+1=lh(\vec{\V})$. The array $\vec{\V}$ typically comes from an $X$-validated sts construction or a hybrid $K^c$-construction (described later).

Let $\R$ be honest as witnessed by $(\vec{\V}, p)$. Suppose $\T$ is a normal tree on $\R$ such that $\pi^{\T,b}$ exists and $\delta$ is a Woodin cardinal of $\pi^{\T,b}(\mathcal{P})$ and $\delta^*$ is a Woodin cardinal of $\mathcal{P}$. Then we have:
\begin{itemize}
\item cof$(\delta^*)  < \kappa$ because there is a hod pair $(\Q,\Lambda)\in \mathcal{F}$ and $\gamma$ such that $\Q\models ``\xi$ is Woodin" and $\delta^* = \pi_{\Q,\infty}^{\Lambda}(\xi)$.
\item if $\delta > \textrm{sup}(\pi^{\T,b}[\delta^\mathcal{H}])$, then by \cite[Lemma 8.11]{sargsyan2019exact}, cof$(\delta) = \textrm{cof}(\sf{ord}$$(\mathcal{P}))$. By our assumption and Lemma \ref{lem:small}, $\textrm{cof}(\sf{ord}$$(\mathcal{P})) < \gamma$, so cof$(\delta) <\gamma$.
\end{itemize}

Suppose $\M$ is the $\eta$-th model appearing in the $X$-validated sts construction over $Y$ based on $\R$ ($\M=\N_{\eta-1}$) and $\T^*\in \M$ is a $\sf{uvs}$ iteration of  $\R$ such that the indexing scheme requires that we index a branch of $\T^*$ at $\sf{ord}$$(\M)$. We need to show that there is a branch $b$ of $\T^*$ such that $(\M, b)$ is $X$-validated. Because of \rprop{unique branches above mu} and \cite[Proposition 9.11]{sargsyan2019exact}, there can be at most one such branch. The proof of this is given in \cite[Proposition 9.5]{sargsyan2019exact}. We outline the main points of the proof in the following. Suppose $c$ is a certified branch; we need to see that $c$ is $(X,\vec{\V})$-embeddable. If $U$ is an $X$-good hull such that $(\R, \vec{\V},p, \T^*,c)\in U$. Let $\V'=\pi_U^{-1}(\mathcal{C}_n(\V_\eta))$ and $k:\R_U\rightarrow \mathcal{C}_n(\V_\eta)$ be such that $\pi_U\rest \V'=k\circ \pi^{p_U}$. We now suppose that there is a cofinal branch $d$ of $\T^*_U$ such that for some $\b\leq\eta$ there is $m:\M^{\T^*_U}_d\rightarrow \V_\b$ and $\Q(d, \T^*_U)$-exists. Let $\M=\Q(d, \T^*_U)$ and $\N=\Q(c_U, \T_U^*)$. Both $\M$ and $\N$ are $X$-approved. Let $\S_0=\textrm{m}^+(\T^*_U)$. If we could conclude that $\M=\N$ then we would get that $c_U=d$, and that would finish the proof. To conclude that $\M=\N$, \cite[Proposition 9.5]{sargsyan2019exact} argues that $\S_0$ is not infinitely descending. Otherwise, there is a sequence $(p_i, \S_i: i<\omega)$ witnessing that $\S_0$ is infinitely descending such that for some $\b<\eta$ and for some $i_0<\omega$ for every $i<j\in (i_0, \omega)$ there are weak $n_i$-embeddings $m_i:\S_i\rightarrow \mathcal{C}_{n_i}(\V_\b)$ such that $m_i=m_j\circ \pi^{p_i}$. The existence of a sequence as in the claim above gives us a contradiction, as the sequence must have a well-founded branch. The uniqueness proof is similar.

Because $\T^*$ is $\sf{uvs}$, we have a normal iteration $\T\in \M$ with last model $\S$ such that $\pi^{\T}$ is defined and a normal iteration $\U$ based on $\S^b$ such that $\T^\frown \U=\T^*$. At this point, we assume that $(ii)_\eta$ fails, so we have that $\M$ is $X$-validated and therefore, $\T$ is $X$-validated. Also, we can assume that $\U$ is not based on $\S|\xi$ where $\xi=\sup(\pi^\T[\d^\mH])$, as otherwise the desired branch of $\U$ is given by $\Sigma$. 

We now show that $\U$ has a branch $b$ such that $(\M, b)$ is $X$-validated. Given an $X$-good hull $U$ such that $\{\M, \T, \S, \U\}\subseteq U$, let $b_U=\Sigma_W(\pi^{-1}_U(\U))$ where $W$ is any extension of $X$ such that $\pi^{-1}_U(\S^b)=\Q^X_W$. First we claim that for all $U$ as above,

\begin{claim}\label{claim:b_U_in} $b_U\in M_U$. 
\end{claim}
\begin{proof} Fix then a $U$ as above. Suppose first that $\Q(b_U, \U_U)$ doesn't exist. As we are assuming $\U$ is not based on $\S|\xi$, the remark above gives that $\cf(\d(\U)) < \gamma$. Because $M_U$ is $\gamma$-closed it follows that $b_U\in M_U$.

Suppose next that $\Q(b_U, \U_U)$ exists. It is easy to see that letting $\Sigma_U$ be the $\pi_U$-pullback of $\Sigma$, then Lp$^{\Omega,\Sigma_U}(a)\in M_U$\footnote{This is true because Lp$^{j(\Omega),\Sigma}(A)\in V$ and $\Sigma = j(\Sigma)^j$.}, where $a = \pi_U^{-1}(A)$ and $A$ is a transitive set that codes $\{\M,\T,\S,\U\}$. 

Let $Y=U\cap \mH$. Clearly $Y$ is an extension of $X$ and because $\M$ is $X$-validated, we must have $W^*$ an extension of $X\cup Y$ such that $\S_U^b=\Q^X_{W^*}$. Notice that because $\Sigma_{W^*}$ is computable from $\Sigma_U$ and because $Lp^{\Omega, \Sigma_U}(a)\in M_U$, we must have that $\Q(b_U, \U_U)\in M_U$. Hence, $b_U\in M_U$.
\end{proof}

Suppose first that $\cf(lh(\U))>\omega$. In this case, let $U$ be as above and set $c=\pi_U(b_U)$. Then $c$ is the unique well-founded branch of $\U$ and hence, for any $X$-good hull $Z$ such that  $U\cup \{(\M, c), U\}\in Z$, $c_Z=b_Z$. Hence, $(\M, c)$ is $X$-validated (see  \rprop{one hull witness for premice}). 

Suppose then $lh(\U)=\omega$. We now claim that there is an $X$-good hull $Z$ such that for all $X$-good hull $Y$ such that $Z\cup \{\M, Z\}\in Y$, $\pi_{X, Y}(b_X)=b_Y$. Assuming not we get a  continuous chain $(X_\a: \a<\kappa)$ such that
\begin{enumerate}
\item $\M, \U\in X_0$,
\item for all $\a<\kappa$, $X_{\a+1}$ is an $X$-good hull,
\item for all $\a<\kappa$, $X_{\a}\cup \{X_{\a}\}\in X_{\a+1}$,
\item for all $\a<\kappa$, $\pi_{X_{\a+1}, X_{\a+2}}(b_{X_{\a+1}})\not =b_{X_{\a+1}}$.
\end{enumerate}
Let $\nu\in (\gamma, \kappa)$ be an inaccessible cardinal such that $X_\nu\cap \kappa=\nu$. Fix $\a<\nu$ such that 
\begin{center}
$\sup(b_{X_\nu}\cap  rng(\pi_{X_\a, X_\nu}))=\d(\U_{\X_\nu})$.
\end{center}
 As $\cf(lh(\U_{X_\nu}))=\omega$ this is easy to achieve. For $\b\in [\a, \nu)$ let $c_\b=\pi_{X_\a, X_\nu}^{-1}[b_{X_\nu}]$. Let for $\b\in [\a, \nu]$, $W_\b$ be such that $\S^b_{X_\b}=\Q^X_{W_\b}$. It follows that $c_\b$ is according to $\pi_{X_\b, X_\nu}$-pullback of $\Sigma_{W_\nu}$. Because $\Sigma_{W_\b}$ depends only $\S^b_{X_\b}$, we have that $c_\b=b_{X_\b}$\footnote{The  $\pi_{X_\b, X_\nu}$-pullback of $\Sigma_{W_\nu}$ is a strategy of the form $\Sigma_Y$ where $\Q^X_Y=\S^b_{X_\b}$ by Lemma \ref{lem:condensing_set}.} It follows that for all $\b<\gg\in [\a, \nu)$, $\pi_{X_\b, X_\gg}(b_{X_\b})=b_{X_\gg}$. This contradiction proves the claim.
 
 Fix now a $Z$ as above. Set $c=\pi_Z(b_Z)$. The the above property of $Z$ guarantees that $(\M, c)$ is $X$-validated. Indeed, fix an $X$-good hull $U$ such that $\M, c\in U$. Let $Y$ be an $X$-good hull such that $Z\cup U \cup \{Z, U\}\in Y$. Then $\pi_{U, Y}(c_U)=\pi_{Z,Y}(b_Z)=b_Y$. It follows that $c_U=\Sigma_W^{\pi_{U, Y}}(\pi_U^{-1}(\U))$ where $W$ is such that $\S^b_Y=\Q^X_W$. Hence, $c_U=b_U$.

This completes the outline of the proof that the $X$-validated sts construction over an honest $\R$ cannot break down because $(iii)_\eta$ holds for some $\eta$.

\subsubsection{Failure of $(iv)_\eta$}\label{failure of iv}

Now we sketch the proof that $(iv)_\eta$ fails. Suppose that the first time the $X$-validated sts construction over $X$ breaks down because $(iv)_\eta$ holds, where $\eta$ is the least such. This means that setting $\W=\N_{\eta-1}$ 
\begin{enumerate}
\item $\W$ is $X$-validated,
\item there is a $\sf{nuvs}$ $\T\in \W$ that is according to $S^{\W}$ and is such that one of the following holds:
\begin{enumerate}
\item there is a cofinal branch $b\in \W$ such that $\Q(b, \T)$ exists and is authenticated in $\W$ but $(\W, b)$ is not $X$-validated, or
\item there is a $\Q$-structure $\Q\in \W$ that is authenticated but there is no branch $b\in \W$ such that $\Q(b, \T)=\Q$.
\end{enumerate}
\end{enumerate}
We first show that case 2.a holds as it illustrates the main arguments and utilizes most of the concepts introduced above. Let $\b$ be such that $\W|\b$ authenticates $b$. Thus $\W|\b$ is a model of $\sf{ZFC}$ in which there is a limit of Woodin cardinals $\nu$ and the derived model of $\W|\b$ at $\nu$ has a strategy for $\Q(b, \T)$ that is $\W|\b$-authenticated. 

Fix now a $X$-good hull $U$ such that $(\R, \W, \T) \in U$ and $\T_U^\frown \{b_U\}$ is not a correctly guided $X$-realizable iteration of $\R_U$. Because $\W$ is $X$-validated, we can assume that $\T_U$ is correctly guided $X$-realizable iteration. It must then be that $\Q(b_U, \T_U)$ is not $X$-approved. 

We show that $\N=_{def}\Q(b_U, \T_U)$ is $X$-approved of depth 1. The proof of depth n is the same, we will leave the rest to the reader. To start with, notice that since $\T_U$ itself is correctly guided $X$-realizable, we have that $\S=\textrm{m}^+(\T_U)$ is weakly $X$-suitable. To prove that $\N$ is $X$-approved of depth 1 we need to show that if $\U\in \N$ is according to $S^{\N}$ then $\U$ is $X$-realizable. 

Fix then $\alpha\in R^\U$ and $\X = \M^\U_\alpha$. First we show that there is $Z$, an extension of $X$ such that $\Q^X_Z=\X^b$. Because $\T_U^\frown \{b_U\}$ is authenticated inside $\W_U|\b_U$, we must have an iteration $\Y$ of $\R_U$ according to $S^{\W_U}$ with last model $\R_1$ such that there is an embedding $k:\X^b\rightarrow \R_1^b$ with the property that $\pi^{\Y, b}=k\circ\pi^{\U_{\leq \X}, b}\circ \pi^{\T_U,b}_{b_U}$. Because $\Y$ is $X$-realizable, we must have $Y$ an extension of $X$ such that $\R_1^b=\Q^X_Y$. Composing $k$ with $\tau^X_Y$ we have that $\X^b=\Q^X_Z$ for some $Z$.

The rest is similar. If $\U^*$ is the longest initial segment of $\U_{\geq \X}$ that is based on $\X^b$ then there are $\Y$ and $k$ as above such that $\U^*$ is according to $k$-pullback of $S^{\W_U}_{\R_1^b}$. But because $\W_U$ is $X$-approved, $S^{\W_U}_{\R_1^b}$ is a fragment of $\Sigma_Y$ where $Y$ is as above. Hence, $\U^*$ is according to $\Sigma_Z$ for $Z$ extending $X$ as above (see \ref{lem:condensing_set}). This finishes proof of 2.a.

Next we assume that case 2.b holds. Let then $U$ be a $X$-good hull such that $(\R, \W, \T, \Q)\in U$. Because $\T$ is $X$-validated, we have that the $\pi_U$-realizable branch of $d$ of $\T_U$ is cofinal. Suppose then $\Q(d, \T_U)$ exists. Then because it is $X$-approved, we must have that $\Q(d, \T_U)=\Q_U$ (see \cite[Proposition 9.5]{sargsyan2019exact}). It follows that $d\in \W_U$, and so $\pi_U(d)$ is our desired branch.

We claim that $\Q(d, \T_U)$ exists. Suppose not. We regard $\T_U$ as a tree on $\W_U$ based on $\R_U$. First set $\N=\M^{\T_U}_d$ and $j=\pi^{\T_U}_d$. Note that $d\cap D^{\T_U}=\emptyset$ and $\pi^{\T_U}_d(\delta^{\R_U}) = \delta(\T_U)$; so $\N\models ``\delta(\T_U)$ is a Woodin cardinal". We have that $j(\Q_U)\in \N$ and is authenticated in $\N$. Let $\gg=j(\b_U)$. Then $\N|\gg$ has Woodin cardinals bigger than $\delta(j(\T_U))$. Let $\d$ be the least one that is $>\d(j(\T_U))$. We can now iterate $\N$ below $\d$ and above $j(\Q(b_U,\T_U))$ to make $\Q_U$ generic for the extender algebra at the image of $\d$. This iteration produced $i:\N\rightarrow \N_1$ such that $\cp(i)>\d(j(\T_U))$. Letting $h\subseteq Coll(\omega, i(\d))$ be $\N_1$-generic such that $\Q_U\in \N_1[h]$, we can find $l:\Q_U\rightarrow i(j(\Q_U))=j(\Q_U)$ such that
\begin{itemize}
\item $l\in \N_1[h]$, and 
\item $l\rest (\textrm{m}^+(\T_U))^b = \pi^{j(\T_U),b}$.
\end{itemize}
As $\N_1[h]\models ``j(\Q_U)$ is authenticated and has an authenticated strategy", $\N_1[h]\models ``\Q_U$ has an  authenticated iteration strategy", and hence $\Q_U$ is definable in  $\N_1[h]$ from objects in $\N_1$. It follows that $\Q_U\in \N_1$, implying that $\N_1\models ``\d(\T_U)$ is not a Woodin cardinal". Hence, $\N\models ``\d(\T_U)$ is not a Woodin cardinal". Therefore, $\Q(d, \T_U)$ exists. This completes the proof of case 2.b, and the proof that $(iv)_\eta$ fails.

The proofs above give us the following corollary. 

\begin{corollary}\label{cor:not_small}
Suppose the $X$-validated sts construction above breaks down because of $(iv)_\alpha$ for some $\alpha$, then
\begin{enumerate}
\item $\vec{\mathcal{V}}$ is not small (so $\mathcal{V}_\eta = \R$), and
\item letting $(\T,b)\in \N_{\alpha-1}$ witnessing the construction fails because of $(iv)_\alpha$, then $\N_{\alpha-1} \models ``\delta^\R$ is not a Woodin cardinal".
\end{enumerate}
\end{corollary}

The following corollary is also useful. We will use it in later sections. See \cite[Propositions 10.6--10.7]{sargsyan2019exact} for the corresponding versions of these two corollaries.

\begin{corollary}\label{cor:validated_suitable}
Suppose $\vec{\V}$ is a small array with the $X$-realizability property. Then either 
\begin{enumerate}
\item $\V_\eta$ has a $X$-validated iteration strategy

or
\item there is a $X$-validated $\sf{nuvs}$ iteration $p$ of $\V_\eta$ such that $\m^+(p)$ is $Z$-suitable\footnote{See \rdef{z-full}.}. 

\end{enumerate}
\end{corollary}

\begin{proof}
We outline the argument here:  the argument earlier in this section shows that if $p$ is an $X$-validated $\sf{uvs}$ of $\M_\xi$ of limit length then there is a unique branch $b$ of $p$ such that $p^\frown\{b\}$ is $X$-validated. Therefore, since picking $X$-validated branches is not defining an iteration strategy for $\M_\xi$, we must have an $\sf{nuvs}$ $X$-validated iteration $p$ of $\M_\xi$ which does not have a $X$-validated branch.\footnote{Note that there may not be any $\Q$-structure for $p$.} We now claim that $\textrm{m}^+(p)$ is a $X$-suitable hod premouse. Indeed, suppose there is some $X$-validated sts premouse $\Q$ extending $\R=_{def}\textrm{m}^+(p)$ such that $\Q$ is a $\Q$-structure for $p$. Let then $U$ be an $X$-good hull such that$\{\vec{\V}, p, \Q\}\in U$. It is not hard to see, using the fact that $\vec{\V}$ is a small array with the $X$-realizability property, that there is $\b\leq \xi$, a branch $b$ of $p_U$ such that $\Q(b, p_U)$ exists and a weak $l$-embedding $k:\M_b^{p_U}\rightarrow \mathcal{C}_l(\M'_\b)$ for an appropriate $l$. It follows that $\Q(b, p_U)$ is $X$-approved and hence, $\Q(b, p_U)=\Q_U$. Because $\Q_U\in M_U$, we have that $b\in M_U$. Then $c=_{def}\pi_U(b)$ is a (cofinal) branch of $p$ such that $p^\frown\{c\}$ is $X$-validated.  
\end{proof}

We end this section by defining the following. 
\begin{definition}\label{dfn:stop_prematurely}
We say that the $X$-validated sts construction over $Y$ based on $\R$ \textit{stops prematurely} if $\Upsilon$ is the least such that the following hold for $\M_\Upsilon$:

\begin{enumerate}[(i)]

\item There is an increasing sequence $(\delta_n: n<\omega.2)$ of Woodin cardinals above $\delta^\P$ such that $\delta^\P$ is the least $<\delta_0$-strong and $(\delta_n : n < \omega.2)$ are the only Woodin cardinals above $\delta_0$.
\item There are no extenders $E$ on the $\M_\Upsilon$-sequence such that there is some $n$ such that cr$(E) \leq \delta_i < \rm{lh}(E)$.

\item $\M_\Upsilon$ is an $X$-validated sts hod premouse over $Y$ based on $\R$.

\item $\M_\Upsilon$ is $E$-active with top extender $F$ such that cr$(F)>\delta_n$ for all $n<\omega.2$.
\item $\rho_{\omega}(\M_\Upsilon) \geq \sf{ord}$$(\R)$.
\end{enumerate}
$\myqedhere$
\end{definition}
We note that (ii) easily follows from (i), but for clarity, we make it explicit. In the later sections, we will obtain a contradiction (using stacking mice techniques) from the assumption that the construction does not stop prematurely and $\Upsilon = \kappa^{+++}$. The subsequent several sections will show that the construction does not fail and stop prematurely. From this, we then show that there must be a model of $\sf{LSA}$. 


\section{Hybrid $K^c$-constructions} \label{sec:hybrid_kc}

We continue with the notations of the previous section. Now we describe the hybrid $K^c$-construction. We first describe the ``bottom structure" that the hybrid $K^c$-construction is built on top of. We say $\mathcal{S}$ is \textit{$K^c$-appropriate}\index{$K^c$-appropriate model} if $\S = \P$ or $\mathcal{S} = \R$ where the following hold for $\R$:
\begin{itemize}
\item $\P\lhd \R$ and $o(\R)$ is Woodin in $\R$, where $o(\R) = \textrm{sup} \{\xi : E^\R_\xi \neq \emptyset \wedge \textrm{crit}(E^\R_\xi)=\delta^\P\}$.
\item $\R$ is a sound, $X$-validated sts premouse over $(\R|o(\R))^\sharp$ such that $\J_1[\R]\models ``o(\R)$ is not Woodin". 
\item There exists a (unique) $X$-validated strategy for $\R$.
\end{itemize}  

We build in $V[G]$ a sequence $(\M_\xi,\N_\xi: \xi \leq \Upsilon)$ of levels of our $K^c$-construction such that $\N_0=\M_0 = \R$ for a $K^c$-appropriate $\R$, $\N_\xi = \mathcal{C}(\M_\xi)$ for all $\xi\leq \Upsilon$ and $\Upsilon\leq\kappa^{+++}$. It will be clear from the construction that $\N_\xi,\M_\xi\in V$ for all $\xi$.  

We maintain during the construction that $\M_\xi$ (and hence $\M^Y_\xi$) is \textit{small}, i.e. either $\M_\xi$ has no Woodin cardinal $>\delta^\P$ or else letting $\delta$ be  the least such Woodin, then $\rho(\M_\xi) \leq\delta$ and $\M_\xi$ defines a failure of Woodinness of $\delta$. Note that this implies that 
\begin{center}
$\J_1[\N_\xi] \models $`` there are no Woodin cardinals $> \delta^\P$".
\end{center}
If $\M_\xi$ is not small or that we fail to construct an $X$-validated strategy for $\M_\xi$, then we let $\Upsilon = \xi$ and stop the construction.

Suppose we have constructed $\M_\xi, \N_\xi$ for some $\xi$. Let $\gamma_\xi = o(\N_\xi)$ be the supremum of indices of extenders on the $\N_\xi$-sequence with critical point $\delta^\P$ if there are such extenders; otherwise, let $\gamma_\xi = \sf{ord}$$(\P)$. Suppose $\gamma_\xi < \sf{ord}$$(\N_\xi)$ and let $\gamma_\xi\leq \lambda_\xi \leq \sf{ord}$$(\N_\xi)$ be such that $\rho_{\omega}(\N_\xi)\geq\lambda_\xi$. Suppose there is a stack $\vec{\T}\in \N_\xi$ based on $\N_\xi|\lambda_\xi$ according to the internal strategy $S^{\N_\xi}$ such that $S^{\N_\xi}(\vec{\T})$ is undefined.\footnote{By the notations earlier in the book, $\N_\xi|\lambda_\xi$ is a complete hod initial segment of $\N_\xi$.} Suppose also $\vec{\T}$ is such that the theory developed above (Chapter 3) dictates that a cofinal branch $b$ for $\vec{\T}$ needs to be added to $\N_\xi$ and $\N_\xi$ is so that $(\N_\xi, B_b)$ is amenable. We call such a tuple $(\N_\xi,\lambda_\xi,\vec{\T})$ \textit{branch-ready}.\index{branch-ready tuple}



For a branch-ready tuple $(\N_\xi,\lambda_\xi,\vec{\T})$ in our hybrid $K^c$-construction, we need to see that $(\vec{\T},b)$ is $X$-validated and that $(\N_\xi,B_b)$ is an $X$-validated hod premouse; this is accomplished by constructing an external $X$-validated strategy $\Lambda$ for $\N_\xi|\lambda_\xi$ and let $b = \Lambda(\T)$, see below for a more detailed discussion. Furthermore, we need to construct an external $X$-validated strategy for $(\N_\xi,B_b)$.  So we maintain that $S^{\N_\xi}$ \footnote{Recall, this is the sequence of branch predicates that codes up some internal strategy of $\N_\xi$.} above $\lambda_\xi$ is according to: 
\begin{enumerate}[(a)]
\item either the  strategy $\Sigma_{\xi}$ of $\N_\xi|\lambda_\xi$, where $\Sigma_{\xi}$ is the canonical $Q$-structure guided, $X$-validated strategy of $\N_\xi|\lambda_\xi$ if $\sf{ord}$$(\N_\xi) < \sf{ord}$$(\M_2^{\Sigma_{\xi},\sharp})$;
\item or else the canonical $\Sigma_{\xi}$-strategy $\Lambda_{\xi}$ of $\N_\xi|\lambda_\xi = \M_2^{\Sigma_{\xi},\sharp}(\N_\xi|\epsilon)$, where $\Sigma_{\xi}$ is the canonical $Q$-structure guided,  $X$-validated strategy of $\N_\xi|\epsilon$.
\end{enumerate}
We let $\Psi_\xi$ denote $\Sigma_{\xi}$ in case (a), and $\Lambda_{\xi}$ in case (b).


We will discuss the construction of the strategy in (a), or (b) in the next section. At this point, we assume it exists and just want to extend the internal strategy of $\N_\xi$ one more step. The key thing we want to maintain here is that the indexed branch $b$ for $\N_\xi$ is according to the strategies $\Psi_\xi$. Roughly, what we need to do to construct such a strategy is as follows. 

In the following, we write $\forall^* Y$ to mean ``for some club $\mathfrak{C}$, $Y\in \mathfrak{C}\cap \mathfrak{S}_{\phi,\Omega}$". What we show in the next sections is that $\forall^* Y$ such that $Y$ is $X$-good, we can construct an $X$-realizable strategy $\Psi_{Y,\xi}$ as in case (a), (b). $\Psi_{\xi}$ is then determined from the strategies $\Psi_{Y, \xi}$ by the procedure described in the previous section. $b=\Psi_\xi(\vec{\T})$ if there is an $X$-good $Z$ such that letting $b_Z = \Psi_{Z,\xi}(\vec{\T}_Z)$, then for any $X$-good hull $Y$ such that $\N_\xi, Z\in Y$, 
\begin{center}
$\pi_{Z, Y}[b_Z] \subseteq b_Y$ and $\tau^X_{Z}[b_Z] \subseteq \tau^X_Y[b_Y] \subseteq b$. 
\end{center}
We say that $b$ described above is \textit{suitable} for $(\N_\xi,\lambda_\xi,\vec{\T})$. We define the above notions in a similar manner for $\M_\xi$; for brevity, we also use the symbols $\gamma_\xi, \lambda_\xi$ for $\M_\xi$ when no confusion arises from the context.

\begin{remark}\label{rem:strategy}

If such a strategy does not exist, we will stop the hybrid $K^c$-construction over $\R$ at stage $\xi$. In this case, there will be an $X$-validated iteration $p$ witnessing this. We will then switch to the sts $X$-validated construction over $\M(p)^\sharp$.

The argument in the previous section also shows that if such a strategy exists and is constructed according to the aforementioned procedure, then it must agree with $S^{\N_\xi}$ (similarly for $\M_\xi$). $\myqedhere$
\end{remark}

\begin{remark}\label{rem:g_organized}
The reason we have case (b) is because we want our hod mice to be g-organized in the sense of \cite{trang2013}. g-organization ensures that $S$-constructions go through as discussed in Chapter 6. $\myqedhere$
\end{remark}

The procedure above allows us to define the object Lp$^{\Psi_{\xi}}(\N_\xi)$ in the case $\gamma_\xi < \lambda_\xi=\sf{ord}$$(\N_\xi)$ as follows. In the definition below, following \cite{Trang2015PFA}, we let $\rm{Lp}^{\Psi_{Y,\xi},\Omega}(\N_\xi^Y)$ be the union of sound, $\Psi_{Y,\xi}$-mouse $\R$ over $\N_\xi^Y$ such that $\rho(\R)\leq \sf{ord}$$(\N^Y_\xi)$ with unique $X$-realizable iteration strategy (above $\sf{ord}$$(\N^Y_\xi)$) in $\Omega$. 
\begin{definition}\label{def:Lp}
Suppose $\gamma_\xi < \lambda_\xi=\sf{ord}$$(\N_\xi)$. We let Lp$^{\Psi_{\xi}}(\N_\xi)$ be the union of $\N_\xi \lhd \M$ such that $\rho_\omega(\M)\leq \sf{ord}$$(\N_\xi)$, $\M$ is $\sf{ord}$$(\N_\xi)$-sound, and $\forall^* Y$, $Y$ is  $X$-good and contains all relevant objects, $\pi_Y^{-1}(\M)\lhd \rm{Lp}^{\Psi_{Y,\xi},\Omega}(\N_\xi^Y)$.  Sometimes, we write Lp$^{\Sigma_{\N_\xi}}(\N_\xi)$ for Lp$^{\Psi_{\xi}}(\N_\xi)$.


We also define $\Sigma_{\M_\xi}$ and Lp$^{\Sigma_{\M_\xi}}(\M_\xi)$ in a similar manner. $\myqedhere$

\end{definition}

\begin{remark}
By \cite[Lemma 3.78]{Trang2015PFA} and Lemma \ref{lem:condensing_set}, $\forall^* Y \pi_Y^{-1}(\rm{Lp}$$^{\Sigma_{\N_\xi}}(\N_\xi)) = \rm{Lp}^{\Psi_{Y,\xi},\Omega}(\N^Y_\xi)$. We also remind the reader, for reasons mentioned before, levels of $\rm{Lp}^{\Psi_{Y,\xi},\Omega}(\N^Y_\xi)$, etc are g-organized in the sense of \cite{trang2013}. $\myqedhere$

\end{remark}

\begin{definition}[Relevant extender]\label{def:relevant_extenders}\index{relevant extender}
Suppose $F$ is on the $\N_\xi$-extender sequence for some $\xi\leq \Upsilon$. We say that $F$ is \textbf{relevant} if $F = G\cap \N_\xi$\footnote{If $\cp(F) = \delta^\P$, we confuse $F$ with its amenable code for $G\cap \N_\xi$ and in the case $\cp(F)>\delta^\P$, we think of $F$ as a ``map" as in \cite{Zeman}.} for $G$ a  correctly backgrounded extender. $\myqedhere$
\end{definition}


We let $\N_\xi^+=\J_\gamma[\N_\xi]$ for $\gamma$ being least such that 
\begin{enumerate}
\item either $\J_\gamma[\N_\xi]$ is not sound or $\rho(\J_\gamma[\N_\xi])<\rho(\N_\xi)$,

\item or else $\J_\gamma[\N_\xi]$ satisfies $\sf{ZFC}$$^-$ and there is a correctly backgrounded extender $F$ that coheres the $\J_\gamma[\N_\xi]$-sequence, 
\item or else there are $\lambda_\xi\geq \gamma_\xi$, $\vec{\T}\in \J_\gamma[\N_\xi]$ such that $(\J_\gamma[\N_\xi],\lambda_\xi,\vec{\T})$ is branch-ready.
\end{enumerate}

\begin{definition}[Hybrid $K^c$-construction]\label{def:N_xi}\index{hybrid $K^c$-construction} The models $\M_\xi, \N_\xi$ are defined as follows: for all $\xi\leq \Upsilon$,
\begin{enumerate}[(a)]
\item if $\xi =0$, then $\N_\xi = \M_\xi = \R$ for a $K^c$-appropriate $\R$;
\item if $\xi$ is limit, let $\M_\xi$ be liminf$_{\xi^* < \xi} \N_{\xi^*}$;
\item if $\xi = \xi^*+1$, the following hold:
\begin{enumerate}[(i)]
\item if $\N_{\xi^*}$ is passive and there is a correctly backgrounded extender $F$ that coheres the $\N_{\xi^*}$-sequence, then let $\M_\xi = (\N_{\xi^*},F)$.\footnote{See Definition \ref{def:correctly_backgrounded_ext}. The uniqueness of the extender $F$ with crit$(F)>\delta^\P$ follows from a standard bicephalous argument, cf. \cite[Section 7]{FSIT}, and tools developed in the previous Chapters. If crit$(F) = \delta^\P$, the uniqueness of $F$ follows from the proof of Lemma \ref{lem:next_certified_extender}.}
\item if $\N_{\xi^*}$ is passive and there is some $\vec{\T}\in \N_{\xi^*}$ and some $\lambda$ such that $\vec{\T}$ is based on $\N_{\xi^*}|\lambda$ and $(\N_{\xi^*},\lambda,\vec{\T})$ is branch-ready, then letting $b$ be given by the procedure above, we set $\M_\xi = (\N_{\xi^*},B_b)$. 
\item  if $\M_{\xi^*}$ is passive and cases (i) and (ii) do not hold, then we set $\M_{\xi} = \N_{\xi^*}^+$.
\end{enumerate}
\item $\N_\xi = \mathcal{C}(\M_\xi)$.
\end{enumerate}
$\myqedhere$
\end{definition}

\begin{remark}\label{rem:strat_activation_level}



If $\N_{\xi^*}$ is as in (c)(i), we say that $\N_{\xi^*}$ is \textit{extender-ready}. We give priority to adding extenders, that is, if $\N_\xi$ is both an extender-ready level and a branch-ready level, then we are in case (c)(i).

If $\N_\xi$ is weakly suitable or that $\Sigma_{\M_\xi}$ is not defined and $p$ witnesses this, then as mentioned above, we let $\Upsilon = \xi$ and stop the construction. We then start a new $X$-validated sts construction over $\N_\xi$ in the first case and over $\M(p)^\sharp$ in the second case.$\myqedhere$

\end{remark}

Let $Y\in \mathfrak{S}_{\phi,\Omega}$, so $Y\cap \P$ is an honest extension of $X$. Let $\pi_Y$ be the uncollapse map and $\N^Y_\xi = \pi^{-1}_Y(\N_\xi)$. We recall that $\Lambda$ is the \textit{$X$-realizable strategy of $\N^Y_\xi$}\index{$X$-realizable strategy} if whenever $\vec{\T}$ is according to $\Lambda$, $i:\N^Y_\xi\rightarrow \Q$ is the iteration map according to $\Lambda$, where $\Q = \M^{\vec{\T}}_\alpha$ and $\alpha\in R^{\vec{\T}}$, then there is some $Z$, honest extension of $X$, such that $\Q^b = \Q^X_Z$, and the map $k: \Q\rightarrow \N_\xi$ defined as: for $f\in \N^Y_\xi$, $a\in (\delta^\Q)^{<\omega}$,
\begin{center}
$k(i(f)(a)) = \pi_Y(f)(\pi^{\Lambda_{\vec{\T},\Q}}_{\Q,\infty}(a))$
\end{center}
is well-defined, elementary\footnote{By this, we mean if $\vec{\T}$ is a $k$-maximal stack then $k$ is a weak $k$-embedding in the sense of \cite{FSIT}.}, $k\circ i = \pi_Y$, and $k\rest \delta^\Q = \pi^{\Lambda_{\vec{\T},\Q}}_{\Q,\infty}\rest \delta^\Q$. Similar definitions are given for $\M_\xi, \M^Y_\xi$.

We maintain as part of the construction the following for $\xi < \Upsilon$:\index{$(1)_\xi-(3)_\xi$}
\begin{enumerate}
\item[$(1)_\xi$] 
\begin{enumerate}[(a)]
\item $\M_\xi$ is $X$-validated.
\item $\forall^*Y$, there is an $X$-realizable strategy for $\M^Y_\xi$, called $\Sigma_{\M^Y_\xi}$.\footnote{It will be clear that $\Sigma_{\M^Y_\xi}$ is unique. In essence, $(1)_\xi(b)$ is equivalent to the statement that the natural $X$-realizable strategy defined in Section \ref{sec:iterability} is total; in particular, all iterates according to this strategy are $X$-approved.} So $\Sigma_{\M_\xi}$ is an $X$-validated strategy  for $\M_\xi$. 
\item $\forall^*Y$, $\Sigma_{\M^Y_\xi}$ has (locally) strong branch condensation and is (locally) strongly $\Omega$-fullness preserving.
\end{enumerate}
\item[$(2)_\xi$] $\rho_{\omega}(\M_\xi)\geq \sf{ord}$$(\P)$. In other words, $\sf{ord}$$(\P)$ is $(\delta^\P)^+$ in $\M_\xi$ and in $\N_\xi$.

\item[$(3)_\xi$] $\M_\xi$ is solid and universal. So $\N_\xi$ is sound.
\end{enumerate}

See $(\dag\dag)_\xi$ in Section \ref{sec:iterability} for how these statements are precisely maintained. $\Upsilon$ is the least such that one of the conditions fail at $\Upsilon$ or that the hybrid $K^c$-construction stops prematurely. By Section \ref{sec:iterability}, if one of $(1)_\Upsilon, (2)_\Upsilon, (3)_\Upsilon$ fails, then in fact $(1)_\Upsilon(b)$ fails.

In the next section, we will obtain a contradiction (using stacking mice techniques) from the assumption that $\Upsilon = \kappa^{+++}$. 

\begin{remark}\label{rem:K-c}

The extender sequence of $\N_\xi$ utilizes two indexing schemes: the cutpoint indexing scheme (for extenders with critical point $\delta^\P$) and the Jensen indexing scheme (for extenders with critical point $>\delta^\P$). This follows from the definition of correctly backgrounded extenders for relevant extenders. The Jensen indexing scheme could be replaced by the Mitchell-Steel indexing scheme, but we choose not to do so out of convenience; we want to quote direct results from \cite{JSSS} and \cite{ANS} as well as using results of Chapter 11.$\myqedhere$


\end{remark}

\section{Stacking mice}\label{sec:stacking_mice}

Suppose the constructions described above do not stop prematurely and therefore result in a model $\M_{\Upsilon}$ such that $\sf{ord}$$(\M_{\Upsilon})=\kappa^{+++}$ (see Lemma \ref{lem:not_all_the_way}). Let $\N = \N_{\Upsilon} = \mathcal{C}(\M_\Upsilon)$. It is clear that $\N = \M_\Upsilon$.\footnote{One way to see that is to recall that $\M_\Upsilon = \textrm{liminf}_{\xi\rightarrow \Upsilon} \M_\xi$. Since $\Upsilon = \kappa^{+++}$, and $\sf{PFA}$ holds, $\sf{ord}$$(\M_\Upsilon)=\kappa^{+++}$ is a limit cardinal in $\M_{\Upsilon}$ and $\M_{\Upsilon}\models \sf{ZFC}$. So $\rho(\M_\Upsilon) = \kappa^{+++}$.} So $\N$ is an $X$-validated hod premouse or an $X$-validated sts premouse. There are two possible cases on how we reach such an $\N$.  In the first case, we must have alternated the two constructions (the $X$-validated sts construction and the hybrid $K^c$-construction) until we reach a $K^c$-appropriate $\R$ and the rest of the construction is the hybrid $K^c$-construction with $\N_0 = \R$ and $\N_\Upsilon = \N$\footnote{It could be that $\R = \N_\Upsilon$.}; in the second case, we reached a suitable $\R$ and $\N$ is obtained by the $X$-validated sts construction over $\R$.

Let $\delta^\N > \delta^\P$ be the unique $\gamma$ such that $\N\vDash ``\delta^\P$ is strong to $\gamma$ and $\gamma$ is Woodin" if it exists; otherwise we let $\delta^\N= 0$. Note that by the remarks above, which is a consequence of our smallness assumption $(\dag)$, $\delta^\N$ is a strong cutpoint of $\N$. Following \cite{JSSS}, we define the following \textit{stack of hod mice above $\N$}. The following definition takes place in $V[G]$ but it is easily seen that $S(\N)\in V$ (see Lemma \ref{lem:stack_facts}).
\begin{definition}\label{def:stack}\index{stack of hod mice}
Let $\delta$ denote $\delta^\N$. Let $S(\N)$\index{$S(\N)$} be the stack of sound $X$-validated hod premice $\M$ if $\delta=0$ or else $X$-validated sts-premice $\M$ extending $\N$ such that $\rho_\omega(\M)=\sf{ord}$$(\N) = \kappa^{+++}$ and for every $\M^*$ embeddable into $\M$ via $\pi_{\M^*}$ such that $|\M^*|<\kappa$, $X\cup \{X, \P^-, \P,\N,\delta \}\subset \rm{rng}(\pi_{\M^*})$, $\rm{rng}(\pi_{\M^*})\cap \P$ is an honest extension of $X$, $\M^*$ is $(\omega_1+1)$-iterable above $\pi^{-1}(\delta)$. Furthermore, in the case $\delta=0$, the strategy for $\M^*$ is $X$-realizable and in the case $\delta > 0$, the strategy witnesses $\M^*$ is an $X$-approved sts mouse. $\myqedhere$
\end{definition}

$\N$ is the $\kappa^{+++}$-th model in the hybrid $K^c$-construction or the $X$-validated sts construction, and hence is passive. In the first case, one can show easily that items $(1)-(3)$ hold for $\N$. In the above, if $\delta=0$, then $\M^*$ has an $X$-realizable strategy $\Lambda_{\M^*}$ such that $\Lambda_{\M^*}\rest\rm{HC}\in \Omega$, and  $\Lambda_{\M^*}$ is locally $\Omega$-fullness preserving and has local strong branch condensation (see next section). Furthermore, by the fact that $X$ is a condensing set and Section \ref{sec:sts_constructions}, $\Lambda_{\M^*}$ witnesses $\M^*$ is an $X$-approved hod mouse. In particular, if $E$ is on the $\M$-sequence such that cr$(E)=\delta^\P$ and lh$(E)\geq o(\N)$, then for every such $\M^*$ as above such that $E\in \rm{rng}(\pi_{\M^*})$, letting $\nu$ be the length of $\pi^{-1}_{\M^*}(E)$, then for any $a\in[\nu]^{<\omega}$, $A\in \powerset(\delta^\P)^{|a|}\cap \P$ such that $(a,A)\in E\cap \rm{rng}(\pi_{\M^*})$, then $\pi^{\Lambda_{\M^*}}_{\M^*||\nu,\infty}(\pi_{\M^*}^{-1}(a)) \in A$. In the case $\delta>0$, we demand as part of Definition \ref{def:stack} that $\M^*$ is iterable above $\pi^{-1}(\delta)$ as an $X$-approved sts mouse; note also that $\delta$ is a strong cutpoint of $\M$. The following facts about $S(\N)$ more or less follow immediately from results in \cite{JSSS}. 
\begin{lemma}\label{lem:stack_facts}
Suppose $\Upsilon = \kappa^{+++}$ and $\N = \N_{\Upsilon}$. 
\begin{enumerate}[(i)]
\item For $\M_0,\M_1\in S(\N)$, either $\M_0\unlhd \M_1$ or $\M_1\unlhd \M_0$. 
\item For all $\M\unlhd S(\N)$, there is some $\R\lhd S(\N)$ such that $\M\lhd \R$. In particular, $S(\N)\vDash \sf{ZFC}^-$.
\item cof$(\sf{ord}$$(S(\N)))\geq \kappa^{+++}$.
\end{enumerate}
\end{lemma}
\begin{proof}
(i) and (ii) are analogs of \cite[Lemma 3.1]{JSSS} and \cite[Lemma 3.3]{JSSS} respectively and follow straightforwardly from the condensation lemma, Theorem \ref{thm:condensation_lemma}. The point is that if $\delta^\N=0$, then the theory developed above allows us to perform comparisons (and shows that no strategy disagreement can occur); otherwise, the construction above $\delta$ is an $X$-validated sts construction over a fixed suitable $\R$, so the comparison is again an extender comparison (and is above $\sf{ord}$$(\R) = \delta$).\footnote{See Corollary \ref{cor:solidity_universality} for a similar argument with more details. The point is that the comparisons involve two $X$-approved sts mice, so no strategy disagreements appear.} Therefore, by an easy application of Theorem \ref{thm:condensation_lemma} (see also \cite[Lemma 1.3]{JSSS}), letting $\pi: H\rightarrow H_{\kappa^{+4}}$ be elementary and such that $\{\M_1,\M_2, X,\R\}\in H$, $H\cap \kappa^{+++}\in \kappa^{+++}$, then $\pi^{-1}(\M_0) \lhd \N$ and $\pi^{-1}(\M_1)\lhd \N$. By elementarity of $\pi$, (i) holds for $\M_0,\M_1$. The proof of (ii) follows from \cite[Lemma 3.3]{JSSS} and the discussion above.

(iii) follows from the proof of \cite[Theorem 3.4]{JSSS} with obvious modifications, noting that by our assumption, $\kappa^{+}, \kappa^{++}$ are $\kappa$-closed, and $2^{\kappa^{++}} = \kappa^{+++}$ in $V[G]$. We note that $\kappa$ plays the role of $\omega$, $\kappa^{+++}$ plays the role of $\kappa$ in that proof and all hulls taken are closed under $\kappa$-sequences in this case. 
\end{proof}

\section{Iterability of lsa-small, non-lsa type levels}\label{sec:iterability}
First, we verify that (1)-(3) holds for $\xi = 0$ and $\M_0 = \N_0 = \P$. By Lemma \ref{lem:small}, no $\P^- \lhd \M \lhd \P$ projects across $\delta^\P$; also $\P \vDash \sf{ZFC}^-$, and hence $\rho_{\omega}(\P) = \sf{ord}$$(\P)$. 

\begin{lemma}\label{lemma:base_case}
(1)-(3) hold for $\xi = 0$. 
\end{lemma}
\begin{proof}
Fix $Y$ as in the statement of (1); so $\M^Y_0 = \N^Y_0$. Let $\delta_Y = \delta^{\N^Y_0}$ and $\Sigma^Y_0 = \Sigma_{\N^Y_0}$ be $\Sigma_Y$. By definition, $\Sigma^Y_0$ has branch condensation as it is the join of strategies with those properties. Furthermore, note that $\Sigma^Y_0$ acts on $\N^Y_0$ in the following way.  Let $(\Q,\vec{\mathcal{T}}) \in I(\N_0^Y,\Sigma^Y_0)$ and let $i:\N_0^Y \rightarrow \Q$ be the iteration map and $\Sigma_{\Q,\vec{\T}}$ be the $\vec{\T}$-tail of $\Sigma_0^Y$.

Suppose $x\in \Q$, then there is some $f\in \N_0^Y$ and $a\in i(\delta_Y)^{<\omega}$ such that
\begin{center}
$x = i(f)(a)$.
\end{center}

Let $k:\Q \rightarrow \N_0$ be defined as follows:
\begin{center}
$k(i(f)(a) = \pi_Y(f)(\pi^{\Sigma_{\Q,\vec{\T}}}_{\Q(i(\delta_Y)),\infty}(a))$,
\end{center}
for any $f\in \N^Y_0$ and any $a\in i(\delta_Y)^{<\omega}$. Note that since $X$ is a condensing set and $i\circ \pi_{X,Y}\rest \delta_X$ is according to $\Sigma_0^X$, rng$(k)$ is an honest extension of $X$. By Lemma \ref{lem:condensing_set}, $k$ is well-defined, $\Sigma_1$-elementary (and cofinal), $k\circ i = \pi_Y$, and $k\rest \delta^\Q = \pi^{\Sigma_{\vec{\T},\Q}}_{\Q|\delta^\Q,\infty}\rest \delta^\Q$. It is clear that this is the only way to define $k$; the uniqueness of $\Sigma^Y_0$ also follows.

We remark that local strong branch condensation is just branch condensation in this case. Now to see that $\Sigma^Y_0$ is $\Omega$-fullness preserving, it suffices to show $\Q$ is $\Omega$-full. But this follows from the definition of condensing sets and the fact that $Y$ and rng$(k)$ are honest extensions of $X$. Also, we get local strong $\Omega$-fullness preservation.

We have shown (1). (2) holds by the remark immediately before the lemma and (3) follows from (2) and (1) by Remark \ref{rem:K-c}. The usual proof of universality and solidity goes through with the iterability proved in (1).
\end{proof}
Now we inductively in $\xi$ prove the following: 
\begin{enumerate}
\item[$(\dag\dag)_\xi$] suppose $(1)_{\xi^*}-(3)_{\xi^*}$ hold for all $\xi^* < \xi$, then $(1)_\xi(a)$ holds and suppose $(1)_\xi(b)$ holds, then  $(1)_\xi(c), (2)_\xi, (3)_\xi$ hold.\index{$(\dag\dag)_\xi$}
\end{enumerate}

We continue with proving $(\dag\dag)_\xi$. So suppose $(\dag\dag)_{\xi^*}$ holds for all $\xi^* <\xi$. Now, we verify $(\dag\dag)_\xi$. Let $Y\in V$ be an honest extension of $X$; as before, we assume also $Y=Y^*\cap \M_\xi$ for some $Y^*\prec H_{\kappa^{+4}}^V$. We recall that $\M_\xi$ (and hence $\M^Y_\xi$) is \textit{small}, i.e. either $\M_\xi$ has no Woodin cardinal $>\delta^\P$ or else letting $\delta$ be  the least such Woodin, then $\rho(\M_\xi) \leq\delta$ and $\M_\xi$ defines a failure of Woodinness of $\delta$. Note that this implies that 
\begin{center}
$\J_1[\N_\xi] \models $`` there are no Woodin cardinals $> \delta^\P$".
\end{center}

If $\M_\xi$ is not small and that $\M_\xi$ is of lsa type, i.e. $\M_\xi = \textrm{m}^+(\M_\xi|\delta) \models ``\delta$ is Woodin and $\delta^\P$ is $<\delta$-strong", then we stop the hybrid $K^c$-construction, set $\Upsilon = \xi$, and switch to the $X$-validated sts construction over $\M_\xi$.

If $(1)_\xi(b)$ fails, we stop the hybrid $K^c$-construction and let $\Upsilon = \xi$. In this case Corollary \ref{cor:validated_suitable} shows that there is an $X$-validated $\sf{nuvs}$ $p$ of $\M_\Upsilon$ such that $\R = \textrm{m}^+(p)$ is $X$-suitable. We then continue with our $X$-validated sts construction over $\R$ or over some transitive $W$ containing $\R$. 

So we assume $\M_\xi$ is small and $(1)_\xi(b)$. We verify the other clauses. We now define the strategy $\Sigma_{\M^Y_\xi}$ for $\M^Y_\xi$; we sometimes write $\Sigma^Y_\xi$ for $\Sigma_{\M^Y_\xi}$.\footnote{Technically, we construct $\Sigma_{\M^Y_\xi}$ in $V[G]$ but  $\Sigma_{\M^Y_\xi}\cap V \in V$ and  $\Sigma_{\M^Y_\xi}$ does not depend on the choice of $G$. This will be clear from the construction of  $\Sigma_{\M^Y_\xi}$. So in effect, we are constructing an invariant name $\dot{\Sigma}$ in $V$ whose interpretation in $V[G]$ is  $\Sigma_{\M^Y_\xi}$ for any $G$. This justifies our notation $\Sigma^Y_\xi$.} We write $x^Y$ for $\pi^{-1}_Y(x)$ for $x\in \M_\xi\cap \rm{rng}(\pi_Y)$.


\begin{definition}[Normal form]\label{def:normal_form}\index{stacks of normal form}
An iteration $((\P_\alpha,\vec{\T}_\alpha) \ | \ \alpha < \eta)$ on $\P_0 = \M^Y_\xi$ is said to be in \textbf{normal form} if the following hold:
\begin{enumerate}[(i)]
\item $\vec{\T}_\alpha$ is a stack of normal trees with base model $\P_\alpha$ and last model $\P_{\alpha+1}$.
\item If $\lambda \leq \eta$ is limit, $\P_\lambda = \textrm{lim}_{\alpha<\lambda}\P_\alpha$.
\item Either $\vec{\mathcal{T}}_\alpha$  uses no extenders in the top block of $\P_\alpha$ or its images or  $\P_{\alpha+1} = \textrm{Ult}(\P_\alpha,E)$ for some extender $E$ on the $\P_\alpha$-sequence with cr$(E)=\delta^{\P_\alpha}$ or else $\vec{\mathcal{T}}_\alpha$ is completely above $\delta^{\P_\alpha}$.
\item If $\eta=\alpha+1$ for some $\alpha$, then for all $\beta<\alpha$, $\vec{\T_\beta}$ does not drop.
\end{enumerate}
$\myqedhere$
\end{definition}

We define  $\Sigma^Y_\xi$ for stacks in normal form. We say that a stack $((\P_\alpha,\vec{\T}_\alpha) \ | \ \alpha < \eta)$ in normal form, where $\P_0 = \M^Y_\xi$, is according to  $\Sigma^Y_\xi$ if: letting $\tau_0=\pi_Y\rest\P_0$, $i_{\gamma,\tau}:\P_\gamma\rightarrow \P_\tau$ be iteration maps, and $\kappa^{\P_\gamma} = i_{0,\gamma}(\kappa^{\P_0})$, where $\kappa^{\P_0} = \delta^{\P_0^b} = \delta^\P$, then 
\begin{enumerate}[(A)]
\item $\P_\alpha$ is $X$-approved and  there are maps $\tau_\alpha:\P_\alpha\rightarrow \M_\xi$ for all $\alpha < \eta$;
\item for all $\gamma\leq\alpha<\eta$, if $i_{\gamma,\alpha}$ exists then $\tau_\gamma = \tau_\alpha \circ i_{\gamma,\alpha}$;
\item for all $\alpha<\eta$, letting $\Lambda_\alpha$ be the $\tau_\alpha$-pullback strategy and $\pi_{\P_\alpha|\kappa^{\P_\alpha},\infty}^{\Lambda_\alpha} = \tau_\alpha\rest\P_\alpha|\kappa^{\P_\alpha}$; furthermore, $\P_\alpha^b = \Q^X_Z$ for some $Z$ extending $X$;
\item if $\eta=\alpha+1$ and $\vec{\T}_\alpha$ drops, then there is a (unique) branch $b$ of $\vec{\T}_{\alpha}$, some $\xi'<\xi$, and a weak-deg$(b)$-embedding\footnote{deg$(b)$ is the degree of soundness of model corresponding to the last drop along $b$.} $\tau_{\eta}:\M^{\vec{\T}_\alpha}_b\rightarrow \M_{\xi'}$. Otherwise, there is a (unique) branch $b$ and map $\tau_\eta:\M^{\vec{\T}_\alpha}_b \rightarrow \M_\xi$ such that $\tau_\alpha=\tau_\eta\circ i_{\alpha,\eta}$.
\end{enumerate}  

It is clear how to extend $\Sigma^Y_\xi$ to all stacks of normal trees. This is because all stacks of normal trees on $\N^Y_\xi$ can be decomposed into stacks in normal form. We will need to define maps $\tau_\alpha$ in the definition of $\Sigma^Y_\xi$ in such a way that makes $\Sigma^Y_\xi$ an $X$-realizable strategy.

\begin{lemma}\label{lem:certified_extender}
Suppose $Y\prec H_{\kappa^{+4}}$ is a countable such that $Y\cap\P$ is an extension of $X$. Suppose $i: \M^Y_\xi \rightarrow \R$ and $\sigma:\R\rightarrow \M_\xi$ are such that $\pi_Y = \sigma\circ i$, and rge$(\sigma)\cap \P$ is an honest extension of $Y\cap P$. Let $\Lambda$ be the $\sigma$-pullback strategy on $\R$. Then:
\begin{enumerate}[(a)]
\item If $j:\R\rightarrow \S$ is a $\Lambda$-iteration based on $\R^b$ and suppose $\kappa^\S= \sup j[\kappa^\R] = j(\kappa^\R)$, then letting $\tau:\S\rightarrow \M_\xi$ be the map: $\tau(j(f)(a)) = \sigma(f)(\pi^{\Lambda_\S}_{\S|\kappa^\S,\infty}(a))$, where $f\in \R$, $a \in [\kappa^\S]^{<\omega}$, and $\Lambda_\S$ is the tail of $\Lambda$. Then $\tau$ is well-defined, elementary, and $\pi^{\Lambda_\S}_{\S|\kappa^\S,\infty} = \tau\rest (\S|\kappa^\S)$. As before (and later on in this chapter), $\kappa^\R = \delta^{\R_b}$ etc. is the cutpoint cardinal that begins the top block of $\R$.
\item Suppose $F$ is an extender on the $\R$-sequence with cr$(F) = \kappa^\R = i(\delta^\P)$.  Then $F$ is $\sigma$-certified over $(\R',\Lambda_{\R'})$, where $\R' = \R||lh(F)$. This means for $a\in lh(F)^{<\omega}$, $A\subseteq \kappa^\R$ in $\R$, $(a,A)\in F$ if and only if $\pi^{\Lambda_{\R'}}_{\R',\infty}(a)\in \sigma(A)$.
\end{enumerate}
\end{lemma}
\begin{remark}
The extender $F$ in the lemma is said to be \textbf{certified} for short instead of ``$\sigma$-certified over $(\R||lh(F),\Lambda_{\R||lh(F)})$".
\end{remark}
\begin{figure}
\centering
\resizebox{0.4\textwidth}{!}{
\begin{tikzpicture}[node distance=3cm, auto]
 \node (A) {$\M_{\xi}^Y$};
\node (B) [node distance=2cm,below of=A] {$\R$};
  \node (C) [right of=A] {$\M_{\xi}^Z$};
 \node (E) [node distance=2cm, above of=C] {$\M_{\xi}$};
 \node (F) [right of=B]{$\S$};
 \node (G)[right of=C]{$\W$};
  \draw[->] (A) to node [swap] {$i$} (B);
  \draw[->] (B) to node [swap] {$\sigma_Z$} (C);
  \draw[->] (A) to node {$\pi_{Y}$} (E);
  \draw[->] (A) to node {$\pi_{Y,Z}$} (C);
  \draw[->] (C) to node {$i_H$} (G);
  \draw[->] (C) to node {$\pi_Z$} (E);
  \draw[->] (B) to node {$i_F$} (F);
  \draw[->] (F) to node [swap]{$\tau_Z$} (G);
  \draw[->] (G) to node {$\psi$} (E);
  \draw[dashed,->, bend right=33] (F) to node  {$\tau$} (E);
\end{tikzpicture}
}
\caption{Hypothesis of Lemma \ref{lem:certified_extender}}
\label{LemmaExtendersOfIteratesAreCertified}
\end{figure}

\begin{proof}
(a) follows from Lemma \ref{lem:condensing_set} and the fact that the iteration map $j$ is continuous at $\delta^\R$.

For (b), first, note that $i$ is continuous at $(\delta^+)^{\N^Y_{\xi}}$ and is cofinal in $((\kappa^\R)^+)^\R$. This is because $\pi_{Y}$ is continuous at $(\delta^+)^{\N^Y_{\xi}}$ and is cofinal in $(\delta^+)^{\P}$. Finally, $F$ is total over $\R$; this follows from the continuity of $i$.

Now, let $\S = \rm{Ult}$$(\R,F)$, $i_F$ be the ultrapower map. Let $Y\prec Z$ be countable such that $Z\cap \P$ is an honest extension of $X$ and such that  rng$(\sigma)\subseteq \rm{rng}$$(\pi_Z)$ . Let $\sigma_Z = \pi_Z^{-1}\circ \sigma$. Let $H = \sigma_Z(F)$ \footnote{If $F$ is the top extender of $\R$, then by $\sigma_Z(F)$, we mean $\sigma_Z[F]$.} and $i_H: \M^Z_\xi \rightarrow \rm{Ult}$$(\M^Z_\xi,H)=_{\rm{def}} \W$ be the ultrapower map. Let $\tau_Z:\S\rightarrow \W$ be the copy map and $\psi:\rm{Ult}$$(\M^Z_\xi,H)\rightarrow \M_\xi$ be the map
\begin{center}
$\psi(i_H(f)(a)) = \pi_Z(f)(\pi^{\Psi}_{\N^Z_\xi||lh(H),\infty}(a))$,
\end{center}
where $\Psi$ is $(\Sigma^Z_\xi)_{\M^Z_\xi||lh(H)}$ Since $H$ is $\pi_Z$-certified over $(\N^Z_\xi||lh(H),\Psi)$, $\psi$ is well-defined, elementary, and $\psi\circ i_H = \pi_Z$. Now, 
\begin{center}
$\sigma = \psi\circ \tau_Z\circ i_F$,
\end{center}
so letting $\Lambda_\S$ be the $\psi\circ \tau_Z$-pullback strategy for $\S$, then by strategy coherence for hod mice, $\Lambda_\S$ agrees with $\Lambda_\R$ on $\R||lh(F)$ (see Remark \ref{rem:how_to_realize}(\ref{rem:countable_reflection})). Now let $\tau: \S\rightarrow \M_\xi$ be defined as follows: for all $a\in [lh(F)]^{<\omega}$ and $f\in \R$,
\begin{center}
$\tau(i^\R_F(f)(a)) = \sigma(f)(\pi^{\Lambda_\R}_{\R||lh(F),\infty}(a))$.
\end{center}
By Lemma \ref{lem:condensing_set}, $\tau$ is well-defined, elementary, and agrees with $\pi^{\Lambda_\S}_{\S,\infty}$ up to $\delta^\S$ and with $\pi^{\Lambda_\R}_{\R||lh(F),\infty}$ up to $\R||lh(F)$. This proves part (b).
\end{proof}

The following remarks summarize how we can inductively define maps $\tau_\alpha$ and hence define $\Sigma^Y_\xi$ on stacks in normal form.
\begin{remark}\label{rem:how_to_realize}
\begin{enumerate}[(i)]
\item \label{rem:realization} If $\vec{\T}_\alpha = \langle E \rangle$ for cr$(E) = \kappa^{\P_\alpha}$, then 
\begin{center}$\tau_{\alpha+1}(i^{\P_\alpha}_E(f)(a)) = \tau_\alpha(f)(\pi_{\P_\alpha||lh(E),\infty}^{\Lambda_\alpha}(a))$.\end{center} Lemma \ref{lem:certified_extender}(b) shows that $\tau_{\alpha+1}$ is well-defined, elementary,\footnote{If $\tau_\alpha$ is a weak $k$-embedding for some $k$, as is typical of realization maps, then so is $\tau_{\alpha+1}$.}, agrees with $\tau_\alpha$ up to $\kappa^{\P_\alpha}$ and with $\pi_{\P_\alpha||lh(E),\infty}^{\Lambda_\alpha}$ up to $lh(E)$.

\begin{figure}
\centering
\resizebox{0.5\textwidth}{!}{
\begin{tikzpicture}[node distance=3cm, auto]
  \node (A) {$\bar{\P}_0$};
  \node (B) [node distance=0.5cm,above of=A] {$\bar{\W}_0$};
  \node (C) [right of=A] {$\bar{\P}_{\alpha+1}$};
  \node (D) [node distance=0.5cm, above of=C] {$\bar{\W}_{\alpha+1}$};
  \node (E) [node distance=2cm, above of=D] {$\bar{\M}^Z_{\xi},\Phi$};
  \node (F) [node distance=0.5cm,above of=E] {$\bar{\W}^*$};
  \node (G) [node distance=4cm, above of=B] {$\P_0=\M^Y_{\xi}$};
  \node (H) [node distance=0.5cm,above of=G] {$\W_0=\W_Y$};
  \node (I) [right of=G] {$\P_\alpha$};
  \node (J) [right of=H] {$\W_\alpha$};
  \node (K) [right of=I] {$\P_{\alpha+1}$};
  \node (L) [right of=J] {$\ \ \ \ \W_{\alpha+1}=\W$};
  \node (M) [node distance=2cm, above of=I] {$\M^Z_{\xi}$};
  \node (N) [node distance=0.5cm, above of=M] {$\W^*= \W_Z$};
  \draw[->,bend right=45] (I) to node [swap] {$\tau_\alpha$} (M);
  \draw[->] (B) to node [swap] {$\pi$} (G);
  \draw[dashed,->] (A) to node {} (C);
  \draw[dashed,->] (A) to node {$\bar{\tau}_0$} (E);
  \draw[dashed,->, bend right=15] (F) to node {$\epsilon$} (H);
  \draw[->] (G) to node  {$i_{0,\alpha}=i$} (I);
  \draw[dashed,->] (H) to node [swap] {$\tau_0$} (N);
  \draw[->] (I) to node  {$i_E=i_{\alpha,\alpha+1}$} (K);
  \draw[->] (K) to node [swap] {$\tau_{\alpha+1}=\tau$} (M);
 \draw[dashed, ->,bend right=45] (C) to node {$\bar{\tau} = \bar{\tau}_{\alpha+1}$} (E);

\end{tikzpicture}
}
\caption{Sketch of Remark \ref{rem:how_to_realize}(\ref{rem:countable_reflection})}
\label{fig:Sketch}
\end{figure}

\item\label{rem:countable_reflection} With the exact same situation as (\ref{rem:realization}) and suppose cof$^V(\sf{ord}$$(\M_\xi))<\kappa$,\footnote{Cofinally many $\xi'$ has the property that cof$^V(o(\M_{\xi'}))<\kappa$. In our case, $\xi = \xi^*+1$, this holds because $\N_{\xi^*}\vDash \forall\xi \square_{\xi,2}$ by Chapter 11 and the hypothesis of Theorem \ref{thm:square_lsa}.} we claim that the $\S=_{\textrm{def}}\P_{\alpha+1}$-tail of $\Psi =_{\textrm{def}}\Sigma^Y_{\xi-1}$ agrees with $\Psi_\S$, the $\tau_{\alpha+1}$-pullback strategy of $\S$. This is \textit{strategy coherence at $\alpha+1$}. Suppose not. Write $\tau$ for $\tau_{\alpha+1}$ and $i$ for $i_{0,\alpha}$. This is basically combining the proof of Theorem 2.7.6 in \cite{ATHM} and Lemma \ref{lem:condensing_set} (see Figure \ref{fig:Sketch}). We briefly sketch it here. Let $Y\prec Z$ and $Z\in V$ be countable (in $V[G]$), such that $Z\cap \P$ is an honest extension of $X$. Furthermore, we assume $Y\cap \kappa\in \kappa$, $Y^{<|Y|}\subset Y$, and letting $\iota = \textrm{cof}^V(\sf{ord}$$(\M_\xi))$, then $\iota < |Y|$.\footnote{This is possible because $\kappa$ is strongly inaccessible.}  Let $\W_Y$ be a $\Psi$-hod mouse with cof($\lambda^{\W_Y})=\omega$ and $\W_Z = Ult(\W_Y,F)$, where $F$ is the $(\textrm{crit}(\pi_Y),\sf{ord}$$(\M^Z_{\xi}))$-extender derived from $\pi_{Y,Z}$. So letting $j=i_{\alpha,\alpha+1}\circ i$, $j$ extends to $j^+: \W_Y\rightarrow \W$ and $\tau$ extends to $\tau^+:\W\rightarrow \W_Z$, where $\M^Z_{\xi}\lhd \W_Z$ (this is because $\sf{ord}$$(\M^Y_{\xi})$ is a cutpoint in $\W_Y$ and $\pi_{Y,Z}$ is cofinal in $\sf{ord}$$(\M^Z_{\xi})$). Let $\pi:M\rightarrow H^V_{\kappa^{+4}}$ be the inverse of the transitive collapse of some elementary substructure of $H^V_{\kappa^{+4}}$  in $V$ containing all relevant objects such that $X\subset \textrm{ran}(\pi)$ and $|M| < |Y|$. For any $a\in H^V_{\kappa^{+4}}\cap \textrm{ran}(\pi)$, let $\bar{a} = \pi^{-1}(a)$. Let $\bar{g}\subseteq Col(\omega,\bar{\kappa})$ be $M$-generic with $g\in V$ and $\bar{\S}$, $\bar{j}$, $\bar{\tau}$ be the objects in $M[\bar{g}]$ witnessing the failure of the claim in $M[\bar{g}]$. Since $|M|<|Y|$ and $Y^{<|Y|}\subseteq Y$, there is a map $\epsilon: \bar{W_Z} \rightarrow W_Y$ such that $\pi\rest \bar{\W_Y} = \epsilon\circ \bar{\pi_Y}\rest \bar{\W_Y}$. Let $\Phi$ be the $\epsilon$-pullback of $\Psi$. By the proof of Theorem 2.7.6 in \cite{ATHM}, working in $M[\bar{g}]$, the uB code for $\bar{\Psi}$ gets moved to the uB code for its $\bar{\S}$-tail and also to the uB code for the $\bar{\tau}$-pullback of $\Phi$; by Lemma \ref{lem:condensing_set}, this is also the $\tau\circ \pi = \bar{\tau}\circ \pi$-pullback of $\Sigma^Z_{\xi-1}$. This is a contradiction.
\item If $\vec{\T}_\alpha$ is below $\delta^{\P_\alpha}$ then it is according to $\Lambda_\alpha$ and so $\tau_{\alpha+1}$ is given by the inductive assumption on $\Lambda_\alpha$. Strategy coherence at $\alpha+1$ is maintained here. See Lemma \ref{lem:certified_extender}(a).
\item If $\vec{\T}_\alpha$ is above $\delta^{\P_\alpha}$  then $\vec{\T}_\alpha$ is correctly guided.\footnote{Recall this means that for $\beta < lh(\vec{\T}_\alpha)$, letting $c = [0,\beta]_{\vec{\T}_\alpha}$, then $\Q(c,\vec{\T}_\alpha\rest \beta)$ exists and is $X$-validated.} The map $\tau_{\alpha+1}$ is given by the K$^c$-construction theorem (cf. \cite[Theorem 3.2]{ANS}) and our smallness assumption on the hod mice that we are constructing; in fact, using the argument in Section \ref{sec:sts_constructions} and $(1)_\xi(b)$, we get that the branch giving rise to $\tau_{\alpha+1}$ is the unique branch $b$ such that $\Q(b,\vec{\T}_\alpha)$ exists and is $X$-validated. 

\item Suppose $\lambda < \eta$ is limit. Let for $\alpha<\lambda$
\begin{center}
$\tau_\lambda(i_{\alpha,\lambda}(x)) = \tau_{\alpha}(x)$.
\end{center} 
So we get $\tau_\lambda:\P_\lambda \rightarrow \M^Y_{\xi}$ is such that for all $\alpha<\lambda$, $\tau_\alpha = \tau_\lambda\circ i_{\alpha,\lambda}$. Using the above argument, we get strategy coherence at $\lambda$. Finally, we verify that letting $\pi:\P_\lambda|\delta^{\P_\lambda}\rightarrow \M_\xi$ be the iteration maps by the $\tau_\lambda$-pullback strategy $\Lambda_\lambda$, $\pi = \tau_\lambda\rest \delta^{\P_\lambda}$. Let $\nu<\delta^{\P_\lambda}$. We note that $\Lambda_\lambda$ is the $\Lambda_\alpha$-tail by strategy coherence at $\lambda$. Let $i_{\alpha,\lambda}(\nu^*)=\nu$ for some $\alpha<\lambda$ and $\nu^*<\delta^{\P_\alpha}$. Then
\begin{center}
$\tau_\lambda(\nu) = \tau_\alpha(i_{\alpha,\lambda}(\nu^*)) = \pi_{\P_\alpha|\kappa^{\P_\alpha},\infty}(\nu^*) = \pi(i_{\alpha,\lambda}(\nu^*)) = \pi(\nu)$.
\end{center}
$ \myqedhere$
\end{enumerate}
\end{remark}

The following lemma gives some useful consequences regarding uniqueness of strategies.
\begin{lemma}\label{lem:uniqueness}
\begin{enumerate}[(i)]
\item Suppose $\pi:\Q\rightarrow \P$ is elementary such that $\rm{rng}(\pi)$ is an extension of $X$. Suppose $i:\Q\rightarrow \R$ is such that $i\rest \delta^\P$ is according to the $\pi$-pullback strategy $\Sigma^\pi$ and $\tau_0,\tau_1:\R\rightarrow \P$ are such that $\tau_0 \circ i =  \tau_1\circ i = \pi$. Then the $\tau_0$-pullback strategy $\Sigma^{\tau_0}$ is the same as the $\tau_1$-pullback strategy $\Sigma^{\tau_1}$. 
\item Suppose $Y$ is countable, elementary in $\M_\xi$ and $Y\cap \P$ is an (honest) extension of $X$. Suppose $n$ is such that $\omega\rho_{\M^Y_\xi}^{n+1} < \omega\rho_{\M^Y_\xi}^{n}$. Let $\Psi = \Sigma^Y_\xi$.
\begin{enumerate}[(a)]
\item Suppose $(\vec{\T},\R)\in I(\M^Y_\xi,\Psi)$ is such that $\pi^{\vec{\T},b}$ exists and $\tau: \M^Y_\xi \rightarrow \S$ is $\Sigma^{(n)}_0$ and cardinal preserving and $\S\unlhd \R$. Suppose $(\vec{\U},\Q)\in I(\M^Y_\xi,\Psi)$ is such that $\pi^{\vec{\U},b}$ exists, $\Q^b=\S^b$, and $\tau\rest \P_Y = \pi^{\vec{\U},b}\rest \P_Y$, then $\Psi_{\vec{\T},\S}^\tau = \Psi$.
\item Suppose $(\vec{\T},\R)$ is such that $\pi^{\vec{\T},b}$ exists and is according to $\Psi$.  Suppose $\mathcal{U}$ is a normal tree of limit length on $\R(\beta)$ according to $\Psi_{\vec{\T},\R}$, where $\R(\beta) \lhd_c^{hod} \R$.\footnote{Recall this means that $\R(\beta)$ is a complete layer of $\R$ and $\R(\beta)\neq \R$.} Suppose $c$ is a cofinal branch of $\U$ (considered as a tree on $\R$) and there is a map $\tau_c:\M^\U_c\rightarrow \M_\xi$ such that $\pi_Y\rest \P_Y =  \tau_c\circ \pi^\U_c\circ \pi^{\vec{\T},b}$. Then $c = \Psi_{\vec{\T},\R}(\U)$.
\end{enumerate}
\end{enumerate}
\end{lemma}
\begin{proof}
(i) follows straightforwardly from Lemma \ref{lem:condensing_set} (iv). The main point is that, letting $\Lambda_i$ be the $\tau_i$-pullback strategy $\Sigma^{\tau_i}$ (for $i=0,1$), then letting $\sigma_i:\R\rightarrow \P$ be
\begin{center}
$\sigma(i(f)(a)) = \pi(f)(\pi^{\Lambda_i}_{\R,\infty}(a))$
\end{center}
for $f\in \Q$ and $a\in \delta^\R$. Then $\sigma_i[\R]$ is an honest extension of $X$.

(ii)(b) follows easily from (i) and Remark \ref{rem:how_to_realize}(ii). For (ii)(a) (see Figure \ref{fig:low_levels}), suppose $\Psi_{\vec{\T},\S}^\tau \neq \Psi$, then by results of Section 4.7, there is a (minimal) low-level disagreement, i.e. there is $(\vec{\W},\R_0, \W^*)$ such that:
\begin{itemize}
\item $\vec{\W}$ is according to both strategies.
\item $\R_0$ is the last model of $\vec{\W}$.
\item $\W^*$ is a tree of limit length on $\R_0(\beta)$ for some $\R_0(\beta) \lhd_c^{hod} \R_0$.
\end{itemize}
Let $b = \Psi(\vec{\W}^\smallfrown \W^*)$ and $c = \Psi_{\vec{\T},\S}^\tau(\vec{\W}^\smallfrown \W^*)$. Let $\sigma:\Q\rightarrow \M_\xi$ be the realization map; hence 
\begin{equation}\label{eqn:one}
\pi_Y\rest \Q_Y = \sigma\circ \pi^{\vec{\U},b} = \sigma\circ \tau\rest\Q_Y. \footnote{Recall $\Q_Y = \pi_Y^{-1}(\P)$.}
\end{equation}
Let $\tau_b: \M^{\W^*,b}_b\rightarrow \P$ and $\tau_c:\M^{\W^*,b}_c\rightarrow \P$ be the natural realization maps. We have:
\begin{equation}\label{eqn:two}
\sigma \circ \tau\rest \Q_Y = \tau_c\circ \pi^{\W^*}_c\circ \pi^{\vec{\W}}
\end{equation} 
and 
\begin{equation}\label{eqn:three}
\pi_Y\rest \Q_Y = \tau_b\circ \pi^{\W^*}_b\circ \pi^{\vec{\W}}.
\end{equation}
By (i) and Equations \ref{eqn:one}, \ref{eqn:two}, \ref{eqn:three}, $b = c$. Contradiction.
\end{proof}

\begin{figure}
\centering
\resizebox{0.4\textwidth}{!}{%
\begin{tikzpicture}[node distance=3cm, auto]
 \node (A) {$\M^{Y,b}_{\xi}$};
\node (B) [node distance=2cm,above of=A] {$\S^b$};
  \node (C) [right of=A] {$\R_0$};
 \node (E) [node distance=4cm, above of=C] {$\M_{\xi}$};
 \node (F) [node distance=1.5cm, right of=B]{$\M^{\W^*}_b$};
 \node (G)[node distance=3cm, right of=F]{$\M^{\W^*}_c$};
  \draw[->] (A) to node [swap] {$\tau$} (B);
  \draw[->, bend left=45] (A) to node {$\pi^{\vec{\U},b}$} (B);
  \draw[->] (B) to node {$\sigma$} (E);
  \draw[->] (A) to node [swap]{$\vec{\W}$} (C);
  \draw[->] (C) to node {$\pi^{\W^*}_b$} (F);
  \draw[->] (C) to node [swap]{$\pi^{\W^*}_c$} (G);
  \draw[->] (F) to node [swap]{$\tau_b$} (E);
  \draw[->] (G) to node [swap]{$\tau_c$} (E);
 \end{tikzpicture}
 }
 \caption{Lemma \ref{lem:uniqueness} (ii)(a)}
 \label{fig:low_levels}
\end{figure}

\begin{lemma}\label{lem:locally_strong_branch_condensation}
Suppose $Y$ is countable, elementary in $\M_\xi$ and $Y\cap \P$ is an (honest) extension of $X$. Then $\Sigma^Y_\xi$ has locally strong branch condensation, and is $\Omega$-fullness preserving.
\end{lemma}
\begin{proof} 
$\Omega$-fullness preservation follows from the construction of $\Sigma^Y_\xi$ and the fact that $X$ is a condensing set (see Lemma \ref{lem:condensing_set}). We first prove branch condensation (see Figure \ref{fig:branch_condensation}). Suppose not. Let $\N=\M^Y_\xi$ and $\Psi=\Sigma^Y_\xi$ and suppose the following hold: there are stacks $\vec{\T}^\smallfrown \U$ and $\vec{\W}$ on $\N$ such that
\begin{itemize}
\item $\vec{\T}$ is via $\Psi$ with end model $\R$.
\item $\U$ is according to $\Psi_\R$, $\vec{\W}$ is according to $\Psi$, and $i = \pi^{\vec{\W}}:\Q_Y \rightarrow \Q$ is the iteration map. 
\item There are cofinal branches $b, c$ of $\U$ and $\pi:\M^\U_b\rightarrow \Q$ such that
\begin{enumerate}
\item $i = \pi \circ \pi^\U_b \circ i^{\vec{\T}}$.
\item $c = \Psi(\vec{\T}^\smallfrown \U)$.
\item $b \neq c$.
\end{enumerate}
\end{itemize}
\begin{figure}

\centering
\resizebox{0.8\textwidth}{!}{
\begin{tikzpicture}[node distance=3cm, auto]
 \node (A) {$\N$};
\node (B) [right of=A] {$\R$};
\node (L) [node distance=2cm, above of=B]{$\mathcal{Q}$};
  \node (C) [right of=B] {$\M^+(\U)$};
 \node (D) [node distance=1cm, above of=C] {$\M^\U_b, \Psi_0$};
 \node (E) [node distance=1cm, below of=C]{$\M^\U_c, \Psi_1$};
 \node (F)[node distance=4cm, above of=C]{$\M_\xi$};
 \node (G)[node distance=4cm, below of=C]{$\M_\xi$};
 \node (H)[right of = C]{$\mathcal{Y}$};
 \node (I)[right of =H]{};
 \node (J)[node distance=1cm, above of=I]{$\M^{\vec{\U}^*}_{b^*}$};
 \node (K)[node distance=1cm, below of=I]{$\M^{\vec{\U}^*}_{c^*}$};
  \draw[->] (A) to node [swap] {$\vec{\W}, i$} (L);
  \draw[->] (A) to node [swap] {$\vec{\T}$} (B);
  \draw[->] (L) to node {$\sigma$} (F);
  \draw[->, bend left=45] (A) to node [swap] {$\pi_Y$} (F);
  \draw[->] (A) to node {$\pi_Y$} (G);
  \draw[->] (B) to node {$\U,b$} (D);
  \draw[->] (B) to node [swap]{$\U,c$} (E);
  \draw[->] (D) to node [swap]{$\pi$} (L);
  \draw[->] (C) to node {$\vec{\W}^*$} (H);
  \draw[->] (H) to node {$\vec{\U}^*,b^*$} (J);
  \draw[->] (H) to node [swap]{$\vec{\U}^*,c^*$} (K);
  \draw[->] (J) to node [swap]{$\tau_{b^*}$}(F);
  \draw[->] (K) to node {$\tau_{c^*}$}(G);  

 \end{tikzpicture}
 }
 \caption{Branch condensation}
 \label{fig:branch_condensation}
\end{figure}

Let $\Psi_0$ be the $\pi$-pullback strategy of $\Psi_{\vec{\W},\Q}$ and $\Psi_1$ be $\Psi_{\vec{\T}^\smallfrown \U ^\smallfrown c}$. Recall $\textrm{m}^+(\U) = \M(\U)^\sharp$. We first show:
\begin{equation}\label{eqn:lower_level_disagreement}
\Lambda_0=_{\rm{def}} (\Psi_0)^{sts}_{\textrm{m}^+(\U)} = (\Psi_1)^{sts}_{\textrm{m}^+(\U)}=_{\rm{def}} \Lambda_1.
\end{equation}

In the case there is $\Q\unlhd \textrm{m}^+(\U)$ which is a $Q$-structure for $\delta(\U)$ then $(\Psi_0)^{sts}_{\textrm{m}^+(\U)} = (\Psi_0)_{\textrm{m}^+(\U)}$ and similarly for $\Psi_1$. We assume this is not the case; otherwise, the argument is similar and simpler. 

Let $\sigma:\Q^b\rightarrow \P$ be the $\pi_Y$-realization map, so that 
\begin{center}
$\pi_Y\rest \Q_Y = \sigma\circ \pi \circ\pi^{\vec{\T}^\smallfrown \U,b}$.
\end{center}
In the above, we note that $\pi^{\vec{\T}^\smallfrown \U,b}$ exists and is the same as $\pi_c^{\vec{\T}^\smallfrown \U,b}$ and this map does not depend on the choice of the cofinal branch; so $\pi^{\vec{\T}^\smallfrown \U,b}$ is also $\pi_b^{\vec{\T}^\smallfrown \U,b} = \pi^{\U,b}\circ \pi^{\vec{\T}}$.

By results of Section 4.7, if (\ref{eqn:lower_level_disagreement}) fails, then there is a minimal disagreement $(\vec{\W}^*,\mathcal{Y})\in B(\textrm{m}^+(\U),\Lambda_0)\cap B(\textrm{m}^+(\U),\Lambda_1)$. Note that $\mathcal{Y}$ is of successor type and $(\Lambda_0)_{\vec{\W}^*,\mathcal{Y}(\alpha)} = (\Lambda_1)_{\vec{\W}^*,\mathcal{Y}(\alpha)}$ for all $\mathcal{Y}(\alpha)\lhd_c^{hod} \mathcal{Y}$. Furthermore, there is a stack $\vec{\U}^*$ on $\mathcal{Y}$ such that there are distinct branches $b^*=(\Lambda_0)_{\vec{\W}^*,\mathcal{Y}} \neq c^* = (\Lambda_1)_{\vec{\W}^*,\mathcal{Y}}$. Note that 
\begin{center}
$\pi^{\U,b}_b\circ\pi^{\vec{\T}}\rest\Q_Y=\pi^{\U}_c\circ\pi^{\vec{\T}}\rest\Q_Y$.
\end{center}
Note further that there are $\tau_{b^*}:\M^{\vec{\U}^*,b}_{b^*} \rightarrow \P$ and $\tau_{c^*}:\M^{\vec{\U}^*,b}_{c^*} \rightarrow \P$ such that 
\begin{equation}\label{eqn:four}
\pi_Y\rest \Q_Y = \tau_{b^*}\circ \pi^{\vec{\U}^*}_{b^*}\circ \pi^{\vec{\W}^*,b}\circ \pi^{\U,b}_b\circ \pi^{\vec{\T}}\rest \Q_Y,
\end{equation}
and
\begin{equation}\label{eqn:five}
\pi_Y\rest \P_Y = \sigma\circ i = \tau_{c^*}\circ \pi^{\vec{\U}^*}_{c^*}\circ \pi^{\vec{\W}^*,b}\circ \pi^{\U,b}_c\circ \pi^{\vec{\T}}\rest \P_Y.
\end{equation}
This is because
\begin{center}
$\pi_{b^*}^{\vec{\U}^*}\circ \pi^{\vec{\W}^*,b}=\pi_{c^*}^{\vec{\U}^*}\circ \pi^{\vec{\W}^*,b}$.
\end{center}
Equations (\ref{eqn:four}), (\ref{eqn:five}) give us 
\begin{center}
$\tau_{b^*}\circ \pi^{\vec{\U}^*}_{b^*}\circ \pi^{\vec{\W}^*,b}\rest (\M^+(\U))^b = \tau_{c^*}\circ \pi^{\vec{\U}^*}_{c^*}\circ \pi^{\vec{\W}^*,b}\rest  (\M^+(\U))^b$.
\end{center}
Lemma \ref{lem:uniqueness} then implies that $b^* = c^*$. This is a contradiction.

So (\ref{eqn:lower_level_disagreement}) holds. By our assumption, $\Q(c,\U)\unlhd \M^\U_c$ and is a $X$-validated $\Lambda_1$-mouse and $\Q(b,\U)\unlhd \M^\U_b$ and is a $X$-validated $\Lambda_0$-mouse. Results of Chapter \ref{lsa internal theory chapter} and earlier sections of this chapter imply that $\Q(b,\U)=\Q(c,\U)$ and hence $b = c$. Contradiction.

The argument above shows branch condensation. The other clause of strong branch condensation follows from a very similar argument, so we leave it to the reader.
\end{proof}

\begin{lemma}\label{lem:fullness_preserving}
$\Sigma^Y_\xi$ is locally strongly $\Omega$-fullness preserving.
\end{lemma}
\begin{proof}
$\Omega$-fullness preservation follows from the previous lemma. We now prove the other clause of locally strongly $\Omega$-fullness preservation (see Figure \ref{fig:strong_fullness_preserving}). Let $\N=\M^Y_\xi$ and $\Psi = \Sigma^Y_\xi$. Suppose $(\vec{\T},\S)\in I(\N,\Psi)$ (so $\pi^{\vec{\T},b}$ exists). Suppose $\S^b \lhd \W \unlhd \S$ is such that for some $n$, $\W$ is $n$-sound and, 
\begin{center}
$o(\S^b)\leq \omega\rho^{n+1}_\W <  \omega\rho^n_\W$.
\end{center}
Suppose $\tau:\R\rightarrow \W$ is cardinal preserving, is $\Sigma_0^{(n)}$, and $\omega\rho^n_{\R} > $ cr$(\tau) \geq \omega\rho^{n+1}_\R=\omega\rho^{n+1}_\W$. We want to show the $\tau$-pullback of the strategy $\Sigma_{\vec{\T},\W}$ is $\Omega$-fullness preserving.

Note that $\tau\rest \R^b = \rm{id}$ and $\R^b = \W^b$. This implies rng$(\pi^{\vec{\T},b})\subseteq \rm{rng}(\tau)$. Let $\sigma:\W^b= \S^b \rightarrow \P$ be the $\pi_Y$-realization map, so that $\pi_Y\rest\N^b = \sigma\circ \pi^{\vec{\T},b}$. Note that rng$(\sigma)$ is an honest extension of $X$. 

We now show $\Sigma^\tau_{\vec{\T},\W}$ is $\Omega$-fullness preserving. To see this, let $(\W^*,\vec{\U})\in I(\W,\Sigma^\tau_{\vec{\T},\W})$ be such that $\pi^{\vec{\U},b}: \R^b\rightarrow (\R^*)^b$ exists and let $\tau\vec{\U}$ be the copy tree on $\W$ with last model $\W^*$. So $\pi^{\tau\vec{\U},b}:\W^b\rightarrow (\W^*)^b$ exists. Let $\psi: (\R^*)^b\rightarrow (\W^*)^b$ be the copy map and $\sigma^*:(\W^*)^b\rightarrow \P$ be given by the construction of $\Psi$, so $\sigma^*\circ \pi^{\tau\vec{\U},b} = \sigma$ and $\sigma^*\circ \pi^{\tau\vec{U},b} \circ \sigma = \pi_Y\rest \N^b$.

Note that $\psi = \rm{id}$ and rng$(\sigma^*)$ is an honest extension of $X$. So $(\W^*)^b$ is $\Omega$-full. This is our desired conclusion.
\end{proof}

\begin{figure}
\centering
\resizebox{0.4\textwidth}{!}{
\begin{tikzpicture}[node distance=3cm, auto]
 \node (A) {$\N^{b}$};
\node (B) [right of=A] {$\R^b$};
  \node (C) [above of=B] {$\W^b$};
 \node (D) [right of=B] {$(\R^*)^b$};
 \node (E) [right of=C]{$(\W^*)^b$};
 \node (F)[above of=E]{$\P$};

  \draw[->] (A) to node [swap] {$\pi^{\vec{\T},b}$} (C);
  \draw[->] (B) to node {$\tau$} (C);
  \draw[->] (B) to node {$\pi^{\vec{\U},b}$} (D);
  \draw[->] (C) to node [swap]{$\pi^{\tau\vec{\U},b}$} (E);
  \draw[->] (E) to node {$\sigma^*$} (F);
  \draw[->] (D) to node [swap]{$\psi$} (E);
  \draw[->] (C) to node [swap]{$\sigma$} (F);
  \draw[->, bend left=45] (A) to node [swap]{$\pi_Y$} (F);
 \end{tikzpicture}
 }
 \caption{Strong $\Omega$-fullness preservation}
 \label{fig:strong_fullness_preserving}
\end{figure}

An easy corollary of the above Lemmata is the following.
\begin{corollary}\label{lem:pullback_equal}
Suppose $Y\prec Z\prec \M_\xi$ are countable (in $V[G]$), and such that $Y\cap \P, Z\cap \P$ are honest extensions of $X$, $Y,Z\in V$, and $Y=Y^* \cap \M_\xi$, $Z = Z^*\cap \M_\xi$ for some $Y^*\prec Z^*\prec H_{\kappa^{+4}}^V$. Let $\pi_{Y,Z} = \pi_Z^{-1}\circ \pi_Y$. Then $\Sigma^{Y}_\xi = (\Sigma^{Z}_\xi)^{\pi_{Y,Z}}$.
\end{corollary}
\begin{proof}
Let $\delta_Y=\pi_Y^{-1}(\delta^\P)$ and $\delta_Z=\pi_Z^{-1}(\delta^\P)$. By our assumption on $Y$ and $Z$, we have:
\begin{center}
$\pi_Z\rest \delta_Z = \pi^{\Sigma^Z_\xi}_{\N^Z_\xi,\infty}\rest \delta_Z$,
\end{center}
and
\begin{center}
$\pi_Y\rest \delta_Y = \pi^{\Sigma^Y_\xi}_{\N^Y_\xi,\infty}\rest \delta_Y = \pi_Z\circ \pi_{Y,Z}\rest \delta_Y$.
\end{center}

Using the above equations and the proof of Lemma \ref{lem:uniqueness} (especially the idea that if two strategies disagree, then there is a lower-level disagreement), we obtain the desired conclusion.

\end{proof}

\begin{corollary}\label{cor:commuting_positional}
$\Sigma^Y_\xi$ is positional and commuting. 
\end{corollary}
\begin{proof}
This follows from Lemmata \ref{lem:locally_strong_branch_condensation}, \ref{lem:fullness_preserving}, and results of Section 4.7.
\end{proof}



We have verified $(1)_\xi(c)$ holds, assuming $(1)_\xi(b)$. Let $Y$ be as above, i.e. $Y\prec H^V_{\kappa^{+4}}$ is such that $|Y| < \kappa$, $Y\cap \P$ is an honest extension of $X$. We discuss how to lift $\Psi=\Sigma^Y_\xi$ to a (necessarily unique) $(\kappa^{+},\kappa^{+})$-strategy $\Psi^+$ with branch condensation and show $Code(\Psi)\in \Omega$.

Recall $\Psi$ is an $(\omega_1,\omega_1)$-strategy for $\M^Y_{\xi}$ with branch condensation, is positional and $\Omega$-fullness preserving. Furthermore, $\Psi\cap V\in V$ and is independent of the choice of generic $G$. $\Psi\cap V$ can be uniquely extended to an $(\kappa^{+},\kappa^{+})$ strategy with branch condensation and is positional in $V$. We also call this extension $\Psi$. We outline how this extension works. We define $\Psi(\T)$ for $\T$, a normal tree of length $<\kappa^+$. Suppose cof$(lh(\T))\geq \omega_2$, then letting $\xi = \textrm{cof}(lh(\T))$, we can construe $\vec{C} = ([0,\alpha]_T : \alpha < lh(\T) \wedge \alpha \textrm{ is a limit ordinal})$ as a coherent sequence. Applying $\neg \square(cof(\xi))$ to $\vec{C}$, we get a club $D\subset lh(\T)$ that threads the sequence $\vec{C}$. $D$ gives a cofinal branch $b$ through $\T$.  This branch is necessarily the unique well-founded branch of $\T$. We define $\Psi(\T) = b$. Suppose cof$(lh(\T))< \omega_2$, the arguments in   \cite{PFA} or \cite[Lemma 3.62]{Trang2015PFA} show that there is $W\prec H_{\kappa^{+4}}$ with the properties: 
\begin{enumerate}[(a)]
\item $|W| < \kappa$, $W\cap \kappa\in \kappa$;
\item $W^{<|W|} \subset W$;
\item $\{\M^Y_\xi, \P, \T, \Psi\}\in W$;
\item for any $W\prec W_0 \prec W_1$ with properties (a)-(c), letting $\pi_{W_i}:W_i\rightarrow H_{\kappa^{+4}}$, $\T_{W_i} = \pi_{W_i}^{-1}(\T)$, and $b_i = \Psi(\T_{W_i})$ for $i\in \{0,1\}$, then $\pi_{W_0}[b_0]\downarrow \subseteq \pi_{W_1}[b_1] \downarrow$, where $\pi_{W_i}[b_i]\downarrow$ is the downward closure of $\pi_{W_i}[b_i]$ in $\T$.
\end{enumerate}
The $W$ as above is called $\T$-stable and we define $\Psi(\T) = b$ where $b$ is the downward closure of $\pi_W[\Psi(\T_W)]$ in $\T$. It is clear that the definition of $\Psi(\T)$ does not depend on the choice of $\T$-stable $W$.

We briefly give a sketch as to how to obtain a $(\kappa^{+},\kappa^{+})$-strategy $\Psi^+$ extending $\Psi$ with branch condensation and is positional in $V[G]$. In $V[G]$, suppose $\mathcal{T}$ is of limit length $< \kappa^{+}$ and is according to $\Psi^+$. We show how to define $\Psi^+(\mathcal{T})$ (stacks of normal trees can be handled similarly). In $V$, let $A\subseteq \kappa^{}$ code $H_{\kappa^{}}$ and a (nice) $Col(\omega,<\kappa)$-name $\dot{\T}\in H_{\kappa^{+}}$ for $\T$. Let 
\begin{center}
$M_A = L_{\kappa^{+}}^{\Lambda}[A,\M_2^{\Psi,\sharp}]$
\end{center}
where $\Lambda$ is the unique $(\kappa^{+},\kappa^{+})$-strategy for $\mathfrak{M}=_{\rm{def}}\M_2^{\Psi,\sharp}$, the minimal $E$-active $\Psi$-mouse with two Woodin cardinals. We note that the existence of $\M_2^{\Psi,\sharp}$ follows from \cite[Section 3.2]{Trang2015PFA}. By $\neg \square(\kappa^{+})$, 
\begin{center}
$M_A \vDash$ there are no largest cardinals.
\end{center}
In particular $(\kappa^+)^{M_A}<\kappa^{+}$, so in $M_A$, which is closed under $\Lambda$, we can use $\Lambda$ to perform a generic genericity iteration to make $A$-generically generic (see \cite{ATHM} or \cite{trang2013} for more on generic genericity iterations). Let $\Q\in M_A$ be the result of such an iteration. There is a $\Q$-generic $h\subseteq Col(\omega,\delta_0^\Q)$ such that $H_{\kappa},G,\dot{\T}\in \Q[h]$, where $\delta_0^\Q$ is the first Woodin cardinal of $\Q$. Since $\Q$ is closed under $\Psi$; we can generically interpret $\Psi$ on any generic extensions of $\Q$ (as done in \cite{trang2013} or in Chapter \ref{lsa internal theory chapter}).\footnote{If $\Psi$ is a strategy, we could have simply let $\mathfrak{M}=\M_1^{\Psi,\sharp}$; but if $\Psi$ is a short-tree strategy, then one seems to need $\M_2^{\Psi,\sharp}$ to apply results in Chapter \ref{lsa internal theory chapter}. Relevant results in \cite{trang2013} can be applied to $\M_2^{\Psi,\sharp}$ as well.} This allows us to define $\Psi^+(\T)$ as the branch chosen by the interpretation of $\Psi$ applied to $\T$ in $\Q[h]$. The well-definition and uniqueness of $\Psi^+$ follow from hull arguments in \cite[Section 3.2]{Trang2015PFA}.\footnote{Let $M, M^*$ be such that $\T\in M\cap M^*$; let $\tau,\tau^*$ be nice Col$(\omega,<\kappa)$-terms for $M, M^*$ respectively. In $V[G]$, let $W[G]$ contain all relevant objects and $W\prec H_{\kappa^{+4}}$ is good. Let $\bar{a}=\pi_W^{-1}(a)$ for all $a\in W[G]$. Then letting $b_0,b_1$ be the branches of $\bar{\mathcal{U}}$ given by applying \cite[Lemma 4.8]{trang2013} in $L^{\Lambda}[tr.cl.(\bar{\tau}),<_1,\mathfrak{M}], L^{\Lambda}[tr.cl.(\bar{\tau^*}),<_2,\mathfrak{M}]$ (built inside $M_W[G]$), where $<_1$ is a well-ordering of $\bar{\tau}$ and $<_2$ is a well-ordering of $\bar{\tau^*}$. Then $b_0 = b_1$ as both are according to $\Psi$, since $(\mathfrak{M},\Lambda)$ generically interprets $\Psi$ in $V[G]$.}

Using $\Psi^+$ and assuming $\Psi$ is a strategy, we can define the stack of $\Theta$-g-organized mice over $\mathbb{R}$, Lp$^{^\gTheta\Psi^+}(\mathbb{R}, Code(\Psi))$, in $V[G]$ (cf. \cite[Definition 4.23]{trang2013}),\footnote{\cite{Trang2015PFA} shows that $Code(\Psi)$ is self-scaled in the sense of \cite[Definition 4.22]{trang2013} if $\Psi$ is a strategy.} and show that there is a maximal initial segment $\M\unlhd \rm{Lp}$$^{^\gTheta\Psi^+}(\mathbb{R}, Code(\Psi))$ such that $\M$ is constructibly closed and $\M\vDash \sf{AD}^+ + \textsf{SMC}$ $+ \ \Theta = \theta_{\Psi}$. This implies $Code(\Psi)\in \Omega$.

\begin{remark}\label{rem:extension}
The arguments given above show that we can further extend $\Psi^+$ to a $(\kappa^{+4},\kappa^{+4})$-strategy.

If $\Psi$ is not total (so $(1)_\xi(b)$ fails) and that $\M_\xi$ is of $\sharp$-lsa-type, then we stop the hybrid $K^c$-construction and continue the $X$-validated sts construction above $\mathcal{C}(\M_\xi) = \N_\xi$. The idea is that we'll wait until we reach a level $\M'_\gamma$ (if exists) of the $X$-validated sts construction extending $\N_\xi$ such that some $\R\unlhd \M'_\gamma$ is a $\Q$-structure for $\delta^{\N_\xi}$ and then $\forall^* Y$, we can construction the canonical $X$-realizable strategy $(\Sigma^Y_\xi)$ of $\R^Y$ and show that it is in $\Omega$. If we reach a level $\M'_\gamma$ such that there is no $X$-validated strategy for $\M'_\gamma$ as witnessed by $p$, then we need to continue with an $X$-validated sts construction over $\M(p)^\sharp$.$\myqedhere$
\end{remark}

\begin{definition}[Certified-extender-ready levels]\label{def:extender_ready}\index{certified-extender-ready levels}
For $\xi < \Upsilon$, $\N_\xi$ is \textbf{certified-extender-ready}\index{certified-extender-ready levels} if for a $V$-club $\mathcal{C}_\xi$ of $Y\prec \N_\xi$ such that $Y\cap \P$ is an honest extension of $X$ and $Y\in V$ is countable in $V[G]$, letting $\N^Y_\xi = \pi_Y^{-1}(\N_\xi)$, $\Psi=\Sigma^Y_\xi$, and $\gamma^Y_\xi$ be the supremum of the indices of extenders on the $\N^Y_\xi$-sequence with critical point $\delta_Y=_{\rm{def}}$$ \pi_Y^{-1}(\delta^\P)$ (we let $\gamma^Y_\xi=((\delta_Y)^+)^{\N^Y_\xi}$ if $\N^Y_\xi$ has no such extenders on its sequence), we have that $\Psi$ is a strategy\footnote{This means $(1)_\xi(b)$ holds and hence $(\N^Y_\xi,\Psi)$ is not a sts hod pair.} and no $\M\unlhd$ Lp$^{\Psi,\Omega}(\N^Y_\xi)$ projects across $\gamma^Y_\xi$. For $Y\in \mathcal{C}_\xi$, we also say $\N^Y_\xi$ is \textbf{$\pi_Y$-certified extender-ready}. $\myqedhere$
\end{definition}
\begin{remark}\label{rem:ext_ready}
Extender-ready levels are those $\N_\xi$'s that are eligible to be extended to a hod premouse $(\N_\xi,F)$ where $F$ has critical point $\delta^\P$. Let $Y,\M$ be as in the above definition, it is easy to see that $\M$ also does not project across $\sf{ord}$$(\N^Y_\xi)$.

$\myqedhere$
\end{remark}

The lemma below shows that the collection of correctly-backgrounded extenders with critical point $\delta^\P$ is sufficiently rich. For instance, if $\P_Y=\pi_Y^{-1}(\P)$, and $\N_\xi^Y = \rm{Lp}$$^{\Sigma_{\P_Y},\Omega}(\P_Y)$, then $\N^Y_\xi$ is extender-ready (Corollary \ref{cor:solidity_universality} shows that no level of $\N^Y_\xi$ projects below $\sf{ord}$$(\P_Y)$ and Theorem \ref{thm:condensation_lemma} and Corollary \ref{cor:solidity_universality} show that every level of Lp$^{\Sigma_{\P_Y},\Omega}(\N^Y_\xi)$ is sound). Lemma \ref{lem:next_certified_extender} shows that if $\N_\xi$ is extender-ready then for every $Y\in \mathcal{C}_\xi$, there is a correctly backgrounded extender $E$ with critical point $\delta_Y$ such that $(\N^Y_\xi,E)$ is a hod premouse.

\begin{lemma}
 \label{lem:next_certified_extender}
Suppose $\N_{\xi^*}$ is extender-ready where $\xi = \xi^*+1$ and suppose $(1)_{\xi^*}-(3)_{\xi^*}$ hold. Fix $Y\prec \N_{\xi^*}$ in $\mathcal{C}_{\xi^*}$. Let $\N=\N^Y_{\xi^*}$, $ \delta^Y=\pi_Y^{-1}(\delta^\P)$, and $\Psi=\Sigma^Y_{\xi^*}$ be the $X$-realizable strategy for $\N$. Then there is an extender $E_Y$ with $\cp(E_Y) = \delta^Y$ such that $E_Y$ is $\pi_Y$-certified over $(\N,\Psi)$.
\end{lemma}
\begin{proof} Let $\gamma=\sf{ord}$$(\N)$. Let $E = E_Y$ be the following extender over $\N$: for $a\in[\gamma]^{<\omega}$ and $A\in \powerset(\delta^Y)^{|a|}\cap \N$,
\begin{center}
$(a,A)\in E \Leftrightarrow \pi^{\Psi}_{\N,\infty}(a)\in \pi_Y(A)$.
\end{center}

Fix a $Y\prec Z\in \mathcal{C}_{\xi^*}$ such that $Z=Z'\cap H_{\kappa^{+4}}^V$ and $Z'\prec H_{\kappa^{+4}}^V$ contains all relevant objects. Let $\pi_{Z'}:M_{Z'}\rightarrow Z'$ be the uncollapse map and $\iota = \cp(\pi_{Z'}) = Z'\cap \kappa$. Naturally, $\pi_{Z'}$ extends to act on all of $M_{Z'}[G\rest \iota]$ and induces an elementary embedding from $M_{Z'}[G\rest \iota]$ into $H_{\kappa^{+4}}[G]$; we also denote the extension map $\pi_{Z'}$. Let $\pi=\pi^\Psi_{\N,\infty}$ and $\pi'=(\pi^{\Psi}_{\N,\infty})^{M_{Z'}}$. By our assumption, $\Psi$ is $\Omega$-fullness preserving, commuting, and has branch condensation; furthermore, $\pi\rest \N|\delta^Y = \pi_Y\rest \N|\delta^Y$ and $\pi'\rest \N|\delta^Y = \pi_{Y,Z}\rest \N|\delta^Y$.

It is easy to see that $E$ is the extender $E'$ defined as follows: for $a\in \gamma^{<\omega}$ and $A\in \powerset(\delta^Y)^{|a|}\cap \N$,
\begin{center}
 $(a,A)\in E' \Leftrightarrow \pi'(a)\in \pi_{Y,Z}(A)$.
\end{center}

We need to see that $(\N,E)$ is a hod premouse.

\underline{Amenability}: Let $\eta<\gamma$ and $\xi < (\delta^{Y,+})^\N$, we show: $E\cap (\eta^{<\omega}\times \N|\xi)\in \N$.

Let $\mathcal{A} = (A_\alpha \ | \ \alpha<\delta^Y)$ enumerate $\N|\xi\cap \powerset(\delta^{Y,<\omega})$. Let
\begin{center}
 $B = \pi_{Y}(\mathcal{A})\cap (\pi(\eta)\times \pi(\eta))$.
\end{center}
Then $B\in \N_\xi|\delta^\P$ and so is $OD^{\Omega}$. Now for all $a\in \eta^{<\omega}$, for all $\alpha<\delta^Y$,
\begin{eqnarray*}
 (a,A_\alpha)\in E  & \Leftrightarrow & \pi(a)\in B_{\pi(\alpha)}.
\end{eqnarray*}
This shows $E\cap (\eta^{<\omega}\times \N|\xi)$ is $OD^{\Omega}_{\Psi}$. By $\textsf{SMC}$ and the fact that $\N$ is $\pi_Y$-certified extender-ready, $E\cap (\eta^{<\omega}\times \N|\xi)\in \N$.

\underline{Normality}: Let $c\in \gamma^{<\omega}$, $f:[\delta_Y]^{|c|}\rightarrow \delta_Y$ be such that $f\in \N^b$ and $\forall^*_{E_c} u \ f(u) < \textrm{max}(u)$ or equivalently $\pi_{Y}(f)(\pi(c)) < \textrm{max}(\pi(c))$. We want to find a $\xi < \textrm{max}(c)$ such that 
\begin{center}
$\pi_{Y}(f)(\pi(c)) = \pi(\xi) = c_{\pi(\xi)}(\pi(c))$, 
\end{center}
where $c_\xi$ is the constant function with range $\{\xi\}$.

Let  $\M$ be a $\Psi$-iterate of $\N$ such that $\pi^\Psi_{\N,\M} = \pi_{\N,\M}$ exists, $\tau_\M:\M\rightarrow \N_\xi$ be the $\pi_Y$-realization map given by the definition of $\Psi$, and let $\Psi_\M$ be the $\tau_\M$-pullback of $\Sigma_{\N_\xi^*}$. Let $E_\M$ be the extender that is $\tau_\M$-certified over $(\M,\Psi_\M)$, that is:
\begin{center}
$(a,A)\in E_\M \Leftrightarrow \pi^{\Psi_\M}_{\M,\infty}(a)\in \tau_\M(A)$.
\end{center}
It is easy to see, using Lemma \ref{lem:derived_model} that $\pi_{\N,\M}[E_\N]\subseteq E_\M$ and $\Psi_\M$ extends the tail strategy induced by $\Psi$ and $\pi_{\N,\M}$. 

We can find $\M$ such that $\pi_Y(f)(\pi(c))\in \rng(\tau_\M)$. Let $\M^* = \textrm{Ult}(\M,E_\M)$ and we note that $\Psi_{\M^*||lh(E_\M)}=\Psi_{\M||lh(E_\M)}$; call this strategy $\Lambda$. Note that $\pi(c) = \pi^\Lambda_{\M||lh(E_\M),\infty}\circ \pi^{\Psi}_{\N||lh(E),\M}(c)$. We have then that 
\begin{center}
$L(\Omega,\P) \models ``\pi_Y(f)(\pi(c))\in \rng(\pi^\Lambda_{\M||lh(E_\M),\infty})"$. 
\end{center}
By Lemma \ref{lem:derived_model}, 
\begin{center}
$L(\Omega,\P)\models ``\pi_Y(f)(c)\in \rng(\pi^\Psi_{N||lh(E),\infty})"$. 
\end{center}
This is what we want.

\underline{Coherence:} We now show:
\begin{enumerate}
\item $\textrm{Ult}_0(\N,E)|\gamma = \N$.
\item Let $\nu = \textrm{max}\{(\delta_Y^+)^\N,\gamma^Y_{\xi^*}\}$.\footnote{See Definition \ref{def:extender_ready}.} Then $\nu$ is a cutpoint of $\textrm{Ult}_0(\N,E)$ and $\gamma = ((\nu)^+)^{\textrm{Ult}_0(\N,E)}$.
\end{enumerate}

For 1), let $\tilde{\tau}:\textrm{Ult}_0(\N,E)\rightarrow \N_{\xi^*}$ be the natural map:
\begin{center}
$\tilde{\tau}(i_E(f)(a)) = \pi_{Y}(f)(\pi^\Psi_{\N,\infty}(a))$
\end{center}
for all $f\in \P_X$ and $a\in \gamma^{<\omega}$. It's clear from the fact that $\Psi$ is $X$-realizable and rge$(\pi_Y)$ is an honest extension of $X$ that $\tilde{\tau}\rest \gamma = \pi^\Psi_{\N,\infty}\rest \gamma$. This implies $\textrm{Ult}_0(\N,E)|\gamma$ is isomorphic to $\pi_{Y}[\N]$ and hence isomorphic to $\N$.
\begin{figure}
\centering
\resizebox{0.6\textwidth}{!}{
\begin{tikzpicture}[node distance=3cm, auto]
  \node (A) {$\N$};
  \node (B) [right of=A] {$\N^Z_{\xi^*}$};
  \draw[->] (A) to node {$\pi_{Y,Z}$} (B);  
  \node (C) [node distance=1.5cm, below of=B] {$\M$};
  \draw[->] (A) to node {$i$} (C); 
 \node (D) [node distance=0.8cm, right of=B] {$\tau(F)$};
  \node (E)[node distance=0.8cm, right of=C] {$F$};
 \draw[dashed,->] (E) to node  {} (D);
  \draw[->] (C) to node {$\tau$} (B);
  \node (F) [node distance=5cm, right of=B] {$\textrm{Ult}(\N^Z_{\xi^*},\tau(F))$};
  \draw[->] (D) to node {$u$} (F);
  \node (G) [node distance=5cm, right of=C] {$\textrm{Ult}(\M,F)$};
  \draw[->] (E) to node {$t$} (G);
  \draw[->] (G) to node {$k$} (F);
  \node (H) [node distance=1.5cm, right of=B] {};
  \node (I) [node distance=2cm, above of=H] {$\N_{\xi^*}$};
  \draw[->] (B) to node [swap]{$\pi_{Z}$} (I);
  \draw[->] (F) to node {$\sigma$} (I);
  \draw[->] (A) to node {$\pi_Y$} (I);
\end{tikzpicture}
}
\caption{Coherence}
\label{fig:proof_of_coherence}
\end{figure}

For 2), suppose not. Using the fact that $\N$ is extender-ready, we first observe that, 
\begin{equation}\label{eqn:nu_largest_cardinal}
\N\vDash \forall \nu\leq \alpha < \gamma \ (|\alpha|\leq \nu).
\end{equation}

Let $F$ be on the sequence of $\textrm{Ult}_0(\N,E)$ such that 
\begin{enumerate}[(i)]
\item crit$(F)=\delta_Y$.
\item lh$(F)\geq \nu$.
\item lh$(F)$ is the least such that (i) and (ii) hold.
\end{enumerate}
We have then that lh$(F)\geq \gamma$ by the definition of $\nu$ and the fact that $\textrm{Ult}_0(\N,E)|\gamma = \N$.\footnote{If $\nu$ is not a cutpoint of $\textrm{Ult}_0(\N,E)$, then there is some extender $H$ on the sequence of $\textrm{Ult}_0(\N,E)$ such that cr$(H)\leq \nu < \rm{lh}$$(H)$. This easily implies that there is some extender $F$ on the sequence of $\textrm{Ult}_0(\N,E)$ such that $\cp(F)=\delta_Y$ and lh$(E)\geq \nu$.}

Let $\tilde{\tau}$ and $Z$ be defined as above. Recall $Y\prec Z\in\mathcal{C}_\xi$.  Let $\M = \textrm{Ult}_0(\N,E)$, $i$ be the corresponding ultrapower map. Let $\tau: \M\rightarrow \N^Z_{\xi^*}$ be the natural map so that $\tilde{\tau}=\pi_Z\circ \tau$. Let $t:\M \rightarrow \textrm{Ult}(\M,F)$ be the ultrapower map by $F$ and $u:\N^Z_{\xi}\rightarrow \textrm{Ult}(\N^Z_{\xi},\tau(F))$ be the ultrapower map by $\tau(F)$. Let $k:\textrm{Ult}(\M,F)\rightarrow  \textrm{Ult}(\N^Z_{\xi},\tau(F))$ be the natural map and $\sigma:\textrm{Ult}(\N^Z_{\xi},\tau(F))\rightarrow \N_{\xi}$ be the realization map. The existence of $\sigma$ comes from the fact that $\tau(F)$ is $\pi_{Z}$-certified over $(\N^Z_{\xi^*}||\textrm{lh}(\tau(F)),(\Sigma^Z_{\xi^*})_{\N^Z_{\xi^*}||\textrm{lh}(\tau(F))})$.

\begin{claim} 
lh$(F) = \gamma$.
\end{claim}
\begin{proof}
Note that $\nu$ is a cutpoint in $\textrm{Ult}(\M,F)$ and is the least such $> \delta_Y$. So by (\ref{eqn:nu_largest_cardinal}),
\begin{center}
lh$(F) = (\nu^+)^{\textrm{Ult}(\M,F)}$
\end{center}
Suppose lh$(F) > \gamma$. Let $\Q \lhd \M||\textrm{lh}(F)$ be least such that
\begin{center}
$\N\lhd \Q \wedge \Q \vDash |\gamma| = \nu$.
\end{center}
Note that $\Q$ is a level of $Lp^{\Psi,\Omega}(\N)$. This is by $\sf{SMC}$ and the fact that 
\begin{center}$\pi_Z\circ\tau\rest\N = \tilde{\tau}\rest \N =  \pi^\Psi_{\N,\infty}$ . 
\end{center}
This contradicts the assumption that $(\N,\Psi)$ is extender-ready.
\end{proof}
Now we show $F$ is $\pi_{Y}$-certified over $(\N,\Psi)$. This would give $E = F\in \textrm{Ult}(\N,E)$. Contradiction.

Let $\Lambda=\Sigma^Z_{\xi^*}$. First note that
\begin{enumerate}[(a)]
\item $\pi^\Psi_{\N,\infty} = \pi^{\Lambda_{\tau(\N)}}_{\tau(\N),\infty} \circ \pi^\Psi_{\N,\N^Z_{\xi^*}}$.\footnote{$\N^Z_{\xi^*}$ is not literally a $\Psi$-iterate of $\N$, but $\N$ iterates into a hod initial segment of $\N^Z_{\xi^*}$. By $\pi^\Psi_{\N,\N^Z_{\xi^*}}$, we mean $(\pi^\Psi_{\N,\infty})^{M_{Z}}$.}
\item $\tau\rest\N = \pi^\Psi_{\N,\N^Z_{\xi^*}}$.
\item $\sigma\rest\tau(\N) = \pi^{\Lambda_{\tau(\N)}}_{\tau(\N),\infty}$.

\end{enumerate}
Let $c\in [\sf{ord}$$(\N)]^{<\omega}$, $A\in \P_Y$, we have:
\begin{eqnarray*}
A\in F_c & \Leftrightarrow & c \in t(A) \\&
\Leftrightarrow & c \in t(i(A)) \ \ \  \\ &
\Leftrightarrow & k(c) \in k(t(i(A)) \\ &
\Leftrightarrow & \tau(c) \in u\circ \tau (i(A)) \\ &
\Leftrightarrow & \tau(c) \in u(\pi_{Y,Z}(A)) \\ &
\Leftrightarrow & \sigma(\tau(c)) \in \pi_{Y}(A)  \\ &
\Leftrightarrow & \pi^{\Psi}_{\N,\infty}(c)\in \pi_{Y}(A).
\end{eqnarray*}
The second equivalence holds becase $i(A)\cap \delta_Y = A$. The third equivalence uses Corollary \ref{cor:strategy_condensation}, noting that $\rng(\sigma)$ is an honest extension of $\rng(\sigma\circ k)$. The fourth equivalence uses the fact that $\tau(c) = k(c)$ and $k\circ t = u\circ \tau$. The fifth equivalence is true because $\pi_{Y,Z}(A) = \tau(i(A))$. The sixth equivalence is true because $\pi_{Y}(A) = \sigma(u(\pi_{Y,Z}(A)))$. The last equivalence follows from equations (a)--(c). This finishes the proof of the lemma.
\end{proof}

Lemma \ref{lem:next_certified_extender} implies that if $\N_\xi$ is extender-ready then $\M_{\xi^*+1} = (\N_{\xi^*}, E)$ where using the notation of Lemma \ref{lem:next_certified_extender}
\begin{center}
$(a,A)\in E \Leftrightarrow \forall^* Y\in \mathcal{C}_\xi ((a,A)\in Y \rightarrow \pi_Y^{-1}(a,A)\in E_Y)$.
\end{center}

It also follows from Lemma \ref{lem:next_certified_extender} that $(1)_\xi(a)$ holds if $(1)_{\xi^*}-(3)_{\xi^*}$ hold. We continue by proving another condensation lemma for relevant extenders with critical point $\delta^\P$. This condensation property does not seem to follow from Theorem \ref{thm:condensation_lemma}.
\begin{lemma}\label{lem:condensation_relevant_extenders}
\begin{enumerate}[(a)]
\item Suppose $\M_\xi$ is of the form $(\M_\xi^-,F_\xi)$, where $\cp(F_\xi) = \delta^\P$. Suppose $\pi:\M=(\M^-,\tilde{F}) \rightarrow \M_\xi$ is $\Sigma_0$ and cofinal, or $\Sigma_2$, with $\cp(\pi)>\sf{ord}$$(\P)$ and suppose further that $\M^-\unlhd \M^-_\xi$. Furthermore, let $Y$ be a good hull that contains all relevant objects, let $\pi_Y:M_Y[G]\rightarrow H_{\kappa^{+4}}[G]$ be the uncollapse map, and let $\M^Y=\pi^{-1}_Y(\M)$. Let $\Psi$ be the $\pi_Y$-pullback strategy for $\M^Y$ and suppose that $\Psi_{\M^{Y,-}} = (\Sigma^Y_\xi)_{\M^{Y,-}} = (\Sigma^Y_\xi)_{\M^{Y,-}}^{\pi^Y}$. Then $\tilde{F}$ is on the sequence of $\M_\xi$ and $\Psi = (\Sigma^Y_\xi)_{\M^Y} = (\Sigma^Y_\xi)_{\M^Y}^{\pi^Y}$, where $\Sigma^Y_\xi$ is the strategy for $\N^Y_\xi$ defined above.

\item More generally, suppose $\M_\xi$ is as above and $\pi: \M=(\M^-,\tilde{F}) \rightarrow \M_\xi$ is $\Sigma_0$ and cofinal, or $\Sigma_2$, with $\cp(\pi)>\sf{ord}$$(\P)$. Suppose $Y,\pi_Y, \Psi$ are as above and $\tilde{F}^Y = \pi_Y^{-1}(\tilde{F})$, then $\tilde{F}^Y$ is $\pi_Y$-certified over $(\M^{Y,-},\Psi_{\M^{Y,-}})$.
\end{enumerate}

\end{lemma}
\begin{proof}
The preservation of $\pi$ guarantees that $\M$ is a hod premouse. Recall that $\sf{ord}$$(\P)$ is the cardinal successor of $\delta^\P$ in both $\M_\xi$ and $\M$ and the models agree up to $\P$. 

We first prove (a). Let $Y$ be as in the hypothesis. Let $\tilde{F}^Y = \pi_Y^{-1}(\tilde{F})$, $(\P^Y,\delta^Y, F_\xi^Y)=\pi_Y^{-1}((\P,\delta^\P, F_\xi))$, and $\pi^Y = \pi_Y^{-1}(\pi)$. We work with $\M^Y$ and $\N^Y_\xi$ and first show that $\tilde{F}^Y$ is on the sequence of $\N^Y_\xi$. Let $\Lambda = \Psi_{\M^{Y,-}} = (\Sigma^Y_\xi)_{\M^{Y,-}}$. 
\begin{claim}\label{claim:certified}
For $a\in [\sf{ord}$$(\M^Y)]^{<\omega}$ and $A\subset [\delta^Y]^{|a|}$ in $\P^Y$, $(a,A)\in \tilde{F}^Y$ if and only if $\pi^{\Lambda}_{\M^{Y,-},\infty}(a) \in \pi_Y(A)$.
\end{claim}
\begin{proof}
First, note that $\tilde{F}^Y$ is total over $\N^Y_\xi$ and hence it makes sense to apply $\tilde{F}^Y$ to $\N^Y_\xi$. Also, Ult$(\N_\xi^Y,\tilde{F}^Y)$ embeds into Ult$(\N^Y_\xi,F_\xi^Y)$ via the natural map $\tau$:
\begin{center}
$\tau(i_{\tilde{F}^Y}(f)(b)) = i_{F^Y_\xi}(f)(\pi^Y(b))$.
\end{center}
Note that
\begin{center}
$\tau\rest \M^Y||lh(\tilde{F}^Y) = \pi^Y\rest \M^Y||lh(\tilde{F}^Y)$. 
\end{center}
Now,
\begin{eqnarray*}
(a,A)\in \tilde{F}^Y & \Leftrightarrow & (\tau(a)=\pi^Y(a), A)\in F^Y_\xi  \\ & \Leftrightarrow & \pi^{\Lambda}_{\M^{Y,-},\infty}(\pi^Y(a))\in \pi_Y(A) 
\\ &
\Leftrightarrow & \pi^{\Lambda}_{\M^{Y,-},\infty}(a)\in \pi_Y(A). 
\end{eqnarray*}
The first equivalence holds because $\pi^Y(A)= \tau(A)=A$. The second equivalence holds by the definition of $F^Y_\xi$ and our assumption on $\Lambda$. The last equivalence follows from Lemma \ref{lem:derived_model}. This finishes the proof of the claim. 
\end{proof}

The claim and Lemma \ref{lem:next_certified_extender} imply that $\tilde{F}^Y$ is on the $\N_\xi^Y$-sequence. By elementarity, $\tilde{F}$ is on the $\N_\xi$-sequence. 

$\Psi = (\Sigma^Y_\xi)_{\M^Y} = (\Sigma^Y_\xi)_{\M^Y}^{\pi^Y}$ follows from Lemma \ref{lem:uniqueness} and the proof of Lemma \ref{lem:locally_strong_branch_condensation} (the main points are $\pi_Y\circ\pi^Y\rest \P^Y = \pi_Y\rest \P^Y$, and the fact that if the strategies disagree then we can find a lower-level disagreement just as in the proof of Lemma \ref{lem:locally_strong_branch_condensation}). This proves (a).

The proof of (b) is very similar, note that we have  the following equivalences
\begin{eqnarray*}
(a,A)\in \tilde{F}^Y & \Leftrightarrow & (\pi^Y(a), A)\in F^Y_\xi  \\ & \Leftrightarrow & \pi^{\Sigma^Y_\xi}_{\M_\xi^{Y,-},\infty}(\pi^Y(a))\in \pi_Y(A) 
\\ &
\Leftrightarrow & \pi^{\Psi}_{\M^{Y,-},\infty}(a)\in \pi_Y(A). 
\end{eqnarray*}
The last equivalence easily follows from Lemma \ref{lem:derived_model} and shows that $\tilde{F}^Y$ is $\pi_Y$-certified over $(\M^{Y,-},\Psi_{\M^{Y,-}})$.

\end{proof}

\begin{corollary}\label{cor:solidity_universality}
\begin{enumerate}
\item Let $\N=\M_\xi$ be the $\xi$-th model in the hybrid $K^c$-construction. Suppose $(1)_\xi$ hold. Then $\rho(\N) \geq \sf{ord}$$(\N^b)=\sf{ord}$$(\P)$ and $\N$ is $k(\N)+1$-solid and $k(\N)+1$-universal.

\item Suppose $\M_\xi$ is the $\xi$-th model in the $X$-validated sts construction over a weakly suitable $\R$ that is the result of a hybrid $K^c$-construction. Suppose $\rho(\M_\xi) < \delta$, where $\delta$ is the lsa Woodin cardinal of $\R$. Let $\gamma = \textrm{max}(\rho(\M_\xi),\delta^\P)$ and let $\N$ be the transitive collapse of Hull$^{\M_\xi}(\gamma\cup \{p(\M_\xi\}))$.\footnote{$\N$ is obtained by decoding the $\Sigma_1$-hull of the $k(\M_\xi)$-reduct with parameters in $\gamma\cup \{p(\M_\xi)\}$.} If there is an $X$-validated iteration strategy $\Psi$ of $\N$ (i.e. if $(1)_\xi(b)$ holds), then $\rho(\N) > \delta^\P$.  
\end{enumerate}
\end{corollary}
\begin{proof}
We prove (1) first. We prove $\rho(\N)\geq \sf{ord}$$(\N^b)$ and $\N$ is $n$-solid and $n$-universal, where $n = k(\N)+1$. Without loss of generality we assume $n=1$. The case $n>1$ is similar (one just has to work with the $n-1$-reduct). 

\begin{claim}\label{claim:projectum}
$\rho_1(\N)\geq \sf{ord}$$(\N^b)$.
\end{claim}
\begin{proof}
Suppose not. Let $Y\prec H^V_{\kappa^{+4}}$ be $X$-good such that $\N\in Y$. Let $\delta_Y=\pi_Y^{-1}(\delta^\P)$, $\N_Y = \pi^{-1}_Y(\N)$, and $\Psi=\Sigma^Y_\xi$. Let $\Q = \textrm{Ult}_0(\N_Y,\nu)$ where $\nu$ is the order $0$ total measure with critical point $\delta_Y$ and $i_\nu: \N_Y\rightarrow \Q$ be the canonical embedding. Let $q = i_\nu(p)$ where $p = p_1(\N_Y)$. Hence
\begin{enumerate}[(i)]
 \item $\N_Y^b$ is a cutpoint initial segment of $\Q$ and $\sf{ord}$$(\N_Y^b)$ is the cardinal successor of $\delta_Y$ in $\Q$.
\item We can regard $\Q$ as a hod premouse over $(\N^b,\Psi_{\N_Y^b})$ with strategy $\Sigma_\Q\in \Omega$ that is commuting and is $\Omega$-fullness preserving.\footnote{We can take $\Sigma_\Q$ be the $\Q$-tail of $\Psi$. By Lemma \ref{lem:fullness_preserving}, $\Sigma_\Q$ is $\Omega$-fullness preserving. By Corollary \ref{cor:commuting_positional} and results of Section 4.7, $\Sigma_\Q$ is positional and commuting.}
\item There is some $A\subseteq \delta_Y$ such that $A$ is $\Sigma_1$-definable over $\Q$ from $q$ and $A\notin \N_Y^b$.
\end{enumerate}
We say that a triple $(\Q,\Sigma_\Q,q)$ satisfying (i)-(iii) is \textit{minimal} if there is no iteration $\vec{\T}$ according to $\Sigma_\Q$ with iteration map $i:\Q\rightarrow \R$ and some $r < i(q)$ (in the reverse lexicographic order) such that $(\R,\Sigma_{\R,\vec{\T}}, i(\N^b),r)$ satisfies (i)-(iii).

Fix two minimal triples $(\R,\Sigma_{\R},r)$ and $(\S,\Sigma_\S,s)$. We can then compare them above $\N_Y^b$. Letting $i:\R\rightarrow \W$ and $j:\S\rightarrow \W$ be iteration maps. Note that $i(r)=j(s)$ and so 
\begin{center}
 $\textrm{Th}_{\Sigma_1}^{\R}(\delta_Y\cup \{r\})=\textrm{Th}_{\Sigma_1}^{\S}(\delta_Y\cup \{s\})$.
\end{center}
This means $\textrm{Th}_{\Sigma_1}^{\R}(\delta_Y\cup \{r\})$ is $OD^\Omega_{\N_Y^b,\Psi_{\N_Y^b}}$ for any minimal $(\R,\Sigma_\R,r)$. By $\textsf{MC}(\Psi_{\N_Y^b})$, 
\begin{center}
$\textrm{Th}_{\Sigma_1}^{\R}(\delta_Y\cup \{r\})\in \N_Y^b$.\footnote{Note that we take $Y$ so that $\N_Y^b = \rm{Lp}$$^{\Psi_{\N_Y|\delta_Y},\Omega}(\N_Y|\delta_Y)$. } 
\end{center}
This contradicts (iii).

\end{proof}
The claim and Theorem \ref{thm:condensation_lemma} (which is built on the results of Section 4.9) imply that $\N_Y$ and hence $\N$ is $1$-solid and $1$-universal. The point is that relevant phalanx comparisons of the form $(H, \N_Y,\alpha)$ where $H$ is the transitive collapse of Hull$_1^{\N_Y}(\alpha\cup \{p-(\alpha+1)\})$ for $\alpha \in p$ are such that  $\alpha > \sf{ord}$$((\N_Y)^b)$ (by the claim), no strategy disagreements can occur (see Lemma \ref{lem:no_strat_disagreement}) and do not use extenders with critical point $\delta_Y$ (by Lemma \ref{lem:condensation_relevant_extenders}). The last item holds because Lemma \ref{lem:condensation_relevant_extenders} shows that extenders with critical point $\delta_Y$ on the $H$-sequence are certified and the proof of \ref{lem:next_certified_extender} shows that extenders with critical point $\delta_Y$ and its images on the sequence of iterates of $H, \N_Y$ are certified.\footnote{See the proof of Claim \ref{claim:soundness} for a similar argument.} These comparisons terminate successfully by the usual arguments. By similar arguments, we get the conclusion for all $n\in\omega$.


For (2), letting $Y$ be $X$-good and using the notations as above, we describe the $X$-realizable iteration strategy $\Psi_Y$ for $\N_Y$ witnessing $\N_Y$ is an $X$-approved sts mouse; then the $X$-validated strategy $\Psi$ for $\N$ is defined from the $\Psi_Y$'s as before.  First, note that by \cite{ANS}, for such a $Y$, $\M_\xi^Y$ has a $\tau$-realizable strategy above $\gamma$ for some map $\tau: \M^Y_\xi\rightarrow \M_\xi$.\footnote{$\tau$ is a minimal map relative to some enumeration $\vec{e}$ of $\M^Y_\xi$ in order type $\omega$.} The usual proof of solidity/universality then shows that $\N$ is $\gamma$-sound and if $\gamma = \rho(\M_\xi)$, then $\N$ is sound. In the case $\N$ is sound, then it is just $\N_\xi$.

Let $(\S,\delta')$ be the image of $(\R,\delta)$ under the collapse map $\pi^{-1}$ and let $(\S_Y,\delta'_Y)$ be the image of $(\S,\delta')$ under $\pi_Y^{-1}$. Note that $\rho(\N)<\delta'$. Now we outline the description of $\Psi_Y$. $\Psi_Y$ on stacks based on $\S_Y^b$ has been defined in great detail before (using the fact that $X$ is a condensing set),\footnote{Extenders on these stacks have critical point $\leq \delta^{\S^b_Y}$ and their images. Note also that $\R^b = \S^b$.} so we focus on stacks $\vec{\T}$ above $\S_Y^b$. Suppose $\T$ is above $\S^b_Y$, based on $\S_Y$, and is correctly guided. Then by \cite{ANS}, there is a maximal branch $b$ and a realizing map $\sigma: \M^{\T}_b\rightarrow \M_\xi$ such that $\sigma\circ i^{\T}_b = \pi\circ \pi_Y\rest \N_Y$. But note that there is a $\Q\unlhd \M^\T_b$ such that $\Q = \Q(\T,b)$. This comes from the fact that $\rho(\N)<\delta'$ and that $\N$ is $\gamma$-sound; so since $\gamma < \delta'$, $\J_1[\N]\models$ ``there are no Woodin cardinals $>\gamma$", and hence $\Q$ exists. Therefore, the branch $b$ is the canonical $\Q$-structure guided branch for $\T$. The case where $\T$ is above $\S$ is similar. We can then easily define $\Psi_Y$ on arbitrary stacks on $\N_Y$. An argument similar to the proof of Claim \ref{claim:projectum} then shows that $\rho(\N) > \delta^\P$.

\end{proof}
Corollary \ref{cor:solidity_universality} verifies $(2)_\xi, (3)_\xi$ hold, given that $(1)_\xi$ holds. 

Now suppose for some $\xi$, $\M_\xi$ and $\N$ are as in Corollary \ref{cor:solidity_universality}(2). By Corollary \ref{cor:solidity_universality}(2), $\N = \N_\xi$. Then as in \ref{cor:solidity_universality}, $\forall^* Y$, $Y$ is $X$-good,  $\N^{Y}_\xi$ is iterable via the $X$-realizable strategy. This induces an $X$-validated strategy $\Lambda$ for $\N_\xi$\footnote{We assumed this strategy  exists.}; so $\S =\N_\xi$ is $K^c$-appropriate. We can then start a hybrid $K^c$-construction over $\S$, producing models $(\M'_\xi,\N'_\xi: \xi \leq \Upsilon')$, maintaining $(1)-(3)$ along the way. At some $\xi\leq \Upsilon'$, if $(1)_\xi(b)$ fails then this implies that there is a $\sf{nuvs}$ $p$ witnessing such a failure; so $(\vec{\V} =_{def}  \{\M'_\alpha, \N'_\alpha: \alpha < \xi\}\cup\{\M'_\xi\},p)$ witnesses that $\textrm{m}^+(p)$ is honest weakly $X$-suitable. Corollary \ref{cor:validated_suitable} then shows that $\textrm{m}^+(p)$ is suitable. At this point, we will continue with the $X$-validated sts construction over $\textrm{m}^+(p)$ until it stops prematurely or it produces $\M_{\kappa^{+++}}$.

Suppose $\M_\xi$ is the $\xi$-th model in an $X$-validated sts construction over some $\R$ and does not define a failure of Woodinness of $\delta$, the lsa Woodin of $\R$, then we continue with the $X$-validated sts construction over $\R = \textrm{m}^+(\M_\xi|\delta)$. There are several cases. The first case occurs when we reach an $X$-validated model $\M_\epsilon$ for some $\epsilon > \xi$ such that $\epsilon$ is the least such that $\rho(\M_\epsilon) < \delta$ and $\rho(\M_\epsilon)\geq \sf{ord}$$(\P)$. So letting $\N_\epsilon = \mathcal{C}(\M_\epsilon)$, then as in Section \ref{sec:sts_constructions} and Corollary \ref{cor:solidity_universality}, we can show $\N_\epsilon$ is $X$-validated and if it has an $X$-validated iteration strategy, then $\N_\epsilon$ is sound and by the assumption on $\M_\epsilon$, we have that $\J_1[\N_\epsilon]\models$ ``there are no Woodin cardinals $> \delta^\P$". In this case, letting $\S = \N_\epsilon$, then $\S$ is $K^c$-appropriate as before and we continue with our hybrid $K^c$-construction over $\S$ until we reach a level $\M'_\xi$ with no $X$-validated iteration strategy, so letting $p$ witness this, we then continue with the sts $X$-validated construction over $\textrm{m}^+(p)$, which is suitable, as above. If $\M_\epsilon$ does not have an $X$-validated strategy, then letting $p$ witness this, we continue with the $X$-validated sts construction above $\textrm{m}^+(p)$ as before. The cases where the construction reaches $\M_{\kappa^{+++}}$ or stops prematurely, are handled in the next section.



\section{$K^c$ stops prematurely and a model of $\sf{LSA}$}\label{sec:main}


Suppose the construction lasts $\kappa^{+++}$ steps; as in the previous subsection, let $\N = \N_{\kappa^{+++}}$ and $\S = \S(\N)$.
\begin{lemma}\label{lem:small_cof}
cof$(\S) < \kappa^{+++}$.
\end{lemma}
\begin{proof}
Let $\lambda = \kappa^{+++}$. Note that $\S\in V$. If $\N$ is lsa-small,\footnote{In this case, $\N$ is the result of a hybrid $K^c$-construction over some $K^c$-appropriate $\R$ or the result of alternating the hybrid $K^c$-constructions and the $X$-validated sts construction in a manner described in the previous section.} then as shown in Chapter 11, $\S \vDash \square_{\lambda,2}$. Suppose $\N$ is the result of the $X$-validated sts construction over some suitable $\R$, then $\N$ is not lsa-small but it is an sts mouse over $\R$ and so is $\S$. In this case, the standard proof (\cite{schimmerling2004characterization}) shows $\S\models \square_{\lambda,2}$.\footnote{In fact $\S\models \square_{\lambda}$. This is because $\sf{ord}$$(\R)$ is a strong cut point of $\S$ and all relevant comparisons are above $\R$ and in fact extender comparisons.} Working in $V$, $\neg\square(3,\kappa^{+4})$ implies then that o$(\S) < \kappa^{+4}$ and $\neg \square(3,\lambda)$ now implies that cof$(\sf{ord}$$(\S)) < \lambda$ since otherwise, the canonical $\square_{\lambda,2}$-sequence $\vec{C}$ of $\S$ (as defined in Chapter 11) has a thread $D$. The thread $D$ will produce a sound hod mouse (or sts mouse) $\M$ such that $\sf{ord}$$(\M) \geq \sf{ord}$$(\S)$ and $\rho_{\omega}(\M)\leq \lambda$. This contradicts (ii) of Lemma \ref{lem:stack_facts}. 
\end{proof}

Lemma \ref{lem:small_cof} contradicts (iii) of Lemma \ref{lem:stack_facts}. Now we assume the construction stops prematurely. Without loss of generality, we assume the $X$-validated sts construction over some $\Q$, where $\Q$ is produced by a hybrid $K^c$-construction,\footnote{The other case where $\Q = \textrm{m}^+(p)$ for some $p$ is similar.} reaches a model $\N_\Upsilon$ which is a sts hod premouse that satisfies:
\begin{enumerate}[(i)]
\item There is a unique Woodin cardinal $\delta_0>\delta^\P$ such that $\delta^\P$ is the least $<\delta_0$-strong.
\item There are $\omega.2$ many Woodin cardinals above $\delta_0$, say these Woodin cardinals are $(\delta_n : 1\leq n<\omega.2)$.
\item There is an extender $F$ with crt$(F)>\sup_n \delta_n$ such that $\N_\Upsilon = (\N_\Upsilon^-,F)$ for some extender $F$. In fact, $\N_\Upsilon = (\N_\Upsilon^-)^\sharp$.
\item $\N_\Upsilon$ is a sts hod premouse over $\Q =_{\rm{def}}$$(\N_\Upsilon|\delta_0)^\sharp$, $\Q$ is of lsa type with $\delta^\P$ is $<\delta_0$-strong in $\Q$.
\item $\rho_\omega(\N_\Upsilon)\geq \sf{ord}$$(\Q)$. 
\end{enumerate}



Let $\lambda = \rm{sup}$$_n \delta_n$ and for every $\Q\lhd \M \unlhd \N_\Upsilon$, let $\Sigma^{\M}$ be the internal sts strategy of $\Q$ as defined in $\M$. 

\begin{lemma}\label{lem:not_all_the_way}
Suppose the construction stops prematurely. Then $\Upsilon<\kappa^{+++}$.
\end{lemma}
\begin{proof}
If the construction stops prematurely, then $\N_{\Upsilon}$ is $E$-active. This clearly implies that $\Upsilon<\kappa^{+++}$ because if $\Upsilon=\kappa^{+++}$, then $\N_{\Upsilon}$ is the lim inf of $\N_{\alpha}$ for $\alpha<\Upsilon$ and hence is passive.
\end{proof}

Now suppose $\N_\Upsilon$ satisfies Definition \ref{omega woodins over lsa}, then the results of Section 8.2 show that the derived model of $\N_\Upsilon$ (at the supremum of its Woodin cardinals) satisfies $\sf{LSA}$. Suppose this is not the case. We would like to produce an active $\omega.2$ Woodin lsa mouse as in Definition \ref{omega woodins over lsa} from $\N_\Upsilon$. 
\begin{lemma}\label{lem:hull}
Then there is a countable substructure of some $\M^*\unlhd \N_\Upsilon$ satisfying Definition \ref{omega woodins over lsa}.
\end{lemma}
\begin{proof}
Recall that we have $\rho_\omega(\N_\Upsilon) \geq \sf{ord}$$(\Q)$. Let $W\prec \N_\Upsilon$ be such that $\P\cup \{\P\}\subset W$ and $W\cap \delta_0\in \delta_0$.\footnote{Such a $W$ can be found easily because $\rho_\omega(\N_\Upsilon)\geq \delta_0$ and $\delta_0$ is definably inaccessible over $
\N_\Upsilon$.}  Let $\pi^*: \M^*\rightarrow W$ be the uncollapse map. 


Let $\M$ be the transitive collapse of Hull$^{\M^*}(\P\cup p(\M^*))$ and $\pi': \M \rightarrow \M^*$ be the uncollapse map. Let $\pi = \pi^*\circ \pi'$. First, note that we have the following:
\begin{center}
$\rho_\omega(\N_\Upsilon) \geq \sf{ord}$$(\P)$ and $\rho_\omega(\M)\leq \sf{ord}$$(\P)$ (so in fact, $\rho_\omega(\M)=\sf{ord}$$(\P)$).
\end{center}
Now we claim that
\begin{claim}\label{claim:soundness}
\begin{enumerate}[(i)]
\item $\M^* \lhd \N_\Upsilon$.
\item If $Y$ is $X$-good such that $\{\M,\M^*,\N_\Upsilon\}\in Y$, letting $\pi_Y$ be the uncollapse map and $x^Y=\pi^{-1}(x)$ for $x$ in range $\pi_Y$ or $x=\M$, then $\M^Y$ is iterable via the $X$-realizable strategy.
\item Suppose $Y$ is as in (ii) and $\tau: \N\rightarrow \M^Y$ is either $\Sigma_0$ cofinal or $\Sigma_2$ elementary and cr$(\tau)>o(\P^Y)$, then the comparison $(\M^Y,\N, \cp(\tau))$ against $\M^Y$ does not use extenders with critical point $(\delta^{\P^Y})$. 
\item $\M$ is $\omega$-sound.
\end{enumerate}
\end{claim}
\begin{proof}
For (i), first note that letting $\xi = \textrm{crit}(\pi^*)$ and $\Q^* = (\pi^*)^{-1}(\Q)$, then $\xi = \delta_0^{\M^*}$ and $\pi^*(\xi) = \delta_0^{\N_\Upsilon}$. Now, note that $\Q^*\lhd \Q$. Now let $\R\lhd \Q$ be the largest such that $\R \models ``\xi$ is Woodin". Let $Y$ be $X$-good such that $\{\M^*, \Q^*, \R, \N_\Upsilon,\xi\}\in Y$. Let $((\M^*)^Y, \R^Y,(\Q^*)^Y, \xi^Y) = \pi_Y^{-1}(\M^*,\R,\Q^*,\xi)$. Note that $\xi$ is a strong cut point of both $\M^*, \R$. Since $(\M^*)^Y$ is $\xi^Y$-sound and projects to $\xi^Y$, by \cite{ANS}, $(\M^*)^Y$ has an iteration strategy $\Psi$ above $\xi^Y$ and $\Psi$ can be taken to be a $\tau$-realizable strategy for some $\tau:(\M^*)^Y\rightarrow \M^*$ such that $X\subset \textrm{rng}(\tau)$.\footnote{We fix an enumeration $\vec{e}$ of $(\M^*)^Y$, in type $\omega$ in $V[g]$. Let $\tau: (\M^*)^Y \rightarrow \M^* $ be a $\vec{e}$-minimal embedding $\tau^*$ with the property that $\tau^*\rest (\Q^*)^Y = \pi_Y\rest (\Q^*)^Y$. Such an embedding can be obtained as the left-most branch of the tree that builds approximations to embeddings from $(\M^*)^Y$ into $\M^*$ that agrees with $\pi_Y$ on $(\Q^*)^Y$. By results of \cite{ANS}, we get a $\tau$-realizable strategy for $(\M^*)^Y$ for stacks above $(\Q^*)^Y$.} We can then compare $(\M^*)^Y, \R^Y$ above $\xi^Y$.\footnote{By what has been shown, the comparison does not encounter strategy disagreements. The point is that since $\Psi$ is a $\tau$-realizable strategy and $X\subset \textrm{rng}(\tau)$, it witnesses that $\Psi$-iterates of $(\M^*)^Y$ are $X$-approved. A similar comment applies to the the canonical strategy of $\R^Y$.} Since both models are $\xi$-sound and $\R^Y$ is a $\Q$-structure for $\xi^Y$, $(\M^*)^Y\unlhd \R^Y$. This implies $\M^*\lhd \N_\Upsilon$.


(ii) follows from the fact that $\M^*$ is an initial segment of a model in a hybrid $K^c$-construction that produces $\Q$. Furthermore, since there is a $\Q$-structure $\R\lhd \Q$ for $\xi$ extending $\M^*$, the previous section indeed produces a unique $X$-realizable iteration strategy for $\M^Y$. This strategy is $\pi'\circ \pi_Y$-realizable since it is the $(\pi')^Y$-pullback of the unique $X$-realizable strategy $\Lambda$ of $(\M^*)^Y$, and by \cite{ANS} and the fact that $(\M^*)^Y\unlhd \R^Y$, $\Lambda$ is $\pi$-realizable. 


(iii) follows from an argument as in Corolloary \ref{cor:solidity_universality}; the main point is that by Lemmata \ref{lem:condensation_relevant_extenders} and \ref{lem:next_certified_extender}, letting $\R$ be $\M^Y$ or $\N$ or its iterate, extenders on the sequence of $\R$ with critical point $\delta^{\R^b}$ are certified. More precisely, let $\S$ be the tree on the phalanx side and $\T$ be the tree on the $\M^Y$-side that participate in the comparison. Suppose $\R = \M^\S_\alpha$ and $\R' = \M^\T_\alpha$ and that $\S\rest \alpha$, $\T\rest \alpha$ have not used extenders with critical point $\delta^{\P^Y}$. Suppose $\xi$ is the largest such that $\R||\xi = \S||\xi$ and that $\R|\xi \neq \S|\xi$. Since no strategy disagreement occurs by Lemma \ref{lem:no_strat_disagreement} and the fact that $\R,\S$ are $X$-approved, letting $\Psi$ be the common strategy for $\R||\xi, \S||\xi$, then if $\xi$ indices an extender $E$ with $\cp(E) = \delta^{\P^Y}$ on the $\R$-sequence, then $E$ is $\pi^Y$-certified over $(\R||\xi, \Psi) = (\S||\xi, \Psi)$, this implies $E$ is the extender indexed at $\xi$ on the $\S$-sequence. Contradiction.

For (iv), the point is that in the relevant phalanx comparisons in the proof of solidity and universality, no extenders with critical point $(\delta^{\P})^Y$ are used by (iii), no strategy disagreements occur, and hence these comparisons are successful.
\end{proof}

Suppose without loss of generality, no countable substructures of any $\N\lhd \M$ satisfies Definition \ref{omega woodins over lsa}. We claim that for $Y$ as in (ii) of the above claim, $\M_Y$ does. Again, let $Y$ be as above and it suffices to show $\M^Y$ satisfies Definition \ref{omega woodins over lsa}. Everything is clear except, perhaps, for (1). So let $\Lambda$ be the $\pi' \circ \pi_Y$-realizable strategy for $\M^Y$ and $\Q = (\M^Y|\delta_0^Y)^{\sharp}$. By the argument as in Claim \ref{claim:soundness} and Lemma \ref{lem:locally_strong_branch_condensation}, $\Lambda^{stc}_\Q$ has (locally) strong branch condensation. Similarly to \ref{lem:fullness_preserving}, $\Lambda^{stc}_\Q$ is also (locally) strongly $\Omega$-fullness preserving and hence is (locallly) strongly $\Gamma(\Q,\Lambda^{stc}_\Q)$-fullness preserving. Lastly, the arguments in Section \ref{sec:sts_constructions} (particularly the proof of case 2.b) show that $\Lambda_\Q^{stc}\rest \M^Y \supseteq \Sigma^{\M_Y}$.


\end{proof}

Again, Lemma \ref{lem:hull} and results in Section 8.2 show that the new derived model of $\N$ as in the conclusion of Lemma \ref{lem:hull} (at the sup of its Woodin cardinals) satisfies $\sf{LSA}$. 

Now by boolean comparisons, there is some $(\M,\Sigma)\in V$ satisfying Definition \ref{omega woodins over lsa}. By taking a countable hull of $\M$ if necessary, we may assume $\M$ is countable (in $V$). Let $\M^-$ be the class model obtained by iterating the top extender of $\M$ OR many times and $\M_\infty$ be the result of an $\mathbb{R}$-genericity iteration of $\M^-$ via $\Sigma$. Then (new) derived model $N$ of $\M_\infty$ satisfies $\sf{LSA}$ as shown by Section 8.2. By homogeneity of $Col(\omega, < \kappa)$, there is in $V$ a model $M$ containing $\mathbb{R} \cup \rm{OR}$ such that $M\vDash \sf{LSA}$.
\newline

\noindent \textit{Proof of Theorem \ref{thm:square_lsa}.} The arguments above prove the consistency of $\sf{LSA}$ from the hypothesis of Theorem \ref{thm:square_lsa} plus the simplifying assumption (\ref{eqn:simplifying}). For this argument, since $2^{<\kappa}=\kappa$ and the core model induction is carried out in $V^{Coll(\omega,<\kappa)}$, $\sf{ord}$$(\P) < \kappa^{+}$ as shown in Lemma \ref{lem:small}. The constructions on top of $\P$ go on for at most $\kappa^{+++}$ many steps and the hypothesis we need to carry out the arguments in the previous sections is $\forall \alpha\in [\gamma,\kappa^{+4}] \ \neg \square(3,\alpha)$. 

To eliminate (\ref{eqn:simplifying}), we carry out the core model induction in $V^{Coll(\omega,\kappa)}$ much like that in \cite{Trang2015PFA} to obtain objects like $\Gamma, \P$. We note that in this case, using the argument in \ref{lem:small}, we can show $\sf{ord}$$(\P^-)\leq (2^\kappa)^+$ and $\sf{ord}$$(\P) < (2^\kappa)^{++}$. Now similar to the proof of \cite[Theorem 4.1]{JSSS}, letting $\xi = 2^\kappa$, $\mathbb{P} = Coll(\omega,\kappa)\star Coll(\xi^+,\xi^+)\star Coll(\xi^{++},\xi^{++})\star Coll(\xi^{+++},\xi^{+++})$, we carry out the hybrid $K^c$-constructions and the $X$-validated sts constructions over $\P$ in $V^{\mathbb{P}}$. By homogeneity, the objects constructed are in $V$. We use \rthm{thm:condensing_X}, which in turns was built on \ref{thm:weak_condensing_X}, \ref{lem:key_realizing_lemma}, and \ref{lem:full_X}, to obtain condensing sets and adapt the constructions in the previous sections of this chapter to obtain a model of $\sf{LSA}$ in a straightforward way. Note that in $V^{\mathbb{P}}$, $\xi^{+++}$ is countably closed and $2^{<\xi^{+++}}=\xi^{+++}$; this allows the proof of Lemma \ref{lem:stack_facts} to go through in this case. The constructions above $\P$ in the previous sections go on for at most $\xi^{+++}$ many steps and we need the full hypothesis $\forall \alpha\in [\gamma,(\powerset_4(\kappa)^{+}] \ \neg \square(3,\alpha)$ to carry out the arguments.\footnote{In $V^{\mathbb{P}}$, $\powerset_1^V(\xi)$ is collapsed to $\xi^+$, $\powerset_2^V(\xi)$ is collapsed to $\xi^{++}$, and $\powerset_3^V(\xi)$ is collapsed to $\xi^{+++}$. It seems very plausible that Lemma \ref{lem:stack_facts} can be proven with less than what we assumed, but we have not checked this thoroughly.} We leave the details to the kind reader.\footnote{Some definitions are modified in an obvious way. For instance, a good hull will now have size $\kappa$ and be countably closed; in Definition \ref{def:correctly_backgrounded_ext}, we demand that $|Z|\geq 2^{\kappa^+}$ and $Z^\kappa\subset Z$. } $\hfill \square$

\cleardoublepage
\phantomsection
\addcontentsline{toc}{chapter}{\indexname}
\printindex

\cleardoublepage

\phantomsection

\addcontentsline{toc}{chapter}{Bibliography}

\bibliographystyle{plain}
\bibliography{main}
\end{document}